\documentclass[11pt,reqno]{amsart}

\makeatletter
\def\@tocline#1#2#3#4#5#6#7{\relax
  \ifnum #1>\c@tocdepth % then omit
  \else
    \par \addpenalty\@secpenalty\addvspace{#2}%
    \begingroup \hyphenpenalty\@M
    \@ifempty{#4}{%
      \@tempdima\csname r@tocindent\number#1\endcsname\relax
    }{%
      \@tempdima#4\relax
    }%
    \parindent\z@ \leftskip#3\relax \advance\leftskip\@tempdima\relax
    \rightskip\@pnumwidth plus4em \parfillskip-\@pnumwidth
    #5\leavevmode\hskip-\@tempdima
      \ifcase #1
       \or\or \hskip 1em \or \hskip 2em \else \hskip 3em \fi%
      #6\nobreak\relax
    \dotfill\hbox to\@pnumwidth{\@tocpagenum{#7}}\par
    \nobreak
    \endgroup
  \fi}
\makeatother

\makeatletter
\newcounter{savesection}
\newcounter{apdxsection}
\renewcommand\appendix{\par
  \setcounter{savesection}{\value{section}}%
  \setcounter{section}{\value{apdxsection}}%
  \setcounter{subsection}{0}%
  \gdef\thesection{\@Alph\c@section}%
}
\newcommand\unappendix{\par
  \setcounter{apdxsection}{\value{section}}%
  \setcounter{section}{\value{savesection}}%
  \setcounter{subsection}{0}%
  \gdef\thesection{\@arabic\c@section}%
}
\makeatother

\usepackage{hyperref}
\hypersetup{
    colorlinks = true,
    allcolors = blue,
}

\usepackage{multirow}
\usepackage{fullpage} 
\usepackage{bbm}
\usepackage{tikz-cd}
\usepackage{mathabx,amsmath,amscd, latexsym, amssymb, bbold}
\usepackage{endnotes}
\usepackage{ulem} 
\usepackage[all, knot]{xy}

\usepackage{longtable}
 
\def\arraystretch{1.2} 

\usepackage[utf8]{inputenc}
\usepackage{amsmath,amsfonts}

\newcommand{\frc}[2]{{{#1}/{#2}}}

\def\1{\mathbbm{1}}
\def\a{\alpha}
\def\b{\beta}
\def\e{\epsilon}
\def\ve{\varepsilon}

\def\s{\sigma}
\def\es{{\varepsilon_\s}}

\def\hs{\hat\sigma}

\def\g{\gamma}

\def\tj{\tilde{j}}

\def\fs{\mathfrak{s}}
\def\ft{\mathfrak{t}}
\def\fg{\mathfrak{g}}

\def\BO{{\mathbb{O}}}
\def\BQ{{\mathbb{Q}}}
\def\BC{{\mathbb{C}}}

\def\BZ{{\mathbb{Z}}}
\def\BR{{\mathbb{R}}}
\def\BN{{\mathbb{N}}}

\def\cB{{\mathcal{B}}}

\def\SO{{\operatorname{SO}}}
\def\O{{\operatorname{O}}}

\def\qed{{\ \ \ \mbox{$\square$}}}
\def\Rep{\operatorname{Rep}}
\newcommand\mtx[1]{\begin{bmatrix} #1 \end{bmatrix}}

\newcommand\Inv{\operatorname{Inv}}
\newcommand\rank{\operatorname{rank}}
\newcommand\irr{\operatorname{irr}}
\newcommand\spec{\operatorname{spec}}
\newcommand\GQ{\Gal(\bar\BQ)}

\renewcommand{\deg}{\dim}

\newcommand\rd[3]{\rho_{{#1}_{#2}^{#3}}}

\def\Gal{{\mbox{\rm Gal}}}

\newtheorem{thm}{Theorem}[section]
\newtheorem{cor}[thm]{Corollary}
\newtheorem{prop}[thm]{Proposition}
\newtheorem{lem}[thm]{Lemma}

\newtheorem{defn}[thm]{Definition}

\newtheorem{remark}[thm]{Remark}

\numberwithin{equation}{section}

 \newcommand{\gbeg}[2]{
   \unitlength=1pt
   \grrow = #2
   \grcolumn = 0
   \grcalca = #1
   \grcalcb = #2
   \multiply \grcalca by \factor
   \grwidth = \grcalca
   \multiply \grcalcb by \factor
   \begin{minipage}{\grcalca pt}
   \begin{picture}(\grcalca,\grcalcb)
   \advance \grcalcb by -\factor
}

\newcommand\ord{\operatorname{ord}}
\newcommand\pord{\operatorname{pord}}
\newcommand\nd{{\operatorname{ndeg}}}
\newcommand{\diag}{\operatorname{diag}}
\newcommand{\FPdim}{\operatorname{FPdim}}
\newcommand{\MD}{{\operatorname{MD}}}
\newcommand{\isum}{{\operatorname{isum}}}
\newcommand{\pMD}{{\operatorname{pMD}}}
\newcommand\sgn{{\operatorname{sgn}}}

\newcommand\Z{\BZ}
\newcommand\ol[1]{\overline{#1}}
\newcommand\replace[1]{}

\newcommand\td{\tilde}

\newcommand\id{\operatorname{id}}
\renewcommand\o{\otimes}

\newcommand\Tr{\operatorname{Tr}}

\newcommand\SL{\operatorname{SL}_2(\BZ)}
\newcommand\hSL{{\widehat{\operatorname{SL}_2(\BZ)}}}
\newcommand\SLn[1]{\operatorname{SL}_2(\BZ/#1\BZ)}
\newcommand\qsl[1]{{\operatorname{SL}_2(\BZ / #1 \BZ)}}

\newcommand\inv{^{-1}}

\newcommand\CC{\mathcal C}

\def\ft{\mathfrak{t}}
\def\fs{\mathfrak{s}}

\makeatletter
\def\namelabel#1#2{\@bsphack
  \protected@write\@auxout{}%
         {\string\newlabel{#1.nme}{{#2}{#2}}}%
  \@esphack}

\makeatother

\def\Gal{{\mbox{\rm Gal}}}

\newcommand{\jacobi}[2]{{\left(\frac{#1}{#2}\right)}}

\usepackage{color}
\definecolor{white}{rgb}{1,1,1}
\definecolor{grey}{rgb}{0.5,0.5,0.5}
\newcommand{\white}[1]{\color{white}{#1}}

\def\ee{\mathrm{e}}
\def\ii{\mathrm{i}}

\allowdisplaybreaks
\numberwithin{equation}{section}

\usepackage[symbol]{footmisc}

\title{Reconstruction of modular data\\ from $\SL$ representations}
\author[Ng]{Siu-Hung Ng}
\author[Rowell]{Eric C Rowell}
\author[Wang]{Zhenghan Wang}
\author[Wen]{Xiao-Gang Wen}

\begin{document}

\maketitle

\begin{abstract} Modular data  is the most significant invariant of a modular tensor category.
We pursue an approach to the classification of modular data of modular tensor categories by building the modular $S$ and $T$ matrices directly from irreducible representations of $\qsl{n}$.  We discover and collect
many conditions on the $\qsl{n}$ representations to identify those that correspond to some modular data.  To arrive at concrete matrices from representations, we also develop methods that allow us to select the proper basis of the
$\qsl{n}$ representations so that they have the form of modular data.  We apply this technique to the classification of rank-$6$ modular tensor categories, obtaining a classification up to modular data.  Most of the calculations can be automated using a computer algebraic system, which can be employed to classify modular data of higher rank modular tensor categories.

\end{abstract}

\tableofcontents
\setcounter{tocdepth}{2}

\section{Introduction}

Just as conventional symmetries are described by groups, gapped quantum liquid phases of bosonic matter (\emph{i.e} bosonic topological order) seem to be described by non-degenerate higher braided fusion categories.  It has been conjectured that topological orders are classified by the collection of projective representations of mapping class groups for various topologies of closed space manifolds \cite{W9039}.  In particular, we believe that a gapped phase of quantum matter in two spacial dimensions is classified by a pair $(\mathcal{B},c)$, where $\mathcal{B}$ is a unitary modular tensor category (MTC) and $c$ is a rational number equal to the central charge of $\mathcal{B}$ mod 8.  Physically, $\mathcal{B}$ models the topological excitations (\emph{i.e.} the anyons) in the gapped phase \cite{K062}, and $c$ measures the possible stacking of $E_8$ quantum Hall state, which has central charge $c=8$. Therefore, a classification of unitary MTCs would give rise to a classification of all gapped quantum phases of bosons without symmetry in two spacial dimensions.

MTCs are defined by very complicated data.  The classification of MTCs naturally breaks into two steps: the first step is to classify the modular data (MD), and the second is to classify modular isotopes with a given MD if not unique. The MD $(S,T)$ of an MTC form a projective representation of the mapping class group of the 2-dimensional torus. (In fact, the notion of topological order was first introduced based on modular data $(S,T)$ \cite{W9039}.) We will see that the classification of MDs is much more manageable than the full classification of MTCs.

Modular data $(S,T)$ corresponding to MTCs of rank $r \leq 5$ have been
completely classified \cite{CRank5,RSW0777,HR}.  More recently,
such a classification for MTCs of rank 6 containing a pair of non-self-dual simple objects and a partial classification of general MTCs of rank $6$ has also been obtained \cite{Creamer}.
%, relying heavily upon computational techniques.  
The strategy employed in those classifications begins with a stratification of the Galois group of the extension of $\BQ$ by the entries of the modular $S$ matrix, followed by a case by case analysis on the inferred polynomial
constraints. As the Galois group is isomorphic to an abelian subgroup of
$\mathfrak{S}_r$, this program is tractable, although somewhat tedious. As a
last resort in a few cases, the classification of low-dimensional
representations of $\SLn{n}$ for small $n$ was required as well.  The typical
outcome is that most Galois groups can be eliminated and one eventually finds a finite list of modular data which can then be realized from known
constructions.

In this article we complete the classification of rank $6$ MDs using
the reverse strategy: we build upon the approach in \cite{Eh, CRank5} by constructing the MDs directly from $\SLn{n}$
representations of low dimension.  Since $n$ is bounded in terms of the rank,
expressing irreducible $\qsl{n}$ representations as tensor products of
prime-power level  representations (i.e. $\SLn{p^k}$ for primes $p$) allows us to stratify by representation type and level. Thus, up to basis choice, the $\SL$ representations can be presented as pairs $(s,t)$, where $s$ is
symmetric and $t$ is diagonal.  The construction of \emph{symmetric} representations of $\SL$ is an interesting problem of its own \cite{NWW21, NWW}.  We
note that the number of inequivalent $\qsl{n}$ representations is finite at a given
dimension, since the dimension and $n$ are bounded in terms of the rank.
These facts make our classification possible.  

In the next step of our classification, for each representation $(s,t)$, we
conjugate $s$ by an arbitrary (real orthogonal) matrix that commutes with $t$ to reconstruct the potential MD $(S,T)$ with $S$ symmetric and $T$ diagonal.  We find several methods that allow us to select a finite number of possible real orthogonal matrices from the uncountable set of real orthogonal matrices, so that the resulting $(S,T)$ include all the MDs. Up to reordering the objects in the category, \emph{i.e.}, the rows/columns of the resulting $(S,T)$, these must satisfy the algebraic and number-theoretic constraints of MDs.  Case by case analyses, following a similar pattern, then yield our classification.

The approach to the classification of MDs by building the modular $S$ and $T$ matrices directly from irreducible representations of $\qsl{n}$ is applicable to much more general cases than the rank 6 case in this paper.  
One version of our approaches that is presented in the Appendix can be automated and almost all of the calculations in this approach can be implemented using the GAP computer algebra system.

The content of the paper is as follows: In sections 2 and 3,  we discover and collect
many conditions on the $\qsl{n}$ representations to help
us identifying those that are from some MDs.  To arrive at concrete matrices from representations, we also develop methods that allow us to select the proper basis of the
$\qsl{n}$ representations so that they become the MDs. In sections 4 and 5,  we apply this technique to the classification of rank-$6$ MTCs, obtaining a classification up to MD.  Most of the calculations can be automated using a computer algebraic system, which can be employed to classify MDs of higher rank MTCs.

\section{Modular tensor categories and modular data}

Given a modular tensor category (MTC)\footnote{We use the terminology of MTC as in its original sense \cite{MS89}, which is equivalent to a semi-simple modular category of \cite{TuraevBook}, i.e. a semi-simple modular category.} $\mathcal{B}$, the modular data (MD) of $\mathcal{B}$ consists of the un-normalized $S$- and $T$- matrices of $\mathcal{B}$, hence the MD of an MTC is independent of any normalizations.  Though the MD of an MTC does not determine the MTC uniquely \cite{MigSch}, it is still the most useful and important invariant of an MTC.  Moreover, the MDs of MTCs have enchanting relations with diverse areas from congruence subgroups to vector-valued modular forms to topological phases of matter.

\subsection{Necessary conditions for the modular data of an MTC}

An obvious strategy to classify MDs would be first to find all necessary and sufficient conditions for MDs, and then simply look for solutions.  But it seems very hard to find such a complete characterization of MD.  Instead we will list some necessary conditions and then appeal to other methods to finish a classification.

The following collection of results on modular data which will be useful in the sequel.  Many are well-known and found in, e.g. \cite{BK}.  

\begin{thm}
\label{p:MD}
The modular data $(S,T)$ of an MTC satisfies:
\begin{enumerate}
\item
$S,T$ are symmetric complex matrices, indexed by $i,j=0,\ldots, r-1$.\footnote{The index
also labels the simple objects in the MTC, with $i=0$ corresponding to the unit
object, and $r$ is the \textbf{rank} of the modular data and the MTC.}

\item
$T$ is unitary, diagonal, and $T_{00}=1$.

\item
$S_{00}=1$. Let  $d_i = S_{0i} $ and
$D=\sqrt{\sum_{i=0}^{r-1} d_i d_i^*}$. Then 
\begin{align}
S S^\dag = D^2\id,
\end{align}
and the $d_i\in\mathbb{R}$.
\item
$S_{ij}$ are cyclotomic integers in $\BQ_{\ord(T)}$\footnote{Here $\BQ_n$ denotes the field $\BQ(\zeta_n)$ for a primitive $n$th root of unity $\zeta_n$} \cite{NS10}. 
The ratios
$S_{ij}/S_{0j}$ are cyclotomic integers for all $i,j$ \cite{CG}. Also there is a $j$ such
that $S_{ij}/S_{0j} \in [1,+\infty) $ for all $i$ \cite{ENO}.

\item

Let $\theta_i = T_{ii}$ and $p_\pm = \sum_{i=0}^{r-1} d_i^2 (\theta_i)^{\pm 1}$.

Then
$p_+/p_-$ is a root of unity, and $p_+=D\ee^{\ii 2\pi c/8}$ for some rational number $c$.\footnote{The \textbf{central charge} $c$ of the modular data and of the MTC is only defined modulo $8$.} Moreover, the modular data $(S,T)$ is associated with a projective $\SL$ representation, since:
\begin{align}
% p p^* = D^2, \ \ \ \ redundant
 (ST)^3 = p_+ S^2,\ \ \ \frac{S^2}{D^2} =C,\ \ \ C^2 =\id,
\end{align}
where $C$ is a permutation matrix
satisfying
\begin{align}
 \Tr(C) > 0.
\end{align}
%The root of unity $\ee^{2\pi\ii\frac{c}{8}}$ is called the multiplicative
%central charge of the MTC.  

\item Cauchy Theorem \cite{BNRW}: The set of distinct prime
factors of $\ord(T)$ coincides with the distinct prime factors of
$\mathrm{norm}(D^2)$.\footnote{Here $\mathrm{norm}(x)$ is the product of the distinct Galois
conjugates of the algebraic number $x$.}

\item  Verlinde formula $($cf. \cite{V8860}$):$

\begin{equation}\label{Ver0} 
 N^{ij}_k = \frac{1}{D^2}\sum_{l=0}^{r-1} \frac{S_{li} S_{lj} S_{lk}^*}{ d_l }
\in \BN ,  
\end{equation} 
where $i,j, k=0,1,\ldots,r-1$ and $\BN$ is the set of non-negative integers.\footnote{The $N^{ij}_k$ are called the fusion coefficients.} The $N^{ij}_0$ satisfy
\begin{align}
 N^{ij}_0 = C_{ij},
\end{align}
which defines a
charge conjugation $i \to \bar i$ via
\begin{align}
 N^{\bar i j}_0 = \delta_{i j}.
\end{align}

\item 
\label{FScnd}
Let $n \in \BN_+$.  The $n^\text{th}$ Frobenius-Schur indicator of the $i$-th
simple object
\begin{equation}
 \label{nunFS}
 \nu_n(i)= \sum_{j, k} N_i^{jk} (d_j\theta_j^n) (d_k\theta_k^n)^*
\end{equation}
%\red{XG: The above expression needs to updated for complex $d_i$.}
%\begin{equation}
% \label{nunFS}
% \nu_n(i)= \frac{1}{D^2}\sum_{j, k} N_i^{jk}
%d_{j}\theta_j^n (d_{k}\theta_k^n)^*
%\end{equation}
 is a cyclotomic integer whose conductor divides $n$ and
$\ord(T)$ \cite{NS07a, NS07b}.  The 1st Frobenius-Schur indicator satisfies $\nu_1(i)=\delta_{i,0}$ while the 2nd
Frobenius-Schur indicator $\nu_2(i)$ satisfies $\nu_2(i)=0$ if $i\neq \bar i$,
and $\nu_2(i)=\pm 1$ if $i = \bar i$ (see \cite{Bantay, NS07a, RSW0777}).

%\item 
%For any Galois conjugation $\s$ in $\Gal(\BQ_{\ord(\rho_\a(\ft))})$, there is a
%permutation of the indices, $i \to \hs(i)$, and $\e_\s(i)\in \{1,-1\}$, such
%that
%\begin{align}
%\label{Galact}
%\s \big(\rho_\a(\fs)_{ij}\big) &
%= \e_\s(i)\rho_\a(\fs)_{\hat \s (i),j} 
%= \rho_\a(\fs)_{i,\hat \s (j)}\e_\s(j) 
%\nonumber\\
%\s^2 \big(\rho_\a(\ft)_{i,i}\big) &= \rho_\a(\ft)_{\hat \s (i),\hat \s (i)},
%\end{align}
%for all $i,j$. 

%\red{XG: we may drop the follow items.}
%
%\item The matrix elements of $S,T$ are cyclotomic integers
%\begin{align}
% S_{ij} \in \BO_{\ord(T)},\ \ \ \ \
% T_{ij} \in \BO_{\ord(T)},
%\end{align}
%where $\BO_{\ord(T)}$ is the ring extension of the ring of integers $\BZ$ by the
%$\ord(T)^\mathrm{th}$ roots of unity ({\it i.e.} by $\ee^{\ii 2\pi/ \ord(T)}$).
%
%
%
%\item $\Gal(\BQ(S)/\BQ)$ is isomorphic with an Abelian subgroup of the symmetric
%group $S_r$. $\Gal(\BQ(T)/\BQ(S))=\BZ_2^l$ for an integer $l$ \cite{DLN}.  Let
%$\omega(n)$ denotes the number of distinct prime factors of an integer $n$,
%then $l \le \omega(\ord(T))$ if $8 \nmid\ord(T) $, and $l \le
%\omega(\ord(T))+1$ if $8 \mid \ord(T)$.   We have
%\begin{align}
% \frac{\Gal(\BQ(T)/\BQ)}{ \Gal(\BQ(T)/\BQ(S))} \cong \Gal(\BQ(S)/\BQ).
%\end{align}
%
%\item
%The prime divisors of $\mathrm{norm}(D^2)$ and $\ord(T)$ coincide \cite{BNRW}
%and $D^5/\ord(T)$ is a cyclotomic integers \cite{E02}. Therefore, $D^6/\ord(T)$
%is a cyclotomic integer in $\BQ_{\ord(T)}$.
%
%

\end{enumerate}
\end{thm}

We denote by $\Gal(\BQ_n)$ the Galois group of the cyclotomic field $\BQ_n$.
\begin{remark}
{\rm
The above conditions are for modular data of unitary or
non-unitary MTCs.  In particular, the above conditions are invariant under
Galois conjugations in $\Gal(\BQ_{\ord(T)}/\BQ)$.  Therefore, we can group modular
data into Galois orbits. 

Mathematical definition of Frobenius-Schur indicators of an object in pivotal fusion category was introduced in \cite{NS07a} and the trichotomy of the 2nd Frobenius-Schur indicator of a simple object was also proved therein. If the underlying \emph{pivotal structure} is not spherical, the $d_i$ in the preceding theorem could be complex. We do not need this for the sequel, but it may lead to an interesting generalization.
}
\end{remark}

\subsection{Classification of modular data up to rank=$5$ and candidate list of rank=$6$}

\subsubsection{Rank 1-5 MTCs}

\begin{table}[tb]
\caption{Rank $\leq 5$ modular data
\label{rank5}
}
\centerline{
\small
\begin{tabular}{ |c|c|l|l||c|c|l|l| } \hline $N_c$ & $D^2$ &
$d_0,d_1,\cdots$ & $s_0,s_1,\cdots$ & $N_c$ & $D^2$ & $d_0,d_1,\cdots$
& $s_0,s_1,\cdots$ \\
 \hline $1_{ 1}$ & $1$ & $1$ & $0$ & & & & \\
 \hline $2_{
1}$ & $2$ & $1,1$ & $0, \frac{1}{4}$ & $2_{-1}$ & $2$ & $1,1$ &
$0,-\frac{1}{4}$ \\
 $2_{\frc{14}{5}}$ & $3.6180$ & $1,\frac{1+\sqrt{5}}{2}$ & $0,
\frac{2}{5}$ & $2_{-\frc{14}{5}}$ & $3.6180$ & $1,\frac{1+\sqrt{5}}{2}$ &
$0,-\frac{2}{5}$ \\
 \hline $3_{ 2}$ & $3$ & $1,1,1$ & $0, \frac{1}{3},
\frac{1}{3}$ & $3_{-2}$ & $3$ & $1,1,1$ & $0,-\frac{1}{3},-\frac{1}{3}$\\
$3_{\frc{1}{2}}$ & $4$ & $1,1,\sqrt{2}$ & $0,\frac{1}{2}, \frac{1}{16}$ &
$3_{-\frc{1}{2}}$ & $4$ & $1,1,\sqrt{2}$ & $0,\frac{1}{2},-\frac{1}{16}$ \\
$3_{\frc{3}{2}}$ & $4$ & $1,1,\sqrt{2}$ & $0,\frac{1}{2}, \frac{3}{16}$ &
$3_{-\frc{3}{2}}$ & $4$ & $1,1,\sqrt{2}$ & $0,\frac{1}{2},-\frac{3}{16}$ \\
$3_{\frc{5}{2}}$ & $4$ & $1,1,\sqrt{2}$ & $0,\frac{1}{2}, \frac{5}{16}$ &
$3_{-\frc{5}{2}}$ & $4$ & $1,1,\sqrt{2}$ & $0,\frac{1}{2},-\frac{5}{16}$ \\
$3_{\frc{7}{2}}$ & $4$ & $1,1,\sqrt{2}$ & $0,\frac{1}{2}, \frac{7}{16}$ &
$3_{-\frc{7}{2}}$ & $4$ & $1,1,\sqrt{2}$ & $0,\frac{1}{2},-\frac{7}{16}$ \\
$3_{\frc{8}{7}}$ & $9.2946$ & $1,\xi_7^{2},\xi_7^{3}$ & $0,-\frac{1}{7},
\frac{2}{7}$ & $3_{-\frc{8}{7}}$ & $9.2946$ & $1,\xi_7^{2},\xi_7^{3}$ & $0,
\frac{1}{7},-\frac{2}{7}$ \\
 \hline $4^{a}_{ 0}$ & $4$ & $1,1,1,1$ & $0, 0,
0,\frac{1}{2}$ &  $4^{b}_{ 0}$ & $4$ &
$1,1,1,1$ &  $0, 0, \frac{1}{4},-\frac{1}{4}$ \\
 $4_{ 1}$ & $4$ &
$1,1,1,1$ & $0, \frac{1}{8}, \frac{1}{8},\frac{1}{2}$ & $4_{-1}$ & $4$ &
$1,1,1,1$ & $0,-\frac{1}{8},-\frac{1}{8},\frac{1}{2}$ \\
  $4_{ 2}$
& $4$ &  $1,1,1,1$ &  $0, \frac{1}{4},
\frac{1}{4},\frac{1}{2}$ &  $4_{-2}$ & $4$
&  $1,1,1,1$ &
$0,-\frac{1}{4},-\frac{1}{4},\frac{1}{2}$ \\
 $4_{ 3}$ & $4$ & $1,1,1,1$ & $0,
\frac{3}{8}, \frac{3}{8},\frac{1}{2}$ & $4_{-3}$ & $4$ & $1,1,1,1$ &
$0,-\frac{3}{8},-\frac{3}{8},\frac{1}{2}$ \\
 $4_{4}$ & $4$ & $1,1,1,1$ &
$0,\frac{1}{2},\frac{1}{2},\frac{1}{2}$ &  $4_{\frc{9}{5}}$
& $7.2360$ & $1,1,\frac{1+\sqrt{5}}{2},\frac{1+\sqrt{5}}{2}$ &
$0,-\frac{1}{4}, \frac{3}{20}, \frac{2}{5}$ \\
  $4_{-\frc{9}{5}}$
& $7.2360$ & $1,1,\frac{1+\sqrt{5}}{2},\frac{1+\sqrt{5}}{2}$ &
$0, \frac{1}{4},-\frac{3}{20},-\frac{2}{5}$ & $4_{ \frc{19}{5}}$ & $7.2360$ & $1,1,\frac{1+\sqrt{5}}{2},\frac{1+\sqrt{5}}{2}$
& $0, \frac{1}{4},-\frac{7}{20}, \frac{2}{5}$ \\
$4_{-\frc{19}{5}}$ & $7.2360$ &
$1,1,\frac{1+\sqrt{5}}{2},\frac{1+\sqrt{5}}{2}$ & $0,-\frac{1}{4},
\frac{7}{20},-\frac{2}{5}$ &  $4^{c}_{ 0}$ &
$13.090$ & $1,\frac{1+\sqrt{5}}{2},\frac{1+\sqrt{5}}{2},\frac{1+\sqrt{5}}{2}\frac{1+\sqrt{5}}{2}$ &
$0, \frac{2}{5},-\frac{2}{5}, 0$ \\
  $4_{\frc{12}{5}}$
& $13.090$ & $1,\frac{1+\sqrt{5}}{2},\frac{1+\sqrt{5}}{2},\frac{1+\sqrt{5}}{2}\frac{1+\sqrt{5}}{2}$
& $0,-\frac{2}{5},-\frac{2}{5}, \frac{1}{5}$ &
 $4_{-\frc{12}{5}}$ & $13.090$ &
$1,\frac{1+\sqrt{5}}{2},\frac{1+\sqrt{5}}{2},\frac{1+\sqrt{5}}{2}\frac{1+\sqrt{5}}{2}$ & $0, \frac{2}{5},
\frac{2}{5},-\frac{1}{5}$ \\
 $4_{\frc{10}{3}}$ & $19.234$ &
$1,\xi_9^{2},\xi_9^{3},\xi_9^{4}$ & $0, \frac{1}{3}, \frac{2}{9},-\frac{1}{3}$
& $4_{-\frc{10}{3}}$ & $19.234$ & $1,\xi_9^{2},\xi_9^{3},\xi_9^{4}$ &
$0,-\frac{1}{3},-\frac{2}{9}, \frac{1}{3}$ \\
 \hline $5_{ 0}$ & $5$ &
$1,1,1,1,1$ & $0, \frac{1}{5}, \frac{1}{5},-\frac{1}{5},-\frac{1}{5}$ & $5_{4}$
& $5$ & $1,1,1,1,1$ & $0, \frac{2}{5},
\frac{2}{5},-\frac{2}{5},-\frac{2}{5}$ \\
 $5^{a}_{ 2}$ & $12$ &
$1,1,\xi_6^{2},\xi_6^{2},2$ & $0, 0, \frac{1}{8},-\frac{3}{8}, \frac{1}{3}$ &
$5^{b}_{ 2}$ & $12$ & $1,1,\xi_6^{2},\xi_6^{2},2$ & $0, 0,-\frac{1}{8},
\frac{3}{8}, \frac{1}{3}$ \\
 $5^{b}_{-2}$ & $12$ & $1,1,\xi_6^{2},\xi_6^{2},2$
& $0, 0, \frac{1}{8},-\frac{3}{8},-\frac{1}{3}$ & $5^{a}_{-2}$ & $12$ &
$1,1,\xi_6^{2},\xi_6^{2},2$ & $0, 0,-\frac{1}{8}, \frac{3}{8},-\frac{1}{3}$ \\
$5_{\frc{16}{11}}$ & $34.645$ & $1,\xi_{11}^{2},\xi_{11}^{3},\xi_{11}^{4},\xi_{11}^{5}$ &
$0,-\frac{2}{11}, \frac{2}{11}, \frac{1}{11},-\frac{5}{11}$ &
$5_{-\frc{16}{11}}$ & $34.645$ & $1,\xi_{11}^{2},\xi_{11}^{3},\xi_{11}^{4},\xi_{11}^{5}$ &
$0, \frac{2}{11},-\frac{2}{11},-\frac{1}{11}, \frac{5}{11}$ \\
  $5_{\frc{18}{7}}$ &  $35.339$ &
$1,\xi_7^{3},\xi_7^{3},\xi_{14}^{3},\xi_{14}^{5}$ &
$0,-\frac{1}{7},-\frac{1}{7}, \frac{1}{7}, \frac{3}{7}$ &
$5_{-\frc{18}{7}}$ &  $35.339$ &
$1,\xi_7^{3},\xi_7^{3},\xi_{14}^{3},\xi_{14}^{5}$ &  $0,
\frac{1}{7}, \frac{1}{7},-\frac{1}{7},-\frac{3}{7}$ \\
 \hline \end{tabular}
}
\end{table}

The rank$\leq 5$ unitary MTCs are classified \cite{CRank5,RSW0777,HR}; Table \ref{rank5} lists all 45 rank $\leq 5$ cases, only the quantum dimensions and twists are displayed.  These are labeled by $N_c$, where $N$ is the rank and
$c$ the central charge.  The entries of the table are ordered by the total
quantum dimension $D^2$. Also $d_i$ is the quantum dimension and
$s_i=\mathrm{arg}(T_{ii})$ is the \textbf{topological spin} of the $i^\text{th}$
simple object in the MTC. For $r\in\BQ$ we often use the notation $e(r)=e^{2\pi i r}$, so that the topological spin of the object labeled by $i$ is $e(s_i)$.  The  quantum dimensions are given in terms of
$\xi_n^{m,k}=\frac{\sin(m\pi/n)}{\sin(k\pi/n)}$ and $\xi_n^{m}=\xi_n^{m,1}$.  The fusion coefficients
$N^{ij}_k$ and the $S$-matrices of MTCs can be deduced from the given data in these low rank cases, and we do not list them for brevity's sake.

\subsubsection{Known rank-6 MD of MTCs and their Galois Groups}
Among the known rank 6 modular tensor categories there are 11 distinct fusion rules. We can determine their Galois groups $\Gal(\BQ(S_{ij})/\BQ)$ and the representation type (i.e. dimensions of their irreducible subrepresentations) of their $\SL$ representation, displayed in Table \ref{rank6realizations}. Six are realized as product categories, the other $5$ by prime categories. Note that there are two types that yield the fusion rules of $SO(5)_2$: $(3,2,1)$ is realized by a zesting of $SO(5)_2$, see Theorem \ref{t:321}. \begin{table}[tb]
\caption{Realizations of known rank $6$ modular data, their Galois groups and representation types.}
\label{rank6realizations}

\begin{tabular}{ |c|c|c| } 
\hline 
$\CC$ & $\Gal(\CC)$ & Type  \\
 \hline 
 $PSU(2)_3\boxtimes SU(2)_2$& $\langle(0\;1)(2\;3),(0\;2)(1\;3)(4\;5) \rangle\cong \mathbb{Z}_2\times\mathbb{Z}_2$&$(6)$\\
 \hline
 $PSU(2)_3\boxtimes U(3)_1$ & $\langle (0\;1)(2\;3)(4\;5),(2\;4)(3\;5) \rangle\cong\mathbb{Z}_2\times\mathbb{Z}_2$& $(4,2)$\\
 \hline
 $PSU(2)_3\boxtimes PSU(2)_5$ &  $\langle  (0\;1)(2\;3)(4\;5),(0\;2\;4)(1\;3\;5) \rangle\cong\mathbb{Z}_6$& $(6)$
    \\ \hline   $U(2)_1\boxtimes SU(2)_2$ &  $\langle (0\;1)(2\;3)\rangle\cong\mathbb{Z}_2$& $(6)$
    \\ \hline  $U(2)_1\boxtimes U(3)_1$ &  $\langle (1\;2)(3\;4)\rangle\cong\mathbb{Z}_2$& $(4,2)$
    \\ \hline  $U(2)_1\boxtimes PSU(2)_5$ &  $\langle (0\;1\;2)(3\;4\;5)\rangle\cong\mathbb{Z}_3$& $(6)$
    \\ \hline $SO(5)_2$ &  $\langle(0\;1)(2\;3) \rangle\cong\mathbb{Z}_2$& $(3,3)$, $(3,2,1)$
    \\ \hline $PSU(2)_{11}$ &  $\langle (0\;1\;2\;3\;4\;5)\rangle\cong\mathbb{Z}_6$& $(6)$
    \\ \hline $G(2)_3$ &  $\langle (0\;1),(2\;3\;4) \rangle\cong\mathbb{Z}_6$& $(4,2)$
    \\ \hline $PSO(8)_3$ &  $\langle(0\;1\;2) \rangle\cong\mathbb{Z}_3$& $(4,1,1)$
    \\ \hline $PSO(5)_{\frac{3}{2}}$  &  $\langle (0\;1\;2)(3\;4\;5)\rangle\cong\mathbb{Z}_3$& $(6)$\\ \hline
\end{tabular}\end{table}

\begin{table}[tb]
\caption{Table of rank 6 modular data with $N^{ij}_k \leq 3$ and $D^2\leq 18$.
\label{rank6a}
}
\centerline{
\small
\begin{tabular}{ |c|c|c|l|l|c| } 
\hline 
$N_c$ & $D^2$ & $d_0,d_1,\cdots$ & $s_0,s_1,\cdots$ & comment \\
 \hline 
$6_{ 1}$ &
 $6$ & $1,1,1,1,1,1$ & $0, \frac{1}{12}, \frac{1}{12},-\frac{1}{4},
\frac{1}{3}, \frac{1}{3}$ & $2_{-1}\boxtimes 3_{ 2}$\\
 $6_{-1}$ &  $6$
& $1,1,1,1,1,1$ & $0,-\frac{1}{12},-\frac{1}{12},
\frac{1}{4},-\frac{1}{3},-\frac{1}{3}$ & $2_{ 1}\boxtimes 3_{-2}$\\
 $6_{ 3}$ &
 $6$ & $1,1,1,1,1,1$ & $0, \frac{1}{4}, \frac{1}{3},
\frac{1}{3},-\frac{5}{12},-\frac{5}{12}$ & $2_{ 1}\boxtimes 3_{ 2}$\\
 $6_{-3}$ &
 $6$ & $1,1,1,1,1,1$ & $0,-\frac{1}{4},-\frac{1}{3},-\frac{1}{3},
\frac{5}{12}, \frac{5}{12}$ & $2_{-1}\boxtimes 3_{-2}$\\
 \hline $6_{
\frc{1}{2}}$ & $8$ & $1,1,1,1,\sqrt{2},\sqrt{2}$ & $0,
\frac{1}{4},-\frac{1}{4},\frac{1}{2},-\frac{1}{16}, \frac{3}{16}$ & $2_{
1}\boxtimes 3_{-\frc{1}{2}}$\\
 $6_{-\frc{1}{2}}$ &  $8$ &
$1,1,1,1,\sqrt{2},\sqrt{2}$ & $0, \frac{1}{4},-\frac{1}{4},\frac{1}{2},
\frac{1}{16},-\frac{3}{16}$ & $2_{ 1}\boxtimes 3_{-\frc{3}{2}}$\\
 $6_{
\frc{3}{2}}$ &  $8$ & $1,1,1,1,\sqrt{2},\sqrt{2}$ & $0,
\frac{1}{4},-\frac{1}{4},\frac{1}{2}, \frac{1}{16}, \frac{5}{16}$ & $2_{
1}\boxtimes 3_{\frc{1}{2}}$\\
 $6_{-\frc{3}{2}}$ & $8$ &
$1,1,1,1,\sqrt{2},\sqrt{2}$ & $0,
\frac{1}{4},-\frac{1}{4},\frac{1}{2},-\frac{1}{16},-\frac{5}{16}$ & $2_{
1}\boxtimes 3_{-\frc{5}{2}}$\\
 $6_{\frc{5}{2}}$ & $8$ &
$1,1,1,1,\sqrt{2},\sqrt{2}$ & $0, \frac{1}{4},-\frac{1}{4},\frac{1}{2},
\frac{3}{16}, \frac{7}{16}$ & $2_{ 1}\boxtimes 3_{\frc{3}{2}}$\\
$6_{-\frc{5}{2}}$ & $8$ & $1,1,1,1,\sqrt{2},\sqrt{2}$ & $0,
\frac{1}{4},-\frac{1}{4},\frac{1}{2},-\frac{3}{16},-\frac{7}{16}$ & $2_{
1}\boxtimes 3_{-\frc{7}{2}}$\\
 $6_{\frc{7}{2}}$ & $8$ &
$1,1,1,1,\sqrt{2},\sqrt{2}$ & $0, \frac{1}{4},-\frac{1}{4},\frac{1}{2},
\frac{5}{16},-\frac{7}{16}$ & $2_{ 1}\boxtimes 3_{\frc{5}{2}}$\\
$6_{-\frc{7}{2}}$ & $8$ & $1,1,1,1,\sqrt{2},\sqrt{2}$ & $0,
\frac{1}{4},-\frac{1}{4},\frac{1}{2},-\frac{5}{16}, \frac{7}{16}$ & $2_{
1}\boxtimes 3_{\frc{7}{2}}$\\
 \hline $6_{\frc{4}{5}}$ & $10.854$ &
$1,1,1,\frac{1+\sqrt{5}}{2},\frac{1+\sqrt{5}}{2},\frac{1+\sqrt{5}}{2}$ & $0,-\frac{1}{3},-\frac{1}{3},
\frac{1}{15}, \frac{1}{15}, \frac{2}{5}$ & $2_{\frc{14}{5}}\boxtimes 3_{-2}$\\
$6_{-\frc{4}{5}}$ & $10.854$ &
$1,1,1,\frac{1+\sqrt{5}}{2},\frac{1+\sqrt{5}}{2},\frac{1+\sqrt{5}}{2}$ & $0, \frac{1}{3},
\frac{1}{3},-\frac{1}{15},-\frac{1}{15},-\frac{2}{5}$ &
$2_{-\frc{14}{5}}\boxtimes 3_{ 2}$\\
 $6_{\frc{16}{5}}$ & $10.854$ &
$1,1,1,\frac{1+\sqrt{5}}{2},\frac{1+\sqrt{5}}{2},\frac{1+\sqrt{5}}{2}$ & $0,-\frac{1}{3},-\frac{1}{3},
\frac{4}{15}, \frac{4}{15},-\frac{2}{5}$ & $2_{-\frc{14}{5}}\boxtimes 3_{-2}$\\
$6_{-\frc{16}{5}}$ & $10.854$ &
$1,1,1,\frac{1+\sqrt{5}}{2},\frac{1+\sqrt{5}}{2},\frac{1+\sqrt{5}}{2}$ & $0, \frac{1}{3},
\frac{1}{3},-\frac{4}{15},-\frac{4}{15}, \frac{2}{5}$ & $2_{
\frc{14}{5}}\boxtimes 3_{ 2}$\\
 \hline $6_{\frc{3}{10}}$ & $14.472$
& $1,1,\sqrt{2},\frac{1+\sqrt{5}}{2},\frac{1+\sqrt{5}}{2},\sqrt{2}\frac{1+\sqrt{5}}{2}$ &
$0,\frac{1}{2},-\frac{5}{16},-\frac{1}{10}, \frac{2}{5}, \frac{7}{80}$ & $2_{
\frc{14}{5}}\boxtimes 3_{-\frc{5}{2}}$\\
 $6_{-\frc{3}{10}}$ &
$14.472$ & $1,1,\sqrt{2},\frac{1+\sqrt{5}}{2},\frac{1+\sqrt{5}}{2},\sqrt{2}\frac{1+\sqrt{5}}{2}$ &
$0,\frac{1}{2}, \frac{5}{16}, \frac{1}{10},-\frac{2}{5},-\frac{7}{80}$ &
$2_{-\frc{14}{5}}\boxtimes 3_{\frc{5}{2}}$\\
 $6_{\frc{7}{10}}$ &
$14.472$ & $1,1,\sqrt{2},\frac{1+\sqrt{5}}{2},\frac{1+\sqrt{5}}{2},\sqrt{2}\frac{1+\sqrt{5}}{2}$ &
$0,\frac{1}{2}, \frac{7}{16}, \frac{1}{10},-\frac{2}{5}, \frac{3}{80}$ &
$2_{-\frc{14}{5}}\boxtimes 3_{\frc{7}{2}}$\\
 $6_{-\frc{7}{10}}$ &
$14.472$ & $1,1,\sqrt{2},\frac{1+\sqrt{5}}{2},\frac{1+\sqrt{5}}{2},\sqrt{2}\frac{1+\sqrt{5}}{2}$ &
$0,\frac{1}{2},-\frac{7}{16},-\frac{1}{10}, \frac{2}{5},-\frac{3}{80}$ & $2_{
\frc{14}{5}}\boxtimes 3_{-\frc{7}{2}}$\\
 $6_{\frc{13}{10}}$ &
$14.472$ & $1,1,\sqrt{2},\frac{1+\sqrt{5}}{2},\frac{1+\sqrt{5}}{2},\sqrt{2}\frac{1+\sqrt{5}}{2}$ &
$0,\frac{1}{2},-\frac{3}{16},-\frac{1}{10}, \frac{2}{5}, \frac{17}{80}$ & $2_{
\frc{14}{5}}\boxtimes 3_{-\frc{3}{2}}$\\
 $6_{-\frc{13}{10}}$ &
$14.472$ & $1,1,\sqrt{2},\frac{1+\sqrt{5}}{2},\frac{1+\sqrt{5}}{2},\sqrt{2}\frac{1+\sqrt{5}}{2}$ &
$0,\frac{1}{2}, \frac{3}{16}, \frac{1}{10},-\frac{2}{5},-\frac{17}{80}$ &
$2_{-\frc{14}{5}}\boxtimes 3_{\frc{3}{2}}$\\
 $6_{\frc{17}{10}}$ &
$14.472$ & $1,1,\sqrt{2},\frac{1+\sqrt{5}}{2},\frac{1+\sqrt{5}}{2},\sqrt{2}\frac{1+\sqrt{5}}{2}$ &
$0,\frac{1}{2},-\frac{7}{16}, \frac{1}{10},-\frac{2}{5}, \frac{13}{80}$ &
$2_{-\frc{14}{5}}\boxtimes 3_{-\frc{7}{2}}$\\
 $6_{-\frc{17}{10}}$ &
$14.472$ & $1,1,\sqrt{2},\frac{1+\sqrt{5}}{2},\frac{1+\sqrt{5}}{2},\sqrt{2}\frac{1+\sqrt{5}}{2}$ &
$0,\frac{1}{2}, \frac{7}{16},-\frac{1}{10}, \frac{2}{5},-\frac{13}{80}$ & $2_{
\frc{14}{5}}\boxtimes 3_{\frc{7}{2}}$\\
 $6_{\frc{23}{10}}$ &
$14.472$ & $1,1,\sqrt{2},\frac{1+\sqrt{5}}{2},\frac{1+\sqrt{5}}{2},\sqrt{2}\frac{1+\sqrt{5}}{2}$ &
$0,\frac{1}{2},-\frac{1}{16},-\frac{1}{10}, \frac{2}{5}, \frac{27}{80}$ & $2_{
\frc{14}{5}}\boxtimes 3_{-\frc{1}{2}}$\\
 $6_{-\frc{23}{10}}$ &
$14.472$ & $1,1,\sqrt{2},\frac{1+\sqrt{5}}{2},\frac{1+\sqrt{5}}{2},\sqrt{2}\frac{1+\sqrt{5}}{2}$ &
$0,\frac{1}{2}, \frac{1}{16}, \frac{1}{10},-\frac{2}{5},-\frac{27}{80}$ &
$2_{-\frc{14}{5}}\boxtimes 3_{\frc{1}{2}}$\\
 $6_{\frc{27}{10}}$ &
$14.472$ & $1,1,\sqrt{2},\frac{1+\sqrt{5}}{2},\frac{1+\sqrt{5}}{2},\sqrt{2}\frac{1+\sqrt{5}}{2}$ &
$0,\frac{1}{2},-\frac{5}{16}, \frac{1}{10},-\frac{2}{5}, \frac{23}{80}$ &
$2_{-\frc{14}{5}}\boxtimes 3_{-\frc{5}{2}}$\\
 $6_{-\frc{27}{10}}$ &
$14.472$ & $1,1,\sqrt{2},\frac{1+\sqrt{5}}{2},\frac{1+\sqrt{5}}{2},\sqrt{2}\frac{1+\sqrt{5}}{2}$ &
$0,\frac{1}{2}, \frac{5}{16},-\frac{1}{10}, \frac{2}{5},-\frac{23}{80}$ & $2_{
\frc{14}{5}}\boxtimes 3_{\frc{5}{2}}$\\
 $6_{\frc{33}{10}}$ &
$14.472$ & $1,1,\sqrt{2},\frac{1+\sqrt{5}}{2},\frac{1+\sqrt{5}}{2},\sqrt{2}\frac{1+\sqrt{5}}{2}$ &
$0,\frac{1}{2}, \frac{1}{16},-\frac{1}{10}, \frac{2}{5}, \frac{37}{80}$ & $2_{
\frc{14}{5}}\boxtimes 3_{\frc{1}{2}}$\\
 $6_{-\frc{33}{10}}$ &
$14.472$ & $1,1,\sqrt{2},\frac{1+\sqrt{5}}{2},\frac{1+\sqrt{5}}{2},\sqrt{2}\frac{1+\sqrt{5}}{2}$ &
$0,\frac{1}{2},-\frac{1}{16}, \frac{1}{10},-\frac{2}{5},-\frac{37}{80}$ &
$2_{-\frc{14}{5}}\boxtimes 3_{-\frc{1}{2}}$\\
 $6_{\frc{37}{10}}$ &
$14.472$ & $1,1,\sqrt{2},\frac{1+\sqrt{5}}{2},\frac{1+\sqrt{5}}{2},\sqrt{2}\frac{1+\sqrt{5}}{2}$ &
$0,\frac{1}{2},-\frac{3}{16}, \frac{1}{10},-\frac{2}{5}, \frac{33}{80}$ &
$2_{-\frc{14}{5}}\boxtimes 3_{-\frc{3}{2}}$\\
 $6_{-\frc{37}{10}}$ &
$14.472$ & $1,1,\sqrt{2},\frac{1+\sqrt{5}}{2},\frac{1+\sqrt{5}}{2},\sqrt{2}\frac{1+\sqrt{5}}{2}$ &
$0,\frac{1}{2}, \frac{3}{16},-\frac{1}{10}, \frac{2}{5},-\frac{33}{80}$ & $2_{
\frc{14}{5}}\boxtimes 3_{\frc{3}{2}}$\\
 \hline 
\end{tabular} 
}
\end{table}

\begin{table}[tb]
\caption{Table of rank 6 modular data with $N^{ij}_k \leq 3$ and $D^2>18$.
\label{rank6b}
}
\centerline{
\small
\begin{tabular}{ |c|c|c|l|l|c| } 
\hline 
$N_c$ & $D^2$ & $d_0,d_1,\cdots$ & $s_0,s_1,\cdots$ & comment \\
 \hline 
$6_{\frc{1}{7}}$ &
$18.591$ & $1,1,\xi_7^{2},\xi_7^{2},\xi_7^{3},\xi_7^{3}$ &
$0,-\frac{1}{4},-\frac{1}{7},-\frac{11}{28}, \frac{1}{28}, \frac{2}{7}$ &
$2_{-1}\boxtimes 3_{\frc{8}{7}}$\\
 $6_{-\frc{1}{7}}$ & $18.591$ &
$1,1,\xi_7^{2},\xi_7^{2},\xi_7^{3},\xi_7^{3}$ & $0, \frac{1}{4}, \frac{1}{7},
\frac{11}{28},-\frac{1}{28},-\frac{2}{7}$ & $2_{ 1}\boxtimes 3_{-\frc{8}{7}}$\\
$6_{\frc{15}{7}}$ & $18.591$ &
$1,1,\xi_7^{2},\xi_7^{2},\xi_7^{3},\xi_7^{3}$ & $0, \frac{1}{4},
\frac{3}{28},-\frac{1}{7}, \frac{2}{7},-\frac{13}{28}$ & $2_{ 1}\boxtimes 3_{
\frc{8}{7}}$\\
 $6_{-\frc{15}{7}}$ & $18.591$ &
$1,1,\xi_7^{2},\xi_7^{2},\xi_7^{3},\xi_7^{3}$ & $0,-\frac{1}{4},-\frac{3}{28},
\frac{1}{7},-\frac{2}{7}, \frac{13}{28}$ & $2_{-1}\boxtimes 3_{-\frc{8}{7}}$\\
\hline 
$6^{a}_{ 0}$ & $20$ & $1,1,2,2,\sqrt{5},\sqrt{5}$ & $0, 0,
\frac{1}{5},-\frac{1}{5}, 0,\frac{1}{2}$ & root of $SO(10)_2$ \\
 $6^{b}_{ 0}$ & $20$ &
$1,1,2,2,\sqrt{5},\sqrt{5}$ & $0, 0, \frac{1}{5},-\frac{1}{5},
\frac{1}{4},-\frac{1}{4}$ & root of $SO(10)_2$ \\
\hline 
$6^{a}_{4}$ & $20$ &
$1,1,2,2,\sqrt{5},\sqrt{5}$ & $0, 0, \frac{2}{5},-\frac{2}{5}, 0,\frac{1}{2}$ &
root of $SO(5)_2$\\
 $6^{b}_{4}$ & $20$ & $1,1,2,2,\sqrt{5},\sqrt{5}$ & $0, 0,
\frac{2}{5},-\frac{2}{5}, \frac{1}{4},-\frac{1}{4}$ &  $SO(5)_2$\\
\hline 
$6_{ \frc{58}{35}}$ & $33.632$ &
$1,\frac{1+\sqrt{5}}{2},\xi_7^{2},\xi_7^{3},\frac{1+\sqrt{5}}{2}\xi_7^{2},\frac{1+\sqrt{5}}{2}\xi_7^{3}$ & $0,
\frac{2}{5}, \frac{1}{7},-\frac{2}{7},-\frac{16}{35}, \frac{4}{35}$ & $2_{
\frc{14}{5}}\boxtimes 3_{-\frc{8}{7}}$\\
 $6_{-\frc{58}{35}}$ &
$33.632$ &
$1,\frac{1+\sqrt{5}}{2},\xi_7^{2},\xi_7^{3},\frac{1+\sqrt{5}}{2}\xi_7^{2},\frac{1+\sqrt{5}}{2}\xi_7^{3}$ &
$0,-\frac{2}{5},-\frac{1}{7}, \frac{2}{7}, \frac{16}{35},-\frac{4}{35}$ &
$2_{-\frc{14}{5}}\boxtimes 3_{\frc{8}{7}}$\\
 $6_{\frc{138}{35}}$ 
& $33.632$ &
$1,\frac{1+\sqrt{5}}{2},\xi_7^{2},\xi_7^{3},\frac{1+\sqrt{5}}{2}\xi_7^{2},\frac{1+\sqrt{5}}{2}\xi_7^{3}$ & $0,
\frac{2}{5},-\frac{1}{7}, \frac{2}{7}, \frac{9}{35},-\frac{11}{35}$ & $2_{
\frc{14}{5}}\boxtimes 3_{\frc{8}{7}}$\\
 $6_{-\frc{138}{35}}$ &
$33.632$ &
$1,\frac{1+\sqrt{5}}{2},\xi_7^{2},\xi_7^{3},\frac{1+\sqrt{5}}{2}\xi_7^{2},\frac{1+\sqrt{5}}{2}\xi_7^{3}$ &
$0,-\frac{2}{5}, \frac{1}{7},-\frac{2}{7},-\frac{9}{35}, \frac{11}{35}$ &
$2_{-\frc{14}{5}}\boxtimes 3_{-\frc{8}{7}}$\\
 \hline {$6_{\frc{46}{13}}$}
& {$56.746$} &
{$1,\xi_{13}^{2},\xi_{13}^{3},\xi_{13}^{4},\xi_{13}^{5},\xi_{13}^{6}$}
& {$0, \frac{4}{13}, \frac{2}{13},-\frac{6}{13},
\frac{6}{13},-\frac{1}{13}$} & root of $SU(2)_{11}$\\
 {$6_{-\frc{46}{13}}$} &
{$56.746$} &
{$1,\xi_{13}^{2},\xi_{13}^{3},\xi_{13}^{4},\xi_{13}^{5},\xi_{13}^{6}$}
& {$0,-\frac{4}{13},-\frac{2}{13}, \frac{6}{13},-\frac{6}{13},
\frac{1}{13}$} &  root of $SU(2)_{\widebar{11}}$\\
 \hline {$6_{\frc{8}{3}}$} &
{$74.617$} &
{$1,\xi_{18}^{3},\xi_{18}^{3},\xi_{18}^{3},\xi_{18}^{5},\xi_{18}^{7}$}
& {$0, \frac{1}{9}, \frac{1}{9}, \frac{1}{9}, \frac{1}{3},-\frac{1}{3}$} &
 root of $SO(8)_{\bar 3}$ \\
 {$6_{-\frc{8}{3}}$} & {$74.617$} &
{$1,\xi_{18}^{3},\xi_{18}^{3},\xi_{18}^{3},\xi_{18}^{5},\xi_{18}^{7}$}
& {$0,-\frac{1}{9},-\frac{1}{9},-\frac{1}{9},-\frac{1}{3}, \frac{1}{3}$} &
root of $SO(8)_{3}$ \\
\hline 
 $6_{2}$ &   $100.61$ &  $1,\frac{3+\sqrt{21}}{2},\frac{3+\sqrt{21}}{2},\frac{3+\sqrt{21}}{2},\frac{5+\sqrt{21}}{2},\frac{7+\sqrt{21}}{2}$ &  $0,-\frac{1}{7},-\frac{2}{7}, \frac{3}{7}, 0, \frac{1}{3}$ & root of $G(2)_{\bar 3}$ \\
 $6_{-2}$ &   $100.61$ &  $1,\frac{3+\sqrt{21}}{2},\frac{3+\sqrt{21}}{2},\frac{3+\sqrt{21}}{2},\frac{5+\sqrt{21}}{2},\frac{7+\sqrt{21}}{2}$ &  $0, \frac{1}{7}, \frac{2}{7},-\frac{3}{7}, 0,-\frac{1}{3}$ & root of $G(2)_3$ \\
 \hline 
\end{tabular} 
}
\end{table}

\iffalse
\begin{enumerate}
    \item $PSU(2)_3\boxtimes SU(2)_2$ Galois group $\langle(0\;1)(2\;3),(0\;2)(1\;3)(4\;5) \rangle\cong \mathbb{Z}_2\times\mathbb{Z}_2$, type $(6)$.
    \item  $PSU(2)_3\boxtimes U(3)_1$ Galois group  $\langle (0\;1)(2\;3)(4\;5),(2\;4)(3\;5) \rangle\cong\mathbb{Z}_2\times\mathbb{Z}_2$, type $(4,2)$.
    \item  $PSU(2)_3\boxtimes PSU(2)_5$ Galois group  $\langle  (0\;1)(2\;3)(4\;5),(0\;2\;4)(1\;3\;5) \rangle\cong\mathbb{Z}_6$, type $(6)$.
    \item   $U(2)_1\boxtimes SU(2)_2$ Galois group  $\langle (0\;1)(2\;3)\rangle\cong\mathbb{Z}_2$, type $(6)$.
    \item  $U(2)_1\boxtimes U(3)_1$ Galois group  $\langle (1\;2)(3\;4)\rangle\cong\mathbb{Z}_2$, type $(4,2)$.
    \item  $U(2)_1\boxtimes PSU(2)_5$ Galois group  $\langle (0\;1\;2)(3\;4\;5)\rangle\cong\mathbb{Z}_3$, type $(6)$.

    \item $SO(5)_2$ Galois group  $\langle(0\;1)(2\;3) \rangle\cong\mathbb{Z}_2$, types $(3,3)$ and $(3,2,1)$\footnote{There are two types that yield the same fusion rules, type $(3,2,1)$ is realized by a zesting of $SO(5)_2$, see Theorem \ref{t:321}.}. 
    \item $PSU(2)_{11}$ Galois group  $\langle (0\;1\;2\;3\;4\;5)\rangle\cong\mathbb{Z}_6$, type $(6)$.
    \item $G(2)_3$ Galois group  $\langle (0\;1),(2\;3\;4) \rangle\cong\mathbb{Z}_6$, type $(4,2)$.
    \item $PSO(8)_3$ Galois group  $\langle(0\;1\;2) \rangle\cong\mathbb{Z}_3$, type $(4,1,1)$.
    \item $PSO(5)_{\frac{3}{2}}$  Galois group  $\langle (0\;1\;2)(3\;4\;5)\rangle\cong\mathbb{Z}_3$, type $(6)$.
\end{enumerate}
\fi

The example $PSO(5)_{\frac{3}{2}}$ is noteworthy--it is the smallest example of a MTC the fusion rules of which are never realized as those of a \emph{unitary} MTC.  We also remark that the fusion rules of $SO(5)_2$ are realized by categories with distinct representation types: namely the \emph{zested} version of $SO(5)_2$, see Theorem \ref{t:321}.  In particular, the fusion rules do not determine the representation type.

We also did a computer search for all rank-6 unitary modular data with
$N^{ij}_k \leq 3$.  (Ref. \cite{W150605768} computed all rank-6 unitary modular
data with $N^{ij}_k \leq 2$.) The Tables \ref{rank6a} and \ref{rank6b} list all
 50 of the resulting modular data, we include only the quantum dimensions and twists.  In the last column, $N_c\boxtimes N'_{c'}$
indicates that the rank-6 MTC is the product of two MTCs labeled by $N_c$ and
$N'_{c'}$.  The prime MTCs are all non-Abelian roots of MTCs from Kac-Moody
algebra. (The notion of non-Abelian roots is introduced in Ref.
\cite{LW170107820}.) In this paper, we will show that the Tables \ref{rank6a}
and \ref{rank6b} include all modular data of rank-6 unitary MTCs.

\section{Modular data representations of modular tensor categories}

While the number theoretical properties of MD allow the classification of MTCs up to rank=$4$, the deeper properties of the $\SL$ representations of MD (cf. Definition \ref{d:MDrep})  lead to a more streamlined approach with the potential to achieve a classification up to rank=$10$. The classification of rank=$5$ MTCs is  already a mixture of both Galois theory and representation techniques.  Instead of working on cases labeled by abelian subgroups of $S_r$ for rank=$r$ as in earlier classification, we introduce the notion {\it type} of the MD of an MTC--the list of dimensions of irreducible subrepresentations, so that the cases are indexed by Young diagrams with $r$ boxes.

Every MTC $\mathcal{B}$ leads to a $(2+1)$-TQFT, hence there is a corresponding projective matrix representation $\ol\rho_{\mathcal{B}}$ of $\SL$---the mapping class group of the torus.  We will refer to this representation as the projective $\SL$ representation of the MTC $\mathcal{B}$, and is given by the $S$-, $T$- matrices of $\mathcal{B}$. The 
 linearizations of this projective matrix $\SL$ representation $\ol\rho_\cB$, called $\SL$ representations  of $\mathcal{B}$,  will be elaborated upon in next section.

\subsection{$\SL$ representations of MTC or MD}
\label{SL2Zrep}

%Modular data $S,T$ are closely related to certain $\SL$ representations with
%\emph{a certain choices of basis}, which will be referred to as $\SL$ matrix
%representations.  We  use the term \emph{matrix representation}, instead of
%representation, to stress the choice of basis.  

Since our classification is based on $\SL$ representations, let us first
summarize some important facts about them.  Let $\fs=\begin{bmatrix} 0&-1\\ 1&
0 \end{bmatrix} $, $\ft =\begin{bmatrix} 1&1\\ 0& 1 \end{bmatrix} $ be the standard generators of $\SL$. This admits the presentation:
$$
\SL = \langle \fs, \ft\mid \fs^4 = \id, (\fs \ft)^3 = \fs^2\rangle\,.
$$
The 1-dimensional representations of $\SL$, denoted $\hSL$, form a cyclic group of order 12 under tensor product.  We will take $\chi \in \hSL$ defined by $\chi(\ft) = \zeta_{12}$ to be the generator, where $\zeta_n^k:=\ee^{2 \pi \ii k/n}$. Under this convention, every 1-dimensional representation of $\SL$ is equivalent to $\chi^\a$ for some  integer $\a$, unique modulo 12:
\begin{align}
\chi^\a(\fs) = \ol \zeta_4^\a, \quad \chi^\a(\ft)=\zeta_{12}^\a .
\end{align}

Given a modular tensor category $\cB$ with the modular data $(S,T)$ and central charge $c$,
the assignment
\begin{align}
 \rho_\a(\fs) = \ol\zeta_4^\a S/D, \ \
 \rho_\a(\ft) = \zeta_{12}^\a \ee^{-2\pi \ii \frac{c}{24}} T
\ \ \ (\a \in \BZ_{12}). 
\end{align}
define a (linear) representation of $\SL$, and we call these  representations $\rho_\a$ the {\bf $\SL$ representations of $\cB$} or  the {\bf $\SL$ representations of the modular data $(S,T)$}.  For any $\a, \a' \in \BZ_{12}$,
$$
\rho_\a \cong \chi^{\a-\a'} \o \rho_{\a'}
$$
as $\SL$ representations. Therefore, the $\SL$ representation $\rho_\cB$ of $\cB$ is unique up to a tensor factor of linear characters of $\widehat{\SL}$.

Note that two modular data $(S,T)$ and $(S',T')$ are regarded as \emph{the same} if
they differ only by a permutation of indices:
\begin{align}
 S' = P S P^\top,\ \ \ \
 T' = P T P^\top,
\end{align}
where $P$ is a permutation matrix. Throughout this paper, we simply \emph{identify} $\rho_\a$ and its conjugations by permutation matrices. 
\begin{defn} \label{d:MDrep} {\rm
A unitary matrix representation $\rho$ of $\SL$ is called an \emph{MD representation} if $\rho$ is an $\SL$ representation of some modular tensor category. It is called a \emph{ pseudo-MD (pMD) representation} if $V \rho V$ is an  MD representation for some signed diagonal matrix $V$.}
\end{defn}

\subsection{Type and level of modular data}

\begin{defn}{\rm
Given an MTC $\mathcal{B}$ of rank $r$, an $\SL$ representation $\rho_{\mathcal{B}}$  decomposes into direct sum of irreducible representations of dimensions $\lambda_1, \ldots, \lambda_m$ in non-increasing order.
The \emph{type} of the corresponding MD of $\mathcal{B}$ of rank=$r$ is the Young diagram of $r$ boxes $(\lambda_1,\ldots,\lambda_m)$ with $\sum_{i=1}^m\lambda_i=r$. The type of an MTC simply refers to the type of its MD.}
\end{defn}

The modular representations of the Fibonacci and Ising theories are both irreducible, so they are of types $(2), (3)$, respectively.  The modular representation of the toric code has an image isomorphic to $\qsl{2}$ and is reducible of type $(2,1,1)$. 

We note that for any positive integer $n$, the reduction $\BZ \to \BZ/n \BZ$
defines a surjective group homomorphism  $\pi_n:\SL \to \qsl{n}$. Thus, a
representation of $\qsl{n}$ is also a representation of $\SL$, which will be
called a \emph{congruence} representation of $\SL$ in this paper.  It is
immediate to see that a representation of $\qsl{n}$ is also a $\qsl{mn}$ representation for
any positive integer $m$. The smallest positive integer $n$ such that a
congruence representation $\rho$ of $\SL$  factors through $\pi_n : \SL \to
\qsl{n}$ is called the \emph{level} of $\rho$. It is known that the level
$n=\ord(\rho(\ft))$ (cf.  \cite[Lem. A.1]{DLN}).  Here $\ord(t)$ is the \emph{order} of $t$, i.e., the 
smallest positive integer such that
\begin{align}
 t^{\ord(t)} = \id.
\end{align}

There are many more finite-dimensional noncongruence representations of $\SL$ (cf. \cite{KL}) but they are not associated with any modular tensor category by \cite[Thm. II]{DLN}. Since we only deal with congruence representations of $\SL$,  all the representations of $\SL$ throughout this paper are assumed to be congruence and finite-dimensional over $\BC$. 

An $\SL$ representation $\rho$ of an MTC is also \emph{symmetric}, which means $\rho$ is a unitary matrix representation with $\rho(\fs)$ symmetric and $\rho(\ft)$ diagonal. The following theorem proved in \cite{NWW} provides the theoretic background for the GAP package \cite{NWW21} and our reconstruction process:
\begin{thm}
Every finite-dimensional congruence representation of $\SL$ is equivalent to a symmetric one.
\end{thm}
Therefore, throughout this paper, we always assume our \emph{general} representations of $\SL$ to be congruence and symmetric.

In Appendix \ref{repPP}, we list all the irreducible $\SL$ representations, generated by \cite{NWW21}, of prime-power levels and dimensions $\le 6$. These $\SL$ representations are congruence and symmetric. From these representations, we can
construct all the inequivalent $\SL$ representations with dimensions  $\le 6$. The MD representations of dimensions $\le 6$ can be reconstructed from these symmetric representations with the help of the following theorem.
\begin{thm}\label{t:ortho_eqv}
Let $\rho, \rho': \SL{} \to U_n(\BC)$ be  unitarily equivalent symmetric representations of $\SL{}$ such that $\rho(\ft) = \rho'(\ft)=t$, and define $s=\rho(\fs)$ and $s'=\rho'(\fs)$. Then there exists a (real) orthogonal matrix $U$ such that 
$$
s'=UsU^{\top}\quad \text{and} \quad Ut = tU.
$$ 
\end{thm}
\begin{proof}
Let $Q$ be a unitary matrix such that    
$$
s'=QsQ^\dag\quad \text{and} \quad Qt = tQ.
$$ 
Since $t$ is diagonal and unitary, $t^\dag=\ol t$. Taking the conjugate transpose of the second equality implies
$$
Q^\dag \ol t = \ol t Q^\dag\quad \text{or}\quad \ol Q t = t \ol Q\,.
$$
Let $Q=X_1+i X_2$ for some real matrices $X_1$ and $X_2$. Then we have
$$
(X_1 \pm i X_2) t = t (X_1 \pm  i X_2)
$$
which implies $[X_i, t] = 0$ for $i=1,2$. Similarly, $s'Q = Qs$ implies $\ol s' Q =  Q \ol s $ since both $s$ and $s'$ are symmetric. Therefore, we also have $s' \ol Q= \ol Q s$, which implies
$$
X_i s =  s' X_i \quad \text{for } i =1,2.
$$
Since there are only finitely many roots for the equation $\det(X_1 + x X_2)=0$, one can take $\lambda \in \BR$ such that $X=X_1+ \lambda X_2$ is invertible. Then
$$
X s = s'X \quad \text{and}\quad X t = t X\,.
$$
 Let $X = UP$ be the polar decomposition of $X$ where $U$ is orthogonal and $P$ is the unique positive definite satisfying $P^2 = X^{\top} X$. In fact, $P$ is a polynomial of $P^2$ (cf. \cite[Chap.9. Thm 11.]{HK}). 
Since  $s\inv = \ol s$ and ${s'}\inv = \ol{s'}$, 
$$
P^2  = X^{\top} X   =(s' X \ol s)^\dag (s' X \ol s)=s X^{\top} {s'}^\dag s' X s^{\top} =  s P^2 \ol s
$$  
and 
$$
X^{\top} t = t X^{\top} \,.
$$ 
Therefore,
$$
P^2 s = s P^2 \quad\text{and}\quad P^2 t = t P^2\,.
$$
Since $P$ is a polynomial of $P^2$, we find
$$
P s = s P \quad\text{and}\quad P t = t P\,.
$$
Therefore,
$$
Us = UPsP\inv = X sP\inv=s' XP\inv = s'U 
$$
and
$$
Ut = UPtP\inv = Xt P\inv =t XP\inv=tU\,.\qedhere
$$
\end{proof}

\begin{remark}
{\rm
An $\SL$ representation $\rho$ is said to be \emph{even} (resp. \emph{odd}) if $\rho(\fs^2) = \id$ (resp. $\rho(\fs^2)=-\id$).  If $\rho$ is  symmetric and irreducible, then $\rho(\fs)$ or $i\rho(\fs)$ is a real symmetric matrix, depending on whether $\rho$ is even or odd respectively. A direct sum of irreducible representations of opposite parties is neither even nor odd. In particular, if $\rho$ is an $\SL$ representation of a modular tensor category $\CC$, then  $\rho$ is even or odd if, and only if, $\CC$ is self-dual. 
}
\end{remark}

\subsection{Useful conditions on $\SL$ representations}

The set of all the roots of unity can be totally ordered as follows:   For any roots of unity $x, y$, we say that $x < y$ if one the following conditions hold:
\begin{enumerate}
    \item[(i)] $\ord(x) < \ord(y)$, or 
    \item[(ii)] $\ord(x)= \ord(y)$   and  $\arg(x) < \arg(y)$,
\end{enumerate}
where $\arg(\zeta)$ denotes the unique number $s_\zeta \in [0, 1) \cap \BQ$ such that $e^{2 i \pi  s_\zeta} = \zeta$. 

For any representation $\rho$ of $\SL$, $\rho(\ft)$ has finite order. We
denoted by $\spec(\rho(\ft))$ the increasing ordered set of eigenvalues of
$\rho(\ft)$ with multiplicities. If
$\spec(\rho(\ft))$ is multiplicity free $\rho$ is called \emph{non-degenerate}.  If $\rho'$ is another representation
of $\SL$, $\spec(\rho(\ft)) =\{x_1, \dots, x_m\} $ and
$\spec(\rho'(\ft))=\{y_1, \dots, y_n\}$ can be compared by the lexicographical
order.

Two representations $\rho, \rho'$ of $\SL$ are called \emph{projectively equivalent} if
$$
\rho' \cong \chi^\a \o \rho \text{ for some } \a \in \BZ/12\BZ\,. 
$$
A representations $\rho$ of $\SL$ is said to have a \emph{minimal $\ft$-spectrum} if $\spec(\rho(\ft))$ is minimal among all the representations projectively equivalent to $\rho$, i.e.,
$$
\spec(\rho(\ft))  \le \spec((\chi^\a \o \rho)(\ft)) \text{ for all } \alpha \in \BZ/12 \BZ\,.
$$
\begin{defn}{\rm
Let $t$ be any matrix over $\BC$. The smallest positive integer $n$ such that
$t^n=\a \id$ for some $\a \in \BC$ is called the \emph{projective order} of $t$,
and denoted by $\pord(t):=n$. If such integer does not exist, we define
$\pord(t):=\infty$. }
\end{defn}
%\begin{defn}
%A set of 12 projectively equivalent matrix representations of $\SL$, $\rho_\a$,
%$\a=1,\cdots,12$, is called a set of \textbf{MD representations} for a modular
%data $S,T$, if one of the matrix representations $\rho_\mathrm{MD} \in
%\{\rho_\a\}$ satisfies the following condition: there exist an index $u
%\in\{0,1,\cdots,r-1\}$ and a permutation matrix $P$ such that
%\begin{align}
%\label{STrho}
% \frac{\rho_\mathrm{MD}(\fs)}{\big(\rho_\mathrm{MD}(\fs)\big)_{u,u}} =  PSP^\top,\ \ \ \
% \frac{\rho_\mathrm{MD}(\ft)}{\big(\rho_\mathrm{MD}(\ft)\big)_{u,u}} = PTP^\top.
%\end{align}
%%where $c$ is the central charge.
%\end{defn}
%Here, the index $u$ corresponds to the unit object in the MTC and the MD
%representations $\rho_\a$ are unitary.  Note that, by definition, $\rho(\ft)$
%in an MD representation is always diagonal. Also if one of the 12
%representations satisfies the above condition, then all the 12 representations
%satisfy the above condition.  The reverse is also true: if $\rho_\mathrm{MD}$
%is an MD representation, then the set of the 12 MD representations of the same
%modular data are given by 
%\begin{align}
% \rho_\a(\ft) = \rho_\mathrm{MD}(\ft) \ee^{\ii 2\pi \frac{\a}{12}}, \ \ \
% \rho_\a(\fs) = \rho_\mathrm{MD}(\fs) \ee^{-\ii 2\pi \frac{\a}{4}}, \ \ \
%\a = 1,\cdots,12.
%\end{align}
%Two MD representations $\rho$ and $\rho'$ of a modular data will be considered
%to be equivalent if they differ only by a permutation of indices.

We can organize the irreducible representations of $\SL$ by the level and the
dimension of the representations.  Due to the Chinese remainder theorem, if the
level of a irreducible representation $\rho$ factors as $n=\prod_i
p_i^{k_i}$ where $p_i$ are distinct primes, then $\rho \cong \bigotimes_i \rho_i$
where $\rho_i$ are level $p_i^{k_i}$ representations.  Thus we can construct
all irreducible $\SL$ representations as tensor products of irreducible $\SL$ representations
of prime-power levels, which in turn, yields a construction of all $\SL$
representations $\rho$ via direct sums of the irreducible representations.

Define $\BQ_n=\BQ(\zeta_{n})$ to be the cyclotomic field of order $n$. For any positive integer $n$, we can construct a faithful representation
$D_n: \Gal(\BQ_n) \to \qsl{n}$, which identifies the Galois group $\Gal(\BQ_n) \cong \BZ_n^\times$ with the diagonal subgroup of $\qsl{n}$ \cite[Remark 4.5]{DLN}. More generally,  for any $\s \in \GQ$, $\s(\BQ_n) = \BQ_n$ and so there exists an integer $a$ (unique modulo $n$) such that $\s( \zeta_n )=\zeta_n^a$ and  
\begin{align}
\label{DSLn}
 D_n(\s) :=
\ft^a \fs \ft^b \fs \ft^a \fs^{-1}  
=
\begin{pmatrix}
 a & 0\\
 0 & b\\
\end{pmatrix} \in \qsl{n}\,,
\end{align}
 where $b$ satisfies $ab \equiv 1 \mod n$. If $\rho$ is a  level $n$ representation of $\SL$, the composition 
\begin{equation} \label{eq:Drho}
    D_\rho(\s):=\rho \circ D_n(\s)
\end{equation}
defines a representation of $\GQ$. We may also write $D_{n}(\s)$ as $D_n(a)$.   Such a representation of Galois group captures the Galois conjugation action on $\SL$ representations $\rho_{\MD}$ of  modular data, and plays a very important role in our classification. Many of the following collection of results on $\rho_{\MD}$ were proved in \cite{NS10, DLN}.

\begin{thm}
\label{p:MD1}
Every $\SL$ representation $\rho$ of an MTC $\cB$ is a matrix representation with the standard basis $(e_0, \dots, e_{r-1})$ identified with $\irr(\cB)$. Assuming $e_0 = \1$, $\rho$ satisfies the following conditions:
\begin{enumerate}
\item Let $n = \ord(\rho(\ft))$. For any $\fg \in \SL$, $\rho(\fg)$ is a matrix over $\BQ_n$. 
In particular, $\rho(\fs)_{ij}$ are cyclotomic numbers in $\BQ_n$ for all $i,j$.
\item The modular data $(S,T)$ of $\cB$ is given by 
\begin{align}
\label{STrho1}
 S = \frac{\rho(\fs)}{\rho(\fs)_{00}} , \ \ \ \ \
 T= \frac{\rho(\ft)}{\rho(\ft)_{00}}.
\end{align}
\item In particular, $\rho$ is symmetric,  $\ord(T) = \pord\rho(\ft))$ and (cf. Theorem \ref{p:MD}(4)) 
$$
\frac{\rho(\fs)_{ij}}{\rho(\fs)_{0j}}   \in \BZ[\zeta_{\ord(T)}]\,.
$$
\item The representation $\rho$ is congruence of level $n$
$\ord(T) \mid n \mid 12 \ord(T)$. Thus, $\rho$ is a symmetric and congruence $\SL$ representation. 

\item
 One has $1/\rho(\fs)_{i0} \in \BZ[\zeta_n]$, and the set of distinct prime
factors of $\ord(T)$ coincides with that of the integer 
$\mathrm{norm}(1/\rho(\fs)_{00})$.
%\footnote{Here $\mathrm{norm}(c)$ be the product of the distinct Galois conjugates of the cyclotomic number $c$.}  
\item 
%Galois conjugation on the MD representation $\rho$ has a nice simple form.
Let $\s\in\Gal(\BQ_n)$ be a Galois automorphism. Then (cf. \eqref{DSLn}) 
\begin{align}
 D_{\rho}(\s)_{ij} = \e_\s(i) \delta_{\hs(i),j}, 
\end{align}
where $\e_\s(i) \in \{1,-1\}$ and $\hs$ is a permutation on $\{0, \dots, r-1\}$ determined by
\begin{equation}\label{eq:perm}
  \s \left(\frac{\rho(\fs)_{ij}}{\rho(\fs)_{0j}}\right) =  \frac{\rho(\fs)_{i\hs(j)}}{\rho(\fs)_{0\hs(j)}}\,.  
\end{equation}
Moreover,
\begin{equation} \label{siDrho1}
\s(\rho(\fs)) = D_{\rho}(\s) \rho(\fs) =\rho(\fs)D_{\rho}^\top(\s)\quad \text{and}\quad
\s^2(\rho(\ft)) = D_{\rho}(\s) \rho(\ft) D_{\rho}^\top(\s)\,.
\end{equation}
%This condition plays an important role in our approach.

\item The matrix $\rho(\fs)$ satisfies the Verlinde formula (cf. \cite{V8860}):
\begin{align} 
\label{Ver} 
 N^{ij}_k = \sum_{l=0}^{r-1} \frac{
\rho(\fs)_{li} \rho(\fs)_{lj} \rho(\fs)_{lk}^*}{ \rho(\fs)_{l0} } , \quad i,j=0,1,\ldots,r-1\,.
\end{align}
\item 
\label{FScnd1}
For $m \in \BN_+$, the $m^\text{th}$ Frobenius-Schur indicator of the $i$-th
simple object can also be expressed in terms of $\rho(\fs)$ and $\rho(\ft)$:
\begin{equation}
 \label{nunFS1}
 \nu_m(i)= \sum_{j, k} N_i^{jk}
\rho(\fs)_{j0}\rho(\ft)_{jj}^m \cdot (\rho(\fs)_{k0}\rho(\ft)_{kk}^m)^*\,.
\end{equation}
%\red{XG: The above expression needs to updated for complex $d_i$.}
%\begin{equation}
% \label{nunFS}
% \nu_n(i)= \frac{1}{D^2}\sum_{j, k} N_i^{jk}
%d_{j}\theta_j^n (d_{k}\theta_k^n)^*
%\end{equation}
%\item { If $\rholpha$ decomposes as $\rho_1\oplus\rho_2$ then $\spec(\rho_1(\ft))\cap \spec(\rho_2(\ft))\neq \emptyset$} (the $\ft$-spectrum criteria, cf. \cite{CRank5}).

\end{enumerate}
\end{thm}

\begin{remark} \label{r:pMD}
{\rm It is worth noting that a pMD representation $\rho_\pMD$ shares arithmetic properties with MD representations as $\rho= V \rho_\pMD V$ is an MD representation for some signed diagonal matrix $V$. Therefore, Theorem \ref{p:MD1} (1) and (3-6) also hold for any pMD representation. In particular, for $\s \in \GQ$,  $D_{\rho_{\pMD}}(\s) =V D_\rho(\s)V$, and so 
$$
\s(\rho_{\pMD}(\fs)_{ij}) = \e'_\s(i) \rho_{\pMD}(\fs)_{\hs(i)j} =\e'_\s(j)\rho_\pMD(\fs)_{i\hs(j)}
$$
but the sign function $\e'_\s$ is different from $\e_\s$ in   Theorem \ref{p:MD1} (6) in general.
}
\end{remark}

\subsection{Modular data representations and our classification strategy}
\label{strat}

The MD representation introduced in Definition \ref{d:MDrep} plays an important role in our approach. We now explain the strategy of a systematic construction for low
rank modular data, implementable on a computer.  In Section \ref{s:hand approach} we provide a largely by-hand approach to the classification of rank $6$ MD.

For a given rank, we first construct all the inequivalent $\SL$ representations
$\rho_\mathrm{isum}$ of finite levels, as direct sums of
irreducible $\SL$ representations obtained as tensor products of the prime-power level representations listed in Appendix \ref{repPP}.  Each of these $\SL$ irreducible representations is symmetric, and so is $\rho_{\isum}$.

Although the number of the $\SL$ representations $\rho_\mathrm{isum}$ is
finite, most of these representations are not associated to any MTC.  In next
section, we introduce and collect conditions on MD representations, to reject as much as possible the $\SL$ representations that are not associated to
MTCs.

After we obtain a short list of candidate $\SL$ representations
$\rho_\mathrm{isum}$, we permute the indices using a permutation matrix $P$
\begin{align}
\td\rho = P\rho_\mathrm{isum}P^\top
\end{align}
such that $\arg( \td\rho(\ft)_{ii})$ is ordered for computer implementation or mathematical deduction. 

Suppose $\tilde \rho$ is equivalent to an MD representation $\rho$. Without losing generality, we can further assume $\rho(\ft)=\tilde \rho(\ft)$. It follows from Theorem \ref{t:ortho_eqv} there exists an orthogonal matrix $U$ such that $\rho(\fs) = U \tilde \rho(\fs) U^\top$ and  
$\rho(\ft) = U  \tilde\rho(\ft) U^\top$.
In this case, $U$ is a block-diagonal orthogonal matrix. The size of each block $U_i$ is equal to the multiplicity of the eigenvalue $\tilde\rho(\ft)_{ii}$. We first assume that each of these blocks is of determinant 1. Then
\begin{align}\label{eq:MD}
 \rho_\mathrm{pMD} =  U \td\rho U^\top
\end{align}
is a pseudo-MD representation. Using Theorem \ref{p:MD1}, Remark \ref{r:pMD} and the conditions established in the next section, the existence of such $U$ could either imply contradiction or be determined for all the rank 6 modular data. In the former case, representation $\rho_{\isum}$ will be rejected. Once the matrix $U$ is determined, one can determine the correct signed diagonal matrix by using the Frobenius-Perron dimensions or the Verlinde formula. 

The eigenvectors of the diagonal matrix $\tilde\rho(\ft)$  corresponding to the eigenvalues of multiplicity 1 are of particular importance in the determination of the orthogonal matrix $U$. We simply called the block of $\tilde\rho(\fs)$ corresponding to these eigenvectors the \emph{non-degenerate block}, and denoted by $\tilde\rho(\fs)^\nd$. The following proposition provide a convenient sufficient condition for any $\SL$ representation equivalent to an MD representation.

\begin{prop}\label{p:ndeg}
Let $\tilde\rho$ be any (symmetric) $\SL$ representation. If $\tilde\rho$ is equivalent to an MD representation, then  the entries of $\tilde \rho(\fs)^\nd$ are cyclotomic numbers in $\BQ_{\ord(\tilde \rho)}$.
\end{prop}
\begin{proof}
The statement is an immediate consequence of Theorem \ref{t:ortho_eqv} and Theorem \ref{p:MD1}(1).
\end{proof}

 The proposition can be implemented for computer automation to eliminate many $\rho_{\isum}$.  Theorem \ref{p:MD1} (6) and the property of second Frobenius-Schur indicators are implemented to eliminate $\rho_{\isum}$ or solving the matrix $U$. When the matrix $U$ is determined, the signed diagonal matrix $P_{\sgn}$ can be searched by using the nonnegative integral fusion coefficients (cf. Theorem \ref{p:MD1} (7)). The potential MD representation $\rho_{\MD}$ is then given by
\begin{align}
 \rho_{\MD} =  P_{\sgn} \rho_{\pMD} P_{\sgn}^\top,
\end{align}
Again, $\rho_{\isum}$ will be rejected if no such $P_{\sgn}$ is found.
From the potential MD representations $\rho_\MD$ we can then obtain the
potential modular data $(S,T)$ via \eqref{STrho1}, and they will be verified if  Theorems  \ref{p:MD} and \ref{p:MD1} are satisfied. This
allows us to get a list of $(S,T)$ pairs that include all the modular data. The computer automation for the endeavor is robust particularly when $\rho_{\isum} = \rho_{\isum}^\nd$.

By
comparing the list of $(S,T)$ pairs to known rank-6 MTCs, we obtain a
classification of all modular data via matrix representations of $\SL$.

\subsection{More general properties of $\SL$ representations}

In this subsection, we introduce and collect conditions on $\SL$
representations necessary for them to be MD representations

The decomposition criteria on $\ft$-spectrum \cite{CRank5} of a linear representation of $\SL$ associated with a MTC is one of the major tools. 
\begin{thm} [$\ft$-spectrum  criteria] \label{p:partition}
Let $\rho$ be an MD representation.  If
$$
\rho \cong \rho_1 \oplus \rho_2
$$ 
for some representations $\rho_1, \rho_2$ of $\SL$, then $\spec(\rho_1(\ft))\cap \spec(\rho_2(\ft)) \ne \emptyset$.
\end{thm}

Let $p$ be a prime. We denote by $G_p$ the Galois group $\Gal(\BQ_p)$. The
least dimension of an irreducible representation of $\SL$ of level $p$ is $\frac{p-1}{2}$. Their
$\ft$-spectrum is either $G_p^2 \cdot \zeta_p$ or $G_p^2 \cdot \zeta_p^a$ where $x^2 \equiv a \mod p$ has no integer solution. Note that an integer
$a$ is called a \emph{nonresidue} modulo $p$ if $x^2 \equiv a \mod p$ has no
integral solution. The second least dimension irreducible representation $\rho$
of $\SL$ of level $p$ is $\frac{p+1}{2}$ whose $\ft$-spectrum is either $G_p^2
\cdot e^{2\pi i/p} \cup \{1\}$ or $G_p^2 \cdot e^{2\pi i a/p}\cup \{1\}$ where
$a$ is any nonresidue modulo $p$. In this case, $\rho(\fs)^2 =
\jacobi{-1}{p}\id$ (see for example \cite{Hum}).
\begin{prop}
Let $3< p < q$ be prime such that $pq \equiv 3 \mod 4$. For any modular tensor category $\CC$ such that $\ord(T) = pq$, then 
$rank(\CC) \ne \frac{p+q}{2}+1$. Moreover, if $p > 5$, $rank(\CC) > \frac{p+q}{2}+1$.
\end{prop}
\begin{proof}
Let $\CC$ be a modular tensor category of rank $r \le \frac{p+q}{2}+1$ and $\ord(T) = pq$. There exists an $\SL$ representation $\rho$  of $\CC$ with level $pq$ \cite{DLN}. Suppose $\rho$ has an irreducible subrepresentation $\rho'$ of level $pq$. By the Chinese remainder theorem, the $\rho' \cong \rho_1 \o \rho_2$, where $\rho_1, \rho_2$ are  irreducible representations of $\SL$ of levels $p$ and $q$ respectively. Then 
$$
\frac{p+q}{2}+1 \ge  \deg \rho'  = (\deg \rho_1) (\deg \rho_2) \ge \left(\frac{p-1}{2}\right)\left( \frac{q-1}{2}\right).
$$
The inequality implies  $p=5$ and $q =7$,  and hence $\dim \rho'  = 6$. Therefore, the $\ft$-spectrum of $\rho'$ consists of 6 distinct primitive $35$-th roots of unity, and $rank(\CC) = 6$ or $7$. There exists a modular tensor category of rank 6 with $\ord(T)=35$. However, if $rank(\CC)=7$, then $\rho \cong \rho'\oplus \rho_0$ where $\rho_0$ is trivial. This is not possible by Proposition \ref{p:partition}. In conclusion, if $\rho$ has an irreducible subrepresentation of level $pq$, then $p=5$, $q=7$ and $rank(\CC)=6$. 

Now, we assume $\rho$ has no  irreducible subrepresentation of level $pq$. Then $\rho$ must have irreducible subrepresentations  $\rho_1, \rho_2$ of levels $p$ and $q$ respectively. If $\deg \rho_1 < \frac{p+1}{2}$ or $\deg \rho_2< \frac{q+1}{2}$, then
$$
\rho \cong \rho_1 \oplus \rho_2 \oplus \rho_3
$$
where $\rho_3$ is a subrepresentation of $\rho$ of dimension $\le 2$. If $\rho_3$ has a 1-dimensional component $\rho_4$, then $\rho_4$ must be a trivial. However,  this contradicts Proposition \ref{p:partition}.  Note that irreducible $\SL$ representation of dimension 2 at prime levels only appear for the primes 2, 3 and 5. Therefore,  if $\rho_3$ is irreducible of dimension 2, then $p=5$ and $\rho_3$ is of level 5,  but this contradicts Proposition \ref{p:partition} again. Thus, $\deg  \rho_1 \ge \frac{p+1}{2}$ and $\deg  \rho_2 \ge \frac{q+1}{2}$. Since $rank(\CC) \le  \frac{p+q}{2}+1$, we find $rank(\CC) = \frac{p+q}{2}+1$,  $\deg  \rho_1 = \frac{p+1}{2}$ and $\deg  \rho_2 =\frac{q+1}{2}$.  Now, we would like to show that this also impossible. 

Without loss of generality, we may assume $\jacobi{-1}{p}=1$ and  $\jacobi{-1}{q}=-1$. Then $\rho(\fs)^2$ is a signed diagonal matrix and the multiplicities  $1, -1$ are respectively    $\frac{p+1}{2}$ and  $\frac{q+1}{2}$.  Thus, $|\Tr(\rho(\fs)^2)| =\frac{q-p}{2}$, which  means $\CC$ has   $\frac{p+1}{2} \ge 3$ pairs of simple objects which are  not self-dual. Since $\rho(\ft)$ has only one eigenvalue of multiplicity 2 and all other eigenvalues are of multiplicity 1, $\CC$ has at most 1 pair of simple objects which dual of each other, a contradiction!
\end{proof}

%\red{

%For any positive integer $n$ and $\s \in \Gal(\BQ_n) \cong \BZ_n^\times$, 
%\begin{align}
%\label{DSLn}
% D_{\SLn{n}}(\s) :=
%\ft^a \fs \ft^b \fs \ft^a \fs^{-1}  \mod  n
%=
%\begin{pmatrix}
% a & 0\\
% 0 & b\\
%\end{pmatrix}
%.
%\end{align}
%where $a,b$ are given by $\s( \ee^{ 2 \pi\ii/\ord(\rho_\a(\ft))} )=\ee^{ a 2
%\pi\ii/\ord(\rho_\a(\ft))}$ and $ab \equiv 1 \mod n$.
%So $D_{\mathrm{SL}(2,\BZ_n)}$ form a faithful representation of $\Gal(\BQ_n)$.

%Let $\tilde\rho$ be an $r$-dimensional representation of $\SL$ with $\ord(\tilde\rho(\ft))=n < \infty$. For any integer $a$ coprime to $n$, we define
%$$
%D_{\tilde\rho}(a) := \tilde\rho(\ft^a \fs \ft^b \fs \ft^a \fs\inv), \text{ where } ab \equiv 1 \mod n\,. 
%$$
%For any $\s \in \GQ$, $\s(\zeta_n) = \zeta_n^a$ for some unique integer $a$ modulo n. We define
%$$
%D_{\tilde\rho}(\s) := D_{\tilde\rho}(a)\,. 
%$$
%It is worth noting that if $\tilde\rho$ is a level $n$ representation of $\SL$, then $D_{\tilde\rho}:(\BZ/n\BZ)^\times \to \GL_r(\BC)$ is a representation equivalent to the restriction of $\tilde{\rho}$ on the diagonal subgroup of $\qsl{n}$. Therefore, the composition $\GQ \xrightarrow{\res} (\BZ/n\BZ)^\times \xrightarrow{D_{\tilde\rho}}\GL_r(\BC)$, denoted by the abused notation $D_{\tilde\rho}$, defines a representation of $\GQ.$
%}

Let $\rho$ be an $\SL{}$ representation of a modular tensor category $\CC$ and let $n$ be the level of $\rho$. For any  $\s \in \GQ$,  $D_\rho(\s)$ defined in \eqref{eq:Drho}
is a signed permutation matrix of $\hs$ by \cite[Theorem II]{DLN} (or Theorem \ref{p:MD1} (6)). The permutation $\hs$ on $\irr(\CC)$  is given by \eqref{eq:perm}, and we set
$$
\Inv_\CC(\s):= \{ i \in \irr(\CC) \mid \hs(i)=i\}\,.
$$
If $\gamma$ is complex conjugation, by  \eqref{siDrho1},
$$
D_\rho(\gamma) =\ol{\rho(\fs)} \rho(\fs)\inv =\rho(\fs)^2 = \pm C,
$$
where $C$ is the charge conjugation matrix of $\CC$. Since $\hat{\gamma}(i)=i^*$ for $i \in \irr(\CC)$,
$$
|\Tr(D_\rho(\gamma))| = |\Tr(\rho(\fs^2))| =  \Tr(C) = |\{ i \in \irr(\CC) \mid i^*=i\} = |\Inv(\gamma)|\,.
$$ 
This equality can be generalized to any $\sigma \in \GQ$ as an inequality in the following proposition.
\begin{prop} \label{p:fixpt}
Let $\rho$ be an $\SL{}$ representation of a modular tensor category  $\CC$. For any  $\s \in \GQ$,
$$
|\Tr(D_\rho(\s))| \le |\Inv_\CC(\s)|\,.
$$
Let $s:=\rho(\fs)$, and follow the notation of Theorem \ref{p:MD1}(6). If $s_{ij} \ne 0$ for any $i,j \in \Inv_\CC(\s)$, then $\e_\s(i)= \e_\s(j)$. If there exists  $i \in \Inv_\CC(\s)$ such that  $s_{ij} \ne 0$ for all $j \in \Inv_\CC(\s)$, then
$$
|\Tr(D_\rho(\s))| = |\Inv_\CC(\s)|\,.
$$
In particular, 
$$
\Tr(s^2) = |\{i \in \irr(\CC)\mid i^*=i\}| > 0\,.
$$
\end{prop}
\begin{proof}
By Theorem \ref{p:MD1}(6),  $D_\rho(\s) = \es(i)\delta_{\hs(i),j}$. Therefore, 
$$
|\Tr(D_\rho(\s))| = \left|\sum_{i \in \Inv_\CC(\s)} \es(i)\right| \le |\Inv_\CC(\s)|\,.
$$
If $s_{ij} \ne 0$ for any $i, j \in \Inv_\CC(\s)$, then $\s(s_{ij}) = \es(i) s_{ij} = \es(j) s_{ij}$, and so $\es(i)=\es(j)$. Thus, if there exists $i \in \Inv_\CC(\s)$ such that $s_{ij} \ne 0$ for all $j \in \Inv_\CC(\s)$, then $\es(i)=\es(j)$ for all $j$ and hence the equality
$$
|\Tr(D_\rho(\s))| =  |\Inv_\CC(\s)|\,.
$$
The last statement was proved in the preceding remark and since $\1^*=\1$ this completes the proof of the proposition.
\end{proof}

According to \cite{DLN}, if $\rho$ is an MD representation of an  integral modular tensor category $\CC$, then $\rho(\ft)_{\1, \1} = \zeta$ for some 24-th root of unity $\zeta$ under the identification of the standard basis for $\rho$ and $\irr(\CC)$. The following proposition provides a sufficient condition on the representation type of $\rho$ for $\CC$ to be integral.

\begin{prop}\label{p:5.2}
Let $\tilde \rho$ be any $\SL$ representation. For any $\zeta \in \spec(\td\rho(\ft))$,  denote by $E_{\zeta}(\tilde\rho)$ the eigenspace of $\tilde \rho(\ft)$ for the eigenvalue $\zeta$. Suppose $\tilde\rho$ is equivalent to an MD representation $\rho$ of a modular tensor category $\CC$. Then
\begin{enumerate}
    \item $D_{\tilde\rho}(\s)(E_{\zeta}(\tilde\rho)) \subseteq E_{\zeta}(\tilde\rho)$ for all $\s \in \GQ$ if and only if $\zeta^{24}=1$. 
    \item  If $\1 \in E_{\zeta}(\rho)$ for some $\zeta \in \spec(\rho(\ft))$, and  for each $\s \in \GQ$, there exists $\e_\s = \pm 1$ such that 
$$
D_{\tilde\rho}(\s)|_{E_\zeta(\td\rho)} = \e_\s\, \id_{E_\zeta(\td\rho)}\,,
$$
then $\CC$ is integral. In particular, $\zeta^{24}=1$.
\item If $\1 \in \bigoplus_{\g \in A} E_{\g}(\rho)$ for some subset $A \subseteq \spec(\rho(\ft))$, and
 for any  $\g \in A$, $\s \in \GQ$, there exists $\e_\s(\g) = \pm 1$ such that $D_{\tilde\rho}(\s)|_{E_{\g}(\td\rho)}= \e_\s(\g) \id_{E_{\g}(\td\rho)}$, then $A$ is a set of 24-th roots of unity and $\CC$ is integral. 
\end{enumerate}

\end{prop}
\begin{proof}
Assuming the identification of the standard basis for $\rho$ and $\irr(\CC)$, $E_{\zeta}(\rho)$ is spanned by the objects $X \in \irr(\CC)$ such that $\rho(\ft)_{X,X} =\zeta$. Let $\s \in \GQ$.
It follows from Theorem \ref{p:MD1}(6) that $D_\rho(\s)$ is a signed permutation matrix of a permutation $\hs$ on $\irr(\CC)$, and that $\s^2(\rho(\ft))=D_\rho(\s)\rho(\ft)D_\rho(\s)^{-1}$, or equivalently $\rho(\ft) D_\rho(\s) =  D_\rho(\s) \s^{-2}(\rho(\ft))$.  If $\zeta^{24}=1$, then $\s^2(\zeta) = \zeta$ for all $\s \in \GQ$.
Thus, for any simple object $X \in E_{\zeta}(\rho)$, 
$$
\rho(\ft) D_\rho(\s)(X) =  D_\rho(\s) \s^{-2}(\rho(\ft)) (X) = \s^{-2}(\zeta) D_\rho(\s)(X) = \zeta  D_\rho(\s)(X)\,.
$$
Therefore, $D_\rho(\s) (E_{\zeta}(\rho)) \subseteq E_{\zeta}(\rho)$. Let $\phi: \tilde\rho \to \rho$ be an isomorphism of $\SL$ representations. Then $\phi(E_{\zeta}(\tilde\rho) )= E_{\zeta}(\rho)$, and $\phi D_{\tilde\rho}(\s) =  D_\rho(\s)\phi$ for any $\s \in \GQ$. This implies $D_{\tilde\rho}(\s) (E_{\zeta}(\tilde\rho)) \subseteq E_{\zeta}(\tilde\rho)$. 

Conversely, if $D_{\td\rho}(\s) (E_{\zeta}(\td\rho)) \subseteq E_{\zeta}(\td\rho)$, then $D_\rho(\s) (E_{\zeta}(\rho)) \subseteq E_{\zeta}(\rho)$ by the same reason. Thus, for any $X \in E_\zeta(\rho)$, $\rho(\ft) D_\rho(\s)(X) =  \zeta D_\rho(\s)(X)$. However, we also have 
$$
\rho(\ft) D_\rho(\s)(X) = D_\rho(\s) \s^{-2}(\rho(\ft)) (X) = \s^{-2}(\zeta)  D_\rho(\s)(X).
$$
Therefore, $\s^{-2}(\zeta) = \zeta$ for all $\s \in \GQ$. This implies $\zeta$ is a 24-th root (cf. \cite[Prop. 6.7 and Lem. A.2]{DLN}). This proves statement (1).

For statement (2), we assume $\1 \in E_{\zeta}(\rho)$, and for each $\s \in \GQ$  there exists $\e_\s = \pm 1$ such that $D_{\tilde\rho}(\s)|_{E_{\zeta}(\tilde\rho)} = \e_\s \id_{E_{\zeta}(\tilde\rho)}$. It follows from (1) that $\zeta^{24}=1$. Moreover,
$D_{\rho}(\s)|_{E_{\zeta}(\rho)} = \e_\s \id_{E_{\zeta}(\rho)}$ and  
$D_\rho(\s)(\1) = \e_\s \1 = \pm \hs(\1)$.
Therefore,  $\hs(\1)=\1$, and hence $\s(\dim(V)) = \dim(V)$ for any $V \in \irr(\CC)$ by Theorem \ref{p:MD1}(6). Thus, $\dim(V)$ are integers for $V \in \irr(\CC)$. It follows from \cite[Rem. 6.3]{DLN} that $\FPdim(V) \in \BZ$, and hence $\CC$ is integral.

Statement (3) follows directly from (2), and this completes the proof of the proposition.
\end{proof}

The following result in \cite{BR} is important for determining whether an $\SL$ representation of small level
is an MD representation.
\begin{thm}\label{t:2346}
Modular tensor categories with $\ord(T) = 2,3,4,6$ are  integral.
\end{thm}
Then case for $\ord(T)=2$ is completely classified in \cite{WW}, and the types of these MTCs are given in the following proposition.
\begin{prop}
Let $\CC$ be a modular tensor category with $\ord(T)=2$. Then $\rank(\CC)=4^n$ for some positive integer $n$, and every $\SL$ representation $\rho$ of $\CC$ is projectively equivalent to
$$
(\rho_2 \oplus 2\chi_0)^{\otimes n} \cong a_n \rho_2 \oplus b_n \rho_1 \oplus c_n \chi_0\,,
$$
where $\rho_1, \rho_2$ are respectively the level 2 irreducible representations of dimension 1 and 2, and
$$
a_n=\frac{4^n-1}{3}, \quad b_n=\frac{2\cdot 4^{n-1}+1}{3} - 2^{n-1}, \quad c_n = \frac{2\cdot 4^{n-1}+1}{3} + 2^{n-1}\,.
$$
\end{prop}
\begin{proof}
By \cite{WW}, $\CC$ is a Deligne product of the pointed modular tensor categories $\CC(\BZ_2^2, q)$ and $\CC(\BZ_2^2, q')$ with the quadratic forms $q, q': \BZ_2^2 \to \{\pm 1\}$ given by 
$$
q(x,y)=(-1)^{xy}, \quad  q'(x,y)=(-1)^{x^2+xy+y^2}\,.
$$
Both modular tensor categories,  up to a linear character,  have a representation of $\SL$ equivalent to
$
\rho_2 \oplus 2 \chi_0\,.
$
Thus, $\rho \cong (\rho_2 \oplus 2 \chi_0)^{\otimes n} = a_n \rho_2 \oplus b_n \rho_1 \oplus c_n \chi_0$.
Since the character of $\rho_2 \oplus 2 \chi_0$ satisfies the polynomial $(x-1)(x-2)(x-4)$. Then the formulae for the coefficients $a_n$, $b_n$ and $c_n$ follow by induction on $n$.
\end{proof}

The following proposition follows immediately from the classification of \cite{CWeakly}, where strictly weakly integral means $\FPdim(\CC)\in\Z$ while $\FPdim(X)\not\in\Z$ for some object $X$.
\begin{prop}\label{p:weakly_int}
Let $\CC$ be a modular tensor category of rank 6. 
\begin{enumerate}
\item[(1)] If $\CC$ is integral, then $\CC$ is pointed and hence $\CC$ is of type $(4, 2)$ and every $\SL$ representation of $\CC$ has level 24.
\item[(2)] If $\CC$ is strictly weakly integral, then $\CC$ is braided equivalent to a Galois conjugate of $U(2)_1 \boxtimes SU(2)_3$, $SO(5)_2$ or its zesting. If $\rho$ is an $\SL$ representations of $\CC$ with a minimal $\ft$-spectrum, then one of the following holds: {\rm (i)} $\CC$ is of type $(6)$ and $\rho$ has level 16, {\rm (ii)} $\CC$ is of type $(3,3)$ and $\rho$ has level 20, or {\rm (iii)} $\CC$ is of type $(3,2,1)$ and $\rho$ has level 10\,.
\item[(3)] In particular, if $\CC$ is weakly integral, then $\dim(\CC) = 6, 8, 20$.
\end{enumerate}
\end{prop}

When a potential modular data is obtained from a representation of $\SL$, one could obtain the $\FPdim(X)$ of each simple object $X$. Those simple objects $X$ with $\FPdim(X)=1$ generate a pointed ribbon subcategory. The next proposition, which should be well-known to experts, describes some relations between the rank of a pointed   
ribbon category and the orders of the twists. 
\begin{prop}
Let $\CC$ be a pointed ribbon category of rank $n$. Then $\ord(T) \mid n$ if $n$ is odd, and $\ord(T) \mid 2n$ if $n$ is even. If, in addition, $\CC$ is symmetric and $\dim(a)>0$ for all $a \in \irr(\CC)$, then either $\ord(T)=1$ or $2$. In the latter case, $n$ must be even and there are exactly $n/2$ simple objects with twist $-1$.
\end{prop}
\begin{proof}
Since $\CC$ pointed, the set $G=\irr(\CC)$ forms an abelian group under the tensor product and the map $q:G \to \BC^\times, q(a) = \theta_a$ defines a quadratic form on $G$. Therefore, $B_q(a,b)=\frac{q(ab)}{q(a)q(b)}$ defines a bicharacter on $G$. In particular, $B_q(a,b)$ is an $n$-th root of unity for any $a,b \in G$. Now, for any positive integer $m$ and $a \in G$, we have
$$
q(a^m) = q(a)q(a^{m-1}) B_q(a, a^{m-1}) = q(a)q(a^{m-1}) B_q(a, a)^{m-1}\,.
$$
Therefore, by induction, we have
$$
q(a^m) = q(a)^m B_q(a,a)^{m(m-1)/2}\,.
$$
In particular, $q(a)^n = B_q(a,a)^{-n(n-1)/2}$. If $n$ is odd, $\frac{n-1}{2} \in \BZ$ and so $q(a)^n=1$. If $n$ is even, then $q(a)^{2n} = 1$. This completes the proof of the first statement.

If, in addition, $\CC$ is symmetric and $\dim(a)=1$ for $a \in G$, then
$$
1 =S_{a,b} = B_q(a^{-1}, b) = B_q(a, b)^{-1} = \frac{q(a)q(b)}{q(ab)}
$$
for any $a, b \in G$. Therefore, $q$ is a character of $G$. Since $q(a^{-1})=q(a)$, $q(a)^2=1$ for all $a \in G$. If $q$ is of order 1, then $q(a)=1$ for all $a \in G$ or $T=\id$. However, if $q$ is of order 2, then the image of $q$ is the group $\{\pm 1\}$ which is of order 2. Therefore, $\ker q$ is of index 2 which means there are exactly $n/2$ simple objects in $G$ with twists are 1. Thus, the second statement follows.
\end{proof}

It is worth noting that last statement of the preceding proposition does not hold for super-Tannakian fusion categories which are not pointed. For example, if we take $Q$ to be the quaternion group of order 8 and $z$ the unique central element of order 2, then the super-Tannakian fusion categories $\Rep(Q,z)$ has 4 simple objects $a$ of dimension 1 with $\theta_a=1$ and a unique simple object $b$ of dimension 2 with $\theta_b= -1$.

For any legitimate fusion rules $N_{ij}^k$, one could obtain the possible $\theta_k = \ee^{2 \pi \ii s_k}$ by solving  a system of linear equations with unknowns $s_k$. The following proposition provides a condition for legitimate $s_k$ of a potential modular data. 
\begin{prop}\label{p:sum_of_twists}
Let $\CC$ be a modular tensor category of rank $n$ and  central charge $c$. If the twists of $\CC$ are $\ee^{2 \pi \ii s_1}, \dots, \ee^{2 \pi \ii s_n}$ for some rational numbers $s_1, s_2, \dots, s_n$, then
$$
12\sum_{k=1}^n s_k - nc/2 \in \BZ
$$
\end{prop}
\begin{proof}
Note that $\ee^{ \pi \ii c/4} = \frac{1}{D} \sum_{k=1}^n d_k^2\, \ee^{2 \pi \ii s_k}$ where $d_k$ denotes the dimension of the simple object $k$ with twist $\ee^{2 \pi \ii s_k}$ and $D = \sqrt{\dim(\CC)}$. Let $(S,T)$ be the modular data of $\CC$. Then
$$
\rho(\fs) = \frac{1}{D} S, \quad \rho(\ft) = \ee^{-2 \pi \ii c/24}\,T
$$
defines an $\SL$ representation of $\CC$. Thus, $\det \circ \rho$ is a 1-dimensional representation of $\SL{}$. Since the group of linear characters of $\SL{}$ is a cyclic group of order 12, $\det\rho(\fg)^{12} = 1$ for all $\fg \in \SL{}$. In particular,
$$
1= \det\rho(\ft)^{12} = (\ee^{2 \pi \ii s_1} \cdots \ee^{2 \pi \ii s_n} \cdot \ee^{-2 \pi \ii n c/24})^{12}= \ee^{2 \pi \ii\cdot  12\sum_{k=1}^n s_k - nc/2}\,.
$$
This implies  $12\sum_{k=1}^n s_k - nc/2 \in \BZ$. 
\end{proof}

The following proposition is proved in \cite{PSYZ}  will also  be useful later.
\begin{prop} \label{p:multiple}
Let $\rho$ an MD linear representation. Then 
$$
\rho \not\cong n \rho_0
$$ 
for any integer $n >1$ and any non-degenerate  representation $\rho_0$ of $\SL$.
\end{prop}

\subsection{Modular tensor categories of type $(d,1,\ldots,1)$} 
A representation $\rho_\isum$ of $\SL$  with such a irreducible decomposition dimensions is, in general, more difficult to determine whether it is equivalent to an MD representation. However, this type of MTC exists. It is desirable to deduce some conditions for such MD representations.

\begin{lem} \label{l:2}
Let  $\rho$ an MD representation. If $\rho \cong \rho_0 \oplus \rho_1 \oplus \dots \oplus \rho_\ell$ for some 1-dimensional representations $\rho_1, \dots, \rho_\ell$ of $\SL$, then $\spec(\rho_i(\ft)) \subset \spec(\rho_0(\ft))$ for all $i > 0$. In particular, 
if $\rho_0(\ft)$ has  exactly one eigenvalue which is a 12-th root of unity, then $\rho_1, \ldots, \rho_\ell$ are all equivalent, and $\rho\cong \rho_0 \oplus \ell \rho_1$. 
\end{lem}
\begin{proof} By the $\ft$-spectrum criteria, $\spec(\rho_j(\ft)) \subset\spec(\rho_0(\ft))$ for some $j > 0$.  Suppose there exists $j>0$  such that $\spec(\rho_j(\ft)) \not\subset\spec(\rho_0(\ft))$. Let  $J = \{j \in \{0,\dots, \ell\}\mid \spec(\rho_j(\ft))  \not\subset \spec(\rho_0(\ft))\}$. Then,  the decomposition 
$$
\rho\cong \left( \sum_{j \in J} \rho_j\right) \oplus \left( \sum_{j \not\in J} \rho_j\right) 
$$
does not satisfies the $\ft$-spectrum criteria. Therefore, $\spec(\rho_j(\ft)) \subset \spec(\rho_0(\ft))$ for all $j$.

If, in particular, $\spec(\rho(t))$ contains exactly one 12-th root of unity $\zeta$, then $\spec\rho_i(t) = \{\zeta\}$ for all $i >0$. Hence $\rho_1\cong \rho_i$ for $i > 1$, and the last assertion follows.
\end{proof}
\begin{cor}\label{l:3}
Let $\rho$ be an $\SL$ representation of a modular tensor category   $\CC$. Suppose that $\rho \cong \rho_0 \oplus \rho_1 \oplus \dots \oplus \rho_\ell$ for some 1-dimensional representations $\rho_1, \dots, \rho_\ell$ and some non-degenerate  irreducible representation $\rho_0$ of $\SL$ such that $\spec(\rho_0(\ft))$ has a unique 12-th root of unity. Then 
$\CC$  admits an MD representation $\rho' \cong \rho'_0 \oplus \ell \chi_0$, where $\chi_0$ is the trivial representation and $\rho'_0$ is projectively equivalent to $\rho_0$ with $1 \in \spec(\rho'_0(\ft))$.

If $\ell \not\in \{1, 2\deg \rho_0 -1\}$,  then $\CC$ is self-dual,  and $\rho'_0$ is even. If $\ell \in \{1, 2\deg \rho_0 -1\}$ and $\CC$ is not self-dual, then $\rho'_0$ is odd, and the set of non-self-dual objects is given by $\{i \in \irr(\CC) \mid \rho'(\ft)_{ii} =1\}$.
\end{cor}
\begin{proof}
By Lemma \ref{l:3}, $\rho \cong \rho_0 \oplus \ell \rho_1$. Since $\dim \rho_1 =1$, $\rho' = \rho_1^* \o \rho$ is another $\SL$ representation of $\CC$. Moreover, $\rho' \cong \rho'_0 \oplus \ell \chi_0$, where $\rho'_0 = \rho^*_1\o \rho_0$ which is also non-degenerate. 

Suppose $\rho'_0(\fs^2) = -\id$. By Proposition \ref{p:fixpt}, the number of self-dual objects in $\irr(\CC)$    is given by
$$
|\Tr(\rho'(\fs^2))|=|\ell -\deg\rho'_0| > 0
$$
since $\1$ is self-dual simple object. If $\ell > \deg \rho_0$,  then $|\Tr(\rho'(\fs^2))| = \ell-\deg \rho_0$ and so number of non-self-dual objects in $\irr(\CC)$ is $2\deg \rho_0$. 
The non-degeneracy of $\rho'_0$ implies that   $\rho'(\ft)_{ii} =1$ for any non-self-dual $i \in \irr(\CC)$. Therefore, $2\deg \rho_0= \ell+1$ or $\ell = 2\deg \rho_0-1$. 

On the other hand, if $\ell < \deg \rho_0$,  then $|\Tr(\rho'(\fs^2))| = \deg \rho_1-\ell$ and so number of non-self-dual objects in $\irr(\CC)$ is $2 \ell$. Since  $\rho'(\ft)_{ii} =1$ for any non-self-dual simple object $i$,  $2\ell = \ell+1$ or $\ell = 1$. 

Thus, if $\ell \ne 1$ or $2\deg \rho_1-1$, then $\rho'_0(\fs^2) =\id$  and so $\CC$ is self-dual. On the other hand, if $\ell \in \{1, 2\deg \rho_0 -1\}$ and $\CC$ is not self-dual, then $\rho'_0(\fs^2)=-\id$ and the above discussion shows that the non-self-dual objects $i\in \irr(\CC)$ are exactly those $i$ satisfying $\rho'(\ft)_{ii}=1$. 
\end{proof}

Now, we can prove a sufficient condition for any MD representation of prime level $p>3$ and of type $(\frac{p+1}{2}, 1,\dots, 1)$. 

\begin{prop}\label{p:d111} Let $\CC$ be an MTC of type $(d, 1,\dots,1)$ such that $\ord(T)$ is a prime $p>3$, where  $d = \frac{p+1}{2}$. Then  $\CC$ is of type $(d,1)$, and hence $\rank \CC = d+1$.  Moreover, $\Inv_\CC(\s)=\emptyset$ for any generator $\s \in \Gal(\BQ_p/\BQ)$. Furthermore, if $p \equiv 1 \mod 4$, then $\CC$ is self-dual; otherwise $\CC$ is not self-dual.
\end{prop}
\begin{proof}
 By \cite{DLN}, there is an $\SL$ representation $\rho$ of $\CC$, which has level $p$. Then, every subrepresentation  of $\rho$ must have a level dividing $p$.  Since $\CC$ is of type $(d, 1,\dots,1)$, $\rho$ has a irreducible subsrepresentation $\rho_0$ of dimension $d$ and level $p$. By the classification of irreducible representation $\qsl{p}$, $\rho_0(\fs^2) =\jacobi{-1}{p} \id$, $\rho_0$ is non-degenerate and $1$  is the unique 12-th root of unity in $\spec(\rho_0(\ft))$.  By Corollary \ref{l:3},
 $$
 \rho \cong \rho_0 \oplus \ell \chi_0.
 $$
 Thus, if $p \equiv 1 \mod 4$, then $\rho_0$ is even and hence $\CC$ is self-dual. However, if $p \equiv 3 \mod 4$, then $\rho_0$ is odd and so $\CC$ is not self-dual. 
 
One can derive from \cite{NW} that
$$
\rho_0(\fs) = \frac{\jacobi{a}{p}}{\sqrt{p^*}} \left[\begin{array}{c|c}
1 & \sqrt{2}\,\,\,  \cdots\,\,\,  \sqrt{2} \\\hline
\sqrt{2}  &\multirow{3}{*}{$2 \cos\left(\frac{4 \pi a i j}{p}\right)$} \\
\vdots &\\
\sqrt{2}  & 
\end{array}
\right],\quad \rho_0(\ft)=\left[\begin{array}{lll} \zeta_p^{a\cdot 0} &&  \\
&  \ddots & \\
 & & \zeta_p^{a(d-1)^2}
 \end{array}
\right] 
$$
where $1 \le i, j \le d-1$, $p^*=\jacobi{-1}{p} p$, and $a$ an integer coprime to $p$. One may assume  $\rho(\ft) = \diag(1, \dots, 1, \zeta_p^a, \dots,\zeta_p^{a (d-1)^2})$. By Theorem \ref{t:ortho_eqv}, there exists $W \in \O_{d+\ell}(\BR)$ such that $\rho =W (\ell \chi_0 \oplus \rho_0) W^\top$. Note that $W = V U$ for some signed diagonal matrix $V$ and 
$$
U = \left[\begin{array}{c|c}
f & 0\\ \hline
0 & I_{d-1}
\end{array}\right], \quad \text{where } f \in \SO_{\ell+1}(\BR)\,,
$$
and $\rho_\pMD = U (\ell \chi_0 \oplus \rho_0) U^\top$ is a pseud0-MD representation, where $I_{d-1}$ denotes the identity matrix of dimension $d-1$. 

By direct computation, 
\begin{equation}\label{eq:conjugate}
\rho_{\pMD}(\fs) = U \left[\begin{array}{c|c}
I_\ell  & 0 \\\hline
0 & \rho_0(\fs)
\end{array}
\right] U^{\top}  = \left[\begin{array}{c|c}
I_{\ell+1}+f_{*, \ell+1} f_{*, \ell+1} ^{\top} (x-1)   & 
x\sqrt{2}  f_{*, \ell+1} r_{d-1}\\ \hline
x\sqrt{2}\,    r_{d-1}^\top f_{*, \ell+1}^\top & 2 x \cos\left(\frac{4 \pi a i j}{p}\right)
\end{array}\right]
\end{equation}
where $f_{*,\ell+1}=[f_{1, \ell+1} , \cdots,  f_{\ell+1, \ell+1}]^\top$, $r_{d-1} = [1, \cdots, 1] \in \BR^{d-1}$, and $x = \jacobi{a}{p}/\sqrt{p^*}$. 

Let $\s$ be the generator of $\Gal(\BQ_p/\BQ)$.  For any $j \in \{1,\dots, d-1\}$, there exists $\tj \in \{1,\dots, d-1\}$ such that
$$
\s\left(2 \cos(2 \pi j/p)\right) = 2 \cos(2 \pi \tj /p)\,.
$$
Since $\sqrt{p^*} \in \BQ_p$, $\s(\sqrt{p^*})=-\sqrt{p^*}$, and so
$$
\s\left(2 x \cos(4 \pi aij/p)\right) = -2 x \cos(4 \pi a i \tj/p)\,.
$$
for any $i,j \in  \{1,\dots, d-1\}$.
If one identifies $\irr(\CC)$ with $\{1, \dots, d+\ell\}$, then we have  $\hs(\ell+1+j) = \ell+1+\tj$ for each $j \in \{1,\dots, d-1\}$. In particular, $\hs$ has no fixed point in $\{\ell+2, \dots, \ell+d\}$.

 By \eqref{eq:conjugate} and Remark \ref{r:pMD}, 
$$
\s(x\sqrt{2} f_{i, \ell+1})= -x\sqrt{2} f_{i, \ell+1}  \quad \text{for all }i \in \{1,\dots, \ell+1\}\,.
$$
Since $x, x\sqrt{2} f_{i, \ell+1} \in \BQ_p$ and $\s(x)=-x$,  $\sqrt{2} f_{i, \ell+1} \in \BQ_p$  and $\s(\sqrt{2} f_{i, \ell+1}) = \sqrt{2} f_{i, \ell+1}$. 
Therefore,  $\sqrt{2} f_{i, \ell+1} \in \BQ$ for all $i \in \{1,\dots, \ell+1\}$, and hence $f_{i, \ell+1}f_{j, \ell+1} \in \BQ$ for all $i,j \in \{1,\dots, \ell+1\}$.

We claim that $0 < f_{i, \ell+1}^2 < 1$ for all  $i \in \{1, \dots,\ell+1\}$. If $f_{i, \ell+1} = 0$ for some $i$, then each row of $\rho_{pMD}(\fs)$ has a zero entry by \eqref{eq:conjugate}. Therefore, $f_{i, \ell+1}\ne 0$ for all $i$. Since  $f_{*, \ell+1}$ as unit length,  if $f^2_{i, \ell+1} = 1$, then $f_{k, \ell+1}=0$ for all $k \ne i \le \ell+1$, a contradiction.  This proves the claim.

Now we can show that $\Inv_\CC(\s) = \emptyset$. It suffices to show that $\hs$ has no fixed point in $\{1, \dots, \ell+1\}$. Suppose the $i$-th column of $s:=\rho_{pMD}(\fs)$ is fixed by $\hs$ for some $i \in \{1, \dots, \ell+1\}$. Then $\s(s_{ii}) = \e'_\s(i) s_{ii}$, where $\e'_\s(i) = \pm 1$. Since
$s_{ii}=1+f_{i, \ell+1}^2(x-1)$, the preceding equality implies
$$
\e'_\s(i) (1+f_{i, \ell+1}^2(x-1)) = 1+f_{i, \ell+1}^2(-x-1)\,.
$$
Since $f_{i,\ell+1}^2 \le 1$ is rational, the equation forces $f^2_{i,\ell+1}=1, \e'_\s(i)=-1$ or $f^2_{i,\ell+1}=0, \e'_\s(i)=1$. Both are not possible as $0 < f_{i, \ell+1}^2 < 1$.   Therefore, $\hs$ has no fixed point in $\irr(\CC)$.

Let $\s(\zeta_p) = \zeta_p^v$. Then $\Tr(D_{\rho_0}(\s))=\Tr(\rho_0(\ft^v \fs\ft^u \fs\ft^v \fs\inv))=-1$ (cf. \cite{Hum}),  where $uv \equiv 1 \mod p$. It follows from Proposition \ref{p:fixpt}, $|\Inv_\CC(\s)| \ge  |\Tr(D_\rho(\s))| = \ell-1$.  Therefore, $\ell =1$.
\end{proof}
\subsection{MD representations with multiplicities}
In this subsection, we investigate the MD representations $\rho \cong \rho_1 \oplus \rho_2$ such that $\rho_1$, $\rho_2$ are non-degenerate, symmetric, and their $\ft$-spectrums have nonempty intersection.

\begin{thm}\label{t:solution1}
Let $\rho_1$, $\rho_2$ be non-degenerate symmetric representations of $\SL$ such that the intersection of their $\ft$-spectra is of size $l\geq 1$. Let $\dim\rho_1 = l+k$ and $\dim \rho_2=l+m$ and suppose   $k,  m  \ge 1$. Let $\rho_1(\fs)=[\psi_{ij}]$, $\rho_1(\ft)=\diag(\a_1,\dots, \a_{k+l})$, $\rho_2(\fs)=[\eta_{ij}]$ and $\rho_2(\ft)=\diag( \b_1, \cdots, \b_{m+l})$ with $\a_i = \b_i$ for $i=1,\dots, l$. Suppose  $\rho_1 \oplus \rho_2$ is equivalent to an $\SL$ representation $\rho$ of a modular tensor category $\CC$. Then 
\begin{enumerate}
   \item[(i)] there exists a signed diagonal matrix 
    $V$ and $2 \times 2$ orthogonal matrices $U_i =\mtx{a_i & -b_i \\ b_i& a_i }$ with  $a_i \ge 0$ ($i =1, \dots, l$) such that 
    $$
   \rho(\fs) = V\left[
  \begin{array}{c|c|c}
 A & B^{\top} & C^{\top} \\ \hline
 B & \psi' & 0 \\\hline
 C & 0 &\eta'
\end{array}
\right]V  \text{ and } \rho(\ft) = \diag(\a_1 I_2, \dots, \a_l I_2, \a_{l+1},\dots, \a_{l+k}, \b_{l+1},\dots, \b_{l+m} ),
    $$
where $A, B$ and $C$ are block matrices with
$$
A_{ij} =  U_i \mtx{\psi_{ij} & 0 \\ 0 & \eta_{ij}} U_j^{\top}, \quad B_{i'j} = [\psi_{l+i',j} \,\, 0 ] U_j^{\top} \quad \text{and} \quad C_{i''j} =  [0\,\, \eta_{l+i'', j}]   U_j^{\top},
$$    
$1 \le i,j \le l$, $1 \le i' \le k$ and $1 \le i'' \le m$, and $\psi'$, $\eta'$ are respectively the $k \times k$ and the $m \times m$  bottom diagonal blocks of $\rho_1(\fs)$ and $\rho_2(\fs)$, i.e.,
$$
\rho_1(\fs)=\left[\begin{array}{c|c} 
* & * \\ \hline
* & \psi' 
\end{array}\right] 
\quad \text{and}  \quad
\rho_2(\fs)=\left[
\begin{array}{c|c}
* & * \\ \hline
* & \eta'
\end{array}\right]\,.
$$
    \item[(ii)]  Let $(e_1, \dots, e_{2l+m+k})$ be the standard basis for $\rho$ which is identified with  $\irr(\CC)$. Then the unit object $\1$ of $\CC$ is $e_{2u-1}$ or $e_{2u}$ for some $u \le l$ such that 
    \begin{enumerate}
    \item[(a)] $\psi_{uu} + \eta_{uu}\ne 0$ and $\frac{\psi_{uu} - \eta_{uu}}{\psi_{uu} + \eta_{uu}}$ is a real unit in $\BZ[\zeta_N]$  where $N$ is the order of
     $
     T= \a_u^{-1} \rho(\ft) \,;    
     $
    \item[(b)] the list of elements
    $${\small
     \frac{4}{(\psi_{uu}+\eta_{uu})^2},\, 
     \frac{\sqrt{2}\psi_{u, l+1}}{\psi_{uu} + \eta_{uu}},
     \dots, \frac{\sqrt{2}\psi_{u, l+k}}{\psi_{uu} + \eta_{uu}},\,  
     \frac{\sqrt{2}\eta_{u, l+1}}{\psi_{uu} + \eta_{uu}},\dots, 
     \frac{\sqrt{2}\eta_{u, l+m}}{\psi_{uu} + \eta_{uu}}  }
      $$
      are nonzero real  numbers in  $\BZ[\zeta_N]$; 
     \item[(c)]    $\frac{\sqrt{2} \psi_{i, l+i'}}{\psi_{u, l+i'}}, \frac{\sqrt{2} \eta_{i, l+i''}}{\eta_{u, l+i''}}  \in \BZ[\zeta_N] \text{ for } l < i,\,  1 \le i' \le k,\, 1 \le i'' \le  m$\,. 
    \end{enumerate}  
      In this case, $\e_u b_u=  a_u = \frac{1}{\sqrt{2}}$   for some  $\e_u=\pm 1$,  and the modular data of $\CC$ is given by
    \begin{equation}\label{eq:S-mat}
    S = \frac{2}{\psi_{u,u}+\eta_{u,u}} \rho(\fs) \quad \text{and}  \quad T=\a_u^{-1} \rho(\ft) \,.
    \end{equation}
    In particular,
    $
     \dim(\CC)  =  \frac{4}{|\psi_{uu}+\eta_{uu}|^2}$, $\frac{\psi_{uu}-\eta_{uu}}{\psi_{uu}+\eta_{uu}} \in \{\pm\dim(e_{2u-1}), \pm \dim(e_{2u})\}$,
    and the dimensions of  $e_{2l+1}, \dots, e_{2l+k+m}$, up to some signs, are respectively given by
    $$
    \frac{\sqrt{2}\psi_{u, l+1}}{\psi_{uu} + \eta_{uu}},
     \dots, \frac{\sqrt{2}\psi_{u, l+k}}{\psi_{uu} + \eta_{uu}},\,  
     \frac{\sqrt{2}\eta_{u, l+1}}{\psi_{uu} + \eta_{uu}},\dots, 
     \frac{\sqrt{2}\eta_{u, l+m}}{\psi_{uu} + \eta_{uu}} \,.
    $$
 
    \item[(iii)] If $\rho_1$ and $\rho_2$ are irreducible, then $\rho_1$ and $\rho_2$ must have the same parity and $\CC$ is self-dual.
\end{enumerate}
\end{thm}
\begin{proof} We first obtain a representation $\tilde\rho$ by conjugating $\rho_1 \oplus \rho_2$ with a permutation matrix so that
$$
\tilde \rho(\ft) = \diag(\a_1 I_2, \dots, \a_l I_2, \a_{l+1}, \dots, \a_{l+k}, \b_{l+1}, \dots,  \b_{l+k})\text{ and } 
   \rho(\fs) = \left[
  \begin{array}{c|c|c}
 \tilde A & \tilde B^{\top} & \tilde C^{\top} \\ \hline
 \tilde B & \psi' & 0 \\\hline
 \tilde C & 0 &\eta'
\end{array}
\right]
    $$
where $I_2$ is the $2\times 2$ identity matrix, and $\tilde A$, $\tilde B$, $\tilde C$ are block matrices given by
   $$
   \tilde A_{ij} = \mtx{\psi_{ij} & 0 \\ 0 & \eta_{ij}},\quad \tilde B_{i'j} = [\psi_{l+i',j} \,\, 0 ] \quad \text{and} \quad \tilde C_{i''j} =  [0 \,\, \eta_{l+i'', j}]
$$ 
with  $1 \le i,j \le l$, $1 \le i' \le k$ and $1 \le i'' \le m$,
Suppose there exists an MD representation $\rho$ of a modular tensor category $\CC$ such that 
$\rho \cong \rho_1 \oplus \rho_2$. Then $\rho \cong \tilde \rho$ and we may assume $\rho(\ft) = \tilde \rho(\ft)$ by conjugating a permutation matrix to $\rho$.  According to Theorem \ref{t:ortho_eqv}, 
there exists a block diagonal orthogonal matrix $U$ of the form
$$
U = \diag(U_1, \dots, U_l,  \g_{2l+1}, \dots, \g_{2l+m+k}) 
$$
such that $\rho(\fs) = U \tilde \rho(\fs) U^{\top}$ and $\rho(\ft) =  \tilde \rho(\ft) $, 
where $\g_j = \pm 1$ and $U_i$ is a $2\times 2$ orthogonal matrix for $i=1, \dots, l$ and $j=2l+1, \dots, 2l+k+m$. We can always write $U_i = V_i \mtx{a_i & -b_i\\ b_i & a_i}$ where $a_i^2+b_i^2=1$, $a_i \ge 0$ and $V_i$ a signed diagonal matrix. Now, we set $V= \diag(V_1, \dots, V_l, \g_{2l+1}, \dots, \g_{2l+k+m})$. Then statement (i) follows.

The standard basis $(e_1, \dots, e_{2l+k+m})$ is now identified with $\irr(\CC)$. Since only the first $2l$ rows of $\rho(\fs)$ may not contain any zero entries,  the unit object $\1$  can only be $e_x$ with $ 1 \le x \le 2l$. Let $u=\lceil x/2 \rceil$, the least integer $\ge u/2$. Then, 
$$T= \a_u^{-1} \diag(\a_1 I_2, \dots, \a_l I_2, \a_{l+1}, \dots, \a_{k+l}, \b_{l+1}, \dots, \b_{l+m})$$ 
and the $(2u-1)$-th and $2u$-th rows of $\rho(\fs)$ are given by
 
$$
    A_{u,i} = {\small\mtx{a_u a_j \psi_{uj}+ b_u b_j  \eta_{uj} & a_u b_j \psi_{uj}- b_u  a_j \eta_{uj}\\
  b_u a_j \psi_{uj}- a_u b_j  \eta_{uj}  & b_u b_j \psi_{uj}+  a_u a_j \eta_{uj}
    }}, \quad (B^{\top})_{u,i'} = \psi_{u, l+i'}\mtx{a_u \\ b_u},\quad (C^{\top})_{u, i''} = \eta_{u, l+i''}\mtx{-b_u \\ a_u}\,.
    $$
Since $e_x=\1$ and $x \in \{2u-1, 2u\}$, $a_u, b_u, \psi_{u, l+i'}$ and $\eta_{u, l+i''}$ are non-zero for $1 \le i' \le k$ and $1 \le  i'' \le m$.

Now, we assume $x=2u-1$. Then, by \cite{NS10},
$$
\frac{\rho(\fs)_{2u,2l+i'}}{\rho(\fs)_{2u-1,2l+i'}} = \frac{b_u \psi_{u, l+i'}}{a_u \psi_{u, l+i'} a_u} = \frac{b_u}{a_u}\quad \text{and} \quad \frac{\rho(\fs)_{2u,2l+k+i''}}{\rho(\fs)_{2u-1,2l+k+i''}} = \frac{-a_u \eta_{u, l+i''}}{b_u \eta_{u, l+i''}} = \frac{-a_u}{b_u} \in  \BZ[\zeta_N]
$$
where $N = \ord(T)$. Therefore, $\frac{a_u}{b_u}$ is a unit in $\BZ[\zeta_N]$. According to \cite{NWZ}, both $\spec(\rho_1(\ft))$ and   $\spec(\rho_2(\ft))$ are closed under the action of $\s^2$ for any $\s \in \GQ$. Therefore, the subsets 
$$\{\a_{l+1}, \dots, \a_{l+k}\}  \subset \spec(\rho_1(\ft)) \quad\text{and}\quad  \{\b_{l+1}, \dots, \b_{l+m}\}\subset \spec(\rho_2(\ft))$$ are closed  under $\s^2$ for all $\s \in \GQ$. Thus, $\{2l+1,\dots, 2l+k\}$ and $\{2l+k+1,\dots, 2l+k+m\}$ are both closed under the action of $\hs$ for $\s \in \GQ$. In particular, for $\s \in \GQ$, $\hs(2l+1) = 2l+i'$ for some positive integer $i' \le k$. Hence,  
 $$
 \s\left(\frac{b_u}{a_u}\right) = \s\left(\frac{\rho(\fs)_{2u,2l+1}}{\rho(\fs)_{2u-1,2l+1}}\right) = \frac{\rho(\fs)_{2u,\hs(2l+1)}}{\rho(\fs)_{2u-1,\hs(2l+1)}} = 
 \frac{\rho(\fs)_{2u,2l+i'}}{\rho(\fs)_{2u-1,2l+i'}} =  \frac{b_u}{a_u}\,. 
 $$
 So, $b_u/a_u \in \BQ$ and hence $b_u/a_u = \pm 1 = \e_u$. Since $a_u^2+b_u^2=1$, we have $a_u= \frac{1}{\sqrt{2}}$. This implies that 
$$
    A_{u,u} = \frac{1}{2}{\small\mtx{ \psi_{uu}+   \eta_{uu} & \e_u (\psi_{uu}-  \eta_{uu})\\
  \e_u (\psi_{uu}-  \eta_{uu})  &  \psi_{uu}+  \eta_{uu}
    }}, \quad (B^{\top})_{u,i'} = \frac{1}{\sqrt{2}}\psi_{u, l+i'}\mtx{1 \\ \e_u},\quad (C^{\top})_{u, i''} = \frac{1}{\sqrt{2}}\eta_{u, l+i''}\mtx{-\e_u \\ 1}\,.
    $$ 
 In particular, $\frac{\zeta_4^i}{D} = \frac{\psi_{uu}+\eta_{uu}}{2}$. Therefore,  
 $\psi_{uu}+\eta_{uu} \ne 0$ and so the $S$-matrix \eqref{eq:S-mat} of $\CC$ is then obtained. The global dimension 
 $$
 \dim(\CC) =  \frac{\pm 4}{(\psi_{uu}+\eta_{uu})^2} \in \BR^\times  \cap \BZ[\zeta_N]\,.
 $$
 It follows from \eqref{eq:S-mat} and $b_u = \e_u a_u =\frac{\pm 1}{\sqrt{2}}$  that the dimensions of $e_{2u}, e_{2l+1},\dots, e_{2l+k+m}$, up to some signs, are respectively given by
 $$
 \frac{\psi_{uu} - \eta_{uu}}{\psi_{uu}+\eta_{uu}},\,  \frac{\sqrt{2}\psi_{u, l+1}}{\psi_{uu} + \eta_{uu}},
     \dots, \frac{\sqrt{2}\psi_{u, l+k}}{\psi_{uu} + \eta_{uu}},\,  
     \frac{\sqrt{2}\eta_{u, l+1}}{\psi_{uu} + \eta_{uu}},\dots, 
     \frac{\sqrt{2}\eta_{u, l+m}}{\psi_{uu} + \eta_{uu}},\,
 $$
 in particular, they are non-zero real numbers in $\BZ[\zeta_N]$. It follows from \cite{NS10} that $\frac{\rho(\fs)_{y, z}}{ \rho(\fs)_{2u-1, z}} \in \BZ[\zeta_N]$ for any $y,z=1, \dots, 2l+k+m$. For $y=z=2u$, we find
 $$
 \frac{\psi_{uu}+\eta_{uu}}{\psi_{uu}-\eta_{uu}} \in \BZ[\zeta_N], 
 $$
 and so $\frac{\psi_{uu}-\eta_{uu}}{\psi_{uu}+\eta_{uu}}$ is a real unit in $\BZ[\zeta_N]$. For $y, z > 2l$, we find
 $$
 \frac{\sqrt{2}\psi_{i, l+i'}}{\psi_{u, l+i'}} \quad \text{and}\quad 
 \frac{\sqrt{2}\eta_{i, l+i''}}{\eta_{u, l+i''}} \in \BZ[\zeta_N].
 $$
 for $i > l,\, 1\le i'\le k,\,  1\le i''\le m$. This completes the case for $x=2u-1$.

 One can follow the same argument for  the case when $x=2u$. However, the conclusions are identical to the case $x =2u-1$. Therefore, the proof of statement (ii) is completed. \\
 
 \noindent (iii). Assume the contrary. Then $\rho_1$, $\rho_2$ are irreducible representations with opposite parities. Thus, $|\Tr(\rho(\fs)^2)|=|k-m|$, which is the number of self-dual objects in $\irr(\CC)$. Since $\rho(\ft)$ has $m+k$ eigenvalues of multiplicity 1, the number of self-dual objects in $\irr(\CC)$  is at least $m+k$ which is greater than $|k-m|$, a contradiction. The proof of statement (iii) is completed.
\end{proof}

As a consequence of the preceding theorem,  two non-degenerate irreducible representations with opposite parities will never satisfy the conditions of the theorem. However, we can solve the modular data if the $\ft$-spectrum of $\rho_2$ is subset of that of $\rho_1$.

\begin{thm}\label{t:solution2}
Let $\rho_1$, $\rho_2$ be non-degenerate symmetric representations of $\SL$such that 
$$
\spec(\rho_2(\ft)) \subsetneq \spec(\rho_1(\ft)).
$$
 Let $l+k = \dim\rho_1$ and $l=\dim \rho_2$, $\rho_1(\fs)=[\psi_{ij}]$, $\rho_1(\ft)=\diag(\a_1,\dots, \a_{k+l})$, $\rho_2(\fs)=[\eta_{ij}]$,  $\rho_2(\ft)=\diag(\a_{1},\dots, \a_{l})$. Suppose  $ \rho_1 \oplus \rho_2$ is equivalent to an $\SL$ representation $\rho$ of a modular tensor category $\CC$. Then 
\begin{enumerate}
   \item[(i)] 
there exists a signed diagonal matrix 
    $V$
    and $2 \times 2$ orthogonal matrices $U_i =\mtx{a_i & -b_i \\ b_i& a_i }$ with $a_i^2+b_i^2=1$ and $a_i \ge 0$ ($i =1, \dots, l$) such that 
    $$
\rho(\fs) = V\left[\begin{array}{c|c}
 A &  B^\top \\ \hline
 B & \psi'\end{array}\right]V \quad \text{and}\quad \rho(\ft) = \diag(\a_1I_2,\dots, \a_l I_2, \a_{l+1}, \dots, \a_{k+l}),
$$ 
where $\psi'$ is the $k \times k$ lower right corner block of $\rho_1(\fs)$ and $A, B$ are block matrices given by
$$
A_{ij}=U_i \left[\begin{array}{cc}
\psi_{ij} & 0 \\ 
0  & \eta_{ij}
\end{array}\right] U_j^\top, \quad 
B_{i'j}=[\psi_{l+i',j} \,\,  0]  U_j^\top, \quad 
$$
for $1 \le i,j \le l$ and $1 \le j' \le k$.

    \item[(ii)] Suppose $\rho_1$ and $\rho_2$ have opposite parities. We identify the standard basis  $(e_1, \dots, e_{2l+k})$ of $\rho$ with $\irr(\CC)$. Then 
    \begin{enumerate}
    \item[(a)] $e_{2i-1}$ and $e_{2i}$ form a dual pair for $i=1, \dots, l$.
    \item[(b)] The unit object $\1$ can only be $e_{2l+u}$ with $1 \le u \le k$ such that $\psi_{l+u, l+u} \ne 0$ and 
    $$
    \dim(\CC)=|\psi_{l+u,l+u}|^{-2},\quad \dim(e_{2i-1}) = \dim(e_{2i}) = \frac{\pm \psi_{i, l+u}}{\sqrt{2} \psi_{l+u, l+u}},\quad  \dim(e_{j})=\frac{\pm \psi_{j, l+u}}{\psi_{l+u, l+u}}
    $$
     for $i=1, \dots, l$ and $j = l+1, \dots, l+k$\,.
   In particular, they are elements of $\BZ[\zeta_N]\cap \BR^\times$ where $N$ is the     order of  $T=\a_{l+u}^{-1}\rho(\ft)$, and the $S$-matrix of $\CC$ is   given by 
    \begin{equation} \label{eq:s-mat2}
     S = \psi_{l+u,l+u}^{-1}\,V' \left[
    \begin{array}{c|c}
     A'
    & 
     B'^\top \\ \hline
     B' & \psi'
    \end{array}
    \right]V'
    \end{equation}
     for some signed diagonal matrix $V'$ and  block matrices $A'$, $B'$  given by
    $$
    A'_{ij} = 
    \mtx{\frac{\psi_{i, j}+\e_i\e_j\eta_{i,j}}{2} &\frac{\psi_{i, j}-\e_i\e_j\eta_{i,j}}{2}\\
    \frac{\psi_{i, j}-\e_i\e_j\eta_{i,j}}{2} & \frac{\psi_{i, j}+\e_i\e_j\eta_{i,j}}{2}
    } \quad \text{and}\quad
    B'_{i',j} = \frac{\psi_{l+i', j} }{\sqrt{2}}[1  \,\,  1]
    $$
    where $\e_j =\pm 1$, $1\le i, j \le l$ and $1 \le i' \le k$.
    \end{enumerate}
    \end{enumerate}
 \end{thm} 
\begin{proof} By conjugating a permutation matrix to $\rho_1 \oplus \rho_2$, we can obtain an equivalent  representation $\tilde\rho$ given by
$$
\tilde \rho(\fs)=\left[\begin{array}{c|c}
\tilde A & \tilde B^\top \\ \hline
\tilde B  & \psi'
\end{array}\right] \quad\text{and}\quad  \tilde \rho(\ft) = \diag(\a_1I_2,\dots, \a_l I_2, \a_{l+1}, \dots, \a_{k+l}),
$$ 
where $\psi'$ is the $k \times k$ bottom diagonal block of $\rho_1(\fs)$, and $\tilde A$, $\tilde B$ are block matrices with
$$
\tilde A_{ij} = \left[\begin{array}{cc}
\psi_{ij} & 0 \\ 
0  & \eta_{ij}
\end{array}\right], \quad \tilde B_{i', j} 
=[\psi_{l+i',j} \,\,  0 ]
$$
for $1 \le i,j \le l$ and $1 \le j' \le k$. By Theorem \ref{t:ortho_eqv}, there exists an orthogonal matrix $U= \diag(U_1, \dots, U_l, \g_{2l+1}, \dots, \g_{2l+1})$ such that $\rho(\fs)= U \tilde \rho(\fs) U^\top$ and $\rho(\ft)=  \tilde \rho(\ft) $ where $\g_j =\pm 1$ and $U_i$ is a $2\times 2$ orthogonal matrix for $i = 1, \dots, l$ and $j = 2l+1, \dots, 2l+k$. As before, we write $U_i = V_i \mtx{a_i & -b_i\\ b_i & a_i}$ where $a_i^2+b_i^2=1$, $a_i \ge 0$ and $V_i$ a signed diagonal matrix. Now, we set $V= \diag(V_1, \dots, V_l, \g_{2l+1}, \dots, \g_{2l+k})$, and statement (i) follows immediately.

(ii). Now we assume $\rho_1$ and $\rho_2$ are of opposite parities. Then $|\Tr(\rho(\fs)^2)|=k$ and so there are exactly $k$ self-dual simple objects in $\irr(\CC)$ and $l$ dual pairs. Since $e_{2i-1}$ and $e_{2i}$ give rise to the same eigenvalue of $\rho(\ft)$ for $i=1 , \dots, l$, and $\rho(\ft)_{2i, 2i} \ne \rho(\ft)_{j, j}$ for $j \not\in \{2i-1, 2i\}$, they must form a dual pair. Since the unit object $\1$ is self-dual, $\1 = e_{2l+u}$ for some positive integer $u \le k$, and so $1/\sqrt{\dim(\CC)}$, up to a 4-th root, is $\rho(\fs)_{2l+u, 2l+u} = \psi_{l+u, l+u}$. In particular, $\psi_{l+u, l+u} \ne 0$, $\dim(\CC) = |\psi_{l+u, l+u}|^{-2}$  and  $\psi_{l+u, l+u}^{-2} \in \BZ[\zeta_N]\cap \BR^\times$, where  $N$ is the order of $T = \a_{l+u}^{-1}\rho(\ft)$. By (i),  
$$
S= \psi_{l+u,l+u}^{-1}\,V \left[
    \begin{array}{c|c}
     A
    & 
     B^\top \\ \hline
     B & \psi'
    \end{array}
    \right]V,
$$
 where $A, B$ are block matrices given by
    $$
    A_{ij} = 
    \mtx{\frac{a_i a_j \psi_{i, j}+ b_i b_j\eta_{i,j}}{2} &\frac{a_i b_j \psi_{i, j}-a_j b_i\eta_{i,j}}{2}\\
    \frac{a_j b_i \psi_{i, j}-a_i b_j\eta_{i,j}}{2} & \frac{b_i b_j \psi_{i, j}+a_i a_j\eta_{i,j}}{2}
    } \quad \text{and}\quad
    B_{i', j} =\psi_{l+i',j} [a_j \,\, b_j]\,. 
    $$
    Thus, the dimensions of $e_{2j-1}$ and $e_{2j}$ are respectively given by
    $$
    \frac{\psi_{l+u, j}a_j}{\psi_{l+u, l+u}} \quad \text{and}\quad \frac{-\psi_{l+u, j} b_j}{\psi_{l+u, l+u}}
    $$
    which implies $\pm a_j = b_j$. Since $a_j^2+b_j^2=1$ and $a_j \ge 0$, we have $a_j = \frac{1}{\sqrt{2}}$ and $b_j =\frac{\e_j}{\sqrt{2}}$ for some $\e_j = \pm 1$ ($j=1, \dots, l$). Therefore,
    $
    A_{ij} = 
    \mtx{\frac{ \psi_{i, j}+ \e_i \e_j \eta_{i,j}}{2} &\frac{\e_j  \psi_{i, j}-\e_i\eta_{i,j}}{2}\\
    \frac{\e_i \psi_{i, j}-\e_j\eta_{i,j}}{2} & \frac{\e_i\e_j \psi_{i, j}+\eta_{i,j}}{2}
    }$ and
    $B_{i', j} =\frac{\psi_{l+i',j}}{\sqrt{2}} [1 \,\,\e_j]$. Let $E_j = \mtx{1 & 0 \\ 0 & \e_j}$ for $j = 1, \dots, l$. Then 
    $$
    A_{ij} = E_i A'_{ij} E_j \quad \text{and}\quad B_{i'j} =  B'_{i'j} E_j 
    $$    
    and the expression \eqref{eq:s-mat2} of the $S$-matrix follows immediately by setting $V' = V E$ where $E=\diag(E_1, \dots, E_l, 1,\dots, 1)$. Moreover, $\dim(e_{2j-1}) = \dim(e_{2j}) =\frac{\pm \psi_{l+u, j}}{\sqrt{2}\psi_{l+u, l+u}}$ for $j=1, \dots, l$, and $\dim(e_{2l+i'}) = \frac{\pm \psi_{l+i', l+u}}{\psi_{l+u, l+u}}$ for $1 \le i' \le k$. It follows from \cite{NS07b} that they are elements of $\BZ[\zeta_N] \cap \BR^\times$. This completes the proof of statement (ii).
\end{proof}

\section{Classification of modular data of rank=6: admissible types}\label{s:hand approach}
In this section, we prove that admissible types of MDs that can be realized by some rank=$6$ MTCs include $(4,1,1), (4,2), (3,3)$, and $(3,2,1)$.

\subsection{Classification of modular data of type (4,1,1)} Recall that  $SO(8)_3\cong PSO(8)_3\boxtimes SO(8)_1$ as modular tensor categories, which defines the notation $PSO(8)_3$. Alternatively, the modular data of $PSO(8)_3$ can be obtained from $SU(3)_6$ via boson condensation \cite{Schop}. We will prove in this section that the Galois conjugates of the modular data of $PSO(8)_3$ are characterized by the MTCs of type (4,1,1).
\begin{thm} \label{t:411}
Let $\CC$ be a rank 6 modular tensor category of type $(4,1,1)$. Then the modular data of $\CC$ is a Galois conjugate of $PSO(8)_3$.
\end{thm}

Let $\CC$ be an MTC of type $(4,1,1)$, and $\rho$ an $\SL$ representation of $\CC$. Then $\rho$ admits an irreducible decomposition $\rho_0 \oplus \rho_1 \oplus \rho_2$ in which $\dim \rho_0,  \dim \rho_1, \dim \rho_2$ respectively 4,1,1. By tensoring a suitable 1-dimensional representation of $\SL$, we will assume $\rho_0$ has a minimal $\ft$-spectrum.

In particular, all the 4-dimensional irreducible representations of level 6 are even. Now, can prove
\begin{lem}\label{l:411}  $\CC$ is self-dual, $\rho_0$ must be even of level $9$, and $\rho \cong \rho_0 \oplus 2\chi_0$.
\end{lem}
\begin{proof}
From Appendix \ref{repPP}, 4-dimensional irreducible representations of $\SL$ with minimal $\ft$-spectrums appear at the levels 5, 6, 7, 8, 9, 10, 12, 15, 20, 24 and 40.  The $\ft$-spectrums of those 4-dimensional irreducible representations of levels 5, 8, 10, 15, 20, 24 and 40 do not contain any 12-th root of unity. It follows from Lemma \ref{l:2} that $\rho_0$ cannot be of any of these levels.

It remains to show that the level of $\rho_0$ cannot be 6,  7 or 12. Suppose $\rho_0$ has level  7.  Then $\CC$ is of type (4,1,1), which contradicts Proposition \ref{p:d111}. Therefore,  the level of $\rho_0$ cannot be 7.  

 Suppose $\rho_0$ has level  6 or 12.  Then $\rho_0$ is a tensor product of two 2-dimensional representations, namely $\rho_0$ is projectively equivalent to $\rd{2}{2}{1,0} \o \rd{2}{3}{1,0}$ or $\rd{2}{4}{1,0} \o \rd{2}{3}{1,0}$ (cf. Appendix \ref{repPP}). However, $\rd{2}{2}{1,0}$ and $\rd{2}{4}{1,0}$ are projectively equivalent and so $\rd{2}{2}{1,0} \o \rd{2}{3}{1,0}$ and  $\rd{2}{4}{1,0} \o \rd{2}{3}{1,0}$. So $\rho_0$ is projectively equivalent to $\rd{2}{2}{1,0} \o \rd{2}{3}{1,0}$, which has a  minimal $\ft$-spectrums  $\{1, -1, \zeta_3, -\zeta_3\}$. Therefore,  $\rho_0 \cong \rd{2}{2}{1,0} \o \rd{2}{3}{1,0}$.   
 
By Lemma \ref{l:2}, the levels of $\rho_1$ and $\rho_2$ are divisors of 6, and so is the level of $\rho$. Therefore, $\ord(T)|6$ and hence $\CC$ is integral by Theorem \ref{t:2346}. It follows from Proposition \ref{p:weakly_int} that $\CC$ is of type (4,2), a contradiction. Therefore,  the level of $\rho_0$ cannot be 6 or 12.

As a consequence, $\rho_0$ must have level 9, and  $\rho  \cong \rho_0 \oplus 2\chi_0$ by Lemma \ref{l:2} since 1 is the unique eigenvalue of $\rho(\ft)$ with order dividing 12. It follows from Corollary \ref{l:3} that $\rho_0(\fs^2)=\id$ and $\CC$ is self-dual.
\end{proof}

\subsubsection{Solving modular data of type (4,1,1)} By Appendix \ref{repPP}, there is only one Galois orbit of 4-dimensional irreducible representations of level 9 which is even. This Galois orbit has two projectively equivalent classes given by $\rd{4}{9,1}{1,0}$ and  $\rd{4}{9,1}{8,0}$ which are complex conjugate of each other. First, we consider $\rho_0=\rd{4}{9,1}{1,0}$.

Let $z_1 = c_9^2$, $z_2=c^4_9$ and $z_3=c_9^1$.   Then
$$
\rho_0(\fs) =\frac{1}{3} \left[
\begin{array}{cccc}
 0& -\sqrt{3}& -\sqrt{3}& -\sqrt{3}\\
-\sqrt{3}& z_1 & z_2 & z_3\\
 -\sqrt{3}& z_2 & z_3 & z_1\\ 
  -\sqrt{3}& z_3 & z_1 & z_2
\end{array}\right] , \quad \rho_0(\ft) = \diag(1, \zeta_9, \zeta_9^4, \zeta_9^7)\,.
$$

Let $\tilde\rho= 2\chi_0  \oplus  \rho_0$ and set $s :=\rho(\fs)$ and 
 $t :=\rho(\ft)$. By reordering $\irr(\CC)$, one can assume 
 $$
 \rho(\ft) = \tilde\rho (\ft) = \diag(1,1,1, \zeta_9, \zeta_9^4, \zeta_9^7).
 $$
By Theorem \ref{t:ortho_eqv}, there exists $U \in O_6(\BR)$ such that $\rho= U \tilde \rho U^\top$. Then $U =f \oplus V$ for some signed diagonal matrix $V=\diag(\ve_1,\ve_2,\ve_3)$ and $f \in O_3(\BR)$ where $f \oplus V$ denotes the block direct sum of $f$ and $V$. We may further assume $\ve_3=1$, and we get
\begin{eqnarray*}
s=U\tilde\rho(\fs) U^\top 
&=& {\small \left[\begin{array}{cccccc}
f^2_{11} + f^2_{12} &  f_{11} f_{21}+f_{12}f_{22} & f_{11} f_{31}+f_{12}f_{32}& \frac{\ve_1 f_{13}}{-\sqrt{3}} &  \frac{  f_{13}}{-\sqrt{3}} & \frac{ f_{13}}{-\sqrt{3}} \\
  f_{11} f_{21}+f_{12}f_{22} & f^2_{21} + f^2_{22} &  f_{21}f_{31}+f_{22} f_{32} & \frac{\ve_1 f_{23}}{-\sqrt{3}} &  \frac{\ve_2 f_{23}}{-\sqrt{3}} & \frac{ f_{23}}{-\sqrt{3}} \\
 f_{11} f_{31}+f_{12}f_{32} &  f_{21}f_{31}+f_{22} f_{32}&  f^2_{31}+f^2_{32} &  \frac{\ve_1 f_{33}}{-\sqrt{3}} &  \frac{\ve_2 f_{33}}{-\sqrt{3}} & \frac{ f_{33}}{-\sqrt{3}} \\
  \frac{\ve_1 f_{13}}{-\sqrt{3}} & \frac{\ve_1 f_{23}}{-\sqrt{3}} & \frac{\ve_1 f_{33}}{-\sqrt{3}} & \frac{z_1}{3}  & \frac{\ve_1 \ve_2 z_2}{3}  & \frac{\ve_1  z_3}{3}  \\
  \frac{\ve_2 f_{13}}{-\sqrt{3}} & \frac{\ve_2 f_{23}}{-\sqrt{3}} & \frac{\ve_2 f_{33}}{-\sqrt{3}} & \frac{\ve_1 \ve_2 z_2}{3}  & \frac{ z_3}{3}  & \frac{\ve_2  z_1}{3} \\
   \frac{ f_{13}}{-\sqrt{3}} & \frac{ f_{23}}{-\sqrt{3}} & \frac{ f_{33}}{-\sqrt{3}} & \frac{\ve_1  z_3}{3}   &  \frac{\ve_2  z_1}{3}  &  \frac{z_2}{3} \\
\end{array}\right] }\\
&=&{\small \left[\begin{array}{cccccc}
1- f^2_{13} &  -f_{13}f_{23} & -f_{13}f_{33} & \frac{\ve_1 f_{13}}{-\sqrt{3}} &  \frac{\ve_2 f_{13}}{-\sqrt{3}} & \frac{ f_{13}}{-\sqrt{3}} \\
 -f_{13}f_{23} & 1- f^2_{23} &  -f_{23}f_{33} & \frac{\ve_1 f_{23}}{-\sqrt{3}} &  \frac{\ve_2 f_{23}}{-\sqrt{3}} & \frac{ f_{23}}{-\sqrt{3}} \\
 -f_{13}f_{33} & -f_{23}f_{34} &  1- f^2_{33} &  \frac{\ve_1 f_{33}}{-\sqrt{3}} &  \frac{\ve_2 f_{33}}{-\sqrt{3}} & \frac{ f_{33}}{-\sqrt{3}} \\
  \frac{\ve_1 f_{13}}{-\sqrt{3}} & \frac{\ve_1 f_{23}}{-\sqrt{3}} & \frac{\ve_1 f_{33}}{-\sqrt{3}} & \frac{z_1}{3}  & \frac{\ve_1 \ve_2 z_2}{3}  & \frac{\ve_1  z_3}{3}  \\
  \frac{\ve_2 f_{13}}{-\sqrt{3}} & \frac{\ve_2 f_{23}}{-\sqrt{3}} & \frac{\ve_2 f_{33}}{-\sqrt{3}} & \frac{\ve_1 \ve_2 z_2}{3}  & \frac{ z_3}{3}  & \frac{\ve_2  z_1}{3} \\
   \frac{ f_{13}}{-\sqrt{3}} & \frac{ f_{23}}{-\sqrt{3}} & \frac{ f_{33}}{-\sqrt{3}} & \frac{\ve_1  z_3}{3}   &  \frac{\ve_2  z_1}{3}  &  \frac{z_2}{3} \\
\end{array}\right]}.
\end{eqnarray*}

We now apply the Galois symmetry \cite[Theorem II]{DLN} of $\rho$ to determine $f$ and $\ve_1, \ve_2$  (cf. Theorem \ref{p:MD1} (6)). Since $\ord(t)=9$, then $s$ is a matrix over $\BQ_9$. The Galois group $\Gal(\BQ_9/\BQ)$ is generated by $\s$ defined by $\s: \zeta_9 \mapsto \zeta_9^2$, and $\hs$ denotes the corresponding  permutation on $\irr(\CC)=\{1, \dots,  6\}$. The $i$-th diagonal entry of $t$ will be denoted by $t_i$. Under the action of $\s^2$, 
$$
t_4 \mapsto t_5, \quad t_5 \mapsto t_6 ,\quad \text{and}\quad  t_6 \mapsto t_4 \,.
$$
We find $\hs(4)=5$, $\hs(5)=6$ and $\hs(6)=4$. Recall that $\s(s_{ij}) = \e_\s(i) s_{\hs(i) j}$ where $\e_\s(i) = \pm 1$.  Applying $\s$ to those $s_{ij}$ with $i, j \in \{4,5,6\}$, we have
$$
\s(z_1) = \e_\s(4) \ve_1\ve_2 z_2, \quad \s(\ve_1 \ve_2  z_2) = \e_\s(5) \ve_1 z_3 \quad  \text{and}\quad \s(\ve_1 z_3) = \e_\s(6) z_1\,.
$$
 Since 
$\s(z_1) = z_2$, $\s(z_2) = z_3$ and $\s(z_3) = z_1$, we find
$$
\e_\s(4) = \ve_1 \ve_2, \quad \e_\s(5) = \ve_2 \quad \text{ and }\quad \e_\s(6) = \ve_1\,.
$$
Now, we apply $\s$ to those $s_{ij}$ with $i\in \{1,2,3\}$ and  $j \in \{4,5,6\}$.  We have   $\s(\frac{f_{i3}}{\sqrt{3}}) = \frac{f_{i3}}{\sqrt{3}}$ , and hence $\frac{f_{i3}}{\sqrt{3}}  \in \BQ$  for $i =1,2,3$. This implies that $f_{i3}f_{j3} \in \BQ$ for any $i, j \in \{1,2,3\}$. Therefore,  the first 3 rows of $s$ have rational entries, and hence  $\hs$ fixes 1,2,3. Now, we can conclude that $\hs = (4,5,6)$.

Since $\CC$ is not integral by Proposition \ref{p:weakly_int}, None of $1, 2$ or $3$ cannot be the isomorphism class of the unit object $\1$ or the simple object $\iota$ for the Frobenius-Perron dimensions. Therefore, $\dim(\CC)$ and $\FPdim(\CC)$ are Galois conjugates, and $\FPdim(\CC)$ is the largest conjugate of $\dim(\CC)$. The global dimension $\dim(\CC)$ can be $9 z^{-2}_1$, $9 z^{-2}_3$ or $9 z_2^{-2}$ depending which of the classes 4,5,6 corresponds $\1$. Since they are conjugates and $-z_2 > z_3 > z_1>0$,   $\FPdim(\CC)=9 z_1^{-2} $.

Let $(S,T)$ be the modular data of $\CC$. Note that  $z_1, z_2, z_3$  are units, and they are roots of the irreducible polynomial $x^3 -3 x+1$. No matter which of 4,5,6 is the isomorphism class of $\1$, for $i \in \{1,2,3\}$ and $j \in \{4,5,6\}$,
$$
S_{ij} = \pm \frac{\sqrt{3} f_{i3}}{z_k} 
$$
for some $k \in \{1,2,3\}$. Since $S_{ij}$ is a cyclotomic integer, so is $\sqrt{3} f_{i3}$. Thus, $\sqrt{3} f_{i3}$ is an integer and they satisfy 
$$
(\sqrt{3} f_{13})^2 +(\sqrt{3} f_{23})^2 +(\sqrt{3} f_{33})^2 =3\,.
$$
Therefore, $\sqrt{3} f_{i3} = \pm 1$ or equivalently  $f_{i3}=\pm \frac{1}{\sqrt{3}}$ for $i = 1,2,3$. 
Now, we can compute the modular data for the cases  $\1=4,5$ or $6$:\\
\noindent
(i) Suppose $4$ is the isomorphism class of $\1$. Then $D=3/z_1$ and 
$$
S = \left[\begin{array}{cccccc}
3 \frac{1-f^2_{13}}{z_1} &  3\frac{f_{13}f_{23}}{-z_1}  & 3\frac{f_{13}f_{33}}{-z_1} & \frac{\ve_1 \sqrt{3} f_{13}}{-z_1} &  \frac{\ve_2 \sqrt{3} f_{13}}{-z_1} & \frac{\sqrt{3} f_{13}}{-z_1} \\
3 \frac{f_{13}f_{23}}{-z_1}  & 3\frac{1-f^2_{23}}{z_1} &  3\frac{f_{23}f_{33}}{-z_1}  & \frac{\ve_1 \sqrt{3} f_{23}}{-z_1} &  \frac{\ve_2 \sqrt{3} f_{23}}{-z_1} & \frac{\sqrt{3} f_{23}}{-z_1}  \\
 3\frac{f_{13}f_{33}}{-z_1}  & 3\frac{f_{23}f_{33}}{-z_1}  & 3 \frac{1-f^2_{33}}{z_1} &  \frac{\ve_1 \sqrt{3} f_{33}}{-z_1}  &  \frac{\ve_2 \sqrt{3} f_{33}}{-z_1} & \frac{\sqrt{3} f_{33}}{-z_1} \\
  \frac{\ve_1 \sqrt{3} f_{13}}{-z_1} & \frac{\ve_1 \sqrt{3} f_{23}}{-z_1} & \frac{\ve_1 \sqrt{3} f_{33}}{-z_1 } & 1  & \frac{\ve_1 \ve_2 z_2}{z_1} & \frac{\ve_1  z_3}{z_1} \\
  \frac{\ve_2 \sqrt{3} f_{13}}{-z_1} & \frac{\ve_2 \sqrt{3} f_{23}}{-z_1 } & \frac{\ve_2 \sqrt{3} f_{33}}{-z_1 } & \frac{\ve_1 \ve_2 z_2}{z_1} &  \frac{z_3}{z_1} &  \ve_2  \\
   \frac{\sqrt{3} f_{13}}{-z_1 } & \frac{\sqrt{3} f_{23}}{-z_1 } & \frac{\sqrt{3} f_{33}}{-z_1 } & \frac{\ve_1  z_3}{z_1}   &  \ve_2  &  \frac{z_2}{z_1}\\
\end{array}\right].
$$
Note that 
$$
\sum_{i=1}^6 \left(\frac{S_{i,4}}{S_{4,4}}\right)^2 = \frac{9}{z_1^2} =\FPdim(\CC) \,.
$$
Therefore, $4$ is also the isomorphism class of $\iota$. In particular, $\CC$ is pseudounitary and the entries of 4th row of $S$ must be positive. Since $\frac{z_2}{z_1} <0$ and $\frac{z_3}{z_1} >0$, we have $\ve_1 =1$, $\ve_2=-1$ and $f_{i3}  < 0$ for $i = 1,2,3$.  
This implies $\sqrt{3} f_{i3}=-1$ for $i=1, 2,3$ and
$$
S = \left[\begin{array}{cccccc}
2 z_1^{-1} & -z_1^{-1} & -z_1^{-1} &  z_1^{-1} & -z_1^{-1} & z_1^{-1} \\
 -z_1^{-1} & 2 z_1^{-1} &  -z_1^{-1} &   z_1^{-1} & -z_1^{-1} & z_1^{-1}   \\
 -z_1^{-1} & -z_1^{-1} &  2 z_1^{-1}  & z_1^{-1} & -z_1^{-1} & z_1^{-1} \\
 z_1^{-1} & z_1^{-1} & z_1^{-1}  & 1  & \frac{-z_2}{z_1}  & \frac{z_3}{z_1} \\
  -z_1^{-1} & -z_1^{-1} & -z_1^{-1} & \frac{-z_2}{z_1} &  \frac{z_3}{z_1} & -1 \\
    z_1^{-1} & z_1^{-1} & z_1^{-1} &  \frac{z_3}{z_1}   &  -1  &   \frac{z_2}{z_1}   \\
\end{array}\right]\quad\text{and}\quad
T= \diag(\zeta_9^8, \zeta_9^8, \zeta_9^8, 1, \zeta_3, \zeta_3^2 )\,.
$$
\noindent
(ii) Suppose $5$ is the isomorphism class of $\1$. Then $D=3/z_3$ and hence
$$
S = \left[\begin{array}{cccccc}
3 \frac{1-f^2_{13}}{z_3} &  3\frac{f_{13}f_{23}}{-z_3}  & 3\frac{f_{13}f_{33}}{-z_3} & \frac{\ve_1 \sqrt{3} f_{13}}{-z_3} &  \frac{\ve_2 \sqrt{3} f_{13}}{-z_3} & \frac{\sqrt{3} f_{13}}{-z_3} \\
3 \frac{f_{13}f_{23}}{-z_3}  & 3\frac{1-f^2_{23}}{z_3} &  3\frac{f_{23}f_{33}}{-z_3}  & \frac{\ve_1 \sqrt{3} f_{23}}{-z_3} &  \frac{\ve_2 \sqrt{3} f_{23}}{-z_3} & \frac{\sqrt{3} f_{23}}{-z_3}  \\
 3\frac{f_{13}f_{33}}{-z_3}  & 3\frac{f_{23}f_{33}}{-z_3}  & 3 \frac{1-f^2_{33}}{z_3} &  \frac{\ve_1 \sqrt{3} f_{33}}{-z_3}  &  \frac{\ve_2 \sqrt{3} f_{33}}{-z_3} & \frac{\sqrt{3} f_{33}}{-z_3} \\
  \frac{\ve_1 \sqrt{3} f_{13}}{-z_3} & \frac{\ve_1 \sqrt{3} f_{23}}{-z_3} & \frac{\ve_1 \sqrt{3} f_{33}}{-z_3} & \frac{z_1}{z_3}  & \frac{\ve_1 \ve_2 z_2}{z_3} & \ve_1  \\
  \frac{\ve_2 \sqrt{3} f_{13}}{-z_3} & \frac{\ve_2 \sqrt{3} f_{23}}{-z_3} & \frac{\ve_2 \sqrt{3} f_{33}}{-z_3} & \frac{\ve_1 \ve_2 z_2}{z_3} &  1 & \frac{\ve_2 z_1}{z_3}  \\
   \frac{\sqrt{3} f_{13}}{-z_3} & \frac{\sqrt{3} f_{23}}{-z_3} & \frac{\sqrt{3} f_{33}}{-z_3} & 
   \ve_1     &  \frac{\ve_2 z_1}{z_3}  &  \frac{z_2}{z_3}\\
\end{array}\right].
$$
Now, one can check directly that
$$
\sum_{i=1}^6 (\frac{S_{i4}}{S_{54}})^2 = \frac{9}{z_2^2} \quad \text{and}\quad \sum_{i=1}^6 (\frac{S_{i6}}{S_{56}})^2 = \frac{9}{z_1^2}\,,
$$
which implies $6$ is the isomorphism class of $\iota$. Thus, all the entries of the 6th row of $S$ have the same sign. Since $z_2/z_3 <0$ and $z_1/z_3 >0$, we obtain that $\ve_1 =\ve_2 = -1$ and $f_{i3} > 0$ for $i =1,2,3$.  Therefore, 
$$
\sqrt{3} f_{13} =  \sqrt{3}f_{23}=\sqrt{3} f_{33}  =1
$$
and hence
$$
S = \left[\begin{array}{cccccc}
2z_3^{-1}   &  -z_3^{-1}  &  -z_3^{-1} & z_3^{-1}   & z_3^{-1}    & -z_3^{-1}\\
 -z_3^{-1}& 2z_3^{-1}   &  -z_3^{-1} & z_3^{-1}   & z_3^{-1}    & -z_3^{-1}  \\
 -z_3^{-1} & -z_3^{-1}  &  2z_3^{-1}   & z_3^{-1}   & z_3^{-1}    & -z_3^{-1} \\
  z_3^{-1}   & z_3^{-1}    & z_3^{-1}   & \frac{z_1}{z_3}  & \frac{z_2}{z_3} & -1\\
  z_3^{-1}  & z_3^{-1}  & z_3^{-1}  & \frac{z_2}{z_3} & 1 & \frac{-z_1}{z_3}   \\
   -z_3^{-1}   & -z_3^{-1}    & -z_3^{-1}  & -1  &  \frac{-z_1}{z_3}    &  \frac{z_2}{z_3} \\
\end{array}\right]
\quad\text{and}\quad
T= \diag(\zeta_9^5, \zeta_9^5, \zeta_9^5, \zeta_3^2, 1, \zeta_3 )\,.
$$
\noindent
(iii) Suppose $6$ is the isomorphism class of $\1$. Then $D=3/z_2$ and 
$$
S = \left[\begin{array}{cccccc}
3\frac{1-f^2_{13}}{z_2} &  3\frac{f_{13}f_{23}}{-z_2}  & 3\frac{f_{13}f_{33}}{-z_2} & \frac{\ve_1 \sqrt{3} f_{13}}{-z_2} &  \frac{\ve_2 \sqrt{3} f_{13}}{-z_2} & \frac{ \sqrt{3}f_{13}}{-z_2} \\
3 \frac{f_{13}f_{23}}{-z_2}  & 3\frac{1-f^2_{23}}{z_2} &  3\frac{f_{23}f_{33}}{-z_2}  & \frac{\ve_1 \sqrt{3} f_{23}}{-z_2} &  \frac{\ve_2 \sqrt{3} f_{23}}{-z_2} & \frac{\sqrt{3} f_{23}}{-z_2}  \\
 3 \frac{f_{13}f_{33}}{-z_2}  & 3\frac{f_{23}f_{33}}{-z_2}  &  3\frac{1-f^2_{33}}{z_2} &  \frac{\ve_1 \sqrt{3} f_{33}}{-z_2}  &  \frac{\ve_2 \sqrt{3} f_{33}}{-z_2} & \frac{ \sqrt{3} f_{33}}{-z_2} \\
  \frac{\ve_1 \sqrt{3} f_{13}}{-z_2} & \frac{\ve_1\sqrt{3} f_{23}}{-z_2} & \frac{\ve_1\sqrt{3}  f_{33}}{-z_2 } & \frac{z_1}{z_2}  &  \ve_1 \ve_2 & \frac{\ve_1 z_3}{z_2}  \\
  \frac{\ve_2 \sqrt{3} f_{13}}{-z_2} & \frac{\ve_2\sqrt{3} f_{23}}{-z_2 } & \frac{\ve_2 \sqrt{3} f_{33}}{-z_2 } &\ve_1 \ve_2 &  \frac{z_3}{z_2}   & \frac{\ve_2 z_1}{z_2}  \\
   \frac{\sqrt{3} f_{13}}{-z_2 } & \frac{\sqrt{3} f_{23}}{-z_2 } & \frac{\sqrt{3} f_{33}}{-z_2 } & \frac{\ve_1 z_3}{z_2}   &  \frac{\ve_2 z_1}{z_2}   &  1\\
\end{array}\right].
$$
Now, 
$$
\sum_{i=1}^6 (\frac{S_{i4}}{S_{64}})^2 = \frac{9}{z_3^2} \quad \text{and}\quad \sum_{i=1}^6 (\frac{S_{i5}}{S_{65}})^2 = \frac{9}{z_1^2}\,,
$$
which implies $5$ is the isomorphism class of $\iota$. Thus,  all the entries of the 4th row have the same signs. Since $z_3/z_2<0$ and $z_1/z_2 <0$,  $\ve_1 = -1$, $\ve_2=1$ and $f_{i3} > 0$ for $i =1,2,3$.  Therefore, 
$$
 \sqrt{3} f_{13} =  \sqrt{3}f_{23}=\sqrt{3} f_{33}  =1
$$
and hence
$$
S = \left[\begin{array}{cccccc}
2z_2^{-1} & -z_2^{-1} & -z_2^{-1}&-z_2^{-1}  &z_2^{-1}   & z_2^{-1}\\
-z_2^{-1} &2z_2^{-1}  & -z_2^{-1}&-z_2^{-1}  &z_2^{-1}   &z_2^{-1} \\
-z_2^{-1}&-z_2^{-1} & 2z_2^{-1}  &-z_2^{-1}  &z_2^{-1}   &z_2^{-1}\\
 -z_2^{-1}  &-z_2^{-1}   &-z_2^{-1}  & \frac{z_1}{z_2}  & -1& \frac{-z_3}{z_2}  \\
 z_2^{-1} &z_2^{-1} &z_2^{-1} & -1 &  \frac{z_3}{z_2}& \frac{z_1}{z_2}   \\
  z_2^{-1}  &z_2^{-1}   &z_2^{-1} & \frac{-z_3}{z_2}    &  \frac{z_1}{z_2}    &  1\\
\end{array}\right] \quad\text{and}\quad
T= \diag(\zeta_9^2, \zeta_9^2, \zeta_9^2, \zeta_3,\zeta_3^2, 1 )\,.
$$

Now, we compute the modular data for $\rho_0 = \rd{4}{9,1}{8,0}$, which is the complex conjugate of $\rd{4}{9,1}{1,0}(\fs)$. Since $\rd{4}{9,1}{1,0}(\fs))=\rd{4}{9,1}{8,0}(\fs)$, modular data  are complex conjugations of those obtained for $\rho_0 = \rd{4}{9,1}{1,0}$. They are:\\
\noindent
(iv)
$$
S = \left[\begin{array}{cccccc}
2 z_1^{-1} & -z_1^{-1} & -z_1^{-1} &  z_1^{-1} & -z_1^{-1} & z_1^{-1} \\
 -z_1^{-1} & 2 z_1^{-1} &  -z_1^{-1} &   z_1^{-1} & -z_1^{-1} & z_1^{-1}   \\
 -z_1^{-1} & -z_1^{-1} &  2 z_1^{-1}  & z_1^{-1} & -z_1^{-1} & z_1^{-1} \\
 z_1^{-1} & z_1^{-1} & z_1^{-1}  & 1  & \frac{-z_2}{z_1}  & \frac{z_3}{z_1} \\
  -z_1^{-1} & -z_1^{-1} & -z_1^{-1} & \frac{-z_2}{z_1} &  \frac{z_3}{z_1} & -1 \\
    z_1^{-1} & z_1^{-1} & z_1^{-1} &  \frac{z_3}{z_1}   &  -1  &   \frac{z_2}{z_1}   \\
\end{array}\right]\quad\text{and}\quad
T= \diag(\zeta_9, \zeta_9, \zeta_9, 1, \zeta_3^2, \zeta_3 )\,.
$$
\noindent
(v)
$$
S = \left[\begin{array}{cccccc}
2z_3^{-1}   &  -z_3^{-1}  &  -z_3^{-1} & z_3^{-1}   & z_3^{-1}    & -z_3^{-1}\\
 -z_3^{-1}& 2z_3^{-1}   &  -z_3^{-1} & z_3^{-1}   & z_3^{-1}    & -z_3^{-1}  \\
 -z_3^{-1} & -z_3^{-1}  &  2z_3^{-1}   & z_3^{-1}   & z_3^{-1}    & -z_3^{-1} \\
  z_3^{-1}   & z_3^{-1}    & z_3^{-1}   & \frac{z_1}{z_3}  & \frac{z_2}{z_3} & -1\\
  z_3^{-1}  & z_3^{-1}  & z_3^{-1}  & \frac{z_2}{z_3} & 1 & \frac{-z_1}{z_3}   \\
   -z_3^{-1}   & -z_3^{-1}    & -z_3^{-1}  & -1  &  \frac{-z_1}{z_3}    &  \frac{z_2}{z_3} \\
\end{array}\right]
\quad\text{and}\quad
T= \diag(\zeta_9^4, \zeta_9^4, \zeta_9^4, \zeta_3, 1, \zeta_3^2 )\,.
$$
\noindent
(vi)
$$
S = \left[\begin{array}{cccccc}
2z_2^{-1} & -z_2^{-1} & -z_2^{-1}&-z_2^{-1}  &z_2^{-1}   & z_2^{-1}\\
-z_2^{-1} &2z_2^{-1}  & -z_2^{-1}&-z_2^{-1}  &z_2^{-1}   &z_2^{-1} \\
-z_2^{-1}&-z_2^{-1} & 2z_2^{-1}  &-z_2^{-1}  &z_2^{-1}   &z_2^{-1}\\
 -z_2^{-1}  &-z_2^{-1}   &-z_2^{-1}  & \frac{z_1}{z_2}  & -1& \frac{-z_3}{z_2}  \\
 z_2^{-1} &z_2^{-1} &z_2^{-1} & -1 &  \frac{z_3}{z_2}& \frac{z_1}{z_2}   \\
  z_2^{-1}  &z_2^{-1}   &z_2^{-1} & \frac{-z_3}{z_2}    &  \frac{z_1}{z_2}    &  1\\
\end{array}\right] \quad\text{and}\quad
T= \diag(\zeta_9^7, \zeta_9^7, \zeta_9^7, \zeta_3^2,\zeta_3, 1 )\,.
$$

\subsubsection{Proof of Theorem \ref{t:411}}
Since modular data of Type (4,1,1) have been completely solved in the last subsection. The modular data of $PSO(8)_3$ coincides with (i) up to a permutation. Let $\s \in \Gal(\BQ_9)$ be the generator defined by $\s: \zeta_9  \mapsto\zeta_9^2$. Applying $\s$ to the modular data (i)-(vi), One can check directly 
$$
\text{(i)} \xrightarrow{\s}  \text{(vi)} \xrightarrow{\s} \text{(ii)}  \xrightarrow{\s} \text{(iv)} \xrightarrow{\s} \text{(iii)} \xrightarrow{\s} \text{(v)} \xrightarrow{\s} \text{(i)}        
$$ 
up to permutations of the objects. This completes the proof of Theorem \ref{t:411}. \qed

%%%old existence cases

\subsection{Classification of modular data of type (4,2)}

In this section, we will complete the classification of modular data of type (4,2)  in the following theorem.
\begin{thm} \label{t:42}
Let $\CC$ be a rank 6 modular tensor category of type $(4,2)$. Then the modular data of $\CC$ can only be a  Galois conjugate of the modular data of the following modular tensor categories:
\begin{enumerate}
    \item[(1)] $\CC(\BZ_6, q)$ with $q(1) = \zeta_{12}$;
    \item[(2)]  $\CC(\BZ_3, q) \boxtimes PSU(2)_3$ with $q(1) = \zeta_{3}$;
    \item[(3)] $G(2)_3$\,.
\end{enumerate}
\end{thm}
We will use the following level 5 irreducible representations $\rd{2}{5}{1}$, $\rd{4}{5,1}{1}$ and $\rd{4}{5,2}{1}$ when necessary.
\begin{equation} \label{eq:rd251}
\rd{2}{5}{1}(\fs)= \frac{1}{s_5^1}\left[
\arraycolsep=3pt
\begin{array}{cr}
 1 & \varphi  \\
 \varphi  & {-1} 
\end{array}
\right], \quad   \rd{2}{5}{1}(\ft)=\diag(\zeta_5, \zeta_5^4)\,. 
\end{equation}
Note that $\rd{2}{5}{1}$ is defined over $\BQ_5$. Let $\s \in \GQ$ such that $\s(\zeta_5)=\zeta_5^2$.  Then $\rd{2}{5}{2}:=\s \circ \rd{2}{5}{1}$. $\rd{2}{5}{i}$, $i=1,2$, form a complete set of inequivalent 2-dimensional representations of level 5. The following irreducible representations also  form a complete set of inequivalent 4-dimensional representations of level 5:
\begin{equation} \label{eq:rd451}
\rd{4}{5,1}{1}(\fs)=\frac{s_5^3}{5}\left[
\arraycolsep=3pt\def\arraystretch{.9}
\begin{array}{cccc}
 -\varphi ^2 & \varphi^{-1} & \sqrt{3} \varphi  & \sqrt{3} \\
 \varphi^{-1} & \varphi ^2 & \sqrt{3} & -\sqrt{3} \varphi  \\
 \sqrt{3} \varphi  & \sqrt{3} & \varphi^{-1} & \varphi ^2 \\
 \sqrt{3} & -\sqrt{3} \varphi  & \varphi ^2 & -\varphi^{-1} \\
\end{array}
\right], \quad \rd{4}{5,1}{1}(\ft)=\diag(\zeta_5, \zeta_5^2, \zeta_5^3, \zeta_5^4)\,. 
\end{equation}
\begin{equation} \label{eq:rd452}
\rd{4}{5,2}{1}(\fs)=\frac{1}{\sqrt{5}}\left[
\arraycolsep=3pt\def\arraystretch{.9}
\begin{array}{cccc}
 1 & -1 & \varphi^{-1} & \varphi  \\
 -1 & 1 & \varphi  & \varphi^{-1} \\
 \varphi^{-1} & \varphi  & -1 & 1 \\
 \varphi  & \varphi^{-1} & 1 & -1 \\
\end{array}
\right], \quad \rd{4}{5,2}{1}(\ft)=\diag(\zeta_5, \zeta_5^2, \zeta_5^3, \zeta_5^4)\,. 
\end{equation}

We will need to establish a few lemmas to complete the proof of this theorem. 
Let $\CC$ be a modular tensor category of type (4,2) and $\rho$ an $\SL$ representation of $\CC$. Then 
$$
\rho \cong \rho_1 \oplus \rho_2
$$
for some irreducible representations $\rho_1, \rho_2$ of dimensions 4 and 2 respectively. By tensoring with a suitable $\chi^i \in \hSL$, we may assume that the $\ft$-spectrum of $\rho_1$ is minimal. Therefore,  $\rho_1$ has a prime power level or $\rho_1$ is a tensor product of two 2-dimensional irreducible representations of distinct prime power levels. 

According to Appendix A, $\rho_1$ can only have the prime power levels 5, 7, 8, 9 or the composite levels 6, 10, 15, 24, 40.  Note that a 4-dimensional irreducible representation of level 12 is projectively equivalent to an irreducible representation of level  6 as shown in the proof of Lemma \ref{l:411}. We will prove that only the levels 7, 15 and 24 are possible. 

It follows from Appendix \ref{repPP} that the eigenvalues of  $\rho_1(\ft)$ and $\rho_2(\ft)$ are multiplicity free. By the $\ft$-spectrum criteria, $\spec(\rho_1(\ft)) \cap \spec(\rho_2(\ft))=\{\tilde \theta_0\}$ or $\spec(\rho_2(\ft)) \subset \spec(\rho_1(\ft)$. These situations have been studied in Theorems  \ref{t:solution1} and \ref{t:solution2}. Now, we can begin to prove the level of $\rho_1$ cannot 5, 8, or 9. 
\begin{lem}\label{l:9.3}
The level of $\rho_1$ cannot be $5$.
\end{lem}
\begin{proof}
Suppose $\rho_1$ is of level 5. Since there are exactly two inequivalent irreducible representations of level 5 and dimension 4, which are given by  $\rd{4}{5,1}{1}$ and $\rd{4}{5,2}{1}$,  $\rho_1$ must be equivalent one of them. In particular, the spectrum of $\rho_1(\ft)$ consists of all the primitive $5$-th root of unity. By the $\ft$-spectrum criteria, $\rho_2$ can only be equivalent to $\rd{2}{5}{1}$ or $\rd{2}{5}{2}$, which are the inequivalent irreducible representations of level 5 and dimension 2. Therefore, $\rho$ is of level 5 and hence $\rho(\fs)$ is a matrix over $\BQ_5$. Let $\s \in \GQ$  such that $\s(\zeta_5) = \zeta_5^2$. Then $\rd{2}{5}{2} = \s\circ \rd{2}{5}{1}$.

Note that $\tau \circ \rd{4}{5,i}{1}\cong \rd{4}{5,i}{1}$ for all $\tau \in \GQ$ and $i=1,2$.  Thus, if $\rho_1 \oplus \rd{2}{5}{1}$ is not equivalent to any MD representation, then so is $\s\circ (\rho_1 \oplus \rd{2}{5}{1}) \cong \rho_1 \oplus \rd{2}{5}{2}$. Therefore, it suffices to show that $\rd{4}{5,1}{1} \oplus \rd{2}{5}{1}$ and $\rd{4}{5,2}{1} \oplus \rd{2}{5}{1}$ are not equivalent to any MD representation.

(i) Suppose $\rho_1=\rd{4}{5,1}{1}$ and $\rho_2=\rd{2}{5}{1}$. Using the representations $\rd{4}{5,1}{1}$ and   $\rd{2}{5}{1}$ presented in \eqref{eq:rd451} and \eqref{eq:rd251}, we have
\begin{eqnarray*}
(\rho_1 \oplus \rho_2)(\fs)& = &
\frac{s_5^3}{5}\left[
\arraycolsep=3pt\def\arraystretch{.9}
\begin{array}{cccc}
 -\varphi ^2 & \varphi^{-1} & \sqrt{3} \varphi  & \sqrt{3} \\
 \varphi^{-1} & \varphi ^2 & \sqrt{3} & -\sqrt{3} \varphi  \\
 \sqrt{3} \varphi  & \sqrt{3} & \varphi^{-1} & \varphi ^2 \\
 \sqrt{3} & -\sqrt{3} \varphi  & \varphi ^2 & -\varphi^{-1} \\
\end{array}
\right] \oplus 
 \frac{1}{s_5^1}\left[
\arraycolsep=3pt
\begin{array}{cr}
 1 & \varphi  \\
 \varphi  & {-1} 
\end{array}
\right]\\
(\rho_1 \oplus \rho_2)(\ft) &=& \diag(\zeta_5, \zeta_5^4, \zeta_5^2, \zeta_5^3, \zeta_5, \zeta_5^4 )\,. 
\end{eqnarray*}
By Theorem \ref{t:solution2} (1),  There exists a block diagonal orthogonal matrix 
$$U = 
\left[\begin{array}{cr}
a & -  b  \\
 b &   a 
\end{array}
\right]\oplus  
\left[
\begin{array}{cr}
c & - d  \\
 d &  c 
\end{array}
\right] \oplus 
I_2 \quad\text{with } a^2+b^2=1,\, c^2+d^2=1,
$$
such that $\rho(\ft) = \diag(\zeta_5, \zeta_5, \zeta_5^4, \zeta_5^4, \zeta_5^2, \zeta_5^3)$ and $\rho(\fs)$ is a conjugation of $s'$ by a signed diagonal matrix, where $s'$ is given by
$$
s'=\frac{s_5^3}{5}
\left[
\begin{array}{cccccc}
 * & * & * & * & -\sqrt{3} b     \varphi  & -\sqrt{3} b   \\
  * & * & * & *  & \sqrt{3} a     \varphi  & \sqrt{3} a   \\
  * & * & * & * & -\sqrt{3} d     & \sqrt{3} d   \varphi  \\
  * & * & * & * & \sqrt{3} c     & -\sqrt{3} c   \varphi  \\
 -\sqrt{3} b     \varphi  & \sqrt{3} a     \varphi  & -\sqrt{3} d     & \sqrt{3} c     & \varphi^{-1} &   \varphi ^2 \\
 -\sqrt{3} b   & \sqrt{3} a   & \sqrt{3} d   \varphi  & -\sqrt{3} c   \varphi  &   \varphi ^2 & -\varphi^{-1} \\
\end{array}
\right]\,.
$$
It follows from the action of $\s^2$ on $\rho(\ft)$, we find $\hs(5)=6$. Since 
$$
\s(s_5^3/5)=\frac{s_5^1}{5} = - \frac{s_5^3}{5}\varphi \quad \text{and} \quad \s(\varphi) = -\varphi^{-1}\,,
$$
the action of $\s$ on $s'_{55}$  implies $\e'_\s(5)=1$. Hence, by the action of $\s$ on the 5-th column, we have
$$
\s(\sqrt{3} x) = \sqrt{3}x \text{ for } x=a,b,c,d\,. 
$$
Therefore, $\sqrt{3}a, \sqrt{3} b, \sqrt{3} c, \sqrt{3} d \in \BQ$ as $\s|_{\BQ_5}$ generates $\Gal(\BQ_5/\BQ)$. If $5$ (resp.  $6$) corresponds to the unit object $\1$, then $s'/s'_{55}$ (resp. $s'/s'_{66}$) is a matrix $\BZ[\zeta_5]$. Since $\varphi$ is a unit in $\BZ[\zeta_5]$,  $\sqrt{3}a, \sqrt{3} b, \sqrt{3} c, \sqrt{3} d \in \BZ[\zeta_5]$  and hence $\sqrt{3}a, \sqrt{3} b, \sqrt{3} c, \sqrt{3} d \in \BZ\setminus\{0\}$. However, this contradicts that $(\sqrt{3}a)^2+(\sqrt{3} b)^2 = 3$. Therefore, 5 and 6 cannot be $\1$.  

Suppose 1 is the isomorphism class of $\1$. Then $s'_{i,5}/s'_{1,5} \in \BZ[\zeta_5]$ for all $i$. In particular,  $\frac{1}{\sqrt{3}b}$, $a/b \in   \BZ[\zeta_5]$. So, $\frac{1}{\sqrt{3}b}, a/b \in \BZ$. Let $m, n \in \BZ$ such that $a =mb$ and $1 = \sqrt{3} b n$. The equality $a^2+b^2=1$ implies
$(m^2+1)3b^2 = 3$ and so  $m^2+1 = 3n^2$. However, $3 \nmid (m^2+1)$ for any integer $m$. Therefore, 1 cannot the unit object. By the same reason, 2, 3, and 4 are not the isomorphism class of $\1$. This ultimate contradiction implies that
$\rd{4}{5,1}{1}\oplus \rd{2}{5}{1}$ is not equivalent any MD representation. 

(ii) Now we assume $\rho_1=\rd{4}{5,2}{1}$ and $\rho_2=\rd{2}{5,1}{1}$. It follows from \eqref{eq:rd251}  and \eqref{eq:rd452} that
\begin{eqnarray*}
(\rho_1 \oplus \rho_2)(\fs)& = &
\frac{1}{\sqrt{5}}\left[
\arraycolsep=3pt\def\arraystretch{.9}
\begin{array}{cccc}
 1 & -1 & \varphi^{-1} & \varphi  \\
 -1 & 1 & \varphi  & \varphi^{-1} \\
 \varphi^{-1} & \varphi  & -1 & 1 \\
 \varphi  & \varphi^{-1} & 1 & -1 \\
\end{array}
\right] \oplus 
 \frac{1}{s_5^1}\left[
\arraycolsep=3pt
\begin{array}{cr}
 1 & \varphi  \\
 \varphi  & {-1} 
\end{array}
\right]\\
(\rho_1 \oplus \rho_2)(\ft) &=& \diag(\zeta_5, \zeta_5^4, \zeta_5^2, \zeta_5^3, \zeta_5, \zeta_5^4 )\,. 
\end{eqnarray*}
Note that $\rho_1$, $\rho_2$ have opposite parities. We reorder the simple objects as in Theorem \ref{t:solution2} so that 
$\rho(\ft)=\diag(\zeta_5, \zeta_5, \zeta_5^4, \zeta_5^4, \zeta_5^2, \zeta_5^3)$. The unit object can only be $e_5$ or $e_6$. In either case, we find $\dim(\CC)=5$, and $\dim(e_1)=\dim(e_2)=\frac{\pm \rho_1(\fs)_{13}}{\sqrt{2}\rho_1(\fs)_{33}} = \frac{\pm \varphi^{-1}}{\sqrt{2}} \not\in \BQ_5$. This contradicts Theorem \ref{p:MD} (4). Therefore, $\rd{4}{5,2}{1} \oplus \rd{2}{5,1}{1}$ is not equivalent to any MD representation. This completes the proof of this lemma.
\end{proof}
\begin{lem} \label{l:9.4}
The level $\rho_1$ cannot be 8. 
\end{lem}
\begin{proof}  Suppose $\rho_1$ has level 8.
Since there is only one projectively equivalent class of irreducible representations of level 8 and dimension 4.  One can assume $\rho_1 = \rd{4}{8}{1,0}$ (cf. Appendix \ref{repPP}). In particular, $\rho_1$ is odd, and $\spec(\rho_1(\ft))$ consists of all the primitive 8-th roots of unity.

 By the $\ft$-spectrum criteria, $\spec(\rho_2(\ft))$ must be a set of primitive 8-th roots of unity, and hence $\rho_2$ has level 8.  Therefore, $\rho_2$ must be projectively equivalent $\rd{2}{8}{1,0}$, or $\rho_2 \cong \rd{2}{8}{1,\ell}$, where $\ell = 0,3,6,9$. Note that $\rho_1$ is equivalent to its complex conjugation while $\{\rd{2}{8}{1,0}, \rd{2}{8}{1,6}\}$ and  $\{\rd{2}{8}{1,3}, \rd{2}{8}{1,9}\}$ are complex conjugation pairs. It suffices to show that $\rho_2$ is not equivalent to  (i) $\rd{2}{8}{1,0}$ or (ii) $\rd{2}{8}{1,3}$. 

(i) Suppose $\rho_2 \cong \rd{2}{8}{1,0}$. Then  $\spec(\rho_2(\ft)) \subset \spec(\rho_1(\ft))$ and $\rho_1, \rho_2$ have opposite parities. Their direct sum $\tilde\rho = \rd{4}{8}{1,0} \oplus \rd{2}{8}{1,0}$ is given by
$$
\tilde{\rho}(\fs)= {\textstyle\frac{\ii}{\sqrt{8}}} 
\left[
{\small
\begin{array}{cccc}
 1 & \sqrt{3} & \sqrt{3} & 1 \\
 \sqrt{3} & 1 & -1 & -\sqrt{3} \\
 \sqrt{3} & -1 & -1 & \sqrt{3} \\
 1 & -\sqrt{3} & \sqrt{3} & -1 \\
\end{array}}
\right] \oplus 
{\textstyle\frac{1}{\sqrt{2}}} 
\left[
{\small
\begin{array}{cc}
 -1 & 1 \\
 1 & 1 \\
\end{array}}
\right]
\text{ and }  \tilde\rho(\ft) = \diag(\zeta_8,\zeta_8^3, \zeta_8^5, \zeta_8^7, \zeta_8, \zeta_8^3)\,.
$$
However,
$$
\tilde\rho(\fs)^{\nd} = \frac{\ii}{\sqrt{8}} 
\left[
\begin{array}{cc}
 -1 & \sqrt{3} \\
 \sqrt{3} & 1 
\end{array}
\right]
$$
is not a matrix over $\BQ_8$, a contradiction to Proposition \ref{p:ndeg}. Therefore, $\rho_2 \not\cong \rd{2}{8}{1,0}$. 

(ii) Now, we assume $\rho_2 \cong \rd{2}{8}{1,3}$.  Then $\rho_1, \rho_2$ have the same party, and $\tilde\rho = \rd{4}{8}{1,0} \oplus \rd{2}{8}{1,3}$ is given by
$$
\tilde{\rho}(\fs)= {\textstyle\frac{\ii}{\sqrt{8}}} 
\left[
{\small
\begin{array}{cccc}
 1 & \sqrt{3} & \sqrt{3} & 1 \\
 \sqrt{3} & 1 & -1 & -\sqrt{3} \\
 \sqrt{3} & -1 & -1 & \sqrt{3} \\
 1 & -\sqrt{3} & \sqrt{3} & -1 \\
\end{array}}
\right] \oplus 
{\textstyle\frac{\ii}{\sqrt{2}}} 
\left[
{\small
\begin{array}{cc}
 -1 & 1 \\
 1 & 1 \\
\end{array}}
\right]
\text{ and }  \tilde\rho(\ft) = \diag(\zeta_8,\zeta_8^3, \zeta_8^5, \zeta_8^7, \zeta_8^3, \zeta_8^5)\,.
$$
However,
$$
\tilde\rho(\fs)^{\nd} = \frac{\ii}{\sqrt{8}} 
\left[
\begin{array}{cc}
 1 & \sqrt{3} \\
 \sqrt{3} & -1 
\end{array}
\right]
$$
is not a matrix over $\BQ_8$, a contradiction to Proposition \ref{p:ndeg}. Therefore, $\rho_2 \not\cong \rd{2}{8}{1,3}$. 
\end{proof}
\begin{lem} \label{l:9.5}
The level of $\rho_1$ cannot be 9.
\end{lem}
\begin{proof}
There are 4 projectively inequivalent 4-dimensional irreducible $\SL$ representations of level 9, which are given by $\rd{4}{9,1}{1,0}$, $\rd{4}{9,1}{8,0}$, $\rd{4}{9,2}{1,0}$, and $\rd{4}{9,2}{8,0}$ (cf. Appendix \ref{repPP}).
$\rd{4}{9,1}{1,0}$, $\rd{4}{9,1}{8,0}$ are complex conjugate of each other and so are $\rd{4}{9,2}{1,0}$, and $\rd{4}{9,2}{8,0}$. Therefore, it suffices to show that $\rho_1$ cannot be equivalent to (i) $\rd{4}{9,1}{1,0}$ or (ii) $\rd{4}{9,2}{1,0}$.

(i) Suppose $\rho_1 \cong \rd{4}{9,1}{1,0}$, which is odd. By the $\ft$-spectrum criteria, $\rho_2$ can only be projectively equivalent to $\rd{2}{2}{1,0}$ or $\rd{2}{3}{1,0}$, and this implies
$\rho_2 \cong \rd{2}{2}{1,0}$, $\rd{2}{3}{1,0}$ or $\rd{2}{3}{1,8}$. In any of these cases, $\spec( \rho_1(\ft)) \cap \spec(\rho_2(\ft) ) = \{1\}$. Therefore, by Theorem \ref{t:solution1} (iii), $\rho_2$ is also odd, which means $\rho_2 \not\cong \rd{2}{2}{1,0}$ as it is even.

Now $\rho_2 \cong \rd{2}{3}{1,0}$ or $\rd{2}{3}{1,8}$. Note that
\begin{eqnarray*}
\rd{2}{3}{1,0}(\fs) = \frac{\ii}{\sqrt{3}} \left[ 
\arraycolsep=1pt\def\arraystretch{.8} {\small
\begin {array}{cr} -1&\sqrt {2}\\ \noalign{\medskip}\sqrt {2}&
1\end {array}} \right], &\quad& \rd{2}{3}{1,0}(\ft)=\diag(1,\zeta_3), \\
\rd{2}{3}{1,8}(\fs) = \frac{\ii}{\sqrt{3}} \left[ 
\arraycolsep=1pt\def\arraystretch{.8}{\small
\begin {array}{cr} 1&\sqrt {2}\\ \noalign{\medskip}\sqrt {2}& 
-1\end {array}} \right], &\quad&
\rd{2}{3}{1,8}(\ft)=\diag(1,\zeta_3^2)\,.
\end{eqnarray*}
By Theorem \ref{t:solution1} (ii), the unit object $\1$ of $\CC$ is an eigenvector of $\rho(\ft)$ of eigenvalue 1, and $\dim(\CC) = 4/|\rho_1(\fs)_{11}+ \rho_2(\fs)_{11}|^2 = 12$. By the Cauchy Theorem of modular categories, $2 \mid \ord(T) \mid \ord(\rho(\ft) = 9$, a contradiction. Therefore, $\rho_1 \not\cong \rd{4}{9,1}{1,0}$. 

(ii) Now, we assume $\rho_1 \cong \rd{4}{9,2}{1,0}$, which is even. Using similar argument as in Case (i),  $\rho_2 \cong \rd{2}{2}{1,0}$ by  the $\ft$-spectrum criteria and Theorem \ref{t:solution1} (iii).  In this case, $\spec( \rho_1(\ft)) \cap \spec(\rho_2(\ft) ) = \{1\}$ and $\rho$ has level 18.  Theorem \ref{t:solution1} (ii), the unit object of $\CC$ is an eigenvector of $\rho(\ft)$ of eigenvalue 1, and $\dim(\CC)=4/|\rho_1(\fs)_{11}+ \rho_2(\fs)_{11}|^2=16$. By the Cauchy Theorem of modular categories, $\ord(T)$ is a 2-power, but this contradicts Theorem \ref{p:MD1} (4). Therefore, $\rho_1 \not\cong \rd{4}{9,2}{1,0}$.
\end{proof}
\begin{lem}
If $\rho_1$ projectively equivalent to an irreducible representation of prime power level, then the modular data of $\CC$ is a Galois conjugate of that of 
$G(2)_3$.
\end{lem}
\begin{proof}
By Lemmas \ref{l:9.3}, \ref{l:9.4}, \ref{l:9.5} and Appendix \ref{repPP}, $\rho_1$ can only be projective equivalent a level 7 irreducible representation. By the $\ft$-spectrum criteria, $\rho_1 \cong \rd{4}{7}{1}$ or  its complex conjugate $\rd{4}{7}{6}$, which They are  defined over $\BQ_{56}$. 

If there exists some modular data $(S,T)$ whose associated $\SL$ representation $\rho \cong \rd{4}{7}{1} \oplus \rho_2$ for some irreducible 2-dimensional representation $\rho_2$, one can obtain the modular data derived from the MD representation which admits the decomposition  $\rd{4}{7}{6}  \oplus \ol\rho_2$ by the complex conjugation of $(S,T)$.

(I) Assume $\rho_1 \cong \rd{4}{7}{1}$, which is odd. It follows the $\ft$-spectrum criteria that $\rho_2$ must be equivalent to a level 2 or level 3 irreducible representation. In any of these cases, $\spec( \rho_1(\ft)) \cap \spec(\rho_2(\ft) ) = \{1\}$. There is only one 2-dimensional irreducible representation of level 2 which is even. By Theorem \ref{t:solution1} (iii), $\rho_2 \cong \rd{2}{3}{1,0}$ or $\rd{2}{3}{1,8}$, which is odd. Since 
$$
\rd{2}{3}{1,8}  \cong \ol{\rd{2}{3}{1,0}} = \rd{2}{3}{2,0}\,. 
$$
We will solve the modular data for (i) $\rho \cong \rd{4}{7}{1} \oplus \rd{2}{3}{1,0}$ and   (ii) $\rho \cong  \rd{4}{7}{1} \oplus \rd{2}{3}{2,0}$.    

(i) Let $\tilde\rho =\rd{4}{7}{1} \oplus \rd{2}{3}{1,0}$. Then $\tilde\rho(\ft)=\diag(1,\zeta_7, \zeta^2_7, \zeta^4_7,  1, \zeta_3)$ and
$$
\tilde\rho(\fs) =  
\frac{\ii}{\sqrt{7}}\left[
\arraycolsep=1pt\def\arraystretch{.8} {\small
\begin{array}{cccc}
 -1 & \sqrt{2} & \sqrt{2} &\sqrt{2}\\
 \sqrt{2} & \gamma_1 & \gamma_2 & \gamma_3 \\
 \sqrt{2} & \gamma_2 & \gamma_3 & \gamma_1 \\
 \sqrt{2} & \gamma_3 & \gamma_1 & \gamma_2 
 \end{array}}\right] \oplus
 \frac{\ii}{\sqrt{3}} \left[
 \arraycolsep=1pt\def\arraystretch{.8} {\small
\begin{array}{cc}
 -1 &  \sqrt{2} \\
\sqrt{2} & 1 
\end{array}}
\right]
$$
where $\g_1 = -c_7^2$, $\g_2 = -c_7^1$ and $\g_3=-c_7^3$. We reorder $\irr(\CC)$ so that 
$\rho(\ft)= \diag(1,1, \zeta_7, \zeta^2_7, \zeta^4_7, \zeta_3)$, and identify $\irr(\CC)$ with the standard basis of $\BC^6$. By Theorem \ref{t:solution1}, 
$$
\rho(\fs)=\frac{i}{\sqrt{7}}\left[
\arraycolsep=1.5pt\def\arraystretch{1} {\small
\begin{array}{cccccc}
 \frac{-\sqrt{21}-3}{6}  & \frac{(\sqrt{21}-3) \ve_1}{6} & \ve_2 & \ve_3 & \ve_4 & -\sqrt{\frac{7}{3}} \ve_1 \ve_5 \\
 \frac{(\sqrt{21}-3) \ve_1}{6} & \frac{-\sqrt{21}-3}{6} & \ve_1\ve_2 & \ve_1\ve_3 & \ve_1\ve_4 & \sqrt{\frac{7}{3}}\ve_5 \\
 \ve_2 & \ve_1\ve_2 & \g_1 & \g_2 \ve_2 \ve_3 & \g_3 \ve_2 \ve_4 & 0 \\
 \ve_3 & \ve_1\ve_3 & \g_2 \ve_2 \ve_3 & \g_3 & \g_1 \ve_3 \ve_4 & 0 \\
 \ve_4 & \ve_1\ve_4 & \g_3  \ve_2 \ve_4 & \g_1 \ve_3 \ve_4 & \g_2 & 0 \\
 -\sqrt{\frac{7}{3}} \ve_1 \ve_5 & \sqrt{\frac{7}{3}}\ve_5 & 0 & 0 & 0 & \sqrt{\frac{7}{3}} \\
\end{array}}
\right]\,.
$$
for some $\ve_i =\pm 1$, and so $D = 2\left(\frac{1}{\sqrt{3}} + \frac{1}{\sqrt{7}}\right)^{-1}$ or $\dim(\CC) =\frac{21}{2} \left(5-\sqrt{21}\right)$. Since $\frac{21}{2} \left(5+\sqrt{21}\right)$ is a Galois conjugate of $\dim(\CC)$ and 
$$
\dim(\CC) < \frac{21}{2} \left(5+\sqrt{21}\right)  \le \FPdim(\CC),
$$
the objects $\1$ and $\iota$ are distinct. By Theorem \ref{t:solution1} (ii),  $\{\1, \iota\} = \{e_1, e_2\}$, and the modular data of $\CC$ is given by
$$
S =  \left[
\arraycolsep=1.5pt\def\arraystretch{1.2} {\small
\begin{array}{cccccc}
 1 & -d_1 \ve_1 & -d_2 \ve_2 & -d_2 \ve_3 & -d_2 \ve_4 &  d_3 \ve_1 \ve_5\\\
 -d_1 \ve_1 & 1 & - d_2\ve_1 \ve_2 & - d_2\ve_1 \ve_3 & - d_2\ve_1 \ve_4 & -d_3  \ve_5\\
 -d_2 \ve_2 & - d_2\ve_1 \ve_2 & -d_2\g_1 & -d_2\g_2 \ve_2 \ve_3 & -d_2\g_3 \ve_2 \ve_4& 0 \\
 -d_2  \ve_3& - d_2\ve_1 \ve_3 & -d_2\g_2 \ve_2 \ve_3 & -d_2\g_3 & -d_2\g_1 \ve_3 \ve_4 & 0 \\
 -d_2 \ve_4 & - d_2\ve_1 \ve_4 & -d_2\g_3 \ve_2 \ve_4 & -d_2\g_1 & -d_2\g_2 & 0 \\
  d_3 \ve_1\ve_5 & -d_3  \ve_5& 0 & 0 & 0 & -d_3 \\
\end{array}}
\right]  \quad\text{and}\quad T =\rho(\ft)\,, 
$$
where  $d_1=\frac{1}{2} \left(5-\sqrt{21}\right)$, $d_2=\frac{1}{2} \left(\sqrt{21}-3\right)$, $d_3=\frac{1}{2} \left(7-\sqrt{21}\right)$.

If $\1=e_1$, then $\iota=e_2$ and so $S_{2,*} = [d_1, 1, d_2, d_2, d_2, d_3]$. This forces $\ve_1=\ve_5=-1, \ve_2=\ve_3=\ve_4=1$. Thus,
$$
S =  \left[
\arraycolsep=2pt\def\arraystretch{1.2} {\small
\begin{array}{cccccc}
 1 & d_1 & -d_2   & -d_2    & -d_2    &  d_3 \\
 d_1 & 1 &  d_2   &  d_2    &  d_2    &  d_3  \\
 -d_2   &  d_2   & -d_2\g_1 & -d_2\g_2      & -d_2\g_3     & 0 \\
 -d_2    &  d_2    & -d_2\g_2      & -d_2\g_3 & -d_2\g_1       & 0 \\
 -d_2    &  d_2    & -d_2\g_3      & -d_2\g_1 & -d_2\g_2 & 0 \\
  d_3  & d_3  & 0 & 0 & 0 & -d_3 \\
\end{array}}
\right]  \quad\text{and}\quad T =\diag(1,1, \zeta_7, \zeta^2_7, \zeta^4_7, \zeta_3)\,. 
$$

If $\1=e_2$, then $\iota=e_1$ and so $S_{1,*} = [1, d_1, d_2, d_2, d_2, d_3]$. This forces $\ve_i=-1$ for $i=1,\dots, 5$, and so resulting $S$-matrix is equivalent to the preceding one interchanges  the indexes of $e_1$ and $e_2$.

(ii)  Let $\tilde\rho =\rd{4}{7}{1} \oplus \rd{2}{3}{2,0}$. Then $\tilde\rho(\ft)=\diag(1,\zeta_7, \zeta^2_7, \zeta^4_7,  1, \zeta^2_3)$ and
$$
\tilde\rho(\fs) =  
\frac{\ii}{\sqrt{7}}\left[
\arraycolsep=1.5pt\def\arraystretch{.8} {\small
\begin{array}{cccc}
 -1 & \sqrt{2} & \sqrt{2} &\sqrt{2}\\
 \sqrt{2} & \gamma_1 & \gamma_2 & \gamma_3 \\
 \sqrt{2} & \gamma_2 & \gamma_3 & \gamma_1 \\
 \sqrt{2} & \gamma_3 & \gamma_1 & \gamma_2 
 \end{array}}\right] \oplus
 \frac{\ii}{\sqrt{3}} \left[
 \arraycolsep=2pt\def\arraystretch{.8} {\small
\begin{array}{cc}
 1 &  \sqrt{2} \\
\sqrt{2} & -1 
\end{array}}
\right]\,.
$$
Note that $\tilde \rho$ is defined over $\BQ_{168}$. Let $\s \in \Gal(\BQ_{168}/\BQ)$ such that $\s(\zeta_{168})=\zeta_{168}^{113}$. Then $\s|_{\BQ_{56}}= \id$ and $\s(\zeta_3) = \zeta_3^2$. One can see easily that $$
\s\circ (\rd{4}{7}{1} \oplus \rd{2}{3}{1,0}) =\rd{4}{7}{1} \oplus \rd{2}{3}{2,0}\,.
$$
Thus the modular data $(S', T')$ for the MD representation equivalent to $\tilde \rho$ is the Galois conjugate by $\s$ of the modular data $(S,T)$ obtained in (i). Therefore,
$$
S' =  \left[
\arraycolsep=2pt\def\arraystretch{1.2} {\small
\begin{array}{cccccc}
 1 & d'_1 & -d'_2   & -d'_2    & -d'_2    &  d'_3 \\
 d'_1 & 1 &  d'_2   &  d'_2    &  d'_2    &  d'_3  \\
 -d'_2   &  d'_2   & -d'_2\g_1 & -d'_2\g_2      & -d'_2\g_3     & 0 \\
 -d'_2    &  d'_2    & -d'_2\g_2      & -d'_2\g_3 & -d'_2\g_1       & 0 \\
 -d'_2    &  d'_2    & -d'_2\g_3      & -d'_2\g_1 & -d'_2\g_2 & 0 \\
  d'_3  & d'_3  & 0 & 0 & 0 & -d'_3 \\
\end{array}}
\right]  \quad\text{and}\quad T' =\diag(1,1, \zeta_7, \zeta^2_7, \zeta^4_7, \zeta^2_3)\,,
$$
where $d'_1 = \s(d_1) =\frac{1}{2} \left(5+\sqrt{21}\right)$, $d'_2=\s(d_2)= -\frac{1}{2} \left(3+\sqrt{21}\right)$, $d'_3 = \s(d_3) =\frac{1}{2} \left(7+\sqrt{21}\right)$. Since $S'_{1,j} >0$, $e_1=\1 = \iota$, and so $\CC$ is pseudounitary and $\dim(\CC) = \s( \frac{21}{2} \left(5-\sqrt{21}\right)) = \frac{21}{2} \left(5+\sqrt{21}\right)$. 
The modular data of $G(2)_3$ is also $(S',T')$.

(II) Now, we assume $\rho_1=\rd{4}{7}{6}$ and proceed to solve the modular data for (i) $\rho \cong \rd{4}{7}{6} \oplus \rd{2}{3}{1,0}$ and (ii) $\rho \cong \rd{4}{7}{6} \oplus \rd{2}{3}{2,0}$. Note that both of them are defined over $\BQ_{168}$. 

(i) Let $\tilde\rho = \rd{4}{7}{6} \oplus \rd{2}{3}{1,0}$. Then $\ol{\rd{4}{7}{1} \oplus \rd{2}{3}{2,0}} = \tilde\rho$. Thus the modular data $(S'', T'')$ of the MD representations equivalent to $\tilde \rho$ is $(\ol S', \ol T')$, which is given by
$$
S'' = S' = \left[
\arraycolsep=2pt\def\arraystretch{1.2} {\small
\begin{array}{cccccc}
 1 & d'_1 & -d'_2   & -d'_2    & -d'_2    &  d'_3 \\
 d'_1 & 1 &  d'_2   &  d'_2    &  d'_2    &  d'_3  \\
 -d'_2   &  d'_2   & -d'_2\g_1 & -d'_2\g_2      & -d'_2\g_3     & 0 \\
 -d'_2    &  d'_2    & -d'_2\g_2      & -d'_2\g_3 & -d'_2\g_1       & 0 \\
 -d'_2    &  d'_2    & -d'_2\g_3      & -d'_2\g_1 & -d'_2\g_2 & 0 \\
  d'_3  & d'_3  & 0 & 0 & 0 & -d'_3 \\
\end{array}}
\right]  \quad\text{and}\quad T'' = \ol T' =\diag(1,1, \zeta^6_7, \zeta^5_7, \zeta^3_7, \zeta_3)\,,
$$
In particular, the MTC $\CC$ is also pseudounitary with $\dim(\CC)=\frac{21}{2} \left(5+\sqrt{21}\right)$.

(ii) Finally, we consider $\tilde\rho = \rd{4}{7}{6} \oplus \rd{2}{3}{2,0}$ which is the complex conjugate of $\rd{4}{7}{1} \oplus \rd{2}{3}{1,0}$. Thus the modular data $(S''', T''')$ of the MD representations equivalent to $\tilde \rho$ is $(\ol S, \ol T)$, which is given by
$$
S'''=S =  \left[
\arraycolsep=2pt\def\arraystretch{1.2} {\small
\begin{array}{cccccc}
 1 & d_1 & -d_2   & -d_2    & -d_2    &  d_3 \\
 d_1 & 1 &  d_2   &  d_2    &  d_2    &  d_3  \\
 -d_2   &  d_2   & -d_2\g_1 & -d_2\g_2      & -d_2\g_3     & 0 \\
 -d_2    &  d_2    & -d_2\g_2      & -d_2\g_3 & -d_2\g_1       & 0 \\
 -d_2    &  d_2    & -d_2\g_3      & -d_2\g_1 & -d_2\g_2 & 0 \\
  d_3  & d_3  & 0 & 0 & 0 & -d_3 \\
\end{array}}
\right]  \quad\text{and}\quad T'''=\ol T =\diag(1,1, \zeta^6_7, \zeta^5_7, \zeta^3_7, \zeta^2_3)\,. 
$$
Therefore, the MTC $\CC$ is not pseudounitary and $\dim(\CC)=\frac{21}{2} \left(5-\sqrt{21}\right)$.
\end{proof}
\begin{lem}
The level of $\rho_1$ cannot be 6, 10 or 40.
\end{lem}
\begin{proof} 

(i)  Suppose $\rho_1$ is of level 6. Then $\rho_1 \cong \psi \o \eta$ for some 2-dimensional irreducible representations $\psi$ and $\eta$ of level 2 and 3 respectively. 
There is only one 2-dimensional irreducible representation, up to projective equivalence, of levels 2 and 3. Since the $\ft$-spectrum of $\rho_1$ is minimal,  $\rho_1 \cong \rd{2}{2}{1,0} \o\rd{2}{3}{1,0}$ . In particular, $\spec(\rho_1(\ft)) = \{1, -1, \zeta_3, -\zeta_3\}$. By the $\ft$-spectrum criteria, $\rho_2$ can only be equivalent to $\rho_2 \cong \rd{2}{2}{1,i}$, $i \in \{0, 4, 6, 10\}$, or $\rd{2}{3}{1,j}$, $j$ even.  Therefore, $\ord(\rho_2(\ft)) \mid 6$ and so $\ord(\rho(\ft))=6$. This implies $\ord(T) \mid 6$ and so $\CC$ is integral by Theorem \ref{t:2346}. However,  this contradicts Proposition \ref{p:weakly_int}. Therefore, the level of $\rho_1$ cannot be  6. 

(ii) Suppose $\rho_1$ is of level 40. Then $\rho_1$ is projectively equivalent to $\rd{2}{8}{1,0} \o \rd{2}{5}{1}$ or $\rd{2}{8}{1,0} \o \rd{2}{5}{2}$ (cf. Appendix \ref{repPP}).  In particular, $\spec(\rho_1(\ft))$ is a set of primitive $40$-th roots of unity. However, there does not exist any 2-dimensional representation $\rho_2$ which satisfies the $\ft$-spectrum criteria. Therefore, the level $\rho_1$ cannot be 40.

(iii)  Suppose $\rho_1$ is of level 10. Then $\rho_1$ is projectively equivalent to $\rd{2}{2}{1,0} \o \rd{2}{5}{1}$ or $\rd{2}{2}{1,0} \o \rd{2}{5}{2}$. Since $\rd{2}{2}{1,0}$ is equivalent to any of it Galois conjugates,  $\rd{2}{2}{1,0} \o \rd{2}{5}{1}$ or $\rd{2}{2}{1,0} \o \rd{2}{5}{2}$ are  Galois conjugate. So, it suffices to show that $\rho_1 \cong \rd{2}{2}{1,0} \o \rd{2}{5}{1}$ is not possible. 

Assume $\rho_1 \cong \rd{2}{2}{1,0} \o \rd{2}{5}{1}$. Then $\spec(\rho_1(\ft)) = \{\zeta_5,\zeta_5^4, -\zeta_5, -\zeta^4_5\}$. By the $\ft$-spectrum criteria, $\rho_2 \cong \rd{2}{5}{1}$ or $\chi^6\o \rd{2}{5}{1}$.  Since $\chi^6 \o \rd{2}{2}{1,0} \cong \rd{2}{2}{1,0}$, $\rho_1 \oplus \rd{2}{5}{1}$  and $\rho_1 \oplus \chi^6\o \rd{2}{5}{1}$ are projectively equivalent. Therefore,  $\rho$ is projectively equivalent to $\tilde\rho= (\rd{2}{2}{1,0} \o \rd{2}{5}{1}) \oplus \rd{2}{5}{1}$ and we can simply assume $\rho \cong \tilde \rho$. As in Lemma \ref{l:9.3}, we the use the following equivalent form of $\rd{2}{5}{1}$:
$$
\rd{2}{5}{1}(\fs)=\frac{1}{s_5^1}\left[
\arraycolsep=3pt
\begin{array}{cr}
 1 & \varphi  \\
 \varphi  & {-1} 
\end{array}
\right],\quad \rd{2}{5}{1}(\ft)  =\diag(\zeta_5, \zeta_5^4)\,.
$$
Thus, $\tilde\rho(\ft) = \diag(\zeta_5, \zeta_5^4, -\zeta_5, -\zeta_5^4, \zeta_5, \zeta_5^4)$ and 
$$
\tilde\rho(\fs)=\frac{1}{2 s_5^1 }\left[
\arraycolsep=2pt\def\arraystretch{1.2} {\small
\begin{array}{cccccc}
 -1 & -\varphi  & -\sqrt{3} & -\sqrt{3} \varphi  \\
 -\varphi  & 1 &  -\sqrt{3} \varphi  & \sqrt{3} \\
 -\sqrt{3} &  -\sqrt{3} \varphi  & 1 & \varphi  \\
 - \sqrt{3} \varphi  &  \sqrt{3} &  \varphi  & -1 
\end{array}}
\right]\oplus \frac{1}{s_5^1}\left[
\arraycolsep=2pt\def\arraystretch{1.2} {\small
\begin{array}{cr}
 1 & \varphi  \\
 \varphi  & {-1} 
\end{array}}
\right]\,.
$$
By  Theorem \ref{t:solution2} (i), if we reorder $\irr(\CC)$ so that $\rho(\ft) = \diag(\zeta_5\, I_2, \zeta_5^4\, I_2, -\zeta_5, -\zeta_5^4)$, then  $\rho(\fs) = V s' V$ for some signed diagonal matrix $V$ and
$$
s'=\frac{1}{2 s_5^1} \left[
\arraycolsep=3pt\def\arraystretch{1.2} {\small
\begin{array}{cccccc}
 * & * & *&* & -\sqrt{3} a  & -\sqrt{3} a \varphi  \\
  * & * & *&*  & -\sqrt{3} b  & -\sqrt{3} b \varphi  \\
  * & * & *&* & -\sqrt{3} c \varphi   & \sqrt{3} c \\
  * & * & *&* & -\sqrt{3} d \varphi  & \sqrt{3} d \\
 -\sqrt{3} a & -\sqrt{3} b  & -\sqrt{3} c \varphi  & -\sqrt{3} d \varphi   & 1 & \varphi   \\
 -\sqrt{3} a \varphi  & -\sqrt{3} b \varphi  & \sqrt{3} c & \sqrt{3} d & \varphi & -1 \\
\end{array}}
\right]\,.
$$
where $a,b,c,d \in \BR$ satisfying $a^2+b^2=1$ and $c^2+d^2=1$.

Note that $\varphi$ is a unit in $\BZ[\zeta_5]$, and the automorphism $\s$ defined by $\s(\zeta_{10})=\zeta_{10}^7$ generates $\Gal(\BQ_{10})$. By the action of $\s^2$ on $\rho(\ft)$, we see $\hs(5)=6$.  Since 
$$
\s(s'_{5,5}) = \frac{1}{2 s_5^2} =  \frac{\varphi}{2 s_5^1} =s'_{56}\,,
$$
$\s(s'_{i,5}) = s'_{i,6}$ for $i=1,\dots,6$. This implies $\sqrt{3} a, \sqrt{3} b, \sqrt{3} c, \sqrt{3} c$ are fixed by $\s$ and so they are rational. 

The unit object cannot be $e_5$, for otherwise $\sqrt{3}a, \sqrt{3} b \in \BZ$ and they satisfy the equation $(\sqrt{3}a)^2+ (\sqrt{3} b)^2 = 3$, which is not possible. Similarly, $e_6 \ne \1$. So, the unit object $\1 \in \{e_1 , e_2, e_3, e_4\}$. 

Assume $\1 =e_1$. Then $a \ne 0$, $b/a \in \BZ$ and $\frac{1}{\sqrt{3}a} \in \BZ$. However, this will imply $3 \mid (1 + (b/a)^2)$ which is not possible. Therefore, $\1 \ne e_1$. Since $\varphi$ is a unit in $\BZ[\zeta_5]$, if $\1 \not \in \{e_2, e_3, e_4\}$ for similar reason. Now, we find $1 \not \in  \{e_1, \dots, e_6\}$, a contradiction. Therefore, the level of $\rho_1$ cannot be 10.
\end{proof}

\begin{lem}
If the level of $\rho_1$ is 24, then $\CC$ is equivalent to $\CC(\BZ_6, q)$ for some non-degenerate quadratic $q:\BZ_6 \to \BC^\times$. 
\end{lem}
\begin{proof}
Since $\rho_1$ is of level 24, $\rho_1$ is projectively equivalent to $\rd{2}{3}{1,0} \o \rd{2}{8}{1,0}$ according to Appendix \ref{repPP}. Therefore, we can simply assume $\rho_1\cong \rd{2}{3}{1,0} \o \rd{2}{8}{1,0}$ as it has a minimal $\ft$-spectrum.  Then,  $\rho_1$ is odd and 
$$
\rho_1(\fs) = \frac{\ii}{\sqrt{6}}\left[
\begin{array}{cccc}
 1 & -1 & -\sqrt{2} & \sqrt{2} \\
 -1 & -1 & \sqrt{2} & \sqrt{2} \\
 -\sqrt{2} & \sqrt{2} & -1 & 1 \\
 \sqrt{2} & \sqrt{2} & 1 & 1 \\
\end{array}
\right], \quad\rho_1(\ft)=\diag(\zeta_8, \zeta_8^3, \zeta_{24}^{11}, \zeta_{24}^{17})\,.
$$
By the $\ft$-spectrum criteria, $\rho_2 \cong \rd{2}{8}{1,j}$, $j \in \{0, 1,3,4,7,9\}$, and
$$
\rd{2}{8}{1,j}(\fs) = \frac{(-i)^j}{\sqrt{2}}\left[
\begin{array}{cc}
-1 & 1 \\
1 & 1
\end{array}
\right],\quad \quad\rd{2}{8}{1,j}(\ft)=\diag(\zeta_{24}^{3+2j}, \zeta_{24}^{9+2j})\,.
$$
For $j=1,3,7,9$, $|\spec(\rho_1(\ft)) \cap \spec(\rho_2(\ft))|=1$ and so Theorem \ref{t:solution1} can be applied.

For $j=1,9$, $\spec(\rho_1(\ft)) \cap \spec(\rho_2(\ft)) = \{\zeta_{24}^{9+2j}\}$, and for $j=3,7$, $\spec(\rho_1(\ft)) \cap \spec(\rho_2(\ft)) = \{\zeta_{24}^{3+2j}\}$.  If $\rho \cong \rho_1 \oplus \rd{2}{8}{1,j}$ is an MD representation of an MTC $\CC$, for $j=1,3,7,9$, then by Theorem  \ref{t:solution1}, $\ord(T) = 12$ and
$$
D=\sqrt{\dim(\CC)}=\frac{2}{\frac{1}{\sqrt{2}}(1 \pm \frac{1}{\sqrt{3}})}=\sqrt{6}(\sqrt{3} \mp 1).
$$
Note that each row of $\rho_1(\fs)$ has an off diagonal entry of the form $\frac{\pm i}{\sqrt{6}}$ and so $\frac{D }{\sqrt{6}}/\sqrt{2}$ is the dimension of an object up to a sign. However, 
$$
\frac{D }{\sqrt{6} \sqrt{2}} = \frac{ \sqrt{3}\pm 1}{\sqrt{2}}  \not\in \BQ_{12}\,.
$$
Therefore, $\rho_1 \oplus \rd{2}{8}{1,j}$ is not equivalent to any MD representation for $j =1,3,7,9$.

Now, we can conclude that $\rho \cong \rho_1 \oplus \rho_2$ where $\rho_2 \cong \rd{2}{8}{1,j}$ for some $j=0,4$. In particular, $\rho_1$ and $\rho_2$ have opposite parties and $\spec(\rho_2(\ft))\subset \spec(\rho_1(\ft))$. By Theorem \ref{t:solution2} (ii), the unit object $\1$ is an eigenvector of $\rho(\ft)$ with eigenvalue $\zeta \in \spec(\rho_1(\ft))\setminus \spec(\rho_2(\ft))$. Let $E_j$ be the subspace of $\BC^6$ spanned by the eigenvectors of $\tilde\rho_j=\rho_1(\ft) \oplus \rd{2}{8}{1,j}(\ft)$ with eigenvalues in $\spec(\rho_1(\ft))\setminus \spec(\rd{2}{8}{1,j}(\ft))$ for $j=0,4$. One can compute that for $\s \in \Gal(\BQ_{24}/\BQ)$,  $D_{\tilde\rho_j}(\sigma)|_{E_j} =\id$ or $-\id$. By Proposition \ref{p:5.2}, $\CC$ is integral. It follows from \cite{CWeakly} that $\CC$ is a pointed modular tensor category, which is equivalent to to $\CC(\BZ_6, q)$ for some non-degenerate quadratic $q:\BZ_6 \to \BC^\times$. 
\end{proof}
\begin{lem} \label{l:4.10}
If the level of $\rho_1$ is 15, then  the modular data of $\CC$ is a Galois conjugate of that of $\CC(\BZ_3, q) \boxtimes PSU(2)_3$, where $q:\BZ_3 \to \BC^\times$ is a quadratic form given by  $q(1)=\zeta_3$.
\end{lem}
\begin{proof}
Since $\rho_1$ has a minimal $\ft$-spectrum, it must be equivalent to a tensor product of two 2-dimensional irreducible representations of levels 3 and 5. According to Appendix \ref{repPP}, $\rho_1 \cong \rd{2}{3}{1,0} \o \rd{2}{5}{i}$,  $i=1,2$. By the $\ft$-spectrum criteria, $\rho_2 \cong \chi^j \o \rd{2}{5}{i}$ with $j = 0,4$. Thus, $\rho$ is equivalent to 
$$
\tilde\rho_{i,j}  =(\rd{2}{3}{1,0} \o \rd{2}{5}{i})  \oplus (\chi^j \o \rd{2}{5}{i}) , \quad i=1,2, \, j=0,4\,.
$$
Note that $\tilde\rho_{i,j}$ is defined over $\BQ_{120}$ for $i,j$. Let $\s_a \in \Gal(\BQ_{120}/\BQ)$ such that $\s_a (\zeta_{120})=\zeta^a_{120}$. Then,  $\s_{97}\circ \tilde\rho_{1,j} = \tilde\rho_{2,j}$ for $j=0,4$.

Since 
$
\s_{41}\circ ((\rd{2}{3}{1,0} \o \rd{2}{5}{1})  \oplus \rd{2}{5}{1})\cong  (\ol{\rd{2}{3}{1,0}} \o \rd{2}{5}{1})  \oplus \rd{2}{5}{1} \cong (\rd{2}{3}{1,8} \o \rd{2}{5}{1})  \oplus \rd{2}{5}{1}\,,
$
we  have $\chi^4 \o \s_{41}\circ \tilde\rho_{1,0} \cong \tilde\rho_{1,4}$.
Therefore, $\tilde\rho_{i,j}$ is projectively equivalent to a Galois conjugate of $\tilde\rho_{1,0}$. Hence, it suffices to consider $\tilde\rho = \tilde\rho_{1,0}$, or equivalently $\rho_1 \cong \rd{2}{3}{1,0} \o \rd{2}{5}{1}$ and $\rho_2 \cong \rd{2}{5}{1}$

Now, the MD representation $\rho$ of $\CC$ is equivalent to $\rho_1 \oplus \rho_2$, where
$\rho_1$ is even, $\rho_2$ is odd and  $\spec(\rho_1(\ft)) \subset \spec(\rho_2(\ft))$. Moreover, 
$$
\tilde\rho(\fs) =\frac{\ii}{\sqrt{3} s_5^1}
\left[
\arraycolsep=2pt\def\arraystretch{1.2} {\small
\begin{array}{cccc}
 -1 & -\varphi  & \sqrt{2} & \sqrt{2} \varphi  \\
 -\varphi  & 1 & \sqrt{2} \varphi  & -\sqrt{2} \\
 \sqrt{2} & \sqrt{2} \varphi  & 1 & \varphi  \\
 \sqrt{2} \varphi  & -\sqrt{2} & \varphi  & -1 \\
\end{array}}
\right] \oplus \frac{1}{s_5^1}\left[
\arraycolsep=2pt\def\arraystretch{1.2} {\small
\begin{array}{cr}
 1 & \varphi  \\
 \varphi  & {-1} 
\end{array}}
\right], \quad \tilde\rho(\ft)=\diag(\zeta_5, \zeta^4_5, \zeta^8_{15}, \zeta^2_{15}, \zeta_5, \zeta^4_5)\,.
$$

By Theorem \ref{t:solution2}, $\dim(\CC)=12 \sin^2(2 \pi /5) = 3(2+\varphi)$. Reorder $\irr(\CC)$ so that 
$$
\rho(\ft)= \diag(\zeta_5, \zeta^4_5, \zeta_5, \zeta^4_5, \zeta^8_{15}, \zeta^2_{15}).
$$
Again, by Theorem \ref{t:solution2} (ii), there exist $\g_i, \kappa_i, \ve_i \in \{\pm 1\}$ such that
$$
\rho(\fs)=\frac{-1}{D}  \left[
\arraycolsep=.9pt\def\arraystretch{1.1} {\small
\begin{array}{cccccc}
 \frac{ 1+\ii \sqrt{3}}{2} & \frac{ \left(1-\ii \sqrt{3}\right) \kappa _1}{2} & \frac{\g_3 \varphi  \left(1+\ii \sqrt{3} \ve_1 \ve_2 \kappa _1 \kappa _2\right)}{2} & \frac{\g_3 \varphi  \left(\kappa _2-\ii \sqrt{3} \ve_1 \ve_2 \kappa _1\right)}{2} & -\g_1 & -\g_2 \varphi \\
 \frac{ \left(1-\ii \sqrt{3}\right) \kappa _1}{2} & \frac{ 1+\ii \sqrt{3}}{2} & \frac{\g_3 \varphi  \left(\kappa _1-\ii \sqrt{3} \ve_1 \ve_2 \kappa _2\right)}{2} & \frac{ \g_3 \varphi  \left(\kappa _1 \kappa _2+\ii \sqrt{3} \ve_1 \ve_2\right)}{2} & -\g_1 \kappa _1 & -\g_2 \kappa_1 \varphi \\
 \frac{ \g_3 \varphi  \left(1+\ii \sqrt{3} \ve_1 \ve_2 \kappa _1 \kappa _2\right)}{2} & \frac{\g_3 \varphi  \left(\kappa _1-\ii \sqrt{3} \ve_1 \ve_2 \kappa _2\right)}{2} & \frac{- \left(1+\ii\sqrt{3}\right)}{2} & \frac{\left(-1+\ii\sqrt{3}\right) \kappa _2}{2} & -\g_1 \g_3 \varphi  & \g_2 \g_3 \\
 \frac{\g_3 \varphi  \left(\kappa _2-\ii \sqrt{3} \ve_1 \ve_2 \kappa _1\right)}{2} & \frac{ \g_3 \varphi  \left(\kappa _1 \kappa _2+\ii \sqrt{3} \ve_1 \ve_2\right)}{2} & \frac{\left(-1+\ii\sqrt{3}\right) \kappa _2}{2} & \frac{- \left(1+\ii\sqrt{3}\right)}{2} & -\g_1 \g_3 \kappa _2 \varphi  & \g_2 \g_3 \kappa _2 \\
 -\g_1 & -\g_1 \kappa _1 & -\g_1 \g_3 \varphi & -\g_1 \g_3 \kappa _2 \varphi & -1 & -\g_1 \g_2 \varphi \\
 -\g_2 \varphi & -\g_2 \kappa_1 \varphi  & \g_2 \g_3 & \g_2 \g_3 \kappa _2 & -\g_1 \g_2 \varphi & 1 \\
\end{array}}
\right]\,. 
$$

We will use the equalities $\frac{ -1+\ii \sqrt{3}}{2}=\zeta_3$ and  $\frac{1+\ii \sqrt{3}}{2}=-\ol \zeta_3$ to simplify $S$-matrix, but we need to determine which of the standard basis elements is the unit object. According to Theorem \ref{t:solution2} (ii), $\1 \in \{e_5,e_6\}$.

(i) Suppose $e_6 = \1$.  Then $T=\diag(\zeta_{15},\zeta_3^2, \zeta_{15}, \zeta_3^2,\zeta_5^2, 1)$. 
Then $\dim(e_5)^2 =\varphi^2 >1$ and so $e_6 = \iota$. Thus, all the entries of 6-th rows of $\rho(\fs)$ has the same signed, we find $\g_2 =\g_3=-1, \g_1= \kappa_1=\kappa_2=1$. Thus,  
$$
S=\left[
\arraycolsep=2pt\def\arraystretch{1.3} {\small
\begin{array}{cccccc}
 -\ol\zeta_3 & -\zeta_3 & \varphi \ol\zeta_3 & \varphi \zeta_3 & -1 & \varphi \\
 -\zeta_3 & -\ol\zeta_3 & \varphi \zeta_3 & \varphi \ol\zeta_3 & -1  & \varphi \\
 \varphi \ol\zeta_3 & \varphi \zeta_3 & \ol\zeta_3 & \zeta_3 & \varphi  & 1 \\
 \varphi \zeta_3 & \varphi \ol\zeta_3 & \zeta_3 & \ol\zeta_3 & \varphi  & 1 \\
 -1 & -1  & \varphi & \varphi & -1 &  \varphi \\
  \varphi &   \varphi  & 1 & 1 & \varphi & 1 \\
\end{array}}
\right]\,. 
$$
By reordering $\irr(\CC)$, we find $T=\diag(1,\zeta_5^2, \zeta^2_3, \zeta_{15},  \zeta^2_3, \zeta_{15}) = T_1 \o T_2$ and
$$
S = \left[
\arraycolsep=2pt\def\arraystretch{1.3} {\small
\begin{array}{cccccc}
1 & \varphi & 1  &\varphi  & 1 & \varphi \\
 \varphi & -1 & \varphi  & -1 & \varphi & -1 \\
 1 & \varphi  & \ol\zeta_3  & 
 \varphi \ol \zeta_3 & \zeta_3  & \varphi \zeta_3\\
 \varphi  & -1 &   \varphi \ol\zeta_3 
&- \ol\zeta_3 &   \varphi \zeta_3 & -\zeta_3  \\
 1 & \varphi  & \zeta_3     &  
 \varphi \zeta_3 & \ol\zeta_3 & \varphi \ol \zeta_3 \\ 
 \varphi  & -1 &   \varphi \zeta_3 
 & -\zeta_3 &   \varphi \ol\zeta_3  & -\ol\zeta_3   \\
\end{array} }
\right] = S_1 \o S_2, 
$$
where $(S_1, T_1)$ is the modular data of $\CC(\BZ_3, \ol q)$ and  $(S_2, T_2)$ given by  $S_2 = \mtx{1 & \varphi \\ \varphi & -1}$ and $T_2=\diag(1, \zeta_5^2)$ is the modular data of $PSU(2)_3$.  In particular, $(S,T)$ is a Galois conjugate of the modular data of $\CC(\BZ/3\BZ, q_1) \boxtimes PSU(2)_3$.

(ii) Now, we assume $e_5 = \1$. Then $T= \diag(\zeta_3^2, \zeta_{15}^4, \zeta_3^2, \zeta_{15}^4, 1, \zeta_5^3)$ and $\dim(e_6)^2=\varphi^2>1$, and so $e_5 = \iota$. Then $\g_1 =\g_2=\g_3 =\kappa_1=\kappa_2=1$, and we obtain 
$$
S=\left[
\arraycolsep=2pt\def\arraystretch{1.3} {\small
\begin{array}{cccccc}
 \ol\zeta_3 & \zeta_3 & \ol\zeta_3 \varphi  & \zeta_3 \varphi  & 1 & \varphi  \\
 \zeta_3 & \ol\zeta_3 & \zeta_3 \varphi  & \ol\zeta_3 \varphi  & 1 & \varphi  \\
 \ol\zeta_3 \varphi  & \zeta_3 \varphi  & -\ol\zeta_3 & -\zeta_3 & \varphi  & -1 \\
 \zeta_3 \varphi  & \ol\zeta_3 \varphi  & -\zeta_3 & -\ol\zeta_3 & \varphi  & -1 \\
 1 & 1 & \varphi  & \varphi  & 1 & \varphi  \\
 \varphi  & \varphi  & -1 & -1 & \varphi  & -1 \\
\end{array}
}
\right]\,.
$$
By reordering $\irr(\CC)$, we find $T=\diag(1,\zeta_5^3, \zeta^2_3, \zeta_{15}^4,  \zeta^2_3, \zeta_{15}^4) = T_1 \o \ol T_2$ and 
$$
S = \left[
\arraycolsep=2pt\def\arraystretch{1.3} {\small
\begin{array}{cccccc}
1 & \varphi & 1  &\varphi  & 1 & \varphi \\
 \varphi & -1 & \varphi  & -1 & \varphi & -1 \\
 1 & \varphi  & \ol\zeta_3  & 
 \varphi \ol \zeta_3 & \zeta_3  & \varphi \zeta_3\\
 \varphi  & -1 &   \varphi \ol\zeta_3 
&- \ol\zeta_3 &   \varphi \zeta_3 & -\zeta_3  \\
 1 & \varphi  & \zeta_3     &  
 \varphi \zeta_3 & \ol\zeta_3 & \varphi \ol \zeta_3 \\ 
 \varphi  & -1 &   \varphi \zeta_3 
 & -\zeta_3 &   \varphi \ol\zeta_3  & -\ol\zeta_3   \\
\end{array} }
\right] = S_1 \o S_2, 
$$
Since $(S_2, \ol T_2)$ is the complex conjugate of modular data of $PSU(2)_3$. Therefore, $(S,T)$ is a Galois conjugate of the modular data of $\CC(\BZ_3, q) \boxtimes PSU(2)_3$.  This completes the proof of statement.

As a consequence, for any $i,j$, $\tilde\rho_{i,j}$ is  equivalent to  $\SL$ representations of some modular tensor categories Galois conjugate to $\CC(\BZ_3,q) \boxtimes PSU(2)_3$. 
\end{proof}

\begin{proof}[Proof of Theorem \ref{t:42}] The result of Theorem \ref{t:42} is a consequence of Lemmas  \ref{l:9.3} to \ref{l:4.10}.
\end{proof}

\subsection{Classification of modular data of type (3,3)}

\begin{thm}\label{t:33}
The modular data of any type $(3,3)$ modular tensor category is a Galois conjugate of that of $SO(5)_2$.
\end{thm}

Let $\CC$ be a modular tensor category of type (3,3) and $\rho$ an $\SL$  representation of $\CC$. Then 
$$
\rho \cong \rho_1 \oplus \rho_2
$$
for some 3-dimensional irreducible representations $\rho_1, \rho_2$. If $\rho_1$, $\rho_2$ have opposite parities, then $\Tr(\rho(\fs))=0$ which contradicts to Proposition \ref{p:fixpt}. Therefore, they have the same parity.  We may assume that $\rho_1$ has a minimal $\ft$-spectrum and show that for $\rho_1$ cannot be projectively equivalent of any 3-dimensional irreducible representation of levels 3, 7, 8 or 16. 
\begin{lem}
Neither $\rho_1$ nor $\rho_2$ is projectively equivalent to a 3-dimensional irreducible representation   of  level $3$, $7$, $8$ or $16$. 
\end{lem} \label{l:4.12}
\begin{proof}
Suppose $\rho_1$ is a 3-dimensional irreducible representation of level 3, 7, 8 or 16 with a minimal $\ft$-spectrum.

(i) {\bf $\rho_1$ cannot be of level $7$}: Suppose  $\rho_1$ is of level $7$. Then,  by the  $\ft$-spectrum criteria and Appendix \ref{repPP},  $\rho_2 \cong \rho_1$ but this contradicts Proposition \ref{p:multiple}. 

(ii) {\bf $\rho_1$ cannot be of level $3$}: Suppose $\rho_1$ is the level 3. Then $\rho_1 \cong \rd{3}{3}{1,0}$. If $\rho_2$ is projectively equivalent to $\rd{3}{3}{1,0}$, then $\rho_2 \cong \rd{3}{3}{1,0}$ by the $\ft$-spectrum criteria, but this contradicts  Proposition \ref{p:multiple}. Therefore, $\rho_2$ is not projectively equivalent to $\rd{3}{3}{1,0}$.

 By the $\ft$-spectrum criteria, $\rho_2$ is not projectively equivalent to any level 16 irreducible representations. If $\rho_2$ is projectively equivalent a level 8 irreducible representation, then $\rho_2 \cong \chi^j\o \psi$ for any level 8 representations in Appendix \ref{repPP}.  Since $\rho_1$ is even,  $j \equiv 0 \mod 4$, and so $\psi$  must be even.   This implies $\psi =\rd{3}{8}{3,3},\, \rd{3}{8}{1,3},\, \rd{3}{8}{3,9},\, \rd{3}{8}{1,9}$, but none of them satisfies the $\ft$-spectrum criteria.  Therefore, $\rho_2$ can only be projectively equivalent to some $\psi$ of level 5 or 4 in Appendix \ref{repPP}. Thus, by the $\ft$-spectrum criteria, $\rho_2 \cong \chi^j \o \rd{3}{4}{1,0}$ or $\chi^j \o \rd{3}{5}{i}$ for $j=0,4,8$ and $i=1,3$. In any of these cases, $|\spec(\rho_1(\ft)) \cap \spec(\rho_2(\ft))|=1$ and $\ord(\rho(\ft)) =12$ or $15$. It follows from Theorem \ref{t:solution1} (ii) that if $\spec(\rho_1(\ft)) \cap \spec(\rho_2(\ft)) = \{\rho_1(\ft)_{u,u}\}$, then $\frac{\sqrt{2} \rho_1(\fs)_{jj}}{\rho_1(\fs)_{uj}} \in \BQ_{12}$ or $\BQ_{15}$ for $u \ne j$.  However, $\frac{\sqrt{2} \rho_1(\fs)_{jj}}{\rho_1(\fs)_{uj}} = \frac{-1}{\sqrt{2}} \not\in \BQ_{12}$ or $\BQ_{15}$.
 Therefore, $\rho_2$ cannot be projectively equivalent to any irreducible of level 4 or 5. This completes the proof that $\rho_1$ cannot be of level $3$.

(iii) {\bf $\rho_1$ cannot be of level $8$}: Let
$$
A=\frac{\ii}{2}\left[
\arraycolsep=2pt\def\arraystretch{1.3} {\small
\begin{array}{ccc}
 0 & \sqrt{2} & \sqrt{2} \\
 \sqrt{2} & -1 & 1 \\
 \sqrt{2} & 1 & -1
 \end{array}}
\right].
$$
Then, by Appendix \ref{repPP}, 
$$
\rd{3}{8}{1,0}(\fs) = A
\quad \text{and} \quad\rd{3}{8}{1,0}(\ft)=\diag(1, \zeta_8, \zeta_8^5 )
$$
which is odd and has a minimal $\ft$-spectrum. Since all other 4-dimensional level 8 irreducible representations are projectively equivalent to a Galois conjugate of $\rd{3}{8}{1,0}$, it suffices to show that $\rho_1\not\cong\rd{3}{8}{1,0}$. 

Assume contrary. Then $\rho_1 \cong\rd{3}{8}{1,0}$, and hence $\rho_2$ must be odd. It follows from (i) and (ii), $\rho_2$ cannot  be projectively equivalent to any irreducible representation of level 3 or 7. 
By the $\ft$-spectrum criteria and the parity constraint, $\rho_2$ cannot be projectively equivalent to any irreducible representations of level 5. Therefore, $\rho_2$ can only be projectively equivalent to an irreducible representation of level 4, 8 or 16. By the $\ft$-spectrum criteria, $\rho_2$ is of level 4, 8 or 16.

Suppose $\rho_2$ has level 4 or 8.  Since $\rho_2$ is odd, $\rho_2 \cong \rd{3}{4}{1,3}, \rd{3}{4}{1,9}, \rd{3}{8}{1,0}, \rd{3}{8}{3,0}, \rd{3}{8}{1,6}, \rd{3}{8}{3,6}$.   However, $D_{\rho_1 \oplus \rho_2}(\s) = \pm \id$ for all $\s \in \Gal(\BQ_8/\BQ)$. By Proposition \ref{p:5.2}, $\CC$ is integral which contradicts Proposition \ref{p:weakly_int}. Therefore, the level of $\rho_2$ is neither 4 nor 8. 

Suppose $\rho_2$ is an odd irreducible representation  of level 16. By the $\ft$-spectrum criteria, $\rho_2 \cong \rd{3}{16}{1,0}\,, \rd{3}{16}{5,6}\,, \rd{3}{16}{1,6}\,, \rd{3}{16}{5,0}$, and they are respectively isomorphic to the following representations:
\begin{enumerate}
    \item[(1)] $\fs \mapsto A$,\quad  $\ft \mapsto \diag(\zeta_8, \zeta_{16}, \zeta_{16}^{9} )$ ;
    \item[(2)] $\fs \mapsto -A$,\quad  $\ft \mapsto \diag(\zeta_8, \zeta_{16}^5, \zeta_{16}^{13} )$ ;
    \item[(3)] $\fs \mapsto -A$,\quad  $\ft \mapsto \diag(\zeta_8^5, \zeta_{16}, \zeta_{16}^9 )$ ;
    \item[(4)] $\fs \mapsto A$,\quad  $\ft \mapsto \diag((\zeta_8^5, \zeta_{16}^5, \zeta_{16}^{13} )$ .
\end{enumerate}
In any of these cases, $\spec(\rho_1(\ft)) \cap \spec(\rho_2(\ft)) = \{\zeta_8\}$ or $\{\zeta_8^5\}$. It follows from Theorem \ref{t:solution1} that $D =4$ and the Frobenius-Perron dimensions of the simple objects of $\CC$ are 
$$
 2, \sqrt{2}, 1, 1, 2, 2\,.
$$
In particular, $\CC$ is weakly integral, which contradicts Proposition \ref{p:weakly_int} (ii). Thus, $\rho_2$ is not of level 16 either.  As a consequence, $\rho_1$ cannot be of level 8.

(iv) {\bf $\rho_1$ cannot be of level $16$}: Assume contrary. Then  $\rho_1 \cong \rd{3}{16}{1,0},\, \rd{3}{16}{3,9},\, \rd{3}{16}{5,6},\, \rd{3}{16}{7,3}$, which are projectively inequivalent and have a minimal $\ft$-spectrum. Moreover,
$$
\rd{3}{16}{1,0} (\fs) = A
\quad \text{and} \quad\eta(\ft)=\diag(\zeta_8, \zeta_{16}, \zeta_{16}^9 )
$$
which is odd. Since all the 3-dimensional level 16 irreducible representations are projectively equivalent to a Galois conjugates of $\rd{3}{16}{1,0}$, its suffices consider the case $\rho_1\cong \rd{3}{16}{1,0}$.

By the $\ft$-spectrum criteria, $\rho_2$ cannot be projectively equivalent to any irreducible representation of level 4 or 5.  By (i), (ii) and (iii),  $\rho_2$ cannot be projectively equivalent to any irreducible representation of level 3, 7, 8. Therefore, $\rho_2$ can only be projectively equivalent to an irreducible representation of level 16. The $\ft$-spectrum criteria forces $\rho_2$ to be an irreducible  representation of level 16. Since $\rho_2$ is odd, by Proposition \ref{p:multiple}, $\rho_2 \cong \rd{3}{16}{1,6}$ or $\rd{3}{16}{5,6}$, which are respectively isomorphic to the following irreducible representations:
\begin{enumerate}
    \item[(1)] $ \fs \mapsto -A$,\quad  $\ft \mapsto \diag(\zeta_8^5, \zeta_{16}, \zeta_{16}^{9} )$;
    \item[(2)] $ \fs \mapsto -A$,\quad  $\ft \mapsto \diag(\zeta_8, \zeta_{16}^5, \zeta_{16}^{13} )$.
\end{enumerate}

For Case (1), $\spec(\rho_1(\ft)) \cap \spec(\rho_2(\ft)) = \{\zeta_{16}, \zeta_{16}^9\}$ but 
$$
\rho_1(\fs)_{ii}  + \rho_2(\fs)_{ii} = A_{ii}-A_{ii} = 0
$$
for $i = 2,3$. Therefore,  $\rho\cong  \rd{3}{16}{1,0} \oplus \rd{3}{16}{1,6}$ is impossible by Theorem \ref{t:solution1}.

For Case (2),  $\spec(\rho_1(\ft)) \cap \spec(\rho_2(\ft)) = \{\zeta_8\}$ and $\rho_1(\fs)_{11}+\rho_2(\fs)_{11} = 0$. It follows from Theorem \ref{t:solution1} that $\rho\cong  \rd{3}{16}{1,0} \oplus \rd{3}{16}{5,6}$ is also not possible. 
\end{proof}

\begin{lem} \label{l:4.13}
If $\rho_1$ is of level 5,  then $\rho_2$ cannot be projectively equivalent to any level 5 irreducible representation.
\end{lem}
\begin{proof}
Suppose $\rho_2$ is projectively equivalent to some level 5 irreducible representation. Then, by the $\ft$-spectrum criteria, $\rho_2$ is a level 5 irreducible representation. Since there are only two inequivalent level 5 irreducible representation, $\rho_1 \not \cong \rho_2$ by Proposition \ref{p:multiple}.  Then  $\spec(\rho_1(\ft))\cap \spec(\rho_2(\ft)) = \{1\}$. It follows from Appendix \ref{repPP} that
$$
\rho_1(\fs)_{11}+\rho_2(\fs)_{11} = 0\,.
$$
By Theorem \ref{t:solution1}(i), $\rho_1 \oplus \rho_2$ is not equivalent to any MD representation.  
 Therefore, $\rho_2$ cannot be projectively equivalent to any level 5 irreducible representation. 
\end{proof}

It follows from Lemmas \ref{l:4.12} and \ref{l:4.13} that the MD representation $\rho$ of $\CC$ of type (3,3) must have the irreducible decomposition $\rho_1 \oplus \rho_2$ where $\rho_1$ and $\rho_2$ are 3-dimensional and of levels 5 and 4.

\subsubsection{Solving modular data of type (3,3) level (5,4)} There are only two inequivalent level 5 irreducible representations $\rd{3}{5}{1}$ and $\rd{3}{5}{3}$. Note that $\s \circ \rd{3}{5}{1}= \rd{3}{5}{3}$ where $\s \in \GQ$ such that $\s(\zeta_5) = \zeta_5^3$.
One may assume $\rho_1 \cong \rd{3}{5}{1}$ which is even, and has a minimal $\ft$-spectrum. 
 
 By the $\ft$-spectrum criteria and the parity constraint, $\rho_2\cong \rd{3}{4}{1,0}$, and so   $$\spec(\rho_1(\ft)) \cap \spec(\rho_2(\ft)) = \{1\}.$$
 By Theorem \ref{t:solution1},  $D = 2/\frac{1}{\sqrt{5}} =2\sqrt{5}$ or $\dim(\CC) = 20$. Moreover, if $\irr(\CC)$ is reordered so that $\rho(\ft) = \diag(1,1, \zeta_5, \zeta_5^4, \ii, \-\ii )$,
then
$$
\rho(\fs)=  \frac{1}{2\sqrt{5}}\left[
\arraycolsep=2pt\def\arraystretch{1.3} {\small
\begin{array}{cccccc}
 1 & \kappa  & -2 \g_1 & -2 \g_2 & -\sqrt{5} \g_3 \kappa  & -\sqrt{5} \g_4 \kappa  \\
 \kappa  & 1 & -2 \g_1 \kappa  & -2 \g_2 \kappa  & \sqrt{5} \g_3 & \sqrt{5} \g_4 \\
 -2 \g_1 & -2 \g_1 \kappa  & -1-\sqrt{5} & (-1+\sqrt{5})\g_1\g_2 & 0 & 0 \\
 -2 \g_2 & -2 \g_2 \kappa  & (-1+\sqrt{5})\g_1\g_2 & -1-\sqrt{5} & 0 & 0 \\
 -\sqrt{5} \g_3 \kappa  & \sqrt{5} \g_3 & 0 & 0 & -\sqrt{5} & \sqrt{5} \g_3 \g_4 \\
 -\sqrt{5} \g_4 \kappa  & \sqrt{5} \g_4 & 0 & 0 & \sqrt{5} \g_3 \g_4 & -\sqrt{5} \\
\end{array}}
\right]
$$
for some  $\kappa, \g_i \in \{\pm 1\}$. One can conclude from $S$ that $\CC$ is pseudounitary, and so we can assume $\1 = \iota = e_1$. This implies $\kappa=1$, $\g_i=-1$ for $i=1,\dots, 4$. 
Thus, the modular data of $\CC$ is given by
$$
S=\left[
\arraycolsep=4pt\def\arraystretch{1.3} {\small
\begin{array}{cccccc}
 1 & 1  & 2 & 2  & \sqrt{5}   & \sqrt{5}   \\
 1  & 1 & 2 & 2  & -\sqrt{5}  & -\sqrt{5} \\
 2  & 2  & -1-\sqrt{5} & -1+\sqrt{5} & 0 & 0 \\
 2  & 2  & -1+\sqrt{5} & -1-\sqrt{5} & 0 & 0 \\
 \sqrt{5}  & -\sqrt{5}  & 0 & 0 & -\sqrt{5} & \sqrt{5} \\
 \sqrt{5}  & -\sqrt{5}  & 0 & 0 & \sqrt{5} & -\sqrt{5} \\
\end{array}}
\right] \quad \text{and}\quad  T= \diag(1,1, \zeta_5, \zeta_5^4, \ii, -\ii )\,.
$$
However, if  $\1 \ne \iota$, then one may assume $e_1=\1$ and $e_2=\iota$. Then the resulting modular data is $(PSP, T)$ where $P$ is the permutation matrix of the transposition $(1,2)$. In this sense, the two modular data corresponding to different spherical structures are the \emph{same}.

For $\rho_1 =\rd{3}{5}{3}$,  the corresponding modular data is $(\s(S), \s(T))$, where $\s \in \GQ$ such that $\s(\zeta_5)=\zeta_5^3$ and $\s(\ii)=\ii$. Precisely, 
$$
\s(S)=\left[
\arraycolsep=4pt\def\arraystretch{1.3} {\small
\begin{array}{cccccc}
 1 & 1  & 2 & 2  & -\sqrt{5}   & -\sqrt{5}   \\
 1  & 1 & 2 & 2  & \sqrt{5}  & \sqrt{5} \\
 2  & 2  & -1+\sqrt{5} & -1-\sqrt{5} & 0 & 0 \\
 2  & 2  & -1-\sqrt{5} & -1+\sqrt{5} & 0 & 0 \\
 -\sqrt{5}  & \sqrt{5}  & 0 & 0 & \sqrt{5} & -\sqrt{5} \\
 -\sqrt{5}  & \sqrt{5}  & 0 & 0 & -\sqrt{5} & \sqrt{5} \\
\end{array}}
\right] \quad \text{and}\quad  \s(T)= \diag(1,1, \zeta_5^3, \zeta_5^2, \ii, -\ii )\,.
$$
In this case, the $e_2 = \iota$. One can use the other spherical structure of $\CC$ so that $\1 =\iota = e_1$. The resulting modular data is $(P\s(S)P, \s(T))$, which is the \emph{same} as the modular data $(\s(S), \s(T))$, and is the  modular data of $SO(5)_2$.  This completes the proof of Theorem \ref{t:33}.

\subsection{Classification of Modular Data of type $(3 ,2,1)$}
We now classify modular tensor categories with $\SL$ representations decomposing as a direct sum of irreducible representations of dimension $3,2$ and $1$.
The main theorem of this section is:
\begin{thm}\label{t:321}
The modular data of any type $(3,2,1)$ modular tensor category is a Galois conjugate of a non-trivial braided zesting of $SO(5)_2$. 
\end{thm}
The zesting procedure is found in \cite{DGPRZ}.  An alternative approach is to consider the classification of metaplectic modular tensor categories in \cite{ACRZ}: this shows that the categories above can be obtained by gauging the particle-hole symmetry on a pointed modular tensor category of the form $\CC(\Z_5,q)$.  In \cite{GRR} it is shown that of the 4 modular tensor categories obtained in this way, 2 are $SO(5)_2$ and its (unitary) Galois conjugate and the other two are the non-trivial zesting of $SO(5)_2$ and its (unitary) Galois conjugate.

Let $\rho=\chi_1\oplus(\rho_2\otimes\chi_2)\oplus(\rho_3\otimes\chi_3)$ be the irreducible decomposition a modular representation with $\rho_i$ irreducible of dimension $i$ of prime power level and $\chi_i$ a character.  This description is possible by the Chinese Remainder Theorem and the fact that $2$ and $3$ are prime. As before, we may assume $\chi_3 = 1$ and require $\rho_3$ has a minimal $\ft$-spectrum.

We consider cases in turn, describing the level triples for $(\rho_3,\rho_2,\chi_1)$.  The $\ft$-spectrum criteria immediately implies that the level of $\rho_3$ cannot be $7$.  Similarly the level of $\rho_3$ cannot be $16$: looking at the eigenvalues of the level $16$ irreducible $3$-dimensional representation we see that  $\chi_1(\ft)\not\in \spec(\rho_3(\ft))$, and hence  $\spec((\rho_2\o \chi_2)(\ft)) \cap \spec(\rho_3(\ft)) \ne \emptyset$. This implies $\rho_2 \o \chi_2$ has level 8 but then   $\chi_1(\ft) \not\in \spec(\rho_3(\ft) \oplus  (\rho_2 \o \chi_2)(\ft))$, which contradicts the $\ft$-spectrum criteria.

Suppose the level of $\rho_3$ is $8$. Then $\rho_3 \cong \rd{2}{8}{1,0}$ or $\rd{2}{8}{3,0}$, and hence $\rho_3$ is odd. Note that   $\spec(\rd{2}{8}{1,0})=\{1,\zeta_8,-\zeta_8\}$ and $\spec(\rd{2}{8}{3,0})=\{1,\zeta_8^3,-\zeta_8^3\}$. The level of $\rho_2$ cannot be 5, by inspection of the corresponding eigenvalues.  If the level of $\rho_2$ is $2$ then the $\ft$-spectrum criteria implies that $(\chi_2)^2 = (\chi_1)^2=1$.  But now $\rho(\fs^2)$ has trace $0$, contradicting Proposition \ref{p:MD}.  Thus the level of $\rho_2$ is either $8$ or $3$.  Applying the  $\ft$-spectrum criteria yields the following possible levels in this case: $(8,8,1),(8,3,3)$ or $(8,3,1)$. In particular, if the level of $\rho_2$ is $8$ we cannot have levels $(8,8,2)$ or $(8,8,4)$ as the $\ft$-spectrum criteria fails in these cases.  In all three cases we see that $\rho_2 \o \chi_2$ must be odd for otherwise $\Tr(\rho(\fs^2))=0$. Hence, the corresponding category would be non-self-dual. 

Now suppose that the level of $\rho_3$ is $5$. Then $\rho_3$ is even.  The $\ft$-spectrum criteria implies the level of $\rho_2$ cannot be $8$.  Inspecting the remaining possibilities we find the following possible level triples: $(5,5,1),(5,3,1),(5,3,3),(5,2,1)$ or $(5,2,2)$.  The parities imply that the corresponding category would be non-self-dual in the first three cases and self-dual for the last two.

Next if the level of $\rho_3$ is $4$, then $\rho_3 \cong \rd{3}{4}{1,3}$ which  is odd, and has the minimal $\ft$-spectrum $\{1, -1, \ii\}$ according Appendix \ref{repPP}. The $\ft$-spectrum criteria show that the level of $\rho_2$ cannot be $8$ or $5$.  If $\rho_2$ had level $2$ then the order of $\rho(\ft)$ would be $4$, yielding a pointed {integral} category (by Theorem \ref{t:2346}) with $T$-matrix of order $4$, which contradicts Proposition \ref{p:weakly_int}.  Thus $\rho_2$ has level $3$ and we find $(4,3,1),(4,3,2),(4,3,3)$ and $(4,3,4)$ as possible level triples.

Finally, if the level of $\rho_3$ were $3$ then the  $\ft$-spectrum criteria implies that the order of $\rho(\ft)$ is a divisor of $6$ and hence pointed {integral} by Theorem \ref{t:2346}.  This contradicts Proposition \ref{p:weakly_int}.

Below we provide the details of the cases of levels $(4,3,2),(5,2,2)$ and $(5,2,1)$ explicitly.  The remaining cases
$\{(8,8,1),(8,3,3),(8,3,1),(5,5,1),(5,3,1),(5,3,3),(4,3,1),(4,3,3),(4,3,4)\}$ can be similarly addressed (and indeed are easier).  We can eliminate all of these cases computationally as well, see Section \ref{repL}.

\subsubsection{Case $(4,3,2)$}
Suppose that the levels of $\rho_3,\rho_2$ and $\chi_1$ are $4,3$ and $2$, respectively.  Without loss of generality we may assume that $\rho_3 \cong \rd{3}{4}{1,3}$ and  $\rho_2\cong \rd{3}{2}{1,0}$, which are odd, and respectively have the minimal $\ft$-spectrums $\{1,-1, \ii\}$ and $\{1,\zeta_3\}$ according to Appendix \ref{repPP}.    Let us determine what $\chi_2$ can be. Note that $\chi_1$ is even.  Now if $\chi_2$ is odd, then $\Tr(\rho(\fs^2))=0$, which is impossible.  The $\ft$-spectrum criteria implies that $\chi_2(\ft)\neq \pm \zeta_3$, and a relabeling eliminates $\zeta_3^{-1}$.  If $\chi_2(\ft) \in \{-1, -\ol\zeta_3\}$, then $\rho$ is projectively equivalent to the complex conjugate the $(4,3,1)$ case.  So we will assume that $\chi_2(\ft)=1, \zeta_3$ or $\rho_2 \o \chi_2 \cong \rd{3}{2}{1,0}$ or $\rd{3}{2}{1,8}$. In either case, $\rho$ is defined over $\BQ_{24}$. If $\s \in \GQ$ such that $\s(\zeta_3)  = \zeta_3^2$ and $\s(\zeta_8)=\zeta_8$, then we have
$$
\rho_3 \oplus \rd{3}{2}{1,8} \oplus \chi_1 \cong \s\circ(\rho_3 \oplus \rd{3}{2}{1,0} \oplus \chi_1)\,.
$$
It suffices to consider $\rho \cong \tilde\rho:=\rho_3 \oplus \rd{3}{2}{1,0} \oplus \chi_1$. By Appendix \ref{repPP}, we have
$$
\tilde\rho(\fs)= \frac {\ii}{2}\left[ 
\arraycolsep=4pt\def\arraystretch{1.3} {\small
\begin {array}{ccc} -1 &1 & \sqrt{2} \\ 
1 & -1 & \sqrt{2}\\
\sqrt{2} & \sqrt{2} & 0
\end {array}} \right] \oplus 
\frac {\ii}{\sqrt {3}} \left[ 
\arraycolsep=4pt\def\arraystretch{1.3} {\small
\begin{array}{cc} -1 & \sqrt{2} \\ 
\sqrt{2} & 1
\end {array} }
\right] 
\oplus [-1]\quad \text{and} \quad \tilde\rho(\ft)=\diag(1,-1,\ii,1,\zeta_3,-1)
$$
Reordering $\irr(\CC)$ so that $\rho(\ft)=\diag(1,1, -1,-1, \ii, \zeta_3)$. By Theorem \ref{t:solution1}, the unit object $\1$ must be an eigenvector $\rho(\ft)$ with eigenvalue $1$ and so $D = 2/(\frac{1}{2}+\frac{1}{\sqrt{3}})=8\sqrt{3}-12$ or $\dim(\CC)=48 (7 - 4 \sqrt{3})$. Moreover, $T= \diag(1,1, -1,-1, \ii, \zeta_3)$ and 
$\rho(\fs)=\frac{-\ii (3 + 2 \sqrt{3})}{12} S$, where 
$$
S=
\left[
\arraycolsep=4pt\def\arraystretch{1.3} {\small
\begin{array}{cccccc}
 1 & -\frac{\left(2 \sqrt{3}-3\right) \kappa }{2 \sqrt{3}+3} & -\frac{3 \sqrt{2} a  }{2 \sqrt{3}+3} & -\frac{3 \sqrt{2} b  }{2 \sqrt{3}+3} & -\frac{6 \g_1}{2 \sqrt{3}+3} & \frac{4 \sqrt{3} \g_2 \kappa }{2 \sqrt{3}+3} \\
 -\frac{\left(2 \sqrt{3}-3\right) \kappa }{2 \sqrt{3}+3} & 1 & -\frac{3 \sqrt{2} a   \kappa }{2 \sqrt{3}+3} & -\frac{3 \sqrt{2} b   \kappa }{2 \sqrt{3}+3} & -\frac{6 \g_1 \kappa }{2 \sqrt{3}+3} & -\frac{4 \sqrt{3} \g_2}{2 \sqrt{3}+3} \\
 -\frac{3 \sqrt{2} a  }{2 \sqrt{3}+3} & -\frac{3 \sqrt{2} a   \kappa }{2 \sqrt{3}+3} & -\frac{6 \left(-1+(1+2 i) b^2\right)}{2 \sqrt{3}+3} & \frac{(6+12 i) a b}{2 \sqrt{3}+3} & -\frac{6 \sqrt{2} a   \g_1}{2 \sqrt{3}+3} & 0 \\
 -\frac{3 \sqrt{2} b  }{2 \sqrt{3}+3} & -\frac{3 \sqrt{2} b   \kappa }{2 \sqrt{3}+3} & \frac{(6+12 i) a b}{2 \sqrt{3}+3} & \frac{12 i \left(-1+\left(1-\frac{i}{2}\right) b^2\right)}{2 \sqrt{3}+3} & -\frac{6 \sqrt{2} b   \g_1}{2 \sqrt{3}+3} & 0 \\
 -\frac{6 \g_1}{2 \sqrt{3}+3} & -\frac{6 \g_1 \kappa }{2 \sqrt{3}+3} & -\frac{6 \sqrt{2} a   \g_1}{2 \sqrt{3}+3} & -\frac{6 \sqrt{2} b   \g_1}{2 \sqrt{3}+3} & 0 & 0 \\
 \frac{4 \sqrt{3} \g_2 \kappa }{2 \sqrt{3}+3} & -\frac{4 \sqrt{3} \g_2}{2 \sqrt{3}+3} & 0 & 0 & 0 & -\frac{4 \sqrt{3}}{2 \sqrt{3}+3} \\
\end{array}
}
\right]
$$
for some $\kappa, \g_i \in \{\pm 1\}$ and $a, b \in \BR$ such that $a^2+b^2=1$. Since $1, \iota \in \{e_1, e_2\}$, $\kappa=-1$, and $a, b \ne 0$. Since $\dim(\iota) = \frac{\left(2 \sqrt{3}-3\right)}{2 \sqrt{3}+3} = 7-4\sqrt{3} < 1$, $\FPdim(\CC)=48 (7 + 4 \sqrt{3})$ and $\iota \ne \1$. We may simply assume $e_1=\1$ and $e_2 = \iota$. Then $\g_1 =1$, $\g_2=-1$ and $a  , b >0$.  Since $\Tr(\rho(\fs^2))=-4$, there is exactly one dual pair of simple objects, and they can only be $e_3, e_4$. Therefore, $a=b = \frac{\pm 1}{\sqrt{2}}$ and so $a   = b   = \frac{1}{\sqrt{2}}$.
Thus,  
$$
S= \left[
\arraycolsep=4pt\def\arraystretch{1} {\small
\begin{array}{cccccc}
 1 & 1-2 d & -d & -d & -2 d & 2-2 d \\
 1-2 d & 1 & d & d & 2 d & 2-2 d \\
 -d & d & (1-2 \ii) d & (1+2 \ii) d & -2 d & 0 \\
 -d & d & (1+2 \ii) d & (1-2 \ii) d & -2 d & 0 \\
 -2 d & 2 d & -2 d & -2 d & 0 & 0 \\
 2-2 d & 2-2 d & 0 & 0 & 0 & 2 d-2 \\
\end{array}
}
\right]
$$
where $d = 2 \sqrt{3}-3$. Remarkably, the Verlinde formula yields a consistent set of fusion rules.  For example the object with twist $\ii$ has the fusion matrix: 
$$
N_5=\left[
\arraycolsep=4pt\def\arraystretch{1} {\small
\begin{array}{cccccc}
 0 & 0 & 0 & 0 & 1 & 0 \\
 0 & 4 & 2 & 2 & 3 & 4 \\
 0 & 2 & 2 & 0 & 1 & 2 \\
 0 & 2 & 0 & 2 & 1 & 2 \\
 1 & 3 & 1 & 1 & 4 & 4 \\
 0 & 4 & 2 & 2 & 4 & 4 \\
\end{array}}
\right]\,.
$$
However, the second FS-indicator for this object is $\nu_2(e_5)=\frac{1}{\dim(\CC)}\sum_{j,k}N_{j,k}^3d_jd_k\left(\frac{\theta_j}{\theta_k}\right)^2=2$, a contradiction.

\subsubsection{Case $(5,2,1)$}
Consider the case of levels $(5,2,1)$. Then $\rho \cong \rd{3}{5}{1} \oplus\rd{2}{2}{1,0} \oplus \chi^0$ or $\rd{3}{5}{3} \oplus\rd{2}{2}{1,0} \oplus \chi^0$ according to Appendix \ref{repPP}.  Since the latter is a Galois conjugate of the former one, it suffices to solve the first case. Let $\tilde\rho = \rd{2}{2}{1,0} \oplus \chi^0 \oplus  \rd{3}{5}{3}$ in which  $\tilde\rho(\ft)=\diag\left(1,-1, 1,1, \zeta_5,\zeta_5^4\right)$.
By permuting the first two basis elements, we may assume that $t=\rho(\ft)=\diag(-1,1,1,1,\zeta_5,\zeta_5^4)$.  Conjugating by a block diagonal matrix of the form $(r_1)\oplus F \oplus (r_2)\oplus(r_3)$ where $F=\left[ \begin {array}{ccc} f_{1,1} & f_{1,2}& f_{1,3}\\f_{2,1} & f_{2,2}& f_{2,3}\\f_{3,1} & f_{3,2}& f_{3,3}
\end{array}\right]$ is real orthogonal matrix (cf. Prop. \ref{t:ortho_eqv}) and 
$r_i=\pm 1$. One may assume $r_1=1$, and we find that $\pm S/D=s=\rho(\mathfrak{s})$ has the form:

$$\left[ \begin {array}{cccccc} {\frac{1}{2}}&{-\frac {f_{{1,1}}\sqrt {3
}}{2}}&{-\frac {f_{{2,1}}\sqrt {3}}{2}}&{-\frac {f_{{3,1}}\sqrt {3}}{2}}
&0&0\\ 
\noalign{\medskip}{-\frac {f_{{1,1}}\sqrt {3}}{2}}&
&&&{\frac {f_{{1,3}}\sqrt 
{10}r_2}{5}}&{\frac {f_{{1,3}}\sqrt {10}r_3}{5}}
\\ 
\noalign{\medskip}{-\frac {f_{{2,1}}\sqrt {3}}{2}}&&A&&{\frac {f_{{2,3}}\sqrt {10}
r_2}{5}}&{\frac {f_{{2,3}}\sqrt {10}r_3}{5}}
\\ 
\noalign{\medskip}{-\frac {f_{{3,1}}\sqrt {3}}{2}}&&&&{\frac {f_{{3,3}}\sqrt {10}
r_2}{5}}&{\frac {f_{{3,3}}\sqrt {10}r_3}{5}}
\\ \noalign{\medskip}0&{\frac {f_{{1,3}}\sqrt {10}r_2}{5}}&{
\frac {f_{{2,3}}\sqrt {10}r_2}{5}}&{\frac {f_{{3,3}}\sqrt {10}r_2}{5}}&-{\frac {  \sqrt {5}+5  }{10}}&{
\frac {r_2\,\sqrt {5} \left( \sqrt {5}-1 \right) r_3}{10}}
\\ \noalign{\medskip}0&{\frac {f_{{1,3}}\sqrt {10}r_3}{5}}&{
\frac {f_{{2,3}}\sqrt {10}r_3}{5}}&{\frac {f_{{3,3}}\sqrt {10} r_3}{5}}&{\frac {r_2\,\sqrt {5} \left( \sqrt {5}-1 \right)r_3}{10}}&-{\frac {\sqrt {5} \left( \sqrt {5}+1 \right) }{10}}
\end {array} \right]$$ 
where $$A=\small
\left[ \begin {array}{ccc} -{\frac {{f_{{1,1}}}^{2}}{2}}+{f_{{1,2}}}^
{2}+{\frac {{f_{{1,3}}}^{2}\sqrt {5}}{5}}&*&*\\ \noalign{\medskip}-{\frac {f_{{1,1}}f_{{2,1}
}}{2}}+f_{{1,2}}f_{{2,2}}+{\frac {f_{{1,3}}\sqrt {5}f_{{2,3}}}{5}}&-{
\frac {{f_{{2,1}}}^{2}}{2}}+{f_{{2,2}}}^{2}+{\frac {{f_{{2,3}}}^{2}
\sqrt {5}}{5}}&*\\ \noalign{\medskip}f_{{1,2}}f_{{3,2}}+{\frac {f_{{1,3}}\sqrt {5}f_
{{3,3}}}{5}}-{\frac {f
_{{1,1}}f_{{3,1}}}{2}}&f_{{2,2}}f_{{3,2}}+{
\frac {f_{{2,3}}\sqrt {5}f_{{3,3}}}{5}}-{\frac {f_{{2,1}}f_{{3,1}}}{2}}&
{f_{{3,2}}}^{2}+{\frac {{f_{{3,3}}}^{2}\sqrt {5}}{5}}-{\frac {{f_{{3,1}}}^{2}}{2}}\end {array}
 \right] .
$$

First we observe that the FP-dimensions and categorical dimensions (which may coincide) must appear as multiples of one of the columns $2,3$ or $4$.  Moreover, since our category is non-integral by Proposition \ref{p:weakly_int}, the Galois orbit of the dimension column has size $2$. The FP-dimension column of $s$ must have all the same sign, which implies that $r_2=r_3$.   

Let $\s \in \Gal(\BQ_5/\BQ)$ be the automorphism defined by $\zeta_5\rightarrow \zeta_5^3$.  
By Galois symmetry we have: $\hat\sigma(1)=1$, $\hat\sigma(5)=6$. {Therefore,} $\hat\sigma$ has order $2$. Reordering the rows of $F$ if necessary (which permutes the corresponding rows/columns of $s$) we may assume that $\hat\sigma(2)=2$ and $\hat\sigma(3)=4$, so that the FP-dimensions and categorical dimensions correspond to either columns $3$ or $4$ (or one of each).
 
 We will make frequent use of the fact that $\sigma(s_{ij})=\epsilon_\sigma(i)s_{\hs(i),j}=\epsilon_\sigma(j)s_{i,\hs(j)}$ where $\epsilon_\sigma(i)$ is a sign.

Now $1/2=\sigma(s_{1,1})=\epsilon_\sigma(1)/2$ so that $\epsilon_\sigma(1)=1$.  By a similar computation $\sigma(s_{1,2})=\epsilon_\sigma(1)s_{1,2}=\epsilon_\sigma(2)s_{1,2}$, so that  $\epsilon_\sigma(2)=1$.
From $\sigma(s_{5,5})=\frac{\sqrt{5}-5}{10}=\epsilon_\sigma(5)s_{5,6}$ we find that $\epsilon_\sigma(5)=-1$.
Now we compute two ways: $\sigma(s_{2,5})=\epsilon_\sigma(2)s_{2,5}=s_{2,5}=\epsilon_\sigma(5)s_{2,6}={-s_{2,6}=-s_{2,5}}$, which implies $s_{2,5}=0$ so that $f_{1,3}=0$.
Now $\sigma(s_{3,5})=\epsilon_\sigma(3)s_{4,5}=\epsilon_\sigma(5)s_{3,6}=-s_{3,6}$ implies $f_{3,3}=\pm f_{2,3}$ so that $(f_{2,3})^2=\frac{1}{2}$. 
Applying a similar calculation we see that $\sigma(s_{1,3})=s_{1,3}=\epsilon_\sigma(3)s_{1,4}$ implies $f_{2,1}=\pm f_{3,1}$.  Setting $z=f_{1,1}$ and $y=f_{2,1}$  orthogonality yields the following:
$$F=\left[ \begin {array}{ccc} z & \delta_1\sqrt{2}y& 0\\y & \frac{-\delta_1z}{\sqrt{2}}& \frac{-\delta_2\delta_3}{\sqrt{2}}\\\delta_2y & \frac{-\delta_1\delta_2 z}{\sqrt{2}}& \frac{\delta_3}{\sqrt{2}}
\end{array}\right].$$  One important consequence is that there are only 2 rows of $\rho(\fs)$ that have strictly non-zero entries: the 3rd and the 4th.  

Next we find that $\sigma(s_{i,1})=s_{i,1}$ since $\epsilon_\sigma(1)=1$ and $\hs(1)=1$.  Thus $f_{i,1}\sqrt{3}\in\BQ$. Note that $z^2+2y^2=1$ where $z,y\in\frac{1}{\sqrt{3}}\BQ$, and one of $s_{2,1}/s_{3,1}=\pm s_{2,1}/s_{4,1}$ is of the form $S_{X,Y}/d_X$, i.e., an eigenvalue of a fusion matrix.  In particular $z/y=\gamma$ is a (rational) algebraic integer, i.e., $\gamma\in\Z$.  From this we find that $\gamma^2+2=1/y^2\in\BZ$ so that $0< y^2\leq 1/3$, and so $1/3\leq z^2\leq 1$.

Let us compute the values of the submatrix $A$ above.  We have:
$$A=\left[ \begin {array}{ccc} -z^2/2+{2} y^2&*&*\\
-3yz/2&\frac{1}{2}(z^2-y^2+\frac{1}{\sqrt{5}})&*\\
-\delta_23yz/2&\frac{\delta_2}{2}(z^2-y^2-\frac{1}{\sqrt{5}})&\frac{1}{2}(z^2-y^2+\frac{1}{\sqrt{5}})\end{array}\right].$$

%Now one of $s_{4,3}/s_{4,4}=s_{3,4}/s_{3,3}$ is the (FP-)dimension of a simple object, up to a sign.  Setting $\alpha=z^2-y^2$ we find that $\left|\frac{\alpha-\frac{1}{\sqrt{5}}}{\alpha+\frac{1}{\sqrt{5}}}\right|\geq 1$ as it is an $FP$-dimension.  But $0\leq \alpha\leq 1$ which implies that $\alpha=0$, i.e., $z=\pm y$. 

Since the unit object can only correspond to either row 3 or 4 and $s_{32} = \pm s_{42}$, $s_{22}/s_{32}$ is an algebraic integer in $\BQ(\sqrt{5})$. Note that 
 $$
 \frac{s_{22}}{s_{32}} = \frac{\gamma}{3}- \frac{4}{3 \gamma} =  \frac{\gamma^2-4}{3 \gamma} \in \BQ\,.
 $$
 Therefore, $\frac{\gamma^2-4}{3 \gamma} \in \BZ$ and so $\gamma \mid 4$. Thus, $\gamma^2 = 1, 4$ or $16$. However, if $\gamma^2 =4$ or $16$, $y= \frac{\pm 1}{\sqrt{\gamma^2+2}} \not\in\frac{1}{\sqrt{3}}\BQ$.  Thus, $\gamma^2=  1$ or $z= \pm y$.

This implies that $y=\pm \frac{1}{\sqrt{3}}$, from which we compute:
 $f_{2,2}=\pm\frac{1}{\sqrt{6}}$, $f_{3,2}=\pm\frac{1}{\sqrt{6}}$, $f_{1,1}=\pm\frac{1}{\sqrt{3}}$, and $f_{1,2}=\pm\frac{2}{\sqrt{6}}$.

\iffalse

As the FP-dimension corresponds to either row $3$ or $4$, we see that $\frac{\pm 1}{D}$ is the entry $s_{3,3}$ or $s_{4,4}$, both of which are equal to $\frac{1}{2\sqrt{5}}$ by the above.  \fi

Now we may assume $F=\left[ \begin {array}{ccc}  -1/\sqrt{3} & 2x_1/\sqrt{6} &0\\
x_2/\sqrt{3} &x_3/\sqrt{6} &x_4/\sqrt{2}\\
x_5/\sqrt{3} & x_6/\sqrt{6} & x_7/\sqrt{2}
\end {array}
 \right]$ where the $x_i=\pm 1$ after an overall rescaling by $\pm 1$.  Orthogonality of $F$ implies several additional conditions on the $x_i$, so that all are determined by the values of $x_2,x_4,x_5$ and $x_7$.

 Substituting into $s$ above, rescaling by $\pm D$ and permuting the rows/columns so that the two non-zero rows appear first, we have:
 
 $$S=\left[ \begin {array}{cccccc} 1&x_{{4}}x_{{7}}&\sqrt {5}x_{{5}}&
\sqrt {5}x_{{5}}&2\,x_{{7}}r_3&2\,x_{{7}}r_3
\\ \noalign{\medskip}x_{{4}}x_{{7}}&1&\sqrt {5}x_{{2}}&\sqrt {5}x_{{2}
}&2\,x_{{4}}r_3&2\,x_{{4}}r_3\\ \noalign{\medskip}\sqrt {5}x
_{{5}}&\sqrt {5}x_{{2}}&\sqrt {5}&-\sqrt {5}&0&0\\ \noalign{\medskip}
\sqrt {5}x_{{5}}&\sqrt {5}x_{{2}}&-\sqrt {5}&\sqrt {5}&0&0
\\ \noalign{\medskip}2\,x_{{7}}r_3&2\,x_{{4}}r_3&0&0&-\sqrt
{5}-1&\sqrt {5}-1\\ \noalign{\medskip}2\,x_{{7}}r_3&2\,x_{{4}}
r_3&0&0&\sqrt {5}-1&-\sqrt {5}-1\end {array} \right].$$
 
 Thus we see that the dimensions and FP-dimensions must be, up to sign choices, among ${1,2,\sqrt{5}}$.  In particular, any such category must be weakly integral, and there is an invertible object of order 2. Therefore, are two spherical structures on $\CC$ which make $\1=\iota$ or $\1 \ne \iota$. We may assume $\1$ corresponds to the first row. For the first case, we find $x_4x_7=x_5=x_7r_3=1$ and $x_2=-1$.

 Thus we obtain the following $S$-matrix:
$$S= \left[\begin{array}{cccccc} 1 & 1& \sqrt{5} & \sqrt{5} & 2 & 2\\ 1 & 1& -\sqrt{5} & -\sqrt{5} & 2 & 2\\ \sqrt{5} & \sqrt{5} &\sqrt{5} & -\sqrt{5} & 0&0\\ \sqrt{5} & \sqrt{5} &-\sqrt{5} & \sqrt{5} & 0&0\\ 2 & 2& 0 & 0& -\sqrt{5} -1 &\sqrt{5} -1\\ 2 & 2& 0 & 0& \sqrt{5} -1 &-\sqrt{5} -1\\\end {array}
 \right] $$
 
  For the second case, one can obtain the same $S$-matrix except the first two rows/columns are interchanged, but the $T$-matrix is unchanged. Therefore, we have only one modular data for either case.
 
 Applying $\s$ to $(S,T)$, we obtain the modular data for  $\rho \cong \rd{3}{5}{3} \oplus \rd{2}{2}{1,0}\oplus \chi^0$ with the $T$-matrix given  by $\s(T)=\diag(1,1,1,-1,\zeta^3_5,\zeta_5^{2})$. 
Both of these  modular data $(S,T)$ and $(\s(S), \s(T))$ are  modular data of non-trivial braided zesting of MTCs (see \cite{DGPRZ})  of type (3,3).  Notice that the MTCs of type (3,3) have $T$-matrix of order $20$.

\subsubsection{Case $(5,2,2)$} 
It suffice to consider the case with $\rho\cong \tilde\rho := \rd{3}{5}{1} \oplus \rd{2}{2}{1,0} \oplus \chi^6$. Then 
$$
\tilde\rho(\fs)= \frac {1}{\sqrt{5}}\left[ 
\arraycolsep=2pt\def\arraystretch{1.3} {\textstyle
\begin {array}{ccc} 1 & -\sqrt{2} &-\sqrt{2} \\ 
 -\sqrt{2} & -\varphi & \varphi^{-1}\\
-\sqrt{2} & \varphi^{-1}  &  -\varphi
\end {array}} \right] \oplus 
\frac {1}{2} \left[ 
\arraycolsep=2pt\def\arraystretch{1.3} {\small
\begin{array}{cc} -1 & -\sqrt{3} \\ 
-\sqrt{3} & 1
\end {array} }
\right] 
\oplus [-1]\quad \text{and} \quad \tilde\rho(\ft)= \diag(1,\zeta_5,\zeta_5^4, 1, -1,-1)\,.
$$
Permute $\irr(\CC)$ so that $\rho(\ft)=\diag(-1,-1,1,1,\zeta_5,\zeta_5^4)$. By Theorem \ref{t:solution1}, the objects $\1, \iota \in\{e_3,e_4\}$,  $D=2/(\frac{1}{2}-\frac{1}{\sqrt{5}}) = 20+8\sqrt{5}$, and
$$
s:=\rho(\fs)= \left[
\arraycolsep=2pt\def\arraystretch{1.3} {\textstyle
\begin{array}{cccccc}
 \frac{3 b^2}{2}-1 & -\frac{1}{2} (3 a b) & \frac{1}{2} \sqrt{\frac{3}{2}} b & \frac{1}{2} \sqrt{\frac{3}{2}} b \kappa  & 0 & 0 \\
 -\frac{1}{2} (3 a b) & \frac{1}{2} \left(1-3 b^2\right) & -\frac{1}{2} \sqrt{\frac{3}{2}} a & -\frac{1}{2} \sqrt{\frac{3}{2}} a \kappa  & 0 & 0 \\
 \frac{1}{2} \sqrt{\frac{3}{2}} b & -\frac{1}{2} \sqrt{\frac{3}{2}} a & \frac{1}{20} \left(2 \sqrt{5}-5\right) & -\frac{1}{20} \left(2 \sqrt{5}+5\right) \kappa  & -\frac{\gamma _1 \kappa }{\sqrt{5}} & -\frac{\gamma _2 \kappa }{\sqrt{5}} \\
 \frac{1}{2} \sqrt{\frac{3}{2}} b \kappa  & -\frac{1}{2} \sqrt{\frac{3}{2}} a \kappa  & -\frac{1}{20} \left(2 \sqrt{5}+5\right) \kappa  & \frac{1}{20} \left(2 \sqrt{5}-5\right) & \frac{\gamma _1}{\sqrt{5}} & \frac{\gamma _2}{\sqrt{5}} \\
 0 & 0 & -\frac{\gamma _1 \kappa }{\sqrt{5}} & \frac{\gamma _1}{\sqrt{5}} & \frac{1}{10} \left(-\sqrt{5}-5\right) & \frac{2 \gamma _1 \gamma _2}{\sqrt{5}+5} \\
 0 & 0 & -\frac{\gamma _2 \kappa }{\sqrt{5}} & \frac{\gamma _2}{\sqrt{5}} & \frac{2 \gamma _1 \gamma _2}{\sqrt{5}+5} & \frac{1}{10} \left(-\sqrt{5}-5\right) \\
\end{array}
}
\right]
$$
for some $\kappa, \g_1, \g_2 \in \{\pm 1 \}$ and $a, b \in \BR$ such that $a^2+b^2=1$. Since $\frac{5+2\sqrt{5}}{5-2\sqrt{5}} > 1$, $\iota =\1$. We may simply assume $e_4=\1$. Then 
$\kappa=1$, $\g_1=\g_2 =-1$ and $a>0$ and $b<0$.

By Proposition \ref{p:weakly_int}, $\CC$ is not integral. Let $\s \in \Gal(\BQ_5/\BQ)$ be a generator. Then $\hs(3)=4$ and $\e_\s(3)=1$ since $\s(s_{3,5})=s_{4,5}$. Therefore, $\s$ fixes $s_{3,1}, s{3,2}$, and so $\sqrt{\frac{3}{2}}a,\sqrt{\frac{3}{2}}b \in \BQ$. Now,  $s_{2,1}= s_{2,1} \in \BQ $ since $ab \in \BQ$. By Theorem \ref{p:MD1}, 
$\frac{s_{2,1}}{s{4,1}}$ and $\frac{s_{1,2}}{s{4,2}}$ are in $\BZ[\zeta_5] \cap \BQ$, $\sqrt{6}a, \sqrt{6}b \in \BZ$ and $(\sqrt{6}a)^2+(\sqrt{6}b)^2 = 6$.  But the Diophantine equation $X^2+Y^2=6$ has no integral  solutions, so we conclude that $\tilde\rho$ has no realization.

\subsection{Classification of modular data of type (6)}

In this subsection, we discuss the possible rank-6 MDs of type $(6)$
({\it i.e.} MDs from dimension-6 irreducible $\SL$ symmetric representations).  This part of the classification relies upon computer computations.
\begin{thm} \label{t:60}
Let $\CC$ be a rank 6 modular tensor category of type $(6)$ with $\dim(\CC)=D^2 \not\in \mathbb{Z}$. Then the modular data of $\CC$ can be obtained, {up to a choice of (spherical) pivotal structure,} as a Galois conjugate of the modular data of the following modular tensor
categories:
\begin{enumerate}
    \item[(i)]  $PSU(2)_{11}$ (entry 10 in Appendix \ref{SsL}); 
    \item[(ii)]  $PSU(2)_3 \boxtimes PSU(2)_5$ (entry 20 in Appendix \ref{SsL}); 
    \item[(iii)]  $SU(2)_1\boxtimes PSU(2)_5$ (entry 24 in Appendix \ref{SsL});
    \item[(iv)] $PSU(2)_3\boxtimes SU(2)_2$ (entry 36 in Appendix \ref{SsL}).
    \item[(v)] $PSU(2)_3\boxtimes E(8)_2$ (entry 28 in Appendix \ref{SsL}).
    \item[(vi)] $PSO(5)_{3/2}$ (non-unitary, entry 9 in Appendix \ref{SsL});
\end{enumerate}
\end{thm}

It is worth noting that (i), (ii) and (vi) have a unique pivotal structure, up to equivalence (cf. \cite{BNRW}). The categories (i) and (ii) are \emph{transitive} \cite{NWZ}, and they are completely determined by their modular data.

Recall that a symmetric  $\SL$ representation $\rho$ is defined to be an unitary representation which has diagonal
$\rho(\ft)$ and symmetric $\rho(\fs)$.  Every finite-dimensional representation of $\qsl{n}$ is equivalent to a symmetric one. Two  symmetric $\SL$ representations are
equivalent if and only if they are related by a conjugation of a real
orthogonal matrix (see Theorem \ref{t:ortho_eqv}).  There are 70 inequivalent 6-dimensional symmetric irreducible $\SL$ representations of prime-power levels (cf. Appendix \ref{repPP}).  Up to tensoring one of the 12 1-dimensional representations, other
6-dimensional irreducible representations are tensor products of one of the 11 2-dimensional  and one of the 33 3-dimensional irreducible symmetric representations of {distinct} prime-power levels.

Since there are only a finite number of $\SL$ representations, up to equivalence, for any given dimension, we can
examine {representatives} of each of those symmetric representations by computer and reject those
representations that do not satisfy the following necessary conditions (for a
 symmetric $\SL$ representation equivalent to an MD representation):
\begin{enumerate}

\item If all the eigenvalues of
$\rho(\ft)$ are distinct (non-degenerate) then $\rho(\fs)$ has a row that contains no zero.  Note that when $\rho(\ft)$ has
non-degenerate spectrum, the matrix $\rho(\fs)$ differs from that of an MD representation
only by a conjugation by signed diagonal matrix.  In this case, $\rho(\fs)$
must have a row that contains no zero (i.e. the row corresponding to the unit object).

\item Let  $\rho(\fs)^\text{ndeg}$ (or $M^\text{ndeg}$) be the non-degenerate block of $\rho(\fs)$ (or $M$), (i.e., corresponding to the multiplicity 1 eigenvalues of the diagonal matrix $\rho(\ft)$, see section \ref{strat}). Then the conductor of $\rho(\fs)^\text{ndeg}$  divides $\ord(\rho(\ft))$ {(cf. Proposition \ref{p:ndeg})}.  If the $\rho(\ft)$-spectrum is non-degenerate then we may drop the $\text{ndeg}$ superscript.

\item $\s(\rho(\fs)^\text{ndeg}) = \big( \rho^a(\ft) \rho(\fs) \rho^b(\ft)
\rho(\fs) \rho^a(\ft) \big)^\text{ndeg}$ for any $\s \in \GQ$, where $\s(\zeta_n)
= \zeta_n^a$ for an unique integer $a$ modulo $n$.  Here $n =\ord(\rho(\ft))$
and $b$ satisfies $ab \equiv 1 \mod n$ {(cf. Theorem \ref{p:MD1}).}  Again, this is because
$\rho(\fs)^\text{ndeg}$ can only differ from that of an MD representation by a
conjugation of signed diagonal matrix. 
%\blue{$\s(\rho(\fs)) = \big( \rho^a(\ft) \rho(\fs) \rho^b(\ft)
%\rho(\fs) \rho^a(\ft) \big)$ holds in a part of matrix, the non-degenerate block.}

\end{enumerate}

Since the weakly integral rank-6 MD of MTCs are classified, we can exclude
symmetric $\SL$ representations that must produce such MDs.  Thus we also reject
the representations that satisfy the following conditions, both of which imply weak integrality:
\begin{enumerate}

\item $\pord(\rho(\ft)) \in \{2,3,4,6\}$. In fact, this implies the category is pointed, see Proposition \ref{p:weakly_int}(i).
\item
The squares of the matrix entries of $\rho(\fs)$ in each row containing no zeros are all rational numbers, and 
 $\rho$ is non-degenerate.  Indeed, in this case $1/D^2,(d_i/D)^2$ and $(d_i\FPdim(X_i)/D)^2$ are rational, where column $i$ is the unique strictly positive (or negative) column. (This condition only
rejects one case. See entry 566 in Supplementary material Section \ref{Section4}.)

\end{enumerate}
For details, see Supplementary material Section \ref{Section4}.

We remark that there are 6-dimensional irreducible $\SL$
representations where $\rho(\ft)$ are degenerate, for example, the representation $6_5^1$ in Appendix A.
Such a representation is rejected since the conductor
of $\rho(\fs)^\text{ndeg}$ is 40 which does not divides
$ \ord(\rho(\ft)) = 5$ (see also entry 582 in Supplementary material Section \ref{Section4}). 

All the passing symmetric $\SL$  representations can be grouped into orbits
generated by Galois conjugations and tensoring 1-dimensional representations.
There are 7 such orbits. A representative for each orbit is listed in Section
\ref{repL}, which have $(\text{dims};\text{levels})=(6;9)$, $(6;13)$, $(6;15)$, $(6;16)$, $(6;35)$,
$(6;56)$, $(6;80)$.  

Fortunately, we find that all these $\SL$ representations have non-degenerate
$\rho(\ft)$, so they can only possibly differ from an MD representation by a
conjugation of signed diagonal matrix, if they indeed are associated with MDs.
We can then search through the finite number of signed diagonal conjugations,
and find the $(S,T)$ matrices that satisfy the conditions listed in Theorems
\ref{p:MD} and \ref{p:MD1}.  The results are given in Section \ref{SsL}, where
$(S,T)$ matrices are found from $\SL$ representations that have $(\text{dims};\text{levels})=(6;9)$,
$(6;13)$, $(6;16)$, $(6;35)$, $(6;56)$, $(6;80)$.  Those computer assisted
calculations are described in detail in the Appendix.

\section{Classification of modular data of rank=6: non-admissible types}

In this section, we complete the classification of rank=$6$ MDs by eliminating the remaining types.

\begin{thm}
There are no rank=$6$ MTCs of types $(3,1,1,1), (2,2,2), (2,2,1,1), (2,1,1,1,1), (5,1)$, or $(1,1,1,1,1,1)$.
\end{thm}

Obviously, type $Vec$ is the only MTC of type (1). However, no MTCs of rank $n > 1$ is of type $(1,\dots, 1)$, as the associated $\SL$ representations $\rho \cong n\chi^i$ for some integer $i$ by Corollary \ref{l:3}. In particular, $\rho(\fs)$ has zeros in each row if $n >1$. 

\subsection{Nonexistence of type $(3,1,1,1)$}
\begin{prop}
There does not exist any modular tensor category of type (3,1,1,1)\,.
\end{prop}
\begin{proof}
Assume  contrary. Let $\CC$ be a modular tensor category of type $(3,1,1,1)$ and $\rho$ an $\SL$  representation of  $\CC$. 
Then 
$$
\rho \cong \rho_0 \oplus \chi_1 \oplus \chi_2 \oplus \chi_3\,.
$$
where $\rho_0$ is irreducible of dimension 3 and  $\chi_i, i=1,2,3$,  are 1-dimensional representations. By Lemma \ref{l:2}, $\spec(\chi_i(\ft)) \subset \spec(\rho_0(\ft))$ for $i=1,2,3$. One may assume $\rho_0$ has a minimal $\ft$-spectrum. Then $\rho_0$ must have a prime power level. By Appendix \ref{repPP}, the level of $\rho_0$ can only be 3, 4, 5, 7, 8 or 16. The $\ft$-spectrum of any 3-dimensional irreducible representations of level 7 or 16 does not contain any 12-th root of unity. Therefore, the level of $\rho_0$ can only be $3,4, 5$, or $8$. It suffices to show that none of these levels is possible.

If $\rho_0$ were of level 3 or 4, then $\ord(\rho(\ft))=3$ or $4$, by Lemma \ref{l:2}. This implies $\ord(T)=2, 3$ or $4$ and hence $\CC$ is integral by Theorem \ref{t:2346}. By Proposition \ref{p:weakly_int}, $\CC$ must be of type (4,2), a contradiction. Therefore, $\rho_0$ can only be of level 5 or 8. 

If $\rho_0$ were of level 5, then $\ord(\rho(\ft))=5$ by Lemma \ref{l:2}. Hence, $\ord(T)=5$ which is not possible by Proposition \ref{p:d111}.   

If the level of $\rho_0$ were 8, then $\rho_0 \cong \rd{3}{8}{1,0}$ or $\rd{3}{8}{3,0}$  as they are the 3-dimensional irreducible representations of level 8 with a minimal $\ft$-spectrum. In either case, $\spec(\rho_0(\ft))$ has exactly one 12-th root of unity, which is 1, and $\rho_0$ is odd.  Therefore,  $\rho \cong \rho_0 \oplus 3 \chi^0$ by Corollary \ref{l:3}.  This implies $\Tr(\rho(\fs^2))=0$, which is impossible for any MD representation.
\end{proof}

\subsection{Nonexistence of types (2,2,2), (2,2,1,1) and (2,1,1,1,1)}

We will prove the following theorem which leads to the nonexistence of modular tensor categories of these types.
\begin{thm} \label{t:h2}
Let $\CC$ a be modular tensor category with $\rank \CC > 2$, and $\rho$ an $\SL$ representation of $\CC$. If all the irreducible subrepresentations of $\rho$ have dimensions $\le 2$, then $\ord(T)=1, 2,3, 4$, or $6$ and therefore $\CC$ is  integral. 
\end{thm}
\begin{proof}
If every irreducible subrepresentation of $\rho$ is  1-dimensional, then $\CC$ is of type $(1,\dots, 1)$ which can only be trivial by the beginning remark of this section. In particular, $\ord(T)=1$ and $\CC$ is integral.

Now, we assume $\rho$ admits a 2-dimensional irreducible subrepresentation $\rho_0$. By tensoring a 1-dimensional representation to $\rho$, we may assume the level of $\rho_0$ to be $2,3,5$, or $8$. 

Suppose $\rho_0$ is of level 5. Then each irreducible subrepresentations $\rho_0'$ of  $\rho$ which is not isomorphic to $\rho_0$ satisfies $\spec(\rho_0'(\ft)) \cap \spec(\rho_0(\ft)) =\emptyset$ by Appendix \ref{repPP}. This implies $\rho \cong \ell \rho_0$, but this is impossible by Proposition \ref{p:multiple}. Therefore, $\rho_0$ cannot have level 5.

Assume $\rho_0$ is of level 8. Note that the $\ft$-spectrum of any 2-dimensional level 8 irreducible representation consists of primitive 8-th roots of unity.  By the  $\ft$-spectra criterion and Appendix \ref{repPP}, all the irreducible  subrepresentations of $\rho$ are of dimension 2 and level 8. In particular,   $\ord(T)=\pord(\rho(\ft))=4$. 

If $\rho_0$ is of level 2 or 3, it follows from the preceding discussion that all the 2-dimensional irreducible subrepresentations of $\rho$ are of level 2 or 3. By Lemma \ref{l:2}, $\ord(\rho(\ft))=2,3$ or $6$ and so $\ord(T) =2,3$ or $6$. 

The last assertion follows from Theorem \ref{t:2346}. 
\end{proof}
\begin{cor}
There is no modular tensor category of types $(2,2,2)$, $(2,2,1,1)$ or $(2,1,1,1,1)$.
\end{cor}
\begin{proof}
Suppose there exists a modular tensor category $\CC$ of any of these types. By Theorem \ref{t:h2}, $\CC$ is integral, but this contradicts Proposition \ref{p:weakly_int} which shows $\CC$ is of type $(4,2)$.  
\end{proof}

\subsection{Nonexistence of type $(5,1)$}

Suppose that $\CC$ is a modular tensor category of type $(5,1)$, and $\rho$ an $\SL$ representation  of $\CC$. Then $\CC$ is not integral by Proposition \ref{p:weakly_int}, and   $\rho \cong \rho_0 \oplus\rho_1$ where $\rho_0, \rho_1$ are irreducible of dimension $5$ and $1$ respectively.  By tensoring a 1-dimensional representation of $\SL$, one may assume $\rho_0$ is of prime power level. By Appendix \ref{repPP}, the level of $\rho_0$ can only be 11 or 5. 

In the former case the $\ft$-spectrum consists primitive 11-th roots of unity.   Since $\rho_1(\ft)$ is a $12$ root of unity, the $\ft$-spectrum criteria shows this is impossible. 

Now if $\rho_2$ has  level $5$ and $\rho_2 \cong \rd{5}{5}{1}$. This implies $\rho_1 \cong \chi^0$. Let $\tilde\rho = \chi^0 \oplus \rd{5}{5}{1}$. Then $\tilde\rho(\ft) =\diag(1,1,\zeta_5, \zeta_5^2, \zeta_5^3, \zeta_5^4)$, and
$$
\tilde\rho(\fs)=
[1] \oplus \left[
\begin{array}{cccccc}
 1 & 0 & 0 & 0 & 0 & 0 \\
 0 & -\frac{1}{5} & \frac{\sqrt{6}}{5} & \frac{\sqrt{6}}{5} & \frac{\sqrt{6}}{5} & \frac{\sqrt{6}}{5} \\
 0 & \frac{\sqrt{6}}{5} & \frac{3-\sqrt{5}}{10} & -\frac{1}{5}-\frac{1}{\sqrt{5}} & \frac{1}{\sqrt{5}}-\frac{1}{5} & \frac{3+\sqrt{5}}{10} \\
 0 & \frac{\sqrt{6}}{5} & -\frac{1}{5}-\frac{1}{\sqrt{5}} & \frac{3+\sqrt{5}}{10} & \frac{3-\sqrt{5}}{10} & \frac{1}{\sqrt{5}}-\frac{1}{5} \\
 0 & \frac{\sqrt{6}}{5} & \frac{1}{\sqrt{5}}-\frac{1}{5} & \frac{3-\sqrt{5}}{10} & \frac{3+\sqrt{5}}{10} & -\frac{1}{5}-\frac{1}{\sqrt{5}} \\
 0 & \frac{\sqrt{6}}{5} & \frac{3+\sqrt{5}}{10} & \frac{1}{\sqrt{5}}-\frac{1}{5} & -\frac{1}{5}-\frac{1}{\sqrt{5}} & \frac{3-\sqrt{5}}{10} \\
\end{array}
\right]\,.
$$
There exists a  real orthogonal matrix $U=\diag(f,\ve_1,\ve_2,\ve_3,\ve_4)$ such that $\rho(\fs)=U\tilde{\rho}(\fs)U^\top$ and $\rho(\ft) =\tilde{\rho}(\ft)$ , where $f\in O_2(\BR)$ and $\ve_i = \pm 1$. 

The group $\Gal(\BQ_5/\BQ)$ is generated by $\sigma$ defined by $\s(\zeta_5)=\zeta_5^2$, and 
$$
D_{\tilde{\rho}}(\s) = I_2 \oplus J_4  \quad \text{ where } J_4=[\delta_{i, 5-j}]_{1 \le i,j \le 4}\,.
$$
So $\hs$ fixes 1 and 2.  Since  $\CC$ is non-integral, the row corresponding to $\1$ must be one of the last 4. Since $\rho(\fs)^\nd$ and $\tilde\rho(\fs)^\nd$ are the same up to some signs, $D = \frac{10}{3\pm \sqrt{5}}$ which has norm 25.

Observe that each row of $\tilde\rho(\fs)^\nd$ has the entries $-\frac{1}{5} \pm \frac{1}{\sqrt{5}}$. Therefore, $(-\frac{1}{5} \pm \frac{1}{\sqrt{5}})/\frac{3\mp \sqrt{5}}{10}=1\pm \sqrt{5}$ are dimensions of some objects up to a sign. However, their norms are -4 which is not a divisor of 25, a contradiction.   So, we conclude that such a category cannot exist.
\section{Summary and Future Directions}

 We have developed tools for classifying modular data directly from representations of $\SL$, and have applied them to provide a classification of rank $6$ modular data.  Sufficiently many of these tools have been implemented as computer algorithms to yield a purely computational approach to the rank $6$ classification.  A purely ``by hand" approach to higher ranks is too involved for the currently theory, but the computational approach can be implemented in higher ranks.  It should be noted that in this work we used the classification of weakly integral modular data \cite{CWeakly} of rank up to $7$ to simplify the computer calculations.  For higher ranks this will require further work. \medskip \\

\noindent
{\bf Acknowledgements:}

The authors are partially supported by an NSF FRG grant: Z.W. by DMS-1664351, E.C.R. by DMS-1664359, X.-G. W. by DMS-1664412 and S.-H. N. by DMS-1664418. Z.W. is also partially supported by CCF 2006463 and ARO MURI contract W911NF-20-1-0082.

\appendix

\section{List of $\SL$ irreducible representations of prime-power levels}
\label{repPP}

In this section, we list all the $\SL$ irreducible representations of dimension
1 -- 6, whose level ($l=\ord(\rho(\ft))$) is a power of single prime number.
In the list, $\rho(\mathfrak{t})$ is presented in term of topological spins
$(\tilde s_{1},\tilde s_{2},\cdots)$ ($\tilde s_{i} =
\arg(\rho_a(\mathfrak{t})_{ii})$).  

Note that $\rho(\mathfrak{s})$ is symmetric and $\rho(\mathfrak{s})_{ij}$'s are
either all real or all imaginary.  When $\rho(\mathfrak{s})_{ij}$'s are all
real, $\rho(\mathfrak{s})$ is presented as $(\rho_{11}, \rho_{12}, \rho_{13},
\rho_{14}, \cdots;\ \  \rho_{22}, \rho_{23}, \rho_{24}, \cdots)$. In this case,
$\rho(\mathfrak{s})^2=\id$ and the representation $\rho$ is said to be even.
When $\rho(\mathfrak{s})_{ij}$'s are all imaginary, $\rho(\mathfrak{s})$ is
presented as $\ii(-\ii\rho_{11}, -\ii\rho_{12}, -\ii\rho_{13}, -\ii\rho_{14},
\cdots;$ $ -\ii\rho_{22}, -\ii\rho_{23}, -\ii\rho_{24}, \cdots)$, or 
as $(s_n^m)^{-1}(s_n^m \rho_{11}, s_n^m \rho_{12}, s_n^m \rho_{13}, s_n^m
\rho_{14},$ $\cdots;$ $ s_n^m \rho_{22}, s_n^m \rho_{23}, s_n^m \rho_{24},
\cdots)$, where $s_n^m := \zeta_n^m-\zeta_n^{-m}$. In this case,
$\rho(\mathfrak{s})^2=-\id$ and the representation $\rho$ is said to be odd.
In any case, the numbers inside the bracket $(\cdots)$ are all real.  We can
tell a representation to be even or odd by the absence or the presence of $\ii$
or $(s_n^m)^{-1}$in front of the bracket $(\cdots)$.

We note that two symmetric representations are equivalent up to a permutation
of indices, and a conjugation of signed diagonal matrix.  To choose the
ordering in indices, we introduce arrays $O_i =[\text{DenominatorOf}(\tilde
s_{i}),\tilde s_{i}, \rho_{ii}]$.  The order of two arrays is determined by
first compare the lengths of the two arrays. If the lengths are equal, we then
compare the first elements of the two arrays.  If the first elements are equal,
we then compare the second elements of the two arrays, {\it etc}.  To compare
two cyclotomic numbers, here we used the ordering of cyclotomic numbers
provided by GAP computer algebraic system.  We order the indices to make
$O_1\leq O_2 \leq O_3 \cdots$.  The  conjugation of signed diagonal matrix is
chosen to make $-\rho(\mathfrak{s})_{1j} < \rho(\mathfrak{s})_{1j} $ for
$j=2,3,\cdots$.  If $\rho(\mathfrak{s})_{1j}=0$, we will try to make
$-\rho(\mathfrak{s})_{2j} < \rho(\mathfrak{s})_{2j} $, {\it etc}.

All the prime-power-level irreducible representations are labeled by index
$d_{l,k}^{a,m}$, where $d$ is the dimension and $l$ is the level of the
representation.  The irreducible representations of a given $d,l$ can be
grouped into several orbits, generated by Galois conjugations and tensoring of
1-dimensional representations that do not change the level $l$.  $k$ in
$d_{l,k}^{a,m}$ labels those different orbits.  If there is only 1 orbit for a
given $d,l$, the index $k$ will be dropped.  

The irreducible representation labeled by $d_{l,k}^{a,m}$ is generated from the
irreducible representation labeled by $d_{l,k}^{1,0}$ via the following Galois
conjugations and tensoring of 1-dimensional representations
\begin{align}
 \rho_{d_{l,k}^{a,m}}(\ft) &= \s_a\big(\rho_{d_{l,k}^{1,0}}(\ft) \big)\ee^{2\pi \ii\frac m{12}} 
\nonumber\\
 \rho_{d_{l,k}^{a,m}}(\fs) &= \s_a\big(\rho_{d_{l,k}^{1,0}}(\fs) \big)\ee^{-2\pi \ii\frac m{4}} 
\end{align}
where the Galois conjugation $\s_a $ is in $ \Gal(\BQ_n)$ with $n$ be the least
common multiple of $\ord(\rho_{d_{l,k}^{1,0}}(\ft))$ and the conductor of
$\rho_{d_{l,k}^{1,0}}(\fs)$. The  Galois conjugation $\s_a $ is labeled by an
integer $a$, which is given by
\begin{align}
 \s_a\big(\ee^{2\pi \ii /n} \big) = 
 \ee^{2\pi \ii a/n} .
\end{align}
Also $m \in \BZ_{12}$ is such that  $ \ord(\rho_{d_{l,k}^{1,0}}(\ft)\ee^{2\pi
\ii\frac m{12}})= \ord(\rho_{d_{l,k}^{1,0}}(\ft))$. Due to this condition, when
$l$ is not divisible by 2 and 3, $m$ can only be $0$.  In this case, we will
drop $m$.  Here we choose $d_{l,k}^{1,0}$ to be the representation in the orbit
with minimal $[\tilde s_1, \tilde s_2,\cdots] $.

%We like to remark that two representations labeled by $d_{l,k}^{a,m}$ and
%$d_{l,k}^{a',m'}$ may be equivalent.  
The numbers of distinct irreducible
representations with prime-power level (PPL) in each dimension are given by
\setlength\tabcolsep{5pt}
\begin{align}
\begin{tabular}{|r|r|r|r|r|r|r|r|r|r|r|r|r|}
\hline
\text{dim}: & 1 & 2 & 3 & 4 & 5 & 6 & 7 & 8 & 9 & 10 & 11 & 12 \\
\hline
\text{\# of irreps with PPL} & 6 & 11 & 33 & 18 & 3 & 70 & 3  & 10  & 4 & 7 & 3 & 176  \\ 
\hline
\text{\# of irreps} & 12 & 54 & 136 & 180 & 36 & 720 & 36 & 456 & 476 & 222 & 36 & 3214 \\ 
\hline
\end{tabular}
\end{align}
In the above we also list the numbers of distinct irreducible representations,
which are tensor products of the irreducible representations with prime-power
levels.

In the following tables, we list all irreducible representations with
prime-power levels for rank $2,3,4,5$.  For rank 6, to save space, we only list
all irreducible representations with prime-power levels that have a form
$\rho_{d_{l,k}^{1,0}}$. Other irreducible representations, with prime-power
levels and the same dimension, can be obtained from those listed ones via
Galois conjugations and tensoring 1-dimensional representations.  In
Supplementary Material Section \ref{Section1}, we list all distinct irreducible
representations of prime-power levels.  In the table $c_n^m := \zeta_n^m +
\zeta_n^{-m}$ and $s_n^m := \zeta_n^m - \zeta_n^{-m}$.

%number of dim 1 irreps: 12 number of dim 2 irreps: 54 number of dim 3 irreps:
%136 number of dim 4 irreps: 180 number of dim 5 irreps: 36 number of dim 6
%irreps: 720 number of dim 7 irreps: 36 number of dim 8 irreps: 456 number of
%dim 9 irreps: 476 number of dim 10 irreps: 222 number of dim 11 irreps: 36
%number of dim 12 irreps: 3214

\def\arraystretch{1.6} \setlength\tabcolsep{3pt}
\begin{longtable}{|l|l|l|}
\hline
$d_{l,k}^{a,m}$ & \# & $\rho(\mathfrak{t})$, \ $\rho(\mathfrak{s})$\\
\hline
$1_{1}^{1}$ & 1 & $( 0 )
$, ~ ($1$)
\\ 
\hline
$1_{2}^{1,0}$ & 2 & $( \frac{1}{2} )
$, ~ ($-1$)
\\ 
\hline
$1_{3}^{1,0}$ & 3 & $( \frac{1}{3} )
$, ~ ($1$)
\\ 
\hline
$1_{3}^{1,4}$ & 4 & $( \frac{2}{3} )
$, ~ ($1$)
\\ 
\hline
$1_{4}^{1,0}$ & 5 & $( \frac{1}{4} )
$, ~ $\mathrm{i}$($1$)
\\ 
\hline
$1_{4}^{1,6}$ & 6 & $( \frac{3}{4} )
$, ~ $\mathrm{i}$($-1$)
\\ 
\hline
\end{longtable}

\begin{longtable}{|l|l|l|}
\hline
$d_{l,k}^{a,m}$ & \# & $\rho(\mathfrak{t})$, \ $\rho(\mathfrak{s})$\\
\hline
$2_{2}^{1,0}$ & 1 & $( 0,
\frac{1}{2} )
$, ~ ($-\frac{1}{2}$,
$-\frac{\sqrt{3}}{2}$;
$\frac{1}{2}$)
\\ 
\hline
$2_{3}^{1,0}$ & 2 & $( 0,
\frac{1}{3} )
$, ~ $(s_{3}^{1})^{-1}$($1$,
$-\sqrt{2}$;\ \ 
$-1$)
\\ 
\hline
$2_{3}^{1,8}$ & 3 & $( 0,
\frac{2}{3} )
$, ~ $(s_{3}^{1})^{-1}$($-1$,
$-\sqrt{2}$;\ \ 
$1$)
\\ 
\hline
$2_{3}^{1,4}$ & 4 & $( \frac{1}{3},
\frac{2}{3} )
$, ~ $(s_{3}^{1})^{-1}$($1$,
$-\sqrt{2}$;\ \ 
$-1$)
\\ 
\hline
$2_{4}^{1,0}$ & 5 & $( \frac{1}{4},
\frac{3}{4} )
$, ~ $\mathrm{i}$($-\frac{1}{2}$,
$\frac{\sqrt{3}}{2}$;\ \ 
$\frac{1}{2}$)
\\ 
\hline
$2_{5}^{1}$ & 6 & $( \frac{1}{5},
\frac{4}{5} )
$, ~ $(s_{5}^{1})^{-1}$($1$,
$-\frac{1+\sqrt{5}}{2}$;\ \ 
$-1$)
\\ 
\hline
$2_{5}^{2}$ & 7 & $( \frac{2}{5},
\frac{3}{5} )
$, ~ $(s_{5}^{2})^{-1}$($1$,
$\frac{1-\sqrt{5}}{2}$;\ \ 
$-1$)
\\ 
\hline
$2_{8}^{1,0}$ & 8 & $( \frac{1}{8},
\frac{3}{8} )
$, ~ ($-\frac{\sqrt{2}}{2}$,
$\frac{\sqrt{2}}{2}$;
$\frac{\sqrt{2}}{2}$)
\\ 
\hline
$2_{8}^{1,9}$ & 9 & $( \frac{1}{8},
\frac{7}{8} )
$, ~ $\mathrm{i}$($-\frac{\sqrt{2}}{2}$,
$\frac{\sqrt{2}}{2}$;\ \ 
$\frac{\sqrt{2}}{2}$)
\\ 
\hline
$2_{8}^{1,3}$ & 10 & $( \frac{3}{8},
\frac{5}{8} )
$, ~ $\mathrm{i}$($-\frac{\sqrt{2}}{2}$,
$\frac{\sqrt{2}}{2}$;\ \ 
$\frac{\sqrt{2}}{2}$)
\\ 
\hline
$2_{8}^{1,6}$ & 11 & $( \frac{5}{8},
\frac{7}{8} )
$, ~ ($\frac{\sqrt{2}}{2}$,
$\frac{\sqrt{2}}{2}$;
$-\frac{\sqrt{2}}{2}$)
\\ 
\hline
\end{longtable}

\begin{longtable}{|l|l|l|}
\hline
$d_{l,k}^{a,m}$ & \# & $\rho(\mathfrak{t})$, \ $\rho(\mathfrak{s})$\\
\hline
$3_{3}^{1,0}$ & 1 & $( 0,
\frac{1}{3},
\frac{2}{3} )
$, ~ ($-\frac{1}{3}$,
$\frac{2}{3}$,
$\frac{2}{3}$;
$-\frac{1}{3}$,
$\frac{2}{3}$;
$-\frac{1}{3}$)
\\ 
\hline
$3_{4}^{1,0}$ & 2 & $( 0,
\frac{1}{4},
\frac{3}{4} )
$, ~ ($0$,
$\frac{\sqrt{2}}{2}$,
$\frac{\sqrt{2}}{2}$;
$-\frac{1}{2}$,
$\frac{1}{2}$;
$-\frac{1}{2}$)
\\ 
\hline
$3_{4}^{1,3}$ & 3 & $( 0,
\frac{1}{2},
\frac{1}{4} )
$, ~ $\mathrm{i}$($-\frac{1}{2}$,
$\frac{1}{2}$,
$\frac{\sqrt{2}}{2}$;\ \ 
$-\frac{1}{2}$,
$\frac{\sqrt{2}}{2}$;\ \ 
$0$)
\\ 
\hline
$3_{4}^{1,9}$ & 4 & $( 0,
\frac{1}{2},
\frac{3}{4} )
$, ~ $\mathrm{i}$($\frac{1}{2}$,
$\frac{1}{2}$,
$\frac{\sqrt{2}}{2}$;\ \ 
$\frac{1}{2}$,
$-\frac{\sqrt{2}}{2}$;\ \ 
$0$)
\\ 
\hline
$3_{4}^{1,6}$ & 5 & $( \frac{1}{2},
\frac{1}{4},
\frac{3}{4} )
$, ~ ($0$,
$\frac{\sqrt{2}}{2}$,
$\frac{\sqrt{2}}{2}$;
$\frac{1}{2}$,
$-\frac{1}{2}$;
$\frac{1}{2}$)
\\ 
\hline
$3_{5}^{1}$ & 6 & $( 0,
\frac{1}{5},
\frac{4}{5} )
$, ~ ($\frac{\sqrt{5}}{5}$,
$-\frac{\sqrt{10}}{5}$,
$-\frac{\sqrt{10}}{5}$;
$-\frac{5+\sqrt{5}}{10}$,
$\frac{5-\sqrt{5}}{10}$;
$-\frac{5+\sqrt{5}}{10}$)
\\ 
\hline
$3_{5}^{3}$ & 7 & $( 0,
\frac{2}{5},
\frac{3}{5} )
$, ~ ($-\frac{\sqrt{5}}{5}$,
$-\frac{\sqrt{10}}{5}$,
$-\frac{\sqrt{10}}{5}$;
$-\frac{5-\sqrt{5}}{10}$,
$\frac{5+\sqrt{5}}{10}$;
$-\frac{5-\sqrt{5}}{10}$)
\\
\hline
\end{longtable}

\begin{longtable}{|l|l|p{4.5in}|}
\hline
$d_{l,k}^{a,m}$ & \# & $\rho(\mathfrak{t})$, \ $\rho(\mathfrak{s})$\\
\hline
$4_{5,1}^{1}$ & 1 & $( \frac{1}{5},
\frac{2}{5},
\frac{3}{5},
\frac{4}{5} )
$, ~ $(s_{5}^{2})^{-1}$($-\frac{5+\sqrt{5}}{10}$,
$-\frac{\sqrt{15}}{5}$,
$\frac{3-3\sqrt{5}}{2\sqrt{15}}$,
$\frac{5-3\sqrt{5}}{10}$;\ \ 
$-\frac{5-3\sqrt{5}}{10}$,
$\frac{5+\sqrt{5}}{10}$,
$-\frac{3-3\sqrt{5}}{2\sqrt{15}}$;\ \ 
$\frac{5-3\sqrt{5}}{10}$,
$-\frac{\sqrt{15}}{5}$;\ \ 
$\frac{5+\sqrt{5}}{10}$)
\\ 
\hline
$4_{5,2}^{1}$ & 2 & $( \frac{1}{5},
\frac{2}{5},
\frac{3}{5},
\frac{4}{5} )
$, ~ ($\frac{\sqrt{5}}{5}$,
$-\frac{5-\sqrt{5}}{10}$,
$-\frac{5+\sqrt{5}}{10}$,
$\frac{\sqrt{5}}{5}$;
$-\frac{\sqrt{5}}{5}$,
$\frac{\sqrt{5}}{5}$,
$\frac{5+\sqrt{5}}{10}$;
$-\frac{\sqrt{5}}{5}$,
$\frac{5-\sqrt{5}}{10}$;
$\frac{\sqrt{5}}{5}$)
\\ 
\hline
$4_{7}^{1}$ & 3 & $( 0,
\frac{1}{7},
\frac{2}{7},
\frac{4}{7} )
$, ~ $\mathrm{i}$($-\frac{\sqrt{7}}{7}$,
$\frac{\sqrt{14}}{7}$,
$\frac{\sqrt{14}}{7}$,
$\frac{\sqrt{14}}{7}$;\ \ 
$-\frac{c_{7}^{2}}{\sqrt{7}}$,
$-\frac{c_{7}^{1}}{\sqrt{7}}$,
$-\frac{c_{7}^{3}}{\sqrt{7}}$;\ \ 
$-\frac{c_{7}^{3}}{\sqrt{7}}$,
$-\frac{c_{7}^{2}}{\sqrt{7}}$;\ \ 
$-\frac{c_{7}^{1}}{\sqrt{7}}$)
\\ 
\hline
$4_{7}^{3}$ & 4 & $( 0,
\frac{3}{7},
\frac{5}{7},
\frac{6}{7} )
$, ~ $\mathrm{i}$($\frac{\sqrt{7}}{7}$,
$\frac{\sqrt{14}}{7}$,
$\frac{\sqrt{14}}{7}$,
$\frac{\sqrt{14}}{7}$;\ \ 
$\frac{c_{7}^{1}}{\sqrt{7}}$,
$\frac{c_{7}^{2}}{\sqrt{7}}$,
$\frac{c_{7}^{3}}{\sqrt{7}}$;\ \ 
$\frac{c_{7}^{3}}{\sqrt{7}}$,
$\frac{c_{7}^{1}}{\sqrt{7}}$;\ \ 
$\frac{c_{7}^{2}}{\sqrt{7}}$)
\\ 
\hline
$4_{8}^{1,0}$ & 5 & $( \frac{1}{8},
\frac{3}{8},
\frac{5}{8},
\frac{7}{8} )
$, ~ $\mathrm{i}$($\frac{\sqrt{2}}{4}$,
$\frac{\sqrt{6}}{4}$,
$\frac{\sqrt{6}}{4}$,
$\frac{\sqrt{2}}{4}$;\ \ 
$\frac{\sqrt{2}}{4}$,
$-\frac{\sqrt{2}}{4}$,
$-\frac{\sqrt{6}}{4}$;\ \ 
$-\frac{\sqrt{2}}{4}$,
$\frac{\sqrt{6}}{4}$;\ \ 
$-\frac{\sqrt{2}}{4}$)
\\
\hline
\end{longtable}

\begin{longtable}{|l|l|p{4.5in}|}
\hline
$d_{l,k}^{a,m}$ & \# & $\rho(\mathfrak{t})$, \ $\rho(\mathfrak{s})$\\
\hline
$5_{5}^{1}$ & 1 & $( 0,
\frac{1}{5},
\frac{2}{5},
\frac{3}{5},
\frac{4}{5} )
$, ~ ($-\frac{1}{5}$,
$\frac{\sqrt{6}}{5}$,
$\frac{\sqrt{6}}{5}$,
$\frac{\sqrt{6}}{5}$,
$\frac{\sqrt{6}}{5}$;
$\frac{3-\sqrt{5}}{10}$,
$-\frac{1+\sqrt{5}}{5}$,
$-\frac{1-\sqrt{5}}{5}$,
$\frac{3+\sqrt{5}}{10}$;
$\frac{3+\sqrt{5}}{10}$,
$\frac{3-\sqrt{5}}{10}$,
$-\frac{1-\sqrt{5}}{5}$;
$\frac{3+\sqrt{5}}{10}$,
$-\frac{1+\sqrt{5}}{5}$;
$\frac{3-\sqrt{5}}{10}$)
\\ 
\hline
$5_{11}^{1}$ & 2 & $( \frac{1}{11},
\frac{3}{11},
\frac{4}{11},
\frac{5}{11},
\frac{9}{11} )
$, ~ ($-\frac{c_{44}^{3}}{\sqrt{11}}$,
$-\frac{c_{44}^{7}}{\sqrt{11}}$,
$-\frac{c_{44}^{5}}{\sqrt{11}}$,
$-\frac{c_{44}^{1}}{\sqrt{11}}$,
$-\frac{c_{44}^{9}}{\sqrt{11}}$;
$\frac{c_{44}^{9}}{\sqrt{11}}$,
$-\frac{c_{44}^{3}}{\sqrt{11}}$,
$\frac{c_{44}^{5}}{\sqrt{11}}$,
$\frac{c_{44}^{1}}{\sqrt{11}}$;
$\frac{c_{44}^{1}}{\sqrt{11}}$,
$\frac{c_{44}^{9}}{\sqrt{11}}$,
$\frac{c_{44}^{7}}{\sqrt{11}}$;
$\frac{c_{44}^{7}}{\sqrt{11}}$,
$-\frac{c_{44}^{3}}{\sqrt{11}}$;
$\frac{c_{44}^{5}}{\sqrt{11}}$)
\\ 
\hline
$5_{11}^{2}$ & 3 & $( \frac{2}{11},
\frac{6}{11},
\frac{7}{11},
\frac{8}{11},
\frac{10}{11} )
$, ~ ($\frac{c_{44}^{5}}{\sqrt{11}}$,
$\frac{c_{44}^{3}}{\sqrt{11}}$,
$-\frac{c_{44}^{7}}{\sqrt{11}}$,
$-\frac{c_{44}^{1}}{\sqrt{11}}$,
$-\frac{c_{44}^{9}}{\sqrt{11}}$;
$\frac{c_{44}^{7}}{\sqrt{11}}$,
$\frac{c_{44}^{9}}{\sqrt{11}}$,
$\frac{c_{44}^{5}}{\sqrt{11}}$,
$\frac{c_{44}^{1}}{\sqrt{11}}$;
$\frac{c_{44}^{1}}{\sqrt{11}}$,
$-\frac{c_{44}^{3}}{\sqrt{11}}$,
$\frac{c_{44}^{5}}{\sqrt{11}}$;
$\frac{c_{44}^{9}}{\sqrt{11}}$,
$\frac{c_{44}^{7}}{\sqrt{11}}$;
$-\frac{c_{44}^{3}}{\sqrt{11}}$)
\\ 
\hline
\end{longtable}

\def\arraystretch{1.6} \setlength\tabcolsep{3pt}

\begin{longtable}{|l|l|p{4.5in}|}
\hline
$d_{l,k}^{a,m}$ & \# & $\rho(\mathfrak{t})$, \ $\rho(\mathfrak{s})$\\
\hline
$6_{5}^{1}$ & 1 & $( 0,
0,
\frac{1}{5},
\frac{2}{5},
\frac{3}{5},
\frac{4}{5} )
$, ~ $(s_{5}^{2})^{-1}$($\frac{\sqrt{5}}{5}$,
$-\frac{5-\sqrt{5}}{10}$,
$-\frac{\sqrt{10}}{5}$,
$-\frac{5+\sqrt{5}}{10}$,
$\frac{5-3\sqrt{5}}{10}$,
$\frac{1-\sqrt{5}}{\sqrt{10}}$;\ \ 
$-\frac{\sqrt{5}}{5}$,
$\frac{-1+\sqrt{5}}{\sqrt{10}}$,
$\frac{5-3\sqrt{5}}{10}$,
$\frac{5+\sqrt{5}}{10}$,
$-\frac{\sqrt{10}}{5}$;\ \ 
$\frac{5-\sqrt{5}}{10}$,
$\frac{1-\sqrt{5}}{\sqrt{10}}$,
$-\frac{\sqrt{10}}{5}$,
$-\frac{\sqrt{5}}{5}$;\ \ 
$-\frac{\sqrt{5}}{5}$,
$\frac{5-\sqrt{5}}{10}$,
$\frac{\sqrt{10}}{5}$;\ \ 
$\frac{\sqrt{5}}{5}$,
$\frac{1-\sqrt{5}}{\sqrt{10}}$;\ \ 
$-\frac{5-\sqrt{5}}{10}$)
\\
\hline
\end{longtable}

\section{A list of all candidate $\SL$ representations of MTCs }

We will follow the strategy out lined in Section \ref{strat}.  We first try of
obtain a list that includes all $\SL$ representations associated with MTCs.
Certainly, one such list is the list of all $\SL$ representations of finite
levels.  But such a list is very inefficient since most representations in the
list are not associated with MTC's.  So this section, we collect the conditions
that a representation of MTC must satisfy, to obtain a shorter list.

\subsection{The conditions on $\SL$ representations }

Some of the conditions on $\SL$ representations are obtained from the necessary
conditions on modular data Propositions \ref{p:MDcond} and \ref{p:MD1}, and
others are discussed in the main text of this paper.  Let us first translate
the conditions on the $(S,T)$ matrices to condition on an MD representations
$\rho_\a$:
\begin{prop}
\label{p:MDcond}
Given a modular data $S,T$ of rank $r$, let $\rho_\a$ be any one of its 12 MD
representations. Then $\rho_\a$ has the following properties:
\begin{enumerate}
\item $\rho_\a$ is an $\SL$ representation of level $\ord(\rho_\a(\ft))$, and $\ord(T) \mid \ord(\rho_\a(\ft)) \mid 12 \ord(T)$\,.

\item The conductor of the elements of $\rho_\a(\fs)$ 
divides $\ord(\rho_\a(\ft))$.

\item 
If $\rho_\a$ is a direct sum of two $\SL$ representations
\begin{align}
 \rho_\a \cong \rho \oplus \rho',
\end{align}
then the eigenvalues
of $\rho(\ft)$ and $\rho'(\ft)$ must overlap.  This implies that if $\rho_\a =
\rho \oplus \chi_1 \oplus \dots \oplus \chi_\ell$ for some 1-dimensional
representations $\chi_1, \dots, \chi_\ell$,  then $\chi_1, \cdots \chi_\ell$
are the same 1-dimensional representation.

\item
Suppose that $\rho_\a \cong \rho \oplus \ell \chi$ for an irreducible
representation $\rho$ with non-degenerate $\rho(\ft)$, and an 1-dimensional
representation $\chi$.  If $\ell \ne 2\dim(\rho) -1$ or $\ell > 1$, then
$(\rho(\fs) \chi(\fs)^{-1})^2 = \id$. 

\item $\rho_\a$ satisfies
\begin{align}
\rho_\a \not\cong n \rho
\end{align}
for any integer $n >1$ and any  representation $\rho$ such that
$\rho(\ft)$ is non-degenerate.

\item
If $\rho_\a(s)^2=\pm \id$ ({\it i.e.} if the modular data or MTC is self dual),
$\pord(\rho_\a(\ft)) $ is a prime and satisfies $\pord(\rho_\a(\ft)) = 1$ mod
4, then the representation $\rho_\a$ cannot be a direct sum of a
$d$-dimensional irreducible $\SL$ representation and two or more 1-dimensional
$\SL$ representations with $d=(p+1)/2$.  

\item Let $3< p < q$ be prime such that $pq \equiv 3 \mod 4$ and
$\pord(\rho_\a(\ft))=pq$, then the rank $r \ne \frac{p+q}{2}+1$.  Moreover, if
$p > 5$, rank $r > \frac{p+q}{2}+1$.

\item  The number of self dual objects is greater than 0. Thus
\begin{align}
 \Tr(\rho_\a(\fs)^2) \neq 0 .
\end{align}
Since $\Tr(\rho_\a(\fs)^2) \neq 0$, let us introduce
\begin{align}
 C = \frac{\Tr(\rho_\a(\fs)^2)}{|\Tr(\rho_\a(\fs)^2)|}   \rho_\a(\fs)^2.
\end{align}
The above $C$ is the charge conjugation operator of MTC,
{\it i.e.} $C$ is a permutation matrix of order 2.  In particular,
$\Tr(C)$ is the number of  self dual objects. Also, for each eigenvalue $\tilde
\theta$ of $\rho_\a(\ft)$, 
\begin{align}
\Tr_{\tilde \theta}(C) \geq 0,
\end{align}
where $\Tr_{\tilde \theta}$ is the trace in the degenerate
subspace of $\rho_\a(\ft)$ with eigenvalue $\tilde \theta$.

\item 
For any Galois conjugation $\s$ in $\Gal(\BQ_{\ord(\rho_\a(\ft))})$, there is a
permutation of the indices, $i \to \hs(i)$, and $\e_\s(i)\in \{1,-1\}$, such
that
\begin{align}
\label{Galact}
\s \big(\rho_\a(\fs)_{i,j}\big) &
= \e_\s(i)\rho_\a(\fs)_{\hat \s (i),j} 
= \rho_\a(\fs)_{i,\hat \s (j)}\e_\s(j) 
\\
\s^2 \big(\rho_\a(\ft)_{i,i}\big) &= \rho_\a(\ft)_{\hat \s (i),\hat \s (i)},
\end{align}
for all $i,j$. 

\item 
By \cite[Theorem II]{DLN}, $D_{\rho_\a}(\s)$ defined in \eqref{eq:Drho} must be a
\textbf{signed permutation}
\begin{align*}
 (D_{\rho_\a}(\s))_{i,j} = \e_\s(i) \delta_{\hs(i),j}.
\end{align*}
and satisfies
\begin{align}
\label{siDrho}
\s(\rho_\a(\fs)) &= D_{\rho_\a}(\s) \rho_\a(\fs) =\rho_\a(\fs)D_{\rho_\a}^\top(\s),
\nonumber\\
\s^2(\rho_\a(\ft)) &= D_{\rho_\a}(\s) \rho_\a(\ft) D_{\rho_\a}^\top(\s)
\end{align}

\item 
There exists a $u$ such that $\rho_\a(\fs)_{uu} \neq 0$ and
\begin{align}
\label{SSN}
&  \rho_\a(\fs)_{ui} \neq 0 \in \BR, \ \ \ \
 \frac{ \rho_\a(\fs)_{ij} }{ \rho_\a(\fs)_{uu} },\ 
 \frac{ \rho_\a(\fs)_{ij} }{ \rho_\a(\fs)_{uj} }
\in \BO_{\ord(T)}, \ \ \ \
 \frac{ \rho_\a(\fs)_{ij} }{ \rho_\a(\fs)_{i'j'} }
\in \BQ_{\ord(T)},
\nonumber\\
& N^{ij}_k = \sum_{l=0}^{r-1} \frac{
\rho_\a(\fs)_{li} \rho_\a(\fs)_{lj} \rho_\a(\fs^{-1})_{lk}}{ \rho_\a(\fs)_{lu} } \in\BN. 
\nonumber\\
& \forall \ i,j,k = 0,1,\ldots,r-1.
\end{align}
($u$ corresponds the unit object of MTC).

\item
Let $n \in \BN_+$.  The $n^\text{th}$ Frobenius-Schur indicator of the $i$-th
simple object
\begin{align}
 \label{nunFSapp}
 \nu_n(i)&= 
\sum_{j, k=0}^{r-1} N_i^{jk} \rho_\a(\fs)_{ju}\theta_j^n [\rho_\a(\fs)_{ku}\theta_k^n]^*
=\sum_{j, k=0}^{r-1} N_i^{jk} 
\rho_\a(\ft^n\fs)_{ju} \rho_\a(\ft^{-n}\fs^{-1})_{ku}
\nonumber\\
&=\sum_{j,k,l=0}^{r-1} \frac{
\rho_\a(\fs)_{lj} \rho_\a(\fs)_{lk} \rho^*_\a(\fs)_{li}}{ \rho_\a(\fs)_{lu} }
\rho_\a(\ft^n\fs)_{ju} \rho_\a(\ft^{-n}\fs^{-1})_{ku}
\nonumber\\
&=\sum_{l=0}^{r-1} \frac{
\rho_\a(\fs\ft^n\fs)_{lu} \rho_\a(\fs\ft^{-n}\fs^{-1})_{lu} \rho_\a(\fs^{-1})_{li}}{ \rho_\a(\fs)_{lu} }
\end{align}

 is a cyclotomic integer whose conductor divides $n$ and $\ord(T)$.  The 1st
Frobenius-Schur indicator satisfies $\nu_1(i)=\delta_{iu}$ while the 2nd
Frobenius-Schur indicator $\nu_2(i)$ satisfies 
$\nu_2(i)=\pm \rho_\a(\fs^2)_{ii}$
(see \cite{Bantay, NS07b, RSW0777}).
%$\nu_2(i)=0$ if $i\neq \bar i$, and $\nu_2(i)=\pm 1$ if $i = \bar i$ 

\item If we further assume the modular data or the MTC to be
non-integral, then $\pord(\tilde\rho_\a(\ft)) = \ord(T) \notin \{2,3,4,6\}$.
In particular, $\ord(\rho_\a(\ft)) \notin \{2,3,4,6\}$.

\end{enumerate}
\end{prop}

In Section \ref{SL2Zrep} and Appendix \ref{repPP}, we have explicitly
constructed all irreducible unitary representations of $\SL$ (up to unitary
equivalence).  However, this only gives the $\SL$ representations in some
arbitrary basis, not in the basis yielding MD representations (\emph{i.e.}
satisfying \eqref{STrho1}).   We can improve the situation by choosing a basis to make
$\rho(\ft)$ diagonal and $\rho(\fs)$ symmetric.  Since we are going to use
several types of basews, let us define these choices:
\begin{defn}
An unitary $\SL$ representations $\tilde\rho$ is called a \textbf{general}
$\SL$ matrix representations if $\tilde\rho(\ft)$ is diagonal \footnote{We
will consider only $\SL$ matrix representations with diagonal $\tilde\rho(\ft)$
in this paper.}.  A general $\SL$ matrix representation $\tilde \rho$ is called
\textbf{symmetric} if $\tilde\rho(\fs)$ is symmetric.  
%A general $\SL$ representation $\tilde \rho$ of a modular data is called
%\textbf{pseodu-MD} if $\tilde\rho(\fs)$ is symmetric and $U$ has a form:
%\begin{align}
% U = V_\mathrm{sd} P V'_\mathrm{sd}
%\end{align}  
%where $P$ is a permutation matrix and $ V_\mathrm{sd},  V_\mathrm{sd}'$ are
%diagonal matrices with $\pm 1$ diagonal entries.  
An general $\SL$ matrix representation $\tilde\rho$ is called
\textbf{irrep-sum} if $\tilde\rho(\fs),\tilde\rho(\ft)$ are matrix-direct sum
of irreducible $\SL$ representations.  An $\SL$ matrix representations
$\tilde\rho$ is called an $\SL$ representation \textbf{of modular data} $S,T$,
if $\tilde\rho$ is unitary equivalent to an MD
representation of the modular data, {\it i.e.},
\begin{align}
\label{STrhoU}
\tilde\rho(\fs) = \ee^{-2\pi \ii \frac{\a}{4}}\frac{1}{D}\, US U^\dag, \ \ \ \
\tilde\rho(\ft) = UTU^\dag \ee^{2\pi \ii (\frac{-c}{24} + \frac{\a}{12})},
\end{align}
for some unitary matrix $U$ and $\a \in \BZ_{12}$, where $c$ is the central
charge.\footnote{Note that $D^2$ is always positive and $D$ in \eqref{STrhoU} is the
positive root of $D^2$, even for non-unitary cases.}
\end{defn}

Through our explicit construction, we observe that all irreducible unitary
representations of $\SL$ are unitary equivalent to symmetric matrix
representations of $\SL$, at least for dimension equal or less than 12. 
%In Appendix \ref{repPP}, we list all the irreducible symmetric matrix
%representations of $\SL$ up to dimension 6.

We note that different choices of orthogonal basis give rise to different
matrix representations of $\SL$. The modular data $S,T$ is obtained from some
particular choices of the basis.  Some properties on the MD representations of
a modular data do not depend on the choices of basis in the eigenspaces of
$\tilde\rho(\ft)$ (induced by the block-diagonal unitary transformation $U$ in
\eqref{STrhoU} that leaves $\tilde\rho(\ft)$ invariant). Those properties remain
valid for any general $\SL$ representations $\tilde\rho$ of the modular
data.  In the following, we collect the basis-independent conditions on the
$\SL$ matrix representations of modular data. Those conditions have been
discussed in the main text.
\begin{prop}
\label{p:gencond1}
Let $\tilde\rho$ be a general $\SL$ matrix representations of a modular data or
a MTC.  Then $\tilde\rho$ must satisfy the following conditions:
\begin{enumerate}
\item 
If $\tilde\rho$ is a direct sum of two $\SL$ representations
\begin{align}
 \tilde\rho \cong \rho \oplus \rho',
\end{align}
%such that $\rho(\ft)$ and $\rho'(\ft)$ are diagonal, 
then the diagonals entries
of $\rho(\ft)$ and $\rho'(\ft)$ must overlap.  

\item
Suppose that $\tilde\rho \cong \rho \oplus \ell \chi$ for an irreducible
representation $\rho$ with $\rho(\ft)$ non-degenerate, and a character $\chi$.  If $\ell \neq 1$ and $\ell \ne 2\dim(\rho) -1$, then
$(\rho(\fs) \chi(\fs)^{-1})^2 = \id$. 

\item
If $\tilde\rho(\fs)^2=\pm \id$, and $\pord(\tilde\rho(\ft)) = 1$ mod 4 is a
prime, then the representation $\td\rho$ cannot be a direct sum of a
$d$-dimensional irreducible $\SL$ representation and two or more 1-dimensional
$\SL$ representations with $d=(\pord(\tilde\rho(\ft)) +1)/2$.  

\item $\tilde\rho$ satisfies
\begin{align}
\tilde\rho \not\cong n \rho
\end{align}
for any integer $n >1$ and any  representation $\rho$ such that $\rho(\ft)$ is non-degenerate.

\item
Let $3< p < q$ be prime such that $pq \equiv 3 \mod 4$ and
$\pord(\rho(\ft))=pq$, then the rank $r \ne \frac{p+q}{2}+1$.  Moreover, if
$p > 5$, rank $r > \frac{p+q}{2}+1$.

\item If we further assume 
$D^2$ of the modular data or 
the MTC to be
non-integral, then $\pord(\tilde\rho(\ft)) = \ord(T) \notin \{2,3,4,6\}$.  This
implies that $\ord(\tilde\rho(\ft)) \notin \{2,3,4,6\}$.

\end{enumerate}
\end{prop}

%we note that
%\begin{prop}
%Let $\tilde\rho$ be a general $\SL$ representations of a modular data or a MTC.
%Then there exists a block-diagonal unitary matrix $U$ that connect  $\tilde\rho$
%to a modular data, up to tensoring an 1-dimensional representation.  The
%block-diagonal $U$ acts within the eigenspace of $\tilde\rho(\ft)$.
%\end{prop}

Some other properties of an MD representation do depend on the choice of basis.
To make use of those properties, we can construct some combinations of
$\tilde\rho(\fs)$'s that are invariant under the block-diagonal unitary
transformation $U$.  

The eigenvalues of $\tilde\rho(\ft)$ partition the indices of the basis vectors. To construct the invariant combinations of $\tilde\rho(\fs)$, for any eigenvalue $\tilde\theta$ of $\tilde\rho(\ft)$, let
\begin{align}
\label{Ith}
 I_{\tilde\theta} = \{ i \, \big|\, \tilde\rho(\ft)_{ii}=\tilde\theta\}. 
 \end{align} 
Let $I = I_{\tilde\theta}$, $J = J_{\tilde\theta'}$, $K = K_{\tilde\theta''}$ for some eigenvalues $\tilde\theta$, $\tilde\theta'$, $\tilde\theta''$ of $\tilde\rho(\ft)$.
We see that the following uniform polynomials of $\tilde\rho(\fs)$ are
invariant 
\begin{align}
\label{invcomb}
P_I(\rho(\fs)) =
\Tr \tilde\rho(\fs)_{II} &\equiv \sum_{i\in I} \tilde\rho(\fs)_{ii},
\nonumber\\
P_{IJ}(\rho(\fs)) =
\Tr \tilde\rho(\fs)_{IJ} \tilde\rho(\fs)_{JI} &\equiv 
\sum_{i\in I,j\in J} \tilde\rho(\fs)_{i,j} \tilde\rho(\fs)_{ji},
\\
P_{IJK}(\rho(\fs)) =
\Tr \tilde\rho(\fs)_{IJ} \tilde\rho(\fs)_{JK} \tilde\rho(\fs)_{KI} 
&\equiv \sum_{i\in I,j\in J,k\in K} 
\tilde\rho(\fs)_{i,j} 
\tilde\rho(\fs)_{j,k}
\tilde\rho(\fs)_{k,i}
.
\nonumber 
\end{align}
Certainly we can construction many other invariant uniform polynomials in the
similar way.  Using those invariant uniform polynomials, we have the following
results
\begin{prop}
\label{p:gencond2}
Let $\tilde\rho$ be a general $\SL$ representations of a modular data or a
MTC.  Then following statements hold:
\begin{enumerate}

\item $\tilde\rho(\fs)$ satisfies
\begin{align}
 \Tr(\tilde\rho(\fs)^2)  \in \BZ \setminus\{0\} .
\end{align}
Let 
\begin{align}
\label{Crho}
 C = \frac{\Tr(\tilde\rho(\fs)^2)}{|\Tr(\tilde\rho(\fs)^2)|}   \tilde\rho(\fs)^2.
\end{align}
For all $I$, 
\begin{align}
\label{PIC}
 P_I(C) \geq 0.
\end{align}

\item 
\label{cndI}
The conductor of $P_\mathrm{odd}(\tilde\rho(\fs))$ divides
$\ord(\tilde\rho(\ft))$ for all the invariant uniform polynomials
$P_\mathrm{odd}$ with odd powers of $\tilde\rho(\fs)$ (such as $P_I$ and
$P_{IJK}$ in \eqref{invcomb}).  The conductor of
$P_\mathrm{even}(\tilde\rho(\fs))$ divides $\pord(\tilde\rho(\ft))$ for all the
invariant uniform polynomials $P_\mathrm{even}$ with even powers of
$\tilde\rho(\fs)$ (such as $P_{IJ}$ in \eqref{invcomb}).

\item 
 For any Galois conjugation $\s\in\Gal(\BQ_{\ord(\rho(\ft))})$, there is a
permutation on the set $\{I\}$, $I \to \hs(I)$, such that
\begin{align}
\s P_{IJ}(\tilde\rho(\fs)) &
= P_{I\hs(J)}(\tilde\rho(\fs))
=P_{\hs(I)J}(\tilde\rho(\fs))
\nonumber\\
\s^2 \big(\tilde\theta_I\big) &= \tilde\theta_{\hat \s (I)},
\end{align}
for all $I,J$.

\item 
For any invariant uniform polynomials $P$ (such as those in \eqref{invcomb})
\begin{align}
\s P\big(\tilde\rho(\fs)\big) = 
 P\big(\s \tilde\rho(\fs)\big) = 
P \big(\tilde\rho(\ft)^a \tilde\rho(\fs) \tilde\rho(\ft)^b \tilde\rho(\fs) \tilde\rho(\ft)^a\big)
\end{align}
where $\s\in\Gal(\BQ_{\ord(\td\rho(\ft))})$,
and $a,b$ are given by $\s( \ee^{\ii 2
\pi/\ord(\td\rho(\ft))} ) = \ee^{a \ii 2 \pi/\ord(\td\rho(\ft))}$ and $ab \equiv 1$ mod
$\ord(\td\rho(\ft))$.  
\end{enumerate}
\end{prop}

%Let $\tilde \rho$ be a symmetric $\SL$ representation of a modular data. It is
%immediate to see that the tensor product $\tilde \rho_\a$ of any 1-dimensional
%$\SL$ representation $\chi_\a$ with $\tilde \rho$ is again a symmetric $\SL$
%representation of the same modular data. 
%For any general $\SL$ representation of a modular data, we will simply call
%the eigenvectors of $\tilde\rho(\ft)$ \emph{non-degenerate} if the
%corresponding eigenvalues are non-degenerate.  The subspace spanned by the
%\emph{non-degenerate} eigenvectors of $\tilde \rho(\ft)$ is the
%\emph{non-degenerate subspace} of the representation $\tilde\rho$.  Obviously,
%the non-degenerate eigenvectors of $\tilde \rho(\ft)$ form a canonical basis
%of non-degenerate subspace of $\tilde\rho$, up to $\pm$ signs.
%we like to choose a more special basis in the eigenspace of
%$\tilde(\rho(\ft)$, so that the basis is closer to the basis that lead to the
%MD representation.  For example, we can choose a basis to make
%$\tilde\rho(\fs)$ symmetric ({\it i.e.} to make $\tilde\rho$ a symmetric
%representation).

Instead of constructing invariants, there is another way to make use of the
properties of an MD representation that depend on the choices of basis.  We can
choose a more special basis, so that the basis is closer to the basis that
leads to the MD representation.  For example, we can choose a basis to make
$\tilde\rho(\fs)$ symmetric ({\it i.e.} to make $\tilde\rho$ a symmetric
representation).

Now consider a symmetric $\SL$ matrix representation $\tilde \rho$ of a modular
data or a MTC.  We find that the restriction of the unitary $U$ in \eqref{STrhoU} on the non-degenerate subspace (see 
 Theorem \ref{t:ortho_eqv}) must be diagonal with diagonal elements $U_{ii}\in \{
1,-1\}$.  Therefore, within the non-degenerate subspace,
$\tilde\rho(\fs)$ of a symmetric representation differs from $\rho(\fs)$ of an MD
representation  only by a diagonal unitary transformation $U$ with diagonal
elements $\pm 1$, i.e., a \textbf{signed diagonal}
matrix.  In this case some properties of MD representation apply to the blocks of the symmetric
representation within the non-degenerate subspace.  This allows us to obtain
\begin{prop}
\label{p:gal_sym}
Let $\tilde\rho$ be a symmetric $\SL$ representations equivalent to an MD representation.
 Let 
\begin{align}
I_\mathrm{ndeg} &:=\{i \mid \tilde\rho(\ft)_{i,i} \text{ is a
non-degenerate eigenvalue}\},
\end{align}
Then there exists an orthogonal $U$ such that $U \tilde \rho U^\top$ is a pMD representation, and the following statements hold:

\begin{enumerate}
\item 
\label{cndndegI}
The conductor of $(U \tilde\rho(\fs) U^\top)_{i,j} $ divides
$\ord(\tilde\rho(\ft))$ for all $i,j$.  This implies that the conductor of
$(\tilde\rho(\fs))_{i,j} $ divides $\ord(\tilde\rho(\ft))$ for all $i,j \in
I_\mathrm{ndeg}$.

\item 
 For any Galois conjugation $\s$ in $\Gal(\BQ_{\ord(\td\rho(\ft))})$, there is a
permutation $i \to \hs(i)$, such that
\begin{align}
\s \big((U \tilde\rho(\fs)U^\top)_{i,j}\big) &
= \e_\s(i) (U \tilde\rho(\fs) U^\top)_{\hat \s (i),j} 
= (U \tilde\rho(\fs) U^\top)_{i,\hat \s (j)} \e_\s(j)
\nonumber\\
\s^2 \big( \tilde\rho(\ft)_{i,i}\big) &= 
 \tilde\rho(\ft)_{\hat \s (i),\hat \s (i)},
\end{align}
for all $i,j$, where $\e_\s(i)\in \{1,-1\}$. 
This implies that
\begin{align}
\s \big( \tilde\rho(\fs)_{i,j}\big) & = \tilde\rho(\fs)_{\hat \s (i),j} 
\ \ \mathrm{or}  \ \
\s \big( \tilde\rho(\fs)_{i,j}\big)  = -  \tilde\rho(\fs)_{\hat \s (i),j} 
\nonumber\\
\s \big( \tilde\rho(\fs)_{i,j}\big) & =  \tilde\rho(\fs) _{i,\hat \s (j)} 
\ \ \mathrm{or}  \ \
\s \big( \tilde\rho(\fs)_{i,j}\big)  = - \tilde\rho(\fs) _{i,\hat \s (j)} 
\end{align}
for all $i,j \in I_\mathrm{ndeg}$.  This also implies that $D_{\td\rho}(\s)$
defined in \eqref{eq:Drho} is a signed permutation matrix in the $I_\mathrm{ndeg}$
block, {\it i.e.} $(D_{\td\rho}(\s))_{i,j}$ for $i,j \in I_\mathrm{ndeg}$ are
matrix elements of a signed permutation matrix.

\item 
\label{sirhondegI}
For all $i,j$, 
\begin{align}
\s\big((U\tilde\rho(\fs) U^\top)_{i,j}\big) = 
\big(U \tilde\rho(\ft)^a \tilde\rho(\fs) \tilde\rho(\ft)^b \tilde\rho(\fs) \tilde\rho(\ft)^a U^\top\big)_{i,j}
\end{align}
where $\s\in\Gal(\BQ_{\ord(\td\rho(\ft))})$, and
$a,b$ are given by $\s( \ee^{\ii 2 \pi/\ord(\td\rho(\ft))} ) = \ee^{a \ii 2
\pi/\ord(\td\rho(\ft))}$ and $ab \equiv 1$ mod $\ord(\td\rho(\ft))$.  
This implies that
\begin{align}
\label{sirhondeg}
\s\big((\tilde\rho(\fs) )_{i,j}\big) = 
\big( \tilde\rho(\ft)^a \tilde\rho(\fs) \tilde\rho(\ft)^b \tilde\rho(\fs) \tilde\rho(\ft)^a \big)_{i,j}.
\end{align}
for all $i,j \in I_\mathrm{ndeg}$.

\item
\label{nonzeroI} Both $T$ and $\td\rho(\ft)$ are diagonal, and without loss of
generality, we may assume $\td\rho(\ft)$ is a scalar multiple of $T$. In this
case $U$ in \eqref{STrhoU} is a block diagonal matrix preserving the
eigenspaces of $\td\rho(\ft)$.  Let $I_\mathrm{nonzero} = \{i\}$ be a set of
indices such that the $i^\mathrm{th}$ row of  $U \td\rho(\fs) U^\top$ contains
no zeros for some othorgonal $U$ satisfying $U\td\rho(\ft)U^\top =
\td\rho(\ft)$.  The index for the unit object of MTC must be in
$I_\mathrm{nonzero}$.  Thus $I_\mathrm{nonzero}$ must be nonempty:
\begin{align}
\label{Ineq0}
 I_\mathrm{nonzero} \neq \emptyset.
\end{align}

\item 
Let $I_{\td\theta}$ be a set of indices for an eigenspace $E_{\td\theta}$
of $\tilde\rho(\ft)$
\begin{align}
I_{\td\theta} &:=\{i \mid \tilde\rho(\ft)_{i,i} = \td\theta \}.
\end{align}
Then there exists a $I_{\td\theta}$ such that
\begin{align}
\label{TrCunit}
I_{\td\theta} \cap I_\mathrm{nonzero} \neq \emptyset\ \text{ and }\
\Tr_{E_{\td\theta}} C > 0,
\end{align}
where $C$ is given in \eqref{Crho}.

%\item There exists a $u$ in $I_\mathrm{nonzero}$, such that
%\begin{align}
%\frac{ \rho_\a(\fs)_{u,j} }{ \rho_\a(\fs)_{u,i} } \in \BR, \ \ \
%\forall \ i,j \in I_\mathrm{ndeg}
%.
%\end{align}

\item 
If we further assume the modular data to be non-integral, then there
exists a $I_{\td\theta}$ that has a non-empty overlap with $I_\mathrm{nonzero}$,
such that $D_{\td\rho}(\s)_{I_{\td\theta}} \neq \pm \id$ for some $\s \in
\Gal(\BQ_{\ord(\td\rho(\ft))}/\BQ)$.  Here $D_{\td\rho}(\s)$ is defined in
 \eqref{eq:Drho}:
\begin{align}
\label{Drho}
 D_{\td\rho}(\s)=
\td\rho(\ft)^a \td\rho(\fs) \td\rho(\ft)^b \td\rho(\fs) \td\rho(\ft)^a \td\rho^{-1}(\fs)
\end{align}
where $a,b$ are given by $\s( \ee^{ 2 \pi\ii/\ord(\td\rho(\ft))} )=\ee^{ a 2
\pi\ii/\ord(\td\rho(\ft))}$ and $ab \equiv 1 \mod \ord(\td\rho(\ft))$.  
%Note that $D_{\td\rho}(\s)$ is a signed permutation matrix only within the
%non-degenerate block.  
Also $D_{\td\rho}(\s)_{I_{\td\theta}}$ is the block of $D_{\td\rho}(\s)$ with
indices in $I_{\td\theta}$, {\it i.e.} the matrix elements of
$D_{\td\rho}(\s)_{I_{\td\theta}}$ are given by $(D_{\td\rho}(\s))_{i,j},\ i,j
\in I_{\td\theta}$.

\end{enumerate}
\end{prop}

Proposition \ref{p:gal_sym}(6) is a consequence of Theorem \ref{p:5.2}(3).
Using GAP System for Computational Discrete Algebra, we obtain a list of symmetric irrep-sum $\SL$ matrix
representations that satisfy the conditions in Propositions \ref{p:gencond1}, \ref{p:gencond2}, and
\ref{p:gal_sym}.  The list is given below for rank $r=6$ case (see Appendix \ref{repL}).   The lists for
$r=2,3,\ldots,5$ are also obtained, and are given in the Supplementary Material
Section \ref{Section5} -- \ref{Section8}.  

Some of those symmetric irrep-sum $\SL$
matrix representations are representations of modular data, while others are
not.  However, the list includes all the symmetric irrep-sum $\SL$ matrix
representations of modular data or MTC's which are not integral. 

In Supplementary Material Section \ref{Section3}, we give a list of all rank-6
symmetric irrep-sum representations, and indicate the passing and failing
representations, as well as which conditions of Propositions \ref{p:gencond1},
\ref{p:gencond2}, and \ref{p:gal_sym} that a failing representation fail to
satisfy.

\subsection{List of symmetric irrep-sum representations}
\label{repL}

The following is a list the all rank-6 symmetric irrep-sum representations that
satisfy the conditions in Propositions \ref{p:gencond1}, \ref{p:gencond2}, and
\ref{p:gal_sym}.  The list contains all the rank-6 symmetric irrep-sum
representations that are unitarily equivalent to rank-6 MD representations, plus
some extra ones.

For each symmetric irrep-sum representation, we may generate an orbit by
orthogonal transformations
\begin{align}
 \rho_\text{isum}(\fs) \to U \rho_\text{isum}(\fs) U^\top,\ \ \ \rho_\text{isum}(\ft) \to U \rho_\text{isum}(\ft) U^\top,
\end{align}
tensoring 1-dimensional $\SL$ representations $\chi_\a$, $\a =1,\ldots,12$:
\begin{align}
\rho_\text{isum}(\fs) \to \chi_\a(\fs)\rho_\text{isum}(\fs),\ \ \ \rho_\text{isum}(\ft) \to \chi_\a(\ft)\rho_\text{isum}(\ft) ,
\end{align}
and applying Galois conjugations $\s$ in
$\Gal(\BQ_{\ord(\rho_\text{isum}(\ft))})$:
\begin{align}
\rho_\text{isum}(\fs) \to \s(\rho_\text{isum}(\fs)),\ \ \ \rho_\text{isum}(\ft) \to \s(\rho_\text{isum}(\ft)) .
\end{align}
We will call such an orbit a \textbf{GT orbit}.  The following list includes
only one representative for each GT orbit.  The list can also be regarded as a
list GT orbits.

%For example, from Appendix \ref{repPP}, we see that all the dimension-6
%irreducible representations form 16 GT orbits.  The $2^\mathrm{th}$,
%$6^\mathrm{th}$ GT orbits are rejected since they fail to satisfy Proposition
%\ref{p:gencond2} (\ref{cndI}) ({\it i.e.} the conductor of $
%\rho_\text{isum}(\fs)_{i,i}$ does not divides $\ord(\rho_\text{isum}(\ft))$).  The
%$1^\mathrm{th}$, $3^\mathrm{th}$, $7^\mathrm{th}$, $9^\mathrm{th}$ GT orbits
%are rejected since they fail to satisfy Proposition \ref{p:gal_sym}
%(\ref{cndndegI}) ({\it i.e.} the conductor of $ \rho_\text{isum}(\fs)_{i,j}$ does not
%divides $\ord(\rho_\text{isum}(\ft))$).  The $4^\mathrm{th}$, $12^\mathrm{th}$,
%$13^\mathrm{th}$, $14^\mathrm{th}$ GT orbits are rejected since they fail to
%satisfy Proposition \ref{p:gal_sym} (\ref{nonzeroI}) ({\it i.e.} 
%all the rows of $U
%\rho_\text{isum}(\fs)U^\dag$ contain zeros).  The
%$5^\mathrm{th}$, $15^\mathrm{th}$, $16^\mathrm{th}$ GT orbits are rejected
%since they fail to satisfy \eqref{sirhondeg} in Proposition \ref{p:gal_sym}
%(\ref{sirhondegI}) ({\it i.e.} $\s_a\big((\rho_\text{isum}(\fs) )\big) \neq \big(
%\rho_\text{isum}(\ft)^a \rho_\text{isum}(\fs) \rho_\text{isum}(\ft)^b \rho_\text{isum}(\fs)
%\rho_\text{isum}(\ft)^a \big)$ in non-degenerate block $I_\mathrm{ndeg}$).  The
%$8^\mathrm{th}$, $10^\mathrm{th}$, $11^\mathrm{th}$ GT orbits satisfy all the
%conditions in Propositions \ref{p:gencond1}, \ref{p:gencond2}, and
%\ref{p:gal_sym}. Thus they appear in the following list and are indicated by
%ind = $6_8^{0,1}$, $6_{10}^{0,1}$, and $6_{11}^{0,1}$.

In the list, a representation $\rho_\text{isum}$ is expressed as the direct sum of
irreducible representations $\rho_\text{isum} = \rho_1\oplus \rho_2 \oplus \cdots$,
where $\rho_a(\mathfrak{t})$ is presented as $(\tilde s_{1},\tilde
s_{2},\cdots)$ with $\tilde s_{i} = \arg(\rho_a(\mathfrak{t})_{ii})$, 
and $\rho_a(\mathfrak{s})$ is presented as $(\rho_{11},
\rho_{12}, \rho_{13}, \rho_{14}, \cdots;\ \  \rho_{22},
\rho_{23}, \rho_{24}, \cdots)$.  The direct sum is also given via
an index form, for example, irreps = $2_{2}^{1,0} \hskip -1.5pt \otimes \hskip
-1.5pt 2_{5}^{1,0}\oplus 2_{5}^{1,0}$. It means that the representation
$\rho_\text{isum}$ is a direct sum of two irreducible representations $2_{2}^{1,0}
\hskip -1.5pt \otimes \hskip -1.5pt 2_{5}^{1,0}$ and $2_{5}^{1,0}$.  Here
$2_{2}^{1,0}$, $2_{5}^{1,0}$ are indices of $\SL$ irreducible representations
with prime-power levels.  Those prime-power-level $\SL$ irreducible
representation are listed in Appendix \ref{repPP}, where the meaning of the
indices is explained further.  $2_{2}^{1,0} \hskip -1.5pt \otimes \hskip -1.5pt
2_{5}^{1,0}$ is the irreducible representation obtained by the tensor product
of $2_{2}^{1,0}$ and $2_{5}^{1,0}$.

The dimensions of the representations $\rho_\text{isum}$ are given by dims =
$(r_1, r_2 , \cdots)$, where $r_a$ is the dimension of the irreducible
representation $\rho_a$, satisfying $r_1\geq r_2 \geq \cdots$.  The levels  of
the representations $\rho_a$ are given by levels = $(l_1, l_2 , \cdots )$,
where $l_a =\mathrm{ord}(\rho_a(\ft))$.  We will use (dims;levels) = $(r_1, r_2
, \cdots; l_1,l_2,\cdots)$ to label those representations.  Now we can explain
how the representative of a GT orbit is chosen.  The representative for a GT
orbit is chosen to be the one with minimal $[ [r_1,r_2,\cdots],
\ord(\rho_\text{isum}(\ft)), [l_1,l_2,\cdots] ]$.  Here the order of two lists
is determined by first compare the first elements of the two lists.  If the
first elements are equal, we then compare the second elements, {\it etc}.
The order of cyclotomic numbers are given by GAP.

To describe the entries of $\rho_a(\fs)$, we also introduced the following
notations: 
\begin{align}
\zeta^m_n &=\mathrm{e}^{2\pi \mathrm{i} m/n}
,\ \ \
c^m_n  = \zeta^m_n+\zeta^{-m}_n
,\ \ \
s^m_n  = \zeta^m_n-\zeta^{-m}_n, 
\nonumber\\
\xi^{m,k}_n  &= (\zeta^m_{2n}-\zeta^{-m}_{2n})/(\zeta_{2n}^k-\zeta_{2n}^{-k}),\ \ \ \ 
\xi^{m}_n = \xi^{m,1}_n
.
\end{align}
%We write $\xi^{m,l}_{2n}$ as $\xi^{m,l}_{n\times 2}$ if it actually has a
%conductor $n$.

We find that, for rank 6, there are only 25 GT orbits.  The GT orbits can be
divided into two classes, resolved and unresolved, whose definition will to
given in the next section.  Below each GT orbit, we indicate whether it is
resolved or unresolved.  Among  25 GT orbits, 17 are resolved and 8 are
unresolved.

For the 17 resolved GT orbits, it is easy to compute all the corresponding
pairs of $(S,T)$ matrices that satisfied the conditions in Proposition
\ref{p:MDcond}, which will be done in next section.  Below each
resolved GT orbit, we indicate the number valid $(S,T)$ pairs obtain with such
a computation.  Those valid $(S,T)$ pairs will be listed in Appendix \ref{SsL}.
The 8 unresolved GT orbits are difficult to handle by computer, which are
discussed in the main text. (The main text also discussed most of the resolved
cases.)

\

\noindent1. (dims;levels) =$(3,
2,
1;5,
5,
1
)$,
irreps = $3_{5}^{1}\oplus
2_{5}^{1}\oplus
1_{1}^{1}$,
pord$(\rho_\text{isum}(\mathfrak{t})) = 5$,

\vskip 0.7ex
\hangindent=4em \hangafter=1
 $\rho_\text{isum}(\mathfrak{t})$ =
 $( 0,
\frac{1}{5},
\frac{4}{5} )
\oplus
( \frac{1}{5},
\frac{4}{5} )
\oplus
( 0 )
$,

\vskip 0.7ex
\hangindent=4em \hangafter=1
 $\rho_\text{isum}(\mathfrak{s})$ =
($\sqrt{\frac{1}{5}}$,
$-\sqrt{\frac{2}{5}}$,
$-\sqrt{\frac{2}{5}}$;
$-\frac{5+\sqrt{5}}{10}$,
$\frac{5-\sqrt{5}}{10}$;
$-\frac{5+\sqrt{5}}{10}$)
 $\oplus$
$\mathrm{i}$($-\frac{1}{\sqrt{5}}c^{3}_{20}
$,
$\frac{1}{\sqrt{5}}c^{1}_{20}
$;\ \ 
$\frac{1}{\sqrt{5}}c^{3}_{20}
$)
 $\oplus$
($1$)

Resolved. Number of valid $(S,T)$ pairs = 0.

\vskip 2ex

 \noindent2. (dims;levels) =$(3,
2,
1;8,
8,
1
)$,
irreps = $3_{8}^{1,0}\oplus
2_{8}^{1,9}\oplus
1_{1}^{1}$,
pord$(\rho_\text{isum}(\mathfrak{t})) = 8$,

\vskip 0.7ex
\hangindent=4em \hangafter=1
 $\rho_\text{isum}(\mathfrak{t})$ =
 $( 0,
\frac{1}{8},
\frac{5}{8} )
\oplus
( \frac{1}{8},
\frac{7}{8} )
\oplus
( 0 )
$,

\vskip 0.7ex
\hangindent=4em \hangafter=1
 $\rho_\text{isum}(\mathfrak{s})$ =
$\mathrm{i}$($0$,
$\sqrt{\frac{1}{2}}$,
$\sqrt{\frac{1}{2}}$;\ \ 
$-\frac{1}{2}$,
$\frac{1}{2}$;\ \ 
$-\frac{1}{2}$)
 $\oplus$
$\mathrm{i}$($-\sqrt{\frac{1}{2}}$,
$\sqrt{\frac{1}{2}}$;\ \ 
$\sqrt{\frac{1}{2}}$)
 $\oplus$
($1$)

Resolved. Number of valid $(S,T)$ pairs = 0.

\vskip 2ex

 \noindent3. (dims;levels) =$(3,
2,
1;5,
2,
1
)$,
irreps = $3_{5}^{1}\oplus
2_{2}^{1,0}\oplus
1_{1}^{1}$,
pord$(\rho_\text{isum}(\mathfrak{t})) = 10$,

\vskip 0.7ex
\hangindent=4em \hangafter=1
 $\rho_\text{isum}(\mathfrak{t})$ =
 $( 0,
\frac{1}{5},
\frac{4}{5} )
\oplus
( 0,
\frac{1}{2} )
\oplus
( 0 )
$,

\vskip 0.7ex
\hangindent=4em \hangafter=1
 $\rho_\text{isum}(\mathfrak{s})$ =
($\sqrt{\frac{1}{5}}$,
$-\sqrt{\frac{2}{5}}$,
$-\sqrt{\frac{2}{5}}$;
$-\frac{5+\sqrt{5}}{10}$,
$\frac{5-\sqrt{5}}{10}$;
$-\frac{5+\sqrt{5}}{10}$)
 $\oplus$
($-\frac{1}{2}$,
$-\sqrt{\frac{3}{4}}$;
$\frac{1}{2}$)
 $\oplus$
($1$)

Unresolved. 

\vskip 2ex

\noindent4. (dims;levels) =$(3,
2,
1;5,
2,
2
)$,
irreps = $3_{5}^{1}\oplus
2_{2}^{1,0}\oplus
1_{2}^{1,0}$,
pord$(\rho_\text{isum}(\mathfrak{t})) = 10$,

\vskip 0.7ex
\hangindent=4em \hangafter=1
 $\rho_\text{isum}(\mathfrak{t})$ =
 $( 0,
\frac{1}{5},
\frac{4}{5} )
\oplus
( 0,
\frac{1}{2} )
\oplus
( \frac{1}{2} )
$,

\vskip 0.7ex
\hangindent=4em \hangafter=1
 $\rho_\text{isum}(\mathfrak{s})$ =
($\sqrt{\frac{1}{5}}$,
$-\sqrt{\frac{2}{5}}$,
$-\sqrt{\frac{2}{5}}$;
$-\frac{5+\sqrt{5}}{10}$,
$\frac{5-\sqrt{5}}{10}$;
$-\frac{5+\sqrt{5}}{10}$)
 $\oplus$
($-\frac{1}{2}$,
$-\sqrt{\frac{3}{4}}$;
$\frac{1}{2}$)
 $\oplus$
($-1$)

Unresolved. 

\vskip 2ex

\noindent5. (dims;levels) =$(3,
2,
1;4,
3,
2
)$,
irreps = $3_{4}^{1,3}\oplus
2_{3}^{1,0}\oplus
1_{2}^{1,0}$,
pord$(\rho_\text{isum}(\mathfrak{t})) = 12$,

\vskip 0.7ex
\hangindent=4em \hangafter=1
 $\rho_\text{isum}(\mathfrak{t})$ =
 $( 0,
\frac{1}{2},
\frac{1}{4} )
\oplus
( 0,
\frac{1}{3} )
\oplus
( \frac{1}{2} )
$,

\vskip 0.7ex
\hangindent=4em \hangafter=1
 $\rho_\text{isum}(\mathfrak{s})$ =
$\mathrm{i}$($-\frac{1}{2}$,
$\frac{1}{2}$,
$\sqrt{\frac{1}{2}}$;\ \ 
$-\frac{1}{2}$,
$\sqrt{\frac{1}{2}}$;\ \ 
$0$)
 $\oplus$
$\mathrm{i}$($-\sqrt{\frac{1}{3}}$,
$\sqrt{\frac{2}{3}}$;\ \ 
$\sqrt{\frac{1}{3}}$)
 $\oplus$
($-1$)

Resolved. Number of valid $(S,T)$ pairs = 0.

\vskip 2ex

 \noindent6. (dims;levels) =$(3,
2,
1;4,
3,
4
)$,
irreps = $3_{4}^{1,3}\oplus
2_{3}^{1,0}\oplus
1_{4}^{1,0}$,
pord$(\rho_\text{isum}(\mathfrak{t})) = 12$,

\vskip 0.7ex
\hangindent=4em \hangafter=1
 $\rho_\text{isum}(\mathfrak{t})$ =
 $( 0,
\frac{1}{2},
\frac{1}{4} )
\oplus
( 0,
\frac{1}{3} )
\oplus
( \frac{1}{4} )
$,

\vskip 0.7ex
\hangindent=4em \hangafter=1
 $\rho_\text{isum}(\mathfrak{s})$ =
$\mathrm{i}$($-\frac{1}{2}$,
$\frac{1}{2}$,
$\sqrt{\frac{1}{2}}$;\ \ 
$-\frac{1}{2}$,
$\sqrt{\frac{1}{2}}$;\ \ 
$0$)
 $\oplus$
$\mathrm{i}$($-\sqrt{\frac{1}{3}}$,
$\sqrt{\frac{2}{3}}$;\ \ 
$\sqrt{\frac{1}{3}}$)
 $\oplus$
$\mathrm{i}$($1$)

Unresolved. 

\vskip 2ex

\noindent7. (dims;levels) =$(3,
2,
1;8,
3,
1
)$,
irreps = $3_{8}^{1,0}\oplus
2_{3}^{1,0}\oplus
1_{1}^{1}$,
pord$(\rho_\text{isum}(\mathfrak{t})) = 24$,

\vskip 0.7ex
\hangindent=4em \hangafter=1
 $\rho_\text{isum}(\mathfrak{t})$ =
 $( 0,
\frac{1}{8},
\frac{5}{8} )
\oplus
( 0,
\frac{1}{3} )
\oplus
( 0 )
$,

\vskip 0.7ex
\hangindent=4em \hangafter=1
 $\rho_\text{isum}(\mathfrak{s})$ =
$\mathrm{i}$($0$,
$\sqrt{\frac{1}{2}}$,
$\sqrt{\frac{1}{2}}$;\ \ 
$-\frac{1}{2}$,
$\frac{1}{2}$;\ \ 
$-\frac{1}{2}$)
 $\oplus$
$\mathrm{i}$($-\sqrt{\frac{1}{3}}$,
$\sqrt{\frac{2}{3}}$;\ \ 
$\sqrt{\frac{1}{3}}$)
 $\oplus$
($1$)

Resolved. Number of valid $(S,T)$ pairs = 0.

\vskip 2ex

 \noindent8. (dims;levels) =$(3,
2,
1;8,
3,
3
)$,
irreps = $3_{8}^{1,0}\oplus
2_{3}^{1,0}\oplus
1_{3}^{1,0}$,
pord$(\rho_\text{isum}(\mathfrak{t})) = 24$,

\vskip 0.7ex
\hangindent=4em \hangafter=1
 $\rho_\text{isum}(\mathfrak{t})$ =
 $( 0,
\frac{1}{8},
\frac{5}{8} )
\oplus
( 0,
\frac{1}{3} )
\oplus
( \frac{1}{3} )
$,

\vskip 0.7ex
\hangindent=4em \hangafter=1
 $\rho_\text{isum}(\mathfrak{s})$ =
$\mathrm{i}$($0$,
$\sqrt{\frac{1}{2}}$,
$\sqrt{\frac{1}{2}}$;\ \ 
$-\frac{1}{2}$,
$\frac{1}{2}$;\ \ 
$-\frac{1}{2}$)
 $\oplus$
$\mathrm{i}$($-\sqrt{\frac{1}{3}}$,
$\sqrt{\frac{2}{3}}$;\ \ 
$\sqrt{\frac{1}{3}}$)
 $\oplus$
($1$)

Resolved. Number of valid $(S,T)$ pairs = 0.

\vskip 2ex

 \noindent9. (dims;levels) =$(3,
3;5,
3
)$,
irreps = $3_{5}^{1}\oplus
3_{3}^{1,0}$,
pord$(\rho_\text{isum}(\mathfrak{t})) = 15$,

\vskip 0.7ex
\hangindent=4em \hangafter=1
 $\rho_\text{isum}(\mathfrak{t})$ =
 $( 0,
\frac{1}{5},
\frac{4}{5} )
\oplus
( 0,
\frac{1}{3},
\frac{2}{3} )
$,

\vskip 0.7ex
\hangindent=4em \hangafter=1
 $\rho_\text{isum}(\mathfrak{s})$ =
($\sqrt{\frac{1}{5}}$,
$-\sqrt{\frac{2}{5}}$,
$-\sqrt{\frac{2}{5}}$;
$-\frac{5+\sqrt{5}}{10}$,
$\frac{5-\sqrt{5}}{10}$;
$-\frac{5+\sqrt{5}}{10}$)
 $\oplus$
($-\frac{1}{3}$,
$\frac{2}{3}$,
$\frac{2}{3}$;
$-\frac{1}{3}$,
$\frac{2}{3}$;
$-\frac{1}{3}$)

Resolved. Number of valid $(S,T)$ pairs = 0.

\vskip 2ex

 \noindent10. (dims;levels) =$(3,
3;16,
16
)$,
irreps = $3_{16}^{1,0}\oplus
3_{16}^{1,6}$,
pord$(\rho_\text{isum}(\mathfrak{t})) = 16$,

\vskip 0.7ex
\hangindent=4em \hangafter=1
 $\rho_\text{isum}(\mathfrak{t})$ =
 $( \frac{1}{8},
\frac{1}{16},
\frac{9}{16} )
\oplus
( \frac{5}{8},
\frac{1}{16},
\frac{9}{16} )
$,

\vskip 0.7ex
\hangindent=4em \hangafter=1
 $\rho_\text{isum}(\mathfrak{s})$ =
$\mathrm{i}$($0$,
$\sqrt{\frac{1}{2}}$,
$\sqrt{\frac{1}{2}}$;\ \ 
$-\frac{1}{2}$,
$\frac{1}{2}$;\ \ 
$-\frac{1}{2}$)
 $\oplus$
$\mathrm{i}$($0$,
$\sqrt{\frac{1}{2}}$,
$\sqrt{\frac{1}{2}}$;\ \ 
$\frac{1}{2}$,
$-\frac{1}{2}$;\ \ 
$\frac{1}{2}$)

Unresolved. 

\vskip 2ex

\noindent11. (dims;levels) =$(3,
3;5,
4
)$,
irreps = $3_{5}^{1}\oplus
3_{4}^{1,0}$,
pord$(\rho_\text{isum}(\mathfrak{t})) = 20$,

\vskip 0.7ex
\hangindent=4em \hangafter=1
 $\rho_\text{isum}(\mathfrak{t})$ =
 $( 0,
\frac{1}{5},
\frac{4}{5} )
\oplus
( 0,
\frac{1}{4},
\frac{3}{4} )
$,

\vskip 0.7ex
\hangindent=4em \hangafter=1
 $\rho_\text{isum}(\mathfrak{s})$ =
($\sqrt{\frac{1}{5}}$,
$-\sqrt{\frac{2}{5}}$,
$-\sqrt{\frac{2}{5}}$;
$-\frac{5+\sqrt{5}}{10}$,
$\frac{5-\sqrt{5}}{10}$;
$-\frac{5+\sqrt{5}}{10}$)
 $\oplus$
($0$,
$\sqrt{\frac{1}{2}}$,
$\sqrt{\frac{1}{2}}$;
$-\frac{1}{2}$,
$\frac{1}{2}$;
$-\frac{1}{2}$)

Resolved. Number of valid $(S,T)$ pairs = 2.

\vskip 2ex

 \noindent12. (dims;levels) =$(4,
1,
1;9,
1,
1
)$,
irreps = $4_{9,2}^{1,0}\oplus
1_{1}^{1}\oplus
1_{1}^{1}$,
pord$(\rho_\text{isum}(\mathfrak{t})) = 9$,

\vskip 0.7ex
\hangindent=4em \hangafter=1
 $\rho_\text{isum}(\mathfrak{t})$ =
 $( 0,
\frac{1}{9},
\frac{4}{9},
\frac{7}{9} )
\oplus
( 0 )
\oplus
( 0 )
$,

\vskip 0.7ex
\hangindent=4em \hangafter=1
 $\rho_\text{isum}(\mathfrak{s})$ =
($0$,
$-\sqrt{\frac{1}{3}}$,
$-\sqrt{\frac{1}{3}}$,
$-\sqrt{\frac{1}{3}}$;
$\frac{1}{3}c^{2}_{9}
$,
$\frac{1}{3} c_9^4 $,
$\frac{1}{3}c^{1}_{9}
$;
$\frac{1}{3}c^{1}_{9}
$,
$\frac{1}{3}c^{2}_{9}
$;
$\frac{1}{3} c_9^4 $)
 $\oplus$
($1$)
 $\oplus$
($1$)

Unresolved. 

\vskip 2ex

\noindent13. (dims;levels) =$(4,
2;5,
5
)$,
irreps = $4_{5,1}^{1}\oplus
2_{5}^{1}$,
pord$(\rho_\text{isum}(\mathfrak{t})) = 5$,

\vskip 0.7ex
\hangindent=4em \hangafter=1
 $\rho_\text{isum}(\mathfrak{t})$ =
 $( \frac{1}{5},
\frac{2}{5},
\frac{3}{5},
\frac{4}{5} )
\oplus
( \frac{1}{5},
\frac{4}{5} )
$,

\vskip 0.7ex
\hangindent=4em \hangafter=1
 $\rho_\text{isum}(\mathfrak{s})$ =
$\mathrm{i}$($\frac{1}{5}c^{1}_{20}
+\frac{1}{5}c^{3}_{20}
$,
$\frac{2}{5}c^{2}_{15}
+\frac{1}{5}c^{3}_{15}
$,
$-\frac{1}{5}+\frac{2}{5}c^{1}_{15}
-\frac{1}{5}c^{3}_{15}
$,
$\frac{1}{5}c^{1}_{20}
-\frac{1}{5}c^{3}_{20}
$;\ \ 
$-\frac{1}{5}c^{1}_{20}
+\frac{1}{5}c^{3}_{20}
$,
$-\frac{1}{5}c^{1}_{20}
-\frac{1}{5}c^{3}_{20}
$,
$\frac{1}{5}-\frac{2}{5}c^{1}_{15}
+\frac{1}{5}c^{3}_{15}
$;\ \ 
$\frac{1}{5}c^{1}_{20}
-\frac{1}{5}c^{3}_{20}
$,
$\frac{2}{5}c^{2}_{15}
+\frac{1}{5}c^{3}_{15}
$;\ \ 
$-\frac{1}{5}c^{1}_{20}
-\frac{1}{5}c^{3}_{20}
$)
 $\oplus$
$\mathrm{i}$($-\frac{1}{\sqrt{5}}c^{3}_{20}
$,
$\frac{1}{\sqrt{5}}c^{1}_{20}
$;\ \ 
$\frac{1}{\sqrt{5}}c^{3}_{20}
$)

Unresolved. 

\vskip 2ex

\noindent14. (dims;levels) =$(4,
2;5,
5
;a)$,
irreps = $4_{5,2}^{1}\oplus
2_{5}^{1}$,
pord$(\rho_\text{isum}(\mathfrak{t})) = 5$,

\vskip 0.7ex
\hangindent=4em \hangafter=1
 $\rho_\text{isum}(\mathfrak{t})$ =
 $( \frac{1}{5},
\frac{2}{5},
\frac{3}{5},
\frac{4}{5} )
\oplus
( \frac{1}{5},
\frac{4}{5} )
$,

\vskip 0.7ex
\hangindent=4em \hangafter=1
 $\rho_\text{isum}(\mathfrak{s})$ =
($\sqrt{\frac{1}{5}}$,
$\frac{-5+\sqrt{5}}{10}$,
$-\frac{5+\sqrt{5}}{10}$,
$\sqrt{\frac{1}{5}}$;
$-\sqrt{\frac{1}{5}}$,
$\sqrt{\frac{1}{5}}$,
$\frac{5+\sqrt{5}}{10}$;
$-\sqrt{\frac{1}{5}}$,
$\frac{5-\sqrt{5}}{10}$;
$\sqrt{\frac{1}{5}}$)
 $\oplus$
$\mathrm{i}$($-\frac{1}{\sqrt{5}}c^{3}_{20}
$,
$\frac{1}{\sqrt{5}}c^{1}_{20}
$;\ \ 
$\frac{1}{\sqrt{5}}c^{3}_{20}
$)

Resolved. Number of valid $(S,T)$ pairs = 0.

\vskip 2ex

 \noindent15. (dims;levels) =$(4,
2;10,
5
)$,
irreps = $2_{5}^{1}
\hskip -1.5pt \otimes \hskip -1.5pt
2_{2}^{1,0}\oplus
2_{5}^{1}$,
pord$(\rho_\text{isum}(\mathfrak{t})) = 10$,

\vskip 0.7ex
\hangindent=4em \hangafter=1
 $\rho_\text{isum}(\mathfrak{t})$ =
 $( \frac{1}{5},
\frac{4}{5},
\frac{3}{10},
\frac{7}{10} )
\oplus
( \frac{1}{5},
\frac{4}{5} )
$,

\vskip 0.7ex
\hangindent=4em \hangafter=1
 $\rho_\text{isum}(\mathfrak{s})$ =
$\mathrm{i}$($\frac{1}{2\sqrt{5}}c^{3}_{20}
$,
$\frac{1}{2\sqrt{5}}c^{1}_{20}
$,
$\frac{3}{2\sqrt{15}}c^{1}_{20}
$,
$\frac{3}{2\sqrt{15}}c^{3}_{20}
$;\ \ 
$-\frac{1}{2\sqrt{5}}c^{3}_{20}
$,
$-\frac{3}{2\sqrt{15}}c^{3}_{20}
$,
$\frac{3}{2\sqrt{15}}c^{1}_{20}
$;\ \ 
$\frac{1}{2\sqrt{5}}c^{3}_{20}
$,
$-\frac{1}{2\sqrt{5}}c^{1}_{20}
$;\ \ 
$-\frac{1}{2\sqrt{5}}c^{3}_{20}
$)
 $\oplus$
$\mathrm{i}$($-\frac{1}{\sqrt{5}}c^{3}_{20}
$,
$\frac{1}{\sqrt{5}}c^{1}_{20}
$;\ \ 
$\frac{1}{\sqrt{5}}c^{3}_{20}
$)

Unresolved. 

\vskip 2ex

\noindent16. (dims;levels) =$(4,
2;15,
5
)$,
irreps = $2_{5}^{1}
\hskip -1.5pt \otimes \hskip -1.5pt
2_{3}^{1,0}\oplus
2_{5}^{1}$,
pord$(\rho_\text{isum}(\mathfrak{t})) = 15$,

\vskip 0.7ex
\hangindent=4em \hangafter=1
 $\rho_\text{isum}(\mathfrak{t})$ =
 $( \frac{1}{5},
\frac{4}{5},
\frac{2}{15},
\frac{8}{15} )
\oplus
( \frac{1}{5},
\frac{4}{5} )
$,

\vskip 0.7ex
\hangindent=4em \hangafter=1
 $\rho_\text{isum}(\mathfrak{s})$ =
($-\frac{1}{\sqrt{15}}c^{3}_{20}
$,
$\frac{1}{\sqrt{15}}c^{1}_{20}
$,
$\frac{2}{\sqrt{30}}c^{1}_{20}
$,
$-\frac{2}{\sqrt{30}}c^{3}_{20}
$;
$\frac{1}{\sqrt{15}}c^{3}_{20}
$,
$\frac{2}{\sqrt{30}}c^{3}_{20}
$,
$\frac{2}{\sqrt{30}}c^{1}_{20}
$;
$-\frac{1}{\sqrt{15}}c^{3}_{20}
$,
$-\frac{1}{\sqrt{15}}c^{1}_{20}
$;
$\frac{1}{\sqrt{15}}c^{3}_{20}
$)
 $\oplus$
$\mathrm{i}$($-\frac{1}{\sqrt{5}}c^{3}_{20}
$,
$\frac{1}{\sqrt{5}}c^{1}_{20}
$;\ \ 
$\frac{1}{\sqrt{5}}c^{3}_{20}
$)

Resolved. Number of valid $(S,T)$ pairs = 1.

\vskip 2ex

 \noindent17. (dims;levels) =$(4,
2;7,
3
)$,
irreps = $4_{7}^{1}\oplus
2_{3}^{1,0}$,
pord$(\rho_\text{isum}(\mathfrak{t})) = 21$,

\vskip 0.7ex
\hangindent=4em \hangafter=1
 $\rho_\text{isum}(\mathfrak{t})$ =
 $( 0,
\frac{1}{7},
\frac{2}{7},
\frac{4}{7} )
\oplus
( 0,
\frac{1}{3} )
$,

\vskip 0.7ex
\hangindent=4em \hangafter=1
 $\rho_\text{isum}(\mathfrak{s})$ =
$\mathrm{i}$($-\sqrt{\frac{1}{7}}$,
$\sqrt{\frac{2}{7}}$,
$\sqrt{\frac{2}{7}}$,
$\sqrt{\frac{2}{7}}$;\ \ 
$-\frac{1}{\sqrt{7}}c^{2}_{7}
$,
$-\frac{1}{\sqrt{7}}c^{1}_{7}
$,
$\frac{1}{\sqrt{7}\mathrm{i}}s^{5}_{28}
$;\ \ 
$\frac{1}{\sqrt{7}\mathrm{i}}s^{5}_{28}
$,
$-\frac{1}{\sqrt{7}}c^{2}_{7}
$;\ \ 
$-\frac{1}{\sqrt{7}}c^{1}_{7}
$)
 $\oplus$
$\mathrm{i}$($-\sqrt{\frac{1}{3}}$,
$\sqrt{\frac{2}{3}}$;\ \ 
$\sqrt{\frac{1}{3}}$)

Resolved. Number of valid $(S,T)$ pairs = 1.

\vskip 2ex

 \noindent18. (dims;levels) =$(5,
1;5,
1
)$,
irreps = $5_{5}^{1}\oplus
1_{1}^{1}$,
pord$(\rho_\text{isum}(\mathfrak{t})) = 5$,

\vskip 0.7ex
\hangindent=4em \hangafter=1
 $\rho_\text{isum}(\mathfrak{t})$ =
 $( 0,
\frac{1}{5},
\frac{2}{5},
\frac{3}{5},
\frac{4}{5} )
\oplus
( 0 )
$,

\vskip 0.7ex
\hangindent=4em \hangafter=1
 $\rho_\text{isum}(\mathfrak{s})$ =
($-\frac{1}{5}$,
$\sqrt{\frac{6}{25}}$,
$\sqrt{\frac{6}{25}}$,
$\sqrt{\frac{6}{25}}$,
$\sqrt{\frac{6}{25}}$;
$\frac{3-\sqrt{5}}{10}$,
$-\frac{1+\sqrt{5}}{5}$,
$\frac{-1+\sqrt{5}}{5}$,
$\frac{3+\sqrt{5}}{10}$;
$\frac{3+\sqrt{5}}{10}$,
$\frac{3-\sqrt{5}}{10}$,
$\frac{-1+\sqrt{5}}{5}$;
$\frac{3+\sqrt{5}}{10}$,
$-\frac{1+\sqrt{5}}{5}$;
$\frac{3-\sqrt{5}}{10}$)
 $\oplus$
($1$)

Unresolved. 

\vskip 2ex

\noindent19. (dims;levels) =$(6;9
)$,
irreps = $6_{9,3}^{1,0}$,
pord$(\rho_\text{isum}(\mathfrak{t})) = 9$,

\vskip 0.7ex
\hangindent=4em \hangafter=1
 $\rho_\text{isum}(\mathfrak{t})$ =
 $( \frac{1}{9},
\frac{2}{9},
\frac{4}{9},
\frac{5}{9},
\frac{7}{9},
\frac{8}{9} )
$,

\vskip 0.7ex
\hangindent=4em \hangafter=1
 $\rho_\text{isum}(\mathfrak{s})$ =
($\frac{1}{3}$,
$\frac{1}{3}c^{2}_{9}
$,
$\frac{1}{3}$,
$-\frac{1}{3}c^{1}_{9}
$,
$\frac{1}{3}$,
$\frac{1}{3} c_9^4 $;
$\frac{1}{3}$,
$\frac{1}{3} c_9^4 $,
$-\frac{1}{3}$,
$\frac{1}{3}c^{1}_{9}
$,
$\frac{1}{3}$;
$\frac{1}{3}$,
$-\frac{1}{3}c^{2}_{9}
$,
$\frac{1}{3}$,
$\frac{1}{3}c^{1}_{9}
$;
$\frac{1}{3}$,
$-\frac{1}{3} c_9^4 $,
$-\frac{1}{3}$;
$\frac{1}{3}$,
$\frac{1}{3}c^{2}_{9}
$;
$\frac{1}{3}$)

Resolved. Number of valid $(S,T)$ pairs = 1.

\vskip 2ex

 \noindent20. (dims;levels) =$(6;13
)$,
irreps = $6_{13}^{1}$,
pord$(\rho_\text{isum}(\mathfrak{t})) = 13$,

\vskip 0.7ex
\hangindent=4em \hangafter=1
 $\rho_\text{isum}(\mathfrak{t})$ =
 $( \frac{1}{13},
\frac{3}{13},
\frac{4}{13},
\frac{9}{13},
\frac{10}{13},
\frac{12}{13} )
$,

\vskip 0.7ex
\hangindent=4em \hangafter=1
 $\rho_\text{isum}(\mathfrak{s})$ =
$\mathrm{i}$($-\frac{1}{\sqrt{13}}c^{5}_{52}
$,
$\frac{1}{\sqrt{13}}c^{7}_{52}
$,
$\frac{1}{\sqrt{13}}c^{3}_{52}
$,
$\frac{1}{\sqrt{13}}c^{11}_{52}
$,
$\frac{1}{\sqrt{13}}c^{9}_{52}
$,
$-\frac{1}{\sqrt{13}}c^{1}_{52}
$;\ \ 
$-\frac{1}{\sqrt{13}}c^{11}_{52}
$,
$\frac{1}{\sqrt{13}}c^{1}_{52}
$,
$-\frac{1}{\sqrt{13}}c^{5}_{52}
$,
$\frac{1}{\sqrt{13}}c^{3}_{52}
$,
$\frac{1}{\sqrt{13}}c^{9}_{52}
$;\ \ 
$\frac{1}{\sqrt{13}}c^{7}_{52}
$,
$\frac{1}{\sqrt{13}}c^{9}_{52}
$,
$-\frac{1}{\sqrt{13}}c^{5}_{52}
$,
$\frac{1}{\sqrt{13}}c^{11}_{52}
$;\ \ 
$-\frac{1}{\sqrt{13}}c^{7}_{52}
$,
$-\frac{1}{\sqrt{13}}c^{1}_{52}
$,
$-\frac{1}{\sqrt{13}}c^{3}_{52}
$;\ \ 
$\frac{1}{\sqrt{13}}c^{11}_{52}
$,
$-\frac{1}{\sqrt{13}}c^{7}_{52}
$;\ \ 
$\frac{1}{\sqrt{13}}c^{5}_{52}
$)

Resolved. Number of valid $(S,T)$ pairs = 1.

\vskip 2ex

 \noindent21. (dims;levels) =$(6;15
)$,
irreps = $3_{3}^{1,0}
\hskip -1.5pt \otimes \hskip -1.5pt
2_{5}^{1}$,
pord$(\rho_\text{isum}(\mathfrak{t})) = 15$,

\vskip 0.7ex
\hangindent=4em \hangafter=1
 $\rho_\text{isum}(\mathfrak{t})$ =
 $( \frac{1}{5},
\frac{4}{5},
\frac{2}{15},
\frac{7}{15},
\frac{8}{15},
\frac{13}{15} )
$,

\vskip 0.7ex
\hangindent=4em \hangafter=1
 $\rho_\text{isum}(\mathfrak{s})$ =
$\mathrm{i}$($\frac{1}{3\sqrt{5}}c^{3}_{20}
$,
$\frac{1}{3\sqrt{5}}c^{1}_{20}
$,
$\frac{2}{3\sqrt{5}}c^{1}_{20}
$,
$\frac{2}{3\sqrt{5}}c^{1}_{20}
$,
$\frac{2}{3\sqrt{5}}c^{3}_{20}
$,
$\frac{2}{3\sqrt{5}}c^{3}_{20}
$;\ \ 
$-\frac{1}{3\sqrt{5}}c^{3}_{20}
$,
$-\frac{2}{3\sqrt{5}}c^{3}_{20}
$,
$-\frac{2}{3\sqrt{5}}c^{3}_{20}
$,
$\frac{2}{3\sqrt{5}}c^{1}_{20}
$,
$\frac{2}{3\sqrt{5}}c^{1}_{20}
$;\ \ 
$-\frac{1}{3\sqrt{5}}c^{3}_{20}
$,
$\frac{2}{3\sqrt{5}}c^{3}_{20}
$,
$\frac{1}{3\sqrt{5}}c^{1}_{20}
$,
$-\frac{2}{3\sqrt{5}}c^{1}_{20}
$;\ \ 
$-\frac{1}{3\sqrt{5}}c^{3}_{20}
$,
$-\frac{2}{3\sqrt{5}}c^{1}_{20}
$,
$\frac{1}{3\sqrt{5}}c^{1}_{20}
$;\ \ 
$\frac{1}{3\sqrt{5}}c^{3}_{20}
$,
$-\frac{2}{3\sqrt{5}}c^{3}_{20}
$;\ \ 
$\frac{1}{3\sqrt{5}}c^{3}_{20}
$)

Resolved. Number of valid $(S,T)$ pairs = 0.

\vskip 2ex

 \noindent22. (dims;levels) =$(6;16
)$,
irreps = $6_{16,1}^{1,0}$,
pord$(\rho_\text{isum}(\mathfrak{t})) = 16$,

\vskip 0.7ex
\hangindent=4em \hangafter=1
 $\rho_\text{isum}(\mathfrak{t})$ =
 $( 0,
\frac{1}{4},
\frac{1}{16},
\frac{5}{16},
\frac{9}{16},
\frac{13}{16} )
$,

\vskip 0.7ex
\hangindent=4em \hangafter=1
 $\rho_\text{isum}(\mathfrak{s})$ =
$\mathrm{i}$($0$,
$0$,
$\frac{1}{2}$,
$\frac{1}{2}$,
$\frac{1}{2}$,
$\frac{1}{2}$;\ \ 
$0$,
$\frac{1}{2}$,
$-\frac{1}{2}$,
$\frac{1}{2}$,
$-\frac{1}{2}$;\ \ 
$-\sqrt{\frac{1}{8}}$,
$-\sqrt{\frac{1}{8}}$,
$\sqrt{\frac{1}{8}}$,
$\sqrt{\frac{1}{8}}$;\ \ 
$\sqrt{\frac{1}{8}}$,
$\sqrt{\frac{1}{8}}$,
$-\sqrt{\frac{1}{8}}$;\ \ 
$-\sqrt{\frac{1}{8}}$,
$-\sqrt{\frac{1}{8}}$;\ \ 
$\sqrt{\frac{1}{8}}$)

Resolved. Number of valid $(S,T)$ pairs = 4.

\vskip 2ex

 \noindent23. (dims;levels) =$(6;35
)$,
irreps = $3_{7}^{3}
\hskip -1.5pt \otimes \hskip -1.5pt
2_{5}^{2}$,
pord$(\rho_\text{isum}(\mathfrak{t})) = 35$,

\vskip 0.7ex
\hangindent=4em \hangafter=1
 $\rho_\text{isum}(\mathfrak{t})$ =
 $( \frac{1}{35},
\frac{4}{35},
\frac{9}{35},
\frac{11}{35},
\frac{16}{35},
\frac{29}{35} )
$,

\vskip 0.7ex
\hangindent=4em \hangafter=1
 $\rho_\text{isum}(\mathfrak{s})$ =
$\mathrm{i}$($-\frac{4}{35}c^{1}_{140}
-\frac{3}{35}c^{3}_{140}
-\frac{1}{7}c^{5}_{140}
+\frac{1}{35}c^{7}_{140}
+\frac{1}{35}c^{9}_{140}
+\frac{4}{35}c^{13}_{140}
+\frac{2}{35}c^{15}_{140}
-\frac{3}{35}c^{17}_{140}
+\frac{9}{35}c^{19}_{140}
-\frac{4}{35}c^{21}_{140}
-\frac{2}{7}c^{23}_{140}
$,
$-\frac{1}{\sqrt{35}}c^{4}_{35}
+\frac{1}{\sqrt{35}}c^{11}_{35}
$,
$\frac{1}{\sqrt{35}}c^{1}_{35}
-\frac{1}{\sqrt{35}}c^{6}_{35}
$,
$\frac{2}{\sqrt{35}}c^{3}_{35}
+\frac{1}{\sqrt{35}}c^{4}_{35}
+\frac{1}{\sqrt{35}}c^{10}_{35}
+\frac{1}{\sqrt{35}}c^{11}_{35}
$,
$-\frac{1}{\sqrt{35}\mathrm{i}}s^{3}_{140}
-\frac{1}{\sqrt{35}\mathrm{i}}s^{17}_{140}
$,
$\frac{2}{35}c^{1}_{140}
-\frac{1}{35}c^{3}_{140}
-\frac{1}{7}c^{5}_{140}
-\frac{3}{35}c^{7}_{140}
+\frac{1}{5}c^{9}_{140}
-\frac{2}{35}c^{13}_{140}
-\frac{1}{35}c^{15}_{140}
-\frac{1}{35}c^{17}_{140}
+\frac{3}{35}c^{19}_{140}
+\frac{2}{35}c^{21}_{140}
-\frac{2}{7}c^{23}_{140}
$;\ \ 
$-\frac{1}{\sqrt{35}\mathrm{i}}s^{3}_{140}
-\frac{1}{\sqrt{35}\mathrm{i}}s^{17}_{140}
$,
$\frac{4}{35}c^{1}_{140}
+\frac{3}{35}c^{3}_{140}
+\frac{1}{7}c^{5}_{140}
-\frac{1}{35}c^{7}_{140}
-\frac{1}{35}c^{9}_{140}
-\frac{4}{35}c^{13}_{140}
-\frac{2}{35}c^{15}_{140}
+\frac{3}{35}c^{17}_{140}
-\frac{9}{35}c^{19}_{140}
+\frac{4}{35}c^{21}_{140}
+\frac{2}{7}c^{23}_{140}
$,
$-\frac{1}{\sqrt{35}}c^{1}_{35}
+\frac{1}{\sqrt{35}}c^{6}_{35}
$,
$-\frac{2}{35}c^{1}_{140}
+\frac{1}{35}c^{3}_{140}
+\frac{1}{7}c^{5}_{140}
+\frac{3}{35}c^{7}_{140}
-\frac{1}{5}c^{9}_{140}
+\frac{2}{35}c^{13}_{140}
+\frac{1}{35}c^{15}_{140}
+\frac{1}{35}c^{17}_{140}
-\frac{3}{35}c^{19}_{140}
-\frac{2}{35}c^{21}_{140}
+\frac{2}{7}c^{23}_{140}
$,
$\frac{2}{\sqrt{35}}c^{3}_{35}
+\frac{1}{\sqrt{35}}c^{4}_{35}
+\frac{1}{\sqrt{35}}c^{10}_{35}
+\frac{1}{\sqrt{35}}c^{11}_{35}
$;\ \ 
$\frac{2}{\sqrt{35}}c^{3}_{35}
+\frac{1}{\sqrt{35}}c^{4}_{35}
+\frac{1}{\sqrt{35}}c^{10}_{35}
+\frac{1}{\sqrt{35}}c^{11}_{35}
$,
$-\frac{2}{35}c^{1}_{140}
+\frac{1}{35}c^{3}_{140}
+\frac{1}{7}c^{5}_{140}
+\frac{3}{35}c^{7}_{140}
-\frac{1}{5}c^{9}_{140}
+\frac{2}{35}c^{13}_{140}
+\frac{1}{35}c^{15}_{140}
+\frac{1}{35}c^{17}_{140}
-\frac{3}{35}c^{19}_{140}
-\frac{2}{35}c^{21}_{140}
+\frac{2}{7}c^{23}_{140}
$,
$\frac{1}{\sqrt{35}}c^{4}_{35}
-\frac{1}{\sqrt{35}}c^{11}_{35}
$,
$-\frac{1}{\sqrt{35}\mathrm{i}}s^{3}_{140}
-\frac{1}{\sqrt{35}\mathrm{i}}s^{17}_{140}
$;\ \ 
$\frac{1}{\sqrt{35}\mathrm{i}}s^{3}_{140}
+\frac{1}{\sqrt{35}\mathrm{i}}s^{17}_{140}
$,
$-\frac{4}{35}c^{1}_{140}
-\frac{3}{35}c^{3}_{140}
-\frac{1}{7}c^{5}_{140}
+\frac{1}{35}c^{7}_{140}
+\frac{1}{35}c^{9}_{140}
+\frac{4}{35}c^{13}_{140}
+\frac{2}{35}c^{15}_{140}
-\frac{3}{35}c^{17}_{140}
+\frac{9}{35}c^{19}_{140}
-\frac{4}{35}c^{21}_{140}
-\frac{2}{7}c^{23}_{140}
$,
$\frac{1}{\sqrt{35}}c^{4}_{35}
-\frac{1}{\sqrt{35}}c^{11}_{35}
$;\ \ 
$-\frac{2}{\sqrt{35}}c^{3}_{35}
-\frac{1}{\sqrt{35}}c^{4}_{35}
-\frac{1}{\sqrt{35}}c^{10}_{35}
-\frac{1}{\sqrt{35}}c^{11}_{35}
$,
$-\frac{1}{\sqrt{35}}c^{1}_{35}
+\frac{1}{\sqrt{35}}c^{6}_{35}
$;\ \ 
$\frac{4}{35}c^{1}_{140}
+\frac{3}{35}c^{3}_{140}
+\frac{1}{7}c^{5}_{140}
-\frac{1}{35}c^{7}_{140}
-\frac{1}{35}c^{9}_{140}
-\frac{4}{35}c^{13}_{140}
-\frac{2}{35}c^{15}_{140}
+\frac{3}{35}c^{17}_{140}
-\frac{9}{35}c^{19}_{140}
+\frac{4}{35}c^{21}_{140}
+\frac{2}{7}c^{23}_{140}
$)

Resolved. Number of valid $(S,T)$ pairs = 1.

\vskip 2ex

 \noindent24. (dims;levels) =$(6;56
)$,
irreps = $3_{7}^{1}
\hskip -1.5pt \otimes \hskip -1.5pt
2_{8}^{1,6}$,
pord$(\rho_\text{isum}(\mathfrak{t})) = 28$,

\vskip 0.7ex
\hangindent=4em \hangafter=1
 $\rho_\text{isum}(\mathfrak{t})$ =
 $( \frac{1}{56},
\frac{9}{56},
\frac{11}{56},
\frac{25}{56},
\frac{43}{56},
\frac{51}{56} )
$,

\vskip 0.7ex
\hangindent=4em \hangafter=1
 $\rho_\text{isum}(\mathfrak{s})$ =
($\frac{1}{\sqrt{14}}c^{1}_{28}
$,
$\frac{1}{\sqrt{14}}c^{3}_{28}
$,
$-\frac{1}{\sqrt{14}}c^{5}_{28}
$,
$-\frac{1}{\sqrt{14}}c^{5}_{28}
$,
$\frac{1}{\sqrt{14}}c^{1}_{28}
$,
$\frac{1}{\sqrt{14}}c^{3}_{28}
$;
$-\frac{1}{\sqrt{14}}c^{5}_{28}
$,
$\frac{1}{\sqrt{14}}c^{1}_{28}
$,
$\frac{1}{\sqrt{14}}c^{1}_{28}
$,
$\frac{1}{\sqrt{14}}c^{3}_{28}
$,
$-\frac{1}{\sqrt{14}}c^{5}_{28}
$;
$-\frac{1}{\sqrt{14}}c^{3}_{28}
$,
$\frac{1}{\sqrt{14}}c^{3}_{28}
$,
$\frac{1}{\sqrt{14}}c^{5}_{28}
$,
$-\frac{1}{\sqrt{14}}c^{1}_{28}
$;
$\frac{1}{\sqrt{14}}c^{3}_{28}
$,
$-\frac{1}{\sqrt{14}}c^{5}_{28}
$,
$\frac{1}{\sqrt{14}}c^{1}_{28}
$;
$-\frac{1}{\sqrt{14}}c^{1}_{28}
$,
$-\frac{1}{\sqrt{14}}c^{3}_{28}
$;
$\frac{1}{\sqrt{14}}c^{5}_{28}
$)

Resolved. Number of valid $(S,T)$ pairs = 2.

\vskip 2ex

 \noindent25. (dims;levels) =$(6;80
)$,
irreps = $3_{16}^{3,3}
\hskip -1.5pt \otimes \hskip -1.5pt
2_{5}^{2}$,
pord$(\rho_\text{isum}(\mathfrak{t})) = 80$,

\vskip 0.7ex
\hangindent=4em \hangafter=1
 $\rho_\text{isum}(\mathfrak{t})$ =
 $( \frac{1}{40},
\frac{9}{40},
\frac{3}{80},
\frac{27}{80},
\frac{43}{80},
\frac{67}{80} )
$,

\vskip 0.7ex
\hangindent=4em \hangafter=1
 $\rho_\text{isum}(\mathfrak{s})$ =
$\mathrm{i}$($0$,
$0$,
$\frac{1}{\sqrt{10}}c^{3}_{20}
$,
$\frac{1}{\sqrt{10}}c^{1}_{20}
$,
$\frac{1}{\sqrt{10}}c^{3}_{20}
$,
$\frac{1}{\sqrt{10}}c^{1}_{20}
$;\ \ 
$0$,
$\frac{1}{\sqrt{10}}c^{1}_{20}
$,
$-\frac{1}{\sqrt{10}}c^{3}_{20}
$,
$\frac{1}{\sqrt{10}}c^{1}_{20}
$,
$-\frac{1}{\sqrt{10}}c^{3}_{20}
$;\ \ 
$-\frac{1}{2\sqrt{5}}c^{1}_{20}
$,
$-\frac{1}{2\sqrt{5}}c^{3}_{20}
$,
$\frac{1}{2\sqrt{5}}c^{1}_{20}
$,
$\frac{1}{2\sqrt{5}}c^{3}_{20}
$;\ \ 
$\frac{1}{2\sqrt{5}}c^{1}_{20}
$,
$\frac{1}{2\sqrt{5}}c^{3}_{20}
$,
$-\frac{1}{2\sqrt{5}}c^{1}_{20}
$;\ \ 
$-\frac{1}{2\sqrt{5}}c^{1}_{20}
$,
$-\frac{1}{2\sqrt{5}}c^{3}_{20}
$;\ \ 
$\frac{1}{2\sqrt{5}}c^{1}_{20}
$)

Resolved. Number of valid $(S,T)$ pairs = 2.

\vskip 2ex

\

\section{A list of candidate modular data from resolved
$\SL$ representations }
\label{ndegMD}

\subsection{The notion of resolved $\SL$ matrix representations}

In the above, we have chosen a special basis in the eigenspaces of an $\SL$
matrix representation $\tilde\rho$ to make $\tilde\rho(\fs)$ symmetric.
But such a special basis is still not special enough to make $\tilde\rho$ to be
an MD representation $\rho$.

We can choose a more special basis to make $\tilde\rho(\fs^2)$ a signed
permutation matrix, and $\tilde\rho(\fs)$ symmetric.  We know
that, for an MD representation $\rho$, $\rho(\fs^2)$ is a signed permutation
matrix. So the new special basis makes $\tilde\rho$ closer to the MD
representation $\rho$.

We can choose an even more special basis in the eigenspaces of
$\tilde\rho(\ft)$ to make $\tilde\rho$ even closer to the MD representation
$\rho$, by using the matrix $D_{\td\rho}(\s)$ in \eqref{Drho}.  For an MD
representation $\rho$, $D_{\rho}(\s)$ is suppose to be signed permutations.  So we
will try to choose a basis to transform each $D_{\td\rho}(\s)$ into signed
permutations.  We like to point out that, since both $\td\rho$ and $\rho$ are
symmetric $\SL$ matrix representations that are related by an unitary
transformation, according to Theorem \ref{t:ortho_eqv}, they can be related by
an orthogonal transformation.  

%Some times, the number of orthogonal transformations that
%transform $D_{\td\rho}(\s)$'s into signed permutations is finite, which give
%rise to a finite number of potential MD representations $\rho$'s.  In this
%case, we can directly check which of those $\rho$'s satisfy Proposition
%\ref{p:MD} and may actually correspond to the MD representations.  If the
%number of orthogonal transformations is infinite, we have to use other methods
%to handle those more difficult cases (see Section ???). 

%If $\s_\text{c}$ is the complex conjugation, then $a=b=-1$.  From
%$(\rho(\fs)^{-1} \rho(\ft))^3=\id$, we find $\rho(\ft)^{-1} \rho(\fs)
%\rho(\ft)^{-1} \rho(\fs) \rho(\ft)^{-1}= \rho(\fs)^{-1}$ and \begin{align}
%D_\rho(\s_\text{c}) 
%%=\ol{\rho(\fs)} \rho(\fs)\inv 
%=\rho(\fs)^{-2} 
%=\rho(\fs)^2 
%= \pm C,
%\end{align}
%where $C$ is the charge conjugation matrix of $\CC$.  Since $\s_\text{c}^2(
%\rho(\ft))  = \rho(\ft)$, $C  \rho(\ft) =  \rho(\ft) C$ and we can diagonalize
%$\rho(\ft)$ and $C$ simultaneously.  This allows us to define non-degenerate
%subspace using the pair of eigenvalues of $\rho(\ft)$ and $C$.  

Let us consider a simple case to demonstrate our approach.  If 
$\td\rho(\ft)$ is non-degenerate, then  $D_{\td\rho}(\s)$ will
automatically be a signed permutation matrix. Using signed diagonal matrices
$V_\mathrm{sd}$, we can transform $\td\rho$ to many other symmetric
representations, $\rho$'s: 
\begin{align}
 \rho = V_\mathrm{sd} \td\rho V_\mathrm{sd} ,
\end{align}
where $D_{\rho}(\s)$ remains a signed permutation.  In fact the
signed diagonal matrices $V_\mathrm{sd}$ are the most general orthogonal
matrices that fix $\td\rho(\ft)$ and transform all
$D_{\td\rho}(\s)$'s into (potentially different) signed permutations.  Thus the
resulting symmetric representations, $\rho$'s, include all the symmetric
representations where $D_{\rho}(\s)$'s are signed permutations.
From those $\rho$'s, we can then construct many pairs of $S,T$ matrices via
\eqref{STrho1}, and check which one satisfies the conditions in Proposition
\ref{p:MDcond}.  Those $S,T$ matrices that satisfy those conditions
may very likely correspond to modular data (or MTC's).  If none of the $S,T$
matrices satisfy the conditions, then the representation $\td\rho$ will not be
an $\SL$ representation of any modular data.

When some eigenspaces of $\td\rho(\ft)$ are more than 1-dimensional, then the
$D_{\td\rho}(\s)$ may not be signed permutations.  There may be infinite many
orthogonal matrices that can transform $D_{\td\rho}(\s)$ into signed
permutations, which make the subsequent selection difficult.  In the following,
we will generalize the above notion of non-degenerate representation, to
include some cases where some eigenspaces of $\td\rho(\ft)$ are 2-dimensional
or more.  We will show that, for those special representations, there is only a
finite number of orthogonal matrices that can transform $D_{\td\rho}(\s)$ into
signed permutations. 

To carry through this program, 
%let
%\begin{align}
%\Omega_{\td\rho}^{(2)} = \{ \s \in \Gal(\BQ_{\ord(\td\rho(\ft))}) \mid \s^2 =\id\}\,.
%\end{align}
%Then $\Omega_{\td\rho}^{(2)}$ is an elementary 2-subgroup of
%$\Gal(\BQ_{\ord(\td\rho(\ft))})$. It follows from \eqref{Galact} that $\s \in
%\Omega_{\td\rho}^{(2)}$ if and only if 
%\begin{align}
%[D_{\td\rho}(\s), \td\rho(\fs)] =
%[D_{\td\rho}(\s), \td\rho(\ft)]  =0\,.  
%\end{align}
%Therefore,  $D_{\td\rho}(\s)$ stabilizes the $\td\theta$-eigenspace for $\s \in
%\Omega_{\td\rho}^{(2)}$, and commute with each other.  In particular,
%$D_{\td\rho}$ defines a representation of $\Omega_{\td\rho}^{(2)}$ on the
%$\td\theta$-eigenspace.  We can diagonalize $\td\rho(\ft)$ and
%$\{D_{\td\rho}(\s) \mid  \s \in \Omega_{\td\rho}^{(2)} \} $ simultaneously.  If
%all the simultaneous eigenspaces are 1-dimensional, then $\td\rho(\ft)$ is
%unitary equivalent to a finite number of symmetric representations $\rho$ where
%$D_{\rho}(\s)$ for $\s \in \Omega_{\td\rho}^{(2)} $ are signed permutations.
%We can then, from the finite list of $\rho$ representations, select those that
%satisfy Proposition \ref{p:MD}.
%
%In fact, the above discussion can be generalized one step further.  
let us concentrate on an eigenspace $E_{\td\theta}$ of $\td \rho(\ft)$
corresponding to an eigenvalue $\td\theta$, and let
\begin{align}
\Omega_{\td\rho}(\td\theta) 
= \{ \s \in \Gal(\BQ_{\ord(\td\rho(\ft))}) \mid \s^2(\td\theta) =\td\theta \}
\,.
\end{align}
Then $\Omega_{\td\rho}(\td\theta)$ is a subgroup of
$\Gal(\BQ_{\ord(\td\rho(\ft))})$.  By definition, $D_{\td\rho}(\s)$ stabilizes
the $\td\theta$-eigenspace $E_{\td\theta}$ for $\s \in
\Omega_{\td\rho}(\td\theta)$, and commute with each other.  In particular,
$D_{\td\rho}|_{E_{\td\theta}}$ (restricted on $E_{\td\theta}$) defines a
representation of $\Omega_{\td\rho}(\td\theta)$ on $E_{\td\theta}$.
%For simplicity, we denote a Galois automorphism in
%$\Omega_{\td\rho}(\td\theta)$ as $\s_{\invo}$.

We can diagonalize $\{D_{\td\rho}(\s)|_{E_{\td\theta}} \mid  \s \in
\Omega_{\td\rho}(\td\theta) \} $ simultaneously within $E_{\td\theta}$.  The
degeneracy of the $\td\theta$-eigenspace $E_{\td\theta}$ is fully resolved by
these $D_{\td\rho}(\s)$'s, if the common eigenspace of these
$D_{\td\rho}(\s)|_{E_{\td\theta}}$'s are all 1-dimensional. In terms of the
characters of $\Omega_{\td\rho}(\td\theta)$, the degeneracy of $E_{\td\theta}$
can be fully resolved if each irreducible character of
$\Omega_{\td\rho}(\td\theta)$ has multiplicity at most 1 in the character
decomposition of $E_{\td\theta}$ as a representation of
$\Omega_{\td\rho}(\td\theta)$.
Now we can introduce the notion of resolved representation:
\begin{defn}
\label{ndegDef}
A general $\SL$ matrix representation $\td\rho$ is called
\textbf{resolved} if the degeneracy of each of eigenspace of
$\td\rho(\ft)$ is fully resolved by $D_{\td\rho}(\s)$, $\s \in
\Omega_{\td\rho}(\td\theta)$, as described above.  
%In other words, in the
%direct sum decomposition of $D_{\td\rho}(\s)$, as a represenation of
%$\Omega_{\td\rho}(\td\theta)$,
%an irreducible represenation of
%$\Omega_{\td\rho}(\td\theta)$ can appear at most once.
\end{defn}

%Let us also introduce a notion of pseudo-MD representation:
%\begin{defn}
%A \textbf{pseudo-MD} representation $\rho$ is a symmetric matrix
%representation of $\SL$ such that $D_{\rho}(\s)|_{E_{\td\theta}}$,
%for all $\s \in \Omega_{\td\rho}(\td\theta)$, are signed permutations within
%the $\td\theta$-eigenspace.
%\end{defn}

%The basis of a non-degenerate $D$-diagonal symmetric matrix
%representation is almost fixed: Two unitary equivalent non-degenerate
%$D$-diagonal symmetric matrix representations are equivalent under a signed
%diagonal matrix $V_\mathrm{sd}$, according to Proposition \ref{p:ortho_eqv}.
%Because of this, a non-degenerate $D$-diagonal symmetric matrix representation
%$\td\rho$ is unitary equivalent to only a finite number of pseudo-MD
%representations $\rho$.  Here

%We have constructed a list of symmetric irrep-sum $\SL$ matrix representations
%(see Appendix \ref{repL}), that include all the representations from modular
%data. 
%It turns out that those symmetric irrep-sum representations are
%automatically $D$-diagonal, where the diagonal entries of $D_{\td\rho}(\s)$,
%for $\s \in \Omega_{\td\rho}(\td\theta)$, are $\pm 1$, at least for rank 10 and
%less.

Given a symmetric irrep-sum matrix representation (denoted as
$\rho_\mathrm{isum}$), we can use unitary matrices, $U$'s, to transform it into a
symmetric  
%pseudo-MD matrix 
representation $\rho$ via
\begin{align}
\rho(\ft) = U \rho_\mathrm{isum}(\ft) U^\dag,\ \ \ \
\rho(\fs) = U \rho_\mathrm{isum}(\fs) U^\dag.
\end{align}
where $D_{\rho}(\s)|_{E_{\td\theta}}$, for all $\s \in
\Omega_{\td\rho}(\td\theta)$, are signed permutations within the
$\td\theta$-eigenspace.  If $\rho_\mathrm{isum}$ is resolved, then there is
only a finite number of such representations.  We then can check which of those
representations satisfy Proposition \ref{p:MDcond}. This is how we
compute the potential modular data $S,T$'s from resolved
$\rho_\mathrm{isum}$'s.

To show a resolved $\rho_\mathrm{isum}$ is unitarily equivalent to only a
finite number representations whose $D_{\rho}(\s)|_{E_{\td\theta}}$  are signed
permutations, we note that both $\rho$ and $\rho_\mathrm{isum}$ are symmetric,
and according to Theorem \ref{t:ortho_eqv}, $\rho$ and $\rho_\mathrm{isum}$
are in fact orthogonally equivalent, {\it i.e.} the above $U$ can be chosen to
satisfy $U = U^*$ and $UU^\top =\id$.  If the number of most general orthogonal
matrices $U$ that transform $\rho_\mathrm{isum}$ to $\rho$ is finite, then the
number of representations $\rho$ are finite.

Since the orthogonal $U$ acts within the eigenspace of
$\rho_\mathrm{isum}(\ft)$, to show the number of possible $U$'s are finite, we
can concentrate on a single $\td\theta$-eigenspace $E_{\td\theta}$, and denote
$\s \in \Omega_{\td\rho}(\td\theta)$ as $\s_\mathrm{inv}$.  In the following,
we will consider the cases where $E_{\td\theta}$ is 1-dimensional,
2-dimensional, {\it etc.}.  For each case, we will show the number of possible
$U$'s are finite, and give the possible choices of $U$'s.

%We will only consider
%$D_{\rho_\mathrm{isum}}(\s_\mathrm{inv})|_{E_{\td\theta}}$'s within the
%$\td\theta$-eigenspace.

\subsubsection{Within an 1-dimensional eigenspace of $\rho_\mathrm{isum}(\ft)$}
$D_{ \rho_\mathrm{isum}}(\s_\mathrm{inv})|_{E_{\td\theta}} = \pm 1$ are already
signed permutations. In this case the orthogonal matrix $U$ (within the
1-dimensional eigenspace) has only two choices 
\begin{align}
\label{U1x1}
U=\pm 1, 
\end{align}
which is finite.

\subsubsection{Within a 2-dimensional eigenspace of $\rho_\mathrm{isum}(\ft)$}

In this case, the matrix groups $MG$ generated by 2-by-2 matrices, $D_{
\rho_\mathrm{isum}}(\s_\mathrm{inv})|_{E_{\td\theta}}$, can have several
different forms, for those passing representations. By examine the
computer results, we find that, for unresolved cases, matrix groups $MG$ can be
\begin{align}
 MG &=\Big\{  
\begin{pmatrix} 1 &0\\ 0 & 1\\ \end{pmatrix} 
\Big\}, 
&&\text{ for } \dim(\rho_\mathrm{isum}) \geq 5;
\nonumber\\
 MG &=\Big\{  
\begin{pmatrix} 1 &0\\ 0 & 1\\ \end{pmatrix}, 
-\begin{pmatrix} 1 &0\\ 0 & 1\\ \end{pmatrix} 
\Big\}, 
&&\text{ for } \dim(\rho_\mathrm{isum}) \geq 6.
\end{align}
For resolved cases, we have 
\begin{align}
 MG &=\Big\{  
\begin{pmatrix} 1 &0\\ 0 & 1\\ \end{pmatrix}, 
\begin{pmatrix} 1 &0\\ 0 & -1\\ \end{pmatrix} 
\Big\}, 
&&\text{ for } \dim(\rho_\mathrm{isum}) \geq 4;
\nonumber\\
 MG &=\Big\{  
\begin{pmatrix} 1 &0\\ 0 & 1\\ \end{pmatrix}, 
-\begin{pmatrix} 1 &0\\ 0 & 1\\ \end{pmatrix}, 
\begin{pmatrix} 1 &0\\ 0 & -1\\ \end{pmatrix}, 
\begin{pmatrix} -1 &0\\ 0 & 1\\ \end{pmatrix} 
\Big\}, 
&&\text{ for } \dim(\rho_\mathrm{isum}) \geq 6.
\end{align}
In those two cases,
\begin{align}
 U = \frac1{\sqrt 2} \begin{pmatrix}
 1& 1\\
 1&-1\\
\end{pmatrix}
\ \ \text{ or } \ \
 U = \frac1{\sqrt 2} \begin{pmatrix}
-1& 1\\
 1& 1\\
\end{pmatrix}
\ \ \text{ or } \ \
 U = \begin{pmatrix}
 1& 0\\
 0& 1\\
\end{pmatrix}
\end{align}
will transform all  $D_{
\rho_\mathrm{isum}}(\s_\mathrm{inv})|_{E_{\td\theta}}$'s into signed
permutations.  
In general we have
\begin{thm}
\label{MG2x2}
Let
\begin{align}
 MG_2 &=\Big\{  
\begin{pmatrix} 1 &0\\ 0 & 1\\ \end{pmatrix}, 
\begin{pmatrix} 1 &0\\ 0 & -1\\ \end{pmatrix} 
\Big\}, 
\nonumber\\
 MG_4 &=\Big\{  
\begin{pmatrix} 1 &0\\ 0 & 1\\ \end{pmatrix}, 
-\begin{pmatrix} 1 &0\\ 0 & 1\\ \end{pmatrix}, 
\begin{pmatrix} 1 &0\\ 0 & -1\\ \end{pmatrix}, 
\begin{pmatrix} -1 &0\\ 0 & 1\\ \end{pmatrix} 
\Big\}. 
\end{align}
The most general orthogonal matrices that transform all  matrices in $MG_2$ or
$MG_4$ into signed permutations must have one of the following forms 
\begin{align}
\label{U2x2}
 U =  
\frac{PV_\mathrm{sd}}{\sqrt 2}  
\begin{pmatrix}
 1& 1\\
 1&-1\\
\end{pmatrix} ,
 \text{ or } 
 U =   \frac{PV_\mathrm{sd}}{\sqrt 2} 
 \begin{pmatrix}
-1& 1\\
 1& 1\\
\end{pmatrix} ,
 \text{ or } 
 U =  P V_\mathrm{sd} \begin{pmatrix}
 1& 0\\
 0& 1\\
\end{pmatrix}
\end{align}
where $V_\mathrm{sd}$ are signed diagonal matrices, and $P$ are permutation
matrices.  The number of the orthogonal transformations $U$ is finite.
\end{thm}
\begin{proof}[Proof of Theorem \ref{MG2x2}] 
We only needs to consider the first matrix group $MG_2$, where the matrix group
is isomorphic to the $\mathbb{Z}_2$ group.  There are only four matrix groups
formed by 2-dimensional signed permutations matrices, that are  isomorphic
$\mathbb{Z}_2$.  The four matrix groups are generated by the following four
generators respectively:
\begin{align}
 \begin{pmatrix}
 1& 0\\
 0&-1\\
\end{pmatrix},\ \
 \begin{pmatrix}
-1& 0\\
 0&1\\
\end{pmatrix},\ \
 \begin{pmatrix}
0&1\\
1&0\\
\end{pmatrix},\ \
 \begin{pmatrix}
0&-1\\
-1&0\\
\end{pmatrix}.
\end{align}
An orthogonal transformation $U$ that transforms $MG$ to one of the above
matrix groups must have a from $U=VU_0$, where $V$ transforms $MG_2$ 
into itself,
and $U_0$ is a fixed orthogonal transformation that transforms $MG_2$ to one of
the above matrix groups.
We can choose $U_0$ to have the following form
\begin{align}
 U_0 =  
\frac{P}{\sqrt 2}  
\begin{pmatrix}
 1& 1\\
 1&-1\\
\end{pmatrix} ,
 \text{ or } 
 U_0 =   \frac{P}{\sqrt 2} 
 \begin{pmatrix}
-1& 1\\
 1& 1\\
\end{pmatrix} ,
 \text{ or } 
 U_0 =  P  \begin{pmatrix}
 1& 0\\
 0& 1\\
\end{pmatrix}.
\end{align}
To keep $MG$ unchanged
$V$ must satisfy
\begin{align}
 V  
\begin{pmatrix}
 1& 0\\
 0& -1\\
\end{pmatrix}
=
\begin{pmatrix}
 1& 0\\
 0& -1\\
\end{pmatrix}
V.
\end{align}
We find that $V$ must be diagonal. Thus $V$, as an orthogonal matrix, must be
signed diagonal. This gives us the result \eqref{U2x2}.
\end{proof}

If $\dim(\rho_\mathrm{isum}) \geq 8$, it is possible that the matrix group
of $D_{
\rho_\mathrm{isum}}(\s_\mathrm{inv})|_{E_{\td\theta}}$'s is generated by the
following non-diagonal matrix 
\begin{align}
\pm \begin{pmatrix}
0 & -1\\
1 &  0\\
\end{pmatrix}
\end{align}
This is because the direct sum decomposition of $\rho_\mathrm{isum}$ contains a
dimension-6 irreducible representation $6_1^{0,1}$ in Appendix \ref{repPP},
whose $\rho(\ft)$ has a 2-dimensional eigenspace.  The representation
$6_1^{0,1}$ can give rise to such form of $D_{
\rho_\mathrm{isum}}(\s_\mathrm{inv})|_{E_{\td\theta}}$'s.

The eigenvalues of the matrices are $(\ii,-\ii)$.  
The most general  orthogonal matrices that transform all  $D_{
\rho_\mathrm{isum}}(\s_\mathrm{inv})|_{E_{\td\theta}}$'s into signed
permutations must have a form 
\begin{align}
 U =  P V_\mathrm{sd} \begin{pmatrix}
 1& 0\\
 0& 1\\
\end{pmatrix}
.
\end{align}

If $\dim(\rho_\mathrm{isum}) \geq 8$, it is also possible that $D_{
\rho_\mathrm{isum}}(\s_\mathrm{inv})|_{E_{\td\theta}}$'s 
form the following matrix group:
\begin{align}
 \begin{pmatrix}
1& 0 \\ 
0 & 1 \\ 
\end{pmatrix}
,\ \ 
 \begin{pmatrix}
-\frac{1}{2}& -\sqrt{\frac{3}{4}} \\ 
\sqrt{\frac{3}{4}}& -\frac{1}{2} \\ 
\end{pmatrix}
,\ \ 
 \begin{pmatrix}
-\frac{1}{2}& \sqrt{\frac{3}{4}} \\ 
-\sqrt{\frac{3}{4}}& -\frac{1}{2} \\
\end{pmatrix}
\end{align}
This is because the direct sum decomposition of $\rho_\mathrm{isum}$ contains a
dimension-8 irreducible representation whose $\rho(\ft)$ has a 2-dimensional
eigenspace, which gives rise to the such form of $D_{
\rho_\mathrm{isum}}(\s_\mathrm{inv})|_{E_{\td\theta}}$'s .

The eigenvalues of the later two matrices are $\pm (\ee^{\ii 2\pi/3},\ee^{-\ii
2\pi/3})$.  A permutation of two elements can only have orders 1 or 2. The
corresponding $2 \times 2$ signed permutation matrix  can only have eigenvalues
$1$, $-1$ or $\pm \ii$. Any other eigenvalue is not possible.  Thus, there is
no orthogonal matrix that can transform the above two matrices into signed
permutation.  Such $\rho_\mathrm{isum}$ is not a representation of any modular
data.

\subsubsection{Within a 3-dimensional eigenspace of $\rho_\mathrm{isum}(\ft)$
for rank $\leq 6$}

There is only one such case for rank $\leq 6$.
The $3\times 3$ matrix group $MG$ generated by
$D_{
\rho_\mathrm{isum}}(\s_\mathrm{inv})|_{E_{\td\theta}}$'s
is given by
\begin{align}
\label{MG3x3}
 MG &=\Big\{  
\begin{pmatrix} 1 &0 & 0\\ 0 & 1 & 0\\ 0&0&1\\ \end{pmatrix}, 
\begin{pmatrix} -1 &0 & 0\\ 0 & -1 & 0\\ 0&0&1\\ \end{pmatrix}, 
\begin{pmatrix} 1 &0 & 0\\ 0 & -1 & 0\\ 0&0&1\\ \end{pmatrix}, 
\begin{pmatrix} -1 &0 & 0\\ 0 & 1 & 0\\ 0&0&1\\ \end{pmatrix} 
\Big\}, 
&&\text{ for } \dim(\rho_\mathrm{isum}) = 6.
\end{align}
which is resolved. To find the most general orthogonal matrices that transform the above $3\times 3$
matrices in $MG$ into signed permutation matrices, we first show
\begin{thm}
\label{PPsgn}
If $P$ is a permutation matrix with $P^2 = \id$, then $P$ is a direct sum of
$2\times 2$ and $1\times 1$ matrices.  If $P_\mathrm{sgn}$ is a signed
permutation matrix with $P_\mathrm{sgn}^2 = \id$, then $P_\mathrm{sgn}$ is
a direct sum of $2\times 2$ and $1\times 1$ matrices.
\end{thm}
\begin{proof}[Proof of Theorem \ref{PPsgn}] 
If $P$ is a permutation matrix with $P^2 = \id$, then $P$ must be a pair-wise
permutation, and thus $P$ is a direct sum of $2\times 2$ and $1\times 1$
matrices.  The reduction from signed permutation matrix to permutation matrix
by ignoring the signs is homomorphism of the matrix product.  If
$P_\mathrm{sgn}$ is a signed permutation matrix with $P_\mathrm{sgn}^2 = \id$, then
its reduction gives rise to a permutation matrix $P$ with $P^2 = \id$.  Since
$P$ is a direct sum of $2\times 2$ and $1\times 1$ matrices, $P_\mathrm{sgn}$
is also a direct sum of $2\times 2$ and $1\times 1$ matrices.
\end{proof}
Using the above result, similarly, we can show that
the most general orthogonal matrices that transform all  $D_{
\rho_\mathrm{isum}}(\s_\mathrm{inv})|_{E_{\td\theta}}$'s into signed
permutations must have a form 
\begin{align}
\label{U3x3a}
 U &=  
\frac{PV_\mathrm{sd}}{\sqrt 2}  
\begin{pmatrix}
 1& 1 & 0\\
 1&-1 & 0\\
 0& 0 & 1\\
\end{pmatrix} ,
 \text{ or } \
 U =   \frac{PV_\mathrm{sd}}{\sqrt 2} 
 \begin{pmatrix}
-1& 1 & 0\\
 1& 1 & 0\\
 0& 0 & 1\\
\end{pmatrix} ,
\nonumber\\
\text{or }\  U &=  
\frac{PV_\mathrm{sd}}{\sqrt 2}  
\begin{pmatrix}
 1& 0& 1 \\
 0& 1& 0\\
 1& 0&-1 \\
\end{pmatrix} ,
 \text{ or } \
 U =   \frac{PV_\mathrm{sd}}{\sqrt 2} 
\begin{pmatrix}
-1& 0& 1 \\
 0& 1& 0\\
 1& 0& 1 \\
\end{pmatrix} ,
\nonumber\\
\text{or }\  U &=  
\frac{PV_\mathrm{sd}}{\sqrt 2}  
\begin{pmatrix}
 1& 0& 0 \\
 0& 1& 1\\
 0& 1&-1 \\
\end{pmatrix} ,
 \text{ or } \
 U =   \frac{PV_\mathrm{sd}}{\sqrt 2} 
\begin{pmatrix}
 1& 0& 0 \\
 0&-1& 1\\
 0& 1& 1 \\
\end{pmatrix} ,
\nonumber\\
\text{or }\  U &=  
\frac{PV_\mathrm{sd}}{\sqrt 2}  
\begin{pmatrix}
 1& 0& 0 \\
 0& 1& 0\\
 0& 0& 1 \\
\end{pmatrix} .
\end{align}
where $V_\mathrm{sd}$ are signed diagonal matrices, and $P$ are permutation
matrices.  We note that the non-trivial part of $U$ is a $2\times 2$ block for
index $(1,2)$, $(1,3)$, and $(2,3)$.  The $2\times 2$ block has three
possibilities given in \eqref{U2x2}.  
Such $U$'s transform the diagonal matrices in $MG$ into a direct sum of a $2\times 2$ and an $1\times 1$
matrices. This is a general pattern that apply for
all resolved diagonal matrix group $MG$ generated by $D_{
\rho_\mathrm{isum}}(\s_\mathrm{inv})|_{E_{\td\theta}}$.

The above are all the possibilities that can appear in resolved dimension-6
representations. In the following, we will consider more possibilities, that
appear only for resolved representations with dimension larger than 6.

\subsection{List of $S,T$ matrices from resolved representations}
\label{SsL}

We have constructed a list of irrep-sum symmetric representations
(see Appendix \ref{repL}) that include all the representations of modular data.
Among them, we can select a sublist of resolved
symmetric representations, denoted as $\{ \rho_\mathrm{res} \}$.
We then use the orthogonal matrix $U$ constructed above (see \eqref{U1x1},
\eqref{U2x2} and \eqref{U3x3a})
%, and \eqref{U4x4}) 
to transform the resolved symmetric
representations $ \rho_\mathrm{res}$ to representations, $\rho$'s:
\begin{align}
 \rho(\ft) = U\rho_\mathrm{res}(\ft) U^\top , \ \ \  
 \rho(\fs) = U\rho_\mathrm{res}(\fs) U^\top .
\end{align}
such that the corresponding $D_{\rho}(\s)$ are either zero or signed
permutation in each eigenspace of $\rho(\ft)$.  Since the number of such
representations is finite, we can examine all resulting representations one by
one.

For each $U$, the resulting representation $\rho$ should satisfy Proposition
\ref{p:MDcond}.  In particular, we examine all possible choices of index $u$
that may correspond to the unit object, to see if $\rho$ satisfy the condition
\eqref{SSN}.  If no choices of $u$ can satisfy \eqref{SSN}, then the
representations $\rho$ is rejected.
If some $u$'s satisfy \eqref{SSN}, then for each $u$, we can
construct $S,T$ matrices via \eqref{STrho1}.
We then check if the resulting $S,T$ matrices satisfy the conditions of modular
data summarized in Proposition \ref{p:MDcond} (for details, see Supplementary 
Material Section \ref{Section2}).

In the following, we list all the pairs of $S,T$ matrices that satisfy the
conditions in Proposition \ref{p:MDcond}, and come from the dimension-6 resolved
$\SL$ representations listed in Appendix \ref{repL}.  The list includes all the
unitary and non-unitary modular data with $D^2 \notin \BZ$ from resolved $\SL$
representations.

Each pair of $S,T$ matrices is indexed by $(r_1, r_2, \cdots;
l_1,l_2,\cdots)_k^a$, such as $(3, 3; 5,4)_2^1$.  The first part of index,
$(3, 3; 5,4)$ = (dims;levels), is the index of GT orbit listed in Appendix
\ref{repL}, indicating that the $S,T$ matrices arise from a particular $\SL$
representation in the GT orbit.  This particular $\SL$ representation may
produce several $S,T$ matrices, indexed by the subscript $k$.  We label those
$S,T$ matrices by $(r_1, r_2, \cdots; l_1,l_2,\cdots)_k$.  The $S,T$
matrices from other $\SL$ representations in the GT orbit can obtained by the
$S,T$ matrices labeled by $(r_1, r_2, \cdots; l_1,l_2,\cdots)_k$ via
Galois conjugation $\s_a: \ee^{\ii 2\pi/\ord(T)} \to \ee^{a \ii 2\pi/\ord(T)}$.
The resulting $S,T$ matrices are labeled by $(r_1, r_2, \cdots;
l_1,l_2,\cdots)_k^a$.  Those $a$-indexed $S,T$ matrices  form a Galois orbit.

In the list, $T$ is presented in terms of topological spin
$(s_{1},s_{2},\cdots)$ with $s_{i} = \arg(T_{ii})$.
$S$ is presented as $(S_{00},S_{01}, S_{02}, S_{03}, \cdots;\ S_{11}, S_{12},
S_{13},  \cdots)$.  $d_i = S_{0i}$ are the quantum dimensions.

In the following list, the $S,T$ matrices are grouped into Galois orbits.  To
save space, in general, we only list the $S,T$ matrices whose quantum
dimensions are all positive.  There are 42 of them.  

However, there are Galois orbits that have no $S,T$ matrices whose quantum
dimensions are all positive.  In this case, very often, Galois orbits contain a
pair of $S,T$ matrices that is \textbf{pseudo-unitary}, {\it i.e.} the $S,T$
matrices can be transformed, via a signed diagonal matrix, into a new pair of
$S,T$ matrices whose quantum dimensions are all positive.  The $S,T$ matrices
is such Galois orbits are not listed.

There is only one Galois orbit which contains no unitary nor pseudo-unitary
$S,T$ matrices.  The following list also contain one pair of representative
$S,T$ matrices for such a Galois orbit (the 9$^\text{th}$ in the list).

Our calculation actually produces 174 pairs of $S,T$ matrices, which are given
in Supplementary Material Section \ref{Section3}.  All those 174 pairs of $S,T$
matrices can be obtained form the 43 pairs of $S,T$ matrices in the following
list, via Galois conjugations and conjugations by signed diagonal matrices.
The following 43 pairs of $S,T$ matrices include all with positive quantum
dimensions. The grey entries are the Galois conjugations of the previous black
entry.

%We also list a particular pair of $S,T$ matrices with negative quantum
%dimensions.  The rest of $174-43$ pairs of $S,T$ matrices can be obtained by
%Galois conjugations $\s \in \Gal(\BQ_{\ord(T)})$ of those 43 pairs of $S,T$
%matrices.

\

\noindent1. ind = $(3,
3;5,
4
)_{1}^{1}$:\ \ 
$d_i$ = ($1.0$,
$1.0$,
$2.0$,
$2.0$,
$2.236$,
$2.236$) 

\vskip 0.7ex
\hangindent=3em \hangafter=1
$D^2=$ 20.0 = 
 $20$

\vskip 0.7ex
\hangindent=3em \hangafter=1
$T = ( 0,
0,
\frac{1}{5},
\frac{4}{5},
\frac{1}{4},
\frac{3}{4} )
$,

\vskip 0.7ex
\hangindent=3em \hangafter=1
$S$ = ($ 1$,
$ 1$,
$ 2$,
$ 2$,
$ \sqrt{5}$,
$ \sqrt{5}$;\ \ 
$ 1$,
$ 2$,
$ 2$,
$ -\sqrt{5}$,
$ -\sqrt{5}$;\ \ 
$ -1-\sqrt{5}$,
$ -1+\sqrt{5}$,
$0$,
$0$;\ \ 
$ -1-\sqrt{5}$,
$0$,
$0$;\ \ 
$ -\sqrt{5}$,
$ \sqrt{5}$;\ \ 
$ -\sqrt{5}$)

\vskip 1ex 

 \color{black} \vskip 2ex

\noindent2. ind = $(3,
3;5,
4
)_{2}^{1}$:\ \ 
$d_i$ = ($1.0$,
$1.0$,
$2.0$,
$2.0$,
$2.236$,
$2.236$) 

\vskip 0.7ex
\hangindent=3em \hangafter=1
$D^2=$ 20.0 = 
 $20$

\vskip 0.7ex
\hangindent=3em \hangafter=1
$T = ( 0,
0,
\frac{2}{5},
\frac{3}{5},
\frac{1}{4},
\frac{3}{4} )
$,

\vskip 0.7ex
\hangindent=3em \hangafter=1
$S$ = ($ 1$,
$ 1$,
$ 2$,
$ 2$,
$ \sqrt{5}$,
$ \sqrt{5}$;\ \ 
$ 1$,
$ 2$,
$ 2$,
$ -\sqrt{5}$,
$ -\sqrt{5}$;\ \ 
$ -1+\sqrt{5}$,
$ -1-\sqrt{5}$,
$0$,
$0$;\ \ 
$ -1+\sqrt{5}$,
$0$,
$0$;\ \ 
$ \sqrt{5}$,
$ -\sqrt{5}$;\ \ 
$ \sqrt{5}$)

\vskip 1ex 

 \color{black} \vskip 2ex

\noindent3. ind = $(4,
2;15,
5
)_{1}^{1}$:\ \ 
$d_i$ = ($1.0$,
$1.0$,
$1.0$,
$1.618$,
$1.618$,
$1.618$) 

\vskip 0.7ex
\hangindent=3em \hangafter=1
$D^2=$ 10.854 = 
 $\frac{15+3\sqrt{5}}{2}$

\vskip 0.7ex
\hangindent=3em \hangafter=1
$T = ( 0,
\frac{1}{3},
\frac{1}{3},
\frac{2}{5},
\frac{11}{15},
\frac{11}{15} )
$,

\vskip 0.7ex
\hangindent=3em \hangafter=1
$S$ = ($ 1$,
$ 1$,
$ 1$,
$ \frac{1+\sqrt{5}}{2}$,
$ \frac{1+\sqrt{5}}{2}$,
$ \frac{1+\sqrt{5}}{2}$;\ \ 
$ \zeta_{3}^{1}$,
$ -\zeta_{6}^{1}$,
$ \frac{1+\sqrt{5}}{2}$,
$ \frac{1+\sqrt{5}}{2}\zeta_{3}^{1}$,
$ -\frac{1+\sqrt{5}}{2}\zeta_{6}^{1}$;\ \ 
$ \zeta_{3}^{1}$,
$ \frac{1+\sqrt{5}}{2}$,
$ -\frac{1+\sqrt{5}}{2}\zeta_{6}^{1}$,
$ \frac{1+\sqrt{5}}{2}\zeta_{3}^{1}$;\ \ 
$ -1$,
$ -1$,
$ -1$;\ \ 
$ -\zeta_{3}^{1}$,
$ \zeta_{6}^{1}$;\ \ 
$ -\zeta_{3}^{1}$)

\vskip 1ex 
\color{grey}

\noindent4. ind = $(4,
2;15,
5
)_{1}^{4}$:\ \ 
$d_i$ = ($1.0$,
$1.0$,
$1.0$,
$1.618$,
$1.618$,
$1.618$) 

\vskip 0.7ex
\hangindent=3em \hangafter=1
$D^2=$ 10.854 = 
 $\frac{15+3\sqrt{5}}{2}$

\vskip 0.7ex
\hangindent=3em \hangafter=1
$T = ( 0,
\frac{1}{3},
\frac{1}{3},
\frac{3}{5},
\frac{14}{15},
\frac{14}{15} )
$,

\vskip 0.7ex
\hangindent=3em \hangafter=1
$S$ = ($ 1$,
$ 1$,
$ 1$,
$ \frac{1+\sqrt{5}}{2}$,
$ \frac{1+\sqrt{5}}{2}$,
$ \frac{1+\sqrt{5}}{2}$;\ \ 
$ \zeta_{3}^{1}$,
$ -\zeta_{6}^{1}$,
$ \frac{1+\sqrt{5}}{2}$,
$ \frac{1+\sqrt{5}}{2}\zeta_{3}^{1}$,
$ -\frac{1+\sqrt{5}}{2}\zeta_{6}^{1}$;\ \ 
$ \zeta_{3}^{1}$,
$ \frac{1+\sqrt{5}}{2}$,
$ -\frac{1+\sqrt{5}}{2}\zeta_{6}^{1}$,
$ \frac{1+\sqrt{5}}{2}\zeta_{3}^{1}$;\ \ 
$ -1$,
$ -1$,
$ -1$;\ \ 
$ -\zeta_{3}^{1}$,
$ \zeta_{6}^{1}$;\ \ 
$ -\zeta_{3}^{1}$)

\vskip 1ex 
\color{grey}

\noindent5. ind = $(4,
2;15,
5
)_{1}^{11}$:\ \ 
$d_i$ = ($1.0$,
$1.0$,
$1.0$,
$1.618$,
$1.618$,
$1.618$) 

\vskip 0.7ex
\hangindent=3em \hangafter=1
$D^2=$ 10.854 = 
 $\frac{15+3\sqrt{5}}{2}$

\vskip 0.7ex
\hangindent=3em \hangafter=1
$T = ( 0,
\frac{2}{3},
\frac{2}{3},
\frac{2}{5},
\frac{1}{15},
\frac{1}{15} )
$,

\vskip 0.7ex
\hangindent=3em \hangafter=1
$S$ = ($ 1$,
$ 1$,
$ 1$,
$ \frac{1+\sqrt{5}}{2}$,
$ \frac{1+\sqrt{5}}{2}$,
$ \frac{1+\sqrt{5}}{2}$;\ \ 
$ -\zeta_{6}^{1}$,
$ \zeta_{3}^{1}$,
$ \frac{1+\sqrt{5}}{2}$,
$ -\frac{1+\sqrt{5}}{2}\zeta_{6}^{1}$,
$ \frac{1+\sqrt{5}}{2}\zeta_{3}^{1}$;\ \ 
$ -\zeta_{6}^{1}$,
$ \frac{1+\sqrt{5}}{2}$,
$ \frac{1+\sqrt{5}}{2}\zeta_{3}^{1}$,
$ -\frac{1+\sqrt{5}}{2}\zeta_{6}^{1}$;\ \ 
$ -1$,
$ -1$,
$ -1$;\ \ 
$ \zeta_{6}^{1}$,
$ -\zeta_{3}^{1}$;\ \ 
$ \zeta_{6}^{1}$)

\vskip 1ex 
\color{grey}

\noindent6. ind = $(4,
2;15,
5
)_{1}^{14}$:\ \ 
$d_i$ = ($1.0$,
$1.0$,
$1.0$,
$1.618$,
$1.618$,
$1.618$) 

\vskip 0.7ex
\hangindent=3em \hangafter=1
$D^2=$ 10.854 = 
 $\frac{15+3\sqrt{5}}{2}$

\vskip 0.7ex
\hangindent=3em \hangafter=1
$T = ( 0,
\frac{2}{3},
\frac{2}{3},
\frac{3}{5},
\frac{4}{15},
\frac{4}{15} )
$,

\vskip 0.7ex
\hangindent=3em \hangafter=1
$S$ = ($ 1$,
$ 1$,
$ 1$,
$ \frac{1+\sqrt{5}}{2}$,
$ \frac{1+\sqrt{5}}{2}$,
$ \frac{1+\sqrt{5}}{2}$;\ \ 
$ -\zeta_{6}^{1}$,
$ \zeta_{3}^{1}$,
$ \frac{1+\sqrt{5}}{2}$,
$ -\frac{1+\sqrt{5}}{2}\zeta_{6}^{1}$,
$ \frac{1+\sqrt{5}}{2}\zeta_{3}^{1}$;\ \ 
$ -\zeta_{6}^{1}$,
$ \frac{1+\sqrt{5}}{2}$,
$ \frac{1+\sqrt{5}}{2}\zeta_{3}^{1}$,
$ -\frac{1+\sqrt{5}}{2}\zeta_{6}^{1}$;\ \ 
$ -1$,
$ -1$,
$ -1$;\ \ 
$ \zeta_{6}^{1}$,
$ -\zeta_{3}^{1}$;\ \ 
$ \zeta_{6}^{1}$)

\vskip 1ex 

 \color{black} \vskip 2ex

\noindent7. ind = $(4,
2;7,
3
)_{1}^{1}$:\ \ 
$d_i$ = ($1.0$,
$3.791$,
$3.791$,
$3.791$,
$4.791$,
$5.791$) 

\vskip 0.7ex
\hangindent=3em \hangafter=1
$D^2=$ 100.617 = 
 $\frac{105+21\sqrt{21}}{2}$

\vskip 0.7ex
\hangindent=3em \hangafter=1
$T = ( 0,
\frac{1}{7},
\frac{2}{7},
\frac{4}{7},
0,
\frac{2}{3} )
$,

\vskip 0.7ex
\hangindent=3em \hangafter=1
$S$ = ($ 1$,
$ \frac{3+\sqrt{21}}{2}$,
$ \frac{3+\sqrt{21}}{2}$,
$ \frac{3+\sqrt{21}}{2}$,
$ \frac{5+\sqrt{21}}{2}$,
$ \frac{7+\sqrt{21}}{2}$;\ \ 
$ 2-c^{1}_{21}
-2c^{2}_{21}
+3c^{3}_{21}
+2c^{4}_{21}
-2c^{5}_{21}
$,
$ -c^{2}_{21}
-2c^{3}_{21}
-c^{4}_{21}
+c^{5}_{21}
$,
$ -1+2c^{1}_{21}
+3c^{2}_{21}
-c^{3}_{21}
+2c^{5}_{21}
$,
$ -\frac{3+\sqrt{21}}{2}$,
$0$;\ \ 
$ -1+2c^{1}_{21}
+3c^{2}_{21}
-c^{3}_{21}
+2c^{5}_{21}
$,
$ 2-c^{1}_{21}
-2c^{2}_{21}
+3c^{3}_{21}
+2c^{4}_{21}
-2c^{5}_{21}
$,
$ -\frac{3+\sqrt{21}}{2}$,
$0$;\ \ 
$ -c^{2}_{21}
-2c^{3}_{21}
-c^{4}_{21}
+c^{5}_{21}
$,
$ -\frac{3+\sqrt{21}}{2}$,
$0$;\ \ 
$ 1$,
$ \frac{7+\sqrt{21}}{2}$;\ \ 
$ -\frac{7+\sqrt{21}}{2}$)

\vskip 1ex 
\color{grey}

\noindent8. ind = $(4,
2;7,
3
)_{1}^{5}$:\ \ 
$d_i$ = ($1.0$,
$3.791$,
$3.791$,
$3.791$,
$4.791$,
$5.791$) 

\vskip 0.7ex
\hangindent=3em \hangafter=1
$D^2=$ 100.617 = 
 $\frac{105+21\sqrt{21}}{2}$

\vskip 0.7ex
\hangindent=3em \hangafter=1
$T = ( 0,
\frac{3}{7},
\frac{5}{7},
\frac{6}{7},
0,
\frac{1}{3} )
$,

\vskip 0.7ex
\hangindent=3em \hangafter=1
$S$ = ($ 1$,
$ \frac{3+\sqrt{21}}{2}$,
$ \frac{3+\sqrt{21}}{2}$,
$ \frac{3+\sqrt{21}}{2}$,
$ \frac{5+\sqrt{21}}{2}$,
$ \frac{7+\sqrt{21}}{2}$;\ \ 
$ -c^{2}_{21}
-2c^{3}_{21}
-c^{4}_{21}
+c^{5}_{21}
$,
$ 2-c^{1}_{21}
-2c^{2}_{21}
+3c^{3}_{21}
+2c^{4}_{21}
-2c^{5}_{21}
$,
$ -1+2c^{1}_{21}
+3c^{2}_{21}
-c^{3}_{21}
+2c^{5}_{21}
$,
$ -\frac{3+\sqrt{21}}{2}$,
$0$;\ \ 
$ -1+2c^{1}_{21}
+3c^{2}_{21}
-c^{3}_{21}
+2c^{5}_{21}
$,
$ -c^{2}_{21}
-2c^{3}_{21}
-c^{4}_{21}
+c^{5}_{21}
$,
$ -\frac{3+\sqrt{21}}{2}$,
$0$;\ \ 
$ 2-c^{1}_{21}
-2c^{2}_{21}
+3c^{3}_{21}
+2c^{4}_{21}
-2c^{5}_{21}
$,
$ -\frac{3+\sqrt{21}}{2}$,
$0$;\ \ 
$ 1$,
$ \frac{7+\sqrt{21}}{2}$;\ \ 
$ -\frac{7+\sqrt{21}}{2}$)

\vskip 1ex 

 \color{black} \vskip 2ex 
\color{blue}

\noindent9. ind = $(6;9
)_{1}^{1}$:\ \ 
$d_i$ = ($1.0$,
$0.347$,
$1.0$,
$1.532$,
$-1.0$,
$-1.879$) 

\vskip 0.7ex
\hangindent=3em \hangafter=1
$D^2=$ 9.0 = 
 $9$

\vskip 0.7ex
\hangindent=3em \hangafter=1
$T = ( 0,
\frac{1}{9},
\frac{2}{3},
\frac{4}{9},
\frac{1}{3},
\frac{7}{9} )
$,

\vskip 0.7ex
\hangindent=3em \hangafter=1
$S$ = ($ 1$,
$ c^{2}_{9}
$,
$ 1$,
$ c^{1}_{9}
$,
$ -1$,
$  c_9^4 $;\ \ 
$ 1$,
$ c^{1}_{9}
$,
$ 1$,
$  -c_9^4 $,
$ 1$;\ \ 
$ 1$,
$  c_9^4 $,
$ -1$,
$ c^{2}_{9}
$;\ \ 
$ 1$,
$ -c^{2}_{9}
$,
$ 1$;\ \ 
$ 1$,
$ -c^{1}_{9}
$;\ \ 
$ 1$)

\vskip 1ex 

 \color{black} \vskip 2ex

\noindent10. ind = $(6;13
)_{1}^{1}$:\ \ 
$d_i$ = ($1.0$,
$1.941$,
$2.770$,
$3.438$,
$3.907$,
$4.148$) 

\vskip 0.7ex
\hangindent=3em \hangafter=1
$D^2=$ 56.746 = 
 $21+15c^{1}_{13}
+10c^{2}_{13}
+6c^{3}_{13}
+3c^{4}_{13}
+c^{5}_{13}
$

\vskip 0.7ex
\hangindent=3em \hangafter=1
$T = ( 0,
\frac{4}{13},
\frac{2}{13},
\frac{7}{13},
\frac{6}{13},
\frac{12}{13} )
$,

\vskip 0.7ex
\hangindent=3em \hangafter=1
$S$ = ($ 1$,
$ \xi_{13}^{2}$,
$ \xi_{13}^{3}$,
$ \xi_{13}^{4}$,
$ \xi_{13}^{5}$,
$ \xi_{13}^{6}$;\ \ 
$ -\xi_{13}^{4}$,
$ \xi_{13}^{6}$,
$ -\xi_{13}^{5}$,
$ \xi_{13}^{3}$,
$ -1$;\ \ 
$ \xi_{13}^{4}$,
$ 1$,
$ -\xi_{13}^{2}$,
$ -\xi_{13}^{5}$;\ \ 
$ \xi_{13}^{3}$,
$ -\xi_{13}^{6}$,
$ \xi_{13}^{2}$;\ \ 
$ -1$,
$ \xi_{13}^{4}$;\ \ 
$ -\xi_{13}^{3}$)

\vskip 1ex 
\color{grey}

\noindent11. ind = $(6;13
)_{1}^{12}$:\ \ 
$d_i$ = ($1.0$,
$1.941$,
$2.770$,
$3.438$,
$3.907$,
$4.148$) 

\vskip 0.7ex
\hangindent=3em \hangafter=1
$D^2=$ 56.746 = 
 $21+15c^{1}_{13}
+10c^{2}_{13}
+6c^{3}_{13}
+3c^{4}_{13}
+c^{5}_{13}
$

\vskip 0.7ex
\hangindent=3em \hangafter=1
$T = ( 0,
\frac{9}{13},
\frac{11}{13},
\frac{6}{13},
\frac{7}{13},
\frac{1}{13} )
$,

\vskip 0.7ex
\hangindent=3em \hangafter=1
$S$ = ($ 1$,
$ \xi_{13}^{2}$,
$ \xi_{13}^{3}$,
$ \xi_{13}^{4}$,
$ \xi_{13}^{5}$,
$ \xi_{13}^{6}$;\ \ 
$ -\xi_{13}^{4}$,
$ \xi_{13}^{6}$,
$ -\xi_{13}^{5}$,
$ \xi_{13}^{3}$,
$ -1$;\ \ 
$ \xi_{13}^{4}$,
$ 1$,
$ -\xi_{13}^{2}$,
$ -\xi_{13}^{5}$;\ \ 
$ \xi_{13}^{3}$,
$ -\xi_{13}^{6}$,
$ \xi_{13}^{2}$;\ \ 
$ -1$,
$ \xi_{13}^{4}$;\ \ 
$ -\xi_{13}^{3}$)

\vskip 1ex 

 \color{black} \vskip 2ex

\noindent12. ind = $(6;16
)_{1}^{1}$:\ \ 
$d_i$ = ($1.0$,
$1.0$,
$1.0$,
$1.0$,
$1.414$,
$1.414$) 

\vskip 0.7ex
\hangindent=3em \hangafter=1
$D^2=$ 8.0 = 
 $8$

\vskip 0.7ex
\hangindent=3em \hangafter=1
$T = ( 0,
\frac{1}{2},
\frac{1}{4},
\frac{3}{4},
\frac{1}{16},
\frac{5}{16} )
$,

\vskip 0.7ex
\hangindent=3em \hangafter=1
$S$ = ($ 1$,
$ 1$,
$ 1$,
$ 1$,
$ \sqrt{2}$,
$ \sqrt{2}$;\ \ 
$ 1$,
$ 1$,
$ 1$,
$ -\sqrt{2}$,
$ -\sqrt{2}$;\ \ 
$ -1$,
$ -1$,
$ \sqrt{2}$,
$ -\sqrt{2}$;\ \ 
$ -1$,
$ -\sqrt{2}$,
$ \sqrt{2}$;\ \ 
$0$,
$0$;\ \ 
$0$)

\vskip 1ex 
\color{grey}

\noindent13. ind = $(6;16
)_{1}^{7}$:\ \ 
$d_i$ = ($1.0$,
$1.0$,
$1.0$,
$1.0$,
$1.414$,
$1.414$) 

\vskip 0.7ex
\hangindent=3em \hangafter=1
$D^2=$ 8.0 = 
 $8$

\vskip 0.7ex
\hangindent=3em \hangafter=1
$T = ( 0,
\frac{1}{2},
\frac{1}{4},
\frac{3}{4},
\frac{3}{16},
\frac{7}{16} )
$,

\vskip 0.7ex
\hangindent=3em \hangafter=1
$S$ = ($ 1$,
$ 1$,
$ 1$,
$ 1$,
$ \sqrt{2}$,
$ \sqrt{2}$;\ \ 
$ 1$,
$ 1$,
$ 1$,
$ -\sqrt{2}$,
$ -\sqrt{2}$;\ \ 
$ -1$,
$ -1$,
$ \sqrt{2}$,
$ -\sqrt{2}$;\ \ 
$ -1$,
$ -\sqrt{2}$,
$ \sqrt{2}$;\ \ 
$0$,
$0$;\ \ 
$0$)

\vskip 1ex 
\color{grey}

\noindent14. ind = $(6;16
)_{1}^{9}$:\ \ 
$d_i$ = ($1.0$,
$1.0$,
$1.0$,
$1.0$,
$1.414$,
$1.414$) 

\vskip 0.7ex
\hangindent=3em \hangafter=1
$D^2=$ 8.0 = 
 $8$

\vskip 0.7ex
\hangindent=3em \hangafter=1
$T = ( 0,
\frac{1}{2},
\frac{1}{4},
\frac{3}{4},
\frac{9}{16},
\frac{13}{16} )
$,

\vskip 0.7ex
\hangindent=3em \hangafter=1
$S$ = ($ 1$,
$ 1$,
$ 1$,
$ 1$,
$ \sqrt{2}$,
$ \sqrt{2}$;\ \ 
$ 1$,
$ 1$,
$ 1$,
$ -\sqrt{2}$,
$ -\sqrt{2}$;\ \ 
$ -1$,
$ -1$,
$ \sqrt{2}$,
$ -\sqrt{2}$;\ \ 
$ -1$,
$ -\sqrt{2}$,
$ \sqrt{2}$;\ \ 
$0$,
$0$;\ \ 
$0$)

\vskip 1ex 
\color{grey}

\noindent15. ind = $(6;16
)_{1}^{15}$:\ \ 
$d_i$ = ($1.0$,
$1.0$,
$1.0$,
$1.0$,
$1.414$,
$1.414$) 

\vskip 0.7ex
\hangindent=3em \hangafter=1
$D^2=$ 8.0 = 
 $8$

\vskip 0.7ex
\hangindent=3em \hangafter=1
$T = ( 0,
\frac{1}{2},
\frac{1}{4},
\frac{3}{4},
\frac{11}{16},
\frac{15}{16} )
$,

\vskip 0.7ex
\hangindent=3em \hangafter=1
$S$ = ($ 1$,
$ 1$,
$ 1$,
$ 1$,
$ \sqrt{2}$,
$ \sqrt{2}$;\ \ 
$ 1$,
$ 1$,
$ 1$,
$ -\sqrt{2}$,
$ -\sqrt{2}$;\ \ 
$ -1$,
$ -1$,
$ \sqrt{2}$,
$ -\sqrt{2}$;\ \ 
$ -1$,
$ -\sqrt{2}$,
$ \sqrt{2}$;\ \ 
$0$,
$0$;\ \ 
$0$)

\vskip 1ex 

 \color{black} \vskip 2ex

\noindent16. ind = $(6;16
)_{2}^{1}$:\ \ 
$d_i$ = ($1.0$,
$1.0$,
$1.0$,
$1.0$,
$1.414$,
$1.414$) 

\vskip 0.7ex
\hangindent=3em \hangafter=1
$D^2=$ 8.0 = 
 $8$

\vskip 0.7ex
\hangindent=3em \hangafter=1
$T = ( 0,
\frac{1}{2},
\frac{1}{4},
\frac{3}{4},
\frac{1}{16},
\frac{13}{16} )
$,

\vskip 0.7ex
\hangindent=3em \hangafter=1
$S$ = ($ 1$,
$ 1$,
$ 1$,
$ 1$,
$ \sqrt{2}$,
$ \sqrt{2}$;\ \ 
$ 1$,
$ 1$,
$ 1$,
$ -\sqrt{2}$,
$ -\sqrt{2}$;\ \ 
$ -1$,
$ -1$,
$ -\sqrt{2}$,
$ \sqrt{2}$;\ \ 
$ -1$,
$ \sqrt{2}$,
$ -\sqrt{2}$;\ \ 
$0$,
$0$;\ \ 
$0$)

\vskip 1ex 
\color{grey}

\noindent17. ind = $(6;16
)_{2}^{15}$:\ \ 
$d_i$ = ($1.0$,
$1.0$,
$1.0$,
$1.0$,
$1.414$,
$1.414$) 

\vskip 0.7ex
\hangindent=3em \hangafter=1
$D^2=$ 8.0 = 
 $8$

\vskip 0.7ex
\hangindent=3em \hangafter=1
$T = ( 0,
\frac{1}{2},
\frac{1}{4},
\frac{3}{4},
\frac{3}{16},
\frac{15}{16} )
$,

\vskip 0.7ex
\hangindent=3em \hangafter=1
$S$ = ($ 1$,
$ 1$,
$ 1$,
$ 1$,
$ \sqrt{2}$,
$ \sqrt{2}$;\ \ 
$ 1$,
$ 1$,
$ 1$,
$ -\sqrt{2}$,
$ -\sqrt{2}$;\ \ 
$ -1$,
$ -1$,
$ -\sqrt{2}$,
$ \sqrt{2}$;\ \ 
$ -1$,
$ \sqrt{2}$,
$ -\sqrt{2}$;\ \ 
$0$,
$0$;\ \ 
$0$)

\vskip 1ex 
\color{grey}

\noindent18. ind = $(6;16
)_{2}^{9}$:\ \ 
$d_i$ = ($1.0$,
$1.0$,
$1.0$,
$1.0$,
$1.414$,
$1.414$) 

\vskip 0.7ex
\hangindent=3em \hangafter=1
$D^2=$ 8.0 = 
 $8$

\vskip 0.7ex
\hangindent=3em \hangafter=1
$T = ( 0,
\frac{1}{2},
\frac{1}{4},
\frac{3}{4},
\frac{5}{16},
\frac{9}{16} )
$,

\vskip 0.7ex
\hangindent=3em \hangafter=1
$S$ = ($ 1$,
$ 1$,
$ 1$,
$ 1$,
$ \sqrt{2}$,
$ \sqrt{2}$;\ \ 
$ 1$,
$ 1$,
$ 1$,
$ -\sqrt{2}$,
$ -\sqrt{2}$;\ \ 
$ -1$,
$ -1$,
$ \sqrt{2}$,
$ -\sqrt{2}$;\ \ 
$ -1$,
$ -\sqrt{2}$,
$ \sqrt{2}$;\ \ 
$0$,
$0$;\ \ 
$0$)

\vskip 1ex 
\color{grey}

\noindent19. ind = $(6;16
)_{2}^{7}$:\ \ 
$d_i$ = ($1.0$,
$1.0$,
$1.0$,
$1.0$,
$1.414$,
$1.414$) 

\vskip 0.7ex
\hangindent=3em \hangafter=1
$D^2=$ 8.0 = 
 $8$

\vskip 0.7ex
\hangindent=3em \hangafter=1
$T = ( 0,
\frac{1}{2},
\frac{1}{4},
\frac{3}{4},
\frac{7}{16},
\frac{11}{16} )
$,

\vskip 0.7ex
\hangindent=3em \hangafter=1
$S$ = ($ 1$,
$ 1$,
$ 1$,
$ 1$,
$ \sqrt{2}$,
$ \sqrt{2}$;\ \ 
$ 1$,
$ 1$,
$ 1$,
$ -\sqrt{2}$,
$ -\sqrt{2}$;\ \ 
$ -1$,
$ -1$,
$ \sqrt{2}$,
$ -\sqrt{2}$;\ \ 
$ -1$,
$ -\sqrt{2}$,
$ \sqrt{2}$;\ \ 
$0$,
$0$;\ \ 
$0$)

\vskip 1ex 

 \color{black} \vskip 2ex

\noindent20. ind = $(6;35
)_{1}^{1}$:\ \ 
$d_i$ = ($1.0$,
$1.618$,
$1.801$,
$2.246$,
$2.915$,
$3.635$) 

\vskip 0.7ex
\hangindent=3em \hangafter=1
$D^2=$ 33.632 = 
 $15+3c^{1}_{35}
+2c^{4}_{35}
+6c^{5}_{35}
+3c^{6}_{35}
+3c^{7}_{35}
+2c^{10}_{35}
+2c^{11}_{35}
$

\vskip 0.7ex
\hangindent=3em \hangafter=1
$T = ( 0,
\frac{2}{5},
\frac{1}{7},
\frac{5}{7},
\frac{19}{35},
\frac{4}{35} )
$,

\vskip 0.7ex
\hangindent=3em \hangafter=1
$S$ = ($ 1$,
$ \frac{1+\sqrt{5}}{2}$,
$ \xi_{7}^{2}$,
$ \xi_{7}^{3}$,
$ c^{1}_{35}
+c^{6}_{35}
$,
$ c^{1}_{35}
+c^{4}_{35}
+c^{6}_{35}
+c^{11}_{35}
$;\ \ 
$ -1$,
$ c^{1}_{35}
+c^{6}_{35}
$,
$ c^{1}_{35}
+c^{4}_{35}
+c^{6}_{35}
+c^{11}_{35}
$,
$ -\xi_{7}^{2}$,
$ -\xi_{7}^{3}$;\ \ 
$ -\xi_{7}^{3}$,
$ 1$,
$ -c^{1}_{35}
-c^{4}_{35}
-c^{6}_{35}
-c^{11}_{35}
$,
$ \frac{1+\sqrt{5}}{2}$;\ \ 
$ -\xi_{7}^{2}$,
$ \frac{1+\sqrt{5}}{2}$,
$ -c^{1}_{35}
-c^{6}_{35}
$;\ \ 
$ \xi_{7}^{3}$,
$ -1$;\ \ 
$ \xi_{7}^{2}$)

\vskip 1ex 
\color{grey}

\noindent21. ind = $(6;35
)_{1}^{6}$:\ \ 
$d_i$ = ($1.0$,
$1.618$,
$1.801$,
$2.246$,
$2.915$,
$3.635$) 

\vskip 0.7ex
\hangindent=3em \hangafter=1
$D^2=$ 33.632 = 
 $15+3c^{1}_{35}
+2c^{4}_{35}
+6c^{5}_{35}
+3c^{6}_{35}
+3c^{7}_{35}
+2c^{10}_{35}
+2c^{11}_{35}
$

\vskip 0.7ex
\hangindent=3em \hangafter=1
$T = ( 0,
\frac{2}{5},
\frac{6}{7},
\frac{2}{7},
\frac{9}{35},
\frac{24}{35} )
$,

\vskip 0.7ex
\hangindent=3em \hangafter=1
$S$ = ($ 1$,
$ \frac{1+\sqrt{5}}{2}$,
$ \xi_{7}^{2}$,
$ \xi_{7}^{3}$,
$ c^{1}_{35}
+c^{6}_{35}
$,
$ c^{1}_{35}
+c^{4}_{35}
+c^{6}_{35}
+c^{11}_{35}
$;\ \ 
$ -1$,
$ c^{1}_{35}
+c^{6}_{35}
$,
$ c^{1}_{35}
+c^{4}_{35}
+c^{6}_{35}
+c^{11}_{35}
$,
$ -\xi_{7}^{2}$,
$ -\xi_{7}^{3}$;\ \ 
$ -\xi_{7}^{3}$,
$ 1$,
$ -c^{1}_{35}
-c^{4}_{35}
-c^{6}_{35}
-c^{11}_{35}
$,
$ \frac{1+\sqrt{5}}{2}$;\ \ 
$ -\xi_{7}^{2}$,
$ \frac{1+\sqrt{5}}{2}$,
$ -c^{1}_{35}
-c^{6}_{35}
$;\ \ 
$ \xi_{7}^{3}$,
$ -1$;\ \ 
$ \xi_{7}^{2}$)

\vskip 1ex 
\color{grey}

\noindent22. ind = $(6;35
)_{1}^{29}$:\ \ 
$d_i$ = ($1.0$,
$1.618$,
$1.801$,
$2.246$,
$2.915$,
$3.635$) 

\vskip 0.7ex
\hangindent=3em \hangafter=1
$D^2=$ 33.632 = 
 $15+3c^{1}_{35}
+2c^{4}_{35}
+6c^{5}_{35}
+3c^{6}_{35}
+3c^{7}_{35}
+2c^{10}_{35}
+2c^{11}_{35}
$

\vskip 0.7ex
\hangindent=3em \hangafter=1
$T = ( 0,
\frac{3}{5},
\frac{1}{7},
\frac{5}{7},
\frac{26}{35},
\frac{11}{35} )
$,

\vskip 0.7ex
\hangindent=3em \hangafter=1
$S$ = ($ 1$,
$ \frac{1+\sqrt{5}}{2}$,
$ \xi_{7}^{2}$,
$ \xi_{7}^{3}$,
$ c^{1}_{35}
+c^{6}_{35}
$,
$ c^{1}_{35}
+c^{4}_{35}
+c^{6}_{35}
+c^{11}_{35}
$;\ \ 
$ -1$,
$ c^{1}_{35}
+c^{6}_{35}
$,
$ c^{1}_{35}
+c^{4}_{35}
+c^{6}_{35}
+c^{11}_{35}
$,
$ -\xi_{7}^{2}$,
$ -\xi_{7}^{3}$;\ \ 
$ -\xi_{7}^{3}$,
$ 1$,
$ -c^{1}_{35}
-c^{4}_{35}
-c^{6}_{35}
-c^{11}_{35}
$,
$ \frac{1+\sqrt{5}}{2}$;\ \ 
$ -\xi_{7}^{2}$,
$ \frac{1+\sqrt{5}}{2}$,
$ -c^{1}_{35}
-c^{6}_{35}
$;\ \ 
$ \xi_{7}^{3}$,
$ -1$;\ \ 
$ \xi_{7}^{2}$)

\vskip 1ex 
\color{grey}

\noindent23. ind = $(6;35
)_{1}^{34}$:\ \ 
$d_i$ = ($1.0$,
$1.618$,
$1.801$,
$2.246$,
$2.915$,
$3.635$) 

\vskip 0.7ex
\hangindent=3em \hangafter=1
$D^2=$ 33.632 = 
 $15+3c^{1}_{35}
+2c^{4}_{35}
+6c^{5}_{35}
+3c^{6}_{35}
+3c^{7}_{35}
+2c^{10}_{35}
+2c^{11}_{35}
$

\vskip 0.7ex
\hangindent=3em \hangafter=1
$T = ( 0,
\frac{3}{5},
\frac{6}{7},
\frac{2}{7},
\frac{16}{35},
\frac{31}{35} )
$,

\vskip 0.7ex
\hangindent=3em \hangafter=1
$S$ = ($ 1$,
$ \frac{1+\sqrt{5}}{2}$,
$ \xi_{7}^{2}$,
$ \xi_{7}^{3}$,
$ c^{1}_{35}
+c^{6}_{35}
$,
$ c^{1}_{35}
+c^{4}_{35}
+c^{6}_{35}
+c^{11}_{35}
$;\ \ 
$ -1$,
$ c^{1}_{35}
+c^{6}_{35}
$,
$ c^{1}_{35}
+c^{4}_{35}
+c^{6}_{35}
+c^{11}_{35}
$,
$ -\xi_{7}^{2}$,
$ -\xi_{7}^{3}$;\ \ 
$ -\xi_{7}^{3}$,
$ 1$,
$ -c^{1}_{35}
-c^{4}_{35}
-c^{6}_{35}
-c^{11}_{35}
$,
$ \frac{1+\sqrt{5}}{2}$;\ \ 
$ -\xi_{7}^{2}$,
$ \frac{1+\sqrt{5}}{2}$,
$ -c^{1}_{35}
-c^{6}_{35}
$;\ \ 
$ \xi_{7}^{3}$,
$ -1$;\ \ 
$ \xi_{7}^{2}$)

\vskip 1ex 

 \color{black} \vskip 2ex

\noindent24. ind = $(6;56
)_{1}^{1}$:\ \ 
$d_i$ = ($1.0$,
$1.0$,
$1.801$,
$1.801$,
$2.246$,
$2.246$) 

\vskip 0.7ex
\hangindent=3em \hangafter=1
$D^2=$ 18.591 = 
 $12+6c^{1}_{7}
+2c^{2}_{7}
$

\vskip 0.7ex
\hangindent=3em \hangafter=1
$T = ( 0,
\frac{1}{4},
\frac{1}{7},
\frac{11}{28},
\frac{5}{7},
\frac{27}{28} )
$,

\vskip 0.7ex
\hangindent=3em \hangafter=1
$S$ = ($ 1$,
$ 1$,
$ \xi_{7}^{2}$,
$ \xi_{7}^{2}$,
$ \xi_{7}^{3}$,
$ \xi_{7}^{3}$;\ \ 
$ -1$,
$ \xi_{7}^{2}$,
$ -\xi_{7}^{2}$,
$ \xi_{7}^{3}$,
$ -\xi_{7}^{3}$;\ \ 
$ -\xi_{7}^{3}$,
$ -\xi_{7}^{3}$,
$ 1$,
$ 1$;\ \ 
$ \xi_{7}^{3}$,
$ 1$,
$ -1$;\ \ 
$ -\xi_{7}^{2}$,
$ -\xi_{7}^{2}$;\ \ 
$ \xi_{7}^{2}$)

\vskip 1ex 
\color{grey}

\noindent25. ind = $(6;56
)_{1}^{13}$:\ \ 
$d_i$ = ($1.0$,
$1.0$,
$1.801$,
$1.801$,
$2.246$,
$2.246$) 

\vskip 0.7ex
\hangindent=3em \hangafter=1
$D^2=$ 18.591 = 
 $12+6c^{1}_{7}
+2c^{2}_{7}
$

\vskip 0.7ex
\hangindent=3em \hangafter=1
$T = ( 0,
\frac{1}{4},
\frac{6}{7},
\frac{3}{28},
\frac{2}{7},
\frac{15}{28} )
$,

\vskip 0.7ex
\hangindent=3em \hangafter=1
$S$ = ($ 1$,
$ 1$,
$ \xi_{7}^{2}$,
$ \xi_{7}^{2}$,
$ \xi_{7}^{3}$,
$ \xi_{7}^{3}$;\ \ 
$ -1$,
$ \xi_{7}^{2}$,
$ -\xi_{7}^{2}$,
$ \xi_{7}^{3}$,
$ -\xi_{7}^{3}$;\ \ 
$ -\xi_{7}^{3}$,
$ -\xi_{7}^{3}$,
$ 1$,
$ 1$;\ \ 
$ \xi_{7}^{3}$,
$ 1$,
$ -1$;\ \ 
$ -\xi_{7}^{2}$,
$ -\xi_{7}^{2}$;\ \ 
$ \xi_{7}^{2}$)

\vskip 1ex 
\color{grey}

\noindent26. ind = $(6;56
)_{1}^{15}$:\ \ 
$d_i$ = ($1.0$,
$1.0$,
$1.801$,
$1.801$,
$2.246$,
$2.246$) 

\vskip 0.7ex
\hangindent=3em \hangafter=1
$D^2=$ 18.591 = 
 $12+6c^{1}_{7}
+2c^{2}_{7}
$

\vskip 0.7ex
\hangindent=3em \hangafter=1
$T = ( 0,
\frac{3}{4},
\frac{1}{7},
\frac{25}{28},
\frac{5}{7},
\frac{13}{28} )
$,

\vskip 0.7ex
\hangindent=3em \hangafter=1
$S$ = ($ 1$,
$ 1$,
$ \xi_{7}^{2}$,
$ \xi_{7}^{2}$,
$ \xi_{7}^{3}$,
$ \xi_{7}^{3}$;\ \ 
$ -1$,
$ \xi_{7}^{2}$,
$ -\xi_{7}^{2}$,
$ \xi_{7}^{3}$,
$ -\xi_{7}^{3}$;\ \ 
$ -\xi_{7}^{3}$,
$ -\xi_{7}^{3}$,
$ 1$,
$ 1$;\ \ 
$ \xi_{7}^{3}$,
$ 1$,
$ -1$;\ \ 
$ -\xi_{7}^{2}$,
$ -\xi_{7}^{2}$;\ \ 
$ \xi_{7}^{2}$)

\vskip 1ex 
\color{grey}

\noindent27. ind = $(6;56
)_{1}^{27}$:\ \ 
$d_i$ = ($1.0$,
$1.0$,
$1.801$,
$1.801$,
$2.246$,
$2.246$) 

\vskip 0.7ex
\hangindent=3em \hangafter=1
$D^2=$ 18.591 = 
 $12+6c^{1}_{7}
+2c^{2}_{7}
$

\vskip 0.7ex
\hangindent=3em \hangafter=1
$T = ( 0,
\frac{3}{4},
\frac{6}{7},
\frac{17}{28},
\frac{2}{7},
\frac{1}{28} )
$,

\vskip 0.7ex
\hangindent=3em \hangafter=1
$S$ = ($ 1$,
$ 1$,
$ \xi_{7}^{2}$,
$ \xi_{7}^{2}$,
$ \xi_{7}^{3}$,
$ \xi_{7}^{3}$;\ \ 
$ -1$,
$ \xi_{7}^{2}$,
$ -\xi_{7}^{2}$,
$ \xi_{7}^{3}$,
$ -\xi_{7}^{3}$;\ \ 
$ -\xi_{7}^{3}$,
$ -\xi_{7}^{3}$,
$ 1$,
$ 1$;\ \ 
$ \xi_{7}^{3}$,
$ 1$,
$ -1$;\ \ 
$ -\xi_{7}^{2}$,
$ -\xi_{7}^{2}$;\ \ 
$ \xi_{7}^{2}$)

\vskip 1ex 

 \color{black} \vskip 2ex

\noindent28. ind = $(6;80
)_{1}^{1}$:\ \ 
$d_i$ = ($1.0$,
$1.0$,
$1.414$,
$1.618$,
$1.618$,
$2.288$) 

\vskip 0.7ex
\hangindent=3em \hangafter=1
$D^2=$ 14.472 = 
 $10+2\sqrt{5}$

\vskip 0.7ex
\hangindent=3em \hangafter=1
$T = ( 0,
\frac{1}{2},
\frac{1}{16},
\frac{2}{5},
\frac{9}{10},
\frac{37}{80} )
$,

\vskip 0.7ex
\hangindent=3em \hangafter=1
$S$ = ($ 1$,
$ 1$,
$ \sqrt{2}$,
$ \frac{1+\sqrt{5}}{2}$,
$ \frac{1+\sqrt{5}}{2}$,
$ c^{3}_{40}
+c^{5}_{40}
-c^{7}_{40}
$;\ \ 
$ 1$,
$ -\sqrt{2}$,
$ \frac{1+\sqrt{5}}{2}$,
$ \frac{1+\sqrt{5}}{2}$,
$ -c^{3}_{40}
-c^{5}_{40}
+c^{7}_{40}
$;\ \ 
$0$,
$ c^{3}_{40}
+c^{5}_{40}
-c^{7}_{40}
$,
$ -c^{3}_{40}
-c^{5}_{40}
+c^{7}_{40}
$,
$0$;\ \ 
$ -1$,
$ -1$,
$ -\sqrt{2}$;\ \ 
$ -1$,
$ \sqrt{2}$;\ \ 
$0$)

\vskip 1ex 
\color{grey}

\noindent29. ind = $(6;80
)_{1}^{49}$:\ \ 
$d_i$ = ($1.0$,
$1.0$,
$1.414$,
$1.618$,
$1.618$,
$2.288$) 

\vskip 0.7ex
\hangindent=3em \hangafter=1
$D^2=$ 14.472 = 
 $10+2\sqrt{5}$

\vskip 0.7ex
\hangindent=3em \hangafter=1
$T = ( 0,
\frac{1}{2},
\frac{1}{16},
\frac{3}{5},
\frac{1}{10},
\frac{53}{80} )
$,

\vskip 0.7ex
\hangindent=3em \hangafter=1
$S$ = ($ 1$,
$ 1$,
$ \sqrt{2}$,
$ \frac{1+\sqrt{5}}{2}$,
$ \frac{1+\sqrt{5}}{2}$,
$ c^{3}_{40}
+c^{5}_{40}
-c^{7}_{40}
$;\ \ 
$ 1$,
$ -\sqrt{2}$,
$ \frac{1+\sqrt{5}}{2}$,
$ \frac{1+\sqrt{5}}{2}$,
$ -c^{3}_{40}
-c^{5}_{40}
+c^{7}_{40}
$;\ \ 
$0$,
$ c^{3}_{40}
+c^{5}_{40}
-c^{7}_{40}
$,
$ -c^{3}_{40}
-c^{5}_{40}
+c^{7}_{40}
$,
$0$;\ \ 
$ -1$,
$ -1$,
$ -\sqrt{2}$;\ \ 
$ -1$,
$ \sqrt{2}$;\ \ 
$0$)

\vskip 1ex 
\color{grey}

\noindent30. ind = $(6;80
)_{1}^{71}$:\ \ 
$d_i$ = ($1.0$,
$1.0$,
$1.414$,
$1.618$,
$1.618$,
$2.288$) 

\vskip 0.7ex
\hangindent=3em \hangafter=1
$D^2=$ 14.472 = 
 $10+2\sqrt{5}$

\vskip 0.7ex
\hangindent=3em \hangafter=1
$T = ( 0,
\frac{1}{2},
\frac{7}{16},
\frac{2}{5},
\frac{9}{10},
\frac{67}{80} )
$,

\vskip 0.7ex
\hangindent=3em \hangafter=1
$S$ = ($ 1$,
$ 1$,
$ \sqrt{2}$,
$ \frac{1+\sqrt{5}}{2}$,
$ \frac{1+\sqrt{5}}{2}$,
$ c^{3}_{40}
+c^{5}_{40}
-c^{7}_{40}
$;\ \ 
$ 1$,
$ -\sqrt{2}$,
$ \frac{1+\sqrt{5}}{2}$,
$ \frac{1+\sqrt{5}}{2}$,
$ -c^{3}_{40}
-c^{5}_{40}
+c^{7}_{40}
$;\ \ 
$0$,
$ c^{3}_{40}
+c^{5}_{40}
-c^{7}_{40}
$,
$ -c^{3}_{40}
-c^{5}_{40}
+c^{7}_{40}
$,
$0$;\ \ 
$ -1$,
$ -1$,
$ -\sqrt{2}$;\ \ 
$ -1$,
$ \sqrt{2}$;\ \ 
$0$)

\vskip 1ex 
\color{grey}

\noindent31. ind = $(6;80
)_{1}^{39}$:\ \ 
$d_i$ = ($1.0$,
$1.0$,
$1.414$,
$1.618$,
$1.618$,
$2.288$) 

\vskip 0.7ex
\hangindent=3em \hangafter=1
$D^2=$ 14.472 = 
 $10+2\sqrt{5}$

\vskip 0.7ex
\hangindent=3em \hangafter=1
$T = ( 0,
\frac{1}{2},
\frac{7}{16},
\frac{3}{5},
\frac{1}{10},
\frac{3}{80} )
$,

\vskip 0.7ex
\hangindent=3em \hangafter=1
$S$ = ($ 1$,
$ 1$,
$ \sqrt{2}$,
$ \frac{1+\sqrt{5}}{2}$,
$ \frac{1+\sqrt{5}}{2}$,
$ c^{3}_{40}
+c^{5}_{40}
-c^{7}_{40}
$;\ \ 
$ 1$,
$ -\sqrt{2}$,
$ \frac{1+\sqrt{5}}{2}$,
$ \frac{1+\sqrt{5}}{2}$,
$ -c^{3}_{40}
-c^{5}_{40}
+c^{7}_{40}
$;\ \ 
$0$,
$ c^{3}_{40}
+c^{5}_{40}
-c^{7}_{40}
$,
$ -c^{3}_{40}
-c^{5}_{40}
+c^{7}_{40}
$,
$0$;\ \ 
$ -1$,
$ -1$,
$ -\sqrt{2}$;\ \ 
$ -1$,
$ \sqrt{2}$;\ \ 
$0$)

\vskip 1ex 
\color{grey}

\noindent32. ind = $(6;80
)_{1}^{41}$:\ \ 
$d_i$ = ($1.0$,
$1.0$,
$1.414$,
$1.618$,
$1.618$,
$2.288$) 

\vskip 0.7ex
\hangindent=3em \hangafter=1
$D^2=$ 14.472 = 
 $10+2\sqrt{5}$

\vskip 0.7ex
\hangindent=3em \hangafter=1
$T = ( 0,
\frac{1}{2},
\frac{9}{16},
\frac{2}{5},
\frac{9}{10},
\frac{77}{80} )
$,

\vskip 0.7ex
\hangindent=3em \hangafter=1
$S$ = ($ 1$,
$ 1$,
$ \sqrt{2}$,
$ \frac{1+\sqrt{5}}{2}$,
$ \frac{1+\sqrt{5}}{2}$,
$ c^{3}_{40}
+c^{5}_{40}
-c^{7}_{40}
$;\ \ 
$ 1$,
$ -\sqrt{2}$,
$ \frac{1+\sqrt{5}}{2}$,
$ \frac{1+\sqrt{5}}{2}$,
$ -c^{3}_{40}
-c^{5}_{40}
+c^{7}_{40}
$;\ \ 
$0$,
$ c^{3}_{40}
+c^{5}_{40}
-c^{7}_{40}
$,
$ -c^{3}_{40}
-c^{5}_{40}
+c^{7}_{40}
$,
$0$;\ \ 
$ -1$,
$ -1$,
$ -\sqrt{2}$;\ \ 
$ -1$,
$ \sqrt{2}$;\ \ 
$0$)

\vskip 1ex 
\color{grey}

\noindent33. ind = $(6;80
)_{1}^{9}$:\ \ 
$d_i$ = ($1.0$,
$1.0$,
$1.414$,
$1.618$,
$1.618$,
$2.288$) 

\vskip 0.7ex
\hangindent=3em \hangafter=1
$D^2=$ 14.472 = 
 $10+2\sqrt{5}$

\vskip 0.7ex
\hangindent=3em \hangafter=1
$T = ( 0,
\frac{1}{2},
\frac{9}{16},
\frac{3}{5},
\frac{1}{10},
\frac{13}{80} )
$,

\vskip 0.7ex
\hangindent=3em \hangafter=1
$S$ = ($ 1$,
$ 1$,
$ \sqrt{2}$,
$ \frac{1+\sqrt{5}}{2}$,
$ \frac{1+\sqrt{5}}{2}$,
$ c^{3}_{40}
+c^{5}_{40}
-c^{7}_{40}
$;\ \ 
$ 1$,
$ -\sqrt{2}$,
$ \frac{1+\sqrt{5}}{2}$,
$ \frac{1+\sqrt{5}}{2}$,
$ -c^{3}_{40}
-c^{5}_{40}
+c^{7}_{40}
$;\ \ 
$0$,
$ c^{3}_{40}
+c^{5}_{40}
-c^{7}_{40}
$,
$ -c^{3}_{40}
-c^{5}_{40}
+c^{7}_{40}
$,
$0$;\ \ 
$ -1$,
$ -1$,
$ -\sqrt{2}$;\ \ 
$ -1$,
$ \sqrt{2}$;\ \ 
$0$)

\vskip 1ex 
\color{grey}

\noindent34. ind = $(6;80
)_{1}^{31}$:\ \ 
$d_i$ = ($1.0$,
$1.0$,
$1.414$,
$1.618$,
$1.618$,
$2.288$) 

\vskip 0.7ex
\hangindent=3em \hangafter=1
$D^2=$ 14.472 = 
 $10+2\sqrt{5}$

\vskip 0.7ex
\hangindent=3em \hangafter=1
$T = ( 0,
\frac{1}{2},
\frac{15}{16},
\frac{2}{5},
\frac{9}{10},
\frac{27}{80} )
$,

\vskip 0.7ex
\hangindent=3em \hangafter=1
$S$ = ($ 1$,
$ 1$,
$ \sqrt{2}$,
$ \frac{1+\sqrt{5}}{2}$,
$ \frac{1+\sqrt{5}}{2}$,
$ c^{3}_{40}
+c^{5}_{40}
-c^{7}_{40}
$;\ \ 
$ 1$,
$ -\sqrt{2}$,
$ \frac{1+\sqrt{5}}{2}$,
$ \frac{1+\sqrt{5}}{2}$,
$ -c^{3}_{40}
-c^{5}_{40}
+c^{7}_{40}
$;\ \ 
$0$,
$ c^{3}_{40}
+c^{5}_{40}
-c^{7}_{40}
$,
$ -c^{3}_{40}
-c^{5}_{40}
+c^{7}_{40}
$,
$0$;\ \ 
$ -1$,
$ -1$,
$ -\sqrt{2}$;\ \ 
$ -1$,
$ \sqrt{2}$;\ \ 
$0$)

\vskip 1ex 
\color{grey}

\noindent35. ind = $(6;80
)_{1}^{79}$:\ \ 
$d_i$ = ($1.0$,
$1.0$,
$1.414$,
$1.618$,
$1.618$,
$2.288$) 

\vskip 0.7ex
\hangindent=3em \hangafter=1
$D^2=$ 14.472 = 
 $10+2\sqrt{5}$

\vskip 0.7ex
\hangindent=3em \hangafter=1
$T = ( 0,
\frac{1}{2},
\frac{15}{16},
\frac{3}{5},
\frac{1}{10},
\frac{43}{80} )
$,

\vskip 0.7ex
\hangindent=3em \hangafter=1
$S$ = ($ 1$,
$ 1$,
$ \sqrt{2}$,
$ \frac{1+\sqrt{5}}{2}$,
$ \frac{1+\sqrt{5}}{2}$,
$ c^{3}_{40}
+c^{5}_{40}
-c^{7}_{40}
$;\ \ 
$ 1$,
$ -\sqrt{2}$,
$ \frac{1+\sqrt{5}}{2}$,
$ \frac{1+\sqrt{5}}{2}$,
$ -c^{3}_{40}
-c^{5}_{40}
+c^{7}_{40}
$;\ \ 
$0$,
$ c^{3}_{40}
+c^{5}_{40}
-c^{7}_{40}
$,
$ -c^{3}_{40}
-c^{5}_{40}
+c^{7}_{40}
$,
$0$;\ \ 
$ -1$,
$ -1$,
$ -\sqrt{2}$;\ \ 
$ -1$,
$ \sqrt{2}$;\ \ 
$0$)

\vskip 1ex 

 \color{black} \vskip 2ex

\noindent36. ind = $(6;80
)_{2}^{1}$:\ \ 
$d_i$ = ($1.0$,
$1.0$,
$1.414$,
$1.618$,
$1.618$,
$2.288$) 

\vskip 0.7ex
\hangindent=3em \hangafter=1
$D^2=$ 14.472 = 
 $10+2\sqrt{5}$

\vskip 0.7ex
\hangindent=3em \hangafter=1
$T = ( 0,
\frac{1}{2},
\frac{3}{16},
\frac{2}{5},
\frac{9}{10},
\frac{47}{80} )
$,

\vskip 0.7ex
\hangindent=3em \hangafter=1
$S$ = ($ 1$,
$ 1$,
$ \sqrt{2}$,
$ \frac{1+\sqrt{5}}{2}$,
$ \frac{1+\sqrt{5}}{2}$,
$ c^{3}_{40}
+c^{5}_{40}
-c^{7}_{40}
$;\ \ 
$ 1$,
$ -\sqrt{2}$,
$ \frac{1+\sqrt{5}}{2}$,
$ \frac{1+\sqrt{5}}{2}$,
$ -c^{3}_{40}
-c^{5}_{40}
+c^{7}_{40}
$;\ \ 
$0$,
$ c^{3}_{40}
+c^{5}_{40}
-c^{7}_{40}
$,
$ -c^{3}_{40}
-c^{5}_{40}
+c^{7}_{40}
$,
$0$;\ \ 
$ -1$,
$ -1$,
$ -\sqrt{2}$;\ \ 
$ -1$,
$ \sqrt{2}$;\ \ 
$0$)

\vskip 1ex 
\color{grey}

\noindent37. ind = $(6;80
)_{2}^{49}$:\ \ 
$d_i$ = ($1.0$,
$1.0$,
$1.414$,
$1.618$,
$1.618$,
$2.288$) 

\vskip 0.7ex
\hangindent=3em \hangafter=1
$D^2=$ 14.472 = 
 $10+2\sqrt{5}$

\vskip 0.7ex
\hangindent=3em \hangafter=1
$T = ( 0,
\frac{1}{2},
\frac{3}{16},
\frac{3}{5},
\frac{1}{10},
\frac{63}{80} )
$,

\vskip 0.7ex
\hangindent=3em \hangafter=1
$S$ = ($ 1$,
$ 1$,
$ \sqrt{2}$,
$ \frac{1+\sqrt{5}}{2}$,
$ \frac{1+\sqrt{5}}{2}$,
$ c^{3}_{40}
+c^{5}_{40}
-c^{7}_{40}
$;\ \ 
$ 1$,
$ -\sqrt{2}$,
$ \frac{1+\sqrt{5}}{2}$,
$ \frac{1+\sqrt{5}}{2}$,
$ -c^{3}_{40}
-c^{5}_{40}
+c^{7}_{40}
$;\ \ 
$0$,
$ c^{3}_{40}
+c^{5}_{40}
-c^{7}_{40}
$,
$ -c^{3}_{40}
-c^{5}_{40}
+c^{7}_{40}
$,
$0$;\ \ 
$ -1$,
$ -1$,
$ -\sqrt{2}$;\ \ 
$ -1$,
$ \sqrt{2}$;\ \ 
$0$)

\vskip 1ex 
\color{grey}

\noindent38. ind = $(6;80
)_{2}^{71}$:\ \ 
$d_i$ = ($1.0$,
$1.0$,
$1.414$,
$1.618$,
$1.618$,
$2.288$) 

\vskip 0.7ex
\hangindent=3em \hangafter=1
$D^2=$ 14.472 = 
 $10+2\sqrt{5}$

\vskip 0.7ex
\hangindent=3em \hangafter=1
$T = ( 0,
\frac{1}{2},
\frac{5}{16},
\frac{2}{5},
\frac{9}{10},
\frac{57}{80} )
$,

\vskip 0.7ex
\hangindent=3em \hangafter=1
$S$ = ($ 1$,
$ 1$,
$ \sqrt{2}$,
$ \frac{1+\sqrt{5}}{2}$,
$ \frac{1+\sqrt{5}}{2}$,
$ c^{3}_{40}
+c^{5}_{40}
-c^{7}_{40}
$;\ \ 
$ 1$,
$ -\sqrt{2}$,
$ \frac{1+\sqrt{5}}{2}$,
$ \frac{1+\sqrt{5}}{2}$,
$ -c^{3}_{40}
-c^{5}_{40}
+c^{7}_{40}
$;\ \ 
$0$,
$ c^{3}_{40}
+c^{5}_{40}
-c^{7}_{40}
$,
$ -c^{3}_{40}
-c^{5}_{40}
+c^{7}_{40}
$,
$0$;\ \ 
$ -1$,
$ -1$,
$ -\sqrt{2}$;\ \ 
$ -1$,
$ \sqrt{2}$;\ \ 
$0$)

\vskip 1ex 
\color{grey}

\noindent39. ind = $(6;80
)_{2}^{39}$:\ \ 
$d_i$ = ($1.0$,
$1.0$,
$1.414$,
$1.618$,
$1.618$,
$2.288$) 

\vskip 0.7ex
\hangindent=3em \hangafter=1
$D^2=$ 14.472 = 
 $10+2\sqrt{5}$

\vskip 0.7ex
\hangindent=3em \hangafter=1
$T = ( 0,
\frac{1}{2},
\frac{5}{16},
\frac{3}{5},
\frac{1}{10},
\frac{73}{80} )
$,

\vskip 0.7ex
\hangindent=3em \hangafter=1
$S$ = ($ 1$,
$ 1$,
$ \sqrt{2}$,
$ \frac{1+\sqrt{5}}{2}$,
$ \frac{1+\sqrt{5}}{2}$,
$ c^{3}_{40}
+c^{5}_{40}
-c^{7}_{40}
$;\ \ 
$ 1$,
$ -\sqrt{2}$,
$ \frac{1+\sqrt{5}}{2}$,
$ \frac{1+\sqrt{5}}{2}$,
$ -c^{3}_{40}
-c^{5}_{40}
+c^{7}_{40}
$;\ \ 
$0$,
$ c^{3}_{40}
+c^{5}_{40}
-c^{7}_{40}
$,
$ -c^{3}_{40}
-c^{5}_{40}
+c^{7}_{40}
$,
$0$;\ \ 
$ -1$,
$ -1$,
$ -\sqrt{2}$;\ \ 
$ -1$,
$ \sqrt{2}$;\ \ 
$0$)

\vskip 1ex 
\color{grey}

\noindent40. ind = $(6;80
)_{2}^{41}$:\ \ 
$d_i$ = ($1.0$,
$1.0$,
$1.414$,
$1.618$,
$1.618$,
$2.288$) 

\vskip 0.7ex
\hangindent=3em \hangafter=1
$D^2=$ 14.472 = 
 $10+2\sqrt{5}$

\vskip 0.7ex
\hangindent=3em \hangafter=1
$T = ( 0,
\frac{1}{2},
\frac{11}{16},
\frac{2}{5},
\frac{9}{10},
\frac{7}{80} )
$,

\vskip 0.7ex
\hangindent=3em \hangafter=1
$S$ = ($ 1$,
$ 1$,
$ \sqrt{2}$,
$ \frac{1+\sqrt{5}}{2}$,
$ \frac{1+\sqrt{5}}{2}$,
$ c^{3}_{40}
+c^{5}_{40}
-c^{7}_{40}
$;\ \ 
$ 1$,
$ -\sqrt{2}$,
$ \frac{1+\sqrt{5}}{2}$,
$ \frac{1+\sqrt{5}}{2}$,
$ -c^{3}_{40}
-c^{5}_{40}
+c^{7}_{40}
$;\ \ 
$0$,
$ c^{3}_{40}
+c^{5}_{40}
-c^{7}_{40}
$,
$ -c^{3}_{40}
-c^{5}_{40}
+c^{7}_{40}
$,
$0$;\ \ 
$ -1$,
$ -1$,
$ -\sqrt{2}$;\ \ 
$ -1$,
$ \sqrt{2}$;\ \ 
$0$)

\vskip 1ex 
\color{grey}

\noindent41. ind = $(6;80
)_{2}^{9}$:\ \ 
$d_i$ = ($1.0$,
$1.0$,
$1.414$,
$1.618$,
$1.618$,
$2.288$) 

\vskip 0.7ex
\hangindent=3em \hangafter=1
$D^2=$ 14.472 = 
 $10+2\sqrt{5}$

\vskip 0.7ex
\hangindent=3em \hangafter=1
$T = ( 0,
\frac{1}{2},
\frac{11}{16},
\frac{3}{5},
\frac{1}{10},
\frac{23}{80} )
$,

\vskip 0.7ex
\hangindent=3em \hangafter=1
$S$ = ($ 1$,
$ 1$,
$ \sqrt{2}$,
$ \frac{1+\sqrt{5}}{2}$,
$ \frac{1+\sqrt{5}}{2}$,
$ c^{3}_{40}
+c^{5}_{40}
-c^{7}_{40}
$;\ \ 
$ 1$,
$ -\sqrt{2}$,
$ \frac{1+\sqrt{5}}{2}$,
$ \frac{1+\sqrt{5}}{2}$,
$ -c^{3}_{40}
-c^{5}_{40}
+c^{7}_{40}
$;\ \ 
$0$,
$ c^{3}_{40}
+c^{5}_{40}
-c^{7}_{40}
$,
$ -c^{3}_{40}
-c^{5}_{40}
+c^{7}_{40}
$,
$0$;\ \ 
$ -1$,
$ -1$,
$ -\sqrt{2}$;\ \ 
$ -1$,
$ \sqrt{2}$;\ \ 
$0$)

\vskip 1ex 
\color{grey}

\noindent42. ind = $(6;80
)_{2}^{31}$:\ \ 
$d_i$ = ($1.0$,
$1.0$,
$1.414$,
$1.618$,
$1.618$,
$2.288$) 

\vskip 0.7ex
\hangindent=3em \hangafter=1
$D^2=$ 14.472 = 
 $10+2\sqrt{5}$

\vskip 0.7ex
\hangindent=3em \hangafter=1
$T = ( 0,
\frac{1}{2},
\frac{13}{16},
\frac{2}{5},
\frac{9}{10},
\frac{17}{80} )
$,

\vskip 0.7ex
\hangindent=3em \hangafter=1
$S$ = ($ 1$,
$ 1$,
$ \sqrt{2}$,
$ \frac{1+\sqrt{5}}{2}$,
$ \frac{1+\sqrt{5}}{2}$,
$ c^{3}_{40}
+c^{5}_{40}
-c^{7}_{40}
$;\ \ 
$ 1$,
$ -\sqrt{2}$,
$ \frac{1+\sqrt{5}}{2}$,
$ \frac{1+\sqrt{5}}{2}$,
$ -c^{3}_{40}
-c^{5}_{40}
+c^{7}_{40}
$;\ \ 
$0$,
$ c^{3}_{40}
+c^{5}_{40}
-c^{7}_{40}
$,
$ -c^{3}_{40}
-c^{5}_{40}
+c^{7}_{40}
$,
$0$;\ \ 
$ -1$,
$ -1$,
$ -\sqrt{2}$;\ \ 
$ -1$,
$ \sqrt{2}$;\ \ 
$0$)

\vskip 1ex 
\color{grey}

\noindent43. ind = $(6;80
)_{2}^{79}$:\ \ 
$d_i$ = ($1.0$,
$1.0$,
$1.414$,
$1.618$,
$1.618$,
$2.288$) 

\vskip 0.7ex
\hangindent=3em \hangafter=1
$D^2=$ 14.472 = 
 $10+2\sqrt{5}$

\vskip 0.7ex
\hangindent=3em \hangafter=1
$T = ( 0,
\frac{1}{2},
\frac{13}{16},
\frac{3}{5},
\frac{1}{10},
\frac{33}{80} )
$,

\vskip 0.7ex
\hangindent=3em \hangafter=1
$S$ = ($ 1$,
$ 1$,
$ \sqrt{2}$,
$ \frac{1+\sqrt{5}}{2}$,
$ \frac{1+\sqrt{5}}{2}$,
$ c^{3}_{40}
+c^{5}_{40}
-c^{7}_{40}
$;\ \ 
$ 1$,
$ -\sqrt{2}$,
$ \frac{1+\sqrt{5}}{2}$,
$ \frac{1+\sqrt{5}}{2}$,
$ -c^{3}_{40}
-c^{5}_{40}
+c^{7}_{40}
$;\ \ 
$0$,
$ c^{3}_{40}
+c^{5}_{40}
-c^{7}_{40}
$,
$ -c^{3}_{40}
-c^{5}_{40}
+c^{7}_{40}
$,
$0$;\ \ 
$ -1$,
$ -1$,
$ -\sqrt{2}$;\ \ 
$ -1$,
$ \sqrt{2}$;\ \ 
$0$)

\vskip 1ex 

 \color{black} \vskip 2ex

\

%A pair of $S,T$ matrices with negative quantum dimensions:
%
%\
%
%\noindent $19_{1}^{1}$: dims = $6$, levels = $( 9 ) $, $\rho(\mathfrak{t})$ =
%$( \frac{1}{9}, \frac{2}{9}, \frac{4}{9}, \frac{5}{9}, \frac{7}{9},
%\frac{8}{9} ) $
%
%\vskip 0.7ex \hangindent=5em \hangafter=1 $d_i$ = ($1.$, $0.34729$, $-1.$,
%$1.$, $1.53208$, $-1.8793$) 
%
%\vskip 0.7ex \hangindent=5em \hangafter=1 $D^2=$ 9. = $9$
%
%\vskip 0.7ex \hangindent=5em \hangafter=1 $T = ( 0, \frac{1}{9}, \frac{1}{3},
%\frac{2}{3}, \frac{4}{9}, \frac{7}{9} ) $,
%
%\vskip 0.7ex \hangindent=5em \hangafter=1 $S$ = ($ 1$, $ c^{2}_{9} $, $ -1$, $
%1$, $ c^{1}_{9} $, $ -\xi_{9}^{2,1}$;\ \ $ 1$, $ \xi_{9}^{2,1}$, $ c^{1}_{9}
%$, $ 1$, $ 1$;\ \ $ 1$, $ -1$, $ -c^{2}_{9} $, $ -c^{1}_{9} $;\ \ $ 1$, $
%-\xi_{9}^{2,1}$, $ c^{2}_{9} $;\ \ $ 1$, $ 1$;\ \ $ 1$)
%
%%Number of MD's: 43
%
%\

The above list include all rank-6 modular data with non-integral $D^2$ and
coming from resolved $\SL$ representations.  The list misses two known modular
data
% and two known modular data
%with $D^2 = 20$, $\ord(T)=10$.  The modular data with $D^2 = 20$, $\ord(T)=10$
%and $s_i = (0,0,0,\frac12, \frac15,\frac45), (0,0,0,\frac12, \frac25,\frac35)$
%must come from the GT orbit $3$.  
with non-integral $D^2 = 74.617$, whose topological spins are $s_i =
(0,\frac19,\frac19,\frac19,\frac13,\frac23)$ or $s_i=
(0,\frac89,\frac89,\frac89,\frac13,\frac23) $. From those $s_i$'s, we find that
they must come from the unresolved GT orbit $(4, 1,1;9,1,1)$.  In the main text
of this paper, we showed that the unresolved $\SL$ representations can only
produce such modular data (and its conjugations by Galois action and signed
diagonal matrices).  The unresolved cases are handled in the main text of
the paper, which leads to a complete classification of all rank-6 modular data.

%\begin{align}
% p_\pm = \sum_i d_i^2 \theta_i^{\pm} = D \ee^{\pm \ii 2\pi c/8}
%\end{align}
%
%\begin{align}
% \rho_\a(\fs) = \ee^{-2\pi \ii \frac{\a}{4}} S/D, \ \
% \rho_\a(\ft) = \ee^{2\pi \ii \frac{\a}{12}} \ee^{-2\pi \ii \frac{c}{24}} T,
%\ \ \ \a \in \BZ_{12}. 
%\end{align}
%
%\begin{align}
% \rho_\a(\fs)_{uu} = \ee^{-2\pi \ii \frac{\a}{4}}/D, \ \ \ \
%\tilde \theta_u = \ee^{2\pi \ii \frac{\a}{12}} \ee^{-2\pi \ii \frac{c}{24}}.
%\end{align}
%
%\begin{align}
% \sum_i 
%\frac{\rho_\a(\fs)_{iu}^2} 
%{\rho_\a(\fs)_{uu}^2 }
%\tilde \theta_i  \tilde \theta_u^{-1}
%= 
%\frac{\ee^{-2\pi \ii \frac{\a}{4}}}{\rho_\a(\fs)_{uu}}
% \ee^{ \ii 2\pi c/8}
%= 
%\frac{\tilde \theta_u^{-3}}{\rho_\a(\fs)_{uu}}
%\end{align}
%
%\begin{align}
% \sum_i \rho_\a(\fs)_{iu}^2 \tilde \theta_i  \tilde \theta_u^2
%=  \rho_\a(\fs)_{uu} 
%\end{align}
%
%
%\begin{align}
% (\rho_\a(\ft) \rho_\a(\fs) \rho_\a(\ft) \rho_\a(\fs) \rho_\a(\ft) )_{uu}
%=
% \rho_\a(\ft\fs\ft \fs \ft )_{uu}
%=  \rho_\a(\fs)_{uu} 
%\end{align}

\bibliographystyle{plain}
\bibliography{ref.bib}

\vfill
\pagebreak

%\usepackage[top=1.0in, bottom=.75in, left=.75in, right=.75in]{geometry}

%\parindent = 1em
%\numberwithin{equation}{section}

\part*{
\Large Supplementary materials
}
\label{suppl}

\

\unappendix
\setcounter{section}{0}

\centerline{\Large
Reconstruction of modular data from $\SL$ representations
}

\

\centerline{Siu-Hung Ng, Eric Rowell, Zhenghan Wang, Xiao-Gang Wen}

\

%\begin{abstract} 
%
%\end{abstract}

%\setcounter{tocdepth}{2} \tableofcontents 

\section{Tables of irreducible representations of
prime-power levels}
\label{Section1}

In this section, we list all the $\SL$ irreducible representations of dimension
1 -- 9, whose level ($l=\ord(\rho(\ft))$) is a power of single prime number.
In the list, $\rho(\mathfrak{t})$ is presented in term of topological spins
$(\tilde s_{1},\tilde s_{2},\cdots)$ ($\tilde s_{i} =
\arg(\rho_a(\mathfrak{t})_{ii})$).

Note that $\rho(\mathfrak{s})$ is symmetric and $\rho(\mathfrak{s})_{ij}$'s are
either all real or all imaginary.  When $\rho(\mathfrak{s})_{ij}$'s are all
real, $\rho(\mathfrak{s})$ is presented as $(\rho_{11}$, $ \rho_{12}$, $
\rho_{13}$, $ \rho_{14}$, $ \cdots;\ \  \rho_{22}$, $ \rho_{23}$, $ \rho_{24}$,
$ \cdots)$.  We introduce
\begin{align}
\zeta^m_n &= \mathrm{e}^{2\pi \mathrm{i} m/n}, \ \ \
c^m_n = \zeta^m_n+\zeta^{-m}_n, \ \ \
s^m_n = \zeta^m_n-\zeta^{-m}_n,
\nonumber\\
\xi^{m,k}_n &= (\zeta^m_{2n}-\zeta^{-m}_{2n})/(\zeta_{2n}^k-\zeta_{2n}^{-k}), 
\ \ \ \ \
\xi^{m}_n = \xi_n^{m,1}, 
\end{align}
to describe the matrix
elements of $\rho(\mathfrak{s})$.

When $\rho(\mathfrak{s})_{ij}$'s are all imaginary, $\rho(\mathfrak{s})$ is
presented as $\ii(-\ii\rho_{11}$, $ -\ii\rho_{12}$, $ -\ii\rho_{13}$, $ -\ii\rho_{14}$, $
\cdots;\ \  -\ii\rho_{22}$, $ -\ii\rho_{23}$, $ -\ii\rho_{24}$, $ \cdots)$, or as
$(s_n^m)^{-1}(s_n^m\rho_{11}$, $ s_n^m\rho_{12}$, $ s_n^m\rho_{13}$, $ s_n^m\rho_{14}$, $
\cdots;\ \  s_n^m\rho_{22}$, $ s_n^m\rho_{23}$, $ s_n^m\rho_{24}$, $ \cdots)$.  In any
case, the numbers inside the bracket $(\cdots)$ are all real.  We can tell a
representation to be even or odd by the absence or the presence of $\ii$ in
front of the bracket $(\cdots)$.

We note that two symmetric representations are equivalent up to a permutation
of indices, and a conjugation of signed diagonal matrix.  To choose the
ordering in indices, we introduce arrays $O_i =[\text{DenominatorOf}(\tilde
s_{i}),\tilde s_{i}, \rho_{ii}]$.  The order of two arrays is determined by
first compare the lengths of the two arrays. If the lengths are equal, we then
compare the first elements of the two arrays.  If the first elements are equal,
we then compare the second elements of the two arrays, {\it etc}.  To compare
two cyclotomic numbers, here we used the ordering of cyclotomic numbers
provided by GAP computer algebraic system.  We order the indices to make
$O_1\leq O_2 \leq O_3 \cdots$.  The  conjugation of signed diagonal matrix is
chosen to make $-\rho(\mathfrak{s})_{1j} < \rho(\mathfrak{s})_{1j} $ for
$j=2,3,\cdots$.  If $\rho(\mathfrak{s})_{1j}=0$, we will try to make
$-\rho(\mathfrak{s})_{2j} < \rho(\mathfrak{s})_{2j} $, {\it etc}.

All the prime-power-level irreducible representations are labeled by index
$d_{l,k}^{a,m}$, where $d$ is the dimension and $l$ is the level of the
representation.  The irreducible representations of a given $d,l$ can be
grouped into several orbits, generated by Galois conjugations and tensoring of
1-dimensional representations that do not change the level $l$.  $k$ in
$d_{l,k}^{a,m}$ labels those different orbits.  If there is only 1 orbit for a
given $d,l$, the index $k$ will be dropped.  

The irreducible representation labeled by $d_{l,k}^{a,m}$ is generated from the
irreducible representation labeled by $d_{l,k}^{1,0}$ via the following Galois
conjugations and tensoring of 1-dimensional representations
\begin{align}
 \rho_{d_{l,k}^{a,m}}(\ft) &= \s_a\big(\rho_{d_{l,k}^{1,0}}(\ft) \big)\ee^{2\pi \ii\frac m{12}} 
\nonumber\\
 \rho_{d_{l,k}^{a,m}}(\fs) &= \s_a\big(\rho_{d_{l,k}^{1,0}}(\fs) \big)\ee^{-2\pi \ii\frac m{4}} 
\end{align}
where the Galois conjugation $\s_a $ is in $ \mathrm{Gal}(\BQ_n)$ with $n$ be the least
common multiple of $\ord(\rho_{d_{l,k}^{1,0}}(\ft))$ and the conductor of
$\rho_{d_{l,k}^{1,0}}(\fs)$, and the integer $a$ is given by
\begin{align}
 \s_a\big(\ee^{2\pi \ii /n} \big) = 
 \ee^{2\pi \ii a/n} .
\end{align}
Also $m \in \BZ_{12}$ is such that  $ \ord(\rho_{d_{l,k}^{1,0}}(\ft)\ee^{2\pi
\ii\frac m{12}})= \ord(\rho_{d_{l,k}^{1,0}}(\ft))$. Due to this condition, when
$l$ is not divisible by 2 and 3, $m$ can only be $0$.  In this case, we will
drop $m$.  Here $d_{l,k}^{1,0}$ is the representation in the orbit with minimal
$[\tilde s_1, \tilde s_2,\cdots] $.

We like to remark that two representations labeled by $d_{l,k}^{a,m}$ and
$d_{l,k}^{a',m'}$ may be equivalent.  The numbers of distinct irreducible
representations with prime-power level (PPL) in each dimension are given by
\begin{align}
\begin{tabular}{|r|r|r|r|r|r|r|r|r|r|r|r|r|}
\hline
\text{dim}: & 1 & 2 & 3 & 4 & 5 & 6 & 7 & 8 & 9 & 10 & 11 & 12 \\
\hline
\text{\# of irreps with PPL} & 6 & 11 & 33 & 18 & 3 & 70 & 3  & 10  & 4 & 7 & 3 & 176  \\ 
\hline
\text{\# of irreps} & 12 & 54 & 136 & 180 & 36 & 720 & 36 & 456 & 476 & 222 & 36 & 3214 \\ 
\hline
\end{tabular}
\end{align}
In the above we also list the numbers of distinct irreducible representations,
which are tensor products of the irreducible representations with prime-power
levels.

\def\arraystretch{1.6} \setlength\tabcolsep{3pt}

\begin{longtable}{|l|l|p{4.5in}|}
\hline
$d_{l,k}^{a,m}$ & \# & $\rho(\mathfrak{t})$, \ $\rho(\mathfrak{s})$\\
\hline
$6_{5}^{1}$ & 1 & $( 0,
0,
\frac{1}{5},
\frac{2}{5},
\frac{3}{5},
\frac{4}{5} )
$, ~ $(s_{5}^{2})^{-1}$($\frac{\sqrt{5}}{5}$,
$-\frac{5-\sqrt{5}}{10}$,
$-\frac{\sqrt{10}}{5}$,
$-\frac{5+\sqrt{5}}{10}$,
$\frac{5-3\sqrt{5}}{10}$,
$\frac{1-\sqrt{5}}{\sqrt{10}}$;\ \ 
$-\frac{\sqrt{5}}{5}$,
$\frac{-1+\sqrt{5}}{\sqrt{10}}$,
$\frac{5-3\sqrt{5}}{10}$,
$\frac{5+\sqrt{5}}{10}$,
$-\frac{\sqrt{10}}{5}$;\ \ 
$\frac{5-\sqrt{5}}{10}$,
$\frac{1-\sqrt{5}}{\sqrt{10}}$,
$-\frac{\sqrt{10}}{5}$,
$-\frac{\sqrt{5}}{5}$;\ \ 
$-\frac{\sqrt{5}}{5}$,
$\frac{5-\sqrt{5}}{10}$,
$\frac{\sqrt{10}}{5}$;\ \ 
$\frac{\sqrt{5}}{5}$,
$\frac{1-\sqrt{5}}{\sqrt{10}}$;\ \ 
$-\frac{5-\sqrt{5}}{10}$)
\\
\hline
\end{longtable}

\section{Compute rank-6 $S,T$ matrices from resolved irrep-sum
representations}
\label{Section2}
\label{transU}

In Appendix B.2, we list all 25 irrep-sum representations, $\rho_\text{isum}$,
that may potentially produce MD representations $\rho$
\begin{align}
 \rho = U  \rho_\text{isum} U^\top,\ \ \ 
U = V_\mathrm{sgn} U_0 P
\end{align}
after proper orthogonal transformations $U$.  Here $V_\mathrm{sgn}$ is a signed
diagonal matrix and $P$ is a permutation matrix.  We choose a $P$ so that the
topological spins of $\td\rho(\ft) = P\rho_\text{isum}(\ft)P$, $\td s_i \in
[0,1)$, form an ordered list, which are listed in the following table.  (The
rational topological spins, $\td s_i$, are ordered first according to the
denominators of $\td s_i$'s.  If the denominators are the same, then they are
ordered according to the value of $\td s_i$.) $U_0$ is block diagonal in this
basis, that keeps $\td\rho(\ft)$ unchanged:
\begin{align}
 U_0 \td\rho(\ft) U_0^\top =
 \td\rho(\ft)  .
\end{align}

In the following list, we examine all the irrep-sum representations listed in
Appendix B.2.  We determine which irrep-sum representations are resolved, as
indicated in the list.  For resolved representations, we list all $U_0$'s that
transform $D_{\td \rho}(\s)$ into $D_{\rho}(\s)$, such that $D_{\rho}(\s)$,
when restricted within each eigenspace of $\rho(\ft)$, are either signed
permutations or zero. 

Only $U_0$'s that are not related by  signed diagonal matrices,
$V_\mathrm{sgn}$, are listed.  So the most general transformations that
transform $D_{\td \rho}(\s)$ into $D_{\rho}(\s)$ described above
are given by $V_\mathrm{sgn} U_0$ (up to a permutation of
indices).

For each $U_0$, we examine all possible choices of index $u$ that may
correspond to the unit object, and construct $S,T$ matrices via 
\begin{align}
S_{i,j} = \frac{\rho(\fs)_{i,j} }{\rho(\fs)_{u,u}},\ \ \ \ 
T_{i,j} = \frac{\rho(\ft)_{i,j} }{\rho(\ft)_{u,u}},  
\end{align}
as well as the fusion rule $N^{ij}_k$ via
\begin{align} 
\label{Ver1} 
 N^{i,j}_k = \sum_{l=0}^{r-1} \frac{
\rho(\fs)_{l,i} \rho(\fs)_{l,j} \rho^*(\fs)_{l,k}}{ \rho(\fs)_{l,u} } . 
\end{align}
We also consider all possible choices of the  signed diagonal matrices
$V_\mathrm{sgn}$, that make the $N^{ij}_k$ non-negative.  We then check if the
resulting $S,T$ matrices satisfy the conditions of modular data summarized in
Proposition B.1.

In fact, we check the following conditions for each choice of $u$, and
summarize the results in the tables in the list. We use  error code ``0'' to
indicates that the condition is satisfied.  We use ``-'' (corresponding to
internal error code ``9'') to indicates that the condition is not checked
(usually because it does not apply).  A non-zero error code indicates that the
condition is not satisfied. The meanings of the non-zero error codes are
explained below:
\begin{enumerate}

\item ``$D_\rho$ conditions'':
The MD representations $\rho$ associated with a modular data
must satisfy (see Proposition B.1)
\begin{align}
\s(\rho(\fs)) &= D_{\rho}(\s) \rho(\fs) =\rho(\fs)D_{\rho}^\top(\s),
\nonumber\\
\s^2(\rho(\ft)) &= D_{\rho}(\s) \rho(\ft) D_{\rho}^\top(\s)
\end{align}
The error code is
given by $e_1 + 2e_2 +4e_3$.\\
a)  $e_3=0$ if $D_{\rho}(\s)$ is a signed permutation matrix.
$e_3=1$ otherwise.\\
b)  $e_2=0$ if $\s^2(\rho(\ft)) = D_{\rho}(\s) \rho(\ft) D_{\rho}^\top(\s)$.
$e_2=1$ otherwise.\\
c)  $e_1=0$ if $\s(\rho(\fs)) = D_{\rho}(\s) \rho(\fs) $.
$e_1=1$ otherwise.\\
We test the above conditions for generators $\s$ of the Galois group.
Note that the error code depends on the
the Galois conjugations $\s$.
In the following table, we list the maximum error code.

\item ``$[\rho(\fs)\rho(\ft)]^3 = \rho(\fs)^2 = \tilde C $'': The error code is
given by $e_1 + 2e_2 +4e_3$.\\
a)  $e_3=0$ if $[\rho(\fs)\rho(\ft)]^3 = \rho(\fs)^2$ and
$\rho(\fs)^4 = \text{id}$.
$e_3=1$ otherwise.\\
b) $e_2=0$ if $\text{Tr} (\tilde C) \neq 0$.
$e_2=1$ otherwise.\\
c) $e_1=0$ if $\tilde C$ is a signed permutation.
$e_1=1$ otherwise.

\item
``$\rho(\fs)_{iu} \rho^*(\fs)_{ju} \in \BR $'': 
If $u$ corresponds to the unit object, then
$\rho(\fs)_{iu}
\rho^*(\fs)_{ju}  $ must be real for all $i,j$.

\item ``$\rho(\fs)_{iu} \neq 0$'': If $u$ corresponds to the unit object, then
$\rho_{iu}(\fs) \neq 0$ for all $i$.  The error code is ``2'' if
$\rho(\fs)_{uu}=0$.  The error code is ``1'' if $\rho(\fs)_{uu}\neq 0$, but
$\rho(\fs)_{iu}= 0$ for some $i$'s.  Once this condition is satisfied, we can
check the following conditions.

\item ``cnd$(S)$, cnd$(\rho(\fs))$'': Conductor of $S$ matrix, cnd$(S)$, must
divides $\ord(T)$.  Conductor of $\rho(\fs)$, cnd$(\rho(\fs))$, must divides
$\ord(\rho(\ft))$.  The error code is given by $e_1 + 2e_2$.\\
a) $e_2=0$ if cnd$(\rho(\fs))\big|\ord(\rho(\ft))$.
$e_2=1$ otherwise.\\
b) $e_1=0$ if cnd$(S)\big|\ord(T)$.
$e_1=1$ otherwise.

\item ``norm$(D^2)$ factors'': The distinct prime factors of norm$(D^2)$ must
coincide with the  distinct prime factors of $\ord(T)$.  Here norm$(c)$ is the
product of distinct Galois conjugations of $c$.  The error code is ``2'' if
norm$(D^2)$ is not an integer.  The error code is ``1'' if norm$(D^2)$ is an
integer, but the distinct prime factors of norm$(D^2)$ do not coincide with the
distinct prime factors of $\ord(T)$

\item
``$1/\rho(\fs)_{iu} = $ cyc-int'': $1/\rho(\fs)_{iu}$ must be a cyclotomic
integer for all $i$'s.

\item ``norm$(1/\rho(\fs)_{iu})$ factors'': 
The error code is $e_1+3e_2$.\
a) $e_2=0$ if norm$(1/\rho(\fs)_{uu})$ is an integer and the distinct prime
factors of norm$(1/\rho(\fs)_{uu})$ coincide with the  distinct prime factors
of $\ord(T)$.
$e_2=1$ if norm$(1/\rho(\fs)_{uu})$ is an integer, but the distinct prime
factors of norm$(1/\rho(\fs)_{uu})$ do not coincide with the  distinct prime factors
of $\ord(T)$.
$e_2=2$ if norm$(1/\rho(\fs)_{uu})$ is not an integer.\\
b) $e_1=0$ if norm$(1/\rho(\fs)_{iu})$ is an integer and the distinct prime
factors of norm$(1/\rho(\fs)_{iu})$ is included in the  distinct prime factors
of $\ord(T)$.
$e_1=1$ if norm$(1/\rho(\fs)_{iu})$ is an integer, but the distinct prime
factors of norm$(1/\rho(\fs)_{iu})$ is not included in the  distinct prime factors
of $\ord(T)$.
$e_1=2$ if norm$(1/\rho(\fs)_{iu})$ is not an integer.
\\
Note that the error code depends on $i$. The maximum error code is displayed in
the following table.

\item
``$\frac{S_{ij}}{S_{uj}} = $ cyc-int'': $\frac{S_{ij}}{S_{uj}}$ must be
cyclotomic integers for all $i,j$'s.

\item
``$N^{ij}_k \in \BN$'': $N^{ij}_k$ must be non-negative integers.  If some
$N^{ij}_k$ are not rational numbers, an error code  ``4'' is given.  In this
case, we assume $V_\mathrm{sgn}=\text{id}$, and continue to check the following
conditions.  If $N^{ij}_k$ are rational numbers but some are not integers, an
error code ``3'' or ``2'' is given.  ``3'' means that we cannot find signed
diagonal matrices $V_\mathrm{sgn}$ to make all $N^{ij}_k$ non-negative.  In
this case, we assume $V_\mathrm{sgn}=\text{id}$, and continue to check the
following conditions.  ``2'' means that we can find one or more signed diagonal
matrices, $V_\mathrm{sgn}$, to make all $N^{ij}_k$ non-negative.  We will use
those $V_\mathrm{sgn}$ to continue to check the following conditions.  If
$N^{ij}_k$ are integers, an error code ``1'' or ``0'' is given.  ``1'' means
that we cannot find signed diagonal matrices $V_\mathrm{sgn}$ to make all
$N^{ij}_k$ non-negative.  In this case, we assume $V_\mathrm{sgn}=\text{id}$,
and continue to check the following conditions.  ``0'' means that we can find
one or more signed diagonal matrices, $V_\mathrm{sgn}$, to make all $N^{ij}_k$
non-negative.  We will use those $V_\mathrm{sgn}$ to continue to check the
following conditions. 

\item
``$\exists\ j$ that $\frac{S_{ij}}{S_{uj}} \geq 1$ '': The error code is ``3''
if $\frac{S_{ij}}{S_{uj}}$ is not real.  The error code is ``2'' if
$\frac{S_{ij}}{S_{uj}}$ is real and $-1<\frac{S_{ij}}{S_{uj}}<1$.  The error
code is ``1'' if $\frac{S_{ij}}{S_{uj}}$ is real and $\frac{S_{ij}}{S_{uj}}\leq
-1$.  The error code is ``0'' if $\frac{S_{ij}}{S_{uj}}$ is real and
$\frac{S_{ij}}{S_{uj}}\geq 1$.  Such an error code depends on $i,j$.  We first
maximize the error code against $i$. Then we minimize the maximized error code
against $j$.  Such maximized-minimized error code is displayed in the following
tables.

\item
``FS indicator'': $S,T$ matrices must satisfy the FS indicator condition in
Proposition B.1 (12) eqn. (B.10).  An error code $n>0$ is given if the
$n^\text{th}$ FS indicator condition is violated.  Error code 0 implies that
the FS indicator condition is satisfied.

\item
``$C$ = perm-mat'': $C=S^2$ must be a permutation matrix.

\end{enumerate}

We remark that in the following list, we only list the error codes for each
choice of $U_0$ and $u$.  In fact, for each choice of $U_0$ and $u$, there are
also many choices the signed diagonal matrices $V_\mathrm{sgn}$, which lead to
different $S,T$ matrices and different error codes.  In the following list, we
only display one set of error codes with minimal sum, among all different
choices of $V_\mathrm{sgn}$.

Since the minimized error codes are displayed, if a displayed error code is
non-zero for a give $(U_0,u)$, then for any choices of signed diagonal matrices
$V_\mathrm{sgn}$, some error codes must be non-zero.  If all the displayed
error codes are zero for a give $u$, then usually there are several signed
diagonal matrices $V_\mathrm{sgn}$, that will give rise to several pairs of
$S,T$ matrices that have no error.  The number of the inequivalent error-free
$S,T$ pairs, from different choices $u$'s and $V_\mathrm{sgn}$'s, is given
below each table and the total number of the inequivalent $S,T$ pairs is given
at the end of each list entry.  If the total number of the inequivalent $S,T$
pairs is zero, the corresponding irrep-sum representation will fail to produce
any modular data.  If the total number of the inequivalent $S,T$ pairs is not
zero, the corresponding $S,T$ pairs will be listed in next section and in the
Appendix C.2.

We note that the conditions (1) and (2) do not depend on the choices of unit
index $u$.  The conditions (1 -- 9) do not depend on the choices of the signed
diagonal matrix $V_\text{sd}$.

\

\newcommand{\grey}[1]{\color{grey} #1 \color{black} }

{\tiny 
\input{MDGT2ST6}
}

\

\section{A list of rank-6 $S,T$ matrices}
\label{Section3}

In the last section, we obtained all the $S,T$ matrices coming from certain
resolved dimension-6 $\SL$ representations (one particular representation from
each GT orbit).  Starting from those $S,T$ matrices, after applying Galois
conjugations, we obtain the following list of 174 $S,T$ matrices.  The $S,T$
matrices are group into orbits generated by Galois conjugations.  Each Galois
orbit is led by a black entry (with $d_i>0$) or a blue entry (with some
$d_i<0$), followed by grey entries.  The notations are explained in Appendix
C.2, where a shorten version of this list is presented.

An pair of $S,T$ matrices is called unitary if all its quantum dimensions are
positive: $d_i>0$, and is called pseudo-unitary if there is a signed diagonal
matrix $V_\mathrm{sgn}$ that can transform the pair into a valid pair $S',T'$
with positive quantum dimensions and satisfying the conditions Proposition B.1
on the modular data. (In fact, we only check the 8 conditions listed in the
last section).  Below each entry in the following list, we indicate if the
corresponding pair of $S,T$ is pseudo-unitary or not pseudo-unitary.  If it is
pseudo-unitary, we indicate the index of its corresponding unitary $S,T$
matrices.

\

\noindent1. ind = $(3 , 
3;5,
4
)_{1}^{1}$:\ \ 
$d_i$ = ($1.0$,
$1.0$,
$2.0$,
$2.0$,
$2.236$,
$2.236$) 

\vskip 0.7ex
\hangindent=3em \hangafter=1
$D^2=$ 20.0 = 
 $20$

\vskip 0.7ex
\hangindent=3em \hangafter=1
$T = ( 0,
0,
\frac{1}{5},
\frac{4}{5},
\frac{1}{4},
\frac{3}{4} )
$,

\vskip 0.7ex
\hangindent=3em \hangafter=1
$S$ = ($ 1$,
$ 1$,
$ 2$,
$ 2$,
$ \sqrt{5}$,
$ \sqrt{5}$;\ \ 
$ 1$,
$ 2$,
$ 2$,
$ -\sqrt{5}$,
$ -\sqrt{5}$;\ \ 
$ -1-\sqrt{5}$,
$ -1+\sqrt{5}$,
$0$,
$0$;\ \ 
$ -1-\sqrt{5}$,
$0$,
$0$;\ \ 
$ -\sqrt{5}$,
$ \sqrt{5}$;\ \ 
$ -\sqrt{5}$)

\vskip 1ex 
\color{grey}

\noindent2. ind = $(3 , 
3;5,
4
)_{1}^{3}$:\ \ 
$d_i$ = ($1.0$,
$1.0$,
$2.0$,
$2.0$,
$-2.236$,
$-2.236$) 

\vskip 0.7ex
\hangindent=3em \hangafter=1
$D^2=$ 20.0 = 
 $20$

\vskip 0.7ex
\hangindent=3em \hangafter=1
$T = ( 0,
0,
\frac{2}{5},
\frac{3}{5},
\frac{1}{4},
\frac{3}{4} )
$,

\vskip 0.7ex
\hangindent=3em \hangafter=1
$S$ = ($ 1$,
$ 1$,
$ 2$,
$ 2$,
$ -\sqrt{5}$,
$ -\sqrt{5}$;\ \ 
$ 1$,
$ 2$,
$ 2$,
$ \sqrt{5}$,
$ \sqrt{5}$;\ \ 
$ -1+\sqrt{5}$,
$ -1-\sqrt{5}$,
$0$,
$0$;\ \ 
$ -1+\sqrt{5}$,
$0$,
$0$;\ \ 
$ \sqrt{5}$,
$ -\sqrt{5}$;\ \ 
$ \sqrt{5}$)

Pseudo-unitary $\sim$  
$(3 , 
3;5,
4
)_{2}^{1}$

\vskip 1ex 

 \color{black} \vskip 2ex

\noindent3. ind = $(3 , 
3;5,
4
)_{2}^{1}$:\ \ 
$d_i$ = ($1.0$,
$1.0$,
$2.0$,
$2.0$,
$2.236$,
$2.236$) 

\vskip 0.7ex
\hangindent=3em \hangafter=1
$D^2=$ 20.0 = 
 $20$

\vskip 0.7ex
\hangindent=3em \hangafter=1
$T = ( 0,
0,
\frac{2}{5},
\frac{3}{5},
\frac{1}{4},
\frac{3}{4} )
$,

\vskip 0.7ex
\hangindent=3em \hangafter=1
$S$ = ($ 1$,
$ 1$,
$ 2$,
$ 2$,
$ \sqrt{5}$,
$ \sqrt{5}$;\ \ 
$ 1$,
$ 2$,
$ 2$,
$ -\sqrt{5}$,
$ -\sqrt{5}$;\ \ 
$ -1+\sqrt{5}$,
$ -1-\sqrt{5}$,
$0$,
$0$;\ \ 
$ -1+\sqrt{5}$,
$0$,
$0$;\ \ 
$ \sqrt{5}$,
$ -\sqrt{5}$;\ \ 
$ \sqrt{5}$)

\vskip 1ex 
\color{grey}

\noindent4. ind = $(3 , 
3;5,
4
)_{2}^{3}$:\ \ 
$d_i$ = ($1.0$,
$1.0$,
$2.0$,
$2.0$,
$-2.236$,
$-2.236$) 

\vskip 0.7ex
\hangindent=3em \hangafter=1
$D^2=$ 20.0 = 
 $20$

\vskip 0.7ex
\hangindent=3em \hangafter=1
$T = ( 0,
0,
\frac{1}{5},
\frac{4}{5},
\frac{1}{4},
\frac{3}{4} )
$,

\vskip 0.7ex
\hangindent=3em \hangafter=1
$S$ = ($ 1$,
$ 1$,
$ 2$,
$ 2$,
$ -\sqrt{5}$,
$ -\sqrt{5}$;\ \ 
$ 1$,
$ 2$,
$ 2$,
$ \sqrt{5}$,
$ \sqrt{5}$;\ \ 
$ -1-\sqrt{5}$,
$ -1+\sqrt{5}$,
$0$,
$0$;\ \ 
$ -1-\sqrt{5}$,
$0$,
$0$;\ \ 
$ -\sqrt{5}$,
$ \sqrt{5}$;\ \ 
$ -\sqrt{5}$)

Pseudo-unitary $\sim$  
$(3 , 
3;5,
4
)_{1}^{1}$

\vskip 1ex 

 \color{black} \vskip 2ex

\noindent5. ind = $(4 , 
2;15,
5
)_{1}^{1}$:\ \ 
$d_i$ = ($1.0$,
$1.0$,
$1.0$,
$1.618$,
$1.618$,
$1.618$) 

\vskip 0.7ex
\hangindent=3em \hangafter=1
$D^2=$ 10.854 = 
 $\frac{15+3\sqrt{5}}{2}$

\vskip 0.7ex
\hangindent=3em \hangafter=1
$T = ( 0,
\frac{1}{3},
\frac{1}{3},
\frac{2}{5},
\frac{11}{15},
\frac{11}{15} )
$,

\vskip 0.7ex
\hangindent=3em \hangafter=1
$S$ = ($ 1$,
$ 1$,
$ 1$,
$ \frac{1+\sqrt{5}}{2}$,
$ \frac{1+\sqrt{5}}{2}$,
$ \frac{1+\sqrt{5}}{2}$;\ \ 
$ \zeta_{3}^{1}$,
$ -\zeta_{6}^{1}$,
$ \frac{1+\sqrt{5}}{2}$,
$ \frac{1+\sqrt{5}}{2}\zeta_{3}^{1}$,
$ -\frac{1+\sqrt{5}}{2}\zeta_{6}^{1}$;\ \ 
$ \zeta_{3}^{1}$,
$ \frac{1+\sqrt{5}}{2}$,
$ -\frac{1+\sqrt{5}}{2}\zeta_{6}^{1}$,
$ \frac{1+\sqrt{5}}{2}\zeta_{3}^{1}$;\ \ 
$ -1$,
$ -1$,
$ -1$;\ \ 
$ -\zeta_{3}^{1}$,
$ \zeta_{6}^{1}$;\ \ 
$ -\zeta_{3}^{1}$)

\vskip 1ex 
\color{grey}

\noindent6. ind = $(4 , 
2;15,
5
)_{1}^{4}$:\ \ 
$d_i$ = ($1.0$,
$1.0$,
$1.0$,
$1.618$,
$1.618$,
$1.618$) 

\vskip 0.7ex
\hangindent=3em \hangafter=1
$D^2=$ 10.854 = 
 $\frac{15+3\sqrt{5}}{2}$

\vskip 0.7ex
\hangindent=3em \hangafter=1
$T = ( 0,
\frac{1}{3},
\frac{1}{3},
\frac{3}{5},
\frac{14}{15},
\frac{14}{15} )
$,

\vskip 0.7ex
\hangindent=3em \hangafter=1
$S$ = ($ 1$,
$ 1$,
$ 1$,
$ \frac{1+\sqrt{5}}{2}$,
$ \frac{1+\sqrt{5}}{2}$,
$ \frac{1+\sqrt{5}}{2}$;\ \ 
$ \zeta_{3}^{1}$,
$ -\zeta_{6}^{1}$,
$ \frac{1+\sqrt{5}}{2}$,
$ \frac{1+\sqrt{5}}{2}\zeta_{3}^{1}$,
$ -\frac{1+\sqrt{5}}{2}\zeta_{6}^{1}$;\ \ 
$ \zeta_{3}^{1}$,
$ \frac{1+\sqrt{5}}{2}$,
$ -\frac{1+\sqrt{5}}{2}\zeta_{6}^{1}$,
$ \frac{1+\sqrt{5}}{2}\zeta_{3}^{1}$;\ \ 
$ -1$,
$ -1$,
$ -1$;\ \ 
$ -\zeta_{3}^{1}$,
$ \zeta_{6}^{1}$;\ \ 
$ -\zeta_{3}^{1}$)

\vskip 1ex 
\color{grey}

\noindent7. ind = $(4 , 
2;15,
5
)_{1}^{11}$:\ \ 
$d_i$ = ($1.0$,
$1.0$,
$1.0$,
$1.618$,
$1.618$,
$1.618$) 

\vskip 0.7ex
\hangindent=3em \hangafter=1
$D^2=$ 10.854 = 
 $\frac{15+3\sqrt{5}}{2}$

\vskip 0.7ex
\hangindent=3em \hangafter=1
$T = ( 0,
\frac{2}{3},
\frac{2}{3},
\frac{2}{5},
\frac{1}{15},
\frac{1}{15} )
$,

\vskip 0.7ex
\hangindent=3em \hangafter=1
$S$ = ($ 1$,
$ 1$,
$ 1$,
$ \frac{1+\sqrt{5}}{2}$,
$ \frac{1+\sqrt{5}}{2}$,
$ \frac{1+\sqrt{5}}{2}$;\ \ 
$ -\zeta_{6}^{1}$,
$ \zeta_{3}^{1}$,
$ \frac{1+\sqrt{5}}{2}$,
$ -\frac{1+\sqrt{5}}{2}\zeta_{6}^{1}$,
$ \frac{1+\sqrt{5}}{2}\zeta_{3}^{1}$;\ \ 
$ -\zeta_{6}^{1}$,
$ \frac{1+\sqrt{5}}{2}$,
$ \frac{1+\sqrt{5}}{2}\zeta_{3}^{1}$,
$ -\frac{1+\sqrt{5}}{2}\zeta_{6}^{1}$;\ \ 
$ -1$,
$ -1$,
$ -1$;\ \ 
$ \zeta_{6}^{1}$,
$ -\zeta_{3}^{1}$;\ \ 
$ \zeta_{6}^{1}$)

\vskip 1ex 
\color{grey}

\noindent8. ind = $(4 , 
2;15,
5
)_{1}^{14}$:\ \ 
$d_i$ = ($1.0$,
$1.0$,
$1.0$,
$1.618$,
$1.618$,
$1.618$) 

\vskip 0.7ex
\hangindent=3em \hangafter=1
$D^2=$ 10.854 = 
 $\frac{15+3\sqrt{5}}{2}$

\vskip 0.7ex
\hangindent=3em \hangafter=1
$T = ( 0,
\frac{2}{3},
\frac{2}{3},
\frac{3}{5},
\frac{4}{15},
\frac{4}{15} )
$,

\vskip 0.7ex
\hangindent=3em \hangafter=1
$S$ = ($ 1$,
$ 1$,
$ 1$,
$ \frac{1+\sqrt{5}}{2}$,
$ \frac{1+\sqrt{5}}{2}$,
$ \frac{1+\sqrt{5}}{2}$;\ \ 
$ -\zeta_{6}^{1}$,
$ \zeta_{3}^{1}$,
$ \frac{1+\sqrt{5}}{2}$,
$ -\frac{1+\sqrt{5}}{2}\zeta_{6}^{1}$,
$ \frac{1+\sqrt{5}}{2}\zeta_{3}^{1}$;\ \ 
$ -\zeta_{6}^{1}$,
$ \frac{1+\sqrt{5}}{2}$,
$ \frac{1+\sqrt{5}}{2}\zeta_{3}^{1}$,
$ -\frac{1+\sqrt{5}}{2}\zeta_{6}^{1}$;\ \ 
$ -1$,
$ -1$,
$ -1$;\ \ 
$ \zeta_{6}^{1}$,
$ -\zeta_{3}^{1}$;\ \ 
$ \zeta_{6}^{1}$)

\vskip 1ex 
\color{grey}

\noindent9. ind = $(4 , 
2;15,
5
)_{1}^{13}$:\ \ 
$d_i$ = ($1.0$,
$1.0$,
$1.0$,
$-0.618$,
$-0.618$,
$-0.618$) 

\vskip 0.7ex
\hangindent=3em \hangafter=1
$D^2=$ 4.145 = 
 $\frac{15-3\sqrt{5}}{2}$

\vskip 0.7ex
\hangindent=3em \hangafter=1
$T = ( 0,
\frac{1}{3},
\frac{1}{3},
\frac{1}{5},
\frac{8}{15},
\frac{8}{15} )
$,

\vskip 0.7ex
\hangindent=3em \hangafter=1
$S$ = ($ 1$,
$ 1$,
$ 1$,
$ \frac{1-\sqrt{5}}{2}$,
$ \frac{1-\sqrt{5}}{2}$,
$ \frac{1-\sqrt{5}}{2}$;\ \ 
$ \zeta_{3}^{1}$,
$ -\zeta_{6}^{1}$,
$ \frac{1-\sqrt{5}}{2}$,
$ \frac{-1+\sqrt{5}}{2}\zeta_{6}^{1}$,
$ \frac{1-\sqrt{5}}{2}\zeta_{3}^{1}$;\ \ 
$ \zeta_{3}^{1}$,
$ \frac{1-\sqrt{5}}{2}$,
$ \frac{1-\sqrt{5}}{2}\zeta_{3}^{1}$,
$ \frac{-1+\sqrt{5}}{2}\zeta_{6}^{1}$;\ \ 
$ -1$,
$ -1$,
$ -1$;\ \ 
$ -\zeta_{3}^{1}$,
$ \zeta_{6}^{1}$;\ \ 
$ -\zeta_{3}^{1}$)

Not pseudo-unitary. 

\vskip 1ex 
\color{grey}

\noindent10. ind = $(4 , 
2;15,
5
)_{1}^{7}$:\ \ 
$d_i$ = ($1.0$,
$1.0$,
$1.0$,
$-0.618$,
$-0.618$,
$-0.618$) 

\vskip 0.7ex
\hangindent=3em \hangafter=1
$D^2=$ 4.145 = 
 $\frac{15-3\sqrt{5}}{2}$

\vskip 0.7ex
\hangindent=3em \hangafter=1
$T = ( 0,
\frac{1}{3},
\frac{1}{3},
\frac{4}{5},
\frac{2}{15},
\frac{2}{15} )
$,

\vskip 0.7ex
\hangindent=3em \hangafter=1
$S$ = ($ 1$,
$ 1$,
$ 1$,
$ \frac{1-\sqrt{5}}{2}$,
$ \frac{1-\sqrt{5}}{2}$,
$ \frac{1-\sqrt{5}}{2}$;\ \ 
$ \zeta_{3}^{1}$,
$ -\zeta_{6}^{1}$,
$ \frac{1-\sqrt{5}}{2}$,
$ \frac{-1+\sqrt{5}}{2}\zeta_{6}^{1}$,
$ \frac{1-\sqrt{5}}{2}\zeta_{3}^{1}$;\ \ 
$ \zeta_{3}^{1}$,
$ \frac{1-\sqrt{5}}{2}$,
$ \frac{1-\sqrt{5}}{2}\zeta_{3}^{1}$,
$ \frac{-1+\sqrt{5}}{2}\zeta_{6}^{1}$;\ \ 
$ -1$,
$ -1$,
$ -1$;\ \ 
$ -\zeta_{3}^{1}$,
$ \zeta_{6}^{1}$;\ \ 
$ -\zeta_{3}^{1}$)

Not pseudo-unitary. 

\vskip 1ex 
\color{grey}

\noindent11. ind = $(4 , 
2;15,
5
)_{1}^{8}$:\ \ 
$d_i$ = ($1.0$,
$1.0$,
$1.0$,
$-0.618$,
$-0.618$,
$-0.618$) 

\vskip 0.7ex
\hangindent=3em \hangafter=1
$D^2=$ 4.145 = 
 $\frac{15-3\sqrt{5}}{2}$

\vskip 0.7ex
\hangindent=3em \hangafter=1
$T = ( 0,
\frac{2}{3},
\frac{2}{3},
\frac{1}{5},
\frac{13}{15},
\frac{13}{15} )
$,

\vskip 0.7ex
\hangindent=3em \hangafter=1
$S$ = ($ 1$,
$ 1$,
$ 1$,
$ \frac{1-\sqrt{5}}{2}$,
$ \frac{1-\sqrt{5}}{2}$,
$ \frac{1-\sqrt{5}}{2}$;\ \ 
$ -\zeta_{6}^{1}$,
$ \zeta_{3}^{1}$,
$ \frac{1-\sqrt{5}}{2}$,
$ \frac{1-\sqrt{5}}{2}\zeta_{3}^{1}$,
$ \frac{-1+\sqrt{5}}{2}\zeta_{6}^{1}$;\ \ 
$ -\zeta_{6}^{1}$,
$ \frac{1-\sqrt{5}}{2}$,
$ \frac{-1+\sqrt{5}}{2}\zeta_{6}^{1}$,
$ \frac{1-\sqrt{5}}{2}\zeta_{3}^{1}$;\ \ 
$ -1$,
$ -1$,
$ -1$;\ \ 
$ \zeta_{6}^{1}$,
$ -\zeta_{3}^{1}$;\ \ 
$ \zeta_{6}^{1}$)

Not pseudo-unitary. 

\vskip 1ex 
\color{grey}

\noindent12. ind = $(4 , 
2;15,
5
)_{1}^{2}$:\ \ 
$d_i$ = ($1.0$,
$1.0$,
$1.0$,
$-0.618$,
$-0.618$,
$-0.618$) 

\vskip 0.7ex
\hangindent=3em \hangafter=1
$D^2=$ 4.145 = 
 $\frac{15-3\sqrt{5}}{2}$

\vskip 0.7ex
\hangindent=3em \hangafter=1
$T = ( 0,
\frac{2}{3},
\frac{2}{3},
\frac{4}{5},
\frac{7}{15},
\frac{7}{15} )
$,

\vskip 0.7ex
\hangindent=3em \hangafter=1
$S$ = ($ 1$,
$ 1$,
$ 1$,
$ \frac{1-\sqrt{5}}{2}$,
$ \frac{1-\sqrt{5}}{2}$,
$ \frac{1-\sqrt{5}}{2}$;\ \ 
$ -\zeta_{6}^{1}$,
$ \zeta_{3}^{1}$,
$ \frac{1-\sqrt{5}}{2}$,
$ \frac{1-\sqrt{5}}{2}\zeta_{3}^{1}$,
$ \frac{-1+\sqrt{5}}{2}\zeta_{6}^{1}$;\ \ 
$ -\zeta_{6}^{1}$,
$ \frac{1-\sqrt{5}}{2}$,
$ \frac{-1+\sqrt{5}}{2}\zeta_{6}^{1}$,
$ \frac{1-\sqrt{5}}{2}\zeta_{3}^{1}$;\ \ 
$ -1$,
$ -1$,
$ -1$;\ \ 
$ \zeta_{6}^{1}$,
$ -\zeta_{3}^{1}$;\ \ 
$ \zeta_{6}^{1}$)

Not pseudo-unitary. 

\vskip 1ex 

 \color{black} \vskip 2ex

\noindent13. ind = $(4 , 
2;7,
3
)_{1}^{1}$:\ \ 
$d_i$ = ($1.0$,
$3.791$,
$3.791$,
$3.791$,
$4.791$,
$5.791$) 

\vskip 0.7ex
\hangindent=3em \hangafter=1
$D^2=$ 100.617 = 
 $\frac{105+21\sqrt{21}}{2}$

\vskip 0.7ex
\hangindent=3em \hangafter=1
$T = ( 0,
\frac{1}{7},
\frac{2}{7},
\frac{4}{7},
0,
\frac{2}{3} )
$,

\vskip 0.7ex
\hangindent=3em \hangafter=1
$S$ = ($ 1$,
$ \frac{3+\sqrt{21}}{2}$,
$ \frac{3+\sqrt{21}}{2}$,
$ \frac{3+\sqrt{21}}{2}$,
$ \frac{5+\sqrt{21}}{2}$,
$ \frac{7+\sqrt{21}}{2}$;\ \ 
$ 2-c^{1}_{21}
-2c^{2}_{21}
+3c^{3}_{21}
+2c^{4}_{21}
-2c^{5}_{21}
$,
$ -c^{2}_{21}
-2c^{3}_{21}
-c^{4}_{21}
+c^{5}_{21}
$,
$ -1+2c^{1}_{21}
+3c^{2}_{21}
-c^{3}_{21}
+2c^{5}_{21}
$,
$ -\frac{3+\sqrt{21}}{2}$,
$0$;\ \ 
$ -1+2c^{1}_{21}
+3c^{2}_{21}
-c^{3}_{21}
+2c^{5}_{21}
$,
$ 2-c^{1}_{21}
-2c^{2}_{21}
+3c^{3}_{21}
+2c^{4}_{21}
-2c^{5}_{21}
$,
$ -\frac{3+\sqrt{21}}{2}$,
$0$;\ \ 
$ -c^{2}_{21}
-2c^{3}_{21}
-c^{4}_{21}
+c^{5}_{21}
$,
$ -\frac{3+\sqrt{21}}{2}$,
$0$;\ \ 
$ 1$,
$ \frac{7+\sqrt{21}}{2}$;\ \ 
$ -\frac{7+\sqrt{21}}{2}$)

\vskip 1ex 
\color{grey}

\noindent14. ind = $(4 , 
2;7,
3
)_{1}^{5}$:\ \ 
$d_i$ = ($1.0$,
$3.791$,
$3.791$,
$3.791$,
$4.791$,
$5.791$) 

\vskip 0.7ex
\hangindent=3em \hangafter=1
$D^2=$ 100.617 = 
 $\frac{105+21\sqrt{21}}{2}$

\vskip 0.7ex
\hangindent=3em \hangafter=1
$T = ( 0,
\frac{3}{7},
\frac{5}{7},
\frac{6}{7},
0,
\frac{1}{3} )
$,

\vskip 0.7ex
\hangindent=3em \hangafter=1
$S$ = ($ 1$,
$ \frac{3+\sqrt{21}}{2}$,
$ \frac{3+\sqrt{21}}{2}$,
$ \frac{3+\sqrt{21}}{2}$,
$ \frac{5+\sqrt{21}}{2}$,
$ \frac{7+\sqrt{21}}{2}$;\ \ 
$ -c^{2}_{21}
-2c^{3}_{21}
-c^{4}_{21}
+c^{5}_{21}
$,
$ 2-c^{1}_{21}
-2c^{2}_{21}
+3c^{3}_{21}
+2c^{4}_{21}
-2c^{5}_{21}
$,
$ -1+2c^{1}_{21}
+3c^{2}_{21}
-c^{3}_{21}
+2c^{5}_{21}
$,
$ -\frac{3+\sqrt{21}}{2}$,
$0$;\ \ 
$ -1+2c^{1}_{21}
+3c^{2}_{21}
-c^{3}_{21}
+2c^{5}_{21}
$,
$ -c^{2}_{21}
-2c^{3}_{21}
-c^{4}_{21}
+c^{5}_{21}
$,
$ -\frac{3+\sqrt{21}}{2}$,
$0$;\ \ 
$ 2-c^{1}_{21}
-2c^{2}_{21}
+3c^{3}_{21}
+2c^{4}_{21}
-2c^{5}_{21}
$,
$ -\frac{3+\sqrt{21}}{2}$,
$0$;\ \ 
$ 1$,
$ \frac{7+\sqrt{21}}{2}$;\ \ 
$ -\frac{7+\sqrt{21}}{2}$)

\vskip 1ex 
\color{grey}

\noindent15. ind = $(4 , 
2;7,
3
)_{1}^{2}$:\ \ 
$d_i$ = ($1.0$,
$0.208$,
$1.208$,
$-0.791$,
$-0.791$,
$-0.791$) 

\vskip 0.7ex
\hangindent=3em \hangafter=1
$D^2=$ 4.382 = 
 $\frac{105-21\sqrt{21}}{2}$

\vskip 0.7ex
\hangindent=3em \hangafter=1
$T = ( 0,
0,
\frac{1}{3},
\frac{1}{7},
\frac{2}{7},
\frac{4}{7} )
$,

\vskip 0.7ex
\hangindent=3em \hangafter=1
$S$ = ($ 1$,
$ \frac{5-\sqrt{21}}{2}$,
$ \frac{7-\sqrt{21}}{2}$,
$ \frac{3-\sqrt{21}}{2}$,
$ \frac{3-\sqrt{21}}{2}$,
$ \frac{3-\sqrt{21}}{2}$;\ \ 
$ 1$,
$ \frac{7-\sqrt{21}}{2}$,
$ \frac{-3+\sqrt{21}}{2}$,
$ \frac{-3+\sqrt{21}}{2}$,
$ \frac{-3+\sqrt{21}}{2}$;\ \ 
$ \frac{-7+\sqrt{21}}{2}$,
$0$,
$0$,
$0$;\ \ 
$ 1+c^{1}_{21}
-c^{2}_{21}
-2c^{4}_{21}
-c^{5}_{21}
$,
$ c^{2}_{21}
-c^{3}_{21}
+c^{4}_{21}
-c^{5}_{21}
$,
$ 1-2c^{1}_{21}
+c^{3}_{21}
+c^{5}_{21}
$;\ \ 
$ 1-2c^{1}_{21}
+c^{3}_{21}
+c^{5}_{21}
$,
$ 1+c^{1}_{21}
-c^{2}_{21}
-2c^{4}_{21}
-c^{5}_{21}
$;\ \ 
$ c^{2}_{21}
-c^{3}_{21}
+c^{4}_{21}
-c^{5}_{21}
$)

Not pseudo-unitary. 

\vskip 1ex 
\color{grey}

\noindent16. ind = $(4 , 
2;7,
3
)_{1}^{10}$:\ \ 
$d_i$ = ($1.0$,
$0.208$,
$1.208$,
$-0.791$,
$-0.791$,
$-0.791$) 

\vskip 0.7ex
\hangindent=3em \hangafter=1
$D^2=$ 4.382 = 
 $\frac{105-21\sqrt{21}}{2}$

\vskip 0.7ex
\hangindent=3em \hangafter=1
$T = ( 0,
0,
\frac{2}{3},
\frac{3}{7},
\frac{5}{7},
\frac{6}{7} )
$,

\vskip 0.7ex
\hangindent=3em \hangafter=1
$S$ = ($ 1$,
$ \frac{5-\sqrt{21}}{2}$,
$ \frac{7-\sqrt{21}}{2}$,
$ \frac{3-\sqrt{21}}{2}$,
$ \frac{3-\sqrt{21}}{2}$,
$ \frac{3-\sqrt{21}}{2}$;\ \ 
$ 1$,
$ \frac{7-\sqrt{21}}{2}$,
$ \frac{-3+\sqrt{21}}{2}$,
$ \frac{-3+\sqrt{21}}{2}$,
$ \frac{-3+\sqrt{21}}{2}$;\ \ 
$ \frac{-7+\sqrt{21}}{2}$,
$0$,
$0$,
$0$;\ \ 
$ c^{2}_{21}
-c^{3}_{21}
+c^{4}_{21}
-c^{5}_{21}
$,
$ 1+c^{1}_{21}
-c^{2}_{21}
-2c^{4}_{21}
-c^{5}_{21}
$,
$ 1-2c^{1}_{21}
+c^{3}_{21}
+c^{5}_{21}
$;\ \ 
$ 1-2c^{1}_{21}
+c^{3}_{21}
+c^{5}_{21}
$,
$ c^{2}_{21}
-c^{3}_{21}
+c^{4}_{21}
-c^{5}_{21}
$;\ \ 
$ 1+c^{1}_{21}
-c^{2}_{21}
-2c^{4}_{21}
-c^{5}_{21}
$)

Not pseudo-unitary. 

\vskip 1ex 

 \color{black} \vskip 2ex 
\color{blue}

\noindent17. ind = $(6;9
)_{1}^{1}$:\ \ 
$d_i$ = ($1.0$,
$0.347$,
$1.0$,
$1.532$,
$-1.0$,
$-1.879$) 

\vskip 0.7ex
\hangindent=3em \hangafter=1
$D^2=$ 9.0 = 
 $9$

\vskip 0.7ex
\hangindent=3em \hangafter=1
$T = ( 0,
\frac{1}{9},
\frac{2}{3},
\frac{4}{9},
\frac{1}{3},
\frac{7}{9} )
$,

\vskip 0.7ex
\hangindent=3em \hangafter=1
$S$ = ($ 1$,
$ c^{2}_{9}
$,
$ 1$,
$ c^{1}_{9}
$,
$ -1$,
$  c_9^4 $;\ \ 
$ 1$,
$ c^{1}_{9}
$,
$ 1$,
$  -c_9^4 $,
$ 1$;\ \ 
$ 1$,
$  c_9^4 $,
$ -1$,
$ c^{2}_{9}
$;\ \ 
$ 1$,
$ -c^{2}_{9}
$,
$ 1$;\ \ 
$ 1$,
$ -c^{1}_{9}
$;\ \ 
$ 1$)

Not pseudo-unitary. 

\vskip 1ex 
\color{grey}

\noindent18. ind = $(6;9
)_{1}^{2}$:\ \ 
$d_i$ = ($1.0$,
$0.347$,
$1.0$,
$1.532$,
$-1.0$,
$-1.879$) 

\vskip 0.7ex
\hangindent=3em \hangafter=1
$D^2=$ 9.0 = 
 $9$

\vskip 0.7ex
\hangindent=3em \hangafter=1
$T = ( 0,
\frac{8}{9},
\frac{1}{3},
\frac{5}{9},
\frac{2}{3},
\frac{2}{9} )
$,

\vskip 0.7ex
\hangindent=3em \hangafter=1
$S$ = ($ 1$,
$ c^{2}_{9}
$,
$ 1$,
$ c^{1}_{9}
$,
$ -1$,
$  c_9^4 $;\ \ 
$ 1$,
$ c^{1}_{9}
$,
$ 1$,
$  -c_9^4 $,
$ 1$;\ \ 
$ 1$,
$  c_9^4 $,
$ -1$,
$ c^{2}_{9}
$;\ \ 
$ 1$,
$ -c^{2}_{9}
$,
$ 1$;\ \ 
$ 1$,
$ -c^{1}_{9}
$;\ \ 
$ 1$)

Not pseudo-unitary. 

\vskip 1ex 

 \color{black} \vskip 2ex

\noindent19. ind = $(6;13
)_{1}^{1}$:\ \ 
$d_i$ = ($1.0$,
$1.941$,
$2.770$,
$3.438$,
$3.907$,
$4.148$) 

\vskip 0.7ex
\hangindent=3em \hangafter=1
$D^2=$ 56.746 = 
 $21+15c^{1}_{13}
+10c^{2}_{13}
+6c^{3}_{13}
+3c^{4}_{13}
+c^{5}_{13}
$

\vskip 0.7ex
\hangindent=3em \hangafter=1
$T = ( 0,
\frac{4}{13},
\frac{2}{13},
\frac{7}{13},
\frac{6}{13},
\frac{12}{13} )
$,

\vskip 0.7ex
\hangindent=3em \hangafter=1
$S$ = ($ 1$,
$ \xi_{13}^{2}$,
$ \xi_{13}^{3}$,
$ \xi_{13}^{4}$,
$ \xi_{13}^{5}$,
$ \xi_{13}^{6}$;\ \ 
$ -\xi_{13}^{4}$,
$ \xi_{13}^{6}$,
$ -\xi_{13}^{5}$,
$ \xi_{13}^{3}$,
$ -1$;\ \ 
$ \xi_{13}^{4}$,
$ 1$,
$ -\xi_{13}^{2}$,
$ -\xi_{13}^{5}$;\ \ 
$ \xi_{13}^{3}$,
$ -\xi_{13}^{6}$,
$ \xi_{13}^{2}$;\ \ 
$ -1$,
$ \xi_{13}^{4}$;\ \ 
$ -\xi_{13}^{3}$)

\vskip 1ex 
\color{grey}

\noindent20. ind = $(6;13
)_{1}^{12}$:\ \ 
$d_i$ = ($1.0$,
$1.941$,
$2.770$,
$3.438$,
$3.907$,
$4.148$) 

\vskip 0.7ex
\hangindent=3em \hangafter=1
$D^2=$ 56.746 = 
 $21+15c^{1}_{13}
+10c^{2}_{13}
+6c^{3}_{13}
+3c^{4}_{13}
+c^{5}_{13}
$

\vskip 0.7ex
\hangindent=3em \hangafter=1
$T = ( 0,
\frac{9}{13},
\frac{11}{13},
\frac{6}{13},
\frac{7}{13},
\frac{1}{13} )
$,

\vskip 0.7ex
\hangindent=3em \hangafter=1
$S$ = ($ 1$,
$ \xi_{13}^{2}$,
$ \xi_{13}^{3}$,
$ \xi_{13}^{4}$,
$ \xi_{13}^{5}$,
$ \xi_{13}^{6}$;\ \ 
$ -\xi_{13}^{4}$,
$ \xi_{13}^{6}$,
$ -\xi_{13}^{5}$,
$ \xi_{13}^{3}$,
$ -1$;\ \ 
$ \xi_{13}^{4}$,
$ 1$,
$ -\xi_{13}^{2}$,
$ -\xi_{13}^{5}$;\ \ 
$ \xi_{13}^{3}$,
$ -\xi_{13}^{6}$,
$ \xi_{13}^{2}$;\ \ 
$ -1$,
$ \xi_{13}^{4}$;\ \ 
$ -\xi_{13}^{3}$)

\vskip 1ex 
\color{grey}

\noindent21. ind = $(6;13
)_{1}^{11}$:\ \ 
$d_i$ = ($1.0$,
$1.426$,
$2.136$,
$-0.514$,
$-1.770$,
$-2.11$) 

\vskip 0.7ex
\hangindent=3em \hangafter=1
$D^2=$ 15.48 = 
 $15-6c^{1}_{13}
+9c^{2}_{13}
-5c^{3}_{13}
+4c^{4}_{13}
-3c^{5}_{13}
$

\vskip 0.7ex
\hangindent=3em \hangafter=1
$T = ( 0,
\frac{1}{13},
\frac{9}{13},
\frac{2}{13},
\frac{5}{13},
\frac{12}{13} )
$,

\vskip 0.7ex
\hangindent=3em \hangafter=1
$S$ = ($ 1$,
$ \xi_{13}^{3,2}$,
$ \xi_{13}^{6,2}$,
$ -\xi_{13}^{1,2}$,
$ -c^{1}_{13}
$,
$ -\xi_{13}^{5,2}$;\ \ 
$ -1$,
$ c^{1}_{13}
$,
$ -\xi_{13}^{5,2}$,
$ \xi_{13}^{6,2}$,
$ \xi_{13}^{1,2}$;\ \ 
$ -\xi_{13}^{5,2}$,
$ -\xi_{13}^{3,2}$,
$ -\xi_{13}^{1,2}$,
$ 1$;\ \ 
$ -\xi_{13}^{6,2}$,
$ -1$,
$ -c^{1}_{13}
$;\ \ 
$ \xi_{13}^{5,2}$,
$ -\xi_{13}^{3,2}$;\ \ 
$ \xi_{13}^{6,2}$)

Not pseudo-unitary. 

\vskip 1ex 
\color{grey}

\noindent22. ind = $(6;13
)_{1}^{6}$:\ \ 
$d_i$ = ($1.0$,
$0.468$,
$0.829$,
$-0.241$,
$-0.667$,
$-0.941$) 

\vskip 0.7ex
\hangindent=3em \hangafter=1
$D^2=$ 3.297 = 
 $6-5c^{1}_{13}
-12c^{2}_{13}
-15c^{3}_{13}
-14c^{4}_{13}
-9c^{5}_{13}
$

\vskip 0.7ex
\hangindent=3em \hangafter=1
$T = ( 0,
\frac{3}{13},
\frac{10}{13},
\frac{11}{13},
\frac{7}{13},
\frac{12}{13} )
$,

\vskip 0.7ex
\hangindent=3em \hangafter=1
$S$ = ($ 1$,
$ \xi_{13}^{2,6}$,
$ \xi_{13}^{4,6}$,
$ -c^{3}_{13}
$,
$ -\xi_{13}^{3,6}$,
$ -\xi_{13}^{5,6}$;\ \ 
$ -\xi_{13}^{5,6}$,
$ \xi_{13}^{3,6}$,
$ -\xi_{13}^{4,6}$,
$ -c^{3}_{13}
$,
$ 1$;\ \ 
$ -1$,
$ -\xi_{13}^{5,6}$,
$ \xi_{13}^{2,6}$,
$ c^{3}_{13}
$;\ \ 
$ -\xi_{13}^{2,6}$,
$ -1$,
$ -\xi_{13}^{3,6}$;\ \ 
$ \xi_{13}^{5,6}$,
$ -\xi_{13}^{4,6}$;\ \ 
$ \xi_{13}^{2,6}$)

Not pseudo-unitary. 

\vskip 1ex 
\color{grey}

\noindent23. ind = $(6;13
)_{1}^{9}$:\ \ 
$d_i$ = ($1.0$,
$0.290$,
$0.564$,
$0.805$,
$-1.136$,
$-1.206$) 

\vskip 0.7ex
\hangindent=3em \hangafter=1
$D^2=$ 4.798 = 
 $20+5c^{1}_{13}
-c^{2}_{13}
+2c^{3}_{13}
+14c^{4}_{13}
+9c^{5}_{13}
$

\vskip 0.7ex
\hangindent=3em \hangafter=1
$T = ( 0,
\frac{5}{13},
\frac{4}{13},
\frac{11}{13},
\frac{10}{13},
\frac{2}{13} )
$,

\vskip 0.7ex
\hangindent=3em \hangafter=1
$S$ = ($ 1$,
$ \xi_{13}^{1,4}$,
$ \xi_{13}^{2,4}$,
$ \xi_{13}^{3,4}$,
$ -c^{2}_{13}
$,
$ -\xi_{13}^{6,4}$;\ \ 
$ \xi_{13}^{3,4}$,
$ \xi_{13}^{6,4}$,
$ 1$,
$ \xi_{13}^{2,4}$,
$ c^{2}_{13}
$;\ \ 
$ -\xi_{13}^{1,4}$,
$ -c^{2}_{13}
$,
$ -1$,
$ \xi_{13}^{3,4}$;\ \ 
$ \xi_{13}^{1,4}$,
$ \xi_{13}^{6,4}$,
$ -\xi_{13}^{2,4}$;\ \ 
$ -\xi_{13}^{3,4}$,
$ \xi_{13}^{1,4}$;\ \ 
$ -1$)

Not pseudo-unitary. 

\vskip 1ex 
\color{grey}

\noindent24. ind = $(6;13
)_{1}^{10}$:\ \ 
$d_i$ = ($1.0$,
$0.360$,
$1.241$,
$1.497$,
$-0.700$,
$-1.410$) 

\vskip 0.7ex
\hangindent=3em \hangafter=1
$D^2=$ 7.390 = 
 $11-7c^{1}_{13}
-9c^{2}_{13}
+5c^{3}_{13}
-4c^{4}_{13}
-10c^{5}_{13}
$

\vskip 0.7ex
\hangindent=3em \hangafter=1
$T = ( 0,
\frac{5}{13},
\frac{7}{13},
\frac{1}{13},
\frac{8}{13},
\frac{3}{13} )
$,

\vskip 0.7ex
\hangindent=3em \hangafter=1
$S$ = ($ 1$,
$ \xi_{13}^{1,3}$,
$ \xi_{13}^{4,3}$,
$ -c^{5}_{13}
$,
$ -\xi_{13}^{2,3}$,
$ -\xi_{13}^{5,3}$;\ \ 
$ \xi_{13}^{4,3}$,
$ 1$,
$ \xi_{13}^{2,3}$,
$ \xi_{13}^{5,3}$,
$ -c^{5}_{13}
$;\ \ 
$ \xi_{13}^{1,3}$,
$ -\xi_{13}^{5,3}$,
$ c^{5}_{13}
$,
$ \xi_{13}^{2,3}$;\ \ 
$ -\xi_{13}^{1,3}$,
$ \xi_{13}^{4,3}$,
$ -1$;\ \ 
$ -1$,
$ \xi_{13}^{1,3}$;\ \ 
$ -\xi_{13}^{4,3}$)

Not pseudo-unitary. 

\vskip 1ex 
\color{grey}

\noindent25. ind = $(6;13
)_{1}^{8}$:\ \ 
$d_i$ = ($1.0$,
$0.709$,
$0.880$,
$-0.255$,
$-0.497$,
$-1.61$) 

\vskip 0.7ex
\hangindent=3em \hangafter=1
$D^2=$ 3.717 = 
 $18-2c^{1}_{13}
+3c^{2}_{13}
+7c^{3}_{13}
-3c^{4}_{13}
+12c^{5}_{13}
$

\vskip 0.7ex
\hangindent=3em \hangafter=1
$T = ( 0,
\frac{6}{13},
\frac{5}{13},
\frac{9}{13},
\frac{3}{13},
\frac{4}{13} )
$,

\vskip 0.7ex
\hangindent=3em \hangafter=1
$S$ = ($ 1$,
$ -c^{4}_{13}
$,
$ \xi_{13}^{4,5}$,
$ -\xi_{13}^{1,5}$,
$ -\xi_{13}^{2,5}$,
$ -\xi_{13}^{6,5}$;\ \ 
$ \xi_{13}^{6,5}$,
$ -1$,
$ -\xi_{13}^{2,5}$,
$ \xi_{13}^{4,5}$,
$ \xi_{13}^{1,5}$;\ \ 
$ \xi_{13}^{2,5}$,
$ -\xi_{13}^{6,5}$,
$ \xi_{13}^{1,5}$,
$ -c^{4}_{13}
$;\ \ 
$ -1$,
$ c^{4}_{13}
$,
$ -\xi_{13}^{4,5}$;\ \ 
$ -\xi_{13}^{6,5}$,
$ 1$;\ \ 
$ -\xi_{13}^{2,5}$)

Not pseudo-unitary. 

\vskip 1ex 
\color{grey}

\noindent26. ind = $(6;13
)_{1}^{5}$:\ \ 
$d_i$ = ($1.0$,
$0.709$,
$0.880$,
$-0.255$,
$-0.497$,
$-1.61$) 

\vskip 0.7ex
\hangindent=3em \hangafter=1
$D^2=$ 3.717 = 
 $18-2c^{1}_{13}
+3c^{2}_{13}
+7c^{3}_{13}
-3c^{4}_{13}
+12c^{5}_{13}
$

\vskip 0.7ex
\hangindent=3em \hangafter=1
$T = ( 0,
\frac{7}{13},
\frac{8}{13},
\frac{4}{13},
\frac{10}{13},
\frac{9}{13} )
$,

\vskip 0.7ex
\hangindent=3em \hangafter=1
$S$ = ($ 1$,
$ -c^{4}_{13}
$,
$ \xi_{13}^{4,5}$,
$ -\xi_{13}^{1,5}$,
$ -\xi_{13}^{2,5}$,
$ -\xi_{13}^{6,5}$;\ \ 
$ \xi_{13}^{6,5}$,
$ -1$,
$ -\xi_{13}^{2,5}$,
$ \xi_{13}^{4,5}$,
$ \xi_{13}^{1,5}$;\ \ 
$ \xi_{13}^{2,5}$,
$ -\xi_{13}^{6,5}$,
$ \xi_{13}^{1,5}$,
$ -c^{4}_{13}
$;\ \ 
$ -1$,
$ c^{4}_{13}
$,
$ -\xi_{13}^{4,5}$;\ \ 
$ -\xi_{13}^{6,5}$,
$ 1$;\ \ 
$ -\xi_{13}^{2,5}$)

Not pseudo-unitary. 

\vskip 1ex 
\color{grey}

\noindent27. ind = $(6;13
)_{1}^{3}$:\ \ 
$d_i$ = ($1.0$,
$0.360$,
$1.241$,
$1.497$,
$-0.700$,
$-1.410$) 

\vskip 0.7ex
\hangindent=3em \hangafter=1
$D^2=$ 7.390 = 
 $11-7c^{1}_{13}
-9c^{2}_{13}
+5c^{3}_{13}
-4c^{4}_{13}
-10c^{5}_{13}
$

\vskip 0.7ex
\hangindent=3em \hangafter=1
$T = ( 0,
\frac{8}{13},
\frac{6}{13},
\frac{12}{13},
\frac{5}{13},
\frac{10}{13} )
$,

\vskip 0.7ex
\hangindent=3em \hangafter=1
$S$ = ($ 1$,
$ \xi_{13}^{1,3}$,
$ \xi_{13}^{4,3}$,
$ -c^{5}_{13}
$,
$ -\xi_{13}^{2,3}$,
$ -\xi_{13}^{5,3}$;\ \ 
$ \xi_{13}^{4,3}$,
$ 1$,
$ \xi_{13}^{2,3}$,
$ \xi_{13}^{5,3}$,
$ -c^{5}_{13}
$;\ \ 
$ \xi_{13}^{1,3}$,
$ -\xi_{13}^{5,3}$,
$ c^{5}_{13}
$,
$ \xi_{13}^{2,3}$;\ \ 
$ -\xi_{13}^{1,3}$,
$ \xi_{13}^{4,3}$,
$ -1$;\ \ 
$ -1$,
$ \xi_{13}^{1,3}$;\ \ 
$ -\xi_{13}^{4,3}$)

Not pseudo-unitary. 

\vskip 1ex 
\color{grey}

\noindent28. ind = $(6;13
)_{1}^{4}$:\ \ 
$d_i$ = ($1.0$,
$0.290$,
$0.564$,
$0.805$,
$-1.136$,
$-1.206$) 

\vskip 0.7ex
\hangindent=3em \hangafter=1
$D^2=$ 4.798 = 
 $20+5c^{1}_{13}
-c^{2}_{13}
+2c^{3}_{13}
+14c^{4}_{13}
+9c^{5}_{13}
$

\vskip 0.7ex
\hangindent=3em \hangafter=1
$T = ( 0,
\frac{8}{13},
\frac{9}{13},
\frac{2}{13},
\frac{3}{13},
\frac{11}{13} )
$,

\vskip 0.7ex
\hangindent=3em \hangafter=1
$S$ = ($ 1$,
$ \xi_{13}^{1,4}$,
$ \xi_{13}^{2,4}$,
$ \xi_{13}^{3,4}$,
$ -c^{2}_{13}
$,
$ -\xi_{13}^{6,4}$;\ \ 
$ \xi_{13}^{3,4}$,
$ \xi_{13}^{6,4}$,
$ 1$,
$ \xi_{13}^{2,4}$,
$ c^{2}_{13}
$;\ \ 
$ -\xi_{13}^{1,4}$,
$ -c^{2}_{13}
$,
$ -1$,
$ \xi_{13}^{3,4}$;\ \ 
$ \xi_{13}^{1,4}$,
$ \xi_{13}^{6,4}$,
$ -\xi_{13}^{2,4}$;\ \ 
$ -\xi_{13}^{3,4}$,
$ \xi_{13}^{1,4}$;\ \ 
$ -1$)

Not pseudo-unitary. 

\vskip 1ex 
\color{grey}

\noindent29. ind = $(6;13
)_{1}^{7}$:\ \ 
$d_i$ = ($1.0$,
$0.468$,
$0.829$,
$-0.241$,
$-0.667$,
$-0.941$) 

\vskip 0.7ex
\hangindent=3em \hangafter=1
$D^2=$ 3.297 = 
 $6-5c^{1}_{13}
-12c^{2}_{13}
-15c^{3}_{13}
-14c^{4}_{13}
-9c^{5}_{13}
$

\vskip 0.7ex
\hangindent=3em \hangafter=1
$T = ( 0,
\frac{10}{13},
\frac{3}{13},
\frac{2}{13},
\frac{6}{13},
\frac{1}{13} )
$,

\vskip 0.7ex
\hangindent=3em \hangafter=1
$S$ = ($ 1$,
$ \xi_{13}^{2,6}$,
$ \xi_{13}^{4,6}$,
$ -c^{3}_{13}
$,
$ -\xi_{13}^{3,6}$,
$ -\xi_{13}^{5,6}$;\ \ 
$ -\xi_{13}^{5,6}$,
$ \xi_{13}^{3,6}$,
$ -\xi_{13}^{4,6}$,
$ -c^{3}_{13}
$,
$ 1$;\ \ 
$ -1$,
$ -\xi_{13}^{5,6}$,
$ \xi_{13}^{2,6}$,
$ c^{3}_{13}
$;\ \ 
$ -\xi_{13}^{2,6}$,
$ -1$,
$ -\xi_{13}^{3,6}$;\ \ 
$ \xi_{13}^{5,6}$,
$ -\xi_{13}^{4,6}$;\ \ 
$ \xi_{13}^{2,6}$)

Not pseudo-unitary. 

\vskip 1ex 
\color{grey}

\noindent30. ind = $(6;13
)_{1}^{2}$:\ \ 
$d_i$ = ($1.0$,
$1.426$,
$2.136$,
$-0.514$,
$-1.770$,
$-2.11$) 

\vskip 0.7ex
\hangindent=3em \hangafter=1
$D^2=$ 15.48 = 
 $15-6c^{1}_{13}
+9c^{2}_{13}
-5c^{3}_{13}
+4c^{4}_{13}
-3c^{5}_{13}
$

\vskip 0.7ex
\hangindent=3em \hangafter=1
$T = ( 0,
\frac{12}{13},
\frac{4}{13},
\frac{11}{13},
\frac{8}{13},
\frac{1}{13} )
$,

\vskip 0.7ex
\hangindent=3em \hangafter=1
$S$ = ($ 1$,
$ \xi_{13}^{3,2}$,
$ \xi_{13}^{6,2}$,
$ -\xi_{13}^{1,2}$,
$ -c^{1}_{13}
$,
$ -\xi_{13}^{5,2}$;\ \ 
$ -1$,
$ c^{1}_{13}
$,
$ -\xi_{13}^{5,2}$,
$ \xi_{13}^{6,2}$,
$ \xi_{13}^{1,2}$;\ \ 
$ -\xi_{13}^{5,2}$,
$ -\xi_{13}^{3,2}$,
$ -\xi_{13}^{1,2}$,
$ 1$;\ \ 
$ -\xi_{13}^{6,2}$,
$ -1$,
$ -c^{1}_{13}
$;\ \ 
$ \xi_{13}^{5,2}$,
$ -\xi_{13}^{3,2}$;\ \ 
$ \xi_{13}^{6,2}$)

Not pseudo-unitary. 

\vskip 1ex 

 \color{black} \vskip 2ex

\noindent31. ind = $(6;16
)_{1}^{1}$:\ \ 
$d_i$ = ($1.0$,
$1.0$,
$1.0$,
$1.0$,
$1.414$,
$1.414$) 

\vskip 0.7ex
\hangindent=3em \hangafter=1
$D^2=$ 8.0 = 
 $8$

\vskip 0.7ex
\hangindent=3em \hangafter=1
$T = ( 0,
\frac{1}{2},
\frac{1}{4},
\frac{3}{4},
\frac{1}{16},
\frac{5}{16} )
$,

\vskip 0.7ex
\hangindent=3em \hangafter=1
$S$ = ($ 1$,
$ 1$,
$ 1$,
$ 1$,
$ \sqrt{2}$,
$ \sqrt{2}$;\ \ 
$ 1$,
$ 1$,
$ 1$,
$ -\sqrt{2}$,
$ -\sqrt{2}$;\ \ 
$ -1$,
$ -1$,
$ \sqrt{2}$,
$ -\sqrt{2}$;\ \ 
$ -1$,
$ -\sqrt{2}$,
$ \sqrt{2}$;\ \ 
$0$,
$0$;\ \ 
$0$)

\vskip 1ex 
\color{grey}

\noindent32. ind = $(6;16
)_{1}^{7}$:\ \ 
$d_i$ = ($1.0$,
$1.0$,
$1.0$,
$1.0$,
$1.414$,
$1.414$) 

\vskip 0.7ex
\hangindent=3em \hangafter=1
$D^2=$ 8.0 = 
 $8$

\vskip 0.7ex
\hangindent=3em \hangafter=1
$T = ( 0,
\frac{1}{2},
\frac{1}{4},
\frac{3}{4},
\frac{3}{16},
\frac{7}{16} )
$,

\vskip 0.7ex
\hangindent=3em \hangafter=1
$S$ = ($ 1$,
$ 1$,
$ 1$,
$ 1$,
$ \sqrt{2}$,
$ \sqrt{2}$;\ \ 
$ 1$,
$ 1$,
$ 1$,
$ -\sqrt{2}$,
$ -\sqrt{2}$;\ \ 
$ -1$,
$ -1$,
$ \sqrt{2}$,
$ -\sqrt{2}$;\ \ 
$ -1$,
$ -\sqrt{2}$,
$ \sqrt{2}$;\ \ 
$0$,
$0$;\ \ 
$0$)

\vskip 1ex 
\color{grey}

\noindent33. ind = $(6;16
)_{1}^{9}$:\ \ 
$d_i$ = ($1.0$,
$1.0$,
$1.0$,
$1.0$,
$1.414$,
$1.414$) 

\vskip 0.7ex
\hangindent=3em \hangafter=1
$D^2=$ 8.0 = 
 $8$

\vskip 0.7ex
\hangindent=3em \hangafter=1
$T = ( 0,
\frac{1}{2},
\frac{1}{4},
\frac{3}{4},
\frac{9}{16},
\frac{13}{16} )
$,

\vskip 0.7ex
\hangindent=3em \hangafter=1
$S$ = ($ 1$,
$ 1$,
$ 1$,
$ 1$,
$ \sqrt{2}$,
$ \sqrt{2}$;\ \ 
$ 1$,
$ 1$,
$ 1$,
$ -\sqrt{2}$,
$ -\sqrt{2}$;\ \ 
$ -1$,
$ -1$,
$ \sqrt{2}$,
$ -\sqrt{2}$;\ \ 
$ -1$,
$ -\sqrt{2}$,
$ \sqrt{2}$;\ \ 
$0$,
$0$;\ \ 
$0$)

\vskip 1ex 
\color{grey}

\noindent34. ind = $(6;16
)_{1}^{15}$:\ \ 
$d_i$ = ($1.0$,
$1.0$,
$1.0$,
$1.0$,
$1.414$,
$1.414$) 

\vskip 0.7ex
\hangindent=3em \hangafter=1
$D^2=$ 8.0 = 
 $8$

\vskip 0.7ex
\hangindent=3em \hangafter=1
$T = ( 0,
\frac{1}{2},
\frac{1}{4},
\frac{3}{4},
\frac{11}{16},
\frac{15}{16} )
$,

\vskip 0.7ex
\hangindent=3em \hangafter=1
$S$ = ($ 1$,
$ 1$,
$ 1$,
$ 1$,
$ \sqrt{2}$,
$ \sqrt{2}$;\ \ 
$ 1$,
$ 1$,
$ 1$,
$ -\sqrt{2}$,
$ -\sqrt{2}$;\ \ 
$ -1$,
$ -1$,
$ \sqrt{2}$,
$ -\sqrt{2}$;\ \ 
$ -1$,
$ -\sqrt{2}$,
$ \sqrt{2}$;\ \ 
$0$,
$0$;\ \ 
$0$)

\vskip 1ex 
\color{grey}

\noindent35. ind = $(6;16
)_{1}^{13}$:\ \ 
$d_i$ = ($1.0$,
$1.0$,
$1.0$,
$1.0$,
$-1.414$,
$-1.414$) 

\vskip 0.7ex
\hangindent=3em \hangafter=1
$D^2=$ 8.0 = 
 $8$

\vskip 0.7ex
\hangindent=3em \hangafter=1
$T = ( 0,
\frac{1}{2},
\frac{1}{4},
\frac{3}{4},
\frac{1}{16},
\frac{13}{16} )
$,

\vskip 0.7ex
\hangindent=3em \hangafter=1
$S$ = ($ 1$,
$ 1$,
$ 1$,
$ 1$,
$ -\sqrt{2}$,
$ -\sqrt{2}$;\ \ 
$ 1$,
$ 1$,
$ 1$,
$ \sqrt{2}$,
$ \sqrt{2}$;\ \ 
$ -1$,
$ -1$,
$ \sqrt{2}$,
$ -\sqrt{2}$;\ \ 
$ -1$,
$ -\sqrt{2}$,
$ \sqrt{2}$;\ \ 
$0$,
$0$;\ \ 
$0$)

Pseudo-unitary $\sim$  
$(6;16
)_{2}^{1}$

\vskip 1ex 
\color{grey}

\noindent36. ind = $(6;16
)_{1}^{3}$:\ \ 
$d_i$ = ($1.0$,
$1.0$,
$1.0$,
$1.0$,
$-1.414$,
$-1.414$) 

\vskip 0.7ex
\hangindent=3em \hangafter=1
$D^2=$ 8.0 = 
 $8$

\vskip 0.7ex
\hangindent=3em \hangafter=1
$T = ( 0,
\frac{1}{2},
\frac{1}{4},
\frac{3}{4},
\frac{3}{16},
\frac{15}{16} )
$,

\vskip 0.7ex
\hangindent=3em \hangafter=1
$S$ = ($ 1$,
$ 1$,
$ 1$,
$ 1$,
$ -\sqrt{2}$,
$ -\sqrt{2}$;\ \ 
$ 1$,
$ 1$,
$ 1$,
$ \sqrt{2}$,
$ \sqrt{2}$;\ \ 
$ -1$,
$ -1$,
$ \sqrt{2}$,
$ -\sqrt{2}$;\ \ 
$ -1$,
$ -\sqrt{2}$,
$ \sqrt{2}$;\ \ 
$0$,
$0$;\ \ 
$0$)

Pseudo-unitary $\sim$  
$(6;16
)_{2}^{15}$

\vskip 1ex 
\color{grey}

\noindent37. ind = $(6;16
)_{1}^{5}$:\ \ 
$d_i$ = ($1.0$,
$1.0$,
$1.0$,
$1.0$,
$-1.414$,
$-1.414$) 

\vskip 0.7ex
\hangindent=3em \hangafter=1
$D^2=$ 8.0 = 
 $8$

\vskip 0.7ex
\hangindent=3em \hangafter=1
$T = ( 0,
\frac{1}{2},
\frac{1}{4},
\frac{3}{4},
\frac{5}{16},
\frac{9}{16} )
$,

\vskip 0.7ex
\hangindent=3em \hangafter=1
$S$ = ($ 1$,
$ 1$,
$ 1$,
$ 1$,
$ -\sqrt{2}$,
$ -\sqrt{2}$;\ \ 
$ 1$,
$ 1$,
$ 1$,
$ \sqrt{2}$,
$ \sqrt{2}$;\ \ 
$ -1$,
$ -1$,
$ -\sqrt{2}$,
$ \sqrt{2}$;\ \ 
$ -1$,
$ \sqrt{2}$,
$ -\sqrt{2}$;\ \ 
$0$,
$0$;\ \ 
$0$)

Pseudo-unitary $\sim$  
$(6;16
)_{2}^{9}$

\vskip 1ex 
\color{grey}

\noindent38. ind = $(6;16
)_{1}^{11}$:\ \ 
$d_i$ = ($1.0$,
$1.0$,
$1.0$,
$1.0$,
$-1.414$,
$-1.414$) 

\vskip 0.7ex
\hangindent=3em \hangafter=1
$D^2=$ 8.0 = 
 $8$

\vskip 0.7ex
\hangindent=3em \hangafter=1
$T = ( 0,
\frac{1}{2},
\frac{1}{4},
\frac{3}{4},
\frac{7}{16},
\frac{11}{16} )
$,

\vskip 0.7ex
\hangindent=3em \hangafter=1
$S$ = ($ 1$,
$ 1$,
$ 1$,
$ 1$,
$ -\sqrt{2}$,
$ -\sqrt{2}$;\ \ 
$ 1$,
$ 1$,
$ 1$,
$ \sqrt{2}$,
$ \sqrt{2}$;\ \ 
$ -1$,
$ -1$,
$ -\sqrt{2}$,
$ \sqrt{2}$;\ \ 
$ -1$,
$ \sqrt{2}$,
$ -\sqrt{2}$;\ \ 
$0$,
$0$;\ \ 
$0$)

Pseudo-unitary $\sim$  
$(6;16
)_{2}^{7}$

\vskip 1ex 

 \color{black} \vskip 2ex

\noindent39. ind = $(6;16
)_{2}^{1}$:\ \ 
$d_i$ = ($1.0$,
$1.0$,
$1.0$,
$1.0$,
$1.414$,
$1.414$) 

\vskip 0.7ex
\hangindent=3em \hangafter=1
$D^2=$ 8.0 = 
 $8$

\vskip 0.7ex
\hangindent=3em \hangafter=1
$T = ( 0,
\frac{1}{2},
\frac{1}{4},
\frac{3}{4},
\frac{1}{16},
\frac{13}{16} )
$,

\vskip 0.7ex
\hangindent=3em \hangafter=1
$S$ = ($ 1$,
$ 1$,
$ 1$,
$ 1$,
$ \sqrt{2}$,
$ \sqrt{2}$;\ \ 
$ 1$,
$ 1$,
$ 1$,
$ -\sqrt{2}$,
$ -\sqrt{2}$;\ \ 
$ -1$,
$ -1$,
$ -\sqrt{2}$,
$ \sqrt{2}$;\ \ 
$ -1$,
$ \sqrt{2}$,
$ -\sqrt{2}$;\ \ 
$0$,
$0$;\ \ 
$0$)

\vskip 1ex 
\color{grey}

\noindent40. ind = $(6;16
)_{2}^{15}$:\ \ 
$d_i$ = ($1.0$,
$1.0$,
$1.0$,
$1.0$,
$1.414$,
$1.414$) 

\vskip 0.7ex
\hangindent=3em \hangafter=1
$D^2=$ 8.0 = 
 $8$

\vskip 0.7ex
\hangindent=3em \hangafter=1
$T = ( 0,
\frac{1}{2},
\frac{1}{4},
\frac{3}{4},
\frac{3}{16},
\frac{15}{16} )
$,

\vskip 0.7ex
\hangindent=3em \hangafter=1
$S$ = ($ 1$,
$ 1$,
$ 1$,
$ 1$,
$ \sqrt{2}$,
$ \sqrt{2}$;\ \ 
$ 1$,
$ 1$,
$ 1$,
$ -\sqrt{2}$,
$ -\sqrt{2}$;\ \ 
$ -1$,
$ -1$,
$ -\sqrt{2}$,
$ \sqrt{2}$;\ \ 
$ -1$,
$ \sqrt{2}$,
$ -\sqrt{2}$;\ \ 
$0$,
$0$;\ \ 
$0$)

\vskip 1ex 
\color{grey}

\noindent41. ind = $(6;16
)_{2}^{9}$:\ \ 
$d_i$ = ($1.0$,
$1.0$,
$1.0$,
$1.0$,
$1.414$,
$1.414$) 

\vskip 0.7ex
\hangindent=3em \hangafter=1
$D^2=$ 8.0 = 
 $8$

\vskip 0.7ex
\hangindent=3em \hangafter=1
$T = ( 0,
\frac{1}{2},
\frac{1}{4},
\frac{3}{4},
\frac{5}{16},
\frac{9}{16} )
$,

\vskip 0.7ex
\hangindent=3em \hangafter=1
$S$ = ($ 1$,
$ 1$,
$ 1$,
$ 1$,
$ \sqrt{2}$,
$ \sqrt{2}$;\ \ 
$ 1$,
$ 1$,
$ 1$,
$ -\sqrt{2}$,
$ -\sqrt{2}$;\ \ 
$ -1$,
$ -1$,
$ \sqrt{2}$,
$ -\sqrt{2}$;\ \ 
$ -1$,
$ -\sqrt{2}$,
$ \sqrt{2}$;\ \ 
$0$,
$0$;\ \ 
$0$)

\vskip 1ex 
\color{grey}

\noindent42. ind = $(6;16
)_{2}^{7}$:\ \ 
$d_i$ = ($1.0$,
$1.0$,
$1.0$,
$1.0$,
$1.414$,
$1.414$) 

\vskip 0.7ex
\hangindent=3em \hangafter=1
$D^2=$ 8.0 = 
 $8$

\vskip 0.7ex
\hangindent=3em \hangafter=1
$T = ( 0,
\frac{1}{2},
\frac{1}{4},
\frac{3}{4},
\frac{7}{16},
\frac{11}{16} )
$,

\vskip 0.7ex
\hangindent=3em \hangafter=1
$S$ = ($ 1$,
$ 1$,
$ 1$,
$ 1$,
$ \sqrt{2}$,
$ \sqrt{2}$;\ \ 
$ 1$,
$ 1$,
$ 1$,
$ -\sqrt{2}$,
$ -\sqrt{2}$;\ \ 
$ -1$,
$ -1$,
$ \sqrt{2}$,
$ -\sqrt{2}$;\ \ 
$ -1$,
$ -\sqrt{2}$,
$ \sqrt{2}$;\ \ 
$0$,
$0$;\ \ 
$0$)

\vskip 1ex 
\color{grey}

\noindent43. ind = $(6;16
)_{2}^{5}$:\ \ 
$d_i$ = ($1.0$,
$1.0$,
$1.0$,
$1.0$,
$-1.414$,
$-1.414$) 

\vskip 0.7ex
\hangindent=3em \hangafter=1
$D^2=$ 8.0 = 
 $8$

\vskip 0.7ex
\hangindent=3em \hangafter=1
$T = ( 0,
\frac{1}{2},
\frac{1}{4},
\frac{3}{4},
\frac{1}{16},
\frac{5}{16} )
$,

\vskip 0.7ex
\hangindent=3em \hangafter=1
$S$ = ($ 1$,
$ 1$,
$ 1$,
$ 1$,
$ -\sqrt{2}$,
$ -\sqrt{2}$;\ \ 
$ 1$,
$ 1$,
$ 1$,
$ \sqrt{2}$,
$ \sqrt{2}$;\ \ 
$ -1$,
$ -1$,
$ -\sqrt{2}$,
$ \sqrt{2}$;\ \ 
$ -1$,
$ \sqrt{2}$,
$ -\sqrt{2}$;\ \ 
$0$,
$0$;\ \ 
$0$)

Pseudo-unitary $\sim$  
$(6;16
)_{1}^{1}$

\vskip 1ex 
\color{grey}

\noindent44. ind = $(6;16
)_{2}^{3}$:\ \ 
$d_i$ = ($1.0$,
$1.0$,
$1.0$,
$1.0$,
$-1.414$,
$-1.414$) 

\vskip 0.7ex
\hangindent=3em \hangafter=1
$D^2=$ 8.0 = 
 $8$

\vskip 0.7ex
\hangindent=3em \hangafter=1
$T = ( 0,
\frac{1}{2},
\frac{1}{4},
\frac{3}{4},
\frac{3}{16},
\frac{7}{16} )
$,

\vskip 0.7ex
\hangindent=3em \hangafter=1
$S$ = ($ 1$,
$ 1$,
$ 1$,
$ 1$,
$ -\sqrt{2}$,
$ -\sqrt{2}$;\ \ 
$ 1$,
$ 1$,
$ 1$,
$ \sqrt{2}$,
$ \sqrt{2}$;\ \ 
$ -1$,
$ -1$,
$ -\sqrt{2}$,
$ \sqrt{2}$;\ \ 
$ -1$,
$ \sqrt{2}$,
$ -\sqrt{2}$;\ \ 
$0$,
$0$;\ \ 
$0$)

Pseudo-unitary $\sim$  
$(6;16
)_{1}^{7}$

\vskip 1ex 
\color{grey}

\noindent45. ind = $(6;16
)_{2}^{13}$:\ \ 
$d_i$ = ($1.0$,
$1.0$,
$1.0$,
$1.0$,
$-1.414$,
$-1.414$) 

\vskip 0.7ex
\hangindent=3em \hangafter=1
$D^2=$ 8.0 = 
 $8$

\vskip 0.7ex
\hangindent=3em \hangafter=1
$T = ( 0,
\frac{1}{2},
\frac{1}{4},
\frac{3}{4},
\frac{9}{16},
\frac{13}{16} )
$,

\vskip 0.7ex
\hangindent=3em \hangafter=1
$S$ = ($ 1$,
$ 1$,
$ 1$,
$ 1$,
$ -\sqrt{2}$,
$ -\sqrt{2}$;\ \ 
$ 1$,
$ 1$,
$ 1$,
$ \sqrt{2}$,
$ \sqrt{2}$;\ \ 
$ -1$,
$ -1$,
$ -\sqrt{2}$,
$ \sqrt{2}$;\ \ 
$ -1$,
$ \sqrt{2}$,
$ -\sqrt{2}$;\ \ 
$0$,
$0$;\ \ 
$0$)

Pseudo-unitary $\sim$  
$(6;16
)_{1}^{9}$

\vskip 1ex 
\color{grey}

\noindent46. ind = $(6;16
)_{2}^{11}$:\ \ 
$d_i$ = ($1.0$,
$1.0$,
$1.0$,
$1.0$,
$-1.414$,
$-1.414$) 

\vskip 0.7ex
\hangindent=3em \hangafter=1
$D^2=$ 8.0 = 
 $8$

\vskip 0.7ex
\hangindent=3em \hangafter=1
$T = ( 0,
\frac{1}{2},
\frac{1}{4},
\frac{3}{4},
\frac{11}{16},
\frac{15}{16} )
$,

\vskip 0.7ex
\hangindent=3em \hangafter=1
$S$ = ($ 1$,
$ 1$,
$ 1$,
$ 1$,
$ -\sqrt{2}$,
$ -\sqrt{2}$;\ \ 
$ 1$,
$ 1$,
$ 1$,
$ \sqrt{2}$,
$ \sqrt{2}$;\ \ 
$ -1$,
$ -1$,
$ -\sqrt{2}$,
$ \sqrt{2}$;\ \ 
$ -1$,
$ \sqrt{2}$,
$ -\sqrt{2}$;\ \ 
$0$,
$0$;\ \ 
$0$)

Pseudo-unitary $\sim$  
$(6;16
)_{1}^{15}$

\vskip 1ex 

 \color{black} \vskip 2ex 
\color{blue}

\noindent47. ind = $(6;16
)_{3}^{1}$:\ \ 
$d_i$ = ($1.0$,
$1.0$,
$1.414$,
$-1.0$,
$-1.0$,
$-1.414$) 

\vskip 0.7ex
\hangindent=3em \hangafter=1
$D^2=$ 8.0 = 
 $8$

\vskip 0.7ex
\hangindent=3em \hangafter=1
$T = ( 0,
\frac{1}{2},
\frac{1}{16},
\frac{1}{4},
\frac{3}{4},
\frac{5}{16} )
$,

\vskip 0.7ex
\hangindent=3em \hangafter=1
$S$ = ($ 1$,
$ 1$,
$ \sqrt{2}$,
$ -1$,
$ -1$,
$ -\sqrt{2}$;\ \ 
$ 1$,
$ -\sqrt{2}$,
$ -1$,
$ -1$,
$ \sqrt{2}$;\ \ 
$0$,
$ -\sqrt{2}$,
$ \sqrt{2}$,
$0$;\ \ 
$ -1$,
$ -1$,
$ -\sqrt{2}$;\ \ 
$ -1$,
$ \sqrt{2}$;\ \ 
$0$)

Pseudo-unitary $\sim$  
$(6;16
)_{1}^{1}$

\vskip 1ex 
\color{grey}

\noindent48. ind = $(6;16
)_{3}^{13}$:\ \ 
$d_i$ = ($1.0$,
$1.0$,
$1.414$,
$-1.0$,
$-1.0$,
$-1.414$) 

\vskip 0.7ex
\hangindent=3em \hangafter=1
$D^2=$ 8.0 = 
 $8$

\vskip 0.7ex
\hangindent=3em \hangafter=1
$T = ( 0,
\frac{1}{2},
\frac{1}{16},
\frac{1}{4},
\frac{3}{4},
\frac{13}{16} )
$,

\vskip 0.7ex
\hangindent=3em \hangafter=1
$S$ = ($ 1$,
$ 1$,
$ \sqrt{2}$,
$ -1$,
$ -1$,
$ -\sqrt{2}$;\ \ 
$ 1$,
$ -\sqrt{2}$,
$ -1$,
$ -1$,
$ \sqrt{2}$;\ \ 
$0$,
$ \sqrt{2}$,
$ -\sqrt{2}$,
$0$;\ \ 
$ -1$,
$ -1$,
$ \sqrt{2}$;\ \ 
$ -1$,
$ -\sqrt{2}$;\ \ 
$0$)

Pseudo-unitary $\sim$  
$(6;16
)_{2}^{1}$

\vskip 1ex 
\color{grey}

\noindent49. ind = $(6;16
)_{3}^{7}$:\ \ 
$d_i$ = ($1.0$,
$1.0$,
$1.414$,
$-1.0$,
$-1.0$,
$-1.414$) 

\vskip 0.7ex
\hangindent=3em \hangafter=1
$D^2=$ 8.0 = 
 $8$

\vskip 0.7ex
\hangindent=3em \hangafter=1
$T = ( 0,
\frac{1}{2},
\frac{7}{16},
\frac{1}{4},
\frac{3}{4},
\frac{3}{16} )
$,

\vskip 0.7ex
\hangindent=3em \hangafter=1
$S$ = ($ 1$,
$ 1$,
$ \sqrt{2}$,
$ -1$,
$ -1$,
$ -\sqrt{2}$;\ \ 
$ 1$,
$ -\sqrt{2}$,
$ -1$,
$ -1$,
$ \sqrt{2}$;\ \ 
$0$,
$ \sqrt{2}$,
$ -\sqrt{2}$,
$0$;\ \ 
$ -1$,
$ -1$,
$ \sqrt{2}$;\ \ 
$ -1$,
$ -\sqrt{2}$;\ \ 
$0$)

Pseudo-unitary $\sim$  
$(6;16
)_{1}^{7}$

\vskip 1ex 
\color{grey}

\noindent50. ind = $(6;16
)_{3}^{11}$:\ \ 
$d_i$ = ($1.0$,
$1.0$,
$1.414$,
$-1.0$,
$-1.0$,
$-1.414$) 

\vskip 0.7ex
\hangindent=3em \hangafter=1
$D^2=$ 8.0 = 
 $8$

\vskip 0.7ex
\hangindent=3em \hangafter=1
$T = ( 0,
\frac{1}{2},
\frac{7}{16},
\frac{1}{4},
\frac{3}{4},
\frac{11}{16} )
$,

\vskip 0.7ex
\hangindent=3em \hangafter=1
$S$ = ($ 1$,
$ 1$,
$ \sqrt{2}$,
$ -1$,
$ -1$,
$ -\sqrt{2}$;\ \ 
$ 1$,
$ -\sqrt{2}$,
$ -1$,
$ -1$,
$ \sqrt{2}$;\ \ 
$0$,
$ -\sqrt{2}$,
$ \sqrt{2}$,
$0$;\ \ 
$ -1$,
$ -1$,
$ -\sqrt{2}$;\ \ 
$ -1$,
$ \sqrt{2}$;\ \ 
$0$)

Pseudo-unitary $\sim$  
$(6;16
)_{2}^{7}$

\vskip 1ex 
\color{grey}

\noindent51. ind = $(6;16
)_{3}^{5}$:\ \ 
$d_i$ = ($1.0$,
$1.0$,
$1.414$,
$-1.0$,
$-1.0$,
$-1.414$) 

\vskip 0.7ex
\hangindent=3em \hangafter=1
$D^2=$ 8.0 = 
 $8$

\vskip 0.7ex
\hangindent=3em \hangafter=1
$T = ( 0,
\frac{1}{2},
\frac{9}{16},
\frac{1}{4},
\frac{3}{4},
\frac{5}{16} )
$,

\vskip 0.7ex
\hangindent=3em \hangafter=1
$S$ = ($ 1$,
$ 1$,
$ \sqrt{2}$,
$ -1$,
$ -1$,
$ -\sqrt{2}$;\ \ 
$ 1$,
$ -\sqrt{2}$,
$ -1$,
$ -1$,
$ \sqrt{2}$;\ \ 
$0$,
$ \sqrt{2}$,
$ -\sqrt{2}$,
$0$;\ \ 
$ -1$,
$ -1$,
$ \sqrt{2}$;\ \ 
$ -1$,
$ -\sqrt{2}$;\ \ 
$0$)

Pseudo-unitary $\sim$  
$(6;16
)_{2}^{9}$

\vskip 1ex 
\color{grey}

\noindent52. ind = $(6;16
)_{3}^{9}$:\ \ 
$d_i$ = ($1.0$,
$1.0$,
$1.414$,
$-1.0$,
$-1.0$,
$-1.414$) 

\vskip 0.7ex
\hangindent=3em \hangafter=1
$D^2=$ 8.0 = 
 $8$

\vskip 0.7ex
\hangindent=3em \hangafter=1
$T = ( 0,
\frac{1}{2},
\frac{9}{16},
\frac{1}{4},
\frac{3}{4},
\frac{13}{16} )
$,

\vskip 0.7ex
\hangindent=3em \hangafter=1
$S$ = ($ 1$,
$ 1$,
$ \sqrt{2}$,
$ -1$,
$ -1$,
$ -\sqrt{2}$;\ \ 
$ 1$,
$ -\sqrt{2}$,
$ -1$,
$ -1$,
$ \sqrt{2}$;\ \ 
$0$,
$ -\sqrt{2}$,
$ \sqrt{2}$,
$0$;\ \ 
$ -1$,
$ -1$,
$ -\sqrt{2}$;\ \ 
$ -1$,
$ \sqrt{2}$;\ \ 
$0$)

Pseudo-unitary $\sim$  
$(6;16
)_{1}^{9}$

\vskip 1ex 
\color{grey}

\noindent53. ind = $(6;16
)_{3}^{3}$:\ \ 
$d_i$ = ($1.0$,
$1.0$,
$1.414$,
$-1.0$,
$-1.0$,
$-1.414$) 

\vskip 0.7ex
\hangindent=3em \hangafter=1
$D^2=$ 8.0 = 
 $8$

\vskip 0.7ex
\hangindent=3em \hangafter=1
$T = ( 0,
\frac{1}{2},
\frac{15}{16},
\frac{1}{4},
\frac{3}{4},
\frac{3}{16} )
$,

\vskip 0.7ex
\hangindent=3em \hangafter=1
$S$ = ($ 1$,
$ 1$,
$ \sqrt{2}$,
$ -1$,
$ -1$,
$ -\sqrt{2}$;\ \ 
$ 1$,
$ -\sqrt{2}$,
$ -1$,
$ -1$,
$ \sqrt{2}$;\ \ 
$0$,
$ -\sqrt{2}$,
$ \sqrt{2}$,
$0$;\ \ 
$ -1$,
$ -1$,
$ -\sqrt{2}$;\ \ 
$ -1$,
$ \sqrt{2}$;\ \ 
$0$)

Pseudo-unitary $\sim$  
$(6;16
)_{2}^{15}$

\vskip 1ex 
\color{grey}

\noindent54. ind = $(6;16
)_{3}^{15}$:\ \ 
$d_i$ = ($1.0$,
$1.0$,
$1.414$,
$-1.0$,
$-1.0$,
$-1.414$) 

\vskip 0.7ex
\hangindent=3em \hangafter=1
$D^2=$ 8.0 = 
 $8$

\vskip 0.7ex
\hangindent=3em \hangafter=1
$T = ( 0,
\frac{1}{2},
\frac{15}{16},
\frac{1}{4},
\frac{3}{4},
\frac{11}{16} )
$,

\vskip 0.7ex
\hangindent=3em \hangafter=1
$S$ = ($ 1$,
$ 1$,
$ \sqrt{2}$,
$ -1$,
$ -1$,
$ -\sqrt{2}$;\ \ 
$ 1$,
$ -\sqrt{2}$,
$ -1$,
$ -1$,
$ \sqrt{2}$;\ \ 
$0$,
$ \sqrt{2}$,
$ -\sqrt{2}$,
$0$;\ \ 
$ -1$,
$ -1$,
$ \sqrt{2}$;\ \ 
$ -1$,
$ -\sqrt{2}$;\ \ 
$0$)

Pseudo-unitary $\sim$  
$(6;16
)_{1}^{15}$

\vskip 1ex 

 \color{black} \vskip 2ex 
\color{blue}

\noindent55. ind = $(6;16
)_{4}^{1}$:\ \ 
$d_i$ = ($1.0$,
$1.0$,
$1.414$,
$-1.0$,
$-1.0$,
$-1.414$) 

\vskip 0.7ex
\hangindent=3em \hangafter=1
$D^2=$ 8.0 = 
 $8$

\vskip 0.7ex
\hangindent=3em \hangafter=1
$T = ( 0,
\frac{1}{2},
\frac{3}{16},
\frac{1}{4},
\frac{3}{4},
\frac{7}{16} )
$,

\vskip 0.7ex
\hangindent=3em \hangafter=1
$S$ = ($ 1$,
$ 1$,
$ \sqrt{2}$,
$ -1$,
$ -1$,
$ -\sqrt{2}$;\ \ 
$ 1$,
$ -\sqrt{2}$,
$ -1$,
$ -1$,
$ \sqrt{2}$;\ \ 
$0$,
$ -\sqrt{2}$,
$ \sqrt{2}$,
$0$;\ \ 
$ -1$,
$ -1$,
$ -\sqrt{2}$;\ \ 
$ -1$,
$ \sqrt{2}$;\ \ 
$0$)

Pseudo-unitary $\sim$  
$(6;16
)_{1}^{7}$

\vskip 1ex 
\color{grey}

\noindent56. ind = $(6;16
)_{4}^{5}$:\ \ 
$d_i$ = ($1.0$,
$1.0$,
$1.414$,
$-1.0$,
$-1.0$,
$-1.414$) 

\vskip 0.7ex
\hangindent=3em \hangafter=1
$D^2=$ 8.0 = 
 $8$

\vskip 0.7ex
\hangindent=3em \hangafter=1
$T = ( 0,
\frac{1}{2},
\frac{3}{16},
\frac{1}{4},
\frac{3}{4},
\frac{15}{16} )
$,

\vskip 0.7ex
\hangindent=3em \hangafter=1
$S$ = ($ 1$,
$ 1$,
$ \sqrt{2}$,
$ -1$,
$ -1$,
$ -\sqrt{2}$;\ \ 
$ 1$,
$ -\sqrt{2}$,
$ -1$,
$ -1$,
$ \sqrt{2}$;\ \ 
$0$,
$ \sqrt{2}$,
$ -\sqrt{2}$,
$0$;\ \ 
$ -1$,
$ -1$,
$ \sqrt{2}$;\ \ 
$ -1$,
$ -\sqrt{2}$;\ \ 
$0$)

Pseudo-unitary $\sim$  
$(6;16
)_{2}^{15}$

\vskip 1ex 
\color{grey}

\noindent57. ind = $(6;16
)_{4}^{7}$:\ \ 
$d_i$ = ($1.0$,
$1.0$,
$1.414$,
$-1.0$,
$-1.0$,
$-1.414$) 

\vskip 0.7ex
\hangindent=3em \hangafter=1
$D^2=$ 8.0 = 
 $8$

\vskip 0.7ex
\hangindent=3em \hangafter=1
$T = ( 0,
\frac{1}{2},
\frac{5}{16},
\frac{1}{4},
\frac{3}{4},
\frac{1}{16} )
$,

\vskip 0.7ex
\hangindent=3em \hangafter=1
$S$ = ($ 1$,
$ 1$,
$ \sqrt{2}$,
$ -1$,
$ -1$,
$ -\sqrt{2}$;\ \ 
$ 1$,
$ -\sqrt{2}$,
$ -1$,
$ -1$,
$ \sqrt{2}$;\ \ 
$0$,
$ \sqrt{2}$,
$ -\sqrt{2}$,
$0$;\ \ 
$ -1$,
$ -1$,
$ \sqrt{2}$;\ \ 
$ -1$,
$ -\sqrt{2}$;\ \ 
$0$)

Pseudo-unitary $\sim$  
$(6;16
)_{1}^{1}$

\vskip 1ex 
\color{grey}

\noindent58. ind = $(6;16
)_{4}^{3}$:\ \ 
$d_i$ = ($1.0$,
$1.0$,
$1.414$,
$-1.0$,
$-1.0$,
$-1.414$) 

\vskip 0.7ex
\hangindent=3em \hangafter=1
$D^2=$ 8.0 = 
 $8$

\vskip 0.7ex
\hangindent=3em \hangafter=1
$T = ( 0,
\frac{1}{2},
\frac{5}{16},
\frac{1}{4},
\frac{3}{4},
\frac{9}{16} )
$,

\vskip 0.7ex
\hangindent=3em \hangafter=1
$S$ = ($ 1$,
$ 1$,
$ \sqrt{2}$,
$ -1$,
$ -1$,
$ -\sqrt{2}$;\ \ 
$ 1$,
$ -\sqrt{2}$,
$ -1$,
$ -1$,
$ \sqrt{2}$;\ \ 
$0$,
$ -\sqrt{2}$,
$ \sqrt{2}$,
$0$;\ \ 
$ -1$,
$ -1$,
$ -\sqrt{2}$;\ \ 
$ -1$,
$ \sqrt{2}$;\ \ 
$0$)

Pseudo-unitary $\sim$  
$(6;16
)_{2}^{9}$

\vskip 1ex 
\color{grey}

\noindent59. ind = $(6;16
)_{4}^{13}$:\ \ 
$d_i$ = ($1.0$,
$1.0$,
$1.414$,
$-1.0$,
$-1.0$,
$-1.414$) 

\vskip 0.7ex
\hangindent=3em \hangafter=1
$D^2=$ 8.0 = 
 $8$

\vskip 0.7ex
\hangindent=3em \hangafter=1
$T = ( 0,
\frac{1}{2},
\frac{11}{16},
\frac{1}{4},
\frac{3}{4},
\frac{7}{16} )
$,

\vskip 0.7ex
\hangindent=3em \hangafter=1
$S$ = ($ 1$,
$ 1$,
$ \sqrt{2}$,
$ -1$,
$ -1$,
$ -\sqrt{2}$;\ \ 
$ 1$,
$ -\sqrt{2}$,
$ -1$,
$ -1$,
$ \sqrt{2}$;\ \ 
$0$,
$ \sqrt{2}$,
$ -\sqrt{2}$,
$0$;\ \ 
$ -1$,
$ -1$,
$ \sqrt{2}$;\ \ 
$ -1$,
$ -\sqrt{2}$;\ \ 
$0$)

Pseudo-unitary $\sim$  
$(6;16
)_{2}^{7}$

\vskip 1ex 
\color{grey}

\noindent60. ind = $(6;16
)_{4}^{9}$:\ \ 
$d_i$ = ($1.0$,
$1.0$,
$1.414$,
$-1.0$,
$-1.0$,
$-1.414$) 

\vskip 0.7ex
\hangindent=3em \hangafter=1
$D^2=$ 8.0 = 
 $8$

\vskip 0.7ex
\hangindent=3em \hangafter=1
$T = ( 0,
\frac{1}{2},
\frac{11}{16},
\frac{1}{4},
\frac{3}{4},
\frac{15}{16} )
$,

\vskip 0.7ex
\hangindent=3em \hangafter=1
$S$ = ($ 1$,
$ 1$,
$ \sqrt{2}$,
$ -1$,
$ -1$,
$ -\sqrt{2}$;\ \ 
$ 1$,
$ -\sqrt{2}$,
$ -1$,
$ -1$,
$ \sqrt{2}$;\ \ 
$0$,
$ -\sqrt{2}$,
$ \sqrt{2}$,
$0$;\ \ 
$ -1$,
$ -1$,
$ -\sqrt{2}$;\ \ 
$ -1$,
$ \sqrt{2}$;\ \ 
$0$)

Pseudo-unitary $\sim$  
$(6;16
)_{1}^{15}$

\vskip 1ex 
\color{grey}

\noindent61. ind = $(6;16
)_{4}^{11}$:\ \ 
$d_i$ = ($1.0$,
$1.0$,
$1.414$,
$-1.0$,
$-1.0$,
$-1.414$) 

\vskip 0.7ex
\hangindent=3em \hangafter=1
$D^2=$ 8.0 = 
 $8$

\vskip 0.7ex
\hangindent=3em \hangafter=1
$T = ( 0,
\frac{1}{2},
\frac{13}{16},
\frac{1}{4},
\frac{3}{4},
\frac{1}{16} )
$,

\vskip 0.7ex
\hangindent=3em \hangafter=1
$S$ = ($ 1$,
$ 1$,
$ \sqrt{2}$,
$ -1$,
$ -1$,
$ -\sqrt{2}$;\ \ 
$ 1$,
$ -\sqrt{2}$,
$ -1$,
$ -1$,
$ \sqrt{2}$;\ \ 
$0$,
$ -\sqrt{2}$,
$ \sqrt{2}$,
$0$;\ \ 
$ -1$,
$ -1$,
$ -\sqrt{2}$;\ \ 
$ -1$,
$ \sqrt{2}$;\ \ 
$0$)

Pseudo-unitary $\sim$  
$(6;16
)_{2}^{1}$

\vskip 1ex 
\color{grey}

\noindent62. ind = $(6;16
)_{4}^{15}$:\ \ 
$d_i$ = ($1.0$,
$1.0$,
$1.414$,
$-1.0$,
$-1.0$,
$-1.414$) 

\vskip 0.7ex
\hangindent=3em \hangafter=1
$D^2=$ 8.0 = 
 $8$

\vskip 0.7ex
\hangindent=3em \hangafter=1
$T = ( 0,
\frac{1}{2},
\frac{13}{16},
\frac{1}{4},
\frac{3}{4},
\frac{9}{16} )
$,

\vskip 0.7ex
\hangindent=3em \hangafter=1
$S$ = ($ 1$,
$ 1$,
$ \sqrt{2}$,
$ -1$,
$ -1$,
$ -\sqrt{2}$;\ \ 
$ 1$,
$ -\sqrt{2}$,
$ -1$,
$ -1$,
$ \sqrt{2}$;\ \ 
$0$,
$ \sqrt{2}$,
$ -\sqrt{2}$,
$0$;\ \ 
$ -1$,
$ -1$,
$ \sqrt{2}$;\ \ 
$ -1$,
$ -\sqrt{2}$;\ \ 
$0$)

Pseudo-unitary $\sim$  
$(6;16
)_{1}^{9}$

\vskip 1ex 

 \color{black} \vskip 2ex

\noindent63. ind = $(6;35
)_{1}^{1}$:\ \ 
$d_i$ = ($1.0$,
$1.618$,
$1.801$,
$2.246$,
$2.915$,
$3.635$) 

\vskip 0.7ex
\hangindent=3em \hangafter=1
$D^2=$ 33.632 = 
 $15+3c^{1}_{35}
+2c^{4}_{35}
+6c^{5}_{35}
+3c^{6}_{35}
+3c^{7}_{35}
+2c^{10}_{35}
+2c^{11}_{35}
$

\vskip 0.7ex
\hangindent=3em \hangafter=1
$T = ( 0,
\frac{2}{5},
\frac{1}{7},
\frac{5}{7},
\frac{19}{35},
\frac{4}{35} )
$,

\vskip 0.7ex
\hangindent=3em \hangafter=1
$S$ = ($ 1$,
$ \frac{1+\sqrt{5}}{2}$,
$ \xi_{7}^{2}$,
$ \xi_{7}^{3}$,
$ c^{1}_{35}
+c^{6}_{35}
$,
$ c^{1}_{35}
+c^{4}_{35}
+c^{6}_{35}
+c^{11}_{35}
$;\ \ 
$ -1$,
$ c^{1}_{35}
+c^{6}_{35}
$,
$ c^{1}_{35}
+c^{4}_{35}
+c^{6}_{35}
+c^{11}_{35}
$,
$ -\xi_{7}^{2}$,
$ -\xi_{7}^{3}$;\ \ 
$ -\xi_{7}^{3}$,
$ 1$,
$ -c^{1}_{35}
-c^{4}_{35}
-c^{6}_{35}
-c^{11}_{35}
$,
$ \frac{1+\sqrt{5}}{2}$;\ \ 
$ -\xi_{7}^{2}$,
$ \frac{1+\sqrt{5}}{2}$,
$ -c^{1}_{35}
-c^{6}_{35}
$;\ \ 
$ \xi_{7}^{3}$,
$ -1$;\ \ 
$ \xi_{7}^{2}$)

\vskip 1ex 
\color{grey}

\noindent64. ind = $(6;35
)_{1}^{6}$:\ \ 
$d_i$ = ($1.0$,
$1.618$,
$1.801$,
$2.246$,
$2.915$,
$3.635$) 

\vskip 0.7ex
\hangindent=3em \hangafter=1
$D^2=$ 33.632 = 
 $15+3c^{1}_{35}
+2c^{4}_{35}
+6c^{5}_{35}
+3c^{6}_{35}
+3c^{7}_{35}
+2c^{10}_{35}
+2c^{11}_{35}
$

\vskip 0.7ex
\hangindent=3em \hangafter=1
$T = ( 0,
\frac{2}{5},
\frac{6}{7},
\frac{2}{7},
\frac{9}{35},
\frac{24}{35} )
$,

\vskip 0.7ex
\hangindent=3em \hangafter=1
$S$ = ($ 1$,
$ \frac{1+\sqrt{5}}{2}$,
$ \xi_{7}^{2}$,
$ \xi_{7}^{3}$,
$ c^{1}_{35}
+c^{6}_{35}
$,
$ c^{1}_{35}
+c^{4}_{35}
+c^{6}_{35}
+c^{11}_{35}
$;\ \ 
$ -1$,
$ c^{1}_{35}
+c^{6}_{35}
$,
$ c^{1}_{35}
+c^{4}_{35}
+c^{6}_{35}
+c^{11}_{35}
$,
$ -\xi_{7}^{2}$,
$ -\xi_{7}^{3}$;\ \ 
$ -\xi_{7}^{3}$,
$ 1$,
$ -c^{1}_{35}
-c^{4}_{35}
-c^{6}_{35}
-c^{11}_{35}
$,
$ \frac{1+\sqrt{5}}{2}$;\ \ 
$ -\xi_{7}^{2}$,
$ \frac{1+\sqrt{5}}{2}$,
$ -c^{1}_{35}
-c^{6}_{35}
$;\ \ 
$ \xi_{7}^{3}$,
$ -1$;\ \ 
$ \xi_{7}^{2}$)

\vskip 1ex 
\color{grey}

\noindent65. ind = $(6;35
)_{1}^{29}$:\ \ 
$d_i$ = ($1.0$,
$1.618$,
$1.801$,
$2.246$,
$2.915$,
$3.635$) 

\vskip 0.7ex
\hangindent=3em \hangafter=1
$D^2=$ 33.632 = 
 $15+3c^{1}_{35}
+2c^{4}_{35}
+6c^{5}_{35}
+3c^{6}_{35}
+3c^{7}_{35}
+2c^{10}_{35}
+2c^{11}_{35}
$

\vskip 0.7ex
\hangindent=3em \hangafter=1
$T = ( 0,
\frac{3}{5},
\frac{1}{7},
\frac{5}{7},
\frac{26}{35},
\frac{11}{35} )
$,

\vskip 0.7ex
\hangindent=3em \hangafter=1
$S$ = ($ 1$,
$ \frac{1+\sqrt{5}}{2}$,
$ \xi_{7}^{2}$,
$ \xi_{7}^{3}$,
$ c^{1}_{35}
+c^{6}_{35}
$,
$ c^{1}_{35}
+c^{4}_{35}
+c^{6}_{35}
+c^{11}_{35}
$;\ \ 
$ -1$,
$ c^{1}_{35}
+c^{6}_{35}
$,
$ c^{1}_{35}
+c^{4}_{35}
+c^{6}_{35}
+c^{11}_{35}
$,
$ -\xi_{7}^{2}$,
$ -\xi_{7}^{3}$;\ \ 
$ -\xi_{7}^{3}$,
$ 1$,
$ -c^{1}_{35}
-c^{4}_{35}
-c^{6}_{35}
-c^{11}_{35}
$,
$ \frac{1+\sqrt{5}}{2}$;\ \ 
$ -\xi_{7}^{2}$,
$ \frac{1+\sqrt{5}}{2}$,
$ -c^{1}_{35}
-c^{6}_{35}
$;\ \ 
$ \xi_{7}^{3}$,
$ -1$;\ \ 
$ \xi_{7}^{2}$)

\vskip 1ex 
\color{grey}

\noindent66. ind = $(6;35
)_{1}^{34}$:\ \ 
$d_i$ = ($1.0$,
$1.618$,
$1.801$,
$2.246$,
$2.915$,
$3.635$) 

\vskip 0.7ex
\hangindent=3em \hangafter=1
$D^2=$ 33.632 = 
 $15+3c^{1}_{35}
+2c^{4}_{35}
+6c^{5}_{35}
+3c^{6}_{35}
+3c^{7}_{35}
+2c^{10}_{35}
+2c^{11}_{35}
$

\vskip 0.7ex
\hangindent=3em \hangafter=1
$T = ( 0,
\frac{3}{5},
\frac{6}{7},
\frac{2}{7},
\frac{16}{35},
\frac{31}{35} )
$,

\vskip 0.7ex
\hangindent=3em \hangafter=1
$S$ = ($ 1$,
$ \frac{1+\sqrt{5}}{2}$,
$ \xi_{7}^{2}$,
$ \xi_{7}^{3}$,
$ c^{1}_{35}
+c^{6}_{35}
$,
$ c^{1}_{35}
+c^{4}_{35}
+c^{6}_{35}
+c^{11}_{35}
$;\ \ 
$ -1$,
$ c^{1}_{35}
+c^{6}_{35}
$,
$ c^{1}_{35}
+c^{4}_{35}
+c^{6}_{35}
+c^{11}_{35}
$,
$ -\xi_{7}^{2}$,
$ -\xi_{7}^{3}$;\ \ 
$ -\xi_{7}^{3}$,
$ 1$,
$ -c^{1}_{35}
-c^{4}_{35}
-c^{6}_{35}
-c^{11}_{35}
$,
$ \frac{1+\sqrt{5}}{2}$;\ \ 
$ -\xi_{7}^{2}$,
$ \frac{1+\sqrt{5}}{2}$,
$ -c^{1}_{35}
-c^{6}_{35}
$;\ \ 
$ \xi_{7}^{3}$,
$ -1$;\ \ 
$ \xi_{7}^{2}$)

\vskip 1ex 
\color{grey}

\noindent67. ind = $(6;35
)_{1}^{8}$:\ \ 
$d_i$ = ($1.0$,
$1.801$,
$2.246$,
$-0.618$,
$-1.113$,
$-1.388$) 

\vskip 0.7ex
\hangindent=3em \hangafter=1
$D^2=$ 12.846 = 
 $15-3c^{1}_{35}
-2c^{4}_{35}
+9c^{5}_{35}
-3c^{6}_{35}
-3c^{7}_{35}
+3c^{10}_{35}
-2c^{11}_{35}
$

\vskip 0.7ex
\hangindent=3em \hangafter=1
$T = ( 0,
\frac{1}{7},
\frac{5}{7},
\frac{1}{5},
\frac{12}{35},
\frac{32}{35} )
$,

\vskip 0.7ex
\hangindent=3em \hangafter=1
$S$ = ($ 1$,
$ \xi_{7}^{2}$,
$ \xi_{7}^{3}$,
$ \frac{1-\sqrt{5}}{2}$,
$ 1-c^{1}_{35}
+c^{5}_{35}
-c^{6}_{35}
+c^{10}_{35}
$,
$ 1-c^{1}_{35}
-c^{4}_{35}
+c^{5}_{35}
-c^{6}_{35}
-c^{11}_{35}
$;\ \ 
$ -\xi_{7}^{3}$,
$ 1$,
$ 1-c^{1}_{35}
+c^{5}_{35}
-c^{6}_{35}
+c^{10}_{35}
$,
$ -1+c^{1}_{35}
+c^{4}_{35}
-c^{5}_{35}
+c^{6}_{35}
+c^{11}_{35}
$,
$ \frac{1-\sqrt{5}}{2}$;\ \ 
$ -\xi_{7}^{2}$,
$ 1-c^{1}_{35}
-c^{4}_{35}
+c^{5}_{35}
-c^{6}_{35}
-c^{11}_{35}
$,
$ \frac{1-\sqrt{5}}{2}$,
$ -1+c^{1}_{35}
-c^{5}_{35}
+c^{6}_{35}
-c^{10}_{35}
$;\ \ 
$ -1$,
$ -\xi_{7}^{2}$,
$ -\xi_{7}^{3}$;\ \ 
$ \xi_{7}^{3}$,
$ -1$;\ \ 
$ \xi_{7}^{2}$)

Not pseudo-unitary. 

\vskip 1ex 
\color{grey}

\noindent68. ind = $(6;35
)_{1}^{22}$:\ \ 
$d_i$ = ($1.0$,
$1.801$,
$2.246$,
$-0.618$,
$-1.113$,
$-1.388$) 

\vskip 0.7ex
\hangindent=3em \hangafter=1
$D^2=$ 12.846 = 
 $15-3c^{1}_{35}
-2c^{4}_{35}
+9c^{5}_{35}
-3c^{6}_{35}
-3c^{7}_{35}
+3c^{10}_{35}
-2c^{11}_{35}
$

\vskip 0.7ex
\hangindent=3em \hangafter=1
$T = ( 0,
\frac{1}{7},
\frac{5}{7},
\frac{4}{5},
\frac{33}{35},
\frac{18}{35} )
$,

\vskip 0.7ex
\hangindent=3em \hangafter=1
$S$ = ($ 1$,
$ \xi_{7}^{2}$,
$ \xi_{7}^{3}$,
$ \frac{1-\sqrt{5}}{2}$,
$ 1-c^{1}_{35}
+c^{5}_{35}
-c^{6}_{35}
+c^{10}_{35}
$,
$ 1-c^{1}_{35}
-c^{4}_{35}
+c^{5}_{35}
-c^{6}_{35}
-c^{11}_{35}
$;\ \ 
$ -\xi_{7}^{3}$,
$ 1$,
$ 1-c^{1}_{35}
+c^{5}_{35}
-c^{6}_{35}
+c^{10}_{35}
$,
$ -1+c^{1}_{35}
+c^{4}_{35}
-c^{5}_{35}
+c^{6}_{35}
+c^{11}_{35}
$,
$ \frac{1-\sqrt{5}}{2}$;\ \ 
$ -\xi_{7}^{2}$,
$ 1-c^{1}_{35}
-c^{4}_{35}
+c^{5}_{35}
-c^{6}_{35}
-c^{11}_{35}
$,
$ \frac{1-\sqrt{5}}{2}$,
$ -1+c^{1}_{35}
-c^{5}_{35}
+c^{6}_{35}
-c^{10}_{35}
$;\ \ 
$ -1$,
$ -\xi_{7}^{2}$,
$ -\xi_{7}^{3}$;\ \ 
$ \xi_{7}^{3}$,
$ -1$;\ \ 
$ \xi_{7}^{2}$)

Not pseudo-unitary. 

\vskip 1ex 
\color{grey}

\noindent69. ind = $(6;35
)_{1}^{24}$:\ \ 
$d_i$ = ($1.0$,
$0.445$,
$0.720$,
$1.618$,
$-0.801$,
$-1.297$) 

\vskip 0.7ex
\hangindent=3em \hangafter=1
$D^2=$ 6.661 = 
 $11-2c^{1}_{35}
+c^{4}_{35}
-4c^{5}_{35}
-2c^{6}_{35}
+5c^{7}_{35}
-6c^{10}_{35}
+c^{11}_{35}
$

\vskip 0.7ex
\hangindent=3em \hangafter=1
$T = ( 0,
\frac{3}{7},
\frac{1}{35},
\frac{3}{5},
\frac{1}{7},
\frac{26}{35} )
$,

\vskip 0.7ex
\hangindent=3em \hangafter=1
$S$ = ($ 1$,
$ -c^{2}_{7}
$,
$ c^{4}_{35}
+c^{11}_{35}
$,
$ \frac{1+\sqrt{5}}{2}$,
$ -\xi_{7}^{2,3}$,
$ 1-c^{1}_{35}
-c^{6}_{35}
+c^{7}_{35}
$;\ \ 
$ \xi_{7}^{2,3}$,
$ -1+c^{1}_{35}
+c^{6}_{35}
-c^{7}_{35}
$,
$ c^{4}_{35}
+c^{11}_{35}
$,
$ 1$,
$ \frac{1+\sqrt{5}}{2}$;\ \ 
$ -\xi_{7}^{2,3}$,
$ c^{2}_{7}
$,
$ \frac{1+\sqrt{5}}{2}$,
$ -1$;\ \ 
$ -1$,
$ 1-c^{1}_{35}
-c^{6}_{35}
+c^{7}_{35}
$,
$ \xi_{7}^{2,3}$;\ \ 
$ c^{2}_{7}
$,
$ -c^{4}_{35}
-c^{11}_{35}
$;\ \ 
$ -c^{2}_{7}
$)

Not pseudo-unitary. 

\vskip 1ex 
\color{grey}

\noindent70. ind = $(6;35
)_{1}^{9}$:\ \ 
$d_i$ = ($1.0$,
$0.554$,
$0.897$,
$1.618$,
$-1.246$,
$-2.17$) 

\vskip 0.7ex
\hangindent=3em \hangafter=1
$D^2=$ 10.358 = 
 $\frac{35}{\sqrt{35}\mathrm{i}}s^{1}_{35}
+\frac{35}{\sqrt{35}\mathrm{i}}s^{4}_{35}
-\frac{35}{\sqrt{35}\mathrm{i}}s^{6}_{35}
+\frac{35}{\sqrt{35}\mathrm{i}}s^{11}_{35}
$

\vskip 0.7ex
\hangindent=3em \hangafter=1
$T = ( 0,
\frac{3}{7},
\frac{1}{35},
\frac{3}{5},
\frac{2}{7},
\frac{31}{35} )
$,

\vskip 0.7ex
\hangindent=3em \hangafter=1
$S$ = ($ 1$,
$ \xi_{7}^{1,2}$,
$ 1-c^{4}_{35}
+c^{7}_{35}
-c^{11}_{35}
$,
$ \frac{1+\sqrt{5}}{2}$,
$ -c^{1}_{7}
$,
$ 1-c^{1}_{35}
-c^{4}_{35}
-c^{6}_{35}
+c^{7}_{35}
-c^{11}_{35}
$;\ \ 
$ c^{1}_{7}
$,
$ -1+c^{1}_{35}
+c^{4}_{35}
+c^{6}_{35}
-c^{7}_{35}
+c^{11}_{35}
$,
$ 1-c^{4}_{35}
+c^{7}_{35}
-c^{11}_{35}
$,
$ 1$,
$ \frac{1+\sqrt{5}}{2}$;\ \ 
$ -c^{1}_{7}
$,
$ -\xi_{7}^{1,2}$,
$ \frac{1+\sqrt{5}}{2}$,
$ -1$;\ \ 
$ -1$,
$ 1-c^{1}_{35}
-c^{4}_{35}
-c^{6}_{35}
+c^{7}_{35}
-c^{11}_{35}
$,
$ c^{1}_{7}
$;\ \ 
$ -\xi_{7}^{1,2}$,
$ -1+c^{4}_{35}
-c^{7}_{35}
+c^{11}_{35}
$;\ \ 
$ \xi_{7}^{1,2}$)

Not pseudo-unitary. 

\vskip 1ex 
\color{grey}

\noindent71. ind = $(6;35
)_{1}^{2}$:\ \ 
$d_i$ = ($1.0$,
$0.554$,
$0.770$,
$-0.342$,
$-0.618$,
$-1.246$) 

\vskip 0.7ex
\hangindent=3em \hangafter=1
$D^2=$ 3.956 = 
 $9+c^{1}_{35}
+3c^{4}_{35}
-3c^{5}_{35}
+c^{6}_{35}
-6c^{7}_{35}
+6c^{10}_{35}
+3c^{11}_{35}
$

\vskip 0.7ex
\hangindent=3em \hangafter=1
$T = ( 0,
\frac{3}{7},
\frac{3}{35},
\frac{8}{35},
\frac{4}{5},
\frac{2}{7} )
$,

\vskip 0.7ex
\hangindent=3em \hangafter=1
$S$ = ($ 1$,
$ \xi_{7}^{1,2}$,
$ -1+c^{1}_{35}
+c^{4}_{35}
-c^{5}_{35}
+c^{6}_{35}
-c^{7}_{35}
+c^{11}_{35}
$,
$ c^{4}_{35}
-c^{7}_{35}
+c^{10}_{35}
+c^{11}_{35}
$,
$ \frac{1-\sqrt{5}}{2}$,
$ -c^{1}_{7}
$;\ \ 
$ c^{1}_{7}
$,
$ \frac{1-\sqrt{5}}{2}$,
$ 1-c^{1}_{35}
-c^{4}_{35}
+c^{5}_{35}
-c^{6}_{35}
+c^{7}_{35}
-c^{11}_{35}
$,
$ c^{4}_{35}
-c^{7}_{35}
+c^{10}_{35}
+c^{11}_{35}
$,
$ 1$;\ \ 
$ \xi_{7}^{1,2}$,
$ -1$,
$ c^{1}_{7}
$,
$ -c^{4}_{35}
+c^{7}_{35}
-c^{10}_{35}
-c^{11}_{35}
$;\ \ 
$ -c^{1}_{7}
$,
$ -\xi_{7}^{1,2}$,
$ \frac{1-\sqrt{5}}{2}$;\ \ 
$ -1$,
$ -1+c^{1}_{35}
+c^{4}_{35}
-c^{5}_{35}
+c^{6}_{35}
-c^{7}_{35}
+c^{11}_{35}
$;\ \ 
$ -\xi_{7}^{1,2}$)

Not pseudo-unitary. 

\vskip 1ex 
\color{grey}

\noindent72. ind = $(6;35
)_{1}^{3}$:\ \ 
$d_i$ = ($1.0$,
$0.445$,
$0.495$,
$-0.275$,
$-0.618$,
$-0.801$) 

\vskip 0.7ex
\hangindent=3em \hangafter=1
$D^2=$ 2.544 = 
 $4+2c^{1}_{35}
-c^{4}_{35}
-6c^{5}_{35}
+2c^{6}_{35}
-5c^{7}_{35}
-9c^{10}_{35}
-c^{11}_{35}
$

\vskip 0.7ex
\hangindent=3em \hangafter=1
$T = ( 0,
\frac{3}{7},
\frac{12}{35},
\frac{22}{35},
\frac{1}{5},
\frac{1}{7} )
$,

\vskip 0.7ex
\hangindent=3em \hangafter=1
$S$ = ($ 1$,
$ -c^{2}_{7}
$,
$ -1+c^{1}_{35}
-c^{5}_{35}
+c^{6}_{35}
-c^{7}_{35}
-c^{10}_{35}
$,
$ -c^{4}_{35}
-c^{10}_{35}
-c^{11}_{35}
$,
$ \frac{1-\sqrt{5}}{2}$,
$ -\xi_{7}^{2,3}$;\ \ 
$ \xi_{7}^{2,3}$,
$ \frac{1-\sqrt{5}}{2}$,
$ 1-c^{1}_{35}
+c^{5}_{35}
-c^{6}_{35}
+c^{7}_{35}
+c^{10}_{35}
$,
$ -c^{4}_{35}
-c^{10}_{35}
-c^{11}_{35}
$,
$ 1$;\ \ 
$ -c^{2}_{7}
$,
$ -1$,
$ \xi_{7}^{2,3}$,
$ c^{4}_{35}
+c^{10}_{35}
+c^{11}_{35}
$;\ \ 
$ -\xi_{7}^{2,3}$,
$ c^{2}_{7}
$,
$ \frac{1-\sqrt{5}}{2}$;\ \ 
$ -1$,
$ -1+c^{1}_{35}
-c^{5}_{35}
+c^{6}_{35}
-c^{7}_{35}
-c^{10}_{35}
$;\ \ 
$ c^{2}_{7}
$)

Not pseudo-unitary. 

\vskip 1ex 
\color{grey}

\noindent73. ind = $(6;35
)_{1}^{23}$:\ \ 
$d_i$ = ($1.0$,
$0.554$,
$0.770$,
$-0.342$,
$-0.618$,
$-1.246$) 

\vskip 0.7ex
\hangindent=3em \hangafter=1
$D^2=$ 3.956 = 
 $9+c^{1}_{35}
+3c^{4}_{35}
-3c^{5}_{35}
+c^{6}_{35}
-6c^{7}_{35}
+6c^{10}_{35}
+3c^{11}_{35}
$

\vskip 0.7ex
\hangindent=3em \hangafter=1
$T = ( 0,
\frac{3}{7},
\frac{17}{35},
\frac{22}{35},
\frac{1}{5},
\frac{2}{7} )
$,

\vskip 0.7ex
\hangindent=3em \hangafter=1
$S$ = ($ 1$,
$ \xi_{7}^{1,2}$,
$ -1+c^{1}_{35}
+c^{4}_{35}
-c^{5}_{35}
+c^{6}_{35}
-c^{7}_{35}
+c^{11}_{35}
$,
$ c^{4}_{35}
-c^{7}_{35}
+c^{10}_{35}
+c^{11}_{35}
$,
$ \frac{1-\sqrt{5}}{2}$,
$ -c^{1}_{7}
$;\ \ 
$ c^{1}_{7}
$,
$ \frac{1-\sqrt{5}}{2}$,
$ 1-c^{1}_{35}
-c^{4}_{35}
+c^{5}_{35}
-c^{6}_{35}
+c^{7}_{35}
-c^{11}_{35}
$,
$ c^{4}_{35}
-c^{7}_{35}
+c^{10}_{35}
+c^{11}_{35}
$,
$ 1$;\ \ 
$ \xi_{7}^{1,2}$,
$ -1$,
$ c^{1}_{7}
$,
$ -c^{4}_{35}
+c^{7}_{35}
-c^{10}_{35}
-c^{11}_{35}
$;\ \ 
$ -c^{1}_{7}
$,
$ -\xi_{7}^{1,2}$,
$ \frac{1-\sqrt{5}}{2}$;\ \ 
$ -1$,
$ -1+c^{1}_{35}
+c^{4}_{35}
-c^{5}_{35}
+c^{6}_{35}
-c^{7}_{35}
+c^{11}_{35}
$;\ \ 
$ -\xi_{7}^{1,2}$)

Not pseudo-unitary. 

\vskip 1ex 
\color{grey}

\noindent74. ind = $(6;35
)_{1}^{31}$:\ \ 
$d_i$ = ($1.0$,
$0.445$,
$0.720$,
$1.618$,
$-0.801$,
$-1.297$) 

\vskip 0.7ex
\hangindent=3em \hangafter=1
$D^2=$ 6.661 = 
 $11-2c^{1}_{35}
+c^{4}_{35}
-4c^{5}_{35}
-2c^{6}_{35}
+5c^{7}_{35}
-6c^{10}_{35}
+c^{11}_{35}
$

\vskip 0.7ex
\hangindent=3em \hangafter=1
$T = ( 0,
\frac{3}{7},
\frac{29}{35},
\frac{2}{5},
\frac{1}{7},
\frac{19}{35} )
$,

\vskip 0.7ex
\hangindent=3em \hangafter=1
$S$ = ($ 1$,
$ -c^{2}_{7}
$,
$ c^{4}_{35}
+c^{11}_{35}
$,
$ \frac{1+\sqrt{5}}{2}$,
$ -\xi_{7}^{2,3}$,
$ 1-c^{1}_{35}
-c^{6}_{35}
+c^{7}_{35}
$;\ \ 
$ \xi_{7}^{2,3}$,
$ -1+c^{1}_{35}
+c^{6}_{35}
-c^{7}_{35}
$,
$ c^{4}_{35}
+c^{11}_{35}
$,
$ 1$,
$ \frac{1+\sqrt{5}}{2}$;\ \ 
$ -\xi_{7}^{2,3}$,
$ c^{2}_{7}
$,
$ \frac{1+\sqrt{5}}{2}$,
$ -1$;\ \ 
$ -1$,
$ 1-c^{1}_{35}
-c^{6}_{35}
+c^{7}_{35}
$,
$ \xi_{7}^{2,3}$;\ \ 
$ c^{2}_{7}
$,
$ -c^{4}_{35}
-c^{11}_{35}
$;\ \ 
$ -c^{2}_{7}
$)

Not pseudo-unitary. 

\vskip 1ex 
\color{grey}

\noindent75. ind = $(6;35
)_{1}^{16}$:\ \ 
$d_i$ = ($1.0$,
$0.554$,
$0.897$,
$1.618$,
$-1.246$,
$-2.17$) 

\vskip 0.7ex
\hangindent=3em \hangafter=1
$D^2=$ 10.358 = 
 $\frac{35}{\sqrt{35}\mathrm{i}}s^{1}_{35}
+\frac{35}{\sqrt{35}\mathrm{i}}s^{4}_{35}
-\frac{35}{\sqrt{35}\mathrm{i}}s^{6}_{35}
+\frac{35}{\sqrt{35}\mathrm{i}}s^{11}_{35}
$

\vskip 0.7ex
\hangindent=3em \hangafter=1
$T = ( 0,
\frac{3}{7},
\frac{29}{35},
\frac{2}{5},
\frac{2}{7},
\frac{24}{35} )
$,

\vskip 0.7ex
\hangindent=3em \hangafter=1
$S$ = ($ 1$,
$ \xi_{7}^{1,2}$,
$ 1-c^{4}_{35}
+c^{7}_{35}
-c^{11}_{35}
$,
$ \frac{1+\sqrt{5}}{2}$,
$ -c^{1}_{7}
$,
$ 1-c^{1}_{35}
-c^{4}_{35}
-c^{6}_{35}
+c^{7}_{35}
-c^{11}_{35}
$;\ \ 
$ c^{1}_{7}
$,
$ -1+c^{1}_{35}
+c^{4}_{35}
+c^{6}_{35}
-c^{7}_{35}
+c^{11}_{35}
$,
$ 1-c^{4}_{35}
+c^{7}_{35}
-c^{11}_{35}
$,
$ 1$,
$ \frac{1+\sqrt{5}}{2}$;\ \ 
$ -c^{1}_{7}
$,
$ -\xi_{7}^{1,2}$,
$ \frac{1+\sqrt{5}}{2}$,
$ -1$;\ \ 
$ -1$,
$ 1-c^{1}_{35}
-c^{4}_{35}
-c^{6}_{35}
+c^{7}_{35}
-c^{11}_{35}
$,
$ c^{1}_{7}
$;\ \ 
$ -\xi_{7}^{1,2}$,
$ -1+c^{4}_{35}
-c^{7}_{35}
+c^{11}_{35}
$;\ \ 
$ \xi_{7}^{1,2}$)

Not pseudo-unitary. 

\vskip 1ex 
\color{grey}

\noindent76. ind = $(6;35
)_{1}^{17}$:\ \ 
$d_i$ = ($1.0$,
$0.445$,
$0.495$,
$-0.275$,
$-0.618$,
$-0.801$) 

\vskip 0.7ex
\hangindent=3em \hangafter=1
$D^2=$ 2.544 = 
 $4+2c^{1}_{35}
-c^{4}_{35}
-6c^{5}_{35}
+2c^{6}_{35}
-5c^{7}_{35}
-9c^{10}_{35}
-c^{11}_{35}
$

\vskip 0.7ex
\hangindent=3em \hangafter=1
$T = ( 0,
\frac{3}{7},
\frac{33}{35},
\frac{8}{35},
\frac{4}{5},
\frac{1}{7} )
$,

\vskip 0.7ex
\hangindent=3em \hangafter=1
$S$ = ($ 1$,
$ -c^{2}_{7}
$,
$ -1+c^{1}_{35}
-c^{5}_{35}
+c^{6}_{35}
-c^{7}_{35}
-c^{10}_{35}
$,
$ -c^{4}_{35}
-c^{10}_{35}
-c^{11}_{35}
$,
$ \frac{1-\sqrt{5}}{2}$,
$ -\xi_{7}^{2,3}$;\ \ 
$ \xi_{7}^{2,3}$,
$ \frac{1-\sqrt{5}}{2}$,
$ 1-c^{1}_{35}
+c^{5}_{35}
-c^{6}_{35}
+c^{7}_{35}
+c^{10}_{35}
$,
$ -c^{4}_{35}
-c^{10}_{35}
-c^{11}_{35}
$,
$ 1$;\ \ 
$ -c^{2}_{7}
$,
$ -1$,
$ \xi_{7}^{2,3}$,
$ c^{4}_{35}
+c^{10}_{35}
+c^{11}_{35}
$;\ \ 
$ -\xi_{7}^{2,3}$,
$ c^{2}_{7}
$,
$ \frac{1-\sqrt{5}}{2}$;\ \ 
$ -1$,
$ -1+c^{1}_{35}
-c^{5}_{35}
+c^{6}_{35}
-c^{7}_{35}
-c^{10}_{35}
$;\ \ 
$ c^{2}_{7}
$)

Not pseudo-unitary. 

\vskip 1ex 
\color{grey}

\noindent77. ind = $(6;35
)_{1}^{18}$:\ \ 
$d_i$ = ($1.0$,
$0.445$,
$0.495$,
$-0.275$,
$-0.618$,
$-0.801$) 

\vskip 0.7ex
\hangindent=3em \hangafter=1
$D^2=$ 2.544 = 
 $4+2c^{1}_{35}
-c^{4}_{35}
-6c^{5}_{35}
+2c^{6}_{35}
-5c^{7}_{35}
-9c^{10}_{35}
-c^{11}_{35}
$

\vskip 0.7ex
\hangindent=3em \hangafter=1
$T = ( 0,
\frac{4}{7},
\frac{2}{35},
\frac{27}{35},
\frac{1}{5},
\frac{6}{7} )
$,

\vskip 0.7ex
\hangindent=3em \hangafter=1
$S$ = ($ 1$,
$ -c^{2}_{7}
$,
$ -1+c^{1}_{35}
-c^{5}_{35}
+c^{6}_{35}
-c^{7}_{35}
-c^{10}_{35}
$,
$ -c^{4}_{35}
-c^{10}_{35}
-c^{11}_{35}
$,
$ \frac{1-\sqrt{5}}{2}$,
$ -\xi_{7}^{2,3}$;\ \ 
$ \xi_{7}^{2,3}$,
$ \frac{1-\sqrt{5}}{2}$,
$ 1-c^{1}_{35}
+c^{5}_{35}
-c^{6}_{35}
+c^{7}_{35}
+c^{10}_{35}
$,
$ -c^{4}_{35}
-c^{10}_{35}
-c^{11}_{35}
$,
$ 1$;\ \ 
$ -c^{2}_{7}
$,
$ -1$,
$ \xi_{7}^{2,3}$,
$ c^{4}_{35}
+c^{10}_{35}
+c^{11}_{35}
$;\ \ 
$ -\xi_{7}^{2,3}$,
$ c^{2}_{7}
$,
$ \frac{1-\sqrt{5}}{2}$;\ \ 
$ -1$,
$ -1+c^{1}_{35}
-c^{5}_{35}
+c^{6}_{35}
-c^{7}_{35}
-c^{10}_{35}
$;\ \ 
$ c^{2}_{7}
$)

Not pseudo-unitary. 

\vskip 1ex 
\color{grey}

\noindent78. ind = $(6;35
)_{1}^{19}$:\ \ 
$d_i$ = ($1.0$,
$0.554$,
$0.897$,
$1.618$,
$-1.246$,
$-2.17$) 

\vskip 0.7ex
\hangindent=3em \hangafter=1
$D^2=$ 10.358 = 
 $\frac{35}{\sqrt{35}\mathrm{i}}s^{1}_{35}
+\frac{35}{\sqrt{35}\mathrm{i}}s^{4}_{35}
-\frac{35}{\sqrt{35}\mathrm{i}}s^{6}_{35}
+\frac{35}{\sqrt{35}\mathrm{i}}s^{11}_{35}
$

\vskip 0.7ex
\hangindent=3em \hangafter=1
$T = ( 0,
\frac{4}{7},
\frac{6}{35},
\frac{3}{5},
\frac{5}{7},
\frac{11}{35} )
$,

\vskip 0.7ex
\hangindent=3em \hangafter=1
$S$ = ($ 1$,
$ \xi_{7}^{1,2}$,
$ 1-c^{4}_{35}
+c^{7}_{35}
-c^{11}_{35}
$,
$ \frac{1+\sqrt{5}}{2}$,
$ -c^{1}_{7}
$,
$ 1-c^{1}_{35}
-c^{4}_{35}
-c^{6}_{35}
+c^{7}_{35}
-c^{11}_{35}
$;\ \ 
$ c^{1}_{7}
$,
$ -1+c^{1}_{35}
+c^{4}_{35}
+c^{6}_{35}
-c^{7}_{35}
+c^{11}_{35}
$,
$ 1-c^{4}_{35}
+c^{7}_{35}
-c^{11}_{35}
$,
$ 1$,
$ \frac{1+\sqrt{5}}{2}$;\ \ 
$ -c^{1}_{7}
$,
$ -\xi_{7}^{1,2}$,
$ \frac{1+\sqrt{5}}{2}$,
$ -1$;\ \ 
$ -1$,
$ 1-c^{1}_{35}
-c^{4}_{35}
-c^{6}_{35}
+c^{7}_{35}
-c^{11}_{35}
$,
$ c^{1}_{7}
$;\ \ 
$ -\xi_{7}^{1,2}$,
$ -1+c^{4}_{35}
-c^{7}_{35}
+c^{11}_{35}
$;\ \ 
$ \xi_{7}^{1,2}$)

Not pseudo-unitary. 

\vskip 1ex 
\color{grey}

\noindent79. ind = $(6;35
)_{1}^{4}$:\ \ 
$d_i$ = ($1.0$,
$0.445$,
$0.720$,
$1.618$,
$-0.801$,
$-1.297$) 

\vskip 0.7ex
\hangindent=3em \hangafter=1
$D^2=$ 6.661 = 
 $11-2c^{1}_{35}
+c^{4}_{35}
-4c^{5}_{35}
-2c^{6}_{35}
+5c^{7}_{35}
-6c^{10}_{35}
+c^{11}_{35}
$

\vskip 0.7ex
\hangindent=3em \hangafter=1
$T = ( 0,
\frac{4}{7},
\frac{6}{35},
\frac{3}{5},
\frac{6}{7},
\frac{16}{35} )
$,

\vskip 0.7ex
\hangindent=3em \hangafter=1
$S$ = ($ 1$,
$ -c^{2}_{7}
$,
$ c^{4}_{35}
+c^{11}_{35}
$,
$ \frac{1+\sqrt{5}}{2}$,
$ -\xi_{7}^{2,3}$,
$ 1-c^{1}_{35}
-c^{6}_{35}
+c^{7}_{35}
$;\ \ 
$ \xi_{7}^{2,3}$,
$ -1+c^{1}_{35}
+c^{6}_{35}
-c^{7}_{35}
$,
$ c^{4}_{35}
+c^{11}_{35}
$,
$ 1$,
$ \frac{1+\sqrt{5}}{2}$;\ \ 
$ -\xi_{7}^{2,3}$,
$ c^{2}_{7}
$,
$ \frac{1+\sqrt{5}}{2}$,
$ -1$;\ \ 
$ -1$,
$ 1-c^{1}_{35}
-c^{6}_{35}
+c^{7}_{35}
$,
$ \xi_{7}^{2,3}$;\ \ 
$ c^{2}_{7}
$,
$ -c^{4}_{35}
-c^{11}_{35}
$;\ \ 
$ -c^{2}_{7}
$)

Not pseudo-unitary. 

\vskip 1ex 
\color{grey}

\noindent80. ind = $(6;35
)_{1}^{12}$:\ \ 
$d_i$ = ($1.0$,
$0.554$,
$0.770$,
$-0.342$,
$-0.618$,
$-1.246$) 

\vskip 0.7ex
\hangindent=3em \hangafter=1
$D^2=$ 3.956 = 
 $9+c^{1}_{35}
+3c^{4}_{35}
-3c^{5}_{35}
+c^{6}_{35}
-6c^{7}_{35}
+6c^{10}_{35}
+3c^{11}_{35}
$

\vskip 0.7ex
\hangindent=3em \hangafter=1
$T = ( 0,
\frac{4}{7},
\frac{18}{35},
\frac{13}{35},
\frac{4}{5},
\frac{5}{7} )
$,

\vskip 0.7ex
\hangindent=3em \hangafter=1
$S$ = ($ 1$,
$ \xi_{7}^{1,2}$,
$ -1+c^{1}_{35}
+c^{4}_{35}
-c^{5}_{35}
+c^{6}_{35}
-c^{7}_{35}
+c^{11}_{35}
$,
$ c^{4}_{35}
-c^{7}_{35}
+c^{10}_{35}
+c^{11}_{35}
$,
$ \frac{1-\sqrt{5}}{2}$,
$ -c^{1}_{7}
$;\ \ 
$ c^{1}_{7}
$,
$ \frac{1-\sqrt{5}}{2}$,
$ 1-c^{1}_{35}
-c^{4}_{35}
+c^{5}_{35}
-c^{6}_{35}
+c^{7}_{35}
-c^{11}_{35}
$,
$ c^{4}_{35}
-c^{7}_{35}
+c^{10}_{35}
+c^{11}_{35}
$,
$ 1$;\ \ 
$ \xi_{7}^{1,2}$,
$ -1$,
$ c^{1}_{7}
$,
$ -c^{4}_{35}
+c^{7}_{35}
-c^{10}_{35}
-c^{11}_{35}
$;\ \ 
$ -c^{1}_{7}
$,
$ -\xi_{7}^{1,2}$,
$ \frac{1-\sqrt{5}}{2}$;\ \ 
$ -1$,
$ -1+c^{1}_{35}
+c^{4}_{35}
-c^{5}_{35}
+c^{6}_{35}
-c^{7}_{35}
+c^{11}_{35}
$;\ \ 
$ -\xi_{7}^{1,2}$)

Not pseudo-unitary. 

\vskip 1ex 
\color{grey}

\noindent81. ind = $(6;35
)_{1}^{32}$:\ \ 
$d_i$ = ($1.0$,
$0.445$,
$0.495$,
$-0.275$,
$-0.618$,
$-0.801$) 

\vskip 0.7ex
\hangindent=3em \hangafter=1
$D^2=$ 2.544 = 
 $4+2c^{1}_{35}
-c^{4}_{35}
-6c^{5}_{35}
+2c^{6}_{35}
-5c^{7}_{35}
-9c^{10}_{35}
-c^{11}_{35}
$

\vskip 0.7ex
\hangindent=3em \hangafter=1
$T = ( 0,
\frac{4}{7},
\frac{23}{35},
\frac{13}{35},
\frac{4}{5},
\frac{6}{7} )
$,

\vskip 0.7ex
\hangindent=3em \hangafter=1
$S$ = ($ 1$,
$ -c^{2}_{7}
$,
$ -1+c^{1}_{35}
-c^{5}_{35}
+c^{6}_{35}
-c^{7}_{35}
-c^{10}_{35}
$,
$ -c^{4}_{35}
-c^{10}_{35}
-c^{11}_{35}
$,
$ \frac{1-\sqrt{5}}{2}$,
$ -\xi_{7}^{2,3}$;\ \ 
$ \xi_{7}^{2,3}$,
$ \frac{1-\sqrt{5}}{2}$,
$ 1-c^{1}_{35}
+c^{5}_{35}
-c^{6}_{35}
+c^{7}_{35}
+c^{10}_{35}
$,
$ -c^{4}_{35}
-c^{10}_{35}
-c^{11}_{35}
$,
$ 1$;\ \ 
$ -c^{2}_{7}
$,
$ -1$,
$ \xi_{7}^{2,3}$,
$ c^{4}_{35}
+c^{10}_{35}
+c^{11}_{35}
$;\ \ 
$ -\xi_{7}^{2,3}$,
$ c^{2}_{7}
$,
$ \frac{1-\sqrt{5}}{2}$;\ \ 
$ -1$,
$ -1+c^{1}_{35}
-c^{5}_{35}
+c^{6}_{35}
-c^{7}_{35}
-c^{10}_{35}
$;\ \ 
$ c^{2}_{7}
$)

Not pseudo-unitary. 

\vskip 1ex 
\color{grey}

\noindent82. ind = $(6;35
)_{1}^{33}$:\ \ 
$d_i$ = ($1.0$,
$0.554$,
$0.770$,
$-0.342$,
$-0.618$,
$-1.246$) 

\vskip 0.7ex
\hangindent=3em \hangafter=1
$D^2=$ 3.956 = 
 $9+c^{1}_{35}
+3c^{4}_{35}
-3c^{5}_{35}
+c^{6}_{35}
-6c^{7}_{35}
+6c^{10}_{35}
+3c^{11}_{35}
$

\vskip 0.7ex
\hangindent=3em \hangafter=1
$T = ( 0,
\frac{4}{7},
\frac{32}{35},
\frac{27}{35},
\frac{1}{5},
\frac{5}{7} )
$,

\vskip 0.7ex
\hangindent=3em \hangafter=1
$S$ = ($ 1$,
$ \xi_{7}^{1,2}$,
$ -1+c^{1}_{35}
+c^{4}_{35}
-c^{5}_{35}
+c^{6}_{35}
-c^{7}_{35}
+c^{11}_{35}
$,
$ c^{4}_{35}
-c^{7}_{35}
+c^{10}_{35}
+c^{11}_{35}
$,
$ \frac{1-\sqrt{5}}{2}$,
$ -c^{1}_{7}
$;\ \ 
$ c^{1}_{7}
$,
$ \frac{1-\sqrt{5}}{2}$,
$ 1-c^{1}_{35}
-c^{4}_{35}
+c^{5}_{35}
-c^{6}_{35}
+c^{7}_{35}
-c^{11}_{35}
$,
$ c^{4}_{35}
-c^{7}_{35}
+c^{10}_{35}
+c^{11}_{35}
$,
$ 1$;\ \ 
$ \xi_{7}^{1,2}$,
$ -1$,
$ c^{1}_{7}
$,
$ -c^{4}_{35}
+c^{7}_{35}
-c^{10}_{35}
-c^{11}_{35}
$;\ \ 
$ -c^{1}_{7}
$,
$ -\xi_{7}^{1,2}$,
$ \frac{1-\sqrt{5}}{2}$;\ \ 
$ -1$,
$ -1+c^{1}_{35}
+c^{4}_{35}
-c^{5}_{35}
+c^{6}_{35}
-c^{7}_{35}
+c^{11}_{35}
$;\ \ 
$ -\xi_{7}^{1,2}$)

Not pseudo-unitary. 

\vskip 1ex 
\color{grey}

\noindent83. ind = $(6;35
)_{1}^{26}$:\ \ 
$d_i$ = ($1.0$,
$0.554$,
$0.897$,
$1.618$,
$-1.246$,
$-2.17$) 

\vskip 0.7ex
\hangindent=3em \hangafter=1
$D^2=$ 10.358 = 
 $\frac{35}{\sqrt{35}\mathrm{i}}s^{1}_{35}
+\frac{35}{\sqrt{35}\mathrm{i}}s^{4}_{35}
-\frac{35}{\sqrt{35}\mathrm{i}}s^{6}_{35}
+\frac{35}{\sqrt{35}\mathrm{i}}s^{11}_{35}
$

\vskip 0.7ex
\hangindent=3em \hangafter=1
$T = ( 0,
\frac{4}{7},
\frac{34}{35},
\frac{2}{5},
\frac{5}{7},
\frac{4}{35} )
$,

\vskip 0.7ex
\hangindent=3em \hangafter=1
$S$ = ($ 1$,
$ \xi_{7}^{1,2}$,
$ 1-c^{4}_{35}
+c^{7}_{35}
-c^{11}_{35}
$,
$ \frac{1+\sqrt{5}}{2}$,
$ -c^{1}_{7}
$,
$ 1-c^{1}_{35}
-c^{4}_{35}
-c^{6}_{35}
+c^{7}_{35}
-c^{11}_{35}
$;\ \ 
$ c^{1}_{7}
$,
$ -1+c^{1}_{35}
+c^{4}_{35}
+c^{6}_{35}
-c^{7}_{35}
+c^{11}_{35}
$,
$ 1-c^{4}_{35}
+c^{7}_{35}
-c^{11}_{35}
$,
$ 1$,
$ \frac{1+\sqrt{5}}{2}$;\ \ 
$ -c^{1}_{7}
$,
$ -\xi_{7}^{1,2}$,
$ \frac{1+\sqrt{5}}{2}$,
$ -1$;\ \ 
$ -1$,
$ 1-c^{1}_{35}
-c^{4}_{35}
-c^{6}_{35}
+c^{7}_{35}
-c^{11}_{35}
$,
$ c^{1}_{7}
$;\ \ 
$ -\xi_{7}^{1,2}$,
$ -1+c^{4}_{35}
-c^{7}_{35}
+c^{11}_{35}
$;\ \ 
$ \xi_{7}^{1,2}$)

Not pseudo-unitary. 

\vskip 1ex 
\color{grey}

\noindent84. ind = $(6;35
)_{1}^{11}$:\ \ 
$d_i$ = ($1.0$,
$0.445$,
$0.720$,
$1.618$,
$-0.801$,
$-1.297$) 

\vskip 0.7ex
\hangindent=3em \hangafter=1
$D^2=$ 6.661 = 
 $11-2c^{1}_{35}
+c^{4}_{35}
-4c^{5}_{35}
-2c^{6}_{35}
+5c^{7}_{35}
-6c^{10}_{35}
+c^{11}_{35}
$

\vskip 0.7ex
\hangindent=3em \hangafter=1
$T = ( 0,
\frac{4}{7},
\frac{34}{35},
\frac{2}{5},
\frac{6}{7},
\frac{9}{35} )
$,

\vskip 0.7ex
\hangindent=3em \hangafter=1
$S$ = ($ 1$,
$ -c^{2}_{7}
$,
$ c^{4}_{35}
+c^{11}_{35}
$,
$ \frac{1+\sqrt{5}}{2}$,
$ -\xi_{7}^{2,3}$,
$ 1-c^{1}_{35}
-c^{6}_{35}
+c^{7}_{35}
$;\ \ 
$ \xi_{7}^{2,3}$,
$ -1+c^{1}_{35}
+c^{6}_{35}
-c^{7}_{35}
$,
$ c^{4}_{35}
+c^{11}_{35}
$,
$ 1$,
$ \frac{1+\sqrt{5}}{2}$;\ \ 
$ -\xi_{7}^{2,3}$,
$ c^{2}_{7}
$,
$ \frac{1+\sqrt{5}}{2}$,
$ -1$;\ \ 
$ -1$,
$ 1-c^{1}_{35}
-c^{6}_{35}
+c^{7}_{35}
$,
$ \xi_{7}^{2,3}$;\ \ 
$ c^{2}_{7}
$,
$ -c^{4}_{35}
-c^{11}_{35}
$;\ \ 
$ -c^{2}_{7}
$)

Not pseudo-unitary. 

\vskip 1ex 
\color{grey}

\noindent85. ind = $(6;35
)_{1}^{13}$:\ \ 
$d_i$ = ($1.0$,
$1.801$,
$2.246$,
$-0.618$,
$-1.113$,
$-1.388$) 

\vskip 0.7ex
\hangindent=3em \hangafter=1
$D^2=$ 12.846 = 
 $15-3c^{1}_{35}
-2c^{4}_{35}
+9c^{5}_{35}
-3c^{6}_{35}
-3c^{7}_{35}
+3c^{10}_{35}
-2c^{11}_{35}
$

\vskip 0.7ex
\hangindent=3em \hangafter=1
$T = ( 0,
\frac{6}{7},
\frac{2}{7},
\frac{1}{5},
\frac{2}{35},
\frac{17}{35} )
$,

\vskip 0.7ex
\hangindent=3em \hangafter=1
$S$ = ($ 1$,
$ \xi_{7}^{2}$,
$ \xi_{7}^{3}$,
$ \frac{1-\sqrt{5}}{2}$,
$ 1-c^{1}_{35}
+c^{5}_{35}
-c^{6}_{35}
+c^{10}_{35}
$,
$ 1-c^{1}_{35}
-c^{4}_{35}
+c^{5}_{35}
-c^{6}_{35}
-c^{11}_{35}
$;\ \ 
$ -\xi_{7}^{3}$,
$ 1$,
$ 1-c^{1}_{35}
+c^{5}_{35}
-c^{6}_{35}
+c^{10}_{35}
$,
$ -1+c^{1}_{35}
+c^{4}_{35}
-c^{5}_{35}
+c^{6}_{35}
+c^{11}_{35}
$,
$ \frac{1-\sqrt{5}}{2}$;\ \ 
$ -\xi_{7}^{2}$,
$ 1-c^{1}_{35}
-c^{4}_{35}
+c^{5}_{35}
-c^{6}_{35}
-c^{11}_{35}
$,
$ \frac{1-\sqrt{5}}{2}$,
$ -1+c^{1}_{35}
-c^{5}_{35}
+c^{6}_{35}
-c^{10}_{35}
$;\ \ 
$ -1$,
$ -\xi_{7}^{2}$,
$ -\xi_{7}^{3}$;\ \ 
$ \xi_{7}^{3}$,
$ -1$;\ \ 
$ \xi_{7}^{2}$)

Not pseudo-unitary. 

\vskip 1ex 
\color{grey}

\noindent86. ind = $(6;35
)_{1}^{27}$:\ \ 
$d_i$ = ($1.0$,
$1.801$,
$2.246$,
$-0.618$,
$-1.113$,
$-1.388$) 

\vskip 0.7ex
\hangindent=3em \hangafter=1
$D^2=$ 12.846 = 
 $15-3c^{1}_{35}
-2c^{4}_{35}
+9c^{5}_{35}
-3c^{6}_{35}
-3c^{7}_{35}
+3c^{10}_{35}
-2c^{11}_{35}
$

\vskip 0.7ex
\hangindent=3em \hangafter=1
$T = ( 0,
\frac{6}{7},
\frac{2}{7},
\frac{4}{5},
\frac{23}{35},
\frac{3}{35} )
$,

\vskip 0.7ex
\hangindent=3em \hangafter=1
$S$ = ($ 1$,
$ \xi_{7}^{2}$,
$ \xi_{7}^{3}$,
$ \frac{1-\sqrt{5}}{2}$,
$ 1-c^{1}_{35}
+c^{5}_{35}
-c^{6}_{35}
+c^{10}_{35}
$,
$ 1-c^{1}_{35}
-c^{4}_{35}
+c^{5}_{35}
-c^{6}_{35}
-c^{11}_{35}
$;\ \ 
$ -\xi_{7}^{3}$,
$ 1$,
$ 1-c^{1}_{35}
+c^{5}_{35}
-c^{6}_{35}
+c^{10}_{35}
$,
$ -1+c^{1}_{35}
+c^{4}_{35}
-c^{5}_{35}
+c^{6}_{35}
+c^{11}_{35}
$,
$ \frac{1-\sqrt{5}}{2}$;\ \ 
$ -\xi_{7}^{2}$,
$ 1-c^{1}_{35}
-c^{4}_{35}
+c^{5}_{35}
-c^{6}_{35}
-c^{11}_{35}
$,
$ \frac{1-\sqrt{5}}{2}$,
$ -1+c^{1}_{35}
-c^{5}_{35}
+c^{6}_{35}
-c^{10}_{35}
$;\ \ 
$ -1$,
$ -\xi_{7}^{2}$,
$ -\xi_{7}^{3}$;\ \ 
$ \xi_{7}^{3}$,
$ -1$;\ \ 
$ \xi_{7}^{2}$)

Not pseudo-unitary. 

\vskip 1ex 

 \color{black} \vskip 2ex

\noindent87. ind = $(6;56
)_{1}^{1}$:\ \ 
$d_i$ = ($1.0$,
$1.0$,
$1.801$,
$1.801$,
$2.246$,
$2.246$) 

\vskip 0.7ex
\hangindent=3em \hangafter=1
$D^2=$ 18.591 = 
 $12+6c^{1}_{7}
+2c^{2}_{7}
$

\vskip 0.7ex
\hangindent=3em \hangafter=1
$T = ( 0,
\frac{1}{4},
\frac{1}{7},
\frac{11}{28},
\frac{5}{7},
\frac{27}{28} )
$,

\vskip 0.7ex
\hangindent=3em \hangafter=1
$S$ = ($ 1$,
$ 1$,
$ \xi_{7}^{2}$,
$ \xi_{7}^{2}$,
$ \xi_{7}^{3}$,
$ \xi_{7}^{3}$;\ \ 
$ -1$,
$ \xi_{7}^{2}$,
$ -\xi_{7}^{2}$,
$ \xi_{7}^{3}$,
$ -\xi_{7}^{3}$;\ \ 
$ -\xi_{7}^{3}$,
$ -\xi_{7}^{3}$,
$ 1$,
$ 1$;\ \ 
$ \xi_{7}^{3}$,
$ 1$,
$ -1$;\ \ 
$ -\xi_{7}^{2}$,
$ -\xi_{7}^{2}$;\ \ 
$ \xi_{7}^{2}$)

\vskip 1ex 
\color{grey}

\noindent88. ind = $(6;56
)_{1}^{13}$:\ \ 
$d_i$ = ($1.0$,
$1.0$,
$1.801$,
$1.801$,
$2.246$,
$2.246$) 

\vskip 0.7ex
\hangindent=3em \hangafter=1
$D^2=$ 18.591 = 
 $12+6c^{1}_{7}
+2c^{2}_{7}
$

\vskip 0.7ex
\hangindent=3em \hangafter=1
$T = ( 0,
\frac{1}{4},
\frac{6}{7},
\frac{3}{28},
\frac{2}{7},
\frac{15}{28} )
$,

\vskip 0.7ex
\hangindent=3em \hangafter=1
$S$ = ($ 1$,
$ 1$,
$ \xi_{7}^{2}$,
$ \xi_{7}^{2}$,
$ \xi_{7}^{3}$,
$ \xi_{7}^{3}$;\ \ 
$ -1$,
$ \xi_{7}^{2}$,
$ -\xi_{7}^{2}$,
$ \xi_{7}^{3}$,
$ -\xi_{7}^{3}$;\ \ 
$ -\xi_{7}^{3}$,
$ -\xi_{7}^{3}$,
$ 1$,
$ 1$;\ \ 
$ \xi_{7}^{3}$,
$ 1$,
$ -1$;\ \ 
$ -\xi_{7}^{2}$,
$ -\xi_{7}^{2}$;\ \ 
$ \xi_{7}^{2}$)

\vskip 1ex 
\color{grey}

\noindent89. ind = $(6;56
)_{1}^{15}$:\ \ 
$d_i$ = ($1.0$,
$1.0$,
$1.801$,
$1.801$,
$2.246$,
$2.246$) 

\vskip 0.7ex
\hangindent=3em \hangafter=1
$D^2=$ 18.591 = 
 $12+6c^{1}_{7}
+2c^{2}_{7}
$

\vskip 0.7ex
\hangindent=3em \hangafter=1
$T = ( 0,
\frac{3}{4},
\frac{1}{7},
\frac{25}{28},
\frac{5}{7},
\frac{13}{28} )
$,

\vskip 0.7ex
\hangindent=3em \hangafter=1
$S$ = ($ 1$,
$ 1$,
$ \xi_{7}^{2}$,
$ \xi_{7}^{2}$,
$ \xi_{7}^{3}$,
$ \xi_{7}^{3}$;\ \ 
$ -1$,
$ \xi_{7}^{2}$,
$ -\xi_{7}^{2}$,
$ \xi_{7}^{3}$,
$ -\xi_{7}^{3}$;\ \ 
$ -\xi_{7}^{3}$,
$ -\xi_{7}^{3}$,
$ 1$,
$ 1$;\ \ 
$ \xi_{7}^{3}$,
$ 1$,
$ -1$;\ \ 
$ -\xi_{7}^{2}$,
$ -\xi_{7}^{2}$;\ \ 
$ \xi_{7}^{2}$)

\vskip 1ex 
\color{grey}

\noindent90. ind = $(6;56
)_{1}^{27}$:\ \ 
$d_i$ = ($1.0$,
$1.0$,
$1.801$,
$1.801$,
$2.246$,
$2.246$) 

\vskip 0.7ex
\hangindent=3em \hangafter=1
$D^2=$ 18.591 = 
 $12+6c^{1}_{7}
+2c^{2}_{7}
$

\vskip 0.7ex
\hangindent=3em \hangafter=1
$T = ( 0,
\frac{3}{4},
\frac{6}{7},
\frac{17}{28},
\frac{2}{7},
\frac{1}{28} )
$,

\vskip 0.7ex
\hangindent=3em \hangafter=1
$S$ = ($ 1$,
$ 1$,
$ \xi_{7}^{2}$,
$ \xi_{7}^{2}$,
$ \xi_{7}^{3}$,
$ \xi_{7}^{3}$;\ \ 
$ -1$,
$ \xi_{7}^{2}$,
$ -\xi_{7}^{2}$,
$ \xi_{7}^{3}$,
$ -\xi_{7}^{3}$;\ \ 
$ -\xi_{7}^{3}$,
$ -\xi_{7}^{3}$,
$ 1$,
$ 1$;\ \ 
$ \xi_{7}^{3}$,
$ 1$,
$ -1$;\ \ 
$ -\xi_{7}^{2}$,
$ -\xi_{7}^{2}$;\ \ 
$ \xi_{7}^{2}$)

\vskip 1ex 
\color{grey}

\noindent91. ind = $(6;56
)_{1}^{3}$:\ \ 
$d_i$ = ($1.0$,
$0.445$,
$0.445$,
$1.0$,
$-0.801$,
$-0.801$) 

\vskip 0.7ex
\hangindent=3em \hangafter=1
$D^2=$ 3.682 = 
 $6-4c^{1}_{7}
-6c^{2}_{7}
$

\vskip 0.7ex
\hangindent=3em \hangafter=1
$T = ( 0,
\frac{3}{7},
\frac{5}{28},
\frac{3}{4},
\frac{1}{7},
\frac{25}{28} )
$,

\vskip 0.7ex
\hangindent=3em \hangafter=1
$S$ = ($ 1$,
$ -c^{2}_{7}
$,
$ -c^{2}_{7}
$,
$ 1$,
$ -\xi_{7}^{2,3}$,
$ -\xi_{7}^{2,3}$;\ \ 
$ \xi_{7}^{2,3}$,
$ \xi_{7}^{2,3}$,
$ -c^{2}_{7}
$,
$ 1$,
$ 1$;\ \ 
$ -\xi_{7}^{2,3}$,
$ c^{2}_{7}
$,
$ 1$,
$ -1$;\ \ 
$ -1$,
$ -\xi_{7}^{2,3}$,
$ \xi_{7}^{2,3}$;\ \ 
$ c^{2}_{7}
$,
$ c^{2}_{7}
$;\ \ 
$ -c^{2}_{7}
$)

Not pseudo-unitary. 

\vskip 1ex 
\color{grey}

\noindent92. ind = $(6;56
)_{1}^{23}$:\ \ 
$d_i$ = ($1.0$,
$0.554$,
$0.554$,
$1.0$,
$-1.246$,
$-1.246$) 

\vskip 0.7ex
\hangindent=3em \hangafter=1
$D^2=$ 5.725 = 
 $10-2c^{1}_{7}
+4c^{2}_{7}
$

\vskip 0.7ex
\hangindent=3em \hangafter=1
$T = ( 0,
\frac{3}{7},
\frac{5}{28},
\frac{3}{4},
\frac{2}{7},
\frac{1}{28} )
$,

\vskip 0.7ex
\hangindent=3em \hangafter=1
$S$ = ($ 1$,
$ \xi_{7}^{1,2}$,
$ \xi_{7}^{1,2}$,
$ 1$,
$ -c^{1}_{7}
$,
$ -c^{1}_{7}
$;\ \ 
$ c^{1}_{7}
$,
$ c^{1}_{7}
$,
$ \xi_{7}^{1,2}$,
$ 1$,
$ 1$;\ \ 
$ -c^{1}_{7}
$,
$ -\xi_{7}^{1,2}$,
$ 1$,
$ -1$;\ \ 
$ -1$,
$ -c^{1}_{7}
$,
$ c^{1}_{7}
$;\ \ 
$ -\xi_{7}^{1,2}$,
$ -\xi_{7}^{1,2}$;\ \ 
$ \xi_{7}^{1,2}$)

Not pseudo-unitary. 

\vskip 1ex 
\color{grey}

\noindent93. ind = $(6;56
)_{1}^{17}$:\ \ 
$d_i$ = ($1.0$,
$0.445$,
$0.445$,
$1.0$,
$-0.801$,
$-0.801$) 

\vskip 0.7ex
\hangindent=3em \hangafter=1
$D^2=$ 3.682 = 
 $6-4c^{1}_{7}
-6c^{2}_{7}
$

\vskip 0.7ex
\hangindent=3em \hangafter=1
$T = ( 0,
\frac{3}{7},
\frac{19}{28},
\frac{1}{4},
\frac{1}{7},
\frac{11}{28} )
$,

\vskip 0.7ex
\hangindent=3em \hangafter=1
$S$ = ($ 1$,
$ -c^{2}_{7}
$,
$ -c^{2}_{7}
$,
$ 1$,
$ -\xi_{7}^{2,3}$,
$ -\xi_{7}^{2,3}$;\ \ 
$ \xi_{7}^{2,3}$,
$ \xi_{7}^{2,3}$,
$ -c^{2}_{7}
$,
$ 1$,
$ 1$;\ \ 
$ -\xi_{7}^{2,3}$,
$ c^{2}_{7}
$,
$ 1$,
$ -1$;\ \ 
$ -1$,
$ -\xi_{7}^{2,3}$,
$ \xi_{7}^{2,3}$;\ \ 
$ c^{2}_{7}
$,
$ c^{2}_{7}
$;\ \ 
$ -c^{2}_{7}
$)

Not pseudo-unitary. 

\vskip 1ex 
\color{grey}

\noindent94. ind = $(6;56
)_{1}^{9}$:\ \ 
$d_i$ = ($1.0$,
$0.554$,
$0.554$,
$1.0$,
$-1.246$,
$-1.246$) 

\vskip 0.7ex
\hangindent=3em \hangafter=1
$D^2=$ 5.725 = 
 $10-2c^{1}_{7}
+4c^{2}_{7}
$

\vskip 0.7ex
\hangindent=3em \hangafter=1
$T = ( 0,
\frac{3}{7},
\frac{19}{28},
\frac{1}{4},
\frac{2}{7},
\frac{15}{28} )
$,

\vskip 0.7ex
\hangindent=3em \hangafter=1
$S$ = ($ 1$,
$ \xi_{7}^{1,2}$,
$ \xi_{7}^{1,2}$,
$ 1$,
$ -c^{1}_{7}
$,
$ -c^{1}_{7}
$;\ \ 
$ c^{1}_{7}
$,
$ c^{1}_{7}
$,
$ \xi_{7}^{1,2}$,
$ 1$,
$ 1$;\ \ 
$ -c^{1}_{7}
$,
$ -\xi_{7}^{1,2}$,
$ 1$,
$ -1$;\ \ 
$ -1$,
$ -c^{1}_{7}
$,
$ c^{1}_{7}
$;\ \ 
$ -\xi_{7}^{1,2}$,
$ -\xi_{7}^{1,2}$;\ \ 
$ \xi_{7}^{1,2}$)

Not pseudo-unitary. 

\vskip 1ex 
\color{grey}

\noindent95. ind = $(6;56
)_{1}^{19}$:\ \ 
$d_i$ = ($1.0$,
$0.554$,
$0.554$,
$1.0$,
$-1.246$,
$-1.246$) 

\vskip 0.7ex
\hangindent=3em \hangafter=1
$D^2=$ 5.725 = 
 $10-2c^{1}_{7}
+4c^{2}_{7}
$

\vskip 0.7ex
\hangindent=3em \hangafter=1
$T = ( 0,
\frac{4}{7},
\frac{9}{28},
\frac{3}{4},
\frac{5}{7},
\frac{13}{28} )
$,

\vskip 0.7ex
\hangindent=3em \hangafter=1
$S$ = ($ 1$,
$ \xi_{7}^{1,2}$,
$ \xi_{7}^{1,2}$,
$ 1$,
$ -c^{1}_{7}
$,
$ -c^{1}_{7}
$;\ \ 
$ c^{1}_{7}
$,
$ c^{1}_{7}
$,
$ \xi_{7}^{1,2}$,
$ 1$,
$ 1$;\ \ 
$ -c^{1}_{7}
$,
$ -\xi_{7}^{1,2}$,
$ 1$,
$ -1$;\ \ 
$ -1$,
$ -c^{1}_{7}
$,
$ c^{1}_{7}
$;\ \ 
$ -\xi_{7}^{1,2}$,
$ -\xi_{7}^{1,2}$;\ \ 
$ \xi_{7}^{1,2}$)

Not pseudo-unitary. 

\vskip 1ex 
\color{grey}

\noindent96. ind = $(6;56
)_{1}^{11}$:\ \ 
$d_i$ = ($1.0$,
$0.445$,
$0.445$,
$1.0$,
$-0.801$,
$-0.801$) 

\vskip 0.7ex
\hangindent=3em \hangafter=1
$D^2=$ 3.682 = 
 $6-4c^{1}_{7}
-6c^{2}_{7}
$

\vskip 0.7ex
\hangindent=3em \hangafter=1
$T = ( 0,
\frac{4}{7},
\frac{9}{28},
\frac{3}{4},
\frac{6}{7},
\frac{17}{28} )
$,

\vskip 0.7ex
\hangindent=3em \hangafter=1
$S$ = ($ 1$,
$ -c^{2}_{7}
$,
$ -c^{2}_{7}
$,
$ 1$,
$ -\xi_{7}^{2,3}$,
$ -\xi_{7}^{2,3}$;\ \ 
$ \xi_{7}^{2,3}$,
$ \xi_{7}^{2,3}$,
$ -c^{2}_{7}
$,
$ 1$,
$ 1$;\ \ 
$ -\xi_{7}^{2,3}$,
$ c^{2}_{7}
$,
$ 1$,
$ -1$;\ \ 
$ -1$,
$ -\xi_{7}^{2,3}$,
$ \xi_{7}^{2,3}$;\ \ 
$ c^{2}_{7}
$,
$ c^{2}_{7}
$;\ \ 
$ -c^{2}_{7}
$)

Not pseudo-unitary. 

\vskip 1ex 
\color{grey}

\noindent97. ind = $(6;56
)_{1}^{5}$:\ \ 
$d_i$ = ($1.0$,
$0.554$,
$0.554$,
$1.0$,
$-1.246$,
$-1.246$) 

\vskip 0.7ex
\hangindent=3em \hangafter=1
$D^2=$ 5.725 = 
 $10-2c^{1}_{7}
+4c^{2}_{7}
$

\vskip 0.7ex
\hangindent=3em \hangafter=1
$T = ( 0,
\frac{4}{7},
\frac{23}{28},
\frac{1}{4},
\frac{5}{7},
\frac{27}{28} )
$,

\vskip 0.7ex
\hangindent=3em \hangafter=1
$S$ = ($ 1$,
$ \xi_{7}^{1,2}$,
$ \xi_{7}^{1,2}$,
$ 1$,
$ -c^{1}_{7}
$,
$ -c^{1}_{7}
$;\ \ 
$ c^{1}_{7}
$,
$ c^{1}_{7}
$,
$ \xi_{7}^{1,2}$,
$ 1$,
$ 1$;\ \ 
$ -c^{1}_{7}
$,
$ -\xi_{7}^{1,2}$,
$ 1$,
$ -1$;\ \ 
$ -1$,
$ -c^{1}_{7}
$,
$ c^{1}_{7}
$;\ \ 
$ -\xi_{7}^{1,2}$,
$ -\xi_{7}^{1,2}$;\ \ 
$ \xi_{7}^{1,2}$)

Not pseudo-unitary. 

\vskip 1ex 
\color{grey}

\noindent98. ind = $(6;56
)_{1}^{25}$:\ \ 
$d_i$ = ($1.0$,
$0.445$,
$0.445$,
$1.0$,
$-0.801$,
$-0.801$) 

\vskip 0.7ex
\hangindent=3em \hangafter=1
$D^2=$ 3.682 = 
 $6-4c^{1}_{7}
-6c^{2}_{7}
$

\vskip 0.7ex
\hangindent=3em \hangafter=1
$T = ( 0,
\frac{4}{7},
\frac{23}{28},
\frac{1}{4},
\frac{6}{7},
\frac{3}{28} )
$,

\vskip 0.7ex
\hangindent=3em \hangafter=1
$S$ = ($ 1$,
$ -c^{2}_{7}
$,
$ -c^{2}_{7}
$,
$ 1$,
$ -\xi_{7}^{2,3}$,
$ -\xi_{7}^{2,3}$;\ \ 
$ \xi_{7}^{2,3}$,
$ \xi_{7}^{2,3}$,
$ -c^{2}_{7}
$,
$ 1$,
$ 1$;\ \ 
$ -\xi_{7}^{2,3}$,
$ c^{2}_{7}
$,
$ 1$,
$ -1$;\ \ 
$ -1$,
$ -\xi_{7}^{2,3}$,
$ \xi_{7}^{2,3}$;\ \ 
$ c^{2}_{7}
$,
$ c^{2}_{7}
$;\ \ 
$ -c^{2}_{7}
$)

Not pseudo-unitary. 

\vskip 1ex 

 \color{black} \vskip 2ex 
\color{blue}

\noindent99. ind = $(6;56
)_{2}^{1}$:\ \ 
$d_i$ = ($1.0$,
$1.801$,
$2.246$,
$-1.0$,
$-1.801$,
$-2.246$) 

\vskip 0.7ex
\hangindent=3em \hangafter=1
$D^2=$ 18.591 = 
 $12+6c^{1}_{7}
+2c^{2}_{7}
$

\vskip 0.7ex
\hangindent=3em \hangafter=1
$T = ( 0,
\frac{1}{7},
\frac{5}{7},
\frac{1}{4},
\frac{11}{28},
\frac{27}{28} )
$,

\vskip 0.7ex
\hangindent=3em \hangafter=1
$S$ = ($ 1$,
$ \xi_{7}^{2}$,
$ \xi_{7}^{3}$,
$ -1$,
$ -\xi_{7}^{2}$,
$ -\xi_{7}^{3}$;\ \ 
$ -\xi_{7}^{3}$,
$ 1$,
$ -\xi_{7}^{2}$,
$ \xi_{7}^{3}$,
$ -1$;\ \ 
$ -\xi_{7}^{2}$,
$ -\xi_{7}^{3}$,
$ -1$,
$ \xi_{7}^{2}$;\ \ 
$ -1$,
$ -\xi_{7}^{2}$,
$ -\xi_{7}^{3}$;\ \ 
$ \xi_{7}^{3}$,
$ -1$;\ \ 
$ \xi_{7}^{2}$)

Pseudo-unitary $\sim$  
$(6;56
)_{1}^{1}$

\vskip 1ex 
\color{grey}

\noindent100. ind = $(6;56
)_{2}^{15}$:\ \ 
$d_i$ = ($1.0$,
$1.801$,
$2.246$,
$-1.0$,
$-1.801$,
$-2.246$) 

\vskip 0.7ex
\hangindent=3em \hangafter=1
$D^2=$ 18.591 = 
 $12+6c^{1}_{7}
+2c^{2}_{7}
$

\vskip 0.7ex
\hangindent=3em \hangafter=1
$T = ( 0,
\frac{1}{7},
\frac{5}{7},
\frac{3}{4},
\frac{25}{28},
\frac{13}{28} )
$,

\vskip 0.7ex
\hangindent=3em \hangafter=1
$S$ = ($ 1$,
$ \xi_{7}^{2}$,
$ \xi_{7}^{3}$,
$ -1$,
$ -\xi_{7}^{2}$,
$ -\xi_{7}^{3}$;\ \ 
$ -\xi_{7}^{3}$,
$ 1$,
$ -\xi_{7}^{2}$,
$ \xi_{7}^{3}$,
$ -1$;\ \ 
$ -\xi_{7}^{2}$,
$ -\xi_{7}^{3}$,
$ -1$,
$ \xi_{7}^{2}$;\ \ 
$ -1$,
$ -\xi_{7}^{2}$,
$ -\xi_{7}^{3}$;\ \ 
$ \xi_{7}^{3}$,
$ -1$;\ \ 
$ \xi_{7}^{2}$)

Pseudo-unitary $\sim$  
$(6;56
)_{1}^{15}$

\vskip 1ex 
\color{grey}

\noindent101. ind = $(6;56
)_{2}^{23}$:\ \ 
$d_i$ = ($1.0$,
$0.554$,
$1.246$,
$-0.554$,
$-1.0$,
$-1.246$) 

\vskip 0.7ex
\hangindent=3em \hangafter=1
$D^2=$ 5.725 = 
 $10-2c^{1}_{7}
+4c^{2}_{7}
$

\vskip 0.7ex
\hangindent=3em \hangafter=1
$T = ( 0,
\frac{3}{7},
\frac{1}{28},
\frac{5}{28},
\frac{3}{4},
\frac{2}{7} )
$,

\vskip 0.7ex
\hangindent=3em \hangafter=1
$S$ = ($ 1$,
$ \xi_{7}^{1,2}$,
$ c^{1}_{7}
$,
$ -\xi_{7}^{1,2}$,
$ -1$,
$ -c^{1}_{7}
$;\ \ 
$ c^{1}_{7}
$,
$ -1$,
$ -c^{1}_{7}
$,
$ -\xi_{7}^{1,2}$,
$ 1$;\ \ 
$ \xi_{7}^{1,2}$,
$ -1$,
$ c^{1}_{7}
$,
$ \xi_{7}^{1,2}$;\ \ 
$ -c^{1}_{7}
$,
$ -\xi_{7}^{1,2}$,
$ -1$;\ \ 
$ -1$,
$ c^{1}_{7}
$;\ \ 
$ -\xi_{7}^{1,2}$)

Not pseudo-unitary. 

\vskip 1ex 
\color{grey}

\noindent102. ind = $(6;56
)_{2}^{17}$:\ \ 
$d_i$ = ($1.0$,
$0.445$,
$0.801$,
$-0.445$,
$-0.801$,
$-1.0$) 

\vskip 0.7ex
\hangindent=3em \hangafter=1
$D^2=$ 3.682 = 
 $6-4c^{1}_{7}
-6c^{2}_{7}
$

\vskip 0.7ex
\hangindent=3em \hangafter=1
$T = ( 0,
\frac{3}{7},
\frac{11}{28},
\frac{19}{28},
\frac{1}{7},
\frac{1}{4} )
$,

\vskip 0.7ex
\hangindent=3em \hangafter=1
$S$ = ($ 1$,
$ -c^{2}_{7}
$,
$ \xi_{7}^{2,3}$,
$ c^{2}_{7}
$,
$ -\xi_{7}^{2,3}$,
$ -1$;\ \ 
$ \xi_{7}^{2,3}$,
$ -1$,
$ -\xi_{7}^{2,3}$,
$ 1$,
$ c^{2}_{7}
$;\ \ 
$ -c^{2}_{7}
$,
$ -1$,
$ -c^{2}_{7}
$,
$ \xi_{7}^{2,3}$;\ \ 
$ -\xi_{7}^{2,3}$,
$ -1$,
$ c^{2}_{7}
$;\ \ 
$ c^{2}_{7}
$,
$ \xi_{7}^{2,3}$;\ \ 
$ -1$)

Not pseudo-unitary. 

\vskip 1ex 
\color{grey}

\noindent103. ind = $(6;56
)_{2}^{9}$:\ \ 
$d_i$ = ($1.0$,
$0.554$,
$1.246$,
$-0.554$,
$-1.0$,
$-1.246$) 

\vskip 0.7ex
\hangindent=3em \hangafter=1
$D^2=$ 5.725 = 
 $10-2c^{1}_{7}
+4c^{2}_{7}
$

\vskip 0.7ex
\hangindent=3em \hangafter=1
$T = ( 0,
\frac{3}{7},
\frac{15}{28},
\frac{19}{28},
\frac{1}{4},
\frac{2}{7} )
$,

\vskip 0.7ex
\hangindent=3em \hangafter=1
$S$ = ($ 1$,
$ \xi_{7}^{1,2}$,
$ c^{1}_{7}
$,
$ -\xi_{7}^{1,2}$,
$ -1$,
$ -c^{1}_{7}
$;\ \ 
$ c^{1}_{7}
$,
$ -1$,
$ -c^{1}_{7}
$,
$ -\xi_{7}^{1,2}$,
$ 1$;\ \ 
$ \xi_{7}^{1,2}$,
$ -1$,
$ c^{1}_{7}
$,
$ \xi_{7}^{1,2}$;\ \ 
$ -c^{1}_{7}
$,
$ -\xi_{7}^{1,2}$,
$ -1$;\ \ 
$ -1$,
$ c^{1}_{7}
$;\ \ 
$ -\xi_{7}^{1,2}$)

Not pseudo-unitary. 

\vskip 1ex 
\color{grey}

\noindent104. ind = $(6;56
)_{2}^{3}$:\ \ 
$d_i$ = ($1.0$,
$0.445$,
$0.801$,
$-0.445$,
$-0.801$,
$-1.0$) 

\vskip 0.7ex
\hangindent=3em \hangafter=1
$D^2=$ 3.682 = 
 $6-4c^{1}_{7}
-6c^{2}_{7}
$

\vskip 0.7ex
\hangindent=3em \hangafter=1
$T = ( 0,
\frac{3}{7},
\frac{25}{28},
\frac{5}{28},
\frac{1}{7},
\frac{3}{4} )
$,

\vskip 0.7ex
\hangindent=3em \hangafter=1
$S$ = ($ 1$,
$ -c^{2}_{7}
$,
$ \xi_{7}^{2,3}$,
$ c^{2}_{7}
$,
$ -\xi_{7}^{2,3}$,
$ -1$;\ \ 
$ \xi_{7}^{2,3}$,
$ -1$,
$ -\xi_{7}^{2,3}$,
$ 1$,
$ c^{2}_{7}
$;\ \ 
$ -c^{2}_{7}
$,
$ -1$,
$ -c^{2}_{7}
$,
$ \xi_{7}^{2,3}$;\ \ 
$ -\xi_{7}^{2,3}$,
$ -1$,
$ c^{2}_{7}
$;\ \ 
$ c^{2}_{7}
$,
$ \xi_{7}^{2,3}$;\ \ 
$ -1$)

Not pseudo-unitary. 

\vskip 1ex 
\color{grey}

\noindent105. ind = $(6;56
)_{2}^{25}$:\ \ 
$d_i$ = ($1.0$,
$0.445$,
$0.801$,
$-0.445$,
$-0.801$,
$-1.0$) 

\vskip 0.7ex
\hangindent=3em \hangafter=1
$D^2=$ 3.682 = 
 $6-4c^{1}_{7}
-6c^{2}_{7}
$

\vskip 0.7ex
\hangindent=3em \hangafter=1
$T = ( 0,
\frac{4}{7},
\frac{3}{28},
\frac{23}{28},
\frac{6}{7},
\frac{1}{4} )
$,

\vskip 0.7ex
\hangindent=3em \hangafter=1
$S$ = ($ 1$,
$ -c^{2}_{7}
$,
$ \xi_{7}^{2,3}$,
$ c^{2}_{7}
$,
$ -\xi_{7}^{2,3}$,
$ -1$;\ \ 
$ \xi_{7}^{2,3}$,
$ -1$,
$ -\xi_{7}^{2,3}$,
$ 1$,
$ c^{2}_{7}
$;\ \ 
$ -c^{2}_{7}
$,
$ -1$,
$ -c^{2}_{7}
$,
$ \xi_{7}^{2,3}$;\ \ 
$ -\xi_{7}^{2,3}$,
$ -1$,
$ c^{2}_{7}
$;\ \ 
$ c^{2}_{7}
$,
$ \xi_{7}^{2,3}$;\ \ 
$ -1$)

Not pseudo-unitary. 

\vskip 1ex 
\color{grey}

\noindent106. ind = $(6;56
)_{2}^{19}$:\ \ 
$d_i$ = ($1.0$,
$0.554$,
$1.246$,
$-0.554$,
$-1.0$,
$-1.246$) 

\vskip 0.7ex
\hangindent=3em \hangafter=1
$D^2=$ 5.725 = 
 $10-2c^{1}_{7}
+4c^{2}_{7}
$

\vskip 0.7ex
\hangindent=3em \hangafter=1
$T = ( 0,
\frac{4}{7},
\frac{13}{28},
\frac{9}{28},
\frac{3}{4},
\frac{5}{7} )
$,

\vskip 0.7ex
\hangindent=3em \hangafter=1
$S$ = ($ 1$,
$ \xi_{7}^{1,2}$,
$ c^{1}_{7}
$,
$ -\xi_{7}^{1,2}$,
$ -1$,
$ -c^{1}_{7}
$;\ \ 
$ c^{1}_{7}
$,
$ -1$,
$ -c^{1}_{7}
$,
$ -\xi_{7}^{1,2}$,
$ 1$;\ \ 
$ \xi_{7}^{1,2}$,
$ -1$,
$ c^{1}_{7}
$,
$ \xi_{7}^{1,2}$;\ \ 
$ -c^{1}_{7}
$,
$ -\xi_{7}^{1,2}$,
$ -1$;\ \ 
$ -1$,
$ c^{1}_{7}
$;\ \ 
$ -\xi_{7}^{1,2}$)

Not pseudo-unitary. 

\vskip 1ex 
\color{grey}

\noindent107. ind = $(6;56
)_{2}^{11}$:\ \ 
$d_i$ = ($1.0$,
$0.445$,
$0.801$,
$-0.445$,
$-0.801$,
$-1.0$) 

\vskip 0.7ex
\hangindent=3em \hangafter=1
$D^2=$ 3.682 = 
 $6-4c^{1}_{7}
-6c^{2}_{7}
$

\vskip 0.7ex
\hangindent=3em \hangafter=1
$T = ( 0,
\frac{4}{7},
\frac{17}{28},
\frac{9}{28},
\frac{6}{7},
\frac{3}{4} )
$,

\vskip 0.7ex
\hangindent=3em \hangafter=1
$S$ = ($ 1$,
$ -c^{2}_{7}
$,
$ \xi_{7}^{2,3}$,
$ c^{2}_{7}
$,
$ -\xi_{7}^{2,3}$,
$ -1$;\ \ 
$ \xi_{7}^{2,3}$,
$ -1$,
$ -\xi_{7}^{2,3}$,
$ 1$,
$ c^{2}_{7}
$;\ \ 
$ -c^{2}_{7}
$,
$ -1$,
$ -c^{2}_{7}
$,
$ \xi_{7}^{2,3}$;\ \ 
$ -\xi_{7}^{2,3}$,
$ -1$,
$ c^{2}_{7}
$;\ \ 
$ c^{2}_{7}
$,
$ \xi_{7}^{2,3}$;\ \ 
$ -1$)

Not pseudo-unitary. 

\vskip 1ex 
\color{grey}

\noindent108. ind = $(6;56
)_{2}^{5}$:\ \ 
$d_i$ = ($1.0$,
$0.554$,
$1.246$,
$-0.554$,
$-1.0$,
$-1.246$) 

\vskip 0.7ex
\hangindent=3em \hangafter=1
$D^2=$ 5.725 = 
 $10-2c^{1}_{7}
+4c^{2}_{7}
$

\vskip 0.7ex
\hangindent=3em \hangafter=1
$T = ( 0,
\frac{4}{7},
\frac{27}{28},
\frac{23}{28},
\frac{1}{4},
\frac{5}{7} )
$,

\vskip 0.7ex
\hangindent=3em \hangafter=1
$S$ = ($ 1$,
$ \xi_{7}^{1,2}$,
$ c^{1}_{7}
$,
$ -\xi_{7}^{1,2}$,
$ -1$,
$ -c^{1}_{7}
$;\ \ 
$ c^{1}_{7}
$,
$ -1$,
$ -c^{1}_{7}
$,
$ -\xi_{7}^{1,2}$,
$ 1$;\ \ 
$ \xi_{7}^{1,2}$,
$ -1$,
$ c^{1}_{7}
$,
$ \xi_{7}^{1,2}$;\ \ 
$ -c^{1}_{7}
$,
$ -\xi_{7}^{1,2}$,
$ -1$;\ \ 
$ -1$,
$ c^{1}_{7}
$;\ \ 
$ -\xi_{7}^{1,2}$)

Not pseudo-unitary. 

\vskip 1ex 
\color{grey}

\noindent109. ind = $(6;56
)_{2}^{13}$:\ \ 
$d_i$ = ($1.0$,
$1.801$,
$2.246$,
$-1.0$,
$-1.801$,
$-2.246$) 

\vskip 0.7ex
\hangindent=3em \hangafter=1
$D^2=$ 18.591 = 
 $12+6c^{1}_{7}
+2c^{2}_{7}
$

\vskip 0.7ex
\hangindent=3em \hangafter=1
$T = ( 0,
\frac{6}{7},
\frac{2}{7},
\frac{1}{4},
\frac{3}{28},
\frac{15}{28} )
$,

\vskip 0.7ex
\hangindent=3em \hangafter=1
$S$ = ($ 1$,
$ \xi_{7}^{2}$,
$ \xi_{7}^{3}$,
$ -1$,
$ -\xi_{7}^{2}$,
$ -\xi_{7}^{3}$;\ \ 
$ -\xi_{7}^{3}$,
$ 1$,
$ -\xi_{7}^{2}$,
$ \xi_{7}^{3}$,
$ -1$;\ \ 
$ -\xi_{7}^{2}$,
$ -\xi_{7}^{3}$,
$ -1$,
$ \xi_{7}^{2}$;\ \ 
$ -1$,
$ -\xi_{7}^{2}$,
$ -\xi_{7}^{3}$;\ \ 
$ \xi_{7}^{3}$,
$ -1$;\ \ 
$ \xi_{7}^{2}$)

Pseudo-unitary $\sim$  
$(6;56
)_{1}^{13}$

\vskip 1ex 
\color{grey}

\noindent110. ind = $(6;56
)_{2}^{27}$:\ \ 
$d_i$ = ($1.0$,
$1.801$,
$2.246$,
$-1.0$,
$-1.801$,
$-2.246$) 

\vskip 0.7ex
\hangindent=3em \hangafter=1
$D^2=$ 18.591 = 
 $12+6c^{1}_{7}
+2c^{2}_{7}
$

\vskip 0.7ex
\hangindent=3em \hangafter=1
$T = ( 0,
\frac{6}{7},
\frac{2}{7},
\frac{3}{4},
\frac{17}{28},
\frac{1}{28} )
$,

\vskip 0.7ex
\hangindent=3em \hangafter=1
$S$ = ($ 1$,
$ \xi_{7}^{2}$,
$ \xi_{7}^{3}$,
$ -1$,
$ -\xi_{7}^{2}$,
$ -\xi_{7}^{3}$;\ \ 
$ -\xi_{7}^{3}$,
$ 1$,
$ -\xi_{7}^{2}$,
$ \xi_{7}^{3}$,
$ -1$;\ \ 
$ -\xi_{7}^{2}$,
$ -\xi_{7}^{3}$,
$ -1$,
$ \xi_{7}^{2}$;\ \ 
$ -1$,
$ -\xi_{7}^{2}$,
$ -\xi_{7}^{3}$;\ \ 
$ \xi_{7}^{3}$,
$ -1$;\ \ 
$ \xi_{7}^{2}$)

Pseudo-unitary $\sim$  
$(6;56
)_{1}^{27}$

\vskip 1ex 

 \color{black} \vskip 2ex

\noindent111. ind = $(6;80
)_{1}^{1}$:\ \ 
$d_i$ = ($1.0$,
$1.0$,
$1.414$,
$1.618$,
$1.618$,
$2.288$) 

\vskip 0.7ex
\hangindent=3em \hangafter=1
$D^2=$ 14.472 = 
 $10+2\sqrt{5}$

\vskip 0.7ex
\hangindent=3em \hangafter=1
$T = ( 0,
\frac{1}{2},
\frac{1}{16},
\frac{2}{5},
\frac{9}{10},
\frac{37}{80} )
$,

\vskip 0.7ex
\hangindent=3em \hangafter=1
$S$ = ($ 1$,
$ 1$,
$ \sqrt{2}$,
$ \frac{1+\sqrt{5}}{2}$,
$ \frac{1+\sqrt{5}}{2}$,
$ c^{3}_{40}
+c^{5}_{40}
-c^{7}_{40}
$;\ \ 
$ 1$,
$ -\sqrt{2}$,
$ \frac{1+\sqrt{5}}{2}$,
$ \frac{1+\sqrt{5}}{2}$,
$ -c^{3}_{40}
-c^{5}_{40}
+c^{7}_{40}
$;\ \ 
$0$,
$ c^{3}_{40}
+c^{5}_{40}
-c^{7}_{40}
$,
$ -c^{3}_{40}
-c^{5}_{40}
+c^{7}_{40}
$,
$0$;\ \ 
$ -1$,
$ -1$,
$ -\sqrt{2}$;\ \ 
$ -1$,
$ \sqrt{2}$;\ \ 
$0$)

\vskip 1ex 
\color{grey}

\noindent112. ind = $(6;80
)_{1}^{49}$:\ \ 
$d_i$ = ($1.0$,
$1.0$,
$1.414$,
$1.618$,
$1.618$,
$2.288$) 

\vskip 0.7ex
\hangindent=3em \hangafter=1
$D^2=$ 14.472 = 
 $10+2\sqrt{5}$

\vskip 0.7ex
\hangindent=3em \hangafter=1
$T = ( 0,
\frac{1}{2},
\frac{1}{16},
\frac{3}{5},
\frac{1}{10},
\frac{53}{80} )
$,

\vskip 0.7ex
\hangindent=3em \hangafter=1
$S$ = ($ 1$,
$ 1$,
$ \sqrt{2}$,
$ \frac{1+\sqrt{5}}{2}$,
$ \frac{1+\sqrt{5}}{2}$,
$ c^{3}_{40}
+c^{5}_{40}
-c^{7}_{40}
$;\ \ 
$ 1$,
$ -\sqrt{2}$,
$ \frac{1+\sqrt{5}}{2}$,
$ \frac{1+\sqrt{5}}{2}$,
$ -c^{3}_{40}
-c^{5}_{40}
+c^{7}_{40}
$;\ \ 
$0$,
$ c^{3}_{40}
+c^{5}_{40}
-c^{7}_{40}
$,
$ -c^{3}_{40}
-c^{5}_{40}
+c^{7}_{40}
$,
$0$;\ \ 
$ -1$,
$ -1$,
$ -\sqrt{2}$;\ \ 
$ -1$,
$ \sqrt{2}$;\ \ 
$0$)

\vskip 1ex 
\color{grey}

\noindent113. ind = $(6;80
)_{1}^{71}$:\ \ 
$d_i$ = ($1.0$,
$1.0$,
$1.414$,
$1.618$,
$1.618$,
$2.288$) 

\vskip 0.7ex
\hangindent=3em \hangafter=1
$D^2=$ 14.472 = 
 $10+2\sqrt{5}$

\vskip 0.7ex
\hangindent=3em \hangafter=1
$T = ( 0,
\frac{1}{2},
\frac{7}{16},
\frac{2}{5},
\frac{9}{10},
\frac{67}{80} )
$,

\vskip 0.7ex
\hangindent=3em \hangafter=1
$S$ = ($ 1$,
$ 1$,
$ \sqrt{2}$,
$ \frac{1+\sqrt{5}}{2}$,
$ \frac{1+\sqrt{5}}{2}$,
$ c^{3}_{40}
+c^{5}_{40}
-c^{7}_{40}
$;\ \ 
$ 1$,
$ -\sqrt{2}$,
$ \frac{1+\sqrt{5}}{2}$,
$ \frac{1+\sqrt{5}}{2}$,
$ -c^{3}_{40}
-c^{5}_{40}
+c^{7}_{40}
$;\ \ 
$0$,
$ c^{3}_{40}
+c^{5}_{40}
-c^{7}_{40}
$,
$ -c^{3}_{40}
-c^{5}_{40}
+c^{7}_{40}
$,
$0$;\ \ 
$ -1$,
$ -1$,
$ -\sqrt{2}$;\ \ 
$ -1$,
$ \sqrt{2}$;\ \ 
$0$)

\vskip 1ex 
\color{grey}

\noindent114. ind = $(6;80
)_{1}^{39}$:\ \ 
$d_i$ = ($1.0$,
$1.0$,
$1.414$,
$1.618$,
$1.618$,
$2.288$) 

\vskip 0.7ex
\hangindent=3em \hangafter=1
$D^2=$ 14.472 = 
 $10+2\sqrt{5}$

\vskip 0.7ex
\hangindent=3em \hangafter=1
$T = ( 0,
\frac{1}{2},
\frac{7}{16},
\frac{3}{5},
\frac{1}{10},
\frac{3}{80} )
$,

\vskip 0.7ex
\hangindent=3em \hangafter=1
$S$ = ($ 1$,
$ 1$,
$ \sqrt{2}$,
$ \frac{1+\sqrt{5}}{2}$,
$ \frac{1+\sqrt{5}}{2}$,
$ c^{3}_{40}
+c^{5}_{40}
-c^{7}_{40}
$;\ \ 
$ 1$,
$ -\sqrt{2}$,
$ \frac{1+\sqrt{5}}{2}$,
$ \frac{1+\sqrt{5}}{2}$,
$ -c^{3}_{40}
-c^{5}_{40}
+c^{7}_{40}
$;\ \ 
$0$,
$ c^{3}_{40}
+c^{5}_{40}
-c^{7}_{40}
$,
$ -c^{3}_{40}
-c^{5}_{40}
+c^{7}_{40}
$,
$0$;\ \ 
$ -1$,
$ -1$,
$ -\sqrt{2}$;\ \ 
$ -1$,
$ \sqrt{2}$;\ \ 
$0$)

\vskip 1ex 
\color{grey}

\noindent115. ind = $(6;80
)_{1}^{41}$:\ \ 
$d_i$ = ($1.0$,
$1.0$,
$1.414$,
$1.618$,
$1.618$,
$2.288$) 

\vskip 0.7ex
\hangindent=3em \hangafter=1
$D^2=$ 14.472 = 
 $10+2\sqrt{5}$

\vskip 0.7ex
\hangindent=3em \hangafter=1
$T = ( 0,
\frac{1}{2},
\frac{9}{16},
\frac{2}{5},
\frac{9}{10},
\frac{77}{80} )
$,

\vskip 0.7ex
\hangindent=3em \hangafter=1
$S$ = ($ 1$,
$ 1$,
$ \sqrt{2}$,
$ \frac{1+\sqrt{5}}{2}$,
$ \frac{1+\sqrt{5}}{2}$,
$ c^{3}_{40}
+c^{5}_{40}
-c^{7}_{40}
$;\ \ 
$ 1$,
$ -\sqrt{2}$,
$ \frac{1+\sqrt{5}}{2}$,
$ \frac{1+\sqrt{5}}{2}$,
$ -c^{3}_{40}
-c^{5}_{40}
+c^{7}_{40}
$;\ \ 
$0$,
$ c^{3}_{40}
+c^{5}_{40}
-c^{7}_{40}
$,
$ -c^{3}_{40}
-c^{5}_{40}
+c^{7}_{40}
$,
$0$;\ \ 
$ -1$,
$ -1$,
$ -\sqrt{2}$;\ \ 
$ -1$,
$ \sqrt{2}$;\ \ 
$0$)

\vskip 1ex 
\color{grey}

\noindent116. ind = $(6;80
)_{1}^{9}$:\ \ 
$d_i$ = ($1.0$,
$1.0$,
$1.414$,
$1.618$,
$1.618$,
$2.288$) 

\vskip 0.7ex
\hangindent=3em \hangafter=1
$D^2=$ 14.472 = 
 $10+2\sqrt{5}$

\vskip 0.7ex
\hangindent=3em \hangafter=1
$T = ( 0,
\frac{1}{2},
\frac{9}{16},
\frac{3}{5},
\frac{1}{10},
\frac{13}{80} )
$,

\vskip 0.7ex
\hangindent=3em \hangafter=1
$S$ = ($ 1$,
$ 1$,
$ \sqrt{2}$,
$ \frac{1+\sqrt{5}}{2}$,
$ \frac{1+\sqrt{5}}{2}$,
$ c^{3}_{40}
+c^{5}_{40}
-c^{7}_{40}
$;\ \ 
$ 1$,
$ -\sqrt{2}$,
$ \frac{1+\sqrt{5}}{2}$,
$ \frac{1+\sqrt{5}}{2}$,
$ -c^{3}_{40}
-c^{5}_{40}
+c^{7}_{40}
$;\ \ 
$0$,
$ c^{3}_{40}
+c^{5}_{40}
-c^{7}_{40}
$,
$ -c^{3}_{40}
-c^{5}_{40}
+c^{7}_{40}
$,
$0$;\ \ 
$ -1$,
$ -1$,
$ -\sqrt{2}$;\ \ 
$ -1$,
$ \sqrt{2}$;\ \ 
$0$)

\vskip 1ex 
\color{grey}

\noindent117. ind = $(6;80
)_{1}^{31}$:\ \ 
$d_i$ = ($1.0$,
$1.0$,
$1.414$,
$1.618$,
$1.618$,
$2.288$) 

\vskip 0.7ex
\hangindent=3em \hangafter=1
$D^2=$ 14.472 = 
 $10+2\sqrt{5}$

\vskip 0.7ex
\hangindent=3em \hangafter=1
$T = ( 0,
\frac{1}{2},
\frac{15}{16},
\frac{2}{5},
\frac{9}{10},
\frac{27}{80} )
$,

\vskip 0.7ex
\hangindent=3em \hangafter=1
$S$ = ($ 1$,
$ 1$,
$ \sqrt{2}$,
$ \frac{1+\sqrt{5}}{2}$,
$ \frac{1+\sqrt{5}}{2}$,
$ c^{3}_{40}
+c^{5}_{40}
-c^{7}_{40}
$;\ \ 
$ 1$,
$ -\sqrt{2}$,
$ \frac{1+\sqrt{5}}{2}$,
$ \frac{1+\sqrt{5}}{2}$,
$ -c^{3}_{40}
-c^{5}_{40}
+c^{7}_{40}
$;\ \ 
$0$,
$ c^{3}_{40}
+c^{5}_{40}
-c^{7}_{40}
$,
$ -c^{3}_{40}
-c^{5}_{40}
+c^{7}_{40}
$,
$0$;\ \ 
$ -1$,
$ -1$,
$ -\sqrt{2}$;\ \ 
$ -1$,
$ \sqrt{2}$;\ \ 
$0$)

\vskip 1ex 
\color{grey}

\noindent118. ind = $(6;80
)_{1}^{79}$:\ \ 
$d_i$ = ($1.0$,
$1.0$,
$1.414$,
$1.618$,
$1.618$,
$2.288$) 

\vskip 0.7ex
\hangindent=3em \hangafter=1
$D^2=$ 14.472 = 
 $10+2\sqrt{5}$

\vskip 0.7ex
\hangindent=3em \hangafter=1
$T = ( 0,
\frac{1}{2},
\frac{15}{16},
\frac{3}{5},
\frac{1}{10},
\frac{43}{80} )
$,

\vskip 0.7ex
\hangindent=3em \hangafter=1
$S$ = ($ 1$,
$ 1$,
$ \sqrt{2}$,
$ \frac{1+\sqrt{5}}{2}$,
$ \frac{1+\sqrt{5}}{2}$,
$ c^{3}_{40}
+c^{5}_{40}
-c^{7}_{40}
$;\ \ 
$ 1$,
$ -\sqrt{2}$,
$ \frac{1+\sqrt{5}}{2}$,
$ \frac{1+\sqrt{5}}{2}$,
$ -c^{3}_{40}
-c^{5}_{40}
+c^{7}_{40}
$;\ \ 
$0$,
$ c^{3}_{40}
+c^{5}_{40}
-c^{7}_{40}
$,
$ -c^{3}_{40}
-c^{5}_{40}
+c^{7}_{40}
$,
$0$;\ \ 
$ -1$,
$ -1$,
$ -\sqrt{2}$;\ \ 
$ -1$,
$ \sqrt{2}$;\ \ 
$0$)

\vskip 1ex 
\color{grey}

\noindent119. ind = $(6;80
)_{1}^{13}$:\ \ 
$d_i$ = ($1.0$,
$0.874$,
$1.0$,
$-0.618$,
$-0.618$,
$-1.414$) 

\vskip 0.7ex
\hangindent=3em \hangafter=1
$D^2=$ 5.527 = 
 $10-2\sqrt{5}$

\vskip 0.7ex
\hangindent=3em \hangafter=1
$T = ( 0,
\frac{1}{80},
\frac{1}{2},
\frac{1}{5},
\frac{7}{10},
\frac{13}{16} )
$,

\vskip 0.7ex
\hangindent=3em \hangafter=1
$S$ = ($ 1$,
$ c^{3}_{40}
-c^{7}_{40}
$,
$ 1$,
$ \frac{1-\sqrt{5}}{2}$,
$ \frac{1-\sqrt{5}}{2}$,
$ -\sqrt{2}$;\ \ 
$0$,
$ -c^{3}_{40}
+c^{7}_{40}
$,
$ \sqrt{2}$,
$ -\sqrt{2}$,
$0$;\ \ 
$ 1$,
$ \frac{1-\sqrt{5}}{2}$,
$ \frac{1-\sqrt{5}}{2}$,
$ \sqrt{2}$;\ \ 
$ -1$,
$ -1$,
$ c^{3}_{40}
-c^{7}_{40}
$;\ \ 
$ -1$,
$ -c^{3}_{40}
+c^{7}_{40}
$;\ \ 
$0$)

Not pseudo-unitary. 

\vskip 1ex 
\color{grey}

\noindent120. ind = $(6;80
)_{1}^{37}$:\ \ 
$d_i$ = ($1.0$,
$0.874$,
$1.0$,
$-0.618$,
$-0.618$,
$-1.414$) 

\vskip 0.7ex
\hangindent=3em \hangafter=1
$D^2=$ 5.527 = 
 $10-2\sqrt{5}$

\vskip 0.7ex
\hangindent=3em \hangafter=1
$T = ( 0,
\frac{9}{80},
\frac{1}{2},
\frac{4}{5},
\frac{3}{10},
\frac{5}{16} )
$,

\vskip 0.7ex
\hangindent=3em \hangafter=1
$S$ = ($ 1$,
$ c^{3}_{40}
-c^{7}_{40}
$,
$ 1$,
$ \frac{1-\sqrt{5}}{2}$,
$ \frac{1-\sqrt{5}}{2}$,
$ -\sqrt{2}$;\ \ 
$0$,
$ -c^{3}_{40}
+c^{7}_{40}
$,
$ \sqrt{2}$,
$ -\sqrt{2}$,
$0$;\ \ 
$ 1$,
$ \frac{1-\sqrt{5}}{2}$,
$ \frac{1-\sqrt{5}}{2}$,
$ \sqrt{2}$;\ \ 
$ -1$,
$ -1$,
$ c^{3}_{40}
-c^{7}_{40}
$;\ \ 
$ -1$,
$ -c^{3}_{40}
+c^{7}_{40}
$;\ \ 
$0$)

Not pseudo-unitary. 

\vskip 1ex 
\color{grey}

\noindent121. ind = $(6;80
)_{1}^{3}$:\ \ 
$d_i$ = ($1.0$,
$0.874$,
$1.0$,
$-0.618$,
$-0.618$,
$-1.414$) 

\vskip 0.7ex
\hangindent=3em \hangafter=1
$D^2=$ 5.527 = 
 $10-2\sqrt{5}$

\vskip 0.7ex
\hangindent=3em \hangafter=1
$T = ( 0,
\frac{31}{80},
\frac{1}{2},
\frac{1}{5},
\frac{7}{10},
\frac{3}{16} )
$,

\vskip 0.7ex
\hangindent=3em \hangafter=1
$S$ = ($ 1$,
$ c^{3}_{40}
-c^{7}_{40}
$,
$ 1$,
$ \frac{1-\sqrt{5}}{2}$,
$ \frac{1-\sqrt{5}}{2}$,
$ -\sqrt{2}$;\ \ 
$0$,
$ -c^{3}_{40}
+c^{7}_{40}
$,
$ \sqrt{2}$,
$ -\sqrt{2}$,
$0$;\ \ 
$ 1$,
$ \frac{1-\sqrt{5}}{2}$,
$ \frac{1-\sqrt{5}}{2}$,
$ \sqrt{2}$;\ \ 
$ -1$,
$ -1$,
$ c^{3}_{40}
-c^{7}_{40}
$;\ \ 
$ -1$,
$ -c^{3}_{40}
+c^{7}_{40}
$;\ \ 
$0$)

Not pseudo-unitary. 

\vskip 1ex 
\color{grey}

\noindent122. ind = $(6;80
)_{1}^{27}$:\ \ 
$d_i$ = ($1.0$,
$0.874$,
$1.0$,
$-0.618$,
$-0.618$,
$-1.414$) 

\vskip 0.7ex
\hangindent=3em \hangafter=1
$D^2=$ 5.527 = 
 $10-2\sqrt{5}$

\vskip 0.7ex
\hangindent=3em \hangafter=1
$T = ( 0,
\frac{39}{80},
\frac{1}{2},
\frac{4}{5},
\frac{3}{10},
\frac{11}{16} )
$,

\vskip 0.7ex
\hangindent=3em \hangafter=1
$S$ = ($ 1$,
$ c^{3}_{40}
-c^{7}_{40}
$,
$ 1$,
$ \frac{1-\sqrt{5}}{2}$,
$ \frac{1-\sqrt{5}}{2}$,
$ -\sqrt{2}$;\ \ 
$0$,
$ -c^{3}_{40}
+c^{7}_{40}
$,
$ \sqrt{2}$,
$ -\sqrt{2}$,
$0$;\ \ 
$ 1$,
$ \frac{1-\sqrt{5}}{2}$,
$ \frac{1-\sqrt{5}}{2}$,
$ \sqrt{2}$;\ \ 
$ -1$,
$ -1$,
$ c^{3}_{40}
-c^{7}_{40}
$;\ \ 
$ -1$,
$ -c^{3}_{40}
+c^{7}_{40}
$;\ \ 
$0$)

Not pseudo-unitary. 

\vskip 1ex 
\color{grey}

\noindent123. ind = $(6;80
)_{1}^{33}$:\ \ 
$d_i$ = ($1.0$,
$1.0$,
$1.414$,
$-0.618$,
$-0.618$,
$-0.874$) 

\vskip 0.7ex
\hangindent=3em \hangafter=1
$D^2=$ 5.527 = 
 $10-2\sqrt{5}$

\vskip 0.7ex
\hangindent=3em \hangafter=1
$T = ( 0,
\frac{1}{2},
\frac{1}{16},
\frac{1}{5},
\frac{7}{10},
\frac{21}{80} )
$,

\vskip 0.7ex
\hangindent=3em \hangafter=1
$S$ = ($ 1$,
$ 1$,
$ \sqrt{2}$,
$ \frac{1-\sqrt{5}}{2}$,
$ \frac{1-\sqrt{5}}{2}$,
$ -c^{3}_{40}
+c^{7}_{40}
$;\ \ 
$ 1$,
$ -\sqrt{2}$,
$ \frac{1-\sqrt{5}}{2}$,
$ \frac{1-\sqrt{5}}{2}$,
$ c^{3}_{40}
-c^{7}_{40}
$;\ \ 
$0$,
$ -c^{3}_{40}
+c^{7}_{40}
$,
$ c^{3}_{40}
-c^{7}_{40}
$,
$0$;\ \ 
$ -1$,
$ -1$,
$ -\sqrt{2}$;\ \ 
$ -1$,
$ \sqrt{2}$;\ \ 
$0$)

Not pseudo-unitary. 

\vskip 1ex 
\color{grey}

\noindent124. ind = $(6;80
)_{1}^{17}$:\ \ 
$d_i$ = ($1.0$,
$1.0$,
$1.414$,
$-0.618$,
$-0.618$,
$-0.874$) 

\vskip 0.7ex
\hangindent=3em \hangafter=1
$D^2=$ 5.527 = 
 $10-2\sqrt{5}$

\vskip 0.7ex
\hangindent=3em \hangafter=1
$T = ( 0,
\frac{1}{2},
\frac{1}{16},
\frac{4}{5},
\frac{3}{10},
\frac{69}{80} )
$,

\vskip 0.7ex
\hangindent=3em \hangafter=1
$S$ = ($ 1$,
$ 1$,
$ \sqrt{2}$,
$ \frac{1-\sqrt{5}}{2}$,
$ \frac{1-\sqrt{5}}{2}$,
$ -c^{3}_{40}
+c^{7}_{40}
$;\ \ 
$ 1$,
$ -\sqrt{2}$,
$ \frac{1-\sqrt{5}}{2}$,
$ \frac{1-\sqrt{5}}{2}$,
$ c^{3}_{40}
-c^{7}_{40}
$;\ \ 
$0$,
$ -c^{3}_{40}
+c^{7}_{40}
$,
$ c^{3}_{40}
-c^{7}_{40}
$,
$0$;\ \ 
$ -1$,
$ -1$,
$ -\sqrt{2}$;\ \ 
$ -1$,
$ \sqrt{2}$;\ \ 
$0$)

Not pseudo-unitary. 

\vskip 1ex 
\color{grey}

\noindent125. ind = $(6;80
)_{1}^{51}$:\ \ 
$d_i$ = ($1.0$,
$1.0$,
$1.618$,
$1.618$,
$-1.414$,
$-2.288$) 

\vskip 0.7ex
\hangindent=3em \hangafter=1
$D^2=$ 14.472 = 
 $10+2\sqrt{5}$

\vskip 0.7ex
\hangindent=3em \hangafter=1
$T = ( 0,
\frac{1}{2},
\frac{2}{5},
\frac{9}{10},
\frac{3}{16},
\frac{47}{80} )
$,

\vskip 0.7ex
\hangindent=3em \hangafter=1
$S$ = ($ 1$,
$ 1$,
$ \frac{1+\sqrt{5}}{2}$,
$ \frac{1+\sqrt{5}}{2}$,
$ -\sqrt{2}$,
$ -c^{3}_{40}
-c^{5}_{40}
+c^{7}_{40}
$;\ \ 
$ 1$,
$ \frac{1+\sqrt{5}}{2}$,
$ \frac{1+\sqrt{5}}{2}$,
$ \sqrt{2}$,
$ c^{3}_{40}
+c^{5}_{40}
-c^{7}_{40}
$;\ \ 
$ -1$,
$ -1$,
$ -c^{3}_{40}
-c^{5}_{40}
+c^{7}_{40}
$,
$ \sqrt{2}$;\ \ 
$ -1$,
$ c^{3}_{40}
+c^{5}_{40}
-c^{7}_{40}
$,
$ -\sqrt{2}$;\ \ 
$0$,
$0$;\ \ 
$0$)

Pseudo-unitary $\sim$  
$(6;80
)_{2}^{1}$

\vskip 1ex 
\color{grey}

\noindent126. ind = $(6;80
)_{1}^{21}$:\ \ 
$d_i$ = ($1.0$,
$1.0$,
$1.618$,
$1.618$,
$-1.414$,
$-2.288$) 

\vskip 0.7ex
\hangindent=3em \hangafter=1
$D^2=$ 14.472 = 
 $10+2\sqrt{5}$

\vskip 0.7ex
\hangindent=3em \hangafter=1
$T = ( 0,
\frac{1}{2},
\frac{2}{5},
\frac{9}{10},
\frac{5}{16},
\frac{57}{80} )
$,

\vskip 0.7ex
\hangindent=3em \hangafter=1
$S$ = ($ 1$,
$ 1$,
$ \frac{1+\sqrt{5}}{2}$,
$ \frac{1+\sqrt{5}}{2}$,
$ -\sqrt{2}$,
$ -c^{3}_{40}
-c^{5}_{40}
+c^{7}_{40}
$;\ \ 
$ 1$,
$ \frac{1+\sqrt{5}}{2}$,
$ \frac{1+\sqrt{5}}{2}$,
$ \sqrt{2}$,
$ c^{3}_{40}
+c^{5}_{40}
-c^{7}_{40}
$;\ \ 
$ -1$,
$ -1$,
$ -c^{3}_{40}
-c^{5}_{40}
+c^{7}_{40}
$,
$ \sqrt{2}$;\ \ 
$ -1$,
$ c^{3}_{40}
+c^{5}_{40}
-c^{7}_{40}
$,
$ -\sqrt{2}$;\ \ 
$0$,
$0$;\ \ 
$0$)

Pseudo-unitary $\sim$  
$(6;80
)_{2}^{71}$

\vskip 1ex 
\color{grey}

\noindent127. ind = $(6;80
)_{1}^{11}$:\ \ 
$d_i$ = ($1.0$,
$1.0$,
$1.618$,
$1.618$,
$-1.414$,
$-2.288$) 

\vskip 0.7ex
\hangindent=3em \hangafter=1
$D^2=$ 14.472 = 
 $10+2\sqrt{5}$

\vskip 0.7ex
\hangindent=3em \hangafter=1
$T = ( 0,
\frac{1}{2},
\frac{2}{5},
\frac{9}{10},
\frac{11}{16},
\frac{7}{80} )
$,

\vskip 0.7ex
\hangindent=3em \hangafter=1
$S$ = ($ 1$,
$ 1$,
$ \frac{1+\sqrt{5}}{2}$,
$ \frac{1+\sqrt{5}}{2}$,
$ -\sqrt{2}$,
$ -c^{3}_{40}
-c^{5}_{40}
+c^{7}_{40}
$;\ \ 
$ 1$,
$ \frac{1+\sqrt{5}}{2}$,
$ \frac{1+\sqrt{5}}{2}$,
$ \sqrt{2}$,
$ c^{3}_{40}
+c^{5}_{40}
-c^{7}_{40}
$;\ \ 
$ -1$,
$ -1$,
$ -c^{3}_{40}
-c^{5}_{40}
+c^{7}_{40}
$,
$ \sqrt{2}$;\ \ 
$ -1$,
$ c^{3}_{40}
+c^{5}_{40}
-c^{7}_{40}
$,
$ -\sqrt{2}$;\ \ 
$0$,
$0$;\ \ 
$0$)

Pseudo-unitary $\sim$  
$(6;80
)_{2}^{41}$

\vskip 1ex 
\color{grey}

\noindent128. ind = $(6;80
)_{1}^{61}$:\ \ 
$d_i$ = ($1.0$,
$1.0$,
$1.618$,
$1.618$,
$-1.414$,
$-2.288$) 

\vskip 0.7ex
\hangindent=3em \hangafter=1
$D^2=$ 14.472 = 
 $10+2\sqrt{5}$

\vskip 0.7ex
\hangindent=3em \hangafter=1
$T = ( 0,
\frac{1}{2},
\frac{2}{5},
\frac{9}{10},
\frac{13}{16},
\frac{17}{80} )
$,

\vskip 0.7ex
\hangindent=3em \hangafter=1
$S$ = ($ 1$,
$ 1$,
$ \frac{1+\sqrt{5}}{2}$,
$ \frac{1+\sqrt{5}}{2}$,
$ -\sqrt{2}$,
$ -c^{3}_{40}
-c^{5}_{40}
+c^{7}_{40}
$;\ \ 
$ 1$,
$ \frac{1+\sqrt{5}}{2}$,
$ \frac{1+\sqrt{5}}{2}$,
$ \sqrt{2}$,
$ c^{3}_{40}
+c^{5}_{40}
-c^{7}_{40}
$;\ \ 
$ -1$,
$ -1$,
$ -c^{3}_{40}
-c^{5}_{40}
+c^{7}_{40}
$,
$ \sqrt{2}$;\ \ 
$ -1$,
$ c^{3}_{40}
+c^{5}_{40}
-c^{7}_{40}
$,
$ -\sqrt{2}$;\ \ 
$0$,
$0$;\ \ 
$0$)

Pseudo-unitary $\sim$  
$(6;80
)_{2}^{31}$

\vskip 1ex 
\color{grey}

\noindent129. ind = $(6;80
)_{1}^{23}$:\ \ 
$d_i$ = ($1.0$,
$1.0$,
$1.414$,
$-0.618$,
$-0.618$,
$-0.874$) 

\vskip 0.7ex
\hangindent=3em \hangafter=1
$D^2=$ 5.527 = 
 $10-2\sqrt{5}$

\vskip 0.7ex
\hangindent=3em \hangafter=1
$T = ( 0,
\frac{1}{2},
\frac{7}{16},
\frac{1}{5},
\frac{7}{10},
\frac{51}{80} )
$,

\vskip 0.7ex
\hangindent=3em \hangafter=1
$S$ = ($ 1$,
$ 1$,
$ \sqrt{2}$,
$ \frac{1-\sqrt{5}}{2}$,
$ \frac{1-\sqrt{5}}{2}$,
$ -c^{3}_{40}
+c^{7}_{40}
$;\ \ 
$ 1$,
$ -\sqrt{2}$,
$ \frac{1-\sqrt{5}}{2}$,
$ \frac{1-\sqrt{5}}{2}$,
$ c^{3}_{40}
-c^{7}_{40}
$;\ \ 
$0$,
$ -c^{3}_{40}
+c^{7}_{40}
$,
$ c^{3}_{40}
-c^{7}_{40}
$,
$0$;\ \ 
$ -1$,
$ -1$,
$ -\sqrt{2}$;\ \ 
$ -1$,
$ \sqrt{2}$;\ \ 
$0$)

Not pseudo-unitary. 

\vskip 1ex 
\color{grey}

\noindent130. ind = $(6;80
)_{1}^{7}$:\ \ 
$d_i$ = ($1.0$,
$1.0$,
$1.414$,
$-0.618$,
$-0.618$,
$-0.874$) 

\vskip 0.7ex
\hangindent=3em \hangafter=1
$D^2=$ 5.527 = 
 $10-2\sqrt{5}$

\vskip 0.7ex
\hangindent=3em \hangafter=1
$T = ( 0,
\frac{1}{2},
\frac{7}{16},
\frac{4}{5},
\frac{3}{10},
\frac{19}{80} )
$,

\vskip 0.7ex
\hangindent=3em \hangafter=1
$S$ = ($ 1$,
$ 1$,
$ \sqrt{2}$,
$ \frac{1-\sqrt{5}}{2}$,
$ \frac{1-\sqrt{5}}{2}$,
$ -c^{3}_{40}
+c^{7}_{40}
$;\ \ 
$ 1$,
$ -\sqrt{2}$,
$ \frac{1-\sqrt{5}}{2}$,
$ \frac{1-\sqrt{5}}{2}$,
$ c^{3}_{40}
-c^{7}_{40}
$;\ \ 
$0$,
$ -c^{3}_{40}
+c^{7}_{40}
$,
$ c^{3}_{40}
-c^{7}_{40}
$,
$0$;\ \ 
$ -1$,
$ -1$,
$ -\sqrt{2}$;\ \ 
$ -1$,
$ \sqrt{2}$;\ \ 
$0$)

Not pseudo-unitary. 

\vskip 1ex 
\color{grey}

\noindent131. ind = $(6;80
)_{1}^{73}$:\ \ 
$d_i$ = ($1.0$,
$1.0$,
$1.414$,
$-0.618$,
$-0.618$,
$-0.874$) 

\vskip 0.7ex
\hangindent=3em \hangafter=1
$D^2=$ 5.527 = 
 $10-2\sqrt{5}$

\vskip 0.7ex
\hangindent=3em \hangafter=1
$T = ( 0,
\frac{1}{2},
\frac{9}{16},
\frac{1}{5},
\frac{7}{10},
\frac{61}{80} )
$,

\vskip 0.7ex
\hangindent=3em \hangafter=1
$S$ = ($ 1$,
$ 1$,
$ \sqrt{2}$,
$ \frac{1-\sqrt{5}}{2}$,
$ \frac{1-\sqrt{5}}{2}$,
$ -c^{3}_{40}
+c^{7}_{40}
$;\ \ 
$ 1$,
$ -\sqrt{2}$,
$ \frac{1-\sqrt{5}}{2}$,
$ \frac{1-\sqrt{5}}{2}$,
$ c^{3}_{40}
-c^{7}_{40}
$;\ \ 
$0$,
$ -c^{3}_{40}
+c^{7}_{40}
$,
$ c^{3}_{40}
-c^{7}_{40}
$,
$0$;\ \ 
$ -1$,
$ -1$,
$ -\sqrt{2}$;\ \ 
$ -1$,
$ \sqrt{2}$;\ \ 
$0$)

Not pseudo-unitary. 

\vskip 1ex 
\color{grey}

\noindent132. ind = $(6;80
)_{1}^{57}$:\ \ 
$d_i$ = ($1.0$,
$1.0$,
$1.414$,
$-0.618$,
$-0.618$,
$-0.874$) 

\vskip 0.7ex
\hangindent=3em \hangafter=1
$D^2=$ 5.527 = 
 $10-2\sqrt{5}$

\vskip 0.7ex
\hangindent=3em \hangafter=1
$T = ( 0,
\frac{1}{2},
\frac{9}{16},
\frac{4}{5},
\frac{3}{10},
\frac{29}{80} )
$,

\vskip 0.7ex
\hangindent=3em \hangafter=1
$S$ = ($ 1$,
$ 1$,
$ \sqrt{2}$,
$ \frac{1-\sqrt{5}}{2}$,
$ \frac{1-\sqrt{5}}{2}$,
$ -c^{3}_{40}
+c^{7}_{40}
$;\ \ 
$ 1$,
$ -\sqrt{2}$,
$ \frac{1-\sqrt{5}}{2}$,
$ \frac{1-\sqrt{5}}{2}$,
$ c^{3}_{40}
-c^{7}_{40}
$;\ \ 
$0$,
$ -c^{3}_{40}
+c^{7}_{40}
$,
$ c^{3}_{40}
-c^{7}_{40}
$,
$0$;\ \ 
$ -1$,
$ -1$,
$ -\sqrt{2}$;\ \ 
$ -1$,
$ \sqrt{2}$;\ \ 
$0$)

Not pseudo-unitary. 

\vskip 1ex 
\color{grey}

\noindent133. ind = $(6;80
)_{1}^{19}$:\ \ 
$d_i$ = ($1.0$,
$1.0$,
$1.618$,
$1.618$,
$-1.414$,
$-2.288$) 

\vskip 0.7ex
\hangindent=3em \hangafter=1
$D^2=$ 14.472 = 
 $10+2\sqrt{5}$

\vskip 0.7ex
\hangindent=3em \hangafter=1
$T = ( 0,
\frac{1}{2},
\frac{3}{5},
\frac{1}{10},
\frac{3}{16},
\frac{63}{80} )
$,

\vskip 0.7ex
\hangindent=3em \hangafter=1
$S$ = ($ 1$,
$ 1$,
$ \frac{1+\sqrt{5}}{2}$,
$ \frac{1+\sqrt{5}}{2}$,
$ -\sqrt{2}$,
$ -c^{3}_{40}
-c^{5}_{40}
+c^{7}_{40}
$;\ \ 
$ 1$,
$ \frac{1+\sqrt{5}}{2}$,
$ \frac{1+\sqrt{5}}{2}$,
$ \sqrt{2}$,
$ c^{3}_{40}
+c^{5}_{40}
-c^{7}_{40}
$;\ \ 
$ -1$,
$ -1$,
$ -c^{3}_{40}
-c^{5}_{40}
+c^{7}_{40}
$,
$ \sqrt{2}$;\ \ 
$ -1$,
$ c^{3}_{40}
+c^{5}_{40}
-c^{7}_{40}
$,
$ -\sqrt{2}$;\ \ 
$0$,
$0$;\ \ 
$0$)

Pseudo-unitary $\sim$  
$(6;80
)_{2}^{49}$

\vskip 1ex 
\color{grey}

\noindent134. ind = $(6;80
)_{1}^{69}$:\ \ 
$d_i$ = ($1.0$,
$1.0$,
$1.618$,
$1.618$,
$-1.414$,
$-2.288$) 

\vskip 0.7ex
\hangindent=3em \hangafter=1
$D^2=$ 14.472 = 
 $10+2\sqrt{5}$

\vskip 0.7ex
\hangindent=3em \hangafter=1
$T = ( 0,
\frac{1}{2},
\frac{3}{5},
\frac{1}{10},
\frac{5}{16},
\frac{73}{80} )
$,

\vskip 0.7ex
\hangindent=3em \hangafter=1
$S$ = ($ 1$,
$ 1$,
$ \frac{1+\sqrt{5}}{2}$,
$ \frac{1+\sqrt{5}}{2}$,
$ -\sqrt{2}$,
$ -c^{3}_{40}
-c^{5}_{40}
+c^{7}_{40}
$;\ \ 
$ 1$,
$ \frac{1+\sqrt{5}}{2}$,
$ \frac{1+\sqrt{5}}{2}$,
$ \sqrt{2}$,
$ c^{3}_{40}
+c^{5}_{40}
-c^{7}_{40}
$;\ \ 
$ -1$,
$ -1$,
$ -c^{3}_{40}
-c^{5}_{40}
+c^{7}_{40}
$,
$ \sqrt{2}$;\ \ 
$ -1$,
$ c^{3}_{40}
+c^{5}_{40}
-c^{7}_{40}
$,
$ -\sqrt{2}$;\ \ 
$0$,
$0$;\ \ 
$0$)

Pseudo-unitary $\sim$  
$(6;80
)_{2}^{39}$

\vskip 1ex 
\color{grey}

\noindent135. ind = $(6;80
)_{1}^{59}$:\ \ 
$d_i$ = ($1.0$,
$1.0$,
$1.618$,
$1.618$,
$-1.414$,
$-2.288$) 

\vskip 0.7ex
\hangindent=3em \hangafter=1
$D^2=$ 14.472 = 
 $10+2\sqrt{5}$

\vskip 0.7ex
\hangindent=3em \hangafter=1
$T = ( 0,
\frac{1}{2},
\frac{3}{5},
\frac{1}{10},
\frac{11}{16},
\frac{23}{80} )
$,

\vskip 0.7ex
\hangindent=3em \hangafter=1
$S$ = ($ 1$,
$ 1$,
$ \frac{1+\sqrt{5}}{2}$,
$ \frac{1+\sqrt{5}}{2}$,
$ -\sqrt{2}$,
$ -c^{3}_{40}
-c^{5}_{40}
+c^{7}_{40}
$;\ \ 
$ 1$,
$ \frac{1+\sqrt{5}}{2}$,
$ \frac{1+\sqrt{5}}{2}$,
$ \sqrt{2}$,
$ c^{3}_{40}
+c^{5}_{40}
-c^{7}_{40}
$;\ \ 
$ -1$,
$ -1$,
$ -c^{3}_{40}
-c^{5}_{40}
+c^{7}_{40}
$,
$ \sqrt{2}$;\ \ 
$ -1$,
$ c^{3}_{40}
+c^{5}_{40}
-c^{7}_{40}
$,
$ -\sqrt{2}$;\ \ 
$0$,
$0$;\ \ 
$0$)

Pseudo-unitary $\sim$  
$(6;80
)_{2}^{9}$

\vskip 1ex 
\color{grey}

\noindent136. ind = $(6;80
)_{1}^{29}$:\ \ 
$d_i$ = ($1.0$,
$1.0$,
$1.618$,
$1.618$,
$-1.414$,
$-2.288$) 

\vskip 0.7ex
\hangindent=3em \hangafter=1
$D^2=$ 14.472 = 
 $10+2\sqrt{5}$

\vskip 0.7ex
\hangindent=3em \hangafter=1
$T = ( 0,
\frac{1}{2},
\frac{3}{5},
\frac{1}{10},
\frac{13}{16},
\frac{33}{80} )
$,

\vskip 0.7ex
\hangindent=3em \hangafter=1
$S$ = ($ 1$,
$ 1$,
$ \frac{1+\sqrt{5}}{2}$,
$ \frac{1+\sqrt{5}}{2}$,
$ -\sqrt{2}$,
$ -c^{3}_{40}
-c^{5}_{40}
+c^{7}_{40}
$;\ \ 
$ 1$,
$ \frac{1+\sqrt{5}}{2}$,
$ \frac{1+\sqrt{5}}{2}$,
$ \sqrt{2}$,
$ c^{3}_{40}
+c^{5}_{40}
-c^{7}_{40}
$;\ \ 
$ -1$,
$ -1$,
$ -c^{3}_{40}
-c^{5}_{40}
+c^{7}_{40}
$,
$ \sqrt{2}$;\ \ 
$ -1$,
$ c^{3}_{40}
+c^{5}_{40}
-c^{7}_{40}
$,
$ -\sqrt{2}$;\ \ 
$0$,
$0$;\ \ 
$0$)

Pseudo-unitary $\sim$  
$(6;80
)_{2}^{79}$

\vskip 1ex 
\color{grey}

\noindent137. ind = $(6;80
)_{1}^{63}$:\ \ 
$d_i$ = ($1.0$,
$1.0$,
$1.414$,
$-0.618$,
$-0.618$,
$-0.874$) 

\vskip 0.7ex
\hangindent=3em \hangafter=1
$D^2=$ 5.527 = 
 $10-2\sqrt{5}$

\vskip 0.7ex
\hangindent=3em \hangafter=1
$T = ( 0,
\frac{1}{2},
\frac{15}{16},
\frac{1}{5},
\frac{7}{10},
\frac{11}{80} )
$,

\vskip 0.7ex
\hangindent=3em \hangafter=1
$S$ = ($ 1$,
$ 1$,
$ \sqrt{2}$,
$ \frac{1-\sqrt{5}}{2}$,
$ \frac{1-\sqrt{5}}{2}$,
$ -c^{3}_{40}
+c^{7}_{40}
$;\ \ 
$ 1$,
$ -\sqrt{2}$,
$ \frac{1-\sqrt{5}}{2}$,
$ \frac{1-\sqrt{5}}{2}$,
$ c^{3}_{40}
-c^{7}_{40}
$;\ \ 
$0$,
$ -c^{3}_{40}
+c^{7}_{40}
$,
$ c^{3}_{40}
-c^{7}_{40}
$,
$0$;\ \ 
$ -1$,
$ -1$,
$ -\sqrt{2}$;\ \ 
$ -1$,
$ \sqrt{2}$;\ \ 
$0$)

Not pseudo-unitary. 

\vskip 1ex 
\color{grey}

\noindent138. ind = $(6;80
)_{1}^{47}$:\ \ 
$d_i$ = ($1.0$,
$1.0$,
$1.414$,
$-0.618$,
$-0.618$,
$-0.874$) 

\vskip 0.7ex
\hangindent=3em \hangafter=1
$D^2=$ 5.527 = 
 $10-2\sqrt{5}$

\vskip 0.7ex
\hangindent=3em \hangafter=1
$T = ( 0,
\frac{1}{2},
\frac{15}{16},
\frac{4}{5},
\frac{3}{10},
\frac{59}{80} )
$,

\vskip 0.7ex
\hangindent=3em \hangafter=1
$S$ = ($ 1$,
$ 1$,
$ \sqrt{2}$,
$ \frac{1-\sqrt{5}}{2}$,
$ \frac{1-\sqrt{5}}{2}$,
$ -c^{3}_{40}
+c^{7}_{40}
$;\ \ 
$ 1$,
$ -\sqrt{2}$,
$ \frac{1-\sqrt{5}}{2}$,
$ \frac{1-\sqrt{5}}{2}$,
$ c^{3}_{40}
-c^{7}_{40}
$;\ \ 
$0$,
$ -c^{3}_{40}
+c^{7}_{40}
$,
$ c^{3}_{40}
-c^{7}_{40}
$,
$0$;\ \ 
$ -1$,
$ -1$,
$ -\sqrt{2}$;\ \ 
$ -1$,
$ \sqrt{2}$;\ \ 
$0$)

Not pseudo-unitary. 

\vskip 1ex 
\color{grey}

\noindent139. ind = $(6;80
)_{1}^{53}$:\ \ 
$d_i$ = ($1.0$,
$0.874$,
$1.0$,
$-0.618$,
$-0.618$,
$-1.414$) 

\vskip 0.7ex
\hangindent=3em \hangafter=1
$D^2=$ 5.527 = 
 $10-2\sqrt{5}$

\vskip 0.7ex
\hangindent=3em \hangafter=1
$T = ( 0,
\frac{41}{80},
\frac{1}{2},
\frac{1}{5},
\frac{7}{10},
\frac{5}{16} )
$,

\vskip 0.7ex
\hangindent=3em \hangafter=1
$S$ = ($ 1$,
$ c^{3}_{40}
-c^{7}_{40}
$,
$ 1$,
$ \frac{1-\sqrt{5}}{2}$,
$ \frac{1-\sqrt{5}}{2}$,
$ -\sqrt{2}$;\ \ 
$0$,
$ -c^{3}_{40}
+c^{7}_{40}
$,
$ \sqrt{2}$,
$ -\sqrt{2}$,
$0$;\ \ 
$ 1$,
$ \frac{1-\sqrt{5}}{2}$,
$ \frac{1-\sqrt{5}}{2}$,
$ \sqrt{2}$;\ \ 
$ -1$,
$ -1$,
$ c^{3}_{40}
-c^{7}_{40}
$;\ \ 
$ -1$,
$ -c^{3}_{40}
+c^{7}_{40}
$;\ \ 
$0$)

Not pseudo-unitary. 

\vskip 1ex 
\color{grey}

\noindent140. ind = $(6;80
)_{1}^{77}$:\ \ 
$d_i$ = ($1.0$,
$0.874$,
$1.0$,
$-0.618$,
$-0.618$,
$-1.414$) 

\vskip 0.7ex
\hangindent=3em \hangafter=1
$D^2=$ 5.527 = 
 $10-2\sqrt{5}$

\vskip 0.7ex
\hangindent=3em \hangafter=1
$T = ( 0,
\frac{49}{80},
\frac{1}{2},
\frac{4}{5},
\frac{3}{10},
\frac{13}{16} )
$,

\vskip 0.7ex
\hangindent=3em \hangafter=1
$S$ = ($ 1$,
$ c^{3}_{40}
-c^{7}_{40}
$,
$ 1$,
$ \frac{1-\sqrt{5}}{2}$,
$ \frac{1-\sqrt{5}}{2}$,
$ -\sqrt{2}$;\ \ 
$0$,
$ -c^{3}_{40}
+c^{7}_{40}
$,
$ \sqrt{2}$,
$ -\sqrt{2}$,
$0$;\ \ 
$ 1$,
$ \frac{1-\sqrt{5}}{2}$,
$ \frac{1-\sqrt{5}}{2}$,
$ \sqrt{2}$;\ \ 
$ -1$,
$ -1$,
$ c^{3}_{40}
-c^{7}_{40}
$;\ \ 
$ -1$,
$ -c^{3}_{40}
+c^{7}_{40}
$;\ \ 
$0$)

Not pseudo-unitary. 

\vskip 1ex 
\color{grey}

\noindent141. ind = $(6;80
)_{1}^{43}$:\ \ 
$d_i$ = ($1.0$,
$0.874$,
$1.0$,
$-0.618$,
$-0.618$,
$-1.414$) 

\vskip 0.7ex
\hangindent=3em \hangafter=1
$D^2=$ 5.527 = 
 $10-2\sqrt{5}$

\vskip 0.7ex
\hangindent=3em \hangafter=1
$T = ( 0,
\frac{71}{80},
\frac{1}{2},
\frac{1}{5},
\frac{7}{10},
\frac{11}{16} )
$,

\vskip 0.7ex
\hangindent=3em \hangafter=1
$S$ = ($ 1$,
$ c^{3}_{40}
-c^{7}_{40}
$,
$ 1$,
$ \frac{1-\sqrt{5}}{2}$,
$ \frac{1-\sqrt{5}}{2}$,
$ -\sqrt{2}$;\ \ 
$0$,
$ -c^{3}_{40}
+c^{7}_{40}
$,
$ \sqrt{2}$,
$ -\sqrt{2}$,
$0$;\ \ 
$ 1$,
$ \frac{1-\sqrt{5}}{2}$,
$ \frac{1-\sqrt{5}}{2}$,
$ \sqrt{2}$;\ \ 
$ -1$,
$ -1$,
$ c^{3}_{40}
-c^{7}_{40}
$;\ \ 
$ -1$,
$ -c^{3}_{40}
+c^{7}_{40}
$;\ \ 
$0$)

Not pseudo-unitary. 

\vskip 1ex 
\color{grey}

\noindent142. ind = $(6;80
)_{1}^{67}$:\ \ 
$d_i$ = ($1.0$,
$0.874$,
$1.0$,
$-0.618$,
$-0.618$,
$-1.414$) 

\vskip 0.7ex
\hangindent=3em \hangafter=1
$D^2=$ 5.527 = 
 $10-2\sqrt{5}$

\vskip 0.7ex
\hangindent=3em \hangafter=1
$T = ( 0,
\frac{79}{80},
\frac{1}{2},
\frac{4}{5},
\frac{3}{10},
\frac{3}{16} )
$,

\vskip 0.7ex
\hangindent=3em \hangafter=1
$S$ = ($ 1$,
$ c^{3}_{40}
-c^{7}_{40}
$,
$ 1$,
$ \frac{1-\sqrt{5}}{2}$,
$ \frac{1-\sqrt{5}}{2}$,
$ -\sqrt{2}$;\ \ 
$0$,
$ -c^{3}_{40}
+c^{7}_{40}
$,
$ \sqrt{2}$,
$ -\sqrt{2}$,
$0$;\ \ 
$ 1$,
$ \frac{1-\sqrt{5}}{2}$,
$ \frac{1-\sqrt{5}}{2}$,
$ \sqrt{2}$;\ \ 
$ -1$,
$ -1$,
$ c^{3}_{40}
-c^{7}_{40}
$;\ \ 
$ -1$,
$ -c^{3}_{40}
+c^{7}_{40}
$;\ \ 
$0$)

Not pseudo-unitary. 

\vskip 1ex 

 \color{black} \vskip 2ex

\noindent143. ind = $(6;80
)_{2}^{1}$:\ \ 
$d_i$ = ($1.0$,
$1.0$,
$1.414$,
$1.618$,
$1.618$,
$2.288$) 

\vskip 0.7ex
\hangindent=3em \hangafter=1
$D^2=$ 14.472 = 
 $10+2\sqrt{5}$

\vskip 0.7ex
\hangindent=3em \hangafter=1
$T = ( 0,
\frac{1}{2},
\frac{3}{16},
\frac{2}{5},
\frac{9}{10},
\frac{47}{80} )
$,

\vskip 0.7ex
\hangindent=3em \hangafter=1
$S$ = ($ 1$,
$ 1$,
$ \sqrt{2}$,
$ \frac{1+\sqrt{5}}{2}$,
$ \frac{1+\sqrt{5}}{2}$,
$ c^{3}_{40}
+c^{5}_{40}
-c^{7}_{40}
$;\ \ 
$ 1$,
$ -\sqrt{2}$,
$ \frac{1+\sqrt{5}}{2}$,
$ \frac{1+\sqrt{5}}{2}$,
$ -c^{3}_{40}
-c^{5}_{40}
+c^{7}_{40}
$;\ \ 
$0$,
$ c^{3}_{40}
+c^{5}_{40}
-c^{7}_{40}
$,
$ -c^{3}_{40}
-c^{5}_{40}
+c^{7}_{40}
$,
$0$;\ \ 
$ -1$,
$ -1$,
$ -\sqrt{2}$;\ \ 
$ -1$,
$ \sqrt{2}$;\ \ 
$0$)

\vskip 1ex 
\color{grey}

\noindent144. ind = $(6;80
)_{2}^{49}$:\ \ 
$d_i$ = ($1.0$,
$1.0$,
$1.414$,
$1.618$,
$1.618$,
$2.288$) 

\vskip 0.7ex
\hangindent=3em \hangafter=1
$D^2=$ 14.472 = 
 $10+2\sqrt{5}$

\vskip 0.7ex
\hangindent=3em \hangafter=1
$T = ( 0,
\frac{1}{2},
\frac{3}{16},
\frac{3}{5},
\frac{1}{10},
\frac{63}{80} )
$,

\vskip 0.7ex
\hangindent=3em \hangafter=1
$S$ = ($ 1$,
$ 1$,
$ \sqrt{2}$,
$ \frac{1+\sqrt{5}}{2}$,
$ \frac{1+\sqrt{5}}{2}$,
$ c^{3}_{40}
+c^{5}_{40}
-c^{7}_{40}
$;\ \ 
$ 1$,
$ -\sqrt{2}$,
$ \frac{1+\sqrt{5}}{2}$,
$ \frac{1+\sqrt{5}}{2}$,
$ -c^{3}_{40}
-c^{5}_{40}
+c^{7}_{40}
$;\ \ 
$0$,
$ c^{3}_{40}
+c^{5}_{40}
-c^{7}_{40}
$,
$ -c^{3}_{40}
-c^{5}_{40}
+c^{7}_{40}
$,
$0$;\ \ 
$ -1$,
$ -1$,
$ -\sqrt{2}$;\ \ 
$ -1$,
$ \sqrt{2}$;\ \ 
$0$)

\vskip 1ex 
\color{grey}

\noindent145. ind = $(6;80
)_{2}^{71}$:\ \ 
$d_i$ = ($1.0$,
$1.0$,
$1.414$,
$1.618$,
$1.618$,
$2.288$) 

\vskip 0.7ex
\hangindent=3em \hangafter=1
$D^2=$ 14.472 = 
 $10+2\sqrt{5}$

\vskip 0.7ex
\hangindent=3em \hangafter=1
$T = ( 0,
\frac{1}{2},
\frac{5}{16},
\frac{2}{5},
\frac{9}{10},
\frac{57}{80} )
$,

\vskip 0.7ex
\hangindent=3em \hangafter=1
$S$ = ($ 1$,
$ 1$,
$ \sqrt{2}$,
$ \frac{1+\sqrt{5}}{2}$,
$ \frac{1+\sqrt{5}}{2}$,
$ c^{3}_{40}
+c^{5}_{40}
-c^{7}_{40}
$;\ \ 
$ 1$,
$ -\sqrt{2}$,
$ \frac{1+\sqrt{5}}{2}$,
$ \frac{1+\sqrt{5}}{2}$,
$ -c^{3}_{40}
-c^{5}_{40}
+c^{7}_{40}
$;\ \ 
$0$,
$ c^{3}_{40}
+c^{5}_{40}
-c^{7}_{40}
$,
$ -c^{3}_{40}
-c^{5}_{40}
+c^{7}_{40}
$,
$0$;\ \ 
$ -1$,
$ -1$,
$ -\sqrt{2}$;\ \ 
$ -1$,
$ \sqrt{2}$;\ \ 
$0$)

\vskip 1ex 
\color{grey}

\noindent146. ind = $(6;80
)_{2}^{39}$:\ \ 
$d_i$ = ($1.0$,
$1.0$,
$1.414$,
$1.618$,
$1.618$,
$2.288$) 

\vskip 0.7ex
\hangindent=3em \hangafter=1
$D^2=$ 14.472 = 
 $10+2\sqrt{5}$

\vskip 0.7ex
\hangindent=3em \hangafter=1
$T = ( 0,
\frac{1}{2},
\frac{5}{16},
\frac{3}{5},
\frac{1}{10},
\frac{73}{80} )
$,

\vskip 0.7ex
\hangindent=3em \hangafter=1
$S$ = ($ 1$,
$ 1$,
$ \sqrt{2}$,
$ \frac{1+\sqrt{5}}{2}$,
$ \frac{1+\sqrt{5}}{2}$,
$ c^{3}_{40}
+c^{5}_{40}
-c^{7}_{40}
$;\ \ 
$ 1$,
$ -\sqrt{2}$,
$ \frac{1+\sqrt{5}}{2}$,
$ \frac{1+\sqrt{5}}{2}$,
$ -c^{3}_{40}
-c^{5}_{40}
+c^{7}_{40}
$;\ \ 
$0$,
$ c^{3}_{40}
+c^{5}_{40}
-c^{7}_{40}
$,
$ -c^{3}_{40}
-c^{5}_{40}
+c^{7}_{40}
$,
$0$;\ \ 
$ -1$,
$ -1$,
$ -\sqrt{2}$;\ \ 
$ -1$,
$ \sqrt{2}$;\ \ 
$0$)

\vskip 1ex 
\color{grey}

\noindent147. ind = $(6;80
)_{2}^{41}$:\ \ 
$d_i$ = ($1.0$,
$1.0$,
$1.414$,
$1.618$,
$1.618$,
$2.288$) 

\vskip 0.7ex
\hangindent=3em \hangafter=1
$D^2=$ 14.472 = 
 $10+2\sqrt{5}$

\vskip 0.7ex
\hangindent=3em \hangafter=1
$T = ( 0,
\frac{1}{2},
\frac{11}{16},
\frac{2}{5},
\frac{9}{10},
\frac{7}{80} )
$,

\vskip 0.7ex
\hangindent=3em \hangafter=1
$S$ = ($ 1$,
$ 1$,
$ \sqrt{2}$,
$ \frac{1+\sqrt{5}}{2}$,
$ \frac{1+\sqrt{5}}{2}$,
$ c^{3}_{40}
+c^{5}_{40}
-c^{7}_{40}
$;\ \ 
$ 1$,
$ -\sqrt{2}$,
$ \frac{1+\sqrt{5}}{2}$,
$ \frac{1+\sqrt{5}}{2}$,
$ -c^{3}_{40}
-c^{5}_{40}
+c^{7}_{40}
$;\ \ 
$0$,
$ c^{3}_{40}
+c^{5}_{40}
-c^{7}_{40}
$,
$ -c^{3}_{40}
-c^{5}_{40}
+c^{7}_{40}
$,
$0$;\ \ 
$ -1$,
$ -1$,
$ -\sqrt{2}$;\ \ 
$ -1$,
$ \sqrt{2}$;\ \ 
$0$)

\vskip 1ex 
\color{grey}

\noindent148. ind = $(6;80
)_{2}^{9}$:\ \ 
$d_i$ = ($1.0$,
$1.0$,
$1.414$,
$1.618$,
$1.618$,
$2.288$) 

\vskip 0.7ex
\hangindent=3em \hangafter=1
$D^2=$ 14.472 = 
 $10+2\sqrt{5}$

\vskip 0.7ex
\hangindent=3em \hangafter=1
$T = ( 0,
\frac{1}{2},
\frac{11}{16},
\frac{3}{5},
\frac{1}{10},
\frac{23}{80} )
$,

\vskip 0.7ex
\hangindent=3em \hangafter=1
$S$ = ($ 1$,
$ 1$,
$ \sqrt{2}$,
$ \frac{1+\sqrt{5}}{2}$,
$ \frac{1+\sqrt{5}}{2}$,
$ c^{3}_{40}
+c^{5}_{40}
-c^{7}_{40}
$;\ \ 
$ 1$,
$ -\sqrt{2}$,
$ \frac{1+\sqrt{5}}{2}$,
$ \frac{1+\sqrt{5}}{2}$,
$ -c^{3}_{40}
-c^{5}_{40}
+c^{7}_{40}
$;\ \ 
$0$,
$ c^{3}_{40}
+c^{5}_{40}
-c^{7}_{40}
$,
$ -c^{3}_{40}
-c^{5}_{40}
+c^{7}_{40}
$,
$0$;\ \ 
$ -1$,
$ -1$,
$ -\sqrt{2}$;\ \ 
$ -1$,
$ \sqrt{2}$;\ \ 
$0$)

\vskip 1ex 
\color{grey}

\noindent149. ind = $(6;80
)_{2}^{31}$:\ \ 
$d_i$ = ($1.0$,
$1.0$,
$1.414$,
$1.618$,
$1.618$,
$2.288$) 

\vskip 0.7ex
\hangindent=3em \hangafter=1
$D^2=$ 14.472 = 
 $10+2\sqrt{5}$

\vskip 0.7ex
\hangindent=3em \hangafter=1
$T = ( 0,
\frac{1}{2},
\frac{13}{16},
\frac{2}{5},
\frac{9}{10},
\frac{17}{80} )
$,

\vskip 0.7ex
\hangindent=3em \hangafter=1
$S$ = ($ 1$,
$ 1$,
$ \sqrt{2}$,
$ \frac{1+\sqrt{5}}{2}$,
$ \frac{1+\sqrt{5}}{2}$,
$ c^{3}_{40}
+c^{5}_{40}
-c^{7}_{40}
$;\ \ 
$ 1$,
$ -\sqrt{2}$,
$ \frac{1+\sqrt{5}}{2}$,
$ \frac{1+\sqrt{5}}{2}$,
$ -c^{3}_{40}
-c^{5}_{40}
+c^{7}_{40}
$;\ \ 
$0$,
$ c^{3}_{40}
+c^{5}_{40}
-c^{7}_{40}
$,
$ -c^{3}_{40}
-c^{5}_{40}
+c^{7}_{40}
$,
$0$;\ \ 
$ -1$,
$ -1$,
$ -\sqrt{2}$;\ \ 
$ -1$,
$ \sqrt{2}$;\ \ 
$0$)

\vskip 1ex 
\color{grey}

\noindent150. ind = $(6;80
)_{2}^{79}$:\ \ 
$d_i$ = ($1.0$,
$1.0$,
$1.414$,
$1.618$,
$1.618$,
$2.288$) 

\vskip 0.7ex
\hangindent=3em \hangafter=1
$D^2=$ 14.472 = 
 $10+2\sqrt{5}$

\vskip 0.7ex
\hangindent=3em \hangafter=1
$T = ( 0,
\frac{1}{2},
\frac{13}{16},
\frac{3}{5},
\frac{1}{10},
\frac{33}{80} )
$,

\vskip 0.7ex
\hangindent=3em \hangafter=1
$S$ = ($ 1$,
$ 1$,
$ \sqrt{2}$,
$ \frac{1+\sqrt{5}}{2}$,
$ \frac{1+\sqrt{5}}{2}$,
$ c^{3}_{40}
+c^{5}_{40}
-c^{7}_{40}
$;\ \ 
$ 1$,
$ -\sqrt{2}$,
$ \frac{1+\sqrt{5}}{2}$,
$ \frac{1+\sqrt{5}}{2}$,
$ -c^{3}_{40}
-c^{5}_{40}
+c^{7}_{40}
$;\ \ 
$0$,
$ c^{3}_{40}
+c^{5}_{40}
-c^{7}_{40}
$,
$ -c^{3}_{40}
-c^{5}_{40}
+c^{7}_{40}
$,
$0$;\ \ 
$ -1$,
$ -1$,
$ -\sqrt{2}$;\ \ 
$ -1$,
$ \sqrt{2}$;\ \ 
$0$)

\vskip 1ex 
\color{grey}

\noindent151. ind = $(6;80
)_{2}^{53}$:\ \ 
$d_i$ = ($1.0$,
$0.874$,
$1.0$,
$-0.618$,
$-0.618$,
$-1.414$) 

\vskip 0.7ex
\hangindent=3em \hangafter=1
$D^2=$ 5.527 = 
 $10-2\sqrt{5}$

\vskip 0.7ex
\hangindent=3em \hangafter=1
$T = ( 0,
\frac{11}{80},
\frac{1}{2},
\frac{1}{5},
\frac{7}{10},
\frac{15}{16} )
$,

\vskip 0.7ex
\hangindent=3em \hangafter=1
$S$ = ($ 1$,
$ c^{3}_{40}
-c^{7}_{40}
$,
$ 1$,
$ \frac{1-\sqrt{5}}{2}$,
$ \frac{1-\sqrt{5}}{2}$,
$ -\sqrt{2}$;\ \ 
$0$,
$ -c^{3}_{40}
+c^{7}_{40}
$,
$ \sqrt{2}$,
$ -\sqrt{2}$,
$0$;\ \ 
$ 1$,
$ \frac{1-\sqrt{5}}{2}$,
$ \frac{1-\sqrt{5}}{2}$,
$ \sqrt{2}$;\ \ 
$ -1$,
$ -1$,
$ c^{3}_{40}
-c^{7}_{40}
$;\ \ 
$ -1$,
$ -c^{3}_{40}
+c^{7}_{40}
$;\ \ 
$0$)

Not pseudo-unitary. 

\vskip 1ex 
\color{grey}

\noindent152. ind = $(6;80
)_{2}^{77}$:\ \ 
$d_i$ = ($1.0$,
$0.874$,
$1.0$,
$-0.618$,
$-0.618$,
$-1.414$) 

\vskip 0.7ex
\hangindent=3em \hangafter=1
$D^2=$ 5.527 = 
 $10-2\sqrt{5}$

\vskip 0.7ex
\hangindent=3em \hangafter=1
$T = ( 0,
\frac{19}{80},
\frac{1}{2},
\frac{4}{5},
\frac{3}{10},
\frac{7}{16} )
$,

\vskip 0.7ex
\hangindent=3em \hangafter=1
$S$ = ($ 1$,
$ c^{3}_{40}
-c^{7}_{40}
$,
$ 1$,
$ \frac{1-\sqrt{5}}{2}$,
$ \frac{1-\sqrt{5}}{2}$,
$ -\sqrt{2}$;\ \ 
$0$,
$ -c^{3}_{40}
+c^{7}_{40}
$,
$ \sqrt{2}$,
$ -\sqrt{2}$,
$0$;\ \ 
$ 1$,
$ \frac{1-\sqrt{5}}{2}$,
$ \frac{1-\sqrt{5}}{2}$,
$ \sqrt{2}$;\ \ 
$ -1$,
$ -1$,
$ c^{3}_{40}
-c^{7}_{40}
$;\ \ 
$ -1$,
$ -c^{3}_{40}
+c^{7}_{40}
$;\ \ 
$0$)

Not pseudo-unitary. 

\vskip 1ex 
\color{grey}

\noindent153. ind = $(6;80
)_{2}^{43}$:\ \ 
$d_i$ = ($1.0$,
$0.874$,
$1.0$,
$-0.618$,
$-0.618$,
$-1.414$) 

\vskip 0.7ex
\hangindent=3em \hangafter=1
$D^2=$ 5.527 = 
 $10-2\sqrt{5}$

\vskip 0.7ex
\hangindent=3em \hangafter=1
$T = ( 0,
\frac{21}{80},
\frac{1}{2},
\frac{1}{5},
\frac{7}{10},
\frac{1}{16} )
$,

\vskip 0.7ex
\hangindent=3em \hangafter=1
$S$ = ($ 1$,
$ c^{3}_{40}
-c^{7}_{40}
$,
$ 1$,
$ \frac{1-\sqrt{5}}{2}$,
$ \frac{1-\sqrt{5}}{2}$,
$ -\sqrt{2}$;\ \ 
$0$,
$ -c^{3}_{40}
+c^{7}_{40}
$,
$ \sqrt{2}$,
$ -\sqrt{2}$,
$0$;\ \ 
$ 1$,
$ \frac{1-\sqrt{5}}{2}$,
$ \frac{1-\sqrt{5}}{2}$,
$ \sqrt{2}$;\ \ 
$ -1$,
$ -1$,
$ c^{3}_{40}
-c^{7}_{40}
$;\ \ 
$ -1$,
$ -c^{3}_{40}
+c^{7}_{40}
$;\ \ 
$0$)

Not pseudo-unitary. 

\vskip 1ex 
\color{grey}

\noindent154. ind = $(6;80
)_{2}^{67}$:\ \ 
$d_i$ = ($1.0$,
$0.874$,
$1.0$,
$-0.618$,
$-0.618$,
$-1.414$) 

\vskip 0.7ex
\hangindent=3em \hangafter=1
$D^2=$ 5.527 = 
 $10-2\sqrt{5}$

\vskip 0.7ex
\hangindent=3em \hangafter=1
$T = ( 0,
\frac{29}{80},
\frac{1}{2},
\frac{4}{5},
\frac{3}{10},
\frac{9}{16} )
$,

\vskip 0.7ex
\hangindent=3em \hangafter=1
$S$ = ($ 1$,
$ c^{3}_{40}
-c^{7}_{40}
$,
$ 1$,
$ \frac{1-\sqrt{5}}{2}$,
$ \frac{1-\sqrt{5}}{2}$,
$ -\sqrt{2}$;\ \ 
$0$,
$ -c^{3}_{40}
+c^{7}_{40}
$,
$ \sqrt{2}$,
$ -\sqrt{2}$,
$0$;\ \ 
$ 1$,
$ \frac{1-\sqrt{5}}{2}$,
$ \frac{1-\sqrt{5}}{2}$,
$ \sqrt{2}$;\ \ 
$ -1$,
$ -1$,
$ c^{3}_{40}
-c^{7}_{40}
$;\ \ 
$ -1$,
$ -c^{3}_{40}
+c^{7}_{40}
$;\ \ 
$0$)

Not pseudo-unitary. 

\vskip 1ex 
\color{grey}

\noindent155. ind = $(6;80
)_{2}^{33}$:\ \ 
$d_i$ = ($1.0$,
$1.0$,
$1.414$,
$-0.618$,
$-0.618$,
$-0.874$) 

\vskip 0.7ex
\hangindent=3em \hangafter=1
$D^2=$ 5.527 = 
 $10-2\sqrt{5}$

\vskip 0.7ex
\hangindent=3em \hangafter=1
$T = ( 0,
\frac{1}{2},
\frac{3}{16},
\frac{1}{5},
\frac{7}{10},
\frac{31}{80} )
$,

\vskip 0.7ex
\hangindent=3em \hangafter=1
$S$ = ($ 1$,
$ 1$,
$ \sqrt{2}$,
$ \frac{1-\sqrt{5}}{2}$,
$ \frac{1-\sqrt{5}}{2}$,
$ -c^{3}_{40}
+c^{7}_{40}
$;\ \ 
$ 1$,
$ -\sqrt{2}$,
$ \frac{1-\sqrt{5}}{2}$,
$ \frac{1-\sqrt{5}}{2}$,
$ c^{3}_{40}
-c^{7}_{40}
$;\ \ 
$0$,
$ -c^{3}_{40}
+c^{7}_{40}
$,
$ c^{3}_{40}
-c^{7}_{40}
$,
$0$;\ \ 
$ -1$,
$ -1$,
$ -\sqrt{2}$;\ \ 
$ -1$,
$ \sqrt{2}$;\ \ 
$0$)

Not pseudo-unitary. 

\vskip 1ex 
\color{grey}

\noindent156. ind = $(6;80
)_{2}^{17}$:\ \ 
$d_i$ = ($1.0$,
$1.0$,
$1.414$,
$-0.618$,
$-0.618$,
$-0.874$) 

\vskip 0.7ex
\hangindent=3em \hangafter=1
$D^2=$ 5.527 = 
 $10-2\sqrt{5}$

\vskip 0.7ex
\hangindent=3em \hangafter=1
$T = ( 0,
\frac{1}{2},
\frac{3}{16},
\frac{4}{5},
\frac{3}{10},
\frac{79}{80} )
$,

\vskip 0.7ex
\hangindent=3em \hangafter=1
$S$ = ($ 1$,
$ 1$,
$ \sqrt{2}$,
$ \frac{1-\sqrt{5}}{2}$,
$ \frac{1-\sqrt{5}}{2}$,
$ -c^{3}_{40}
+c^{7}_{40}
$;\ \ 
$ 1$,
$ -\sqrt{2}$,
$ \frac{1-\sqrt{5}}{2}$,
$ \frac{1-\sqrt{5}}{2}$,
$ c^{3}_{40}
-c^{7}_{40}
$;\ \ 
$0$,
$ -c^{3}_{40}
+c^{7}_{40}
$,
$ c^{3}_{40}
-c^{7}_{40}
$,
$0$;\ \ 
$ -1$,
$ -1$,
$ -\sqrt{2}$;\ \ 
$ -1$,
$ \sqrt{2}$;\ \ 
$0$)

Not pseudo-unitary. 

\vskip 1ex 
\color{grey}

\noindent157. ind = $(6;80
)_{2}^{23}$:\ \ 
$d_i$ = ($1.0$,
$1.0$,
$1.414$,
$-0.618$,
$-0.618$,
$-0.874$) 

\vskip 0.7ex
\hangindent=3em \hangafter=1
$D^2=$ 5.527 = 
 $10-2\sqrt{5}$

\vskip 0.7ex
\hangindent=3em \hangafter=1
$T = ( 0,
\frac{1}{2},
\frac{5}{16},
\frac{1}{5},
\frac{7}{10},
\frac{41}{80} )
$,

\vskip 0.7ex
\hangindent=3em \hangafter=1
$S$ = ($ 1$,
$ 1$,
$ \sqrt{2}$,
$ \frac{1-\sqrt{5}}{2}$,
$ \frac{1-\sqrt{5}}{2}$,
$ -c^{3}_{40}
+c^{7}_{40}
$;\ \ 
$ 1$,
$ -\sqrt{2}$,
$ \frac{1-\sqrt{5}}{2}$,
$ \frac{1-\sqrt{5}}{2}$,
$ c^{3}_{40}
-c^{7}_{40}
$;\ \ 
$0$,
$ -c^{3}_{40}
+c^{7}_{40}
$,
$ c^{3}_{40}
-c^{7}_{40}
$,
$0$;\ \ 
$ -1$,
$ -1$,
$ -\sqrt{2}$;\ \ 
$ -1$,
$ \sqrt{2}$;\ \ 
$0$)

Not pseudo-unitary. 

\vskip 1ex 
\color{grey}

\noindent158. ind = $(6;80
)_{2}^{7}$:\ \ 
$d_i$ = ($1.0$,
$1.0$,
$1.414$,
$-0.618$,
$-0.618$,
$-0.874$) 

\vskip 0.7ex
\hangindent=3em \hangafter=1
$D^2=$ 5.527 = 
 $10-2\sqrt{5}$

\vskip 0.7ex
\hangindent=3em \hangafter=1
$T = ( 0,
\frac{1}{2},
\frac{5}{16},
\frac{4}{5},
\frac{3}{10},
\frac{9}{80} )
$,

\vskip 0.7ex
\hangindent=3em \hangafter=1
$S$ = ($ 1$,
$ 1$,
$ \sqrt{2}$,
$ \frac{1-\sqrt{5}}{2}$,
$ \frac{1-\sqrt{5}}{2}$,
$ -c^{3}_{40}
+c^{7}_{40}
$;\ \ 
$ 1$,
$ -\sqrt{2}$,
$ \frac{1-\sqrt{5}}{2}$,
$ \frac{1-\sqrt{5}}{2}$,
$ c^{3}_{40}
-c^{7}_{40}
$;\ \ 
$0$,
$ -c^{3}_{40}
+c^{7}_{40}
$,
$ c^{3}_{40}
-c^{7}_{40}
$,
$0$;\ \ 
$ -1$,
$ -1$,
$ -\sqrt{2}$;\ \ 
$ -1$,
$ \sqrt{2}$;\ \ 
$0$)

Not pseudo-unitary. 

\vskip 1ex 
\color{grey}

\noindent159. ind = $(6;80
)_{2}^{11}$:\ \ 
$d_i$ = ($1.0$,
$1.0$,
$1.618$,
$1.618$,
$-1.414$,
$-2.288$) 

\vskip 0.7ex
\hangindent=3em \hangafter=1
$D^2=$ 14.472 = 
 $10+2\sqrt{5}$

\vskip 0.7ex
\hangindent=3em \hangafter=1
$T = ( 0,
\frac{1}{2},
\frac{2}{5},
\frac{9}{10},
\frac{1}{16},
\frac{37}{80} )
$,

\vskip 0.7ex
\hangindent=3em \hangafter=1
$S$ = ($ 1$,
$ 1$,
$ \frac{1+\sqrt{5}}{2}$,
$ \frac{1+\sqrt{5}}{2}$,
$ -\sqrt{2}$,
$ -c^{3}_{40}
-c^{5}_{40}
+c^{7}_{40}
$;\ \ 
$ 1$,
$ \frac{1+\sqrt{5}}{2}$,
$ \frac{1+\sqrt{5}}{2}$,
$ \sqrt{2}$,
$ c^{3}_{40}
+c^{5}_{40}
-c^{7}_{40}
$;\ \ 
$ -1$,
$ -1$,
$ -c^{3}_{40}
-c^{5}_{40}
+c^{7}_{40}
$,
$ \sqrt{2}$;\ \ 
$ -1$,
$ c^{3}_{40}
+c^{5}_{40}
-c^{7}_{40}
$,
$ -\sqrt{2}$;\ \ 
$0$,
$0$;\ \ 
$0$)

Pseudo-unitary $\sim$  
$(6;80
)_{1}^{1}$

\vskip 1ex 
\color{grey}

\noindent160. ind = $(6;80
)_{2}^{61}$:\ \ 
$d_i$ = ($1.0$,
$1.0$,
$1.618$,
$1.618$,
$-1.414$,
$-2.288$) 

\vskip 0.7ex
\hangindent=3em \hangafter=1
$D^2=$ 14.472 = 
 $10+2\sqrt{5}$

\vskip 0.7ex
\hangindent=3em \hangafter=1
$T = ( 0,
\frac{1}{2},
\frac{2}{5},
\frac{9}{10},
\frac{7}{16},
\frac{67}{80} )
$,

\vskip 0.7ex
\hangindent=3em \hangafter=1
$S$ = ($ 1$,
$ 1$,
$ \frac{1+\sqrt{5}}{2}$,
$ \frac{1+\sqrt{5}}{2}$,
$ -\sqrt{2}$,
$ -c^{3}_{40}
-c^{5}_{40}
+c^{7}_{40}
$;\ \ 
$ 1$,
$ \frac{1+\sqrt{5}}{2}$,
$ \frac{1+\sqrt{5}}{2}$,
$ \sqrt{2}$,
$ c^{3}_{40}
+c^{5}_{40}
-c^{7}_{40}
$;\ \ 
$ -1$,
$ -1$,
$ -c^{3}_{40}
-c^{5}_{40}
+c^{7}_{40}
$,
$ \sqrt{2}$;\ \ 
$ -1$,
$ c^{3}_{40}
+c^{5}_{40}
-c^{7}_{40}
$,
$ -\sqrt{2}$;\ \ 
$0$,
$0$;\ \ 
$0$)

Pseudo-unitary $\sim$  
$(6;80
)_{1}^{71}$

\vskip 1ex 
\color{grey}

\noindent161. ind = $(6;80
)_{2}^{51}$:\ \ 
$d_i$ = ($1.0$,
$1.0$,
$1.618$,
$1.618$,
$-1.414$,
$-2.288$) 

\vskip 0.7ex
\hangindent=3em \hangafter=1
$D^2=$ 14.472 = 
 $10+2\sqrt{5}$

\vskip 0.7ex
\hangindent=3em \hangafter=1
$T = ( 0,
\frac{1}{2},
\frac{2}{5},
\frac{9}{10},
\frac{9}{16},
\frac{77}{80} )
$,

\vskip 0.7ex
\hangindent=3em \hangafter=1
$S$ = ($ 1$,
$ 1$,
$ \frac{1+\sqrt{5}}{2}$,
$ \frac{1+\sqrt{5}}{2}$,
$ -\sqrt{2}$,
$ -c^{3}_{40}
-c^{5}_{40}
+c^{7}_{40}
$;\ \ 
$ 1$,
$ \frac{1+\sqrt{5}}{2}$,
$ \frac{1+\sqrt{5}}{2}$,
$ \sqrt{2}$,
$ c^{3}_{40}
+c^{5}_{40}
-c^{7}_{40}
$;\ \ 
$ -1$,
$ -1$,
$ -c^{3}_{40}
-c^{5}_{40}
+c^{7}_{40}
$,
$ \sqrt{2}$;\ \ 
$ -1$,
$ c^{3}_{40}
+c^{5}_{40}
-c^{7}_{40}
$,
$ -\sqrt{2}$;\ \ 
$0$,
$0$;\ \ 
$0$)

Pseudo-unitary $\sim$  
$(6;80
)_{1}^{41}$

\vskip 1ex 
\color{grey}

\noindent162. ind = $(6;80
)_{2}^{21}$:\ \ 
$d_i$ = ($1.0$,
$1.0$,
$1.618$,
$1.618$,
$-1.414$,
$-2.288$) 

\vskip 0.7ex
\hangindent=3em \hangafter=1
$D^2=$ 14.472 = 
 $10+2\sqrt{5}$

\vskip 0.7ex
\hangindent=3em \hangafter=1
$T = ( 0,
\frac{1}{2},
\frac{2}{5},
\frac{9}{10},
\frac{15}{16},
\frac{27}{80} )
$,

\vskip 0.7ex
\hangindent=3em \hangafter=1
$S$ = ($ 1$,
$ 1$,
$ \frac{1+\sqrt{5}}{2}$,
$ \frac{1+\sqrt{5}}{2}$,
$ -\sqrt{2}$,
$ -c^{3}_{40}
-c^{5}_{40}
+c^{7}_{40}
$;\ \ 
$ 1$,
$ \frac{1+\sqrt{5}}{2}$,
$ \frac{1+\sqrt{5}}{2}$,
$ \sqrt{2}$,
$ c^{3}_{40}
+c^{5}_{40}
-c^{7}_{40}
$;\ \ 
$ -1$,
$ -1$,
$ -c^{3}_{40}
-c^{5}_{40}
+c^{7}_{40}
$,
$ \sqrt{2}$;\ \ 
$ -1$,
$ c^{3}_{40}
+c^{5}_{40}
-c^{7}_{40}
$,
$ -\sqrt{2}$;\ \ 
$0$,
$0$;\ \ 
$0$)

Pseudo-unitary $\sim$  
$(6;80
)_{1}^{31}$

\vskip 1ex 
\color{grey}

\noindent163. ind = $(6;80
)_{2}^{59}$:\ \ 
$d_i$ = ($1.0$,
$1.0$,
$1.618$,
$1.618$,
$-1.414$,
$-2.288$) 

\vskip 0.7ex
\hangindent=3em \hangafter=1
$D^2=$ 14.472 = 
 $10+2\sqrt{5}$

\vskip 0.7ex
\hangindent=3em \hangafter=1
$T = ( 0,
\frac{1}{2},
\frac{3}{5},
\frac{1}{10},
\frac{1}{16},
\frac{53}{80} )
$,

\vskip 0.7ex
\hangindent=3em \hangafter=1
$S$ = ($ 1$,
$ 1$,
$ \frac{1+\sqrt{5}}{2}$,
$ \frac{1+\sqrt{5}}{2}$,
$ -\sqrt{2}$,
$ -c^{3}_{40}
-c^{5}_{40}
+c^{7}_{40}
$;\ \ 
$ 1$,
$ \frac{1+\sqrt{5}}{2}$,
$ \frac{1+\sqrt{5}}{2}$,
$ \sqrt{2}$,
$ c^{3}_{40}
+c^{5}_{40}
-c^{7}_{40}
$;\ \ 
$ -1$,
$ -1$,
$ -c^{3}_{40}
-c^{5}_{40}
+c^{7}_{40}
$,
$ \sqrt{2}$;\ \ 
$ -1$,
$ c^{3}_{40}
+c^{5}_{40}
-c^{7}_{40}
$,
$ -\sqrt{2}$;\ \ 
$0$,
$0$;\ \ 
$0$)

Pseudo-unitary $\sim$  
$(6;80
)_{1}^{49}$

\vskip 1ex 
\color{grey}

\noindent164. ind = $(6;80
)_{2}^{29}$:\ \ 
$d_i$ = ($1.0$,
$1.0$,
$1.618$,
$1.618$,
$-1.414$,
$-2.288$) 

\vskip 0.7ex
\hangindent=3em \hangafter=1
$D^2=$ 14.472 = 
 $10+2\sqrt{5}$

\vskip 0.7ex
\hangindent=3em \hangafter=1
$T = ( 0,
\frac{1}{2},
\frac{3}{5},
\frac{1}{10},
\frac{7}{16},
\frac{3}{80} )
$,

\vskip 0.7ex
\hangindent=3em \hangafter=1
$S$ = ($ 1$,
$ 1$,
$ \frac{1+\sqrt{5}}{2}$,
$ \frac{1+\sqrt{5}}{2}$,
$ -\sqrt{2}$,
$ -c^{3}_{40}
-c^{5}_{40}
+c^{7}_{40}
$;\ \ 
$ 1$,
$ \frac{1+\sqrt{5}}{2}$,
$ \frac{1+\sqrt{5}}{2}$,
$ \sqrt{2}$,
$ c^{3}_{40}
+c^{5}_{40}
-c^{7}_{40}
$;\ \ 
$ -1$,
$ -1$,
$ -c^{3}_{40}
-c^{5}_{40}
+c^{7}_{40}
$,
$ \sqrt{2}$;\ \ 
$ -1$,
$ c^{3}_{40}
+c^{5}_{40}
-c^{7}_{40}
$,
$ -\sqrt{2}$;\ \ 
$0$,
$0$;\ \ 
$0$)

Pseudo-unitary $\sim$  
$(6;80
)_{1}^{39}$

\vskip 1ex 
\color{grey}

\noindent165. ind = $(6;80
)_{2}^{19}$:\ \ 
$d_i$ = ($1.0$,
$1.0$,
$1.618$,
$1.618$,
$-1.414$,
$-2.288$) 

\vskip 0.7ex
\hangindent=3em \hangafter=1
$D^2=$ 14.472 = 
 $10+2\sqrt{5}$

\vskip 0.7ex
\hangindent=3em \hangafter=1
$T = ( 0,
\frac{1}{2},
\frac{3}{5},
\frac{1}{10},
\frac{9}{16},
\frac{13}{80} )
$,

\vskip 0.7ex
\hangindent=3em \hangafter=1
$S$ = ($ 1$,
$ 1$,
$ \frac{1+\sqrt{5}}{2}$,
$ \frac{1+\sqrt{5}}{2}$,
$ -\sqrt{2}$,
$ -c^{3}_{40}
-c^{5}_{40}
+c^{7}_{40}
$;\ \ 
$ 1$,
$ \frac{1+\sqrt{5}}{2}$,
$ \frac{1+\sqrt{5}}{2}$,
$ \sqrt{2}$,
$ c^{3}_{40}
+c^{5}_{40}
-c^{7}_{40}
$;\ \ 
$ -1$,
$ -1$,
$ -c^{3}_{40}
-c^{5}_{40}
+c^{7}_{40}
$,
$ \sqrt{2}$;\ \ 
$ -1$,
$ c^{3}_{40}
+c^{5}_{40}
-c^{7}_{40}
$,
$ -\sqrt{2}$;\ \ 
$0$,
$0$;\ \ 
$0$)

Pseudo-unitary $\sim$  
$(6;80
)_{1}^{9}$

\vskip 1ex 
\color{grey}

\noindent166. ind = $(6;80
)_{2}^{69}$:\ \ 
$d_i$ = ($1.0$,
$1.0$,
$1.618$,
$1.618$,
$-1.414$,
$-2.288$) 

\vskip 0.7ex
\hangindent=3em \hangafter=1
$D^2=$ 14.472 = 
 $10+2\sqrt{5}$

\vskip 0.7ex
\hangindent=3em \hangafter=1
$T = ( 0,
\frac{1}{2},
\frac{3}{5},
\frac{1}{10},
\frac{15}{16},
\frac{43}{80} )
$,

\vskip 0.7ex
\hangindent=3em \hangafter=1
$S$ = ($ 1$,
$ 1$,
$ \frac{1+\sqrt{5}}{2}$,
$ \frac{1+\sqrt{5}}{2}$,
$ -\sqrt{2}$,
$ -c^{3}_{40}
-c^{5}_{40}
+c^{7}_{40}
$;\ \ 
$ 1$,
$ \frac{1+\sqrt{5}}{2}$,
$ \frac{1+\sqrt{5}}{2}$,
$ \sqrt{2}$,
$ c^{3}_{40}
+c^{5}_{40}
-c^{7}_{40}
$;\ \ 
$ -1$,
$ -1$,
$ -c^{3}_{40}
-c^{5}_{40}
+c^{7}_{40}
$,
$ \sqrt{2}$;\ \ 
$ -1$,
$ c^{3}_{40}
+c^{5}_{40}
-c^{7}_{40}
$,
$ -\sqrt{2}$;\ \ 
$0$,
$0$;\ \ 
$0$)

Pseudo-unitary $\sim$  
$(6;80
)_{1}^{79}$

\vskip 1ex 
\color{grey}

\noindent167. ind = $(6;80
)_{2}^{73}$:\ \ 
$d_i$ = ($1.0$,
$1.0$,
$1.414$,
$-0.618$,
$-0.618$,
$-0.874$) 

\vskip 0.7ex
\hangindent=3em \hangafter=1
$D^2=$ 5.527 = 
 $10-2\sqrt{5}$

\vskip 0.7ex
\hangindent=3em \hangafter=1
$T = ( 0,
\frac{1}{2},
\frac{11}{16},
\frac{1}{5},
\frac{7}{10},
\frac{71}{80} )
$,

\vskip 0.7ex
\hangindent=3em \hangafter=1
$S$ = ($ 1$,
$ 1$,
$ \sqrt{2}$,
$ \frac{1-\sqrt{5}}{2}$,
$ \frac{1-\sqrt{5}}{2}$,
$ -c^{3}_{40}
+c^{7}_{40}
$;\ \ 
$ 1$,
$ -\sqrt{2}$,
$ \frac{1-\sqrt{5}}{2}$,
$ \frac{1-\sqrt{5}}{2}$,
$ c^{3}_{40}
-c^{7}_{40}
$;\ \ 
$0$,
$ -c^{3}_{40}
+c^{7}_{40}
$,
$ c^{3}_{40}
-c^{7}_{40}
$,
$0$;\ \ 
$ -1$,
$ -1$,
$ -\sqrt{2}$;\ \ 
$ -1$,
$ \sqrt{2}$;\ \ 
$0$)

Not pseudo-unitary. 

\vskip 1ex 
\color{grey}

\noindent168. ind = $(6;80
)_{2}^{57}$:\ \ 
$d_i$ = ($1.0$,
$1.0$,
$1.414$,
$-0.618$,
$-0.618$,
$-0.874$) 

\vskip 0.7ex
\hangindent=3em \hangafter=1
$D^2=$ 5.527 = 
 $10-2\sqrt{5}$

\vskip 0.7ex
\hangindent=3em \hangafter=1
$T = ( 0,
\frac{1}{2},
\frac{11}{16},
\frac{4}{5},
\frac{3}{10},
\frac{39}{80} )
$,

\vskip 0.7ex
\hangindent=3em \hangafter=1
$S$ = ($ 1$,
$ 1$,
$ \sqrt{2}$,
$ \frac{1-\sqrt{5}}{2}$,
$ \frac{1-\sqrt{5}}{2}$,
$ -c^{3}_{40}
+c^{7}_{40}
$;\ \ 
$ 1$,
$ -\sqrt{2}$,
$ \frac{1-\sqrt{5}}{2}$,
$ \frac{1-\sqrt{5}}{2}$,
$ c^{3}_{40}
-c^{7}_{40}
$;\ \ 
$0$,
$ -c^{3}_{40}
+c^{7}_{40}
$,
$ c^{3}_{40}
-c^{7}_{40}
$,
$0$;\ \ 
$ -1$,
$ -1$,
$ -\sqrt{2}$;\ \ 
$ -1$,
$ \sqrt{2}$;\ \ 
$0$)

Not pseudo-unitary. 

\vskip 1ex 
\color{grey}

\noindent169. ind = $(6;80
)_{2}^{63}$:\ \ 
$d_i$ = ($1.0$,
$1.0$,
$1.414$,
$-0.618$,
$-0.618$,
$-0.874$) 

\vskip 0.7ex
\hangindent=3em \hangafter=1
$D^2=$ 5.527 = 
 $10-2\sqrt{5}$

\vskip 0.7ex
\hangindent=3em \hangafter=1
$T = ( 0,
\frac{1}{2},
\frac{13}{16},
\frac{1}{5},
\frac{7}{10},
\frac{1}{80} )
$,

\vskip 0.7ex
\hangindent=3em \hangafter=1
$S$ = ($ 1$,
$ 1$,
$ \sqrt{2}$,
$ \frac{1-\sqrt{5}}{2}$,
$ \frac{1-\sqrt{5}}{2}$,
$ -c^{3}_{40}
+c^{7}_{40}
$;\ \ 
$ 1$,
$ -\sqrt{2}$,
$ \frac{1-\sqrt{5}}{2}$,
$ \frac{1-\sqrt{5}}{2}$,
$ c^{3}_{40}
-c^{7}_{40}
$;\ \ 
$0$,
$ -c^{3}_{40}
+c^{7}_{40}
$,
$ c^{3}_{40}
-c^{7}_{40}
$,
$0$;\ \ 
$ -1$,
$ -1$,
$ -\sqrt{2}$;\ \ 
$ -1$,
$ \sqrt{2}$;\ \ 
$0$)

Not pseudo-unitary. 

\vskip 1ex 
\color{grey}

\noindent170. ind = $(6;80
)_{2}^{47}$:\ \ 
$d_i$ = ($1.0$,
$1.0$,
$1.414$,
$-0.618$,
$-0.618$,
$-0.874$) 

\vskip 0.7ex
\hangindent=3em \hangafter=1
$D^2=$ 5.527 = 
 $10-2\sqrt{5}$

\vskip 0.7ex
\hangindent=3em \hangafter=1
$T = ( 0,
\frac{1}{2},
\frac{13}{16},
\frac{4}{5},
\frac{3}{10},
\frac{49}{80} )
$,

\vskip 0.7ex
\hangindent=3em \hangafter=1
$S$ = ($ 1$,
$ 1$,
$ \sqrt{2}$,
$ \frac{1-\sqrt{5}}{2}$,
$ \frac{1-\sqrt{5}}{2}$,
$ -c^{3}_{40}
+c^{7}_{40}
$;\ \ 
$ 1$,
$ -\sqrt{2}$,
$ \frac{1-\sqrt{5}}{2}$,
$ \frac{1-\sqrt{5}}{2}$,
$ c^{3}_{40}
-c^{7}_{40}
$;\ \ 
$0$,
$ -c^{3}_{40}
+c^{7}_{40}
$,
$ c^{3}_{40}
-c^{7}_{40}
$,
$0$;\ \ 
$ -1$,
$ -1$,
$ -\sqrt{2}$;\ \ 
$ -1$,
$ \sqrt{2}$;\ \ 
$0$)

Not pseudo-unitary. 

\vskip 1ex 
\color{grey}

\noindent171. ind = $(6;80
)_{2}^{13}$:\ \ 
$d_i$ = ($1.0$,
$0.874$,
$1.0$,
$-0.618$,
$-0.618$,
$-1.414$) 

\vskip 0.7ex
\hangindent=3em \hangafter=1
$D^2=$ 5.527 = 
 $10-2\sqrt{5}$

\vskip 0.7ex
\hangindent=3em \hangafter=1
$T = ( 0,
\frac{51}{80},
\frac{1}{2},
\frac{1}{5},
\frac{7}{10},
\frac{7}{16} )
$,

\vskip 0.7ex
\hangindent=3em \hangafter=1
$S$ = ($ 1$,
$ c^{3}_{40}
-c^{7}_{40}
$,
$ 1$,
$ \frac{1-\sqrt{5}}{2}$,
$ \frac{1-\sqrt{5}}{2}$,
$ -\sqrt{2}$;\ \ 
$0$,
$ -c^{3}_{40}
+c^{7}_{40}
$,
$ \sqrt{2}$,
$ -\sqrt{2}$,
$0$;\ \ 
$ 1$,
$ \frac{1-\sqrt{5}}{2}$,
$ \frac{1-\sqrt{5}}{2}$,
$ \sqrt{2}$;\ \ 
$ -1$,
$ -1$,
$ c^{3}_{40}
-c^{7}_{40}
$;\ \ 
$ -1$,
$ -c^{3}_{40}
+c^{7}_{40}
$;\ \ 
$0$)

Not pseudo-unitary. 

\vskip 1ex 
\color{grey}

\noindent172. ind = $(6;80
)_{2}^{37}$:\ \ 
$d_i$ = ($1.0$,
$0.874$,
$1.0$,
$-0.618$,
$-0.618$,
$-1.414$) 

\vskip 0.7ex
\hangindent=3em \hangafter=1
$D^2=$ 5.527 = 
 $10-2\sqrt{5}$

\vskip 0.7ex
\hangindent=3em \hangafter=1
$T = ( 0,
\frac{59}{80},
\frac{1}{2},
\frac{4}{5},
\frac{3}{10},
\frac{15}{16} )
$,

\vskip 0.7ex
\hangindent=3em \hangafter=1
$S$ = ($ 1$,
$ c^{3}_{40}
-c^{7}_{40}
$,
$ 1$,
$ \frac{1-\sqrt{5}}{2}$,
$ \frac{1-\sqrt{5}}{2}$,
$ -\sqrt{2}$;\ \ 
$0$,
$ -c^{3}_{40}
+c^{7}_{40}
$,
$ \sqrt{2}$,
$ -\sqrt{2}$,
$0$;\ \ 
$ 1$,
$ \frac{1-\sqrt{5}}{2}$,
$ \frac{1-\sqrt{5}}{2}$,
$ \sqrt{2}$;\ \ 
$ -1$,
$ -1$,
$ c^{3}_{40}
-c^{7}_{40}
$;\ \ 
$ -1$,
$ -c^{3}_{40}
+c^{7}_{40}
$;\ \ 
$0$)

Not pseudo-unitary. 

\vskip 1ex 
\color{grey}

\noindent173. ind = $(6;80
)_{2}^{3}$:\ \ 
$d_i$ = ($1.0$,
$0.874$,
$1.0$,
$-0.618$,
$-0.618$,
$-1.414$) 

\vskip 0.7ex
\hangindent=3em \hangafter=1
$D^2=$ 5.527 = 
 $10-2\sqrt{5}$

\vskip 0.7ex
\hangindent=3em \hangafter=1
$T = ( 0,
\frac{61}{80},
\frac{1}{2},
\frac{1}{5},
\frac{7}{10},
\frac{9}{16} )
$,

\vskip 0.7ex
\hangindent=3em \hangafter=1
$S$ = ($ 1$,
$ c^{3}_{40}
-c^{7}_{40}
$,
$ 1$,
$ \frac{1-\sqrt{5}}{2}$,
$ \frac{1-\sqrt{5}}{2}$,
$ -\sqrt{2}$;\ \ 
$0$,
$ -c^{3}_{40}
+c^{7}_{40}
$,
$ \sqrt{2}$,
$ -\sqrt{2}$,
$0$;\ \ 
$ 1$,
$ \frac{1-\sqrt{5}}{2}$,
$ \frac{1-\sqrt{5}}{2}$,
$ \sqrt{2}$;\ \ 
$ -1$,
$ -1$,
$ c^{3}_{40}
-c^{7}_{40}
$;\ \ 
$ -1$,
$ -c^{3}_{40}
+c^{7}_{40}
$;\ \ 
$0$)

Not pseudo-unitary. 

\vskip 1ex 
\color{grey}

\noindent174. ind = $(6;80
)_{2}^{27}$:\ \ 
$d_i$ = ($1.0$,
$0.874$,
$1.0$,
$-0.618$,
$-0.618$,
$-1.414$) 

\vskip 0.7ex
\hangindent=3em \hangafter=1
$D^2=$ 5.527 = 
 $10-2\sqrt{5}$

\vskip 0.7ex
\hangindent=3em \hangafter=1
$T = ( 0,
\frac{69}{80},
\frac{1}{2},
\frac{4}{5},
\frac{3}{10},
\frac{1}{16} )
$,

\vskip 0.7ex
\hangindent=3em \hangafter=1
$S$ = ($ 1$,
$ c^{3}_{40}
-c^{7}_{40}
$,
$ 1$,
$ \frac{1-\sqrt{5}}{2}$,
$ \frac{1-\sqrt{5}}{2}$,
$ -\sqrt{2}$;\ \ 
$0$,
$ -c^{3}_{40}
+c^{7}_{40}
$,
$ \sqrt{2}$,
$ -\sqrt{2}$,
$0$;\ \ 
$ 1$,
$ \frac{1-\sqrt{5}}{2}$,
$ \frac{1-\sqrt{5}}{2}$,
$ \sqrt{2}$;\ \ 
$ -1$,
$ -1$,
$ c^{3}_{40}
-c^{7}_{40}
$;\ \ 
$ -1$,
$ -c^{3}_{40}
+c^{7}_{40}
$;\ \ 
$0$)

Not pseudo-unitary. 

\vskip 1ex 

 \color{black} \vskip 2ex

\

\section{A list of (almost) all dimension-6 $\SL$ representations}
\label{Section4}

In this section, we give a list of all dimension-6 irrep-sum representations of
$\SL$: $\rho=\rho_1\oplus \rho_2\oplus \cdots$ where $\rho_a$'s are irreducible
representations. We use the same notations as in Appendix B.2.  Below each
representation, we indicate if it passes all the conditions on Propositions
B.3, B.4, and B.5.  If a representation fails to satisfy some of those
conditions, we indicate which condition it fails to satisfy.
 
We note that a representation $\rho$ satisfies Propositions B.3, B.4, and B.5
iff $\rho\otimes \chi$ satisfies Propositions B.3, B.4, and B.5, where $\chi$
is any 1-dimensional representations.  The representations connected by
tensoring 1-dimensional representations form a so called T-orbit.  We also note
that different representations in the same T-orbit produce the same set of
modular data.  So we need to consider only one representation in each T-orbit
and skip other representations.

To determine which representations to skip, we consider a irreducible
representation $\rho_{a}$ in the direct sum decomposition of $\rho$ that has a
minimal dimension.  (If two irreducible representations have the same minimal
dimension, we just pick one of them.) We denote this minimal-dimension
irreducible representation as $\rho_m$.  Such an irreducible representation is
a tensor product of  irreducible representations of prime-power levels.  Those
prime-power level irreducible representations are labeled by $d_{l,k}^{a,m}$
(see Appendix A).  Thus
\begin{align}
\label{rhom}
 \rho_m = d_{l,k}^{a,m} \otimes \td d_{\td l,\td k}^{\td a,\td m} \otimes \cdots
\end{align}
We note that $\rho_m$ differ from
\begin{align}
\label{trhom}
 \td \rho_m = d_{l,k}^{a,0} \otimes \td d_{\td l,\td k}^{\td a,0} \otimes \cdots
\end{align}
via tensoring a 1-dimensional representations. Thus we can skip all the
representations, except those whose minimal-dimension irreducible
representation has a form \eqref{trhom}.  In fact, to be more precise, we can
skip all the representations, except those whose minimal-dimension irreducible
representation has a form \eqref{rhom}, with $m+\td m + \cdots = 0$ mod 12.  In
this following table, we do not list those skipped representations.

%we introduce the \textbf{index} of
%a irrep-sum representation $\rho=\rho_1\oplus\rho_2\oplus \cdots$ as $[l,
%[l_1,l_2,\cdots], [\td s_0, \td s_1, \cdots,\td s_{r-1}] ]$.  Here we have
%ordered the irreducible representations, $\rho_1,\rho_2, \cdots$, such that
%dim($\rho_1$) $\geq$ dim($\rho_2$) $\geq \cdots$.  $l_a$ is the level of
%$\rho_a$ and $l$ is the lest common multiple of $l_a$.  $\td s_i \in [0,1)$ and
%$\ee^{\ii 2\pi \td s_i}$ are the eigenvalues of $\rho(\ft)$. $\td s_i$'s in
%$[\td s_0, \td s_1, \cdots,\td s_{r-1}]$ are ordered according to the
%denominators of the rational number $\td s_i$.  If the denominators are the
%same, then they are ordered according to the value of $\td s_i$.  We also
%introduce the \textbf{minimal index} of a representation, as the smallest index
%among the representations in the same T-orbit.  We use the following scheme to
%compare two indices: we first compare the two arrays by their first terms to
%determine which array is smaller.  If the first terms are the same, then we
%compare the second terms, {\it etc}.  This way, we can skip the representations
%whose minimal index is less then its index.

To reduce the length of the list further, we also omit the following three
types of representations from the list: (1) the representations that can be
written as $\rho =\rho'\oplus \rho''$ with non-overlapping spectra of
$\rho'(\ft)$ and $\rho''(\ft)$.  (2) the representations that can be written as
$\rho =n\rho'$ for any integer $n > 1$ and any representation $\rho'$ such that
all the eigenvalues of $\rho'(\ft)$ are non-degenerate (see Proposition B.3 (4)
eqn.  (B.13) ).  (3) the representations whose $\pord(\rho(\ft)) \in
\{2,3,4,6\}$ (see Proposition B.3 (6) eqn. (B.13) ).  (4) the representations
whose $m$-index is zero.  A representation is a direct sum of several
irreducible representations.  We consider one of the irreducible
representations with highest dimension, which has a form $
(d_1)_{l_1,k_1}^{a_1,m_1} \otimes (d_2)_{l_2,k_2}^{a_2,m_2} \otimes \cdots $.
The $m$-index is $m=m_1+m_2+\cdots$ mod 12.  The first two types contain only
failed representations and the third type contain only failed representations
and integral modular data.

The entries of the list, $\rho=\rho_1\oplus\rho_2\oplus \cdots$, are ordered
according to $[[r_1,r_2,\cdots]$, $ l$, $ [l_1,l_2,\cdots]$, $ [\td s_0, \td
s_1, \cdots,\td s_{r-1}] ]$.  Here $r_a$ is the dimension of $\rho_a$.  For
more details and notations, see Appendix B.2.

\

\noindent 1: (dims,levels) = $(3 , 
1 , 
1 , 
1;5,
1,
1,
1
)$,
irreps = $3_{5}^{1}\oplus
1_{1}^{1}\oplus
1_{1}^{1}\oplus
1_{1}^{1}$,
pord$(\rho_\text{isum}(\mathfrak{t})) = 5$,

\vskip 0.7ex
\hangindent=5.5em \hangafter=1
{\white .}\hskip 1em $\rho_\text{isum}(\mathfrak{t})$ =
 $( 0,
\frac{1}{5},
\frac{4}{5} )
\oplus
( 0 )
\oplus
( 0 )
\oplus
( 0 )
$,

\vskip 0.7ex
\hangindent=5.5em \hangafter=1
{\white .}\hskip 1em $\rho_\text{isum}(\mathfrak{s})$ =
($\sqrt{\frac{1}{5}}$,
$-\sqrt{\frac{2}{5}}$,
$-\sqrt{\frac{2}{5}}$;
$-\frac{5+\sqrt{5}}{10}$,
$\frac{5-\sqrt{5}}{10}$;
$-\frac{5+\sqrt{5}}{10}$)
 $\oplus$
($1$)
 $\oplus$
($1$)
 $\oplus$
($1$)

Fail:
dims is $(d,1,1,\cdots)$, $d=(\mathrm{ord}(T)+1)/2$, self-dual, ...
 Prop. B.3 (3)

 \ \color{black}

\noindent 2: (dims,levels) = $(3 , 
1 , 
1 , 
1;5,
1,
1,
1
)$,
irreps = $3_{5}^{3}\oplus
1_{1}^{1}\oplus
1_{1}^{1}\oplus
1_{1}^{1}$,
pord$(\rho_\text{isum}(\mathfrak{t})) = 5$,

\vskip 0.7ex
\hangindent=5.5em \hangafter=1
{\white .}\hskip 1em $\rho_\text{isum}(\mathfrak{t})$ =
 $( 0,
\frac{2}{5},
\frac{3}{5} )
\oplus
( 0 )
\oplus
( 0 )
\oplus
( 0 )
$,

\vskip 0.7ex
\hangindent=5.5em \hangafter=1
{\white .}\hskip 1em $\rho_\text{isum}(\mathfrak{s})$ =
($-\sqrt{\frac{1}{5}}$,
$-\sqrt{\frac{2}{5}}$,
$-\sqrt{\frac{2}{5}}$;
$\frac{-5+\sqrt{5}}{10}$,
$\frac{5+\sqrt{5}}{10}$;
$\frac{-5+\sqrt{5}}{10}$)
 $\oplus$
($1$)
 $\oplus$
($1$)
 $\oplus$
($1$)

Fail:
dims is $(d,1,1,\cdots)$, $d=(\mathrm{ord}(T)+1)/2$, self-dual, ...
 Prop. B.3 (3)

 \ \color{black}

\noindent 3: (dims,levels) = $(3 , 
1 , 
1 , 
1;8,
1,
1,
1
)$,
irreps = $3_{8}^{1,0}\oplus
1_{1}^{1}\oplus
1_{1}^{1}\oplus
1_{1}^{1}$,
pord$(\rho_\text{isum}(\mathfrak{t})) = 8$,

\vskip 0.7ex
\hangindent=5.5em \hangafter=1
{\white .}\hskip 1em $\rho_\text{isum}(\mathfrak{t})$ =
 $( 0,
\frac{1}{8},
\frac{5}{8} )
\oplus
( 0 )
\oplus
( 0 )
\oplus
( 0 )
$,

\vskip 0.7ex
\hangindent=5.5em \hangafter=1
{\white .}\hskip 1em $\rho_\text{isum}(\mathfrak{s})$ =
$\mathrm{i}$($0$,
$\sqrt{\frac{1}{2}}$,
$\sqrt{\frac{1}{2}}$;\ \ 
$-\frac{1}{2}$,
$\frac{1}{2}$;\ \ 
$-\frac{1}{2}$)
 $\oplus$
($1$)
 $\oplus$
($1$)
 $\oplus$
($1$)

Fail:
for $\rho = \rho_1+l\chi, ...,
 (\rho_1(\mathfrak s)/\chi(\mathfrak s))^2\neq$ id. Prop. B.3 (2)

 \ \color{black}

\noindent 4: (dims,levels) = $(3 , 
1 , 
1 , 
1;8,
1,
1,
1
)$,
irreps = $3_{8}^{3,0}\oplus
1_{1}^{1}\oplus
1_{1}^{1}\oplus
1_{1}^{1}$,
pord$(\rho_\text{isum}(\mathfrak{t})) = 8$,

\vskip 0.7ex
\hangindent=5.5em \hangafter=1
{\white .}\hskip 1em $\rho_\text{isum}(\mathfrak{t})$ =
 $( 0,
\frac{3}{8},
\frac{7}{8} )
\oplus
( 0 )
\oplus
( 0 )
\oplus
( 0 )
$,

\vskip 0.7ex
\hangindent=5.5em \hangafter=1
{\white .}\hskip 1em $\rho_\text{isum}(\mathfrak{s})$ =
$\mathrm{i}$($0$,
$\sqrt{\frac{1}{2}}$,
$\sqrt{\frac{1}{2}}$;\ \ 
$\frac{1}{2}$,
$-\frac{1}{2}$;\ \ 
$\frac{1}{2}$)
 $\oplus$
($1$)
 $\oplus$
($1$)
 $\oplus$
($1$)

Fail:
for $\rho = \rho_1+l\chi, ...,
 (\rho_1(\mathfrak s)/\chi(\mathfrak s))^2\neq$ id. Prop. B.3 (2)

 \ \color{black}

\noindent 5: (dims,levels) = $(3 , 
1 , 
1 , 
1;10,
2,
2,
2
)$,
irreps = $3_{5}^{3}
\hskip -1.5pt \otimes \hskip -1.5pt
1_{2}^{1,0}\oplus
1_{2}^{1,0}\oplus
1_{2}^{1,0}\oplus
1_{2}^{1,0}$,
pord$(\rho_\text{isum}(\mathfrak{t})) = 5$,

\vskip 0.7ex
\hangindent=5.5em \hangafter=1
{\white .}\hskip 1em $\rho_\text{isum}(\mathfrak{t})$ =
 $( \frac{1}{2},
\frac{1}{10},
\frac{9}{10} )
\oplus
( \frac{1}{2} )
\oplus
( \frac{1}{2} )
\oplus
( \frac{1}{2} )
$,

\vskip 0.7ex
\hangindent=5.5em \hangafter=1
{\white .}\hskip 1em $\rho_\text{isum}(\mathfrak{s})$ =
($\sqrt{\frac{1}{5}}$,
$-\sqrt{\frac{2}{5}}$,
$-\sqrt{\frac{2}{5}}$;
$\frac{5-\sqrt{5}}{10}$,
$-\frac{5+\sqrt{5}}{10}$;
$\frac{5-\sqrt{5}}{10}$)
 $\oplus$
($-1$)
 $\oplus$
($-1$)
 $\oplus$
($-1$)

Fail:
dims is $(d,1,1,\cdots)$, $d=(\mathrm{ord}(T)+1)/2$, self-dual, ...
 Prop. B.3 (3)

 \ \color{black}

\noindent 6: (dims,levels) = $(3 , 
1 , 
1 , 
1;10,
2,
2,
2
)$,
irreps = $3_{5}^{1}
\hskip -1.5pt \otimes \hskip -1.5pt
1_{2}^{1,0}\oplus
1_{2}^{1,0}\oplus
1_{2}^{1,0}\oplus
1_{2}^{1,0}$,
pord$(\rho_\text{isum}(\mathfrak{t})) = 5$,

\vskip 0.7ex
\hangindent=5.5em \hangafter=1
{\white .}\hskip 1em $\rho_\text{isum}(\mathfrak{t})$ =
 $( \frac{1}{2},
\frac{3}{10},
\frac{7}{10} )
\oplus
( \frac{1}{2} )
\oplus
( \frac{1}{2} )
\oplus
( \frac{1}{2} )
$,

\vskip 0.7ex
\hangindent=5.5em \hangafter=1
{\white .}\hskip 1em $\rho_\text{isum}(\mathfrak{s})$ =
($-\sqrt{\frac{1}{5}}$,
$-\sqrt{\frac{2}{5}}$,
$-\sqrt{\frac{2}{5}}$;
$\frac{5+\sqrt{5}}{10}$,
$\frac{-5+\sqrt{5}}{10}$;
$\frac{5+\sqrt{5}}{10}$)
 $\oplus$
($-1$)
 $\oplus$
($-1$)
 $\oplus$
($-1$)

Fail:
dims is $(d,1,1,\cdots)$, $d=(\mathrm{ord}(T)+1)/2$, self-dual, ...
 Prop. B.3 (3)

 \ \color{black}

\noindent 7: (dims,levels) = $(3 , 
1 , 
1 , 
1;15,
3,
3,
3
)$,
irreps = $3_{5}^{1}
\hskip -1.5pt \otimes \hskip -1.5pt
1_{3}^{1,0}\oplus
1_{3}^{1,0}\oplus
1_{3}^{1,0}\oplus
1_{3}^{1,0}$,
pord$(\rho_\text{isum}(\mathfrak{t})) = 5$,

\vskip 0.7ex
\hangindent=5.5em \hangafter=1
{\white .}\hskip 1em $\rho_\text{isum}(\mathfrak{t})$ =
 $( \frac{1}{3},
\frac{2}{15},
\frac{8}{15} )
\oplus
( \frac{1}{3} )
\oplus
( \frac{1}{3} )
\oplus
( \frac{1}{3} )
$,

\vskip 0.7ex
\hangindent=5.5em \hangafter=1
{\white .}\hskip 1em $\rho_\text{isum}(\mathfrak{s})$ =
($\sqrt{\frac{1}{5}}$,
$-\sqrt{\frac{2}{5}}$,
$-\sqrt{\frac{2}{5}}$;
$-\frac{5+\sqrt{5}}{10}$,
$\frac{5-\sqrt{5}}{10}$;
$-\frac{5+\sqrt{5}}{10}$)
 $\oplus$
($1$)
 $\oplus$
($1$)
 $\oplus$
($1$)

Fail:
dims is $(d,1,1,\cdots)$, $d=(\mathrm{ord}(T)+1)/2$, self-dual, ...
 Prop. B.3 (3)

 \ \color{black}

\noindent 8: (dims,levels) = $(3 , 
1 , 
1 , 
1;15,
3,
3,
3
)$,
irreps = $3_{5}^{3}
\hskip -1.5pt \otimes \hskip -1.5pt
1_{3}^{1,0}\oplus
1_{3}^{1,0}\oplus
1_{3}^{1,0}\oplus
1_{3}^{1,0}$,
pord$(\rho_\text{isum}(\mathfrak{t})) = 5$,

\vskip 0.7ex
\hangindent=5.5em \hangafter=1
{\white .}\hskip 1em $\rho_\text{isum}(\mathfrak{t})$ =
 $( \frac{1}{3},
\frac{11}{15},
\frac{14}{15} )
\oplus
( \frac{1}{3} )
\oplus
( \frac{1}{3} )
\oplus
( \frac{1}{3} )
$,

\vskip 0.7ex
\hangindent=5.5em \hangafter=1
{\white .}\hskip 1em $\rho_\text{isum}(\mathfrak{s})$ =
($-\sqrt{\frac{1}{5}}$,
$-\sqrt{\frac{2}{5}}$,
$-\sqrt{\frac{2}{5}}$;
$\frac{-5+\sqrt{5}}{10}$,
$\frac{5+\sqrt{5}}{10}$;
$\frac{-5+\sqrt{5}}{10}$)
 $\oplus$
($1$)
 $\oplus$
($1$)
 $\oplus$
($1$)

Fail:
dims is $(d,1,1,\cdots)$, $d=(\mathrm{ord}(T)+1)/2$, self-dual, ...
 Prop. B.3 (3)

 \ \color{black}

\noindent 9: (dims,levels) = $(3 , 
1 , 
1 , 
1;20,
4,
4,
4
)$,
irreps = $3_{5}^{1}
\hskip -1.5pt \otimes \hskip -1.5pt
1_{4}^{1,0}\oplus
1_{4}^{1,0}\oplus
1_{4}^{1,0}\oplus
1_{4}^{1,0}$,
pord$(\rho_\text{isum}(\mathfrak{t})) = 5$,

\vskip 0.7ex
\hangindent=5.5em \hangafter=1
{\white .}\hskip 1em $\rho_\text{isum}(\mathfrak{t})$ =
 $( \frac{1}{4},
\frac{1}{20},
\frac{9}{20} )
\oplus
( \frac{1}{4} )
\oplus
( \frac{1}{4} )
\oplus
( \frac{1}{4} )
$,

\vskip 0.7ex
\hangindent=5.5em \hangafter=1
{\white .}\hskip 1em $\rho_\text{isum}(\mathfrak{s})$ =
$\mathrm{i}$($\sqrt{\frac{1}{5}}$,
$\sqrt{\frac{2}{5}}$,
$\sqrt{\frac{2}{5}}$;\ \ 
$-\frac{5+\sqrt{5}}{10}$,
$\frac{5-\sqrt{5}}{10}$;\ \ 
$-\frac{5+\sqrt{5}}{10}$)
 $\oplus$
$\mathrm{i}$($1$)
 $\oplus$
$\mathrm{i}$($1$)
 $\oplus$
$\mathrm{i}$($1$)

Fail:
dims is $(d,1,1,\cdots)$, $d=(\mathrm{ord}(T)+1)/2$, self-dual, ...
 Prop. B.3 (3)

 \ \color{black}

\noindent 10: (dims,levels) = $(3 , 
1 , 
1 , 
1;20,
4,
4,
4
)$,
irreps = $3_{5}^{3}
\hskip -1.5pt \otimes \hskip -1.5pt
1_{4}^{1,0}\oplus
1_{4}^{1,0}\oplus
1_{4}^{1,0}\oplus
1_{4}^{1,0}$,
pord$(\rho_\text{isum}(\mathfrak{t})) = 5$,

\vskip 0.7ex
\hangindent=5.5em \hangafter=1
{\white .}\hskip 1em $\rho_\text{isum}(\mathfrak{t})$ =
 $( \frac{1}{4},
\frac{13}{20},
\frac{17}{20} )
\oplus
( \frac{1}{4} )
\oplus
( \frac{1}{4} )
\oplus
( \frac{1}{4} )
$,

\vskip 0.7ex
\hangindent=5.5em \hangafter=1
{\white .}\hskip 1em $\rho_\text{isum}(\mathfrak{s})$ =
$\mathrm{i}$($-\sqrt{\frac{1}{5}}$,
$\sqrt{\frac{2}{5}}$,
$\sqrt{\frac{2}{5}}$;\ \ 
$\frac{-5+\sqrt{5}}{10}$,
$\frac{5+\sqrt{5}}{10}$;\ \ 
$\frac{-5+\sqrt{5}}{10}$)
 $\oplus$
$\mathrm{i}$($1$)
 $\oplus$
$\mathrm{i}$($1$)
 $\oplus$
$\mathrm{i}$($1$)

Fail:
dims is $(d,1,1,\cdots)$, $d=(\mathrm{ord}(T)+1)/2$, self-dual, ...
 Prop. B.3 (3)

 \ \color{black}

\noindent 11: (dims,levels) = $(3 , 
1 , 
1 , 
1;24,
3,
3,
3
)$,
irreps = $3_{8}^{3,0}
\hskip -1.5pt \otimes \hskip -1.5pt
1_{3}^{1,0}\oplus
1_{3}^{1,0}\oplus
1_{3}^{1,0}\oplus
1_{3}^{1,0}$,
pord$(\rho_\text{isum}(\mathfrak{t})) = 8$,

\vskip 0.7ex
\hangindent=5.5em \hangafter=1
{\white .}\hskip 1em $\rho_\text{isum}(\mathfrak{t})$ =
 $( \frac{1}{3},
\frac{5}{24},
\frac{17}{24} )
\oplus
( \frac{1}{3} )
\oplus
( \frac{1}{3} )
\oplus
( \frac{1}{3} )
$,

\vskip 0.7ex
\hangindent=5.5em \hangafter=1
{\white .}\hskip 1em $\rho_\text{isum}(\mathfrak{s})$ =
$\mathrm{i}$($0$,
$\sqrt{\frac{1}{2}}$,
$\sqrt{\frac{1}{2}}$;\ \ 
$\frac{1}{2}$,
$-\frac{1}{2}$;\ \ 
$\frac{1}{2}$)
 $\oplus$
($1$)
 $\oplus$
($1$)
 $\oplus$
($1$)

Fail:
for $\rho = \rho_1+l\chi, ...,
 (\rho_1(\mathfrak s)/\chi(\mathfrak s))^2\neq$ id. Prop. B.3 (2)

 \ \color{black}

\noindent 12: (dims,levels) = $(3 , 
1 , 
1 , 
1;24,
3,
3,
3
)$,
irreps = $3_{8}^{1,0}
\hskip -1.5pt \otimes \hskip -1.5pt
1_{3}^{1,0}\oplus
1_{3}^{1,0}\oplus
1_{3}^{1,0}\oplus
1_{3}^{1,0}$,
pord$(\rho_\text{isum}(\mathfrak{t})) = 8$,

\vskip 0.7ex
\hangindent=5.5em \hangafter=1
{\white .}\hskip 1em $\rho_\text{isum}(\mathfrak{t})$ =
 $( \frac{1}{3},
\frac{11}{24},
\frac{23}{24} )
\oplus
( \frac{1}{3} )
\oplus
( \frac{1}{3} )
\oplus
( \frac{1}{3} )
$,

\vskip 0.7ex
\hangindent=5.5em \hangafter=1
{\white .}\hskip 1em $\rho_\text{isum}(\mathfrak{s})$ =
$\mathrm{i}$($0$,
$\sqrt{\frac{1}{2}}$,
$\sqrt{\frac{1}{2}}$;\ \ 
$-\frac{1}{2}$,
$\frac{1}{2}$;\ \ 
$-\frac{1}{2}$)
 $\oplus$
($1$)
 $\oplus$
($1$)
 $\oplus$
($1$)

Fail:
for $\rho = \rho_1+l\chi, ...,
 (\rho_1(\mathfrak s)/\chi(\mathfrak s))^2\neq$ id. Prop. B.3 (2)

 \ \color{black}

\noindent 13: (dims,levels) = $(3 , 
1 , 
1 , 
1;30,
6,
6,
6
)$,
irreps = $3_{5}^{1}
\hskip -1.5pt \otimes \hskip -1.5pt
1_{3}^{1,0}
\hskip -1.5pt \otimes \hskip -1.5pt
1_{2}^{1,0}\oplus
1_{3}^{1,0}
\hskip -1.5pt \otimes \hskip -1.5pt
1_{2}^{1,0}\oplus
1_{3}^{1,0}
\hskip -1.5pt \otimes \hskip -1.5pt
1_{2}^{1,0}\oplus
1_{3}^{1,0}
\hskip -1.5pt \otimes \hskip -1.5pt
1_{2}^{1,0}$,
pord$(\rho_\text{isum}(\mathfrak{t})) = 5$,

\vskip 0.7ex
\hangindent=5.5em \hangafter=1
{\white .}\hskip 1em $\rho_\text{isum}(\mathfrak{t})$ =
 $( \frac{5}{6},
\frac{1}{30},
\frac{19}{30} )
\oplus
( \frac{5}{6} )
\oplus
( \frac{5}{6} )
\oplus
( \frac{5}{6} )
$,

\vskip 0.7ex
\hangindent=5.5em \hangafter=1
{\white .}\hskip 1em $\rho_\text{isum}(\mathfrak{s})$ =
($-\sqrt{\frac{1}{5}}$,
$-\sqrt{\frac{2}{5}}$,
$-\sqrt{\frac{2}{5}}$;
$\frac{5+\sqrt{5}}{10}$,
$\frac{-5+\sqrt{5}}{10}$;
$\frac{5+\sqrt{5}}{10}$)
 $\oplus$
($-1$)
 $\oplus$
($-1$)
 $\oplus$
($-1$)

Fail:
dims is $(d,1,1,\cdots)$, $d=(\mathrm{ord}(T)+1)/2$, self-dual, ...
 Prop. B.3 (3)

 \ \color{black}

\noindent 14: (dims,levels) = $(3 , 
1 , 
1 , 
1;30,
6,
6,
6
)$,
irreps = $3_{5}^{3}
\hskip -1.5pt \otimes \hskip -1.5pt
1_{3}^{1,0}
\hskip -1.5pt \otimes \hskip -1.5pt
1_{2}^{1,0}\oplus
1_{3}^{1,0}
\hskip -1.5pt \otimes \hskip -1.5pt
1_{2}^{1,0}\oplus
1_{3}^{1,0}
\hskip -1.5pt \otimes \hskip -1.5pt
1_{2}^{1,0}\oplus
1_{3}^{1,0}
\hskip -1.5pt \otimes \hskip -1.5pt
1_{2}^{1,0}$,
pord$(\rho_\text{isum}(\mathfrak{t})) = 5$,

\vskip 0.7ex
\hangindent=5.5em \hangafter=1
{\white .}\hskip 1em $\rho_\text{isum}(\mathfrak{t})$ =
 $( \frac{5}{6},
\frac{7}{30},
\frac{13}{30} )
\oplus
( \frac{5}{6} )
\oplus
( \frac{5}{6} )
\oplus
( \frac{5}{6} )
$,

\vskip 0.7ex
\hangindent=5.5em \hangafter=1
{\white .}\hskip 1em $\rho_\text{isum}(\mathfrak{s})$ =
($\sqrt{\frac{1}{5}}$,
$-\sqrt{\frac{2}{5}}$,
$-\sqrt{\frac{2}{5}}$;
$\frac{5-\sqrt{5}}{10}$,
$-\frac{5+\sqrt{5}}{10}$;
$\frac{5-\sqrt{5}}{10}$)
 $\oplus$
($-1$)
 $\oplus$
($-1$)
 $\oplus$
($-1$)

Fail:
dims is $(d,1,1,\cdots)$, $d=(\mathrm{ord}(T)+1)/2$, self-dual, ...
 Prop. B.3 (3)

 \ \color{black}

\noindent 15: (dims,levels) = $(3 , 
1 , 
1 , 
1;60,
12,
12,
12
)$,
irreps = $3_{5}^{3}
\hskip -1.5pt \otimes \hskip -1.5pt
1_{4}^{1,0}
\hskip -1.5pt \otimes \hskip -1.5pt
1_{3}^{1,0}\oplus
1_{4}^{1,0}
\hskip -1.5pt \otimes \hskip -1.5pt
1_{3}^{1,0}\oplus
1_{4}^{1,0}
\hskip -1.5pt \otimes \hskip -1.5pt
1_{3}^{1,0}\oplus
1_{4}^{1,0}
\hskip -1.5pt \otimes \hskip -1.5pt
1_{3}^{1,0}$,
pord$(\rho_\text{isum}(\mathfrak{t})) = 5$,

\vskip 0.7ex
\hangindent=5.5em \hangafter=1
{\white .}\hskip 1em $\rho_\text{isum}(\mathfrak{t})$ =
 $( \frac{7}{12},
\frac{11}{60},
\frac{59}{60} )
\oplus
( \frac{7}{12} )
\oplus
( \frac{7}{12} )
\oplus
( \frac{7}{12} )
$,

\vskip 0.7ex
\hangindent=5.5em \hangafter=1
{\white .}\hskip 1em $\rho_\text{isum}(\mathfrak{s})$ =
$\mathrm{i}$($-\sqrt{\frac{1}{5}}$,
$\sqrt{\frac{2}{5}}$,
$\sqrt{\frac{2}{5}}$;\ \ 
$\frac{-5+\sqrt{5}}{10}$,
$\frac{5+\sqrt{5}}{10}$;\ \ 
$\frac{-5+\sqrt{5}}{10}$)
 $\oplus$
$\mathrm{i}$($1$)
 $\oplus$
$\mathrm{i}$($1$)
 $\oplus$
$\mathrm{i}$($1$)

Fail:
dims is $(d,1,1,\cdots)$, $d=(\mathrm{ord}(T)+1)/2$, self-dual, ...
 Prop. B.3 (3)

 \ \color{black}

\noindent 16: (dims,levels) = $(3 , 
1 , 
1 , 
1;60,
12,
12,
12
)$,
irreps = $3_{5}^{1}
\hskip -1.5pt \otimes \hskip -1.5pt
1_{4}^{1,0}
\hskip -1.5pt \otimes \hskip -1.5pt
1_{3}^{1,0}\oplus
1_{4}^{1,0}
\hskip -1.5pt \otimes \hskip -1.5pt
1_{3}^{1,0}\oplus
1_{4}^{1,0}
\hskip -1.5pt \otimes \hskip -1.5pt
1_{3}^{1,0}\oplus
1_{4}^{1,0}
\hskip -1.5pt \otimes \hskip -1.5pt
1_{3}^{1,0}$,
pord$(\rho_\text{isum}(\mathfrak{t})) = 5$,

\vskip 0.7ex
\hangindent=5.5em \hangafter=1
{\white .}\hskip 1em $\rho_\text{isum}(\mathfrak{t})$ =
 $( \frac{7}{12},
\frac{23}{60},
\frac{47}{60} )
\oplus
( \frac{7}{12} )
\oplus
( \frac{7}{12} )
\oplus
( \frac{7}{12} )
$,

\vskip 0.7ex
\hangindent=5.5em \hangafter=1
{\white .}\hskip 1em $\rho_\text{isum}(\mathfrak{s})$ =
$\mathrm{i}$($\sqrt{\frac{1}{5}}$,
$\sqrt{\frac{2}{5}}$,
$\sqrt{\frac{2}{5}}$;\ \ 
$-\frac{5+\sqrt{5}}{10}$,
$\frac{5-\sqrt{5}}{10}$;\ \ 
$-\frac{5+\sqrt{5}}{10}$)
 $\oplus$
$\mathrm{i}$($1$)
 $\oplus$
$\mathrm{i}$($1$)
 $\oplus$
$\mathrm{i}$($1$)

Fail:
dims is $(d,1,1,\cdots)$, $d=(\mathrm{ord}(T)+1)/2$, self-dual, ...
 Prop. B.3 (3)

 \ \color{black}

\noindent 17: (dims,levels) = $(3 , 
2 , 
1;4,
3,
1
)$,
irreps = $3_{4}^{1,0}\oplus
2_{3}^{1,0}\oplus
1_{1}^{1}$,
pord$(\rho_\text{isum}(\mathfrak{t})) = 12$,

\vskip 0.7ex
\hangindent=5.5em \hangafter=1
{\white .}\hskip 1em $\rho_\text{isum}(\mathfrak{t})$ =
 $( 0,
\frac{1}{4},
\frac{3}{4} )
\oplus
( 0,
\frac{1}{3} )
\oplus
( 0 )
$,

\vskip 0.7ex
\hangindent=5.5em \hangafter=1
{\white .}\hskip 1em $\rho_\text{isum}(\mathfrak{s})$ =
($0$,
$\sqrt{\frac{1}{2}}$,
$\sqrt{\frac{1}{2}}$;
$-\frac{1}{2}$,
$\frac{1}{2}$;
$-\frac{1}{2}$)
 $\oplus$
$\mathrm{i}$($-\sqrt{\frac{1}{3}}$,
$\sqrt{\frac{2}{3}}$;\ \ 
$\sqrt{\frac{1}{3}}$)
 $\oplus$
($1$)

Fail:
Tr$_I(C) = -1 <$  0 for I = [ 1/3 ]. Prop. B.4 (1) eqn. (B.18)

 \ \color{black}

\noindent 18: (dims,levels) = $(3 , 
2 , 
1;4,
3,
1
)$,
irreps = $3_{4}^{1,0}\oplus
2_{3}^{1,8}\oplus
1_{1}^{1}$,
pord$(\rho_\text{isum}(\mathfrak{t})) = 12$,

\vskip 0.7ex
\hangindent=5.5em \hangafter=1
{\white .}\hskip 1em $\rho_\text{isum}(\mathfrak{t})$ =
 $( 0,
\frac{1}{4},
\frac{3}{4} )
\oplus
( 0,
\frac{2}{3} )
\oplus
( 0 )
$,

\vskip 0.7ex
\hangindent=5.5em \hangafter=1
{\white .}\hskip 1em $\rho_\text{isum}(\mathfrak{s})$ =
($0$,
$\sqrt{\frac{1}{2}}$,
$\sqrt{\frac{1}{2}}$;
$-\frac{1}{2}$,
$\frac{1}{2}$;
$-\frac{1}{2}$)
 $\oplus$
$\mathrm{i}$($\sqrt{\frac{1}{3}}$,
$\sqrt{\frac{2}{3}}$;\ \ 
$-\sqrt{\frac{1}{3}}$)
 $\oplus$
($1$)

Fail:
Tr$_I(C) = -1 <$  0 for I = [ 2/3 ]. Prop. B.4 (1) eqn. (B.18)

 \ \color{black}

\noindent 19: (dims,levels) = $(3 , 
2 , 
1;4,
3,
3
)$,
irreps = $3_{4}^{1,0}\oplus
2_{3}^{1,0}\oplus
1_{3}^{1,0}$,
pord$(\rho_\text{isum}(\mathfrak{t})) = 12$,

\vskip 0.7ex
\hangindent=5.5em \hangafter=1
{\white .}\hskip 1em $\rho_\text{isum}(\mathfrak{t})$ =
 $( 0,
\frac{1}{4},
\frac{3}{4} )
\oplus
( 0,
\frac{1}{3} )
\oplus
( \frac{1}{3} )
$,

\vskip 0.7ex
\hangindent=5.5em \hangafter=1
{\white .}\hskip 1em $\rho_\text{isum}(\mathfrak{s})$ =
($0$,
$\sqrt{\frac{1}{2}}$,
$\sqrt{\frac{1}{2}}$;
$-\frac{1}{2}$,
$\frac{1}{2}$;
$-\frac{1}{2}$)
 $\oplus$
$\mathrm{i}$($-\sqrt{\frac{1}{3}}$,
$\sqrt{\frac{2}{3}}$;\ \ 
$\sqrt{\frac{1}{3}}$)
 $\oplus$
($1$)

Fail:
all $\theta$-eigenspaces that can contain unit
 have Tr$_{E_\theta}(C) \leq 0 $. Prop. B.5 (5) eqn. (B.29)

 \ \color{black}

\noindent 20: (dims,levels) = $(3 , 
2 , 
1;4,
3,
3
)$,
irreps = $3_{4}^{1,0}\oplus
2_{3}^{1,8}\oplus
1_{3}^{1,4}$,
pord$(\rho_\text{isum}(\mathfrak{t})) = 12$,

\vskip 0.7ex
\hangindent=5.5em \hangafter=1
{\white .}\hskip 1em $\rho_\text{isum}(\mathfrak{t})$ =
 $( 0,
\frac{1}{4},
\frac{3}{4} )
\oplus
( 0,
\frac{2}{3} )
\oplus
( \frac{2}{3} )
$,

\vskip 0.7ex
\hangindent=5.5em \hangafter=1
{\white .}\hskip 1em $\rho_\text{isum}(\mathfrak{s})$ =
($0$,
$\sqrt{\frac{1}{2}}$,
$\sqrt{\frac{1}{2}}$;
$-\frac{1}{2}$,
$\frac{1}{2}$;
$-\frac{1}{2}$)
 $\oplus$
$\mathrm{i}$($\sqrt{\frac{1}{3}}$,
$\sqrt{\frac{2}{3}}$;\ \ 
$-\sqrt{\frac{1}{3}}$)
 $\oplus$
($1$)

Fail:
all $\theta$-eigenspaces that can contain unit
 have Tr$_{E_\theta}(C) \leq 0 $. Prop. B.5 (5) eqn. (B.29)

 \ \color{black}

\noindent 21: (dims,levels) = $(3 , 
2 , 
1;4,
3,
4
)$,
irreps = $3_{4}^{1,0}\oplus
2_{3}^{1,0}\oplus
1_{4}^{1,0}$,
pord$(\rho_\text{isum}(\mathfrak{t})) = 12$,

\vskip 0.7ex
\hangindent=5.5em \hangafter=1
{\white .}\hskip 1em $\rho_\text{isum}(\mathfrak{t})$ =
 $( 0,
\frac{1}{4},
\frac{3}{4} )
\oplus
( 0,
\frac{1}{3} )
\oplus
( \frac{1}{4} )
$,

\vskip 0.7ex
\hangindent=5.5em \hangafter=1
{\white .}\hskip 1em $\rho_\text{isum}(\mathfrak{s})$ =
($0$,
$\sqrt{\frac{1}{2}}$,
$\sqrt{\frac{1}{2}}$;
$-\frac{1}{2}$,
$\frac{1}{2}$;
$-\frac{1}{2}$)
 $\oplus$
$\mathrm{i}$($-\sqrt{\frac{1}{3}}$,
$\sqrt{\frac{2}{3}}$;\ \ 
$\sqrt{\frac{1}{3}}$)
 $\oplus$
$\mathrm{i}$($1$)

Fail:
number of self dual objects $|$Tr($\rho(\mathfrak s^2)$)$|$ = 0. Prop. B.4 (1)\
 eqn. (B.16)

 \ \color{black}

\noindent 22: (dims,levels) = $(3 , 
2 , 
1;4,
3,
4
)$,
irreps = $3_{4}^{1,0}\oplus
2_{3}^{1,0}\oplus
1_{4}^{1,6}$,
pord$(\rho_\text{isum}(\mathfrak{t})) = 12$,

\vskip 0.7ex
\hangindent=5.5em \hangafter=1
{\white .}\hskip 1em $\rho_\text{isum}(\mathfrak{t})$ =
 $( 0,
\frac{1}{4},
\frac{3}{4} )
\oplus
( 0,
\frac{1}{3} )
\oplus
( \frac{3}{4} )
$,

\vskip 0.7ex
\hangindent=5.5em \hangafter=1
{\white .}\hskip 1em $\rho_\text{isum}(\mathfrak{s})$ =
($0$,
$\sqrt{\frac{1}{2}}$,
$\sqrt{\frac{1}{2}}$;
$-\frac{1}{2}$,
$\frac{1}{2}$;
$-\frac{1}{2}$)
 $\oplus$
$\mathrm{i}$($-\sqrt{\frac{1}{3}}$,
$\sqrt{\frac{2}{3}}$;\ \ 
$\sqrt{\frac{1}{3}}$)
 $\oplus$
$\mathrm{i}$($-1$)

Fail:
number of self dual objects $|$Tr($\rho(\mathfrak s^2)$)$|$ = 0. Prop. B.4 (1)\
 eqn. (B.16)

 \ \color{black}

\noindent 23: (dims,levels) = $(3 , 
2 , 
1;4,
3,
4
)$,
irreps = $3_{4}^{1,0}\oplus
2_{3}^{1,8}\oplus
1_{4}^{1,0}$,
pord$(\rho_\text{isum}(\mathfrak{t})) = 12$,

\vskip 0.7ex
\hangindent=5.5em \hangafter=1
{\white .}\hskip 1em $\rho_\text{isum}(\mathfrak{t})$ =
 $( 0,
\frac{1}{4},
\frac{3}{4} )
\oplus
( 0,
\frac{2}{3} )
\oplus
( \frac{1}{4} )
$,

\vskip 0.7ex
\hangindent=5.5em \hangafter=1
{\white .}\hskip 1em $\rho_\text{isum}(\mathfrak{s})$ =
($0$,
$\sqrt{\frac{1}{2}}$,
$\sqrt{\frac{1}{2}}$;
$-\frac{1}{2}$,
$\frac{1}{2}$;
$-\frac{1}{2}$)
 $\oplus$
$\mathrm{i}$($\sqrt{\frac{1}{3}}$,
$\sqrt{\frac{2}{3}}$;\ \ 
$-\sqrt{\frac{1}{3}}$)
 $\oplus$
$\mathrm{i}$($1$)

Fail:
number of self dual objects $|$Tr($\rho(\mathfrak s^2)$)$|$ = 0. Prop. B.4 (1)\
 eqn. (B.16)

 \ \color{black}

\noindent 24: (dims,levels) = $(3 , 
2 , 
1;4,
3,
4
)$,
irreps = $3_{4}^{1,0}\oplus
2_{3}^{1,8}\oplus
1_{4}^{1,6}$,
pord$(\rho_\text{isum}(\mathfrak{t})) = 12$,

\vskip 0.7ex
\hangindent=5.5em \hangafter=1
{\white .}\hskip 1em $\rho_\text{isum}(\mathfrak{t})$ =
 $( 0,
\frac{1}{4},
\frac{3}{4} )
\oplus
( 0,
\frac{2}{3} )
\oplus
( \frac{3}{4} )
$,

\vskip 0.7ex
\hangindent=5.5em \hangafter=1
{\white .}\hskip 1em $\rho_\text{isum}(\mathfrak{s})$ =
($0$,
$\sqrt{\frac{1}{2}}$,
$\sqrt{\frac{1}{2}}$;
$-\frac{1}{2}$,
$\frac{1}{2}$;
$-\frac{1}{2}$)
 $\oplus$
$\mathrm{i}$($\sqrt{\frac{1}{3}}$,
$\sqrt{\frac{2}{3}}$;\ \ 
$-\sqrt{\frac{1}{3}}$)
 $\oplus$
$\mathrm{i}$($-1$)

Fail:
number of self dual objects $|$Tr($\rho(\mathfrak s^2)$)$|$ = 0. Prop. B.4 (1)\
 eqn. (B.16)

 \ \color{black}

 \color{blue}

\noindent 25: (dims,levels) = $(3 , 
2 , 
1;4,
12,
1
)$,
irreps = $3_{4}^{1,0}\oplus
2_{3}^{1,0}
\hskip -1.5pt \otimes \hskip -1.5pt
1_{4}^{1,0}\oplus
1_{1}^{1}$,
pord$(\rho_\text{isum}(\mathfrak{t})) = 12$,

\vskip 0.7ex
\hangindent=5.5em \hangafter=1
{\white .}\hskip 1em $\rho_\text{isum}(\mathfrak{t})$ =
 $( 0,
\frac{1}{4},
\frac{3}{4} )
\oplus
( \frac{1}{4},
\frac{7}{12} )
\oplus
( 0 )
$,

\vskip 0.7ex
\hangindent=5.5em \hangafter=1
{\white .}\hskip 1em $\rho_\text{isum}(\mathfrak{s})$ =
($0$,
$\sqrt{\frac{1}{2}}$,
$\sqrt{\frac{1}{2}}$;
$-\frac{1}{2}$,
$\frac{1}{2}$;
$-\frac{1}{2}$)
 $\oplus$
($\sqrt{\frac{1}{3}}$,
$\sqrt{\frac{2}{3}}$;
$-\sqrt{\frac{1}{3}}$)
 $\oplus$
($1$)

Pass. 

 \ \color{black}

 \color{blue}

\noindent 26: (dims,levels) = $(3 , 
2 , 
1;4,
12,
1
)$,
irreps = $3_{4}^{1,0}\oplus
2_{3}^{1,8}
\hskip -1.5pt \otimes \hskip -1.5pt
1_{4}^{1,0}\oplus
1_{1}^{1}$,
pord$(\rho_\text{isum}(\mathfrak{t})) = 12$,

\vskip 0.7ex
\hangindent=5.5em \hangafter=1
{\white .}\hskip 1em $\rho_\text{isum}(\mathfrak{t})$ =
 $( 0,
\frac{1}{4},
\frac{3}{4} )
\oplus
( \frac{1}{4},
\frac{11}{12} )
\oplus
( 0 )
$,

\vskip 0.7ex
\hangindent=5.5em \hangafter=1
{\white .}\hskip 1em $\rho_\text{isum}(\mathfrak{s})$ =
($0$,
$\sqrt{\frac{1}{2}}$,
$\sqrt{\frac{1}{2}}$;
$-\frac{1}{2}$,
$\frac{1}{2}$;
$-\frac{1}{2}$)
 $\oplus$
($-\sqrt{\frac{1}{3}}$,
$\sqrt{\frac{2}{3}}$;
$\sqrt{\frac{1}{3}}$)
 $\oplus$
($1$)

Pass. 

 \ \color{black}

 \color{blue}

\noindent 27: (dims,levels) = $(3 , 
2 , 
1;4,
12,
1
)$,
irreps = $3_{4}^{1,0}\oplus
2_{3}^{1,0}
\hskip -1.5pt \otimes \hskip -1.5pt
1_{4}^{1,6}\oplus
1_{1}^{1}$,
pord$(\rho_\text{isum}(\mathfrak{t})) = 12$,

\vskip 0.7ex
\hangindent=5.5em \hangafter=1
{\white .}\hskip 1em $\rho_\text{isum}(\mathfrak{t})$ =
 $( 0,
\frac{1}{4},
\frac{3}{4} )
\oplus
( \frac{3}{4},
\frac{1}{12} )
\oplus
( 0 )
$,

\vskip 0.7ex
\hangindent=5.5em \hangafter=1
{\white .}\hskip 1em $\rho_\text{isum}(\mathfrak{s})$ =
($0$,
$\sqrt{\frac{1}{2}}$,
$\sqrt{\frac{1}{2}}$;
$-\frac{1}{2}$,
$\frac{1}{2}$;
$-\frac{1}{2}$)
 $\oplus$
($-\sqrt{\frac{1}{3}}$,
$\sqrt{\frac{2}{3}}$;
$\sqrt{\frac{1}{3}}$)
 $\oplus$
($1$)

Pass. 

 \ \color{black}

 \color{blue}

\noindent 28: (dims,levels) = $(3 , 
2 , 
1;4,
12,
1
)$,
irreps = $3_{4}^{1,0}\oplus
2_{3}^{1,8}
\hskip -1.5pt \otimes \hskip -1.5pt
1_{4}^{1,6}\oplus
1_{1}^{1}$,
pord$(\rho_\text{isum}(\mathfrak{t})) = 12$,

\vskip 0.7ex
\hangindent=5.5em \hangafter=1
{\white .}\hskip 1em $\rho_\text{isum}(\mathfrak{t})$ =
 $( 0,
\frac{1}{4},
\frac{3}{4} )
\oplus
( \frac{3}{4},
\frac{5}{12} )
\oplus
( 0 )
$,

\vskip 0.7ex
\hangindent=5.5em \hangafter=1
{\white .}\hskip 1em $\rho_\text{isum}(\mathfrak{s})$ =
($0$,
$\sqrt{\frac{1}{2}}$,
$\sqrt{\frac{1}{2}}$;
$-\frac{1}{2}$,
$\frac{1}{2}$;
$-\frac{1}{2}$)
 $\oplus$
($\sqrt{\frac{1}{3}}$,
$\sqrt{\frac{2}{3}}$;
$-\sqrt{\frac{1}{3}}$)
 $\oplus$
($1$)

Pass. 

 \ \color{black}

\noindent 29: (dims,levels) = $(3 , 
2 , 
1;4,
12,
4
)$,
irreps = $3_{4}^{1,0}\oplus
2_{3}^{1,0}
\hskip -1.5pt \otimes \hskip -1.5pt
1_{4}^{1,0}\oplus
1_{4}^{1,0}$,
pord$(\rho_\text{isum}(\mathfrak{t})) = 12$,

\vskip 0.7ex
\hangindent=5.5em \hangafter=1
{\white .}\hskip 1em $\rho_\text{isum}(\mathfrak{t})$ =
 $( 0,
\frac{1}{4},
\frac{3}{4} )
\oplus
( \frac{1}{4},
\frac{7}{12} )
\oplus
( \frac{1}{4} )
$,

\vskip 0.7ex
\hangindent=5.5em \hangafter=1
{\white .}\hskip 1em $\rho_\text{isum}(\mathfrak{s})$ =
($0$,
$\sqrt{\frac{1}{2}}$,
$\sqrt{\frac{1}{2}}$;
$-\frac{1}{2}$,
$\frac{1}{2}$;
$-\frac{1}{2}$)
 $\oplus$
($\sqrt{\frac{1}{3}}$,
$\sqrt{\frac{2}{3}}$;
$-\sqrt{\frac{1}{3}}$)
 $\oplus$
$\mathrm{i}$($1$)

Fail:
cnd($\rho(\mathfrak s)_\mathrm{ndeg}$) = 24 does not divide
 ord($\rho(\mathfrak t)$)=12. Prop. B.4 (2)

 \ \color{black}

 \color{blue}

\noindent 30: (dims,levels) = $(3 , 
2 , 
1;4,
12,
4
)$,
irreps = $3_{4}^{1,0}\oplus
2_{3}^{1,0}
\hskip -1.5pt \otimes \hskip -1.5pt
1_{4}^{1,0}\oplus
1_{4}^{1,6}$,
pord$(\rho_\text{isum}(\mathfrak{t})) = 12$,

\vskip 0.7ex
\hangindent=5.5em \hangafter=1
{\white .}\hskip 1em $\rho_\text{isum}(\mathfrak{t})$ =
 $( 0,
\frac{1}{4},
\frac{3}{4} )
\oplus
( \frac{1}{4},
\frac{7}{12} )
\oplus
( \frac{3}{4} )
$,

\vskip 0.7ex
\hangindent=5.5em \hangafter=1
{\white .}\hskip 1em $\rho_\text{isum}(\mathfrak{s})$ =
($0$,
$\sqrt{\frac{1}{2}}$,
$\sqrt{\frac{1}{2}}$;
$-\frac{1}{2}$,
$\frac{1}{2}$;
$-\frac{1}{2}$)
 $\oplus$
($\sqrt{\frac{1}{3}}$,
$\sqrt{\frac{2}{3}}$;
$-\sqrt{\frac{1}{3}}$)
 $\oplus$
$\mathrm{i}$($-1$)

Pass. 

 \ \color{black}

\noindent 31: (dims,levels) = $(3 , 
2 , 
1;4,
12,
4
)$,
irreps = $3_{4}^{1,0}\oplus
2_{3}^{1,8}
\hskip -1.5pt \otimes \hskip -1.5pt
1_{4}^{1,0}\oplus
1_{4}^{1,0}$,
pord$(\rho_\text{isum}(\mathfrak{t})) = 12$,

\vskip 0.7ex
\hangindent=5.5em \hangafter=1
{\white .}\hskip 1em $\rho_\text{isum}(\mathfrak{t})$ =
 $( 0,
\frac{1}{4},
\frac{3}{4} )
\oplus
( \frac{1}{4},
\frac{11}{12} )
\oplus
( \frac{1}{4} )
$,

\vskip 0.7ex
\hangindent=5.5em \hangafter=1
{\white .}\hskip 1em $\rho_\text{isum}(\mathfrak{s})$ =
($0$,
$\sqrt{\frac{1}{2}}$,
$\sqrt{\frac{1}{2}}$;
$-\frac{1}{2}$,
$\frac{1}{2}$;
$-\frac{1}{2}$)
 $\oplus$
($-\sqrt{\frac{1}{3}}$,
$\sqrt{\frac{2}{3}}$;
$\sqrt{\frac{1}{3}}$)
 $\oplus$
$\mathrm{i}$($1$)

Fail:
cnd($\rho(\mathfrak s)_\mathrm{ndeg}$) = 24 does not divide
 ord($\rho(\mathfrak t)$)=12. Prop. B.4 (2)

 \ \color{black}

 \color{blue}

\noindent 32: (dims,levels) = $(3 , 
2 , 
1;4,
12,
4
)$,
irreps = $3_{4}^{1,0}\oplus
2_{3}^{1,8}
\hskip -1.5pt \otimes \hskip -1.5pt
1_{4}^{1,0}\oplus
1_{4}^{1,6}$,
pord$(\rho_\text{isum}(\mathfrak{t})) = 12$,

\vskip 0.7ex
\hangindent=5.5em \hangafter=1
{\white .}\hskip 1em $\rho_\text{isum}(\mathfrak{t})$ =
 $( 0,
\frac{1}{4},
\frac{3}{4} )
\oplus
( \frac{1}{4},
\frac{11}{12} )
\oplus
( \frac{3}{4} )
$,

\vskip 0.7ex
\hangindent=5.5em \hangafter=1
{\white .}\hskip 1em $\rho_\text{isum}(\mathfrak{s})$ =
($0$,
$\sqrt{\frac{1}{2}}$,
$\sqrt{\frac{1}{2}}$;
$-\frac{1}{2}$,
$\frac{1}{2}$;
$-\frac{1}{2}$)
 $\oplus$
($-\sqrt{\frac{1}{3}}$,
$\sqrt{\frac{2}{3}}$;
$\sqrt{\frac{1}{3}}$)
 $\oplus$
$\mathrm{i}$($-1$)

Pass. 

 \ \color{black}

 \color{blue}

\noindent 33: (dims,levels) = $(3 , 
2 , 
1;4,
12,
4
)$,
irreps = $3_{4}^{1,0}\oplus
2_{3}^{1,0}
\hskip -1.5pt \otimes \hskip -1.5pt
1_{4}^{1,6}\oplus
1_{4}^{1,0}$,
pord$(\rho_\text{isum}(\mathfrak{t})) = 12$,

\vskip 0.7ex
\hangindent=5.5em \hangafter=1
{\white .}\hskip 1em $\rho_\text{isum}(\mathfrak{t})$ =
 $( 0,
\frac{1}{4},
\frac{3}{4} )
\oplus
( \frac{3}{4},
\frac{1}{12} )
\oplus
( \frac{1}{4} )
$,

\vskip 0.7ex
\hangindent=5.5em \hangafter=1
{\white .}\hskip 1em $\rho_\text{isum}(\mathfrak{s})$ =
($0$,
$\sqrt{\frac{1}{2}}$,
$\sqrt{\frac{1}{2}}$;
$-\frac{1}{2}$,
$\frac{1}{2}$;
$-\frac{1}{2}$)
 $\oplus$
($-\sqrt{\frac{1}{3}}$,
$\sqrt{\frac{2}{3}}$;
$\sqrt{\frac{1}{3}}$)
 $\oplus$
$\mathrm{i}$($1$)

Pass. 

 \ \color{black}

\noindent 34: (dims,levels) = $(3 , 
2 , 
1;4,
12,
4
)$,
irreps = $3_{4}^{1,0}\oplus
2_{3}^{1,0}
\hskip -1.5pt \otimes \hskip -1.5pt
1_{4}^{1,6}\oplus
1_{4}^{1,6}$,
pord$(\rho_\text{isum}(\mathfrak{t})) = 12$,

\vskip 0.7ex
\hangindent=5.5em \hangafter=1
{\white .}\hskip 1em $\rho_\text{isum}(\mathfrak{t})$ =
 $( 0,
\frac{1}{4},
\frac{3}{4} )
\oplus
( \frac{3}{4},
\frac{1}{12} )
\oplus
( \frac{3}{4} )
$,

\vskip 0.7ex
\hangindent=5.5em \hangafter=1
{\white .}\hskip 1em $\rho_\text{isum}(\mathfrak{s})$ =
($0$,
$\sqrt{\frac{1}{2}}$,
$\sqrt{\frac{1}{2}}$;
$-\frac{1}{2}$,
$\frac{1}{2}$;
$-\frac{1}{2}$)
 $\oplus$
($-\sqrt{\frac{1}{3}}$,
$\sqrt{\frac{2}{3}}$;
$\sqrt{\frac{1}{3}}$)
 $\oplus$
$\mathrm{i}$($-1$)

Fail:
cnd($\rho(\mathfrak s)_\mathrm{ndeg}$) = 24 does not divide
 ord($\rho(\mathfrak t)$)=12. Prop. B.4 (2)

 \ \color{black}

 \color{blue}

\noindent 35: (dims,levels) = $(3 , 
2 , 
1;4,
12,
4
)$,
irreps = $3_{4}^{1,0}\oplus
2_{3}^{1,8}
\hskip -1.5pt \otimes \hskip -1.5pt
1_{4}^{1,6}\oplus
1_{4}^{1,0}$,
pord$(\rho_\text{isum}(\mathfrak{t})) = 12$,

\vskip 0.7ex
\hangindent=5.5em \hangafter=1
{\white .}\hskip 1em $\rho_\text{isum}(\mathfrak{t})$ =
 $( 0,
\frac{1}{4},
\frac{3}{4} )
\oplus
( \frac{3}{4},
\frac{5}{12} )
\oplus
( \frac{1}{4} )
$,

\vskip 0.7ex
\hangindent=5.5em \hangafter=1
{\white .}\hskip 1em $\rho_\text{isum}(\mathfrak{s})$ =
($0$,
$\sqrt{\frac{1}{2}}$,
$\sqrt{\frac{1}{2}}$;
$-\frac{1}{2}$,
$\frac{1}{2}$;
$-\frac{1}{2}$)
 $\oplus$
($\sqrt{\frac{1}{3}}$,
$\sqrt{\frac{2}{3}}$;
$-\sqrt{\frac{1}{3}}$)
 $\oplus$
$\mathrm{i}$($1$)

Pass. 

 \ \color{black}

\noindent 36: (dims,levels) = $(3 , 
2 , 
1;4,
12,
4
)$,
irreps = $3_{4}^{1,0}\oplus
2_{3}^{1,8}
\hskip -1.5pt \otimes \hskip -1.5pt
1_{4}^{1,6}\oplus
1_{4}^{1,6}$,
pord$(\rho_\text{isum}(\mathfrak{t})) = 12$,

\vskip 0.7ex
\hangindent=5.5em \hangafter=1
{\white .}\hskip 1em $\rho_\text{isum}(\mathfrak{t})$ =
 $( 0,
\frac{1}{4},
\frac{3}{4} )
\oplus
( \frac{3}{4},
\frac{5}{12} )
\oplus
( \frac{3}{4} )
$,

\vskip 0.7ex
\hangindent=5.5em \hangafter=1
{\white .}\hskip 1em $\rho_\text{isum}(\mathfrak{s})$ =
($0$,
$\sqrt{\frac{1}{2}}$,
$\sqrt{\frac{1}{2}}$;
$-\frac{1}{2}$,
$\frac{1}{2}$;
$-\frac{1}{2}$)
 $\oplus$
($\sqrt{\frac{1}{3}}$,
$\sqrt{\frac{2}{3}}$;
$-\sqrt{\frac{1}{3}}$)
 $\oplus$
$\mathrm{i}$($-1$)

Fail:
cnd($\rho(\mathfrak s)_\mathrm{ndeg}$) = 24 does not divide
 ord($\rho(\mathfrak t)$)=12. Prop. B.4 (2)

 \ \color{black}

\noindent 37: (dims,levels) = $(3 , 
2 , 
1;4,
12,
12
)$,
irreps = $3_{4}^{1,0}\oplus
2_{3}^{1,0}
\hskip -1.5pt \otimes \hskip -1.5pt
1_{4}^{1,0}\oplus
1_{4}^{1,0}
\hskip -1.5pt \otimes \hskip -1.5pt
1_{3}^{1,0}$,
pord$(\rho_\text{isum}(\mathfrak{t})) = 12$,

\vskip 0.7ex
\hangindent=5.5em \hangafter=1
{\white .}\hskip 1em $\rho_\text{isum}(\mathfrak{t})$ =
 $( 0,
\frac{1}{4},
\frac{3}{4} )
\oplus
( \frac{1}{4},
\frac{7}{12} )
\oplus
( \frac{7}{12} )
$,

\vskip 0.7ex
\hangindent=5.5em \hangafter=1
{\white .}\hskip 1em $\rho_\text{isum}(\mathfrak{s})$ =
($0$,
$\sqrt{\frac{1}{2}}$,
$\sqrt{\frac{1}{2}}$;
$-\frac{1}{2}$,
$\frac{1}{2}$;
$-\frac{1}{2}$)
 $\oplus$
($\sqrt{\frac{1}{3}}$,
$\sqrt{\frac{2}{3}}$;
$-\sqrt{\frac{1}{3}}$)
 $\oplus$
$\mathrm{i}$($1$)

Fail:
cnd($\rho(\mathfrak s)_\mathrm{ndeg}$) = 8 does not divide
 ord($\rho(\mathfrak t)$)=12. Prop. B.4 (2)

 \ \color{black}

\noindent 38: (dims,levels) = $(3 , 
2 , 
1;4,
12,
12
)$,
irreps = $3_{4}^{1,0}\oplus
2_{3}^{1,8}
\hskip -1.5pt \otimes \hskip -1.5pt
1_{4}^{1,0}\oplus
1_{4}^{1,0}
\hskip -1.5pt \otimes \hskip -1.5pt
1_{3}^{1,4}$,
pord$(\rho_\text{isum}(\mathfrak{t})) = 12$,

\vskip 0.7ex
\hangindent=5.5em \hangafter=1
{\white .}\hskip 1em $\rho_\text{isum}(\mathfrak{t})$ =
 $( 0,
\frac{1}{4},
\frac{3}{4} )
\oplus
( \frac{1}{4},
\frac{11}{12} )
\oplus
( \frac{11}{12} )
$,

\vskip 0.7ex
\hangindent=5.5em \hangafter=1
{\white .}\hskip 1em $\rho_\text{isum}(\mathfrak{s})$ =
($0$,
$\sqrt{\frac{1}{2}}$,
$\sqrt{\frac{1}{2}}$;
$-\frac{1}{2}$,
$\frac{1}{2}$;
$-\frac{1}{2}$)
 $\oplus$
($-\sqrt{\frac{1}{3}}$,
$\sqrt{\frac{2}{3}}$;
$\sqrt{\frac{1}{3}}$)
 $\oplus$
$\mathrm{i}$($1$)

Fail:
cnd($\rho(\mathfrak s)_\mathrm{ndeg}$) = 8 does not divide
 ord($\rho(\mathfrak t)$)=12. Prop. B.4 (2)

 \ \color{black}

\noindent 39: (dims,levels) = $(3 , 
2 , 
1;4,
12,
12
)$,
irreps = $3_{4}^{1,0}\oplus
2_{3}^{1,0}
\hskip -1.5pt \otimes \hskip -1.5pt
1_{4}^{1,6}\oplus
1_{4}^{1,6}
\hskip -1.5pt \otimes \hskip -1.5pt
1_{3}^{1,0}$,
pord$(\rho_\text{isum}(\mathfrak{t})) = 12$,

\vskip 0.7ex
\hangindent=5.5em \hangafter=1
{\white .}\hskip 1em $\rho_\text{isum}(\mathfrak{t})$ =
 $( 0,
\frac{1}{4},
\frac{3}{4} )
\oplus
( \frac{3}{4},
\frac{1}{12} )
\oplus
( \frac{1}{12} )
$,

\vskip 0.7ex
\hangindent=5.5em \hangafter=1
{\white .}\hskip 1em $\rho_\text{isum}(\mathfrak{s})$ =
($0$,
$\sqrt{\frac{1}{2}}$,
$\sqrt{\frac{1}{2}}$;
$-\frac{1}{2}$,
$\frac{1}{2}$;
$-\frac{1}{2}$)
 $\oplus$
($-\sqrt{\frac{1}{3}}$,
$\sqrt{\frac{2}{3}}$;
$\sqrt{\frac{1}{3}}$)
 $\oplus$
$\mathrm{i}$($-1$)

Fail:
cnd($\rho(\mathfrak s)_\mathrm{ndeg}$) = 8 does not divide
 ord($\rho(\mathfrak t)$)=12. Prop. B.4 (2)

 \ \color{black}

\noindent 40: (dims,levels) = $(3 , 
2 , 
1;4,
12,
12
)$,
irreps = $3_{4}^{1,0}\oplus
2_{3}^{1,8}
\hskip -1.5pt \otimes \hskip -1.5pt
1_{4}^{1,6}\oplus
1_{4}^{1,6}
\hskip -1.5pt \otimes \hskip -1.5pt
1_{3}^{1,4}$,
pord$(\rho_\text{isum}(\mathfrak{t})) = 12$,

\vskip 0.7ex
\hangindent=5.5em \hangafter=1
{\white .}\hskip 1em $\rho_\text{isum}(\mathfrak{t})$ =
 $( 0,
\frac{1}{4},
\frac{3}{4} )
\oplus
( \frac{3}{4},
\frac{5}{12} )
\oplus
( \frac{5}{12} )
$,

\vskip 0.7ex
\hangindent=5.5em \hangafter=1
{\white .}\hskip 1em $\rho_\text{isum}(\mathfrak{s})$ =
($0$,
$\sqrt{\frac{1}{2}}$,
$\sqrt{\frac{1}{2}}$;
$-\frac{1}{2}$,
$\frac{1}{2}$;
$-\frac{1}{2}$)
 $\oplus$
($\sqrt{\frac{1}{3}}$,
$\sqrt{\frac{2}{3}}$;
$-\sqrt{\frac{1}{3}}$)
 $\oplus$
$\mathrm{i}$($-1$)

Fail:
cnd($\rho(\mathfrak s)_\mathrm{ndeg}$) = 8 does not divide
 ord($\rho(\mathfrak t)$)=12. Prop. B.4 (2)

 \ \color{black}

 \color{blue}

\noindent 41: (dims,levels) = $(3 , 
2 , 
1;5,
2,
1
)$,
irreps = $3_{5}^{1}\oplus
2_{2}^{1,0}\oplus
1_{1}^{1}$,
pord$(\rho_\text{isum}(\mathfrak{t})) = 10$,

\vskip 0.7ex
\hangindent=5.5em \hangafter=1
{\white .}\hskip 1em $\rho_\text{isum}(\mathfrak{t})$ =
 $( 0,
\frac{1}{5},
\frac{4}{5} )
\oplus
( 0,
\frac{1}{2} )
\oplus
( 0 )
$,

\vskip 0.7ex
\hangindent=5.5em \hangafter=1
{\white .}\hskip 1em $\rho_\text{isum}(\mathfrak{s})$ =
($\sqrt{\frac{1}{5}}$,
$-\sqrt{\frac{2}{5}}$,
$-\sqrt{\frac{2}{5}}$;
$-\frac{5+\sqrt{5}}{10}$,
$\frac{5-\sqrt{5}}{10}$;
$-\frac{5+\sqrt{5}}{10}$)
 $\oplus$
($-\frac{1}{2}$,
$-\sqrt{\frac{3}{4}}$;
$\frac{1}{2}$)
 $\oplus$
($1$)

Pass. 

 \ \color{black}

 \color{blue}

\noindent 42: (dims,levels) = $(3 , 
2 , 
1;5,
2,
1
)$,
irreps = $3_{5}^{3}\oplus
2_{2}^{1,0}\oplus
1_{1}^{1}$,
pord$(\rho_\text{isum}(\mathfrak{t})) = 10$,

\vskip 0.7ex
\hangindent=5.5em \hangafter=1
{\white .}\hskip 1em $\rho_\text{isum}(\mathfrak{t})$ =
 $( 0,
\frac{2}{5},
\frac{3}{5} )
\oplus
( 0,
\frac{1}{2} )
\oplus
( 0 )
$,

\vskip 0.7ex
\hangindent=5.5em \hangafter=1
{\white .}\hskip 1em $\rho_\text{isum}(\mathfrak{s})$ =
($-\sqrt{\frac{1}{5}}$,
$-\sqrt{\frac{2}{5}}$,
$-\sqrt{\frac{2}{5}}$;
$\frac{-5+\sqrt{5}}{10}$,
$\frac{5+\sqrt{5}}{10}$;
$\frac{-5+\sqrt{5}}{10}$)
 $\oplus$
($-\frac{1}{2}$,
$-\sqrt{\frac{3}{4}}$;
$\frac{1}{2}$)
 $\oplus$
($1$)

Pass. 

 \ \color{black}

 \color{blue}

\noindent 43: (dims,levels) = $(3 , 
2 , 
1;5,
2,
2
)$,
irreps = $3_{5}^{1}\oplus
2_{2}^{1,0}\oplus
1_{2}^{1,0}$,
pord$(\rho_\text{isum}(\mathfrak{t})) = 10$,

\vskip 0.7ex
\hangindent=5.5em \hangafter=1
{\white .}\hskip 1em $\rho_\text{isum}(\mathfrak{t})$ =
 $( 0,
\frac{1}{5},
\frac{4}{5} )
\oplus
( 0,
\frac{1}{2} )
\oplus
( \frac{1}{2} )
$,

\vskip 0.7ex
\hangindent=5.5em \hangafter=1
{\white .}\hskip 1em $\rho_\text{isum}(\mathfrak{s})$ =
($\sqrt{\frac{1}{5}}$,
$-\sqrt{\frac{2}{5}}$,
$-\sqrt{\frac{2}{5}}$;
$-\frac{5+\sqrt{5}}{10}$,
$\frac{5-\sqrt{5}}{10}$;
$-\frac{5+\sqrt{5}}{10}$)
 $\oplus$
($-\frac{1}{2}$,
$-\sqrt{\frac{3}{4}}$;
$\frac{1}{2}$)
 $\oplus$
($-1$)

Pass. 

 \ \color{black}

 \color{blue}

\noindent 44: (dims,levels) = $(3 , 
2 , 
1;5,
2,
2
)$,
irreps = $3_{5}^{3}\oplus
2_{2}^{1,0}\oplus
1_{2}^{1,0}$,
pord$(\rho_\text{isum}(\mathfrak{t})) = 10$,

\vskip 0.7ex
\hangindent=5.5em \hangafter=1
{\white .}\hskip 1em $\rho_\text{isum}(\mathfrak{t})$ =
 $( 0,
\frac{2}{5},
\frac{3}{5} )
\oplus
( 0,
\frac{1}{2} )
\oplus
( \frac{1}{2} )
$,

\vskip 0.7ex
\hangindent=5.5em \hangafter=1
{\white .}\hskip 1em $\rho_\text{isum}(\mathfrak{s})$ =
($-\sqrt{\frac{1}{5}}$,
$-\sqrt{\frac{2}{5}}$,
$-\sqrt{\frac{2}{5}}$;
$\frac{-5+\sqrt{5}}{10}$,
$\frac{5+\sqrt{5}}{10}$;
$\frac{-5+\sqrt{5}}{10}$)
 $\oplus$
($-\frac{1}{2}$,
$-\sqrt{\frac{3}{4}}$;
$\frac{1}{2}$)
 $\oplus$
($-1$)

Pass. 

 \ \color{black}

\noindent 45: (dims,levels) = $(3 , 
2 , 
1;5,
3,
1
)$,
irreps = $3_{5}^{1}\oplus
2_{3}^{1,0}\oplus
1_{1}^{1}$,
pord$(\rho_\text{isum}(\mathfrak{t})) = 15$,

\vskip 0.7ex
\hangindent=5.5em \hangafter=1
{\white .}\hskip 1em $\rho_\text{isum}(\mathfrak{t})$ =
 $( 0,
\frac{1}{5},
\frac{4}{5} )
\oplus
( 0,
\frac{1}{3} )
\oplus
( 0 )
$,

\vskip 0.7ex
\hangindent=5.5em \hangafter=1
{\white .}\hskip 1em $\rho_\text{isum}(\mathfrak{s})$ =
($\sqrt{\frac{1}{5}}$,
$-\sqrt{\frac{2}{5}}$,
$-\sqrt{\frac{2}{5}}$;
$-\frac{5+\sqrt{5}}{10}$,
$\frac{5-\sqrt{5}}{10}$;
$-\frac{5+\sqrt{5}}{10}$)
 $\oplus$
$\mathrm{i}$($-\sqrt{\frac{1}{3}}$,
$\sqrt{\frac{2}{3}}$;\ \ 
$\sqrt{\frac{1}{3}}$)
 $\oplus$
($1$)

Fail:
Tr$_I(C) = -1 <$  0 for I = [ 1/3 ]. Prop. B.4 (1) eqn. (B.18)

 \ \color{black}

\noindent 46: (dims,levels) = $(3 , 
2 , 
1;5,
3,
1
)$,
irreps = $3_{5}^{1}\oplus
2_{3}^{1,8}\oplus
1_{1}^{1}$,
pord$(\rho_\text{isum}(\mathfrak{t})) = 15$,

\vskip 0.7ex
\hangindent=5.5em \hangafter=1
{\white .}\hskip 1em $\rho_\text{isum}(\mathfrak{t})$ =
 $( 0,
\frac{1}{5},
\frac{4}{5} )
\oplus
( 0,
\frac{2}{3} )
\oplus
( 0 )
$,

\vskip 0.7ex
\hangindent=5.5em \hangafter=1
{\white .}\hskip 1em $\rho_\text{isum}(\mathfrak{s})$ =
($\sqrt{\frac{1}{5}}$,
$-\sqrt{\frac{2}{5}}$,
$-\sqrt{\frac{2}{5}}$;
$-\frac{5+\sqrt{5}}{10}$,
$\frac{5-\sqrt{5}}{10}$;
$-\frac{5+\sqrt{5}}{10}$)
 $\oplus$
$\mathrm{i}$($\sqrt{\frac{1}{3}}$,
$\sqrt{\frac{2}{3}}$;\ \ 
$-\sqrt{\frac{1}{3}}$)
 $\oplus$
($1$)

Fail:
Tr$_I(C) = -1 <$  0 for I = [ 2/3 ]. Prop. B.4 (1) eqn. (B.18)

 \ \color{black}

\noindent 47: (dims,levels) = $(3 , 
2 , 
1;5,
3,
1
)$,
irreps = $3_{5}^{3}\oplus
2_{3}^{1,0}\oplus
1_{1}^{1}$,
pord$(\rho_\text{isum}(\mathfrak{t})) = 15$,

\vskip 0.7ex
\hangindent=5.5em \hangafter=1
{\white .}\hskip 1em $\rho_\text{isum}(\mathfrak{t})$ =
 $( 0,
\frac{2}{5},
\frac{3}{5} )
\oplus
( 0,
\frac{1}{3} )
\oplus
( 0 )
$,

\vskip 0.7ex
\hangindent=5.5em \hangafter=1
{\white .}\hskip 1em $\rho_\text{isum}(\mathfrak{s})$ =
($-\sqrt{\frac{1}{5}}$,
$-\sqrt{\frac{2}{5}}$,
$-\sqrt{\frac{2}{5}}$;
$\frac{-5+\sqrt{5}}{10}$,
$\frac{5+\sqrt{5}}{10}$;
$\frac{-5+\sqrt{5}}{10}$)
 $\oplus$
$\mathrm{i}$($-\sqrt{\frac{1}{3}}$,
$\sqrt{\frac{2}{3}}$;\ \ 
$\sqrt{\frac{1}{3}}$)
 $\oplus$
($1$)

Fail:
Tr$_I(C) = -1 <$  0 for I = [ 1/3 ]. Prop. B.4 (1) eqn. (B.18)

 \ \color{black}

\noindent 48: (dims,levels) = $(3 , 
2 , 
1;5,
3,
1
)$,
irreps = $3_{5}^{3}\oplus
2_{3}^{1,8}\oplus
1_{1}^{1}$,
pord$(\rho_\text{isum}(\mathfrak{t})) = 15$,

\vskip 0.7ex
\hangindent=5.5em \hangafter=1
{\white .}\hskip 1em $\rho_\text{isum}(\mathfrak{t})$ =
 $( 0,
\frac{2}{5},
\frac{3}{5} )
\oplus
( 0,
\frac{2}{3} )
\oplus
( 0 )
$,

\vskip 0.7ex
\hangindent=5.5em \hangafter=1
{\white .}\hskip 1em $\rho_\text{isum}(\mathfrak{s})$ =
($-\sqrt{\frac{1}{5}}$,
$-\sqrt{\frac{2}{5}}$,
$-\sqrt{\frac{2}{5}}$;
$\frac{-5+\sqrt{5}}{10}$,
$\frac{5+\sqrt{5}}{10}$;
$\frac{-5+\sqrt{5}}{10}$)
 $\oplus$
$\mathrm{i}$($\sqrt{\frac{1}{3}}$,
$\sqrt{\frac{2}{3}}$;\ \ 
$-\sqrt{\frac{1}{3}}$)
 $\oplus$
($1$)

Fail:
Tr$_I(C) = -1 <$  0 for I = [ 2/3 ]. Prop. B.4 (1) eqn. (B.18)

 \ \color{black}

\noindent 49: (dims,levels) = $(3 , 
2 , 
1;5,
3,
3
)$,
irreps = $3_{5}^{1}\oplus
2_{3}^{1,0}\oplus
1_{3}^{1,0}$,
pord$(\rho_\text{isum}(\mathfrak{t})) = 15$,

\vskip 0.7ex
\hangindent=5.5em \hangafter=1
{\white .}\hskip 1em $\rho_\text{isum}(\mathfrak{t})$ =
 $( 0,
\frac{1}{5},
\frac{4}{5} )
\oplus
( 0,
\frac{1}{3} )
\oplus
( \frac{1}{3} )
$,

\vskip 0.7ex
\hangindent=5.5em \hangafter=1
{\white .}\hskip 1em $\rho_\text{isum}(\mathfrak{s})$ =
($\sqrt{\frac{1}{5}}$,
$-\sqrt{\frac{2}{5}}$,
$-\sqrt{\frac{2}{5}}$;
$-\frac{5+\sqrt{5}}{10}$,
$\frac{5-\sqrt{5}}{10}$;
$-\frac{5+\sqrt{5}}{10}$)
 $\oplus$
$\mathrm{i}$($-\sqrt{\frac{1}{3}}$,
$\sqrt{\frac{2}{3}}$;\ \ 
$\sqrt{\frac{1}{3}}$)
 $\oplus$
($1$)

Fail:
all $\theta$-eigenspaces that can contain unit
 have Tr$_{E_\theta}(C) \leq 0 $. Prop. B.5 (5) eqn. (B.29)

 \ \color{black}

\noindent 50: (dims,levels) = $(3 , 
2 , 
1;5,
3,
3
)$,
irreps = $3_{5}^{1}\oplus
2_{3}^{1,8}\oplus
1_{3}^{1,4}$,
pord$(\rho_\text{isum}(\mathfrak{t})) = 15$,

\vskip 0.7ex
\hangindent=5.5em \hangafter=1
{\white .}\hskip 1em $\rho_\text{isum}(\mathfrak{t})$ =
 $( 0,
\frac{1}{5},
\frac{4}{5} )
\oplus
( 0,
\frac{2}{3} )
\oplus
( \frac{2}{3} )
$,

\vskip 0.7ex
\hangindent=5.5em \hangafter=1
{\white .}\hskip 1em $\rho_\text{isum}(\mathfrak{s})$ =
($\sqrt{\frac{1}{5}}$,
$-\sqrt{\frac{2}{5}}$,
$-\sqrt{\frac{2}{5}}$;
$-\frac{5+\sqrt{5}}{10}$,
$\frac{5-\sqrt{5}}{10}$;
$-\frac{5+\sqrt{5}}{10}$)
 $\oplus$
$\mathrm{i}$($\sqrt{\frac{1}{3}}$,
$\sqrt{\frac{2}{3}}$;\ \ 
$-\sqrt{\frac{1}{3}}$)
 $\oplus$
($1$)

Fail:
all $\theta$-eigenspaces that can contain unit
 have Tr$_{E_\theta}(C) \leq 0 $. Prop. B.5 (5) eqn. (B.29)

 \ \color{black}

\noindent 51: (dims,levels) = $(3 , 
2 , 
1;5,
3,
3
)$,
irreps = $3_{5}^{3}\oplus
2_{3}^{1,0}\oplus
1_{3}^{1,0}$,
pord$(\rho_\text{isum}(\mathfrak{t})) = 15$,

\vskip 0.7ex
\hangindent=5.5em \hangafter=1
{\white .}\hskip 1em $\rho_\text{isum}(\mathfrak{t})$ =
 $( 0,
\frac{2}{5},
\frac{3}{5} )
\oplus
( 0,
\frac{1}{3} )
\oplus
( \frac{1}{3} )
$,

\vskip 0.7ex
\hangindent=5.5em \hangafter=1
{\white .}\hskip 1em $\rho_\text{isum}(\mathfrak{s})$ =
($-\sqrt{\frac{1}{5}}$,
$-\sqrt{\frac{2}{5}}$,
$-\sqrt{\frac{2}{5}}$;
$\frac{-5+\sqrt{5}}{10}$,
$\frac{5+\sqrt{5}}{10}$;
$\frac{-5+\sqrt{5}}{10}$)
 $\oplus$
$\mathrm{i}$($-\sqrt{\frac{1}{3}}$,
$\sqrt{\frac{2}{3}}$;\ \ 
$\sqrt{\frac{1}{3}}$)
 $\oplus$
($1$)

Fail:
all $\theta$-eigenspaces that can contain unit
 have Tr$_{E_\theta}(C) \leq 0 $. Prop. B.5 (5) eqn. (B.29)

 \ \color{black}

\noindent 52: (dims,levels) = $(3 , 
2 , 
1;5,
3,
3
)$,
irreps = $3_{5}^{3}\oplus
2_{3}^{1,8}\oplus
1_{3}^{1,4}$,
pord$(\rho_\text{isum}(\mathfrak{t})) = 15$,

\vskip 0.7ex
\hangindent=5.5em \hangafter=1
{\white .}\hskip 1em $\rho_\text{isum}(\mathfrak{t})$ =
 $( 0,
\frac{2}{5},
\frac{3}{5} )
\oplus
( 0,
\frac{2}{3} )
\oplus
( \frac{2}{3} )
$,

\vskip 0.7ex
\hangindent=5.5em \hangafter=1
{\white .}\hskip 1em $\rho_\text{isum}(\mathfrak{s})$ =
($-\sqrt{\frac{1}{5}}$,
$-\sqrt{\frac{2}{5}}$,
$-\sqrt{\frac{2}{5}}$;
$\frac{-5+\sqrt{5}}{10}$,
$\frac{5+\sqrt{5}}{10}$;
$\frac{-5+\sqrt{5}}{10}$)
 $\oplus$
$\mathrm{i}$($\sqrt{\frac{1}{3}}$,
$\sqrt{\frac{2}{3}}$;\ \ 
$-\sqrt{\frac{1}{3}}$)
 $\oplus$
($1$)

Fail:
all $\theta$-eigenspaces that can contain unit
 have Tr$_{E_\theta}(C) \leq 0 $. Prop. B.5 (5) eqn. (B.29)

 \ \color{black}

 \color{blue}

\noindent 53: (dims,levels) = $(3 , 
2 , 
1;5,
5,
1
)$,
irreps = $3_{5}^{1}\oplus
2_{5}^{1}\oplus
1_{1}^{1}$,
pord$(\rho_\text{isum}(\mathfrak{t})) = 5$,

\vskip 0.7ex
\hangindent=5.5em \hangafter=1
{\white .}\hskip 1em $\rho_\text{isum}(\mathfrak{t})$ =
 $( 0,
\frac{1}{5},
\frac{4}{5} )
\oplus
( \frac{1}{5},
\frac{4}{5} )
\oplus
( 0 )
$,

\vskip 0.7ex
\hangindent=5.5em \hangafter=1
{\white .}\hskip 1em $\rho_\text{isum}(\mathfrak{s})$ =
($\sqrt{\frac{1}{5}}$,
$-\sqrt{\frac{2}{5}}$,
$-\sqrt{\frac{2}{5}}$;
$-\frac{5+\sqrt{5}}{10}$,
$\frac{5-\sqrt{5}}{10}$;
$-\frac{5+\sqrt{5}}{10}$)
 $\oplus$
$\mathrm{i}$($-\frac{1}{\sqrt{5}}c^{3}_{20}
$,
$\frac{1}{\sqrt{5}}c^{1}_{20}
$;\ \ 
$\frac{1}{\sqrt{5}}c^{3}_{20}
$)
 $\oplus$
($1$)

Pass. 

 \ \color{black}

 \color{blue}

\noindent 54: (dims,levels) = $(3 , 
2 , 
1;5,
5,
1
)$,
irreps = $3_{5}^{3}\oplus
2_{5}^{2}\oplus
1_{1}^{1}$,
pord$(\rho_\text{isum}(\mathfrak{t})) = 5$,

\vskip 0.7ex
\hangindent=5.5em \hangafter=1
{\white .}\hskip 1em $\rho_\text{isum}(\mathfrak{t})$ =
 $( 0,
\frac{2}{5},
\frac{3}{5} )
\oplus
( \frac{2}{5},
\frac{3}{5} )
\oplus
( 0 )
$,

\vskip 0.7ex
\hangindent=5.5em \hangafter=1
{\white .}\hskip 1em $\rho_\text{isum}(\mathfrak{s})$ =
($-\sqrt{\frac{1}{5}}$,
$-\sqrt{\frac{2}{5}}$,
$-\sqrt{\frac{2}{5}}$;
$\frac{-5+\sqrt{5}}{10}$,
$\frac{5+\sqrt{5}}{10}$;
$\frac{-5+\sqrt{5}}{10}$)
 $\oplus$
$\mathrm{i}$($-\frac{1}{\sqrt{5}}c^{1}_{20}
$,
$\frac{1}{\sqrt{5}}c^{3}_{20}
$;\ \ 
$\frac{1}{\sqrt{5}}c^{1}_{20}
$)
 $\oplus$
($1$)

Pass. 

 \ \color{black}

\noindent 55: (dims,levels) = $(3 , 
2 , 
1;8,
2,
1
)$,
irreps = $3_{8}^{1,0}\oplus
2_{2}^{1,0}\oplus
1_{1}^{1}$,
pord$(\rho_\text{isum}(\mathfrak{t})) = 8$,

\vskip 0.7ex
\hangindent=5.5em \hangafter=1
{\white .}\hskip 1em $\rho_\text{isum}(\mathfrak{t})$ =
 $( 0,
\frac{1}{8},
\frac{5}{8} )
\oplus
( 0,
\frac{1}{2} )
\oplus
( 0 )
$,

\vskip 0.7ex
\hangindent=5.5em \hangafter=1
{\white .}\hskip 1em $\rho_\text{isum}(\mathfrak{s})$ =
$\mathrm{i}$($0$,
$\sqrt{\frac{1}{2}}$,
$\sqrt{\frac{1}{2}}$;\ \ 
$-\frac{1}{2}$,
$\frac{1}{2}$;\ \ 
$-\frac{1}{2}$)
 $\oplus$
($-\frac{1}{2}$,
$-\sqrt{\frac{3}{4}}$;
$\frac{1}{2}$)
 $\oplus$
($1$)

Fail:
number of self dual objects $|$Tr($\rho(\mathfrak s^2)$)$|$ = 0. Prop. B.4 (1)\
 eqn. (B.16)

 \ \color{black}

\noindent 56: (dims,levels) = $(3 , 
2 , 
1;8,
2,
1
)$,
irreps = $3_{8}^{3,0}\oplus
2_{2}^{1,0}\oplus
1_{1}^{1}$,
pord$(\rho_\text{isum}(\mathfrak{t})) = 8$,

\vskip 0.7ex
\hangindent=5.5em \hangafter=1
{\white .}\hskip 1em $\rho_\text{isum}(\mathfrak{t})$ =
 $( 0,
\frac{3}{8},
\frac{7}{8} )
\oplus
( 0,
\frac{1}{2} )
\oplus
( 0 )
$,

\vskip 0.7ex
\hangindent=5.5em \hangafter=1
{\white .}\hskip 1em $\rho_\text{isum}(\mathfrak{s})$ =
$\mathrm{i}$($0$,
$\sqrt{\frac{1}{2}}$,
$\sqrt{\frac{1}{2}}$;\ \ 
$\frac{1}{2}$,
$-\frac{1}{2}$;\ \ 
$\frac{1}{2}$)
 $\oplus$
($-\frac{1}{2}$,
$-\sqrt{\frac{3}{4}}$;
$\frac{1}{2}$)
 $\oplus$
($1$)

Fail:
number of self dual objects $|$Tr($\rho(\mathfrak s^2)$)$|$ = 0. Prop. B.4 (1)\
 eqn. (B.16)

 \ \color{black}

\noindent 57: (dims,levels) = $(3 , 
2 , 
1;8,
2,
2
)$,
irreps = $3_{8}^{1,0}\oplus
2_{2}^{1,0}\oplus
1_{2}^{1,0}$,
pord$(\rho_\text{isum}(\mathfrak{t})) = 8$,

\vskip 0.7ex
\hangindent=5.5em \hangafter=1
{\white .}\hskip 1em $\rho_\text{isum}(\mathfrak{t})$ =
 $( 0,
\frac{1}{8},
\frac{5}{8} )
\oplus
( 0,
\frac{1}{2} )
\oplus
( \frac{1}{2} )
$,

\vskip 0.7ex
\hangindent=5.5em \hangafter=1
{\white .}\hskip 1em $\rho_\text{isum}(\mathfrak{s})$ =
$\mathrm{i}$($0$,
$\sqrt{\frac{1}{2}}$,
$\sqrt{\frac{1}{2}}$;\ \ 
$-\frac{1}{2}$,
$\frac{1}{2}$;\ \ 
$-\frac{1}{2}$)
 $\oplus$
($-\frac{1}{2}$,
$-\sqrt{\frac{3}{4}}$;
$\frac{1}{2}$)
 $\oplus$
($-1$)

Fail:
number of self dual objects $|$Tr($\rho(\mathfrak s^2)$)$|$ = 0. Prop. B.4 (1)\
 eqn. (B.16)

 \ \color{black}

\noindent 58: (dims,levels) = $(3 , 
2 , 
1;8,
2,
2
)$,
irreps = $3_{8}^{3,0}\oplus
2_{2}^{1,0}\oplus
1_{2}^{1,0}$,
pord$(\rho_\text{isum}(\mathfrak{t})) = 8$,

\vskip 0.7ex
\hangindent=5.5em \hangafter=1
{\white .}\hskip 1em $\rho_\text{isum}(\mathfrak{t})$ =
 $( 0,
\frac{3}{8},
\frac{7}{8} )
\oplus
( 0,
\frac{1}{2} )
\oplus
( \frac{1}{2} )
$,

\vskip 0.7ex
\hangindent=5.5em \hangafter=1
{\white .}\hskip 1em $\rho_\text{isum}(\mathfrak{s})$ =
$\mathrm{i}$($0$,
$\sqrt{\frac{1}{2}}$,
$\sqrt{\frac{1}{2}}$;\ \ 
$\frac{1}{2}$,
$-\frac{1}{2}$;\ \ 
$\frac{1}{2}$)
 $\oplus$
($-\frac{1}{2}$,
$-\sqrt{\frac{3}{4}}$;
$\frac{1}{2}$)
 $\oplus$
($-1$)

Fail:
number of self dual objects $|$Tr($\rho(\mathfrak s^2)$)$|$ = 0. Prop. B.4 (1)\
 eqn. (B.16)

 \ \color{black}

 \color{blue}

\noindent 59: (dims,levels) = $(3 , 
2 , 
1;8,
3,
1
)$,
irreps = $3_{8}^{1,0}\oplus
2_{3}^{1,0}\oplus
1_{1}^{1}$,
pord$(\rho_\text{isum}(\mathfrak{t})) = 24$,

\vskip 0.7ex
\hangindent=5.5em \hangafter=1
{\white .}\hskip 1em $\rho_\text{isum}(\mathfrak{t})$ =
 $( 0,
\frac{1}{8},
\frac{5}{8} )
\oplus
( 0,
\frac{1}{3} )
\oplus
( 0 )
$,

\vskip 0.7ex
\hangindent=5.5em \hangafter=1
{\white .}\hskip 1em $\rho_\text{isum}(\mathfrak{s})$ =
$\mathrm{i}$($0$,
$\sqrt{\frac{1}{2}}$,
$\sqrt{\frac{1}{2}}$;\ \ 
$-\frac{1}{2}$,
$\frac{1}{2}$;\ \ 
$-\frac{1}{2}$)
 $\oplus$
$\mathrm{i}$($-\sqrt{\frac{1}{3}}$,
$\sqrt{\frac{2}{3}}$;\ \ 
$\sqrt{\frac{1}{3}}$)
 $\oplus$
($1$)

Pass. 

 \ \color{black}

 \color{blue}

\noindent 60: (dims,levels) = $(3 , 
2 , 
1;8,
3,
1
)$,
irreps = $3_{8}^{1,0}\oplus
2_{3}^{1,8}\oplus
1_{1}^{1}$,
pord$(\rho_\text{isum}(\mathfrak{t})) = 24$,

\vskip 0.7ex
\hangindent=5.5em \hangafter=1
{\white .}\hskip 1em $\rho_\text{isum}(\mathfrak{t})$ =
 $( 0,
\frac{1}{8},
\frac{5}{8} )
\oplus
( 0,
\frac{2}{3} )
\oplus
( 0 )
$,

\vskip 0.7ex
\hangindent=5.5em \hangafter=1
{\white .}\hskip 1em $\rho_\text{isum}(\mathfrak{s})$ =
$\mathrm{i}$($0$,
$\sqrt{\frac{1}{2}}$,
$\sqrt{\frac{1}{2}}$;\ \ 
$-\frac{1}{2}$,
$\frac{1}{2}$;\ \ 
$-\frac{1}{2}$)
 $\oplus$
$\mathrm{i}$($\sqrt{\frac{1}{3}}$,
$\sqrt{\frac{2}{3}}$;\ \ 
$-\sqrt{\frac{1}{3}}$)
 $\oplus$
($1$)

Pass. 

 \ \color{black}

 \color{blue}

\noindent 61: (dims,levels) = $(3 , 
2 , 
1;8,
3,
1
)$,
irreps = $3_{8}^{3,0}\oplus
2_{3}^{1,0}\oplus
1_{1}^{1}$,
pord$(\rho_\text{isum}(\mathfrak{t})) = 24$,

\vskip 0.7ex
\hangindent=5.5em \hangafter=1
{\white .}\hskip 1em $\rho_\text{isum}(\mathfrak{t})$ =
 $( 0,
\frac{3}{8},
\frac{7}{8} )
\oplus
( 0,
\frac{1}{3} )
\oplus
( 0 )
$,

\vskip 0.7ex
\hangindent=5.5em \hangafter=1
{\white .}\hskip 1em $\rho_\text{isum}(\mathfrak{s})$ =
$\mathrm{i}$($0$,
$\sqrt{\frac{1}{2}}$,
$\sqrt{\frac{1}{2}}$;\ \ 
$\frac{1}{2}$,
$-\frac{1}{2}$;\ \ 
$\frac{1}{2}$)
 $\oplus$
$\mathrm{i}$($-\sqrt{\frac{1}{3}}$,
$\sqrt{\frac{2}{3}}$;\ \ 
$\sqrt{\frac{1}{3}}$)
 $\oplus$
($1$)

Pass. 

 \ \color{black}

 \color{blue}

\noindent 62: (dims,levels) = $(3 , 
2 , 
1;8,
3,
1
)$,
irreps = $3_{8}^{3,0}\oplus
2_{3}^{1,8}\oplus
1_{1}^{1}$,
pord$(\rho_\text{isum}(\mathfrak{t})) = 24$,

\vskip 0.7ex
\hangindent=5.5em \hangafter=1
{\white .}\hskip 1em $\rho_\text{isum}(\mathfrak{t})$ =
 $( 0,
\frac{3}{8},
\frac{7}{8} )
\oplus
( 0,
\frac{2}{3} )
\oplus
( 0 )
$,

\vskip 0.7ex
\hangindent=5.5em \hangafter=1
{\white .}\hskip 1em $\rho_\text{isum}(\mathfrak{s})$ =
$\mathrm{i}$($0$,
$\sqrt{\frac{1}{2}}$,
$\sqrt{\frac{1}{2}}$;\ \ 
$\frac{1}{2}$,
$-\frac{1}{2}$;\ \ 
$\frac{1}{2}$)
 $\oplus$
$\mathrm{i}$($\sqrt{\frac{1}{3}}$,
$\sqrt{\frac{2}{3}}$;\ \ 
$-\sqrt{\frac{1}{3}}$)
 $\oplus$
($1$)

Pass. 

 \ \color{black}

 \color{blue}

\noindent 63: (dims,levels) = $(3 , 
2 , 
1;8,
3,
3
)$,
irreps = $3_{8}^{1,0}\oplus
2_{3}^{1,0}\oplus
1_{3}^{1,0}$,
pord$(\rho_\text{isum}(\mathfrak{t})) = 24$,

\vskip 0.7ex
\hangindent=5.5em \hangafter=1
{\white .}\hskip 1em $\rho_\text{isum}(\mathfrak{t})$ =
 $( 0,
\frac{1}{8},
\frac{5}{8} )
\oplus
( 0,
\frac{1}{3} )
\oplus
( \frac{1}{3} )
$,

\vskip 0.7ex
\hangindent=5.5em \hangafter=1
{\white .}\hskip 1em $\rho_\text{isum}(\mathfrak{s})$ =
$\mathrm{i}$($0$,
$\sqrt{\frac{1}{2}}$,
$\sqrt{\frac{1}{2}}$;\ \ 
$-\frac{1}{2}$,
$\frac{1}{2}$;\ \ 
$-\frac{1}{2}$)
 $\oplus$
$\mathrm{i}$($-\sqrt{\frac{1}{3}}$,
$\sqrt{\frac{2}{3}}$;\ \ 
$\sqrt{\frac{1}{3}}$)
 $\oplus$
($1$)

Pass. 

 \ \color{black}

 \color{blue}

\noindent 64: (dims,levels) = $(3 , 
2 , 
1;8,
3,
3
)$,
irreps = $3_{8}^{1,0}\oplus
2_{3}^{1,8}\oplus
1_{3}^{1,4}$,
pord$(\rho_\text{isum}(\mathfrak{t})) = 24$,

\vskip 0.7ex
\hangindent=5.5em \hangafter=1
{\white .}\hskip 1em $\rho_\text{isum}(\mathfrak{t})$ =
 $( 0,
\frac{1}{8},
\frac{5}{8} )
\oplus
( 0,
\frac{2}{3} )
\oplus
( \frac{2}{3} )
$,

\vskip 0.7ex
\hangindent=5.5em \hangafter=1
{\white .}\hskip 1em $\rho_\text{isum}(\mathfrak{s})$ =
$\mathrm{i}$($0$,
$\sqrt{\frac{1}{2}}$,
$\sqrt{\frac{1}{2}}$;\ \ 
$-\frac{1}{2}$,
$\frac{1}{2}$;\ \ 
$-\frac{1}{2}$)
 $\oplus$
$\mathrm{i}$($\sqrt{\frac{1}{3}}$,
$\sqrt{\frac{2}{3}}$;\ \ 
$-\sqrt{\frac{1}{3}}$)
 $\oplus$
($1$)

Pass. 

 \ \color{black}

 \color{blue}

\noindent 65: (dims,levels) = $(3 , 
2 , 
1;8,
3,
3
)$,
irreps = $3_{8}^{3,0}\oplus
2_{3}^{1,0}\oplus
1_{3}^{1,0}$,
pord$(\rho_\text{isum}(\mathfrak{t})) = 24$,

\vskip 0.7ex
\hangindent=5.5em \hangafter=1
{\white .}\hskip 1em $\rho_\text{isum}(\mathfrak{t})$ =
 $( 0,
\frac{3}{8},
\frac{7}{8} )
\oplus
( 0,
\frac{1}{3} )
\oplus
( \frac{1}{3} )
$,

\vskip 0.7ex
\hangindent=5.5em \hangafter=1
{\white .}\hskip 1em $\rho_\text{isum}(\mathfrak{s})$ =
$\mathrm{i}$($0$,
$\sqrt{\frac{1}{2}}$,
$\sqrt{\frac{1}{2}}$;\ \ 
$\frac{1}{2}$,
$-\frac{1}{2}$;\ \ 
$\frac{1}{2}$)
 $\oplus$
$\mathrm{i}$($-\sqrt{\frac{1}{3}}$,
$\sqrt{\frac{2}{3}}$;\ \ 
$\sqrt{\frac{1}{3}}$)
 $\oplus$
($1$)

Pass. 

 \ \color{black}

 \color{blue}

\noindent 66: (dims,levels) = $(3 , 
2 , 
1;8,
3,
3
)$,
irreps = $3_{8}^{3,0}\oplus
2_{3}^{1,8}\oplus
1_{3}^{1,4}$,
pord$(\rho_\text{isum}(\mathfrak{t})) = 24$,

\vskip 0.7ex
\hangindent=5.5em \hangafter=1
{\white .}\hskip 1em $\rho_\text{isum}(\mathfrak{t})$ =
 $( 0,
\frac{3}{8},
\frac{7}{8} )
\oplus
( 0,
\frac{2}{3} )
\oplus
( \frac{2}{3} )
$,

\vskip 0.7ex
\hangindent=5.5em \hangafter=1
{\white .}\hskip 1em $\rho_\text{isum}(\mathfrak{s})$ =
$\mathrm{i}$($0$,
$\sqrt{\frac{1}{2}}$,
$\sqrt{\frac{1}{2}}$;\ \ 
$\frac{1}{2}$,
$-\frac{1}{2}$;\ \ 
$\frac{1}{2}$)
 $\oplus$
$\mathrm{i}$($\sqrt{\frac{1}{3}}$,
$\sqrt{\frac{2}{3}}$;\ \ 
$-\sqrt{\frac{1}{3}}$)
 $\oplus$
($1$)

Pass. 

 \ \color{black}

\noindent 67: (dims,levels) = $(3 , 
2 , 
1;8,
8,
1
)$,
irreps = $3_{8}^{1,0}\oplus
2_{8}^{1,0}\oplus
1_{1}^{1}$,
pord$(\rho_\text{isum}(\mathfrak{t})) = 8$,

\vskip 0.7ex
\hangindent=5.5em \hangafter=1
{\white .}\hskip 1em $\rho_\text{isum}(\mathfrak{t})$ =
 $( 0,
\frac{1}{8},
\frac{5}{8} )
\oplus
( \frac{1}{8},
\frac{3}{8} )
\oplus
( 0 )
$,

\vskip 0.7ex
\hangindent=5.5em \hangafter=1
{\white .}\hskip 1em $\rho_\text{isum}(\mathfrak{s})$ =
$\mathrm{i}$($0$,
$\sqrt{\frac{1}{2}}$,
$\sqrt{\frac{1}{2}}$;\ \ 
$-\frac{1}{2}$,
$\frac{1}{2}$;\ \ 
$-\frac{1}{2}$)
 $\oplus$
($-\sqrt{\frac{1}{2}}$,
$\sqrt{\frac{1}{2}}$;
$\sqrt{\frac{1}{2}}$)
 $\oplus$
($1$)

Fail:
number of self dual objects $|$Tr($\rho(\mathfrak s^2)$)$|$ = 0. Prop. B.4 (1)\
 eqn. (B.16)

 \ \color{black}

 \color{blue}

\noindent 68: (dims,levels) = $(3 , 
2 , 
1;8,
8,
1
)$,
irreps = $3_{8}^{1,0}\oplus
2_{8}^{1,9}\oplus
1_{1}^{1}$,
pord$(\rho_\text{isum}(\mathfrak{t})) = 8$,

\vskip 0.7ex
\hangindent=5.5em \hangafter=1
{\white .}\hskip 1em $\rho_\text{isum}(\mathfrak{t})$ =
 $( 0,
\frac{1}{8},
\frac{5}{8} )
\oplus
( \frac{1}{8},
\frac{7}{8} )
\oplus
( 0 )
$,

\vskip 0.7ex
\hangindent=5.5em \hangafter=1
{\white .}\hskip 1em $\rho_\text{isum}(\mathfrak{s})$ =
$\mathrm{i}$($0$,
$\sqrt{\frac{1}{2}}$,
$\sqrt{\frac{1}{2}}$;\ \ 
$-\frac{1}{2}$,
$\frac{1}{2}$;\ \ 
$-\frac{1}{2}$)
 $\oplus$
$\mathrm{i}$($-\sqrt{\frac{1}{2}}$,
$\sqrt{\frac{1}{2}}$;\ \ 
$\sqrt{\frac{1}{2}}$)
 $\oplus$
($1$)

Pass. 

 \ \color{black}

 \color{blue}

\noindent 69: (dims,levels) = $(3 , 
2 , 
1;8,
8,
1
)$,
irreps = $3_{8}^{1,0}\oplus
2_{8}^{1,3}\oplus
1_{1}^{1}$,
pord$(\rho_\text{isum}(\mathfrak{t})) = 8$,

\vskip 0.7ex
\hangindent=5.5em \hangafter=1
{\white .}\hskip 1em $\rho_\text{isum}(\mathfrak{t})$ =
 $( 0,
\frac{1}{8},
\frac{5}{8} )
\oplus
( \frac{3}{8},
\frac{5}{8} )
\oplus
( 0 )
$,

\vskip 0.7ex
\hangindent=5.5em \hangafter=1
{\white .}\hskip 1em $\rho_\text{isum}(\mathfrak{s})$ =
$\mathrm{i}$($0$,
$\sqrt{\frac{1}{2}}$,
$\sqrt{\frac{1}{2}}$;\ \ 
$-\frac{1}{2}$,
$\frac{1}{2}$;\ \ 
$-\frac{1}{2}$)
 $\oplus$
$\mathrm{i}$($-\sqrt{\frac{1}{2}}$,
$\sqrt{\frac{1}{2}}$;\ \ 
$\sqrt{\frac{1}{2}}$)
 $\oplus$
($1$)

Pass. 

 \ \color{black}

\noindent 70: (dims,levels) = $(3 , 
2 , 
1;8,
8,
1
)$,
irreps = $3_{8}^{1,0}\oplus
2_{8}^{1,6}\oplus
1_{1}^{1}$,
pord$(\rho_\text{isum}(\mathfrak{t})) = 8$,

\vskip 0.7ex
\hangindent=5.5em \hangafter=1
{\white .}\hskip 1em $\rho_\text{isum}(\mathfrak{t})$ =
 $( 0,
\frac{1}{8},
\frac{5}{8} )
\oplus
( \frac{5}{8},
\frac{7}{8} )
\oplus
( 0 )
$,

\vskip 0.7ex
\hangindent=5.5em \hangafter=1
{\white .}\hskip 1em $\rho_\text{isum}(\mathfrak{s})$ =
$\mathrm{i}$($0$,
$\sqrt{\frac{1}{2}}$,
$\sqrt{\frac{1}{2}}$;\ \ 
$-\frac{1}{2}$,
$\frac{1}{2}$;\ \ 
$-\frac{1}{2}$)
 $\oplus$
($\sqrt{\frac{1}{2}}$,
$\sqrt{\frac{1}{2}}$;
$-\sqrt{\frac{1}{2}}$)
 $\oplus$
($1$)

Fail:
number of self dual objects $|$Tr($\rho(\mathfrak s^2)$)$|$ = 0. Prop. B.4 (1)\
 eqn. (B.16)

 \ \color{black}

\noindent 71: (dims,levels) = $(3 , 
2 , 
1;8,
8,
1
)$,
irreps = $3_{8}^{3,0}\oplus
2_{8}^{1,0}\oplus
1_{1}^{1}$,
pord$(\rho_\text{isum}(\mathfrak{t})) = 8$,

\vskip 0.7ex
\hangindent=5.5em \hangafter=1
{\white .}\hskip 1em $\rho_\text{isum}(\mathfrak{t})$ =
 $( 0,
\frac{3}{8},
\frac{7}{8} )
\oplus
( \frac{1}{8},
\frac{3}{8} )
\oplus
( 0 )
$,

\vskip 0.7ex
\hangindent=5.5em \hangafter=1
{\white .}\hskip 1em $\rho_\text{isum}(\mathfrak{s})$ =
$\mathrm{i}$($0$,
$\sqrt{\frac{1}{2}}$,
$\sqrt{\frac{1}{2}}$;\ \ 
$\frac{1}{2}$,
$-\frac{1}{2}$;\ \ 
$\frac{1}{2}$)
 $\oplus$
($-\sqrt{\frac{1}{2}}$,
$\sqrt{\frac{1}{2}}$;
$\sqrt{\frac{1}{2}}$)
 $\oplus$
($1$)

Fail:
number of self dual objects $|$Tr($\rho(\mathfrak s^2)$)$|$ = 0. Prop. B.4 (1)\
 eqn. (B.16)

 \ \color{black}

 \color{blue}

\noindent 72: (dims,levels) = $(3 , 
2 , 
1;8,
8,
1
)$,
irreps = $3_{8}^{3,0}\oplus
2_{8}^{1,9}\oplus
1_{1}^{1}$,
pord$(\rho_\text{isum}(\mathfrak{t})) = 8$,

\vskip 0.7ex
\hangindent=5.5em \hangafter=1
{\white .}\hskip 1em $\rho_\text{isum}(\mathfrak{t})$ =
 $( 0,
\frac{3}{8},
\frac{7}{8} )
\oplus
( \frac{1}{8},
\frac{7}{8} )
\oplus
( 0 )
$,

\vskip 0.7ex
\hangindent=5.5em \hangafter=1
{\white .}\hskip 1em $\rho_\text{isum}(\mathfrak{s})$ =
$\mathrm{i}$($0$,
$\sqrt{\frac{1}{2}}$,
$\sqrt{\frac{1}{2}}$;\ \ 
$\frac{1}{2}$,
$-\frac{1}{2}$;\ \ 
$\frac{1}{2}$)
 $\oplus$
$\mathrm{i}$($-\sqrt{\frac{1}{2}}$,
$\sqrt{\frac{1}{2}}$;\ \ 
$\sqrt{\frac{1}{2}}$)
 $\oplus$
($1$)

Pass. 

 \ \color{black}

 \color{blue}

\noindent 73: (dims,levels) = $(3 , 
2 , 
1;8,
8,
1
)$,
irreps = $3_{8}^{3,0}\oplus
2_{8}^{1,3}\oplus
1_{1}^{1}$,
pord$(\rho_\text{isum}(\mathfrak{t})) = 8$,

\vskip 0.7ex
\hangindent=5.5em \hangafter=1
{\white .}\hskip 1em $\rho_\text{isum}(\mathfrak{t})$ =
 $( 0,
\frac{3}{8},
\frac{7}{8} )
\oplus
( \frac{3}{8},
\frac{5}{8} )
\oplus
( 0 )
$,

\vskip 0.7ex
\hangindent=5.5em \hangafter=1
{\white .}\hskip 1em $\rho_\text{isum}(\mathfrak{s})$ =
$\mathrm{i}$($0$,
$\sqrt{\frac{1}{2}}$,
$\sqrt{\frac{1}{2}}$;\ \ 
$\frac{1}{2}$,
$-\frac{1}{2}$;\ \ 
$\frac{1}{2}$)
 $\oplus$
$\mathrm{i}$($-\sqrt{\frac{1}{2}}$,
$\sqrt{\frac{1}{2}}$;\ \ 
$\sqrt{\frac{1}{2}}$)
 $\oplus$
($1$)

Pass. 

 \ \color{black}

\noindent 74: (dims,levels) = $(3 , 
2 , 
1;8,
8,
1
)$,
irreps = $3_{8}^{3,0}\oplus
2_{8}^{1,6}\oplus
1_{1}^{1}$,
pord$(\rho_\text{isum}(\mathfrak{t})) = 8$,

\vskip 0.7ex
\hangindent=5.5em \hangafter=1
{\white .}\hskip 1em $\rho_\text{isum}(\mathfrak{t})$ =
 $( 0,
\frac{3}{8},
\frac{7}{8} )
\oplus
( \frac{5}{8},
\frac{7}{8} )
\oplus
( 0 )
$,

\vskip 0.7ex
\hangindent=5.5em \hangafter=1
{\white .}\hskip 1em $\rho_\text{isum}(\mathfrak{s})$ =
$\mathrm{i}$($0$,
$\sqrt{\frac{1}{2}}$,
$\sqrt{\frac{1}{2}}$;\ \ 
$\frac{1}{2}$,
$-\frac{1}{2}$;\ \ 
$\frac{1}{2}$)
 $\oplus$
($\sqrt{\frac{1}{2}}$,
$\sqrt{\frac{1}{2}}$;
$-\sqrt{\frac{1}{2}}$)
 $\oplus$
($1$)

Fail:
number of self dual objects $|$Tr($\rho(\mathfrak s^2)$)$|$ = 0. Prop. B.4 (1)\
 eqn. (B.16)

 \ \color{black}

 \color{blue}

\noindent 75: (dims,levels) = $(3 , 
2 , 
1;10,
2,
1
)$,
irreps = $3_{5}^{3}
\hskip -1.5pt \otimes \hskip -1.5pt
1_{2}^{1,0}\oplus
2_{2}^{1,0}\oplus
1_{1}^{1}$,
pord$(\rho_\text{isum}(\mathfrak{t})) = 10$,

\vskip 0.7ex
\hangindent=5.5em \hangafter=1
{\white .}\hskip 1em $\rho_\text{isum}(\mathfrak{t})$ =
 $( \frac{1}{2},
\frac{1}{10},
\frac{9}{10} )
\oplus
( 0,
\frac{1}{2} )
\oplus
( 0 )
$,

\vskip 0.7ex
\hangindent=5.5em \hangafter=1
{\white .}\hskip 1em $\rho_\text{isum}(\mathfrak{s})$ =
($\sqrt{\frac{1}{5}}$,
$-\sqrt{\frac{2}{5}}$,
$-\sqrt{\frac{2}{5}}$;
$\frac{5-\sqrt{5}}{10}$,
$-\frac{5+\sqrt{5}}{10}$;
$\frac{5-\sqrt{5}}{10}$)
 $\oplus$
($-\frac{1}{2}$,
$-\sqrt{\frac{3}{4}}$;
$\frac{1}{2}$)
 $\oplus$
($1$)

Pass. 

 \ \color{black}

 \color{blue}

\noindent 76: (dims,levels) = $(3 , 
2 , 
1;10,
2,
1
)$,
irreps = $3_{5}^{1}
\hskip -1.5pt \otimes \hskip -1.5pt
1_{2}^{1,0}\oplus
2_{2}^{1,0}\oplus
1_{1}^{1}$,
pord$(\rho_\text{isum}(\mathfrak{t})) = 10$,

\vskip 0.7ex
\hangindent=5.5em \hangafter=1
{\white .}\hskip 1em $\rho_\text{isum}(\mathfrak{t})$ =
 $( \frac{1}{2},
\frac{3}{10},
\frac{7}{10} )
\oplus
( 0,
\frac{1}{2} )
\oplus
( 0 )
$,

\vskip 0.7ex
\hangindent=5.5em \hangafter=1
{\white .}\hskip 1em $\rho_\text{isum}(\mathfrak{s})$ =
($-\sqrt{\frac{1}{5}}$,
$-\sqrt{\frac{2}{5}}$,
$-\sqrt{\frac{2}{5}}$;
$\frac{5+\sqrt{5}}{10}$,
$\frac{-5+\sqrt{5}}{10}$;
$\frac{5+\sqrt{5}}{10}$)
 $\oplus$
($-\frac{1}{2}$,
$-\sqrt{\frac{3}{4}}$;
$\frac{1}{2}$)
 $\oplus$
($1$)

Pass. 

 \ \color{black}

 \color{blue}

\noindent 77: (dims,levels) = $(3 , 
2 , 
1;10,
2,
2
)$,
irreps = $3_{5}^{3}
\hskip -1.5pt \otimes \hskip -1.5pt
1_{2}^{1,0}\oplus
2_{2}^{1,0}\oplus
1_{2}^{1,0}$,
pord$(\rho_\text{isum}(\mathfrak{t})) = 10$,

\vskip 0.7ex
\hangindent=5.5em \hangafter=1
{\white .}\hskip 1em $\rho_\text{isum}(\mathfrak{t})$ =
 $( \frac{1}{2},
\frac{1}{10},
\frac{9}{10} )
\oplus
( 0,
\frac{1}{2} )
\oplus
( \frac{1}{2} )
$,

\vskip 0.7ex
\hangindent=5.5em \hangafter=1
{\white .}\hskip 1em $\rho_\text{isum}(\mathfrak{s})$ =
($\sqrt{\frac{1}{5}}$,
$-\sqrt{\frac{2}{5}}$,
$-\sqrt{\frac{2}{5}}$;
$\frac{5-\sqrt{5}}{10}$,
$-\frac{5+\sqrt{5}}{10}$;
$\frac{5-\sqrt{5}}{10}$)
 $\oplus$
($-\frac{1}{2}$,
$-\sqrt{\frac{3}{4}}$;
$\frac{1}{2}$)
 $\oplus$
($-1$)

Pass. 

 \ \color{black}

 \color{blue}

\noindent 78: (dims,levels) = $(3 , 
2 , 
1;10,
2,
2
)$,
irreps = $3_{5}^{1}
\hskip -1.5pt \otimes \hskip -1.5pt
1_{2}^{1,0}\oplus
2_{2}^{1,0}\oplus
1_{2}^{1,0}$,
pord$(\rho_\text{isum}(\mathfrak{t})) = 10$,

\vskip 0.7ex
\hangindent=5.5em \hangafter=1
{\white .}\hskip 1em $\rho_\text{isum}(\mathfrak{t})$ =
 $( \frac{1}{2},
\frac{3}{10},
\frac{7}{10} )
\oplus
( 0,
\frac{1}{2} )
\oplus
( \frac{1}{2} )
$,

\vskip 0.7ex
\hangindent=5.5em \hangafter=1
{\white .}\hskip 1em $\rho_\text{isum}(\mathfrak{s})$ =
($-\sqrt{\frac{1}{5}}$,
$-\sqrt{\frac{2}{5}}$,
$-\sqrt{\frac{2}{5}}$;
$\frac{5+\sqrt{5}}{10}$,
$\frac{-5+\sqrt{5}}{10}$;
$\frac{5+\sqrt{5}}{10}$)
 $\oplus$
($-\frac{1}{2}$,
$-\sqrt{\frac{3}{4}}$;
$\frac{1}{2}$)
 $\oplus$
($-1$)

Pass. 

 \ \color{black}

\noindent 79: (dims,levels) = $(3 , 
2 , 
1;10,
6,
2
)$,
irreps = $3_{5}^{3}
\hskip -1.5pt \otimes \hskip -1.5pt
1_{2}^{1,0}\oplus
2_{3}^{1,8}
\hskip -1.5pt \otimes \hskip -1.5pt
1_{2}^{1,0}\oplus
1_{2}^{1,0}$,
pord$(\rho_\text{isum}(\mathfrak{t})) = 15$,

\vskip 0.7ex
\hangindent=5.5em \hangafter=1
{\white .}\hskip 1em $\rho_\text{isum}(\mathfrak{t})$ =
 $( \frac{1}{2},
\frac{1}{10},
\frac{9}{10} )
\oplus
( \frac{1}{2},
\frac{1}{6} )
\oplus
( \frac{1}{2} )
$,

\vskip 0.7ex
\hangindent=5.5em \hangafter=1
{\white .}\hskip 1em $\rho_\text{isum}(\mathfrak{s})$ =
($\sqrt{\frac{1}{5}}$,
$-\sqrt{\frac{2}{5}}$,
$-\sqrt{\frac{2}{5}}$;
$\frac{5-\sqrt{5}}{10}$,
$-\frac{5+\sqrt{5}}{10}$;
$\frac{5-\sqrt{5}}{10}$)
 $\oplus$
$\mathrm{i}$($-\sqrt{\frac{1}{3}}$,
$\sqrt{\frac{2}{3}}$;\ \ 
$\sqrt{\frac{1}{3}}$)
 $\oplus$
($-1$)

Fail:
Tr$_I(C) = -1 <$  0 for I = [ 1/6 ]. Prop. B.4 (1) eqn. (B.18)

 \ \color{black}

\noindent 80: (dims,levels) = $(3 , 
2 , 
1;10,
6,
2
)$,
irreps = $3_{5}^{3}
\hskip -1.5pt \otimes \hskip -1.5pt
1_{2}^{1,0}\oplus
2_{3}^{1,0}
\hskip -1.5pt \otimes \hskip -1.5pt
1_{2}^{1,0}\oplus
1_{2}^{1,0}$,
pord$(\rho_\text{isum}(\mathfrak{t})) = 15$,

\vskip 0.7ex
\hangindent=5.5em \hangafter=1
{\white .}\hskip 1em $\rho_\text{isum}(\mathfrak{t})$ =
 $( \frac{1}{2},
\frac{1}{10},
\frac{9}{10} )
\oplus
( \frac{1}{2},
\frac{5}{6} )
\oplus
( \frac{1}{2} )
$,

\vskip 0.7ex
\hangindent=5.5em \hangafter=1
{\white .}\hskip 1em $\rho_\text{isum}(\mathfrak{s})$ =
($\sqrt{\frac{1}{5}}$,
$-\sqrt{\frac{2}{5}}$,
$-\sqrt{\frac{2}{5}}$;
$\frac{5-\sqrt{5}}{10}$,
$-\frac{5+\sqrt{5}}{10}$;
$\frac{5-\sqrt{5}}{10}$)
 $\oplus$
$\mathrm{i}$($\sqrt{\frac{1}{3}}$,
$\sqrt{\frac{2}{3}}$;\ \ 
$-\sqrt{\frac{1}{3}}$)
 $\oplus$
($-1$)

Fail:
Tr$_I(C) = -1 <$  0 for I = [ 5/6 ]. Prop. B.4 (1) eqn. (B.18)

 \ \color{black}

\noindent 81: (dims,levels) = $(3 , 
2 , 
1;10,
6,
2
)$,
irreps = $3_{5}^{1}
\hskip -1.5pt \otimes \hskip -1.5pt
1_{2}^{1,0}\oplus
2_{3}^{1,8}
\hskip -1.5pt \otimes \hskip -1.5pt
1_{2}^{1,0}\oplus
1_{2}^{1,0}$,
pord$(\rho_\text{isum}(\mathfrak{t})) = 15$,

\vskip 0.7ex
\hangindent=5.5em \hangafter=1
{\white .}\hskip 1em $\rho_\text{isum}(\mathfrak{t})$ =
 $( \frac{1}{2},
\frac{3}{10},
\frac{7}{10} )
\oplus
( \frac{1}{2},
\frac{1}{6} )
\oplus
( \frac{1}{2} )
$,

\vskip 0.7ex
\hangindent=5.5em \hangafter=1
{\white .}\hskip 1em $\rho_\text{isum}(\mathfrak{s})$ =
($-\sqrt{\frac{1}{5}}$,
$-\sqrt{\frac{2}{5}}$,
$-\sqrt{\frac{2}{5}}$;
$\frac{5+\sqrt{5}}{10}$,
$\frac{-5+\sqrt{5}}{10}$;
$\frac{5+\sqrt{5}}{10}$)
 $\oplus$
$\mathrm{i}$($-\sqrt{\frac{1}{3}}$,
$\sqrt{\frac{2}{3}}$;\ \ 
$\sqrt{\frac{1}{3}}$)
 $\oplus$
($-1$)

Fail:
Tr$_I(C) = -1 <$  0 for I = [ 1/6 ]. Prop. B.4 (1) eqn. (B.18)

 \ \color{black}

\noindent 82: (dims,levels) = $(3 , 
2 , 
1;10,
6,
2
)$,
irreps = $3_{5}^{1}
\hskip -1.5pt \otimes \hskip -1.5pt
1_{2}^{1,0}\oplus
2_{3}^{1,0}
\hskip -1.5pt \otimes \hskip -1.5pt
1_{2}^{1,0}\oplus
1_{2}^{1,0}$,
pord$(\rho_\text{isum}(\mathfrak{t})) = 15$,

\vskip 0.7ex
\hangindent=5.5em \hangafter=1
{\white .}\hskip 1em $\rho_\text{isum}(\mathfrak{t})$ =
 $( \frac{1}{2},
\frac{3}{10},
\frac{7}{10} )
\oplus
( \frac{1}{2},
\frac{5}{6} )
\oplus
( \frac{1}{2} )
$,

\vskip 0.7ex
\hangindent=5.5em \hangafter=1
{\white .}\hskip 1em $\rho_\text{isum}(\mathfrak{s})$ =
($-\sqrt{\frac{1}{5}}$,
$-\sqrt{\frac{2}{5}}$,
$-\sqrt{\frac{2}{5}}$;
$\frac{5+\sqrt{5}}{10}$,
$\frac{-5+\sqrt{5}}{10}$;
$\frac{5+\sqrt{5}}{10}$)
 $\oplus$
$\mathrm{i}$($\sqrt{\frac{1}{3}}$,
$\sqrt{\frac{2}{3}}$;\ \ 
$-\sqrt{\frac{1}{3}}$)
 $\oplus$
($-1$)

Fail:
Tr$_I(C) = -1 <$  0 for I = [ 5/6 ]. Prop. B.4 (1) eqn. (B.18)

 \ \color{black}

\noindent 83: (dims,levels) = $(3 , 
2 , 
1;10,
6,
6
)$,
irreps = $3_{5}^{3}
\hskip -1.5pt \otimes \hskip -1.5pt
1_{2}^{1,0}\oplus
2_{3}^{1,8}
\hskip -1.5pt \otimes \hskip -1.5pt
1_{2}^{1,0}\oplus
1_{3}^{1,4}
\hskip -1.5pt \otimes \hskip -1.5pt
1_{2}^{1,0}$,
pord$(\rho_\text{isum}(\mathfrak{t})) = 15$,

\vskip 0.7ex
\hangindent=5.5em \hangafter=1
{\white .}\hskip 1em $\rho_\text{isum}(\mathfrak{t})$ =
 $( \frac{1}{2},
\frac{1}{10},
\frac{9}{10} )
\oplus
( \frac{1}{2},
\frac{1}{6} )
\oplus
( \frac{1}{6} )
$,

\vskip 0.7ex
\hangindent=5.5em \hangafter=1
{\white .}\hskip 1em $\rho_\text{isum}(\mathfrak{s})$ =
($\sqrt{\frac{1}{5}}$,
$-\sqrt{\frac{2}{5}}$,
$-\sqrt{\frac{2}{5}}$;
$\frac{5-\sqrt{5}}{10}$,
$-\frac{5+\sqrt{5}}{10}$;
$\frac{5-\sqrt{5}}{10}$)
 $\oplus$
$\mathrm{i}$($-\sqrt{\frac{1}{3}}$,
$\sqrt{\frac{2}{3}}$;\ \ 
$\sqrt{\frac{1}{3}}$)
 $\oplus$
($-1$)

Fail:
all $\theta$-eigenspaces that can contain unit
 have Tr$_{E_\theta}(C) \leq 0 $. Prop. B.5 (5) eqn. (B.29)

 \ \color{black}

\noindent 84: (dims,levels) = $(3 , 
2 , 
1;10,
6,
6
)$,
irreps = $3_{5}^{3}
\hskip -1.5pt \otimes \hskip -1.5pt
1_{2}^{1,0}\oplus
2_{3}^{1,0}
\hskip -1.5pt \otimes \hskip -1.5pt
1_{2}^{1,0}\oplus
1_{3}^{1,0}
\hskip -1.5pt \otimes \hskip -1.5pt
1_{2}^{1,0}$,
pord$(\rho_\text{isum}(\mathfrak{t})) = 15$,

\vskip 0.7ex
\hangindent=5.5em \hangafter=1
{\white .}\hskip 1em $\rho_\text{isum}(\mathfrak{t})$ =
 $( \frac{1}{2},
\frac{1}{10},
\frac{9}{10} )
\oplus
( \frac{1}{2},
\frac{5}{6} )
\oplus
( \frac{5}{6} )
$,

\vskip 0.7ex
\hangindent=5.5em \hangafter=1
{\white .}\hskip 1em $\rho_\text{isum}(\mathfrak{s})$ =
($\sqrt{\frac{1}{5}}$,
$-\sqrt{\frac{2}{5}}$,
$-\sqrt{\frac{2}{5}}$;
$\frac{5-\sqrt{5}}{10}$,
$-\frac{5+\sqrt{5}}{10}$;
$\frac{5-\sqrt{5}}{10}$)
 $\oplus$
$\mathrm{i}$($\sqrt{\frac{1}{3}}$,
$\sqrt{\frac{2}{3}}$;\ \ 
$-\sqrt{\frac{1}{3}}$)
 $\oplus$
($-1$)

Fail:
all $\theta$-eigenspaces that can contain unit
 have Tr$_{E_\theta}(C) \leq 0 $. Prop. B.5 (5) eqn. (B.29)

 \ \color{black}

\noindent 85: (dims,levels) = $(3 , 
2 , 
1;10,
6,
6
)$,
irreps = $3_{5}^{1}
\hskip -1.5pt \otimes \hskip -1.5pt
1_{2}^{1,0}\oplus
2_{3}^{1,8}
\hskip -1.5pt \otimes \hskip -1.5pt
1_{2}^{1,0}\oplus
1_{3}^{1,4}
\hskip -1.5pt \otimes \hskip -1.5pt
1_{2}^{1,0}$,
pord$(\rho_\text{isum}(\mathfrak{t})) = 15$,

\vskip 0.7ex
\hangindent=5.5em \hangafter=1
{\white .}\hskip 1em $\rho_\text{isum}(\mathfrak{t})$ =
 $( \frac{1}{2},
\frac{3}{10},
\frac{7}{10} )
\oplus
( \frac{1}{2},
\frac{1}{6} )
\oplus
( \frac{1}{6} )
$,

\vskip 0.7ex
\hangindent=5.5em \hangafter=1
{\white .}\hskip 1em $\rho_\text{isum}(\mathfrak{s})$ =
($-\sqrt{\frac{1}{5}}$,
$-\sqrt{\frac{2}{5}}$,
$-\sqrt{\frac{2}{5}}$;
$\frac{5+\sqrt{5}}{10}$,
$\frac{-5+\sqrt{5}}{10}$;
$\frac{5+\sqrt{5}}{10}$)
 $\oplus$
$\mathrm{i}$($-\sqrt{\frac{1}{3}}$,
$\sqrt{\frac{2}{3}}$;\ \ 
$\sqrt{\frac{1}{3}}$)
 $\oplus$
($-1$)

Fail:
all $\theta$-eigenspaces that can contain unit
 have Tr$_{E_\theta}(C) \leq 0 $. Prop. B.5 (5) eqn. (B.29)

 \ \color{black}

\noindent 86: (dims,levels) = $(3 , 
2 , 
1;10,
6,
6
)$,
irreps = $3_{5}^{1}
\hskip -1.5pt \otimes \hskip -1.5pt
1_{2}^{1,0}\oplus
2_{3}^{1,0}
\hskip -1.5pt \otimes \hskip -1.5pt
1_{2}^{1,0}\oplus
1_{3}^{1,0}
\hskip -1.5pt \otimes \hskip -1.5pt
1_{2}^{1,0}$,
pord$(\rho_\text{isum}(\mathfrak{t})) = 15$,

\vskip 0.7ex
\hangindent=5.5em \hangafter=1
{\white .}\hskip 1em $\rho_\text{isum}(\mathfrak{t})$ =
 $( \frac{1}{2},
\frac{3}{10},
\frac{7}{10} )
\oplus
( \frac{1}{2},
\frac{5}{6} )
\oplus
( \frac{5}{6} )
$,

\vskip 0.7ex
\hangindent=5.5em \hangafter=1
{\white .}\hskip 1em $\rho_\text{isum}(\mathfrak{s})$ =
($-\sqrt{\frac{1}{5}}$,
$-\sqrt{\frac{2}{5}}$,
$-\sqrt{\frac{2}{5}}$;
$\frac{5+\sqrt{5}}{10}$,
$\frac{-5+\sqrt{5}}{10}$;
$\frac{5+\sqrt{5}}{10}$)
 $\oplus$
$\mathrm{i}$($\sqrt{\frac{1}{3}}$,
$\sqrt{\frac{2}{3}}$;\ \ 
$-\sqrt{\frac{1}{3}}$)
 $\oplus$
($-1$)

Fail:
all $\theta$-eigenspaces that can contain unit
 have Tr$_{E_\theta}(C) \leq 0 $. Prop. B.5 (5) eqn. (B.29)

 \ \color{black}

 \color{blue}

\noindent 87: (dims,levels) = $(3 , 
2 , 
1;10,
10,
2
)$,
irreps = $3_{5}^{3}
\hskip -1.5pt \otimes \hskip -1.5pt
1_{2}^{1,0}\oplus
2_{5}^{2}
\hskip -1.5pt \otimes \hskip -1.5pt
1_{2}^{1,0}\oplus
1_{2}^{1,0}$,
pord$(\rho_\text{isum}(\mathfrak{t})) = 5$,

\vskip 0.7ex
\hangindent=5.5em \hangafter=1
{\white .}\hskip 1em $\rho_\text{isum}(\mathfrak{t})$ =
 $( \frac{1}{2},
\frac{1}{10},
\frac{9}{10} )
\oplus
( \frac{1}{10},
\frac{9}{10} )
\oplus
( \frac{1}{2} )
$,

\vskip 0.7ex
\hangindent=5.5em \hangafter=1
{\white .}\hskip 1em $\rho_\text{isum}(\mathfrak{s})$ =
($\sqrt{\frac{1}{5}}$,
$-\sqrt{\frac{2}{5}}$,
$-\sqrt{\frac{2}{5}}$;
$\frac{5-\sqrt{5}}{10}$,
$-\frac{5+\sqrt{5}}{10}$;
$\frac{5-\sqrt{5}}{10}$)
 $\oplus$
$\mathrm{i}$($-\frac{1}{\sqrt{5}}c^{1}_{20}
$,
$\frac{1}{\sqrt{5}}c^{3}_{20}
$;\ \ 
$\frac{1}{\sqrt{5}}c^{1}_{20}
$)
 $\oplus$
($-1$)

Pass. 

 \ \color{black}

 \color{blue}

\noindent 88: (dims,levels) = $(3 , 
2 , 
1;10,
10,
2
)$,
irreps = $3_{5}^{1}
\hskip -1.5pt \otimes \hskip -1.5pt
1_{2}^{1,0}\oplus
2_{5}^{1}
\hskip -1.5pt \otimes \hskip -1.5pt
1_{2}^{1,0}\oplus
1_{2}^{1,0}$,
pord$(\rho_\text{isum}(\mathfrak{t})) = 5$,

\vskip 0.7ex
\hangindent=5.5em \hangafter=1
{\white .}\hskip 1em $\rho_\text{isum}(\mathfrak{t})$ =
 $( \frac{1}{2},
\frac{3}{10},
\frac{7}{10} )
\oplus
( \frac{3}{10},
\frac{7}{10} )
\oplus
( \frac{1}{2} )
$,

\vskip 0.7ex
\hangindent=5.5em \hangafter=1
{\white .}\hskip 1em $\rho_\text{isum}(\mathfrak{s})$ =
($-\sqrt{\frac{1}{5}}$,
$-\sqrt{\frac{2}{5}}$,
$-\sqrt{\frac{2}{5}}$;
$\frac{5+\sqrt{5}}{10}$,
$\frac{-5+\sqrt{5}}{10}$;
$\frac{5+\sqrt{5}}{10}$)
 $\oplus$
$\mathrm{i}$($-\frac{1}{\sqrt{5}}c^{3}_{20}
$,
$\frac{1}{\sqrt{5}}c^{1}_{20}
$;\ \ 
$\frac{1}{\sqrt{5}}c^{3}_{20}
$)
 $\oplus$
($-1$)

Pass. 

 \ \color{black}

\noindent 89: (dims,levels) = $(3 , 
2 , 
1;12,
3,
1
)$,
irreps = $3_{4}^{1,0}
\hskip -1.5pt \otimes \hskip -1.5pt
1_{3}^{1,0}\oplus
2_{3}^{1,0}\oplus
1_{1}^{1}$,
pord$(\rho_\text{isum}(\mathfrak{t})) = 12$,

\vskip 0.7ex
\hangindent=5.5em \hangafter=1
{\white .}\hskip 1em $\rho_\text{isum}(\mathfrak{t})$ =
 $( \frac{1}{3},
\frac{1}{12},
\frac{7}{12} )
\oplus
( 0,
\frac{1}{3} )
\oplus
( 0 )
$,

\vskip 0.7ex
\hangindent=5.5em \hangafter=1
{\white .}\hskip 1em $\rho_\text{isum}(\mathfrak{s})$ =
($0$,
$\sqrt{\frac{1}{2}}$,
$\sqrt{\frac{1}{2}}$;
$-\frac{1}{2}$,
$\frac{1}{2}$;
$-\frac{1}{2}$)
 $\oplus$
$\mathrm{i}$($-\sqrt{\frac{1}{3}}$,
$\sqrt{\frac{2}{3}}$;\ \ 
$\sqrt{\frac{1}{3}}$)
 $\oplus$
($1$)

Fail:
all $\theta$-eigenspaces that can contain unit
 have Tr$_{E_\theta}(C) \leq 0 $. Prop. B.5 (5) eqn. (B.29)

 \ \color{black}

\noindent 90: (dims,levels) = $(3 , 
2 , 
1;12,
3,
3
)$,
irreps = $3_{4}^{1,0}
\hskip -1.5pt \otimes \hskip -1.5pt
1_{3}^{1,0}\oplus
2_{3}^{1,0}\oplus
1_{3}^{1,0}$,
pord$(\rho_\text{isum}(\mathfrak{t})) = 12$,

\vskip 0.7ex
\hangindent=5.5em \hangafter=1
{\white .}\hskip 1em $\rho_\text{isum}(\mathfrak{t})$ =
 $( \frac{1}{3},
\frac{1}{12},
\frac{7}{12} )
\oplus
( 0,
\frac{1}{3} )
\oplus
( \frac{1}{3} )
$,

\vskip 0.7ex
\hangindent=5.5em \hangafter=1
{\white .}\hskip 1em $\rho_\text{isum}(\mathfrak{s})$ =
($0$,
$\sqrt{\frac{1}{2}}$,
$\sqrt{\frac{1}{2}}$;
$-\frac{1}{2}$,
$\frac{1}{2}$;
$-\frac{1}{2}$)
 $\oplus$
$\mathrm{i}$($-\sqrt{\frac{1}{3}}$,
$\sqrt{\frac{2}{3}}$;\ \ 
$\sqrt{\frac{1}{3}}$)
 $\oplus$
($1$)

Fail:
Tr$_I(C) = -1 <$  0 for I = [ 0 ]. Prop. B.4 (1) eqn. (B.18)

 \ \color{black}

\noindent 91: (dims,levels) = $(3 , 
2 , 
1;12,
3,
3
)$,
irreps = $3_{4}^{1,0}
\hskip -1.5pt \otimes \hskip -1.5pt
1_{3}^{1,0}\oplus
2_{3}^{1,4}\oplus
1_{3}^{1,0}$,
pord$(\rho_\text{isum}(\mathfrak{t})) = 12$,

\vskip 0.7ex
\hangindent=5.5em \hangafter=1
{\white .}\hskip 1em $\rho_\text{isum}(\mathfrak{t})$ =
 $( \frac{1}{3},
\frac{1}{12},
\frac{7}{12} )
\oplus
( \frac{1}{3},
\frac{2}{3} )
\oplus
( \frac{1}{3} )
$,

\vskip 0.7ex
\hangindent=5.5em \hangafter=1
{\white .}\hskip 1em $\rho_\text{isum}(\mathfrak{s})$ =
($0$,
$\sqrt{\frac{1}{2}}$,
$\sqrt{\frac{1}{2}}$;
$-\frac{1}{2}$,
$\frac{1}{2}$;
$-\frac{1}{2}$)
 $\oplus$
$\mathrm{i}$($-\sqrt{\frac{1}{3}}$,
$\sqrt{\frac{2}{3}}$;\ \ 
$\sqrt{\frac{1}{3}}$)
 $\oplus$
($1$)

Fail:
Tr$_I(C) = -1 <$  0 for I = [ 2/3 ]. Prop. B.4 (1) eqn. (B.18)

 \ \color{black}

\noindent 92: (dims,levels) = $(3 , 
2 , 
1;12,
3,
3
)$,
irreps = $3_{4}^{1,0}
\hskip -1.5pt \otimes \hskip -1.5pt
1_{3}^{1,0}\oplus
2_{3}^{1,4}\oplus
1_{3}^{1,4}$,
pord$(\rho_\text{isum}(\mathfrak{t})) = 12$,

\vskip 0.7ex
\hangindent=5.5em \hangafter=1
{\white .}\hskip 1em $\rho_\text{isum}(\mathfrak{t})$ =
 $( \frac{1}{3},
\frac{1}{12},
\frac{7}{12} )
\oplus
( \frac{1}{3},
\frac{2}{3} )
\oplus
( \frac{2}{3} )
$,

\vskip 0.7ex
\hangindent=5.5em \hangafter=1
{\white .}\hskip 1em $\rho_\text{isum}(\mathfrak{s})$ =
($0$,
$\sqrt{\frac{1}{2}}$,
$\sqrt{\frac{1}{2}}$;
$-\frac{1}{2}$,
$\frac{1}{2}$;
$-\frac{1}{2}$)
 $\oplus$
$\mathrm{i}$($-\sqrt{\frac{1}{3}}$,
$\sqrt{\frac{2}{3}}$;\ \ 
$\sqrt{\frac{1}{3}}$)
 $\oplus$
($1$)

Fail:
all $\theta$-eigenspaces that can contain unit
 have Tr$_{E_\theta}(C) \leq 0 $. Prop. B.5 (5) eqn. (B.29)

 \ \color{black}

\noindent 93: (dims,levels) = $(3 , 
2 , 
1;12,
3,
12
)$,
irreps = $3_{4}^{1,0}
\hskip -1.5pt \otimes \hskip -1.5pt
1_{3}^{1,0}\oplus
2_{3}^{1,0}\oplus
1_{4}^{1,6}
\hskip -1.5pt \otimes \hskip -1.5pt
1_{3}^{1,0}$,
pord$(\rho_\text{isum}(\mathfrak{t})) = 12$,

\vskip 0.7ex
\hangindent=5.5em \hangafter=1
{\white .}\hskip 1em $\rho_\text{isum}(\mathfrak{t})$ =
 $( \frac{1}{3},
\frac{1}{12},
\frac{7}{12} )
\oplus
( 0,
\frac{1}{3} )
\oplus
( \frac{1}{12} )
$,

\vskip 0.7ex
\hangindent=5.5em \hangafter=1
{\white .}\hskip 1em $\rho_\text{isum}(\mathfrak{s})$ =
($0$,
$\sqrt{\frac{1}{2}}$,
$\sqrt{\frac{1}{2}}$;
$-\frac{1}{2}$,
$\frac{1}{2}$;
$-\frac{1}{2}$)
 $\oplus$
$\mathrm{i}$($-\sqrt{\frac{1}{3}}$,
$\sqrt{\frac{2}{3}}$;\ \ 
$\sqrt{\frac{1}{3}}$)
 $\oplus$
$\mathrm{i}$($-1$)

Fail:
number of self dual objects $|$Tr($\rho(\mathfrak s^2)$)$|$ = 0. Prop. B.4 (1)\
 eqn. (B.16)

 \ \color{black}

\noindent 94: (dims,levels) = $(3 , 
2 , 
1;12,
3,
12
)$,
irreps = $3_{4}^{1,0}
\hskip -1.5pt \otimes \hskip -1.5pt
1_{3}^{1,0}\oplus
2_{3}^{1,0}\oplus
1_{4}^{1,0}
\hskip -1.5pt \otimes \hskip -1.5pt
1_{3}^{1,0}$,
pord$(\rho_\text{isum}(\mathfrak{t})) = 12$,

\vskip 0.7ex
\hangindent=5.5em \hangafter=1
{\white .}\hskip 1em $\rho_\text{isum}(\mathfrak{t})$ =
 $( \frac{1}{3},
\frac{1}{12},
\frac{7}{12} )
\oplus
( 0,
\frac{1}{3} )
\oplus
( \frac{7}{12} )
$,

\vskip 0.7ex
\hangindent=5.5em \hangafter=1
{\white .}\hskip 1em $\rho_\text{isum}(\mathfrak{s})$ =
($0$,
$\sqrt{\frac{1}{2}}$,
$\sqrt{\frac{1}{2}}$;
$-\frac{1}{2}$,
$\frac{1}{2}$;
$-\frac{1}{2}$)
 $\oplus$
$\mathrm{i}$($-\sqrt{\frac{1}{3}}$,
$\sqrt{\frac{2}{3}}$;\ \ 
$\sqrt{\frac{1}{3}}$)
 $\oplus$
$\mathrm{i}$($1$)

Fail:
number of self dual objects $|$Tr($\rho(\mathfrak s^2)$)$|$ = 0. Prop. B.4 (1)\
 eqn. (B.16)

 \ \color{black}

\noindent 95: (dims,levels) = $(3 , 
2 , 
1;12,
3,
12
)$,
irreps = $3_{4}^{1,0}
\hskip -1.5pt \otimes \hskip -1.5pt
1_{3}^{1,0}\oplus
2_{3}^{1,4}\oplus
1_{4}^{1,6}
\hskip -1.5pt \otimes \hskip -1.5pt
1_{3}^{1,0}$,
pord$(\rho_\text{isum}(\mathfrak{t})) = 12$,

\vskip 0.7ex
\hangindent=5.5em \hangafter=1
{\white .}\hskip 1em $\rho_\text{isum}(\mathfrak{t})$ =
 $( \frac{1}{3},
\frac{1}{12},
\frac{7}{12} )
\oplus
( \frac{1}{3},
\frac{2}{3} )
\oplus
( \frac{1}{12} )
$,

\vskip 0.7ex
\hangindent=5.5em \hangafter=1
{\white .}\hskip 1em $\rho_\text{isum}(\mathfrak{s})$ =
($0$,
$\sqrt{\frac{1}{2}}$,
$\sqrt{\frac{1}{2}}$;
$-\frac{1}{2}$,
$\frac{1}{2}$;
$-\frac{1}{2}$)
 $\oplus$
$\mathrm{i}$($-\sqrt{\frac{1}{3}}$,
$\sqrt{\frac{2}{3}}$;\ \ 
$\sqrt{\frac{1}{3}}$)
 $\oplus$
$\mathrm{i}$($-1$)

Fail:
number of self dual objects $|$Tr($\rho(\mathfrak s^2)$)$|$ = 0. Prop. B.4 (1)\
 eqn. (B.16)

 \ \color{black}

\noindent 96: (dims,levels) = $(3 , 
2 , 
1;12,
3,
12
)$,
irreps = $3_{4}^{1,0}
\hskip -1.5pt \otimes \hskip -1.5pt
1_{3}^{1,0}\oplus
2_{3}^{1,4}\oplus
1_{4}^{1,0}
\hskip -1.5pt \otimes \hskip -1.5pt
1_{3}^{1,0}$,
pord$(\rho_\text{isum}(\mathfrak{t})) = 12$,

\vskip 0.7ex
\hangindent=5.5em \hangafter=1
{\white .}\hskip 1em $\rho_\text{isum}(\mathfrak{t})$ =
 $( \frac{1}{3},
\frac{1}{12},
\frac{7}{12} )
\oplus
( \frac{1}{3},
\frac{2}{3} )
\oplus
( \frac{7}{12} )
$,

\vskip 0.7ex
\hangindent=5.5em \hangafter=1
{\white .}\hskip 1em $\rho_\text{isum}(\mathfrak{s})$ =
($0$,
$\sqrt{\frac{1}{2}}$,
$\sqrt{\frac{1}{2}}$;
$-\frac{1}{2}$,
$\frac{1}{2}$;
$-\frac{1}{2}$)
 $\oplus$
$\mathrm{i}$($-\sqrt{\frac{1}{3}}$,
$\sqrt{\frac{2}{3}}$;\ \ 
$\sqrt{\frac{1}{3}}$)
 $\oplus$
$\mathrm{i}$($1$)

Fail:
number of self dual objects $|$Tr($\rho(\mathfrak s^2)$)$|$ = 0. Prop. B.4 (1)\
 eqn. (B.16)

 \ \color{black}

 \color{blue}

\noindent 97: (dims,levels) = $(3 , 
2 , 
1;12,
12,
3
)$,
irreps = $3_{4}^{1,0}
\hskip -1.5pt \otimes \hskip -1.5pt
1_{3}^{1,0}\oplus
2_{3}^{1,4}
\hskip -1.5pt \otimes \hskip -1.5pt
1_{4}^{1,6}\oplus
1_{3}^{1,0}$,
pord$(\rho_\text{isum}(\mathfrak{t})) = 12$,

\vskip 0.7ex
\hangindent=5.5em \hangafter=1
{\white .}\hskip 1em $\rho_\text{isum}(\mathfrak{t})$ =
 $( \frac{1}{3},
\frac{1}{12},
\frac{7}{12} )
\oplus
( \frac{1}{12},
\frac{5}{12} )
\oplus
( \frac{1}{3} )
$,

\vskip 0.7ex
\hangindent=5.5em \hangafter=1
{\white .}\hskip 1em $\rho_\text{isum}(\mathfrak{s})$ =
($0$,
$\sqrt{\frac{1}{2}}$,
$\sqrt{\frac{1}{2}}$;
$-\frac{1}{2}$,
$\frac{1}{2}$;
$-\frac{1}{2}$)
 $\oplus$
($-\sqrt{\frac{1}{3}}$,
$\sqrt{\frac{2}{3}}$;
$\sqrt{\frac{1}{3}}$)
 $\oplus$
($1$)

Pass. 

 \ \color{black}

 \color{blue}

\noindent 98: (dims,levels) = $(3 , 
2 , 
1;12,
12,
3
)$,
irreps = $3_{4}^{1,0}
\hskip -1.5pt \otimes \hskip -1.5pt
1_{3}^{1,0}\oplus
2_{3}^{1,0}
\hskip -1.5pt \otimes \hskip -1.5pt
1_{4}^{1,0}\oplus
1_{3}^{1,0}$,
pord$(\rho_\text{isum}(\mathfrak{t})) = 12$,

\vskip 0.7ex
\hangindent=5.5em \hangafter=1
{\white .}\hskip 1em $\rho_\text{isum}(\mathfrak{t})$ =
 $( \frac{1}{3},
\frac{1}{12},
\frac{7}{12} )
\oplus
( \frac{1}{4},
\frac{7}{12} )
\oplus
( \frac{1}{3} )
$,

\vskip 0.7ex
\hangindent=5.5em \hangafter=1
{\white .}\hskip 1em $\rho_\text{isum}(\mathfrak{s})$ =
($0$,
$\sqrt{\frac{1}{2}}$,
$\sqrt{\frac{1}{2}}$;
$-\frac{1}{2}$,
$\frac{1}{2}$;
$-\frac{1}{2}$)
 $\oplus$
($\sqrt{\frac{1}{3}}$,
$\sqrt{\frac{2}{3}}$;
$-\sqrt{\frac{1}{3}}$)
 $\oplus$
($1$)

Pass. 

 \ \color{black}

 \color{blue}

\noindent 99: (dims,levels) = $(3 , 
2 , 
1;12,
12,
3
)$,
irreps = $3_{4}^{1,0}
\hskip -1.5pt \otimes \hskip -1.5pt
1_{3}^{1,0}\oplus
2_{3}^{1,4}
\hskip -1.5pt \otimes \hskip -1.5pt
1_{4}^{1,0}\oplus
1_{3}^{1,0}$,
pord$(\rho_\text{isum}(\mathfrak{t})) = 12$,

\vskip 0.7ex
\hangindent=5.5em \hangafter=1
{\white .}\hskip 1em $\rho_\text{isum}(\mathfrak{t})$ =
 $( \frac{1}{3},
\frac{1}{12},
\frac{7}{12} )
\oplus
( \frac{7}{12},
\frac{11}{12} )
\oplus
( \frac{1}{3} )
$,

\vskip 0.7ex
\hangindent=5.5em \hangafter=1
{\white .}\hskip 1em $\rho_\text{isum}(\mathfrak{s})$ =
($0$,
$\sqrt{\frac{1}{2}}$,
$\sqrt{\frac{1}{2}}$;
$-\frac{1}{2}$,
$\frac{1}{2}$;
$-\frac{1}{2}$)
 $\oplus$
($\sqrt{\frac{1}{3}}$,
$\sqrt{\frac{2}{3}}$;
$-\sqrt{\frac{1}{3}}$)
 $\oplus$
($1$)

Pass. 

 \ \color{black}

 \color{blue}

\noindent 100: (dims,levels) = $(3 , 
2 , 
1;12,
12,
3
)$,
irreps = $3_{4}^{1,0}
\hskip -1.5pt \otimes \hskip -1.5pt
1_{3}^{1,0}\oplus
2_{3}^{1,0}
\hskip -1.5pt \otimes \hskip -1.5pt
1_{4}^{1,6}\oplus
1_{3}^{1,0}$,
pord$(\rho_\text{isum}(\mathfrak{t})) = 12$,

\vskip 0.7ex
\hangindent=5.5em \hangafter=1
{\white .}\hskip 1em $\rho_\text{isum}(\mathfrak{t})$ =
 $( \frac{1}{3},
\frac{1}{12},
\frac{7}{12} )
\oplus
( \frac{3}{4},
\frac{1}{12} )
\oplus
( \frac{1}{3} )
$,

\vskip 0.7ex
\hangindent=5.5em \hangafter=1
{\white .}\hskip 1em $\rho_\text{isum}(\mathfrak{s})$ =
($0$,
$\sqrt{\frac{1}{2}}$,
$\sqrt{\frac{1}{2}}$;
$-\frac{1}{2}$,
$\frac{1}{2}$;
$-\frac{1}{2}$)
 $\oplus$
($-\sqrt{\frac{1}{3}}$,
$\sqrt{\frac{2}{3}}$;
$\sqrt{\frac{1}{3}}$)
 $\oplus$
($1$)

Pass. 

 \ \color{black}

\noindent 101: (dims,levels) = $(3 , 
2 , 
1;12,
12,
4
)$,
irreps = $3_{4}^{1,0}
\hskip -1.5pt \otimes \hskip -1.5pt
1_{3}^{1,0}\oplus
2_{3}^{1,0}
\hskip -1.5pt \otimes \hskip -1.5pt
1_{4}^{1,0}\oplus
1_{4}^{1,0}$,
pord$(\rho_\text{isum}(\mathfrak{t})) = 12$,

\vskip 0.7ex
\hangindent=5.5em \hangafter=1
{\white .}\hskip 1em $\rho_\text{isum}(\mathfrak{t})$ =
 $( \frac{1}{3},
\frac{1}{12},
\frac{7}{12} )
\oplus
( \frac{1}{4},
\frac{7}{12} )
\oplus
( \frac{1}{4} )
$,

\vskip 0.7ex
\hangindent=5.5em \hangafter=1
{\white .}\hskip 1em $\rho_\text{isum}(\mathfrak{s})$ =
($0$,
$\sqrt{\frac{1}{2}}$,
$\sqrt{\frac{1}{2}}$;
$-\frac{1}{2}$,
$\frac{1}{2}$;
$-\frac{1}{2}$)
 $\oplus$
($\sqrt{\frac{1}{3}}$,
$\sqrt{\frac{2}{3}}$;
$-\sqrt{\frac{1}{3}}$)
 $\oplus$
$\mathrm{i}$($1$)

Fail:
cnd($\rho(\mathfrak s)_\mathrm{ndeg}$) = 8 does not divide
 ord($\rho(\mathfrak t)$)=12. Prop. B.4 (2)

 \ \color{black}

\noindent 102: (dims,levels) = $(3 , 
2 , 
1;12,
12,
4
)$,
irreps = $3_{4}^{1,0}
\hskip -1.5pt \otimes \hskip -1.5pt
1_{3}^{1,0}\oplus
2_{3}^{1,0}
\hskip -1.5pt \otimes \hskip -1.5pt
1_{4}^{1,6}\oplus
1_{4}^{1,6}$,
pord$(\rho_\text{isum}(\mathfrak{t})) = 12$,

\vskip 0.7ex
\hangindent=5.5em \hangafter=1
{\white .}\hskip 1em $\rho_\text{isum}(\mathfrak{t})$ =
 $( \frac{1}{3},
\frac{1}{12},
\frac{7}{12} )
\oplus
( \frac{3}{4},
\frac{1}{12} )
\oplus
( \frac{3}{4} )
$,

\vskip 0.7ex
\hangindent=5.5em \hangafter=1
{\white .}\hskip 1em $\rho_\text{isum}(\mathfrak{s})$ =
($0$,
$\sqrt{\frac{1}{2}}$,
$\sqrt{\frac{1}{2}}$;
$-\frac{1}{2}$,
$\frac{1}{2}$;
$-\frac{1}{2}$)
 $\oplus$
($-\sqrt{\frac{1}{3}}$,
$\sqrt{\frac{2}{3}}$;
$\sqrt{\frac{1}{3}}$)
 $\oplus$
$\mathrm{i}$($-1$)

Fail:
cnd($\rho(\mathfrak s)_\mathrm{ndeg}$) = 8 does not divide
 ord($\rho(\mathfrak t)$)=12. Prop. B.4 (2)

 \ \color{black}

\noindent 103: (dims,levels) = $(3 , 
2 , 
1;12,
12,
12
)$,
irreps = $3_{4}^{1,0}
\hskip -1.5pt \otimes \hskip -1.5pt
1_{3}^{1,0}\oplus
2_{3}^{1,4}
\hskip -1.5pt \otimes \hskip -1.5pt
1_{4}^{1,6}\oplus
1_{4}^{1,6}
\hskip -1.5pt \otimes \hskip -1.5pt
1_{3}^{1,0}$,
pord$(\rho_\text{isum}(\mathfrak{t})) = 12$,

\vskip 0.7ex
\hangindent=5.5em \hangafter=1
{\white .}\hskip 1em $\rho_\text{isum}(\mathfrak{t})$ =
 $( \frac{1}{3},
\frac{1}{12},
\frac{7}{12} )
\oplus
( \frac{1}{12},
\frac{5}{12} )
\oplus
( \frac{1}{12} )
$,

\vskip 0.7ex
\hangindent=5.5em \hangafter=1
{\white .}\hskip 1em $\rho_\text{isum}(\mathfrak{s})$ =
($0$,
$\sqrt{\frac{1}{2}}$,
$\sqrt{\frac{1}{2}}$;
$-\frac{1}{2}$,
$\frac{1}{2}$;
$-\frac{1}{2}$)
 $\oplus$
($-\sqrt{\frac{1}{3}}$,
$\sqrt{\frac{2}{3}}$;
$\sqrt{\frac{1}{3}}$)
 $\oplus$
$\mathrm{i}$($-1$)

Fail:
cnd($\rho(\mathfrak s)_\mathrm{ndeg}$) = 24 does not divide
 ord($\rho(\mathfrak t)$)=12. Prop. B.4 (2)

 \ \color{black}

\noindent 104: (dims,levels) = $(3 , 
2 , 
1;12,
12,
12
)$,
irreps = $3_{4}^{1,0}
\hskip -1.5pt \otimes \hskip -1.5pt
1_{3}^{1,0}\oplus
2_{3}^{1,4}
\hskip -1.5pt \otimes \hskip -1.5pt
1_{4}^{1,6}\oplus
1_{4}^{1,6}
\hskip -1.5pt \otimes \hskip -1.5pt
1_{3}^{1,4}$,
pord$(\rho_\text{isum}(\mathfrak{t})) = 12$,

\vskip 0.7ex
\hangindent=5.5em \hangafter=1
{\white .}\hskip 1em $\rho_\text{isum}(\mathfrak{t})$ =
 $( \frac{1}{3},
\frac{1}{12},
\frac{7}{12} )
\oplus
( \frac{1}{12},
\frac{5}{12} )
\oplus
( \frac{5}{12} )
$,

\vskip 0.7ex
\hangindent=5.5em \hangafter=1
{\white .}\hskip 1em $\rho_\text{isum}(\mathfrak{s})$ =
($0$,
$\sqrt{\frac{1}{2}}$,
$\sqrt{\frac{1}{2}}$;
$-\frac{1}{2}$,
$\frac{1}{2}$;
$-\frac{1}{2}$)
 $\oplus$
($-\sqrt{\frac{1}{3}}$,
$\sqrt{\frac{2}{3}}$;
$\sqrt{\frac{1}{3}}$)
 $\oplus$
$\mathrm{i}$($-1$)

Fail:
cnd($\rho(\mathfrak s)_\mathrm{ndeg}$) = 8 does not divide
 ord($\rho(\mathfrak t)$)=12. Prop. B.4 (2)

 \ \color{black}

 \color{blue}

\noindent 105: (dims,levels) = $(3 , 
2 , 
1;12,
12,
12
)$,
irreps = $3_{4}^{1,0}
\hskip -1.5pt \otimes \hskip -1.5pt
1_{3}^{1,0}\oplus
2_{3}^{1,4}
\hskip -1.5pt \otimes \hskip -1.5pt
1_{4}^{1,6}\oplus
1_{4}^{1,0}
\hskip -1.5pt \otimes \hskip -1.5pt
1_{3}^{1,0}$,
pord$(\rho_\text{isum}(\mathfrak{t})) = 12$,

\vskip 0.7ex
\hangindent=5.5em \hangafter=1
{\white .}\hskip 1em $\rho_\text{isum}(\mathfrak{t})$ =
 $( \frac{1}{3},
\frac{1}{12},
\frac{7}{12} )
\oplus
( \frac{1}{12},
\frac{5}{12} )
\oplus
( \frac{7}{12} )
$,

\vskip 0.7ex
\hangindent=5.5em \hangafter=1
{\white .}\hskip 1em $\rho_\text{isum}(\mathfrak{s})$ =
($0$,
$\sqrt{\frac{1}{2}}$,
$\sqrt{\frac{1}{2}}$;
$-\frac{1}{2}$,
$\frac{1}{2}$;
$-\frac{1}{2}$)
 $\oplus$
($-\sqrt{\frac{1}{3}}$,
$\sqrt{\frac{2}{3}}$;
$\sqrt{\frac{1}{3}}$)
 $\oplus$
$\mathrm{i}$($1$)

Pass. 

 \ \color{black}

 \color{blue}

\noindent 106: (dims,levels) = $(3 , 
2 , 
1;12,
12,
12
)$,
irreps = $3_{4}^{1,0}
\hskip -1.5pt \otimes \hskip -1.5pt
1_{3}^{1,0}\oplus
2_{3}^{1,0}
\hskip -1.5pt \otimes \hskip -1.5pt
1_{4}^{1,0}\oplus
1_{4}^{1,6}
\hskip -1.5pt \otimes \hskip -1.5pt
1_{3}^{1,0}$,
pord$(\rho_\text{isum}(\mathfrak{t})) = 12$,

\vskip 0.7ex
\hangindent=5.5em \hangafter=1
{\white .}\hskip 1em $\rho_\text{isum}(\mathfrak{t})$ =
 $( \frac{1}{3},
\frac{1}{12},
\frac{7}{12} )
\oplus
( \frac{1}{4},
\frac{7}{12} )
\oplus
( \frac{1}{12} )
$,

\vskip 0.7ex
\hangindent=5.5em \hangafter=1
{\white .}\hskip 1em $\rho_\text{isum}(\mathfrak{s})$ =
($0$,
$\sqrt{\frac{1}{2}}$,
$\sqrt{\frac{1}{2}}$;
$-\frac{1}{2}$,
$\frac{1}{2}$;
$-\frac{1}{2}$)
 $\oplus$
($\sqrt{\frac{1}{3}}$,
$\sqrt{\frac{2}{3}}$;
$-\sqrt{\frac{1}{3}}$)
 $\oplus$
$\mathrm{i}$($-1$)

Pass. 

 \ \color{black}

\noindent 107: (dims,levels) = $(3 , 
2 , 
1;12,
12,
12
)$,
irreps = $3_{4}^{1,0}
\hskip -1.5pt \otimes \hskip -1.5pt
1_{3}^{1,0}\oplus
2_{3}^{1,0}
\hskip -1.5pt \otimes \hskip -1.5pt
1_{4}^{1,0}\oplus
1_{4}^{1,0}
\hskip -1.5pt \otimes \hskip -1.5pt
1_{3}^{1,0}$,
pord$(\rho_\text{isum}(\mathfrak{t})) = 12$,

\vskip 0.7ex
\hangindent=5.5em \hangafter=1
{\white .}\hskip 1em $\rho_\text{isum}(\mathfrak{t})$ =
 $( \frac{1}{3},
\frac{1}{12},
\frac{7}{12} )
\oplus
( \frac{1}{4},
\frac{7}{12} )
\oplus
( \frac{7}{12} )
$,

\vskip 0.7ex
\hangindent=5.5em \hangafter=1
{\white .}\hskip 1em $\rho_\text{isum}(\mathfrak{s})$ =
($0$,
$\sqrt{\frac{1}{2}}$,
$\sqrt{\frac{1}{2}}$;
$-\frac{1}{2}$,
$\frac{1}{2}$;
$-\frac{1}{2}$)
 $\oplus$
($\sqrt{\frac{1}{3}}$,
$\sqrt{\frac{2}{3}}$;
$-\sqrt{\frac{1}{3}}$)
 $\oplus$
$\mathrm{i}$($1$)

Fail:
cnd($\rho(\mathfrak s)_\mathrm{ndeg}$) = 24 does not divide
 ord($\rho(\mathfrak t)$)=12. Prop. B.4 (2)

 \ \color{black}

 \color{blue}

\noindent 108: (dims,levels) = $(3 , 
2 , 
1;12,
12,
12
)$,
irreps = $3_{4}^{1,0}
\hskip -1.5pt \otimes \hskip -1.5pt
1_{3}^{1,0}\oplus
2_{3}^{1,4}
\hskip -1.5pt \otimes \hskip -1.5pt
1_{4}^{1,0}\oplus
1_{4}^{1,6}
\hskip -1.5pt \otimes \hskip -1.5pt
1_{3}^{1,0}$,
pord$(\rho_\text{isum}(\mathfrak{t})) = 12$,

\vskip 0.7ex
\hangindent=5.5em \hangafter=1
{\white .}\hskip 1em $\rho_\text{isum}(\mathfrak{t})$ =
 $( \frac{1}{3},
\frac{1}{12},
\frac{7}{12} )
\oplus
( \frac{7}{12},
\frac{11}{12} )
\oplus
( \frac{1}{12} )
$,

\vskip 0.7ex
\hangindent=5.5em \hangafter=1
{\white .}\hskip 1em $\rho_\text{isum}(\mathfrak{s})$ =
($0$,
$\sqrt{\frac{1}{2}}$,
$\sqrt{\frac{1}{2}}$;
$-\frac{1}{2}$,
$\frac{1}{2}$;
$-\frac{1}{2}$)
 $\oplus$
($\sqrt{\frac{1}{3}}$,
$\sqrt{\frac{2}{3}}$;
$-\sqrt{\frac{1}{3}}$)
 $\oplus$
$\mathrm{i}$($-1$)

Pass. 

 \ \color{black}

\noindent 109: (dims,levels) = $(3 , 
2 , 
1;12,
12,
12
)$,
irreps = $3_{4}^{1,0}
\hskip -1.5pt \otimes \hskip -1.5pt
1_{3}^{1,0}\oplus
2_{3}^{1,4}
\hskip -1.5pt \otimes \hskip -1.5pt
1_{4}^{1,0}\oplus
1_{4}^{1,0}
\hskip -1.5pt \otimes \hskip -1.5pt
1_{3}^{1,0}$,
pord$(\rho_\text{isum}(\mathfrak{t})) = 12$,

\vskip 0.7ex
\hangindent=5.5em \hangafter=1
{\white .}\hskip 1em $\rho_\text{isum}(\mathfrak{t})$ =
 $( \frac{1}{3},
\frac{1}{12},
\frac{7}{12} )
\oplus
( \frac{7}{12},
\frac{11}{12} )
\oplus
( \frac{7}{12} )
$,

\vskip 0.7ex
\hangindent=5.5em \hangafter=1
{\white .}\hskip 1em $\rho_\text{isum}(\mathfrak{s})$ =
($0$,
$\sqrt{\frac{1}{2}}$,
$\sqrt{\frac{1}{2}}$;
$-\frac{1}{2}$,
$\frac{1}{2}$;
$-\frac{1}{2}$)
 $\oplus$
($\sqrt{\frac{1}{3}}$,
$\sqrt{\frac{2}{3}}$;
$-\sqrt{\frac{1}{3}}$)
 $\oplus$
$\mathrm{i}$($1$)

Fail:
cnd($\rho(\mathfrak s)_\mathrm{ndeg}$) = 24 does not divide
 ord($\rho(\mathfrak t)$)=12. Prop. B.4 (2)

 \ \color{black}

\noindent 110: (dims,levels) = $(3 , 
2 , 
1;12,
12,
12
)$,
irreps = $3_{4}^{1,0}
\hskip -1.5pt \otimes \hskip -1.5pt
1_{3}^{1,0}\oplus
2_{3}^{1,4}
\hskip -1.5pt \otimes \hskip -1.5pt
1_{4}^{1,0}\oplus
1_{4}^{1,0}
\hskip -1.5pt \otimes \hskip -1.5pt
1_{3}^{1,4}$,
pord$(\rho_\text{isum}(\mathfrak{t})) = 12$,

\vskip 0.7ex
\hangindent=5.5em \hangafter=1
{\white .}\hskip 1em $\rho_\text{isum}(\mathfrak{t})$ =
 $( \frac{1}{3},
\frac{1}{12},
\frac{7}{12} )
\oplus
( \frac{7}{12},
\frac{11}{12} )
\oplus
( \frac{11}{12} )
$,

\vskip 0.7ex
\hangindent=5.5em \hangafter=1
{\white .}\hskip 1em $\rho_\text{isum}(\mathfrak{s})$ =
($0$,
$\sqrt{\frac{1}{2}}$,
$\sqrt{\frac{1}{2}}$;
$-\frac{1}{2}$,
$\frac{1}{2}$;
$-\frac{1}{2}$)
 $\oplus$
($\sqrt{\frac{1}{3}}$,
$\sqrt{\frac{2}{3}}$;
$-\sqrt{\frac{1}{3}}$)
 $\oplus$
$\mathrm{i}$($1$)

Fail:
cnd($\rho(\mathfrak s)_\mathrm{ndeg}$) = 8 does not divide
 ord($\rho(\mathfrak t)$)=12. Prop. B.4 (2)

 \ \color{black}

\noindent 111: (dims,levels) = $(3 , 
2 , 
1;12,
12,
12
)$,
irreps = $3_{4}^{1,0}
\hskip -1.5pt \otimes \hskip -1.5pt
1_{3}^{1,0}\oplus
2_{3}^{1,0}
\hskip -1.5pt \otimes \hskip -1.5pt
1_{4}^{1,6}\oplus
1_{4}^{1,6}
\hskip -1.5pt \otimes \hskip -1.5pt
1_{3}^{1,0}$,
pord$(\rho_\text{isum}(\mathfrak{t})) = 12$,

\vskip 0.7ex
\hangindent=5.5em \hangafter=1
{\white .}\hskip 1em $\rho_\text{isum}(\mathfrak{t})$ =
 $( \frac{1}{3},
\frac{1}{12},
\frac{7}{12} )
\oplus
( \frac{3}{4},
\frac{1}{12} )
\oplus
( \frac{1}{12} )
$,

\vskip 0.7ex
\hangindent=5.5em \hangafter=1
{\white .}\hskip 1em $\rho_\text{isum}(\mathfrak{s})$ =
($0$,
$\sqrt{\frac{1}{2}}$,
$\sqrt{\frac{1}{2}}$;
$-\frac{1}{2}$,
$\frac{1}{2}$;
$-\frac{1}{2}$)
 $\oplus$
($-\sqrt{\frac{1}{3}}$,
$\sqrt{\frac{2}{3}}$;
$\sqrt{\frac{1}{3}}$)
 $\oplus$
$\mathrm{i}$($-1$)

Fail:
cnd($\rho(\mathfrak s)_\mathrm{ndeg}$) = 24 does not divide
 ord($\rho(\mathfrak t)$)=12. Prop. B.4 (2)

 \ \color{black}

 \color{blue}

\noindent 112: (dims,levels) = $(3 , 
2 , 
1;12,
12,
12
)$,
irreps = $3_{4}^{1,0}
\hskip -1.5pt \otimes \hskip -1.5pt
1_{3}^{1,0}\oplus
2_{3}^{1,0}
\hskip -1.5pt \otimes \hskip -1.5pt
1_{4}^{1,6}\oplus
1_{4}^{1,0}
\hskip -1.5pt \otimes \hskip -1.5pt
1_{3}^{1,0}$,
pord$(\rho_\text{isum}(\mathfrak{t})) = 12$,

\vskip 0.7ex
\hangindent=5.5em \hangafter=1
{\white .}\hskip 1em $\rho_\text{isum}(\mathfrak{t})$ =
 $( \frac{1}{3},
\frac{1}{12},
\frac{7}{12} )
\oplus
( \frac{3}{4},
\frac{1}{12} )
\oplus
( \frac{7}{12} )
$,

\vskip 0.7ex
\hangindent=5.5em \hangafter=1
{\white .}\hskip 1em $\rho_\text{isum}(\mathfrak{s})$ =
($0$,
$\sqrt{\frac{1}{2}}$,
$\sqrt{\frac{1}{2}}$;
$-\frac{1}{2}$,
$\frac{1}{2}$;
$-\frac{1}{2}$)
 $\oplus$
($-\sqrt{\frac{1}{3}}$,
$\sqrt{\frac{2}{3}}$;
$\sqrt{\frac{1}{3}}$)
 $\oplus$
$\mathrm{i}$($1$)

Pass. 

 \ \color{black}

\noindent 113: (dims,levels) = $(3 , 
2 , 
1;15,
3,
1
)$,
irreps = $3_{5}^{1}
\hskip -1.5pt \otimes \hskip -1.5pt
1_{3}^{1,0}\oplus
2_{3}^{1,0}\oplus
1_{1}^{1}$,
pord$(\rho_\text{isum}(\mathfrak{t})) = 15$,

\vskip 0.7ex
\hangindent=5.5em \hangafter=1
{\white .}\hskip 1em $\rho_\text{isum}(\mathfrak{t})$ =
 $( \frac{1}{3},
\frac{2}{15},
\frac{8}{15} )
\oplus
( 0,
\frac{1}{3} )
\oplus
( 0 )
$,

\vskip 0.7ex
\hangindent=5.5em \hangafter=1
{\white .}\hskip 1em $\rho_\text{isum}(\mathfrak{s})$ =
($\sqrt{\frac{1}{5}}$,
$-\sqrt{\frac{2}{5}}$,
$-\sqrt{\frac{2}{5}}$;
$-\frac{5+\sqrt{5}}{10}$,
$\frac{5-\sqrt{5}}{10}$;
$-\frac{5+\sqrt{5}}{10}$)
 $\oplus$
$\mathrm{i}$($-\sqrt{\frac{1}{3}}$,
$\sqrt{\frac{2}{3}}$;\ \ 
$\sqrt{\frac{1}{3}}$)
 $\oplus$
($1$)

Fail:
all $\theta$-eigenspaces that can contain unit
 have Tr$_{E_\theta}(C) \leq 0 $. Prop. B.5 (5) eqn. (B.29)

 \ \color{black}

\noindent 114: (dims,levels) = $(3 , 
2 , 
1;15,
3,
1
)$,
irreps = $3_{5}^{3}
\hskip -1.5pt \otimes \hskip -1.5pt
1_{3}^{1,0}\oplus
2_{3}^{1,0}\oplus
1_{1}^{1}$,
pord$(\rho_\text{isum}(\mathfrak{t})) = 15$,

\vskip 0.7ex
\hangindent=5.5em \hangafter=1
{\white .}\hskip 1em $\rho_\text{isum}(\mathfrak{t})$ =
 $( \frac{1}{3},
\frac{11}{15},
\frac{14}{15} )
\oplus
( 0,
\frac{1}{3} )
\oplus
( 0 )
$,

\vskip 0.7ex
\hangindent=5.5em \hangafter=1
{\white .}\hskip 1em $\rho_\text{isum}(\mathfrak{s})$ =
($-\sqrt{\frac{1}{5}}$,
$-\sqrt{\frac{2}{5}}$,
$-\sqrt{\frac{2}{5}}$;
$\frac{-5+\sqrt{5}}{10}$,
$\frac{5+\sqrt{5}}{10}$;
$\frac{-5+\sqrt{5}}{10}$)
 $\oplus$
$\mathrm{i}$($-\sqrt{\frac{1}{3}}$,
$\sqrt{\frac{2}{3}}$;\ \ 
$\sqrt{\frac{1}{3}}$)
 $\oplus$
($1$)

Fail:
all $\theta$-eigenspaces that can contain unit
 have Tr$_{E_\theta}(C) \leq 0 $. Prop. B.5 (5) eqn. (B.29)

 \ \color{black}

\noindent 115: (dims,levels) = $(3 , 
2 , 
1;15,
3,
3
)$,
irreps = $3_{5}^{1}
\hskip -1.5pt \otimes \hskip -1.5pt
1_{3}^{1,0}\oplus
2_{3}^{1,0}\oplus
1_{3}^{1,0}$,
pord$(\rho_\text{isum}(\mathfrak{t})) = 15$,

\vskip 0.7ex
\hangindent=5.5em \hangafter=1
{\white .}\hskip 1em $\rho_\text{isum}(\mathfrak{t})$ =
 $( \frac{1}{3},
\frac{2}{15},
\frac{8}{15} )
\oplus
( 0,
\frac{1}{3} )
\oplus
( \frac{1}{3} )
$,

\vskip 0.7ex
\hangindent=5.5em \hangafter=1
{\white .}\hskip 1em $\rho_\text{isum}(\mathfrak{s})$ =
($\sqrt{\frac{1}{5}}$,
$-\sqrt{\frac{2}{5}}$,
$-\sqrt{\frac{2}{5}}$;
$-\frac{5+\sqrt{5}}{10}$,
$\frac{5-\sqrt{5}}{10}$;
$-\frac{5+\sqrt{5}}{10}$)
 $\oplus$
$\mathrm{i}$($-\sqrt{\frac{1}{3}}$,
$\sqrt{\frac{2}{3}}$;\ \ 
$\sqrt{\frac{1}{3}}$)
 $\oplus$
($1$)

Fail:
Tr$_I(C) = -1 <$  0 for I = [ 0 ]. Prop. B.4 (1) eqn. (B.18)

 \ \color{black}

\noindent 116: (dims,levels) = $(3 , 
2 , 
1;15,
3,
3
)$,
irreps = $3_{5}^{1}
\hskip -1.5pt \otimes \hskip -1.5pt
1_{3}^{1,0}\oplus
2_{3}^{1,4}\oplus
1_{3}^{1,0}$,
pord$(\rho_\text{isum}(\mathfrak{t})) = 15$,

\vskip 0.7ex
\hangindent=5.5em \hangafter=1
{\white .}\hskip 1em $\rho_\text{isum}(\mathfrak{t})$ =
 $( \frac{1}{3},
\frac{2}{15},
\frac{8}{15} )
\oplus
( \frac{1}{3},
\frac{2}{3} )
\oplus
( \frac{1}{3} )
$,

\vskip 0.7ex
\hangindent=5.5em \hangafter=1
{\white .}\hskip 1em $\rho_\text{isum}(\mathfrak{s})$ =
($\sqrt{\frac{1}{5}}$,
$-\sqrt{\frac{2}{5}}$,
$-\sqrt{\frac{2}{5}}$;
$-\frac{5+\sqrt{5}}{10}$,
$\frac{5-\sqrt{5}}{10}$;
$-\frac{5+\sqrt{5}}{10}$)
 $\oplus$
$\mathrm{i}$($-\sqrt{\frac{1}{3}}$,
$\sqrt{\frac{2}{3}}$;\ \ 
$\sqrt{\frac{1}{3}}$)
 $\oplus$
($1$)

Fail:
Tr$_I(C) = -1 <$  0 for I = [ 2/3 ]. Prop. B.4 (1) eqn. (B.18)

 \ \color{black}

\noindent 117: (dims,levels) = $(3 , 
2 , 
1;15,
3,
3
)$,
irreps = $3_{5}^{1}
\hskip -1.5pt \otimes \hskip -1.5pt
1_{3}^{1,0}\oplus
2_{3}^{1,4}\oplus
1_{3}^{1,4}$,
pord$(\rho_\text{isum}(\mathfrak{t})) = 15$,

\vskip 0.7ex
\hangindent=5.5em \hangafter=1
{\white .}\hskip 1em $\rho_\text{isum}(\mathfrak{t})$ =
 $( \frac{1}{3},
\frac{2}{15},
\frac{8}{15} )
\oplus
( \frac{1}{3},
\frac{2}{3} )
\oplus
( \frac{2}{3} )
$,

\vskip 0.7ex
\hangindent=5.5em \hangafter=1
{\white .}\hskip 1em $\rho_\text{isum}(\mathfrak{s})$ =
($\sqrt{\frac{1}{5}}$,
$-\sqrt{\frac{2}{5}}$,
$-\sqrt{\frac{2}{5}}$;
$-\frac{5+\sqrt{5}}{10}$,
$\frac{5-\sqrt{5}}{10}$;
$-\frac{5+\sqrt{5}}{10}$)
 $\oplus$
$\mathrm{i}$($-\sqrt{\frac{1}{3}}$,
$\sqrt{\frac{2}{3}}$;\ \ 
$\sqrt{\frac{1}{3}}$)
 $\oplus$
($1$)

Fail:
all $\theta$-eigenspaces that can contain unit
 have Tr$_{E_\theta}(C) \leq 0 $. Prop. B.5 (5) eqn. (B.29)

 \ \color{black}

\noindent 118: (dims,levels) = $(3 , 
2 , 
1;15,
3,
3
)$,
irreps = $3_{5}^{3}
\hskip -1.5pt \otimes \hskip -1.5pt
1_{3}^{1,0}\oplus
2_{3}^{1,0}\oplus
1_{3}^{1,0}$,
pord$(\rho_\text{isum}(\mathfrak{t})) = 15$,

\vskip 0.7ex
\hangindent=5.5em \hangafter=1
{\white .}\hskip 1em $\rho_\text{isum}(\mathfrak{t})$ =
 $( \frac{1}{3},
\frac{11}{15},
\frac{14}{15} )
\oplus
( 0,
\frac{1}{3} )
\oplus
( \frac{1}{3} )
$,

\vskip 0.7ex
\hangindent=5.5em \hangafter=1
{\white .}\hskip 1em $\rho_\text{isum}(\mathfrak{s})$ =
($-\sqrt{\frac{1}{5}}$,
$-\sqrt{\frac{2}{5}}$,
$-\sqrt{\frac{2}{5}}$;
$\frac{-5+\sqrt{5}}{10}$,
$\frac{5+\sqrt{5}}{10}$;
$\frac{-5+\sqrt{5}}{10}$)
 $\oplus$
$\mathrm{i}$($-\sqrt{\frac{1}{3}}$,
$\sqrt{\frac{2}{3}}$;\ \ 
$\sqrt{\frac{1}{3}}$)
 $\oplus$
($1$)

Fail:
Tr$_I(C) = -1 <$  0 for I = [ 0 ]. Prop. B.4 (1) eqn. (B.18)

 \ \color{black}

\noindent 119: (dims,levels) = $(3 , 
2 , 
1;15,
3,
3
)$,
irreps = $3_{5}^{3}
\hskip -1.5pt \otimes \hskip -1.5pt
1_{3}^{1,0}\oplus
2_{3}^{1,4}\oplus
1_{3}^{1,0}$,
pord$(\rho_\text{isum}(\mathfrak{t})) = 15$,

\vskip 0.7ex
\hangindent=5.5em \hangafter=1
{\white .}\hskip 1em $\rho_\text{isum}(\mathfrak{t})$ =
 $( \frac{1}{3},
\frac{11}{15},
\frac{14}{15} )
\oplus
( \frac{1}{3},
\frac{2}{3} )
\oplus
( \frac{1}{3} )
$,

\vskip 0.7ex
\hangindent=5.5em \hangafter=1
{\white .}\hskip 1em $\rho_\text{isum}(\mathfrak{s})$ =
($-\sqrt{\frac{1}{5}}$,
$-\sqrt{\frac{2}{5}}$,
$-\sqrt{\frac{2}{5}}$;
$\frac{-5+\sqrt{5}}{10}$,
$\frac{5+\sqrt{5}}{10}$;
$\frac{-5+\sqrt{5}}{10}$)
 $\oplus$
$\mathrm{i}$($-\sqrt{\frac{1}{3}}$,
$\sqrt{\frac{2}{3}}$;\ \ 
$\sqrt{\frac{1}{3}}$)
 $\oplus$
($1$)

Fail:
Tr$_I(C) = -1 <$  0 for I = [ 2/3 ]. Prop. B.4 (1) eqn. (B.18)

 \ \color{black}

\noindent 120: (dims,levels) = $(3 , 
2 , 
1;15,
3,
3
)$,
irreps = $3_{5}^{3}
\hskip -1.5pt \otimes \hskip -1.5pt
1_{3}^{1,0}\oplus
2_{3}^{1,4}\oplus
1_{3}^{1,4}$,
pord$(\rho_\text{isum}(\mathfrak{t})) = 15$,

\vskip 0.7ex
\hangindent=5.5em \hangafter=1
{\white .}\hskip 1em $\rho_\text{isum}(\mathfrak{t})$ =
 $( \frac{1}{3},
\frac{11}{15},
\frac{14}{15} )
\oplus
( \frac{1}{3},
\frac{2}{3} )
\oplus
( \frac{2}{3} )
$,

\vskip 0.7ex
\hangindent=5.5em \hangafter=1
{\white .}\hskip 1em $\rho_\text{isum}(\mathfrak{s})$ =
($-\sqrt{\frac{1}{5}}$,
$-\sqrt{\frac{2}{5}}$,
$-\sqrt{\frac{2}{5}}$;
$\frac{-5+\sqrt{5}}{10}$,
$\frac{5+\sqrt{5}}{10}$;
$\frac{-5+\sqrt{5}}{10}$)
 $\oplus$
$\mathrm{i}$($-\sqrt{\frac{1}{3}}$,
$\sqrt{\frac{2}{3}}$;\ \ 
$\sqrt{\frac{1}{3}}$)
 $\oplus$
($1$)

Fail:
all $\theta$-eigenspaces that can contain unit
 have Tr$_{E_\theta}(C) \leq 0 $. Prop. B.5 (5) eqn. (B.29)

 \ \color{black}

 \color{blue}

\noindent 121: (dims,levels) = $(3 , 
2 , 
1;15,
6,
3
)$,
irreps = $3_{5}^{1}
\hskip -1.5pt \otimes \hskip -1.5pt
1_{3}^{1,0}\oplus
2_{2}^{1,0}
\hskip -1.5pt \otimes \hskip -1.5pt
1_{3}^{1,0}\oplus
1_{3}^{1,0}$,
pord$(\rho_\text{isum}(\mathfrak{t})) = 10$,

\vskip 0.7ex
\hangindent=5.5em \hangafter=1
{\white .}\hskip 1em $\rho_\text{isum}(\mathfrak{t})$ =
 $( \frac{1}{3},
\frac{2}{15},
\frac{8}{15} )
\oplus
( \frac{1}{3},
\frac{5}{6} )
\oplus
( \frac{1}{3} )
$,

\vskip 0.7ex
\hangindent=5.5em \hangafter=1
{\white .}\hskip 1em $\rho_\text{isum}(\mathfrak{s})$ =
($\sqrt{\frac{1}{5}}$,
$-\sqrt{\frac{2}{5}}$,
$-\sqrt{\frac{2}{5}}$;
$-\frac{5+\sqrt{5}}{10}$,
$\frac{5-\sqrt{5}}{10}$;
$-\frac{5+\sqrt{5}}{10}$)
 $\oplus$
($-\frac{1}{2}$,
$-\sqrt{\frac{3}{4}}$;
$\frac{1}{2}$)
 $\oplus$
($1$)

Pass. 

 \ \color{black}

 \color{blue}

\noindent 122: (dims,levels) = $(3 , 
2 , 
1;15,
6,
3
)$,
irreps = $3_{5}^{3}
\hskip -1.5pt \otimes \hskip -1.5pt
1_{3}^{1,0}\oplus
2_{2}^{1,0}
\hskip -1.5pt \otimes \hskip -1.5pt
1_{3}^{1,0}\oplus
1_{3}^{1,0}$,
pord$(\rho_\text{isum}(\mathfrak{t})) = 10$,

\vskip 0.7ex
\hangindent=5.5em \hangafter=1
{\white .}\hskip 1em $\rho_\text{isum}(\mathfrak{t})$ =
 $( \frac{1}{3},
\frac{11}{15},
\frac{14}{15} )
\oplus
( \frac{1}{3},
\frac{5}{6} )
\oplus
( \frac{1}{3} )
$,

\vskip 0.7ex
\hangindent=5.5em \hangafter=1
{\white .}\hskip 1em $\rho_\text{isum}(\mathfrak{s})$ =
($-\sqrt{\frac{1}{5}}$,
$-\sqrt{\frac{2}{5}}$,
$-\sqrt{\frac{2}{5}}$;
$\frac{-5+\sqrt{5}}{10}$,
$\frac{5+\sqrt{5}}{10}$;
$\frac{-5+\sqrt{5}}{10}$)
 $\oplus$
($-\frac{1}{2}$,
$-\sqrt{\frac{3}{4}}$;
$\frac{1}{2}$)
 $\oplus$
($1$)

Pass. 

 \ \color{black}

 \color{blue}

\noindent 123: (dims,levels) = $(3 , 
2 , 
1;15,
6,
6
)$,
irreps = $3_{5}^{1}
\hskip -1.5pt \otimes \hskip -1.5pt
1_{3}^{1,0}\oplus
2_{2}^{1,0}
\hskip -1.5pt \otimes \hskip -1.5pt
1_{3}^{1,0}\oplus
1_{3}^{1,0}
\hskip -1.5pt \otimes \hskip -1.5pt
1_{2}^{1,0}$,
pord$(\rho_\text{isum}(\mathfrak{t})) = 10$,

\vskip 0.7ex
\hangindent=5.5em \hangafter=1
{\white .}\hskip 1em $\rho_\text{isum}(\mathfrak{t})$ =
 $( \frac{1}{3},
\frac{2}{15},
\frac{8}{15} )
\oplus
( \frac{1}{3},
\frac{5}{6} )
\oplus
( \frac{5}{6} )
$,

\vskip 0.7ex
\hangindent=5.5em \hangafter=1
{\white .}\hskip 1em $\rho_\text{isum}(\mathfrak{s})$ =
($\sqrt{\frac{1}{5}}$,
$-\sqrt{\frac{2}{5}}$,
$-\sqrt{\frac{2}{5}}$;
$-\frac{5+\sqrt{5}}{10}$,
$\frac{5-\sqrt{5}}{10}$;
$-\frac{5+\sqrt{5}}{10}$)
 $\oplus$
($-\frac{1}{2}$,
$-\sqrt{\frac{3}{4}}$;
$\frac{1}{2}$)
 $\oplus$
($-1$)

Pass. 

 \ \color{black}

 \color{blue}

\noindent 124: (dims,levels) = $(3 , 
2 , 
1;15,
6,
6
)$,
irreps = $3_{5}^{3}
\hskip -1.5pt \otimes \hskip -1.5pt
1_{3}^{1,0}\oplus
2_{2}^{1,0}
\hskip -1.5pt \otimes \hskip -1.5pt
1_{3}^{1,0}\oplus
1_{3}^{1,0}
\hskip -1.5pt \otimes \hskip -1.5pt
1_{2}^{1,0}$,
pord$(\rho_\text{isum}(\mathfrak{t})) = 10$,

\vskip 0.7ex
\hangindent=5.5em \hangafter=1
{\white .}\hskip 1em $\rho_\text{isum}(\mathfrak{t})$ =
 $( \frac{1}{3},
\frac{11}{15},
\frac{14}{15} )
\oplus
( \frac{1}{3},
\frac{5}{6} )
\oplus
( \frac{5}{6} )
$,

\vskip 0.7ex
\hangindent=5.5em \hangafter=1
{\white .}\hskip 1em $\rho_\text{isum}(\mathfrak{s})$ =
($-\sqrt{\frac{1}{5}}$,
$-\sqrt{\frac{2}{5}}$,
$-\sqrt{\frac{2}{5}}$;
$\frac{-5+\sqrt{5}}{10}$,
$\frac{5+\sqrt{5}}{10}$;
$\frac{-5+\sqrt{5}}{10}$)
 $\oplus$
($-\frac{1}{2}$,
$-\sqrt{\frac{3}{4}}$;
$\frac{1}{2}$)
 $\oplus$
($-1$)

Pass. 

 \ \color{black}

 \color{blue}

\noindent 125: (dims,levels) = $(3 , 
2 , 
1;15,
15,
3
)$,
irreps = $3_{5}^{1}
\hskip -1.5pt \otimes \hskip -1.5pt
1_{3}^{1,0}\oplus
2_{5}^{1}
\hskip -1.5pt \otimes \hskip -1.5pt
1_{3}^{1,0}\oplus
1_{3}^{1,0}$,
pord$(\rho_\text{isum}(\mathfrak{t})) = 5$,

\vskip 0.7ex
\hangindent=5.5em \hangafter=1
{\white .}\hskip 1em $\rho_\text{isum}(\mathfrak{t})$ =
 $( \frac{1}{3},
\frac{2}{15},
\frac{8}{15} )
\oplus
( \frac{2}{15},
\frac{8}{15} )
\oplus
( \frac{1}{3} )
$,

\vskip 0.7ex
\hangindent=5.5em \hangafter=1
{\white .}\hskip 1em $\rho_\text{isum}(\mathfrak{s})$ =
($\sqrt{\frac{1}{5}}$,
$-\sqrt{\frac{2}{5}}$,
$-\sqrt{\frac{2}{5}}$;
$-\frac{5+\sqrt{5}}{10}$,
$\frac{5-\sqrt{5}}{10}$;
$-\frac{5+\sqrt{5}}{10}$)
 $\oplus$
$\mathrm{i}$($\frac{1}{\sqrt{5}}c^{3}_{20}
$,
$\frac{1}{\sqrt{5}}c^{1}_{20}
$;\ \ 
$-\frac{1}{\sqrt{5}}c^{3}_{20}
$)
 $\oplus$
($1$)

Pass. 

 \ \color{black}

 \color{blue}

\noindent 126: (dims,levels) = $(3 , 
2 , 
1;15,
15,
3
)$,
irreps = $3_{5}^{3}
\hskip -1.5pt \otimes \hskip -1.5pt
1_{3}^{1,0}\oplus
2_{5}^{2}
\hskip -1.5pt \otimes \hskip -1.5pt
1_{3}^{1,0}\oplus
1_{3}^{1,0}$,
pord$(\rho_\text{isum}(\mathfrak{t})) = 5$,

\vskip 0.7ex
\hangindent=5.5em \hangafter=1
{\white .}\hskip 1em $\rho_\text{isum}(\mathfrak{t})$ =
 $( \frac{1}{3},
\frac{11}{15},
\frac{14}{15} )
\oplus
( \frac{11}{15},
\frac{14}{15} )
\oplus
( \frac{1}{3} )
$,

\vskip 0.7ex
\hangindent=5.5em \hangafter=1
{\white .}\hskip 1em $\rho_\text{isum}(\mathfrak{s})$ =
($-\sqrt{\frac{1}{5}}$,
$-\sqrt{\frac{2}{5}}$,
$-\sqrt{\frac{2}{5}}$;
$\frac{-5+\sqrt{5}}{10}$,
$\frac{5+\sqrt{5}}{10}$;
$\frac{-5+\sqrt{5}}{10}$)
 $\oplus$
$\mathrm{i}$($-\frac{1}{\sqrt{5}}c^{1}_{20}
$,
$\frac{1}{\sqrt{5}}c^{3}_{20}
$;\ \ 
$\frac{1}{\sqrt{5}}c^{1}_{20}
$)
 $\oplus$
($1$)

Pass. 

 \ \color{black}

 \color{blue}

\noindent 127: (dims,levels) = $(3 , 
2 , 
1;20,
4,
4
)$,
irreps = $3_{5}^{1}
\hskip -1.5pt \otimes \hskip -1.5pt
1_{4}^{1,0}\oplus
2_{4}^{1,0}\oplus
1_{4}^{1,0}$,
pord$(\rho_\text{isum}(\mathfrak{t})) = 10$,

\vskip 0.7ex
\hangindent=5.5em \hangafter=1
{\white .}\hskip 1em $\rho_\text{isum}(\mathfrak{t})$ =
 $( \frac{1}{4},
\frac{1}{20},
\frac{9}{20} )
\oplus
( \frac{1}{4},
\frac{3}{4} )
\oplus
( \frac{1}{4} )
$,

\vskip 0.7ex
\hangindent=5.5em \hangafter=1
{\white .}\hskip 1em $\rho_\text{isum}(\mathfrak{s})$ =
$\mathrm{i}$($\sqrt{\frac{1}{5}}$,
$\sqrt{\frac{2}{5}}$,
$\sqrt{\frac{2}{5}}$;\ \ 
$-\frac{5+\sqrt{5}}{10}$,
$\frac{5-\sqrt{5}}{10}$;\ \ 
$-\frac{5+\sqrt{5}}{10}$)
 $\oplus$
$\mathrm{i}$($-\frac{1}{2}$,
$\sqrt{\frac{3}{4}}$;\ \ 
$\frac{1}{2}$)
 $\oplus$
$\mathrm{i}$($1$)

Pass. 

 \ \color{black}

 \color{blue}

\noindent 128: (dims,levels) = $(3 , 
2 , 
1;20,
4,
4
)$,
irreps = $3_{5}^{1}
\hskip -1.5pt \otimes \hskip -1.5pt
1_{4}^{1,0}\oplus
2_{4}^{1,0}\oplus
1_{4}^{1,6}$,
pord$(\rho_\text{isum}(\mathfrak{t})) = 10$,

\vskip 0.7ex
\hangindent=5.5em \hangafter=1
{\white .}\hskip 1em $\rho_\text{isum}(\mathfrak{t})$ =
 $( \frac{1}{4},
\frac{1}{20},
\frac{9}{20} )
\oplus
( \frac{1}{4},
\frac{3}{4} )
\oplus
( \frac{3}{4} )
$,

\vskip 0.7ex
\hangindent=5.5em \hangafter=1
{\white .}\hskip 1em $\rho_\text{isum}(\mathfrak{s})$ =
$\mathrm{i}$($\sqrt{\frac{1}{5}}$,
$\sqrt{\frac{2}{5}}$,
$\sqrt{\frac{2}{5}}$;\ \ 
$-\frac{5+\sqrt{5}}{10}$,
$\frac{5-\sqrt{5}}{10}$;\ \ 
$-\frac{5+\sqrt{5}}{10}$)
 $\oplus$
$\mathrm{i}$($-\frac{1}{2}$,
$\sqrt{\frac{3}{4}}$;\ \ 
$\frac{1}{2}$)
 $\oplus$
$\mathrm{i}$($-1$)

Pass. 

 \ \color{black}

 \color{blue}

\noindent 129: (dims,levels) = $(3 , 
2 , 
1;20,
4,
4
)$,
irreps = $3_{5}^{3}
\hskip -1.5pt \otimes \hskip -1.5pt
1_{4}^{1,0}\oplus
2_{4}^{1,0}\oplus
1_{4}^{1,0}$,
pord$(\rho_\text{isum}(\mathfrak{t})) = 10$,

\vskip 0.7ex
\hangindent=5.5em \hangafter=1
{\white .}\hskip 1em $\rho_\text{isum}(\mathfrak{t})$ =
 $( \frac{1}{4},
\frac{13}{20},
\frac{17}{20} )
\oplus
( \frac{1}{4},
\frac{3}{4} )
\oplus
( \frac{1}{4} )
$,

\vskip 0.7ex
\hangindent=5.5em \hangafter=1
{\white .}\hskip 1em $\rho_\text{isum}(\mathfrak{s})$ =
$\mathrm{i}$($-\sqrt{\frac{1}{5}}$,
$\sqrt{\frac{2}{5}}$,
$\sqrt{\frac{2}{5}}$;\ \ 
$\frac{-5+\sqrt{5}}{10}$,
$\frac{5+\sqrt{5}}{10}$;\ \ 
$\frac{-5+\sqrt{5}}{10}$)
 $\oplus$
$\mathrm{i}$($-\frac{1}{2}$,
$\sqrt{\frac{3}{4}}$;\ \ 
$\frac{1}{2}$)
 $\oplus$
$\mathrm{i}$($1$)

Pass. 

 \ \color{black}

 \color{blue}

\noindent 130: (dims,levels) = $(3 , 
2 , 
1;20,
4,
4
)$,
irreps = $3_{5}^{3}
\hskip -1.5pt \otimes \hskip -1.5pt
1_{4}^{1,0}\oplus
2_{4}^{1,0}\oplus
1_{4}^{1,6}$,
pord$(\rho_\text{isum}(\mathfrak{t})) = 10$,

\vskip 0.7ex
\hangindent=5.5em \hangafter=1
{\white .}\hskip 1em $\rho_\text{isum}(\mathfrak{t})$ =
 $( \frac{1}{4},
\frac{13}{20},
\frac{17}{20} )
\oplus
( \frac{1}{4},
\frac{3}{4} )
\oplus
( \frac{3}{4} )
$,

\vskip 0.7ex
\hangindent=5.5em \hangafter=1
{\white .}\hskip 1em $\rho_\text{isum}(\mathfrak{s})$ =
$\mathrm{i}$($-\sqrt{\frac{1}{5}}$,
$\sqrt{\frac{2}{5}}$,
$\sqrt{\frac{2}{5}}$;\ \ 
$\frac{-5+\sqrt{5}}{10}$,
$\frac{5+\sqrt{5}}{10}$;\ \ 
$\frac{-5+\sqrt{5}}{10}$)
 $\oplus$
$\mathrm{i}$($-\frac{1}{2}$,
$\sqrt{\frac{3}{4}}$;\ \ 
$\frac{1}{2}$)
 $\oplus$
$\mathrm{i}$($-1$)

Pass. 

 \ \color{black}

\noindent 131: (dims,levels) = $(3 , 
2 , 
1;20,
12,
4
)$,
irreps = $3_{5}^{1}
\hskip -1.5pt \otimes \hskip -1.5pt
1_{4}^{1,0}\oplus
2_{3}^{1,0}
\hskip -1.5pt \otimes \hskip -1.5pt
1_{4}^{1,0}\oplus
1_{4}^{1,0}$,
pord$(\rho_\text{isum}(\mathfrak{t})) = 15$,

\vskip 0.7ex
\hangindent=5.5em \hangafter=1
{\white .}\hskip 1em $\rho_\text{isum}(\mathfrak{t})$ =
 $( \frac{1}{4},
\frac{1}{20},
\frac{9}{20} )
\oplus
( \frac{1}{4},
\frac{7}{12} )
\oplus
( \frac{1}{4} )
$,

\vskip 0.7ex
\hangindent=5.5em \hangafter=1
{\white .}\hskip 1em $\rho_\text{isum}(\mathfrak{s})$ =
$\mathrm{i}$($\sqrt{\frac{1}{5}}$,
$\sqrt{\frac{2}{5}}$,
$\sqrt{\frac{2}{5}}$;\ \ 
$-\frac{5+\sqrt{5}}{10}$,
$\frac{5-\sqrt{5}}{10}$;\ \ 
$-\frac{5+\sqrt{5}}{10}$)
 $\oplus$
($\sqrt{\frac{1}{3}}$,
$\sqrt{\frac{2}{3}}$;
$-\sqrt{\frac{1}{3}}$)
 $\oplus$
$\mathrm{i}$($1$)

Fail:
Tr$_I(C) = -1 <$  0 for I = [ 7/12 ]. Prop. B.4 (1) eqn. (B.18)

 \ \color{black}

\noindent 132: (dims,levels) = $(3 , 
2 , 
1;20,
12,
4
)$,
irreps = $3_{5}^{1}
\hskip -1.5pt \otimes \hskip -1.5pt
1_{4}^{1,0}\oplus
2_{3}^{1,8}
\hskip -1.5pt \otimes \hskip -1.5pt
1_{4}^{1,0}\oplus
1_{4}^{1,0}$,
pord$(\rho_\text{isum}(\mathfrak{t})) = 15$,

\vskip 0.7ex
\hangindent=5.5em \hangafter=1
{\white .}\hskip 1em $\rho_\text{isum}(\mathfrak{t})$ =
 $( \frac{1}{4},
\frac{1}{20},
\frac{9}{20} )
\oplus
( \frac{1}{4},
\frac{11}{12} )
\oplus
( \frac{1}{4} )
$,

\vskip 0.7ex
\hangindent=5.5em \hangafter=1
{\white .}\hskip 1em $\rho_\text{isum}(\mathfrak{s})$ =
$\mathrm{i}$($\sqrt{\frac{1}{5}}$,
$\sqrt{\frac{2}{5}}$,
$\sqrt{\frac{2}{5}}$;\ \ 
$-\frac{5+\sqrt{5}}{10}$,
$\frac{5-\sqrt{5}}{10}$;\ \ 
$-\frac{5+\sqrt{5}}{10}$)
 $\oplus$
($-\sqrt{\frac{1}{3}}$,
$\sqrt{\frac{2}{3}}$;
$\sqrt{\frac{1}{3}}$)
 $\oplus$
$\mathrm{i}$($1$)

Fail:
Tr$_I(C) = -1 <$  0 for I = [ 11/12 ]. Prop. B.4 (1) eqn. (B.18)

 \ \color{black}

\noindent 133: (dims,levels) = $(3 , 
2 , 
1;20,
12,
4
)$,
irreps = $3_{5}^{3}
\hskip -1.5pt \otimes \hskip -1.5pt
1_{4}^{1,0}\oplus
2_{3}^{1,0}
\hskip -1.5pt \otimes \hskip -1.5pt
1_{4}^{1,0}\oplus
1_{4}^{1,0}$,
pord$(\rho_\text{isum}(\mathfrak{t})) = 15$,

\vskip 0.7ex
\hangindent=5.5em \hangafter=1
{\white .}\hskip 1em $\rho_\text{isum}(\mathfrak{t})$ =
 $( \frac{1}{4},
\frac{13}{20},
\frac{17}{20} )
\oplus
( \frac{1}{4},
\frac{7}{12} )
\oplus
( \frac{1}{4} )
$,

\vskip 0.7ex
\hangindent=5.5em \hangafter=1
{\white .}\hskip 1em $\rho_\text{isum}(\mathfrak{s})$ =
$\mathrm{i}$($-\sqrt{\frac{1}{5}}$,
$\sqrt{\frac{2}{5}}$,
$\sqrt{\frac{2}{5}}$;\ \ 
$\frac{-5+\sqrt{5}}{10}$,
$\frac{5+\sqrt{5}}{10}$;\ \ 
$\frac{-5+\sqrt{5}}{10}$)
 $\oplus$
($\sqrt{\frac{1}{3}}$,
$\sqrt{\frac{2}{3}}$;
$-\sqrt{\frac{1}{3}}$)
 $\oplus$
$\mathrm{i}$($1$)

Fail:
Tr$_I(C) = -1 <$  0 for I = [ 7/12 ]. Prop. B.4 (1) eqn. (B.18)

 \ \color{black}

\noindent 134: (dims,levels) = $(3 , 
2 , 
1;20,
12,
4
)$,
irreps = $3_{5}^{3}
\hskip -1.5pt \otimes \hskip -1.5pt
1_{4}^{1,0}\oplus
2_{3}^{1,8}
\hskip -1.5pt \otimes \hskip -1.5pt
1_{4}^{1,0}\oplus
1_{4}^{1,0}$,
pord$(\rho_\text{isum}(\mathfrak{t})) = 15$,

\vskip 0.7ex
\hangindent=5.5em \hangafter=1
{\white .}\hskip 1em $\rho_\text{isum}(\mathfrak{t})$ =
 $( \frac{1}{4},
\frac{13}{20},
\frac{17}{20} )
\oplus
( \frac{1}{4},
\frac{11}{12} )
\oplus
( \frac{1}{4} )
$,

\vskip 0.7ex
\hangindent=5.5em \hangafter=1
{\white .}\hskip 1em $\rho_\text{isum}(\mathfrak{s})$ =
$\mathrm{i}$($-\sqrt{\frac{1}{5}}$,
$\sqrt{\frac{2}{5}}$,
$\sqrt{\frac{2}{5}}$;\ \ 
$\frac{-5+\sqrt{5}}{10}$,
$\frac{5+\sqrt{5}}{10}$;\ \ 
$\frac{-5+\sqrt{5}}{10}$)
 $\oplus$
($-\sqrt{\frac{1}{3}}$,
$\sqrt{\frac{2}{3}}$;
$\sqrt{\frac{1}{3}}$)
 $\oplus$
$\mathrm{i}$($1$)

Fail:
Tr$_I(C) = -1 <$  0 for I = [ 11/12 ]. Prop. B.4 (1) eqn. (B.18)

 \ \color{black}

\noindent 135: (dims,levels) = $(3 , 
2 , 
1;20,
12,
12
)$,
irreps = $3_{5}^{1}
\hskip -1.5pt \otimes \hskip -1.5pt
1_{4}^{1,0}\oplus
2_{3}^{1,0}
\hskip -1.5pt \otimes \hskip -1.5pt
1_{4}^{1,0}\oplus
1_{4}^{1,0}
\hskip -1.5pt \otimes \hskip -1.5pt
1_{3}^{1,0}$,
pord$(\rho_\text{isum}(\mathfrak{t})) = 15$,

\vskip 0.7ex
\hangindent=5.5em \hangafter=1
{\white .}\hskip 1em $\rho_\text{isum}(\mathfrak{t})$ =
 $( \frac{1}{4},
\frac{1}{20},
\frac{9}{20} )
\oplus
( \frac{1}{4},
\frac{7}{12} )
\oplus
( \frac{7}{12} )
$,

\vskip 0.7ex
\hangindent=5.5em \hangafter=1
{\white .}\hskip 1em $\rho_\text{isum}(\mathfrak{s})$ =
$\mathrm{i}$($\sqrt{\frac{1}{5}}$,
$\sqrt{\frac{2}{5}}$,
$\sqrt{\frac{2}{5}}$;\ \ 
$-\frac{5+\sqrt{5}}{10}$,
$\frac{5-\sqrt{5}}{10}$;\ \ 
$-\frac{5+\sqrt{5}}{10}$)
 $\oplus$
($\sqrt{\frac{1}{3}}$,
$\sqrt{\frac{2}{3}}$;
$-\sqrt{\frac{1}{3}}$)
 $\oplus$
$\mathrm{i}$($1$)

Fail:
all $\theta$-eigenspaces that can contain unit
 have Tr$_{E_\theta}(C) \leq 0 $. Prop. B.5 (5) eqn. (B.29)

 \ \color{black}

\noindent 136: (dims,levels) = $(3 , 
2 , 
1;20,
12,
12
)$,
irreps = $3_{5}^{1}
\hskip -1.5pt \otimes \hskip -1.5pt
1_{4}^{1,0}\oplus
2_{3}^{1,8}
\hskip -1.5pt \otimes \hskip -1.5pt
1_{4}^{1,0}\oplus
1_{4}^{1,0}
\hskip -1.5pt \otimes \hskip -1.5pt
1_{3}^{1,4}$,
pord$(\rho_\text{isum}(\mathfrak{t})) = 15$,

\vskip 0.7ex
\hangindent=5.5em \hangafter=1
{\white .}\hskip 1em $\rho_\text{isum}(\mathfrak{t})$ =
 $( \frac{1}{4},
\frac{1}{20},
\frac{9}{20} )
\oplus
( \frac{1}{4},
\frac{11}{12} )
\oplus
( \frac{11}{12} )
$,

\vskip 0.7ex
\hangindent=5.5em \hangafter=1
{\white .}\hskip 1em $\rho_\text{isum}(\mathfrak{s})$ =
$\mathrm{i}$($\sqrt{\frac{1}{5}}$,
$\sqrt{\frac{2}{5}}$,
$\sqrt{\frac{2}{5}}$;\ \ 
$-\frac{5+\sqrt{5}}{10}$,
$\frac{5-\sqrt{5}}{10}$;\ \ 
$-\frac{5+\sqrt{5}}{10}$)
 $\oplus$
($-\sqrt{\frac{1}{3}}$,
$\sqrt{\frac{2}{3}}$;
$\sqrt{\frac{1}{3}}$)
 $\oplus$
$\mathrm{i}$($1$)

Fail:
all $\theta$-eigenspaces that can contain unit
 have Tr$_{E_\theta}(C) \leq 0 $. Prop. B.5 (5) eqn. (B.29)

 \ \color{black}

\noindent 137: (dims,levels) = $(3 , 
2 , 
1;20,
12,
12
)$,
irreps = $3_{5}^{3}
\hskip -1.5pt \otimes \hskip -1.5pt
1_{4}^{1,0}\oplus
2_{3}^{1,0}
\hskip -1.5pt \otimes \hskip -1.5pt
1_{4}^{1,0}\oplus
1_{4}^{1,0}
\hskip -1.5pt \otimes \hskip -1.5pt
1_{3}^{1,0}$,
pord$(\rho_\text{isum}(\mathfrak{t})) = 15$,

\vskip 0.7ex
\hangindent=5.5em \hangafter=1
{\white .}\hskip 1em $\rho_\text{isum}(\mathfrak{t})$ =
 $( \frac{1}{4},
\frac{13}{20},
\frac{17}{20} )
\oplus
( \frac{1}{4},
\frac{7}{12} )
\oplus
( \frac{7}{12} )
$,

\vskip 0.7ex
\hangindent=5.5em \hangafter=1
{\white .}\hskip 1em $\rho_\text{isum}(\mathfrak{s})$ =
$\mathrm{i}$($-\sqrt{\frac{1}{5}}$,
$\sqrt{\frac{2}{5}}$,
$\sqrt{\frac{2}{5}}$;\ \ 
$\frac{-5+\sqrt{5}}{10}$,
$\frac{5+\sqrt{5}}{10}$;\ \ 
$\frac{-5+\sqrt{5}}{10}$)
 $\oplus$
($\sqrt{\frac{1}{3}}$,
$\sqrt{\frac{2}{3}}$;
$-\sqrt{\frac{1}{3}}$)
 $\oplus$
$\mathrm{i}$($1$)

Fail:
all $\theta$-eigenspaces that can contain unit
 have Tr$_{E_\theta}(C) \leq 0 $. Prop. B.5 (5) eqn. (B.29)

 \ \color{black}

\noindent 138: (dims,levels) = $(3 , 
2 , 
1;20,
12,
12
)$,
irreps = $3_{5}^{3}
\hskip -1.5pt \otimes \hskip -1.5pt
1_{4}^{1,0}\oplus
2_{3}^{1,8}
\hskip -1.5pt \otimes \hskip -1.5pt
1_{4}^{1,0}\oplus
1_{4}^{1,0}
\hskip -1.5pt \otimes \hskip -1.5pt
1_{3}^{1,4}$,
pord$(\rho_\text{isum}(\mathfrak{t})) = 15$,

\vskip 0.7ex
\hangindent=5.5em \hangafter=1
{\white .}\hskip 1em $\rho_\text{isum}(\mathfrak{t})$ =
 $( \frac{1}{4},
\frac{13}{20},
\frac{17}{20} )
\oplus
( \frac{1}{4},
\frac{11}{12} )
\oplus
( \frac{11}{12} )
$,

\vskip 0.7ex
\hangindent=5.5em \hangafter=1
{\white .}\hskip 1em $\rho_\text{isum}(\mathfrak{s})$ =
$\mathrm{i}$($-\sqrt{\frac{1}{5}}$,
$\sqrt{\frac{2}{5}}$,
$\sqrt{\frac{2}{5}}$;\ \ 
$\frac{-5+\sqrt{5}}{10}$,
$\frac{5+\sqrt{5}}{10}$;\ \ 
$\frac{-5+\sqrt{5}}{10}$)
 $\oplus$
($-\sqrt{\frac{1}{3}}$,
$\sqrt{\frac{2}{3}}$;
$\sqrt{\frac{1}{3}}$)
 $\oplus$
$\mathrm{i}$($1$)

Fail:
all $\theta$-eigenspaces that can contain unit
 have Tr$_{E_\theta}(C) \leq 0 $. Prop. B.5 (5) eqn. (B.29)

 \ \color{black}

 \color{blue}

\noindent 139: (dims,levels) = $(3 , 
2 , 
1;20,
20,
4
)$,
irreps = $3_{5}^{1}
\hskip -1.5pt \otimes \hskip -1.5pt
1_{4}^{1,0}\oplus
2_{5}^{1}
\hskip -1.5pt \otimes \hskip -1.5pt
1_{4}^{1,0}\oplus
1_{4}^{1,0}$,
pord$(\rho_\text{isum}(\mathfrak{t})) = 5$,

\vskip 0.7ex
\hangindent=5.5em \hangafter=1
{\white .}\hskip 1em $\rho_\text{isum}(\mathfrak{t})$ =
 $( \frac{1}{4},
\frac{1}{20},
\frac{9}{20} )
\oplus
( \frac{1}{20},
\frac{9}{20} )
\oplus
( \frac{1}{4} )
$,

\vskip 0.7ex
\hangindent=5.5em \hangafter=1
{\white .}\hskip 1em $\rho_\text{isum}(\mathfrak{s})$ =
$\mathrm{i}$($\sqrt{\frac{1}{5}}$,
$\sqrt{\frac{2}{5}}$,
$\sqrt{\frac{2}{5}}$;\ \ 
$-\frac{5+\sqrt{5}}{10}$,
$\frac{5-\sqrt{5}}{10}$;\ \ 
$-\frac{5+\sqrt{5}}{10}$)
 $\oplus$
($-\frac{1}{\sqrt{5}}c^{3}_{20}
$,
$\frac{1}{\sqrt{5}}c^{1}_{20}
$;
$\frac{1}{\sqrt{5}}c^{3}_{20}
$)
 $\oplus$
$\mathrm{i}$($1$)

Pass. 

 \ \color{black}

 \color{blue}

\noindent 140: (dims,levels) = $(3 , 
2 , 
1;20,
20,
4
)$,
irreps = $3_{5}^{3}
\hskip -1.5pt \otimes \hskip -1.5pt
1_{4}^{1,0}\oplus
2_{5}^{2}
\hskip -1.5pt \otimes \hskip -1.5pt
1_{4}^{1,0}\oplus
1_{4}^{1,0}$,
pord$(\rho_\text{isum}(\mathfrak{t})) = 5$,

\vskip 0.7ex
\hangindent=5.5em \hangafter=1
{\white .}\hskip 1em $\rho_\text{isum}(\mathfrak{t})$ =
 $( \frac{1}{4},
\frac{13}{20},
\frac{17}{20} )
\oplus
( \frac{13}{20},
\frac{17}{20} )
\oplus
( \frac{1}{4} )
$,

\vskip 0.7ex
\hangindent=5.5em \hangafter=1
{\white .}\hskip 1em $\rho_\text{isum}(\mathfrak{s})$ =
$\mathrm{i}$($-\sqrt{\frac{1}{5}}$,
$\sqrt{\frac{2}{5}}$,
$\sqrt{\frac{2}{5}}$;\ \ 
$\frac{-5+\sqrt{5}}{10}$,
$\frac{5+\sqrt{5}}{10}$;\ \ 
$\frac{-5+\sqrt{5}}{10}$)
 $\oplus$
($\frac{1}{\sqrt{5}}c^{1}_{20}
$,
$\frac{1}{\sqrt{5}}c^{3}_{20}
$;
$-\frac{1}{\sqrt{5}}c^{1}_{20}
$)
 $\oplus$
$\mathrm{i}$($1$)

Pass. 

 \ \color{black}

 \color{blue}

\noindent 141: (dims,levels) = $(3 , 
2 , 
1;24,
3,
1
)$,
irreps = $3_{8}^{3,0}
\hskip -1.5pt \otimes \hskip -1.5pt
1_{3}^{1,0}\oplus
2_{3}^{1,0}\oplus
1_{1}^{1}$,
pord$(\rho_\text{isum}(\mathfrak{t})) = 24$,

\vskip 0.7ex
\hangindent=5.5em \hangafter=1
{\white .}\hskip 1em $\rho_\text{isum}(\mathfrak{t})$ =
 $( \frac{1}{3},
\frac{5}{24},
\frac{17}{24} )
\oplus
( 0,
\frac{1}{3} )
\oplus
( 0 )
$,

\vskip 0.7ex
\hangindent=5.5em \hangafter=1
{\white .}\hskip 1em $\rho_\text{isum}(\mathfrak{s})$ =
$\mathrm{i}$($0$,
$\sqrt{\frac{1}{2}}$,
$\sqrt{\frac{1}{2}}$;\ \ 
$\frac{1}{2}$,
$-\frac{1}{2}$;\ \ 
$\frac{1}{2}$)
 $\oplus$
$\mathrm{i}$($-\sqrt{\frac{1}{3}}$,
$\sqrt{\frac{2}{3}}$;\ \ 
$\sqrt{\frac{1}{3}}$)
 $\oplus$
($1$)

Pass. 

 \ \color{black}

 \color{blue}

\noindent 142: (dims,levels) = $(3 , 
2 , 
1;24,
3,
1
)$,
irreps = $3_{8}^{1,0}
\hskip -1.5pt \otimes \hskip -1.5pt
1_{3}^{1,0}\oplus
2_{3}^{1,0}\oplus
1_{1}^{1}$,
pord$(\rho_\text{isum}(\mathfrak{t})) = 24$,

\vskip 0.7ex
\hangindent=5.5em \hangafter=1
{\white .}\hskip 1em $\rho_\text{isum}(\mathfrak{t})$ =
 $( \frac{1}{3},
\frac{11}{24},
\frac{23}{24} )
\oplus
( 0,
\frac{1}{3} )
\oplus
( 0 )
$,

\vskip 0.7ex
\hangindent=5.5em \hangafter=1
{\white .}\hskip 1em $\rho_\text{isum}(\mathfrak{s})$ =
$\mathrm{i}$($0$,
$\sqrt{\frac{1}{2}}$,
$\sqrt{\frac{1}{2}}$;\ \ 
$-\frac{1}{2}$,
$\frac{1}{2}$;\ \ 
$-\frac{1}{2}$)
 $\oplus$
$\mathrm{i}$($-\sqrt{\frac{1}{3}}$,
$\sqrt{\frac{2}{3}}$;\ \ 
$\sqrt{\frac{1}{3}}$)
 $\oplus$
($1$)

Pass. 

 \ \color{black}

 \color{blue}

\noindent 143: (dims,levels) = $(3 , 
2 , 
1;24,
3,
3
)$,
irreps = $3_{8}^{3,0}
\hskip -1.5pt \otimes \hskip -1.5pt
1_{3}^{1,0}\oplus
2_{3}^{1,0}\oplus
1_{3}^{1,0}$,
pord$(\rho_\text{isum}(\mathfrak{t})) = 24$,

\vskip 0.7ex
\hangindent=5.5em \hangafter=1
{\white .}\hskip 1em $\rho_\text{isum}(\mathfrak{t})$ =
 $( \frac{1}{3},
\frac{5}{24},
\frac{17}{24} )
\oplus
( 0,
\frac{1}{3} )
\oplus
( \frac{1}{3} )
$,

\vskip 0.7ex
\hangindent=5.5em \hangafter=1
{\white .}\hskip 1em $\rho_\text{isum}(\mathfrak{s})$ =
$\mathrm{i}$($0$,
$\sqrt{\frac{1}{2}}$,
$\sqrt{\frac{1}{2}}$;\ \ 
$\frac{1}{2}$,
$-\frac{1}{2}$;\ \ 
$\frac{1}{2}$)
 $\oplus$
$\mathrm{i}$($-\sqrt{\frac{1}{3}}$,
$\sqrt{\frac{2}{3}}$;\ \ 
$\sqrt{\frac{1}{3}}$)
 $\oplus$
($1$)

Pass. 

 \ \color{black}

 \color{blue}

\noindent 144: (dims,levels) = $(3 , 
2 , 
1;24,
3,
3
)$,
irreps = $3_{8}^{3,0}
\hskip -1.5pt \otimes \hskip -1.5pt
1_{3}^{1,0}\oplus
2_{3}^{1,4}\oplus
1_{3}^{1,0}$,
pord$(\rho_\text{isum}(\mathfrak{t})) = 24$,

\vskip 0.7ex
\hangindent=5.5em \hangafter=1
{\white .}\hskip 1em $\rho_\text{isum}(\mathfrak{t})$ =
 $( \frac{1}{3},
\frac{5}{24},
\frac{17}{24} )
\oplus
( \frac{1}{3},
\frac{2}{3} )
\oplus
( \frac{1}{3} )
$,

\vskip 0.7ex
\hangindent=5.5em \hangafter=1
{\white .}\hskip 1em $\rho_\text{isum}(\mathfrak{s})$ =
$\mathrm{i}$($0$,
$\sqrt{\frac{1}{2}}$,
$\sqrt{\frac{1}{2}}$;\ \ 
$\frac{1}{2}$,
$-\frac{1}{2}$;\ \ 
$\frac{1}{2}$)
 $\oplus$
$\mathrm{i}$($-\sqrt{\frac{1}{3}}$,
$\sqrt{\frac{2}{3}}$;\ \ 
$\sqrt{\frac{1}{3}}$)
 $\oplus$
($1$)

Pass. 

 \ \color{black}

 \color{blue}

\noindent 145: (dims,levels) = $(3 , 
2 , 
1;24,
3,
3
)$,
irreps = $3_{8}^{3,0}
\hskip -1.5pt \otimes \hskip -1.5pt
1_{3}^{1,0}\oplus
2_{3}^{1,4}\oplus
1_{3}^{1,4}$,
pord$(\rho_\text{isum}(\mathfrak{t})) = 24$,

\vskip 0.7ex
\hangindent=5.5em \hangafter=1
{\white .}\hskip 1em $\rho_\text{isum}(\mathfrak{t})$ =
 $( \frac{1}{3},
\frac{5}{24},
\frac{17}{24} )
\oplus
( \frac{1}{3},
\frac{2}{3} )
\oplus
( \frac{2}{3} )
$,

\vskip 0.7ex
\hangindent=5.5em \hangafter=1
{\white .}\hskip 1em $\rho_\text{isum}(\mathfrak{s})$ =
$\mathrm{i}$($0$,
$\sqrt{\frac{1}{2}}$,
$\sqrt{\frac{1}{2}}$;\ \ 
$\frac{1}{2}$,
$-\frac{1}{2}$;\ \ 
$\frac{1}{2}$)
 $\oplus$
$\mathrm{i}$($-\sqrt{\frac{1}{3}}$,
$\sqrt{\frac{2}{3}}$;\ \ 
$\sqrt{\frac{1}{3}}$)
 $\oplus$
($1$)

Pass. 

 \ \color{black}

 \color{blue}

\noindent 146: (dims,levels) = $(3 , 
2 , 
1;24,
3,
3
)$,
irreps = $3_{8}^{1,0}
\hskip -1.5pt \otimes \hskip -1.5pt
1_{3}^{1,0}\oplus
2_{3}^{1,0}\oplus
1_{3}^{1,0}$,
pord$(\rho_\text{isum}(\mathfrak{t})) = 24$,

\vskip 0.7ex
\hangindent=5.5em \hangafter=1
{\white .}\hskip 1em $\rho_\text{isum}(\mathfrak{t})$ =
 $( \frac{1}{3},
\frac{11}{24},
\frac{23}{24} )
\oplus
( 0,
\frac{1}{3} )
\oplus
( \frac{1}{3} )
$,

\vskip 0.7ex
\hangindent=5.5em \hangafter=1
{\white .}\hskip 1em $\rho_\text{isum}(\mathfrak{s})$ =
$\mathrm{i}$($0$,
$\sqrt{\frac{1}{2}}$,
$\sqrt{\frac{1}{2}}$;\ \ 
$-\frac{1}{2}$,
$\frac{1}{2}$;\ \ 
$-\frac{1}{2}$)
 $\oplus$
$\mathrm{i}$($-\sqrt{\frac{1}{3}}$,
$\sqrt{\frac{2}{3}}$;\ \ 
$\sqrt{\frac{1}{3}}$)
 $\oplus$
($1$)

Pass. 

 \ \color{black}

 \color{blue}

\noindent 147: (dims,levels) = $(3 , 
2 , 
1;24,
3,
3
)$,
irreps = $3_{8}^{1,0}
\hskip -1.5pt \otimes \hskip -1.5pt
1_{3}^{1,0}\oplus
2_{3}^{1,4}\oplus
1_{3}^{1,0}$,
pord$(\rho_\text{isum}(\mathfrak{t})) = 24$,

\vskip 0.7ex
\hangindent=5.5em \hangafter=1
{\white .}\hskip 1em $\rho_\text{isum}(\mathfrak{t})$ =
 $( \frac{1}{3},
\frac{11}{24},
\frac{23}{24} )
\oplus
( \frac{1}{3},
\frac{2}{3} )
\oplus
( \frac{1}{3} )
$,

\vskip 0.7ex
\hangindent=5.5em \hangafter=1
{\white .}\hskip 1em $\rho_\text{isum}(\mathfrak{s})$ =
$\mathrm{i}$($0$,
$\sqrt{\frac{1}{2}}$,
$\sqrt{\frac{1}{2}}$;\ \ 
$-\frac{1}{2}$,
$\frac{1}{2}$;\ \ 
$-\frac{1}{2}$)
 $\oplus$
$\mathrm{i}$($-\sqrt{\frac{1}{3}}$,
$\sqrt{\frac{2}{3}}$;\ \ 
$\sqrt{\frac{1}{3}}$)
 $\oplus$
($1$)

Pass. 

 \ \color{black}

 \color{blue}

\noindent 148: (dims,levels) = $(3 , 
2 , 
1;24,
3,
3
)$,
irreps = $3_{8}^{1,0}
\hskip -1.5pt \otimes \hskip -1.5pt
1_{3}^{1,0}\oplus
2_{3}^{1,4}\oplus
1_{3}^{1,4}$,
pord$(\rho_\text{isum}(\mathfrak{t})) = 24$,

\vskip 0.7ex
\hangindent=5.5em \hangafter=1
{\white .}\hskip 1em $\rho_\text{isum}(\mathfrak{t})$ =
 $( \frac{1}{3},
\frac{11}{24},
\frac{23}{24} )
\oplus
( \frac{1}{3},
\frac{2}{3} )
\oplus
( \frac{2}{3} )
$,

\vskip 0.7ex
\hangindent=5.5em \hangafter=1
{\white .}\hskip 1em $\rho_\text{isum}(\mathfrak{s})$ =
$\mathrm{i}$($0$,
$\sqrt{\frac{1}{2}}$,
$\sqrt{\frac{1}{2}}$;\ \ 
$-\frac{1}{2}$,
$\frac{1}{2}$;\ \ 
$-\frac{1}{2}$)
 $\oplus$
$\mathrm{i}$($-\sqrt{\frac{1}{3}}$,
$\sqrt{\frac{2}{3}}$;\ \ 
$\sqrt{\frac{1}{3}}$)
 $\oplus$
($1$)

Pass. 

 \ \color{black}

\noindent 149: (dims,levels) = $(3 , 
2 , 
1;24,
6,
3
)$,
irreps = $3_{8}^{3,0}
\hskip -1.5pt \otimes \hskip -1.5pt
1_{3}^{1,0}\oplus
2_{2}^{1,0}
\hskip -1.5pt \otimes \hskip -1.5pt
1_{3}^{1,0}\oplus
1_{3}^{1,0}$,
pord$(\rho_\text{isum}(\mathfrak{t})) = 8$,

\vskip 0.7ex
\hangindent=5.5em \hangafter=1
{\white .}\hskip 1em $\rho_\text{isum}(\mathfrak{t})$ =
 $( \frac{1}{3},
\frac{5}{24},
\frac{17}{24} )
\oplus
( \frac{1}{3},
\frac{5}{6} )
\oplus
( \frac{1}{3} )
$,

\vskip 0.7ex
\hangindent=5.5em \hangafter=1
{\white .}\hskip 1em $\rho_\text{isum}(\mathfrak{s})$ =
$\mathrm{i}$($0$,
$\sqrt{\frac{1}{2}}$,
$\sqrt{\frac{1}{2}}$;\ \ 
$\frac{1}{2}$,
$-\frac{1}{2}$;\ \ 
$\frac{1}{2}$)
 $\oplus$
($-\frac{1}{2}$,
$-\sqrt{\frac{3}{4}}$;
$\frac{1}{2}$)
 $\oplus$
($1$)

Fail:
number of self dual objects $|$Tr($\rho(\mathfrak s^2)$)$|$ = 0. Prop. B.4 (1)\
 eqn. (B.16)

 \ \color{black}

\noindent 150: (dims,levels) = $(3 , 
2 , 
1;24,
6,
3
)$,
irreps = $3_{8}^{1,0}
\hskip -1.5pt \otimes \hskip -1.5pt
1_{3}^{1,0}\oplus
2_{2}^{1,0}
\hskip -1.5pt \otimes \hskip -1.5pt
1_{3}^{1,0}\oplus
1_{3}^{1,0}$,
pord$(\rho_\text{isum}(\mathfrak{t})) = 8$,

\vskip 0.7ex
\hangindent=5.5em \hangafter=1
{\white .}\hskip 1em $\rho_\text{isum}(\mathfrak{t})$ =
 $( \frac{1}{3},
\frac{11}{24},
\frac{23}{24} )
\oplus
( \frac{1}{3},
\frac{5}{6} )
\oplus
( \frac{1}{3} )
$,

\vskip 0.7ex
\hangindent=5.5em \hangafter=1
{\white .}\hskip 1em $\rho_\text{isum}(\mathfrak{s})$ =
$\mathrm{i}$($0$,
$\sqrt{\frac{1}{2}}$,
$\sqrt{\frac{1}{2}}$;\ \ 
$-\frac{1}{2}$,
$\frac{1}{2}$;\ \ 
$-\frac{1}{2}$)
 $\oplus$
($-\frac{1}{2}$,
$-\sqrt{\frac{3}{4}}$;
$\frac{1}{2}$)
 $\oplus$
($1$)

Fail:
number of self dual objects $|$Tr($\rho(\mathfrak s^2)$)$|$ = 0. Prop. B.4 (1)\
 eqn. (B.16)

 \ \color{black}

\noindent 151: (dims,levels) = $(3 , 
2 , 
1;24,
6,
6
)$,
irreps = $3_{8}^{3,0}
\hskip -1.5pt \otimes \hskip -1.5pt
1_{3}^{1,0}\oplus
2_{2}^{1,0}
\hskip -1.5pt \otimes \hskip -1.5pt
1_{3}^{1,0}\oplus
1_{3}^{1,0}
\hskip -1.5pt \otimes \hskip -1.5pt
1_{2}^{1,0}$,
pord$(\rho_\text{isum}(\mathfrak{t})) = 8$,

\vskip 0.7ex
\hangindent=5.5em \hangafter=1
{\white .}\hskip 1em $\rho_\text{isum}(\mathfrak{t})$ =
 $( \frac{1}{3},
\frac{5}{24},
\frac{17}{24} )
\oplus
( \frac{1}{3},
\frac{5}{6} )
\oplus
( \frac{5}{6} )
$,

\vskip 0.7ex
\hangindent=5.5em \hangafter=1
{\white .}\hskip 1em $\rho_\text{isum}(\mathfrak{s})$ =
$\mathrm{i}$($0$,
$\sqrt{\frac{1}{2}}$,
$\sqrt{\frac{1}{2}}$;\ \ 
$\frac{1}{2}$,
$-\frac{1}{2}$;\ \ 
$\frac{1}{2}$)
 $\oplus$
($-\frac{1}{2}$,
$-\sqrt{\frac{3}{4}}$;
$\frac{1}{2}$)
 $\oplus$
($-1$)

Fail:
number of self dual objects $|$Tr($\rho(\mathfrak s^2)$)$|$ = 0. Prop. B.4 (1)\
 eqn. (B.16)

 \ \color{black}

\noindent 152: (dims,levels) = $(3 , 
2 , 
1;24,
6,
6
)$,
irreps = $3_{8}^{1,0}
\hskip -1.5pt \otimes \hskip -1.5pt
1_{3}^{1,0}\oplus
2_{2}^{1,0}
\hskip -1.5pt \otimes \hskip -1.5pt
1_{3}^{1,0}\oplus
1_{3}^{1,0}
\hskip -1.5pt \otimes \hskip -1.5pt
1_{2}^{1,0}$,
pord$(\rho_\text{isum}(\mathfrak{t})) = 8$,

\vskip 0.7ex
\hangindent=5.5em \hangafter=1
{\white .}\hskip 1em $\rho_\text{isum}(\mathfrak{t})$ =
 $( \frac{1}{3},
\frac{11}{24},
\frac{23}{24} )
\oplus
( \frac{1}{3},
\frac{5}{6} )
\oplus
( \frac{5}{6} )
$,

\vskip 0.7ex
\hangindent=5.5em \hangafter=1
{\white .}\hskip 1em $\rho_\text{isum}(\mathfrak{s})$ =
$\mathrm{i}$($0$,
$\sqrt{\frac{1}{2}}$,
$\sqrt{\frac{1}{2}}$;\ \ 
$-\frac{1}{2}$,
$\frac{1}{2}$;\ \ 
$-\frac{1}{2}$)
 $\oplus$
($-\frac{1}{2}$,
$-\sqrt{\frac{3}{4}}$;
$\frac{1}{2}$)
 $\oplus$
($-1$)

Fail:
number of self dual objects $|$Tr($\rho(\mathfrak s^2)$)$|$ = 0. Prop. B.4 (1)\
 eqn. (B.16)

 \ \color{black}

 \color{blue}

\noindent 153: (dims,levels) = $(3 , 
2 , 
1;24,
24,
3
)$,
irreps = $3_{8}^{3,0}
\hskip -1.5pt \otimes \hskip -1.5pt
1_{3}^{1,0}\oplus
2_{8}^{1,9}
\hskip -1.5pt \otimes \hskip -1.5pt
1_{3}^{1,0}\oplus
1_{3}^{1,0}$,
pord$(\rho_\text{isum}(\mathfrak{t})) = 8$,

\vskip 0.7ex
\hangindent=5.5em \hangafter=1
{\white .}\hskip 1em $\rho_\text{isum}(\mathfrak{t})$ =
 $( \frac{1}{3},
\frac{5}{24},
\frac{17}{24} )
\oplus
( \frac{5}{24},
\frac{11}{24} )
\oplus
( \frac{1}{3} )
$,

\vskip 0.7ex
\hangindent=5.5em \hangafter=1
{\white .}\hskip 1em $\rho_\text{isum}(\mathfrak{s})$ =
$\mathrm{i}$($0$,
$\sqrt{\frac{1}{2}}$,
$\sqrt{\frac{1}{2}}$;\ \ 
$\frac{1}{2}$,
$-\frac{1}{2}$;\ \ 
$\frac{1}{2}$)
 $\oplus$
$\mathrm{i}$($\sqrt{\frac{1}{2}}$,
$\sqrt{\frac{1}{2}}$;\ \ 
$-\sqrt{\frac{1}{2}}$)
 $\oplus$
($1$)

Pass. 

 \ \color{black}

\noindent 154: (dims,levels) = $(3 , 
2 , 
1;24,
24,
3
)$,
irreps = $3_{8}^{3,0}
\hskip -1.5pt \otimes \hskip -1.5pt
1_{3}^{1,0}\oplus
2_{8}^{1,6}
\hskip -1.5pt \otimes \hskip -1.5pt
1_{3}^{1,0}\oplus
1_{3}^{1,0}$,
pord$(\rho_\text{isum}(\mathfrak{t})) = 8$,

\vskip 0.7ex
\hangindent=5.5em \hangafter=1
{\white .}\hskip 1em $\rho_\text{isum}(\mathfrak{t})$ =
 $( \frac{1}{3},
\frac{5}{24},
\frac{17}{24} )
\oplus
( \frac{5}{24},
\frac{23}{24} )
\oplus
( \frac{1}{3} )
$,

\vskip 0.7ex
\hangindent=5.5em \hangafter=1
{\white .}\hskip 1em $\rho_\text{isum}(\mathfrak{s})$ =
$\mathrm{i}$($0$,
$\sqrt{\frac{1}{2}}$,
$\sqrt{\frac{1}{2}}$;\ \ 
$\frac{1}{2}$,
$-\frac{1}{2}$;\ \ 
$\frac{1}{2}$)
 $\oplus$
($-\sqrt{\frac{1}{2}}$,
$\sqrt{\frac{1}{2}}$;
$\sqrt{\frac{1}{2}}$)
 $\oplus$
($1$)

Fail:
number of self dual objects $|$Tr($\rho(\mathfrak s^2)$)$|$ = 0. Prop. B.4 (1)\
 eqn. (B.16)

 \ \color{black}

\noindent 155: (dims,levels) = $(3 , 
2 , 
1;24,
24,
3
)$,
irreps = $3_{8}^{3,0}
\hskip -1.5pt \otimes \hskip -1.5pt
1_{3}^{1,0}\oplus
2_{8}^{1,0}
\hskip -1.5pt \otimes \hskip -1.5pt
1_{3}^{1,0}\oplus
1_{3}^{1,0}$,
pord$(\rho_\text{isum}(\mathfrak{t})) = 8$,

\vskip 0.7ex
\hangindent=5.5em \hangafter=1
{\white .}\hskip 1em $\rho_\text{isum}(\mathfrak{t})$ =
 $( \frac{1}{3},
\frac{5}{24},
\frac{17}{24} )
\oplus
( \frac{11}{24},
\frac{17}{24} )
\oplus
( \frac{1}{3} )
$,

\vskip 0.7ex
\hangindent=5.5em \hangafter=1
{\white .}\hskip 1em $\rho_\text{isum}(\mathfrak{s})$ =
$\mathrm{i}$($0$,
$\sqrt{\frac{1}{2}}$,
$\sqrt{\frac{1}{2}}$;\ \ 
$\frac{1}{2}$,
$-\frac{1}{2}$;\ \ 
$\frac{1}{2}$)
 $\oplus$
($-\sqrt{\frac{1}{2}}$,
$\sqrt{\frac{1}{2}}$;
$\sqrt{\frac{1}{2}}$)
 $\oplus$
($1$)

Fail:
number of self dual objects $|$Tr($\rho(\mathfrak s^2)$)$|$ = 0. Prop. B.4 (1)\
 eqn. (B.16)

 \ \color{black}

 \color{blue}

\noindent 156: (dims,levels) = $(3 , 
2 , 
1;24,
24,
3
)$,
irreps = $3_{8}^{3,0}
\hskip -1.5pt \otimes \hskip -1.5pt
1_{3}^{1,0}\oplus
2_{8}^{1,3}
\hskip -1.5pt \otimes \hskip -1.5pt
1_{3}^{1,0}\oplus
1_{3}^{1,0}$,
pord$(\rho_\text{isum}(\mathfrak{t})) = 8$,

\vskip 0.7ex
\hangindent=5.5em \hangafter=1
{\white .}\hskip 1em $\rho_\text{isum}(\mathfrak{t})$ =
 $( \frac{1}{3},
\frac{5}{24},
\frac{17}{24} )
\oplus
( \frac{17}{24},
\frac{23}{24} )
\oplus
( \frac{1}{3} )
$,

\vskip 0.7ex
\hangindent=5.5em \hangafter=1
{\white .}\hskip 1em $\rho_\text{isum}(\mathfrak{s})$ =
$\mathrm{i}$($0$,
$\sqrt{\frac{1}{2}}$,
$\sqrt{\frac{1}{2}}$;\ \ 
$\frac{1}{2}$,
$-\frac{1}{2}$;\ \ 
$\frac{1}{2}$)
 $\oplus$
$\mathrm{i}$($-\sqrt{\frac{1}{2}}$,
$\sqrt{\frac{1}{2}}$;\ \ 
$\sqrt{\frac{1}{2}}$)
 $\oplus$
($1$)

Pass. 

 \ \color{black}

 \color{blue}

\noindent 157: (dims,levels) = $(3 , 
2 , 
1;24,
24,
3
)$,
irreps = $3_{8}^{1,0}
\hskip -1.5pt \otimes \hskip -1.5pt
1_{3}^{1,0}\oplus
2_{8}^{1,9}
\hskip -1.5pt \otimes \hskip -1.5pt
1_{3}^{1,0}\oplus
1_{3}^{1,0}$,
pord$(\rho_\text{isum}(\mathfrak{t})) = 8$,

\vskip 0.7ex
\hangindent=5.5em \hangafter=1
{\white .}\hskip 1em $\rho_\text{isum}(\mathfrak{t})$ =
 $( \frac{1}{3},
\frac{11}{24},
\frac{23}{24} )
\oplus
( \frac{5}{24},
\frac{11}{24} )
\oplus
( \frac{1}{3} )
$,

\vskip 0.7ex
\hangindent=5.5em \hangafter=1
{\white .}\hskip 1em $\rho_\text{isum}(\mathfrak{s})$ =
$\mathrm{i}$($0$,
$\sqrt{\frac{1}{2}}$,
$\sqrt{\frac{1}{2}}$;\ \ 
$-\frac{1}{2}$,
$\frac{1}{2}$;\ \ 
$-\frac{1}{2}$)
 $\oplus$
$\mathrm{i}$($\sqrt{\frac{1}{2}}$,
$\sqrt{\frac{1}{2}}$;\ \ 
$-\sqrt{\frac{1}{2}}$)
 $\oplus$
($1$)

Pass. 

 \ \color{black}

\noindent 158: (dims,levels) = $(3 , 
2 , 
1;24,
24,
3
)$,
irreps = $3_{8}^{1,0}
\hskip -1.5pt \otimes \hskip -1.5pt
1_{3}^{1,0}\oplus
2_{8}^{1,6}
\hskip -1.5pt \otimes \hskip -1.5pt
1_{3}^{1,0}\oplus
1_{3}^{1,0}$,
pord$(\rho_\text{isum}(\mathfrak{t})) = 8$,

\vskip 0.7ex
\hangindent=5.5em \hangafter=1
{\white .}\hskip 1em $\rho_\text{isum}(\mathfrak{t})$ =
 $( \frac{1}{3},
\frac{11}{24},
\frac{23}{24} )
\oplus
( \frac{5}{24},
\frac{23}{24} )
\oplus
( \frac{1}{3} )
$,

\vskip 0.7ex
\hangindent=5.5em \hangafter=1
{\white .}\hskip 1em $\rho_\text{isum}(\mathfrak{s})$ =
$\mathrm{i}$($0$,
$\sqrt{\frac{1}{2}}$,
$\sqrt{\frac{1}{2}}$;\ \ 
$-\frac{1}{2}$,
$\frac{1}{2}$;\ \ 
$-\frac{1}{2}$)
 $\oplus$
($-\sqrt{\frac{1}{2}}$,
$\sqrt{\frac{1}{2}}$;
$\sqrt{\frac{1}{2}}$)
 $\oplus$
($1$)

Fail:
number of self dual objects $|$Tr($\rho(\mathfrak s^2)$)$|$ = 0. Prop. B.4 (1)\
 eqn. (B.16)

 \ \color{black}

\noindent 159: (dims,levels) = $(3 , 
2 , 
1;24,
24,
3
)$,
irreps = $3_{8}^{1,0}
\hskip -1.5pt \otimes \hskip -1.5pt
1_{3}^{1,0}\oplus
2_{8}^{1,0}
\hskip -1.5pt \otimes \hskip -1.5pt
1_{3}^{1,0}\oplus
1_{3}^{1,0}$,
pord$(\rho_\text{isum}(\mathfrak{t})) = 8$,

\vskip 0.7ex
\hangindent=5.5em \hangafter=1
{\white .}\hskip 1em $\rho_\text{isum}(\mathfrak{t})$ =
 $( \frac{1}{3},
\frac{11}{24},
\frac{23}{24} )
\oplus
( \frac{11}{24},
\frac{17}{24} )
\oplus
( \frac{1}{3} )
$,

\vskip 0.7ex
\hangindent=5.5em \hangafter=1
{\white .}\hskip 1em $\rho_\text{isum}(\mathfrak{s})$ =
$\mathrm{i}$($0$,
$\sqrt{\frac{1}{2}}$,
$\sqrt{\frac{1}{2}}$;\ \ 
$-\frac{1}{2}$,
$\frac{1}{2}$;\ \ 
$-\frac{1}{2}$)
 $\oplus$
($-\sqrt{\frac{1}{2}}$,
$\sqrt{\frac{1}{2}}$;
$\sqrt{\frac{1}{2}}$)
 $\oplus$
($1$)

Fail:
number of self dual objects $|$Tr($\rho(\mathfrak s^2)$)$|$ = 0. Prop. B.4 (1)\
 eqn. (B.16)

 \ \color{black}

 \color{blue}

\noindent 160: (dims,levels) = $(3 , 
2 , 
1;24,
24,
3
)$,
irreps = $3_{8}^{1,0}
\hskip -1.5pt \otimes \hskip -1.5pt
1_{3}^{1,0}\oplus
2_{8}^{1,3}
\hskip -1.5pt \otimes \hskip -1.5pt
1_{3}^{1,0}\oplus
1_{3}^{1,0}$,
pord$(\rho_\text{isum}(\mathfrak{t})) = 8$,

\vskip 0.7ex
\hangindent=5.5em \hangafter=1
{\white .}\hskip 1em $\rho_\text{isum}(\mathfrak{t})$ =
 $( \frac{1}{3},
\frac{11}{24},
\frac{23}{24} )
\oplus
( \frac{17}{24},
\frac{23}{24} )
\oplus
( \frac{1}{3} )
$,

\vskip 0.7ex
\hangindent=5.5em \hangafter=1
{\white .}\hskip 1em $\rho_\text{isum}(\mathfrak{s})$ =
$\mathrm{i}$($0$,
$\sqrt{\frac{1}{2}}$,
$\sqrt{\frac{1}{2}}$;\ \ 
$-\frac{1}{2}$,
$\frac{1}{2}$;\ \ 
$-\frac{1}{2}$)
 $\oplus$
$\mathrm{i}$($-\sqrt{\frac{1}{2}}$,
$\sqrt{\frac{1}{2}}$;\ \ 
$\sqrt{\frac{1}{2}}$)
 $\oplus$
($1$)

Pass. 

 \ \color{black}

\noindent 161: (dims,levels) = $(3 , 
2 , 
1;30,
6,
2
)$,
irreps = $3_{5}^{1}
\hskip -1.5pt \otimes \hskip -1.5pt
1_{3}^{1,0}
\hskip -1.5pt \otimes \hskip -1.5pt
1_{2}^{1,0}\oplus
2_{3}^{1,0}
\hskip -1.5pt \otimes \hskip -1.5pt
1_{2}^{1,0}\oplus
1_{2}^{1,0}$,
pord$(\rho_\text{isum}(\mathfrak{t})) = 15$,

\vskip 0.7ex
\hangindent=5.5em \hangafter=1
{\white .}\hskip 1em $\rho_\text{isum}(\mathfrak{t})$ =
 $( \frac{5}{6},
\frac{1}{30},
\frac{19}{30} )
\oplus
( \frac{1}{2},
\frac{5}{6} )
\oplus
( \frac{1}{2} )
$,

\vskip 0.7ex
\hangindent=5.5em \hangafter=1
{\white .}\hskip 1em $\rho_\text{isum}(\mathfrak{s})$ =
($-\sqrt{\frac{1}{5}}$,
$-\sqrt{\frac{2}{5}}$,
$-\sqrt{\frac{2}{5}}$;
$\frac{5+\sqrt{5}}{10}$,
$\frac{-5+\sqrt{5}}{10}$;
$\frac{5+\sqrt{5}}{10}$)
 $\oplus$
$\mathrm{i}$($\sqrt{\frac{1}{3}}$,
$\sqrt{\frac{2}{3}}$;\ \ 
$-\sqrt{\frac{1}{3}}$)
 $\oplus$
($-1$)

Fail:
all $\theta$-eigenspaces that can contain unit
 have Tr$_{E_\theta}(C) \leq 0 $. Prop. B.5 (5) eqn. (B.29)

 \ \color{black}

\noindent 162: (dims,levels) = $(3 , 
2 , 
1;30,
6,
2
)$,
irreps = $3_{5}^{3}
\hskip -1.5pt \otimes \hskip -1.5pt
1_{3}^{1,0}
\hskip -1.5pt \otimes \hskip -1.5pt
1_{2}^{1,0}\oplus
2_{3}^{1,0}
\hskip -1.5pt \otimes \hskip -1.5pt
1_{2}^{1,0}\oplus
1_{2}^{1,0}$,
pord$(\rho_\text{isum}(\mathfrak{t})) = 15$,

\vskip 0.7ex
\hangindent=5.5em \hangafter=1
{\white .}\hskip 1em $\rho_\text{isum}(\mathfrak{t})$ =
 $( \frac{5}{6},
\frac{7}{30},
\frac{13}{30} )
\oplus
( \frac{1}{2},
\frac{5}{6} )
\oplus
( \frac{1}{2} )
$,

\vskip 0.7ex
\hangindent=5.5em \hangafter=1
{\white .}\hskip 1em $\rho_\text{isum}(\mathfrak{s})$ =
($\sqrt{\frac{1}{5}}$,
$-\sqrt{\frac{2}{5}}$,
$-\sqrt{\frac{2}{5}}$;
$\frac{5-\sqrt{5}}{10}$,
$-\frac{5+\sqrt{5}}{10}$;
$\frac{5-\sqrt{5}}{10}$)
 $\oplus$
$\mathrm{i}$($\sqrt{\frac{1}{3}}$,
$\sqrt{\frac{2}{3}}$;\ \ 
$-\sqrt{\frac{1}{3}}$)
 $\oplus$
($-1$)

Fail:
all $\theta$-eigenspaces that can contain unit
 have Tr$_{E_\theta}(C) \leq 0 $. Prop. B.5 (5) eqn. (B.29)

 \ \color{black}

 \color{blue}

\noindent 163: (dims,levels) = $(3 , 
2 , 
1;30,
6,
3
)$,
irreps = $3_{5}^{1}
\hskip -1.5pt \otimes \hskip -1.5pt
1_{3}^{1,0}
\hskip -1.5pt \otimes \hskip -1.5pt
1_{2}^{1,0}\oplus
2_{2}^{1,0}
\hskip -1.5pt \otimes \hskip -1.5pt
1_{3}^{1,0}\oplus
1_{3}^{1,0}$,
pord$(\rho_\text{isum}(\mathfrak{t})) = 10$,

\vskip 0.7ex
\hangindent=5.5em \hangafter=1
{\white .}\hskip 1em $\rho_\text{isum}(\mathfrak{t})$ =
 $( \frac{5}{6},
\frac{1}{30},
\frac{19}{30} )
\oplus
( \frac{1}{3},
\frac{5}{6} )
\oplus
( \frac{1}{3} )
$,

\vskip 0.7ex
\hangindent=5.5em \hangafter=1
{\white .}\hskip 1em $\rho_\text{isum}(\mathfrak{s})$ =
($-\sqrt{\frac{1}{5}}$,
$-\sqrt{\frac{2}{5}}$,
$-\sqrt{\frac{2}{5}}$;
$\frac{5+\sqrt{5}}{10}$,
$\frac{-5+\sqrt{5}}{10}$;
$\frac{5+\sqrt{5}}{10}$)
 $\oplus$
($-\frac{1}{2}$,
$-\sqrt{\frac{3}{4}}$;
$\frac{1}{2}$)
 $\oplus$
($1$)

Pass. 

 \ \color{black}

 \color{blue}

\noindent 164: (dims,levels) = $(3 , 
2 , 
1;30,
6,
3
)$,
irreps = $3_{5}^{3}
\hskip -1.5pt \otimes \hskip -1.5pt
1_{3}^{1,0}
\hskip -1.5pt \otimes \hskip -1.5pt
1_{2}^{1,0}\oplus
2_{2}^{1,0}
\hskip -1.5pt \otimes \hskip -1.5pt
1_{3}^{1,0}\oplus
1_{3}^{1,0}$,
pord$(\rho_\text{isum}(\mathfrak{t})) = 10$,

\vskip 0.7ex
\hangindent=5.5em \hangafter=1
{\white .}\hskip 1em $\rho_\text{isum}(\mathfrak{t})$ =
 $( \frac{5}{6},
\frac{7}{30},
\frac{13}{30} )
\oplus
( \frac{1}{3},
\frac{5}{6} )
\oplus
( \frac{1}{3} )
$,

\vskip 0.7ex
\hangindent=5.5em \hangafter=1
{\white .}\hskip 1em $\rho_\text{isum}(\mathfrak{s})$ =
($\sqrt{\frac{1}{5}}$,
$-\sqrt{\frac{2}{5}}$,
$-\sqrt{\frac{2}{5}}$;
$\frac{5-\sqrt{5}}{10}$,
$-\frac{5+\sqrt{5}}{10}$;
$\frac{5-\sqrt{5}}{10}$)
 $\oplus$
($-\frac{1}{2}$,
$-\sqrt{\frac{3}{4}}$;
$\frac{1}{2}$)
 $\oplus$
($1$)

Pass. 

 \ \color{black}

\noindent 165: (dims,levels) = $(3 , 
2 , 
1;30,
6,
6
)$,
irreps = $3_{5}^{1}
\hskip -1.5pt \otimes \hskip -1.5pt
1_{3}^{1,0}
\hskip -1.5pt \otimes \hskip -1.5pt
1_{2}^{1,0}\oplus
2_{3}^{1,4}
\hskip -1.5pt \otimes \hskip -1.5pt
1_{2}^{1,0}\oplus
1_{3}^{1,4}
\hskip -1.5pt \otimes \hskip -1.5pt
1_{2}^{1,0}$,
pord$(\rho_\text{isum}(\mathfrak{t})) = 15$,

\vskip 0.7ex
\hangindent=5.5em \hangafter=1
{\white .}\hskip 1em $\rho_\text{isum}(\mathfrak{t})$ =
 $( \frac{5}{6},
\frac{1}{30},
\frac{19}{30} )
\oplus
( \frac{1}{6},
\frac{5}{6} )
\oplus
( \frac{1}{6} )
$,

\vskip 0.7ex
\hangindent=5.5em \hangafter=1
{\white .}\hskip 1em $\rho_\text{isum}(\mathfrak{s})$ =
($-\sqrt{\frac{1}{5}}$,
$-\sqrt{\frac{2}{5}}$,
$-\sqrt{\frac{2}{5}}$;
$\frac{5+\sqrt{5}}{10}$,
$\frac{-5+\sqrt{5}}{10}$;
$\frac{5+\sqrt{5}}{10}$)
 $\oplus$
$\mathrm{i}$($-\sqrt{\frac{1}{3}}$,
$\sqrt{\frac{2}{3}}$;\ \ 
$\sqrt{\frac{1}{3}}$)
 $\oplus$
($-1$)

Fail:
all $\theta$-eigenspaces that can contain unit
 have Tr$_{E_\theta}(C) \leq 0 $. Prop. B.5 (5) eqn. (B.29)

 \ \color{black}

\noindent 166: (dims,levels) = $(3 , 
2 , 
1;30,
6,
6
)$,
irreps = $3_{5}^{1}
\hskip -1.5pt \otimes \hskip -1.5pt
1_{3}^{1,0}
\hskip -1.5pt \otimes \hskip -1.5pt
1_{2}^{1,0}\oplus
2_{3}^{1,4}
\hskip -1.5pt \otimes \hskip -1.5pt
1_{2}^{1,0}\oplus
1_{3}^{1,0}
\hskip -1.5pt \otimes \hskip -1.5pt
1_{2}^{1,0}$,
pord$(\rho_\text{isum}(\mathfrak{t})) = 15$,

\vskip 0.7ex
\hangindent=5.5em \hangafter=1
{\white .}\hskip 1em $\rho_\text{isum}(\mathfrak{t})$ =
 $( \frac{5}{6},
\frac{1}{30},
\frac{19}{30} )
\oplus
( \frac{1}{6},
\frac{5}{6} )
\oplus
( \frac{5}{6} )
$,

\vskip 0.7ex
\hangindent=5.5em \hangafter=1
{\white .}\hskip 1em $\rho_\text{isum}(\mathfrak{s})$ =
($-\sqrt{\frac{1}{5}}$,
$-\sqrt{\frac{2}{5}}$,
$-\sqrt{\frac{2}{5}}$;
$\frac{5+\sqrt{5}}{10}$,
$\frac{-5+\sqrt{5}}{10}$;
$\frac{5+\sqrt{5}}{10}$)
 $\oplus$
$\mathrm{i}$($-\sqrt{\frac{1}{3}}$,
$\sqrt{\frac{2}{3}}$;\ \ 
$\sqrt{\frac{1}{3}}$)
 $\oplus$
($-1$)

Fail:
Tr$_I(C) = -1 <$  0 for I = [ 1/6 ]. Prop. B.4 (1) eqn. (B.18)

 \ \color{black}

 \color{blue}

\noindent 167: (dims,levels) = $(3 , 
2 , 
1;30,
6,
6
)$,
irreps = $3_{5}^{1}
\hskip -1.5pt \otimes \hskip -1.5pt
1_{3}^{1,0}
\hskip -1.5pt \otimes \hskip -1.5pt
1_{2}^{1,0}\oplus
2_{2}^{1,0}
\hskip -1.5pt \otimes \hskip -1.5pt
1_{3}^{1,0}\oplus
1_{3}^{1,0}
\hskip -1.5pt \otimes \hskip -1.5pt
1_{2}^{1,0}$,
pord$(\rho_\text{isum}(\mathfrak{t})) = 10$,

\vskip 0.7ex
\hangindent=5.5em \hangafter=1
{\white .}\hskip 1em $\rho_\text{isum}(\mathfrak{t})$ =
 $( \frac{5}{6},
\frac{1}{30},
\frac{19}{30} )
\oplus
( \frac{1}{3},
\frac{5}{6} )
\oplus
( \frac{5}{6} )
$,

\vskip 0.7ex
\hangindent=5.5em \hangafter=1
{\white .}\hskip 1em $\rho_\text{isum}(\mathfrak{s})$ =
($-\sqrt{\frac{1}{5}}$,
$-\sqrt{\frac{2}{5}}$,
$-\sqrt{\frac{2}{5}}$;
$\frac{5+\sqrt{5}}{10}$,
$\frac{-5+\sqrt{5}}{10}$;
$\frac{5+\sqrt{5}}{10}$)
 $\oplus$
($-\frac{1}{2}$,
$-\sqrt{\frac{3}{4}}$;
$\frac{1}{2}$)
 $\oplus$
($-1$)

Pass. 

 \ \color{black}

\noindent 168: (dims,levels) = $(3 , 
2 , 
1;30,
6,
6
)$,
irreps = $3_{5}^{1}
\hskip -1.5pt \otimes \hskip -1.5pt
1_{3}^{1,0}
\hskip -1.5pt \otimes \hskip -1.5pt
1_{2}^{1,0}\oplus
2_{3}^{1,0}
\hskip -1.5pt \otimes \hskip -1.5pt
1_{2}^{1,0}\oplus
1_{3}^{1,0}
\hskip -1.5pt \otimes \hskip -1.5pt
1_{2}^{1,0}$,
pord$(\rho_\text{isum}(\mathfrak{t})) = 15$,

\vskip 0.7ex
\hangindent=5.5em \hangafter=1
{\white .}\hskip 1em $\rho_\text{isum}(\mathfrak{t})$ =
 $( \frac{5}{6},
\frac{1}{30},
\frac{19}{30} )
\oplus
( \frac{1}{2},
\frac{5}{6} )
\oplus
( \frac{5}{6} )
$,

\vskip 0.7ex
\hangindent=5.5em \hangafter=1
{\white .}\hskip 1em $\rho_\text{isum}(\mathfrak{s})$ =
($-\sqrt{\frac{1}{5}}$,
$-\sqrt{\frac{2}{5}}$,
$-\sqrt{\frac{2}{5}}$;
$\frac{5+\sqrt{5}}{10}$,
$\frac{-5+\sqrt{5}}{10}$;
$\frac{5+\sqrt{5}}{10}$)
 $\oplus$
$\mathrm{i}$($\sqrt{\frac{1}{3}}$,
$\sqrt{\frac{2}{3}}$;\ \ 
$-\sqrt{\frac{1}{3}}$)
 $\oplus$
($-1$)

Fail:
Tr$_I(C) = -1 <$  0 for I = [ 1/2 ]. Prop. B.4 (1) eqn. (B.18)

 \ \color{black}

\noindent 169: (dims,levels) = $(3 , 
2 , 
1;30,
6,
6
)$,
irreps = $3_{5}^{3}
\hskip -1.5pt \otimes \hskip -1.5pt
1_{3}^{1,0}
\hskip -1.5pt \otimes \hskip -1.5pt
1_{2}^{1,0}\oplus
2_{3}^{1,4}
\hskip -1.5pt \otimes \hskip -1.5pt
1_{2}^{1,0}\oplus
1_{3}^{1,4}
\hskip -1.5pt \otimes \hskip -1.5pt
1_{2}^{1,0}$,
pord$(\rho_\text{isum}(\mathfrak{t})) = 15$,

\vskip 0.7ex
\hangindent=5.5em \hangafter=1
{\white .}\hskip 1em $\rho_\text{isum}(\mathfrak{t})$ =
 $( \frac{5}{6},
\frac{7}{30},
\frac{13}{30} )
\oplus
( \frac{1}{6},
\frac{5}{6} )
\oplus
( \frac{1}{6} )
$,

\vskip 0.7ex
\hangindent=5.5em \hangafter=1
{\white .}\hskip 1em $\rho_\text{isum}(\mathfrak{s})$ =
($\sqrt{\frac{1}{5}}$,
$-\sqrt{\frac{2}{5}}$,
$-\sqrt{\frac{2}{5}}$;
$\frac{5-\sqrt{5}}{10}$,
$-\frac{5+\sqrt{5}}{10}$;
$\frac{5-\sqrt{5}}{10}$)
 $\oplus$
$\mathrm{i}$($-\sqrt{\frac{1}{3}}$,
$\sqrt{\frac{2}{3}}$;\ \ 
$\sqrt{\frac{1}{3}}$)
 $\oplus$
($-1$)

Fail:
all $\theta$-eigenspaces that can contain unit
 have Tr$_{E_\theta}(C) \leq 0 $. Prop. B.5 (5) eqn. (B.29)

 \ \color{black}

\noindent 170: (dims,levels) = $(3 , 
2 , 
1;30,
6,
6
)$,
irreps = $3_{5}^{3}
\hskip -1.5pt \otimes \hskip -1.5pt
1_{3}^{1,0}
\hskip -1.5pt \otimes \hskip -1.5pt
1_{2}^{1,0}\oplus
2_{3}^{1,4}
\hskip -1.5pt \otimes \hskip -1.5pt
1_{2}^{1,0}\oplus
1_{3}^{1,0}
\hskip -1.5pt \otimes \hskip -1.5pt
1_{2}^{1,0}$,
pord$(\rho_\text{isum}(\mathfrak{t})) = 15$,

\vskip 0.7ex
\hangindent=5.5em \hangafter=1
{\white .}\hskip 1em $\rho_\text{isum}(\mathfrak{t})$ =
 $( \frac{5}{6},
\frac{7}{30},
\frac{13}{30} )
\oplus
( \frac{1}{6},
\frac{5}{6} )
\oplus
( \frac{5}{6} )
$,

\vskip 0.7ex
\hangindent=5.5em \hangafter=1
{\white .}\hskip 1em $\rho_\text{isum}(\mathfrak{s})$ =
($\sqrt{\frac{1}{5}}$,
$-\sqrt{\frac{2}{5}}$,
$-\sqrt{\frac{2}{5}}$;
$\frac{5-\sqrt{5}}{10}$,
$-\frac{5+\sqrt{5}}{10}$;
$\frac{5-\sqrt{5}}{10}$)
 $\oplus$
$\mathrm{i}$($-\sqrt{\frac{1}{3}}$,
$\sqrt{\frac{2}{3}}$;\ \ 
$\sqrt{\frac{1}{3}}$)
 $\oplus$
($-1$)

Fail:
Tr$_I(C) = -1 <$  0 for I = [ 1/6 ]. Prop. B.4 (1) eqn. (B.18)

 \ \color{black}

 \color{blue}

\noindent 171: (dims,levels) = $(3 , 
2 , 
1;30,
6,
6
)$,
irreps = $3_{5}^{3}
\hskip -1.5pt \otimes \hskip -1.5pt
1_{3}^{1,0}
\hskip -1.5pt \otimes \hskip -1.5pt
1_{2}^{1,0}\oplus
2_{2}^{1,0}
\hskip -1.5pt \otimes \hskip -1.5pt
1_{3}^{1,0}\oplus
1_{3}^{1,0}
\hskip -1.5pt \otimes \hskip -1.5pt
1_{2}^{1,0}$,
pord$(\rho_\text{isum}(\mathfrak{t})) = 10$,

\vskip 0.7ex
\hangindent=5.5em \hangafter=1
{\white .}\hskip 1em $\rho_\text{isum}(\mathfrak{t})$ =
 $( \frac{5}{6},
\frac{7}{30},
\frac{13}{30} )
\oplus
( \frac{1}{3},
\frac{5}{6} )
\oplus
( \frac{5}{6} )
$,

\vskip 0.7ex
\hangindent=5.5em \hangafter=1
{\white .}\hskip 1em $\rho_\text{isum}(\mathfrak{s})$ =
($\sqrt{\frac{1}{5}}$,
$-\sqrt{\frac{2}{5}}$,
$-\sqrt{\frac{2}{5}}$;
$\frac{5-\sqrt{5}}{10}$,
$-\frac{5+\sqrt{5}}{10}$;
$\frac{5-\sqrt{5}}{10}$)
 $\oplus$
($-\frac{1}{2}$,
$-\sqrt{\frac{3}{4}}$;
$\frac{1}{2}$)
 $\oplus$
($-1$)

Pass. 

 \ \color{black}

\noindent 172: (dims,levels) = $(3 , 
2 , 
1;30,
6,
6
)$,
irreps = $3_{5}^{3}
\hskip -1.5pt \otimes \hskip -1.5pt
1_{3}^{1,0}
\hskip -1.5pt \otimes \hskip -1.5pt
1_{2}^{1,0}\oplus
2_{3}^{1,0}
\hskip -1.5pt \otimes \hskip -1.5pt
1_{2}^{1,0}\oplus
1_{3}^{1,0}
\hskip -1.5pt \otimes \hskip -1.5pt
1_{2}^{1,0}$,
pord$(\rho_\text{isum}(\mathfrak{t})) = 15$,

\vskip 0.7ex
\hangindent=5.5em \hangafter=1
{\white .}\hskip 1em $\rho_\text{isum}(\mathfrak{t})$ =
 $( \frac{5}{6},
\frac{7}{30},
\frac{13}{30} )
\oplus
( \frac{1}{2},
\frac{5}{6} )
\oplus
( \frac{5}{6} )
$,

\vskip 0.7ex
\hangindent=5.5em \hangafter=1
{\white .}\hskip 1em $\rho_\text{isum}(\mathfrak{s})$ =
($\sqrt{\frac{1}{5}}$,
$-\sqrt{\frac{2}{5}}$,
$-\sqrt{\frac{2}{5}}$;
$\frac{5-\sqrt{5}}{10}$,
$-\frac{5+\sqrt{5}}{10}$;
$\frac{5-\sqrt{5}}{10}$)
 $\oplus$
$\mathrm{i}$($\sqrt{\frac{1}{3}}$,
$\sqrt{\frac{2}{3}}$;\ \ 
$-\sqrt{\frac{1}{3}}$)
 $\oplus$
($-1$)

Fail:
Tr$_I(C) = -1 <$  0 for I = [ 1/2 ]. Prop. B.4 (1) eqn. (B.18)

 \ \color{black}

 \color{blue}

\noindent 173: (dims,levels) = $(3 , 
2 , 
1;30,
30,
6
)$,
irreps = $3_{5}^{1}
\hskip -1.5pt \otimes \hskip -1.5pt
1_{3}^{1,0}
\hskip -1.5pt \otimes \hskip -1.5pt
1_{2}^{1,0}\oplus
2_{5}^{1}
\hskip -1.5pt \otimes \hskip -1.5pt
1_{3}^{1,0}
\hskip -1.5pt \otimes \hskip -1.5pt
1_{2}^{1,0}\oplus
1_{3}^{1,0}
\hskip -1.5pt \otimes \hskip -1.5pt
1_{2}^{1,0}$,
pord$(\rho_\text{isum}(\mathfrak{t})) = 5$,

\vskip 0.7ex
\hangindent=5.5em \hangafter=1
{\white .}\hskip 1em $\rho_\text{isum}(\mathfrak{t})$ =
 $( \frac{5}{6},
\frac{1}{30},
\frac{19}{30} )
\oplus
( \frac{1}{30},
\frac{19}{30} )
\oplus
( \frac{5}{6} )
$,

\vskip 0.7ex
\hangindent=5.5em \hangafter=1
{\white .}\hskip 1em $\rho_\text{isum}(\mathfrak{s})$ =
($-\sqrt{\frac{1}{5}}$,
$-\sqrt{\frac{2}{5}}$,
$-\sqrt{\frac{2}{5}}$;
$\frac{5+\sqrt{5}}{10}$,
$\frac{-5+\sqrt{5}}{10}$;
$\frac{5+\sqrt{5}}{10}$)
 $\oplus$
$\mathrm{i}$($\frac{1}{\sqrt{5}}c^{3}_{20}
$,
$\frac{1}{\sqrt{5}}c^{1}_{20}
$;\ \ 
$-\frac{1}{\sqrt{5}}c^{3}_{20}
$)
 $\oplus$
($-1$)

Pass. 

 \ \color{black}

 \color{blue}

\noindent 174: (dims,levels) = $(3 , 
2 , 
1;30,
30,
6
)$,
irreps = $3_{5}^{3}
\hskip -1.5pt \otimes \hskip -1.5pt
1_{3}^{1,0}
\hskip -1.5pt \otimes \hskip -1.5pt
1_{2}^{1,0}\oplus
2_{5}^{2}
\hskip -1.5pt \otimes \hskip -1.5pt
1_{3}^{1,0}
\hskip -1.5pt \otimes \hskip -1.5pt
1_{2}^{1,0}\oplus
1_{3}^{1,0}
\hskip -1.5pt \otimes \hskip -1.5pt
1_{2}^{1,0}$,
pord$(\rho_\text{isum}(\mathfrak{t})) = 5$,

\vskip 0.7ex
\hangindent=5.5em \hangafter=1
{\white .}\hskip 1em $\rho_\text{isum}(\mathfrak{t})$ =
 $( \frac{5}{6},
\frac{7}{30},
\frac{13}{30} )
\oplus
( \frac{7}{30},
\frac{13}{30} )
\oplus
( \frac{5}{6} )
$,

\vskip 0.7ex
\hangindent=5.5em \hangafter=1
{\white .}\hskip 1em $\rho_\text{isum}(\mathfrak{s})$ =
($\sqrt{\frac{1}{5}}$,
$-\sqrt{\frac{2}{5}}$,
$-\sqrt{\frac{2}{5}}$;
$\frac{5-\sqrt{5}}{10}$,
$-\frac{5+\sqrt{5}}{10}$;
$\frac{5-\sqrt{5}}{10}$)
 $\oplus$
$\mathrm{i}$($\frac{1}{\sqrt{5}}c^{1}_{20}
$,
$\frac{1}{\sqrt{5}}c^{3}_{20}
$;\ \ 
$-\frac{1}{\sqrt{5}}c^{1}_{20}
$)
 $\oplus$
($-1$)

Pass. 

 \ \color{black}

\noindent 175: (dims,levels) = $(3 , 
2 , 
1;60,
12,
4
)$,
irreps = $3_{5}^{3}
\hskip -1.5pt \otimes \hskip -1.5pt
1_{4}^{1,0}
\hskip -1.5pt \otimes \hskip -1.5pt
1_{3}^{1,0}\oplus
2_{3}^{1,0}
\hskip -1.5pt \otimes \hskip -1.5pt
1_{4}^{1,0}\oplus
1_{4}^{1,0}$,
pord$(\rho_\text{isum}(\mathfrak{t})) = 15$,

\vskip 0.7ex
\hangindent=5.5em \hangafter=1
{\white .}\hskip 1em $\rho_\text{isum}(\mathfrak{t})$ =
 $( \frac{7}{12},
\frac{11}{60},
\frac{59}{60} )
\oplus
( \frac{1}{4},
\frac{7}{12} )
\oplus
( \frac{1}{4} )
$,

\vskip 0.7ex
\hangindent=5.5em \hangafter=1
{\white .}\hskip 1em $\rho_\text{isum}(\mathfrak{s})$ =
$\mathrm{i}$($-\sqrt{\frac{1}{5}}$,
$\sqrt{\frac{2}{5}}$,
$\sqrt{\frac{2}{5}}$;\ \ 
$\frac{-5+\sqrt{5}}{10}$,
$\frac{5+\sqrt{5}}{10}$;\ \ 
$\frac{-5+\sqrt{5}}{10}$)
 $\oplus$
($\sqrt{\frac{1}{3}}$,
$\sqrt{\frac{2}{3}}$;
$-\sqrt{\frac{1}{3}}$)
 $\oplus$
$\mathrm{i}$($1$)

Fail:
all $\theta$-eigenspaces that can contain unit
 have Tr$_{E_\theta}(C) \leq 0 $. Prop. B.5 (5) eqn. (B.29)

 \ \color{black}

\noindent 176: (dims,levels) = $(3 , 
2 , 
1;60,
12,
4
)$,
irreps = $3_{5}^{1}
\hskip -1.5pt \otimes \hskip -1.5pt
1_{4}^{1,0}
\hskip -1.5pt \otimes \hskip -1.5pt
1_{3}^{1,0}\oplus
2_{3}^{1,0}
\hskip -1.5pt \otimes \hskip -1.5pt
1_{4}^{1,0}\oplus
1_{4}^{1,0}$,
pord$(\rho_\text{isum}(\mathfrak{t})) = 15$,

\vskip 0.7ex
\hangindent=5.5em \hangafter=1
{\white .}\hskip 1em $\rho_\text{isum}(\mathfrak{t})$ =
 $( \frac{7}{12},
\frac{23}{60},
\frac{47}{60} )
\oplus
( \frac{1}{4},
\frac{7}{12} )
\oplus
( \frac{1}{4} )
$,

\vskip 0.7ex
\hangindent=5.5em \hangafter=1
{\white .}\hskip 1em $\rho_\text{isum}(\mathfrak{s})$ =
$\mathrm{i}$($\sqrt{\frac{1}{5}}$,
$\sqrt{\frac{2}{5}}$,
$\sqrt{\frac{2}{5}}$;\ \ 
$-\frac{5+\sqrt{5}}{10}$,
$\frac{5-\sqrt{5}}{10}$;\ \ 
$-\frac{5+\sqrt{5}}{10}$)
 $\oplus$
($\sqrt{\frac{1}{3}}$,
$\sqrt{\frac{2}{3}}$;
$-\sqrt{\frac{1}{3}}$)
 $\oplus$
$\mathrm{i}$($1$)

Fail:
all $\theta$-eigenspaces that can contain unit
 have Tr$_{E_\theta}(C) \leq 0 $. Prop. B.5 (5) eqn. (B.29)

 \ \color{black}

 \color{blue}

\noindent 177: (dims,levels) = $(3 , 
2 , 
1;60,
12,
12
)$,
irreps = $3_{5}^{3}
\hskip -1.5pt \otimes \hskip -1.5pt
1_{4}^{1,0}
\hskip -1.5pt \otimes \hskip -1.5pt
1_{3}^{1,0}\oplus
2_{4}^{1,0}
\hskip -1.5pt \otimes \hskip -1.5pt
1_{3}^{1,0}\oplus
1_{4}^{1,6}
\hskip -1.5pt \otimes \hskip -1.5pt
1_{3}^{1,0}$,
pord$(\rho_\text{isum}(\mathfrak{t})) = 10$,

\vskip 0.7ex
\hangindent=5.5em \hangafter=1
{\white .}\hskip 1em $\rho_\text{isum}(\mathfrak{t})$ =
 $( \frac{7}{12},
\frac{11}{60},
\frac{59}{60} )
\oplus
( \frac{1}{12},
\frac{7}{12} )
\oplus
( \frac{1}{12} )
$,

\vskip 0.7ex
\hangindent=5.5em \hangafter=1
{\white .}\hskip 1em $\rho_\text{isum}(\mathfrak{s})$ =
$\mathrm{i}$($-\sqrt{\frac{1}{5}}$,
$\sqrt{\frac{2}{5}}$,
$\sqrt{\frac{2}{5}}$;\ \ 
$\frac{-5+\sqrt{5}}{10}$,
$\frac{5+\sqrt{5}}{10}$;\ \ 
$\frac{-5+\sqrt{5}}{10}$)
 $\oplus$
$\mathrm{i}$($\frac{1}{2}$,
$\sqrt{\frac{3}{4}}$;\ \ 
$-\frac{1}{2}$)
 $\oplus$
$\mathrm{i}$($-1$)

Pass. 

 \ \color{black}

 \color{blue}

\noindent 178: (dims,levels) = $(3 , 
2 , 
1;60,
12,
12
)$,
irreps = $3_{5}^{3}
\hskip -1.5pt \otimes \hskip -1.5pt
1_{4}^{1,0}
\hskip -1.5pt \otimes \hskip -1.5pt
1_{3}^{1,0}\oplus
2_{4}^{1,0}
\hskip -1.5pt \otimes \hskip -1.5pt
1_{3}^{1,0}\oplus
1_{4}^{1,0}
\hskip -1.5pt \otimes \hskip -1.5pt
1_{3}^{1,0}$,
pord$(\rho_\text{isum}(\mathfrak{t})) = 10$,

\vskip 0.7ex
\hangindent=5.5em \hangafter=1
{\white .}\hskip 1em $\rho_\text{isum}(\mathfrak{t})$ =
 $( \frac{7}{12},
\frac{11}{60},
\frac{59}{60} )
\oplus
( \frac{1}{12},
\frac{7}{12} )
\oplus
( \frac{7}{12} )
$,

\vskip 0.7ex
\hangindent=5.5em \hangafter=1
{\white .}\hskip 1em $\rho_\text{isum}(\mathfrak{s})$ =
$\mathrm{i}$($-\sqrt{\frac{1}{5}}$,
$\sqrt{\frac{2}{5}}$,
$\sqrt{\frac{2}{5}}$;\ \ 
$\frac{-5+\sqrt{5}}{10}$,
$\frac{5+\sqrt{5}}{10}$;\ \ 
$\frac{-5+\sqrt{5}}{10}$)
 $\oplus$
$\mathrm{i}$($\frac{1}{2}$,
$\sqrt{\frac{3}{4}}$;\ \ 
$-\frac{1}{2}$)
 $\oplus$
$\mathrm{i}$($1$)

Pass. 

 \ \color{black}

\noindent 179: (dims,levels) = $(3 , 
2 , 
1;60,
12,
12
)$,
irreps = $3_{5}^{3}
\hskip -1.5pt \otimes \hskip -1.5pt
1_{4}^{1,0}
\hskip -1.5pt \otimes \hskip -1.5pt
1_{3}^{1,0}\oplus
2_{3}^{1,0}
\hskip -1.5pt \otimes \hskip -1.5pt
1_{4}^{1,0}\oplus
1_{4}^{1,0}
\hskip -1.5pt \otimes \hskip -1.5pt
1_{3}^{1,0}$,
pord$(\rho_\text{isum}(\mathfrak{t})) = 15$,

\vskip 0.7ex
\hangindent=5.5em \hangafter=1
{\white .}\hskip 1em $\rho_\text{isum}(\mathfrak{t})$ =
 $( \frac{7}{12},
\frac{11}{60},
\frac{59}{60} )
\oplus
( \frac{1}{4},
\frac{7}{12} )
\oplus
( \frac{7}{12} )
$,

\vskip 0.7ex
\hangindent=5.5em \hangafter=1
{\white .}\hskip 1em $\rho_\text{isum}(\mathfrak{s})$ =
$\mathrm{i}$($-\sqrt{\frac{1}{5}}$,
$\sqrt{\frac{2}{5}}$,
$\sqrt{\frac{2}{5}}$;\ \ 
$\frac{-5+\sqrt{5}}{10}$,
$\frac{5+\sqrt{5}}{10}$;\ \ 
$\frac{-5+\sqrt{5}}{10}$)
 $\oplus$
($\sqrt{\frac{1}{3}}$,
$\sqrt{\frac{2}{3}}$;
$-\sqrt{\frac{1}{3}}$)
 $\oplus$
$\mathrm{i}$($1$)

Fail:
Tr$_I(C) = -1 <$  0 for I = [ 1/4 ]. Prop. B.4 (1) eqn. (B.18)

 \ \color{black}

\noindent 180: (dims,levels) = $(3 , 
2 , 
1;60,
12,
12
)$,
irreps = $3_{5}^{3}
\hskip -1.5pt \otimes \hskip -1.5pt
1_{4}^{1,0}
\hskip -1.5pt \otimes \hskip -1.5pt
1_{3}^{1,0}\oplus
2_{3}^{1,4}
\hskip -1.5pt \otimes \hskip -1.5pt
1_{4}^{1,0}\oplus
1_{4}^{1,0}
\hskip -1.5pt \otimes \hskip -1.5pt
1_{3}^{1,0}$,
pord$(\rho_\text{isum}(\mathfrak{t})) = 15$,

\vskip 0.7ex
\hangindent=5.5em \hangafter=1
{\white .}\hskip 1em $\rho_\text{isum}(\mathfrak{t})$ =
 $( \frac{7}{12},
\frac{11}{60},
\frac{59}{60} )
\oplus
( \frac{7}{12},
\frac{11}{12} )
\oplus
( \frac{7}{12} )
$,

\vskip 0.7ex
\hangindent=5.5em \hangafter=1
{\white .}\hskip 1em $\rho_\text{isum}(\mathfrak{s})$ =
$\mathrm{i}$($-\sqrt{\frac{1}{5}}$,
$\sqrt{\frac{2}{5}}$,
$\sqrt{\frac{2}{5}}$;\ \ 
$\frac{-5+\sqrt{5}}{10}$,
$\frac{5+\sqrt{5}}{10}$;\ \ 
$\frac{-5+\sqrt{5}}{10}$)
 $\oplus$
($\sqrt{\frac{1}{3}}$,
$\sqrt{\frac{2}{3}}$;
$-\sqrt{\frac{1}{3}}$)
 $\oplus$
$\mathrm{i}$($1$)

Fail:
Tr$_I(C) = -1 <$  0 for I = [ 11/12 ]. Prop. B.4 (1) eqn. (B.18)

 \ \color{black}

\noindent 181: (dims,levels) = $(3 , 
2 , 
1;60,
12,
12
)$,
irreps = $3_{5}^{3}
\hskip -1.5pt \otimes \hskip -1.5pt
1_{4}^{1,0}
\hskip -1.5pt \otimes \hskip -1.5pt
1_{3}^{1,0}\oplus
2_{3}^{1,4}
\hskip -1.5pt \otimes \hskip -1.5pt
1_{4}^{1,0}\oplus
1_{4}^{1,0}
\hskip -1.5pt \otimes \hskip -1.5pt
1_{3}^{1,4}$,
pord$(\rho_\text{isum}(\mathfrak{t})) = 15$,

\vskip 0.7ex
\hangindent=5.5em \hangafter=1
{\white .}\hskip 1em $\rho_\text{isum}(\mathfrak{t})$ =
 $( \frac{7}{12},
\frac{11}{60},
\frac{59}{60} )
\oplus
( \frac{7}{12},
\frac{11}{12} )
\oplus
( \frac{11}{12} )
$,

\vskip 0.7ex
\hangindent=5.5em \hangafter=1
{\white .}\hskip 1em $\rho_\text{isum}(\mathfrak{s})$ =
$\mathrm{i}$($-\sqrt{\frac{1}{5}}$,
$\sqrt{\frac{2}{5}}$,
$\sqrt{\frac{2}{5}}$;\ \ 
$\frac{-5+\sqrt{5}}{10}$,
$\frac{5+\sqrt{5}}{10}$;\ \ 
$\frac{-5+\sqrt{5}}{10}$)
 $\oplus$
($\sqrt{\frac{1}{3}}$,
$\sqrt{\frac{2}{3}}$;
$-\sqrt{\frac{1}{3}}$)
 $\oplus$
$\mathrm{i}$($1$)

Fail:
all $\theta$-eigenspaces that can contain unit
 have Tr$_{E_\theta}(C) \leq 0 $. Prop. B.5 (5) eqn. (B.29)

 \ \color{black}

 \color{blue}

\noindent 182: (dims,levels) = $(3 , 
2 , 
1;60,
12,
12
)$,
irreps = $3_{5}^{1}
\hskip -1.5pt \otimes \hskip -1.5pt
1_{4}^{1,0}
\hskip -1.5pt \otimes \hskip -1.5pt
1_{3}^{1,0}\oplus
2_{4}^{1,0}
\hskip -1.5pt \otimes \hskip -1.5pt
1_{3}^{1,0}\oplus
1_{4}^{1,6}
\hskip -1.5pt \otimes \hskip -1.5pt
1_{3}^{1,0}$,
pord$(\rho_\text{isum}(\mathfrak{t})) = 10$,

\vskip 0.7ex
\hangindent=5.5em \hangafter=1
{\white .}\hskip 1em $\rho_\text{isum}(\mathfrak{t})$ =
 $( \frac{7}{12},
\frac{23}{60},
\frac{47}{60} )
\oplus
( \frac{1}{12},
\frac{7}{12} )
\oplus
( \frac{1}{12} )
$,

\vskip 0.7ex
\hangindent=5.5em \hangafter=1
{\white .}\hskip 1em $\rho_\text{isum}(\mathfrak{s})$ =
$\mathrm{i}$($\sqrt{\frac{1}{5}}$,
$\sqrt{\frac{2}{5}}$,
$\sqrt{\frac{2}{5}}$;\ \ 
$-\frac{5+\sqrt{5}}{10}$,
$\frac{5-\sqrt{5}}{10}$;\ \ 
$-\frac{5+\sqrt{5}}{10}$)
 $\oplus$
$\mathrm{i}$($\frac{1}{2}$,
$\sqrt{\frac{3}{4}}$;\ \ 
$-\frac{1}{2}$)
 $\oplus$
$\mathrm{i}$($-1$)

Pass. 

 \ \color{black}

 \color{blue}

\noindent 183: (dims,levels) = $(3 , 
2 , 
1;60,
12,
12
)$,
irreps = $3_{5}^{1}
\hskip -1.5pt \otimes \hskip -1.5pt
1_{4}^{1,0}
\hskip -1.5pt \otimes \hskip -1.5pt
1_{3}^{1,0}\oplus
2_{4}^{1,0}
\hskip -1.5pt \otimes \hskip -1.5pt
1_{3}^{1,0}\oplus
1_{4}^{1,0}
\hskip -1.5pt \otimes \hskip -1.5pt
1_{3}^{1,0}$,
pord$(\rho_\text{isum}(\mathfrak{t})) = 10$,

\vskip 0.7ex
\hangindent=5.5em \hangafter=1
{\white .}\hskip 1em $\rho_\text{isum}(\mathfrak{t})$ =
 $( \frac{7}{12},
\frac{23}{60},
\frac{47}{60} )
\oplus
( \frac{1}{12},
\frac{7}{12} )
\oplus
( \frac{7}{12} )
$,

\vskip 0.7ex
\hangindent=5.5em \hangafter=1
{\white .}\hskip 1em $\rho_\text{isum}(\mathfrak{s})$ =
$\mathrm{i}$($\sqrt{\frac{1}{5}}$,
$\sqrt{\frac{2}{5}}$,
$\sqrt{\frac{2}{5}}$;\ \ 
$-\frac{5+\sqrt{5}}{10}$,
$\frac{5-\sqrt{5}}{10}$;\ \ 
$-\frac{5+\sqrt{5}}{10}$)
 $\oplus$
$\mathrm{i}$($\frac{1}{2}$,
$\sqrt{\frac{3}{4}}$;\ \ 
$-\frac{1}{2}$)
 $\oplus$
$\mathrm{i}$($1$)

Pass. 

 \ \color{black}

\noindent 184: (dims,levels) = $(3 , 
2 , 
1;60,
12,
12
)$,
irreps = $3_{5}^{1}
\hskip -1.5pt \otimes \hskip -1.5pt
1_{4}^{1,0}
\hskip -1.5pt \otimes \hskip -1.5pt
1_{3}^{1,0}\oplus
2_{3}^{1,0}
\hskip -1.5pt \otimes \hskip -1.5pt
1_{4}^{1,0}\oplus
1_{4}^{1,0}
\hskip -1.5pt \otimes \hskip -1.5pt
1_{3}^{1,0}$,
pord$(\rho_\text{isum}(\mathfrak{t})) = 15$,

\vskip 0.7ex
\hangindent=5.5em \hangafter=1
{\white .}\hskip 1em $\rho_\text{isum}(\mathfrak{t})$ =
 $( \frac{7}{12},
\frac{23}{60},
\frac{47}{60} )
\oplus
( \frac{1}{4},
\frac{7}{12} )
\oplus
( \frac{7}{12} )
$,

\vskip 0.7ex
\hangindent=5.5em \hangafter=1
{\white .}\hskip 1em $\rho_\text{isum}(\mathfrak{s})$ =
$\mathrm{i}$($\sqrt{\frac{1}{5}}$,
$\sqrt{\frac{2}{5}}$,
$\sqrt{\frac{2}{5}}$;\ \ 
$-\frac{5+\sqrt{5}}{10}$,
$\frac{5-\sqrt{5}}{10}$;\ \ 
$-\frac{5+\sqrt{5}}{10}$)
 $\oplus$
($\sqrt{\frac{1}{3}}$,
$\sqrt{\frac{2}{3}}$;
$-\sqrt{\frac{1}{3}}$)
 $\oplus$
$\mathrm{i}$($1$)

Fail:
Tr$_I(C) = -1 <$  0 for I = [ 1/4 ]. Prop. B.4 (1) eqn. (B.18)

 \ \color{black}

\noindent 185: (dims,levels) = $(3 , 
2 , 
1;60,
12,
12
)$,
irreps = $3_{5}^{1}
\hskip -1.5pt \otimes \hskip -1.5pt
1_{4}^{1,0}
\hskip -1.5pt \otimes \hskip -1.5pt
1_{3}^{1,0}\oplus
2_{3}^{1,4}
\hskip -1.5pt \otimes \hskip -1.5pt
1_{4}^{1,0}\oplus
1_{4}^{1,0}
\hskip -1.5pt \otimes \hskip -1.5pt
1_{3}^{1,0}$,
pord$(\rho_\text{isum}(\mathfrak{t})) = 15$,

\vskip 0.7ex
\hangindent=5.5em \hangafter=1
{\white .}\hskip 1em $\rho_\text{isum}(\mathfrak{t})$ =
 $( \frac{7}{12},
\frac{23}{60},
\frac{47}{60} )
\oplus
( \frac{7}{12},
\frac{11}{12} )
\oplus
( \frac{7}{12} )
$,

\vskip 0.7ex
\hangindent=5.5em \hangafter=1
{\white .}\hskip 1em $\rho_\text{isum}(\mathfrak{s})$ =
$\mathrm{i}$($\sqrt{\frac{1}{5}}$,
$\sqrt{\frac{2}{5}}$,
$\sqrt{\frac{2}{5}}$;\ \ 
$-\frac{5+\sqrt{5}}{10}$,
$\frac{5-\sqrt{5}}{10}$;\ \ 
$-\frac{5+\sqrt{5}}{10}$)
 $\oplus$
($\sqrt{\frac{1}{3}}$,
$\sqrt{\frac{2}{3}}$;
$-\sqrt{\frac{1}{3}}$)
 $\oplus$
$\mathrm{i}$($1$)

Fail:
Tr$_I(C) = -1 <$  0 for I = [ 11/12 ]. Prop. B.4 (1) eqn. (B.18)

 \ \color{black}

\noindent 186: (dims,levels) = $(3 , 
2 , 
1;60,
12,
12
)$,
irreps = $3_{5}^{1}
\hskip -1.5pt \otimes \hskip -1.5pt
1_{4}^{1,0}
\hskip -1.5pt \otimes \hskip -1.5pt
1_{3}^{1,0}\oplus
2_{3}^{1,4}
\hskip -1.5pt \otimes \hskip -1.5pt
1_{4}^{1,0}\oplus
1_{4}^{1,0}
\hskip -1.5pt \otimes \hskip -1.5pt
1_{3}^{1,4}$,
pord$(\rho_\text{isum}(\mathfrak{t})) = 15$,

\vskip 0.7ex
\hangindent=5.5em \hangafter=1
{\white .}\hskip 1em $\rho_\text{isum}(\mathfrak{t})$ =
 $( \frac{7}{12},
\frac{23}{60},
\frac{47}{60} )
\oplus
( \frac{7}{12},
\frac{11}{12} )
\oplus
( \frac{11}{12} )
$,

\vskip 0.7ex
\hangindent=5.5em \hangafter=1
{\white .}\hskip 1em $\rho_\text{isum}(\mathfrak{s})$ =
$\mathrm{i}$($\sqrt{\frac{1}{5}}$,
$\sqrt{\frac{2}{5}}$,
$\sqrt{\frac{2}{5}}$;\ \ 
$-\frac{5+\sqrt{5}}{10}$,
$\frac{5-\sqrt{5}}{10}$;\ \ 
$-\frac{5+\sqrt{5}}{10}$)
 $\oplus$
($\sqrt{\frac{1}{3}}$,
$\sqrt{\frac{2}{3}}$;
$-\sqrt{\frac{1}{3}}$)
 $\oplus$
$\mathrm{i}$($1$)

Fail:
all $\theta$-eigenspaces that can contain unit
 have Tr$_{E_\theta}(C) \leq 0 $. Prop. B.5 (5) eqn. (B.29)

 \ \color{black}

 \color{blue}

\noindent 187: (dims,levels) = $(3 , 
2 , 
1;60,
60,
12
)$,
irreps = $3_{5}^{3}
\hskip -1.5pt \otimes \hskip -1.5pt
1_{4}^{1,0}
\hskip -1.5pt \otimes \hskip -1.5pt
1_{3}^{1,0}\oplus
2_{5}^{2}
\hskip -1.5pt \otimes \hskip -1.5pt
1_{4}^{1,0}
\hskip -1.5pt \otimes \hskip -1.5pt
1_{3}^{1,0}\oplus
1_{4}^{1,0}
\hskip -1.5pt \otimes \hskip -1.5pt
1_{3}^{1,0}$,
pord$(\rho_\text{isum}(\mathfrak{t})) = 5$,

\vskip 0.7ex
\hangindent=5.5em \hangafter=1
{\white .}\hskip 1em $\rho_\text{isum}(\mathfrak{t})$ =
 $( \frac{7}{12},
\frac{11}{60},
\frac{59}{60} )
\oplus
( \frac{11}{60},
\frac{59}{60} )
\oplus
( \frac{7}{12} )
$,

\vskip 0.7ex
\hangindent=5.5em \hangafter=1
{\white .}\hskip 1em $\rho_\text{isum}(\mathfrak{s})$ =
$\mathrm{i}$($-\sqrt{\frac{1}{5}}$,
$\sqrt{\frac{2}{5}}$,
$\sqrt{\frac{2}{5}}$;\ \ 
$\frac{-5+\sqrt{5}}{10}$,
$\frac{5+\sqrt{5}}{10}$;\ \ 
$\frac{-5+\sqrt{5}}{10}$)
 $\oplus$
($-\frac{1}{\sqrt{5}}c^{1}_{20}
$,
$\frac{1}{\sqrt{5}}c^{3}_{20}
$;
$\frac{1}{\sqrt{5}}c^{1}_{20}
$)
 $\oplus$
$\mathrm{i}$($1$)

Pass. 

 \ \color{black}

 \color{blue}

\noindent 188: (dims,levels) = $(3 , 
2 , 
1;60,
60,
12
)$,
irreps = $3_{5}^{1}
\hskip -1.5pt \otimes \hskip -1.5pt
1_{4}^{1,0}
\hskip -1.5pt \otimes \hskip -1.5pt
1_{3}^{1,0}\oplus
2_{5}^{1}
\hskip -1.5pt \otimes \hskip -1.5pt
1_{4}^{1,0}
\hskip -1.5pt \otimes \hskip -1.5pt
1_{3}^{1,0}\oplus
1_{4}^{1,0}
\hskip -1.5pt \otimes \hskip -1.5pt
1_{3}^{1,0}$,
pord$(\rho_\text{isum}(\mathfrak{t})) = 5$,

\vskip 0.7ex
\hangindent=5.5em \hangafter=1
{\white .}\hskip 1em $\rho_\text{isum}(\mathfrak{t})$ =
 $( \frac{7}{12},
\frac{23}{60},
\frac{47}{60} )
\oplus
( \frac{23}{60},
\frac{47}{60} )
\oplus
( \frac{7}{12} )
$,

\vskip 0.7ex
\hangindent=5.5em \hangafter=1
{\white .}\hskip 1em $\rho_\text{isum}(\mathfrak{s})$ =
$\mathrm{i}$($\sqrt{\frac{1}{5}}$,
$\sqrt{\frac{2}{5}}$,
$\sqrt{\frac{2}{5}}$;\ \ 
$-\frac{5+\sqrt{5}}{10}$,
$\frac{5-\sqrt{5}}{10}$;\ \ 
$-\frac{5+\sqrt{5}}{10}$)
 $\oplus$
($-\frac{1}{\sqrt{5}}c^{3}_{20}
$,
$\frac{1}{\sqrt{5}}c^{1}_{20}
$;
$\frac{1}{\sqrt{5}}c^{3}_{20}
$)
 $\oplus$
$\mathrm{i}$($1$)

Pass. 

 \ \color{black}

\noindent 189: (dims,levels) = $(3 , 
3;3,
4
)$,
irreps = $3_{3}^{1,0}\oplus
3_{4}^{1,0}$,
pord$(\rho_\text{isum}(\mathfrak{t})) = 12$,

\vskip 0.7ex
\hangindent=5.5em \hangafter=1
{\white .}\hskip 1em $\rho_\text{isum}(\mathfrak{t})$ =
 $( 0,
\frac{1}{3},
\frac{2}{3} )
\oplus
( 0,
\frac{1}{4},
\frac{3}{4} )
$,

\vskip 0.7ex
\hangindent=5.5em \hangafter=1
{\white .}\hskip 1em $\rho_\text{isum}(\mathfrak{s})$ =
($-\frac{1}{3}$,
$\frac{2}{3}$,
$\frac{2}{3}$;
$-\frac{1}{3}$,
$\frac{2}{3}$;
$-\frac{1}{3}$)
 $\oplus$
($0$,
$\sqrt{\frac{1}{2}}$,
$\sqrt{\frac{1}{2}}$;
$-\frac{1}{2}$,
$\frac{1}{2}$;
$-\frac{1}{2}$)

Fail:
Integral: $D_{\rho}(\sigma)_{\theta} \propto $ id,
 for all $\sigma$ and all $\theta$-eigenspaces that can contain unit. Prop. B.5 (6)

 \ \color{black}

\noindent 190: (dims,levels) = $(3 , 
3;3,
4
)$,
irreps = $3_{3}^{1,0}\oplus
3_{4}^{1,3}$,
pord$(\rho_\text{isum}(\mathfrak{t})) = 12$,

\vskip 0.7ex
\hangindent=5.5em \hangafter=1
{\white .}\hskip 1em $\rho_\text{isum}(\mathfrak{t})$ =
 $( 0,
\frac{1}{3},
\frac{2}{3} )
\oplus
( 0,
\frac{1}{2},
\frac{1}{4} )
$,

\vskip 0.7ex
\hangindent=5.5em \hangafter=1
{\white .}\hskip 1em $\rho_\text{isum}(\mathfrak{s})$ =
($-\frac{1}{3}$,
$\frac{2}{3}$,
$\frac{2}{3}$;
$-\frac{1}{3}$,
$\frac{2}{3}$;
$-\frac{1}{3}$)
 $\oplus$
$\mathrm{i}$($-\frac{1}{2}$,
$\frac{1}{2}$,
$\sqrt{\frac{1}{2}}$;\ \ 
$-\frac{1}{2}$,
$\sqrt{\frac{1}{2}}$;\ \ 
$0$)

Fail:
number of self dual objects $|$Tr($\rho(\mathfrak s^2)$)$|$ = 0. Prop. B.4 (1)\
 eqn. (B.16)

 \ \color{black}

\noindent 191: (dims,levels) = $(3 , 
3;3,
4
)$,
irreps = $3_{3}^{1,0}\oplus
3_{4}^{1,9}$,
pord$(\rho_\text{isum}(\mathfrak{t})) = 12$,

\vskip 0.7ex
\hangindent=5.5em \hangafter=1
{\white .}\hskip 1em $\rho_\text{isum}(\mathfrak{t})$ =
 $( 0,
\frac{1}{3},
\frac{2}{3} )
\oplus
( 0,
\frac{1}{2},
\frac{3}{4} )
$,

\vskip 0.7ex
\hangindent=5.5em \hangafter=1
{\white .}\hskip 1em $\rho_\text{isum}(\mathfrak{s})$ =
($-\frac{1}{3}$,
$\frac{2}{3}$,
$\frac{2}{3}$;
$-\frac{1}{3}$,
$\frac{2}{3}$;
$-\frac{1}{3}$)
 $\oplus$
$\mathrm{i}$($\frac{1}{2}$,
$\frac{1}{2}$,
$\sqrt{\frac{1}{2}}$;\ \ 
$\frac{1}{2}$,
$-\sqrt{\frac{1}{2}}$;\ \ 
$0$)

Fail:
number of self dual objects $|$Tr($\rho(\mathfrak s^2)$)$|$ = 0. Prop. B.4 (1)\
 eqn. (B.16)

 \ \color{black}

 \color{blue}

\noindent 192: (dims,levels) = $(3 , 
3;3,
5
)$,
irreps = $3_{3}^{1,0}\oplus
3_{5}^{1}$,
pord$(\rho_\text{isum}(\mathfrak{t})) = 15$,

\vskip 0.7ex
\hangindent=5.5em \hangafter=1
{\white .}\hskip 1em $\rho_\text{isum}(\mathfrak{t})$ =
 $( 0,
\frac{1}{3},
\frac{2}{3} )
\oplus
( 0,
\frac{1}{5},
\frac{4}{5} )
$,

\vskip 0.7ex
\hangindent=5.5em \hangafter=1
{\white .}\hskip 1em $\rho_\text{isum}(\mathfrak{s})$ =
($-\frac{1}{3}$,
$\frac{2}{3}$,
$\frac{2}{3}$;
$-\frac{1}{3}$,
$\frac{2}{3}$;
$-\frac{1}{3}$)
 $\oplus$
($\sqrt{\frac{1}{5}}$,
$-\sqrt{\frac{2}{5}}$,
$-\sqrt{\frac{2}{5}}$;
$-\frac{5+\sqrt{5}}{10}$,
$\frac{5-\sqrt{5}}{10}$;
$-\frac{5+\sqrt{5}}{10}$)

Pass. 

 \ \color{black}

 \color{blue}

\noindent 193: (dims,levels) = $(3 , 
3;3,
5
)$,
irreps = $3_{3}^{1,0}\oplus
3_{5}^{3}$,
pord$(\rho_\text{isum}(\mathfrak{t})) = 15$,

\vskip 0.7ex
\hangindent=5.5em \hangafter=1
{\white .}\hskip 1em $\rho_\text{isum}(\mathfrak{t})$ =
 $( 0,
\frac{1}{3},
\frac{2}{3} )
\oplus
( 0,
\frac{2}{5},
\frac{3}{5} )
$,

\vskip 0.7ex
\hangindent=5.5em \hangafter=1
{\white .}\hskip 1em $\rho_\text{isum}(\mathfrak{s})$ =
($-\frac{1}{3}$,
$\frac{2}{3}$,
$\frac{2}{3}$;
$-\frac{1}{3}$,
$\frac{2}{3}$;
$-\frac{1}{3}$)
 $\oplus$
($-\sqrt{\frac{1}{5}}$,
$-\sqrt{\frac{2}{5}}$,
$-\sqrt{\frac{2}{5}}$;
$\frac{-5+\sqrt{5}}{10}$,
$\frac{5+\sqrt{5}}{10}$;
$\frac{-5+\sqrt{5}}{10}$)

Pass. 

 \ \color{black}

\noindent 194: (dims,levels) = $(3 , 
3;3,
8
)$,
irreps = $3_{3}^{1,0}\oplus
3_{8}^{1,0}$,
pord$(\rho_\text{isum}(\mathfrak{t})) = 24$,

\vskip 0.7ex
\hangindent=5.5em \hangafter=1
{\white .}\hskip 1em $\rho_\text{isum}(\mathfrak{t})$ =
 $( 0,
\frac{1}{3},
\frac{2}{3} )
\oplus
( 0,
\frac{1}{8},
\frac{5}{8} )
$,

\vskip 0.7ex
\hangindent=5.5em \hangafter=1
{\white .}\hskip 1em $\rho_\text{isum}(\mathfrak{s})$ =
($-\frac{1}{3}$,
$\frac{2}{3}$,
$\frac{2}{3}$;
$-\frac{1}{3}$,
$\frac{2}{3}$;
$-\frac{1}{3}$)
 $\oplus$
$\mathrm{i}$($0$,
$\sqrt{\frac{1}{2}}$,
$\sqrt{\frac{1}{2}}$;\ \ 
$-\frac{1}{2}$,
$\frac{1}{2}$;\ \ 
$-\frac{1}{2}$)

Fail:
number of self dual objects $|$Tr($\rho(\mathfrak s^2)$)$|$ = 0. Prop. B.4 (1)\
 eqn. (B.16)

 \ \color{black}

\noindent 195: (dims,levels) = $(3 , 
3;3,
8
)$,
irreps = $3_{3}^{1,0}\oplus
3_{8}^{3,0}$,
pord$(\rho_\text{isum}(\mathfrak{t})) = 24$,

\vskip 0.7ex
\hangindent=5.5em \hangafter=1
{\white .}\hskip 1em $\rho_\text{isum}(\mathfrak{t})$ =
 $( 0,
\frac{1}{3},
\frac{2}{3} )
\oplus
( 0,
\frac{3}{8},
\frac{7}{8} )
$,

\vskip 0.7ex
\hangindent=5.5em \hangafter=1
{\white .}\hskip 1em $\rho_\text{isum}(\mathfrak{s})$ =
($-\frac{1}{3}$,
$\frac{2}{3}$,
$\frac{2}{3}$;
$-\frac{1}{3}$,
$\frac{2}{3}$;
$-\frac{1}{3}$)
 $\oplus$
$\mathrm{i}$($0$,
$\sqrt{\frac{1}{2}}$,
$\sqrt{\frac{1}{2}}$;\ \ 
$\frac{1}{2}$,
$-\frac{1}{2}$;\ \ 
$\frac{1}{2}$)

Fail:
number of self dual objects $|$Tr($\rho(\mathfrak s^2)$)$|$ = 0. Prop. B.4 (1)\
 eqn. (B.16)

 \ \color{black}

\noindent 196: (dims,levels) = $(3 , 
3;3,
12
)$,
irreps = $3_{3}^{1,0}\oplus
3_{4}^{1,0}
\hskip -1.5pt \otimes \hskip -1.5pt
1_{3}^{1,0}$,
pord$(\rho_\text{isum}(\mathfrak{t})) = 12$,

\vskip 0.7ex
\hangindent=5.5em \hangafter=1
{\white .}\hskip 1em $\rho_\text{isum}(\mathfrak{t})$ =
 $( 0,
\frac{1}{3},
\frac{2}{3} )
\oplus
( \frac{1}{3},
\frac{1}{12},
\frac{7}{12} )
$,

\vskip 0.7ex
\hangindent=5.5em \hangafter=1
{\white .}\hskip 1em $\rho_\text{isum}(\mathfrak{s})$ =
($-\frac{1}{3}$,
$\frac{2}{3}$,
$\frac{2}{3}$;
$-\frac{1}{3}$,
$\frac{2}{3}$;
$-\frac{1}{3}$)
 $\oplus$
($0$,
$\sqrt{\frac{1}{2}}$,
$\sqrt{\frac{1}{2}}$;
$-\frac{1}{2}$,
$\frac{1}{2}$;
$-\frac{1}{2}$)

Fail:
Integral: $D_{\rho}(\sigma)_{\theta} \propto $ id,
 for all $\sigma$ and all $\theta$-eigenspaces that can contain unit. Prop. B.5 (6)

 \ \color{black}

\noindent 197: (dims,levels) = $(3 , 
3;3,
12
)$,
irreps = $3_{3}^{1,0}\oplus
3_{4}^{1,9}
\hskip -1.5pt \otimes \hskip -1.5pt
1_{3}^{1,0}$,
pord$(\rho_\text{isum}(\mathfrak{t})) = 12$,

\vskip 0.7ex
\hangindent=5.5em \hangafter=1
{\white .}\hskip 1em $\rho_\text{isum}(\mathfrak{t})$ =
 $( 0,
\frac{1}{3},
\frac{2}{3} )
\oplus
( \frac{1}{3},
\frac{5}{6},
\frac{1}{12} )
$,

\vskip 0.7ex
\hangindent=5.5em \hangafter=1
{\white .}\hskip 1em $\rho_\text{isum}(\mathfrak{s})$ =
($-\frac{1}{3}$,
$\frac{2}{3}$,
$\frac{2}{3}$;
$-\frac{1}{3}$,
$\frac{2}{3}$;
$-\frac{1}{3}$)
 $\oplus$
$\mathrm{i}$($\frac{1}{2}$,
$\frac{1}{2}$,
$\sqrt{\frac{1}{2}}$;\ \ 
$\frac{1}{2}$,
$-\sqrt{\frac{1}{2}}$;\ \ 
$0$)

Fail:
number of self dual objects $|$Tr($\rho(\mathfrak s^2)$)$|$ = 0. Prop. B.4 (1)\
 eqn. (B.16)

 \ \color{black}

\noindent 198: (dims,levels) = $(3 , 
3;3,
12
)$,
irreps = $3_{3}^{1,0}\oplus
3_{4}^{1,3}
\hskip -1.5pt \otimes \hskip -1.5pt
1_{3}^{1,0}$,
pord$(\rho_\text{isum}(\mathfrak{t})) = 12$,

\vskip 0.7ex
\hangindent=5.5em \hangafter=1
{\white .}\hskip 1em $\rho_\text{isum}(\mathfrak{t})$ =
 $( 0,
\frac{1}{3},
\frac{2}{3} )
\oplus
( \frac{1}{3},
\frac{5}{6},
\frac{7}{12} )
$,

\vskip 0.7ex
\hangindent=5.5em \hangafter=1
{\white .}\hskip 1em $\rho_\text{isum}(\mathfrak{s})$ =
($-\frac{1}{3}$,
$\frac{2}{3}$,
$\frac{2}{3}$;
$-\frac{1}{3}$,
$\frac{2}{3}$;
$-\frac{1}{3}$)
 $\oplus$
$\mathrm{i}$($-\frac{1}{2}$,
$\frac{1}{2}$,
$\sqrt{\frac{1}{2}}$;\ \ 
$-\frac{1}{2}$,
$\sqrt{\frac{1}{2}}$;\ \ 
$0$)

Fail:
number of self dual objects $|$Tr($\rho(\mathfrak s^2)$)$|$ = 0. Prop. B.4 (1)\
 eqn. (B.16)

 \ \color{black}

\noindent 199: (dims,levels) = $(3 , 
3;3,
12
)$,
irreps = $3_{3}^{1,0}\oplus
3_{4}^{1,9}
\hskip -1.5pt \otimes \hskip -1.5pt
1_{3}^{1,4}$,
pord$(\rho_\text{isum}(\mathfrak{t})) = 12$,

\vskip 0.7ex
\hangindent=5.5em \hangafter=1
{\white .}\hskip 1em $\rho_\text{isum}(\mathfrak{t})$ =
 $( 0,
\frac{1}{3},
\frac{2}{3} )
\oplus
( \frac{2}{3},
\frac{1}{6},
\frac{5}{12} )
$,

\vskip 0.7ex
\hangindent=5.5em \hangafter=1
{\white .}\hskip 1em $\rho_\text{isum}(\mathfrak{s})$ =
($-\frac{1}{3}$,
$\frac{2}{3}$,
$\frac{2}{3}$;
$-\frac{1}{3}$,
$\frac{2}{3}$;
$-\frac{1}{3}$)
 $\oplus$
$\mathrm{i}$($\frac{1}{2}$,
$\frac{1}{2}$,
$\sqrt{\frac{1}{2}}$;\ \ 
$\frac{1}{2}$,
$-\sqrt{\frac{1}{2}}$;\ \ 
$0$)

Fail:
number of self dual objects $|$Tr($\rho(\mathfrak s^2)$)$|$ = 0. Prop. B.4 (1)\
 eqn. (B.16)

 \ \color{black}

\noindent 200: (dims,levels) = $(3 , 
3;3,
12
)$,
irreps = $3_{3}^{1,0}\oplus
3_{4}^{1,3}
\hskip -1.5pt \otimes \hskip -1.5pt
1_{3}^{1,4}$,
pord$(\rho_\text{isum}(\mathfrak{t})) = 12$,

\vskip 0.7ex
\hangindent=5.5em \hangafter=1
{\white .}\hskip 1em $\rho_\text{isum}(\mathfrak{t})$ =
 $( 0,
\frac{1}{3},
\frac{2}{3} )
\oplus
( \frac{2}{3},
\frac{1}{6},
\frac{11}{12} )
$,

\vskip 0.7ex
\hangindent=5.5em \hangafter=1
{\white .}\hskip 1em $\rho_\text{isum}(\mathfrak{s})$ =
($-\frac{1}{3}$,
$\frac{2}{3}$,
$\frac{2}{3}$;
$-\frac{1}{3}$,
$\frac{2}{3}$;
$-\frac{1}{3}$)
 $\oplus$
$\mathrm{i}$($-\frac{1}{2}$,
$\frac{1}{2}$,
$\sqrt{\frac{1}{2}}$;\ \ 
$-\frac{1}{2}$,
$\sqrt{\frac{1}{2}}$;\ \ 
$0$)

Fail:
number of self dual objects $|$Tr($\rho(\mathfrak s^2)$)$|$ = 0. Prop. B.4 (1)\
 eqn. (B.16)

 \ \color{black}

\noindent 201: (dims,levels) = $(3 , 
3;3,
12
)$,
irreps = $3_{3}^{1,0}\oplus
3_{4}^{1,0}
\hskip -1.5pt \otimes \hskip -1.5pt
1_{3}^{1,4}$,
pord$(\rho_\text{isum}(\mathfrak{t})) = 12$,

\vskip 0.7ex
\hangindent=5.5em \hangafter=1
{\white .}\hskip 1em $\rho_\text{isum}(\mathfrak{t})$ =
 $( 0,
\frac{1}{3},
\frac{2}{3} )
\oplus
( \frac{2}{3},
\frac{5}{12},
\frac{11}{12} )
$,

\vskip 0.7ex
\hangindent=5.5em \hangafter=1
{\white .}\hskip 1em $\rho_\text{isum}(\mathfrak{s})$ =
($-\frac{1}{3}$,
$\frac{2}{3}$,
$\frac{2}{3}$;
$-\frac{1}{3}$,
$\frac{2}{3}$;
$-\frac{1}{3}$)
 $\oplus$
($0$,
$\sqrt{\frac{1}{2}}$,
$\sqrt{\frac{1}{2}}$;
$-\frac{1}{2}$,
$\frac{1}{2}$;
$-\frac{1}{2}$)

Fail:
Integral: $D_{\rho}(\sigma)_{\theta} \propto $ id,
 for all $\sigma$ and all $\theta$-eigenspaces that can contain unit. Prop. B.5 (6)

 \ \color{black}

 \color{blue}

\noindent 202: (dims,levels) = $(3 , 
3;3,
15
)$,
irreps = $3_{3}^{1,0}\oplus
3_{5}^{1}
\hskip -1.5pt \otimes \hskip -1.5pt
1_{3}^{1,0}$,
pord$(\rho_\text{isum}(\mathfrak{t})) = 15$,

\vskip 0.7ex
\hangindent=5.5em \hangafter=1
{\white .}\hskip 1em $\rho_\text{isum}(\mathfrak{t})$ =
 $( 0,
\frac{1}{3},
\frac{2}{3} )
\oplus
( \frac{1}{3},
\frac{2}{15},
\frac{8}{15} )
$,

\vskip 0.7ex
\hangindent=5.5em \hangafter=1
{\white .}\hskip 1em $\rho_\text{isum}(\mathfrak{s})$ =
($-\frac{1}{3}$,
$\frac{2}{3}$,
$\frac{2}{3}$;
$-\frac{1}{3}$,
$\frac{2}{3}$;
$-\frac{1}{3}$)
 $\oplus$
($\sqrt{\frac{1}{5}}$,
$-\sqrt{\frac{2}{5}}$,
$-\sqrt{\frac{2}{5}}$;
$-\frac{5+\sqrt{5}}{10}$,
$\frac{5-\sqrt{5}}{10}$;
$-\frac{5+\sqrt{5}}{10}$)

Pass. 

 \ \color{black}

 \color{blue}

\noindent 203: (dims,levels) = $(3 , 
3;3,
15
)$,
irreps = $3_{3}^{1,0}\oplus
3_{5}^{3}
\hskip -1.5pt \otimes \hskip -1.5pt
1_{3}^{1,0}$,
pord$(\rho_\text{isum}(\mathfrak{t})) = 15$,

\vskip 0.7ex
\hangindent=5.5em \hangafter=1
{\white .}\hskip 1em $\rho_\text{isum}(\mathfrak{t})$ =
 $( 0,
\frac{1}{3},
\frac{2}{3} )
\oplus
( \frac{1}{3},
\frac{11}{15},
\frac{14}{15} )
$,

\vskip 0.7ex
\hangindent=5.5em \hangafter=1
{\white .}\hskip 1em $\rho_\text{isum}(\mathfrak{s})$ =
($-\frac{1}{3}$,
$\frac{2}{3}$,
$\frac{2}{3}$;
$-\frac{1}{3}$,
$\frac{2}{3}$;
$-\frac{1}{3}$)
 $\oplus$
($-\sqrt{\frac{1}{5}}$,
$-\sqrt{\frac{2}{5}}$,
$-\sqrt{\frac{2}{5}}$;
$\frac{-5+\sqrt{5}}{10}$,
$\frac{5+\sqrt{5}}{10}$;
$\frac{-5+\sqrt{5}}{10}$)

Pass. 

 \ \color{black}

 \color{blue}

\noindent 204: (dims,levels) = $(3 , 
3;3,
15
)$,
irreps = $3_{3}^{1,0}\oplus
3_{5}^{3}
\hskip -1.5pt \otimes \hskip -1.5pt
1_{3}^{1,4}$,
pord$(\rho_\text{isum}(\mathfrak{t})) = 15$,

\vskip 0.7ex
\hangindent=5.5em \hangafter=1
{\white .}\hskip 1em $\rho_\text{isum}(\mathfrak{t})$ =
 $( 0,
\frac{1}{3},
\frac{2}{3} )
\oplus
( \frac{2}{3},
\frac{1}{15},
\frac{4}{15} )
$,

\vskip 0.7ex
\hangindent=5.5em \hangafter=1
{\white .}\hskip 1em $\rho_\text{isum}(\mathfrak{s})$ =
($-\frac{1}{3}$,
$\frac{2}{3}$,
$\frac{2}{3}$;
$-\frac{1}{3}$,
$\frac{2}{3}$;
$-\frac{1}{3}$)
 $\oplus$
($-\sqrt{\frac{1}{5}}$,
$-\sqrt{\frac{2}{5}}$,
$-\sqrt{\frac{2}{5}}$;
$\frac{-5+\sqrt{5}}{10}$,
$\frac{5+\sqrt{5}}{10}$;
$\frac{-5+\sqrt{5}}{10}$)

Pass. 

 \ \color{black}

 \color{blue}

\noindent 205: (dims,levels) = $(3 , 
3;3,
15
)$,
irreps = $3_{3}^{1,0}\oplus
3_{5}^{1}
\hskip -1.5pt \otimes \hskip -1.5pt
1_{3}^{1,4}$,
pord$(\rho_\text{isum}(\mathfrak{t})) = 15$,

\vskip 0.7ex
\hangindent=5.5em \hangafter=1
{\white .}\hskip 1em $\rho_\text{isum}(\mathfrak{t})$ =
 $( 0,
\frac{1}{3},
\frac{2}{3} )
\oplus
( \frac{2}{3},
\frac{7}{15},
\frac{13}{15} )
$,

\vskip 0.7ex
\hangindent=5.5em \hangafter=1
{\white .}\hskip 1em $\rho_\text{isum}(\mathfrak{s})$ =
($-\frac{1}{3}$,
$\frac{2}{3}$,
$\frac{2}{3}$;
$-\frac{1}{3}$,
$\frac{2}{3}$;
$-\frac{1}{3}$)
 $\oplus$
($\sqrt{\frac{1}{5}}$,
$-\sqrt{\frac{2}{5}}$,
$-\sqrt{\frac{2}{5}}$;
$-\frac{5+\sqrt{5}}{10}$,
$\frac{5-\sqrt{5}}{10}$;
$-\frac{5+\sqrt{5}}{10}$)

Pass. 

 \ \color{black}

\noindent 206: (dims,levels) = $(3 , 
3;3,
24
)$,
irreps = $3_{3}^{1,0}\oplus
3_{8}^{3,0}
\hskip -1.5pt \otimes \hskip -1.5pt
1_{3}^{1,0}$,
pord$(\rho_\text{isum}(\mathfrak{t})) = 24$,

\vskip 0.7ex
\hangindent=5.5em \hangafter=1
{\white .}\hskip 1em $\rho_\text{isum}(\mathfrak{t})$ =
 $( 0,
\frac{1}{3},
\frac{2}{3} )
\oplus
( \frac{1}{3},
\frac{5}{24},
\frac{17}{24} )
$,

\vskip 0.7ex
\hangindent=5.5em \hangafter=1
{\white .}\hskip 1em $\rho_\text{isum}(\mathfrak{s})$ =
($-\frac{1}{3}$,
$\frac{2}{3}$,
$\frac{2}{3}$;
$-\frac{1}{3}$,
$\frac{2}{3}$;
$-\frac{1}{3}$)
 $\oplus$
$\mathrm{i}$($0$,
$\sqrt{\frac{1}{2}}$,
$\sqrt{\frac{1}{2}}$;\ \ 
$\frac{1}{2}$,
$-\frac{1}{2}$;\ \ 
$\frac{1}{2}$)

Fail:
number of self dual objects $|$Tr($\rho(\mathfrak s^2)$)$|$ = 0. Prop. B.4 (1)\
 eqn. (B.16)

 \ \color{black}

\noindent 207: (dims,levels) = $(3 , 
3;3,
24
)$,
irreps = $3_{3}^{1,0}\oplus
3_{8}^{1,0}
\hskip -1.5pt \otimes \hskip -1.5pt
1_{3}^{1,0}$,
pord$(\rho_\text{isum}(\mathfrak{t})) = 24$,

\vskip 0.7ex
\hangindent=5.5em \hangafter=1
{\white .}\hskip 1em $\rho_\text{isum}(\mathfrak{t})$ =
 $( 0,
\frac{1}{3},
\frac{2}{3} )
\oplus
( \frac{1}{3},
\frac{11}{24},
\frac{23}{24} )
$,

\vskip 0.7ex
\hangindent=5.5em \hangafter=1
{\white .}\hskip 1em $\rho_\text{isum}(\mathfrak{s})$ =
($-\frac{1}{3}$,
$\frac{2}{3}$,
$\frac{2}{3}$;
$-\frac{1}{3}$,
$\frac{2}{3}$;
$-\frac{1}{3}$)
 $\oplus$
$\mathrm{i}$($0$,
$\sqrt{\frac{1}{2}}$,
$\sqrt{\frac{1}{2}}$;\ \ 
$-\frac{1}{2}$,
$\frac{1}{2}$;\ \ 
$-\frac{1}{2}$)

Fail:
number of self dual objects $|$Tr($\rho(\mathfrak s^2)$)$|$ = 0. Prop. B.4 (1)\
 eqn. (B.16)

 \ \color{black}

\noindent 208: (dims,levels) = $(3 , 
3;3,
24
)$,
irreps = $3_{3}^{1,0}\oplus
3_{8}^{3,0}
\hskip -1.5pt \otimes \hskip -1.5pt
1_{3}^{1,4}$,
pord$(\rho_\text{isum}(\mathfrak{t})) = 24$,

\vskip 0.7ex
\hangindent=5.5em \hangafter=1
{\white .}\hskip 1em $\rho_\text{isum}(\mathfrak{t})$ =
 $( 0,
\frac{1}{3},
\frac{2}{3} )
\oplus
( \frac{2}{3},
\frac{1}{24},
\frac{13}{24} )
$,

\vskip 0.7ex
\hangindent=5.5em \hangafter=1
{\white .}\hskip 1em $\rho_\text{isum}(\mathfrak{s})$ =
($-\frac{1}{3}$,
$\frac{2}{3}$,
$\frac{2}{3}$;
$-\frac{1}{3}$,
$\frac{2}{3}$;
$-\frac{1}{3}$)
 $\oplus$
$\mathrm{i}$($0$,
$\sqrt{\frac{1}{2}}$,
$\sqrt{\frac{1}{2}}$;\ \ 
$\frac{1}{2}$,
$-\frac{1}{2}$;\ \ 
$\frac{1}{2}$)

Fail:
number of self dual objects $|$Tr($\rho(\mathfrak s^2)$)$|$ = 0. Prop. B.4 (1)\
 eqn. (B.16)

 \ \color{black}

\noindent 209: (dims,levels) = $(3 , 
3;3,
24
)$,
irreps = $3_{3}^{1,0}\oplus
3_{8}^{1,0}
\hskip -1.5pt \otimes \hskip -1.5pt
1_{3}^{1,4}$,
pord$(\rho_\text{isum}(\mathfrak{t})) = 24$,

\vskip 0.7ex
\hangindent=5.5em \hangafter=1
{\white .}\hskip 1em $\rho_\text{isum}(\mathfrak{t})$ =
 $( 0,
\frac{1}{3},
\frac{2}{3} )
\oplus
( \frac{2}{3},
\frac{7}{24},
\frac{19}{24} )
$,

\vskip 0.7ex
\hangindent=5.5em \hangafter=1
{\white .}\hskip 1em $\rho_\text{isum}(\mathfrak{s})$ =
($-\frac{1}{3}$,
$\frac{2}{3}$,
$\frac{2}{3}$;
$-\frac{1}{3}$,
$\frac{2}{3}$;
$-\frac{1}{3}$)
 $\oplus$
$\mathrm{i}$($0$,
$\sqrt{\frac{1}{2}}$,
$\sqrt{\frac{1}{2}}$;\ \ 
$-\frac{1}{2}$,
$\frac{1}{2}$;\ \ 
$-\frac{1}{2}$)

Fail:
number of self dual objects $|$Tr($\rho(\mathfrak s^2)$)$|$ = 0. Prop. B.4 (1)\
 eqn. (B.16)

 \ \color{black}

 \color{blue}

\noindent 210: (dims,levels) = $(3 , 
3;4,
5
)$,
irreps = $3_{4}^{1,0}\oplus
3_{5}^{1}$,
pord$(\rho_\text{isum}(\mathfrak{t})) = 20$,

\vskip 0.7ex
\hangindent=5.5em \hangafter=1
{\white .}\hskip 1em $\rho_\text{isum}(\mathfrak{t})$ =
 $( 0,
\frac{1}{4},
\frac{3}{4} )
\oplus
( 0,
\frac{1}{5},
\frac{4}{5} )
$,

\vskip 0.7ex
\hangindent=5.5em \hangafter=1
{\white .}\hskip 1em $\rho_\text{isum}(\mathfrak{s})$ =
($0$,
$\sqrt{\frac{1}{2}}$,
$\sqrt{\frac{1}{2}}$;
$-\frac{1}{2}$,
$\frac{1}{2}$;
$-\frac{1}{2}$)
 $\oplus$
($\sqrt{\frac{1}{5}}$,
$-\sqrt{\frac{2}{5}}$,
$-\sqrt{\frac{2}{5}}$;
$-\frac{5+\sqrt{5}}{10}$,
$\frac{5-\sqrt{5}}{10}$;
$-\frac{5+\sqrt{5}}{10}$)

Pass. 

 \ \color{black}

 \color{blue}

\noindent 211: (dims,levels) = $(3 , 
3;4,
5
)$,
irreps = $3_{4}^{1,0}\oplus
3_{5}^{3}$,
pord$(\rho_\text{isum}(\mathfrak{t})) = 20$,

\vskip 0.7ex
\hangindent=5.5em \hangafter=1
{\white .}\hskip 1em $\rho_\text{isum}(\mathfrak{t})$ =
 $( 0,
\frac{1}{4},
\frac{3}{4} )
\oplus
( 0,
\frac{2}{5},
\frac{3}{5} )
$,

\vskip 0.7ex
\hangindent=5.5em \hangafter=1
{\white .}\hskip 1em $\rho_\text{isum}(\mathfrak{s})$ =
($0$,
$\sqrt{\frac{1}{2}}$,
$\sqrt{\frac{1}{2}}$;
$-\frac{1}{2}$,
$\frac{1}{2}$;
$-\frac{1}{2}$)
 $\oplus$
($-\sqrt{\frac{1}{5}}$,
$-\sqrt{\frac{2}{5}}$,
$-\sqrt{\frac{2}{5}}$;
$\frac{-5+\sqrt{5}}{10}$,
$\frac{5+\sqrt{5}}{10}$;
$\frac{-5+\sqrt{5}}{10}$)

Pass. 

 \ \color{black}

\noindent 212: (dims,levels) = $(3 , 
3;4,
8
)$,
irreps = $3_{4}^{1,0}\oplus
3_{8}^{1,0}$,
pord$(\rho_\text{isum}(\mathfrak{t})) = 8$,

\vskip 0.7ex
\hangindent=5.5em \hangafter=1
{\white .}\hskip 1em $\rho_\text{isum}(\mathfrak{t})$ =
 $( 0,
\frac{1}{4},
\frac{3}{4} )
\oplus
( 0,
\frac{1}{8},
\frac{5}{8} )
$,

\vskip 0.7ex
\hangindent=5.5em \hangafter=1
{\white .}\hskip 1em $\rho_\text{isum}(\mathfrak{s})$ =
($0$,
$\sqrt{\frac{1}{2}}$,
$\sqrt{\frac{1}{2}}$;
$-\frac{1}{2}$,
$\frac{1}{2}$;
$-\frac{1}{2}$)
 $\oplus$
$\mathrm{i}$($0$,
$\sqrt{\frac{1}{2}}$,
$\sqrt{\frac{1}{2}}$;\ \ 
$-\frac{1}{2}$,
$\frac{1}{2}$;\ \ 
$-\frac{1}{2}$)

Fail:
number of self dual objects $|$Tr($\rho(\mathfrak s^2)$)$|$ = 0. Prop. B.4 (1)\
 eqn. (B.16)

 \ \color{black}

\noindent 213: (dims,levels) = $(3 , 
3;4,
8
)$,
irreps = $3_{4}^{1,0}\oplus
3_{8}^{3,0}$,
pord$(\rho_\text{isum}(\mathfrak{t})) = 8$,

\vskip 0.7ex
\hangindent=5.5em \hangafter=1
{\white .}\hskip 1em $\rho_\text{isum}(\mathfrak{t})$ =
 $( 0,
\frac{1}{4},
\frac{3}{4} )
\oplus
( 0,
\frac{3}{8},
\frac{7}{8} )
$,

\vskip 0.7ex
\hangindent=5.5em \hangafter=1
{\white .}\hskip 1em $\rho_\text{isum}(\mathfrak{s})$ =
($0$,
$\sqrt{\frac{1}{2}}$,
$\sqrt{\frac{1}{2}}$;
$-\frac{1}{2}$,
$\frac{1}{2}$;
$-\frac{1}{2}$)
 $\oplus$
$\mathrm{i}$($0$,
$\sqrt{\frac{1}{2}}$,
$\sqrt{\frac{1}{2}}$;\ \ 
$\frac{1}{2}$,
$-\frac{1}{2}$;\ \ 
$\frac{1}{2}$)

Fail:
number of self dual objects $|$Tr($\rho(\mathfrak s^2)$)$|$ = 0. Prop. B.4 (1)\
 eqn. (B.16)

 \ \color{black}

\noindent 214: (dims,levels) = $(3 , 
3;4,
8
)$,
irreps = $3_{4}^{1,0}\oplus
3_{8}^{3,3}$,
pord$(\rho_\text{isum}(\mathfrak{t})) = 8$,

\vskip 0.7ex
\hangindent=5.5em \hangafter=1
{\white .}\hskip 1em $\rho_\text{isum}(\mathfrak{t})$ =
 $( 0,
\frac{1}{4},
\frac{3}{4} )
\oplus
( \frac{1}{4},
\frac{1}{8},
\frac{5}{8} )
$,

\vskip 0.7ex
\hangindent=5.5em \hangafter=1
{\white .}\hskip 1em $\rho_\text{isum}(\mathfrak{s})$ =
($0$,
$\sqrt{\frac{1}{2}}$,
$\sqrt{\frac{1}{2}}$;
$-\frac{1}{2}$,
$\frac{1}{2}$;
$-\frac{1}{2}$)
 $\oplus$
($0$,
$\sqrt{\frac{1}{2}}$,
$\sqrt{\frac{1}{2}}$;
$-\frac{1}{2}$,
$\frac{1}{2}$;
$-\frac{1}{2}$)

Fail:
$\sigma(\rho(\mathfrak s)_\mathrm{ndeg}) \neq
 (\rho(\mathfrak t)^a \rho(\mathfrak s) \rho(\mathfrak t)^b
 \rho(\mathfrak s) \rho(\mathfrak t)^a)_\mathrm{ndeg}$,
 $\sigma = a$ = 3. Prop. B.5 (3) eqn. (B.25)

 \ \color{black}

\noindent 215: (dims,levels) = $(3 , 
3;4,
8
)$,
irreps = $3_{4}^{1,0}\oplus
3_{8}^{1,3}$,
pord$(\rho_\text{isum}(\mathfrak{t})) = 8$,

\vskip 0.7ex
\hangindent=5.5em \hangafter=1
{\white .}\hskip 1em $\rho_\text{isum}(\mathfrak{t})$ =
 $( 0,
\frac{1}{4},
\frac{3}{4} )
\oplus
( \frac{1}{4},
\frac{3}{8},
\frac{7}{8} )
$,

\vskip 0.7ex
\hangindent=5.5em \hangafter=1
{\white .}\hskip 1em $\rho_\text{isum}(\mathfrak{s})$ =
($0$,
$\sqrt{\frac{1}{2}}$,
$\sqrt{\frac{1}{2}}$;
$-\frac{1}{2}$,
$\frac{1}{2}$;
$-\frac{1}{2}$)
 $\oplus$
($0$,
$\sqrt{\frac{1}{2}}$,
$\sqrt{\frac{1}{2}}$;
$\frac{1}{2}$,
$-\frac{1}{2}$;
$\frac{1}{2}$)

Fail:
$\sigma(\rho(\mathfrak s)_\mathrm{ndeg}) \neq
 (\rho(\mathfrak t)^a \rho(\mathfrak s) \rho(\mathfrak t)^b
 \rho(\mathfrak s) \rho(\mathfrak t)^a)_\mathrm{ndeg}$,
 $\sigma = a$ = 3. Prop. B.5 (3) eqn. (B.25)

 \ \color{black}

\noindent 216: (dims,levels) = $(3 , 
3;4,
8
)$,
irreps = $3_{4}^{1,0}\oplus
3_{8}^{3,9}$,
pord$(\rho_\text{isum}(\mathfrak{t})) = 8$,

\vskip 0.7ex
\hangindent=5.5em \hangafter=1
{\white .}\hskip 1em $\rho_\text{isum}(\mathfrak{t})$ =
 $( 0,
\frac{1}{4},
\frac{3}{4} )
\oplus
( \frac{3}{4},
\frac{1}{8},
\frac{5}{8} )
$,

\vskip 0.7ex
\hangindent=5.5em \hangafter=1
{\white .}\hskip 1em $\rho_\text{isum}(\mathfrak{s})$ =
($0$,
$\sqrt{\frac{1}{2}}$,
$\sqrt{\frac{1}{2}}$;
$-\frac{1}{2}$,
$\frac{1}{2}$;
$-\frac{1}{2}$)
 $\oplus$
($0$,
$\sqrt{\frac{1}{2}}$,
$\sqrt{\frac{1}{2}}$;
$\frac{1}{2}$,
$-\frac{1}{2}$;
$\frac{1}{2}$)

Fail:
$\sigma(\rho(\mathfrak s)_\mathrm{ndeg}) \neq
 (\rho(\mathfrak t)^a \rho(\mathfrak s) \rho(\mathfrak t)^b
 \rho(\mathfrak s) \rho(\mathfrak t)^a)_\mathrm{ndeg}$,
 $\sigma = a$ = 3. Prop. B.5 (3) eqn. (B.25)

 \ \color{black}

\noindent 217: (dims,levels) = $(3 , 
3;4,
8
)$,
irreps = $3_{4}^{1,0}\oplus
3_{8}^{1,9}$,
pord$(\rho_\text{isum}(\mathfrak{t})) = 8$,

\vskip 0.7ex
\hangindent=5.5em \hangafter=1
{\white .}\hskip 1em $\rho_\text{isum}(\mathfrak{t})$ =
 $( 0,
\frac{1}{4},
\frac{3}{4} )
\oplus
( \frac{3}{4},
\frac{3}{8},
\frac{7}{8} )
$,

\vskip 0.7ex
\hangindent=5.5em \hangafter=1
{\white .}\hskip 1em $\rho_\text{isum}(\mathfrak{s})$ =
($0$,
$\sqrt{\frac{1}{2}}$,
$\sqrt{\frac{1}{2}}$;
$-\frac{1}{2}$,
$\frac{1}{2}$;
$-\frac{1}{2}$)
 $\oplus$
($0$,
$\sqrt{\frac{1}{2}}$,
$\sqrt{\frac{1}{2}}$;
$-\frac{1}{2}$,
$\frac{1}{2}$;
$-\frac{1}{2}$)

Fail:
$\sigma(\rho(\mathfrak s)_\mathrm{ndeg}) \neq
 (\rho(\mathfrak t)^a \rho(\mathfrak s) \rho(\mathfrak t)^b
 \rho(\mathfrak s) \rho(\mathfrak t)^a)_\mathrm{ndeg}$,
 $\sigma = a$ = 3. Prop. B.5 (3) eqn. (B.25)

 \ \color{black}

\noindent 218: (dims,levels) = $(3 , 
3;4,
12
)$,
irreps = $3_{4}^{1,0}\oplus
3_{3}^{1,0}
\hskip -1.5pt \otimes \hskip -1.5pt
1_{4}^{1,0}$,
pord$(\rho_\text{isum}(\mathfrak{t})) = 12$,

\vskip 0.7ex
\hangindent=5.5em \hangafter=1
{\white .}\hskip 1em $\rho_\text{isum}(\mathfrak{t})$ =
 $( 0,
\frac{1}{4},
\frac{3}{4} )
\oplus
( \frac{1}{4},
\frac{7}{12},
\frac{11}{12} )
$,

\vskip 0.7ex
\hangindent=5.5em \hangafter=1
{\white .}\hskip 1em $\rho_\text{isum}(\mathfrak{s})$ =
($0$,
$\sqrt{\frac{1}{2}}$,
$\sqrt{\frac{1}{2}}$;
$-\frac{1}{2}$,
$\frac{1}{2}$;
$-\frac{1}{2}$)
 $\oplus$
$\mathrm{i}$($-\frac{1}{3}$,
$\frac{2}{3}$,
$\frac{2}{3}$;\ \ 
$-\frac{1}{3}$,
$\frac{2}{3}$;\ \ 
$-\frac{1}{3}$)

Fail:
number of self dual objects $|$Tr($\rho(\mathfrak s^2)$)$|$ = 0. Prop. B.4 (1)\
 eqn. (B.16)

 \ \color{black}

\noindent 219: (dims,levels) = $(3 , 
3;4,
12
)$,
irreps = $3_{4}^{1,0}\oplus
3_{3}^{1,0}
\hskip -1.5pt \otimes \hskip -1.5pt
1_{4}^{1,6}$,
pord$(\rho_\text{isum}(\mathfrak{t})) = 12$,

\vskip 0.7ex
\hangindent=5.5em \hangafter=1
{\white .}\hskip 1em $\rho_\text{isum}(\mathfrak{t})$ =
 $( 0,
\frac{1}{4},
\frac{3}{4} )
\oplus
( \frac{3}{4},
\frac{1}{12},
\frac{5}{12} )
$,

\vskip 0.7ex
\hangindent=5.5em \hangafter=1
{\white .}\hskip 1em $\rho_\text{isum}(\mathfrak{s})$ =
($0$,
$\sqrt{\frac{1}{2}}$,
$\sqrt{\frac{1}{2}}$;
$-\frac{1}{2}$,
$\frac{1}{2}$;
$-\frac{1}{2}$)
 $\oplus$
$\mathrm{i}$($\frac{1}{3}$,
$\frac{2}{3}$,
$\frac{2}{3}$;\ \ 
$\frac{1}{3}$,
$-\frac{2}{3}$;\ \ 
$\frac{1}{3}$)

Fail:
number of self dual objects $|$Tr($\rho(\mathfrak s^2)$)$|$ = 0. Prop. B.4 (1)\
 eqn. (B.16)

 \ \color{black}

\noindent 220: (dims,levels) = $(3 , 
3;4,
20
)$,
irreps = $3_{4}^{1,0}\oplus
3_{5}^{1}
\hskip -1.5pt \otimes \hskip -1.5pt
1_{4}^{1,0}$,
pord$(\rho_\text{isum}(\mathfrak{t})) = 20$,

\vskip 0.7ex
\hangindent=5.5em \hangafter=1
{\white .}\hskip 1em $\rho_\text{isum}(\mathfrak{t})$ =
 $( 0,
\frac{1}{4},
\frac{3}{4} )
\oplus
( \frac{1}{4},
\frac{1}{20},
\frac{9}{20} )
$,

\vskip 0.7ex
\hangindent=5.5em \hangafter=1
{\white .}\hskip 1em $\rho_\text{isum}(\mathfrak{s})$ =
($0$,
$\sqrt{\frac{1}{2}}$,
$\sqrt{\frac{1}{2}}$;
$-\frac{1}{2}$,
$\frac{1}{2}$;
$-\frac{1}{2}$)
 $\oplus$
$\mathrm{i}$($\sqrt{\frac{1}{5}}$,
$\sqrt{\frac{2}{5}}$,
$\sqrt{\frac{2}{5}}$;\ \ 
$-\frac{5+\sqrt{5}}{10}$,
$\frac{5-\sqrt{5}}{10}$;\ \ 
$-\frac{5+\sqrt{5}}{10}$)

Fail:
number of self dual objects $|$Tr($\rho(\mathfrak s^2)$)$|$ = 0. Prop. B.4 (1)\
 eqn. (B.16)

 \ \color{black}

\noindent 221: (dims,levels) = $(3 , 
3;4,
20
)$,
irreps = $3_{4}^{1,0}\oplus
3_{5}^{3}
\hskip -1.5pt \otimes \hskip -1.5pt
1_{4}^{1,0}$,
pord$(\rho_\text{isum}(\mathfrak{t})) = 20$,

\vskip 0.7ex
\hangindent=5.5em \hangafter=1
{\white .}\hskip 1em $\rho_\text{isum}(\mathfrak{t})$ =
 $( 0,
\frac{1}{4},
\frac{3}{4} )
\oplus
( \frac{1}{4},
\frac{13}{20},
\frac{17}{20} )
$,

\vskip 0.7ex
\hangindent=5.5em \hangafter=1
{\white .}\hskip 1em $\rho_\text{isum}(\mathfrak{s})$ =
($0$,
$\sqrt{\frac{1}{2}}$,
$\sqrt{\frac{1}{2}}$;
$-\frac{1}{2}$,
$\frac{1}{2}$;
$-\frac{1}{2}$)
 $\oplus$
$\mathrm{i}$($-\sqrt{\frac{1}{5}}$,
$\sqrt{\frac{2}{5}}$,
$\sqrt{\frac{2}{5}}$;\ \ 
$\frac{-5+\sqrt{5}}{10}$,
$\frac{5+\sqrt{5}}{10}$;\ \ 
$\frac{-5+\sqrt{5}}{10}$)

Fail:
number of self dual objects $|$Tr($\rho(\mathfrak s^2)$)$|$ = 0. Prop. B.4 (1)\
 eqn. (B.16)

 \ \color{black}

\noindent 222: (dims,levels) = $(3 , 
3;4,
20
)$,
irreps = $3_{4}^{1,0}\oplus
3_{5}^{3}
\hskip -1.5pt \otimes \hskip -1.5pt
1_{4}^{1,6}$,
pord$(\rho_\text{isum}(\mathfrak{t})) = 20$,

\vskip 0.7ex
\hangindent=5.5em \hangafter=1
{\white .}\hskip 1em $\rho_\text{isum}(\mathfrak{t})$ =
 $( 0,
\frac{1}{4},
\frac{3}{4} )
\oplus
( \frac{3}{4},
\frac{3}{20},
\frac{7}{20} )
$,

\vskip 0.7ex
\hangindent=5.5em \hangafter=1
{\white .}\hskip 1em $\rho_\text{isum}(\mathfrak{s})$ =
($0$,
$\sqrt{\frac{1}{2}}$,
$\sqrt{\frac{1}{2}}$;
$-\frac{1}{2}$,
$\frac{1}{2}$;
$-\frac{1}{2}$)
 $\oplus$
$\mathrm{i}$($\sqrt{\frac{1}{5}}$,
$\sqrt{\frac{2}{5}}$,
$\sqrt{\frac{2}{5}}$;\ \ 
$\frac{5-\sqrt{5}}{10}$,
$-\frac{5+\sqrt{5}}{10}$;\ \ 
$\frac{5-\sqrt{5}}{10}$)

Fail:
number of self dual objects $|$Tr($\rho(\mathfrak s^2)$)$|$ = 0. Prop. B.4 (1)\
 eqn. (B.16)

 \ \color{black}

\noindent 223: (dims,levels) = $(3 , 
3;4,
20
)$,
irreps = $3_{4}^{1,0}\oplus
3_{5}^{1}
\hskip -1.5pt \otimes \hskip -1.5pt
1_{4}^{1,6}$,
pord$(\rho_\text{isum}(\mathfrak{t})) = 20$,

\vskip 0.7ex
\hangindent=5.5em \hangafter=1
{\white .}\hskip 1em $\rho_\text{isum}(\mathfrak{t})$ =
 $( 0,
\frac{1}{4},
\frac{3}{4} )
\oplus
( \frac{3}{4},
\frac{11}{20},
\frac{19}{20} )
$,

\vskip 0.7ex
\hangindent=5.5em \hangafter=1
{\white .}\hskip 1em $\rho_\text{isum}(\mathfrak{s})$ =
($0$,
$\sqrt{\frac{1}{2}}$,
$\sqrt{\frac{1}{2}}$;
$-\frac{1}{2}$,
$\frac{1}{2}$;
$-\frac{1}{2}$)
 $\oplus$
$\mathrm{i}$($-\sqrt{\frac{1}{5}}$,
$\sqrt{\frac{2}{5}}$,
$\sqrt{\frac{2}{5}}$;\ \ 
$\frac{5+\sqrt{5}}{10}$,
$\frac{-5+\sqrt{5}}{10}$;\ \ 
$\frac{5+\sqrt{5}}{10}$)

Fail:
number of self dual objects $|$Tr($\rho(\mathfrak s^2)$)$|$ = 0. Prop. B.4 (1)\
 eqn. (B.16)

 \ \color{black}

\noindent 224: (dims,levels) = $(3 , 
3;5,
5
)$,
irreps = $3_{5}^{1}\oplus
3_{5}^{3}$,
pord$(\rho_\text{isum}(\mathfrak{t})) = 5$,

\vskip 0.7ex
\hangindent=5.5em \hangafter=1
{\white .}\hskip 1em $\rho_\text{isum}(\mathfrak{t})$ =
 $( 0,
\frac{1}{5},
\frac{4}{5} )
\oplus
( 0,
\frac{2}{5},
\frac{3}{5} )
$,

\vskip 0.7ex
\hangindent=5.5em \hangafter=1
{\white .}\hskip 1em $\rho_\text{isum}(\mathfrak{s})$ =
($\sqrt{\frac{1}{5}}$,
$-\sqrt{\frac{2}{5}}$,
$-\sqrt{\frac{2}{5}}$;
$-\frac{5+\sqrt{5}}{10}$,
$\frac{5-\sqrt{5}}{10}$;
$-\frac{5+\sqrt{5}}{10}$)
 $\oplus$
($-\sqrt{\frac{1}{5}}$,
$-\sqrt{\frac{2}{5}}$,
$-\sqrt{\frac{2}{5}}$;
$\frac{-5+\sqrt{5}}{10}$,
$\frac{5+\sqrt{5}}{10}$;
$\frac{-5+\sqrt{5}}{10}$)

Fail:
Integral: $D_{\rho}(\sigma)_{\theta} \propto $ id,
 for all $\sigma$ and all $\theta$-eigenspaces that can contain unit. Prop. B.5 (6)

 \ \color{black}

\noindent 225: (dims,levels) = $(3 , 
3;5,
8
)$,
irreps = $3_{5}^{1}\oplus
3_{8}^{1,0}$,
pord$(\rho_\text{isum}(\mathfrak{t})) = 40$,

\vskip 0.7ex
\hangindent=5.5em \hangafter=1
{\white .}\hskip 1em $\rho_\text{isum}(\mathfrak{t})$ =
 $( 0,
\frac{1}{5},
\frac{4}{5} )
\oplus
( 0,
\frac{1}{8},
\frac{5}{8} )
$,

\vskip 0.7ex
\hangindent=5.5em \hangafter=1
{\white .}\hskip 1em $\rho_\text{isum}(\mathfrak{s})$ =
($\sqrt{\frac{1}{5}}$,
$-\sqrt{\frac{2}{5}}$,
$-\sqrt{\frac{2}{5}}$;
$-\frac{5+\sqrt{5}}{10}$,
$\frac{5-\sqrt{5}}{10}$;
$-\frac{5+\sqrt{5}}{10}$)
 $\oplus$
$\mathrm{i}$($0$,
$\sqrt{\frac{1}{2}}$,
$\sqrt{\frac{1}{2}}$;\ \ 
$-\frac{1}{2}$,
$\frac{1}{2}$;\ \ 
$-\frac{1}{2}$)

Fail:
number of self dual objects $|$Tr($\rho(\mathfrak s^2)$)$|$ = 0. Prop. B.4 (1)\
 eqn. (B.16)

 \ \color{black}

\noindent 226: (dims,levels) = $(3 , 
3;5,
8
)$,
irreps = $3_{5}^{1}\oplus
3_{8}^{3,0}$,
pord$(\rho_\text{isum}(\mathfrak{t})) = 40$,

\vskip 0.7ex
\hangindent=5.5em \hangafter=1
{\white .}\hskip 1em $\rho_\text{isum}(\mathfrak{t})$ =
 $( 0,
\frac{1}{5},
\frac{4}{5} )
\oplus
( 0,
\frac{3}{8},
\frac{7}{8} )
$,

\vskip 0.7ex
\hangindent=5.5em \hangafter=1
{\white .}\hskip 1em $\rho_\text{isum}(\mathfrak{s})$ =
($\sqrt{\frac{1}{5}}$,
$-\sqrt{\frac{2}{5}}$,
$-\sqrt{\frac{2}{5}}$;
$-\frac{5+\sqrt{5}}{10}$,
$\frac{5-\sqrt{5}}{10}$;
$-\frac{5+\sqrt{5}}{10}$)
 $\oplus$
$\mathrm{i}$($0$,
$\sqrt{\frac{1}{2}}$,
$\sqrt{\frac{1}{2}}$;\ \ 
$\frac{1}{2}$,
$-\frac{1}{2}$;\ \ 
$\frac{1}{2}$)

Fail:
number of self dual objects $|$Tr($\rho(\mathfrak s^2)$)$|$ = 0. Prop. B.4 (1)\
 eqn. (B.16)

 \ \color{black}

\noindent 227: (dims,levels) = $(3 , 
3;5,
8
)$,
irreps = $3_{5}^{3}\oplus
3_{8}^{1,0}$,
pord$(\rho_\text{isum}(\mathfrak{t})) = 40$,

\vskip 0.7ex
\hangindent=5.5em \hangafter=1
{\white .}\hskip 1em $\rho_\text{isum}(\mathfrak{t})$ =
 $( 0,
\frac{2}{5},
\frac{3}{5} )
\oplus
( 0,
\frac{1}{8},
\frac{5}{8} )
$,

\vskip 0.7ex
\hangindent=5.5em \hangafter=1
{\white .}\hskip 1em $\rho_\text{isum}(\mathfrak{s})$ =
($-\sqrt{\frac{1}{5}}$,
$-\sqrt{\frac{2}{5}}$,
$-\sqrt{\frac{2}{5}}$;
$\frac{-5+\sqrt{5}}{10}$,
$\frac{5+\sqrt{5}}{10}$;
$\frac{-5+\sqrt{5}}{10}$)
 $\oplus$
$\mathrm{i}$($0$,
$\sqrt{\frac{1}{2}}$,
$\sqrt{\frac{1}{2}}$;\ \ 
$-\frac{1}{2}$,
$\frac{1}{2}$;\ \ 
$-\frac{1}{2}$)

Fail:
number of self dual objects $|$Tr($\rho(\mathfrak s^2)$)$|$ = 0. Prop. B.4 (1)\
 eqn. (B.16)

 \ \color{black}

\noindent 228: (dims,levels) = $(3 , 
3;5,
8
)$,
irreps = $3_{5}^{3}\oplus
3_{8}^{3,0}$,
pord$(\rho_\text{isum}(\mathfrak{t})) = 40$,

\vskip 0.7ex
\hangindent=5.5em \hangafter=1
{\white .}\hskip 1em $\rho_\text{isum}(\mathfrak{t})$ =
 $( 0,
\frac{2}{5},
\frac{3}{5} )
\oplus
( 0,
\frac{3}{8},
\frac{7}{8} )
$,

\vskip 0.7ex
\hangindent=5.5em \hangafter=1
{\white .}\hskip 1em $\rho_\text{isum}(\mathfrak{s})$ =
($-\sqrt{\frac{1}{5}}$,
$-\sqrt{\frac{2}{5}}$,
$-\sqrt{\frac{2}{5}}$;
$\frac{-5+\sqrt{5}}{10}$,
$\frac{5+\sqrt{5}}{10}$;
$\frac{-5+\sqrt{5}}{10}$)
 $\oplus$
$\mathrm{i}$($0$,
$\sqrt{\frac{1}{2}}$,
$\sqrt{\frac{1}{2}}$;\ \ 
$\frac{1}{2}$,
$-\frac{1}{2}$;\ \ 
$\frac{1}{2}$)

Fail:
number of self dual objects $|$Tr($\rho(\mathfrak s^2)$)$|$ = 0. Prop. B.4 (1)\
 eqn. (B.16)

 \ \color{black}

\noindent 229: (dims,levels) = $(3 , 
3;6,
8
)$,
irreps = $3_{3}^{1,0}
\hskip -1.5pt \otimes \hskip -1.5pt
1_{2}^{1,0}\oplus
3_{8}^{1,6}$,
pord$(\rho_\text{isum}(\mathfrak{t})) = 24$,

\vskip 0.7ex
\hangindent=5.5em \hangafter=1
{\white .}\hskip 1em $\rho_\text{isum}(\mathfrak{t})$ =
 $( \frac{1}{2},
\frac{1}{6},
\frac{5}{6} )
\oplus
( \frac{1}{2},
\frac{1}{8},
\frac{5}{8} )
$,

\vskip 0.7ex
\hangindent=5.5em \hangafter=1
{\white .}\hskip 1em $\rho_\text{isum}(\mathfrak{s})$ =
($\frac{1}{3}$,
$\frac{2}{3}$,
$\frac{2}{3}$;
$\frac{1}{3}$,
$-\frac{2}{3}$;
$\frac{1}{3}$)
 $\oplus$
$\mathrm{i}$($0$,
$\sqrt{\frac{1}{2}}$,
$\sqrt{\frac{1}{2}}$;\ \ 
$\frac{1}{2}$,
$-\frac{1}{2}$;\ \ 
$\frac{1}{2}$)

Fail:
number of self dual objects $|$Tr($\rho(\mathfrak s^2)$)$|$ = 0. Prop. B.4 (1)\
 eqn. (B.16)

 \ \color{black}

\noindent 230: (dims,levels) = $(3 , 
3;6,
8
)$,
irreps = $3_{3}^{1,0}
\hskip -1.5pt \otimes \hskip -1.5pt
1_{2}^{1,0}\oplus
3_{8}^{3,6}$,
pord$(\rho_\text{isum}(\mathfrak{t})) = 24$,

\vskip 0.7ex
\hangindent=5.5em \hangafter=1
{\white .}\hskip 1em $\rho_\text{isum}(\mathfrak{t})$ =
 $( \frac{1}{2},
\frac{1}{6},
\frac{5}{6} )
\oplus
( \frac{1}{2},
\frac{3}{8},
\frac{7}{8} )
$,

\vskip 0.7ex
\hangindent=5.5em \hangafter=1
{\white .}\hskip 1em $\rho_\text{isum}(\mathfrak{s})$ =
($\frac{1}{3}$,
$\frac{2}{3}$,
$\frac{2}{3}$;
$\frac{1}{3}$,
$-\frac{2}{3}$;
$\frac{1}{3}$)
 $\oplus$
$\mathrm{i}$($0$,
$\sqrt{\frac{1}{2}}$,
$\sqrt{\frac{1}{2}}$;\ \ 
$-\frac{1}{2}$,
$\frac{1}{2}$;\ \ 
$-\frac{1}{2}$)

Fail:
number of self dual objects $|$Tr($\rho(\mathfrak s^2)$)$|$ = 0. Prop. B.4 (1)\
 eqn. (B.16)

 \ \color{black}

 \color{blue}

\noindent 231: (dims,levels) = $(3 , 
3;6,
10
)$,
irreps = $3_{3}^{1,0}
\hskip -1.5pt \otimes \hskip -1.5pt
1_{2}^{1,0}\oplus
3_{5}^{3}
\hskip -1.5pt \otimes \hskip -1.5pt
1_{2}^{1,0}$,
pord$(\rho_\text{isum}(\mathfrak{t})) = 15$,

\vskip 0.7ex
\hangindent=5.5em \hangafter=1
{\white .}\hskip 1em $\rho_\text{isum}(\mathfrak{t})$ =
 $( \frac{1}{2},
\frac{1}{6},
\frac{5}{6} )
\oplus
( \frac{1}{2},
\frac{1}{10},
\frac{9}{10} )
$,

\vskip 0.7ex
\hangindent=5.5em \hangafter=1
{\white .}\hskip 1em $\rho_\text{isum}(\mathfrak{s})$ =
($\frac{1}{3}$,
$\frac{2}{3}$,
$\frac{2}{3}$;
$\frac{1}{3}$,
$-\frac{2}{3}$;
$\frac{1}{3}$)
 $\oplus$
($\sqrt{\frac{1}{5}}$,
$-\sqrt{\frac{2}{5}}$,
$-\sqrt{\frac{2}{5}}$;
$\frac{5-\sqrt{5}}{10}$,
$-\frac{5+\sqrt{5}}{10}$;
$\frac{5-\sqrt{5}}{10}$)

Pass. 

 \ \color{black}

 \color{blue}

\noindent 232: (dims,levels) = $(3 , 
3;6,
10
)$,
irreps = $3_{3}^{1,0}
\hskip -1.5pt \otimes \hskip -1.5pt
1_{2}^{1,0}\oplus
3_{5}^{1}
\hskip -1.5pt \otimes \hskip -1.5pt
1_{2}^{1,0}$,
pord$(\rho_\text{isum}(\mathfrak{t})) = 15$,

\vskip 0.7ex
\hangindent=5.5em \hangafter=1
{\white .}\hskip 1em $\rho_\text{isum}(\mathfrak{t})$ =
 $( \frac{1}{2},
\frac{1}{6},
\frac{5}{6} )
\oplus
( \frac{1}{2},
\frac{3}{10},
\frac{7}{10} )
$,

\vskip 0.7ex
\hangindent=5.5em \hangafter=1
{\white .}\hskip 1em $\rho_\text{isum}(\mathfrak{s})$ =
($\frac{1}{3}$,
$\frac{2}{3}$,
$\frac{2}{3}$;
$\frac{1}{3}$,
$-\frac{2}{3}$;
$\frac{1}{3}$)
 $\oplus$
($-\sqrt{\frac{1}{5}}$,
$-\sqrt{\frac{2}{5}}$,
$-\sqrt{\frac{2}{5}}$;
$\frac{5+\sqrt{5}}{10}$,
$\frac{-5+\sqrt{5}}{10}$;
$\frac{5+\sqrt{5}}{10}$)

Pass. 

 \ \color{black}

\noindent 233: (dims,levels) = $(3 , 
3;6,
12
)$,
irreps = $3_{3}^{1,0}
\hskip -1.5pt \otimes \hskip -1.5pt
1_{2}^{1,0}\oplus
3_{4}^{1,6}
\hskip -1.5pt \otimes \hskip -1.5pt
1_{3}^{1,4}$,
pord$(\rho_\text{isum}(\mathfrak{t})) = 12$,

\vskip 0.7ex
\hangindent=5.5em \hangafter=1
{\white .}\hskip 1em $\rho_\text{isum}(\mathfrak{t})$ =
 $( \frac{1}{2},
\frac{1}{6},
\frac{5}{6} )
\oplus
( \frac{1}{6},
\frac{5}{12},
\frac{11}{12} )
$,

\vskip 0.7ex
\hangindent=5.5em \hangafter=1
{\white .}\hskip 1em $\rho_\text{isum}(\mathfrak{s})$ =
($\frac{1}{3}$,
$\frac{2}{3}$,
$\frac{2}{3}$;
$\frac{1}{3}$,
$-\frac{2}{3}$;
$\frac{1}{3}$)
 $\oplus$
($0$,
$\sqrt{\frac{1}{2}}$,
$\sqrt{\frac{1}{2}}$;
$\frac{1}{2}$,
$-\frac{1}{2}$;
$\frac{1}{2}$)

Fail:
Integral: $D_{\rho}(\sigma)_{\theta} \propto $ id,
 for all $\sigma$ and all $\theta$-eigenspaces that can contain unit. Prop. B.5 (6)

 \ \color{black}

\noindent 234: (dims,levels) = $(3 , 
3;6,
12
)$,
irreps = $3_{3}^{1,0}
\hskip -1.5pt \otimes \hskip -1.5pt
1_{2}^{1,0}\oplus
3_{4}^{1,9}
\hskip -1.5pt \otimes \hskip -1.5pt
1_{3}^{1,0}$,
pord$(\rho_\text{isum}(\mathfrak{t})) = 12$,

\vskip 0.7ex
\hangindent=5.5em \hangafter=1
{\white .}\hskip 1em $\rho_\text{isum}(\mathfrak{t})$ =
 $( \frac{1}{2},
\frac{1}{6},
\frac{5}{6} )
\oplus
( \frac{1}{3},
\frac{5}{6},
\frac{1}{12} )
$,

\vskip 0.7ex
\hangindent=5.5em \hangafter=1
{\white .}\hskip 1em $\rho_\text{isum}(\mathfrak{s})$ =
($\frac{1}{3}$,
$\frac{2}{3}$,
$\frac{2}{3}$;
$\frac{1}{3}$,
$-\frac{2}{3}$;
$\frac{1}{3}$)
 $\oplus$
$\mathrm{i}$($\frac{1}{2}$,
$\frac{1}{2}$,
$\sqrt{\frac{1}{2}}$;\ \ 
$\frac{1}{2}$,
$-\sqrt{\frac{1}{2}}$;\ \ 
$0$)

Fail:
number of self dual objects $|$Tr($\rho(\mathfrak s^2)$)$|$ = 0. Prop. B.4 (1)\
 eqn. (B.16)

 \ \color{black}

\noindent 235: (dims,levels) = $(3 , 
3;6,
12
)$,
irreps = $3_{3}^{1,0}
\hskip -1.5pt \otimes \hskip -1.5pt
1_{2}^{1,0}\oplus
3_{4}^{1,3}
\hskip -1.5pt \otimes \hskip -1.5pt
1_{3}^{1,0}$,
pord$(\rho_\text{isum}(\mathfrak{t})) = 12$,

\vskip 0.7ex
\hangindent=5.5em \hangafter=1
{\white .}\hskip 1em $\rho_\text{isum}(\mathfrak{t})$ =
 $( \frac{1}{2},
\frac{1}{6},
\frac{5}{6} )
\oplus
( \frac{1}{3},
\frac{5}{6},
\frac{7}{12} )
$,

\vskip 0.7ex
\hangindent=5.5em \hangafter=1
{\white .}\hskip 1em $\rho_\text{isum}(\mathfrak{s})$ =
($\frac{1}{3}$,
$\frac{2}{3}$,
$\frac{2}{3}$;
$\frac{1}{3}$,
$-\frac{2}{3}$;
$\frac{1}{3}$)
 $\oplus$
$\mathrm{i}$($-\frac{1}{2}$,
$\frac{1}{2}$,
$\sqrt{\frac{1}{2}}$;\ \ 
$-\frac{1}{2}$,
$\sqrt{\frac{1}{2}}$;\ \ 
$0$)

Fail:
number of self dual objects $|$Tr($\rho(\mathfrak s^2)$)$|$ = 0. Prop. B.4 (1)\
 eqn. (B.16)

 \ \color{black}

\noindent 236: (dims,levels) = $(3 , 
3;6,
12
)$,
irreps = $3_{3}^{1,0}
\hskip -1.5pt \otimes \hskip -1.5pt
1_{2}^{1,0}\oplus
3_{4}^{1,9}
\hskip -1.5pt \otimes \hskip -1.5pt
1_{3}^{1,4}$,
pord$(\rho_\text{isum}(\mathfrak{t})) = 12$,

\vskip 0.7ex
\hangindent=5.5em \hangafter=1
{\white .}\hskip 1em $\rho_\text{isum}(\mathfrak{t})$ =
 $( \frac{1}{2},
\frac{1}{6},
\frac{5}{6} )
\oplus
( \frac{2}{3},
\frac{1}{6},
\frac{5}{12} )
$,

\vskip 0.7ex
\hangindent=5.5em \hangafter=1
{\white .}\hskip 1em $\rho_\text{isum}(\mathfrak{s})$ =
($\frac{1}{3}$,
$\frac{2}{3}$,
$\frac{2}{3}$;
$\frac{1}{3}$,
$-\frac{2}{3}$;
$\frac{1}{3}$)
 $\oplus$
$\mathrm{i}$($\frac{1}{2}$,
$\frac{1}{2}$,
$\sqrt{\frac{1}{2}}$;\ \ 
$\frac{1}{2}$,
$-\sqrt{\frac{1}{2}}$;\ \ 
$0$)

Fail:
number of self dual objects $|$Tr($\rho(\mathfrak s^2)$)$|$ = 0. Prop. B.4 (1)\
 eqn. (B.16)

 \ \color{black}

\noindent 237: (dims,levels) = $(3 , 
3;6,
12
)$,
irreps = $3_{3}^{1,0}
\hskip -1.5pt \otimes \hskip -1.5pt
1_{2}^{1,0}\oplus
3_{4}^{1,3}
\hskip -1.5pt \otimes \hskip -1.5pt
1_{3}^{1,4}$,
pord$(\rho_\text{isum}(\mathfrak{t})) = 12$,

\vskip 0.7ex
\hangindent=5.5em \hangafter=1
{\white .}\hskip 1em $\rho_\text{isum}(\mathfrak{t})$ =
 $( \frac{1}{2},
\frac{1}{6},
\frac{5}{6} )
\oplus
( \frac{2}{3},
\frac{1}{6},
\frac{11}{12} )
$,

\vskip 0.7ex
\hangindent=5.5em \hangafter=1
{\white .}\hskip 1em $\rho_\text{isum}(\mathfrak{s})$ =
($\frac{1}{3}$,
$\frac{2}{3}$,
$\frac{2}{3}$;
$\frac{1}{3}$,
$-\frac{2}{3}$;
$\frac{1}{3}$)
 $\oplus$
$\mathrm{i}$($-\frac{1}{2}$,
$\frac{1}{2}$,
$\sqrt{\frac{1}{2}}$;\ \ 
$-\frac{1}{2}$,
$\sqrt{\frac{1}{2}}$;\ \ 
$0$)

Fail:
number of self dual objects $|$Tr($\rho(\mathfrak s^2)$)$|$ = 0. Prop. B.4 (1)\
 eqn. (B.16)

 \ \color{black}

\noindent 238: (dims,levels) = $(3 , 
3;6,
12
)$,
irreps = $3_{3}^{1,0}
\hskip -1.5pt \otimes \hskip -1.5pt
1_{2}^{1,0}\oplus
3_{4}^{1,6}
\hskip -1.5pt \otimes \hskip -1.5pt
1_{3}^{1,0}$,
pord$(\rho_\text{isum}(\mathfrak{t})) = 12$,

\vskip 0.7ex
\hangindent=5.5em \hangafter=1
{\white .}\hskip 1em $\rho_\text{isum}(\mathfrak{t})$ =
 $( \frac{1}{2},
\frac{1}{6},
\frac{5}{6} )
\oplus
( \frac{5}{6},
\frac{1}{12},
\frac{7}{12} )
$,

\vskip 0.7ex
\hangindent=5.5em \hangafter=1
{\white .}\hskip 1em $\rho_\text{isum}(\mathfrak{s})$ =
($\frac{1}{3}$,
$\frac{2}{3}$,
$\frac{2}{3}$;
$\frac{1}{3}$,
$-\frac{2}{3}$;
$\frac{1}{3}$)
 $\oplus$
($0$,
$\sqrt{\frac{1}{2}}$,
$\sqrt{\frac{1}{2}}$;
$\frac{1}{2}$,
$-\frac{1}{2}$;
$\frac{1}{2}$)

Fail:
Integral: $D_{\rho}(\sigma)_{\theta} \propto $ id,
 for all $\sigma$ and all $\theta$-eigenspaces that can contain unit. Prop. B.5 (6)

 \ \color{black}

\noindent 239: (dims,levels) = $(3 , 
3;6,
24
)$,
irreps = $3_{3}^{1,0}
\hskip -1.5pt \otimes \hskip -1.5pt
1_{2}^{1,0}\oplus
3_{8}^{3,6}
\hskip -1.5pt \otimes \hskip -1.5pt
1_{3}^{1,4}$,
pord$(\rho_\text{isum}(\mathfrak{t})) = 24$,

\vskip 0.7ex
\hangindent=5.5em \hangafter=1
{\white .}\hskip 1em $\rho_\text{isum}(\mathfrak{t})$ =
 $( \frac{1}{2},
\frac{1}{6},
\frac{5}{6} )
\oplus
( \frac{1}{6},
\frac{1}{24},
\frac{13}{24} )
$,

\vskip 0.7ex
\hangindent=5.5em \hangafter=1
{\white .}\hskip 1em $\rho_\text{isum}(\mathfrak{s})$ =
($\frac{1}{3}$,
$\frac{2}{3}$,
$\frac{2}{3}$;
$\frac{1}{3}$,
$-\frac{2}{3}$;
$\frac{1}{3}$)
 $\oplus$
$\mathrm{i}$($0$,
$\sqrt{\frac{1}{2}}$,
$\sqrt{\frac{1}{2}}$;\ \ 
$-\frac{1}{2}$,
$\frac{1}{2}$;\ \ 
$-\frac{1}{2}$)

Fail:
number of self dual objects $|$Tr($\rho(\mathfrak s^2)$)$|$ = 0. Prop. B.4 (1)\
 eqn. (B.16)

 \ \color{black}

\noindent 240: (dims,levels) = $(3 , 
3;6,
24
)$,
irreps = $3_{3}^{1,0}
\hskip -1.5pt \otimes \hskip -1.5pt
1_{2}^{1,0}\oplus
3_{8}^{1,6}
\hskip -1.5pt \otimes \hskip -1.5pt
1_{3}^{1,4}$,
pord$(\rho_\text{isum}(\mathfrak{t})) = 24$,

\vskip 0.7ex
\hangindent=5.5em \hangafter=1
{\white .}\hskip 1em $\rho_\text{isum}(\mathfrak{t})$ =
 $( \frac{1}{2},
\frac{1}{6},
\frac{5}{6} )
\oplus
( \frac{1}{6},
\frac{7}{24},
\frac{19}{24} )
$,

\vskip 0.7ex
\hangindent=5.5em \hangafter=1
{\white .}\hskip 1em $\rho_\text{isum}(\mathfrak{s})$ =
($\frac{1}{3}$,
$\frac{2}{3}$,
$\frac{2}{3}$;
$\frac{1}{3}$,
$-\frac{2}{3}$;
$\frac{1}{3}$)
 $\oplus$
$\mathrm{i}$($0$,
$\sqrt{\frac{1}{2}}$,
$\sqrt{\frac{1}{2}}$;\ \ 
$\frac{1}{2}$,
$-\frac{1}{2}$;\ \ 
$\frac{1}{2}$)

Fail:
number of self dual objects $|$Tr($\rho(\mathfrak s^2)$)$|$ = 0. Prop. B.4 (1)\
 eqn. (B.16)

 \ \color{black}

\noindent 241: (dims,levels) = $(3 , 
3;6,
24
)$,
irreps = $3_{3}^{1,0}
\hskip -1.5pt \otimes \hskip -1.5pt
1_{2}^{1,0}\oplus
3_{8}^{3,6}
\hskip -1.5pt \otimes \hskip -1.5pt
1_{3}^{1,0}$,
pord$(\rho_\text{isum}(\mathfrak{t})) = 24$,

\vskip 0.7ex
\hangindent=5.5em \hangafter=1
{\white .}\hskip 1em $\rho_\text{isum}(\mathfrak{t})$ =
 $( \frac{1}{2},
\frac{1}{6},
\frac{5}{6} )
\oplus
( \frac{5}{6},
\frac{5}{24},
\frac{17}{24} )
$,

\vskip 0.7ex
\hangindent=5.5em \hangafter=1
{\white .}\hskip 1em $\rho_\text{isum}(\mathfrak{s})$ =
($\frac{1}{3}$,
$\frac{2}{3}$,
$\frac{2}{3}$;
$\frac{1}{3}$,
$-\frac{2}{3}$;
$\frac{1}{3}$)
 $\oplus$
$\mathrm{i}$($0$,
$\sqrt{\frac{1}{2}}$,
$\sqrt{\frac{1}{2}}$;\ \ 
$-\frac{1}{2}$,
$\frac{1}{2}$;\ \ 
$-\frac{1}{2}$)

Fail:
number of self dual objects $|$Tr($\rho(\mathfrak s^2)$)$|$ = 0. Prop. B.4 (1)\
 eqn. (B.16)

 \ \color{black}

\noindent 242: (dims,levels) = $(3 , 
3;6,
24
)$,
irreps = $3_{3}^{1,0}
\hskip -1.5pt \otimes \hskip -1.5pt
1_{2}^{1,0}\oplus
3_{8}^{1,6}
\hskip -1.5pt \otimes \hskip -1.5pt
1_{3}^{1,0}$,
pord$(\rho_\text{isum}(\mathfrak{t})) = 24$,

\vskip 0.7ex
\hangindent=5.5em \hangafter=1
{\white .}\hskip 1em $\rho_\text{isum}(\mathfrak{t})$ =
 $( \frac{1}{2},
\frac{1}{6},
\frac{5}{6} )
\oplus
( \frac{5}{6},
\frac{11}{24},
\frac{23}{24} )
$,

\vskip 0.7ex
\hangindent=5.5em \hangafter=1
{\white .}\hskip 1em $\rho_\text{isum}(\mathfrak{s})$ =
($\frac{1}{3}$,
$\frac{2}{3}$,
$\frac{2}{3}$;
$\frac{1}{3}$,
$-\frac{2}{3}$;
$\frac{1}{3}$)
 $\oplus$
$\mathrm{i}$($0$,
$\sqrt{\frac{1}{2}}$,
$\sqrt{\frac{1}{2}}$;\ \ 
$\frac{1}{2}$,
$-\frac{1}{2}$;\ \ 
$\frac{1}{2}$)

Fail:
number of self dual objects $|$Tr($\rho(\mathfrak s^2)$)$|$ = 0. Prop. B.4 (1)\
 eqn. (B.16)

 \ \color{black}

 \color{blue}

\noindent 243: (dims,levels) = $(3 , 
3;6,
30
)$,
irreps = $3_{3}^{1,0}
\hskip -1.5pt \otimes \hskip -1.5pt
1_{2}^{1,0}\oplus
3_{5}^{1}
\hskip -1.5pt \otimes \hskip -1.5pt
1_{3}^{1,4}
\hskip -1.5pt \otimes \hskip -1.5pt
1_{2}^{1,0}$,
pord$(\rho_\text{isum}(\mathfrak{t})) = 15$,

\vskip 0.7ex
\hangindent=5.5em \hangafter=1
{\white .}\hskip 1em $\rho_\text{isum}(\mathfrak{t})$ =
 $( \frac{1}{2},
\frac{1}{6},
\frac{5}{6} )
\oplus
( \frac{1}{6},
\frac{11}{30},
\frac{29}{30} )
$,

\vskip 0.7ex
\hangindent=5.5em \hangafter=1
{\white .}\hskip 1em $\rho_\text{isum}(\mathfrak{s})$ =
($\frac{1}{3}$,
$\frac{2}{3}$,
$\frac{2}{3}$;
$\frac{1}{3}$,
$-\frac{2}{3}$;
$\frac{1}{3}$)
 $\oplus$
($-\sqrt{\frac{1}{5}}$,
$-\sqrt{\frac{2}{5}}$,
$-\sqrt{\frac{2}{5}}$;
$\frac{5+\sqrt{5}}{10}$,
$\frac{-5+\sqrt{5}}{10}$;
$\frac{5+\sqrt{5}}{10}$)

Pass. 

 \ \color{black}

 \color{blue}

\noindent 244: (dims,levels) = $(3 , 
3;6,
30
)$,
irreps = $3_{3}^{1,0}
\hskip -1.5pt \otimes \hskip -1.5pt
1_{2}^{1,0}\oplus
3_{5}^{3}
\hskip -1.5pt \otimes \hskip -1.5pt
1_{3}^{1,4}
\hskip -1.5pt \otimes \hskip -1.5pt
1_{2}^{1,0}$,
pord$(\rho_\text{isum}(\mathfrak{t})) = 15$,

\vskip 0.7ex
\hangindent=5.5em \hangafter=1
{\white .}\hskip 1em $\rho_\text{isum}(\mathfrak{t})$ =
 $( \frac{1}{2},
\frac{1}{6},
\frac{5}{6} )
\oplus
( \frac{1}{6},
\frac{17}{30},
\frac{23}{30} )
$,

\vskip 0.7ex
\hangindent=5.5em \hangafter=1
{\white .}\hskip 1em $\rho_\text{isum}(\mathfrak{s})$ =
($\frac{1}{3}$,
$\frac{2}{3}$,
$\frac{2}{3}$;
$\frac{1}{3}$,
$-\frac{2}{3}$;
$\frac{1}{3}$)
 $\oplus$
($\sqrt{\frac{1}{5}}$,
$-\sqrt{\frac{2}{5}}$,
$-\sqrt{\frac{2}{5}}$;
$\frac{5-\sqrt{5}}{10}$,
$-\frac{5+\sqrt{5}}{10}$;
$\frac{5-\sqrt{5}}{10}$)

Pass. 

 \ \color{black}

 \color{blue}

\noindent 245: (dims,levels) = $(3 , 
3;6,
30
)$,
irreps = $3_{3}^{1,0}
\hskip -1.5pt \otimes \hskip -1.5pt
1_{2}^{1,0}\oplus
3_{5}^{1}
\hskip -1.5pt \otimes \hskip -1.5pt
1_{3}^{1,0}
\hskip -1.5pt \otimes \hskip -1.5pt
1_{2}^{1,0}$,
pord$(\rho_\text{isum}(\mathfrak{t})) = 15$,

\vskip 0.7ex
\hangindent=5.5em \hangafter=1
{\white .}\hskip 1em $\rho_\text{isum}(\mathfrak{t})$ =
 $( \frac{1}{2},
\frac{1}{6},
\frac{5}{6} )
\oplus
( \frac{5}{6},
\frac{1}{30},
\frac{19}{30} )
$,

\vskip 0.7ex
\hangindent=5.5em \hangafter=1
{\white .}\hskip 1em $\rho_\text{isum}(\mathfrak{s})$ =
($\frac{1}{3}$,
$\frac{2}{3}$,
$\frac{2}{3}$;
$\frac{1}{3}$,
$-\frac{2}{3}$;
$\frac{1}{3}$)
 $\oplus$
($-\sqrt{\frac{1}{5}}$,
$-\sqrt{\frac{2}{5}}$,
$-\sqrt{\frac{2}{5}}$;
$\frac{5+\sqrt{5}}{10}$,
$\frac{-5+\sqrt{5}}{10}$;
$\frac{5+\sqrt{5}}{10}$)

Pass. 

 \ \color{black}

 \color{blue}

\noindent 246: (dims,levels) = $(3 , 
3;6,
30
)$,
irreps = $3_{3}^{1,0}
\hskip -1.5pt \otimes \hskip -1.5pt
1_{2}^{1,0}\oplus
3_{5}^{3}
\hskip -1.5pt \otimes \hskip -1.5pt
1_{3}^{1,0}
\hskip -1.5pt \otimes \hskip -1.5pt
1_{2}^{1,0}$,
pord$(\rho_\text{isum}(\mathfrak{t})) = 15$,

\vskip 0.7ex
\hangindent=5.5em \hangafter=1
{\white .}\hskip 1em $\rho_\text{isum}(\mathfrak{t})$ =
 $( \frac{1}{2},
\frac{1}{6},
\frac{5}{6} )
\oplus
( \frac{5}{6},
\frac{7}{30},
\frac{13}{30} )
$,

\vskip 0.7ex
\hangindent=5.5em \hangafter=1
{\white .}\hskip 1em $\rho_\text{isum}(\mathfrak{s})$ =
($\frac{1}{3}$,
$\frac{2}{3}$,
$\frac{2}{3}$;
$\frac{1}{3}$,
$-\frac{2}{3}$;
$\frac{1}{3}$)
 $\oplus$
($\sqrt{\frac{1}{5}}$,
$-\sqrt{\frac{2}{5}}$,
$-\sqrt{\frac{2}{5}}$;
$\frac{5-\sqrt{5}}{10}$,
$-\frac{5+\sqrt{5}}{10}$;
$\frac{5-\sqrt{5}}{10}$)

Pass. 

 \ \color{black}

\noindent 247: (dims,levels) = $(3 , 
3;8,
8
)$,
irreps = $3_{8}^{1,0}\oplus
3_{8}^{3,0}$,
pord$(\rho_\text{isum}(\mathfrak{t})) = 8$,

\vskip 0.7ex
\hangindent=5.5em \hangafter=1
{\white .}\hskip 1em $\rho_\text{isum}(\mathfrak{t})$ =
 $( 0,
\frac{1}{8},
\frac{5}{8} )
\oplus
( 0,
\frac{3}{8},
\frac{7}{8} )
$,

\vskip 0.7ex
\hangindent=5.5em \hangafter=1
{\white .}\hskip 1em $\rho_\text{isum}(\mathfrak{s})$ =
$\mathrm{i}$($0$,
$\sqrt{\frac{1}{2}}$,
$\sqrt{\frac{1}{2}}$;\ \ 
$-\frac{1}{2}$,
$\frac{1}{2}$;\ \ 
$-\frac{1}{2}$)
 $\oplus$
$\mathrm{i}$($0$,
$\sqrt{\frac{1}{2}}$,
$\sqrt{\frac{1}{2}}$;\ \ 
$\frac{1}{2}$,
$-\frac{1}{2}$;\ \ 
$\frac{1}{2}$)

Fail:
all rows of $U \rho(\mathfrak s) U^\dagger$
 contain zero for any block-diagonal $U$. Prop. B.5 (4) eqn. (B.27)

 \ \color{black}

\noindent 248: (dims,levels) = $(3 , 
3;8,
8
)$,
irreps = $3_{8}^{1,0}\oplus
3_{8}^{3,3}$,
pord$(\rho_\text{isum}(\mathfrak{t})) = 8$,

\vskip 0.7ex
\hangindent=5.5em \hangafter=1
{\white .}\hskip 1em $\rho_\text{isum}(\mathfrak{t})$ =
 $( 0,
\frac{1}{8},
\frac{5}{8} )
\oplus
( \frac{1}{4},
\frac{1}{8},
\frac{5}{8} )
$,

\vskip 0.7ex
\hangindent=5.5em \hangafter=1
{\white .}\hskip 1em $\rho_\text{isum}(\mathfrak{s})$ =
$\mathrm{i}$($0$,
$\sqrt{\frac{1}{2}}$,
$\sqrt{\frac{1}{2}}$;\ \ 
$-\frac{1}{2}$,
$\frac{1}{2}$;\ \ 
$-\frac{1}{2}$)
 $\oplus$
($0$,
$\sqrt{\frac{1}{2}}$,
$\sqrt{\frac{1}{2}}$;
$-\frac{1}{2}$,
$\frac{1}{2}$;
$-\frac{1}{2}$)

Fail:
number of self dual objects $|$Tr($\rho(\mathfrak s^2)$)$|$ = 0. Prop. B.4 (1)\
 eqn. (B.16)

 \ \color{black}

\noindent 249: (dims,levels) = $(3 , 
3;8,
8
)$,
irreps = $3_{8}^{1,0}\oplus
3_{8}^{1,6}$,
pord$(\rho_\text{isum}(\mathfrak{t})) = 8$,

\vskip 0.7ex
\hangindent=5.5em \hangafter=1
{\white .}\hskip 1em $\rho_\text{isum}(\mathfrak{t})$ =
 $( 0,
\frac{1}{8},
\frac{5}{8} )
\oplus
( \frac{1}{2},
\frac{1}{8},
\frac{5}{8} )
$,

\vskip 0.7ex
\hangindent=5.5em \hangafter=1
{\white .}\hskip 1em $\rho_\text{isum}(\mathfrak{s})$ =
$\mathrm{i}$($0$,
$\sqrt{\frac{1}{2}}$,
$\sqrt{\frac{1}{2}}$;\ \ 
$-\frac{1}{2}$,
$\frac{1}{2}$;\ \ 
$-\frac{1}{2}$)
 $\oplus$
$\mathrm{i}$($0$,
$\sqrt{\frac{1}{2}}$,
$\sqrt{\frac{1}{2}}$;\ \ 
$\frac{1}{2}$,
$-\frac{1}{2}$;\ \ 
$\frac{1}{2}$)

Fail:
Integral: $D_{\rho}(\sigma)_{\theta} \propto $ id,
 for all $\sigma$ and all $\theta$-eigenspaces that can contain unit. Prop. B.5 (6)

 \ \color{black}

\noindent 250: (dims,levels) = $(3 , 
3;8,
8
)$,
irreps = $3_{8}^{1,0}\oplus
3_{8}^{3,9}$,
pord$(\rho_\text{isum}(\mathfrak{t})) = 8$,

\vskip 0.7ex
\hangindent=5.5em \hangafter=1
{\white .}\hskip 1em $\rho_\text{isum}(\mathfrak{t})$ =
 $( 0,
\frac{1}{8},
\frac{5}{8} )
\oplus
( \frac{3}{4},
\frac{1}{8},
\frac{5}{8} )
$,

\vskip 0.7ex
\hangindent=5.5em \hangafter=1
{\white .}\hskip 1em $\rho_\text{isum}(\mathfrak{s})$ =
$\mathrm{i}$($0$,
$\sqrt{\frac{1}{2}}$,
$\sqrt{\frac{1}{2}}$;\ \ 
$-\frac{1}{2}$,
$\frac{1}{2}$;\ \ 
$-\frac{1}{2}$)
 $\oplus$
($0$,
$\sqrt{\frac{1}{2}}$,
$\sqrt{\frac{1}{2}}$;
$\frac{1}{2}$,
$-\frac{1}{2}$;
$\frac{1}{2}$)

Fail:
number of self dual objects $|$Tr($\rho(\mathfrak s^2)$)$|$ = 0. Prop. B.4 (1)\
 eqn. (B.16)

 \ \color{black}

\noindent 251: (dims,levels) = $(3 , 
3;8,
8
)$,
irreps = $3_{8}^{3,0}\oplus
3_{8}^{1,3}$,
pord$(\rho_\text{isum}(\mathfrak{t})) = 8$,

\vskip 0.7ex
\hangindent=5.5em \hangafter=1
{\white .}\hskip 1em $\rho_\text{isum}(\mathfrak{t})$ =
 $( 0,
\frac{3}{8},
\frac{7}{8} )
\oplus
( \frac{1}{4},
\frac{3}{8},
\frac{7}{8} )
$,

\vskip 0.7ex
\hangindent=5.5em \hangafter=1
{\white .}\hskip 1em $\rho_\text{isum}(\mathfrak{s})$ =
$\mathrm{i}$($0$,
$\sqrt{\frac{1}{2}}$,
$\sqrt{\frac{1}{2}}$;\ \ 
$\frac{1}{2}$,
$-\frac{1}{2}$;\ \ 
$\frac{1}{2}$)
 $\oplus$
($0$,
$\sqrt{\frac{1}{2}}$,
$\sqrt{\frac{1}{2}}$;
$\frac{1}{2}$,
$-\frac{1}{2}$;
$\frac{1}{2}$)

Fail:
number of self dual objects $|$Tr($\rho(\mathfrak s^2)$)$|$ = 0. Prop. B.4 (1)\
 eqn. (B.16)

 \ \color{black}

\noindent 252: (dims,levels) = $(3 , 
3;8,
8
)$,
irreps = $3_{8}^{3,0}\oplus
3_{8}^{3,6}$,
pord$(\rho_\text{isum}(\mathfrak{t})) = 8$,

\vskip 0.7ex
\hangindent=5.5em \hangafter=1
{\white .}\hskip 1em $\rho_\text{isum}(\mathfrak{t})$ =
 $( 0,
\frac{3}{8},
\frac{7}{8} )
\oplus
( \frac{1}{2},
\frac{3}{8},
\frac{7}{8} )
$,

\vskip 0.7ex
\hangindent=5.5em \hangafter=1
{\white .}\hskip 1em $\rho_\text{isum}(\mathfrak{s})$ =
$\mathrm{i}$($0$,
$\sqrt{\frac{1}{2}}$,
$\sqrt{\frac{1}{2}}$;\ \ 
$\frac{1}{2}$,
$-\frac{1}{2}$;\ \ 
$\frac{1}{2}$)
 $\oplus$
$\mathrm{i}$($0$,
$\sqrt{\frac{1}{2}}$,
$\sqrt{\frac{1}{2}}$;\ \ 
$-\frac{1}{2}$,
$\frac{1}{2}$;\ \ 
$-\frac{1}{2}$)

Fail:
Integral: $D_{\rho}(\sigma)_{\theta} \propto $ id,
 for all $\sigma$ and all $\theta$-eigenspaces that can contain unit. Prop. B.5 (6)

 \ \color{black}

\noindent 253: (dims,levels) = $(3 , 
3;8,
8
)$,
irreps = $3_{8}^{3,0}\oplus
3_{8}^{1,9}$,
pord$(\rho_\text{isum}(\mathfrak{t})) = 8$,

\vskip 0.7ex
\hangindent=5.5em \hangafter=1
{\white .}\hskip 1em $\rho_\text{isum}(\mathfrak{t})$ =
 $( 0,
\frac{3}{8},
\frac{7}{8} )
\oplus
( \frac{3}{4},
\frac{3}{8},
\frac{7}{8} )
$,

\vskip 0.7ex
\hangindent=5.5em \hangafter=1
{\white .}\hskip 1em $\rho_\text{isum}(\mathfrak{s})$ =
$\mathrm{i}$($0$,
$\sqrt{\frac{1}{2}}$,
$\sqrt{\frac{1}{2}}$;\ \ 
$\frac{1}{2}$,
$-\frac{1}{2}$;\ \ 
$\frac{1}{2}$)
 $\oplus$
($0$,
$\sqrt{\frac{1}{2}}$,
$\sqrt{\frac{1}{2}}$;
$-\frac{1}{2}$,
$\frac{1}{2}$;
$-\frac{1}{2}$)

Fail:
number of self dual objects $|$Tr($\rho(\mathfrak s^2)$)$|$ = 0. Prop. B.4 (1)\
 eqn. (B.16)

 \ \color{black}

\noindent 254: (dims,levels) = $(3 , 
3;8,
16
)$,
irreps = $3_{8}^{1,0}\oplus
3_{16}^{1,0}$,
pord$(\rho_\text{isum}(\mathfrak{t})) = 16$,

\vskip 0.7ex
\hangindent=5.5em \hangafter=1
{\white .}\hskip 1em $\rho_\text{isum}(\mathfrak{t})$ =
 $( 0,
\frac{1}{8},
\frac{5}{8} )
\oplus
( \frac{1}{8},
\frac{1}{16},
\frac{9}{16} )
$,

\vskip 0.7ex
\hangindent=5.5em \hangafter=1
{\white .}\hskip 1em $\rho_\text{isum}(\mathfrak{s})$ =
$\mathrm{i}$($0$,
$\sqrt{\frac{1}{2}}$,
$\sqrt{\frac{1}{2}}$;\ \ 
$-\frac{1}{2}$,
$\frac{1}{2}$;\ \ 
$-\frac{1}{2}$)
 $\oplus$
$\mathrm{i}$($0$,
$\sqrt{\frac{1}{2}}$,
$\sqrt{\frac{1}{2}}$;\ \ 
$-\frac{1}{2}$,
$\frac{1}{2}$;\ \ 
$-\frac{1}{2}$)

Fail:
$\sigma(\rho(\mathfrak s)_\mathrm{ndeg}) \neq
 (\rho(\mathfrak t)^a \rho(\mathfrak s) \rho(\mathfrak t)^b
 \rho(\mathfrak s) \rho(\mathfrak t)^a)_\mathrm{ndeg}$,
 $\sigma = a$ = 3. Prop. B.5 (3) eqn. (B.25)

 \ \color{black}

\noindent 255: (dims,levels) = $(3 , 
3;8,
16
)$,
irreps = $3_{8}^{1,0}\oplus
3_{16}^{7,3}$,
pord$(\rho_\text{isum}(\mathfrak{t})) = 16$,

\vskip 0.7ex
\hangindent=5.5em \hangafter=1
{\white .}\hskip 1em $\rho_\text{isum}(\mathfrak{t})$ =
 $( 0,
\frac{1}{8},
\frac{5}{8} )
\oplus
( \frac{1}{8},
\frac{3}{16},
\frac{11}{16} )
$,

\vskip 0.7ex
\hangindent=5.5em \hangafter=1
{\white .}\hskip 1em $\rho_\text{isum}(\mathfrak{s})$ =
$\mathrm{i}$($0$,
$\sqrt{\frac{1}{2}}$,
$\sqrt{\frac{1}{2}}$;\ \ 
$-\frac{1}{2}$,
$\frac{1}{2}$;\ \ 
$-\frac{1}{2}$)
 $\oplus$
($0$,
$\sqrt{\frac{1}{2}}$,
$\sqrt{\frac{1}{2}}$;
$-\frac{1}{2}$,
$\frac{1}{2}$;
$-\frac{1}{2}$)

Fail:
number of self dual objects $|$Tr($\rho(\mathfrak s^2)$)$|$ = 0. Prop. B.4 (1)\
 eqn. (B.16)

 \ \color{black}

\noindent 256: (dims,levels) = $(3 , 
3;8,
16
)$,
irreps = $3_{8}^{1,0}\oplus
3_{16}^{5,6}$,
pord$(\rho_\text{isum}(\mathfrak{t})) = 16$,

\vskip 0.7ex
\hangindent=5.5em \hangafter=1
{\white .}\hskip 1em $\rho_\text{isum}(\mathfrak{t})$ =
 $( 0,
\frac{1}{8},
\frac{5}{8} )
\oplus
( \frac{1}{8},
\frac{5}{16},
\frac{13}{16} )
$,

\vskip 0.7ex
\hangindent=5.5em \hangafter=1
{\white .}\hskip 1em $\rho_\text{isum}(\mathfrak{s})$ =
$\mathrm{i}$($0$,
$\sqrt{\frac{1}{2}}$,
$\sqrt{\frac{1}{2}}$;\ \ 
$-\frac{1}{2}$,
$\frac{1}{2}$;\ \ 
$-\frac{1}{2}$)
 $\oplus$
$\mathrm{i}$($0$,
$\sqrt{\frac{1}{2}}$,
$\sqrt{\frac{1}{2}}$;\ \ 
$\frac{1}{2}$,
$-\frac{1}{2}$;\ \ 
$\frac{1}{2}$)

Fail:
$\sigma(\rho(\mathfrak s)_\mathrm{ndeg}) \neq
 (\rho(\mathfrak t)^a \rho(\mathfrak s) \rho(\mathfrak t)^b
 \rho(\mathfrak s) \rho(\mathfrak t)^a)_\mathrm{ndeg}$,
 $\sigma = a$ = 3. Prop. B.5 (3) eqn. (B.25)

 \ \color{black}

\noindent 257: (dims,levels) = $(3 , 
3;8,
16
)$,
irreps = $3_{8}^{1,0}\oplus
3_{16}^{3,9}$,
pord$(\rho_\text{isum}(\mathfrak{t})) = 16$,

\vskip 0.7ex
\hangindent=5.5em \hangafter=1
{\white .}\hskip 1em $\rho_\text{isum}(\mathfrak{t})$ =
 $( 0,
\frac{1}{8},
\frac{5}{8} )
\oplus
( \frac{1}{8},
\frac{7}{16},
\frac{15}{16} )
$,

\vskip 0.7ex
\hangindent=5.5em \hangafter=1
{\white .}\hskip 1em $\rho_\text{isum}(\mathfrak{s})$ =
$\mathrm{i}$($0$,
$\sqrt{\frac{1}{2}}$,
$\sqrt{\frac{1}{2}}$;\ \ 
$-\frac{1}{2}$,
$\frac{1}{2}$;\ \ 
$-\frac{1}{2}$)
 $\oplus$
($0$,
$\sqrt{\frac{1}{2}}$,
$\sqrt{\frac{1}{2}}$;
$\frac{1}{2}$,
$-\frac{1}{2}$;
$\frac{1}{2}$)

Fail:
number of self dual objects $|$Tr($\rho(\mathfrak s^2)$)$|$ = 0. Prop. B.4 (1)\
 eqn. (B.16)

 \ \color{black}

\noindent 258: (dims,levels) = $(3 , 
3;8,
16
)$,
irreps = $3_{8}^{1,0}\oplus
3_{16}^{1,6}$,
pord$(\rho_\text{isum}(\mathfrak{t})) = 16$,

\vskip 0.7ex
\hangindent=5.5em \hangafter=1
{\white .}\hskip 1em $\rho_\text{isum}(\mathfrak{t})$ =
 $( 0,
\frac{1}{8},
\frac{5}{8} )
\oplus
( \frac{5}{8},
\frac{1}{16},
\frac{9}{16} )
$,

\vskip 0.7ex
\hangindent=5.5em \hangafter=1
{\white .}\hskip 1em $\rho_\text{isum}(\mathfrak{s})$ =
$\mathrm{i}$($0$,
$\sqrt{\frac{1}{2}}$,
$\sqrt{\frac{1}{2}}$;\ \ 
$-\frac{1}{2}$,
$\frac{1}{2}$;\ \ 
$-\frac{1}{2}$)
 $\oplus$
$\mathrm{i}$($0$,
$\sqrt{\frac{1}{2}}$,
$\sqrt{\frac{1}{2}}$;\ \ 
$\frac{1}{2}$,
$-\frac{1}{2}$;\ \ 
$\frac{1}{2}$)

Fail:
$\sigma(\rho(\mathfrak s)_\mathrm{ndeg}) \neq
 (\rho(\mathfrak t)^a \rho(\mathfrak s) \rho(\mathfrak t)^b
 \rho(\mathfrak s) \rho(\mathfrak t)^a)_\mathrm{ndeg}$,
 $\sigma = a$ = 3. Prop. B.5 (3) eqn. (B.25)

 \ \color{black}

\noindent 259: (dims,levels) = $(3 , 
3;8,
16
)$,
irreps = $3_{8}^{1,0}\oplus
3_{16}^{7,9}$,
pord$(\rho_\text{isum}(\mathfrak{t})) = 16$,

\vskip 0.7ex
\hangindent=5.5em \hangafter=1
{\white .}\hskip 1em $\rho_\text{isum}(\mathfrak{t})$ =
 $( 0,
\frac{1}{8},
\frac{5}{8} )
\oplus
( \frac{5}{8},
\frac{3}{16},
\frac{11}{16} )
$,

\vskip 0.7ex
\hangindent=5.5em \hangafter=1
{\white .}\hskip 1em $\rho_\text{isum}(\mathfrak{s})$ =
$\mathrm{i}$($0$,
$\sqrt{\frac{1}{2}}$,
$\sqrt{\frac{1}{2}}$;\ \ 
$-\frac{1}{2}$,
$\frac{1}{2}$;\ \ 
$-\frac{1}{2}$)
 $\oplus$
($0$,
$\sqrt{\frac{1}{2}}$,
$\sqrt{\frac{1}{2}}$;
$\frac{1}{2}$,
$-\frac{1}{2}$;
$\frac{1}{2}$)

Fail:
number of self dual objects $|$Tr($\rho(\mathfrak s^2)$)$|$ = 0. Prop. B.4 (1)\
 eqn. (B.16)

 \ \color{black}

\noindent 260: (dims,levels) = $(3 , 
3;8,
16
)$,
irreps = $3_{8}^{1,0}\oplus
3_{16}^{5,0}$,
pord$(\rho_\text{isum}(\mathfrak{t})) = 16$,

\vskip 0.7ex
\hangindent=5.5em \hangafter=1
{\white .}\hskip 1em $\rho_\text{isum}(\mathfrak{t})$ =
 $( 0,
\frac{1}{8},
\frac{5}{8} )
\oplus
( \frac{5}{8},
\frac{5}{16},
\frac{13}{16} )
$,

\vskip 0.7ex
\hangindent=5.5em \hangafter=1
{\white .}\hskip 1em $\rho_\text{isum}(\mathfrak{s})$ =
$\mathrm{i}$($0$,
$\sqrt{\frac{1}{2}}$,
$\sqrt{\frac{1}{2}}$;\ \ 
$-\frac{1}{2}$,
$\frac{1}{2}$;\ \ 
$-\frac{1}{2}$)
 $\oplus$
$\mathrm{i}$($0$,
$\sqrt{\frac{1}{2}}$,
$\sqrt{\frac{1}{2}}$;\ \ 
$-\frac{1}{2}$,
$\frac{1}{2}$;\ \ 
$-\frac{1}{2}$)

Fail:
$\sigma(\rho(\mathfrak s)_\mathrm{ndeg}) \neq
 (\rho(\mathfrak t)^a \rho(\mathfrak s) \rho(\mathfrak t)^b
 \rho(\mathfrak s) \rho(\mathfrak t)^a)_\mathrm{ndeg}$,
 $\sigma = a$ = 3. Prop. B.5 (3) eqn. (B.25)

 \ \color{black}

\noindent 261: (dims,levels) = $(3 , 
3;8,
16
)$,
irreps = $3_{8}^{1,0}\oplus
3_{16}^{3,3}$,
pord$(\rho_\text{isum}(\mathfrak{t})) = 16$,

\vskip 0.7ex
\hangindent=5.5em \hangafter=1
{\white .}\hskip 1em $\rho_\text{isum}(\mathfrak{t})$ =
 $( 0,
\frac{1}{8},
\frac{5}{8} )
\oplus
( \frac{5}{8},
\frac{7}{16},
\frac{15}{16} )
$,

\vskip 0.7ex
\hangindent=5.5em \hangafter=1
{\white .}\hskip 1em $\rho_\text{isum}(\mathfrak{s})$ =
$\mathrm{i}$($0$,
$\sqrt{\frac{1}{2}}$,
$\sqrt{\frac{1}{2}}$;\ \ 
$-\frac{1}{2}$,
$\frac{1}{2}$;\ \ 
$-\frac{1}{2}$)
 $\oplus$
($0$,
$\sqrt{\frac{1}{2}}$,
$\sqrt{\frac{1}{2}}$;
$-\frac{1}{2}$,
$\frac{1}{2}$;
$-\frac{1}{2}$)

Fail:
number of self dual objects $|$Tr($\rho(\mathfrak s^2)$)$|$ = 0. Prop. B.4 (1)\
 eqn. (B.16)

 \ \color{black}

\noindent 262: (dims,levels) = $(3 , 
3;8,
16
)$,
irreps = $3_{8}^{3,0}\oplus
3_{16}^{5,9}$,
pord$(\rho_\text{isum}(\mathfrak{t})) = 16$,

\vskip 0.7ex
\hangindent=5.5em \hangafter=1
{\white .}\hskip 1em $\rho_\text{isum}(\mathfrak{t})$ =
 $( 0,
\frac{3}{8},
\frac{7}{8} )
\oplus
( \frac{3}{8},
\frac{1}{16},
\frac{9}{16} )
$,

\vskip 0.7ex
\hangindent=5.5em \hangafter=1
{\white .}\hskip 1em $\rho_\text{isum}(\mathfrak{s})$ =
$\mathrm{i}$($0$,
$\sqrt{\frac{1}{2}}$,
$\sqrt{\frac{1}{2}}$;\ \ 
$\frac{1}{2}$,
$-\frac{1}{2}$;\ \ 
$\frac{1}{2}$)
 $\oplus$
($0$,
$\sqrt{\frac{1}{2}}$,
$\sqrt{\frac{1}{2}}$;
$-\frac{1}{2}$,
$\frac{1}{2}$;
$-\frac{1}{2}$)

Fail:
number of self dual objects $|$Tr($\rho(\mathfrak s^2)$)$|$ = 0. Prop. B.4 (1)\
 eqn. (B.16)

 \ \color{black}

\noindent 263: (dims,levels) = $(3 , 
3;8,
16
)$,
irreps = $3_{8}^{3,0}\oplus
3_{16}^{3,0}$,
pord$(\rho_\text{isum}(\mathfrak{t})) = 16$,

\vskip 0.7ex
\hangindent=5.5em \hangafter=1
{\white .}\hskip 1em $\rho_\text{isum}(\mathfrak{t})$ =
 $( 0,
\frac{3}{8},
\frac{7}{8} )
\oplus
( \frac{3}{8},
\frac{3}{16},
\frac{11}{16} )
$,

\vskip 0.7ex
\hangindent=5.5em \hangafter=1
{\white .}\hskip 1em $\rho_\text{isum}(\mathfrak{s})$ =
$\mathrm{i}$($0$,
$\sqrt{\frac{1}{2}}$,
$\sqrt{\frac{1}{2}}$;\ \ 
$\frac{1}{2}$,
$-\frac{1}{2}$;\ \ 
$\frac{1}{2}$)
 $\oplus$
$\mathrm{i}$($0$,
$\sqrt{\frac{1}{2}}$,
$\sqrt{\frac{1}{2}}$;\ \ 
$\frac{1}{2}$,
$-\frac{1}{2}$;\ \ 
$\frac{1}{2}$)

Fail:
$\sigma(\rho(\mathfrak s)_\mathrm{ndeg}) \neq
 (\rho(\mathfrak t)^a \rho(\mathfrak s) \rho(\mathfrak t)^b
 \rho(\mathfrak s) \rho(\mathfrak t)^a)_\mathrm{ndeg}$,
 $\sigma = a$ = 3. Prop. B.5 (3) eqn. (B.25)

 \ \color{black}

\noindent 264: (dims,levels) = $(3 , 
3;8,
16
)$,
irreps = $3_{8}^{3,0}\oplus
3_{16}^{1,3}$,
pord$(\rho_\text{isum}(\mathfrak{t})) = 16$,

\vskip 0.7ex
\hangindent=5.5em \hangafter=1
{\white .}\hskip 1em $\rho_\text{isum}(\mathfrak{t})$ =
 $( 0,
\frac{3}{8},
\frac{7}{8} )
\oplus
( \frac{3}{8},
\frac{5}{16},
\frac{13}{16} )
$,

\vskip 0.7ex
\hangindent=5.5em \hangafter=1
{\white .}\hskip 1em $\rho_\text{isum}(\mathfrak{s})$ =
$\mathrm{i}$($0$,
$\sqrt{\frac{1}{2}}$,
$\sqrt{\frac{1}{2}}$;\ \ 
$\frac{1}{2}$,
$-\frac{1}{2}$;\ \ 
$\frac{1}{2}$)
 $\oplus$
($0$,
$\sqrt{\frac{1}{2}}$,
$\sqrt{\frac{1}{2}}$;
$\frac{1}{2}$,
$-\frac{1}{2}$;
$\frac{1}{2}$)

Fail:
number of self dual objects $|$Tr($\rho(\mathfrak s^2)$)$|$ = 0. Prop. B.4 (1)\
 eqn. (B.16)

 \ \color{black}

\noindent 265: (dims,levels) = $(3 , 
3;8,
16
)$,
irreps = $3_{8}^{3,0}\oplus
3_{16}^{7,6}$,
pord$(\rho_\text{isum}(\mathfrak{t})) = 16$,

\vskip 0.7ex
\hangindent=5.5em \hangafter=1
{\white .}\hskip 1em $\rho_\text{isum}(\mathfrak{t})$ =
 $( 0,
\frac{3}{8},
\frac{7}{8} )
\oplus
( \frac{3}{8},
\frac{7}{16},
\frac{15}{16} )
$,

\vskip 0.7ex
\hangindent=5.5em \hangafter=1
{\white .}\hskip 1em $\rho_\text{isum}(\mathfrak{s})$ =
$\mathrm{i}$($0$,
$\sqrt{\frac{1}{2}}$,
$\sqrt{\frac{1}{2}}$;\ \ 
$\frac{1}{2}$,
$-\frac{1}{2}$;\ \ 
$\frac{1}{2}$)
 $\oplus$
$\mathrm{i}$($0$,
$\sqrt{\frac{1}{2}}$,
$\sqrt{\frac{1}{2}}$;\ \ 
$-\frac{1}{2}$,
$\frac{1}{2}$;\ \ 
$-\frac{1}{2}$)

Fail:
$\sigma(\rho(\mathfrak s)_\mathrm{ndeg}) \neq
 (\rho(\mathfrak t)^a \rho(\mathfrak s) \rho(\mathfrak t)^b
 \rho(\mathfrak s) \rho(\mathfrak t)^a)_\mathrm{ndeg}$,
 $\sigma = a$ = 3. Prop. B.5 (3) eqn. (B.25)

 \ \color{black}

\noindent 266: (dims,levels) = $(3 , 
3;8,
16
)$,
irreps = $3_{8}^{3,0}\oplus
3_{16}^{5,3}$,
pord$(\rho_\text{isum}(\mathfrak{t})) = 16$,

\vskip 0.7ex
\hangindent=5.5em \hangafter=1
{\white .}\hskip 1em $\rho_\text{isum}(\mathfrak{t})$ =
 $( 0,
\frac{3}{8},
\frac{7}{8} )
\oplus
( \frac{7}{8},
\frac{1}{16},
\frac{9}{16} )
$,

\vskip 0.7ex
\hangindent=5.5em \hangafter=1
{\white .}\hskip 1em $\rho_\text{isum}(\mathfrak{s})$ =
$\mathrm{i}$($0$,
$\sqrt{\frac{1}{2}}$,
$\sqrt{\frac{1}{2}}$;\ \ 
$\frac{1}{2}$,
$-\frac{1}{2}$;\ \ 
$\frac{1}{2}$)
 $\oplus$
($0$,
$\sqrt{\frac{1}{2}}$,
$\sqrt{\frac{1}{2}}$;
$\frac{1}{2}$,
$-\frac{1}{2}$;
$\frac{1}{2}$)

Fail:
number of self dual objects $|$Tr($\rho(\mathfrak s^2)$)$|$ = 0. Prop. B.4 (1)\
 eqn. (B.16)

 \ \color{black}

\noindent 267: (dims,levels) = $(3 , 
3;8,
16
)$,
irreps = $3_{8}^{3,0}\oplus
3_{16}^{3,6}$,
pord$(\rho_\text{isum}(\mathfrak{t})) = 16$,

\vskip 0.7ex
\hangindent=5.5em \hangafter=1
{\white .}\hskip 1em $\rho_\text{isum}(\mathfrak{t})$ =
 $( 0,
\frac{3}{8},
\frac{7}{8} )
\oplus
( \frac{7}{8},
\frac{3}{16},
\frac{11}{16} )
$,

\vskip 0.7ex
\hangindent=5.5em \hangafter=1
{\white .}\hskip 1em $\rho_\text{isum}(\mathfrak{s})$ =
$\mathrm{i}$($0$,
$\sqrt{\frac{1}{2}}$,
$\sqrt{\frac{1}{2}}$;\ \ 
$\frac{1}{2}$,
$-\frac{1}{2}$;\ \ 
$\frac{1}{2}$)
 $\oplus$
$\mathrm{i}$($0$,
$\sqrt{\frac{1}{2}}$,
$\sqrt{\frac{1}{2}}$;\ \ 
$-\frac{1}{2}$,
$\frac{1}{2}$;\ \ 
$-\frac{1}{2}$)

Fail:
$\sigma(\rho(\mathfrak s)_\mathrm{ndeg}) \neq
 (\rho(\mathfrak t)^a \rho(\mathfrak s) \rho(\mathfrak t)^b
 \rho(\mathfrak s) \rho(\mathfrak t)^a)_\mathrm{ndeg}$,
 $\sigma = a$ = 3. Prop. B.5 (3) eqn. (B.25)

 \ \color{black}

\noindent 268: (dims,levels) = $(3 , 
3;8,
16
)$,
irreps = $3_{8}^{3,0}\oplus
3_{16}^{1,9}$,
pord$(\rho_\text{isum}(\mathfrak{t})) = 16$,

\vskip 0.7ex
\hangindent=5.5em \hangafter=1
{\white .}\hskip 1em $\rho_\text{isum}(\mathfrak{t})$ =
 $( 0,
\frac{3}{8},
\frac{7}{8} )
\oplus
( \frac{7}{8},
\frac{5}{16},
\frac{13}{16} )
$,

\vskip 0.7ex
\hangindent=5.5em \hangafter=1
{\white .}\hskip 1em $\rho_\text{isum}(\mathfrak{s})$ =
$\mathrm{i}$($0$,
$\sqrt{\frac{1}{2}}$,
$\sqrt{\frac{1}{2}}$;\ \ 
$\frac{1}{2}$,
$-\frac{1}{2}$;\ \ 
$\frac{1}{2}$)
 $\oplus$
($0$,
$\sqrt{\frac{1}{2}}$,
$\sqrt{\frac{1}{2}}$;
$-\frac{1}{2}$,
$\frac{1}{2}$;
$-\frac{1}{2}$)

Fail:
number of self dual objects $|$Tr($\rho(\mathfrak s^2)$)$|$ = 0. Prop. B.4 (1)\
 eqn. (B.16)

 \ \color{black}

\noindent 269: (dims,levels) = $(3 , 
3;8,
16
)$,
irreps = $3_{8}^{3,0}\oplus
3_{16}^{7,0}$,
pord$(\rho_\text{isum}(\mathfrak{t})) = 16$,

\vskip 0.7ex
\hangindent=5.5em \hangafter=1
{\white .}\hskip 1em $\rho_\text{isum}(\mathfrak{t})$ =
 $( 0,
\frac{3}{8},
\frac{7}{8} )
\oplus
( \frac{7}{8},
\frac{7}{16},
\frac{15}{16} )
$,

\vskip 0.7ex
\hangindent=5.5em \hangafter=1
{\white .}\hskip 1em $\rho_\text{isum}(\mathfrak{s})$ =
$\mathrm{i}$($0$,
$\sqrt{\frac{1}{2}}$,
$\sqrt{\frac{1}{2}}$;\ \ 
$\frac{1}{2}$,
$-\frac{1}{2}$;\ \ 
$\frac{1}{2}$)
 $\oplus$
$\mathrm{i}$($0$,
$\sqrt{\frac{1}{2}}$,
$\sqrt{\frac{1}{2}}$;\ \ 
$\frac{1}{2}$,
$-\frac{1}{2}$;\ \ 
$\frac{1}{2}$)

Fail:
$\sigma(\rho(\mathfrak s)_\mathrm{ndeg}) \neq
 (\rho(\mathfrak t)^a \rho(\mathfrak s) \rho(\mathfrak t)^b
 \rho(\mathfrak s) \rho(\mathfrak t)^a)_\mathrm{ndeg}$,
 $\sigma = a$ = 3. Prop. B.5 (3) eqn. (B.25)

 \ \color{black}

\noindent 270: (dims,levels) = $(3 , 
3;10,
10
)$,
irreps = $3_{5}^{3}
\hskip -1.5pt \otimes \hskip -1.5pt
1_{2}^{1,0}\oplus
3_{5}^{1}
\hskip -1.5pt \otimes \hskip -1.5pt
1_{2}^{1,0}$,
pord$(\rho_\text{isum}(\mathfrak{t})) = 5$,

\vskip 0.7ex
\hangindent=5.5em \hangafter=1
{\white .}\hskip 1em $\rho_\text{isum}(\mathfrak{t})$ =
 $( \frac{1}{2},
\frac{1}{10},
\frac{9}{10} )
\oplus
( \frac{1}{2},
\frac{3}{10},
\frac{7}{10} )
$,

\vskip 0.7ex
\hangindent=5.5em \hangafter=1
{\white .}\hskip 1em $\rho_\text{isum}(\mathfrak{s})$ =
($\sqrt{\frac{1}{5}}$,
$-\sqrt{\frac{2}{5}}$,
$-\sqrt{\frac{2}{5}}$;
$\frac{5-\sqrt{5}}{10}$,
$-\frac{5+\sqrt{5}}{10}$;
$\frac{5-\sqrt{5}}{10}$)
 $\oplus$
($-\sqrt{\frac{1}{5}}$,
$-\sqrt{\frac{2}{5}}$,
$-\sqrt{\frac{2}{5}}$;
$\frac{5+\sqrt{5}}{10}$,
$\frac{-5+\sqrt{5}}{10}$;
$\frac{5+\sqrt{5}}{10}$)

Fail:
Integral: $D_{\rho}(\sigma)_{\theta} \propto $ id,
 for all $\sigma$ and all $\theta$-eigenspaces that can contain unit. Prop. B.5 (6)

 \ \color{black}

\noindent 271: (dims,levels) = $(3 , 
3;12,
12
)$,
irreps = $3_{3}^{1,0}
\hskip -1.5pt \otimes \hskip -1.5pt
1_{4}^{1,0}\oplus
3_{4}^{1,0}
\hskip -1.5pt \otimes \hskip -1.5pt
1_{3}^{1,0}$,
pord$(\rho_\text{isum}(\mathfrak{t})) = 12$,

\vskip 0.7ex
\hangindent=5.5em \hangafter=1
{\white .}\hskip 1em $\rho_\text{isum}(\mathfrak{t})$ =
 $( \frac{1}{4},
\frac{7}{12},
\frac{11}{12} )
\oplus
( \frac{1}{3},
\frac{1}{12},
\frac{7}{12} )
$,

\vskip 0.7ex
\hangindent=5.5em \hangafter=1
{\white .}\hskip 1em $\rho_\text{isum}(\mathfrak{s})$ =
$\mathrm{i}$($-\frac{1}{3}$,
$\frac{2}{3}$,
$\frac{2}{3}$;\ \ 
$-\frac{1}{3}$,
$\frac{2}{3}$;\ \ 
$-\frac{1}{3}$)
 $\oplus$
($0$,
$\sqrt{\frac{1}{2}}$,
$\sqrt{\frac{1}{2}}$;
$-\frac{1}{2}$,
$\frac{1}{2}$;
$-\frac{1}{2}$)

Fail:
number of self dual objects $|$Tr($\rho(\mathfrak s^2)$)$|$ = 0. Prop. B.4 (1)\
 eqn. (B.16)

 \ \color{black}

\noindent 272: (dims,levels) = $(3 , 
3;12,
12
)$,
irreps = $3_{3}^{1,0}
\hskip -1.5pt \otimes \hskip -1.5pt
1_{4}^{1,0}\oplus
3_{4}^{1,3}
\hskip -1.5pt \otimes \hskip -1.5pt
1_{3}^{1,0}$,
pord$(\rho_\text{isum}(\mathfrak{t})) = 12$,

\vskip 0.7ex
\hangindent=5.5em \hangafter=1
{\white .}\hskip 1em $\rho_\text{isum}(\mathfrak{t})$ =
 $( \frac{1}{4},
\frac{7}{12},
\frac{11}{12} )
\oplus
( \frac{1}{3},
\frac{5}{6},
\frac{7}{12} )
$,

\vskip 0.7ex
\hangindent=5.5em \hangafter=1
{\white .}\hskip 1em $\rho_\text{isum}(\mathfrak{s})$ =
$\mathrm{i}$($-\frac{1}{3}$,
$\frac{2}{3}$,
$\frac{2}{3}$;\ \ 
$-\frac{1}{3}$,
$\frac{2}{3}$;\ \ 
$-\frac{1}{3}$)
 $\oplus$
$\mathrm{i}$($-\frac{1}{2}$,
$\frac{1}{2}$,
$\sqrt{\frac{1}{2}}$;\ \ 
$-\frac{1}{2}$,
$\sqrt{\frac{1}{2}}$;\ \ 
$0$)

Fail:
Integral: $D_{\rho}(\sigma)_{\theta} \propto $ id,
 for all $\sigma$ and all $\theta$-eigenspaces that can contain unit. Prop. B.5 (6)

 \ \color{black}

\noindent 273: (dims,levels) = $(3 , 
3;12,
12
)$,
irreps = $3_{3}^{1,0}
\hskip -1.5pt \otimes \hskip -1.5pt
1_{4}^{1,0}\oplus
3_{4}^{1,3}
\hskip -1.5pt \otimes \hskip -1.5pt
1_{3}^{1,4}$,
pord$(\rho_\text{isum}(\mathfrak{t})) = 12$,

\vskip 0.7ex
\hangindent=5.5em \hangafter=1
{\white .}\hskip 1em $\rho_\text{isum}(\mathfrak{t})$ =
 $( \frac{1}{4},
\frac{7}{12},
\frac{11}{12} )
\oplus
( \frac{2}{3},
\frac{1}{6},
\frac{11}{12} )
$,

\vskip 0.7ex
\hangindent=5.5em \hangafter=1
{\white .}\hskip 1em $\rho_\text{isum}(\mathfrak{s})$ =
$\mathrm{i}$($-\frac{1}{3}$,
$\frac{2}{3}$,
$\frac{2}{3}$;\ \ 
$-\frac{1}{3}$,
$\frac{2}{3}$;\ \ 
$-\frac{1}{3}$)
 $\oplus$
$\mathrm{i}$($-\frac{1}{2}$,
$\frac{1}{2}$,
$\sqrt{\frac{1}{2}}$;\ \ 
$-\frac{1}{2}$,
$\sqrt{\frac{1}{2}}$;\ \ 
$0$)

Fail:
Integral: $D_{\rho}(\sigma)_{\theta} \propto $ id,
 for all $\sigma$ and all $\theta$-eigenspaces that can contain unit. Prop. B.5 (6)

 \ \color{black}

\noindent 274: (dims,levels) = $(3 , 
3;12,
12
)$,
irreps = $3_{3}^{1,0}
\hskip -1.5pt \otimes \hskip -1.5pt
1_{4}^{1,0}\oplus
3_{4}^{1,0}
\hskip -1.5pt \otimes \hskip -1.5pt
1_{3}^{1,4}$,
pord$(\rho_\text{isum}(\mathfrak{t})) = 12$,

\vskip 0.7ex
\hangindent=5.5em \hangafter=1
{\white .}\hskip 1em $\rho_\text{isum}(\mathfrak{t})$ =
 $( \frac{1}{4},
\frac{7}{12},
\frac{11}{12} )
\oplus
( \frac{2}{3},
\frac{5}{12},
\frac{11}{12} )
$,

\vskip 0.7ex
\hangindent=5.5em \hangafter=1
{\white .}\hskip 1em $\rho_\text{isum}(\mathfrak{s})$ =
$\mathrm{i}$($-\frac{1}{3}$,
$\frac{2}{3}$,
$\frac{2}{3}$;\ \ 
$-\frac{1}{3}$,
$\frac{2}{3}$;\ \ 
$-\frac{1}{3}$)
 $\oplus$
($0$,
$\sqrt{\frac{1}{2}}$,
$\sqrt{\frac{1}{2}}$;
$-\frac{1}{2}$,
$\frac{1}{2}$;
$-\frac{1}{2}$)

Fail:
number of self dual objects $|$Tr($\rho(\mathfrak s^2)$)$|$ = 0. Prop. B.4 (1)\
 eqn. (B.16)

 \ \color{black}

\noindent 275: (dims,levels) = $(3 , 
3;12,
12
)$,
irreps = $3_{3}^{1,0}
\hskip -1.5pt \otimes \hskip -1.5pt
1_{4}^{1,0}\oplus
3_{4}^{1,6}
\hskip -1.5pt \otimes \hskip -1.5pt
1_{3}^{1,0}$,
pord$(\rho_\text{isum}(\mathfrak{t})) = 12$,

\vskip 0.7ex
\hangindent=5.5em \hangafter=1
{\white .}\hskip 1em $\rho_\text{isum}(\mathfrak{t})$ =
 $( \frac{1}{4},
\frac{7}{12},
\frac{11}{12} )
\oplus
( \frac{5}{6},
\frac{1}{12},
\frac{7}{12} )
$,

\vskip 0.7ex
\hangindent=5.5em \hangafter=1
{\white .}\hskip 1em $\rho_\text{isum}(\mathfrak{s})$ =
$\mathrm{i}$($-\frac{1}{3}$,
$\frac{2}{3}$,
$\frac{2}{3}$;\ \ 
$-\frac{1}{3}$,
$\frac{2}{3}$;\ \ 
$-\frac{1}{3}$)
 $\oplus$
($0$,
$\sqrt{\frac{1}{2}}$,
$\sqrt{\frac{1}{2}}$;
$\frac{1}{2}$,
$-\frac{1}{2}$;
$\frac{1}{2}$)

Fail:
number of self dual objects $|$Tr($\rho(\mathfrak s^2)$)$|$ = 0. Prop. B.4 (1)\
 eqn. (B.16)

 \ \color{black}

\noindent 276: (dims,levels) = $(3 , 
3;12,
12
)$,
irreps = $3_{4}^{1,0}
\hskip -1.5pt \otimes \hskip -1.5pt
1_{3}^{1,0}\oplus
3_{3}^{1,0}
\hskip -1.5pt \otimes \hskip -1.5pt
1_{4}^{1,6}$,
pord$(\rho_\text{isum}(\mathfrak{t})) = 12$,

\vskip 0.7ex
\hangindent=5.5em \hangafter=1
{\white .}\hskip 1em $\rho_\text{isum}(\mathfrak{t})$ =
 $( \frac{1}{3},
\frac{1}{12},
\frac{7}{12} )
\oplus
( \frac{3}{4},
\frac{1}{12},
\frac{5}{12} )
$,

\vskip 0.7ex
\hangindent=5.5em \hangafter=1
{\white .}\hskip 1em $\rho_\text{isum}(\mathfrak{s})$ =
($0$,
$\sqrt{\frac{1}{2}}$,
$\sqrt{\frac{1}{2}}$;
$-\frac{1}{2}$,
$\frac{1}{2}$;
$-\frac{1}{2}$)
 $\oplus$
$\mathrm{i}$($\frac{1}{3}$,
$\frac{2}{3}$,
$\frac{2}{3}$;\ \ 
$\frac{1}{3}$,
$-\frac{2}{3}$;\ \ 
$\frac{1}{3}$)

Fail:
number of self dual objects $|$Tr($\rho(\mathfrak s^2)$)$|$ = 0. Prop. B.4 (1)\
 eqn. (B.16)

 \ \color{black}

 \color{blue}

\noindent 277: (dims,levels) = $(3 , 
3;12,
15
)$,
irreps = $3_{4}^{1,0}
\hskip -1.5pt \otimes \hskip -1.5pt
1_{3}^{1,0}\oplus
3_{5}^{1}
\hskip -1.5pt \otimes \hskip -1.5pt
1_{3}^{1,0}$,
pord$(\rho_\text{isum}(\mathfrak{t})) = 20$,

\vskip 0.7ex
\hangindent=5.5em \hangafter=1
{\white .}\hskip 1em $\rho_\text{isum}(\mathfrak{t})$ =
 $( \frac{1}{3},
\frac{1}{12},
\frac{7}{12} )
\oplus
( \frac{1}{3},
\frac{2}{15},
\frac{8}{15} )
$,

\vskip 0.7ex
\hangindent=5.5em \hangafter=1
{\white .}\hskip 1em $\rho_\text{isum}(\mathfrak{s})$ =
($0$,
$\sqrt{\frac{1}{2}}$,
$\sqrt{\frac{1}{2}}$;
$-\frac{1}{2}$,
$\frac{1}{2}$;
$-\frac{1}{2}$)
 $\oplus$
($\sqrt{\frac{1}{5}}$,
$-\sqrt{\frac{2}{5}}$,
$-\sqrt{\frac{2}{5}}$;
$-\frac{5+\sqrt{5}}{10}$,
$\frac{5-\sqrt{5}}{10}$;
$-\frac{5+\sqrt{5}}{10}$)

Pass. 

 \ \color{black}

 \color{blue}

\noindent 278: (dims,levels) = $(3 , 
3;12,
15
)$,
irreps = $3_{4}^{1,0}
\hskip -1.5pt \otimes \hskip -1.5pt
1_{3}^{1,0}\oplus
3_{5}^{3}
\hskip -1.5pt \otimes \hskip -1.5pt
1_{3}^{1,0}$,
pord$(\rho_\text{isum}(\mathfrak{t})) = 20$,

\vskip 0.7ex
\hangindent=5.5em \hangafter=1
{\white .}\hskip 1em $\rho_\text{isum}(\mathfrak{t})$ =
 $( \frac{1}{3},
\frac{1}{12},
\frac{7}{12} )
\oplus
( \frac{1}{3},
\frac{11}{15},
\frac{14}{15} )
$,

\vskip 0.7ex
\hangindent=5.5em \hangafter=1
{\white .}\hskip 1em $\rho_\text{isum}(\mathfrak{s})$ =
($0$,
$\sqrt{\frac{1}{2}}$,
$\sqrt{\frac{1}{2}}$;
$-\frac{1}{2}$,
$\frac{1}{2}$;
$-\frac{1}{2}$)
 $\oplus$
($-\sqrt{\frac{1}{5}}$,
$-\sqrt{\frac{2}{5}}$,
$-\sqrt{\frac{2}{5}}$;
$\frac{-5+\sqrt{5}}{10}$,
$\frac{5+\sqrt{5}}{10}$;
$\frac{-5+\sqrt{5}}{10}$)

Pass. 

 \ \color{black}

 \color{blue}

\noindent 279: (dims,levels) = $(3 , 
3;12,
20
)$,
irreps = $3_{3}^{1,0}
\hskip -1.5pt \otimes \hskip -1.5pt
1_{4}^{1,0}\oplus
3_{5}^{1}
\hskip -1.5pt \otimes \hskip -1.5pt
1_{4}^{1,0}$,
pord$(\rho_\text{isum}(\mathfrak{t})) = 15$,

\vskip 0.7ex
\hangindent=5.5em \hangafter=1
{\white .}\hskip 1em $\rho_\text{isum}(\mathfrak{t})$ =
 $( \frac{1}{4},
\frac{7}{12},
\frac{11}{12} )
\oplus
( \frac{1}{4},
\frac{1}{20},
\frac{9}{20} )
$,

\vskip 0.7ex
\hangindent=5.5em \hangafter=1
{\white .}\hskip 1em $\rho_\text{isum}(\mathfrak{s})$ =
$\mathrm{i}$($-\frac{1}{3}$,
$\frac{2}{3}$,
$\frac{2}{3}$;\ \ 
$-\frac{1}{3}$,
$\frac{2}{3}$;\ \ 
$-\frac{1}{3}$)
 $\oplus$
$\mathrm{i}$($\sqrt{\frac{1}{5}}$,
$\sqrt{\frac{2}{5}}$,
$\sqrt{\frac{2}{5}}$;\ \ 
$-\frac{5+\sqrt{5}}{10}$,
$\frac{5-\sqrt{5}}{10}$;\ \ 
$-\frac{5+\sqrt{5}}{10}$)

Pass. 

 \ \color{black}

 \color{blue}

\noindent 280: (dims,levels) = $(3 , 
3;12,
20
)$,
irreps = $3_{3}^{1,0}
\hskip -1.5pt \otimes \hskip -1.5pt
1_{4}^{1,0}\oplus
3_{5}^{3}
\hskip -1.5pt \otimes \hskip -1.5pt
1_{4}^{1,0}$,
pord$(\rho_\text{isum}(\mathfrak{t})) = 15$,

\vskip 0.7ex
\hangindent=5.5em \hangafter=1
{\white .}\hskip 1em $\rho_\text{isum}(\mathfrak{t})$ =
 $( \frac{1}{4},
\frac{7}{12},
\frac{11}{12} )
\oplus
( \frac{1}{4},
\frac{13}{20},
\frac{17}{20} )
$,

\vskip 0.7ex
\hangindent=5.5em \hangafter=1
{\white .}\hskip 1em $\rho_\text{isum}(\mathfrak{s})$ =
$\mathrm{i}$($-\frac{1}{3}$,
$\frac{2}{3}$,
$\frac{2}{3}$;\ \ 
$-\frac{1}{3}$,
$\frac{2}{3}$;\ \ 
$-\frac{1}{3}$)
 $\oplus$
$\mathrm{i}$($-\sqrt{\frac{1}{5}}$,
$\sqrt{\frac{2}{5}}$,
$\sqrt{\frac{2}{5}}$;\ \ 
$\frac{-5+\sqrt{5}}{10}$,
$\frac{5+\sqrt{5}}{10}$;\ \ 
$\frac{-5+\sqrt{5}}{10}$)

Pass. 

 \ \color{black}

\noindent 281: (dims,levels) = $(3 , 
3;12,
24
)$,
irreps = $3_{3}^{1,0}
\hskip -1.5pt \otimes \hskip -1.5pt
1_{4}^{1,0}\oplus
3_{8}^{1,3}
\hskip -1.5pt \otimes \hskip -1.5pt
1_{3}^{1,0}$,
pord$(\rho_\text{isum}(\mathfrak{t})) = 24$,

\vskip 0.7ex
\hangindent=5.5em \hangafter=1
{\white .}\hskip 1em $\rho_\text{isum}(\mathfrak{t})$ =
 $( \frac{1}{4},
\frac{7}{12},
\frac{11}{12} )
\oplus
( \frac{7}{12},
\frac{5}{24},
\frac{17}{24} )
$,

\vskip 0.7ex
\hangindent=5.5em \hangafter=1
{\white .}\hskip 1em $\rho_\text{isum}(\mathfrak{s})$ =
$\mathrm{i}$($-\frac{1}{3}$,
$\frac{2}{3}$,
$\frac{2}{3}$;\ \ 
$-\frac{1}{3}$,
$\frac{2}{3}$;\ \ 
$-\frac{1}{3}$)
 $\oplus$
($0$,
$\sqrt{\frac{1}{2}}$,
$\sqrt{\frac{1}{2}}$;
$\frac{1}{2}$,
$-\frac{1}{2}$;
$\frac{1}{2}$)

Fail:
number of self dual objects $|$Tr($\rho(\mathfrak s^2)$)$|$ = 0. Prop. B.4 (1)\
 eqn. (B.16)

 \ \color{black}

\noindent 282: (dims,levels) = $(3 , 
3;12,
24
)$,
irreps = $3_{3}^{1,0}
\hskip -1.5pt \otimes \hskip -1.5pt
1_{4}^{1,0}\oplus
3_{8}^{3,3}
\hskip -1.5pt \otimes \hskip -1.5pt
1_{3}^{1,0}$,
pord$(\rho_\text{isum}(\mathfrak{t})) = 24$,

\vskip 0.7ex
\hangindent=5.5em \hangafter=1
{\white .}\hskip 1em $\rho_\text{isum}(\mathfrak{t})$ =
 $( \frac{1}{4},
\frac{7}{12},
\frac{11}{12} )
\oplus
( \frac{7}{12},
\frac{11}{24},
\frac{23}{24} )
$,

\vskip 0.7ex
\hangindent=5.5em \hangafter=1
{\white .}\hskip 1em $\rho_\text{isum}(\mathfrak{s})$ =
$\mathrm{i}$($-\frac{1}{3}$,
$\frac{2}{3}$,
$\frac{2}{3}$;\ \ 
$-\frac{1}{3}$,
$\frac{2}{3}$;\ \ 
$-\frac{1}{3}$)
 $\oplus$
($0$,
$\sqrt{\frac{1}{2}}$,
$\sqrt{\frac{1}{2}}$;
$-\frac{1}{2}$,
$\frac{1}{2}$;
$-\frac{1}{2}$)

Fail:
number of self dual objects $|$Tr($\rho(\mathfrak s^2)$)$|$ = 0. Prop. B.4 (1)\
 eqn. (B.16)

 \ \color{black}

\noindent 283: (dims,levels) = $(3 , 
3;12,
24
)$,
irreps = $3_{3}^{1,0}
\hskip -1.5pt \otimes \hskip -1.5pt
1_{4}^{1,0}\oplus
3_{8}^{1,3}
\hskip -1.5pt \otimes \hskip -1.5pt
1_{3}^{1,4}$,
pord$(\rho_\text{isum}(\mathfrak{t})) = 24$,

\vskip 0.7ex
\hangindent=5.5em \hangafter=1
{\white .}\hskip 1em $\rho_\text{isum}(\mathfrak{t})$ =
 $( \frac{1}{4},
\frac{7}{12},
\frac{11}{12} )
\oplus
( \frac{11}{12},
\frac{1}{24},
\frac{13}{24} )
$,

\vskip 0.7ex
\hangindent=5.5em \hangafter=1
{\white .}\hskip 1em $\rho_\text{isum}(\mathfrak{s})$ =
$\mathrm{i}$($-\frac{1}{3}$,
$\frac{2}{3}$,
$\frac{2}{3}$;\ \ 
$-\frac{1}{3}$,
$\frac{2}{3}$;\ \ 
$-\frac{1}{3}$)
 $\oplus$
($0$,
$\sqrt{\frac{1}{2}}$,
$\sqrt{\frac{1}{2}}$;
$\frac{1}{2}$,
$-\frac{1}{2}$;
$\frac{1}{2}$)

Fail:
number of self dual objects $|$Tr($\rho(\mathfrak s^2)$)$|$ = 0. Prop. B.4 (1)\
 eqn. (B.16)

 \ \color{black}

\noindent 284: (dims,levels) = $(3 , 
3;12,
24
)$,
irreps = $3_{3}^{1,0}
\hskip -1.5pt \otimes \hskip -1.5pt
1_{4}^{1,0}\oplus
3_{8}^{3,3}
\hskip -1.5pt \otimes \hskip -1.5pt
1_{3}^{1,4}$,
pord$(\rho_\text{isum}(\mathfrak{t})) = 24$,

\vskip 0.7ex
\hangindent=5.5em \hangafter=1
{\white .}\hskip 1em $\rho_\text{isum}(\mathfrak{t})$ =
 $( \frac{1}{4},
\frac{7}{12},
\frac{11}{12} )
\oplus
( \frac{11}{12},
\frac{7}{24},
\frac{19}{24} )
$,

\vskip 0.7ex
\hangindent=5.5em \hangafter=1
{\white .}\hskip 1em $\rho_\text{isum}(\mathfrak{s})$ =
$\mathrm{i}$($-\frac{1}{3}$,
$\frac{2}{3}$,
$\frac{2}{3}$;\ \ 
$-\frac{1}{3}$,
$\frac{2}{3}$;\ \ 
$-\frac{1}{3}$)
 $\oplus$
($0$,
$\sqrt{\frac{1}{2}}$,
$\sqrt{\frac{1}{2}}$;
$-\frac{1}{2}$,
$\frac{1}{2}$;
$-\frac{1}{2}$)

Fail:
number of self dual objects $|$Tr($\rho(\mathfrak s^2)$)$|$ = 0. Prop. B.4 (1)\
 eqn. (B.16)

 \ \color{black}

\noindent 285: (dims,levels) = $(3 , 
3;12,
24
)$,
irreps = $3_{4}^{1,0}
\hskip -1.5pt \otimes \hskip -1.5pt
1_{3}^{1,0}\oplus
3_{8}^{1,9}
\hskip -1.5pt \otimes \hskip -1.5pt
1_{3}^{1,0}$,
pord$(\rho_\text{isum}(\mathfrak{t})) = 8$,

\vskip 0.7ex
\hangindent=5.5em \hangafter=1
{\white .}\hskip 1em $\rho_\text{isum}(\mathfrak{t})$ =
 $( \frac{1}{3},
\frac{1}{12},
\frac{7}{12} )
\oplus
( \frac{1}{12},
\frac{5}{24},
\frac{17}{24} )
$,

\vskip 0.7ex
\hangindent=5.5em \hangafter=1
{\white .}\hskip 1em $\rho_\text{isum}(\mathfrak{s})$ =
($0$,
$\sqrt{\frac{1}{2}}$,
$\sqrt{\frac{1}{2}}$;
$-\frac{1}{2}$,
$\frac{1}{2}$;
$-\frac{1}{2}$)
 $\oplus$
($0$,
$\sqrt{\frac{1}{2}}$,
$\sqrt{\frac{1}{2}}$;
$-\frac{1}{2}$,
$\frac{1}{2}$;
$-\frac{1}{2}$)

Fail:
$\sigma(\rho(\mathfrak s)_\mathrm{ndeg}) \neq
 (\rho(\mathfrak t)^a \rho(\mathfrak s) \rho(\mathfrak t)^b
 \rho(\mathfrak s) \rho(\mathfrak t)^a)_\mathrm{ndeg}$,
 $\sigma = a$ = 5. Prop. B.5 (3) eqn. (B.25)

 \ \color{black}

\noindent 286: (dims,levels) = $(3 , 
3;12,
24
)$,
irreps = $3_{4}^{1,0}
\hskip -1.5pt \otimes \hskip -1.5pt
1_{3}^{1,0}\oplus
3_{8}^{3,9}
\hskip -1.5pt \otimes \hskip -1.5pt
1_{3}^{1,0}$,
pord$(\rho_\text{isum}(\mathfrak{t})) = 8$,

\vskip 0.7ex
\hangindent=5.5em \hangafter=1
{\white .}\hskip 1em $\rho_\text{isum}(\mathfrak{t})$ =
 $( \frac{1}{3},
\frac{1}{12},
\frac{7}{12} )
\oplus
( \frac{1}{12},
\frac{11}{24},
\frac{23}{24} )
$,

\vskip 0.7ex
\hangindent=5.5em \hangafter=1
{\white .}\hskip 1em $\rho_\text{isum}(\mathfrak{s})$ =
($0$,
$\sqrt{\frac{1}{2}}$,
$\sqrt{\frac{1}{2}}$;
$-\frac{1}{2}$,
$\frac{1}{2}$;
$-\frac{1}{2}$)
 $\oplus$
($0$,
$\sqrt{\frac{1}{2}}$,
$\sqrt{\frac{1}{2}}$;
$\frac{1}{2}$,
$-\frac{1}{2}$;
$\frac{1}{2}$)

Fail:
$\sigma(\rho(\mathfrak s)_\mathrm{ndeg}) \neq
 (\rho(\mathfrak t)^a \rho(\mathfrak s) \rho(\mathfrak t)^b
 \rho(\mathfrak s) \rho(\mathfrak t)^a)_\mathrm{ndeg}$,
 $\sigma = a$ = 5. Prop. B.5 (3) eqn. (B.25)

 \ \color{black}

\noindent 287: (dims,levels) = $(3 , 
3;12,
24
)$,
irreps = $3_{4}^{1,0}
\hskip -1.5pt \otimes \hskip -1.5pt
1_{3}^{1,0}\oplus
3_{8}^{3,0}
\hskip -1.5pt \otimes \hskip -1.5pt
1_{3}^{1,0}$,
pord$(\rho_\text{isum}(\mathfrak{t})) = 8$,

\vskip 0.7ex
\hangindent=5.5em \hangafter=1
{\white .}\hskip 1em $\rho_\text{isum}(\mathfrak{t})$ =
 $( \frac{1}{3},
\frac{1}{12},
\frac{7}{12} )
\oplus
( \frac{1}{3},
\frac{5}{24},
\frac{17}{24} )
$,

\vskip 0.7ex
\hangindent=5.5em \hangafter=1
{\white .}\hskip 1em $\rho_\text{isum}(\mathfrak{s})$ =
($0$,
$\sqrt{\frac{1}{2}}$,
$\sqrt{\frac{1}{2}}$;
$-\frac{1}{2}$,
$\frac{1}{2}$;
$-\frac{1}{2}$)
 $\oplus$
$\mathrm{i}$($0$,
$\sqrt{\frac{1}{2}}$,
$\sqrt{\frac{1}{2}}$;\ \ 
$\frac{1}{2}$,
$-\frac{1}{2}$;\ \ 
$\frac{1}{2}$)

Fail:
number of self dual objects $|$Tr($\rho(\mathfrak s^2)$)$|$ = 0. Prop. B.4 (1)\
 eqn. (B.16)

 \ \color{black}

\noindent 288: (dims,levels) = $(3 , 
3;12,
24
)$,
irreps = $3_{4}^{1,0}
\hskip -1.5pt \otimes \hskip -1.5pt
1_{3}^{1,0}\oplus
3_{8}^{1,0}
\hskip -1.5pt \otimes \hskip -1.5pt
1_{3}^{1,0}$,
pord$(\rho_\text{isum}(\mathfrak{t})) = 8$,

\vskip 0.7ex
\hangindent=5.5em \hangafter=1
{\white .}\hskip 1em $\rho_\text{isum}(\mathfrak{t})$ =
 $( \frac{1}{3},
\frac{1}{12},
\frac{7}{12} )
\oplus
( \frac{1}{3},
\frac{11}{24},
\frac{23}{24} )
$,

\vskip 0.7ex
\hangindent=5.5em \hangafter=1
{\white .}\hskip 1em $\rho_\text{isum}(\mathfrak{s})$ =
($0$,
$\sqrt{\frac{1}{2}}$,
$\sqrt{\frac{1}{2}}$;
$-\frac{1}{2}$,
$\frac{1}{2}$;
$-\frac{1}{2}$)
 $\oplus$
$\mathrm{i}$($0$,
$\sqrt{\frac{1}{2}}$,
$\sqrt{\frac{1}{2}}$;\ \ 
$-\frac{1}{2}$,
$\frac{1}{2}$;\ \ 
$-\frac{1}{2}$)

Fail:
number of self dual objects $|$Tr($\rho(\mathfrak s^2)$)$|$ = 0. Prop. B.4 (1)\
 eqn. (B.16)

 \ \color{black}

\noindent 289: (dims,levels) = $(3 , 
3;12,
24
)$,
irreps = $3_{4}^{1,0}
\hskip -1.5pt \otimes \hskip -1.5pt
1_{3}^{1,0}\oplus
3_{8}^{1,3}
\hskip -1.5pt \otimes \hskip -1.5pt
1_{3}^{1,0}$,
pord$(\rho_\text{isum}(\mathfrak{t})) = 8$,

\vskip 0.7ex
\hangindent=5.5em \hangafter=1
{\white .}\hskip 1em $\rho_\text{isum}(\mathfrak{t})$ =
 $( \frac{1}{3},
\frac{1}{12},
\frac{7}{12} )
\oplus
( \frac{7}{12},
\frac{5}{24},
\frac{17}{24} )
$,

\vskip 0.7ex
\hangindent=5.5em \hangafter=1
{\white .}\hskip 1em $\rho_\text{isum}(\mathfrak{s})$ =
($0$,
$\sqrt{\frac{1}{2}}$,
$\sqrt{\frac{1}{2}}$;
$-\frac{1}{2}$,
$\frac{1}{2}$;
$-\frac{1}{2}$)
 $\oplus$
($0$,
$\sqrt{\frac{1}{2}}$,
$\sqrt{\frac{1}{2}}$;
$\frac{1}{2}$,
$-\frac{1}{2}$;
$\frac{1}{2}$)

Fail:
$\sigma(\rho(\mathfrak s)_\mathrm{ndeg}) \neq
 (\rho(\mathfrak t)^a \rho(\mathfrak s) \rho(\mathfrak t)^b
 \rho(\mathfrak s) \rho(\mathfrak t)^a)_\mathrm{ndeg}$,
 $\sigma = a$ = 5. Prop. B.5 (3) eqn. (B.25)

 \ \color{black}

\noindent 290: (dims,levels) = $(3 , 
3;12,
24
)$,
irreps = $3_{4}^{1,0}
\hskip -1.5pt \otimes \hskip -1.5pt
1_{3}^{1,0}\oplus
3_{8}^{3,3}
\hskip -1.5pt \otimes \hskip -1.5pt
1_{3}^{1,0}$,
pord$(\rho_\text{isum}(\mathfrak{t})) = 8$,

\vskip 0.7ex
\hangindent=5.5em \hangafter=1
{\white .}\hskip 1em $\rho_\text{isum}(\mathfrak{t})$ =
 $( \frac{1}{3},
\frac{1}{12},
\frac{7}{12} )
\oplus
( \frac{7}{12},
\frac{11}{24},
\frac{23}{24} )
$,

\vskip 0.7ex
\hangindent=5.5em \hangafter=1
{\white .}\hskip 1em $\rho_\text{isum}(\mathfrak{s})$ =
($0$,
$\sqrt{\frac{1}{2}}$,
$\sqrt{\frac{1}{2}}$;
$-\frac{1}{2}$,
$\frac{1}{2}$;
$-\frac{1}{2}$)
 $\oplus$
($0$,
$\sqrt{\frac{1}{2}}$,
$\sqrt{\frac{1}{2}}$;
$-\frac{1}{2}$,
$\frac{1}{2}$;
$-\frac{1}{2}$)

Fail:
$\sigma(\rho(\mathfrak s)_\mathrm{ndeg}) \neq
 (\rho(\mathfrak t)^a \rho(\mathfrak s) \rho(\mathfrak t)^b
 \rho(\mathfrak s) \rho(\mathfrak t)^a)_\mathrm{ndeg}$,
 $\sigma = a$ = 5. Prop. B.5 (3) eqn. (B.25)

 \ \color{black}

 \color{blue}

\noindent 291: (dims,levels) = $(3 , 
3;12,
60
)$,
irreps = $3_{3}^{1,0}
\hskip -1.5pt \otimes \hskip -1.5pt
1_{4}^{1,0}\oplus
3_{5}^{3}
\hskip -1.5pt \otimes \hskip -1.5pt
1_{4}^{1,0}
\hskip -1.5pt \otimes \hskip -1.5pt
1_{3}^{1,0}$,
pord$(\rho_\text{isum}(\mathfrak{t})) = 15$,

\vskip 0.7ex
\hangindent=5.5em \hangafter=1
{\white .}\hskip 1em $\rho_\text{isum}(\mathfrak{t})$ =
 $( \frac{1}{4},
\frac{7}{12},
\frac{11}{12} )
\oplus
( \frac{7}{12},
\frac{11}{60},
\frac{59}{60} )
$,

\vskip 0.7ex
\hangindent=5.5em \hangafter=1
{\white .}\hskip 1em $\rho_\text{isum}(\mathfrak{s})$ =
$\mathrm{i}$($-\frac{1}{3}$,
$\frac{2}{3}$,
$\frac{2}{3}$;\ \ 
$-\frac{1}{3}$,
$\frac{2}{3}$;\ \ 
$-\frac{1}{3}$)
 $\oplus$
$\mathrm{i}$($-\sqrt{\frac{1}{5}}$,
$\sqrt{\frac{2}{5}}$,
$\sqrt{\frac{2}{5}}$;\ \ 
$\frac{-5+\sqrt{5}}{10}$,
$\frac{5+\sqrt{5}}{10}$;\ \ 
$\frac{-5+\sqrt{5}}{10}$)

Pass. 

 \ \color{black}

 \color{blue}

\noindent 292: (dims,levels) = $(3 , 
3;12,
60
)$,
irreps = $3_{3}^{1,0}
\hskip -1.5pt \otimes \hskip -1.5pt
1_{4}^{1,0}\oplus
3_{5}^{1}
\hskip -1.5pt \otimes \hskip -1.5pt
1_{4}^{1,0}
\hskip -1.5pt \otimes \hskip -1.5pt
1_{3}^{1,0}$,
pord$(\rho_\text{isum}(\mathfrak{t})) = 15$,

\vskip 0.7ex
\hangindent=5.5em \hangafter=1
{\white .}\hskip 1em $\rho_\text{isum}(\mathfrak{t})$ =
 $( \frac{1}{4},
\frac{7}{12},
\frac{11}{12} )
\oplus
( \frac{7}{12},
\frac{23}{60},
\frac{47}{60} )
$,

\vskip 0.7ex
\hangindent=5.5em \hangafter=1
{\white .}\hskip 1em $\rho_\text{isum}(\mathfrak{s})$ =
$\mathrm{i}$($-\frac{1}{3}$,
$\frac{2}{3}$,
$\frac{2}{3}$;\ \ 
$-\frac{1}{3}$,
$\frac{2}{3}$;\ \ 
$-\frac{1}{3}$)
 $\oplus$
$\mathrm{i}$($\sqrt{\frac{1}{5}}$,
$\sqrt{\frac{2}{5}}$,
$\sqrt{\frac{2}{5}}$;\ \ 
$-\frac{5+\sqrt{5}}{10}$,
$\frac{5-\sqrt{5}}{10}$;\ \ 
$-\frac{5+\sqrt{5}}{10}$)

Pass. 

 \ \color{black}

 \color{blue}

\noindent 293: (dims,levels) = $(3 , 
3;12,
60
)$,
irreps = $3_{3}^{1,0}
\hskip -1.5pt \otimes \hskip -1.5pt
1_{4}^{1,0}\oplus
3_{5}^{1}
\hskip -1.5pt \otimes \hskip -1.5pt
1_{4}^{1,0}
\hskip -1.5pt \otimes \hskip -1.5pt
1_{3}^{1,4}$,
pord$(\rho_\text{isum}(\mathfrak{t})) = 15$,

\vskip 0.7ex
\hangindent=5.5em \hangafter=1
{\white .}\hskip 1em $\rho_\text{isum}(\mathfrak{t})$ =
 $( \frac{1}{4},
\frac{7}{12},
\frac{11}{12} )
\oplus
( \frac{11}{12},
\frac{7}{60},
\frac{43}{60} )
$,

\vskip 0.7ex
\hangindent=5.5em \hangafter=1
{\white .}\hskip 1em $\rho_\text{isum}(\mathfrak{s})$ =
$\mathrm{i}$($-\frac{1}{3}$,
$\frac{2}{3}$,
$\frac{2}{3}$;\ \ 
$-\frac{1}{3}$,
$\frac{2}{3}$;\ \ 
$-\frac{1}{3}$)
 $\oplus$
$\mathrm{i}$($\sqrt{\frac{1}{5}}$,
$\sqrt{\frac{2}{5}}$,
$\sqrt{\frac{2}{5}}$;\ \ 
$-\frac{5+\sqrt{5}}{10}$,
$\frac{5-\sqrt{5}}{10}$;\ \ 
$-\frac{5+\sqrt{5}}{10}$)

Pass. 

 \ \color{black}

 \color{blue}

\noindent 294: (dims,levels) = $(3 , 
3;12,
60
)$,
irreps = $3_{3}^{1,0}
\hskip -1.5pt \otimes \hskip -1.5pt
1_{4}^{1,0}\oplus
3_{5}^{3}
\hskip -1.5pt \otimes \hskip -1.5pt
1_{4}^{1,0}
\hskip -1.5pt \otimes \hskip -1.5pt
1_{3}^{1,4}$,
pord$(\rho_\text{isum}(\mathfrak{t})) = 15$,

\vskip 0.7ex
\hangindent=5.5em \hangafter=1
{\white .}\hskip 1em $\rho_\text{isum}(\mathfrak{t})$ =
 $( \frac{1}{4},
\frac{7}{12},
\frac{11}{12} )
\oplus
( \frac{11}{12},
\frac{19}{60},
\frac{31}{60} )
$,

\vskip 0.7ex
\hangindent=5.5em \hangafter=1
{\white .}\hskip 1em $\rho_\text{isum}(\mathfrak{s})$ =
$\mathrm{i}$($-\frac{1}{3}$,
$\frac{2}{3}$,
$\frac{2}{3}$;\ \ 
$-\frac{1}{3}$,
$\frac{2}{3}$;\ \ 
$-\frac{1}{3}$)
 $\oplus$
$\mathrm{i}$($-\sqrt{\frac{1}{5}}$,
$\sqrt{\frac{2}{5}}$,
$\sqrt{\frac{2}{5}}$;\ \ 
$\frac{-5+\sqrt{5}}{10}$,
$\frac{5+\sqrt{5}}{10}$;\ \ 
$\frac{-5+\sqrt{5}}{10}$)

Pass. 

 \ \color{black}

\noindent 295: (dims,levels) = $(3 , 
3;12,
60
)$,
irreps = $3_{4}^{1,0}
\hskip -1.5pt \otimes \hskip -1.5pt
1_{3}^{1,0}\oplus
3_{5}^{1}
\hskip -1.5pt \otimes \hskip -1.5pt
1_{4}^{1,6}
\hskip -1.5pt \otimes \hskip -1.5pt
1_{3}^{1,0}$,
pord$(\rho_\text{isum}(\mathfrak{t})) = 20$,

\vskip 0.7ex
\hangindent=5.5em \hangafter=1
{\white .}\hskip 1em $\rho_\text{isum}(\mathfrak{t})$ =
 $( \frac{1}{3},
\frac{1}{12},
\frac{7}{12} )
\oplus
( \frac{1}{12},
\frac{17}{60},
\frac{53}{60} )
$,

\vskip 0.7ex
\hangindent=5.5em \hangafter=1
{\white .}\hskip 1em $\rho_\text{isum}(\mathfrak{s})$ =
($0$,
$\sqrt{\frac{1}{2}}$,
$\sqrt{\frac{1}{2}}$;
$-\frac{1}{2}$,
$\frac{1}{2}$;
$-\frac{1}{2}$)
 $\oplus$
$\mathrm{i}$($-\sqrt{\frac{1}{5}}$,
$\sqrt{\frac{2}{5}}$,
$\sqrt{\frac{2}{5}}$;\ \ 
$\frac{5+\sqrt{5}}{10}$,
$\frac{-5+\sqrt{5}}{10}$;\ \ 
$\frac{5+\sqrt{5}}{10}$)

Fail:
number of self dual objects $|$Tr($\rho(\mathfrak s^2)$)$|$ = 0. Prop. B.4 (1)\
 eqn. (B.16)

 \ \color{black}

\noindent 296: (dims,levels) = $(3 , 
3;12,
60
)$,
irreps = $3_{4}^{1,0}
\hskip -1.5pt \otimes \hskip -1.5pt
1_{3}^{1,0}\oplus
3_{5}^{3}
\hskip -1.5pt \otimes \hskip -1.5pt
1_{4}^{1,6}
\hskip -1.5pt \otimes \hskip -1.5pt
1_{3}^{1,0}$,
pord$(\rho_\text{isum}(\mathfrak{t})) = 20$,

\vskip 0.7ex
\hangindent=5.5em \hangafter=1
{\white .}\hskip 1em $\rho_\text{isum}(\mathfrak{t})$ =
 $( \frac{1}{3},
\frac{1}{12},
\frac{7}{12} )
\oplus
( \frac{1}{12},
\frac{29}{60},
\frac{41}{60} )
$,

\vskip 0.7ex
\hangindent=5.5em \hangafter=1
{\white .}\hskip 1em $\rho_\text{isum}(\mathfrak{s})$ =
($0$,
$\sqrt{\frac{1}{2}}$,
$\sqrt{\frac{1}{2}}$;
$-\frac{1}{2}$,
$\frac{1}{2}$;
$-\frac{1}{2}$)
 $\oplus$
$\mathrm{i}$($\sqrt{\frac{1}{5}}$,
$\sqrt{\frac{2}{5}}$,
$\sqrt{\frac{2}{5}}$;\ \ 
$\frac{5-\sqrt{5}}{10}$,
$-\frac{5+\sqrt{5}}{10}$;\ \ 
$\frac{5-\sqrt{5}}{10}$)

Fail:
number of self dual objects $|$Tr($\rho(\mathfrak s^2)$)$|$ = 0. Prop. B.4 (1)\
 eqn. (B.16)

 \ \color{black}

\noindent 297: (dims,levels) = $(3 , 
3;12,
60
)$,
irreps = $3_{4}^{1,0}
\hskip -1.5pt \otimes \hskip -1.5pt
1_{3}^{1,0}\oplus
3_{5}^{3}
\hskip -1.5pt \otimes \hskip -1.5pt
1_{4}^{1,0}
\hskip -1.5pt \otimes \hskip -1.5pt
1_{3}^{1,0}$,
pord$(\rho_\text{isum}(\mathfrak{t})) = 20$,

\vskip 0.7ex
\hangindent=5.5em \hangafter=1
{\white .}\hskip 1em $\rho_\text{isum}(\mathfrak{t})$ =
 $( \frac{1}{3},
\frac{1}{12},
\frac{7}{12} )
\oplus
( \frac{7}{12},
\frac{11}{60},
\frac{59}{60} )
$,

\vskip 0.7ex
\hangindent=5.5em \hangafter=1
{\white .}\hskip 1em $\rho_\text{isum}(\mathfrak{s})$ =
($0$,
$\sqrt{\frac{1}{2}}$,
$\sqrt{\frac{1}{2}}$;
$-\frac{1}{2}$,
$\frac{1}{2}$;
$-\frac{1}{2}$)
 $\oplus$
$\mathrm{i}$($-\sqrt{\frac{1}{5}}$,
$\sqrt{\frac{2}{5}}$,
$\sqrt{\frac{2}{5}}$;\ \ 
$\frac{-5+\sqrt{5}}{10}$,
$\frac{5+\sqrt{5}}{10}$;\ \ 
$\frac{-5+\sqrt{5}}{10}$)

Fail:
number of self dual objects $|$Tr($\rho(\mathfrak s^2)$)$|$ = 0. Prop. B.4 (1)\
 eqn. (B.16)

 \ \color{black}

\noindent 298: (dims,levels) = $(3 , 
3;12,
60
)$,
irreps = $3_{4}^{1,0}
\hskip -1.5pt \otimes \hskip -1.5pt
1_{3}^{1,0}\oplus
3_{5}^{1}
\hskip -1.5pt \otimes \hskip -1.5pt
1_{4}^{1,0}
\hskip -1.5pt \otimes \hskip -1.5pt
1_{3}^{1,0}$,
pord$(\rho_\text{isum}(\mathfrak{t})) = 20$,

\vskip 0.7ex
\hangindent=5.5em \hangafter=1
{\white .}\hskip 1em $\rho_\text{isum}(\mathfrak{t})$ =
 $( \frac{1}{3},
\frac{1}{12},
\frac{7}{12} )
\oplus
( \frac{7}{12},
\frac{23}{60},
\frac{47}{60} )
$,

\vskip 0.7ex
\hangindent=5.5em \hangafter=1
{\white .}\hskip 1em $\rho_\text{isum}(\mathfrak{s})$ =
($0$,
$\sqrt{\frac{1}{2}}$,
$\sqrt{\frac{1}{2}}$;
$-\frac{1}{2}$,
$\frac{1}{2}$;
$-\frac{1}{2}$)
 $\oplus$
$\mathrm{i}$($\sqrt{\frac{1}{5}}$,
$\sqrt{\frac{2}{5}}$,
$\sqrt{\frac{2}{5}}$;\ \ 
$-\frac{5+\sqrt{5}}{10}$,
$\frac{5-\sqrt{5}}{10}$;\ \ 
$-\frac{5+\sqrt{5}}{10}$)

Fail:
number of self dual objects $|$Tr($\rho(\mathfrak s^2)$)$|$ = 0. Prop. B.4 (1)\
 eqn. (B.16)

 \ \color{black}

\noindent 299: (dims,levels) = $(3 , 
3;15,
15
)$,
irreps = $3_{5}^{1}
\hskip -1.5pt \otimes \hskip -1.5pt
1_{3}^{1,0}\oplus
3_{5}^{3}
\hskip -1.5pt \otimes \hskip -1.5pt
1_{3}^{1,0}$,
pord$(\rho_\text{isum}(\mathfrak{t})) = 5$,

\vskip 0.7ex
\hangindent=5.5em \hangafter=1
{\white .}\hskip 1em $\rho_\text{isum}(\mathfrak{t})$ =
 $( \frac{1}{3},
\frac{2}{15},
\frac{8}{15} )
\oplus
( \frac{1}{3},
\frac{11}{15},
\frac{14}{15} )
$,

\vskip 0.7ex
\hangindent=5.5em \hangafter=1
{\white .}\hskip 1em $\rho_\text{isum}(\mathfrak{s})$ =
($\sqrt{\frac{1}{5}}$,
$-\sqrt{\frac{2}{5}}$,
$-\sqrt{\frac{2}{5}}$;
$-\frac{5+\sqrt{5}}{10}$,
$\frac{5-\sqrt{5}}{10}$;
$-\frac{5+\sqrt{5}}{10}$)
 $\oplus$
($-\sqrt{\frac{1}{5}}$,
$-\sqrt{\frac{2}{5}}$,
$-\sqrt{\frac{2}{5}}$;
$\frac{-5+\sqrt{5}}{10}$,
$\frac{5+\sqrt{5}}{10}$;
$\frac{-5+\sqrt{5}}{10}$)

Fail:
Integral: $D_{\rho}(\sigma)_{\theta} \propto $ id,
 for all $\sigma$ and all $\theta$-eigenspaces that can contain unit. Prop. B.5 (6)

 \ \color{black}

\noindent 300: (dims,levels) = $(3 , 
3;15,
24
)$,
irreps = $3_{5}^{1}
\hskip -1.5pt \otimes \hskip -1.5pt
1_{3}^{1,0}\oplus
3_{8}^{3,0}
\hskip -1.5pt \otimes \hskip -1.5pt
1_{3}^{1,0}$,
pord$(\rho_\text{isum}(\mathfrak{t})) = 40$,

\vskip 0.7ex
\hangindent=5.5em \hangafter=1
{\white .}\hskip 1em $\rho_\text{isum}(\mathfrak{t})$ =
 $( \frac{1}{3},
\frac{2}{15},
\frac{8}{15} )
\oplus
( \frac{1}{3},
\frac{5}{24},
\frac{17}{24} )
$,

\vskip 0.7ex
\hangindent=5.5em \hangafter=1
{\white .}\hskip 1em $\rho_\text{isum}(\mathfrak{s})$ =
($\sqrt{\frac{1}{5}}$,
$-\sqrt{\frac{2}{5}}$,
$-\sqrt{\frac{2}{5}}$;
$-\frac{5+\sqrt{5}}{10}$,
$\frac{5-\sqrt{5}}{10}$;
$-\frac{5+\sqrt{5}}{10}$)
 $\oplus$
$\mathrm{i}$($0$,
$\sqrt{\frac{1}{2}}$,
$\sqrt{\frac{1}{2}}$;\ \ 
$\frac{1}{2}$,
$-\frac{1}{2}$;\ \ 
$\frac{1}{2}$)

Fail:
number of self dual objects $|$Tr($\rho(\mathfrak s^2)$)$|$ = 0. Prop. B.4 (1)\
 eqn. (B.16)

 \ \color{black}

\noindent 301: (dims,levels) = $(3 , 
3;15,
24
)$,
irreps = $3_{5}^{1}
\hskip -1.5pt \otimes \hskip -1.5pt
1_{3}^{1,0}\oplus
3_{8}^{1,0}
\hskip -1.5pt \otimes \hskip -1.5pt
1_{3}^{1,0}$,
pord$(\rho_\text{isum}(\mathfrak{t})) = 40$,

\vskip 0.7ex
\hangindent=5.5em \hangafter=1
{\white .}\hskip 1em $\rho_\text{isum}(\mathfrak{t})$ =
 $( \frac{1}{3},
\frac{2}{15},
\frac{8}{15} )
\oplus
( \frac{1}{3},
\frac{11}{24},
\frac{23}{24} )
$,

\vskip 0.7ex
\hangindent=5.5em \hangafter=1
{\white .}\hskip 1em $\rho_\text{isum}(\mathfrak{s})$ =
($\sqrt{\frac{1}{5}}$,
$-\sqrt{\frac{2}{5}}$,
$-\sqrt{\frac{2}{5}}$;
$-\frac{5+\sqrt{5}}{10}$,
$\frac{5-\sqrt{5}}{10}$;
$-\frac{5+\sqrt{5}}{10}$)
 $\oplus$
$\mathrm{i}$($0$,
$\sqrt{\frac{1}{2}}$,
$\sqrt{\frac{1}{2}}$;\ \ 
$-\frac{1}{2}$,
$\frac{1}{2}$;\ \ 
$-\frac{1}{2}$)

Fail:
number of self dual objects $|$Tr($\rho(\mathfrak s^2)$)$|$ = 0. Prop. B.4 (1)\
 eqn. (B.16)

 \ \color{black}

\noindent 302: (dims,levels) = $(3 , 
3;15,
24
)$,
irreps = $3_{5}^{3}
\hskip -1.5pt \otimes \hskip -1.5pt
1_{3}^{1,0}\oplus
3_{8}^{3,0}
\hskip -1.5pt \otimes \hskip -1.5pt
1_{3}^{1,0}$,
pord$(\rho_\text{isum}(\mathfrak{t})) = 40$,

\vskip 0.7ex
\hangindent=5.5em \hangafter=1
{\white .}\hskip 1em $\rho_\text{isum}(\mathfrak{t})$ =
 $( \frac{1}{3},
\frac{11}{15},
\frac{14}{15} )
\oplus
( \frac{1}{3},
\frac{5}{24},
\frac{17}{24} )
$,

\vskip 0.7ex
\hangindent=5.5em \hangafter=1
{\white .}\hskip 1em $\rho_\text{isum}(\mathfrak{s})$ =
($-\sqrt{\frac{1}{5}}$,
$-\sqrt{\frac{2}{5}}$,
$-\sqrt{\frac{2}{5}}$;
$\frac{-5+\sqrt{5}}{10}$,
$\frac{5+\sqrt{5}}{10}$;
$\frac{-5+\sqrt{5}}{10}$)
 $\oplus$
$\mathrm{i}$($0$,
$\sqrt{\frac{1}{2}}$,
$\sqrt{\frac{1}{2}}$;\ \ 
$\frac{1}{2}$,
$-\frac{1}{2}$;\ \ 
$\frac{1}{2}$)

Fail:
number of self dual objects $|$Tr($\rho(\mathfrak s^2)$)$|$ = 0. Prop. B.4 (1)\
 eqn. (B.16)

 \ \color{black}

\noindent 303: (dims,levels) = $(3 , 
3;15,
24
)$,
irreps = $3_{5}^{3}
\hskip -1.5pt \otimes \hskip -1.5pt
1_{3}^{1,0}\oplus
3_{8}^{1,0}
\hskip -1.5pt \otimes \hskip -1.5pt
1_{3}^{1,0}$,
pord$(\rho_\text{isum}(\mathfrak{t})) = 40$,

\vskip 0.7ex
\hangindent=5.5em \hangafter=1
{\white .}\hskip 1em $\rho_\text{isum}(\mathfrak{t})$ =
 $( \frac{1}{3},
\frac{11}{15},
\frac{14}{15} )
\oplus
( \frac{1}{3},
\frac{11}{24},
\frac{23}{24} )
$,

\vskip 0.7ex
\hangindent=5.5em \hangafter=1
{\white .}\hskip 1em $\rho_\text{isum}(\mathfrak{s})$ =
($-\sqrt{\frac{1}{5}}$,
$-\sqrt{\frac{2}{5}}$,
$-\sqrt{\frac{2}{5}}$;
$\frac{-5+\sqrt{5}}{10}$,
$\frac{5+\sqrt{5}}{10}$;
$\frac{-5+\sqrt{5}}{10}$)
 $\oplus$
$\mathrm{i}$($0$,
$\sqrt{\frac{1}{2}}$,
$\sqrt{\frac{1}{2}}$;\ \ 
$-\frac{1}{2}$,
$\frac{1}{2}$;\ \ 
$-\frac{1}{2}$)

Fail:
number of self dual objects $|$Tr($\rho(\mathfrak s^2)$)$|$ = 0. Prop. B.4 (1)\
 eqn. (B.16)

 \ \color{black}

\noindent 304: (dims,levels) = $(3 , 
3;16,
16
)$,
irreps = $3_{16}^{1,0}\oplus
3_{16}^{7,3}$,
pord$(\rho_\text{isum}(\mathfrak{t})) = 16$,

\vskip 0.7ex
\hangindent=5.5em \hangafter=1
{\white .}\hskip 1em $\rho_\text{isum}(\mathfrak{t})$ =
 $( \frac{1}{8},
\frac{1}{16},
\frac{9}{16} )
\oplus
( \frac{1}{8},
\frac{3}{16},
\frac{11}{16} )
$,

\vskip 0.7ex
\hangindent=5.5em \hangafter=1
{\white .}\hskip 1em $\rho_\text{isum}(\mathfrak{s})$ =
$\mathrm{i}$($0$,
$\sqrt{\frac{1}{2}}$,
$\sqrt{\frac{1}{2}}$;\ \ 
$-\frac{1}{2}$,
$\frac{1}{2}$;\ \ 
$-\frac{1}{2}$)
 $\oplus$
($0$,
$\sqrt{\frac{1}{2}}$,
$\sqrt{\frac{1}{2}}$;
$-\frac{1}{2}$,
$\frac{1}{2}$;
$-\frac{1}{2}$)

Fail:
number of self dual objects $|$Tr($\rho(\mathfrak s^2)$)$|$ = 0. Prop. B.4 (1)\
 eqn. (B.16)

 \ \color{black}

\noindent 305: (dims,levels) = $(3 , 
3;16,
16
)$,
irreps = $3_{16}^{1,0}\oplus
3_{16}^{5,6}$,
pord$(\rho_\text{isum}(\mathfrak{t})) = 16$,

\vskip 0.7ex
\hangindent=5.5em \hangafter=1
{\white .}\hskip 1em $\rho_\text{isum}(\mathfrak{t})$ =
 $( \frac{1}{8},
\frac{1}{16},
\frac{9}{16} )
\oplus
( \frac{1}{8},
\frac{5}{16},
\frac{13}{16} )
$,

\vskip 0.7ex
\hangindent=5.5em \hangafter=1
{\white .}\hskip 1em $\rho_\text{isum}(\mathfrak{s})$ =
$\mathrm{i}$($0$,
$\sqrt{\frac{1}{2}}$,
$\sqrt{\frac{1}{2}}$;\ \ 
$-\frac{1}{2}$,
$\frac{1}{2}$;\ \ 
$-\frac{1}{2}$)
 $\oplus$
$\mathrm{i}$($0$,
$\sqrt{\frac{1}{2}}$,
$\sqrt{\frac{1}{2}}$;\ \ 
$\frac{1}{2}$,
$-\frac{1}{2}$;\ \ 
$\frac{1}{2}$)

Fail:
all rows of $U \rho(\mathfrak s) U^\dagger$
 contain zero for any block-diagonal $U$. Prop. B.5 (4) eqn. (B.27)

 \ \color{black}

\noindent 306: (dims,levels) = $(3 , 
3;16,
16
)$,
irreps = $3_{16}^{1,0}\oplus
3_{16}^{3,9}$,
pord$(\rho_\text{isum}(\mathfrak{t})) = 16$,

\vskip 0.7ex
\hangindent=5.5em \hangafter=1
{\white .}\hskip 1em $\rho_\text{isum}(\mathfrak{t})$ =
 $( \frac{1}{8},
\frac{1}{16},
\frac{9}{16} )
\oplus
( \frac{1}{8},
\frac{7}{16},
\frac{15}{16} )
$,

\vskip 0.7ex
\hangindent=5.5em \hangafter=1
{\white .}\hskip 1em $\rho_\text{isum}(\mathfrak{s})$ =
$\mathrm{i}$($0$,
$\sqrt{\frac{1}{2}}$,
$\sqrt{\frac{1}{2}}$;\ \ 
$-\frac{1}{2}$,
$\frac{1}{2}$;\ \ 
$-\frac{1}{2}$)
 $\oplus$
($0$,
$\sqrt{\frac{1}{2}}$,
$\sqrt{\frac{1}{2}}$;
$\frac{1}{2}$,
$-\frac{1}{2}$;
$\frac{1}{2}$)

Fail:
number of self dual objects $|$Tr($\rho(\mathfrak s^2)$)$|$ = 0. Prop. B.4 (1)\
 eqn. (B.16)

 \ \color{black}

\noindent 307: (dims,levels) = $(3 , 
3;16,
16
)$,
irreps = $3_{16}^{1,0}\oplus
3_{16}^{5,9}$,
pord$(\rho_\text{isum}(\mathfrak{t})) = 16$,

\vskip 0.7ex
\hangindent=5.5em \hangafter=1
{\white .}\hskip 1em $\rho_\text{isum}(\mathfrak{t})$ =
 $( \frac{1}{8},
\frac{1}{16},
\frac{9}{16} )
\oplus
( \frac{3}{8},
\frac{1}{16},
\frac{9}{16} )
$,

\vskip 0.7ex
\hangindent=5.5em \hangafter=1
{\white .}\hskip 1em $\rho_\text{isum}(\mathfrak{s})$ =
$\mathrm{i}$($0$,
$\sqrt{\frac{1}{2}}$,
$\sqrt{\frac{1}{2}}$;\ \ 
$-\frac{1}{2}$,
$\frac{1}{2}$;\ \ 
$-\frac{1}{2}$)
 $\oplus$
($0$,
$\sqrt{\frac{1}{2}}$,
$\sqrt{\frac{1}{2}}$;
$-\frac{1}{2}$,
$\frac{1}{2}$;
$-\frac{1}{2}$)

Fail:
number of self dual objects $|$Tr($\rho(\mathfrak s^2)$)$|$ = 0. Prop. B.4 (1)\
 eqn. (B.16)

 \ \color{black}

 \color{blue}

\noindent 308: (dims,levels) = $(3 , 
3;16,
16
)$,
irreps = $3_{16}^{1,0}\oplus
3_{16}^{1,6}$,
pord$(\rho_\text{isum}(\mathfrak{t})) = 16$,

\vskip 0.7ex
\hangindent=5.5em \hangafter=1
{\white .}\hskip 1em $\rho_\text{isum}(\mathfrak{t})$ =
 $( \frac{1}{8},
\frac{1}{16},
\frac{9}{16} )
\oplus
( \frac{5}{8},
\frac{1}{16},
\frac{9}{16} )
$,

\vskip 0.7ex
\hangindent=5.5em \hangafter=1
{\white .}\hskip 1em $\rho_\text{isum}(\mathfrak{s})$ =
$\mathrm{i}$($0$,
$\sqrt{\frac{1}{2}}$,
$\sqrt{\frac{1}{2}}$;\ \ 
$-\frac{1}{2}$,
$\frac{1}{2}$;\ \ 
$-\frac{1}{2}$)
 $\oplus$
$\mathrm{i}$($0$,
$\sqrt{\frac{1}{2}}$,
$\sqrt{\frac{1}{2}}$;\ \ 
$\frac{1}{2}$,
$-\frac{1}{2}$;\ \ 
$\frac{1}{2}$)

Pass. 

 \ \color{black}

\noindent 309: (dims,levels) = $(3 , 
3;16,
16
)$,
irreps = $3_{16}^{1,0}\oplus
3_{16}^{5,3}$,
pord$(\rho_\text{isum}(\mathfrak{t})) = 16$,

\vskip 0.7ex
\hangindent=5.5em \hangafter=1
{\white .}\hskip 1em $\rho_\text{isum}(\mathfrak{t})$ =
 $( \frac{1}{8},
\frac{1}{16},
\frac{9}{16} )
\oplus
( \frac{7}{8},
\frac{1}{16},
\frac{9}{16} )
$,

\vskip 0.7ex
\hangindent=5.5em \hangafter=1
{\white .}\hskip 1em $\rho_\text{isum}(\mathfrak{s})$ =
$\mathrm{i}$($0$,
$\sqrt{\frac{1}{2}}$,
$\sqrt{\frac{1}{2}}$;\ \ 
$-\frac{1}{2}$,
$\frac{1}{2}$;\ \ 
$-\frac{1}{2}$)
 $\oplus$
($0$,
$\sqrt{\frac{1}{2}}$,
$\sqrt{\frac{1}{2}}$;
$\frac{1}{2}$,
$-\frac{1}{2}$;
$\frac{1}{2}$)

Fail:
number of self dual objects $|$Tr($\rho(\mathfrak s^2)$)$|$ = 0. Prop. B.4 (1)\
 eqn. (B.16)

 \ \color{black}

\noindent 310: (dims,levels) = $(3 , 
3;16,
16
)$,
irreps = $3_{16}^{3,0}\oplus
3_{16}^{1,3}$,
pord$(\rho_\text{isum}(\mathfrak{t})) = 16$,

\vskip 0.7ex
\hangindent=5.5em \hangafter=1
{\white .}\hskip 1em $\rho_\text{isum}(\mathfrak{t})$ =
 $( \frac{3}{8},
\frac{3}{16},
\frac{11}{16} )
\oplus
( \frac{3}{8},
\frac{5}{16},
\frac{13}{16} )
$,

\vskip 0.7ex
\hangindent=5.5em \hangafter=1
{\white .}\hskip 1em $\rho_\text{isum}(\mathfrak{s})$ =
$\mathrm{i}$($0$,
$\sqrt{\frac{1}{2}}$,
$\sqrt{\frac{1}{2}}$;\ \ 
$\frac{1}{2}$,
$-\frac{1}{2}$;\ \ 
$\frac{1}{2}$)
 $\oplus$
($0$,
$\sqrt{\frac{1}{2}}$,
$\sqrt{\frac{1}{2}}$;
$\frac{1}{2}$,
$-\frac{1}{2}$;
$\frac{1}{2}$)

Fail:
number of self dual objects $|$Tr($\rho(\mathfrak s^2)$)$|$ = 0. Prop. B.4 (1)\
 eqn. (B.16)

 \ \color{black}

\noindent 311: (dims,levels) = $(3 , 
3;16,
16
)$,
irreps = $3_{16}^{3,0}\oplus
3_{16}^{7,6}$,
pord$(\rho_\text{isum}(\mathfrak{t})) = 16$,

\vskip 0.7ex
\hangindent=5.5em \hangafter=1
{\white .}\hskip 1em $\rho_\text{isum}(\mathfrak{t})$ =
 $( \frac{3}{8},
\frac{3}{16},
\frac{11}{16} )
\oplus
( \frac{3}{8},
\frac{7}{16},
\frac{15}{16} )
$,

\vskip 0.7ex
\hangindent=5.5em \hangafter=1
{\white .}\hskip 1em $\rho_\text{isum}(\mathfrak{s})$ =
$\mathrm{i}$($0$,
$\sqrt{\frac{1}{2}}$,
$\sqrt{\frac{1}{2}}$;\ \ 
$\frac{1}{2}$,
$-\frac{1}{2}$;\ \ 
$\frac{1}{2}$)
 $\oplus$
$\mathrm{i}$($0$,
$\sqrt{\frac{1}{2}}$,
$\sqrt{\frac{1}{2}}$;\ \ 
$-\frac{1}{2}$,
$\frac{1}{2}$;\ \ 
$-\frac{1}{2}$)

Fail:
all rows of $U \rho(\mathfrak s) U^\dagger$
 contain zero for any block-diagonal $U$. Prop. B.5 (4) eqn. (B.27)

 \ \color{black}

\noindent 312: (dims,levels) = $(3 , 
3;16,
16
)$,
irreps = $3_{16}^{3,0}\oplus
3_{16}^{7,9}$,
pord$(\rho_\text{isum}(\mathfrak{t})) = 16$,

\vskip 0.7ex
\hangindent=5.5em \hangafter=1
{\white .}\hskip 1em $\rho_\text{isum}(\mathfrak{t})$ =
 $( \frac{3}{8},
\frac{3}{16},
\frac{11}{16} )
\oplus
( \frac{5}{8},
\frac{3}{16},
\frac{11}{16} )
$,

\vskip 0.7ex
\hangindent=5.5em \hangafter=1
{\white .}\hskip 1em $\rho_\text{isum}(\mathfrak{s})$ =
$\mathrm{i}$($0$,
$\sqrt{\frac{1}{2}}$,
$\sqrt{\frac{1}{2}}$;\ \ 
$\frac{1}{2}$,
$-\frac{1}{2}$;\ \ 
$\frac{1}{2}$)
 $\oplus$
($0$,
$\sqrt{\frac{1}{2}}$,
$\sqrt{\frac{1}{2}}$;
$\frac{1}{2}$,
$-\frac{1}{2}$;
$\frac{1}{2}$)

Fail:
number of self dual objects $|$Tr($\rho(\mathfrak s^2)$)$|$ = 0. Prop. B.4 (1)\
 eqn. (B.16)

 \ \color{black}

 \color{blue}

\noindent 313: (dims,levels) = $(3 , 
3;16,
16
)$,
irreps = $3_{16}^{3,0}\oplus
3_{16}^{3,6}$,
pord$(\rho_\text{isum}(\mathfrak{t})) = 16$,

\vskip 0.7ex
\hangindent=5.5em \hangafter=1
{\white .}\hskip 1em $\rho_\text{isum}(\mathfrak{t})$ =
 $( \frac{3}{8},
\frac{3}{16},
\frac{11}{16} )
\oplus
( \frac{7}{8},
\frac{3}{16},
\frac{11}{16} )
$,

\vskip 0.7ex
\hangindent=5.5em \hangafter=1
{\white .}\hskip 1em $\rho_\text{isum}(\mathfrak{s})$ =
$\mathrm{i}$($0$,
$\sqrt{\frac{1}{2}}$,
$\sqrt{\frac{1}{2}}$;\ \ 
$\frac{1}{2}$,
$-\frac{1}{2}$;\ \ 
$\frac{1}{2}$)
 $\oplus$
$\mathrm{i}$($0$,
$\sqrt{\frac{1}{2}}$,
$\sqrt{\frac{1}{2}}$;\ \ 
$-\frac{1}{2}$,
$\frac{1}{2}$;\ \ 
$-\frac{1}{2}$)

Pass. 

 \ \color{black}

\noindent 314: (dims,levels) = $(3 , 
3;16,
16
)$,
irreps = $3_{16}^{5,0}\oplus
3_{16}^{3,3}$,
pord$(\rho_\text{isum}(\mathfrak{t})) = 16$,

\vskip 0.7ex
\hangindent=5.5em \hangafter=1
{\white .}\hskip 1em $\rho_\text{isum}(\mathfrak{t})$ =
 $( \frac{5}{8},
\frac{5}{16},
\frac{13}{16} )
\oplus
( \frac{5}{8},
\frac{7}{16},
\frac{15}{16} )
$,

\vskip 0.7ex
\hangindent=5.5em \hangafter=1
{\white .}\hskip 1em $\rho_\text{isum}(\mathfrak{s})$ =
$\mathrm{i}$($0$,
$\sqrt{\frac{1}{2}}$,
$\sqrt{\frac{1}{2}}$;\ \ 
$-\frac{1}{2}$,
$\frac{1}{2}$;\ \ 
$-\frac{1}{2}$)
 $\oplus$
($0$,
$\sqrt{\frac{1}{2}}$,
$\sqrt{\frac{1}{2}}$;
$-\frac{1}{2}$,
$\frac{1}{2}$;
$-\frac{1}{2}$)

Fail:
number of self dual objects $|$Tr($\rho(\mathfrak s^2)$)$|$ = 0. Prop. B.4 (1)\
 eqn. (B.16)

 \ \color{black}

\noindent 315: (dims,levels) = $(3 , 
3;16,
16
)$,
irreps = $3_{16}^{5,0}\oplus
3_{16}^{1,9}$,
pord$(\rho_\text{isum}(\mathfrak{t})) = 16$,

\vskip 0.7ex
\hangindent=5.5em \hangafter=1
{\white .}\hskip 1em $\rho_\text{isum}(\mathfrak{t})$ =
 $( \frac{5}{8},
\frac{5}{16},
\frac{13}{16} )
\oplus
( \frac{7}{8},
\frac{5}{16},
\frac{13}{16} )
$,

\vskip 0.7ex
\hangindent=5.5em \hangafter=1
{\white .}\hskip 1em $\rho_\text{isum}(\mathfrak{s})$ =
$\mathrm{i}$($0$,
$\sqrt{\frac{1}{2}}$,
$\sqrt{\frac{1}{2}}$;\ \ 
$-\frac{1}{2}$,
$\frac{1}{2}$;\ \ 
$-\frac{1}{2}$)
 $\oplus$
($0$,
$\sqrt{\frac{1}{2}}$,
$\sqrt{\frac{1}{2}}$;
$-\frac{1}{2}$,
$\frac{1}{2}$;
$-\frac{1}{2}$)

Fail:
number of self dual objects $|$Tr($\rho(\mathfrak s^2)$)$|$ = 0. Prop. B.4 (1)\
 eqn. (B.16)

 \ \color{black}

\noindent 316: (dims,levels) = $(3 , 
3;20,
20
)$,
irreps = $3_{5}^{1}
\hskip -1.5pt \otimes \hskip -1.5pt
1_{4}^{1,0}\oplus
3_{5}^{3}
\hskip -1.5pt \otimes \hskip -1.5pt
1_{4}^{1,0}$,
pord$(\rho_\text{isum}(\mathfrak{t})) = 5$,

\vskip 0.7ex
\hangindent=5.5em \hangafter=1
{\white .}\hskip 1em $\rho_\text{isum}(\mathfrak{t})$ =
 $( \frac{1}{4},
\frac{1}{20},
\frac{9}{20} )
\oplus
( \frac{1}{4},
\frac{13}{20},
\frac{17}{20} )
$,

\vskip 0.7ex
\hangindent=5.5em \hangafter=1
{\white .}\hskip 1em $\rho_\text{isum}(\mathfrak{s})$ =
$\mathrm{i}$($\sqrt{\frac{1}{5}}$,
$\sqrt{\frac{2}{5}}$,
$\sqrt{\frac{2}{5}}$;\ \ 
$-\frac{5+\sqrt{5}}{10}$,
$\frac{5-\sqrt{5}}{10}$;\ \ 
$-\frac{5+\sqrt{5}}{10}$)
 $\oplus$
$\mathrm{i}$($-\sqrt{\frac{1}{5}}$,
$\sqrt{\frac{2}{5}}$,
$\sqrt{\frac{2}{5}}$;\ \ 
$\frac{-5+\sqrt{5}}{10}$,
$\frac{5+\sqrt{5}}{10}$;\ \ 
$\frac{-5+\sqrt{5}}{10}$)

Fail:
Integral: $D_{\rho}(\sigma)_{\theta} \propto $ id,
 for all $\sigma$ and all $\theta$-eigenspaces that can contain unit. Prop. B.5 (6)

 \ \color{black}

\noindent 317: (dims,levels) = $(3 , 
3;24,
24
)$,
irreps = $3_{8}^{3,0}
\hskip -1.5pt \otimes \hskip -1.5pt
1_{3}^{1,0}\oplus
3_{8}^{1,0}
\hskip -1.5pt \otimes \hskip -1.5pt
1_{3}^{1,0}$,
pord$(\rho_\text{isum}(\mathfrak{t})) = 8$,

\vskip 0.7ex
\hangindent=5.5em \hangafter=1
{\white .}\hskip 1em $\rho_\text{isum}(\mathfrak{t})$ =
 $( \frac{1}{3},
\frac{5}{24},
\frac{17}{24} )
\oplus
( \frac{1}{3},
\frac{11}{24},
\frac{23}{24} )
$,

\vskip 0.7ex
\hangindent=5.5em \hangafter=1
{\white .}\hskip 1em $\rho_\text{isum}(\mathfrak{s})$ =
$\mathrm{i}$($0$,
$\sqrt{\frac{1}{2}}$,
$\sqrt{\frac{1}{2}}$;\ \ 
$\frac{1}{2}$,
$-\frac{1}{2}$;\ \ 
$\frac{1}{2}$)
 $\oplus$
$\mathrm{i}$($0$,
$\sqrt{\frac{1}{2}}$,
$\sqrt{\frac{1}{2}}$;\ \ 
$-\frac{1}{2}$,
$\frac{1}{2}$;\ \ 
$-\frac{1}{2}$)

Fail:
all rows of $U \rho(\mathfrak s) U^\dagger$
 contain zero for any block-diagonal $U$. Prop. B.5 (4) eqn. (B.27)

 \ \color{black}

\noindent 318: (dims,levels) = $(3 , 
3;24,
24
)$,
irreps = $3_{8}^{3,0}
\hskip -1.5pt \otimes \hskip -1.5pt
1_{3}^{1,0}\oplus
3_{8}^{1,3}
\hskip -1.5pt \otimes \hskip -1.5pt
1_{3}^{1,0}$,
pord$(\rho_\text{isum}(\mathfrak{t})) = 8$,

\vskip 0.7ex
\hangindent=5.5em \hangafter=1
{\white .}\hskip 1em $\rho_\text{isum}(\mathfrak{t})$ =
 $( \frac{1}{3},
\frac{5}{24},
\frac{17}{24} )
\oplus
( \frac{7}{12},
\frac{5}{24},
\frac{17}{24} )
$,

\vskip 0.7ex
\hangindent=5.5em \hangafter=1
{\white .}\hskip 1em $\rho_\text{isum}(\mathfrak{s})$ =
$\mathrm{i}$($0$,
$\sqrt{\frac{1}{2}}$,
$\sqrt{\frac{1}{2}}$;\ \ 
$\frac{1}{2}$,
$-\frac{1}{2}$;\ \ 
$\frac{1}{2}$)
 $\oplus$
($0$,
$\sqrt{\frac{1}{2}}$,
$\sqrt{\frac{1}{2}}$;
$\frac{1}{2}$,
$-\frac{1}{2}$;
$\frac{1}{2}$)

Fail:
number of self dual objects $|$Tr($\rho(\mathfrak s^2)$)$|$ = 0. Prop. B.4 (1)\
 eqn. (B.16)

 \ \color{black}

\noindent 319: (dims,levels) = $(3 , 
3;24,
24
)$,
irreps = $3_{8}^{3,0}
\hskip -1.5pt \otimes \hskip -1.5pt
1_{3}^{1,0}\oplus
3_{8}^{3,6}
\hskip -1.5pt \otimes \hskip -1.5pt
1_{3}^{1,0}$,
pord$(\rho_\text{isum}(\mathfrak{t})) = 8$,

\vskip 0.7ex
\hangindent=5.5em \hangafter=1
{\white .}\hskip 1em $\rho_\text{isum}(\mathfrak{t})$ =
 $( \frac{1}{3},
\frac{5}{24},
\frac{17}{24} )
\oplus
( \frac{5}{6},
\frac{5}{24},
\frac{17}{24} )
$,

\vskip 0.7ex
\hangindent=5.5em \hangafter=1
{\white .}\hskip 1em $\rho_\text{isum}(\mathfrak{s})$ =
$\mathrm{i}$($0$,
$\sqrt{\frac{1}{2}}$,
$\sqrt{\frac{1}{2}}$;\ \ 
$\frac{1}{2}$,
$-\frac{1}{2}$;\ \ 
$\frac{1}{2}$)
 $\oplus$
$\mathrm{i}$($0$,
$\sqrt{\frac{1}{2}}$,
$\sqrt{\frac{1}{2}}$;\ \ 
$-\frac{1}{2}$,
$\frac{1}{2}$;\ \ 
$-\frac{1}{2}$)

Fail:
Integral: $D_{\rho}(\sigma)_{\theta} \propto $ id,
 for all $\sigma$ and all $\theta$-eigenspaces that can contain unit. Prop. B.5 (6)

 \ \color{black}

\noindent 320: (dims,levels) = $(3 , 
3;24,
24
)$,
irreps = $3_{8}^{1,0}
\hskip -1.5pt \otimes \hskip -1.5pt
1_{3}^{1,0}\oplus
3_{8}^{3,3}
\hskip -1.5pt \otimes \hskip -1.5pt
1_{3}^{1,0}$,
pord$(\rho_\text{isum}(\mathfrak{t})) = 8$,

\vskip 0.7ex
\hangindent=5.5em \hangafter=1
{\white .}\hskip 1em $\rho_\text{isum}(\mathfrak{t})$ =
 $( \frac{1}{3},
\frac{11}{24},
\frac{23}{24} )
\oplus
( \frac{7}{12},
\frac{11}{24},
\frac{23}{24} )
$,

\vskip 0.7ex
\hangindent=5.5em \hangafter=1
{\white .}\hskip 1em $\rho_\text{isum}(\mathfrak{s})$ =
$\mathrm{i}$($0$,
$\sqrt{\frac{1}{2}}$,
$\sqrt{\frac{1}{2}}$;\ \ 
$-\frac{1}{2}$,
$\frac{1}{2}$;\ \ 
$-\frac{1}{2}$)
 $\oplus$
($0$,
$\sqrt{\frac{1}{2}}$,
$\sqrt{\frac{1}{2}}$;
$-\frac{1}{2}$,
$\frac{1}{2}$;
$-\frac{1}{2}$)

Fail:
number of self dual objects $|$Tr($\rho(\mathfrak s^2)$)$|$ = 0. Prop. B.4 (1)\
 eqn. (B.16)

 \ \color{black}

\noindent 321: (dims,levels) = $(3 , 
3;24,
24
)$,
irreps = $3_{8}^{1,0}
\hskip -1.5pt \otimes \hskip -1.5pt
1_{3}^{1,0}\oplus
3_{8}^{1,6}
\hskip -1.5pt \otimes \hskip -1.5pt
1_{3}^{1,0}$,
pord$(\rho_\text{isum}(\mathfrak{t})) = 8$,

\vskip 0.7ex
\hangindent=5.5em \hangafter=1
{\white .}\hskip 1em $\rho_\text{isum}(\mathfrak{t})$ =
 $( \frac{1}{3},
\frac{11}{24},
\frac{23}{24} )
\oplus
( \frac{5}{6},
\frac{11}{24},
\frac{23}{24} )
$,

\vskip 0.7ex
\hangindent=5.5em \hangafter=1
{\white .}\hskip 1em $\rho_\text{isum}(\mathfrak{s})$ =
$\mathrm{i}$($0$,
$\sqrt{\frac{1}{2}}$,
$\sqrt{\frac{1}{2}}$;\ \ 
$-\frac{1}{2}$,
$\frac{1}{2}$;\ \ 
$-\frac{1}{2}$)
 $\oplus$
$\mathrm{i}$($0$,
$\sqrt{\frac{1}{2}}$,
$\sqrt{\frac{1}{2}}$;\ \ 
$\frac{1}{2}$,
$-\frac{1}{2}$;\ \ 
$\frac{1}{2}$)

Fail:
Integral: $D_{\rho}(\sigma)_{\theta} \propto $ id,
 for all $\sigma$ and all $\theta$-eigenspaces that can contain unit. Prop. B.5 (6)

 \ \color{black}

\noindent 322: (dims,levels) = $(3 , 
3;24,
48
)$,
irreps = $3_{8}^{3,0}
\hskip -1.5pt \otimes \hskip -1.5pt
1_{3}^{1,0}\oplus
3_{16}^{3,6}
\hskip -1.5pt \otimes \hskip -1.5pt
1_{3}^{1,0}$,
pord$(\rho_\text{isum}(\mathfrak{t})) = 16$,

\vskip 0.7ex
\hangindent=5.5em \hangafter=1
{\white .}\hskip 1em $\rho_\text{isum}(\mathfrak{t})$ =
 $( \frac{1}{3},
\frac{5}{24},
\frac{17}{24} )
\oplus
( \frac{5}{24},
\frac{1}{48},
\frac{25}{48} )
$,

\vskip 0.7ex
\hangindent=5.5em \hangafter=1
{\white .}\hskip 1em $\rho_\text{isum}(\mathfrak{s})$ =
$\mathrm{i}$($0$,
$\sqrt{\frac{1}{2}}$,
$\sqrt{\frac{1}{2}}$;\ \ 
$\frac{1}{2}$,
$-\frac{1}{2}$;\ \ 
$\frac{1}{2}$)
 $\oplus$
$\mathrm{i}$($0$,
$\sqrt{\frac{1}{2}}$,
$\sqrt{\frac{1}{2}}$;\ \ 
$-\frac{1}{2}$,
$\frac{1}{2}$;\ \ 
$-\frac{1}{2}$)

Fail:
$\sigma(\rho(\mathfrak s)_\mathrm{ndeg}) \neq
 (\rho(\mathfrak t)^a \rho(\mathfrak s) \rho(\mathfrak t)^b
 \rho(\mathfrak s) \rho(\mathfrak t)^a)_\mathrm{ndeg}$,
 $\sigma = a$ = 5. Prop. B.5 (3) eqn. (B.25)

 \ \color{black}

\noindent 323: (dims,levels) = $(3 , 
3;24,
48
)$,
irreps = $3_{8}^{3,0}
\hskip -1.5pt \otimes \hskip -1.5pt
1_{3}^{1,0}\oplus
3_{16}^{1,9}
\hskip -1.5pt \otimes \hskip -1.5pt
1_{3}^{1,0}$,
pord$(\rho_\text{isum}(\mathfrak{t})) = 16$,

\vskip 0.7ex
\hangindent=5.5em \hangafter=1
{\white .}\hskip 1em $\rho_\text{isum}(\mathfrak{t})$ =
 $( \frac{1}{3},
\frac{5}{24},
\frac{17}{24} )
\oplus
( \frac{5}{24},
\frac{7}{48},
\frac{31}{48} )
$,

\vskip 0.7ex
\hangindent=5.5em \hangafter=1
{\white .}\hskip 1em $\rho_\text{isum}(\mathfrak{s})$ =
$\mathrm{i}$($0$,
$\sqrt{\frac{1}{2}}$,
$\sqrt{\frac{1}{2}}$;\ \ 
$\frac{1}{2}$,
$-\frac{1}{2}$;\ \ 
$\frac{1}{2}$)
 $\oplus$
($0$,
$\sqrt{\frac{1}{2}}$,
$\sqrt{\frac{1}{2}}$;
$-\frac{1}{2}$,
$\frac{1}{2}$;
$-\frac{1}{2}$)

Fail:
number of self dual objects $|$Tr($\rho(\mathfrak s^2)$)$|$ = 0. Prop. B.4 (1)\
 eqn. (B.16)

 \ \color{black}

\noindent 324: (dims,levels) = $(3 , 
3;24,
48
)$,
irreps = $3_{8}^{3,0}
\hskip -1.5pt \otimes \hskip -1.5pt
1_{3}^{1,0}\oplus
3_{16}^{7,0}
\hskip -1.5pt \otimes \hskip -1.5pt
1_{3}^{1,0}$,
pord$(\rho_\text{isum}(\mathfrak{t})) = 16$,

\vskip 0.7ex
\hangindent=5.5em \hangafter=1
{\white .}\hskip 1em $\rho_\text{isum}(\mathfrak{t})$ =
 $( \frac{1}{3},
\frac{5}{24},
\frac{17}{24} )
\oplus
( \frac{5}{24},
\frac{13}{48},
\frac{37}{48} )
$,

\vskip 0.7ex
\hangindent=5.5em \hangafter=1
{\white .}\hskip 1em $\rho_\text{isum}(\mathfrak{s})$ =
$\mathrm{i}$($0$,
$\sqrt{\frac{1}{2}}$,
$\sqrt{\frac{1}{2}}$;\ \ 
$\frac{1}{2}$,
$-\frac{1}{2}$;\ \ 
$\frac{1}{2}$)
 $\oplus$
$\mathrm{i}$($0$,
$\sqrt{\frac{1}{2}}$,
$\sqrt{\frac{1}{2}}$;\ \ 
$\frac{1}{2}$,
$-\frac{1}{2}$;\ \ 
$\frac{1}{2}$)

Fail:
$\sigma(\rho(\mathfrak s)_\mathrm{ndeg}) \neq
 (\rho(\mathfrak t)^a \rho(\mathfrak s) \rho(\mathfrak t)^b
 \rho(\mathfrak s) \rho(\mathfrak t)^a)_\mathrm{ndeg}$,
 $\sigma = a$ = 5. Prop. B.5 (3) eqn. (B.25)

 \ \color{black}

\noindent 325: (dims,levels) = $(3 , 
3;24,
48
)$,
irreps = $3_{8}^{3,0}
\hskip -1.5pt \otimes \hskip -1.5pt
1_{3}^{1,0}\oplus
3_{16}^{5,3}
\hskip -1.5pt \otimes \hskip -1.5pt
1_{3}^{1,0}$,
pord$(\rho_\text{isum}(\mathfrak{t})) = 16$,

\vskip 0.7ex
\hangindent=5.5em \hangafter=1
{\white .}\hskip 1em $\rho_\text{isum}(\mathfrak{t})$ =
 $( \frac{1}{3},
\frac{5}{24},
\frac{17}{24} )
\oplus
( \frac{5}{24},
\frac{19}{48},
\frac{43}{48} )
$,

\vskip 0.7ex
\hangindent=5.5em \hangafter=1
{\white .}\hskip 1em $\rho_\text{isum}(\mathfrak{s})$ =
$\mathrm{i}$($0$,
$\sqrt{\frac{1}{2}}$,
$\sqrt{\frac{1}{2}}$;\ \ 
$\frac{1}{2}$,
$-\frac{1}{2}$;\ \ 
$\frac{1}{2}$)
 $\oplus$
($0$,
$\sqrt{\frac{1}{2}}$,
$\sqrt{\frac{1}{2}}$;
$\frac{1}{2}$,
$-\frac{1}{2}$;
$\frac{1}{2}$)

Fail:
number of self dual objects $|$Tr($\rho(\mathfrak s^2)$)$|$ = 0. Prop. B.4 (1)\
 eqn. (B.16)

 \ \color{black}

\noindent 326: (dims,levels) = $(3 , 
3;24,
48
)$,
irreps = $3_{8}^{3,0}
\hskip -1.5pt \otimes \hskip -1.5pt
1_{3}^{1,0}\oplus
3_{16}^{3,0}
\hskip -1.5pt \otimes \hskip -1.5pt
1_{3}^{1,0}$,
pord$(\rho_\text{isum}(\mathfrak{t})) = 16$,

\vskip 0.7ex
\hangindent=5.5em \hangafter=1
{\white .}\hskip 1em $\rho_\text{isum}(\mathfrak{t})$ =
 $( \frac{1}{3},
\frac{5}{24},
\frac{17}{24} )
\oplus
( \frac{17}{24},
\frac{1}{48},
\frac{25}{48} )
$,

\vskip 0.7ex
\hangindent=5.5em \hangafter=1
{\white .}\hskip 1em $\rho_\text{isum}(\mathfrak{s})$ =
$\mathrm{i}$($0$,
$\sqrt{\frac{1}{2}}$,
$\sqrt{\frac{1}{2}}$;\ \ 
$\frac{1}{2}$,
$-\frac{1}{2}$;\ \ 
$\frac{1}{2}$)
 $\oplus$
$\mathrm{i}$($0$,
$\sqrt{\frac{1}{2}}$,
$\sqrt{\frac{1}{2}}$;\ \ 
$\frac{1}{2}$,
$-\frac{1}{2}$;\ \ 
$\frac{1}{2}$)

Fail:
$\sigma(\rho(\mathfrak s)_\mathrm{ndeg}) \neq
 (\rho(\mathfrak t)^a \rho(\mathfrak s) \rho(\mathfrak t)^b
 \rho(\mathfrak s) \rho(\mathfrak t)^a)_\mathrm{ndeg}$,
 $\sigma = a$ = 5. Prop. B.5 (3) eqn. (B.25)

 \ \color{black}

\noindent 327: (dims,levels) = $(3 , 
3;24,
48
)$,
irreps = $3_{8}^{3,0}
\hskip -1.5pt \otimes \hskip -1.5pt
1_{3}^{1,0}\oplus
3_{16}^{1,3}
\hskip -1.5pt \otimes \hskip -1.5pt
1_{3}^{1,0}$,
pord$(\rho_\text{isum}(\mathfrak{t})) = 16$,

\vskip 0.7ex
\hangindent=5.5em \hangafter=1
{\white .}\hskip 1em $\rho_\text{isum}(\mathfrak{t})$ =
 $( \frac{1}{3},
\frac{5}{24},
\frac{17}{24} )
\oplus
( \frac{17}{24},
\frac{7}{48},
\frac{31}{48} )
$,

\vskip 0.7ex
\hangindent=5.5em \hangafter=1
{\white .}\hskip 1em $\rho_\text{isum}(\mathfrak{s})$ =
$\mathrm{i}$($0$,
$\sqrt{\frac{1}{2}}$,
$\sqrt{\frac{1}{2}}$;\ \ 
$\frac{1}{2}$,
$-\frac{1}{2}$;\ \ 
$\frac{1}{2}$)
 $\oplus$
($0$,
$\sqrt{\frac{1}{2}}$,
$\sqrt{\frac{1}{2}}$;
$\frac{1}{2}$,
$-\frac{1}{2}$;
$\frac{1}{2}$)

Fail:
number of self dual objects $|$Tr($\rho(\mathfrak s^2)$)$|$ = 0. Prop. B.4 (1)\
 eqn. (B.16)

 \ \color{black}

\noindent 328: (dims,levels) = $(3 , 
3;24,
48
)$,
irreps = $3_{8}^{3,0}
\hskip -1.5pt \otimes \hskip -1.5pt
1_{3}^{1,0}\oplus
3_{16}^{7,6}
\hskip -1.5pt \otimes \hskip -1.5pt
1_{3}^{1,0}$,
pord$(\rho_\text{isum}(\mathfrak{t})) = 16$,

\vskip 0.7ex
\hangindent=5.5em \hangafter=1
{\white .}\hskip 1em $\rho_\text{isum}(\mathfrak{t})$ =
 $( \frac{1}{3},
\frac{5}{24},
\frac{17}{24} )
\oplus
( \frac{17}{24},
\frac{13}{48},
\frac{37}{48} )
$,

\vskip 0.7ex
\hangindent=5.5em \hangafter=1
{\white .}\hskip 1em $\rho_\text{isum}(\mathfrak{s})$ =
$\mathrm{i}$($0$,
$\sqrt{\frac{1}{2}}$,
$\sqrt{\frac{1}{2}}$;\ \ 
$\frac{1}{2}$,
$-\frac{1}{2}$;\ \ 
$\frac{1}{2}$)
 $\oplus$
$\mathrm{i}$($0$,
$\sqrt{\frac{1}{2}}$,
$\sqrt{\frac{1}{2}}$;\ \ 
$-\frac{1}{2}$,
$\frac{1}{2}$;\ \ 
$-\frac{1}{2}$)

Fail:
$\sigma(\rho(\mathfrak s)_\mathrm{ndeg}) \neq
 (\rho(\mathfrak t)^a \rho(\mathfrak s) \rho(\mathfrak t)^b
 \rho(\mathfrak s) \rho(\mathfrak t)^a)_\mathrm{ndeg}$,
 $\sigma = a$ = 5. Prop. B.5 (3) eqn. (B.25)

 \ \color{black}

\noindent 329: (dims,levels) = $(3 , 
3;24,
48
)$,
irreps = $3_{8}^{3,0}
\hskip -1.5pt \otimes \hskip -1.5pt
1_{3}^{1,0}\oplus
3_{16}^{5,9}
\hskip -1.5pt \otimes \hskip -1.5pt
1_{3}^{1,0}$,
pord$(\rho_\text{isum}(\mathfrak{t})) = 16$,

\vskip 0.7ex
\hangindent=5.5em \hangafter=1
{\white .}\hskip 1em $\rho_\text{isum}(\mathfrak{t})$ =
 $( \frac{1}{3},
\frac{5}{24},
\frac{17}{24} )
\oplus
( \frac{17}{24},
\frac{19}{48},
\frac{43}{48} )
$,

\vskip 0.7ex
\hangindent=5.5em \hangafter=1
{\white .}\hskip 1em $\rho_\text{isum}(\mathfrak{s})$ =
$\mathrm{i}$($0$,
$\sqrt{\frac{1}{2}}$,
$\sqrt{\frac{1}{2}}$;\ \ 
$\frac{1}{2}$,
$-\frac{1}{2}$;\ \ 
$\frac{1}{2}$)
 $\oplus$
($0$,
$\sqrt{\frac{1}{2}}$,
$\sqrt{\frac{1}{2}}$;
$-\frac{1}{2}$,
$\frac{1}{2}$;
$-\frac{1}{2}$)

Fail:
number of self dual objects $|$Tr($\rho(\mathfrak s^2)$)$|$ = 0. Prop. B.4 (1)\
 eqn. (B.16)

 \ \color{black}

\noindent 330: (dims,levels) = $(3 , 
3;24,
48
)$,
irreps = $3_{8}^{1,0}
\hskip -1.5pt \otimes \hskip -1.5pt
1_{3}^{1,0}\oplus
3_{16}^{7,3}
\hskip -1.5pt \otimes \hskip -1.5pt
1_{3}^{1,0}$,
pord$(\rho_\text{isum}(\mathfrak{t})) = 16$,

\vskip 0.7ex
\hangindent=5.5em \hangafter=1
{\white .}\hskip 1em $\rho_\text{isum}(\mathfrak{t})$ =
 $( \frac{1}{3},
\frac{11}{24},
\frac{23}{24} )
\oplus
( \frac{11}{24},
\frac{1}{48},
\frac{25}{48} )
$,

\vskip 0.7ex
\hangindent=5.5em \hangafter=1
{\white .}\hskip 1em $\rho_\text{isum}(\mathfrak{s})$ =
$\mathrm{i}$($0$,
$\sqrt{\frac{1}{2}}$,
$\sqrt{\frac{1}{2}}$;\ \ 
$-\frac{1}{2}$,
$\frac{1}{2}$;\ \ 
$-\frac{1}{2}$)
 $\oplus$
($0$,
$\sqrt{\frac{1}{2}}$,
$\sqrt{\frac{1}{2}}$;
$-\frac{1}{2}$,
$\frac{1}{2}$;
$-\frac{1}{2}$)

Fail:
number of self dual objects $|$Tr($\rho(\mathfrak s^2)$)$|$ = 0. Prop. B.4 (1)\
 eqn. (B.16)

 \ \color{black}

\noindent 331: (dims,levels) = $(3 , 
3;24,
48
)$,
irreps = $3_{8}^{1,0}
\hskip -1.5pt \otimes \hskip -1.5pt
1_{3}^{1,0}\oplus
3_{16}^{5,6}
\hskip -1.5pt \otimes \hskip -1.5pt
1_{3}^{1,0}$,
pord$(\rho_\text{isum}(\mathfrak{t})) = 16$,

\vskip 0.7ex
\hangindent=5.5em \hangafter=1
{\white .}\hskip 1em $\rho_\text{isum}(\mathfrak{t})$ =
 $( \frac{1}{3},
\frac{11}{24},
\frac{23}{24} )
\oplus
( \frac{11}{24},
\frac{7}{48},
\frac{31}{48} )
$,

\vskip 0.7ex
\hangindent=5.5em \hangafter=1
{\white .}\hskip 1em $\rho_\text{isum}(\mathfrak{s})$ =
$\mathrm{i}$($0$,
$\sqrt{\frac{1}{2}}$,
$\sqrt{\frac{1}{2}}$;\ \ 
$-\frac{1}{2}$,
$\frac{1}{2}$;\ \ 
$-\frac{1}{2}$)
 $\oplus$
$\mathrm{i}$($0$,
$\sqrt{\frac{1}{2}}$,
$\sqrt{\frac{1}{2}}$;\ \ 
$\frac{1}{2}$,
$-\frac{1}{2}$;\ \ 
$\frac{1}{2}$)

Fail:
$\sigma(\rho(\mathfrak s)_\mathrm{ndeg}) \neq
 (\rho(\mathfrak t)^a \rho(\mathfrak s) \rho(\mathfrak t)^b
 \rho(\mathfrak s) \rho(\mathfrak t)^a)_\mathrm{ndeg}$,
 $\sigma = a$ = 5. Prop. B.5 (3) eqn. (B.25)

 \ \color{black}

\noindent 332: (dims,levels) = $(3 , 
3;24,
48
)$,
irreps = $3_{8}^{1,0}
\hskip -1.5pt \otimes \hskip -1.5pt
1_{3}^{1,0}\oplus
3_{16}^{3,9}
\hskip -1.5pt \otimes \hskip -1.5pt
1_{3}^{1,0}$,
pord$(\rho_\text{isum}(\mathfrak{t})) = 16$,

\vskip 0.7ex
\hangindent=5.5em \hangafter=1
{\white .}\hskip 1em $\rho_\text{isum}(\mathfrak{t})$ =
 $( \frac{1}{3},
\frac{11}{24},
\frac{23}{24} )
\oplus
( \frac{11}{24},
\frac{13}{48},
\frac{37}{48} )
$,

\vskip 0.7ex
\hangindent=5.5em \hangafter=1
{\white .}\hskip 1em $\rho_\text{isum}(\mathfrak{s})$ =
$\mathrm{i}$($0$,
$\sqrt{\frac{1}{2}}$,
$\sqrt{\frac{1}{2}}$;\ \ 
$-\frac{1}{2}$,
$\frac{1}{2}$;\ \ 
$-\frac{1}{2}$)
 $\oplus$
($0$,
$\sqrt{\frac{1}{2}}$,
$\sqrt{\frac{1}{2}}$;
$\frac{1}{2}$,
$-\frac{1}{2}$;
$\frac{1}{2}$)

Fail:
number of self dual objects $|$Tr($\rho(\mathfrak s^2)$)$|$ = 0. Prop. B.4 (1)\
 eqn. (B.16)

 \ \color{black}

\noindent 333: (dims,levels) = $(3 , 
3;24,
48
)$,
irreps = $3_{8}^{1,0}
\hskip -1.5pt \otimes \hskip -1.5pt
1_{3}^{1,0}\oplus
3_{16}^{1,0}
\hskip -1.5pt \otimes \hskip -1.5pt
1_{3}^{1,0}$,
pord$(\rho_\text{isum}(\mathfrak{t})) = 16$,

\vskip 0.7ex
\hangindent=5.5em \hangafter=1
{\white .}\hskip 1em $\rho_\text{isum}(\mathfrak{t})$ =
 $( \frac{1}{3},
\frac{11}{24},
\frac{23}{24} )
\oplus
( \frac{11}{24},
\frac{19}{48},
\frac{43}{48} )
$,

\vskip 0.7ex
\hangindent=5.5em \hangafter=1
{\white .}\hskip 1em $\rho_\text{isum}(\mathfrak{s})$ =
$\mathrm{i}$($0$,
$\sqrt{\frac{1}{2}}$,
$\sqrt{\frac{1}{2}}$;\ \ 
$-\frac{1}{2}$,
$\frac{1}{2}$;\ \ 
$-\frac{1}{2}$)
 $\oplus$
$\mathrm{i}$($0$,
$\sqrt{\frac{1}{2}}$,
$\sqrt{\frac{1}{2}}$;\ \ 
$-\frac{1}{2}$,
$\frac{1}{2}$;\ \ 
$-\frac{1}{2}$)

Fail:
$\sigma(\rho(\mathfrak s)_\mathrm{ndeg}) \neq
 (\rho(\mathfrak t)^a \rho(\mathfrak s) \rho(\mathfrak t)^b
 \rho(\mathfrak s) \rho(\mathfrak t)^a)_\mathrm{ndeg}$,
 $\sigma = a$ = 5. Prop. B.5 (3) eqn. (B.25)

 \ \color{black}

\noindent 334: (dims,levels) = $(3 , 
3;24,
48
)$,
irreps = $3_{8}^{1,0}
\hskip -1.5pt \otimes \hskip -1.5pt
1_{3}^{1,0}\oplus
3_{16}^{7,9}
\hskip -1.5pt \otimes \hskip -1.5pt
1_{3}^{1,0}$,
pord$(\rho_\text{isum}(\mathfrak{t})) = 16$,

\vskip 0.7ex
\hangindent=5.5em \hangafter=1
{\white .}\hskip 1em $\rho_\text{isum}(\mathfrak{t})$ =
 $( \frac{1}{3},
\frac{11}{24},
\frac{23}{24} )
\oplus
( \frac{23}{24},
\frac{1}{48},
\frac{25}{48} )
$,

\vskip 0.7ex
\hangindent=5.5em \hangafter=1
{\white .}\hskip 1em $\rho_\text{isum}(\mathfrak{s})$ =
$\mathrm{i}$($0$,
$\sqrt{\frac{1}{2}}$,
$\sqrt{\frac{1}{2}}$;\ \ 
$-\frac{1}{2}$,
$\frac{1}{2}$;\ \ 
$-\frac{1}{2}$)
 $\oplus$
($0$,
$\sqrt{\frac{1}{2}}$,
$\sqrt{\frac{1}{2}}$;
$\frac{1}{2}$,
$-\frac{1}{2}$;
$\frac{1}{2}$)

Fail:
number of self dual objects $|$Tr($\rho(\mathfrak s^2)$)$|$ = 0. Prop. B.4 (1)\
 eqn. (B.16)

 \ \color{black}

\noindent 335: (dims,levels) = $(3 , 
3;24,
48
)$,
irreps = $3_{8}^{1,0}
\hskip -1.5pt \otimes \hskip -1.5pt
1_{3}^{1,0}\oplus
3_{16}^{5,0}
\hskip -1.5pt \otimes \hskip -1.5pt
1_{3}^{1,0}$,
pord$(\rho_\text{isum}(\mathfrak{t})) = 16$,

\vskip 0.7ex
\hangindent=5.5em \hangafter=1
{\white .}\hskip 1em $\rho_\text{isum}(\mathfrak{t})$ =
 $( \frac{1}{3},
\frac{11}{24},
\frac{23}{24} )
\oplus
( \frac{23}{24},
\frac{7}{48},
\frac{31}{48} )
$,

\vskip 0.7ex
\hangindent=5.5em \hangafter=1
{\white .}\hskip 1em $\rho_\text{isum}(\mathfrak{s})$ =
$\mathrm{i}$($0$,
$\sqrt{\frac{1}{2}}$,
$\sqrt{\frac{1}{2}}$;\ \ 
$-\frac{1}{2}$,
$\frac{1}{2}$;\ \ 
$-\frac{1}{2}$)
 $\oplus$
$\mathrm{i}$($0$,
$\sqrt{\frac{1}{2}}$,
$\sqrt{\frac{1}{2}}$;\ \ 
$-\frac{1}{2}$,
$\frac{1}{2}$;\ \ 
$-\frac{1}{2}$)

Fail:
$\sigma(\rho(\mathfrak s)_\mathrm{ndeg}) \neq
 (\rho(\mathfrak t)^a \rho(\mathfrak s) \rho(\mathfrak t)^b
 \rho(\mathfrak s) \rho(\mathfrak t)^a)_\mathrm{ndeg}$,
 $\sigma = a$ = 5. Prop. B.5 (3) eqn. (B.25)

 \ \color{black}

\noindent 336: (dims,levels) = $(3 , 
3;24,
48
)$,
irreps = $3_{8}^{1,0}
\hskip -1.5pt \otimes \hskip -1.5pt
1_{3}^{1,0}\oplus
3_{16}^{3,3}
\hskip -1.5pt \otimes \hskip -1.5pt
1_{3}^{1,0}$,
pord$(\rho_\text{isum}(\mathfrak{t})) = 16$,

\vskip 0.7ex
\hangindent=5.5em \hangafter=1
{\white .}\hskip 1em $\rho_\text{isum}(\mathfrak{t})$ =
 $( \frac{1}{3},
\frac{11}{24},
\frac{23}{24} )
\oplus
( \frac{23}{24},
\frac{13}{48},
\frac{37}{48} )
$,

\vskip 0.7ex
\hangindent=5.5em \hangafter=1
{\white .}\hskip 1em $\rho_\text{isum}(\mathfrak{s})$ =
$\mathrm{i}$($0$,
$\sqrt{\frac{1}{2}}$,
$\sqrt{\frac{1}{2}}$;\ \ 
$-\frac{1}{2}$,
$\frac{1}{2}$;\ \ 
$-\frac{1}{2}$)
 $\oplus$
($0$,
$\sqrt{\frac{1}{2}}$,
$\sqrt{\frac{1}{2}}$;
$-\frac{1}{2}$,
$\frac{1}{2}$;
$-\frac{1}{2}$)

Fail:
number of self dual objects $|$Tr($\rho(\mathfrak s^2)$)$|$ = 0. Prop. B.4 (1)\
 eqn. (B.16)

 \ \color{black}

\noindent 337: (dims,levels) = $(3 , 
3;24,
48
)$,
irreps = $3_{8}^{1,0}
\hskip -1.5pt \otimes \hskip -1.5pt
1_{3}^{1,0}\oplus
3_{16}^{1,6}
\hskip -1.5pt \otimes \hskip -1.5pt
1_{3}^{1,0}$,
pord$(\rho_\text{isum}(\mathfrak{t})) = 16$,

\vskip 0.7ex
\hangindent=5.5em \hangafter=1
{\white .}\hskip 1em $\rho_\text{isum}(\mathfrak{t})$ =
 $( \frac{1}{3},
\frac{11}{24},
\frac{23}{24} )
\oplus
( \frac{23}{24},
\frac{19}{48},
\frac{43}{48} )
$,

\vskip 0.7ex
\hangindent=5.5em \hangafter=1
{\white .}\hskip 1em $\rho_\text{isum}(\mathfrak{s})$ =
$\mathrm{i}$($0$,
$\sqrt{\frac{1}{2}}$,
$\sqrt{\frac{1}{2}}$;\ \ 
$-\frac{1}{2}$,
$\frac{1}{2}$;\ \ 
$-\frac{1}{2}$)
 $\oplus$
$\mathrm{i}$($0$,
$\sqrt{\frac{1}{2}}$,
$\sqrt{\frac{1}{2}}$;\ \ 
$\frac{1}{2}$,
$-\frac{1}{2}$;\ \ 
$\frac{1}{2}$)

Fail:
$\sigma(\rho(\mathfrak s)_\mathrm{ndeg}) \neq
 (\rho(\mathfrak t)^a \rho(\mathfrak s) \rho(\mathfrak t)^b
 \rho(\mathfrak s) \rho(\mathfrak t)^a)_\mathrm{ndeg}$,
 $\sigma = a$ = 5. Prop. B.5 (3) eqn. (B.25)

 \ \color{black}

\noindent 338: (dims,levels) = $(3 , 
3;30,
30
)$,
irreps = $3_{5}^{1}
\hskip -1.5pt \otimes \hskip -1.5pt
1_{3}^{1,0}
\hskip -1.5pt \otimes \hskip -1.5pt
1_{2}^{1,0}\oplus
3_{5}^{3}
\hskip -1.5pt \otimes \hskip -1.5pt
1_{3}^{1,0}
\hskip -1.5pt \otimes \hskip -1.5pt
1_{2}^{1,0}$,
pord$(\rho_\text{isum}(\mathfrak{t})) = 5$,

\vskip 0.7ex
\hangindent=5.5em \hangafter=1
{\white .}\hskip 1em $\rho_\text{isum}(\mathfrak{t})$ =
 $( \frac{5}{6},
\frac{1}{30},
\frac{19}{30} )
\oplus
( \frac{5}{6},
\frac{7}{30},
\frac{13}{30} )
$,

\vskip 0.7ex
\hangindent=5.5em \hangafter=1
{\white .}\hskip 1em $\rho_\text{isum}(\mathfrak{s})$ =
($-\sqrt{\frac{1}{5}}$,
$-\sqrt{\frac{2}{5}}$,
$-\sqrt{\frac{2}{5}}$;
$\frac{5+\sqrt{5}}{10}$,
$\frac{-5+\sqrt{5}}{10}$;
$\frac{5+\sqrt{5}}{10}$)
 $\oplus$
($\sqrt{\frac{1}{5}}$,
$-\sqrt{\frac{2}{5}}$,
$-\sqrt{\frac{2}{5}}$;
$\frac{5-\sqrt{5}}{10}$,
$-\frac{5+\sqrt{5}}{10}$;
$\frac{5-\sqrt{5}}{10}$)

Fail:
Integral: $D_{\rho}(\sigma)_{\theta} \propto $ id,
 for all $\sigma$ and all $\theta$-eigenspaces that can contain unit. Prop. B.5 (6)

 \ \color{black}

\noindent 339: (dims,levels) = $(3 , 
3;48,
48
)$,
irreps = $3_{16}^{7,0}
\hskip -1.5pt \otimes \hskip -1.5pt
1_{3}^{1,0}\oplus
3_{16}^{5,3}
\hskip -1.5pt \otimes \hskip -1.5pt
1_{3}^{1,0}$,
pord$(\rho_\text{isum}(\mathfrak{t})) = 16$,

\vskip 0.7ex
\hangindent=5.5em \hangafter=1
{\white .}\hskip 1em $\rho_\text{isum}(\mathfrak{t})$ =
 $( \frac{5}{24},
\frac{13}{48},
\frac{37}{48} )
\oplus
( \frac{5}{24},
\frac{19}{48},
\frac{43}{48} )
$,

\vskip 0.7ex
\hangindent=5.5em \hangafter=1
{\white .}\hskip 1em $\rho_\text{isum}(\mathfrak{s})$ =
$\mathrm{i}$($0$,
$\sqrt{\frac{1}{2}}$,
$\sqrt{\frac{1}{2}}$;\ \ 
$\frac{1}{2}$,
$-\frac{1}{2}$;\ \ 
$\frac{1}{2}$)
 $\oplus$
($0$,
$\sqrt{\frac{1}{2}}$,
$\sqrt{\frac{1}{2}}$;
$\frac{1}{2}$,
$-\frac{1}{2}$;
$\frac{1}{2}$)

Fail:
number of self dual objects $|$Tr($\rho(\mathfrak s^2)$)$|$ = 0. Prop. B.4 (1)\
 eqn. (B.16)

 \ \color{black}

\noindent 340: (dims,levels) = $(3 , 
3;48,
48
)$,
irreps = $3_{16}^{7,0}
\hskip -1.5pt \otimes \hskip -1.5pt
1_{3}^{1,0}\oplus
3_{16}^{3,9}
\hskip -1.5pt \otimes \hskip -1.5pt
1_{3}^{1,0}$,
pord$(\rho_\text{isum}(\mathfrak{t})) = 16$,

\vskip 0.7ex
\hangindent=5.5em \hangafter=1
{\white .}\hskip 1em $\rho_\text{isum}(\mathfrak{t})$ =
 $( \frac{5}{24},
\frac{13}{48},
\frac{37}{48} )
\oplus
( \frac{11}{24},
\frac{13}{48},
\frac{37}{48} )
$,

\vskip 0.7ex
\hangindent=5.5em \hangafter=1
{\white .}\hskip 1em $\rho_\text{isum}(\mathfrak{s})$ =
$\mathrm{i}$($0$,
$\sqrt{\frac{1}{2}}$,
$\sqrt{\frac{1}{2}}$;\ \ 
$\frac{1}{2}$,
$-\frac{1}{2}$;\ \ 
$\frac{1}{2}$)
 $\oplus$
($0$,
$\sqrt{\frac{1}{2}}$,
$\sqrt{\frac{1}{2}}$;
$\frac{1}{2}$,
$-\frac{1}{2}$;
$\frac{1}{2}$)

Fail:
number of self dual objects $|$Tr($\rho(\mathfrak s^2)$)$|$ = 0. Prop. B.4 (1)\
 eqn. (B.16)

 \ \color{black}

 \color{blue}

\noindent 341: (dims,levels) = $(3 , 
3;48,
48
)$,
irreps = $3_{16}^{7,0}
\hskip -1.5pt \otimes \hskip -1.5pt
1_{3}^{1,0}\oplus
3_{16}^{7,6}
\hskip -1.5pt \otimes \hskip -1.5pt
1_{3}^{1,0}$,
pord$(\rho_\text{isum}(\mathfrak{t})) = 16$,

\vskip 0.7ex
\hangindent=5.5em \hangafter=1
{\white .}\hskip 1em $\rho_\text{isum}(\mathfrak{t})$ =
 $( \frac{5}{24},
\frac{13}{48},
\frac{37}{48} )
\oplus
( \frac{17}{24},
\frac{13}{48},
\frac{37}{48} )
$,

\vskip 0.7ex
\hangindent=5.5em \hangafter=1
{\white .}\hskip 1em $\rho_\text{isum}(\mathfrak{s})$ =
$\mathrm{i}$($0$,
$\sqrt{\frac{1}{2}}$,
$\sqrt{\frac{1}{2}}$;\ \ 
$\frac{1}{2}$,
$-\frac{1}{2}$;\ \ 
$\frac{1}{2}$)
 $\oplus$
$\mathrm{i}$($0$,
$\sqrt{\frac{1}{2}}$,
$\sqrt{\frac{1}{2}}$;\ \ 
$-\frac{1}{2}$,
$\frac{1}{2}$;\ \ 
$-\frac{1}{2}$)

Pass. 

 \ \color{black}

\noindent 342: (dims,levels) = $(3 , 
3;48,
48
)$,
irreps = $3_{16}^{7,0}
\hskip -1.5pt \otimes \hskip -1.5pt
1_{3}^{1,0}\oplus
3_{16}^{3,3}
\hskip -1.5pt \otimes \hskip -1.5pt
1_{3}^{1,0}$,
pord$(\rho_\text{isum}(\mathfrak{t})) = 16$,

\vskip 0.7ex
\hangindent=5.5em \hangafter=1
{\white .}\hskip 1em $\rho_\text{isum}(\mathfrak{t})$ =
 $( \frac{5}{24},
\frac{13}{48},
\frac{37}{48} )
\oplus
( \frac{23}{24},
\frac{13}{48},
\frac{37}{48} )
$,

\vskip 0.7ex
\hangindent=5.5em \hangafter=1
{\white .}\hskip 1em $\rho_\text{isum}(\mathfrak{s})$ =
$\mathrm{i}$($0$,
$\sqrt{\frac{1}{2}}$,
$\sqrt{\frac{1}{2}}$;\ \ 
$\frac{1}{2}$,
$-\frac{1}{2}$;\ \ 
$\frac{1}{2}$)
 $\oplus$
($0$,
$\sqrt{\frac{1}{2}}$,
$\sqrt{\frac{1}{2}}$;
$-\frac{1}{2}$,
$\frac{1}{2}$;
$-\frac{1}{2}$)

Fail:
number of self dual objects $|$Tr($\rho(\mathfrak s^2)$)$|$ = 0. Prop. B.4 (1)\
 eqn. (B.16)

 \ \color{black}

\noindent 343: (dims,levels) = $(3 , 
3;48,
48
)$,
irreps = $3_{16}^{1,0}
\hskip -1.5pt \otimes \hskip -1.5pt
1_{3}^{1,0}\oplus
3_{16}^{5,9}
\hskip -1.5pt \otimes \hskip -1.5pt
1_{3}^{1,0}$,
pord$(\rho_\text{isum}(\mathfrak{t})) = 16$,

\vskip 0.7ex
\hangindent=5.5em \hangafter=1
{\white .}\hskip 1em $\rho_\text{isum}(\mathfrak{t})$ =
 $( \frac{11}{24},
\frac{19}{48},
\frac{43}{48} )
\oplus
( \frac{17}{24},
\frac{19}{48},
\frac{43}{48} )
$,

\vskip 0.7ex
\hangindent=5.5em \hangafter=1
{\white .}\hskip 1em $\rho_\text{isum}(\mathfrak{s})$ =
$\mathrm{i}$($0$,
$\sqrt{\frac{1}{2}}$,
$\sqrt{\frac{1}{2}}$;\ \ 
$-\frac{1}{2}$,
$\frac{1}{2}$;\ \ 
$-\frac{1}{2}$)
 $\oplus$
($0$,
$\sqrt{\frac{1}{2}}$,
$\sqrt{\frac{1}{2}}$;
$-\frac{1}{2}$,
$\frac{1}{2}$;
$-\frac{1}{2}$)

Fail:
number of self dual objects $|$Tr($\rho(\mathfrak s^2)$)$|$ = 0. Prop. B.4 (1)\
 eqn. (B.16)

 \ \color{black}

 \color{blue}

\noindent 344: (dims,levels) = $(3 , 
3;48,
48
)$,
irreps = $3_{16}^{1,0}
\hskip -1.5pt \otimes \hskip -1.5pt
1_{3}^{1,0}\oplus
3_{16}^{1,6}
\hskip -1.5pt \otimes \hskip -1.5pt
1_{3}^{1,0}$,
pord$(\rho_\text{isum}(\mathfrak{t})) = 16$,

\vskip 0.7ex
\hangindent=5.5em \hangafter=1
{\white .}\hskip 1em $\rho_\text{isum}(\mathfrak{t})$ =
 $( \frac{11}{24},
\frac{19}{48},
\frac{43}{48} )
\oplus
( \frac{23}{24},
\frac{19}{48},
\frac{43}{48} )
$,

\vskip 0.7ex
\hangindent=5.5em \hangafter=1
{\white .}\hskip 1em $\rho_\text{isum}(\mathfrak{s})$ =
$\mathrm{i}$($0$,
$\sqrt{\frac{1}{2}}$,
$\sqrt{\frac{1}{2}}$;\ \ 
$-\frac{1}{2}$,
$\frac{1}{2}$;\ \ 
$-\frac{1}{2}$)
 $\oplus$
$\mathrm{i}$($0$,
$\sqrt{\frac{1}{2}}$,
$\sqrt{\frac{1}{2}}$;\ \ 
$\frac{1}{2}$,
$-\frac{1}{2}$;\ \ 
$\frac{1}{2}$)

Pass. 

 \ \color{black}

\noindent 345: (dims,levels) = $(3 , 
3;48,
48
)$,
irreps = $3_{16}^{3,0}
\hskip -1.5pt \otimes \hskip -1.5pt
1_{3}^{1,0}\oplus
3_{16}^{1,3}
\hskip -1.5pt \otimes \hskip -1.5pt
1_{3}^{1,0}$,
pord$(\rho_\text{isum}(\mathfrak{t})) = 16$,

\vskip 0.7ex
\hangindent=5.5em \hangafter=1
{\white .}\hskip 1em $\rho_\text{isum}(\mathfrak{t})$ =
 $( \frac{17}{24},
\frac{1}{48},
\frac{25}{48} )
\oplus
( \frac{17}{24},
\frac{7}{48},
\frac{31}{48} )
$,

\vskip 0.7ex
\hangindent=5.5em \hangafter=1
{\white .}\hskip 1em $\rho_\text{isum}(\mathfrak{s})$ =
$\mathrm{i}$($0$,
$\sqrt{\frac{1}{2}}$,
$\sqrt{\frac{1}{2}}$;\ \ 
$\frac{1}{2}$,
$-\frac{1}{2}$;\ \ 
$\frac{1}{2}$)
 $\oplus$
($0$,
$\sqrt{\frac{1}{2}}$,
$\sqrt{\frac{1}{2}}$;
$\frac{1}{2}$,
$-\frac{1}{2}$;
$\frac{1}{2}$)

Fail:
number of self dual objects $|$Tr($\rho(\mathfrak s^2)$)$|$ = 0. Prop. B.4 (1)\
 eqn. (B.16)

 \ \color{black}

\noindent 346: (dims,levels) = $(3 , 
3;48,
48
)$,
irreps = $3_{16}^{3,0}
\hskip -1.5pt \otimes \hskip -1.5pt
1_{3}^{1,0}\oplus
3_{16}^{7,6}
\hskip -1.5pt \otimes \hskip -1.5pt
1_{3}^{1,0}$,
pord$(\rho_\text{isum}(\mathfrak{t})) = 16$,

\vskip 0.7ex
\hangindent=5.5em \hangafter=1
{\white .}\hskip 1em $\rho_\text{isum}(\mathfrak{t})$ =
 $( \frac{17}{24},
\frac{1}{48},
\frac{25}{48} )
\oplus
( \frac{17}{24},
\frac{13}{48},
\frac{37}{48} )
$,

\vskip 0.7ex
\hangindent=5.5em \hangafter=1
{\white .}\hskip 1em $\rho_\text{isum}(\mathfrak{s})$ =
$\mathrm{i}$($0$,
$\sqrt{\frac{1}{2}}$,
$\sqrt{\frac{1}{2}}$;\ \ 
$\frac{1}{2}$,
$-\frac{1}{2}$;\ \ 
$\frac{1}{2}$)
 $\oplus$
$\mathrm{i}$($0$,
$\sqrt{\frac{1}{2}}$,
$\sqrt{\frac{1}{2}}$;\ \ 
$-\frac{1}{2}$,
$\frac{1}{2}$;\ \ 
$-\frac{1}{2}$)

Fail:
all rows of $U \rho(\mathfrak s) U^\dagger$
 contain zero for any block-diagonal $U$. Prop. B.5 (4) eqn. (B.27)

 \ \color{black}

\noindent 347: (dims,levels) = $(3 , 
3;48,
48
)$,
irreps = $3_{16}^{3,0}
\hskip -1.5pt \otimes \hskip -1.5pt
1_{3}^{1,0}\oplus
3_{16}^{5,9}
\hskip -1.5pt \otimes \hskip -1.5pt
1_{3}^{1,0}$,
pord$(\rho_\text{isum}(\mathfrak{t})) = 16$,

\vskip 0.7ex
\hangindent=5.5em \hangafter=1
{\white .}\hskip 1em $\rho_\text{isum}(\mathfrak{t})$ =
 $( \frac{17}{24},
\frac{1}{48},
\frac{25}{48} )
\oplus
( \frac{17}{24},
\frac{19}{48},
\frac{43}{48} )
$,

\vskip 0.7ex
\hangindent=5.5em \hangafter=1
{\white .}\hskip 1em $\rho_\text{isum}(\mathfrak{s})$ =
$\mathrm{i}$($0$,
$\sqrt{\frac{1}{2}}$,
$\sqrt{\frac{1}{2}}$;\ \ 
$\frac{1}{2}$,
$-\frac{1}{2}$;\ \ 
$\frac{1}{2}$)
 $\oplus$
($0$,
$\sqrt{\frac{1}{2}}$,
$\sqrt{\frac{1}{2}}$;
$-\frac{1}{2}$,
$\frac{1}{2}$;
$-\frac{1}{2}$)

Fail:
number of self dual objects $|$Tr($\rho(\mathfrak s^2)$)$|$ = 0. Prop. B.4 (1)\
 eqn. (B.16)

 \ \color{black}

\noindent 348: (dims,levels) = $(3 , 
3;48,
48
)$,
irreps = $3_{16}^{3,0}
\hskip -1.5pt \otimes \hskip -1.5pt
1_{3}^{1,0}\oplus
3_{16}^{7,9}
\hskip -1.5pt \otimes \hskip -1.5pt
1_{3}^{1,0}$,
pord$(\rho_\text{isum}(\mathfrak{t})) = 16$,

\vskip 0.7ex
\hangindent=5.5em \hangafter=1
{\white .}\hskip 1em $\rho_\text{isum}(\mathfrak{t})$ =
 $( \frac{17}{24},
\frac{1}{48},
\frac{25}{48} )
\oplus
( \frac{23}{24},
\frac{1}{48},
\frac{25}{48} )
$,

\vskip 0.7ex
\hangindent=5.5em \hangafter=1
{\white .}\hskip 1em $\rho_\text{isum}(\mathfrak{s})$ =
$\mathrm{i}$($0$,
$\sqrt{\frac{1}{2}}$,
$\sqrt{\frac{1}{2}}$;\ \ 
$\frac{1}{2}$,
$-\frac{1}{2}$;\ \ 
$\frac{1}{2}$)
 $\oplus$
($0$,
$\sqrt{\frac{1}{2}}$,
$\sqrt{\frac{1}{2}}$;
$\frac{1}{2}$,
$-\frac{1}{2}$;
$\frac{1}{2}$)

Fail:
number of self dual objects $|$Tr($\rho(\mathfrak s^2)$)$|$ = 0. Prop. B.4 (1)\
 eqn. (B.16)

 \ \color{black}

\noindent 349: (dims,levels) = $(3 , 
3;48,
48
)$,
irreps = $3_{16}^{5,0}
\hskip -1.5pt \otimes \hskip -1.5pt
1_{3}^{1,0}\oplus
3_{16}^{3,3}
\hskip -1.5pt \otimes \hskip -1.5pt
1_{3}^{1,0}$,
pord$(\rho_\text{isum}(\mathfrak{t})) = 16$,

\vskip 0.7ex
\hangindent=5.5em \hangafter=1
{\white .}\hskip 1em $\rho_\text{isum}(\mathfrak{t})$ =
 $( \frac{23}{24},
\frac{7}{48},
\frac{31}{48} )
\oplus
( \frac{23}{24},
\frac{13}{48},
\frac{37}{48} )
$,

\vskip 0.7ex
\hangindent=5.5em \hangafter=1
{\white .}\hskip 1em $\rho_\text{isum}(\mathfrak{s})$ =
$\mathrm{i}$($0$,
$\sqrt{\frac{1}{2}}$,
$\sqrt{\frac{1}{2}}$;\ \ 
$-\frac{1}{2}$,
$\frac{1}{2}$;\ \ 
$-\frac{1}{2}$)
 $\oplus$
($0$,
$\sqrt{\frac{1}{2}}$,
$\sqrt{\frac{1}{2}}$;
$-\frac{1}{2}$,
$\frac{1}{2}$;
$-\frac{1}{2}$)

Fail:
number of self dual objects $|$Tr($\rho(\mathfrak s^2)$)$|$ = 0. Prop. B.4 (1)\
 eqn. (B.16)

 \ \color{black}

\noindent 350: (dims,levels) = $(3 , 
3;48,
48
)$,
irreps = $3_{16}^{5,0}
\hskip -1.5pt \otimes \hskip -1.5pt
1_{3}^{1,0}\oplus
3_{16}^{1,6}
\hskip -1.5pt \otimes \hskip -1.5pt
1_{3}^{1,0}$,
pord$(\rho_\text{isum}(\mathfrak{t})) = 16$,

\vskip 0.7ex
\hangindent=5.5em \hangafter=1
{\white .}\hskip 1em $\rho_\text{isum}(\mathfrak{t})$ =
 $( \frac{23}{24},
\frac{7}{48},
\frac{31}{48} )
\oplus
( \frac{23}{24},
\frac{19}{48},
\frac{43}{48} )
$,

\vskip 0.7ex
\hangindent=5.5em \hangafter=1
{\white .}\hskip 1em $\rho_\text{isum}(\mathfrak{s})$ =
$\mathrm{i}$($0$,
$\sqrt{\frac{1}{2}}$,
$\sqrt{\frac{1}{2}}$;\ \ 
$-\frac{1}{2}$,
$\frac{1}{2}$;\ \ 
$-\frac{1}{2}$)
 $\oplus$
$\mathrm{i}$($0$,
$\sqrt{\frac{1}{2}}$,
$\sqrt{\frac{1}{2}}$;\ \ 
$\frac{1}{2}$,
$-\frac{1}{2}$;\ \ 
$\frac{1}{2}$)

Fail:
all rows of $U \rho(\mathfrak s) U^\dagger$
 contain zero for any block-diagonal $U$. Prop. B.5 (4) eqn. (B.27)

 \ \color{black}

\noindent 351: (dims,levels) = $(3 , 
3;60,
60
)$,
irreps = $3_{5}^{3}
\hskip -1.5pt \otimes \hskip -1.5pt
1_{4}^{1,0}
\hskip -1.5pt \otimes \hskip -1.5pt
1_{3}^{1,0}\oplus
3_{5}^{1}
\hskip -1.5pt \otimes \hskip -1.5pt
1_{4}^{1,0}
\hskip -1.5pt \otimes \hskip -1.5pt
1_{3}^{1,0}$,
pord$(\rho_\text{isum}(\mathfrak{t})) = 5$,

\vskip 0.7ex
\hangindent=5.5em \hangafter=1
{\white .}\hskip 1em $\rho_\text{isum}(\mathfrak{t})$ =
 $( \frac{7}{12},
\frac{11}{60},
\frac{59}{60} )
\oplus
( \frac{7}{12},
\frac{23}{60},
\frac{47}{60} )
$,

\vskip 0.7ex
\hangindent=5.5em \hangafter=1
{\white .}\hskip 1em $\rho_\text{isum}(\mathfrak{s})$ =
$\mathrm{i}$($-\sqrt{\frac{1}{5}}$,
$\sqrt{\frac{2}{5}}$,
$\sqrt{\frac{2}{5}}$;\ \ 
$\frac{-5+\sqrt{5}}{10}$,
$\frac{5+\sqrt{5}}{10}$;\ \ 
$\frac{-5+\sqrt{5}}{10}$)
 $\oplus$
$\mathrm{i}$($\sqrt{\frac{1}{5}}$,
$\sqrt{\frac{2}{5}}$,
$\sqrt{\frac{2}{5}}$;\ \ 
$-\frac{5+\sqrt{5}}{10}$,
$\frac{5-\sqrt{5}}{10}$;\ \ 
$-\frac{5+\sqrt{5}}{10}$)

Fail:
Integral: $D_{\rho}(\sigma)_{\theta} \propto $ id,
 for all $\sigma$ and all $\theta$-eigenspaces that can contain unit. Prop. B.5 (6)

 \ \color{black}

\noindent 352: (dims,levels) = $(4 , 
1 , 
1;7,
1,
1
)$,
irreps = $4_{7}^{1}\oplus
1_{1}^{1}\oplus
1_{1}^{1}$,
pord$(\rho_\text{isum}(\mathfrak{t})) = 7$,

\vskip 0.7ex
\hangindent=5.5em \hangafter=1
{\white .}\hskip 1em $\rho_\text{isum}(\mathfrak{t})$ =
 $( 0,
\frac{1}{7},
\frac{2}{7},
\frac{4}{7} )
\oplus
( 0 )
\oplus
( 0 )
$,

\vskip 0.7ex
\hangindent=5.5em \hangafter=1
{\white .}\hskip 1em $\rho_\text{isum}(\mathfrak{s})$ =
$\mathrm{i}$($-\sqrt{\frac{1}{7}}$,
$\sqrt{\frac{2}{7}}$,
$\sqrt{\frac{2}{7}}$,
$\sqrt{\frac{2}{7}}$;\ \ 
$-\frac{1}{\sqrt{7}}c^{2}_{7}
$,
$-\frac{1}{\sqrt{7}}c^{1}_{7}
$,
$\frac{1}{\sqrt{7}\mathrm{i}}s^{5}_{28}
$;\ \ 
$\frac{1}{\sqrt{7}\mathrm{i}}s^{5}_{28}
$,
$-\frac{1}{\sqrt{7}}c^{2}_{7}
$;\ \ 
$-\frac{1}{\sqrt{7}}c^{1}_{7}
$)
 $\oplus$
($1$)
 $\oplus$
($1$)

Fail:
for $\rho = \rho_1+l\chi, ...,
 (\rho_1(\mathfrak s)/\chi(\mathfrak s))^2\neq$ id. Prop. B.3 (2)

 \ \color{black}

\noindent 353: (dims,levels) = $(4 , 
1 , 
1;7,
1,
1
)$,
irreps = $4_{7}^{3}\oplus
1_{1}^{1}\oplus
1_{1}^{1}$,
pord$(\rho_\text{isum}(\mathfrak{t})) = 7$,

\vskip 0.7ex
\hangindent=5.5em \hangafter=1
{\white .}\hskip 1em $\rho_\text{isum}(\mathfrak{t})$ =
 $( 0,
\frac{3}{7},
\frac{5}{7},
\frac{6}{7} )
\oplus
( 0 )
\oplus
( 0 )
$,

\vskip 0.7ex
\hangindent=5.5em \hangafter=1
{\white .}\hskip 1em $\rho_\text{isum}(\mathfrak{s})$ =
$\mathrm{i}$($\sqrt{\frac{1}{7}}$,
$\sqrt{\frac{2}{7}}$,
$\sqrt{\frac{2}{7}}$,
$\sqrt{\frac{2}{7}}$;\ \ 
$\frac{1}{\sqrt{7}}c^{1}_{7}
$,
$\frac{1}{\sqrt{7}}c^{2}_{7}
$,
$-\frac{1}{\sqrt{7}\mathrm{i}}s^{5}_{28}
$;\ \ 
$-\frac{1}{\sqrt{7}\mathrm{i}}s^{5}_{28}
$,
$\frac{1}{\sqrt{7}}c^{1}_{7}
$;\ \ 
$\frac{1}{\sqrt{7}}c^{2}_{7}
$)
 $\oplus$
($1$)
 $\oplus$
($1$)

Fail:
for $\rho = \rho_1+l\chi, ...,
 (\rho_1(\mathfrak s)/\chi(\mathfrak s))^2\neq$ id. Prop. B.3 (2)

 \ \color{black}

\noindent 354: (dims,levels) = $(4 , 
1 , 
1;9,
1,
1
)$,
irreps = $4_{9,1}^{1,0}\oplus
1_{1}^{1}\oplus
1_{1}^{1}$,
pord$(\rho_\text{isum}(\mathfrak{t})) = 9$,

\vskip 0.7ex
\hangindent=5.5em \hangafter=1
{\white .}\hskip 1em $\rho_\text{isum}(\mathfrak{t})$ =
 $( 0,
\frac{1}{9},
\frac{4}{9},
\frac{7}{9} )
\oplus
( 0 )
\oplus
( 0 )
$,

\vskip 0.7ex
\hangindent=5.5em \hangafter=1
{\white .}\hskip 1em $\rho_\text{isum}(\mathfrak{s})$ =
$\mathrm{i}$($0$,
$\sqrt{\frac{1}{3}}$,
$\sqrt{\frac{1}{3}}$,
$\sqrt{\frac{1}{3}}$;\ \ 
$-\frac{1}{3}c^{1}_{36}
$,
$\frac{1}{3}c^{1}_{36}
-\frac{1}{3}c^{5}_{36}
$,
$\frac{1}{3}c^{5}_{36}
$;\ \ 
$\frac{1}{3}c^{5}_{36}
$,
$-\frac{1}{3}c^{1}_{36}
$;\ \ 
$\frac{1}{3}c^{1}_{36}
-\frac{1}{3}c^{5}_{36}
$)
 $\oplus$
($1$)
 $\oplus$
($1$)

Fail:
for $\rho = \rho_1+l\chi, ...,
 (\rho_1(\mathfrak s)/\chi(\mathfrak s))^2\neq$ id. Prop. B.3 (2)

 \ \color{black}

 \color{blue}

\noindent 355: (dims,levels) = $(4 , 
1 , 
1;9,
1,
1
)$,
irreps = $4_{9,2}^{1,0}\oplus
1_{1}^{1}\oplus
1_{1}^{1}$,
pord$(\rho_\text{isum}(\mathfrak{t})) = 9$,

\vskip 0.7ex
\hangindent=5.5em \hangafter=1
{\white .}\hskip 1em $\rho_\text{isum}(\mathfrak{t})$ =
 $( 0,
\frac{1}{9},
\frac{4}{9},
\frac{7}{9} )
\oplus
( 0 )
\oplus
( 0 )
$,

\vskip 0.7ex
\hangindent=5.5em \hangafter=1
{\white .}\hskip 1em $\rho_\text{isum}(\mathfrak{s})$ =
($0$,
$-\sqrt{\frac{1}{3}}$,
$-\sqrt{\frac{1}{3}}$,
$-\sqrt{\frac{1}{3}}$;
$\frac{1}{3}c^{2}_{9}
$,
$\frac{1}{3} c_9^4 $,
$\frac{1}{3}c^{1}_{9}
$;
$\frac{1}{3}c^{1}_{9}
$,
$\frac{1}{3}c^{2}_{9}
$;
$\frac{1}{3} c_9^4 $)
 $\oplus$
($1$)
 $\oplus$
($1$)

Pass. 

 \ \color{black}

\noindent 356: (dims,levels) = $(4 , 
1 , 
1;9,
1,
1
)$,
irreps = $4_{9,1}^{2,0}\oplus
1_{1}^{1}\oplus
1_{1}^{1}$,
pord$(\rho_\text{isum}(\mathfrak{t})) = 9$,

\vskip 0.7ex
\hangindent=5.5em \hangafter=1
{\white .}\hskip 1em $\rho_\text{isum}(\mathfrak{t})$ =
 $( 0,
\frac{2}{9},
\frac{5}{9},
\frac{8}{9} )
\oplus
( 0 )
\oplus
( 0 )
$,

\vskip 0.7ex
\hangindent=5.5em \hangafter=1
{\white .}\hskip 1em $\rho_\text{isum}(\mathfrak{s})$ =
$\mathrm{i}$($0$,
$\sqrt{\frac{1}{3}}$,
$\sqrt{\frac{1}{3}}$,
$\sqrt{\frac{1}{3}}$;\ \ 
$-\frac{1}{3}c^{1}_{36}
+\frac{1}{3}c^{5}_{36}
$,
$\frac{1}{3}c^{1}_{36}
$,
$-\frac{1}{3}c^{5}_{36}
$;\ \ 
$-\frac{1}{3}c^{5}_{36}
$,
$-\frac{1}{3}c^{1}_{36}
+\frac{1}{3}c^{5}_{36}
$;\ \ 
$\frac{1}{3}c^{1}_{36}
$)
 $\oplus$
($1$)
 $\oplus$
($1$)

Fail:
for $\rho = \rho_1+l\chi, ...,
 (\rho_1(\mathfrak s)/\chi(\mathfrak s))^2\neq$ id. Prop. B.3 (2)

 \ \color{black}

 \color{blue}

\noindent 357: (dims,levels) = $(4 , 
1 , 
1;9,
1,
1
)$,
irreps = $4_{9,2}^{5,0}\oplus
1_{1}^{1}\oplus
1_{1}^{1}$,
pord$(\rho_\text{isum}(\mathfrak{t})) = 9$,

\vskip 0.7ex
\hangindent=5.5em \hangafter=1
{\white .}\hskip 1em $\rho_\text{isum}(\mathfrak{t})$ =
 $( 0,
\frac{2}{9},
\frac{5}{9},
\frac{8}{9} )
\oplus
( 0 )
\oplus
( 0 )
$,

\vskip 0.7ex
\hangindent=5.5em \hangafter=1
{\white .}\hskip 1em $\rho_\text{isum}(\mathfrak{s})$ =
($0$,
$-\sqrt{\frac{1}{3}}$,
$-\sqrt{\frac{1}{3}}$,
$-\sqrt{\frac{1}{3}}$;
$\frac{1}{3} c_9^4 $,
$\frac{1}{3}c^{2}_{9}
$,
$\frac{1}{3}c^{1}_{9}
$;
$\frac{1}{3}c^{1}_{9}
$,
$\frac{1}{3} c_9^4 $;
$\frac{1}{3}c^{2}_{9}
$)
 $\oplus$
($1$)
 $\oplus$
($1$)

Pass. 

 \ \color{black}

\noindent 358: (dims,levels) = $(4 , 
1 , 
1;14,
2,
2
)$,
irreps = $4_{7}^{1}
\hskip -1.5pt \otimes \hskip -1.5pt
1_{2}^{1,0}\oplus
1_{2}^{1,0}\oplus
1_{2}^{1,0}$,
pord$(\rho_\text{isum}(\mathfrak{t})) = 7$,

\vskip 0.7ex
\hangindent=5.5em \hangafter=1
{\white .}\hskip 1em $\rho_\text{isum}(\mathfrak{t})$ =
 $( \frac{1}{2},
\frac{1}{14},
\frac{9}{14},
\frac{11}{14} )
\oplus
( \frac{1}{2} )
\oplus
( \frac{1}{2} )
$,

\vskip 0.7ex
\hangindent=5.5em \hangafter=1
{\white .}\hskip 1em $\rho_\text{isum}(\mathfrak{s})$ =
$\mathrm{i}$($\sqrt{\frac{1}{7}}$,
$\sqrt{\frac{2}{7}}$,
$\sqrt{\frac{2}{7}}$,
$\sqrt{\frac{2}{7}}$;\ \ 
$\frac{1}{\sqrt{7}}c^{1}_{7}
$,
$-\frac{1}{\sqrt{7}\mathrm{i}}s^{5}_{28}
$,
$\frac{1}{\sqrt{7}}c^{2}_{7}
$;\ \ 
$\frac{1}{\sqrt{7}}c^{2}_{7}
$,
$\frac{1}{\sqrt{7}}c^{1}_{7}
$;\ \ 
$-\frac{1}{\sqrt{7}\mathrm{i}}s^{5}_{28}
$)
 $\oplus$
($-1$)
 $\oplus$
($-1$)

Fail:
for $\rho = \rho_1+l\chi, ...,
 (\rho_1(\mathfrak s)/\chi(\mathfrak s))^2\neq$ id. Prop. B.3 (2)

 \ \color{black}

\noindent 359: (dims,levels) = $(4 , 
1 , 
1;14,
2,
2
)$,
irreps = $4_{7}^{3}
\hskip -1.5pt \otimes \hskip -1.5pt
1_{2}^{1,0}\oplus
1_{2}^{1,0}\oplus
1_{2}^{1,0}$,
pord$(\rho_\text{isum}(\mathfrak{t})) = 7$,

\vskip 0.7ex
\hangindent=5.5em \hangafter=1
{\white .}\hskip 1em $\rho_\text{isum}(\mathfrak{t})$ =
 $( \frac{1}{2},
\frac{3}{14},
\frac{5}{14},
\frac{13}{14} )
\oplus
( \frac{1}{2} )
\oplus
( \frac{1}{2} )
$,

\vskip 0.7ex
\hangindent=5.5em \hangafter=1
{\white .}\hskip 1em $\rho_\text{isum}(\mathfrak{s})$ =
$\mathrm{i}$($-\sqrt{\frac{1}{7}}$,
$\sqrt{\frac{2}{7}}$,
$\sqrt{\frac{2}{7}}$,
$\sqrt{\frac{2}{7}}$;\ \ 
$\frac{1}{\sqrt{7}\mathrm{i}}s^{5}_{28}
$,
$-\frac{1}{\sqrt{7}}c^{1}_{7}
$,
$-\frac{1}{\sqrt{7}}c^{2}_{7}
$;\ \ 
$-\frac{1}{\sqrt{7}}c^{2}_{7}
$,
$\frac{1}{\sqrt{7}\mathrm{i}}s^{5}_{28}
$;\ \ 
$-\frac{1}{\sqrt{7}}c^{1}_{7}
$)
 $\oplus$
($-1$)
 $\oplus$
($-1$)

Fail:
for $\rho = \rho_1+l\chi, ...,
 (\rho_1(\mathfrak s)/\chi(\mathfrak s))^2\neq$ id. Prop. B.3 (2)

 \ \color{black}

\noindent 360: (dims,levels) = $(4 , 
1 , 
1;18,
2,
2
)$,
irreps = $4_{9,1}^{2,0}
\hskip -1.5pt \otimes \hskip -1.5pt
1_{2}^{1,0}\oplus
1_{2}^{1,0}\oplus
1_{2}^{1,0}$,
pord$(\rho_\text{isum}(\mathfrak{t})) = 9$,

\vskip 0.7ex
\hangindent=5.5em \hangafter=1
{\white .}\hskip 1em $\rho_\text{isum}(\mathfrak{t})$ =
 $( \frac{1}{2},
\frac{1}{18},
\frac{7}{18},
\frac{13}{18} )
\oplus
( \frac{1}{2} )
\oplus
( \frac{1}{2} )
$,

\vskip 0.7ex
\hangindent=5.5em \hangafter=1
{\white .}\hskip 1em $\rho_\text{isum}(\mathfrak{s})$ =
$\mathrm{i}$($0$,
$\sqrt{\frac{1}{3}}$,
$\sqrt{\frac{1}{3}}$,
$\sqrt{\frac{1}{3}}$;\ \ 
$\frac{1}{3}c^{5}_{36}
$,
$\frac{1}{3}c^{1}_{36}
-\frac{1}{3}c^{5}_{36}
$,
$-\frac{1}{3}c^{1}_{36}
$;\ \ 
$-\frac{1}{3}c^{1}_{36}
$,
$\frac{1}{3}c^{5}_{36}
$;\ \ 
$\frac{1}{3}c^{1}_{36}
-\frac{1}{3}c^{5}_{36}
$)
 $\oplus$
($-1$)
 $\oplus$
($-1$)

Fail:
for $\rho = \rho_1+l\chi, ...,
 (\rho_1(\mathfrak s)/\chi(\mathfrak s))^2\neq$ id. Prop. B.3 (2)

 \ \color{black}

 \color{blue}

\noindent 361: (dims,levels) = $(4 , 
1 , 
1;18,
2,
2
)$,
irreps = $4_{9,2}^{5,0}
\hskip -1.5pt \otimes \hskip -1.5pt
1_{2}^{1,0}\oplus
1_{2}^{1,0}\oplus
1_{2}^{1,0}$,
pord$(\rho_\text{isum}(\mathfrak{t})) = 9$,

\vskip 0.7ex
\hangindent=5.5em \hangafter=1
{\white .}\hskip 1em $\rho_\text{isum}(\mathfrak{t})$ =
 $( \frac{1}{2},
\frac{1}{18},
\frac{7}{18},
\frac{13}{18} )
\oplus
( \frac{1}{2} )
\oplus
( \frac{1}{2} )
$,

\vskip 0.7ex
\hangindent=5.5em \hangafter=1
{\white .}\hskip 1em $\rho_\text{isum}(\mathfrak{s})$ =
($0$,
$-\sqrt{\frac{1}{3}}$,
$-\sqrt{\frac{1}{3}}$,
$-\sqrt{\frac{1}{3}}$;
$-\frac{1}{3}c^{1}_{9}
$,
$-\frac{1}{3} c_9^4 $,
$-\frac{1}{3}c^{2}_{9}
$;
$-\frac{1}{3}c^{2}_{9}
$,
$-\frac{1}{3}c^{1}_{9}
$;
$-\frac{1}{3} c_9^4 $)
 $\oplus$
($-1$)
 $\oplus$
($-1$)

Pass. 

 \ \color{black}

\noindent 362: (dims,levels) = $(4 , 
1 , 
1;18,
2,
2
)$,
irreps = $4_{9,1}^{1,0}
\hskip -1.5pt \otimes \hskip -1.5pt
1_{2}^{1,0}\oplus
1_{2}^{1,0}\oplus
1_{2}^{1,0}$,
pord$(\rho_\text{isum}(\mathfrak{t})) = 9$,

\vskip 0.7ex
\hangindent=5.5em \hangafter=1
{\white .}\hskip 1em $\rho_\text{isum}(\mathfrak{t})$ =
 $( \frac{1}{2},
\frac{5}{18},
\frac{11}{18},
\frac{17}{18} )
\oplus
( \frac{1}{2} )
\oplus
( \frac{1}{2} )
$,

\vskip 0.7ex
\hangindent=5.5em \hangafter=1
{\white .}\hskip 1em $\rho_\text{isum}(\mathfrak{s})$ =
$\mathrm{i}$($0$,
$\sqrt{\frac{1}{3}}$,
$\sqrt{\frac{1}{3}}$,
$\sqrt{\frac{1}{3}}$;\ \ 
$-\frac{1}{3}c^{1}_{36}
+\frac{1}{3}c^{5}_{36}
$,
$-\frac{1}{3}c^{5}_{36}
$,
$\frac{1}{3}c^{1}_{36}
$;\ \ 
$\frac{1}{3}c^{1}_{36}
$,
$-\frac{1}{3}c^{1}_{36}
+\frac{1}{3}c^{5}_{36}
$;\ \ 
$-\frac{1}{3}c^{5}_{36}
$)
 $\oplus$
($-1$)
 $\oplus$
($-1$)

Fail:
for $\rho = \rho_1+l\chi, ...,
 (\rho_1(\mathfrak s)/\chi(\mathfrak s))^2\neq$ id. Prop. B.3 (2)

 \ \color{black}

 \color{blue}

\noindent 363: (dims,levels) = $(4 , 
1 , 
1;18,
2,
2
)$,
irreps = $4_{9,2}^{1,0}
\hskip -1.5pt \otimes \hskip -1.5pt
1_{2}^{1,0}\oplus
1_{2}^{1,0}\oplus
1_{2}^{1,0}$,
pord$(\rho_\text{isum}(\mathfrak{t})) = 9$,

\vskip 0.7ex
\hangindent=5.5em \hangafter=1
{\white .}\hskip 1em $\rho_\text{isum}(\mathfrak{t})$ =
 $( \frac{1}{2},
\frac{5}{18},
\frac{11}{18},
\frac{17}{18} )
\oplus
( \frac{1}{2} )
\oplus
( \frac{1}{2} )
$,

\vskip 0.7ex
\hangindent=5.5em \hangafter=1
{\white .}\hskip 1em $\rho_\text{isum}(\mathfrak{s})$ =
($0$,
$-\sqrt{\frac{1}{3}}$,
$-\sqrt{\frac{1}{3}}$,
$-\sqrt{\frac{1}{3}}$;
$-\frac{1}{3} c_9^4 $,
$-\frac{1}{3}c^{1}_{9}
$,
$-\frac{1}{3}c^{2}_{9}
$;
$-\frac{1}{3}c^{2}_{9}
$,
$-\frac{1}{3} c_9^4 $;
$-\frac{1}{3}c^{1}_{9}
$)
 $\oplus$
($-1$)
 $\oplus$
($-1$)

Pass. 

 \ \color{black}

\noindent 364: (dims,levels) = $(4 , 
1 , 
1;21,
3,
3
)$,
irreps = $4_{7}^{3}
\hskip -1.5pt \otimes \hskip -1.5pt
1_{3}^{1,0}\oplus
1_{3}^{1,0}\oplus
1_{3}^{1,0}$,
pord$(\rho_\text{isum}(\mathfrak{t})) = 7$,

\vskip 0.7ex
\hangindent=5.5em \hangafter=1
{\white .}\hskip 1em $\rho_\text{isum}(\mathfrak{t})$ =
 $( \frac{1}{3},
\frac{1}{21},
\frac{4}{21},
\frac{16}{21} )
\oplus
( \frac{1}{3} )
\oplus
( \frac{1}{3} )
$,

\vskip 0.7ex
\hangindent=5.5em \hangafter=1
{\white .}\hskip 1em $\rho_\text{isum}(\mathfrak{s})$ =
$\mathrm{i}$($\sqrt{\frac{1}{7}}$,
$\sqrt{\frac{2}{7}}$,
$\sqrt{\frac{2}{7}}$,
$\sqrt{\frac{2}{7}}$;\ \ 
$-\frac{1}{\sqrt{7}\mathrm{i}}s^{5}_{28}
$,
$\frac{1}{\sqrt{7}}c^{1}_{7}
$,
$\frac{1}{\sqrt{7}}c^{2}_{7}
$;\ \ 
$\frac{1}{\sqrt{7}}c^{2}_{7}
$,
$-\frac{1}{\sqrt{7}\mathrm{i}}s^{5}_{28}
$;\ \ 
$\frac{1}{\sqrt{7}}c^{1}_{7}
$)
 $\oplus$
($1$)
 $\oplus$
($1$)

Fail:
for $\rho = \rho_1+l\chi, ...,
 (\rho_1(\mathfrak s)/\chi(\mathfrak s))^2\neq$ id. Prop. B.3 (2)

 \ \color{black}

\noindent 365: (dims,levels) = $(4 , 
1 , 
1;21,
3,
3
)$,
irreps = $4_{7}^{1}
\hskip -1.5pt \otimes \hskip -1.5pt
1_{3}^{1,0}\oplus
1_{3}^{1,0}\oplus
1_{3}^{1,0}$,
pord$(\rho_\text{isum}(\mathfrak{t})) = 7$,

\vskip 0.7ex
\hangindent=5.5em \hangafter=1
{\white .}\hskip 1em $\rho_\text{isum}(\mathfrak{t})$ =
 $( \frac{1}{3},
\frac{10}{21},
\frac{13}{21},
\frac{19}{21} )
\oplus
( \frac{1}{3} )
\oplus
( \frac{1}{3} )
$,

\vskip 0.7ex
\hangindent=5.5em \hangafter=1
{\white .}\hskip 1em $\rho_\text{isum}(\mathfrak{s})$ =
$\mathrm{i}$($-\sqrt{\frac{1}{7}}$,
$\sqrt{\frac{2}{7}}$,
$\sqrt{\frac{2}{7}}$,
$\sqrt{\frac{2}{7}}$;\ \ 
$-\frac{1}{\sqrt{7}}c^{2}_{7}
$,
$-\frac{1}{\sqrt{7}}c^{1}_{7}
$,
$\frac{1}{\sqrt{7}\mathrm{i}}s^{5}_{28}
$;\ \ 
$\frac{1}{\sqrt{7}\mathrm{i}}s^{5}_{28}
$,
$-\frac{1}{\sqrt{7}}c^{2}_{7}
$;\ \ 
$-\frac{1}{\sqrt{7}}c^{1}_{7}
$)
 $\oplus$
($1$)
 $\oplus$
($1$)

Fail:
for $\rho = \rho_1+l\chi, ...,
 (\rho_1(\mathfrak s)/\chi(\mathfrak s))^2\neq$ id. Prop. B.3 (2)

 \ \color{black}

\noindent 366: (dims,levels) = $(4 , 
1 , 
1;28,
4,
4
)$,
irreps = $4_{7}^{3}
\hskip -1.5pt \otimes \hskip -1.5pt
1_{4}^{1,0}\oplus
1_{4}^{1,0}\oplus
1_{4}^{1,0}$,
pord$(\rho_\text{isum}(\mathfrak{t})) = 7$,

\vskip 0.7ex
\hangindent=5.5em \hangafter=1
{\white .}\hskip 1em $\rho_\text{isum}(\mathfrak{t})$ =
 $( \frac{1}{4},
\frac{3}{28},
\frac{19}{28},
\frac{27}{28} )
\oplus
( \frac{1}{4} )
\oplus
( \frac{1}{4} )
$,

\vskip 0.7ex
\hangindent=5.5em \hangafter=1
{\white .}\hskip 1em $\rho_\text{isum}(\mathfrak{s})$ =
($-\sqrt{\frac{1}{7}}$,
$\sqrt{\frac{2}{7}}$,
$\sqrt{\frac{2}{7}}$,
$\sqrt{\frac{2}{7}}$;
$-\frac{1}{\sqrt{7}}c^{2}_{7}
$,
$\frac{1}{\sqrt{7}\mathrm{i}}s^{5}_{28}
$,
$-\frac{1}{\sqrt{7}}c^{1}_{7}
$;
$-\frac{1}{\sqrt{7}}c^{1}_{7}
$,
$-\frac{1}{\sqrt{7}}c^{2}_{7}
$;
$\frac{1}{\sqrt{7}\mathrm{i}}s^{5}_{28}
$)
 $\oplus$
$\mathrm{i}$($1$)
 $\oplus$
$\mathrm{i}$($1$)

Fail:
for $\rho = \rho_1+l\chi, ...,
 (\rho_1(\mathfrak s)/\chi(\mathfrak s))^2\neq$ id. Prop. B.3 (2)

 \ \color{black}

\noindent 367: (dims,levels) = $(4 , 
1 , 
1;28,
4,
4
)$,
irreps = $4_{7}^{1}
\hskip -1.5pt \otimes \hskip -1.5pt
1_{4}^{1,0}\oplus
1_{4}^{1,0}\oplus
1_{4}^{1,0}$,
pord$(\rho_\text{isum}(\mathfrak{t})) = 7$,

\vskip 0.7ex
\hangindent=5.5em \hangafter=1
{\white .}\hskip 1em $\rho_\text{isum}(\mathfrak{t})$ =
 $( \frac{1}{4},
\frac{11}{28},
\frac{15}{28},
\frac{23}{28} )
\oplus
( \frac{1}{4} )
\oplus
( \frac{1}{4} )
$,

\vskip 0.7ex
\hangindent=5.5em \hangafter=1
{\white .}\hskip 1em $\rho_\text{isum}(\mathfrak{s})$ =
($\sqrt{\frac{1}{7}}$,
$\sqrt{\frac{2}{7}}$,
$\sqrt{\frac{2}{7}}$,
$\sqrt{\frac{2}{7}}$;
$\frac{1}{\sqrt{7}}c^{2}_{7}
$,
$\frac{1}{\sqrt{7}}c^{1}_{7}
$,
$-\frac{1}{\sqrt{7}\mathrm{i}}s^{5}_{28}
$;
$-\frac{1}{\sqrt{7}\mathrm{i}}s^{5}_{28}
$,
$\frac{1}{\sqrt{7}}c^{2}_{7}
$;
$\frac{1}{\sqrt{7}}c^{1}_{7}
$)
 $\oplus$
$\mathrm{i}$($1$)
 $\oplus$
$\mathrm{i}$($1$)

Fail:
for $\rho = \rho_1+l\chi, ...,
 (\rho_1(\mathfrak s)/\chi(\mathfrak s))^2\neq$ id. Prop. B.3 (2)

 \ \color{black}

 \color{blue}

\noindent 368: (dims,levels) = $(4 , 
1 , 
1;36,
4,
4
)$,
irreps = $4_{9,2}^{1,0}
\hskip -1.5pt \otimes \hskip -1.5pt
1_{4}^{1,0}\oplus
1_{4}^{1,0}\oplus
1_{4}^{1,0}$,
pord$(\rho_\text{isum}(\mathfrak{t})) = 9$,

\vskip 0.7ex
\hangindent=5.5em \hangafter=1
{\white .}\hskip 1em $\rho_\text{isum}(\mathfrak{t})$ =
 $( \frac{1}{4},
\frac{1}{36},
\frac{13}{36},
\frac{25}{36} )
\oplus
( \frac{1}{4} )
\oplus
( \frac{1}{4} )
$,

\vskip 0.7ex
\hangindent=5.5em \hangafter=1
{\white .}\hskip 1em $\rho_\text{isum}(\mathfrak{s})$ =
$\mathrm{i}$($0$,
$\sqrt{\frac{1}{3}}$,
$\sqrt{\frac{1}{3}}$,
$\sqrt{\frac{1}{3}}$;\ \ 
$\frac{1}{3} c_9^4 $,
$\frac{1}{3}c^{1}_{9}
$,
$\frac{1}{3}c^{2}_{9}
$;\ \ 
$\frac{1}{3}c^{2}_{9}
$,
$\frac{1}{3} c_9^4 $;\ \ 
$\frac{1}{3}c^{1}_{9}
$)
 $\oplus$
$\mathrm{i}$($1$)
 $\oplus$
$\mathrm{i}$($1$)

Pass. 

 \ \color{black}

\noindent 369: (dims,levels) = $(4 , 
1 , 
1;36,
4,
4
)$,
irreps = $4_{9,1}^{1,0}
\hskip -1.5pt \otimes \hskip -1.5pt
1_{4}^{1,0}\oplus
1_{4}^{1,0}\oplus
1_{4}^{1,0}$,
pord$(\rho_\text{isum}(\mathfrak{t})) = 9$,

\vskip 0.7ex
\hangindent=5.5em \hangafter=1
{\white .}\hskip 1em $\rho_\text{isum}(\mathfrak{t})$ =
 $( \frac{1}{4},
\frac{1}{36},
\frac{13}{36},
\frac{25}{36} )
\oplus
( \frac{1}{4} )
\oplus
( \frac{1}{4} )
$,

\vskip 0.7ex
\hangindent=5.5em \hangafter=1
{\white .}\hskip 1em $\rho_\text{isum}(\mathfrak{s})$ =
($0$,
$-\sqrt{\frac{1}{3}}$,
$-\sqrt{\frac{1}{3}}$,
$-\sqrt{\frac{1}{3}}$;
$-\frac{1}{3}c^{1}_{36}
+\frac{1}{3}c^{5}_{36}
$,
$-\frac{1}{3}c^{5}_{36}
$,
$\frac{1}{3}c^{1}_{36}
$;
$\frac{1}{3}c^{1}_{36}
$,
$-\frac{1}{3}c^{1}_{36}
+\frac{1}{3}c^{5}_{36}
$;
$-\frac{1}{3}c^{5}_{36}
$)
 $\oplus$
$\mathrm{i}$($1$)
 $\oplus$
$\mathrm{i}$($1$)

Fail:
for $\rho = \rho_1+l\chi, ...,
 (\rho_1(\mathfrak s)/\chi(\mathfrak s))^2\neq$ id. Prop. B.3 (2)

 \ \color{black}

 \color{blue}

\noindent 370: (dims,levels) = $(4 , 
1 , 
1;36,
4,
4
)$,
irreps = $4_{9,2}^{5,0}
\hskip -1.5pt \otimes \hskip -1.5pt
1_{4}^{1,0}\oplus
1_{4}^{1,0}\oplus
1_{4}^{1,0}$,
pord$(\rho_\text{isum}(\mathfrak{t})) = 9$,

\vskip 0.7ex
\hangindent=5.5em \hangafter=1
{\white .}\hskip 1em $\rho_\text{isum}(\mathfrak{t})$ =
 $( \frac{1}{4},
\frac{5}{36},
\frac{17}{36},
\frac{29}{36} )
\oplus
( \frac{1}{4} )
\oplus
( \frac{1}{4} )
$,

\vskip 0.7ex
\hangindent=5.5em \hangafter=1
{\white .}\hskip 1em $\rho_\text{isum}(\mathfrak{s})$ =
$\mathrm{i}$($0$,
$\sqrt{\frac{1}{3}}$,
$\sqrt{\frac{1}{3}}$,
$\sqrt{\frac{1}{3}}$;\ \ 
$\frac{1}{3}c^{2}_{9}
$,
$\frac{1}{3}c^{1}_{9}
$,
$\frac{1}{3} c_9^4 $;\ \ 
$\frac{1}{3} c_9^4 $,
$\frac{1}{3}c^{2}_{9}
$;\ \ 
$\frac{1}{3}c^{1}_{9}
$)
 $\oplus$
$\mathrm{i}$($1$)
 $\oplus$
$\mathrm{i}$($1$)

Pass. 

 \ \color{black}

\noindent 371: (dims,levels) = $(4 , 
1 , 
1;36,
4,
4
)$,
irreps = $4_{9,1}^{2,0}
\hskip -1.5pt \otimes \hskip -1.5pt
1_{4}^{1,0}\oplus
1_{4}^{1,0}\oplus
1_{4}^{1,0}$,
pord$(\rho_\text{isum}(\mathfrak{t})) = 9$,

\vskip 0.7ex
\hangindent=5.5em \hangafter=1
{\white .}\hskip 1em $\rho_\text{isum}(\mathfrak{t})$ =
 $( \frac{1}{4},
\frac{5}{36},
\frac{17}{36},
\frac{29}{36} )
\oplus
( \frac{1}{4} )
\oplus
( \frac{1}{4} )
$,

\vskip 0.7ex
\hangindent=5.5em \hangafter=1
{\white .}\hskip 1em $\rho_\text{isum}(\mathfrak{s})$ =
($0$,
$-\sqrt{\frac{1}{3}}$,
$-\sqrt{\frac{1}{3}}$,
$-\sqrt{\frac{1}{3}}$;
$-\frac{1}{3}c^{1}_{36}
$,
$\frac{1}{3}c^{5}_{36}
$,
$\frac{1}{3}c^{1}_{36}
-\frac{1}{3}c^{5}_{36}
$;
$\frac{1}{3}c^{1}_{36}
-\frac{1}{3}c^{5}_{36}
$,
$-\frac{1}{3}c^{1}_{36}
$;
$\frac{1}{3}c^{5}_{36}
$)
 $\oplus$
$\mathrm{i}$($1$)
 $\oplus$
$\mathrm{i}$($1$)

Fail:
for $\rho = \rho_1+l\chi, ...,
 (\rho_1(\mathfrak s)/\chi(\mathfrak s))^2\neq$ id. Prop. B.3 (2)

 \ \color{black}

\noindent 372: (dims,levels) = $(4 , 
1 , 
1;42,
6,
6
)$,
irreps = $4_{7}^{1}
\hskip -1.5pt \otimes \hskip -1.5pt
1_{3}^{1,0}
\hskip -1.5pt \otimes \hskip -1.5pt
1_{2}^{1,0}\oplus
1_{3}^{1,0}
\hskip -1.5pt \otimes \hskip -1.5pt
1_{2}^{1,0}\oplus
1_{3}^{1,0}
\hskip -1.5pt \otimes \hskip -1.5pt
1_{2}^{1,0}$,
pord$(\rho_\text{isum}(\mathfrak{t})) = 7$,

\vskip 0.7ex
\hangindent=5.5em \hangafter=1
{\white .}\hskip 1em $\rho_\text{isum}(\mathfrak{t})$ =
 $( \frac{5}{6},
\frac{5}{42},
\frac{17}{42},
\frac{41}{42} )
\oplus
( \frac{5}{6} )
\oplus
( \frac{5}{6} )
$,

\vskip 0.7ex
\hangindent=5.5em \hangafter=1
{\white .}\hskip 1em $\rho_\text{isum}(\mathfrak{s})$ =
$\mathrm{i}$($\sqrt{\frac{1}{7}}$,
$\sqrt{\frac{2}{7}}$,
$\sqrt{\frac{2}{7}}$,
$\sqrt{\frac{2}{7}}$;\ \ 
$-\frac{1}{\sqrt{7}\mathrm{i}}s^{5}_{28}
$,
$\frac{1}{\sqrt{7}}c^{2}_{7}
$,
$\frac{1}{\sqrt{7}}c^{1}_{7}
$;\ \ 
$\frac{1}{\sqrt{7}}c^{1}_{7}
$,
$-\frac{1}{\sqrt{7}\mathrm{i}}s^{5}_{28}
$;\ \ 
$\frac{1}{\sqrt{7}}c^{2}_{7}
$)
 $\oplus$
($-1$)
 $\oplus$
($-1$)

Fail:
for $\rho = \rho_1+l\chi, ...,
 (\rho_1(\mathfrak s)/\chi(\mathfrak s))^2\neq$ id. Prop. B.3 (2)

 \ \color{black}

\noindent 373: (dims,levels) = $(4 , 
1 , 
1;42,
6,
6
)$,
irreps = $4_{7}^{3}
\hskip -1.5pt \otimes \hskip -1.5pt
1_{3}^{1,0}
\hskip -1.5pt \otimes \hskip -1.5pt
1_{2}^{1,0}\oplus
1_{3}^{1,0}
\hskip -1.5pt \otimes \hskip -1.5pt
1_{2}^{1,0}\oplus
1_{3}^{1,0}
\hskip -1.5pt \otimes \hskip -1.5pt
1_{2}^{1,0}$,
pord$(\rho_\text{isum}(\mathfrak{t})) = 7$,

\vskip 0.7ex
\hangindent=5.5em \hangafter=1
{\white .}\hskip 1em $\rho_\text{isum}(\mathfrak{t})$ =
 $( \frac{5}{6},
\frac{11}{42},
\frac{23}{42},
\frac{29}{42} )
\oplus
( \frac{5}{6} )
\oplus
( \frac{5}{6} )
$,

\vskip 0.7ex
\hangindent=5.5em \hangafter=1
{\white .}\hskip 1em $\rho_\text{isum}(\mathfrak{s})$ =
$\mathrm{i}$($-\sqrt{\frac{1}{7}}$,
$\sqrt{\frac{2}{7}}$,
$\sqrt{\frac{2}{7}}$,
$\sqrt{\frac{2}{7}}$;\ \ 
$-\frac{1}{\sqrt{7}}c^{1}_{7}
$,
$-\frac{1}{\sqrt{7}}c^{2}_{7}
$,
$\frac{1}{\sqrt{7}\mathrm{i}}s^{5}_{28}
$;\ \ 
$\frac{1}{\sqrt{7}\mathrm{i}}s^{5}_{28}
$,
$-\frac{1}{\sqrt{7}}c^{1}_{7}
$;\ \ 
$-\frac{1}{\sqrt{7}}c^{2}_{7}
$)
 $\oplus$
($-1$)
 $\oplus$
($-1$)

Fail:
for $\rho = \rho_1+l\chi, ...,
 (\rho_1(\mathfrak s)/\chi(\mathfrak s))^2\neq$ id. Prop. B.3 (2)

 \ \color{black}

\noindent 374: (dims,levels) = $(4 , 
1 , 
1;84,
12,
12
)$,
irreps = $4_{7}^{3}
\hskip -1.5pt \otimes \hskip -1.5pt
1_{4}^{1,0}
\hskip -1.5pt \otimes \hskip -1.5pt
1_{3}^{1,0}\oplus
1_{4}^{1,0}
\hskip -1.5pt \otimes \hskip -1.5pt
1_{3}^{1,0}\oplus
1_{4}^{1,0}
\hskip -1.5pt \otimes \hskip -1.5pt
1_{3}^{1,0}$,
pord$(\rho_\text{isum}(\mathfrak{t})) = 7$,

\vskip 0.7ex
\hangindent=5.5em \hangafter=1
{\white .}\hskip 1em $\rho_\text{isum}(\mathfrak{t})$ =
 $( \frac{7}{12},
\frac{1}{84},
\frac{25}{84},
\frac{37}{84} )
\oplus
( \frac{7}{12} )
\oplus
( \frac{7}{12} )
$,

\vskip 0.7ex
\hangindent=5.5em \hangafter=1
{\white .}\hskip 1em $\rho_\text{isum}(\mathfrak{s})$ =
($-\sqrt{\frac{1}{7}}$,
$\sqrt{\frac{2}{7}}$,
$\sqrt{\frac{2}{7}}$,
$\sqrt{\frac{2}{7}}$;
$-\frac{1}{\sqrt{7}}c^{1}_{7}
$,
$-\frac{1}{\sqrt{7}}c^{2}_{7}
$,
$\frac{1}{\sqrt{7}\mathrm{i}}s^{5}_{28}
$;
$\frac{1}{\sqrt{7}\mathrm{i}}s^{5}_{28}
$,
$-\frac{1}{\sqrt{7}}c^{1}_{7}
$;
$-\frac{1}{\sqrt{7}}c^{2}_{7}
$)
 $\oplus$
$\mathrm{i}$($1$)
 $\oplus$
$\mathrm{i}$($1$)

Fail:
for $\rho = \rho_1+l\chi, ...,
 (\rho_1(\mathfrak s)/\chi(\mathfrak s))^2\neq$ id. Prop. B.3 (2)

 \ \color{black}

\noindent 375: (dims,levels) = $(4 , 
1 , 
1;84,
12,
12
)$,
irreps = $4_{7}^{1}
\hskip -1.5pt \otimes \hskip -1.5pt
1_{4}^{1,0}
\hskip -1.5pt \otimes \hskip -1.5pt
1_{3}^{1,0}\oplus
1_{4}^{1,0}
\hskip -1.5pt \otimes \hskip -1.5pt
1_{3}^{1,0}\oplus
1_{4}^{1,0}
\hskip -1.5pt \otimes \hskip -1.5pt
1_{3}^{1,0}$,
pord$(\rho_\text{isum}(\mathfrak{t})) = 7$,

\vskip 0.7ex
\hangindent=5.5em \hangafter=1
{\white .}\hskip 1em $\rho_\text{isum}(\mathfrak{t})$ =
 $( \frac{7}{12},
\frac{13}{84},
\frac{61}{84},
\frac{73}{84} )
\oplus
( \frac{7}{12} )
\oplus
( \frac{7}{12} )
$,

\vskip 0.7ex
\hangindent=5.5em \hangafter=1
{\white .}\hskip 1em $\rho_\text{isum}(\mathfrak{s})$ =
($\sqrt{\frac{1}{7}}$,
$\sqrt{\frac{2}{7}}$,
$\sqrt{\frac{2}{7}}$,
$\sqrt{\frac{2}{7}}$;
$\frac{1}{\sqrt{7}}c^{1}_{7}
$,
$-\frac{1}{\sqrt{7}\mathrm{i}}s^{5}_{28}
$,
$\frac{1}{\sqrt{7}}c^{2}_{7}
$;
$\frac{1}{\sqrt{7}}c^{2}_{7}
$,
$\frac{1}{\sqrt{7}}c^{1}_{7}
$;
$-\frac{1}{\sqrt{7}\mathrm{i}}s^{5}_{28}
$)
 $\oplus$
$\mathrm{i}$($1$)
 $\oplus$
$\mathrm{i}$($1$)

Fail:
for $\rho = \rho_1+l\chi, ...,
 (\rho_1(\mathfrak s)/\chi(\mathfrak s))^2\neq$ id. Prop. B.3 (2)

 \ \color{black}

 \color{blue}

\noindent 376: (dims,levels) = $(4 , 
2;5,
5
)$,
irreps = $4_{5,2}^{1}\oplus
2_{5}^{1}$,
pord$(\rho_\text{isum}(\mathfrak{t})) = 5$,

\vskip 0.7ex
\hangindent=5.5em \hangafter=1
{\white .}\hskip 1em $\rho_\text{isum}(\mathfrak{t})$ =
 $( \frac{1}{5},
\frac{2}{5},
\frac{3}{5},
\frac{4}{5} )
\oplus
( \frac{1}{5},
\frac{4}{5} )
$,

\vskip 0.7ex
\hangindent=5.5em \hangafter=1
{\white .}\hskip 1em $\rho_\text{isum}(\mathfrak{s})$ =
($\sqrt{\frac{1}{5}}$,
$\frac{-5+\sqrt{5}}{10}$,
$-\frac{5+\sqrt{5}}{10}$,
$\sqrt{\frac{1}{5}}$;
$-\sqrt{\frac{1}{5}}$,
$\sqrt{\frac{1}{5}}$,
$\frac{5+\sqrt{5}}{10}$;
$-\sqrt{\frac{1}{5}}$,
$\frac{5-\sqrt{5}}{10}$;
$\sqrt{\frac{1}{5}}$)
 $\oplus$
$\mathrm{i}$($-\frac{1}{\sqrt{5}}c^{3}_{20}
$,
$\frac{1}{\sqrt{5}}c^{1}_{20}
$;\ \ 
$\frac{1}{\sqrt{5}}c^{3}_{20}
$)

Pass. 

 \ \color{black}

 \color{blue}

\noindent 377: (dims,levels) = $(4 , 
2;5,
5
)$,
irreps = $4_{5,1}^{1}\oplus
2_{5}^{1}$,
pord$(\rho_\text{isum}(\mathfrak{t})) = 5$,

\vskip 0.7ex
\hangindent=5.5em \hangafter=1
{\white .}\hskip 1em $\rho_\text{isum}(\mathfrak{t})$ =
 $( \frac{1}{5},
\frac{2}{5},
\frac{3}{5},
\frac{4}{5} )
\oplus
( \frac{1}{5},
\frac{4}{5} )
$,

\vskip 0.7ex
\hangindent=5.5em \hangafter=1
{\white .}\hskip 1em $\rho_\text{isum}(\mathfrak{s})$ =
$\mathrm{i}$($\frac{1}{5}c^{1}_{20}
+\frac{1}{5}c^{3}_{20}
$,
$\frac{2}{5}c^{2}_{15}
+\frac{1}{5}c^{3}_{15}
$,
$-\frac{1}{5}+\frac{2}{5}c^{1}_{15}
-\frac{1}{5}c^{3}_{15}
$,
$\frac{1}{5}c^{1}_{20}
-\frac{1}{5}c^{3}_{20}
$;\ \ 
$-\frac{1}{5}c^{1}_{20}
+\frac{1}{5}c^{3}_{20}
$,
$-\frac{1}{5}c^{1}_{20}
-\frac{1}{5}c^{3}_{20}
$,
$\frac{1}{5}-\frac{2}{5}c^{1}_{15}
+\frac{1}{5}c^{3}_{15}
$;\ \ 
$\frac{1}{5}c^{1}_{20}
-\frac{1}{5}c^{3}_{20}
$,
$\frac{2}{5}c^{2}_{15}
+\frac{1}{5}c^{3}_{15}
$;\ \ 
$-\frac{1}{5}c^{1}_{20}
-\frac{1}{5}c^{3}_{20}
$)
 $\oplus$
$\mathrm{i}$($-\frac{1}{\sqrt{5}}c^{3}_{20}
$,
$\frac{1}{\sqrt{5}}c^{1}_{20}
$;\ \ 
$\frac{1}{\sqrt{5}}c^{3}_{20}
$)

Pass. 

 \ \color{black}

 \color{blue}

\noindent 378: (dims,levels) = $(4 , 
2;5,
5
)$,
irreps = $4_{5,2}^{1}\oplus
2_{5}^{2}$,
pord$(\rho_\text{isum}(\mathfrak{t})) = 5$,

\vskip 0.7ex
\hangindent=5.5em \hangafter=1
{\white .}\hskip 1em $\rho_\text{isum}(\mathfrak{t})$ =
 $( \frac{1}{5},
\frac{2}{5},
\frac{3}{5},
\frac{4}{5} )
\oplus
( \frac{2}{5},
\frac{3}{5} )
$,

\vskip 0.7ex
\hangindent=5.5em \hangafter=1
{\white .}\hskip 1em $\rho_\text{isum}(\mathfrak{s})$ =
($\sqrt{\frac{1}{5}}$,
$\frac{-5+\sqrt{5}}{10}$,
$-\frac{5+\sqrt{5}}{10}$,
$\sqrt{\frac{1}{5}}$;
$-\sqrt{\frac{1}{5}}$,
$\sqrt{\frac{1}{5}}$,
$\frac{5+\sqrt{5}}{10}$;
$-\sqrt{\frac{1}{5}}$,
$\frac{5-\sqrt{5}}{10}$;
$\sqrt{\frac{1}{5}}$)
 $\oplus$
$\mathrm{i}$($-\frac{1}{\sqrt{5}}c^{1}_{20}
$,
$\frac{1}{\sqrt{5}}c^{3}_{20}
$;\ \ 
$\frac{1}{\sqrt{5}}c^{1}_{20}
$)

Pass. 

 \ \color{black}

 \color{blue}

\noindent 379: (dims,levels) = $(4 , 
2;5,
5
)$,
irreps = $4_{5,1}^{1}\oplus
2_{5}^{2}$,
pord$(\rho_\text{isum}(\mathfrak{t})) = 5$,

\vskip 0.7ex
\hangindent=5.5em \hangafter=1
{\white .}\hskip 1em $\rho_\text{isum}(\mathfrak{t})$ =
 $( \frac{1}{5},
\frac{2}{5},
\frac{3}{5},
\frac{4}{5} )
\oplus
( \frac{2}{5},
\frac{3}{5} )
$,

\vskip 0.7ex
\hangindent=5.5em \hangafter=1
{\white .}\hskip 1em $\rho_\text{isum}(\mathfrak{s})$ =
$\mathrm{i}$($\frac{1}{5}c^{1}_{20}
+\frac{1}{5}c^{3}_{20}
$,
$\frac{2}{5}c^{2}_{15}
+\frac{1}{5}c^{3}_{15}
$,
$-\frac{1}{5}+\frac{2}{5}c^{1}_{15}
-\frac{1}{5}c^{3}_{15}
$,
$\frac{1}{5}c^{1}_{20}
-\frac{1}{5}c^{3}_{20}
$;\ \ 
$-\frac{1}{5}c^{1}_{20}
+\frac{1}{5}c^{3}_{20}
$,
$-\frac{1}{5}c^{1}_{20}
-\frac{1}{5}c^{3}_{20}
$,
$\frac{1}{5}-\frac{2}{5}c^{1}_{15}
+\frac{1}{5}c^{3}_{15}
$;\ \ 
$\frac{1}{5}c^{1}_{20}
-\frac{1}{5}c^{3}_{20}
$,
$\frac{2}{5}c^{2}_{15}
+\frac{1}{5}c^{3}_{15}
$;\ \ 
$-\frac{1}{5}c^{1}_{20}
-\frac{1}{5}c^{3}_{20}
$)
 $\oplus$
$\mathrm{i}$($-\frac{1}{\sqrt{5}}c^{1}_{20}
$,
$\frac{1}{\sqrt{5}}c^{3}_{20}
$;\ \ 
$\frac{1}{\sqrt{5}}c^{1}_{20}
$)

Pass. 

 \ \color{black}

\noindent 380: (dims,levels) = $(4 , 
2;7,
2
)$,
irreps = $4_{7}^{1}\oplus
2_{2}^{1,0}$,
pord$(\rho_\text{isum}(\mathfrak{t})) = 14$,

\vskip 0.7ex
\hangindent=5.5em \hangafter=1
{\white .}\hskip 1em $\rho_\text{isum}(\mathfrak{t})$ =
 $( 0,
\frac{1}{7},
\frac{2}{7},
\frac{4}{7} )
\oplus
( 0,
\frac{1}{2} )
$,

\vskip 0.7ex
\hangindent=5.5em \hangafter=1
{\white .}\hskip 1em $\rho_\text{isum}(\mathfrak{s})$ =
$\mathrm{i}$($-\sqrt{\frac{1}{7}}$,
$\sqrt{\frac{2}{7}}$,
$\sqrt{\frac{2}{7}}$,
$\sqrt{\frac{2}{7}}$;\ \ 
$-\frac{1}{\sqrt{7}}c^{2}_{7}
$,
$-\frac{1}{\sqrt{7}}c^{1}_{7}
$,
$\frac{1}{\sqrt{7}\mathrm{i}}s^{5}_{28}
$;\ \ 
$\frac{1}{\sqrt{7}\mathrm{i}}s^{5}_{28}
$,
$-\frac{1}{\sqrt{7}}c^{2}_{7}
$;\ \ 
$-\frac{1}{\sqrt{7}}c^{1}_{7}
$)
 $\oplus$
($-\frac{1}{2}$,
$-\sqrt{\frac{3}{4}}$;
$\frac{1}{2}$)

Fail:
Tr$_I(C) = -1 <$  0 for I = [ 1/2 ]. Prop. B.4 (1) eqn. (B.18)

 \ \color{black}

\noindent 381: (dims,levels) = $(4 , 
2;7,
2
)$,
irreps = $4_{7}^{3}\oplus
2_{2}^{1,0}$,
pord$(\rho_\text{isum}(\mathfrak{t})) = 14$,

\vskip 0.7ex
\hangindent=5.5em \hangafter=1
{\white .}\hskip 1em $\rho_\text{isum}(\mathfrak{t})$ =
 $( 0,
\frac{3}{7},
\frac{5}{7},
\frac{6}{7} )
\oplus
( 0,
\frac{1}{2} )
$,

\vskip 0.7ex
\hangindent=5.5em \hangafter=1
{\white .}\hskip 1em $\rho_\text{isum}(\mathfrak{s})$ =
$\mathrm{i}$($\sqrt{\frac{1}{7}}$,
$\sqrt{\frac{2}{7}}$,
$\sqrt{\frac{2}{7}}$,
$\sqrt{\frac{2}{7}}$;\ \ 
$\frac{1}{\sqrt{7}}c^{1}_{7}
$,
$\frac{1}{\sqrt{7}}c^{2}_{7}
$,
$-\frac{1}{\sqrt{7}\mathrm{i}}s^{5}_{28}
$;\ \ 
$-\frac{1}{\sqrt{7}\mathrm{i}}s^{5}_{28}
$,
$\frac{1}{\sqrt{7}}c^{1}_{7}
$;\ \ 
$\frac{1}{\sqrt{7}}c^{2}_{7}
$)
 $\oplus$
($-\frac{1}{2}$,
$-\sqrt{\frac{3}{4}}$;
$\frac{1}{2}$)

Fail:
Tr$_I(C) = -1 <$  0 for I = [ 1/2 ]. Prop. B.4 (1) eqn. (B.18)

 \ \color{black}

 \color{blue}

\noindent 382: (dims,levels) = $(4 , 
2;7,
3
)$,
irreps = $4_{7}^{1}\oplus
2_{3}^{1,0}$,
pord$(\rho_\text{isum}(\mathfrak{t})) = 21$,

\vskip 0.7ex
\hangindent=5.5em \hangafter=1
{\white .}\hskip 1em $\rho_\text{isum}(\mathfrak{t})$ =
 $( 0,
\frac{1}{7},
\frac{2}{7},
\frac{4}{7} )
\oplus
( 0,
\frac{1}{3} )
$,

\vskip 0.7ex
\hangindent=5.5em \hangafter=1
{\white .}\hskip 1em $\rho_\text{isum}(\mathfrak{s})$ =
$\mathrm{i}$($-\sqrt{\frac{1}{7}}$,
$\sqrt{\frac{2}{7}}$,
$\sqrt{\frac{2}{7}}$,
$\sqrt{\frac{2}{7}}$;\ \ 
$-\frac{1}{\sqrt{7}}c^{2}_{7}
$,
$-\frac{1}{\sqrt{7}}c^{1}_{7}
$,
$\frac{1}{\sqrt{7}\mathrm{i}}s^{5}_{28}
$;\ \ 
$\frac{1}{\sqrt{7}\mathrm{i}}s^{5}_{28}
$,
$-\frac{1}{\sqrt{7}}c^{2}_{7}
$;\ \ 
$-\frac{1}{\sqrt{7}}c^{1}_{7}
$)
 $\oplus$
$\mathrm{i}$($-\sqrt{\frac{1}{3}}$,
$\sqrt{\frac{2}{3}}$;\ \ 
$\sqrt{\frac{1}{3}}$)

Pass. 

 \ \color{black}

 \color{blue}

\noindent 383: (dims,levels) = $(4 , 
2;7,
3
)$,
irreps = $4_{7}^{1}\oplus
2_{3}^{1,8}$,
pord$(\rho_\text{isum}(\mathfrak{t})) = 21$,

\vskip 0.7ex
\hangindent=5.5em \hangafter=1
{\white .}\hskip 1em $\rho_\text{isum}(\mathfrak{t})$ =
 $( 0,
\frac{1}{7},
\frac{2}{7},
\frac{4}{7} )
\oplus
( 0,
\frac{2}{3} )
$,

\vskip 0.7ex
\hangindent=5.5em \hangafter=1
{\white .}\hskip 1em $\rho_\text{isum}(\mathfrak{s})$ =
$\mathrm{i}$($-\sqrt{\frac{1}{7}}$,
$\sqrt{\frac{2}{7}}$,
$\sqrt{\frac{2}{7}}$,
$\sqrt{\frac{2}{7}}$;\ \ 
$-\frac{1}{\sqrt{7}}c^{2}_{7}
$,
$-\frac{1}{\sqrt{7}}c^{1}_{7}
$,
$\frac{1}{\sqrt{7}\mathrm{i}}s^{5}_{28}
$;\ \ 
$\frac{1}{\sqrt{7}\mathrm{i}}s^{5}_{28}
$,
$-\frac{1}{\sqrt{7}}c^{2}_{7}
$;\ \ 
$-\frac{1}{\sqrt{7}}c^{1}_{7}
$)
 $\oplus$
$\mathrm{i}$($\sqrt{\frac{1}{3}}$,
$\sqrt{\frac{2}{3}}$;\ \ 
$-\sqrt{\frac{1}{3}}$)

Pass. 

 \ \color{black}

 \color{blue}

\noindent 384: (dims,levels) = $(4 , 
2;7,
3
)$,
irreps = $4_{7}^{3}\oplus
2_{3}^{1,0}$,
pord$(\rho_\text{isum}(\mathfrak{t})) = 21$,

\vskip 0.7ex
\hangindent=5.5em \hangafter=1
{\white .}\hskip 1em $\rho_\text{isum}(\mathfrak{t})$ =
 $( 0,
\frac{3}{7},
\frac{5}{7},
\frac{6}{7} )
\oplus
( 0,
\frac{1}{3} )
$,

\vskip 0.7ex
\hangindent=5.5em \hangafter=1
{\white .}\hskip 1em $\rho_\text{isum}(\mathfrak{s})$ =
$\mathrm{i}$($\sqrt{\frac{1}{7}}$,
$\sqrt{\frac{2}{7}}$,
$\sqrt{\frac{2}{7}}$,
$\sqrt{\frac{2}{7}}$;\ \ 
$\frac{1}{\sqrt{7}}c^{1}_{7}
$,
$\frac{1}{\sqrt{7}}c^{2}_{7}
$,
$-\frac{1}{\sqrt{7}\mathrm{i}}s^{5}_{28}
$;\ \ 
$-\frac{1}{\sqrt{7}\mathrm{i}}s^{5}_{28}
$,
$\frac{1}{\sqrt{7}}c^{1}_{7}
$;\ \ 
$\frac{1}{\sqrt{7}}c^{2}_{7}
$)
 $\oplus$
$\mathrm{i}$($-\sqrt{\frac{1}{3}}$,
$\sqrt{\frac{2}{3}}$;\ \ 
$\sqrt{\frac{1}{3}}$)

Pass. 

 \ \color{black}

 \color{blue}

\noindent 385: (dims,levels) = $(4 , 
2;7,
3
)$,
irreps = $4_{7}^{3}\oplus
2_{3}^{1,8}$,
pord$(\rho_\text{isum}(\mathfrak{t})) = 21$,

\vskip 0.7ex
\hangindent=5.5em \hangafter=1
{\white .}\hskip 1em $\rho_\text{isum}(\mathfrak{t})$ =
 $( 0,
\frac{3}{7},
\frac{5}{7},
\frac{6}{7} )
\oplus
( 0,
\frac{2}{3} )
$,

\vskip 0.7ex
\hangindent=5.5em \hangafter=1
{\white .}\hskip 1em $\rho_\text{isum}(\mathfrak{s})$ =
$\mathrm{i}$($\sqrt{\frac{1}{7}}$,
$\sqrt{\frac{2}{7}}$,
$\sqrt{\frac{2}{7}}$,
$\sqrt{\frac{2}{7}}$;\ \ 
$\frac{1}{\sqrt{7}}c^{1}_{7}
$,
$\frac{1}{\sqrt{7}}c^{2}_{7}
$,
$-\frac{1}{\sqrt{7}\mathrm{i}}s^{5}_{28}
$;\ \ 
$-\frac{1}{\sqrt{7}\mathrm{i}}s^{5}_{28}
$,
$\frac{1}{\sqrt{7}}c^{1}_{7}
$;\ \ 
$\frac{1}{\sqrt{7}}c^{2}_{7}
$)
 $\oplus$
$\mathrm{i}$($\sqrt{\frac{1}{3}}$,
$\sqrt{\frac{2}{3}}$;\ \ 
$-\sqrt{\frac{1}{3}}$)

Pass. 

 \ \color{black}

\noindent 386: (dims,levels) = $(4 , 
2;9,
2
)$,
irreps = $4_{9,1}^{1,0}\oplus
2_{2}^{1,0}$,
pord$(\rho_\text{isum}(\mathfrak{t})) = 18$,

\vskip 0.7ex
\hangindent=5.5em \hangafter=1
{\white .}\hskip 1em $\rho_\text{isum}(\mathfrak{t})$ =
 $( 0,
\frac{1}{9},
\frac{4}{9},
\frac{7}{9} )
\oplus
( 0,
\frac{1}{2} )
$,

\vskip 0.7ex
\hangindent=5.5em \hangafter=1
{\white .}\hskip 1em $\rho_\text{isum}(\mathfrak{s})$ =
$\mathrm{i}$($0$,
$\sqrt{\frac{1}{3}}$,
$\sqrt{\frac{1}{3}}$,
$\sqrt{\frac{1}{3}}$;\ \ 
$-\frac{1}{3}c^{1}_{36}
$,
$\frac{1}{3}c^{1}_{36}
-\frac{1}{3}c^{5}_{36}
$,
$\frac{1}{3}c^{5}_{36}
$;\ \ 
$\frac{1}{3}c^{5}_{36}
$,
$-\frac{1}{3}c^{1}_{36}
$;\ \ 
$\frac{1}{3}c^{1}_{36}
-\frac{1}{3}c^{5}_{36}
$)
 $\oplus$
($-\frac{1}{2}$,
$-\sqrt{\frac{3}{4}}$;
$\frac{1}{2}$)

Fail:
Tr$_I(C) = -1 <$  0 for I = [ 1/2 ]. Prop. B.4 (1) eqn. (B.18)

 \ \color{black}

\noindent 387: (dims,levels) = $(4 , 
2;9,
2
)$,
irreps = $4_{9,2}^{1,0}\oplus
2_{2}^{1,0}$,
pord$(\rho_\text{isum}(\mathfrak{t})) = 18$,

\vskip 0.7ex
\hangindent=5.5em \hangafter=1
{\white .}\hskip 1em $\rho_\text{isum}(\mathfrak{t})$ =
 $( 0,
\frac{1}{9},
\frac{4}{9},
\frac{7}{9} )
\oplus
( 0,
\frac{1}{2} )
$,

\vskip 0.7ex
\hangindent=5.5em \hangafter=1
{\white .}\hskip 1em $\rho_\text{isum}(\mathfrak{s})$ =
($0$,
$-\sqrt{\frac{1}{3}}$,
$-\sqrt{\frac{1}{3}}$,
$-\sqrt{\frac{1}{3}}$;
$\frac{1}{3}c^{2}_{9}
$,
$\frac{1}{3} c_9^4 $,
$\frac{1}{3}c^{1}_{9}
$;
$\frac{1}{3}c^{1}_{9}
$,
$\frac{1}{3}c^{2}_{9}
$;
$\frac{1}{3} c_9^4 $)
 $\oplus$
($-\frac{1}{2}$,
$-\sqrt{\frac{3}{4}}$;
$\frac{1}{2}$)

Fail:
Integral: $D_{\rho}(\sigma)_{\theta} \propto $ id,
 for all $\sigma$ and all $\theta$-eigenspaces that can contain unit. Prop. B.5 (6)

 \ \color{black}

\noindent 388: (dims,levels) = $(4 , 
2;9,
2
)$,
irreps = $4_{9,1}^{2,0}\oplus
2_{2}^{1,0}$,
pord$(\rho_\text{isum}(\mathfrak{t})) = 18$,

\vskip 0.7ex
\hangindent=5.5em \hangafter=1
{\white .}\hskip 1em $\rho_\text{isum}(\mathfrak{t})$ =
 $( 0,
\frac{2}{9},
\frac{5}{9},
\frac{8}{9} )
\oplus
( 0,
\frac{1}{2} )
$,

\vskip 0.7ex
\hangindent=5.5em \hangafter=1
{\white .}\hskip 1em $\rho_\text{isum}(\mathfrak{s})$ =
$\mathrm{i}$($0$,
$\sqrt{\frac{1}{3}}$,
$\sqrt{\frac{1}{3}}$,
$\sqrt{\frac{1}{3}}$;\ \ 
$-\frac{1}{3}c^{1}_{36}
+\frac{1}{3}c^{5}_{36}
$,
$\frac{1}{3}c^{1}_{36}
$,
$-\frac{1}{3}c^{5}_{36}
$;\ \ 
$-\frac{1}{3}c^{5}_{36}
$,
$-\frac{1}{3}c^{1}_{36}
+\frac{1}{3}c^{5}_{36}
$;\ \ 
$\frac{1}{3}c^{1}_{36}
$)
 $\oplus$
($-\frac{1}{2}$,
$-\sqrt{\frac{3}{4}}$;
$\frac{1}{2}$)

Fail:
Tr$_I(C) = -1 <$  0 for I = [ 1/2 ]. Prop. B.4 (1) eqn. (B.18)

 \ \color{black}

\noindent 389: (dims,levels) = $(4 , 
2;9,
2
)$,
irreps = $4_{9,2}^{5,0}\oplus
2_{2}^{1,0}$,
pord$(\rho_\text{isum}(\mathfrak{t})) = 18$,

\vskip 0.7ex
\hangindent=5.5em \hangafter=1
{\white .}\hskip 1em $\rho_\text{isum}(\mathfrak{t})$ =
 $( 0,
\frac{2}{9},
\frac{5}{9},
\frac{8}{9} )
\oplus
( 0,
\frac{1}{2} )
$,

\vskip 0.7ex
\hangindent=5.5em \hangafter=1
{\white .}\hskip 1em $\rho_\text{isum}(\mathfrak{s})$ =
($0$,
$-\sqrt{\frac{1}{3}}$,
$-\sqrt{\frac{1}{3}}$,
$-\sqrt{\frac{1}{3}}$;
$\frac{1}{3} c_9^4 $,
$\frac{1}{3}c^{2}_{9}
$,
$\frac{1}{3}c^{1}_{9}
$;
$\frac{1}{3}c^{1}_{9}
$,
$\frac{1}{3} c_9^4 $;
$\frac{1}{3}c^{2}_{9}
$)
 $\oplus$
($-\frac{1}{2}$,
$-\sqrt{\frac{3}{4}}$;
$\frac{1}{2}$)

Fail:
Integral: $D_{\rho}(\sigma)_{\theta} \propto $ id,
 for all $\sigma$ and all $\theta$-eigenspaces that can contain unit. Prop. B.5 (6)

 \ \color{black}

\noindent 390: (dims,levels) = $(4 , 
2;9,
3
)$,
irreps = $4_{9,1}^{1,0}\oplus
2_{3}^{1,0}$,
pord$(\rho_\text{isum}(\mathfrak{t})) = 9$,

\vskip 0.7ex
\hangindent=5.5em \hangafter=1
{\white .}\hskip 1em $\rho_\text{isum}(\mathfrak{t})$ =
 $( 0,
\frac{1}{9},
\frac{4}{9},
\frac{7}{9} )
\oplus
( 0,
\frac{1}{3} )
$,

\vskip 0.7ex
\hangindent=5.5em \hangafter=1
{\white .}\hskip 1em $\rho_\text{isum}(\mathfrak{s})$ =
$\mathrm{i}$($0$,
$\sqrt{\frac{1}{3}}$,
$\sqrt{\frac{1}{3}}$,
$\sqrt{\frac{1}{3}}$;\ \ 
$-\frac{1}{3}c^{1}_{36}
$,
$\frac{1}{3}c^{1}_{36}
-\frac{1}{3}c^{5}_{36}
$,
$\frac{1}{3}c^{5}_{36}
$;\ \ 
$\frac{1}{3}c^{5}_{36}
$,
$-\frac{1}{3}c^{1}_{36}
$;\ \ 
$\frac{1}{3}c^{1}_{36}
-\frac{1}{3}c^{5}_{36}
$)
 $\oplus$
$\mathrm{i}$($-\sqrt{\frac{1}{3}}$,
$\sqrt{\frac{2}{3}}$;\ \ 
$\sqrt{\frac{1}{3}}$)

Fail:
Integral: $D_{\rho}(\sigma)_{\theta} \propto $ id,
 for all $\sigma$ and all $\theta$-eigenspaces that can contain unit. Prop. B.5 (6)

 \ \color{black}

\noindent 391: (dims,levels) = $(4 , 
2;9,
3
)$,
irreps = $4_{9,2}^{1,0}\oplus
2_{3}^{1,0}$,
pord$(\rho_\text{isum}(\mathfrak{t})) = 9$,

\vskip 0.7ex
\hangindent=5.5em \hangafter=1
{\white .}\hskip 1em $\rho_\text{isum}(\mathfrak{t})$ =
 $( 0,
\frac{1}{9},
\frac{4}{9},
\frac{7}{9} )
\oplus
( 0,
\frac{1}{3} )
$,

\vskip 0.7ex
\hangindent=5.5em \hangafter=1
{\white .}\hskip 1em $\rho_\text{isum}(\mathfrak{s})$ =
($0$,
$-\sqrt{\frac{1}{3}}$,
$-\sqrt{\frac{1}{3}}$,
$-\sqrt{\frac{1}{3}}$;
$\frac{1}{3}c^{2}_{9}
$,
$\frac{1}{3} c_9^4 $,
$\frac{1}{3}c^{1}_{9}
$;
$\frac{1}{3}c^{1}_{9}
$,
$\frac{1}{3}c^{2}_{9}
$;
$\frac{1}{3} c_9^4 $)
 $\oplus$
$\mathrm{i}$($-\sqrt{\frac{1}{3}}$,
$\sqrt{\frac{2}{3}}$;\ \ 
$\sqrt{\frac{1}{3}}$)

Fail:
Tr$_I(C) = -1 <$  0 for I = [ 1/3 ]. Prop. B.4 (1) eqn. (B.18)

 \ \color{black}

\noindent 392: (dims,levels) = $(4 , 
2;9,
3
)$,
irreps = $4_{9,1}^{1,0}\oplus
2_{3}^{1,8}$,
pord$(\rho_\text{isum}(\mathfrak{t})) = 9$,

\vskip 0.7ex
\hangindent=5.5em \hangafter=1
{\white .}\hskip 1em $\rho_\text{isum}(\mathfrak{t})$ =
 $( 0,
\frac{1}{9},
\frac{4}{9},
\frac{7}{9} )
\oplus
( 0,
\frac{2}{3} )
$,

\vskip 0.7ex
\hangindent=5.5em \hangafter=1
{\white .}\hskip 1em $\rho_\text{isum}(\mathfrak{s})$ =
$\mathrm{i}$($0$,
$\sqrt{\frac{1}{3}}$,
$\sqrt{\frac{1}{3}}$,
$\sqrt{\frac{1}{3}}$;\ \ 
$-\frac{1}{3}c^{1}_{36}
$,
$\frac{1}{3}c^{1}_{36}
-\frac{1}{3}c^{5}_{36}
$,
$\frac{1}{3}c^{5}_{36}
$;\ \ 
$\frac{1}{3}c^{5}_{36}
$,
$-\frac{1}{3}c^{1}_{36}
$;\ \ 
$\frac{1}{3}c^{1}_{36}
-\frac{1}{3}c^{5}_{36}
$)
 $\oplus$
$\mathrm{i}$($\sqrt{\frac{1}{3}}$,
$\sqrt{\frac{2}{3}}$;\ \ 
$-\sqrt{\frac{1}{3}}$)

Fail:
Integral: $D_{\rho}(\sigma)_{\theta} \propto $ id,
 for all $\sigma$ and all $\theta$-eigenspaces that can contain unit. Prop. B.5 (6)

 \ \color{black}

\noindent 393: (dims,levels) = $(4 , 
2;9,
3
)$,
irreps = $4_{9,2}^{1,0}\oplus
2_{3}^{1,8}$,
pord$(\rho_\text{isum}(\mathfrak{t})) = 9$,

\vskip 0.7ex
\hangindent=5.5em \hangafter=1
{\white .}\hskip 1em $\rho_\text{isum}(\mathfrak{t})$ =
 $( 0,
\frac{1}{9},
\frac{4}{9},
\frac{7}{9} )
\oplus
( 0,
\frac{2}{3} )
$,

\vskip 0.7ex
\hangindent=5.5em \hangafter=1
{\white .}\hskip 1em $\rho_\text{isum}(\mathfrak{s})$ =
($0$,
$-\sqrt{\frac{1}{3}}$,
$-\sqrt{\frac{1}{3}}$,
$-\sqrt{\frac{1}{3}}$;
$\frac{1}{3}c^{2}_{9}
$,
$\frac{1}{3} c_9^4 $,
$\frac{1}{3}c^{1}_{9}
$;
$\frac{1}{3}c^{1}_{9}
$,
$\frac{1}{3}c^{2}_{9}
$;
$\frac{1}{3} c_9^4 $)
 $\oplus$
$\mathrm{i}$($\sqrt{\frac{1}{3}}$,
$\sqrt{\frac{2}{3}}$;\ \ 
$-\sqrt{\frac{1}{3}}$)

Fail:
Tr$_I(C) = -1 <$  0 for I = [ 2/3 ]. Prop. B.4 (1) eqn. (B.18)

 \ \color{black}

\noindent 394: (dims,levels) = $(4 , 
2;9,
3
)$,
irreps = $4_{9,1}^{2,0}\oplus
2_{3}^{1,0}$,
pord$(\rho_\text{isum}(\mathfrak{t})) = 9$,

\vskip 0.7ex
\hangindent=5.5em \hangafter=1
{\white .}\hskip 1em $\rho_\text{isum}(\mathfrak{t})$ =
 $( 0,
\frac{2}{9},
\frac{5}{9},
\frac{8}{9} )
\oplus
( 0,
\frac{1}{3} )
$,

\vskip 0.7ex
\hangindent=5.5em \hangafter=1
{\white .}\hskip 1em $\rho_\text{isum}(\mathfrak{s})$ =
$\mathrm{i}$($0$,
$\sqrt{\frac{1}{3}}$,
$\sqrt{\frac{1}{3}}$,
$\sqrt{\frac{1}{3}}$;\ \ 
$-\frac{1}{3}c^{1}_{36}
+\frac{1}{3}c^{5}_{36}
$,
$\frac{1}{3}c^{1}_{36}
$,
$-\frac{1}{3}c^{5}_{36}
$;\ \ 
$-\frac{1}{3}c^{5}_{36}
$,
$-\frac{1}{3}c^{1}_{36}
+\frac{1}{3}c^{5}_{36}
$;\ \ 
$\frac{1}{3}c^{1}_{36}
$)
 $\oplus$
$\mathrm{i}$($-\sqrt{\frac{1}{3}}$,
$\sqrt{\frac{2}{3}}$;\ \ 
$\sqrt{\frac{1}{3}}$)

Fail:
Integral: $D_{\rho}(\sigma)_{\theta} \propto $ id,
 for all $\sigma$ and all $\theta$-eigenspaces that can contain unit. Prop. B.5 (6)

 \ \color{black}

\noindent 395: (dims,levels) = $(4 , 
2;9,
3
)$,
irreps = $4_{9,2}^{5,0}\oplus
2_{3}^{1,0}$,
pord$(\rho_\text{isum}(\mathfrak{t})) = 9$,

\vskip 0.7ex
\hangindent=5.5em \hangafter=1
{\white .}\hskip 1em $\rho_\text{isum}(\mathfrak{t})$ =
 $( 0,
\frac{2}{9},
\frac{5}{9},
\frac{8}{9} )
\oplus
( 0,
\frac{1}{3} )
$,

\vskip 0.7ex
\hangindent=5.5em \hangafter=1
{\white .}\hskip 1em $\rho_\text{isum}(\mathfrak{s})$ =
($0$,
$-\sqrt{\frac{1}{3}}$,
$-\sqrt{\frac{1}{3}}$,
$-\sqrt{\frac{1}{3}}$;
$\frac{1}{3} c_9^4 $,
$\frac{1}{3}c^{2}_{9}
$,
$\frac{1}{3}c^{1}_{9}
$;
$\frac{1}{3}c^{1}_{9}
$,
$\frac{1}{3} c_9^4 $;
$\frac{1}{3}c^{2}_{9}
$)
 $\oplus$
$\mathrm{i}$($-\sqrt{\frac{1}{3}}$,
$\sqrt{\frac{2}{3}}$;\ \ 
$\sqrt{\frac{1}{3}}$)

Fail:
Tr$_I(C) = -1 <$  0 for I = [ 1/3 ]. Prop. B.4 (1) eqn. (B.18)

 \ \color{black}

\noindent 396: (dims,levels) = $(4 , 
2;9,
3
)$,
irreps = $4_{9,1}^{2,0}\oplus
2_{3}^{1,8}$,
pord$(\rho_\text{isum}(\mathfrak{t})) = 9$,

\vskip 0.7ex
\hangindent=5.5em \hangafter=1
{\white .}\hskip 1em $\rho_\text{isum}(\mathfrak{t})$ =
 $( 0,
\frac{2}{9},
\frac{5}{9},
\frac{8}{9} )
\oplus
( 0,
\frac{2}{3} )
$,

\vskip 0.7ex
\hangindent=5.5em \hangafter=1
{\white .}\hskip 1em $\rho_\text{isum}(\mathfrak{s})$ =
$\mathrm{i}$($0$,
$\sqrt{\frac{1}{3}}$,
$\sqrt{\frac{1}{3}}$,
$\sqrt{\frac{1}{3}}$;\ \ 
$-\frac{1}{3}c^{1}_{36}
+\frac{1}{3}c^{5}_{36}
$,
$\frac{1}{3}c^{1}_{36}
$,
$-\frac{1}{3}c^{5}_{36}
$;\ \ 
$-\frac{1}{3}c^{5}_{36}
$,
$-\frac{1}{3}c^{1}_{36}
+\frac{1}{3}c^{5}_{36}
$;\ \ 
$\frac{1}{3}c^{1}_{36}
$)
 $\oplus$
$\mathrm{i}$($\sqrt{\frac{1}{3}}$,
$\sqrt{\frac{2}{3}}$;\ \ 
$-\sqrt{\frac{1}{3}}$)

Fail:
Integral: $D_{\rho}(\sigma)_{\theta} \propto $ id,
 for all $\sigma$ and all $\theta$-eigenspaces that can contain unit. Prop. B.5 (6)

 \ \color{black}

\noindent 397: (dims,levels) = $(4 , 
2;9,
3
)$,
irreps = $4_{9,2}^{5,0}\oplus
2_{3}^{1,8}$,
pord$(\rho_\text{isum}(\mathfrak{t})) = 9$,

\vskip 0.7ex
\hangindent=5.5em \hangafter=1
{\white .}\hskip 1em $\rho_\text{isum}(\mathfrak{t})$ =
 $( 0,
\frac{2}{9},
\frac{5}{9},
\frac{8}{9} )
\oplus
( 0,
\frac{2}{3} )
$,

\vskip 0.7ex
\hangindent=5.5em \hangafter=1
{\white .}\hskip 1em $\rho_\text{isum}(\mathfrak{s})$ =
($0$,
$-\sqrt{\frac{1}{3}}$,
$-\sqrt{\frac{1}{3}}$,
$-\sqrt{\frac{1}{3}}$;
$\frac{1}{3} c_9^4 $,
$\frac{1}{3}c^{2}_{9}
$,
$\frac{1}{3}c^{1}_{9}
$;
$\frac{1}{3}c^{1}_{9}
$,
$\frac{1}{3} c_9^4 $;
$\frac{1}{3}c^{2}_{9}
$)
 $\oplus$
$\mathrm{i}$($\sqrt{\frac{1}{3}}$,
$\sqrt{\frac{2}{3}}$;\ \ 
$-\sqrt{\frac{1}{3}}$)

Fail:
Tr$_I(C) = -1 <$  0 for I = [ 2/3 ]. Prop. B.4 (1) eqn. (B.18)

 \ \color{black}

 \color{blue}

\noindent 398: (dims,levels) = $(4 , 
2;10,
5
)$,
irreps = $2_{5}^{1}
\hskip -1.5pt \otimes \hskip -1.5pt
2_{2}^{1,0}\oplus
2_{5}^{1}$,
pord$(\rho_\text{isum}(\mathfrak{t})) = 10$,

\vskip 0.7ex
\hangindent=5.5em \hangafter=1
{\white .}\hskip 1em $\rho_\text{isum}(\mathfrak{t})$ =
 $( \frac{1}{5},
\frac{4}{5},
\frac{3}{10},
\frac{7}{10} )
\oplus
( \frac{1}{5},
\frac{4}{5} )
$,

\vskip 0.7ex
\hangindent=5.5em \hangafter=1
{\white .}\hskip 1em $\rho_\text{isum}(\mathfrak{s})$ =
$\mathrm{i}$($\frac{1}{2\sqrt{5}}c^{3}_{20}
$,
$\frac{1}{2\sqrt{5}}c^{1}_{20}
$,
$\frac{3}{2\sqrt{15}}c^{1}_{20}
$,
$\frac{3}{2\sqrt{15}}c^{3}_{20}
$;\ \ 
$-\frac{1}{2\sqrt{5}}c^{3}_{20}
$,
$-\frac{3}{2\sqrt{15}}c^{3}_{20}
$,
$\frac{3}{2\sqrt{15}}c^{1}_{20}
$;\ \ 
$\frac{1}{2\sqrt{5}}c^{3}_{20}
$,
$-\frac{1}{2\sqrt{5}}c^{1}_{20}
$;\ \ 
$-\frac{1}{2\sqrt{5}}c^{3}_{20}
$)
 $\oplus$
$\mathrm{i}$($-\frac{1}{\sqrt{5}}c^{3}_{20}
$,
$\frac{1}{\sqrt{5}}c^{1}_{20}
$;\ \ 
$\frac{1}{\sqrt{5}}c^{3}_{20}
$)

Pass. 

 \ \color{black}

 \color{blue}

\noindent 399: (dims,levels) = $(4 , 
2;10,
5
)$,
irreps = $2_{5}^{2}
\hskip -1.5pt \otimes \hskip -1.5pt
2_{2}^{1,0}\oplus
2_{5}^{2}$,
pord$(\rho_\text{isum}(\mathfrak{t})) = 10$,

\vskip 0.7ex
\hangindent=5.5em \hangafter=1
{\white .}\hskip 1em $\rho_\text{isum}(\mathfrak{t})$ =
 $( \frac{2}{5},
\frac{3}{5},
\frac{1}{10},
\frac{9}{10} )
\oplus
( \frac{2}{5},
\frac{3}{5} )
$,

\vskip 0.7ex
\hangindent=5.5em \hangafter=1
{\white .}\hskip 1em $\rho_\text{isum}(\mathfrak{s})$ =
$\mathrm{i}$($\frac{1}{2\sqrt{5}}c^{1}_{20}
$,
$\frac{1}{2\sqrt{5}}c^{3}_{20}
$,
$\frac{3}{2\sqrt{15}}c^{3}_{20}
$,
$\frac{3}{2\sqrt{15}}c^{1}_{20}
$;\ \ 
$-\frac{1}{2\sqrt{5}}c^{1}_{20}
$,
$-\frac{3}{2\sqrt{15}}c^{1}_{20}
$,
$\frac{3}{2\sqrt{15}}c^{3}_{20}
$;\ \ 
$\frac{1}{2\sqrt{5}}c^{1}_{20}
$,
$-\frac{1}{2\sqrt{5}}c^{3}_{20}
$;\ \ 
$-\frac{1}{2\sqrt{5}}c^{1}_{20}
$)
 $\oplus$
$\mathrm{i}$($-\frac{1}{\sqrt{5}}c^{1}_{20}
$,
$\frac{1}{\sqrt{5}}c^{3}_{20}
$;\ \ 
$\frac{1}{\sqrt{5}}c^{1}_{20}
$)

Pass. 

 \ \color{black}

 \color{blue}

\noindent 400: (dims,levels) = $(4 , 
2;10,
10
)$,
irreps = $4_{5,1}^{1}
\hskip -1.5pt \otimes \hskip -1.5pt
1_{2}^{1,0}\oplus
2_{5}^{2}
\hskip -1.5pt \otimes \hskip -1.5pt
1_{2}^{1,0}$,
pord$(\rho_\text{isum}(\mathfrak{t})) = 5$,

\vskip 0.7ex
\hangindent=5.5em \hangafter=1
{\white .}\hskip 1em $\rho_\text{isum}(\mathfrak{t})$ =
 $( \frac{1}{10},
\frac{3}{10},
\frac{7}{10},
\frac{9}{10} )
\oplus
( \frac{1}{10},
\frac{9}{10} )
$,

\vskip 0.7ex
\hangindent=5.5em \hangafter=1
{\white .}\hskip 1em $\rho_\text{isum}(\mathfrak{s})$ =
$\mathrm{i}$($-\frac{1}{5}c^{1}_{20}
+\frac{1}{5}c^{3}_{20}
$,
$\frac{2}{5}c^{2}_{15}
+\frac{1}{5}c^{3}_{15}
$,
$-\frac{1}{5}+\frac{2}{5}c^{1}_{15}
-\frac{1}{5}c^{3}_{15}
$,
$\frac{1}{5}c^{1}_{20}
+\frac{1}{5}c^{3}_{20}
$;\ \ 
$\frac{1}{5}c^{1}_{20}
+\frac{1}{5}c^{3}_{20}
$,
$-\frac{1}{5}c^{1}_{20}
+\frac{1}{5}c^{3}_{20}
$,
$\frac{1}{5}-\frac{2}{5}c^{1}_{15}
+\frac{1}{5}c^{3}_{15}
$;\ \ 
$-\frac{1}{5}c^{1}_{20}
-\frac{1}{5}c^{3}_{20}
$,
$\frac{2}{5}c^{2}_{15}
+\frac{1}{5}c^{3}_{15}
$;\ \ 
$\frac{1}{5}c^{1}_{20}
-\frac{1}{5}c^{3}_{20}
$)
 $\oplus$
$\mathrm{i}$($-\frac{1}{\sqrt{5}}c^{1}_{20}
$,
$\frac{1}{\sqrt{5}}c^{3}_{20}
$;\ \ 
$\frac{1}{\sqrt{5}}c^{1}_{20}
$)

Pass. 

 \ \color{black}

 \color{blue}

\noindent 401: (dims,levels) = $(4 , 
2;10,
10
)$,
irreps = $4_{5,2}^{1}
\hskip -1.5pt \otimes \hskip -1.5pt
1_{2}^{1,0}\oplus
2_{5}^{2}
\hskip -1.5pt \otimes \hskip -1.5pt
1_{2}^{1,0}$,
pord$(\rho_\text{isum}(\mathfrak{t})) = 5$,

\vskip 0.7ex
\hangindent=5.5em \hangafter=1
{\white .}\hskip 1em $\rho_\text{isum}(\mathfrak{t})$ =
 $( \frac{1}{10},
\frac{3}{10},
\frac{7}{10},
\frac{9}{10} )
\oplus
( \frac{1}{10},
\frac{9}{10} )
$,

\vskip 0.7ex
\hangindent=5.5em \hangafter=1
{\white .}\hskip 1em $\rho_\text{isum}(\mathfrak{s})$ =
($\sqrt{\frac{1}{5}}$,
$\frac{-5+\sqrt{5}}{10}$,
$-\frac{5+\sqrt{5}}{10}$,
$\sqrt{\frac{1}{5}}$;
$-\sqrt{\frac{1}{5}}$,
$\sqrt{\frac{1}{5}}$,
$\frac{5+\sqrt{5}}{10}$;
$-\sqrt{\frac{1}{5}}$,
$\frac{5-\sqrt{5}}{10}$;
$\sqrt{\frac{1}{5}}$)
 $\oplus$
$\mathrm{i}$($-\frac{1}{\sqrt{5}}c^{1}_{20}
$,
$\frac{1}{\sqrt{5}}c^{3}_{20}
$;\ \ 
$\frac{1}{\sqrt{5}}c^{1}_{20}
$)

Pass. 

 \ \color{black}

 \color{blue}

\noindent 402: (dims,levels) = $(4 , 
2;10,
10
)$,
irreps = $4_{5,1}^{1}
\hskip -1.5pt \otimes \hskip -1.5pt
1_{2}^{1,0}\oplus
2_{5}^{1}
\hskip -1.5pt \otimes \hskip -1.5pt
1_{2}^{1,0}$,
pord$(\rho_\text{isum}(\mathfrak{t})) = 5$,

\vskip 0.7ex
\hangindent=5.5em \hangafter=1
{\white .}\hskip 1em $\rho_\text{isum}(\mathfrak{t})$ =
 $( \frac{1}{10},
\frac{3}{10},
\frac{7}{10},
\frac{9}{10} )
\oplus
( \frac{3}{10},
\frac{7}{10} )
$,

\vskip 0.7ex
\hangindent=5.5em \hangafter=1
{\white .}\hskip 1em $\rho_\text{isum}(\mathfrak{s})$ =
$\mathrm{i}$($-\frac{1}{5}c^{1}_{20}
+\frac{1}{5}c^{3}_{20}
$,
$\frac{2}{5}c^{2}_{15}
+\frac{1}{5}c^{3}_{15}
$,
$-\frac{1}{5}+\frac{2}{5}c^{1}_{15}
-\frac{1}{5}c^{3}_{15}
$,
$\frac{1}{5}c^{1}_{20}
+\frac{1}{5}c^{3}_{20}
$;\ \ 
$\frac{1}{5}c^{1}_{20}
+\frac{1}{5}c^{3}_{20}
$,
$-\frac{1}{5}c^{1}_{20}
+\frac{1}{5}c^{3}_{20}
$,
$\frac{1}{5}-\frac{2}{5}c^{1}_{15}
+\frac{1}{5}c^{3}_{15}
$;\ \ 
$-\frac{1}{5}c^{1}_{20}
-\frac{1}{5}c^{3}_{20}
$,
$\frac{2}{5}c^{2}_{15}
+\frac{1}{5}c^{3}_{15}
$;\ \ 
$\frac{1}{5}c^{1}_{20}
-\frac{1}{5}c^{3}_{20}
$)
 $\oplus$
$\mathrm{i}$($-\frac{1}{\sqrt{5}}c^{3}_{20}
$,
$\frac{1}{\sqrt{5}}c^{1}_{20}
$;\ \ 
$\frac{1}{\sqrt{5}}c^{3}_{20}
$)

Pass. 

 \ \color{black}

 \color{blue}

\noindent 403: (dims,levels) = $(4 , 
2;10,
10
)$,
irreps = $4_{5,2}^{1}
\hskip -1.5pt \otimes \hskip -1.5pt
1_{2}^{1,0}\oplus
2_{5}^{1}
\hskip -1.5pt \otimes \hskip -1.5pt
1_{2}^{1,0}$,
pord$(\rho_\text{isum}(\mathfrak{t})) = 5$,

\vskip 0.7ex
\hangindent=5.5em \hangafter=1
{\white .}\hskip 1em $\rho_\text{isum}(\mathfrak{t})$ =
 $( \frac{1}{10},
\frac{3}{10},
\frac{7}{10},
\frac{9}{10} )
\oplus
( \frac{3}{10},
\frac{7}{10} )
$,

\vskip 0.7ex
\hangindent=5.5em \hangafter=1
{\white .}\hskip 1em $\rho_\text{isum}(\mathfrak{s})$ =
($\sqrt{\frac{1}{5}}$,
$\frac{-5+\sqrt{5}}{10}$,
$-\frac{5+\sqrt{5}}{10}$,
$\sqrt{\frac{1}{5}}$;
$-\sqrt{\frac{1}{5}}$,
$\sqrt{\frac{1}{5}}$,
$\frac{5+\sqrt{5}}{10}$;
$-\sqrt{\frac{1}{5}}$,
$\frac{5-\sqrt{5}}{10}$;
$\sqrt{\frac{1}{5}}$)
 $\oplus$
$\mathrm{i}$($-\frac{1}{\sqrt{5}}c^{3}_{20}
$,
$\frac{1}{\sqrt{5}}c^{1}_{20}
$;\ \ 
$\frac{1}{\sqrt{5}}c^{3}_{20}
$)

Pass. 

 \ \color{black}

 \color{blue}

\noindent 404: (dims,levels) = $(4 , 
2;10,
10
)$,
irreps = $2_{5}^{1}
\hskip -1.5pt \otimes \hskip -1.5pt
2_{2}^{1,0}\oplus
2_{5}^{1}
\hskip -1.5pt \otimes \hskip -1.5pt
1_{2}^{1,0}$,
pord$(\rho_\text{isum}(\mathfrak{t})) = 10$,

\vskip 0.7ex
\hangindent=5.5em \hangafter=1
{\white .}\hskip 1em $\rho_\text{isum}(\mathfrak{t})$ =
 $( \frac{1}{5},
\frac{4}{5},
\frac{3}{10},
\frac{7}{10} )
\oplus
( \frac{3}{10},
\frac{7}{10} )
$,

\vskip 0.7ex
\hangindent=5.5em \hangafter=1
{\white .}\hskip 1em $\rho_\text{isum}(\mathfrak{s})$ =
$\mathrm{i}$($\frac{1}{2\sqrt{5}}c^{3}_{20}
$,
$\frac{1}{2\sqrt{5}}c^{1}_{20}
$,
$\frac{3}{2\sqrt{15}}c^{1}_{20}
$,
$\frac{3}{2\sqrt{15}}c^{3}_{20}
$;\ \ 
$-\frac{1}{2\sqrt{5}}c^{3}_{20}
$,
$-\frac{3}{2\sqrt{15}}c^{3}_{20}
$,
$\frac{3}{2\sqrt{15}}c^{1}_{20}
$;\ \ 
$\frac{1}{2\sqrt{5}}c^{3}_{20}
$,
$-\frac{1}{2\sqrt{5}}c^{1}_{20}
$;\ \ 
$-\frac{1}{2\sqrt{5}}c^{3}_{20}
$)
 $\oplus$
$\mathrm{i}$($-\frac{1}{\sqrt{5}}c^{3}_{20}
$,
$\frac{1}{\sqrt{5}}c^{1}_{20}
$;\ \ 
$\frac{1}{\sqrt{5}}c^{3}_{20}
$)

Pass. 

 \ \color{black}

 \color{blue}

\noindent 405: (dims,levels) = $(4 , 
2;10,
10
)$,
irreps = $2_{5}^{2}
\hskip -1.5pt \otimes \hskip -1.5pt
2_{2}^{1,0}\oplus
2_{5}^{2}
\hskip -1.5pt \otimes \hskip -1.5pt
1_{2}^{1,0}$,
pord$(\rho_\text{isum}(\mathfrak{t})) = 10$,

\vskip 0.7ex
\hangindent=5.5em \hangafter=1
{\white .}\hskip 1em $\rho_\text{isum}(\mathfrak{t})$ =
 $( \frac{2}{5},
\frac{3}{5},
\frac{1}{10},
\frac{9}{10} )
\oplus
( \frac{1}{10},
\frac{9}{10} )
$,

\vskip 0.7ex
\hangindent=5.5em \hangafter=1
{\white .}\hskip 1em $\rho_\text{isum}(\mathfrak{s})$ =
$\mathrm{i}$($\frac{1}{2\sqrt{5}}c^{1}_{20}
$,
$\frac{1}{2\sqrt{5}}c^{3}_{20}
$,
$\frac{3}{2\sqrt{15}}c^{3}_{20}
$,
$\frac{3}{2\sqrt{15}}c^{1}_{20}
$;\ \ 
$-\frac{1}{2\sqrt{5}}c^{1}_{20}
$,
$-\frac{3}{2\sqrt{15}}c^{1}_{20}
$,
$\frac{3}{2\sqrt{15}}c^{3}_{20}
$;\ \ 
$\frac{1}{2\sqrt{5}}c^{1}_{20}
$,
$-\frac{1}{2\sqrt{5}}c^{3}_{20}
$;\ \ 
$-\frac{1}{2\sqrt{5}}c^{1}_{20}
$)
 $\oplus$
$\mathrm{i}$($-\frac{1}{\sqrt{5}}c^{1}_{20}
$,
$\frac{1}{\sqrt{5}}c^{3}_{20}
$;\ \ 
$\frac{1}{\sqrt{5}}c^{1}_{20}
$)

Pass. 

 \ \color{black}

\noindent 406: (dims,levels) = $(4 , 
2;14,
2
)$,
irreps = $4_{7}^{1}
\hskip -1.5pt \otimes \hskip -1.5pt
1_{2}^{1,0}\oplus
2_{2}^{1,0}$,
pord$(\rho_\text{isum}(\mathfrak{t})) = 14$,

\vskip 0.7ex
\hangindent=5.5em \hangafter=1
{\white .}\hskip 1em $\rho_\text{isum}(\mathfrak{t})$ =
 $( \frac{1}{2},
\frac{1}{14},
\frac{9}{14},
\frac{11}{14} )
\oplus
( 0,
\frac{1}{2} )
$,

\vskip 0.7ex
\hangindent=5.5em \hangafter=1
{\white .}\hskip 1em $\rho_\text{isum}(\mathfrak{s})$ =
$\mathrm{i}$($\sqrt{\frac{1}{7}}$,
$\sqrt{\frac{2}{7}}$,
$\sqrt{\frac{2}{7}}$,
$\sqrt{\frac{2}{7}}$;\ \ 
$\frac{1}{\sqrt{7}}c^{1}_{7}
$,
$-\frac{1}{\sqrt{7}\mathrm{i}}s^{5}_{28}
$,
$\frac{1}{\sqrt{7}}c^{2}_{7}
$;\ \ 
$\frac{1}{\sqrt{7}}c^{2}_{7}
$,
$\frac{1}{\sqrt{7}}c^{1}_{7}
$;\ \ 
$-\frac{1}{\sqrt{7}\mathrm{i}}s^{5}_{28}
$)
 $\oplus$
($-\frac{1}{2}$,
$-\sqrt{\frac{3}{4}}$;
$\frac{1}{2}$)

Fail:
Tr$_I(C) = -1 <$  0 for I = [ 0 ]. Prop. B.4 (1) eqn. (B.18)

 \ \color{black}

\noindent 407: (dims,levels) = $(4 , 
2;14,
2
)$,
irreps = $4_{7}^{3}
\hskip -1.5pt \otimes \hskip -1.5pt
1_{2}^{1,0}\oplus
2_{2}^{1,0}$,
pord$(\rho_\text{isum}(\mathfrak{t})) = 14$,

\vskip 0.7ex
\hangindent=5.5em \hangafter=1
{\white .}\hskip 1em $\rho_\text{isum}(\mathfrak{t})$ =
 $( \frac{1}{2},
\frac{3}{14},
\frac{5}{14},
\frac{13}{14} )
\oplus
( 0,
\frac{1}{2} )
$,

\vskip 0.7ex
\hangindent=5.5em \hangafter=1
{\white .}\hskip 1em $\rho_\text{isum}(\mathfrak{s})$ =
$\mathrm{i}$($-\sqrt{\frac{1}{7}}$,
$\sqrt{\frac{2}{7}}$,
$\sqrt{\frac{2}{7}}$,
$\sqrt{\frac{2}{7}}$;\ \ 
$\frac{1}{\sqrt{7}\mathrm{i}}s^{5}_{28}
$,
$-\frac{1}{\sqrt{7}}c^{1}_{7}
$,
$-\frac{1}{\sqrt{7}}c^{2}_{7}
$;\ \ 
$-\frac{1}{\sqrt{7}}c^{2}_{7}
$,
$\frac{1}{\sqrt{7}\mathrm{i}}s^{5}_{28}
$;\ \ 
$-\frac{1}{\sqrt{7}}c^{1}_{7}
$)
 $\oplus$
($-\frac{1}{2}$,
$-\sqrt{\frac{3}{4}}$;
$\frac{1}{2}$)

Fail:
Tr$_I(C) = -1 <$  0 for I = [ 0 ]. Prop. B.4 (1) eqn. (B.18)

 \ \color{black}

 \color{blue}

\noindent 408: (dims,levels) = $(4 , 
2;14,
6
)$,
irreps = $4_{7}^{1}
\hskip -1.5pt \otimes \hskip -1.5pt
1_{2}^{1,0}\oplus
2_{3}^{1,8}
\hskip -1.5pt \otimes \hskip -1.5pt
1_{2}^{1,0}$,
pord$(\rho_\text{isum}(\mathfrak{t})) = 21$,

\vskip 0.7ex
\hangindent=5.5em \hangafter=1
{\white .}\hskip 1em $\rho_\text{isum}(\mathfrak{t})$ =
 $( \frac{1}{2},
\frac{1}{14},
\frac{9}{14},
\frac{11}{14} )
\oplus
( \frac{1}{2},
\frac{1}{6} )
$,

\vskip 0.7ex
\hangindent=5.5em \hangafter=1
{\white .}\hskip 1em $\rho_\text{isum}(\mathfrak{s})$ =
$\mathrm{i}$($\sqrt{\frac{1}{7}}$,
$\sqrt{\frac{2}{7}}$,
$\sqrt{\frac{2}{7}}$,
$\sqrt{\frac{2}{7}}$;\ \ 
$\frac{1}{\sqrt{7}}c^{1}_{7}
$,
$-\frac{1}{\sqrt{7}\mathrm{i}}s^{5}_{28}
$,
$\frac{1}{\sqrt{7}}c^{2}_{7}
$;\ \ 
$\frac{1}{\sqrt{7}}c^{2}_{7}
$,
$\frac{1}{\sqrt{7}}c^{1}_{7}
$;\ \ 
$-\frac{1}{\sqrt{7}\mathrm{i}}s^{5}_{28}
$)
 $\oplus$
$\mathrm{i}$($-\sqrt{\frac{1}{3}}$,
$\sqrt{\frac{2}{3}}$;\ \ 
$\sqrt{\frac{1}{3}}$)

Pass. 

 \ \color{black}

 \color{blue}

\noindent 409: (dims,levels) = $(4 , 
2;14,
6
)$,
irreps = $4_{7}^{1}
\hskip -1.5pt \otimes \hskip -1.5pt
1_{2}^{1,0}\oplus
2_{3}^{1,0}
\hskip -1.5pt \otimes \hskip -1.5pt
1_{2}^{1,0}$,
pord$(\rho_\text{isum}(\mathfrak{t})) = 21$,

\vskip 0.7ex
\hangindent=5.5em \hangafter=1
{\white .}\hskip 1em $\rho_\text{isum}(\mathfrak{t})$ =
 $( \frac{1}{2},
\frac{1}{14},
\frac{9}{14},
\frac{11}{14} )
\oplus
( \frac{1}{2},
\frac{5}{6} )
$,

\vskip 0.7ex
\hangindent=5.5em \hangafter=1
{\white .}\hskip 1em $\rho_\text{isum}(\mathfrak{s})$ =
$\mathrm{i}$($\sqrt{\frac{1}{7}}$,
$\sqrt{\frac{2}{7}}$,
$\sqrt{\frac{2}{7}}$,
$\sqrt{\frac{2}{7}}$;\ \ 
$\frac{1}{\sqrt{7}}c^{1}_{7}
$,
$-\frac{1}{\sqrt{7}\mathrm{i}}s^{5}_{28}
$,
$\frac{1}{\sqrt{7}}c^{2}_{7}
$;\ \ 
$\frac{1}{\sqrt{7}}c^{2}_{7}
$,
$\frac{1}{\sqrt{7}}c^{1}_{7}
$;\ \ 
$-\frac{1}{\sqrt{7}\mathrm{i}}s^{5}_{28}
$)
 $\oplus$
$\mathrm{i}$($\sqrt{\frac{1}{3}}$,
$\sqrt{\frac{2}{3}}$;\ \ 
$-\sqrt{\frac{1}{3}}$)

Pass. 

 \ \color{black}

 \color{blue}

\noindent 410: (dims,levels) = $(4 , 
2;14,
6
)$,
irreps = $4_{7}^{3}
\hskip -1.5pt \otimes \hskip -1.5pt
1_{2}^{1,0}\oplus
2_{3}^{1,8}
\hskip -1.5pt \otimes \hskip -1.5pt
1_{2}^{1,0}$,
pord$(\rho_\text{isum}(\mathfrak{t})) = 21$,

\vskip 0.7ex
\hangindent=5.5em \hangafter=1
{\white .}\hskip 1em $\rho_\text{isum}(\mathfrak{t})$ =
 $( \frac{1}{2},
\frac{3}{14},
\frac{5}{14},
\frac{13}{14} )
\oplus
( \frac{1}{2},
\frac{1}{6} )
$,

\vskip 0.7ex
\hangindent=5.5em \hangafter=1
{\white .}\hskip 1em $\rho_\text{isum}(\mathfrak{s})$ =
$\mathrm{i}$($-\sqrt{\frac{1}{7}}$,
$\sqrt{\frac{2}{7}}$,
$\sqrt{\frac{2}{7}}$,
$\sqrt{\frac{2}{7}}$;\ \ 
$\frac{1}{\sqrt{7}\mathrm{i}}s^{5}_{28}
$,
$-\frac{1}{\sqrt{7}}c^{1}_{7}
$,
$-\frac{1}{\sqrt{7}}c^{2}_{7}
$;\ \ 
$-\frac{1}{\sqrt{7}}c^{2}_{7}
$,
$\frac{1}{\sqrt{7}\mathrm{i}}s^{5}_{28}
$;\ \ 
$-\frac{1}{\sqrt{7}}c^{1}_{7}
$)
 $\oplus$
$\mathrm{i}$($-\sqrt{\frac{1}{3}}$,
$\sqrt{\frac{2}{3}}$;\ \ 
$\sqrt{\frac{1}{3}}$)

Pass. 

 \ \color{black}

 \color{blue}

\noindent 411: (dims,levels) = $(4 , 
2;14,
6
)$,
irreps = $4_{7}^{3}
\hskip -1.5pt \otimes \hskip -1.5pt
1_{2}^{1,0}\oplus
2_{3}^{1,0}
\hskip -1.5pt \otimes \hskip -1.5pt
1_{2}^{1,0}$,
pord$(\rho_\text{isum}(\mathfrak{t})) = 21$,

\vskip 0.7ex
\hangindent=5.5em \hangafter=1
{\white .}\hskip 1em $\rho_\text{isum}(\mathfrak{t})$ =
 $( \frac{1}{2},
\frac{3}{14},
\frac{5}{14},
\frac{13}{14} )
\oplus
( \frac{1}{2},
\frac{5}{6} )
$,

\vskip 0.7ex
\hangindent=5.5em \hangafter=1
{\white .}\hskip 1em $\rho_\text{isum}(\mathfrak{s})$ =
$\mathrm{i}$($-\sqrt{\frac{1}{7}}$,
$\sqrt{\frac{2}{7}}$,
$\sqrt{\frac{2}{7}}$,
$\sqrt{\frac{2}{7}}$;\ \ 
$\frac{1}{\sqrt{7}\mathrm{i}}s^{5}_{28}
$,
$-\frac{1}{\sqrt{7}}c^{1}_{7}
$,
$-\frac{1}{\sqrt{7}}c^{2}_{7}
$;\ \ 
$-\frac{1}{\sqrt{7}}c^{2}_{7}
$,
$\frac{1}{\sqrt{7}\mathrm{i}}s^{5}_{28}
$;\ \ 
$-\frac{1}{\sqrt{7}}c^{1}_{7}
$)
 $\oplus$
$\mathrm{i}$($\sqrt{\frac{1}{3}}$,
$\sqrt{\frac{2}{3}}$;\ \ 
$-\sqrt{\frac{1}{3}}$)

Pass. 

 \ \color{black}

 \color{blue}

\noindent 412: (dims,levels) = $(4 , 
2;15,
5
)$,
irreps = $2_{5}^{1}
\hskip -1.5pt \otimes \hskip -1.5pt
2_{3}^{1,0}\oplus
2_{5}^{1}$,
pord$(\rho_\text{isum}(\mathfrak{t})) = 15$,

\vskip 0.7ex
\hangindent=5.5em \hangafter=1
{\white .}\hskip 1em $\rho_\text{isum}(\mathfrak{t})$ =
 $( \frac{1}{5},
\frac{4}{5},
\frac{2}{15},
\frac{8}{15} )
\oplus
( \frac{1}{5},
\frac{4}{5} )
$,

\vskip 0.7ex
\hangindent=5.5em \hangafter=1
{\white .}\hskip 1em $\rho_\text{isum}(\mathfrak{s})$ =
($-\frac{1}{\sqrt{15}}c^{3}_{20}
$,
$\frac{1}{\sqrt{15}}c^{1}_{20}
$,
$\frac{2}{\sqrt{30}}c^{1}_{20}
$,
$-\frac{2}{\sqrt{30}}c^{3}_{20}
$;
$\frac{1}{\sqrt{15}}c^{3}_{20}
$,
$\frac{2}{\sqrt{30}}c^{3}_{20}
$,
$\frac{2}{\sqrt{30}}c^{1}_{20}
$;
$-\frac{1}{\sqrt{15}}c^{3}_{20}
$,
$-\frac{1}{\sqrt{15}}c^{1}_{20}
$;
$\frac{1}{\sqrt{15}}c^{3}_{20}
$)
 $\oplus$
$\mathrm{i}$($-\frac{1}{\sqrt{5}}c^{3}_{20}
$,
$\frac{1}{\sqrt{5}}c^{1}_{20}
$;\ \ 
$\frac{1}{\sqrt{5}}c^{3}_{20}
$)

Pass. 

 \ \color{black}

 \color{blue}

\noindent 413: (dims,levels) = $(4 , 
2;15,
5
)$,
irreps = $2_{5}^{2}
\hskip -1.5pt \otimes \hskip -1.5pt
2_{3}^{1,0}\oplus
2_{5}^{2}$,
pord$(\rho_\text{isum}(\mathfrak{t})) = 15$,

\vskip 0.7ex
\hangindent=5.5em \hangafter=1
{\white .}\hskip 1em $\rho_\text{isum}(\mathfrak{t})$ =
 $( \frac{2}{5},
\frac{3}{5},
\frac{11}{15},
\frac{14}{15} )
\oplus
( \frac{2}{5},
\frac{3}{5} )
$,

\vskip 0.7ex
\hangindent=5.5em \hangafter=1
{\white .}\hskip 1em $\rho_\text{isum}(\mathfrak{s})$ =
($-\frac{1}{\sqrt{15}}c^{1}_{20}
$,
$-\frac{1}{\sqrt{15}}c^{3}_{20}
$,
$\frac{2}{\sqrt{30}}c^{1}_{20}
$,
$-\frac{2}{\sqrt{30}}c^{3}_{20}
$;
$\frac{1}{\sqrt{15}}c^{1}_{20}
$,
$\frac{2}{\sqrt{30}}c^{3}_{20}
$,
$\frac{2}{\sqrt{30}}c^{1}_{20}
$;
$\frac{1}{\sqrt{15}}c^{1}_{20}
$,
$-\frac{1}{\sqrt{15}}c^{3}_{20}
$;
$-\frac{1}{\sqrt{15}}c^{1}_{20}
$)
 $\oplus$
$\mathrm{i}$($-\frac{1}{\sqrt{5}}c^{1}_{20}
$,
$\frac{1}{\sqrt{5}}c^{3}_{20}
$;\ \ 
$\frac{1}{\sqrt{5}}c^{1}_{20}
$)

Pass. 

 \ \color{black}

 \color{blue}

\noindent 414: (dims,levels) = $(4 , 
2;15,
15
)$,
irreps = $4_{5,1}^{1}
\hskip -1.5pt \otimes \hskip -1.5pt
1_{3}^{1,0}\oplus
2_{5}^{1}
\hskip -1.5pt \otimes \hskip -1.5pt
1_{3}^{1,0}$,
pord$(\rho_\text{isum}(\mathfrak{t})) = 5$,

\vskip 0.7ex
\hangindent=5.5em \hangafter=1
{\white .}\hskip 1em $\rho_\text{isum}(\mathfrak{t})$ =
 $( \frac{2}{15},
\frac{8}{15},
\frac{11}{15},
\frac{14}{15} )
\oplus
( \frac{2}{15},
\frac{8}{15} )
$,

\vskip 0.7ex
\hangindent=5.5em \hangafter=1
{\white .}\hskip 1em $\rho_\text{isum}(\mathfrak{s})$ =
$\mathrm{i}$($-\frac{1}{5}c^{1}_{20}
-\frac{1}{5}c^{3}_{20}
$,
$\frac{1}{5}c^{1}_{20}
-\frac{1}{5}c^{3}_{20}
$,
$-\frac{1}{5}+\frac{2}{5}c^{1}_{15}
-\frac{1}{5}c^{3}_{15}
$,
$\frac{2}{5}c^{2}_{15}
+\frac{1}{5}c^{3}_{15}
$;\ \ 
$\frac{1}{5}c^{1}_{20}
+\frac{1}{5}c^{3}_{20}
$,
$-\frac{2}{5}c^{2}_{15}
-\frac{1}{5}c^{3}_{15}
$,
$-\frac{1}{5}+\frac{2}{5}c^{1}_{15}
-\frac{1}{5}c^{3}_{15}
$;\ \ 
$-\frac{1}{5}c^{1}_{20}
+\frac{1}{5}c^{3}_{20}
$,
$\frac{1}{5}c^{1}_{20}
+\frac{1}{5}c^{3}_{20}
$;\ \ 
$\frac{1}{5}c^{1}_{20}
-\frac{1}{5}c^{3}_{20}
$)
 $\oplus$
$\mathrm{i}$($\frac{1}{\sqrt{5}}c^{3}_{20}
$,
$\frac{1}{\sqrt{5}}c^{1}_{20}
$;\ \ 
$-\frac{1}{\sqrt{5}}c^{3}_{20}
$)

Pass. 

 \ \color{black}

 \color{blue}

\noindent 415: (dims,levels) = $(4 , 
2;15,
15
)$,
irreps = $4_{5,2}^{1}
\hskip -1.5pt \otimes \hskip -1.5pt
1_{3}^{1,0}\oplus
2_{5}^{1}
\hskip -1.5pt \otimes \hskip -1.5pt
1_{3}^{1,0}$,
pord$(\rho_\text{isum}(\mathfrak{t})) = 5$,

\vskip 0.7ex
\hangindent=5.5em \hangafter=1
{\white .}\hskip 1em $\rho_\text{isum}(\mathfrak{t})$ =
 $( \frac{2}{15},
\frac{8}{15},
\frac{11}{15},
\frac{14}{15} )
\oplus
( \frac{2}{15},
\frac{8}{15} )
$,

\vskip 0.7ex
\hangindent=5.5em \hangafter=1
{\white .}\hskip 1em $\rho_\text{isum}(\mathfrak{s})$ =
($\sqrt{\frac{1}{5}}$,
$\sqrt{\frac{1}{5}}$,
$-\frac{5+\sqrt{5}}{10}$,
$\frac{-5+\sqrt{5}}{10}$;
$\sqrt{\frac{1}{5}}$,
$\frac{5-\sqrt{5}}{10}$,
$\frac{5+\sqrt{5}}{10}$;
$-\sqrt{\frac{1}{5}}$,
$\sqrt{\frac{1}{5}}$;
$-\sqrt{\frac{1}{5}}$)
 $\oplus$
$\mathrm{i}$($\frac{1}{\sqrt{5}}c^{3}_{20}
$,
$\frac{1}{\sqrt{5}}c^{1}_{20}
$;\ \ 
$-\frac{1}{\sqrt{5}}c^{3}_{20}
$)

Pass. 

 \ \color{black}

 \color{blue}

\noindent 416: (dims,levels) = $(4 , 
2;15,
15
)$,
irreps = $4_{5,1}^{1}
\hskip -1.5pt \otimes \hskip -1.5pt
1_{3}^{1,0}\oplus
2_{5}^{2}
\hskip -1.5pt \otimes \hskip -1.5pt
1_{3}^{1,0}$,
pord$(\rho_\text{isum}(\mathfrak{t})) = 5$,

\vskip 0.7ex
\hangindent=5.5em \hangafter=1
{\white .}\hskip 1em $\rho_\text{isum}(\mathfrak{t})$ =
 $( \frac{2}{15},
\frac{8}{15},
\frac{11}{15},
\frac{14}{15} )
\oplus
( \frac{11}{15},
\frac{14}{15} )
$,

\vskip 0.7ex
\hangindent=5.5em \hangafter=1
{\white .}\hskip 1em $\rho_\text{isum}(\mathfrak{s})$ =
$\mathrm{i}$($-\frac{1}{5}c^{1}_{20}
-\frac{1}{5}c^{3}_{20}
$,
$\frac{1}{5}c^{1}_{20}
-\frac{1}{5}c^{3}_{20}
$,
$-\frac{1}{5}+\frac{2}{5}c^{1}_{15}
-\frac{1}{5}c^{3}_{15}
$,
$\frac{2}{5}c^{2}_{15}
+\frac{1}{5}c^{3}_{15}
$;\ \ 
$\frac{1}{5}c^{1}_{20}
+\frac{1}{5}c^{3}_{20}
$,
$-\frac{2}{5}c^{2}_{15}
-\frac{1}{5}c^{3}_{15}
$,
$-\frac{1}{5}+\frac{2}{5}c^{1}_{15}
-\frac{1}{5}c^{3}_{15}
$;\ \ 
$-\frac{1}{5}c^{1}_{20}
+\frac{1}{5}c^{3}_{20}
$,
$\frac{1}{5}c^{1}_{20}
+\frac{1}{5}c^{3}_{20}
$;\ \ 
$\frac{1}{5}c^{1}_{20}
-\frac{1}{5}c^{3}_{20}
$)
 $\oplus$
$\mathrm{i}$($-\frac{1}{\sqrt{5}}c^{1}_{20}
$,
$\frac{1}{\sqrt{5}}c^{3}_{20}
$;\ \ 
$\frac{1}{\sqrt{5}}c^{1}_{20}
$)

Pass. 

 \ \color{black}

 \color{blue}

\noindent 417: (dims,levels) = $(4 , 
2;15,
15
)$,
irreps = $4_{5,2}^{1}
\hskip -1.5pt \otimes \hskip -1.5pt
1_{3}^{1,0}\oplus
2_{5}^{2}
\hskip -1.5pt \otimes \hskip -1.5pt
1_{3}^{1,0}$,
pord$(\rho_\text{isum}(\mathfrak{t})) = 5$,

\vskip 0.7ex
\hangindent=5.5em \hangafter=1
{\white .}\hskip 1em $\rho_\text{isum}(\mathfrak{t})$ =
 $( \frac{2}{15},
\frac{8}{15},
\frac{11}{15},
\frac{14}{15} )
\oplus
( \frac{11}{15},
\frac{14}{15} )
$,

\vskip 0.7ex
\hangindent=5.5em \hangafter=1
{\white .}\hskip 1em $\rho_\text{isum}(\mathfrak{s})$ =
($\sqrt{\frac{1}{5}}$,
$\sqrt{\frac{1}{5}}$,
$-\frac{5+\sqrt{5}}{10}$,
$\frac{-5+\sqrt{5}}{10}$;
$\sqrt{\frac{1}{5}}$,
$\frac{5-\sqrt{5}}{10}$,
$\frac{5+\sqrt{5}}{10}$;
$-\sqrt{\frac{1}{5}}$,
$\sqrt{\frac{1}{5}}$;
$-\sqrt{\frac{1}{5}}$)
 $\oplus$
$\mathrm{i}$($-\frac{1}{\sqrt{5}}c^{1}_{20}
$,
$\frac{1}{\sqrt{5}}c^{3}_{20}
$;\ \ 
$\frac{1}{\sqrt{5}}c^{1}_{20}
$)

Pass. 

 \ \color{black}

 \color{blue}

\noindent 418: (dims,levels) = $(4 , 
2;15,
15
)$,
irreps = $2_{5}^{1}
\hskip -1.5pt \otimes \hskip -1.5pt
2_{3}^{1,0}\oplus
2_{5}^{1}
\hskip -1.5pt \otimes \hskip -1.5pt
1_{3}^{1,0}$,
pord$(\rho_\text{isum}(\mathfrak{t})) = 15$,

\vskip 0.7ex
\hangindent=5.5em \hangafter=1
{\white .}\hskip 1em $\rho_\text{isum}(\mathfrak{t})$ =
 $( \frac{1}{5},
\frac{4}{5},
\frac{2}{15},
\frac{8}{15} )
\oplus
( \frac{2}{15},
\frac{8}{15} )
$,

\vskip 0.7ex
\hangindent=5.5em \hangafter=1
{\white .}\hskip 1em $\rho_\text{isum}(\mathfrak{s})$ =
($-\frac{1}{\sqrt{15}}c^{3}_{20}
$,
$\frac{1}{\sqrt{15}}c^{1}_{20}
$,
$\frac{2}{\sqrt{30}}c^{1}_{20}
$,
$-\frac{2}{\sqrt{30}}c^{3}_{20}
$;
$\frac{1}{\sqrt{15}}c^{3}_{20}
$,
$\frac{2}{\sqrt{30}}c^{3}_{20}
$,
$\frac{2}{\sqrt{30}}c^{1}_{20}
$;
$-\frac{1}{\sqrt{15}}c^{3}_{20}
$,
$-\frac{1}{\sqrt{15}}c^{1}_{20}
$;
$\frac{1}{\sqrt{15}}c^{3}_{20}
$)
 $\oplus$
$\mathrm{i}$($\frac{1}{\sqrt{5}}c^{3}_{20}
$,
$\frac{1}{\sqrt{5}}c^{1}_{20}
$;\ \ 
$-\frac{1}{\sqrt{5}}c^{3}_{20}
$)

Pass. 

 \ \color{black}

 \color{blue}

\noindent 419: (dims,levels) = $(4 , 
2;15,
15
)$,
irreps = $2_{5}^{2}
\hskip -1.5pt \otimes \hskip -1.5pt
2_{3}^{1,0}\oplus
2_{5}^{2}
\hskip -1.5pt \otimes \hskip -1.5pt
1_{3}^{1,0}$,
pord$(\rho_\text{isum}(\mathfrak{t})) = 15$,

\vskip 0.7ex
\hangindent=5.5em \hangafter=1
{\white .}\hskip 1em $\rho_\text{isum}(\mathfrak{t})$ =
 $( \frac{2}{5},
\frac{3}{5},
\frac{11}{15},
\frac{14}{15} )
\oplus
( \frac{11}{15},
\frac{14}{15} )
$,

\vskip 0.7ex
\hangindent=5.5em \hangafter=1
{\white .}\hskip 1em $\rho_\text{isum}(\mathfrak{s})$ =
($-\frac{1}{\sqrt{15}}c^{1}_{20}
$,
$-\frac{1}{\sqrt{15}}c^{3}_{20}
$,
$\frac{2}{\sqrt{30}}c^{1}_{20}
$,
$-\frac{2}{\sqrt{30}}c^{3}_{20}
$;
$\frac{1}{\sqrt{15}}c^{1}_{20}
$,
$\frac{2}{\sqrt{30}}c^{3}_{20}
$,
$\frac{2}{\sqrt{30}}c^{1}_{20}
$;
$\frac{1}{\sqrt{15}}c^{1}_{20}
$,
$-\frac{1}{\sqrt{15}}c^{3}_{20}
$;
$-\frac{1}{\sqrt{15}}c^{1}_{20}
$)
 $\oplus$
$\mathrm{i}$($-\frac{1}{\sqrt{5}}c^{1}_{20}
$,
$\frac{1}{\sqrt{5}}c^{3}_{20}
$;\ \ 
$\frac{1}{\sqrt{5}}c^{1}_{20}
$)

Pass. 

 \ \color{black}

\noindent 420: (dims,levels) = $(4 , 
2;18,
2
)$,
irreps = $4_{9,1}^{2,0}
\hskip -1.5pt \otimes \hskip -1.5pt
1_{2}^{1,0}\oplus
2_{2}^{1,0}$,
pord$(\rho_\text{isum}(\mathfrak{t})) = 18$,

\vskip 0.7ex
\hangindent=5.5em \hangafter=1
{\white .}\hskip 1em $\rho_\text{isum}(\mathfrak{t})$ =
 $( \frac{1}{2},
\frac{1}{18},
\frac{7}{18},
\frac{13}{18} )
\oplus
( 0,
\frac{1}{2} )
$,

\vskip 0.7ex
\hangindent=5.5em \hangafter=1
{\white .}\hskip 1em $\rho_\text{isum}(\mathfrak{s})$ =
$\mathrm{i}$($0$,
$\sqrt{\frac{1}{3}}$,
$\sqrt{\frac{1}{3}}$,
$\sqrt{\frac{1}{3}}$;\ \ 
$\frac{1}{3}c^{5}_{36}
$,
$\frac{1}{3}c^{1}_{36}
-\frac{1}{3}c^{5}_{36}
$,
$-\frac{1}{3}c^{1}_{36}
$;\ \ 
$-\frac{1}{3}c^{1}_{36}
$,
$\frac{1}{3}c^{5}_{36}
$;\ \ 
$\frac{1}{3}c^{1}_{36}
-\frac{1}{3}c^{5}_{36}
$)
 $\oplus$
($-\frac{1}{2}$,
$-\sqrt{\frac{3}{4}}$;
$\frac{1}{2}$)

Fail:
Tr$_I(C) = -1 <$  0 for I = [ 0 ]. Prop. B.4 (1) eqn. (B.18)

 \ \color{black}

\noindent 421: (dims,levels) = $(4 , 
2;18,
2
)$,
irreps = $4_{9,2}^{5,0}
\hskip -1.5pt \otimes \hskip -1.5pt
1_{2}^{1,0}\oplus
2_{2}^{1,0}$,
pord$(\rho_\text{isum}(\mathfrak{t})) = 18$,

\vskip 0.7ex
\hangindent=5.5em \hangafter=1
{\white .}\hskip 1em $\rho_\text{isum}(\mathfrak{t})$ =
 $( \frac{1}{2},
\frac{1}{18},
\frac{7}{18},
\frac{13}{18} )
\oplus
( 0,
\frac{1}{2} )
$,

\vskip 0.7ex
\hangindent=5.5em \hangafter=1
{\white .}\hskip 1em $\rho_\text{isum}(\mathfrak{s})$ =
($0$,
$-\sqrt{\frac{1}{3}}$,
$-\sqrt{\frac{1}{3}}$,
$-\sqrt{\frac{1}{3}}$;
$-\frac{1}{3}c^{1}_{9}
$,
$-\frac{1}{3} c_9^4 $,
$-\frac{1}{3}c^{2}_{9}
$;
$-\frac{1}{3}c^{2}_{9}
$,
$-\frac{1}{3}c^{1}_{9}
$;
$-\frac{1}{3} c_9^4 $)
 $\oplus$
($-\frac{1}{2}$,
$-\sqrt{\frac{3}{4}}$;
$\frac{1}{2}$)

Fail:
Integral: $D_{\rho}(\sigma)_{\theta} \propto $ id,
 for all $\sigma$ and all $\theta$-eigenspaces that can contain unit. Prop. B.5 (6)

 \ \color{black}

\noindent 422: (dims,levels) = $(4 , 
2;18,
2
)$,
irreps = $4_{9,1}^{1,0}
\hskip -1.5pt \otimes \hskip -1.5pt
1_{2}^{1,0}\oplus
2_{2}^{1,0}$,
pord$(\rho_\text{isum}(\mathfrak{t})) = 18$,

\vskip 0.7ex
\hangindent=5.5em \hangafter=1
{\white .}\hskip 1em $\rho_\text{isum}(\mathfrak{t})$ =
 $( \frac{1}{2},
\frac{5}{18},
\frac{11}{18},
\frac{17}{18} )
\oplus
( 0,
\frac{1}{2} )
$,

\vskip 0.7ex
\hangindent=5.5em \hangafter=1
{\white .}\hskip 1em $\rho_\text{isum}(\mathfrak{s})$ =
$\mathrm{i}$($0$,
$\sqrt{\frac{1}{3}}$,
$\sqrt{\frac{1}{3}}$,
$\sqrt{\frac{1}{3}}$;\ \ 
$-\frac{1}{3}c^{1}_{36}
+\frac{1}{3}c^{5}_{36}
$,
$-\frac{1}{3}c^{5}_{36}
$,
$\frac{1}{3}c^{1}_{36}
$;\ \ 
$\frac{1}{3}c^{1}_{36}
$,
$-\frac{1}{3}c^{1}_{36}
+\frac{1}{3}c^{5}_{36}
$;\ \ 
$-\frac{1}{3}c^{5}_{36}
$)
 $\oplus$
($-\frac{1}{2}$,
$-\sqrt{\frac{3}{4}}$;
$\frac{1}{2}$)

Fail:
Tr$_I(C) = -1 <$  0 for I = [ 0 ]. Prop. B.4 (1) eqn. (B.18)

 \ \color{black}

\noindent 423: (dims,levels) = $(4 , 
2;18,
2
)$,
irreps = $4_{9,2}^{1,0}
\hskip -1.5pt \otimes \hskip -1.5pt
1_{2}^{1,0}\oplus
2_{2}^{1,0}$,
pord$(\rho_\text{isum}(\mathfrak{t})) = 18$,

\vskip 0.7ex
\hangindent=5.5em \hangafter=1
{\white .}\hskip 1em $\rho_\text{isum}(\mathfrak{t})$ =
 $( \frac{1}{2},
\frac{5}{18},
\frac{11}{18},
\frac{17}{18} )
\oplus
( 0,
\frac{1}{2} )
$,

\vskip 0.7ex
\hangindent=5.5em \hangafter=1
{\white .}\hskip 1em $\rho_\text{isum}(\mathfrak{s})$ =
($0$,
$-\sqrt{\frac{1}{3}}$,
$-\sqrt{\frac{1}{3}}$,
$-\sqrt{\frac{1}{3}}$;
$-\frac{1}{3} c_9^4 $,
$-\frac{1}{3}c^{1}_{9}
$,
$-\frac{1}{3}c^{2}_{9}
$;
$-\frac{1}{3}c^{2}_{9}
$,
$-\frac{1}{3} c_9^4 $;
$-\frac{1}{3}c^{1}_{9}
$)
 $\oplus$
($-\frac{1}{2}$,
$-\sqrt{\frac{3}{4}}$;
$\frac{1}{2}$)

Fail:
Integral: $D_{\rho}(\sigma)_{\theta} \propto $ id,
 for all $\sigma$ and all $\theta$-eigenspaces that can contain unit. Prop. B.5 (6)

 \ \color{black}

\noindent 424: (dims,levels) = $(4 , 
2;18,
6
)$,
irreps = $4_{9,1}^{2,0}
\hskip -1.5pt \otimes \hskip -1.5pt
1_{2}^{1,0}\oplus
2_{3}^{1,8}
\hskip -1.5pt \otimes \hskip -1.5pt
1_{2}^{1,0}$,
pord$(\rho_\text{isum}(\mathfrak{t})) = 9$,

\vskip 0.7ex
\hangindent=5.5em \hangafter=1
{\white .}\hskip 1em $\rho_\text{isum}(\mathfrak{t})$ =
 $( \frac{1}{2},
\frac{1}{18},
\frac{7}{18},
\frac{13}{18} )
\oplus
( \frac{1}{2},
\frac{1}{6} )
$,

\vskip 0.7ex
\hangindent=5.5em \hangafter=1
{\white .}\hskip 1em $\rho_\text{isum}(\mathfrak{s})$ =
$\mathrm{i}$($0$,
$\sqrt{\frac{1}{3}}$,
$\sqrt{\frac{1}{3}}$,
$\sqrt{\frac{1}{3}}$;\ \ 
$\frac{1}{3}c^{5}_{36}
$,
$\frac{1}{3}c^{1}_{36}
-\frac{1}{3}c^{5}_{36}
$,
$-\frac{1}{3}c^{1}_{36}
$;\ \ 
$-\frac{1}{3}c^{1}_{36}
$,
$\frac{1}{3}c^{5}_{36}
$;\ \ 
$\frac{1}{3}c^{1}_{36}
-\frac{1}{3}c^{5}_{36}
$)
 $\oplus$
$\mathrm{i}$($-\sqrt{\frac{1}{3}}$,
$\sqrt{\frac{2}{3}}$;\ \ 
$\sqrt{\frac{1}{3}}$)

Fail:
Integral: $D_{\rho}(\sigma)_{\theta} \propto $ id,
 for all $\sigma$ and all $\theta$-eigenspaces that can contain unit. Prop. B.5 (6)

 \ \color{black}

\noindent 425: (dims,levels) = $(4 , 
2;18,
6
)$,
irreps = $4_{9,2}^{5,0}
\hskip -1.5pt \otimes \hskip -1.5pt
1_{2}^{1,0}\oplus
2_{3}^{1,8}
\hskip -1.5pt \otimes \hskip -1.5pt
1_{2}^{1,0}$,
pord$(\rho_\text{isum}(\mathfrak{t})) = 9$,

\vskip 0.7ex
\hangindent=5.5em \hangafter=1
{\white .}\hskip 1em $\rho_\text{isum}(\mathfrak{t})$ =
 $( \frac{1}{2},
\frac{1}{18},
\frac{7}{18},
\frac{13}{18} )
\oplus
( \frac{1}{2},
\frac{1}{6} )
$,

\vskip 0.7ex
\hangindent=5.5em \hangafter=1
{\white .}\hskip 1em $\rho_\text{isum}(\mathfrak{s})$ =
($0$,
$-\sqrt{\frac{1}{3}}$,
$-\sqrt{\frac{1}{3}}$,
$-\sqrt{\frac{1}{3}}$;
$-\frac{1}{3}c^{1}_{9}
$,
$-\frac{1}{3} c_9^4 $,
$-\frac{1}{3}c^{2}_{9}
$;
$-\frac{1}{3}c^{2}_{9}
$,
$-\frac{1}{3}c^{1}_{9}
$;
$-\frac{1}{3} c_9^4 $)
 $\oplus$
$\mathrm{i}$($-\sqrt{\frac{1}{3}}$,
$\sqrt{\frac{2}{3}}$;\ \ 
$\sqrt{\frac{1}{3}}$)

Fail:
Tr$_I(C) = -1 <$  0 for I = [ 1/6 ]. Prop. B.4 (1) eqn. (B.18)

 \ \color{black}

\noindent 426: (dims,levels) = $(4 , 
2;18,
6
)$,
irreps = $4_{9,1}^{2,0}
\hskip -1.5pt \otimes \hskip -1.5pt
1_{2}^{1,0}\oplus
2_{3}^{1,0}
\hskip -1.5pt \otimes \hskip -1.5pt
1_{2}^{1,0}$,
pord$(\rho_\text{isum}(\mathfrak{t})) = 9$,

\vskip 0.7ex
\hangindent=5.5em \hangafter=1
{\white .}\hskip 1em $\rho_\text{isum}(\mathfrak{t})$ =
 $( \frac{1}{2},
\frac{1}{18},
\frac{7}{18},
\frac{13}{18} )
\oplus
( \frac{1}{2},
\frac{5}{6} )
$,

\vskip 0.7ex
\hangindent=5.5em \hangafter=1
{\white .}\hskip 1em $\rho_\text{isum}(\mathfrak{s})$ =
$\mathrm{i}$($0$,
$\sqrt{\frac{1}{3}}$,
$\sqrt{\frac{1}{3}}$,
$\sqrt{\frac{1}{3}}$;\ \ 
$\frac{1}{3}c^{5}_{36}
$,
$\frac{1}{3}c^{1}_{36}
-\frac{1}{3}c^{5}_{36}
$,
$-\frac{1}{3}c^{1}_{36}
$;\ \ 
$-\frac{1}{3}c^{1}_{36}
$,
$\frac{1}{3}c^{5}_{36}
$;\ \ 
$\frac{1}{3}c^{1}_{36}
-\frac{1}{3}c^{5}_{36}
$)
 $\oplus$
$\mathrm{i}$($\sqrt{\frac{1}{3}}$,
$\sqrt{\frac{2}{3}}$;\ \ 
$-\sqrt{\frac{1}{3}}$)

Fail:
Integral: $D_{\rho}(\sigma)_{\theta} \propto $ id,
 for all $\sigma$ and all $\theta$-eigenspaces that can contain unit. Prop. B.5 (6)

 \ \color{black}

\noindent 427: (dims,levels) = $(4 , 
2;18,
6
)$,
irreps = $4_{9,2}^{5,0}
\hskip -1.5pt \otimes \hskip -1.5pt
1_{2}^{1,0}\oplus
2_{3}^{1,0}
\hskip -1.5pt \otimes \hskip -1.5pt
1_{2}^{1,0}$,
pord$(\rho_\text{isum}(\mathfrak{t})) = 9$,

\vskip 0.7ex
\hangindent=5.5em \hangafter=1
{\white .}\hskip 1em $\rho_\text{isum}(\mathfrak{t})$ =
 $( \frac{1}{2},
\frac{1}{18},
\frac{7}{18},
\frac{13}{18} )
\oplus
( \frac{1}{2},
\frac{5}{6} )
$,

\vskip 0.7ex
\hangindent=5.5em \hangafter=1
{\white .}\hskip 1em $\rho_\text{isum}(\mathfrak{s})$ =
($0$,
$-\sqrt{\frac{1}{3}}$,
$-\sqrt{\frac{1}{3}}$,
$-\sqrt{\frac{1}{3}}$;
$-\frac{1}{3}c^{1}_{9}
$,
$-\frac{1}{3} c_9^4 $,
$-\frac{1}{3}c^{2}_{9}
$;
$-\frac{1}{3}c^{2}_{9}
$,
$-\frac{1}{3}c^{1}_{9}
$;
$-\frac{1}{3} c_9^4 $)
 $\oplus$
$\mathrm{i}$($\sqrt{\frac{1}{3}}$,
$\sqrt{\frac{2}{3}}$;\ \ 
$-\sqrt{\frac{1}{3}}$)

Fail:
Tr$_I(C) = -1 <$  0 for I = [ 5/6 ]. Prop. B.4 (1) eqn. (B.18)

 \ \color{black}

\noindent 428: (dims,levels) = $(4 , 
2;18,
6
)$,
irreps = $4_{9,1}^{1,0}
\hskip -1.5pt \otimes \hskip -1.5pt
1_{2}^{1,0}\oplus
2_{3}^{1,8}
\hskip -1.5pt \otimes \hskip -1.5pt
1_{2}^{1,0}$,
pord$(\rho_\text{isum}(\mathfrak{t})) = 9$,

\vskip 0.7ex
\hangindent=5.5em \hangafter=1
{\white .}\hskip 1em $\rho_\text{isum}(\mathfrak{t})$ =
 $( \frac{1}{2},
\frac{5}{18},
\frac{11}{18},
\frac{17}{18} )
\oplus
( \frac{1}{2},
\frac{1}{6} )
$,

\vskip 0.7ex
\hangindent=5.5em \hangafter=1
{\white .}\hskip 1em $\rho_\text{isum}(\mathfrak{s})$ =
$\mathrm{i}$($0$,
$\sqrt{\frac{1}{3}}$,
$\sqrt{\frac{1}{3}}$,
$\sqrt{\frac{1}{3}}$;\ \ 
$-\frac{1}{3}c^{1}_{36}
+\frac{1}{3}c^{5}_{36}
$,
$-\frac{1}{3}c^{5}_{36}
$,
$\frac{1}{3}c^{1}_{36}
$;\ \ 
$\frac{1}{3}c^{1}_{36}
$,
$-\frac{1}{3}c^{1}_{36}
+\frac{1}{3}c^{5}_{36}
$;\ \ 
$-\frac{1}{3}c^{5}_{36}
$)
 $\oplus$
$\mathrm{i}$($-\sqrt{\frac{1}{3}}$,
$\sqrt{\frac{2}{3}}$;\ \ 
$\sqrt{\frac{1}{3}}$)

Fail:
Integral: $D_{\rho}(\sigma)_{\theta} \propto $ id,
 for all $\sigma$ and all $\theta$-eigenspaces that can contain unit. Prop. B.5 (6)

 \ \color{black}

\noindent 429: (dims,levels) = $(4 , 
2;18,
6
)$,
irreps = $4_{9,2}^{1,0}
\hskip -1.5pt \otimes \hskip -1.5pt
1_{2}^{1,0}\oplus
2_{3}^{1,8}
\hskip -1.5pt \otimes \hskip -1.5pt
1_{2}^{1,0}$,
pord$(\rho_\text{isum}(\mathfrak{t})) = 9$,

\vskip 0.7ex
\hangindent=5.5em \hangafter=1
{\white .}\hskip 1em $\rho_\text{isum}(\mathfrak{t})$ =
 $( \frac{1}{2},
\frac{5}{18},
\frac{11}{18},
\frac{17}{18} )
\oplus
( \frac{1}{2},
\frac{1}{6} )
$,

\vskip 0.7ex
\hangindent=5.5em \hangafter=1
{\white .}\hskip 1em $\rho_\text{isum}(\mathfrak{s})$ =
($0$,
$-\sqrt{\frac{1}{3}}$,
$-\sqrt{\frac{1}{3}}$,
$-\sqrt{\frac{1}{3}}$;
$-\frac{1}{3} c_9^4 $,
$-\frac{1}{3}c^{1}_{9}
$,
$-\frac{1}{3}c^{2}_{9}
$;
$-\frac{1}{3}c^{2}_{9}
$,
$-\frac{1}{3} c_9^4 $;
$-\frac{1}{3}c^{1}_{9}
$)
 $\oplus$
$\mathrm{i}$($-\sqrt{\frac{1}{3}}$,
$\sqrt{\frac{2}{3}}$;\ \ 
$\sqrt{\frac{1}{3}}$)

Fail:
Tr$_I(C) = -1 <$  0 for I = [ 1/6 ]. Prop. B.4 (1) eqn. (B.18)

 \ \color{black}

\noindent 430: (dims,levels) = $(4 , 
2;18,
6
)$,
irreps = $4_{9,1}^{1,0}
\hskip -1.5pt \otimes \hskip -1.5pt
1_{2}^{1,0}\oplus
2_{3}^{1,0}
\hskip -1.5pt \otimes \hskip -1.5pt
1_{2}^{1,0}$,
pord$(\rho_\text{isum}(\mathfrak{t})) = 9$,

\vskip 0.7ex
\hangindent=5.5em \hangafter=1
{\white .}\hskip 1em $\rho_\text{isum}(\mathfrak{t})$ =
 $( \frac{1}{2},
\frac{5}{18},
\frac{11}{18},
\frac{17}{18} )
\oplus
( \frac{1}{2},
\frac{5}{6} )
$,

\vskip 0.7ex
\hangindent=5.5em \hangafter=1
{\white .}\hskip 1em $\rho_\text{isum}(\mathfrak{s})$ =
$\mathrm{i}$($0$,
$\sqrt{\frac{1}{3}}$,
$\sqrt{\frac{1}{3}}$,
$\sqrt{\frac{1}{3}}$;\ \ 
$-\frac{1}{3}c^{1}_{36}
+\frac{1}{3}c^{5}_{36}
$,
$-\frac{1}{3}c^{5}_{36}
$,
$\frac{1}{3}c^{1}_{36}
$;\ \ 
$\frac{1}{3}c^{1}_{36}
$,
$-\frac{1}{3}c^{1}_{36}
+\frac{1}{3}c^{5}_{36}
$;\ \ 
$-\frac{1}{3}c^{5}_{36}
$)
 $\oplus$
$\mathrm{i}$($\sqrt{\frac{1}{3}}$,
$\sqrt{\frac{2}{3}}$;\ \ 
$-\sqrt{\frac{1}{3}}$)

Fail:
Integral: $D_{\rho}(\sigma)_{\theta} \propto $ id,
 for all $\sigma$ and all $\theta$-eigenspaces that can contain unit. Prop. B.5 (6)

 \ \color{black}

\noindent 431: (dims,levels) = $(4 , 
2;18,
6
)$,
irreps = $4_{9,2}^{1,0}
\hskip -1.5pt \otimes \hskip -1.5pt
1_{2}^{1,0}\oplus
2_{3}^{1,0}
\hskip -1.5pt \otimes \hskip -1.5pt
1_{2}^{1,0}$,
pord$(\rho_\text{isum}(\mathfrak{t})) = 9$,

\vskip 0.7ex
\hangindent=5.5em \hangafter=1
{\white .}\hskip 1em $\rho_\text{isum}(\mathfrak{t})$ =
 $( \frac{1}{2},
\frac{5}{18},
\frac{11}{18},
\frac{17}{18} )
\oplus
( \frac{1}{2},
\frac{5}{6} )
$,

\vskip 0.7ex
\hangindent=5.5em \hangafter=1
{\white .}\hskip 1em $\rho_\text{isum}(\mathfrak{s})$ =
($0$,
$-\sqrt{\frac{1}{3}}$,
$-\sqrt{\frac{1}{3}}$,
$-\sqrt{\frac{1}{3}}$;
$-\frac{1}{3} c_9^4 $,
$-\frac{1}{3}c^{1}_{9}
$,
$-\frac{1}{3}c^{2}_{9}
$;
$-\frac{1}{3}c^{2}_{9}
$,
$-\frac{1}{3} c_9^4 $;
$-\frac{1}{3}c^{1}_{9}
$)
 $\oplus$
$\mathrm{i}$($\sqrt{\frac{1}{3}}$,
$\sqrt{\frac{2}{3}}$;\ \ 
$-\sqrt{\frac{1}{3}}$)

Fail:
Tr$_I(C) = -1 <$  0 for I = [ 5/6 ]. Prop. B.4 (1) eqn. (B.18)

 \ \color{black}

 \color{blue}

\noindent 432: (dims,levels) = $(4 , 
2;20,
20
)$,
irreps = $2_{5}^{1}
\hskip -1.5pt \otimes \hskip -1.5pt
2_{4}^{1,0}\oplus
2_{5}^{1}
\hskip -1.5pt \otimes \hskip -1.5pt
1_{4}^{1,0}$,
pord$(\rho_\text{isum}(\mathfrak{t})) = 10$,

\vskip 0.7ex
\hangindent=5.5em \hangafter=1
{\white .}\hskip 1em $\rho_\text{isum}(\mathfrak{t})$ =
 $( \frac{1}{20},
\frac{9}{20},
\frac{11}{20},
\frac{19}{20} )
\oplus
( \frac{1}{20},
\frac{9}{20} )
$,

\vskip 0.7ex
\hangindent=5.5em \hangafter=1
{\white .}\hskip 1em $\rho_\text{isum}(\mathfrak{s})$ =
($\frac{1}{2\sqrt{5}}c^{3}_{20}
$,
$\frac{1}{2\sqrt{5}}c^{1}_{20}
$,
$-\frac{3}{2\sqrt{15}}c^{3}_{20}
$,
$\frac{3}{2\sqrt{15}}c^{1}_{20}
$;
$-\frac{1}{2\sqrt{5}}c^{3}_{20}
$,
$-\frac{3}{2\sqrt{15}}c^{1}_{20}
$,
$-\frac{3}{2\sqrt{15}}c^{3}_{20}
$;
$-\frac{1}{2\sqrt{5}}c^{3}_{20}
$,
$\frac{1}{2\sqrt{5}}c^{1}_{20}
$;
$\frac{1}{2\sqrt{5}}c^{3}_{20}
$)
 $\oplus$
($-\frac{1}{\sqrt{5}}c^{3}_{20}
$,
$\frac{1}{\sqrt{5}}c^{1}_{20}
$;
$\frac{1}{\sqrt{5}}c^{3}_{20}
$)

Pass. 

 \ \color{black}

 \color{blue}

\noindent 433: (dims,levels) = $(4 , 
2;20,
20
)$,
irreps = $2_{5}^{1}
\hskip -1.5pt \otimes \hskip -1.5pt
2_{4}^{1,0}\oplus
2_{5}^{1}
\hskip -1.5pt \otimes \hskip -1.5pt
1_{4}^{1,6}$,
pord$(\rho_\text{isum}(\mathfrak{t})) = 10$,

\vskip 0.7ex
\hangindent=5.5em \hangafter=1
{\white .}\hskip 1em $\rho_\text{isum}(\mathfrak{t})$ =
 $( \frac{1}{20},
\frac{9}{20},
\frac{11}{20},
\frac{19}{20} )
\oplus
( \frac{11}{20},
\frac{19}{20} )
$,

\vskip 0.7ex
\hangindent=5.5em \hangafter=1
{\white .}\hskip 1em $\rho_\text{isum}(\mathfrak{s})$ =
($\frac{1}{2\sqrt{5}}c^{3}_{20}
$,
$\frac{1}{2\sqrt{5}}c^{1}_{20}
$,
$-\frac{3}{2\sqrt{15}}c^{3}_{20}
$,
$\frac{3}{2\sqrt{15}}c^{1}_{20}
$;
$-\frac{1}{2\sqrt{5}}c^{3}_{20}
$,
$-\frac{3}{2\sqrt{15}}c^{1}_{20}
$,
$-\frac{3}{2\sqrt{15}}c^{3}_{20}
$;
$-\frac{1}{2\sqrt{5}}c^{3}_{20}
$,
$\frac{1}{2\sqrt{5}}c^{1}_{20}
$;
$\frac{1}{2\sqrt{5}}c^{3}_{20}
$)
 $\oplus$
($\frac{1}{\sqrt{5}}c^{3}_{20}
$,
$\frac{1}{\sqrt{5}}c^{1}_{20}
$;
$-\frac{1}{\sqrt{5}}c^{3}_{20}
$)

Pass. 

 \ \color{black}

 \color{blue}

\noindent 434: (dims,levels) = $(4 , 
2;20,
20
)$,
irreps = $4_{5,1}^{1}
\hskip -1.5pt \otimes \hskip -1.5pt
1_{4}^{1,0}\oplus
2_{5}^{1}
\hskip -1.5pt \otimes \hskip -1.5pt
1_{4}^{1,0}$,
pord$(\rho_\text{isum}(\mathfrak{t})) = 5$,

\vskip 0.7ex
\hangindent=5.5em \hangafter=1
{\white .}\hskip 1em $\rho_\text{isum}(\mathfrak{t})$ =
 $( \frac{1}{20},
\frac{9}{20},
\frac{13}{20},
\frac{17}{20} )
\oplus
( \frac{1}{20},
\frac{9}{20} )
$,

\vskip 0.7ex
\hangindent=5.5em \hangafter=1
{\white .}\hskip 1em $\rho_\text{isum}(\mathfrak{s})$ =
($\frac{1}{5}c^{1}_{20}
+\frac{1}{5}c^{3}_{20}
$,
$\frac{1}{5}c^{1}_{20}
-\frac{1}{5}c^{3}_{20}
$,
$-\frac{1}{5}+\frac{2}{5}c^{1}_{15}
-\frac{1}{5}c^{3}_{15}
$,
$\frac{2}{5}c^{2}_{15}
+\frac{1}{5}c^{3}_{15}
$;
$-\frac{1}{5}c^{1}_{20}
-\frac{1}{5}c^{3}_{20}
$,
$\frac{2}{5}c^{2}_{15}
+\frac{1}{5}c^{3}_{15}
$,
$\frac{1}{5}-\frac{2}{5}c^{1}_{15}
+\frac{1}{5}c^{3}_{15}
$;
$\frac{1}{5}c^{1}_{20}
-\frac{1}{5}c^{3}_{20}
$,
$-\frac{1}{5}c^{1}_{20}
-\frac{1}{5}c^{3}_{20}
$;
$-\frac{1}{5}c^{1}_{20}
+\frac{1}{5}c^{3}_{20}
$)
 $\oplus$
($-\frac{1}{\sqrt{5}}c^{3}_{20}
$,
$\frac{1}{\sqrt{5}}c^{1}_{20}
$;
$\frac{1}{\sqrt{5}}c^{3}_{20}
$)

Pass. 

 \ \color{black}

 \color{blue}

\noindent 435: (dims,levels) = $(4 , 
2;20,
20
)$,
irreps = $4_{5,2}^{1}
\hskip -1.5pt \otimes \hskip -1.5pt
1_{4}^{1,0}\oplus
2_{5}^{1}
\hskip -1.5pt \otimes \hskip -1.5pt
1_{4}^{1,0}$,
pord$(\rho_\text{isum}(\mathfrak{t})) = 5$,

\vskip 0.7ex
\hangindent=5.5em \hangafter=1
{\white .}\hskip 1em $\rho_\text{isum}(\mathfrak{t})$ =
 $( \frac{1}{20},
\frac{9}{20},
\frac{13}{20},
\frac{17}{20} )
\oplus
( \frac{1}{20},
\frac{9}{20} )
$,

\vskip 0.7ex
\hangindent=5.5em \hangafter=1
{\white .}\hskip 1em $\rho_\text{isum}(\mathfrak{s})$ =
$\mathrm{i}$($\sqrt{\frac{1}{5}}$,
$\sqrt{\frac{1}{5}}$,
$-\frac{5+\sqrt{5}}{10}$,
$\frac{-5+\sqrt{5}}{10}$;\ \ 
$\sqrt{\frac{1}{5}}$,
$\frac{5-\sqrt{5}}{10}$,
$\frac{5+\sqrt{5}}{10}$;\ \ 
$-\sqrt{\frac{1}{5}}$,
$\sqrt{\frac{1}{5}}$;\ \ 
$-\sqrt{\frac{1}{5}}$)
 $\oplus$
($-\frac{1}{\sqrt{5}}c^{3}_{20}
$,
$\frac{1}{\sqrt{5}}c^{1}_{20}
$;
$\frac{1}{\sqrt{5}}c^{3}_{20}
$)

Pass. 

 \ \color{black}

 \color{blue}

\noindent 436: (dims,levels) = $(4 , 
2;20,
20
)$,
irreps = $4_{5,1}^{1}
\hskip -1.5pt \otimes \hskip -1.5pt
1_{4}^{1,0}\oplus
2_{5}^{2}
\hskip -1.5pt \otimes \hskip -1.5pt
1_{4}^{1,0}$,
pord$(\rho_\text{isum}(\mathfrak{t})) = 5$,

\vskip 0.7ex
\hangindent=5.5em \hangafter=1
{\white .}\hskip 1em $\rho_\text{isum}(\mathfrak{t})$ =
 $( \frac{1}{20},
\frac{9}{20},
\frac{13}{20},
\frac{17}{20} )
\oplus
( \frac{13}{20},
\frac{17}{20} )
$,

\vskip 0.7ex
\hangindent=5.5em \hangafter=1
{\white .}\hskip 1em $\rho_\text{isum}(\mathfrak{s})$ =
($\frac{1}{5}c^{1}_{20}
+\frac{1}{5}c^{3}_{20}
$,
$\frac{1}{5}c^{1}_{20}
-\frac{1}{5}c^{3}_{20}
$,
$-\frac{1}{5}+\frac{2}{5}c^{1}_{15}
-\frac{1}{5}c^{3}_{15}
$,
$\frac{2}{5}c^{2}_{15}
+\frac{1}{5}c^{3}_{15}
$;
$-\frac{1}{5}c^{1}_{20}
-\frac{1}{5}c^{3}_{20}
$,
$\frac{2}{5}c^{2}_{15}
+\frac{1}{5}c^{3}_{15}
$,
$\frac{1}{5}-\frac{2}{5}c^{1}_{15}
+\frac{1}{5}c^{3}_{15}
$;
$\frac{1}{5}c^{1}_{20}
-\frac{1}{5}c^{3}_{20}
$,
$-\frac{1}{5}c^{1}_{20}
-\frac{1}{5}c^{3}_{20}
$;
$-\frac{1}{5}c^{1}_{20}
+\frac{1}{5}c^{3}_{20}
$)
 $\oplus$
($\frac{1}{\sqrt{5}}c^{1}_{20}
$,
$\frac{1}{\sqrt{5}}c^{3}_{20}
$;
$-\frac{1}{\sqrt{5}}c^{1}_{20}
$)

Pass. 

 \ \color{black}

 \color{blue}

\noindent 437: (dims,levels) = $(4 , 
2;20,
20
)$,
irreps = $4_{5,2}^{1}
\hskip -1.5pt \otimes \hskip -1.5pt
1_{4}^{1,0}\oplus
2_{5}^{2}
\hskip -1.5pt \otimes \hskip -1.5pt
1_{4}^{1,0}$,
pord$(\rho_\text{isum}(\mathfrak{t})) = 5$,

\vskip 0.7ex
\hangindent=5.5em \hangafter=1
{\white .}\hskip 1em $\rho_\text{isum}(\mathfrak{t})$ =
 $( \frac{1}{20},
\frac{9}{20},
\frac{13}{20},
\frac{17}{20} )
\oplus
( \frac{13}{20},
\frac{17}{20} )
$,

\vskip 0.7ex
\hangindent=5.5em \hangafter=1
{\white .}\hskip 1em $\rho_\text{isum}(\mathfrak{s})$ =
$\mathrm{i}$($\sqrt{\frac{1}{5}}$,
$\sqrt{\frac{1}{5}}$,
$-\frac{5+\sqrt{5}}{10}$,
$\frac{-5+\sqrt{5}}{10}$;\ \ 
$\sqrt{\frac{1}{5}}$,
$\frac{5-\sqrt{5}}{10}$,
$\frac{5+\sqrt{5}}{10}$;\ \ 
$-\sqrt{\frac{1}{5}}$,
$\sqrt{\frac{1}{5}}$;\ \ 
$-\sqrt{\frac{1}{5}}$)
 $\oplus$
($\frac{1}{\sqrt{5}}c^{1}_{20}
$,
$\frac{1}{\sqrt{5}}c^{3}_{20}
$;
$-\frac{1}{\sqrt{5}}c^{1}_{20}
$)

Pass. 

 \ \color{black}

 \color{blue}

\noindent 438: (dims,levels) = $(4 , 
2;20,
20
)$,
irreps = $2_{5}^{2}
\hskip -1.5pt \otimes \hskip -1.5pt
2_{4}^{1,0}\oplus
2_{5}^{2}
\hskip -1.5pt \otimes \hskip -1.5pt
1_{4}^{1,6}$,
pord$(\rho_\text{isum}(\mathfrak{t})) = 10$,

\vskip 0.7ex
\hangindent=5.5em \hangafter=1
{\white .}\hskip 1em $\rho_\text{isum}(\mathfrak{t})$ =
 $( \frac{3}{20},
\frac{7}{20},
\frac{13}{20},
\frac{17}{20} )
\oplus
( \frac{3}{20},
\frac{7}{20} )
$,

\vskip 0.7ex
\hangindent=5.5em \hangafter=1
{\white .}\hskip 1em $\rho_\text{isum}(\mathfrak{s})$ =
($\frac{1}{2\sqrt{5}}c^{1}_{20}
$,
$\frac{1}{2\sqrt{5}}c^{3}_{20}
$,
$\frac{3}{2\sqrt{15}}c^{1}_{20}
$,
$-\frac{3}{2\sqrt{15}}c^{3}_{20}
$;
$-\frac{1}{2\sqrt{5}}c^{1}_{20}
$,
$\frac{3}{2\sqrt{15}}c^{3}_{20}
$,
$\frac{3}{2\sqrt{15}}c^{1}_{20}
$;
$-\frac{1}{2\sqrt{5}}c^{1}_{20}
$,
$\frac{1}{2\sqrt{5}}c^{3}_{20}
$;
$\frac{1}{2\sqrt{5}}c^{1}_{20}
$)
 $\oplus$
($-\frac{1}{\sqrt{5}}c^{1}_{20}
$,
$\frac{1}{\sqrt{5}}c^{3}_{20}
$;
$\frac{1}{\sqrt{5}}c^{1}_{20}
$)

Pass. 

 \ \color{black}

 \color{blue}

\noindent 439: (dims,levels) = $(4 , 
2;20,
20
)$,
irreps = $2_{5}^{2}
\hskip -1.5pt \otimes \hskip -1.5pt
2_{4}^{1,0}\oplus
2_{5}^{2}
\hskip -1.5pt \otimes \hskip -1.5pt
1_{4}^{1,0}$,
pord$(\rho_\text{isum}(\mathfrak{t})) = 10$,

\vskip 0.7ex
\hangindent=5.5em \hangafter=1
{\white .}\hskip 1em $\rho_\text{isum}(\mathfrak{t})$ =
 $( \frac{3}{20},
\frac{7}{20},
\frac{13}{20},
\frac{17}{20} )
\oplus
( \frac{13}{20},
\frac{17}{20} )
$,

\vskip 0.7ex
\hangindent=5.5em \hangafter=1
{\white .}\hskip 1em $\rho_\text{isum}(\mathfrak{s})$ =
($\frac{1}{2\sqrt{5}}c^{1}_{20}
$,
$\frac{1}{2\sqrt{5}}c^{3}_{20}
$,
$\frac{3}{2\sqrt{15}}c^{1}_{20}
$,
$-\frac{3}{2\sqrt{15}}c^{3}_{20}
$;
$-\frac{1}{2\sqrt{5}}c^{1}_{20}
$,
$\frac{3}{2\sqrt{15}}c^{3}_{20}
$,
$\frac{3}{2\sqrt{15}}c^{1}_{20}
$;
$-\frac{1}{2\sqrt{5}}c^{1}_{20}
$,
$\frac{1}{2\sqrt{5}}c^{3}_{20}
$;
$\frac{1}{2\sqrt{5}}c^{1}_{20}
$)
 $\oplus$
($\frac{1}{\sqrt{5}}c^{1}_{20}
$,
$\frac{1}{\sqrt{5}}c^{3}_{20}
$;
$-\frac{1}{\sqrt{5}}c^{1}_{20}
$)

Pass. 

 \ \color{black}

 \color{blue}

\noindent 440: (dims,levels) = $(4 , 
2;21,
3
)$,
irreps = $4_{7}^{3}
\hskip -1.5pt \otimes \hskip -1.5pt
1_{3}^{1,0}\oplus
2_{3}^{1,0}$,
pord$(\rho_\text{isum}(\mathfrak{t})) = 21$,

\vskip 0.7ex
\hangindent=5.5em \hangafter=1
{\white .}\hskip 1em $\rho_\text{isum}(\mathfrak{t})$ =
 $( \frac{1}{3},
\frac{1}{21},
\frac{4}{21},
\frac{16}{21} )
\oplus
( 0,
\frac{1}{3} )
$,

\vskip 0.7ex
\hangindent=5.5em \hangafter=1
{\white .}\hskip 1em $\rho_\text{isum}(\mathfrak{s})$ =
$\mathrm{i}$($\sqrt{\frac{1}{7}}$,
$\sqrt{\frac{2}{7}}$,
$\sqrt{\frac{2}{7}}$,
$\sqrt{\frac{2}{7}}$;\ \ 
$-\frac{1}{\sqrt{7}\mathrm{i}}s^{5}_{28}
$,
$\frac{1}{\sqrt{7}}c^{1}_{7}
$,
$\frac{1}{\sqrt{7}}c^{2}_{7}
$;\ \ 
$\frac{1}{\sqrt{7}}c^{2}_{7}
$,
$-\frac{1}{\sqrt{7}\mathrm{i}}s^{5}_{28}
$;\ \ 
$\frac{1}{\sqrt{7}}c^{1}_{7}
$)
 $\oplus$
$\mathrm{i}$($-\sqrt{\frac{1}{3}}$,
$\sqrt{\frac{2}{3}}$;\ \ 
$\sqrt{\frac{1}{3}}$)

Pass. 

 \ \color{black}

 \color{blue}

\noindent 441: (dims,levels) = $(4 , 
2;21,
3
)$,
irreps = $4_{7}^{3}
\hskip -1.5pt \otimes \hskip -1.5pt
1_{3}^{1,0}\oplus
2_{3}^{1,4}$,
pord$(\rho_\text{isum}(\mathfrak{t})) = 21$,

\vskip 0.7ex
\hangindent=5.5em \hangafter=1
{\white .}\hskip 1em $\rho_\text{isum}(\mathfrak{t})$ =
 $( \frac{1}{3},
\frac{1}{21},
\frac{4}{21},
\frac{16}{21} )
\oplus
( \frac{1}{3},
\frac{2}{3} )
$,

\vskip 0.7ex
\hangindent=5.5em \hangafter=1
{\white .}\hskip 1em $\rho_\text{isum}(\mathfrak{s})$ =
$\mathrm{i}$($\sqrt{\frac{1}{7}}$,
$\sqrt{\frac{2}{7}}$,
$\sqrt{\frac{2}{7}}$,
$\sqrt{\frac{2}{7}}$;\ \ 
$-\frac{1}{\sqrt{7}\mathrm{i}}s^{5}_{28}
$,
$\frac{1}{\sqrt{7}}c^{1}_{7}
$,
$\frac{1}{\sqrt{7}}c^{2}_{7}
$;\ \ 
$\frac{1}{\sqrt{7}}c^{2}_{7}
$,
$-\frac{1}{\sqrt{7}\mathrm{i}}s^{5}_{28}
$;\ \ 
$\frac{1}{\sqrt{7}}c^{1}_{7}
$)
 $\oplus$
$\mathrm{i}$($-\sqrt{\frac{1}{3}}$,
$\sqrt{\frac{2}{3}}$;\ \ 
$\sqrt{\frac{1}{3}}$)

Pass. 

 \ \color{black}

 \color{blue}

\noindent 442: (dims,levels) = $(4 , 
2;21,
3
)$,
irreps = $4_{7}^{1}
\hskip -1.5pt \otimes \hskip -1.5pt
1_{3}^{1,0}\oplus
2_{3}^{1,0}$,
pord$(\rho_\text{isum}(\mathfrak{t})) = 21$,

\vskip 0.7ex
\hangindent=5.5em \hangafter=1
{\white .}\hskip 1em $\rho_\text{isum}(\mathfrak{t})$ =
 $( \frac{1}{3},
\frac{10}{21},
\frac{13}{21},
\frac{19}{21} )
\oplus
( 0,
\frac{1}{3} )
$,

\vskip 0.7ex
\hangindent=5.5em \hangafter=1
{\white .}\hskip 1em $\rho_\text{isum}(\mathfrak{s})$ =
$\mathrm{i}$($-\sqrt{\frac{1}{7}}$,
$\sqrt{\frac{2}{7}}$,
$\sqrt{\frac{2}{7}}$,
$\sqrt{\frac{2}{7}}$;\ \ 
$-\frac{1}{\sqrt{7}}c^{2}_{7}
$,
$-\frac{1}{\sqrt{7}}c^{1}_{7}
$,
$\frac{1}{\sqrt{7}\mathrm{i}}s^{5}_{28}
$;\ \ 
$\frac{1}{\sqrt{7}\mathrm{i}}s^{5}_{28}
$,
$-\frac{1}{\sqrt{7}}c^{2}_{7}
$;\ \ 
$-\frac{1}{\sqrt{7}}c^{1}_{7}
$)
 $\oplus$
$\mathrm{i}$($-\sqrt{\frac{1}{3}}$,
$\sqrt{\frac{2}{3}}$;\ \ 
$\sqrt{\frac{1}{3}}$)

Pass. 

 \ \color{black}

 \color{blue}

\noindent 443: (dims,levels) = $(4 , 
2;21,
3
)$,
irreps = $4_{7}^{1}
\hskip -1.5pt \otimes \hskip -1.5pt
1_{3}^{1,0}\oplus
2_{3}^{1,4}$,
pord$(\rho_\text{isum}(\mathfrak{t})) = 21$,

\vskip 0.7ex
\hangindent=5.5em \hangafter=1
{\white .}\hskip 1em $\rho_\text{isum}(\mathfrak{t})$ =
 $( \frac{1}{3},
\frac{10}{21},
\frac{13}{21},
\frac{19}{21} )
\oplus
( \frac{1}{3},
\frac{2}{3} )
$,

\vskip 0.7ex
\hangindent=5.5em \hangafter=1
{\white .}\hskip 1em $\rho_\text{isum}(\mathfrak{s})$ =
$\mathrm{i}$($-\sqrt{\frac{1}{7}}$,
$\sqrt{\frac{2}{7}}$,
$\sqrt{\frac{2}{7}}$,
$\sqrt{\frac{2}{7}}$;\ \ 
$-\frac{1}{\sqrt{7}}c^{2}_{7}
$,
$-\frac{1}{\sqrt{7}}c^{1}_{7}
$,
$\frac{1}{\sqrt{7}\mathrm{i}}s^{5}_{28}
$;\ \ 
$\frac{1}{\sqrt{7}\mathrm{i}}s^{5}_{28}
$,
$-\frac{1}{\sqrt{7}}c^{2}_{7}
$;\ \ 
$-\frac{1}{\sqrt{7}}c^{1}_{7}
$)
 $\oplus$
$\mathrm{i}$($-\sqrt{\frac{1}{3}}$,
$\sqrt{\frac{2}{3}}$;\ \ 
$\sqrt{\frac{1}{3}}$)

Pass. 

 \ \color{black}

\noindent 444: (dims,levels) = $(4 , 
2;21,
6
)$,
irreps = $4_{7}^{3}
\hskip -1.5pt \otimes \hskip -1.5pt
1_{3}^{1,0}\oplus
2_{2}^{1,0}
\hskip -1.5pt \otimes \hskip -1.5pt
1_{3}^{1,0}$,
pord$(\rho_\text{isum}(\mathfrak{t})) = 14$,

\vskip 0.7ex
\hangindent=5.5em \hangafter=1
{\white .}\hskip 1em $\rho_\text{isum}(\mathfrak{t})$ =
 $( \frac{1}{3},
\frac{1}{21},
\frac{4}{21},
\frac{16}{21} )
\oplus
( \frac{1}{3},
\frac{5}{6} )
$,

\vskip 0.7ex
\hangindent=5.5em \hangafter=1
{\white .}\hskip 1em $\rho_\text{isum}(\mathfrak{s})$ =
$\mathrm{i}$($\sqrt{\frac{1}{7}}$,
$\sqrt{\frac{2}{7}}$,
$\sqrt{\frac{2}{7}}$,
$\sqrt{\frac{2}{7}}$;\ \ 
$-\frac{1}{\sqrt{7}\mathrm{i}}s^{5}_{28}
$,
$\frac{1}{\sqrt{7}}c^{1}_{7}
$,
$\frac{1}{\sqrt{7}}c^{2}_{7}
$;\ \ 
$\frac{1}{\sqrt{7}}c^{2}_{7}
$,
$-\frac{1}{\sqrt{7}\mathrm{i}}s^{5}_{28}
$;\ \ 
$\frac{1}{\sqrt{7}}c^{1}_{7}
$)
 $\oplus$
($-\frac{1}{2}$,
$-\sqrt{\frac{3}{4}}$;
$\frac{1}{2}$)

Fail:
Tr$_I(C) = -1 <$  0 for I = [ 5/6 ]. Prop. B.4 (1) eqn. (B.18)

 \ \color{black}

\noindent 445: (dims,levels) = $(4 , 
2;21,
6
)$,
irreps = $4_{7}^{1}
\hskip -1.5pt \otimes \hskip -1.5pt
1_{3}^{1,0}\oplus
2_{2}^{1,0}
\hskip -1.5pt \otimes \hskip -1.5pt
1_{3}^{1,0}$,
pord$(\rho_\text{isum}(\mathfrak{t})) = 14$,

\vskip 0.7ex
\hangindent=5.5em \hangafter=1
{\white .}\hskip 1em $\rho_\text{isum}(\mathfrak{t})$ =
 $( \frac{1}{3},
\frac{10}{21},
\frac{13}{21},
\frac{19}{21} )
\oplus
( \frac{1}{3},
\frac{5}{6} )
$,

\vskip 0.7ex
\hangindent=5.5em \hangafter=1
{\white .}\hskip 1em $\rho_\text{isum}(\mathfrak{s})$ =
$\mathrm{i}$($-\sqrt{\frac{1}{7}}$,
$\sqrt{\frac{2}{7}}$,
$\sqrt{\frac{2}{7}}$,
$\sqrt{\frac{2}{7}}$;\ \ 
$-\frac{1}{\sqrt{7}}c^{2}_{7}
$,
$-\frac{1}{\sqrt{7}}c^{1}_{7}
$,
$\frac{1}{\sqrt{7}\mathrm{i}}s^{5}_{28}
$;\ \ 
$\frac{1}{\sqrt{7}\mathrm{i}}s^{5}_{28}
$,
$-\frac{1}{\sqrt{7}}c^{2}_{7}
$;\ \ 
$-\frac{1}{\sqrt{7}}c^{1}_{7}
$)
 $\oplus$
($-\frac{1}{2}$,
$-\sqrt{\frac{3}{4}}$;
$\frac{1}{2}$)

Fail:
Tr$_I(C) = -1 <$  0 for I = [ 5/6 ]. Prop. B.4 (1) eqn. (B.18)

 \ \color{black}

\noindent 446: (dims,levels) = $(4 , 
2;24,
8
)$,
irreps = $2_{8}^{1,0}
\hskip -1.5pt \otimes \hskip -1.5pt
2_{3}^{1,0}\oplus
2_{8}^{1,0}$,
pord$(\rho_\text{isum}(\mathfrak{t})) = 12$,

\vskip 0.7ex
\hangindent=5.5em \hangafter=1
{\white .}\hskip 1em $\rho_\text{isum}(\mathfrak{t})$ =
 $( \frac{1}{8},
\frac{3}{8},
\frac{11}{24},
\frac{17}{24} )
\oplus
( \frac{1}{8},
\frac{3}{8} )
$,

\vskip 0.7ex
\hangindent=5.5em \hangafter=1
{\white .}\hskip 1em $\rho_\text{isum}(\mathfrak{s})$ =
$\mathrm{i}$($\sqrt{\frac{1}{6}}$,
$\sqrt{\frac{1}{6}}$,
$\sqrt{\frac{1}{3}}$,
$\sqrt{\frac{1}{3}}$;\ \ 
$-\sqrt{\frac{1}{6}}$,
$\sqrt{\frac{1}{3}}$,
$-\sqrt{\frac{1}{3}}$;\ \ 
$-\sqrt{\frac{1}{6}}$,
$-\sqrt{\frac{1}{6}}$;\ \ 
$\sqrt{\frac{1}{6}}$)
 $\oplus$
($-\sqrt{\frac{1}{2}}$,
$\sqrt{\frac{1}{2}}$;
$\sqrt{\frac{1}{2}}$)

Fail:
Integral: $D_{\rho}(\sigma)_{\theta} \propto $ id,
 for all $\sigma$ and all $\theta$-eigenspaces that can contain unit. Prop. B.5 (6)

 \ \color{black}

\noindent 447: (dims,levels) = $(4 , 
2;24,
8
)$,
irreps = $2_{8}^{1,0}
\hskip -1.5pt \otimes \hskip -1.5pt
2_{3}^{1,0}\oplus
2_{8}^{1,9}$,
pord$(\rho_\text{isum}(\mathfrak{t})) = 12$,

\vskip 0.7ex
\hangindent=5.5em \hangafter=1
{\white .}\hskip 1em $\rho_\text{isum}(\mathfrak{t})$ =
 $( \frac{1}{8},
\frac{3}{8},
\frac{11}{24},
\frac{17}{24} )
\oplus
( \frac{1}{8},
\frac{7}{8} )
$,

\vskip 0.7ex
\hangindent=5.5em \hangafter=1
{\white .}\hskip 1em $\rho_\text{isum}(\mathfrak{s})$ =
$\mathrm{i}$($\sqrt{\frac{1}{6}}$,
$\sqrt{\frac{1}{6}}$,
$\sqrt{\frac{1}{3}}$,
$\sqrt{\frac{1}{3}}$;\ \ 
$-\sqrt{\frac{1}{6}}$,
$\sqrt{\frac{1}{3}}$,
$-\sqrt{\frac{1}{3}}$;\ \ 
$-\sqrt{\frac{1}{6}}$,
$-\sqrt{\frac{1}{6}}$;\ \ 
$\sqrt{\frac{1}{6}}$)
 $\oplus$
$\mathrm{i}$($-\sqrt{\frac{1}{2}}$,
$\sqrt{\frac{1}{2}}$;\ \ 
$\sqrt{\frac{1}{2}}$)

Fail:
$\sigma(\rho(\mathfrak s)_\mathrm{ndeg}) \neq
 (\rho(\mathfrak t)^a \rho(\mathfrak s) \rho(\mathfrak t)^b
 \rho(\mathfrak s) \rho(\mathfrak t)^a)_\mathrm{ndeg}$,
 $\sigma = a$ = 5. Prop. B.5 (3) eqn. (B.25)

 \ \color{black}

\noindent 448: (dims,levels) = $(4 , 
2;24,
8
)$,
irreps = $2_{8}^{1,0}
\hskip -1.5pt \otimes \hskip -1.5pt
2_{3}^{1,0}\oplus
2_{8}^{1,3}$,
pord$(\rho_\text{isum}(\mathfrak{t})) = 12$,

\vskip 0.7ex
\hangindent=5.5em \hangafter=1
{\white .}\hskip 1em $\rho_\text{isum}(\mathfrak{t})$ =
 $( \frac{1}{8},
\frac{3}{8},
\frac{11}{24},
\frac{17}{24} )
\oplus
( \frac{3}{8},
\frac{5}{8} )
$,

\vskip 0.7ex
\hangindent=5.5em \hangafter=1
{\white .}\hskip 1em $\rho_\text{isum}(\mathfrak{s})$ =
$\mathrm{i}$($\sqrt{\frac{1}{6}}$,
$\sqrt{\frac{1}{6}}$,
$\sqrt{\frac{1}{3}}$,
$\sqrt{\frac{1}{3}}$;\ \ 
$-\sqrt{\frac{1}{6}}$,
$\sqrt{\frac{1}{3}}$,
$-\sqrt{\frac{1}{3}}$;\ \ 
$-\sqrt{\frac{1}{6}}$,
$-\sqrt{\frac{1}{6}}$;\ \ 
$\sqrt{\frac{1}{6}}$)
 $\oplus$
$\mathrm{i}$($-\sqrt{\frac{1}{2}}$,
$\sqrt{\frac{1}{2}}$;\ \ 
$\sqrt{\frac{1}{2}}$)

Fail:
$\sigma(\rho(\mathfrak s)_\mathrm{ndeg}) \neq
 (\rho(\mathfrak t)^a \rho(\mathfrak s) \rho(\mathfrak t)^b
 \rho(\mathfrak s) \rho(\mathfrak t)^a)_\mathrm{ndeg}$,
 $\sigma = a$ = 5. Prop. B.5 (3) eqn. (B.25)

 \ \color{black}

\noindent 449: (dims,levels) = $(4 , 
2;24,
24
)$,
irreps = $2_{8}^{1,0}
\hskip -1.5pt \otimes \hskip -1.5pt
2_{3}^{1,0}\oplus
2_{8}^{1,9}
\hskip -1.5pt \otimes \hskip -1.5pt
1_{3}^{1,0}$,
pord$(\rho_\text{isum}(\mathfrak{t})) = 12$,

\vskip 0.7ex
\hangindent=5.5em \hangafter=1
{\white .}\hskip 1em $\rho_\text{isum}(\mathfrak{t})$ =
 $( \frac{1}{8},
\frac{3}{8},
\frac{11}{24},
\frac{17}{24} )
\oplus
( \frac{5}{24},
\frac{11}{24} )
$,

\vskip 0.7ex
\hangindent=5.5em \hangafter=1
{\white .}\hskip 1em $\rho_\text{isum}(\mathfrak{s})$ =
$\mathrm{i}$($\sqrt{\frac{1}{6}}$,
$\sqrt{\frac{1}{6}}$,
$\sqrt{\frac{1}{3}}$,
$\sqrt{\frac{1}{3}}$;\ \ 
$-\sqrt{\frac{1}{6}}$,
$\sqrt{\frac{1}{3}}$,
$-\sqrt{\frac{1}{3}}$;\ \ 
$-\sqrt{\frac{1}{6}}$,
$-\sqrt{\frac{1}{6}}$;\ \ 
$\sqrt{\frac{1}{6}}$)
 $\oplus$
$\mathrm{i}$($\sqrt{\frac{1}{2}}$,
$\sqrt{\frac{1}{2}}$;\ \ 
$-\sqrt{\frac{1}{2}}$)

Fail:
$\sigma(\rho(\mathfrak s)_\mathrm{ndeg}) \neq
 (\rho(\mathfrak t)^a \rho(\mathfrak s) \rho(\mathfrak t)^b
 \rho(\mathfrak s) \rho(\mathfrak t)^a)_\mathrm{ndeg}$,
 $\sigma = a$ = 5. Prop. B.5 (3) eqn. (B.25)

 \ \color{black}

\noindent 450: (dims,levels) = $(4 , 
2;24,
24
)$,
irreps = $2_{8}^{1,0}
\hskip -1.5pt \otimes \hskip -1.5pt
2_{3}^{1,0}\oplus
2_{8}^{1,0}
\hskip -1.5pt \otimes \hskip -1.5pt
1_{3}^{1,0}$,
pord$(\rho_\text{isum}(\mathfrak{t})) = 12$,

\vskip 0.7ex
\hangindent=5.5em \hangafter=1
{\white .}\hskip 1em $\rho_\text{isum}(\mathfrak{t})$ =
 $( \frac{1}{8},
\frac{3}{8},
\frac{11}{24},
\frac{17}{24} )
\oplus
( \frac{11}{24},
\frac{17}{24} )
$,

\vskip 0.7ex
\hangindent=5.5em \hangafter=1
{\white .}\hskip 1em $\rho_\text{isum}(\mathfrak{s})$ =
$\mathrm{i}$($\sqrt{\frac{1}{6}}$,
$\sqrt{\frac{1}{6}}$,
$\sqrt{\frac{1}{3}}$,
$\sqrt{\frac{1}{3}}$;\ \ 
$-\sqrt{\frac{1}{6}}$,
$\sqrt{\frac{1}{3}}$,
$-\sqrt{\frac{1}{3}}$;\ \ 
$-\sqrt{\frac{1}{6}}$,
$-\sqrt{\frac{1}{6}}$;\ \ 
$\sqrt{\frac{1}{6}}$)
 $\oplus$
($-\sqrt{\frac{1}{2}}$,
$\sqrt{\frac{1}{2}}$;
$\sqrt{\frac{1}{2}}$)

Fail:
Integral: $D_{\rho}(\sigma)_{\theta} \propto $ id,
 for all $\sigma$ and all $\theta$-eigenspaces that can contain unit. Prop. B.5 (6)

 \ \color{black}

\noindent 451: (dims,levels) = $(4 , 
2;24,
24
)$,
irreps = $2_{8}^{1,0}
\hskip -1.5pt \otimes \hskip -1.5pt
2_{3}^{1,0}\oplus
2_{8}^{1,3}
\hskip -1.5pt \otimes \hskip -1.5pt
1_{3}^{1,0}$,
pord$(\rho_\text{isum}(\mathfrak{t})) = 12$,

\vskip 0.7ex
\hangindent=5.5em \hangafter=1
{\white .}\hskip 1em $\rho_\text{isum}(\mathfrak{t})$ =
 $( \frac{1}{8},
\frac{3}{8},
\frac{11}{24},
\frac{17}{24} )
\oplus
( \frac{17}{24},
\frac{23}{24} )
$,

\vskip 0.7ex
\hangindent=5.5em \hangafter=1
{\white .}\hskip 1em $\rho_\text{isum}(\mathfrak{s})$ =
$\mathrm{i}$($\sqrt{\frac{1}{6}}$,
$\sqrt{\frac{1}{6}}$,
$\sqrt{\frac{1}{3}}$,
$\sqrt{\frac{1}{3}}$;\ \ 
$-\sqrt{\frac{1}{6}}$,
$\sqrt{\frac{1}{3}}$,
$-\sqrt{\frac{1}{3}}$;\ \ 
$-\sqrt{\frac{1}{6}}$,
$-\sqrt{\frac{1}{6}}$;\ \ 
$\sqrt{\frac{1}{6}}$)
 $\oplus$
$\mathrm{i}$($-\sqrt{\frac{1}{2}}$,
$\sqrt{\frac{1}{2}}$;\ \ 
$\sqrt{\frac{1}{2}}$)

Fail:
$\sigma(\rho(\mathfrak s)_\mathrm{ndeg}) \neq
 (\rho(\mathfrak t)^a \rho(\mathfrak s) \rho(\mathfrak t)^b
 \rho(\mathfrak s) \rho(\mathfrak t)^a)_\mathrm{ndeg}$,
 $\sigma = a$ = 5. Prop. B.5 (3) eqn. (B.25)

 \ \color{black}

\noindent 452: (dims,levels) = $(4 , 
2;28,
4
)$,
irreps = $4_{7}^{3}
\hskip -1.5pt \otimes \hskip -1.5pt
1_{4}^{1,0}\oplus
2_{4}^{1,0}$,
pord$(\rho_\text{isum}(\mathfrak{t})) = 14$,

\vskip 0.7ex
\hangindent=5.5em \hangafter=1
{\white .}\hskip 1em $\rho_\text{isum}(\mathfrak{t})$ =
 $( \frac{1}{4},
\frac{3}{28},
\frac{19}{28},
\frac{27}{28} )
\oplus
( \frac{1}{4},
\frac{3}{4} )
$,

\vskip 0.7ex
\hangindent=5.5em \hangafter=1
{\white .}\hskip 1em $\rho_\text{isum}(\mathfrak{s})$ =
($-\sqrt{\frac{1}{7}}$,
$\sqrt{\frac{2}{7}}$,
$\sqrt{\frac{2}{7}}$,
$\sqrt{\frac{2}{7}}$;
$-\frac{1}{\sqrt{7}}c^{2}_{7}
$,
$\frac{1}{\sqrt{7}\mathrm{i}}s^{5}_{28}
$,
$-\frac{1}{\sqrt{7}}c^{1}_{7}
$;
$-\frac{1}{\sqrt{7}}c^{1}_{7}
$,
$-\frac{1}{\sqrt{7}}c^{2}_{7}
$;
$\frac{1}{\sqrt{7}\mathrm{i}}s^{5}_{28}
$)
 $\oplus$
$\mathrm{i}$($-\frac{1}{2}$,
$\sqrt{\frac{3}{4}}$;\ \ 
$\frac{1}{2}$)

Fail:
Tr$_I(C) = -1 <$  0 for I = [ 3/4 ]. Prop. B.4 (1) eqn. (B.18)

 \ \color{black}

\noindent 453: (dims,levels) = $(4 , 
2;28,
4
)$,
irreps = $4_{7}^{1}
\hskip -1.5pt \otimes \hskip -1.5pt
1_{4}^{1,0}\oplus
2_{4}^{1,0}$,
pord$(\rho_\text{isum}(\mathfrak{t})) = 14$,

\vskip 0.7ex
\hangindent=5.5em \hangafter=1
{\white .}\hskip 1em $\rho_\text{isum}(\mathfrak{t})$ =
 $( \frac{1}{4},
\frac{11}{28},
\frac{15}{28},
\frac{23}{28} )
\oplus
( \frac{1}{4},
\frac{3}{4} )
$,

\vskip 0.7ex
\hangindent=5.5em \hangafter=1
{\white .}\hskip 1em $\rho_\text{isum}(\mathfrak{s})$ =
($\sqrt{\frac{1}{7}}$,
$\sqrt{\frac{2}{7}}$,
$\sqrt{\frac{2}{7}}$,
$\sqrt{\frac{2}{7}}$;
$\frac{1}{\sqrt{7}}c^{2}_{7}
$,
$\frac{1}{\sqrt{7}}c^{1}_{7}
$,
$-\frac{1}{\sqrt{7}\mathrm{i}}s^{5}_{28}
$;
$-\frac{1}{\sqrt{7}\mathrm{i}}s^{5}_{28}
$,
$\frac{1}{\sqrt{7}}c^{2}_{7}
$;
$\frac{1}{\sqrt{7}}c^{1}_{7}
$)
 $\oplus$
$\mathrm{i}$($-\frac{1}{2}$,
$\sqrt{\frac{3}{4}}$;\ \ 
$\frac{1}{2}$)

Fail:
Tr$_I(C) = -1 <$  0 for I = [ 3/4 ]. Prop. B.4 (1) eqn. (B.18)

 \ \color{black}

 \color{blue}

\noindent 454: (dims,levels) = $(4 , 
2;28,
12
)$,
irreps = $4_{7}^{3}
\hskip -1.5pt \otimes \hskip -1.5pt
1_{4}^{1,0}\oplus
2_{3}^{1,0}
\hskip -1.5pt \otimes \hskip -1.5pt
1_{4}^{1,0}$,
pord$(\rho_\text{isum}(\mathfrak{t})) = 21$,

\vskip 0.7ex
\hangindent=5.5em \hangafter=1
{\white .}\hskip 1em $\rho_\text{isum}(\mathfrak{t})$ =
 $( \frac{1}{4},
\frac{3}{28},
\frac{19}{28},
\frac{27}{28} )
\oplus
( \frac{1}{4},
\frac{7}{12} )
$,

\vskip 0.7ex
\hangindent=5.5em \hangafter=1
{\white .}\hskip 1em $\rho_\text{isum}(\mathfrak{s})$ =
($-\sqrt{\frac{1}{7}}$,
$\sqrt{\frac{2}{7}}$,
$\sqrt{\frac{2}{7}}$,
$\sqrt{\frac{2}{7}}$;
$-\frac{1}{\sqrt{7}}c^{2}_{7}
$,
$\frac{1}{\sqrt{7}\mathrm{i}}s^{5}_{28}
$,
$-\frac{1}{\sqrt{7}}c^{1}_{7}
$;
$-\frac{1}{\sqrt{7}}c^{1}_{7}
$,
$-\frac{1}{\sqrt{7}}c^{2}_{7}
$;
$\frac{1}{\sqrt{7}\mathrm{i}}s^{5}_{28}
$)
 $\oplus$
($\sqrt{\frac{1}{3}}$,
$\sqrt{\frac{2}{3}}$;
$-\sqrt{\frac{1}{3}}$)

Pass. 

 \ \color{black}

 \color{blue}

\noindent 455: (dims,levels) = $(4 , 
2;28,
12
)$,
irreps = $4_{7}^{3}
\hskip -1.5pt \otimes \hskip -1.5pt
1_{4}^{1,0}\oplus
2_{3}^{1,8}
\hskip -1.5pt \otimes \hskip -1.5pt
1_{4}^{1,0}$,
pord$(\rho_\text{isum}(\mathfrak{t})) = 21$,

\vskip 0.7ex
\hangindent=5.5em \hangafter=1
{\white .}\hskip 1em $\rho_\text{isum}(\mathfrak{t})$ =
 $( \frac{1}{4},
\frac{3}{28},
\frac{19}{28},
\frac{27}{28} )
\oplus
( \frac{1}{4},
\frac{11}{12} )
$,

\vskip 0.7ex
\hangindent=5.5em \hangafter=1
{\white .}\hskip 1em $\rho_\text{isum}(\mathfrak{s})$ =
($-\sqrt{\frac{1}{7}}$,
$\sqrt{\frac{2}{7}}$,
$\sqrt{\frac{2}{7}}$,
$\sqrt{\frac{2}{7}}$;
$-\frac{1}{\sqrt{7}}c^{2}_{7}
$,
$\frac{1}{\sqrt{7}\mathrm{i}}s^{5}_{28}
$,
$-\frac{1}{\sqrt{7}}c^{1}_{7}
$;
$-\frac{1}{\sqrt{7}}c^{1}_{7}
$,
$-\frac{1}{\sqrt{7}}c^{2}_{7}
$;
$\frac{1}{\sqrt{7}\mathrm{i}}s^{5}_{28}
$)
 $\oplus$
($-\sqrt{\frac{1}{3}}$,
$\sqrt{\frac{2}{3}}$;
$\sqrt{\frac{1}{3}}$)

Pass. 

 \ \color{black}

 \color{blue}

\noindent 456: (dims,levels) = $(4 , 
2;28,
12
)$,
irreps = $4_{7}^{1}
\hskip -1.5pt \otimes \hskip -1.5pt
1_{4}^{1,0}\oplus
2_{3}^{1,0}
\hskip -1.5pt \otimes \hskip -1.5pt
1_{4}^{1,0}$,
pord$(\rho_\text{isum}(\mathfrak{t})) = 21$,

\vskip 0.7ex
\hangindent=5.5em \hangafter=1
{\white .}\hskip 1em $\rho_\text{isum}(\mathfrak{t})$ =
 $( \frac{1}{4},
\frac{11}{28},
\frac{15}{28},
\frac{23}{28} )
\oplus
( \frac{1}{4},
\frac{7}{12} )
$,

\vskip 0.7ex
\hangindent=5.5em \hangafter=1
{\white .}\hskip 1em $\rho_\text{isum}(\mathfrak{s})$ =
($\sqrt{\frac{1}{7}}$,
$\sqrt{\frac{2}{7}}$,
$\sqrt{\frac{2}{7}}$,
$\sqrt{\frac{2}{7}}$;
$\frac{1}{\sqrt{7}}c^{2}_{7}
$,
$\frac{1}{\sqrt{7}}c^{1}_{7}
$,
$-\frac{1}{\sqrt{7}\mathrm{i}}s^{5}_{28}
$;
$-\frac{1}{\sqrt{7}\mathrm{i}}s^{5}_{28}
$,
$\frac{1}{\sqrt{7}}c^{2}_{7}
$;
$\frac{1}{\sqrt{7}}c^{1}_{7}
$)
 $\oplus$
($\sqrt{\frac{1}{3}}$,
$\sqrt{\frac{2}{3}}$;
$-\sqrt{\frac{1}{3}}$)

Pass. 

 \ \color{black}

 \color{blue}

\noindent 457: (dims,levels) = $(4 , 
2;28,
12
)$,
irreps = $4_{7}^{1}
\hskip -1.5pt \otimes \hskip -1.5pt
1_{4}^{1,0}\oplus
2_{3}^{1,8}
\hskip -1.5pt \otimes \hskip -1.5pt
1_{4}^{1,0}$,
pord$(\rho_\text{isum}(\mathfrak{t})) = 21$,

\vskip 0.7ex
\hangindent=5.5em \hangafter=1
{\white .}\hskip 1em $\rho_\text{isum}(\mathfrak{t})$ =
 $( \frac{1}{4},
\frac{11}{28},
\frac{15}{28},
\frac{23}{28} )
\oplus
( \frac{1}{4},
\frac{11}{12} )
$,

\vskip 0.7ex
\hangindent=5.5em \hangafter=1
{\white .}\hskip 1em $\rho_\text{isum}(\mathfrak{s})$ =
($\sqrt{\frac{1}{7}}$,
$\sqrt{\frac{2}{7}}$,
$\sqrt{\frac{2}{7}}$,
$\sqrt{\frac{2}{7}}$;
$\frac{1}{\sqrt{7}}c^{2}_{7}
$,
$\frac{1}{\sqrt{7}}c^{1}_{7}
$,
$-\frac{1}{\sqrt{7}\mathrm{i}}s^{5}_{28}
$;
$-\frac{1}{\sqrt{7}\mathrm{i}}s^{5}_{28}
$,
$\frac{1}{\sqrt{7}}c^{2}_{7}
$;
$\frac{1}{\sqrt{7}}c^{1}_{7}
$)
 $\oplus$
($-\sqrt{\frac{1}{3}}$,
$\sqrt{\frac{2}{3}}$;
$\sqrt{\frac{1}{3}}$)

Pass. 

 \ \color{black}

 \color{blue}

\noindent 458: (dims,levels) = $(4 , 
2;30,
10
)$,
irreps = $2_{5}^{2}
\hskip -1.5pt \otimes \hskip -1.5pt
2_{3}^{1,0}
\hskip -1.5pt \otimes \hskip -1.5pt
1_{2}^{1,0}\oplus
2_{5}^{2}
\hskip -1.5pt \otimes \hskip -1.5pt
1_{2}^{1,0}$,
pord$(\rho_\text{isum}(\mathfrak{t})) = 15$,

\vskip 0.7ex
\hangindent=5.5em \hangafter=1
{\white .}\hskip 1em $\rho_\text{isum}(\mathfrak{t})$ =
 $( \frac{1}{10},
\frac{9}{10},
\frac{7}{30},
\frac{13}{30} )
\oplus
( \frac{1}{10},
\frac{9}{10} )
$,

\vskip 0.7ex
\hangindent=5.5em \hangafter=1
{\white .}\hskip 1em $\rho_\text{isum}(\mathfrak{s})$ =
($-\frac{1}{\sqrt{15}}c^{1}_{20}
$,
$-\frac{1}{\sqrt{15}}c^{3}_{20}
$,
$-\frac{2}{\sqrt{30}}c^{3}_{20}
$,
$\frac{2}{\sqrt{30}}c^{1}_{20}
$;
$\frac{1}{\sqrt{15}}c^{1}_{20}
$,
$\frac{2}{\sqrt{30}}c^{1}_{20}
$,
$\frac{2}{\sqrt{30}}c^{3}_{20}
$;
$-\frac{1}{\sqrt{15}}c^{1}_{20}
$,
$-\frac{1}{\sqrt{15}}c^{3}_{20}
$;
$\frac{1}{\sqrt{15}}c^{1}_{20}
$)
 $\oplus$
$\mathrm{i}$($-\frac{1}{\sqrt{5}}c^{1}_{20}
$,
$\frac{1}{\sqrt{5}}c^{3}_{20}
$;\ \ 
$\frac{1}{\sqrt{5}}c^{1}_{20}
$)

Pass. 

 \ \color{black}

 \color{blue}

\noindent 459: (dims,levels) = $(4 , 
2;30,
10
)$,
irreps = $2_{5}^{1}
\hskip -1.5pt \otimes \hskip -1.5pt
2_{3}^{1,0}
\hskip -1.5pt \otimes \hskip -1.5pt
1_{2}^{1,0}\oplus
2_{5}^{1}
\hskip -1.5pt \otimes \hskip -1.5pt
1_{2}^{1,0}$,
pord$(\rho_\text{isum}(\mathfrak{t})) = 15$,

\vskip 0.7ex
\hangindent=5.5em \hangafter=1
{\white .}\hskip 1em $\rho_\text{isum}(\mathfrak{t})$ =
 $( \frac{3}{10},
\frac{7}{10},
\frac{1}{30},
\frac{19}{30} )
\oplus
( \frac{3}{10},
\frac{7}{10} )
$,

\vskip 0.7ex
\hangindent=5.5em \hangafter=1
{\white .}\hskip 1em $\rho_\text{isum}(\mathfrak{s})$ =
($-\frac{1}{\sqrt{15}}c^{3}_{20}
$,
$\frac{1}{\sqrt{15}}c^{1}_{20}
$,
$\frac{2}{\sqrt{30}}c^{1}_{20}
$,
$-\frac{2}{\sqrt{30}}c^{3}_{20}
$;
$\frac{1}{\sqrt{15}}c^{3}_{20}
$,
$\frac{2}{\sqrt{30}}c^{3}_{20}
$,
$\frac{2}{\sqrt{30}}c^{1}_{20}
$;
$-\frac{1}{\sqrt{15}}c^{3}_{20}
$,
$-\frac{1}{\sqrt{15}}c^{1}_{20}
$;
$\frac{1}{\sqrt{15}}c^{3}_{20}
$)
 $\oplus$
$\mathrm{i}$($-\frac{1}{\sqrt{5}}c^{3}_{20}
$,
$\frac{1}{\sqrt{5}}c^{1}_{20}
$;\ \ 
$\frac{1}{\sqrt{5}}c^{3}_{20}
$)

Pass. 

 \ \color{black}

 \color{blue}

\noindent 460: (dims,levels) = $(4 , 
2;30,
15
)$,
irreps = $2_{5}^{1}
\hskip -1.5pt \otimes \hskip -1.5pt
2_{2}^{1,0}
\hskip -1.5pt \otimes \hskip -1.5pt
1_{3}^{1,0}\oplus
2_{5}^{1}
\hskip -1.5pt \otimes \hskip -1.5pt
1_{3}^{1,0}$,
pord$(\rho_\text{isum}(\mathfrak{t})) = 10$,

\vskip 0.7ex
\hangindent=5.5em \hangafter=1
{\white .}\hskip 1em $\rho_\text{isum}(\mathfrak{t})$ =
 $( \frac{2}{15},
\frac{8}{15},
\frac{1}{30},
\frac{19}{30} )
\oplus
( \frac{2}{15},
\frac{8}{15} )
$,

\vskip 0.7ex
\hangindent=5.5em \hangafter=1
{\white .}\hskip 1em $\rho_\text{isum}(\mathfrak{s})$ =
$\mathrm{i}$($-\frac{1}{2\sqrt{5}}c^{3}_{20}
$,
$\frac{1}{2\sqrt{5}}c^{1}_{20}
$,
$\frac{3}{2\sqrt{15}}c^{1}_{20}
$,
$\frac{3}{2\sqrt{15}}c^{3}_{20}
$;\ \ 
$\frac{1}{2\sqrt{5}}c^{3}_{20}
$,
$\frac{3}{2\sqrt{15}}c^{3}_{20}
$,
$-\frac{3}{2\sqrt{15}}c^{1}_{20}
$;\ \ 
$-\frac{1}{2\sqrt{5}}c^{3}_{20}
$,
$\frac{1}{2\sqrt{5}}c^{1}_{20}
$;\ \ 
$\frac{1}{2\sqrt{5}}c^{3}_{20}
$)
 $\oplus$
$\mathrm{i}$($\frac{1}{\sqrt{5}}c^{3}_{20}
$,
$\frac{1}{\sqrt{5}}c^{1}_{20}
$;\ \ 
$-\frac{1}{\sqrt{5}}c^{3}_{20}
$)

Pass. 

 \ \color{black}

 \color{blue}

\noindent 461: (dims,levels) = $(4 , 
2;30,
15
)$,
irreps = $2_{5}^{2}
\hskip -1.5pt \otimes \hskip -1.5pt
2_{2}^{1,0}
\hskip -1.5pt \otimes \hskip -1.5pt
1_{3}^{1,0}\oplus
2_{5}^{2}
\hskip -1.5pt \otimes \hskip -1.5pt
1_{3}^{1,0}$,
pord$(\rho_\text{isum}(\mathfrak{t})) = 10$,

\vskip 0.7ex
\hangindent=5.5em \hangafter=1
{\white .}\hskip 1em $\rho_\text{isum}(\mathfrak{t})$ =
 $( \frac{11}{15},
\frac{14}{15},
\frac{7}{30},
\frac{13}{30} )
\oplus
( \frac{11}{15},
\frac{14}{15} )
$,

\vskip 0.7ex
\hangindent=5.5em \hangafter=1
{\white .}\hskip 1em $\rho_\text{isum}(\mathfrak{s})$ =
$\mathrm{i}$($\frac{1}{2\sqrt{5}}c^{1}_{20}
$,
$\frac{1}{2\sqrt{5}}c^{3}_{20}
$,
$\frac{3}{2\sqrt{15}}c^{1}_{20}
$,
$\frac{3}{2\sqrt{15}}c^{3}_{20}
$;\ \ 
$-\frac{1}{2\sqrt{5}}c^{1}_{20}
$,
$\frac{3}{2\sqrt{15}}c^{3}_{20}
$,
$-\frac{3}{2\sqrt{15}}c^{1}_{20}
$;\ \ 
$-\frac{1}{2\sqrt{5}}c^{1}_{20}
$,
$-\frac{1}{2\sqrt{5}}c^{3}_{20}
$;\ \ 
$\frac{1}{2\sqrt{5}}c^{1}_{20}
$)
 $\oplus$
$\mathrm{i}$($-\frac{1}{\sqrt{5}}c^{1}_{20}
$,
$\frac{1}{\sqrt{5}}c^{3}_{20}
$;\ \ 
$\frac{1}{\sqrt{5}}c^{1}_{20}
$)

Pass. 

 \ \color{black}

 \color{blue}

\noindent 462: (dims,levels) = $(4 , 
2;30,
30
)$,
irreps = $4_{5,1}^{1}
\hskip -1.5pt \otimes \hskip -1.5pt
1_{3}^{1,0}
\hskip -1.5pt \otimes \hskip -1.5pt
1_{2}^{1,0}\oplus
2_{5}^{1}
\hskip -1.5pt \otimes \hskip -1.5pt
1_{3}^{1,0}
\hskip -1.5pt \otimes \hskip -1.5pt
1_{2}^{1,0}$,
pord$(\rho_\text{isum}(\mathfrak{t})) = 5$,

\vskip 0.7ex
\hangindent=5.5em \hangafter=1
{\white .}\hskip 1em $\rho_\text{isum}(\mathfrak{t})$ =
 $( \frac{1}{30},
\frac{7}{30},
\frac{13}{30},
\frac{19}{30} )
\oplus
( \frac{1}{30},
\frac{19}{30} )
$,

\vskip 0.7ex
\hangindent=5.5em \hangafter=1
{\white .}\hskip 1em $\rho_\text{isum}(\mathfrak{s})$ =
$\mathrm{i}$($-\frac{1}{5}c^{1}_{20}
-\frac{1}{5}c^{3}_{20}
$,
$\frac{2}{5}c^{2}_{15}
+\frac{1}{5}c^{3}_{15}
$,
$-\frac{1}{5}+\frac{2}{5}c^{1}_{15}
-\frac{1}{5}c^{3}_{15}
$,
$\frac{1}{5}c^{1}_{20}
-\frac{1}{5}c^{3}_{20}
$;\ \ 
$\frac{1}{5}c^{1}_{20}
-\frac{1}{5}c^{3}_{20}
$,
$\frac{1}{5}c^{1}_{20}
+\frac{1}{5}c^{3}_{20}
$,
$-\frac{1}{5}+\frac{2}{5}c^{1}_{15}
-\frac{1}{5}c^{3}_{15}
$;\ \ 
$-\frac{1}{5}c^{1}_{20}
+\frac{1}{5}c^{3}_{20}
$,
$-\frac{2}{5}c^{2}_{15}
-\frac{1}{5}c^{3}_{15}
$;\ \ 
$\frac{1}{5}c^{1}_{20}
+\frac{1}{5}c^{3}_{20}
$)
 $\oplus$
$\mathrm{i}$($\frac{1}{\sqrt{5}}c^{3}_{20}
$,
$\frac{1}{\sqrt{5}}c^{1}_{20}
$;\ \ 
$-\frac{1}{\sqrt{5}}c^{3}_{20}
$)

Pass. 

 \ \color{black}

 \color{blue}

\noindent 463: (dims,levels) = $(4 , 
2;30,
30
)$,
irreps = $4_{5,2}^{1}
\hskip -1.5pt \otimes \hskip -1.5pt
1_{3}^{1,0}
\hskip -1.5pt \otimes \hskip -1.5pt
1_{2}^{1,0}\oplus
2_{5}^{1}
\hskip -1.5pt \otimes \hskip -1.5pt
1_{3}^{1,0}
\hskip -1.5pt \otimes \hskip -1.5pt
1_{2}^{1,0}$,
pord$(\rho_\text{isum}(\mathfrak{t})) = 5$,

\vskip 0.7ex
\hangindent=5.5em \hangafter=1
{\white .}\hskip 1em $\rho_\text{isum}(\mathfrak{t})$ =
 $( \frac{1}{30},
\frac{7}{30},
\frac{13}{30},
\frac{19}{30} )
\oplus
( \frac{1}{30},
\frac{19}{30} )
$,

\vskip 0.7ex
\hangindent=5.5em \hangafter=1
{\white .}\hskip 1em $\rho_\text{isum}(\mathfrak{s})$ =
($-\sqrt{\frac{1}{5}}$,
$\frac{-5+\sqrt{5}}{10}$,
$-\frac{5+\sqrt{5}}{10}$,
$\sqrt{\frac{1}{5}}$;
$\sqrt{\frac{1}{5}}$,
$-\sqrt{\frac{1}{5}}$,
$-\frac{5+\sqrt{5}}{10}$;
$\sqrt{\frac{1}{5}}$,
$\frac{-5+\sqrt{5}}{10}$;
$-\sqrt{\frac{1}{5}}$)
 $\oplus$
$\mathrm{i}$($\frac{1}{\sqrt{5}}c^{3}_{20}
$,
$\frac{1}{\sqrt{5}}c^{1}_{20}
$;\ \ 
$-\frac{1}{\sqrt{5}}c^{3}_{20}
$)

Pass. 

 \ \color{black}

 \color{blue}

\noindent 464: (dims,levels) = $(4 , 
2;30,
30
)$,
irreps = $4_{5,1}^{1}
\hskip -1.5pt \otimes \hskip -1.5pt
1_{3}^{1,0}
\hskip -1.5pt \otimes \hskip -1.5pt
1_{2}^{1,0}\oplus
2_{5}^{2}
\hskip -1.5pt \otimes \hskip -1.5pt
1_{3}^{1,0}
\hskip -1.5pt \otimes \hskip -1.5pt
1_{2}^{1,0}$,
pord$(\rho_\text{isum}(\mathfrak{t})) = 5$,

\vskip 0.7ex
\hangindent=5.5em \hangafter=1
{\white .}\hskip 1em $\rho_\text{isum}(\mathfrak{t})$ =
 $( \frac{1}{30},
\frac{7}{30},
\frac{13}{30},
\frac{19}{30} )
\oplus
( \frac{7}{30},
\frac{13}{30} )
$,

\vskip 0.7ex
\hangindent=5.5em \hangafter=1
{\white .}\hskip 1em $\rho_\text{isum}(\mathfrak{s})$ =
$\mathrm{i}$($-\frac{1}{5}c^{1}_{20}
-\frac{1}{5}c^{3}_{20}
$,
$\frac{2}{5}c^{2}_{15}
+\frac{1}{5}c^{3}_{15}
$,
$-\frac{1}{5}+\frac{2}{5}c^{1}_{15}
-\frac{1}{5}c^{3}_{15}
$,
$\frac{1}{5}c^{1}_{20}
-\frac{1}{5}c^{3}_{20}
$;\ \ 
$\frac{1}{5}c^{1}_{20}
-\frac{1}{5}c^{3}_{20}
$,
$\frac{1}{5}c^{1}_{20}
+\frac{1}{5}c^{3}_{20}
$,
$-\frac{1}{5}+\frac{2}{5}c^{1}_{15}
-\frac{1}{5}c^{3}_{15}
$;\ \ 
$-\frac{1}{5}c^{1}_{20}
+\frac{1}{5}c^{3}_{20}
$,
$-\frac{2}{5}c^{2}_{15}
-\frac{1}{5}c^{3}_{15}
$;\ \ 
$\frac{1}{5}c^{1}_{20}
+\frac{1}{5}c^{3}_{20}
$)
 $\oplus$
$\mathrm{i}$($\frac{1}{\sqrt{5}}c^{1}_{20}
$,
$\frac{1}{\sqrt{5}}c^{3}_{20}
$;\ \ 
$-\frac{1}{\sqrt{5}}c^{1}_{20}
$)

Pass. 

 \ \color{black}

 \color{blue}

\noindent 465: (dims,levels) = $(4 , 
2;30,
30
)$,
irreps = $4_{5,2}^{1}
\hskip -1.5pt \otimes \hskip -1.5pt
1_{3}^{1,0}
\hskip -1.5pt \otimes \hskip -1.5pt
1_{2}^{1,0}\oplus
2_{5}^{2}
\hskip -1.5pt \otimes \hskip -1.5pt
1_{3}^{1,0}
\hskip -1.5pt \otimes \hskip -1.5pt
1_{2}^{1,0}$,
pord$(\rho_\text{isum}(\mathfrak{t})) = 5$,

\vskip 0.7ex
\hangindent=5.5em \hangafter=1
{\white .}\hskip 1em $\rho_\text{isum}(\mathfrak{t})$ =
 $( \frac{1}{30},
\frac{7}{30},
\frac{13}{30},
\frac{19}{30} )
\oplus
( \frac{7}{30},
\frac{13}{30} )
$,

\vskip 0.7ex
\hangindent=5.5em \hangafter=1
{\white .}\hskip 1em $\rho_\text{isum}(\mathfrak{s})$ =
($-\sqrt{\frac{1}{5}}$,
$\frac{-5+\sqrt{5}}{10}$,
$-\frac{5+\sqrt{5}}{10}$,
$\sqrt{\frac{1}{5}}$;
$\sqrt{\frac{1}{5}}$,
$-\sqrt{\frac{1}{5}}$,
$-\frac{5+\sqrt{5}}{10}$;
$\sqrt{\frac{1}{5}}$,
$\frac{-5+\sqrt{5}}{10}$;
$-\sqrt{\frac{1}{5}}$)
 $\oplus$
$\mathrm{i}$($\frac{1}{\sqrt{5}}c^{1}_{20}
$,
$\frac{1}{\sqrt{5}}c^{3}_{20}
$;\ \ 
$-\frac{1}{\sqrt{5}}c^{1}_{20}
$)

Pass. 

 \ \color{black}

 \color{blue}

\noindent 466: (dims,levels) = $(4 , 
2;30,
30
)$,
irreps = $2_{5}^{2}
\hskip -1.5pt \otimes \hskip -1.5pt
2_{3}^{1,0}
\hskip -1.5pt \otimes \hskip -1.5pt
1_{2}^{1,0}\oplus
2_{5}^{2}
\hskip -1.5pt \otimes \hskip -1.5pt
1_{3}^{1,0}
\hskip -1.5pt \otimes \hskip -1.5pt
1_{2}^{1,0}$,
pord$(\rho_\text{isum}(\mathfrak{t})) = 15$,

\vskip 0.7ex
\hangindent=5.5em \hangafter=1
{\white .}\hskip 1em $\rho_\text{isum}(\mathfrak{t})$ =
 $( \frac{1}{10},
\frac{9}{10},
\frac{7}{30},
\frac{13}{30} )
\oplus
( \frac{7}{30},
\frac{13}{30} )
$,

\vskip 0.7ex
\hangindent=5.5em \hangafter=1
{\white .}\hskip 1em $\rho_\text{isum}(\mathfrak{s})$ =
($-\frac{1}{\sqrt{15}}c^{1}_{20}
$,
$-\frac{1}{\sqrt{15}}c^{3}_{20}
$,
$-\frac{2}{\sqrt{30}}c^{3}_{20}
$,
$\frac{2}{\sqrt{30}}c^{1}_{20}
$;
$\frac{1}{\sqrt{15}}c^{1}_{20}
$,
$\frac{2}{\sqrt{30}}c^{1}_{20}
$,
$\frac{2}{\sqrt{30}}c^{3}_{20}
$;
$-\frac{1}{\sqrt{15}}c^{1}_{20}
$,
$-\frac{1}{\sqrt{15}}c^{3}_{20}
$;
$\frac{1}{\sqrt{15}}c^{1}_{20}
$)
 $\oplus$
$\mathrm{i}$($\frac{1}{\sqrt{5}}c^{1}_{20}
$,
$\frac{1}{\sqrt{5}}c^{3}_{20}
$;\ \ 
$-\frac{1}{\sqrt{5}}c^{1}_{20}
$)

Pass. 

 \ \color{black}

 \color{blue}

\noindent 467: (dims,levels) = $(4 , 
2;30,
30
)$,
irreps = $2_{5}^{1}
\hskip -1.5pt \otimes \hskip -1.5pt
2_{2}^{1,0}
\hskip -1.5pt \otimes \hskip -1.5pt
1_{3}^{1,0}\oplus
2_{5}^{1}
\hskip -1.5pt \otimes \hskip -1.5pt
1_{3}^{1,0}
\hskip -1.5pt \otimes \hskip -1.5pt
1_{2}^{1,0}$,
pord$(\rho_\text{isum}(\mathfrak{t})) = 10$,

\vskip 0.7ex
\hangindent=5.5em \hangafter=1
{\white .}\hskip 1em $\rho_\text{isum}(\mathfrak{t})$ =
 $( \frac{2}{15},
\frac{8}{15},
\frac{1}{30},
\frac{19}{30} )
\oplus
( \frac{1}{30},
\frac{19}{30} )
$,

\vskip 0.7ex
\hangindent=5.5em \hangafter=1
{\white .}\hskip 1em $\rho_\text{isum}(\mathfrak{s})$ =
$\mathrm{i}$($-\frac{1}{2\sqrt{5}}c^{3}_{20}
$,
$\frac{1}{2\sqrt{5}}c^{1}_{20}
$,
$\frac{3}{2\sqrt{15}}c^{1}_{20}
$,
$\frac{3}{2\sqrt{15}}c^{3}_{20}
$;\ \ 
$\frac{1}{2\sqrt{5}}c^{3}_{20}
$,
$\frac{3}{2\sqrt{15}}c^{3}_{20}
$,
$-\frac{3}{2\sqrt{15}}c^{1}_{20}
$;\ \ 
$-\frac{1}{2\sqrt{5}}c^{3}_{20}
$,
$\frac{1}{2\sqrt{5}}c^{1}_{20}
$;\ \ 
$\frac{1}{2\sqrt{5}}c^{3}_{20}
$)
 $\oplus$
$\mathrm{i}$($\frac{1}{\sqrt{5}}c^{3}_{20}
$,
$\frac{1}{\sqrt{5}}c^{1}_{20}
$;\ \ 
$-\frac{1}{\sqrt{5}}c^{3}_{20}
$)

Pass. 

 \ \color{black}

 \color{blue}

\noindent 468: (dims,levels) = $(4 , 
2;30,
30
)$,
irreps = $2_{5}^{1}
\hskip -1.5pt \otimes \hskip -1.5pt
2_{3}^{1,0}
\hskip -1.5pt \otimes \hskip -1.5pt
1_{2}^{1,0}\oplus
2_{5}^{1}
\hskip -1.5pt \otimes \hskip -1.5pt
1_{3}^{1,0}
\hskip -1.5pt \otimes \hskip -1.5pt
1_{2}^{1,0}$,
pord$(\rho_\text{isum}(\mathfrak{t})) = 15$,

\vskip 0.7ex
\hangindent=5.5em \hangafter=1
{\white .}\hskip 1em $\rho_\text{isum}(\mathfrak{t})$ =
 $( \frac{3}{10},
\frac{7}{10},
\frac{1}{30},
\frac{19}{30} )
\oplus
( \frac{1}{30},
\frac{19}{30} )
$,

\vskip 0.7ex
\hangindent=5.5em \hangafter=1
{\white .}\hskip 1em $\rho_\text{isum}(\mathfrak{s})$ =
($-\frac{1}{\sqrt{15}}c^{3}_{20}
$,
$\frac{1}{\sqrt{15}}c^{1}_{20}
$,
$\frac{2}{\sqrt{30}}c^{1}_{20}
$,
$-\frac{2}{\sqrt{30}}c^{3}_{20}
$;
$\frac{1}{\sqrt{15}}c^{3}_{20}
$,
$\frac{2}{\sqrt{30}}c^{3}_{20}
$,
$\frac{2}{\sqrt{30}}c^{1}_{20}
$;
$-\frac{1}{\sqrt{15}}c^{3}_{20}
$,
$-\frac{1}{\sqrt{15}}c^{1}_{20}
$;
$\frac{1}{\sqrt{15}}c^{3}_{20}
$)
 $\oplus$
$\mathrm{i}$($\frac{1}{\sqrt{5}}c^{3}_{20}
$,
$\frac{1}{\sqrt{5}}c^{1}_{20}
$;\ \ 
$-\frac{1}{\sqrt{5}}c^{3}_{20}
$)

Pass. 

 \ \color{black}

 \color{blue}

\noindent 469: (dims,levels) = $(4 , 
2;30,
30
)$,
irreps = $2_{5}^{2}
\hskip -1.5pt \otimes \hskip -1.5pt
2_{2}^{1,0}
\hskip -1.5pt \otimes \hskip -1.5pt
1_{3}^{1,0}\oplus
2_{5}^{2}
\hskip -1.5pt \otimes \hskip -1.5pt
1_{3}^{1,0}
\hskip -1.5pt \otimes \hskip -1.5pt
1_{2}^{1,0}$,
pord$(\rho_\text{isum}(\mathfrak{t})) = 10$,

\vskip 0.7ex
\hangindent=5.5em \hangafter=1
{\white .}\hskip 1em $\rho_\text{isum}(\mathfrak{t})$ =
 $( \frac{11}{15},
\frac{14}{15},
\frac{7}{30},
\frac{13}{30} )
\oplus
( \frac{7}{30},
\frac{13}{30} )
$,

\vskip 0.7ex
\hangindent=5.5em \hangafter=1
{\white .}\hskip 1em $\rho_\text{isum}(\mathfrak{s})$ =
$\mathrm{i}$($\frac{1}{2\sqrt{5}}c^{1}_{20}
$,
$\frac{1}{2\sqrt{5}}c^{3}_{20}
$,
$\frac{3}{2\sqrt{15}}c^{1}_{20}
$,
$\frac{3}{2\sqrt{15}}c^{3}_{20}
$;\ \ 
$-\frac{1}{2\sqrt{5}}c^{1}_{20}
$,
$\frac{3}{2\sqrt{15}}c^{3}_{20}
$,
$-\frac{3}{2\sqrt{15}}c^{1}_{20}
$;\ \ 
$-\frac{1}{2\sqrt{5}}c^{1}_{20}
$,
$-\frac{1}{2\sqrt{5}}c^{3}_{20}
$;\ \ 
$\frac{1}{2\sqrt{5}}c^{1}_{20}
$)
 $\oplus$
$\mathrm{i}$($\frac{1}{\sqrt{5}}c^{1}_{20}
$,
$\frac{1}{\sqrt{5}}c^{3}_{20}
$;\ \ 
$-\frac{1}{\sqrt{5}}c^{1}_{20}
$)

Pass. 

 \ \color{black}

\noindent 470: (dims,levels) = $(4 , 
2;36,
4
)$,
irreps = $4_{9,2}^{1,0}
\hskip -1.5pt \otimes \hskip -1.5pt
1_{4}^{1,0}\oplus
2_{4}^{1,0}$,
pord$(\rho_\text{isum}(\mathfrak{t})) = 18$,

\vskip 0.7ex
\hangindent=5.5em \hangafter=1
{\white .}\hskip 1em $\rho_\text{isum}(\mathfrak{t})$ =
 $( \frac{1}{4},
\frac{1}{36},
\frac{13}{36},
\frac{25}{36} )
\oplus
( \frac{1}{4},
\frac{3}{4} )
$,

\vskip 0.7ex
\hangindent=5.5em \hangafter=1
{\white .}\hskip 1em $\rho_\text{isum}(\mathfrak{s})$ =
$\mathrm{i}$($0$,
$\sqrt{\frac{1}{3}}$,
$\sqrt{\frac{1}{3}}$,
$\sqrt{\frac{1}{3}}$;\ \ 
$\frac{1}{3} c_9^4 $,
$\frac{1}{3}c^{1}_{9}
$,
$\frac{1}{3}c^{2}_{9}
$;\ \ 
$\frac{1}{3}c^{2}_{9}
$,
$\frac{1}{3} c_9^4 $;\ \ 
$\frac{1}{3}c^{1}_{9}
$)
 $\oplus$
$\mathrm{i}$($-\frac{1}{2}$,
$\sqrt{\frac{3}{4}}$;\ \ 
$\frac{1}{2}$)

Fail:
Integral: $D_{\rho}(\sigma)_{\theta} \propto $ id,
 for all $\sigma$ and all $\theta$-eigenspaces that can contain unit. Prop. B.5 (6)

 \ \color{black}

\noindent 471: (dims,levels) = $(4 , 
2;36,
4
)$,
irreps = $4_{9,1}^{1,0}
\hskip -1.5pt \otimes \hskip -1.5pt
1_{4}^{1,0}\oplus
2_{4}^{1,0}$,
pord$(\rho_\text{isum}(\mathfrak{t})) = 18$,

\vskip 0.7ex
\hangindent=5.5em \hangafter=1
{\white .}\hskip 1em $\rho_\text{isum}(\mathfrak{t})$ =
 $( \frac{1}{4},
\frac{1}{36},
\frac{13}{36},
\frac{25}{36} )
\oplus
( \frac{1}{4},
\frac{3}{4} )
$,

\vskip 0.7ex
\hangindent=5.5em \hangafter=1
{\white .}\hskip 1em $\rho_\text{isum}(\mathfrak{s})$ =
($0$,
$-\sqrt{\frac{1}{3}}$,
$-\sqrt{\frac{1}{3}}$,
$-\sqrt{\frac{1}{3}}$;
$-\frac{1}{3}c^{1}_{36}
+\frac{1}{3}c^{5}_{36}
$,
$-\frac{1}{3}c^{5}_{36}
$,
$\frac{1}{3}c^{1}_{36}
$;
$\frac{1}{3}c^{1}_{36}
$,
$-\frac{1}{3}c^{1}_{36}
+\frac{1}{3}c^{5}_{36}
$;
$-\frac{1}{3}c^{5}_{36}
$)
 $\oplus$
$\mathrm{i}$($-\frac{1}{2}$,
$\sqrt{\frac{3}{4}}$;\ \ 
$\frac{1}{2}$)

Fail:
Tr$_I(C) = -1 <$  0 for I = [ 3/4 ]. Prop. B.4 (1) eqn. (B.18)

 \ \color{black}

\noindent 472: (dims,levels) = $(4 , 
2;36,
4
)$,
irreps = $4_{9,2}^{5,0}
\hskip -1.5pt \otimes \hskip -1.5pt
1_{4}^{1,0}\oplus
2_{4}^{1,0}$,
pord$(\rho_\text{isum}(\mathfrak{t})) = 18$,

\vskip 0.7ex
\hangindent=5.5em \hangafter=1
{\white .}\hskip 1em $\rho_\text{isum}(\mathfrak{t})$ =
 $( \frac{1}{4},
\frac{5}{36},
\frac{17}{36},
\frac{29}{36} )
\oplus
( \frac{1}{4},
\frac{3}{4} )
$,

\vskip 0.7ex
\hangindent=5.5em \hangafter=1
{\white .}\hskip 1em $\rho_\text{isum}(\mathfrak{s})$ =
$\mathrm{i}$($0$,
$\sqrt{\frac{1}{3}}$,
$\sqrt{\frac{1}{3}}$,
$\sqrt{\frac{1}{3}}$;\ \ 
$\frac{1}{3}c^{2}_{9}
$,
$\frac{1}{3}c^{1}_{9}
$,
$\frac{1}{3} c_9^4 $;\ \ 
$\frac{1}{3} c_9^4 $,
$\frac{1}{3}c^{2}_{9}
$;\ \ 
$\frac{1}{3}c^{1}_{9}
$)
 $\oplus$
$\mathrm{i}$($-\frac{1}{2}$,
$\sqrt{\frac{3}{4}}$;\ \ 
$\frac{1}{2}$)

Fail:
Integral: $D_{\rho}(\sigma)_{\theta} \propto $ id,
 for all $\sigma$ and all $\theta$-eigenspaces that can contain unit. Prop. B.5 (6)

 \ \color{black}

\noindent 473: (dims,levels) = $(4 , 
2;36,
4
)$,
irreps = $4_{9,1}^{2,0}
\hskip -1.5pt \otimes \hskip -1.5pt
1_{4}^{1,0}\oplus
2_{4}^{1,0}$,
pord$(\rho_\text{isum}(\mathfrak{t})) = 18$,

\vskip 0.7ex
\hangindent=5.5em \hangafter=1
{\white .}\hskip 1em $\rho_\text{isum}(\mathfrak{t})$ =
 $( \frac{1}{4},
\frac{5}{36},
\frac{17}{36},
\frac{29}{36} )
\oplus
( \frac{1}{4},
\frac{3}{4} )
$,

\vskip 0.7ex
\hangindent=5.5em \hangafter=1
{\white .}\hskip 1em $\rho_\text{isum}(\mathfrak{s})$ =
($0$,
$-\sqrt{\frac{1}{3}}$,
$-\sqrt{\frac{1}{3}}$,
$-\sqrt{\frac{1}{3}}$;
$-\frac{1}{3}c^{1}_{36}
$,
$\frac{1}{3}c^{5}_{36}
$,
$\frac{1}{3}c^{1}_{36}
-\frac{1}{3}c^{5}_{36}
$;
$\frac{1}{3}c^{1}_{36}
-\frac{1}{3}c^{5}_{36}
$,
$-\frac{1}{3}c^{1}_{36}
$;
$\frac{1}{3}c^{5}_{36}
$)
 $\oplus$
$\mathrm{i}$($-\frac{1}{2}$,
$\sqrt{\frac{3}{4}}$;\ \ 
$\frac{1}{2}$)

Fail:
Tr$_I(C) = -1 <$  0 for I = [ 3/4 ]. Prop. B.4 (1) eqn. (B.18)

 \ \color{black}

\noindent 474: (dims,levels) = $(4 , 
2;36,
12
)$,
irreps = $4_{9,2}^{1,0}
\hskip -1.5pt \otimes \hskip -1.5pt
1_{4}^{1,0}\oplus
2_{3}^{1,0}
\hskip -1.5pt \otimes \hskip -1.5pt
1_{4}^{1,0}$,
pord$(\rho_\text{isum}(\mathfrak{t})) = 9$,

\vskip 0.7ex
\hangindent=5.5em \hangafter=1
{\white .}\hskip 1em $\rho_\text{isum}(\mathfrak{t})$ =
 $( \frac{1}{4},
\frac{1}{36},
\frac{13}{36},
\frac{25}{36} )
\oplus
( \frac{1}{4},
\frac{7}{12} )
$,

\vskip 0.7ex
\hangindent=5.5em \hangafter=1
{\white .}\hskip 1em $\rho_\text{isum}(\mathfrak{s})$ =
$\mathrm{i}$($0$,
$\sqrt{\frac{1}{3}}$,
$\sqrt{\frac{1}{3}}$,
$\sqrt{\frac{1}{3}}$;\ \ 
$\frac{1}{3} c_9^4 $,
$\frac{1}{3}c^{1}_{9}
$,
$\frac{1}{3}c^{2}_{9}
$;\ \ 
$\frac{1}{3}c^{2}_{9}
$,
$\frac{1}{3} c_9^4 $;\ \ 
$\frac{1}{3}c^{1}_{9}
$)
 $\oplus$
($\sqrt{\frac{1}{3}}$,
$\sqrt{\frac{2}{3}}$;
$-\sqrt{\frac{1}{3}}$)

Fail:
Tr$_I(C) = -1 <$  0 for I = [ 7/12 ]. Prop. B.4 (1) eqn. (B.18)

 \ \color{black}

\noindent 475: (dims,levels) = $(4 , 
2;36,
12
)$,
irreps = $4_{9,1}^{1,0}
\hskip -1.5pt \otimes \hskip -1.5pt
1_{4}^{1,0}\oplus
2_{3}^{1,0}
\hskip -1.5pt \otimes \hskip -1.5pt
1_{4}^{1,0}$,
pord$(\rho_\text{isum}(\mathfrak{t})) = 9$,

\vskip 0.7ex
\hangindent=5.5em \hangafter=1
{\white .}\hskip 1em $\rho_\text{isum}(\mathfrak{t})$ =
 $( \frac{1}{4},
\frac{1}{36},
\frac{13}{36},
\frac{25}{36} )
\oplus
( \frac{1}{4},
\frac{7}{12} )
$,

\vskip 0.7ex
\hangindent=5.5em \hangafter=1
{\white .}\hskip 1em $\rho_\text{isum}(\mathfrak{s})$ =
($0$,
$-\sqrt{\frac{1}{3}}$,
$-\sqrt{\frac{1}{3}}$,
$-\sqrt{\frac{1}{3}}$;
$-\frac{1}{3}c^{1}_{36}
+\frac{1}{3}c^{5}_{36}
$,
$-\frac{1}{3}c^{5}_{36}
$,
$\frac{1}{3}c^{1}_{36}
$;
$\frac{1}{3}c^{1}_{36}
$,
$-\frac{1}{3}c^{1}_{36}
+\frac{1}{3}c^{5}_{36}
$;
$-\frac{1}{3}c^{5}_{36}
$)
 $\oplus$
($\sqrt{\frac{1}{3}}$,
$\sqrt{\frac{2}{3}}$;
$-\sqrt{\frac{1}{3}}$)

Fail:
Integral: $D_{\rho}(\sigma)_{\theta} \propto $ id,
 for all $\sigma$ and all $\theta$-eigenspaces that can contain unit. Prop. B.5 (6)

 \ \color{black}

\noindent 476: (dims,levels) = $(4 , 
2;36,
12
)$,
irreps = $4_{9,2}^{1,0}
\hskip -1.5pt \otimes \hskip -1.5pt
1_{4}^{1,0}\oplus
2_{3}^{1,8}
\hskip -1.5pt \otimes \hskip -1.5pt
1_{4}^{1,0}$,
pord$(\rho_\text{isum}(\mathfrak{t})) = 9$,

\vskip 0.7ex
\hangindent=5.5em \hangafter=1
{\white .}\hskip 1em $\rho_\text{isum}(\mathfrak{t})$ =
 $( \frac{1}{4},
\frac{1}{36},
\frac{13}{36},
\frac{25}{36} )
\oplus
( \frac{1}{4},
\frac{11}{12} )
$,

\vskip 0.7ex
\hangindent=5.5em \hangafter=1
{\white .}\hskip 1em $\rho_\text{isum}(\mathfrak{s})$ =
$\mathrm{i}$($0$,
$\sqrt{\frac{1}{3}}$,
$\sqrt{\frac{1}{3}}$,
$\sqrt{\frac{1}{3}}$;\ \ 
$\frac{1}{3} c_9^4 $,
$\frac{1}{3}c^{1}_{9}
$,
$\frac{1}{3}c^{2}_{9}
$;\ \ 
$\frac{1}{3}c^{2}_{9}
$,
$\frac{1}{3} c_9^4 $;\ \ 
$\frac{1}{3}c^{1}_{9}
$)
 $\oplus$
($-\sqrt{\frac{1}{3}}$,
$\sqrt{\frac{2}{3}}$;
$\sqrt{\frac{1}{3}}$)

Fail:
Tr$_I(C) = -1 <$  0 for I = [ 11/12 ]. Prop. B.4 (1) eqn. (B.18)

 \ \color{black}

\noindent 477: (dims,levels) = $(4 , 
2;36,
12
)$,
irreps = $4_{9,1}^{1,0}
\hskip -1.5pt \otimes \hskip -1.5pt
1_{4}^{1,0}\oplus
2_{3}^{1,8}
\hskip -1.5pt \otimes \hskip -1.5pt
1_{4}^{1,0}$,
pord$(\rho_\text{isum}(\mathfrak{t})) = 9$,

\vskip 0.7ex
\hangindent=5.5em \hangafter=1
{\white .}\hskip 1em $\rho_\text{isum}(\mathfrak{t})$ =
 $( \frac{1}{4},
\frac{1}{36},
\frac{13}{36},
\frac{25}{36} )
\oplus
( \frac{1}{4},
\frac{11}{12} )
$,

\vskip 0.7ex
\hangindent=5.5em \hangafter=1
{\white .}\hskip 1em $\rho_\text{isum}(\mathfrak{s})$ =
($0$,
$-\sqrt{\frac{1}{3}}$,
$-\sqrt{\frac{1}{3}}$,
$-\sqrt{\frac{1}{3}}$;
$-\frac{1}{3}c^{1}_{36}
+\frac{1}{3}c^{5}_{36}
$,
$-\frac{1}{3}c^{5}_{36}
$,
$\frac{1}{3}c^{1}_{36}
$;
$\frac{1}{3}c^{1}_{36}
$,
$-\frac{1}{3}c^{1}_{36}
+\frac{1}{3}c^{5}_{36}
$;
$-\frac{1}{3}c^{5}_{36}
$)
 $\oplus$
($-\sqrt{\frac{1}{3}}$,
$\sqrt{\frac{2}{3}}$;
$\sqrt{\frac{1}{3}}$)

Fail:
Integral: $D_{\rho}(\sigma)_{\theta} \propto $ id,
 for all $\sigma$ and all $\theta$-eigenspaces that can contain unit. Prop. B.5 (6)

 \ \color{black}

\noindent 478: (dims,levels) = $(4 , 
2;36,
12
)$,
irreps = $4_{9,2}^{5,0}
\hskip -1.5pt \otimes \hskip -1.5pt
1_{4}^{1,0}\oplus
2_{3}^{1,0}
\hskip -1.5pt \otimes \hskip -1.5pt
1_{4}^{1,0}$,
pord$(\rho_\text{isum}(\mathfrak{t})) = 9$,

\vskip 0.7ex
\hangindent=5.5em \hangafter=1
{\white .}\hskip 1em $\rho_\text{isum}(\mathfrak{t})$ =
 $( \frac{1}{4},
\frac{5}{36},
\frac{17}{36},
\frac{29}{36} )
\oplus
( \frac{1}{4},
\frac{7}{12} )
$,

\vskip 0.7ex
\hangindent=5.5em \hangafter=1
{\white .}\hskip 1em $\rho_\text{isum}(\mathfrak{s})$ =
$\mathrm{i}$($0$,
$\sqrt{\frac{1}{3}}$,
$\sqrt{\frac{1}{3}}$,
$\sqrt{\frac{1}{3}}$;\ \ 
$\frac{1}{3}c^{2}_{9}
$,
$\frac{1}{3}c^{1}_{9}
$,
$\frac{1}{3} c_9^4 $;\ \ 
$\frac{1}{3} c_9^4 $,
$\frac{1}{3}c^{2}_{9}
$;\ \ 
$\frac{1}{3}c^{1}_{9}
$)
 $\oplus$
($\sqrt{\frac{1}{3}}$,
$\sqrt{\frac{2}{3}}$;
$-\sqrt{\frac{1}{3}}$)

Fail:
Tr$_I(C) = -1 <$  0 for I = [ 7/12 ]. Prop. B.4 (1) eqn. (B.18)

 \ \color{black}

\noindent 479: (dims,levels) = $(4 , 
2;36,
12
)$,
irreps = $4_{9,1}^{2,0}
\hskip -1.5pt \otimes \hskip -1.5pt
1_{4}^{1,0}\oplus
2_{3}^{1,0}
\hskip -1.5pt \otimes \hskip -1.5pt
1_{4}^{1,0}$,
pord$(\rho_\text{isum}(\mathfrak{t})) = 9$,

\vskip 0.7ex
\hangindent=5.5em \hangafter=1
{\white .}\hskip 1em $\rho_\text{isum}(\mathfrak{t})$ =
 $( \frac{1}{4},
\frac{5}{36},
\frac{17}{36},
\frac{29}{36} )
\oplus
( \frac{1}{4},
\frac{7}{12} )
$,

\vskip 0.7ex
\hangindent=5.5em \hangafter=1
{\white .}\hskip 1em $\rho_\text{isum}(\mathfrak{s})$ =
($0$,
$-\sqrt{\frac{1}{3}}$,
$-\sqrt{\frac{1}{3}}$,
$-\sqrt{\frac{1}{3}}$;
$-\frac{1}{3}c^{1}_{36}
$,
$\frac{1}{3}c^{5}_{36}
$,
$\frac{1}{3}c^{1}_{36}
-\frac{1}{3}c^{5}_{36}
$;
$\frac{1}{3}c^{1}_{36}
-\frac{1}{3}c^{5}_{36}
$,
$-\frac{1}{3}c^{1}_{36}
$;
$\frac{1}{3}c^{5}_{36}
$)
 $\oplus$
($\sqrt{\frac{1}{3}}$,
$\sqrt{\frac{2}{3}}$;
$-\sqrt{\frac{1}{3}}$)

Fail:
Integral: $D_{\rho}(\sigma)_{\theta} \propto $ id,
 for all $\sigma$ and all $\theta$-eigenspaces that can contain unit. Prop. B.5 (6)

 \ \color{black}

\noindent 480: (dims,levels) = $(4 , 
2;36,
12
)$,
irreps = $4_{9,2}^{5,0}
\hskip -1.5pt \otimes \hskip -1.5pt
1_{4}^{1,0}\oplus
2_{3}^{1,8}
\hskip -1.5pt \otimes \hskip -1.5pt
1_{4}^{1,0}$,
pord$(\rho_\text{isum}(\mathfrak{t})) = 9$,

\vskip 0.7ex
\hangindent=5.5em \hangafter=1
{\white .}\hskip 1em $\rho_\text{isum}(\mathfrak{t})$ =
 $( \frac{1}{4},
\frac{5}{36},
\frac{17}{36},
\frac{29}{36} )
\oplus
( \frac{1}{4},
\frac{11}{12} )
$,

\vskip 0.7ex
\hangindent=5.5em \hangafter=1
{\white .}\hskip 1em $\rho_\text{isum}(\mathfrak{s})$ =
$\mathrm{i}$($0$,
$\sqrt{\frac{1}{3}}$,
$\sqrt{\frac{1}{3}}$,
$\sqrt{\frac{1}{3}}$;\ \ 
$\frac{1}{3}c^{2}_{9}
$,
$\frac{1}{3}c^{1}_{9}
$,
$\frac{1}{3} c_9^4 $;\ \ 
$\frac{1}{3} c_9^4 $,
$\frac{1}{3}c^{2}_{9}
$;\ \ 
$\frac{1}{3}c^{1}_{9}
$)
 $\oplus$
($-\sqrt{\frac{1}{3}}$,
$\sqrt{\frac{2}{3}}$;
$\sqrt{\frac{1}{3}}$)

Fail:
Tr$_I(C) = -1 <$  0 for I = [ 11/12 ]. Prop. B.4 (1) eqn. (B.18)

 \ \color{black}

\noindent 481: (dims,levels) = $(4 , 
2;36,
12
)$,
irreps = $4_{9,1}^{2,0}
\hskip -1.5pt \otimes \hskip -1.5pt
1_{4}^{1,0}\oplus
2_{3}^{1,8}
\hskip -1.5pt \otimes \hskip -1.5pt
1_{4}^{1,0}$,
pord$(\rho_\text{isum}(\mathfrak{t})) = 9$,

\vskip 0.7ex
\hangindent=5.5em \hangafter=1
{\white .}\hskip 1em $\rho_\text{isum}(\mathfrak{t})$ =
 $( \frac{1}{4},
\frac{5}{36},
\frac{17}{36},
\frac{29}{36} )
\oplus
( \frac{1}{4},
\frac{11}{12} )
$,

\vskip 0.7ex
\hangindent=5.5em \hangafter=1
{\white .}\hskip 1em $\rho_\text{isum}(\mathfrak{s})$ =
($0$,
$-\sqrt{\frac{1}{3}}$,
$-\sqrt{\frac{1}{3}}$,
$-\sqrt{\frac{1}{3}}$;
$-\frac{1}{3}c^{1}_{36}
$,
$\frac{1}{3}c^{5}_{36}
$,
$\frac{1}{3}c^{1}_{36}
-\frac{1}{3}c^{5}_{36}
$;
$\frac{1}{3}c^{1}_{36}
-\frac{1}{3}c^{5}_{36}
$,
$-\frac{1}{3}c^{1}_{36}
$;
$\frac{1}{3}c^{5}_{36}
$)
 $\oplus$
($-\sqrt{\frac{1}{3}}$,
$\sqrt{\frac{2}{3}}$;
$\sqrt{\frac{1}{3}}$)

Fail:
Integral: $D_{\rho}(\sigma)_{\theta} \propto $ id,
 for all $\sigma$ and all $\theta$-eigenspaces that can contain unit. Prop. B.5 (6)

 \ \color{black}

 \color{blue}

\noindent 482: (dims,levels) = $(4 , 
2;42,
6
)$,
irreps = $4_{7}^{1}
\hskip -1.5pt \otimes \hskip -1.5pt
1_{3}^{1,0}
\hskip -1.5pt \otimes \hskip -1.5pt
1_{2}^{1,0}\oplus
2_{3}^{1,4}
\hskip -1.5pt \otimes \hskip -1.5pt
1_{2}^{1,0}$,
pord$(\rho_\text{isum}(\mathfrak{t})) = 21$,

\vskip 0.7ex
\hangindent=5.5em \hangafter=1
{\white .}\hskip 1em $\rho_\text{isum}(\mathfrak{t})$ =
 $( \frac{5}{6},
\frac{5}{42},
\frac{17}{42},
\frac{41}{42} )
\oplus
( \frac{1}{6},
\frac{5}{6} )
$,

\vskip 0.7ex
\hangindent=5.5em \hangafter=1
{\white .}\hskip 1em $\rho_\text{isum}(\mathfrak{s})$ =
$\mathrm{i}$($\sqrt{\frac{1}{7}}$,
$\sqrt{\frac{2}{7}}$,
$\sqrt{\frac{2}{7}}$,
$\sqrt{\frac{2}{7}}$;\ \ 
$-\frac{1}{\sqrt{7}\mathrm{i}}s^{5}_{28}
$,
$\frac{1}{\sqrt{7}}c^{2}_{7}
$,
$\frac{1}{\sqrt{7}}c^{1}_{7}
$;\ \ 
$\frac{1}{\sqrt{7}}c^{1}_{7}
$,
$-\frac{1}{\sqrt{7}\mathrm{i}}s^{5}_{28}
$;\ \ 
$\frac{1}{\sqrt{7}}c^{2}_{7}
$)
 $\oplus$
$\mathrm{i}$($-\sqrt{\frac{1}{3}}$,
$\sqrt{\frac{2}{3}}$;\ \ 
$\sqrt{\frac{1}{3}}$)

Pass. 

 \ \color{black}

\noindent 483: (dims,levels) = $(4 , 
2;42,
6
)$,
irreps = $4_{7}^{1}
\hskip -1.5pt \otimes \hskip -1.5pt
1_{3}^{1,0}
\hskip -1.5pt \otimes \hskip -1.5pt
1_{2}^{1,0}\oplus
2_{2}^{1,0}
\hskip -1.5pt \otimes \hskip -1.5pt
1_{3}^{1,0}$,
pord$(\rho_\text{isum}(\mathfrak{t})) = 14$,

\vskip 0.7ex
\hangindent=5.5em \hangafter=1
{\white .}\hskip 1em $\rho_\text{isum}(\mathfrak{t})$ =
 $( \frac{5}{6},
\frac{5}{42},
\frac{17}{42},
\frac{41}{42} )
\oplus
( \frac{1}{3},
\frac{5}{6} )
$,

\vskip 0.7ex
\hangindent=5.5em \hangafter=1
{\white .}\hskip 1em $\rho_\text{isum}(\mathfrak{s})$ =
$\mathrm{i}$($\sqrt{\frac{1}{7}}$,
$\sqrt{\frac{2}{7}}$,
$\sqrt{\frac{2}{7}}$,
$\sqrt{\frac{2}{7}}$;\ \ 
$-\frac{1}{\sqrt{7}\mathrm{i}}s^{5}_{28}
$,
$\frac{1}{\sqrt{7}}c^{2}_{7}
$,
$\frac{1}{\sqrt{7}}c^{1}_{7}
$;\ \ 
$\frac{1}{\sqrt{7}}c^{1}_{7}
$,
$-\frac{1}{\sqrt{7}\mathrm{i}}s^{5}_{28}
$;\ \ 
$\frac{1}{\sqrt{7}}c^{2}_{7}
$)
 $\oplus$
($-\frac{1}{2}$,
$-\sqrt{\frac{3}{4}}$;
$\frac{1}{2}$)

Fail:
Tr$_I(C) = -1 <$  0 for I = [ 1/3 ]. Prop. B.4 (1) eqn. (B.18)

 \ \color{black}

 \color{blue}

\noindent 484: (dims,levels) = $(4 , 
2;42,
6
)$,
irreps = $4_{7}^{1}
\hskip -1.5pt \otimes \hskip -1.5pt
1_{3}^{1,0}
\hskip -1.5pt \otimes \hskip -1.5pt
1_{2}^{1,0}\oplus
2_{3}^{1,0}
\hskip -1.5pt \otimes \hskip -1.5pt
1_{2}^{1,0}$,
pord$(\rho_\text{isum}(\mathfrak{t})) = 21$,

\vskip 0.7ex
\hangindent=5.5em \hangafter=1
{\white .}\hskip 1em $\rho_\text{isum}(\mathfrak{t})$ =
 $( \frac{5}{6},
\frac{5}{42},
\frac{17}{42},
\frac{41}{42} )
\oplus
( \frac{1}{2},
\frac{5}{6} )
$,

\vskip 0.7ex
\hangindent=5.5em \hangafter=1
{\white .}\hskip 1em $\rho_\text{isum}(\mathfrak{s})$ =
$\mathrm{i}$($\sqrt{\frac{1}{7}}$,
$\sqrt{\frac{2}{7}}$,
$\sqrt{\frac{2}{7}}$,
$\sqrt{\frac{2}{7}}$;\ \ 
$-\frac{1}{\sqrt{7}\mathrm{i}}s^{5}_{28}
$,
$\frac{1}{\sqrt{7}}c^{2}_{7}
$,
$\frac{1}{\sqrt{7}}c^{1}_{7}
$;\ \ 
$\frac{1}{\sqrt{7}}c^{1}_{7}
$,
$-\frac{1}{\sqrt{7}\mathrm{i}}s^{5}_{28}
$;\ \ 
$\frac{1}{\sqrt{7}}c^{2}_{7}
$)
 $\oplus$
$\mathrm{i}$($\sqrt{\frac{1}{3}}$,
$\sqrt{\frac{2}{3}}$;\ \ 
$-\sqrt{\frac{1}{3}}$)

Pass. 

 \ \color{black}

 \color{blue}

\noindent 485: (dims,levels) = $(4 , 
2;42,
6
)$,
irreps = $4_{7}^{3}
\hskip -1.5pt \otimes \hskip -1.5pt
1_{3}^{1,0}
\hskip -1.5pt \otimes \hskip -1.5pt
1_{2}^{1,0}\oplus
2_{3}^{1,4}
\hskip -1.5pt \otimes \hskip -1.5pt
1_{2}^{1,0}$,
pord$(\rho_\text{isum}(\mathfrak{t})) = 21$,

\vskip 0.7ex
\hangindent=5.5em \hangafter=1
{\white .}\hskip 1em $\rho_\text{isum}(\mathfrak{t})$ =
 $( \frac{5}{6},
\frac{11}{42},
\frac{23}{42},
\frac{29}{42} )
\oplus
( \frac{1}{6},
\frac{5}{6} )
$,

\vskip 0.7ex
\hangindent=5.5em \hangafter=1
{\white .}\hskip 1em $\rho_\text{isum}(\mathfrak{s})$ =
$\mathrm{i}$($-\sqrt{\frac{1}{7}}$,
$\sqrt{\frac{2}{7}}$,
$\sqrt{\frac{2}{7}}$,
$\sqrt{\frac{2}{7}}$;\ \ 
$-\frac{1}{\sqrt{7}}c^{1}_{7}
$,
$-\frac{1}{\sqrt{7}}c^{2}_{7}
$,
$\frac{1}{\sqrt{7}\mathrm{i}}s^{5}_{28}
$;\ \ 
$\frac{1}{\sqrt{7}\mathrm{i}}s^{5}_{28}
$,
$-\frac{1}{\sqrt{7}}c^{1}_{7}
$;\ \ 
$-\frac{1}{\sqrt{7}}c^{2}_{7}
$)
 $\oplus$
$\mathrm{i}$($-\sqrt{\frac{1}{3}}$,
$\sqrt{\frac{2}{3}}$;\ \ 
$\sqrt{\frac{1}{3}}$)

Pass. 

 \ \color{black}

\noindent 486: (dims,levels) = $(4 , 
2;42,
6
)$,
irreps = $4_{7}^{3}
\hskip -1.5pt \otimes \hskip -1.5pt
1_{3}^{1,0}
\hskip -1.5pt \otimes \hskip -1.5pt
1_{2}^{1,0}\oplus
2_{2}^{1,0}
\hskip -1.5pt \otimes \hskip -1.5pt
1_{3}^{1,0}$,
pord$(\rho_\text{isum}(\mathfrak{t})) = 14$,

\vskip 0.7ex
\hangindent=5.5em \hangafter=1
{\white .}\hskip 1em $\rho_\text{isum}(\mathfrak{t})$ =
 $( \frac{5}{6},
\frac{11}{42},
\frac{23}{42},
\frac{29}{42} )
\oplus
( \frac{1}{3},
\frac{5}{6} )
$,

\vskip 0.7ex
\hangindent=5.5em \hangafter=1
{\white .}\hskip 1em $\rho_\text{isum}(\mathfrak{s})$ =
$\mathrm{i}$($-\sqrt{\frac{1}{7}}$,
$\sqrt{\frac{2}{7}}$,
$\sqrt{\frac{2}{7}}$,
$\sqrt{\frac{2}{7}}$;\ \ 
$-\frac{1}{\sqrt{7}}c^{1}_{7}
$,
$-\frac{1}{\sqrt{7}}c^{2}_{7}
$,
$\frac{1}{\sqrt{7}\mathrm{i}}s^{5}_{28}
$;\ \ 
$\frac{1}{\sqrt{7}\mathrm{i}}s^{5}_{28}
$,
$-\frac{1}{\sqrt{7}}c^{1}_{7}
$;\ \ 
$-\frac{1}{\sqrt{7}}c^{2}_{7}
$)
 $\oplus$
($-\frac{1}{2}$,
$-\sqrt{\frac{3}{4}}$;
$\frac{1}{2}$)

Fail:
Tr$_I(C) = -1 <$  0 for I = [ 1/3 ]. Prop. B.4 (1) eqn. (B.18)

 \ \color{black}

 \color{blue}

\noindent 487: (dims,levels) = $(4 , 
2;42,
6
)$,
irreps = $4_{7}^{3}
\hskip -1.5pt \otimes \hskip -1.5pt
1_{3}^{1,0}
\hskip -1.5pt \otimes \hskip -1.5pt
1_{2}^{1,0}\oplus
2_{3}^{1,0}
\hskip -1.5pt \otimes \hskip -1.5pt
1_{2}^{1,0}$,
pord$(\rho_\text{isum}(\mathfrak{t})) = 21$,

\vskip 0.7ex
\hangindent=5.5em \hangafter=1
{\white .}\hskip 1em $\rho_\text{isum}(\mathfrak{t})$ =
 $( \frac{5}{6},
\frac{11}{42},
\frac{23}{42},
\frac{29}{42} )
\oplus
( \frac{1}{2},
\frac{5}{6} )
$,

\vskip 0.7ex
\hangindent=5.5em \hangafter=1
{\white .}\hskip 1em $\rho_\text{isum}(\mathfrak{s})$ =
$\mathrm{i}$($-\sqrt{\frac{1}{7}}$,
$\sqrt{\frac{2}{7}}$,
$\sqrt{\frac{2}{7}}$,
$\sqrt{\frac{2}{7}}$;\ \ 
$-\frac{1}{\sqrt{7}}c^{1}_{7}
$,
$-\frac{1}{\sqrt{7}}c^{2}_{7}
$,
$\frac{1}{\sqrt{7}\mathrm{i}}s^{5}_{28}
$;\ \ 
$\frac{1}{\sqrt{7}\mathrm{i}}s^{5}_{28}
$,
$-\frac{1}{\sqrt{7}}c^{1}_{7}
$;\ \ 
$-\frac{1}{\sqrt{7}}c^{2}_{7}
$)
 $\oplus$
$\mathrm{i}$($\sqrt{\frac{1}{3}}$,
$\sqrt{\frac{2}{3}}$;\ \ 
$-\sqrt{\frac{1}{3}}$)

Pass. 

 \ \color{black}

 \color{blue}

\noindent 488: (dims,levels) = $(4 , 
2;60,
20
)$,
irreps = $2_{5}^{1}
\hskip -1.5pt \otimes \hskip -1.5pt
2_{3}^{1,0}
\hskip -1.5pt \otimes \hskip -1.5pt
1_{4}^{1,0}\oplus
2_{5}^{1}
\hskip -1.5pt \otimes \hskip -1.5pt
1_{4}^{1,0}$,
pord$(\rho_\text{isum}(\mathfrak{t})) = 15$,

\vskip 0.7ex
\hangindent=5.5em \hangafter=1
{\white .}\hskip 1em $\rho_\text{isum}(\mathfrak{t})$ =
 $( \frac{1}{20},
\frac{9}{20},
\frac{23}{60},
\frac{47}{60} )
\oplus
( \frac{1}{20},
\frac{9}{20} )
$,

\vskip 0.7ex
\hangindent=5.5em \hangafter=1
{\white .}\hskip 1em $\rho_\text{isum}(\mathfrak{s})$ =
$\mathrm{i}$($\frac{1}{\sqrt{15}}c^{3}_{20}
$,
$\frac{1}{\sqrt{15}}c^{1}_{20}
$,
$-\frac{2}{\sqrt{30}}c^{3}_{20}
$,
$\frac{2}{\sqrt{30}}c^{1}_{20}
$;\ \ 
$-\frac{1}{\sqrt{15}}c^{3}_{20}
$,
$-\frac{2}{\sqrt{30}}c^{1}_{20}
$,
$-\frac{2}{\sqrt{30}}c^{3}_{20}
$;\ \ 
$-\frac{1}{\sqrt{15}}c^{3}_{20}
$,
$\frac{1}{\sqrt{15}}c^{1}_{20}
$;\ \ 
$\frac{1}{\sqrt{15}}c^{3}_{20}
$)
 $\oplus$
($-\frac{1}{\sqrt{5}}c^{3}_{20}
$,
$\frac{1}{\sqrt{5}}c^{1}_{20}
$;
$\frac{1}{\sqrt{5}}c^{3}_{20}
$)

Pass. 

 \ \color{black}

 \color{blue}

\noindent 489: (dims,levels) = $(4 , 
2;60,
20
)$,
irreps = $2_{5}^{2}
\hskip -1.5pt \otimes \hskip -1.5pt
2_{3}^{1,0}
\hskip -1.5pt \otimes \hskip -1.5pt
1_{4}^{1,0}\oplus
2_{5}^{2}
\hskip -1.5pt \otimes \hskip -1.5pt
1_{4}^{1,0}$,
pord$(\rho_\text{isum}(\mathfrak{t})) = 15$,

\vskip 0.7ex
\hangindent=5.5em \hangafter=1
{\white .}\hskip 1em $\rho_\text{isum}(\mathfrak{t})$ =
 $( \frac{13}{20},
\frac{17}{20},
\frac{11}{60},
\frac{59}{60} )
\oplus
( \frac{13}{20},
\frac{17}{20} )
$,

\vskip 0.7ex
\hangindent=5.5em \hangafter=1
{\white .}\hskip 1em $\rho_\text{isum}(\mathfrak{s})$ =
$\mathrm{i}$($-\frac{1}{\sqrt{15}}c^{1}_{20}
$,
$\frac{1}{\sqrt{15}}c^{3}_{20}
$,
$-\frac{2}{\sqrt{30}}c^{3}_{20}
$,
$\frac{2}{\sqrt{30}}c^{1}_{20}
$;\ \ 
$\frac{1}{\sqrt{15}}c^{1}_{20}
$,
$-\frac{2}{\sqrt{30}}c^{1}_{20}
$,
$-\frac{2}{\sqrt{30}}c^{3}_{20}
$;\ \ 
$-\frac{1}{\sqrt{15}}c^{1}_{20}
$,
$-\frac{1}{\sqrt{15}}c^{3}_{20}
$;\ \ 
$\frac{1}{\sqrt{15}}c^{1}_{20}
$)
 $\oplus$
($\frac{1}{\sqrt{5}}c^{1}_{20}
$,
$\frac{1}{\sqrt{5}}c^{3}_{20}
$;
$-\frac{1}{\sqrt{5}}c^{1}_{20}
$)

Pass. 

 \ \color{black}

 \color{blue}

\noindent 490: (dims,levels) = $(4 , 
2;60,
60
)$,
irreps = $2_{5}^{1}
\hskip -1.5pt \otimes \hskip -1.5pt
2_{3}^{1,0}
\hskip -1.5pt \otimes \hskip -1.5pt
1_{4}^{1,0}\oplus
2_{5}^{1}
\hskip -1.5pt \otimes \hskip -1.5pt
1_{4}^{1,0}
\hskip -1.5pt \otimes \hskip -1.5pt
1_{3}^{1,0}$,
pord$(\rho_\text{isum}(\mathfrak{t})) = 15$,

\vskip 0.7ex
\hangindent=5.5em \hangafter=1
{\white .}\hskip 1em $\rho_\text{isum}(\mathfrak{t})$ =
 $( \frac{1}{20},
\frac{9}{20},
\frac{23}{60},
\frac{47}{60} )
\oplus
( \frac{23}{60},
\frac{47}{60} )
$,

\vskip 0.7ex
\hangindent=5.5em \hangafter=1
{\white .}\hskip 1em $\rho_\text{isum}(\mathfrak{s})$ =
$\mathrm{i}$($\frac{1}{\sqrt{15}}c^{3}_{20}
$,
$\frac{1}{\sqrt{15}}c^{1}_{20}
$,
$-\frac{2}{\sqrt{30}}c^{3}_{20}
$,
$\frac{2}{\sqrt{30}}c^{1}_{20}
$;\ \ 
$-\frac{1}{\sqrt{15}}c^{3}_{20}
$,
$-\frac{2}{\sqrt{30}}c^{1}_{20}
$,
$-\frac{2}{\sqrt{30}}c^{3}_{20}
$;\ \ 
$-\frac{1}{\sqrt{15}}c^{3}_{20}
$,
$\frac{1}{\sqrt{15}}c^{1}_{20}
$;\ \ 
$\frac{1}{\sqrt{15}}c^{3}_{20}
$)
 $\oplus$
($-\frac{1}{\sqrt{5}}c^{3}_{20}
$,
$\frac{1}{\sqrt{5}}c^{1}_{20}
$;
$\frac{1}{\sqrt{5}}c^{3}_{20}
$)

Pass. 

 \ \color{black}

 \color{blue}

\noindent 491: (dims,levels) = $(4 , 
2;60,
60
)$,
irreps = $4_{5,2}^{1}
\hskip -1.5pt \otimes \hskip -1.5pt
1_{4}^{1,0}
\hskip -1.5pt \otimes \hskip -1.5pt
1_{3}^{1,0}\oplus
2_{5}^{2}
\hskip -1.5pt \otimes \hskip -1.5pt
1_{4}^{1,0}
\hskip -1.5pt \otimes \hskip -1.5pt
1_{3}^{1,0}$,
pord$(\rho_\text{isum}(\mathfrak{t})) = 5$,

\vskip 0.7ex
\hangindent=5.5em \hangafter=1
{\white .}\hskip 1em $\rho_\text{isum}(\mathfrak{t})$ =
 $( \frac{11}{60},
\frac{23}{60},
\frac{47}{60},
\frac{59}{60} )
\oplus
( \frac{11}{60},
\frac{59}{60} )
$,

\vskip 0.7ex
\hangindent=5.5em \hangafter=1
{\white .}\hskip 1em $\rho_\text{isum}(\mathfrak{s})$ =
$\mathrm{i}$($-\sqrt{\frac{1}{5}}$,
$\frac{-5+\sqrt{5}}{10}$,
$-\frac{5+\sqrt{5}}{10}$,
$\sqrt{\frac{1}{5}}$;\ \ 
$\sqrt{\frac{1}{5}}$,
$-\sqrt{\frac{1}{5}}$,
$-\frac{5+\sqrt{5}}{10}$;\ \ 
$\sqrt{\frac{1}{5}}$,
$\frac{-5+\sqrt{5}}{10}$;\ \ 
$-\sqrt{\frac{1}{5}}$)
 $\oplus$
($-\frac{1}{\sqrt{5}}c^{1}_{20}
$,
$\frac{1}{\sqrt{5}}c^{3}_{20}
$;
$\frac{1}{\sqrt{5}}c^{1}_{20}
$)

Pass. 

 \ \color{black}

 \color{blue}

\noindent 492: (dims,levels) = $(4 , 
2;60,
60
)$,
irreps = $4_{5,1}^{1}
\hskip -1.5pt \otimes \hskip -1.5pt
1_{4}^{1,0}
\hskip -1.5pt \otimes \hskip -1.5pt
1_{3}^{1,0}\oplus
2_{5}^{2}
\hskip -1.5pt \otimes \hskip -1.5pt
1_{4}^{1,0}
\hskip -1.5pt \otimes \hskip -1.5pt
1_{3}^{1,0}$,
pord$(\rho_\text{isum}(\mathfrak{t})) = 5$,

\vskip 0.7ex
\hangindent=5.5em \hangafter=1
{\white .}\hskip 1em $\rho_\text{isum}(\mathfrak{t})$ =
 $( \frac{11}{60},
\frac{23}{60},
\frac{47}{60},
\frac{59}{60} )
\oplus
( \frac{11}{60},
\frac{59}{60} )
$,

\vskip 0.7ex
\hangindent=5.5em \hangafter=1
{\white .}\hskip 1em $\rho_\text{isum}(\mathfrak{s})$ =
($-\frac{1}{5}c^{1}_{20}
+\frac{1}{5}c^{3}_{20}
$,
$\frac{2}{5}c^{2}_{15}
+\frac{1}{5}c^{3}_{15}
$,
$-\frac{1}{5}+\frac{2}{5}c^{1}_{15}
-\frac{1}{5}c^{3}_{15}
$,
$\frac{1}{5}c^{1}_{20}
+\frac{1}{5}c^{3}_{20}
$;
$\frac{1}{5}c^{1}_{20}
+\frac{1}{5}c^{3}_{20}
$,
$-\frac{1}{5}c^{1}_{20}
+\frac{1}{5}c^{3}_{20}
$,
$\frac{1}{5}-\frac{2}{5}c^{1}_{15}
+\frac{1}{5}c^{3}_{15}
$;
$-\frac{1}{5}c^{1}_{20}
-\frac{1}{5}c^{3}_{20}
$,
$\frac{2}{5}c^{2}_{15}
+\frac{1}{5}c^{3}_{15}
$;
$\frac{1}{5}c^{1}_{20}
-\frac{1}{5}c^{3}_{20}
$)
 $\oplus$
($-\frac{1}{\sqrt{5}}c^{1}_{20}
$,
$\frac{1}{\sqrt{5}}c^{3}_{20}
$;
$\frac{1}{\sqrt{5}}c^{1}_{20}
$)

Pass. 

 \ \color{black}

 \color{blue}

\noindent 493: (dims,levels) = $(4 , 
2;60,
60
)$,
irreps = $4_{5,2}^{1}
\hskip -1.5pt \otimes \hskip -1.5pt
1_{4}^{1,0}
\hskip -1.5pt \otimes \hskip -1.5pt
1_{3}^{1,0}\oplus
2_{5}^{1}
\hskip -1.5pt \otimes \hskip -1.5pt
1_{4}^{1,0}
\hskip -1.5pt \otimes \hskip -1.5pt
1_{3}^{1,0}$,
pord$(\rho_\text{isum}(\mathfrak{t})) = 5$,

\vskip 0.7ex
\hangindent=5.5em \hangafter=1
{\white .}\hskip 1em $\rho_\text{isum}(\mathfrak{t})$ =
 $( \frac{11}{60},
\frac{23}{60},
\frac{47}{60},
\frac{59}{60} )
\oplus
( \frac{23}{60},
\frac{47}{60} )
$,

\vskip 0.7ex
\hangindent=5.5em \hangafter=1
{\white .}\hskip 1em $\rho_\text{isum}(\mathfrak{s})$ =
$\mathrm{i}$($-\sqrt{\frac{1}{5}}$,
$\frac{-5+\sqrt{5}}{10}$,
$-\frac{5+\sqrt{5}}{10}$,
$\sqrt{\frac{1}{5}}$;\ \ 
$\sqrt{\frac{1}{5}}$,
$-\sqrt{\frac{1}{5}}$,
$-\frac{5+\sqrt{5}}{10}$;\ \ 
$\sqrt{\frac{1}{5}}$,
$\frac{-5+\sqrt{5}}{10}$;\ \ 
$-\sqrt{\frac{1}{5}}$)
 $\oplus$
($-\frac{1}{\sqrt{5}}c^{3}_{20}
$,
$\frac{1}{\sqrt{5}}c^{1}_{20}
$;
$\frac{1}{\sqrt{5}}c^{3}_{20}
$)

Pass. 

 \ \color{black}

 \color{blue}

\noindent 494: (dims,levels) = $(4 , 
2;60,
60
)$,
irreps = $4_{5,1}^{1}
\hskip -1.5pt \otimes \hskip -1.5pt
1_{4}^{1,0}
\hskip -1.5pt \otimes \hskip -1.5pt
1_{3}^{1,0}\oplus
2_{5}^{1}
\hskip -1.5pt \otimes \hskip -1.5pt
1_{4}^{1,0}
\hskip -1.5pt \otimes \hskip -1.5pt
1_{3}^{1,0}$,
pord$(\rho_\text{isum}(\mathfrak{t})) = 5$,

\vskip 0.7ex
\hangindent=5.5em \hangafter=1
{\white .}\hskip 1em $\rho_\text{isum}(\mathfrak{t})$ =
 $( \frac{11}{60},
\frac{23}{60},
\frac{47}{60},
\frac{59}{60} )
\oplus
( \frac{23}{60},
\frac{47}{60} )
$,

\vskip 0.7ex
\hangindent=5.5em \hangafter=1
{\white .}\hskip 1em $\rho_\text{isum}(\mathfrak{s})$ =
($-\frac{1}{5}c^{1}_{20}
+\frac{1}{5}c^{3}_{20}
$,
$\frac{2}{5}c^{2}_{15}
+\frac{1}{5}c^{3}_{15}
$,
$-\frac{1}{5}+\frac{2}{5}c^{1}_{15}
-\frac{1}{5}c^{3}_{15}
$,
$\frac{1}{5}c^{1}_{20}
+\frac{1}{5}c^{3}_{20}
$;
$\frac{1}{5}c^{1}_{20}
+\frac{1}{5}c^{3}_{20}
$,
$-\frac{1}{5}c^{1}_{20}
+\frac{1}{5}c^{3}_{20}
$,
$\frac{1}{5}-\frac{2}{5}c^{1}_{15}
+\frac{1}{5}c^{3}_{15}
$;
$-\frac{1}{5}c^{1}_{20}
-\frac{1}{5}c^{3}_{20}
$,
$\frac{2}{5}c^{2}_{15}
+\frac{1}{5}c^{3}_{15}
$;
$\frac{1}{5}c^{1}_{20}
-\frac{1}{5}c^{3}_{20}
$)
 $\oplus$
($-\frac{1}{\sqrt{5}}c^{3}_{20}
$,
$\frac{1}{\sqrt{5}}c^{1}_{20}
$;
$\frac{1}{\sqrt{5}}c^{3}_{20}
$)

Pass. 

 \ \color{black}

 \color{blue}

\noindent 495: (dims,levels) = $(4 , 
2;60,
60
)$,
irreps = $2_{5}^{2}
\hskip -1.5pt \otimes \hskip -1.5pt
2_{4}^{1,0}
\hskip -1.5pt \otimes \hskip -1.5pt
1_{3}^{1,0}\oplus
2_{5}^{2}
\hskip -1.5pt \otimes \hskip -1.5pt
1_{4}^{1,0}
\hskip -1.5pt \otimes \hskip -1.5pt
1_{3}^{1,0}$,
pord$(\rho_\text{isum}(\mathfrak{t})) = 10$,

\vskip 0.7ex
\hangindent=5.5em \hangafter=1
{\white .}\hskip 1em $\rho_\text{isum}(\mathfrak{t})$ =
 $( \frac{11}{60},
\frac{29}{60},
\frac{41}{60},
\frac{59}{60} )
\oplus
( \frac{11}{60},
\frac{59}{60} )
$,

\vskip 0.7ex
\hangindent=5.5em \hangafter=1
{\white .}\hskip 1em $\rho_\text{isum}(\mathfrak{s})$ =
($\frac{1}{2\sqrt{5}}c^{1}_{20}
$,
$-\frac{3}{2\sqrt{15}}c^{3}_{20}
$,
$\frac{3}{2\sqrt{15}}c^{1}_{20}
$,
$\frac{1}{2\sqrt{5}}c^{3}_{20}
$;
$\frac{1}{2\sqrt{5}}c^{1}_{20}
$,
$\frac{1}{2\sqrt{5}}c^{3}_{20}
$,
$\frac{3}{2\sqrt{15}}c^{1}_{20}
$;
$-\frac{1}{2\sqrt{5}}c^{1}_{20}
$,
$\frac{3}{2\sqrt{15}}c^{3}_{20}
$;
$-\frac{1}{2\sqrt{5}}c^{1}_{20}
$)
 $\oplus$
($-\frac{1}{\sqrt{5}}c^{1}_{20}
$,
$\frac{1}{\sqrt{5}}c^{3}_{20}
$;
$\frac{1}{\sqrt{5}}c^{1}_{20}
$)

Pass. 

 \ \color{black}

 \color{blue}

\noindent 496: (dims,levels) = $(4 , 
2;60,
60
)$,
irreps = $2_{5}^{2}
\hskip -1.5pt \otimes \hskip -1.5pt
2_{4}^{1,0}
\hskip -1.5pt \otimes \hskip -1.5pt
1_{3}^{1,0}\oplus
2_{5}^{2}
\hskip -1.5pt \otimes \hskip -1.5pt
1_{4}^{1,6}
\hskip -1.5pt \otimes \hskip -1.5pt
1_{3}^{1,0}$,
pord$(\rho_\text{isum}(\mathfrak{t})) = 10$,

\vskip 0.7ex
\hangindent=5.5em \hangafter=1
{\white .}\hskip 1em $\rho_\text{isum}(\mathfrak{t})$ =
 $( \frac{11}{60},
\frac{29}{60},
\frac{41}{60},
\frac{59}{60} )
\oplus
( \frac{29}{60},
\frac{41}{60} )
$,

\vskip 0.7ex
\hangindent=5.5em \hangafter=1
{\white .}\hskip 1em $\rho_\text{isum}(\mathfrak{s})$ =
($\frac{1}{2\sqrt{5}}c^{1}_{20}
$,
$-\frac{3}{2\sqrt{15}}c^{3}_{20}
$,
$\frac{3}{2\sqrt{15}}c^{1}_{20}
$,
$\frac{1}{2\sqrt{5}}c^{3}_{20}
$;
$\frac{1}{2\sqrt{5}}c^{1}_{20}
$,
$\frac{1}{2\sqrt{5}}c^{3}_{20}
$,
$\frac{3}{2\sqrt{15}}c^{1}_{20}
$;
$-\frac{1}{2\sqrt{5}}c^{1}_{20}
$,
$\frac{3}{2\sqrt{15}}c^{3}_{20}
$;
$-\frac{1}{2\sqrt{5}}c^{1}_{20}
$)
 $\oplus$
($-\frac{1}{\sqrt{5}}c^{1}_{20}
$,
$\frac{1}{\sqrt{5}}c^{3}_{20}
$;
$\frac{1}{\sqrt{5}}c^{1}_{20}
$)

Pass. 

 \ \color{black}

 \color{blue}

\noindent 497: (dims,levels) = $(4 , 
2;60,
60
)$,
irreps = $2_{5}^{1}
\hskip -1.5pt \otimes \hskip -1.5pt
2_{4}^{1,0}
\hskip -1.5pt \otimes \hskip -1.5pt
1_{3}^{1,0}\oplus
2_{5}^{1}
\hskip -1.5pt \otimes \hskip -1.5pt
1_{4}^{1,6}
\hskip -1.5pt \otimes \hskip -1.5pt
1_{3}^{1,0}$,
pord$(\rho_\text{isum}(\mathfrak{t})) = 10$,

\vskip 0.7ex
\hangindent=5.5em \hangafter=1
{\white .}\hskip 1em $\rho_\text{isum}(\mathfrak{t})$ =
 $( \frac{17}{60},
\frac{23}{60},
\frac{47}{60},
\frac{53}{60} )
\oplus
( \frac{17}{60},
\frac{53}{60} )
$,

\vskip 0.7ex
\hangindent=5.5em \hangafter=1
{\white .}\hskip 1em $\rho_\text{isum}(\mathfrak{s})$ =
($\frac{1}{2\sqrt{5}}c^{3}_{20}
$,
$\frac{3}{2\sqrt{15}}c^{1}_{20}
$,
$-\frac{3}{2\sqrt{15}}c^{3}_{20}
$,
$\frac{1}{2\sqrt{5}}c^{1}_{20}
$;
$\frac{1}{2\sqrt{5}}c^{3}_{20}
$,
$\frac{1}{2\sqrt{5}}c^{1}_{20}
$,
$-\frac{3}{2\sqrt{15}}c^{3}_{20}
$;
$-\frac{1}{2\sqrt{5}}c^{3}_{20}
$,
$-\frac{3}{2\sqrt{15}}c^{1}_{20}
$;
$-\frac{1}{2\sqrt{5}}c^{3}_{20}
$)
 $\oplus$
($-\frac{1}{\sqrt{5}}c^{3}_{20}
$,
$\frac{1}{\sqrt{5}}c^{1}_{20}
$;
$\frac{1}{\sqrt{5}}c^{3}_{20}
$)

Pass. 

 \ \color{black}

 \color{blue}

\noindent 498: (dims,levels) = $(4 , 
2;60,
60
)$,
irreps = $2_{5}^{1}
\hskip -1.5pt \otimes \hskip -1.5pt
2_{4}^{1,0}
\hskip -1.5pt \otimes \hskip -1.5pt
1_{3}^{1,0}\oplus
2_{5}^{1}
\hskip -1.5pt \otimes \hskip -1.5pt
1_{4}^{1,0}
\hskip -1.5pt \otimes \hskip -1.5pt
1_{3}^{1,0}$,
pord$(\rho_\text{isum}(\mathfrak{t})) = 10$,

\vskip 0.7ex
\hangindent=5.5em \hangafter=1
{\white .}\hskip 1em $\rho_\text{isum}(\mathfrak{t})$ =
 $( \frac{17}{60},
\frac{23}{60},
\frac{47}{60},
\frac{53}{60} )
\oplus
( \frac{23}{60},
\frac{47}{60} )
$,

\vskip 0.7ex
\hangindent=5.5em \hangafter=1
{\white .}\hskip 1em $\rho_\text{isum}(\mathfrak{s})$ =
($\frac{1}{2\sqrt{5}}c^{3}_{20}
$,
$\frac{3}{2\sqrt{15}}c^{1}_{20}
$,
$-\frac{3}{2\sqrt{15}}c^{3}_{20}
$,
$\frac{1}{2\sqrt{5}}c^{1}_{20}
$;
$\frac{1}{2\sqrt{5}}c^{3}_{20}
$,
$\frac{1}{2\sqrt{5}}c^{1}_{20}
$,
$-\frac{3}{2\sqrt{15}}c^{3}_{20}
$;
$-\frac{1}{2\sqrt{5}}c^{3}_{20}
$,
$-\frac{3}{2\sqrt{15}}c^{1}_{20}
$;
$-\frac{1}{2\sqrt{5}}c^{3}_{20}
$)
 $\oplus$
($-\frac{1}{\sqrt{5}}c^{3}_{20}
$,
$\frac{1}{\sqrt{5}}c^{1}_{20}
$;
$\frac{1}{\sqrt{5}}c^{3}_{20}
$)

Pass. 

 \ \color{black}

 \color{blue}

\noindent 499: (dims,levels) = $(4 , 
2;60,
60
)$,
irreps = $2_{5}^{2}
\hskip -1.5pt \otimes \hskip -1.5pt
2_{3}^{1,0}
\hskip -1.5pt \otimes \hskip -1.5pt
1_{4}^{1,0}\oplus
2_{5}^{2}
\hskip -1.5pt \otimes \hskip -1.5pt
1_{4}^{1,0}
\hskip -1.5pt \otimes \hskip -1.5pt
1_{3}^{1,0}$,
pord$(\rho_\text{isum}(\mathfrak{t})) = 15$,

\vskip 0.7ex
\hangindent=5.5em \hangafter=1
{\white .}\hskip 1em $\rho_\text{isum}(\mathfrak{t})$ =
 $( \frac{13}{20},
\frac{17}{20},
\frac{11}{60},
\frac{59}{60} )
\oplus
( \frac{11}{60},
\frac{59}{60} )
$,

\vskip 0.7ex
\hangindent=5.5em \hangafter=1
{\white .}\hskip 1em $\rho_\text{isum}(\mathfrak{s})$ =
$\mathrm{i}$($-\frac{1}{\sqrt{15}}c^{1}_{20}
$,
$\frac{1}{\sqrt{15}}c^{3}_{20}
$,
$-\frac{2}{\sqrt{30}}c^{3}_{20}
$,
$\frac{2}{\sqrt{30}}c^{1}_{20}
$;\ \ 
$\frac{1}{\sqrt{15}}c^{1}_{20}
$,
$-\frac{2}{\sqrt{30}}c^{1}_{20}
$,
$-\frac{2}{\sqrt{30}}c^{3}_{20}
$;\ \ 
$-\frac{1}{\sqrt{15}}c^{1}_{20}
$,
$-\frac{1}{\sqrt{15}}c^{3}_{20}
$;\ \ 
$\frac{1}{\sqrt{15}}c^{1}_{20}
$)
 $\oplus$
($-\frac{1}{\sqrt{5}}c^{1}_{20}
$,
$\frac{1}{\sqrt{5}}c^{3}_{20}
$;
$\frac{1}{\sqrt{5}}c^{1}_{20}
$)

Pass. 

 \ \color{black}

\noindent 500: (dims,levels) = $(4 , 
2;84,
12
)$,
irreps = $4_{7}^{3}
\hskip -1.5pt \otimes \hskip -1.5pt
1_{4}^{1,0}
\hskip -1.5pt \otimes \hskip -1.5pt
1_{3}^{1,0}\oplus
2_{4}^{1,0}
\hskip -1.5pt \otimes \hskip -1.5pt
1_{3}^{1,0}$,
pord$(\rho_\text{isum}(\mathfrak{t})) = 14$,

\vskip 0.7ex
\hangindent=5.5em \hangafter=1
{\white .}\hskip 1em $\rho_\text{isum}(\mathfrak{t})$ =
 $( \frac{7}{12},
\frac{1}{84},
\frac{25}{84},
\frac{37}{84} )
\oplus
( \frac{1}{12},
\frac{7}{12} )
$,

\vskip 0.7ex
\hangindent=5.5em \hangafter=1
{\white .}\hskip 1em $\rho_\text{isum}(\mathfrak{s})$ =
($-\sqrt{\frac{1}{7}}$,
$\sqrt{\frac{2}{7}}$,
$\sqrt{\frac{2}{7}}$,
$\sqrt{\frac{2}{7}}$;
$-\frac{1}{\sqrt{7}}c^{1}_{7}
$,
$-\frac{1}{\sqrt{7}}c^{2}_{7}
$,
$\frac{1}{\sqrt{7}\mathrm{i}}s^{5}_{28}
$;
$\frac{1}{\sqrt{7}\mathrm{i}}s^{5}_{28}
$,
$-\frac{1}{\sqrt{7}}c^{1}_{7}
$;
$-\frac{1}{\sqrt{7}}c^{2}_{7}
$)
 $\oplus$
$\mathrm{i}$($\frac{1}{2}$,
$\sqrt{\frac{3}{4}}$;\ \ 
$-\frac{1}{2}$)

Fail:
Tr$_I(C) = -1 <$  0 for I = [ 1/12 ]. Prop. B.4 (1) eqn. (B.18)

 \ \color{black}

 \color{blue}

\noindent 501: (dims,levels) = $(4 , 
2;84,
12
)$,
irreps = $4_{7}^{3}
\hskip -1.5pt \otimes \hskip -1.5pt
1_{4}^{1,0}
\hskip -1.5pt \otimes \hskip -1.5pt
1_{3}^{1,0}\oplus
2_{3}^{1,0}
\hskip -1.5pt \otimes \hskip -1.5pt
1_{4}^{1,0}$,
pord$(\rho_\text{isum}(\mathfrak{t})) = 21$,

\vskip 0.7ex
\hangindent=5.5em \hangafter=1
{\white .}\hskip 1em $\rho_\text{isum}(\mathfrak{t})$ =
 $( \frac{7}{12},
\frac{1}{84},
\frac{25}{84},
\frac{37}{84} )
\oplus
( \frac{1}{4},
\frac{7}{12} )
$,

\vskip 0.7ex
\hangindent=5.5em \hangafter=1
{\white .}\hskip 1em $\rho_\text{isum}(\mathfrak{s})$ =
($-\sqrt{\frac{1}{7}}$,
$\sqrt{\frac{2}{7}}$,
$\sqrt{\frac{2}{7}}$,
$\sqrt{\frac{2}{7}}$;
$-\frac{1}{\sqrt{7}}c^{1}_{7}
$,
$-\frac{1}{\sqrt{7}}c^{2}_{7}
$,
$\frac{1}{\sqrt{7}\mathrm{i}}s^{5}_{28}
$;
$\frac{1}{\sqrt{7}\mathrm{i}}s^{5}_{28}
$,
$-\frac{1}{\sqrt{7}}c^{1}_{7}
$;
$-\frac{1}{\sqrt{7}}c^{2}_{7}
$)
 $\oplus$
($\sqrt{\frac{1}{3}}$,
$\sqrt{\frac{2}{3}}$;
$-\sqrt{\frac{1}{3}}$)

Pass. 

 \ \color{black}

 \color{blue}

\noindent 502: (dims,levels) = $(4 , 
2;84,
12
)$,
irreps = $4_{7}^{3}
\hskip -1.5pt \otimes \hskip -1.5pt
1_{4}^{1,0}
\hskip -1.5pt \otimes \hskip -1.5pt
1_{3}^{1,0}\oplus
2_{3}^{1,4}
\hskip -1.5pt \otimes \hskip -1.5pt
1_{4}^{1,0}$,
pord$(\rho_\text{isum}(\mathfrak{t})) = 21$,

\vskip 0.7ex
\hangindent=5.5em \hangafter=1
{\white .}\hskip 1em $\rho_\text{isum}(\mathfrak{t})$ =
 $( \frac{7}{12},
\frac{1}{84},
\frac{25}{84},
\frac{37}{84} )
\oplus
( \frac{7}{12},
\frac{11}{12} )
$,

\vskip 0.7ex
\hangindent=5.5em \hangafter=1
{\white .}\hskip 1em $\rho_\text{isum}(\mathfrak{s})$ =
($-\sqrt{\frac{1}{7}}$,
$\sqrt{\frac{2}{7}}$,
$\sqrt{\frac{2}{7}}$,
$\sqrt{\frac{2}{7}}$;
$-\frac{1}{\sqrt{7}}c^{1}_{7}
$,
$-\frac{1}{\sqrt{7}}c^{2}_{7}
$,
$\frac{1}{\sqrt{7}\mathrm{i}}s^{5}_{28}
$;
$\frac{1}{\sqrt{7}\mathrm{i}}s^{5}_{28}
$,
$-\frac{1}{\sqrt{7}}c^{1}_{7}
$;
$-\frac{1}{\sqrt{7}}c^{2}_{7}
$)
 $\oplus$
($\sqrt{\frac{1}{3}}$,
$\sqrt{\frac{2}{3}}$;
$-\sqrt{\frac{1}{3}}$)

Pass. 

 \ \color{black}

\noindent 503: (dims,levels) = $(4 , 
2;84,
12
)$,
irreps = $4_{7}^{1}
\hskip -1.5pt \otimes \hskip -1.5pt
1_{4}^{1,0}
\hskip -1.5pt \otimes \hskip -1.5pt
1_{3}^{1,0}\oplus
2_{4}^{1,0}
\hskip -1.5pt \otimes \hskip -1.5pt
1_{3}^{1,0}$,
pord$(\rho_\text{isum}(\mathfrak{t})) = 14$,

\vskip 0.7ex
\hangindent=5.5em \hangafter=1
{\white .}\hskip 1em $\rho_\text{isum}(\mathfrak{t})$ =
 $( \frac{7}{12},
\frac{13}{84},
\frac{61}{84},
\frac{73}{84} )
\oplus
( \frac{1}{12},
\frac{7}{12} )
$,

\vskip 0.7ex
\hangindent=5.5em \hangafter=1
{\white .}\hskip 1em $\rho_\text{isum}(\mathfrak{s})$ =
($\sqrt{\frac{1}{7}}$,
$\sqrt{\frac{2}{7}}$,
$\sqrt{\frac{2}{7}}$,
$\sqrt{\frac{2}{7}}$;
$\frac{1}{\sqrt{7}}c^{1}_{7}
$,
$-\frac{1}{\sqrt{7}\mathrm{i}}s^{5}_{28}
$,
$\frac{1}{\sqrt{7}}c^{2}_{7}
$;
$\frac{1}{\sqrt{7}}c^{2}_{7}
$,
$\frac{1}{\sqrt{7}}c^{1}_{7}
$;
$-\frac{1}{\sqrt{7}\mathrm{i}}s^{5}_{28}
$)
 $\oplus$
$\mathrm{i}$($\frac{1}{2}$,
$\sqrt{\frac{3}{4}}$;\ \ 
$-\frac{1}{2}$)

Fail:
Tr$_I(C) = -1 <$  0 for I = [ 1/12 ]. Prop. B.4 (1) eqn. (B.18)

 \ \color{black}

 \color{blue}

\noindent 504: (dims,levels) = $(4 , 
2;84,
12
)$,
irreps = $4_{7}^{1}
\hskip -1.5pt \otimes \hskip -1.5pt
1_{4}^{1,0}
\hskip -1.5pt \otimes \hskip -1.5pt
1_{3}^{1,0}\oplus
2_{3}^{1,0}
\hskip -1.5pt \otimes \hskip -1.5pt
1_{4}^{1,0}$,
pord$(\rho_\text{isum}(\mathfrak{t})) = 21$,

\vskip 0.7ex
\hangindent=5.5em \hangafter=1
{\white .}\hskip 1em $\rho_\text{isum}(\mathfrak{t})$ =
 $( \frac{7}{12},
\frac{13}{84},
\frac{61}{84},
\frac{73}{84} )
\oplus
( \frac{1}{4},
\frac{7}{12} )
$,

\vskip 0.7ex
\hangindent=5.5em \hangafter=1
{\white .}\hskip 1em $\rho_\text{isum}(\mathfrak{s})$ =
($\sqrt{\frac{1}{7}}$,
$\sqrt{\frac{2}{7}}$,
$\sqrt{\frac{2}{7}}$,
$\sqrt{\frac{2}{7}}$;
$\frac{1}{\sqrt{7}}c^{1}_{7}
$,
$-\frac{1}{\sqrt{7}\mathrm{i}}s^{5}_{28}
$,
$\frac{1}{\sqrt{7}}c^{2}_{7}
$;
$\frac{1}{\sqrt{7}}c^{2}_{7}
$,
$\frac{1}{\sqrt{7}}c^{1}_{7}
$;
$-\frac{1}{\sqrt{7}\mathrm{i}}s^{5}_{28}
$)
 $\oplus$
($\sqrt{\frac{1}{3}}$,
$\sqrt{\frac{2}{3}}$;
$-\sqrt{\frac{1}{3}}$)

Pass. 

 \ \color{black}

 \color{blue}

\noindent 505: (dims,levels) = $(4 , 
2;84,
12
)$,
irreps = $4_{7}^{1}
\hskip -1.5pt \otimes \hskip -1.5pt
1_{4}^{1,0}
\hskip -1.5pt \otimes \hskip -1.5pt
1_{3}^{1,0}\oplus
2_{3}^{1,4}
\hskip -1.5pt \otimes \hskip -1.5pt
1_{4}^{1,0}$,
pord$(\rho_\text{isum}(\mathfrak{t})) = 21$,

\vskip 0.7ex
\hangindent=5.5em \hangafter=1
{\white .}\hskip 1em $\rho_\text{isum}(\mathfrak{t})$ =
 $( \frac{7}{12},
\frac{13}{84},
\frac{61}{84},
\frac{73}{84} )
\oplus
( \frac{7}{12},
\frac{11}{12} )
$,

\vskip 0.7ex
\hangindent=5.5em \hangafter=1
{\white .}\hskip 1em $\rho_\text{isum}(\mathfrak{s})$ =
($\sqrt{\frac{1}{7}}$,
$\sqrt{\frac{2}{7}}$,
$\sqrt{\frac{2}{7}}$,
$\sqrt{\frac{2}{7}}$;
$\frac{1}{\sqrt{7}}c^{1}_{7}
$,
$-\frac{1}{\sqrt{7}\mathrm{i}}s^{5}_{28}
$,
$\frac{1}{\sqrt{7}}c^{2}_{7}
$;
$\frac{1}{\sqrt{7}}c^{2}_{7}
$,
$\frac{1}{\sqrt{7}}c^{1}_{7}
$;
$-\frac{1}{\sqrt{7}\mathrm{i}}s^{5}_{28}
$)
 $\oplus$
($\sqrt{\frac{1}{3}}$,
$\sqrt{\frac{2}{3}}$;
$-\sqrt{\frac{1}{3}}$)

Pass. 

 \ \color{black}

 \color{blue}

\noindent 506: (dims,levels) = $(5 , 
1;5,
1
)$,
irreps = $5_{5}^{1}\oplus
1_{1}^{1}$,
pord$(\rho_\text{isum}(\mathfrak{t})) = 5$,

\vskip 0.7ex
\hangindent=5.5em \hangafter=1
{\white .}\hskip 1em $\rho_\text{isum}(\mathfrak{t})$ =
 $( 0,
\frac{1}{5},
\frac{2}{5},
\frac{3}{5},
\frac{4}{5} )
\oplus
( 0 )
$,

\vskip 0.7ex
\hangindent=5.5em \hangafter=1
{\white .}\hskip 1em $\rho_\text{isum}(\mathfrak{s})$ =
($-\frac{1}{5}$,
$\sqrt{\frac{6}{25}}$,
$\sqrt{\frac{6}{25}}$,
$\sqrt{\frac{6}{25}}$,
$\sqrt{\frac{6}{25}}$;
$\frac{3-\sqrt{5}}{10}$,
$-\frac{1+\sqrt{5}}{5}$,
$\frac{-1+\sqrt{5}}{5}$,
$\frac{3+\sqrt{5}}{10}$;
$\frac{3+\sqrt{5}}{10}$,
$\frac{3-\sqrt{5}}{10}$,
$\frac{-1+\sqrt{5}}{5}$;
$\frac{3+\sqrt{5}}{10}$,
$-\frac{1+\sqrt{5}}{5}$;
$\frac{3-\sqrt{5}}{10}$)
 $\oplus$
($1$)

Pass. 

 \ \color{black}

 \color{blue}

\noindent 507: (dims,levels) = $(5 , 
1;10,
2
)$,
irreps = $5_{5}^{1}
\hskip -1.5pt \otimes \hskip -1.5pt
1_{2}^{1,0}\oplus
1_{2}^{1,0}$,
pord$(\rho_\text{isum}(\mathfrak{t})) = 5$,

\vskip 0.7ex
\hangindent=5.5em \hangafter=1
{\white .}\hskip 1em $\rho_\text{isum}(\mathfrak{t})$ =
 $( \frac{1}{2},
\frac{1}{10},
\frac{3}{10},
\frac{7}{10},
\frac{9}{10} )
\oplus
( \frac{1}{2} )
$,

\vskip 0.7ex
\hangindent=5.5em \hangafter=1
{\white .}\hskip 1em $\rho_\text{isum}(\mathfrak{s})$ =
($\frac{1}{5}$,
$\sqrt{\frac{6}{25}}$,
$\sqrt{\frac{6}{25}}$,
$\sqrt{\frac{6}{25}}$,
$\sqrt{\frac{6}{25}}$;
$-\frac{3+\sqrt{5}}{10}$,
$\frac{1+\sqrt{5}}{5}$,
$\frac{1-\sqrt{5}}{5}$,
$\frac{-3+\sqrt{5}}{10}$;
$\frac{-3+\sqrt{5}}{10}$,
$-\frac{3+\sqrt{5}}{10}$,
$\frac{1-\sqrt{5}}{5}$;
$\frac{-3+\sqrt{5}}{10}$,
$\frac{1+\sqrt{5}}{5}$;
$-\frac{3+\sqrt{5}}{10}$)
 $\oplus$
($-1$)

Pass. 

 \ \color{black}

 \color{blue}

\noindent 508: (dims,levels) = $(5 , 
1;15,
3
)$,
irreps = $5_{5}^{1}
\hskip -1.5pt \otimes \hskip -1.5pt
1_{3}^{1,0}\oplus
1_{3}^{1,0}$,
pord$(\rho_\text{isum}(\mathfrak{t})) = 5$,

\vskip 0.7ex
\hangindent=5.5em \hangafter=1
{\white .}\hskip 1em $\rho_\text{isum}(\mathfrak{t})$ =
 $( \frac{1}{3},
\frac{2}{15},
\frac{8}{15},
\frac{11}{15},
\frac{14}{15} )
\oplus
( \frac{1}{3} )
$,

\vskip 0.7ex
\hangindent=5.5em \hangafter=1
{\white .}\hskip 1em $\rho_\text{isum}(\mathfrak{s})$ =
($-\frac{1}{5}$,
$\sqrt{\frac{6}{25}}$,
$\sqrt{\frac{6}{25}}$,
$\sqrt{\frac{6}{25}}$,
$\sqrt{\frac{6}{25}}$;
$\frac{3-\sqrt{5}}{10}$,
$\frac{3+\sqrt{5}}{10}$,
$\frac{-1+\sqrt{5}}{5}$,
$-\frac{1+\sqrt{5}}{5}$;
$\frac{3-\sqrt{5}}{10}$,
$-\frac{1+\sqrt{5}}{5}$,
$\frac{-1+\sqrt{5}}{5}$;
$\frac{3+\sqrt{5}}{10}$,
$\frac{3-\sqrt{5}}{10}$;
$\frac{3+\sqrt{5}}{10}$)
 $\oplus$
($1$)

Pass. 

 \ \color{black}

 \color{blue}

\noindent 509: (dims,levels) = $(5 , 
1;20,
4
)$,
irreps = $5_{5}^{1}
\hskip -1.5pt \otimes \hskip -1.5pt
1_{4}^{1,0}\oplus
1_{4}^{1,0}$,
pord$(\rho_\text{isum}(\mathfrak{t})) = 5$,

\vskip 0.7ex
\hangindent=5.5em \hangafter=1
{\white .}\hskip 1em $\rho_\text{isum}(\mathfrak{t})$ =
 $( \frac{1}{4},
\frac{1}{20},
\frac{9}{20},
\frac{13}{20},
\frac{17}{20} )
\oplus
( \frac{1}{4} )
$,

\vskip 0.7ex
\hangindent=5.5em \hangafter=1
{\white .}\hskip 1em $\rho_\text{isum}(\mathfrak{s})$ =
$\mathrm{i}$($-\frac{1}{5}$,
$\sqrt{\frac{6}{25}}$,
$\sqrt{\frac{6}{25}}$,
$\sqrt{\frac{6}{25}}$,
$\sqrt{\frac{6}{25}}$;\ \ 
$\frac{3-\sqrt{5}}{10}$,
$\frac{3+\sqrt{5}}{10}$,
$\frac{-1+\sqrt{5}}{5}$,
$-\frac{1+\sqrt{5}}{5}$;\ \ 
$\frac{3-\sqrt{5}}{10}$,
$-\frac{1+\sqrt{5}}{5}$,
$\frac{-1+\sqrt{5}}{5}$;\ \ 
$\frac{3+\sqrt{5}}{10}$,
$\frac{3-\sqrt{5}}{10}$;\ \ 
$\frac{3+\sqrt{5}}{10}$)
 $\oplus$
$\mathrm{i}$($1$)

Pass. 

 \ \color{black}

 \color{blue}

\noindent 510: (dims,levels) = $(5 , 
1;30,
6
)$,
irreps = $5_{5}^{1}
\hskip -1.5pt \otimes \hskip -1.5pt
1_{3}^{1,0}
\hskip -1.5pt \otimes \hskip -1.5pt
1_{2}^{1,0}\oplus
1_{3}^{1,0}
\hskip -1.5pt \otimes \hskip -1.5pt
1_{2}^{1,0}$,
pord$(\rho_\text{isum}(\mathfrak{t})) = 5$,

\vskip 0.7ex
\hangindent=5.5em \hangafter=1
{\white .}\hskip 1em $\rho_\text{isum}(\mathfrak{t})$ =
 $( \frac{5}{6},
\frac{1}{30},
\frac{7}{30},
\frac{13}{30},
\frac{19}{30} )
\oplus
( \frac{5}{6} )
$,

\vskip 0.7ex
\hangindent=5.5em \hangafter=1
{\white .}\hskip 1em $\rho_\text{isum}(\mathfrak{s})$ =
($\frac{1}{5}$,
$\sqrt{\frac{6}{25}}$,
$\sqrt{\frac{6}{25}}$,
$\sqrt{\frac{6}{25}}$,
$\sqrt{\frac{6}{25}}$;
$\frac{-3+\sqrt{5}}{10}$,
$\frac{1+\sqrt{5}}{5}$,
$\frac{1-\sqrt{5}}{5}$,
$-\frac{3+\sqrt{5}}{10}$;
$-\frac{3+\sqrt{5}}{10}$,
$\frac{-3+\sqrt{5}}{10}$,
$\frac{1-\sqrt{5}}{5}$;
$-\frac{3+\sqrt{5}}{10}$,
$\frac{1+\sqrt{5}}{5}$;
$\frac{-3+\sqrt{5}}{10}$)
 $\oplus$
($-1$)

Pass. 

 \ \color{black}

 \color{blue}

\noindent 511: (dims,levels) = $(5 , 
1;60,
12
)$,
irreps = $5_{5}^{1}
\hskip -1.5pt \otimes \hskip -1.5pt
1_{4}^{1,0}
\hskip -1.5pt \otimes \hskip -1.5pt
1_{3}^{1,0}\oplus
1_{4}^{1,0}
\hskip -1.5pt \otimes \hskip -1.5pt
1_{3}^{1,0}$,
pord$(\rho_\text{isum}(\mathfrak{t})) = 5$,

\vskip 0.7ex
\hangindent=5.5em \hangafter=1
{\white .}\hskip 1em $\rho_\text{isum}(\mathfrak{t})$ =
 $( \frac{7}{12},
\frac{11}{60},
\frac{23}{60},
\frac{47}{60},
\frac{59}{60} )
\oplus
( \frac{7}{12} )
$,

\vskip 0.7ex
\hangindent=5.5em \hangafter=1
{\white .}\hskip 1em $\rho_\text{isum}(\mathfrak{s})$ =
$\mathrm{i}$($-\frac{1}{5}$,
$\sqrt{\frac{6}{25}}$,
$\sqrt{\frac{6}{25}}$,
$\sqrt{\frac{6}{25}}$,
$\sqrt{\frac{6}{25}}$;\ \ 
$\frac{3+\sqrt{5}}{10}$,
$-\frac{1+\sqrt{5}}{5}$,
$\frac{-1+\sqrt{5}}{5}$,
$\frac{3-\sqrt{5}}{10}$;\ \ 
$\frac{3-\sqrt{5}}{10}$,
$\frac{3+\sqrt{5}}{10}$,
$\frac{-1+\sqrt{5}}{5}$;\ \ 
$\frac{3-\sqrt{5}}{10}$,
$-\frac{1+\sqrt{5}}{5}$;\ \ 
$\frac{3+\sqrt{5}}{10}$)
 $\oplus$
$\mathrm{i}$($1$)

Pass. 

 \ \color{black}

\noindent 512: (dims,levels) = $(6;5
)$,
irreps = $6_{5}^{1}$,
pord$(\rho_\text{isum}(\mathfrak{t})) = 5$,

\vskip 0.7ex
\hangindent=5.5em \hangafter=1
{\white .}\hskip 1em $\rho_\text{isum}(\mathfrak{t})$ =
 $( 0,
0,
\frac{1}{5},
\frac{2}{5},
\frac{3}{5},
\frac{4}{5} )
$,

\vskip 0.7ex
\hangindent=5.5em \hangafter=1
{\white .}\hskip 1em $\rho_\text{isum}(\mathfrak{s})$ =
$\mathrm{i}$($-\frac{1}{5}c^{1}_{20}
$,
$\frac{1}{5}c^{3}_{20}
$,
$\frac{1}{5}c^{3}_{40}
+\frac{1}{5}c^{7}_{40}
$,
$\frac{1}{5}c^{1}_{20}
+\frac{1}{5}c^{3}_{20}
$,
$\frac{1}{5}c^{1}_{20}
-\frac{1}{5}c^{3}_{20}
$,
$\frac{4}{5\sqrt{10}}c^{1}_{20}
-\frac{2}{5\sqrt{10}}c^{3}_{20}
$;\ \ 
$\frac{1}{5}c^{1}_{20}
$,
$-\frac{4}{5\sqrt{10}}c^{1}_{20}
+\frac{2}{5\sqrt{10}}c^{3}_{20}
$,
$\frac{1}{5}c^{1}_{20}
-\frac{1}{5}c^{3}_{20}
$,
$-\frac{1}{5}c^{1}_{20}
-\frac{1}{5}c^{3}_{20}
$,
$\frac{1}{5}c^{3}_{40}
+\frac{1}{5}c^{7}_{40}
$;\ \ 
$-\frac{1}{5}c^{3}_{20}
$,
$\frac{4}{5\sqrt{10}}c^{1}_{20}
-\frac{2}{5\sqrt{10}}c^{3}_{20}
$,
$\frac{1}{5}c^{3}_{40}
+\frac{1}{5}c^{7}_{40}
$,
$\frac{1}{5}c^{1}_{20}
$;\ \ 
$\frac{1}{5}c^{1}_{20}
$,
$-\frac{1}{5}c^{3}_{20}
$,
$-\frac{1}{5}c^{3}_{40}
-\frac{1}{5}c^{7}_{40}
$;\ \ 
$-\frac{1}{5}c^{1}_{20}
$,
$\frac{4}{5\sqrt{10}}c^{1}_{20}
-\frac{2}{5\sqrt{10}}c^{3}_{20}
$;\ \ 
$\frac{1}{5}c^{3}_{20}
$)

Fail:
cnd($\rho(\mathfrak s)_\mathrm{ndeg}$) = 40 does not divide
 ord($\rho(\mathfrak t)$)=5. Prop. B.4 (2)

 \ \color{black}

\noindent 513: (dims,levels) = $(6;7
)$,
irreps = $6_{7,2}^{1}$,
pord$(\rho_\text{isum}(\mathfrak{t})) = 7$,

\vskip 0.7ex
\hangindent=5.5em \hangafter=1
{\white .}\hskip 1em $\rho_\text{isum}(\mathfrak{t})$ =
 $( \frac{1}{7},
\frac{2}{7},
\frac{3}{7},
\frac{4}{7},
\frac{5}{7},
\frac{6}{7} )
$,

\vskip 0.7ex
\hangindent=5.5em \hangafter=1
{\white .}\hskip 1em $\rho_\text{isum}(\mathfrak{s})$ =
($\frac{2}{7}-\frac{1}{7}c^{2}_{7}
$,
$-\frac{2}{7}+\frac{1}{7}c^{1}_{7}
$,
$\frac{1}{7}c^{3}_{56}
+\frac{2}{7}c^{5}_{56}
-\frac{1}{7}c^{7}_{56}
-\frac{2}{7}c^{9}_{56}
+\frac{1}{7}c^{11}_{56}
$,
$-\frac{3}{7}-\frac{1}{7}c^{1}_{7}
-\frac{1}{7}c^{2}_{7}
$,
$-\frac{1}{7}c^{3}_{56}
+\frac{1}{7}c^{5}_{56}
-\frac{1}{7}c^{9}_{56}
-\frac{1}{7}c^{11}_{56}
$,
$-\frac{2}{7}c^{3}_{56}
-\frac{1}{7}c^{5}_{56}
+\frac{1}{7}c^{7}_{56}
+\frac{1}{7}c^{9}_{56}
-\frac{2}{7}c^{11}_{56}
$;
$\frac{3}{7}+\frac{1}{7}c^{1}_{7}
+\frac{1}{7}c^{2}_{7}
$,
$\frac{2}{7}c^{3}_{56}
+\frac{1}{7}c^{5}_{56}
-\frac{1}{7}c^{7}_{56}
-\frac{1}{7}c^{9}_{56}
+\frac{2}{7}c^{11}_{56}
$,
$\frac{2}{7}-\frac{1}{7}c^{2}_{7}
$,
$\frac{1}{7}c^{3}_{56}
+\frac{2}{7}c^{5}_{56}
-\frac{1}{7}c^{7}_{56}
-\frac{2}{7}c^{9}_{56}
+\frac{1}{7}c^{11}_{56}
$,
$-\frac{1}{7}c^{3}_{56}
+\frac{1}{7}c^{5}_{56}
-\frac{1}{7}c^{9}_{56}
-\frac{1}{7}c^{11}_{56}
$;
$\frac{2}{7}-\frac{1}{7}c^{1}_{7}
$,
$-\frac{1}{7}c^{3}_{56}
+\frac{1}{7}c^{5}_{56}
-\frac{1}{7}c^{9}_{56}
-\frac{1}{7}c^{11}_{56}
$,
$-\frac{2}{7}+\frac{1}{7}c^{2}_{7}
$,
$\frac{3}{7}+\frac{1}{7}c^{1}_{7}
+\frac{1}{7}c^{2}_{7}
$;
$\frac{2}{7}-\frac{1}{7}c^{1}_{7}
$,
$-\frac{2}{7}c^{3}_{56}
-\frac{1}{7}c^{5}_{56}
+\frac{1}{7}c^{7}_{56}
+\frac{1}{7}c^{9}_{56}
-\frac{2}{7}c^{11}_{56}
$,
$-\frac{1}{7}c^{3}_{56}
-\frac{2}{7}c^{5}_{56}
+\frac{1}{7}c^{7}_{56}
+\frac{2}{7}c^{9}_{56}
-\frac{1}{7}c^{11}_{56}
$;
$\frac{3}{7}+\frac{1}{7}c^{1}_{7}
+\frac{1}{7}c^{2}_{7}
$,
$-\frac{2}{7}+\frac{1}{7}c^{1}_{7}
$;
$\frac{2}{7}-\frac{1}{7}c^{2}_{7}
$)

Fail:
cnd($\rho(\mathfrak s)_\mathrm{ndeg}$) = 56 does not divide
 ord($\rho(\mathfrak t)$)=7. Prop. B.4 (2)

 \ \color{black}

\noindent 514: (dims,levels) = $(6;7
)$,
irreps = $6_{7,1}^{1}$,
pord$(\rho_\text{isum}(\mathfrak{t})) = 7$,

\vskip 0.7ex
\hangindent=5.5em \hangafter=1
{\white .}\hskip 1em $\rho_\text{isum}(\mathfrak{t})$ =
 $( \frac{1}{7},
\frac{2}{7},
\frac{3}{7},
\frac{4}{7},
\frac{5}{7},
\frac{6}{7} )
$,

\vskip 0.7ex
\hangindent=5.5em \hangafter=1
{\white .}\hskip 1em $\rho_\text{isum}(\mathfrak{s})$ =
$\mathrm{i}$($\frac{1}{7}c^{2}_{56}
-\frac{1}{7}c^{3}_{56}
+\frac{1}{7}c^{11}_{56}
$,
$\frac{1}{7}c^{5}_{56}
+\frac{1}{7}c^{6}_{56}
+\frac{1}{7}c^{9}_{56}
$,
$\frac{1}{7}c^{3}_{112}
-\frac{1}{7}c^{9}_{112}
+\frac{1}{7}c^{11}_{112}
+\frac{1}{7}c^{23}_{112}
$,
$\frac{2}{7}c^{1}_{56}
-\frac{1}{7}c^{3}_{56}
-\frac{1}{7}c^{5}_{56}
+\frac{1}{7}c^{7}_{56}
+\frac{1}{7}c^{9}_{56}
-\frac{1}{7}c^{10}_{56}
-\frac{1}{7}c^{11}_{56}
$,
$\frac{1}{7}c^{1}_{112}
+\frac{1}{7}c^{3}_{112}
-\frac{1}{7}c^{5}_{112}
-\frac{1}{7}c^{7}_{112}
-\frac{1}{7}c^{9}_{112}
-\frac{1}{7}c^{11}_{112}
+\frac{2}{7}c^{13}_{112}
+\frac{1}{7}c^{15}_{112}
+\frac{2}{7}c^{17}_{112}
+\frac{1}{7}c^{19}_{112}
-\frac{1}{7}c^{21}_{112}
-\frac{1}{7}c^{23}_{112}
$,
$\frac{1}{7}c^{1}_{112}
+\frac{1}{7}c^{5}_{112}
-\frac{1}{7}c^{15}_{112}
+\frac{1}{7}c^{19}_{112}
$;\ \ 
$\frac{2}{7}c^{1}_{56}
-\frac{1}{7}c^{3}_{56}
-\frac{1}{7}c^{5}_{56}
+\frac{1}{7}c^{7}_{56}
+\frac{1}{7}c^{9}_{56}
-\frac{1}{7}c^{10}_{56}
-\frac{1}{7}c^{11}_{56}
$,
$\frac{1}{7}c^{1}_{112}
+\frac{1}{7}c^{5}_{112}
-\frac{1}{7}c^{15}_{112}
+\frac{1}{7}c^{19}_{112}
$,
$\frac{1}{7}c^{2}_{56}
-\frac{1}{7}c^{3}_{56}
+\frac{1}{7}c^{11}_{56}
$,
$-\frac{1}{7}c^{3}_{112}
+\frac{1}{7}c^{9}_{112}
-\frac{1}{7}c^{11}_{112}
-\frac{1}{7}c^{23}_{112}
$,
$-\frac{1}{7}c^{1}_{112}
-\frac{1}{7}c^{3}_{112}
+\frac{1}{7}c^{5}_{112}
+\frac{1}{7}c^{7}_{112}
+\frac{1}{7}c^{9}_{112}
+\frac{1}{7}c^{11}_{112}
-\frac{2}{7}c^{13}_{112}
-\frac{1}{7}c^{15}_{112}
-\frac{2}{7}c^{17}_{112}
-\frac{1}{7}c^{19}_{112}
+\frac{1}{7}c^{21}_{112}
+\frac{1}{7}c^{23}_{112}
$;\ \ 
$-\frac{1}{7}c^{5}_{56}
-\frac{1}{7}c^{6}_{56}
-\frac{1}{7}c^{9}_{56}
$,
$-\frac{1}{7}c^{1}_{112}
-\frac{1}{7}c^{3}_{112}
+\frac{1}{7}c^{5}_{112}
+\frac{1}{7}c^{7}_{112}
+\frac{1}{7}c^{9}_{112}
+\frac{1}{7}c^{11}_{112}
-\frac{2}{7}c^{13}_{112}
-\frac{1}{7}c^{15}_{112}
-\frac{2}{7}c^{17}_{112}
-\frac{1}{7}c^{19}_{112}
+\frac{1}{7}c^{21}_{112}
+\frac{1}{7}c^{23}_{112}
$,
$\frac{1}{7}c^{2}_{56}
-\frac{1}{7}c^{3}_{56}
+\frac{1}{7}c^{11}_{56}
$,
$-\frac{2}{7}c^{1}_{56}
+\frac{1}{7}c^{3}_{56}
+\frac{1}{7}c^{5}_{56}
-\frac{1}{7}c^{7}_{56}
-\frac{1}{7}c^{9}_{56}
+\frac{1}{7}c^{10}_{56}
+\frac{1}{7}c^{11}_{56}
$;\ \ 
$\frac{1}{7}c^{5}_{56}
+\frac{1}{7}c^{6}_{56}
+\frac{1}{7}c^{9}_{56}
$,
$-\frac{1}{7}c^{1}_{112}
-\frac{1}{7}c^{5}_{112}
+\frac{1}{7}c^{15}_{112}
-\frac{1}{7}c^{19}_{112}
$,
$\frac{1}{7}c^{3}_{112}
-\frac{1}{7}c^{9}_{112}
+\frac{1}{7}c^{11}_{112}
+\frac{1}{7}c^{23}_{112}
$;\ \ 
$-\frac{2}{7}c^{1}_{56}
+\frac{1}{7}c^{3}_{56}
+\frac{1}{7}c^{5}_{56}
-\frac{1}{7}c^{7}_{56}
-\frac{1}{7}c^{9}_{56}
+\frac{1}{7}c^{10}_{56}
+\frac{1}{7}c^{11}_{56}
$,
$\frac{1}{7}c^{5}_{56}
+\frac{1}{7}c^{6}_{56}
+\frac{1}{7}c^{9}_{56}
$;\ \ 
$-\frac{1}{7}c^{2}_{56}
+\frac{1}{7}c^{3}_{56}
-\frac{1}{7}c^{11}_{56}
$)

Fail:
cnd( Tr$_I(\rho(\mathfrak s))$ ) =
56 does not divide
 ord($\rho(\mathfrak t)$) =
7, I = [ 1/7 ]. Prop. B.4 (2)

 \ \color{black}

\noindent 515: (dims,levels) = $(6;7
)$,
irreps = $6_{7,1}^{3}$,
pord$(\rho_\text{isum}(\mathfrak{t})) = 7$,

\vskip 0.7ex
\hangindent=5.5em \hangafter=1
{\white .}\hskip 1em $\rho_\text{isum}(\mathfrak{t})$ =
 $( \frac{1}{7},
\frac{2}{7},
\frac{3}{7},
\frac{4}{7},
\frac{5}{7},
\frac{6}{7} )
$,

\vskip 0.7ex
\hangindent=5.5em \hangafter=1
{\white .}\hskip 1em $\rho_\text{isum}(\mathfrak{s})$ =
$\mathrm{i}$($\frac{1}{7}c^{2}_{56}
+\frac{1}{7}c^{3}_{56}
-\frac{1}{7}c^{11}_{56}
$,
$\frac{1}{7}c^{5}_{56}
-\frac{1}{7}c^{6}_{56}
+\frac{1}{7}c^{9}_{56}
$,
$\frac{1}{7}c^{1}_{112}
+\frac{1}{7}c^{5}_{112}
-\frac{1}{7}c^{7}_{112}
-\frac{1}{7}c^{9}_{112}
+\frac{1}{7}c^{15}_{112}
+\frac{2}{7}c^{17}_{112}
+\frac{1}{7}c^{19}_{112}
-\frac{1}{7}c^{23}_{112}
$,
$\frac{2}{7}c^{1}_{56}
-\frac{1}{7}c^{3}_{56}
-\frac{1}{7}c^{5}_{56}
+\frac{1}{7}c^{7}_{56}
+\frac{1}{7}c^{9}_{56}
+\frac{1}{7}c^{10}_{56}
-\frac{1}{7}c^{11}_{56}
$,
$-\frac{1}{7}c^{1}_{112}
+\frac{1}{7}c^{3}_{112}
+\frac{1}{7}c^{11}_{112}
+\frac{1}{7}c^{15}_{112}
$,
$\frac{1}{7}c^{3}_{112}
-\frac{1}{7}c^{5}_{112}
-\frac{1}{7}c^{9}_{112}
-\frac{1}{7}c^{11}_{112}
+\frac{2}{7}c^{13}_{112}
+\frac{1}{7}c^{19}_{112}
-\frac{1}{7}c^{21}_{112}
+\frac{1}{7}c^{23}_{112}
$;\ \ 
$-\frac{2}{7}c^{1}_{56}
+\frac{1}{7}c^{3}_{56}
+\frac{1}{7}c^{5}_{56}
-\frac{1}{7}c^{7}_{56}
-\frac{1}{7}c^{9}_{56}
-\frac{1}{7}c^{10}_{56}
+\frac{1}{7}c^{11}_{56}
$,
$\frac{1}{7}c^{3}_{112}
-\frac{1}{7}c^{5}_{112}
-\frac{1}{7}c^{9}_{112}
-\frac{1}{7}c^{11}_{112}
+\frac{2}{7}c^{13}_{112}
+\frac{1}{7}c^{19}_{112}
-\frac{1}{7}c^{21}_{112}
+\frac{1}{7}c^{23}_{112}
$,
$\frac{1}{7}c^{2}_{56}
+\frac{1}{7}c^{3}_{56}
-\frac{1}{7}c^{11}_{56}
$,
$-\frac{1}{7}c^{1}_{112}
-\frac{1}{7}c^{5}_{112}
+\frac{1}{7}c^{7}_{112}
+\frac{1}{7}c^{9}_{112}
-\frac{1}{7}c^{15}_{112}
-\frac{2}{7}c^{17}_{112}
-\frac{1}{7}c^{19}_{112}
+\frac{1}{7}c^{23}_{112}
$,
$-\frac{1}{7}c^{1}_{112}
+\frac{1}{7}c^{3}_{112}
+\frac{1}{7}c^{11}_{112}
+\frac{1}{7}c^{15}_{112}
$;\ \ 
$\frac{1}{7}c^{5}_{56}
-\frac{1}{7}c^{6}_{56}
+\frac{1}{7}c^{9}_{56}
$,
$\frac{1}{7}c^{1}_{112}
-\frac{1}{7}c^{3}_{112}
-\frac{1}{7}c^{11}_{112}
-\frac{1}{7}c^{15}_{112}
$,
$-\frac{1}{7}c^{2}_{56}
-\frac{1}{7}c^{3}_{56}
+\frac{1}{7}c^{11}_{56}
$,
$-\frac{2}{7}c^{1}_{56}
+\frac{1}{7}c^{3}_{56}
+\frac{1}{7}c^{5}_{56}
-\frac{1}{7}c^{7}_{56}
-\frac{1}{7}c^{9}_{56}
-\frac{1}{7}c^{10}_{56}
+\frac{1}{7}c^{11}_{56}
$;\ \ 
$-\frac{1}{7}c^{5}_{56}
+\frac{1}{7}c^{6}_{56}
-\frac{1}{7}c^{9}_{56}
$,
$\frac{1}{7}c^{3}_{112}
-\frac{1}{7}c^{5}_{112}
-\frac{1}{7}c^{9}_{112}
-\frac{1}{7}c^{11}_{112}
+\frac{2}{7}c^{13}_{112}
+\frac{1}{7}c^{19}_{112}
-\frac{1}{7}c^{21}_{112}
+\frac{1}{7}c^{23}_{112}
$,
$\frac{1}{7}c^{1}_{112}
+\frac{1}{7}c^{5}_{112}
-\frac{1}{7}c^{7}_{112}
-\frac{1}{7}c^{9}_{112}
+\frac{1}{7}c^{15}_{112}
+\frac{2}{7}c^{17}_{112}
+\frac{1}{7}c^{19}_{112}
-\frac{1}{7}c^{23}_{112}
$;\ \ 
$\frac{2}{7}c^{1}_{56}
-\frac{1}{7}c^{3}_{56}
-\frac{1}{7}c^{5}_{56}
+\frac{1}{7}c^{7}_{56}
+\frac{1}{7}c^{9}_{56}
+\frac{1}{7}c^{10}_{56}
-\frac{1}{7}c^{11}_{56}
$,
$-\frac{1}{7}c^{5}_{56}
+\frac{1}{7}c^{6}_{56}
-\frac{1}{7}c^{9}_{56}
$;\ \ 
$-\frac{1}{7}c^{2}_{56}
-\frac{1}{7}c^{3}_{56}
+\frac{1}{7}c^{11}_{56}
$)

Fail:
cnd( Tr$_I(\rho(\mathfrak s))$ ) =
56 does not divide
 ord($\rho(\mathfrak t)$) =
7, I = [ 1/7 ]. Prop. B.4 (2)

 \ \color{black}

\noindent 516: (dims,levels) = $(6;8
)$,
irreps = $6_{8,1}^{1,0}$,
pord$(\rho_\text{isum}(\mathfrak{t})) = 8$,

\vskip 0.7ex
\hangindent=5.5em \hangafter=1
{\white .}\hskip 1em $\rho_\text{isum}(\mathfrak{t})$ =
 $( 0,
\frac{1}{2},
\frac{1}{8},
\frac{3}{8},
\frac{5}{8},
\frac{7}{8} )
$,

\vskip 0.7ex
\hangindent=5.5em \hangafter=1
{\white .}\hskip 1em $\rho_\text{isum}(\mathfrak{s})$ =
($0$,
$0$,
$\frac{1}{2}$,
$\frac{1}{2}$,
$\frac{1}{2}$,
$\frac{1}{2}$;
$0$,
$\frac{1}{2}$,
$-\frac{1}{2}$,
$\frac{1}{2}$,
$-\frac{1}{2}$;
$0$,
$-\frac{1}{2}$,
$0$,
$\frac{1}{2}$;
$0$,
$\frac{1}{2}$,
$0$;
$0$,
$-\frac{1}{2}$;
$0$)

Fail:
all rows of $U \rho(\mathfrak s) U^\dagger$
 contain zero for any block-diagonal $U$. Prop. B.5 (4) eqn. (B.27)

 \ \color{black}

\noindent 517: (dims,levels) = $(6;8
)$,
irreps = $6_{8,2}^{1,0}$,
pord$(\rho_\text{isum}(\mathfrak{t})) = 8$,

\vskip 0.7ex
\hangindent=5.5em \hangafter=1
{\white .}\hskip 1em $\rho_\text{isum}(\mathfrak{t})$ =
 $( 0,
\frac{1}{2},
\frac{1}{4},
\frac{3}{4},
\frac{1}{8},
\frac{3}{8} )
$,

\vskip 0.7ex
\hangindent=5.5em \hangafter=1
{\white .}\hskip 1em $\rho_\text{isum}(\mathfrak{s})$ =
$\mathrm{i}$($-\sqrt{\frac{1}{8}}$,
$\sqrt{\frac{1}{8}}$,
$\sqrt{\frac{1}{8}}$,
$\sqrt{\frac{1}{8}}$,
$\frac{1}{2}$,
$\frac{1}{2}$;\ \ 
$-\sqrt{\frac{1}{8}}$,
$-\sqrt{\frac{1}{8}}$,
$-\sqrt{\frac{1}{8}}$,
$\frac{1}{2}$,
$\frac{1}{2}$;\ \ 
$\sqrt{\frac{1}{8}}$,
$\sqrt{\frac{1}{8}}$,
$-\frac{1}{2}$,
$\frac{1}{2}$;\ \ 
$\sqrt{\frac{1}{8}}$,
$\frac{1}{2}$,
$-\frac{1}{2}$;\ \ 
$0$,
$0$;\ \ 
$0$)

Fail:
$\sigma(\rho(\mathfrak s)_\mathrm{ndeg}) \neq
 (\rho(\mathfrak t)^a \rho(\mathfrak s) \rho(\mathfrak t)^b
 \rho(\mathfrak s) \rho(\mathfrak t)^a)_\mathrm{ndeg}$,
 $\sigma = a$ = 3. Prop. B.5 (3) eqn. (B.25)

 \ \color{black}

\noindent 518: (dims,levels) = $(6;9
)$,
irreps = $6_{9,2}^{1,0}$,
pord$(\rho_\text{isum}(\mathfrak{t})) = 9$,

\vskip 0.7ex
\hangindent=5.5em \hangafter=1
{\white .}\hskip 1em $\rho_\text{isum}(\mathfrak{t})$ =
 $( \frac{1}{9},
\frac{2}{9},
\frac{4}{9},
\frac{5}{9},
\frac{7}{9},
\frac{8}{9} )
$,

\vskip 0.7ex
\hangindent=5.5em \hangafter=1
{\white .}\hskip 1em $\rho_\text{isum}(\mathfrak{s})$ =
($-\frac{1}{3}$,
$\frac{1}{3}c^{1}_{36}
$,
$\frac{1}{3}$,
$\frac{1}{3}c^{5}_{36}
$,
$\frac{1}{3}$,
$-\frac{1}{3}c^{1}_{36}
+\frac{1}{3}c^{5}_{36}
$;
$-\frac{1}{3}$,
$\frac{1}{3}c^{1}_{36}
-\frac{1}{3}c^{5}_{36}
$,
$\frac{1}{3}$,
$\frac{1}{3}c^{5}_{36}
$,
$-\frac{1}{3}$;
$-\frac{1}{3}$,
$\frac{1}{3}c^{1}_{36}
$,
$-\frac{1}{3}$,
$\frac{1}{3}c^{5}_{36}
$;
$-\frac{1}{3}$,
$-\frac{1}{3}c^{1}_{36}
+\frac{1}{3}c^{5}_{36}
$,
$\frac{1}{3}$;
$-\frac{1}{3}$,
$-\frac{1}{3}c^{1}_{36}
$;
$-\frac{1}{3}$)

Fail:
cnd($\rho(\mathfrak s)_\mathrm{ndeg}$) = 36 does not divide
 ord($\rho(\mathfrak t)$)=9. Prop. B.4 (2)

 \ \color{black}

 \color{blue}

\noindent 519: (dims,levels) = $(6;9
)$,
irreps = $6_{9,3}^{1,0}$,
pord$(\rho_\text{isum}(\mathfrak{t})) = 9$,

\vskip 0.7ex
\hangindent=5.5em \hangafter=1
{\white .}\hskip 1em $\rho_\text{isum}(\mathfrak{t})$ =
 $( \frac{1}{9},
\frac{2}{9},
\frac{4}{9},
\frac{5}{9},
\frac{7}{9},
\frac{8}{9} )
$,

\vskip 0.7ex
\hangindent=5.5em \hangafter=1
{\white .}\hskip 1em $\rho_\text{isum}(\mathfrak{s})$ =
($\frac{1}{3}$,
$\frac{1}{3}c^{2}_{9}
$,
$\frac{1}{3}$,
$-\frac{1}{3}c^{1}_{9}
$,
$\frac{1}{3}$,
$\frac{1}{3} c_9^4 $;
$\frac{1}{3}$,
$\frac{1}{3} c_9^4 $,
$-\frac{1}{3}$,
$\frac{1}{3}c^{1}_{9}
$,
$\frac{1}{3}$;
$\frac{1}{3}$,
$-\frac{1}{3}c^{2}_{9}
$,
$\frac{1}{3}$,
$\frac{1}{3}c^{1}_{9}
$;
$\frac{1}{3}$,
$-\frac{1}{3} c_9^4 $,
$-\frac{1}{3}$;
$\frac{1}{3}$,
$\frac{1}{3}c^{2}_{9}
$;
$\frac{1}{3}$)

Pass. 

 \ \color{black}

\noindent 520: (dims,levels) = $(6;9
)$,
irreps = $6_{9,1}^{1,0}$,
pord$(\rho_\text{isum}(\mathfrak{t})) = 9$,

\vskip 0.7ex
\hangindent=5.5em \hangafter=1
{\white .}\hskip 1em $\rho_\text{isum}(\mathfrak{t})$ =
 $( \frac{1}{9},
\frac{2}{9},
\frac{4}{9},
\frac{5}{9},
\frac{7}{9},
\frac{8}{9} )
$,

\vskip 0.7ex
\hangindent=5.5em \hangafter=1
{\white .}\hskip 1em $\rho_\text{isum}(\mathfrak{s})$ =
$\mathrm{i}$($-\frac{1}{3}$,
$\frac{1}{3}c^{7}_{72}
$,
$\frac{1}{3}$,
$-\frac{1}{3}c^{5}_{72}
+\frac{1}{3}c^{7}_{72}
$,
$\frac{1}{3}$,
$-\frac{1}{3}c^{5}_{72}
$;\ \ 
$\frac{1}{3}$,
$\frac{1}{3}c^{5}_{72}
$,
$-\frac{1}{3}$,
$-\frac{1}{3}c^{5}_{72}
+\frac{1}{3}c^{7}_{72}
$,
$\frac{1}{3}$;\ \ 
$-\frac{1}{3}$,
$\frac{1}{3}c^{7}_{72}
$,
$-\frac{1}{3}$,
$-\frac{1}{3}c^{5}_{72}
+\frac{1}{3}c^{7}_{72}
$;\ \ 
$\frac{1}{3}$,
$-\frac{1}{3}c^{5}_{72}
$,
$-\frac{1}{3}$;\ \ 
$-\frac{1}{3}$,
$-\frac{1}{3}c^{7}_{72}
$;\ \ 
$\frac{1}{3}$)

Fail:
cnd( Tr$_I(\rho(\mathfrak s))$ ) =
4 does not divide
 ord($\rho(\mathfrak t)$) =
9, I = [ 1/9 ]. Prop. B.4 (2)

 \ \color{black}

\noindent 521: (dims,levels) = $(6;9
)$,
irreps = $6_{9,1}^{5,0}$,
pord$(\rho_\text{isum}(\mathfrak{t})) = 9$,

\vskip 0.7ex
\hangindent=5.5em \hangafter=1
{\white .}\hskip 1em $\rho_\text{isum}(\mathfrak{t})$ =
 $( \frac{1}{9},
\frac{2}{9},
\frac{4}{9},
\frac{5}{9},
\frac{7}{9},
\frac{8}{9} )
$,

\vskip 0.7ex
\hangindent=5.5em \hangafter=1
{\white .}\hskip 1em $\rho_\text{isum}(\mathfrak{s})$ =
$\mathrm{i}$($\frac{1}{3}$,
$-\frac{1}{3}c^{11}_{72}
$,
$\frac{1}{3}$,
$-\frac{1}{3}c^{1}_{72}
$,
$\frac{1}{3}$,
$-\frac{1}{3}c^{1}_{72}
+\frac{1}{3}c^{11}_{72}
$;\ \ 
$-\frac{1}{3}$,
$-\frac{1}{3}c^{1}_{72}
+\frac{1}{3}c^{11}_{72}
$,
$\frac{1}{3}$,
$\frac{1}{3}c^{1}_{72}
$,
$-\frac{1}{3}$;\ \ 
$\frac{1}{3}$,
$\frac{1}{3}c^{11}_{72}
$,
$\frac{1}{3}$,
$\frac{1}{3}c^{1}_{72}
$;\ \ 
$-\frac{1}{3}$,
$\frac{1}{3}c^{1}_{72}
-\frac{1}{3}c^{11}_{72}
$,
$\frac{1}{3}$;\ \ 
$\frac{1}{3}$,
$-\frac{1}{3}c^{11}_{72}
$;\ \ 
$-\frac{1}{3}$)

Fail:
cnd( Tr$_I(\rho(\mathfrak s))$ ) =
4 does not divide
 ord($\rho(\mathfrak t)$) =
9, I = [ 1/9 ]. Prop. B.4 (2)

 \ \color{black}

\noindent 522: (dims,levels) = $(6;10
)$,
irreps = $3_{5}^{1}
\hskip -1.5pt \otimes \hskip -1.5pt
2_{2}^{1,0}$,
pord$(\rho_\text{isum}(\mathfrak{t})) = 10$,

\vskip 0.7ex
\hangindent=5.5em \hangafter=1
{\white .}\hskip 1em $\rho_\text{isum}(\mathfrak{t})$ =
 $( 0,
\frac{1}{2},
\frac{1}{5},
\frac{4}{5},
\frac{3}{10},
\frac{7}{10} )
$,

\vskip 0.7ex
\hangindent=5.5em \hangafter=1
{\white .}\hskip 1em $\rho_\text{isum}(\mathfrak{s})$ =
($-\sqrt{\frac{1}{20}}$,
$\sqrt{\frac{3}{20}}$,
$-\sqrt{\frac{1}{10}}$,
$-\sqrt{\frac{1}{10}}$,
$\sqrt{\frac{3}{10}}$,
$\sqrt{\frac{3}{10}}$;
$\sqrt{\frac{1}{20}}$,
$\sqrt{\frac{3}{10}}$,
$\sqrt{\frac{3}{10}}$,
$\sqrt{\frac{1}{10}}$,
$\sqrt{\frac{1}{10}}$;
$\frac{5+\sqrt{5}}{20}$,
$\frac{-5+\sqrt{5}}{20}$,
$\frac{3}{2\sqrt{15}}c^{1}_{5}
$,
$-\frac{3}{2\sqrt{15}\mathrm{i}}s^{3}_{20}
$;
$\frac{5+\sqrt{5}}{20}$,
$-\frac{3}{2\sqrt{15}\mathrm{i}}s^{3}_{20}
$,
$\frac{3}{2\sqrt{15}}c^{1}_{5}
$;
$-\frac{5+\sqrt{5}}{20}$,
$\frac{5-\sqrt{5}}{20}$;
$-\frac{5+\sqrt{5}}{20}$)

Fail:
cnd($\rho(\mathfrak s)_\mathrm{ndeg}$) = 120 does not divide
 ord($\rho(\mathfrak t)$)=10. Prop. B.4 (2)

 \ \color{black}

\noindent 523: (dims,levels) = $(6;10
)$,
irreps = $3_{5}^{3}
\hskip -1.5pt \otimes \hskip -1.5pt
2_{2}^{1,0}$,
pord$(\rho_\text{isum}(\mathfrak{t})) = 10$,

\vskip 0.7ex
\hangindent=5.5em \hangafter=1
{\white .}\hskip 1em $\rho_\text{isum}(\mathfrak{t})$ =
 $( 0,
\frac{1}{2},
\frac{2}{5},
\frac{3}{5},
\frac{1}{10},
\frac{9}{10} )
$,

\vskip 0.7ex
\hangindent=5.5em \hangafter=1
{\white .}\hskip 1em $\rho_\text{isum}(\mathfrak{s})$ =
($\sqrt{\frac{1}{20}}$,
$\sqrt{\frac{3}{20}}$,
$-\sqrt{\frac{1}{10}}$,
$-\sqrt{\frac{1}{10}}$,
$\sqrt{\frac{3}{10}}$,
$\sqrt{\frac{3}{10}}$;
$-\sqrt{\frac{1}{20}}$,
$-\sqrt{\frac{3}{10}}$,
$-\sqrt{\frac{3}{10}}$,
$-\sqrt{\frac{1}{10}}$,
$-\sqrt{\frac{1}{10}}$;
$\frac{5-\sqrt{5}}{20}$,
$-\frac{5+\sqrt{5}}{20}$,
$\frac{3}{2\sqrt{15}\mathrm{i}}s^{3}_{20}
$,
$-\frac{3}{2\sqrt{15}}c^{1}_{5}
$;
$\frac{5-\sqrt{5}}{20}$,
$-\frac{3}{2\sqrt{15}}c^{1}_{5}
$,
$\frac{3}{2\sqrt{15}\mathrm{i}}s^{3}_{20}
$;
$\frac{-5+\sqrt{5}}{20}$,
$\frac{5+\sqrt{5}}{20}$;
$\frac{-5+\sqrt{5}}{20}$)

Fail:
cnd($\rho(\mathfrak s)_\mathrm{ndeg}$) = 120 does not divide
 ord($\rho(\mathfrak t)$)=10. Prop. B.4 (2)

 \ \color{black}

\noindent 524: (dims,levels) = $(6;10
)$,
irreps = $6_{5}^{1}
\hskip -1.5pt \otimes \hskip -1.5pt
1_{2}^{1,0}$,
pord$(\rho_\text{isum}(\mathfrak{t})) = 5$,

\vskip 0.7ex
\hangindent=5.5em \hangafter=1
{\white .}\hskip 1em $\rho_\text{isum}(\mathfrak{t})$ =
 $( \frac{1}{2},
\frac{1}{2},
\frac{1}{10},
\frac{3}{10},
\frac{7}{10},
\frac{9}{10} )
$,

\vskip 0.7ex
\hangindent=5.5em \hangafter=1
{\white .}\hskip 1em $\rho_\text{isum}(\mathfrak{s})$ =
$\mathrm{i}$($-\frac{1}{5}c^{1}_{20}
$,
$\frac{1}{5}c^{3}_{20}
$,
$\frac{1}{5}c^{1}_{20}
+\frac{1}{5}c^{3}_{20}
$,
$\frac{1}{5}c^{3}_{40}
+\frac{1}{5}c^{7}_{40}
$,
$\frac{4}{5\sqrt{10}}c^{1}_{20}
-\frac{2}{5\sqrt{10}}c^{3}_{20}
$,
$\frac{1}{5}c^{1}_{20}
-\frac{1}{5}c^{3}_{20}
$;\ \ 
$\frac{1}{5}c^{1}_{20}
$,
$\frac{1}{5}c^{1}_{20}
-\frac{1}{5}c^{3}_{20}
$,
$-\frac{4}{5\sqrt{10}}c^{1}_{20}
+\frac{2}{5\sqrt{10}}c^{3}_{20}
$,
$\frac{1}{5}c^{3}_{40}
+\frac{1}{5}c^{7}_{40}
$,
$-\frac{1}{5}c^{1}_{20}
-\frac{1}{5}c^{3}_{20}
$;\ \ 
$\frac{1}{5}c^{1}_{20}
$,
$\frac{4}{5\sqrt{10}}c^{1}_{20}
-\frac{2}{5\sqrt{10}}c^{3}_{20}
$,
$-\frac{1}{5}c^{3}_{40}
-\frac{1}{5}c^{7}_{40}
$,
$-\frac{1}{5}c^{3}_{20}
$;\ \ 
$-\frac{1}{5}c^{3}_{20}
$,
$\frac{1}{5}c^{1}_{20}
$,
$\frac{1}{5}c^{3}_{40}
+\frac{1}{5}c^{7}_{40}
$;\ \ 
$\frac{1}{5}c^{3}_{20}
$,
$\frac{4}{5\sqrt{10}}c^{1}_{20}
-\frac{2}{5\sqrt{10}}c^{3}_{20}
$;\ \ 
$-\frac{1}{5}c^{1}_{20}
$)

Fail:
cnd($\rho(\mathfrak s)_\mathrm{ndeg}$) = 40 does not divide
 ord($\rho(\mathfrak t)$)=10. Prop. B.4 (2)

 \ \color{black}

\noindent 525: (dims,levels) = $(6;11
)$,
irreps = $6_{11}^{1}$,
pord$(\rho_\text{isum}(\mathfrak{t})) = 11$,

\vskip 0.7ex
\hangindent=5.5em \hangafter=1
{\white .}\hskip 1em $\rho_\text{isum}(\mathfrak{t})$ =
 $( 0,
\frac{1}{11},
\frac{3}{11},
\frac{4}{11},
\frac{5}{11},
\frac{9}{11} )
$,

\vskip 0.7ex
\hangindent=5.5em \hangafter=1
{\white .}\hskip 1em $\rho_\text{isum}(\mathfrak{s})$ =
$\mathrm{i}$($-\sqrt{\frac{1}{11}}$,
$\sqrt{\frac{2}{11}}$,
$\sqrt{\frac{2}{11}}$,
$\sqrt{\frac{2}{11}}$,
$\sqrt{\frac{2}{11}}$,
$\sqrt{\frac{2}{11}}$;\ \ 
$-\frac{1}{\sqrt{11}}c^{2}_{11}
$,
$-\frac{1}{\sqrt{11}}c^{1}_{11}
$,
$-\frac{1}{\sqrt{11}}c^{4}_{11}
$,
$-\frac{1}{\sqrt{11}}c^{3}_{11}
$,
$\frac{1}{\sqrt{11}\mathrm{i}}s^{9}_{44}
$;\ \ 
$\frac{1}{\sqrt{11}\mathrm{i}}s^{9}_{44}
$,
$-\frac{1}{\sqrt{11}}c^{2}_{11}
$,
$-\frac{1}{\sqrt{11}}c^{4}_{11}
$,
$-\frac{1}{\sqrt{11}}c^{3}_{11}
$;\ \ 
$-\frac{1}{\sqrt{11}}c^{3}_{11}
$,
$\frac{1}{\sqrt{11}\mathrm{i}}s^{9}_{44}
$,
$-\frac{1}{\sqrt{11}}c^{1}_{11}
$;\ \ 
$-\frac{1}{\sqrt{11}}c^{1}_{11}
$,
$-\frac{1}{\sqrt{11}}c^{2}_{11}
$;\ \ 
$-\frac{1}{\sqrt{11}}c^{4}_{11}
$)

Fail:
cnd($\rho(\mathfrak s)_\mathrm{ndeg}$) = 88 does not divide
 ord($\rho(\mathfrak t)$)=11. Prop. B.4 (2)

 \ \color{black}

\noindent 526: (dims,levels) = $(6;11
)$,
irreps = $6_{11}^{7}$,
pord$(\rho_\text{isum}(\mathfrak{t})) = 11$,

\vskip 0.7ex
\hangindent=5.5em \hangafter=1
{\white .}\hskip 1em $\rho_\text{isum}(\mathfrak{t})$ =
 $( 0,
\frac{2}{11},
\frac{6}{11},
\frac{7}{11},
\frac{8}{11},
\frac{10}{11} )
$,

\vskip 0.7ex
\hangindent=5.5em \hangafter=1
{\white .}\hskip 1em $\rho_\text{isum}(\mathfrak{s})$ =
$\mathrm{i}$($\sqrt{\frac{1}{11}}$,
$\sqrt{\frac{2}{11}}$,
$\sqrt{\frac{2}{11}}$,
$\sqrt{\frac{2}{11}}$,
$\sqrt{\frac{2}{11}}$,
$\sqrt{\frac{2}{11}}$;\ \ 
$\frac{1}{\sqrt{11}}c^{4}_{11}
$,
$\frac{1}{\sqrt{11}}c^{2}_{11}
$,
$\frac{1}{\sqrt{11}}c^{1}_{11}
$,
$\frac{1}{\sqrt{11}}c^{3}_{11}
$,
$-\frac{1}{\sqrt{11}\mathrm{i}}s^{9}_{44}
$;\ \ 
$\frac{1}{\sqrt{11}}c^{1}_{11}
$,
$-\frac{1}{\sqrt{11}\mathrm{i}}s^{9}_{44}
$,
$\frac{1}{\sqrt{11}}c^{4}_{11}
$,
$\frac{1}{\sqrt{11}}c^{3}_{11}
$;\ \ 
$\frac{1}{\sqrt{11}}c^{3}_{11}
$,
$\frac{1}{\sqrt{11}}c^{2}_{11}
$,
$\frac{1}{\sqrt{11}}c^{4}_{11}
$;\ \ 
$-\frac{1}{\sqrt{11}\mathrm{i}}s^{9}_{44}
$,
$\frac{1}{\sqrt{11}}c^{1}_{11}
$;\ \ 
$\frac{1}{\sqrt{11}}c^{2}_{11}
$)

Fail:
cnd($\rho(\mathfrak s)_\mathrm{ndeg}$) = 88 does not divide
 ord($\rho(\mathfrak t)$)=11. Prop. B.4 (2)

 \ \color{black}

\noindent 527: (dims,levels) = $(6;12
)$,
irreps = $3_{4}^{1,0}
\hskip -1.5pt \otimes \hskip -1.5pt
2_{3}^{1,0}$,
pord$(\rho_\text{isum}(\mathfrak{t})) = 12$,

\vskip 0.7ex
\hangindent=5.5em \hangafter=1
{\white .}\hskip 1em $\rho_\text{isum}(\mathfrak{t})$ =
 $( 0,
\frac{1}{3},
\frac{1}{4},
\frac{3}{4},
\frac{1}{12},
\frac{7}{12} )
$,

\vskip 0.7ex
\hangindent=5.5em \hangafter=1
{\white .}\hskip 1em $\rho_\text{isum}(\mathfrak{s})$ =
$\mathrm{i}$($0$,
$0$,
$\sqrt{\frac{1}{6}}$,
$\sqrt{\frac{1}{6}}$,
$\sqrt{\frac{1}{3}}$,
$\sqrt{\frac{1}{3}}$;\ \ 
$0$,
$\sqrt{\frac{1}{3}}$,
$\sqrt{\frac{1}{3}}$,
$-\sqrt{\frac{1}{6}}$,
$-\sqrt{\frac{1}{6}}$;\ \ 
$\sqrt{\frac{1}{12}}$,
$-\sqrt{\frac{1}{12}}$,
$-\sqrt{\frac{1}{6}}$,
$\sqrt{\frac{1}{6}}$;\ \ 
$\sqrt{\frac{1}{12}}$,
$\sqrt{\frac{1}{6}}$,
$-\sqrt{\frac{1}{6}}$;\ \ 
$-\sqrt{\frac{1}{12}}$,
$\sqrt{\frac{1}{12}}$;\ \ 
$-\sqrt{\frac{1}{12}}$)

Fail:
cnd($\rho(\mathfrak s)_\mathrm{ndeg}$) = 24 does not divide
 ord($\rho(\mathfrak t)$)=12. Prop. B.4 (2)

 \ \color{black}

 \color{blue}

\noindent 528: (dims,levels) = $(6;13
)$,
irreps = $6_{13}^{1}$,
pord$(\rho_\text{isum}(\mathfrak{t})) = 13$,

\vskip 0.7ex
\hangindent=5.5em \hangafter=1
{\white .}\hskip 1em $\rho_\text{isum}(\mathfrak{t})$ =
 $( \frac{1}{13},
\frac{3}{13},
\frac{4}{13},
\frac{9}{13},
\frac{10}{13},
\frac{12}{13} )
$,

\vskip 0.7ex
\hangindent=5.5em \hangafter=1
{\white .}\hskip 1em $\rho_\text{isum}(\mathfrak{s})$ =
$\mathrm{i}$($-\frac{1}{\sqrt{13}}c^{5}_{52}
$,
$\frac{1}{\sqrt{13}}c^{7}_{52}
$,
$\frac{1}{\sqrt{13}}c^{3}_{52}
$,
$\frac{1}{\sqrt{13}}c^{11}_{52}
$,
$\frac{1}{\sqrt{13}}c^{9}_{52}
$,
$-\frac{1}{\sqrt{13}}c^{1}_{52}
$;\ \ 
$-\frac{1}{\sqrt{13}}c^{11}_{52}
$,
$\frac{1}{\sqrt{13}}c^{1}_{52}
$,
$-\frac{1}{\sqrt{13}}c^{5}_{52}
$,
$\frac{1}{\sqrt{13}}c^{3}_{52}
$,
$\frac{1}{\sqrt{13}}c^{9}_{52}
$;\ \ 
$\frac{1}{\sqrt{13}}c^{7}_{52}
$,
$\frac{1}{\sqrt{13}}c^{9}_{52}
$,
$-\frac{1}{\sqrt{13}}c^{5}_{52}
$,
$\frac{1}{\sqrt{13}}c^{11}_{52}
$;\ \ 
$-\frac{1}{\sqrt{13}}c^{7}_{52}
$,
$-\frac{1}{\sqrt{13}}c^{1}_{52}
$,
$-\frac{1}{\sqrt{13}}c^{3}_{52}
$;\ \ 
$\frac{1}{\sqrt{13}}c^{11}_{52}
$,
$-\frac{1}{\sqrt{13}}c^{7}_{52}
$;\ \ 
$\frac{1}{\sqrt{13}}c^{5}_{52}
$)

Pass. 

 \ \color{black}

 \color{blue}

\noindent 529: (dims,levels) = $(6;13
)$,
irreps = $6_{13}^{2}$,
pord$(\rho_\text{isum}(\mathfrak{t})) = 13$,

\vskip 0.7ex
\hangindent=5.5em \hangafter=1
{\white .}\hskip 1em $\rho_\text{isum}(\mathfrak{t})$ =
 $( \frac{2}{13},
\frac{5}{13},
\frac{6}{13},
\frac{7}{13},
\frac{8}{13},
\frac{11}{13} )
$,

\vskip 0.7ex
\hangindent=5.5em \hangafter=1
{\white .}\hskip 1em $\rho_\text{isum}(\mathfrak{s})$ =
$\mathrm{i}$($\frac{1}{\sqrt{13}}c^{3}_{52}
$,
$\frac{1}{\sqrt{13}}c^{9}_{52}
$,
$-\frac{1}{\sqrt{13}}c^{1}_{52}
$,
$-\frac{1}{\sqrt{13}}c^{5}_{52}
$,
$\frac{1}{\sqrt{13}}c^{7}_{52}
$,
$\frac{1}{\sqrt{13}}c^{11}_{52}
$;\ \ 
$-\frac{1}{\sqrt{13}}c^{1}_{52}
$,
$-\frac{1}{\sqrt{13}}c^{3}_{52}
$,
$\frac{1}{\sqrt{13}}c^{11}_{52}
$,
$-\frac{1}{\sqrt{13}}c^{5}_{52}
$,
$-\frac{1}{\sqrt{13}}c^{7}_{52}
$;\ \ 
$-\frac{1}{\sqrt{13}}c^{9}_{52}
$,
$-\frac{1}{\sqrt{13}}c^{7}_{52}
$,
$\frac{1}{\sqrt{13}}c^{11}_{52}
$,
$\frac{1}{\sqrt{13}}c^{5}_{52}
$;\ \ 
$\frac{1}{\sqrt{13}}c^{9}_{52}
$,
$\frac{1}{\sqrt{13}}c^{3}_{52}
$,
$-\frac{1}{\sqrt{13}}c^{1}_{52}
$;\ \ 
$\frac{1}{\sqrt{13}}c^{1}_{52}
$,
$\frac{1}{\sqrt{13}}c^{9}_{52}
$;\ \ 
$-\frac{1}{\sqrt{13}}c^{3}_{52}
$)

Pass. 

 \ \color{black}

\noindent 530: (dims,levels) = $(6;14
)$,
irreps = $6_{7,2}^{1}
\hskip -1.5pt \otimes \hskip -1.5pt
1_{2}^{1,0}$,
pord$(\rho_\text{isum}(\mathfrak{t})) = 7$,

\vskip 0.7ex
\hangindent=5.5em \hangafter=1
{\white .}\hskip 1em $\rho_\text{isum}(\mathfrak{t})$ =
 $( \frac{1}{14},
\frac{3}{14},
\frac{5}{14},
\frac{9}{14},
\frac{11}{14},
\frac{13}{14} )
$,

\vskip 0.7ex
\hangindent=5.5em \hangafter=1
{\white .}\hskip 1em $\rho_\text{isum}(\mathfrak{s})$ =
($-\frac{2}{7}+\frac{1}{7}c^{1}_{7}
$,
$-\frac{2}{7}c^{3}_{56}
-\frac{1}{7}c^{5}_{56}
+\frac{1}{7}c^{7}_{56}
+\frac{1}{7}c^{9}_{56}
-\frac{2}{7}c^{11}_{56}
$,
$\frac{1}{7}c^{3}_{56}
+\frac{2}{7}c^{5}_{56}
-\frac{1}{7}c^{7}_{56}
-\frac{2}{7}c^{9}_{56}
+\frac{1}{7}c^{11}_{56}
$,
$-\frac{3}{7}-\frac{1}{7}c^{1}_{7}
-\frac{1}{7}c^{2}_{7}
$,
$-\frac{2}{7}+\frac{1}{7}c^{2}_{7}
$,
$-\frac{1}{7}c^{3}_{56}
+\frac{1}{7}c^{5}_{56}
-\frac{1}{7}c^{9}_{56}
-\frac{1}{7}c^{11}_{56}
$;
$-\frac{3}{7}-\frac{1}{7}c^{1}_{7}
-\frac{1}{7}c^{2}_{7}
$,
$-\frac{2}{7}+\frac{1}{7}c^{1}_{7}
$,
$\frac{1}{7}c^{3}_{56}
-\frac{1}{7}c^{5}_{56}
+\frac{1}{7}c^{9}_{56}
+\frac{1}{7}c^{11}_{56}
$,
$\frac{1}{7}c^{3}_{56}
+\frac{2}{7}c^{5}_{56}
-\frac{1}{7}c^{7}_{56}
-\frac{2}{7}c^{9}_{56}
+\frac{1}{7}c^{11}_{56}
$,
$\frac{2}{7}-\frac{1}{7}c^{2}_{7}
$;
$-\frac{2}{7}+\frac{1}{7}c^{2}_{7}
$,
$-\frac{2}{7}c^{3}_{56}
-\frac{1}{7}c^{5}_{56}
+\frac{1}{7}c^{7}_{56}
+\frac{1}{7}c^{9}_{56}
-\frac{2}{7}c^{11}_{56}
$,
$\frac{1}{7}c^{3}_{56}
-\frac{1}{7}c^{5}_{56}
+\frac{1}{7}c^{9}_{56}
+\frac{1}{7}c^{11}_{56}
$,
$\frac{3}{7}+\frac{1}{7}c^{1}_{7}
+\frac{1}{7}c^{2}_{7}
$;
$-\frac{2}{7}+\frac{1}{7}c^{2}_{7}
$,
$-\frac{2}{7}+\frac{1}{7}c^{1}_{7}
$,
$-\frac{1}{7}c^{3}_{56}
-\frac{2}{7}c^{5}_{56}
+\frac{1}{7}c^{7}_{56}
+\frac{2}{7}c^{9}_{56}
-\frac{1}{7}c^{11}_{56}
$;
$-\frac{3}{7}-\frac{1}{7}c^{1}_{7}
-\frac{1}{7}c^{2}_{7}
$,
$\frac{2}{7}c^{3}_{56}
+\frac{1}{7}c^{5}_{56}
-\frac{1}{7}c^{7}_{56}
-\frac{1}{7}c^{9}_{56}
+\frac{2}{7}c^{11}_{56}
$;
$-\frac{2}{7}+\frac{1}{7}c^{1}_{7}
$)

Fail:
cnd($\rho(\mathfrak s)_\mathrm{ndeg}$) = 56 does not divide
 ord($\rho(\mathfrak t)$)=14. Prop. B.4 (2)

 \ \color{black}

\noindent 531: (dims,levels) = $(6;14
)$,
irreps = $6_{7,1}^{1}
\hskip -1.5pt \otimes \hskip -1.5pt
1_{2}^{1,0}$,
pord$(\rho_\text{isum}(\mathfrak{t})) = 7$,

\vskip 0.7ex
\hangindent=5.5em \hangafter=1
{\white .}\hskip 1em $\rho_\text{isum}(\mathfrak{t})$ =
 $( \frac{1}{14},
\frac{3}{14},
\frac{5}{14},
\frac{9}{14},
\frac{11}{14},
\frac{13}{14} )
$,

\vskip 0.7ex
\hangindent=5.5em \hangafter=1
{\white .}\hskip 1em $\rho_\text{isum}(\mathfrak{s})$ =
$\mathrm{i}$($-\frac{1}{7}c^{5}_{56}
-\frac{1}{7}c^{6}_{56}
-\frac{1}{7}c^{9}_{56}
$,
$\frac{1}{7}c^{1}_{112}
+\frac{1}{7}c^{5}_{112}
-\frac{1}{7}c^{15}_{112}
+\frac{1}{7}c^{19}_{112}
$,
$\frac{1}{7}c^{3}_{112}
-\frac{1}{7}c^{9}_{112}
+\frac{1}{7}c^{11}_{112}
+\frac{1}{7}c^{23}_{112}
$,
$\frac{2}{7}c^{1}_{56}
-\frac{1}{7}c^{3}_{56}
-\frac{1}{7}c^{5}_{56}
+\frac{1}{7}c^{7}_{56}
+\frac{1}{7}c^{9}_{56}
-\frac{1}{7}c^{10}_{56}
-\frac{1}{7}c^{11}_{56}
$,
$\frac{1}{7}c^{2}_{56}
-\frac{1}{7}c^{3}_{56}
+\frac{1}{7}c^{11}_{56}
$,
$\frac{1}{7}c^{1}_{112}
+\frac{1}{7}c^{3}_{112}
-\frac{1}{7}c^{5}_{112}
-\frac{1}{7}c^{7}_{112}
-\frac{1}{7}c^{9}_{112}
-\frac{1}{7}c^{11}_{112}
+\frac{2}{7}c^{13}_{112}
+\frac{1}{7}c^{15}_{112}
+\frac{2}{7}c^{17}_{112}
+\frac{1}{7}c^{19}_{112}
-\frac{1}{7}c^{21}_{112}
-\frac{1}{7}c^{23}_{112}
$;\ \ 
$\frac{2}{7}c^{1}_{56}
-\frac{1}{7}c^{3}_{56}
-\frac{1}{7}c^{5}_{56}
+\frac{1}{7}c^{7}_{56}
+\frac{1}{7}c^{9}_{56}
-\frac{1}{7}c^{10}_{56}
-\frac{1}{7}c^{11}_{56}
$,
$\frac{1}{7}c^{5}_{56}
+\frac{1}{7}c^{6}_{56}
+\frac{1}{7}c^{9}_{56}
$,
$\frac{1}{7}c^{1}_{112}
+\frac{1}{7}c^{3}_{112}
-\frac{1}{7}c^{5}_{112}
-\frac{1}{7}c^{7}_{112}
-\frac{1}{7}c^{9}_{112}
-\frac{1}{7}c^{11}_{112}
+\frac{2}{7}c^{13}_{112}
+\frac{1}{7}c^{15}_{112}
+\frac{2}{7}c^{17}_{112}
+\frac{1}{7}c^{19}_{112}
-\frac{1}{7}c^{21}_{112}
-\frac{1}{7}c^{23}_{112}
$,
$-\frac{1}{7}c^{3}_{112}
+\frac{1}{7}c^{9}_{112}
-\frac{1}{7}c^{11}_{112}
-\frac{1}{7}c^{23}_{112}
$,
$-\frac{1}{7}c^{2}_{56}
+\frac{1}{7}c^{3}_{56}
-\frac{1}{7}c^{11}_{56}
$;\ \ 
$\frac{1}{7}c^{2}_{56}
-\frac{1}{7}c^{3}_{56}
+\frac{1}{7}c^{11}_{56}
$,
$-\frac{1}{7}c^{1}_{112}
-\frac{1}{7}c^{5}_{112}
+\frac{1}{7}c^{15}_{112}
-\frac{1}{7}c^{19}_{112}
$,
$\frac{1}{7}c^{1}_{112}
+\frac{1}{7}c^{3}_{112}
-\frac{1}{7}c^{5}_{112}
-\frac{1}{7}c^{7}_{112}
-\frac{1}{7}c^{9}_{112}
-\frac{1}{7}c^{11}_{112}
+\frac{2}{7}c^{13}_{112}
+\frac{1}{7}c^{15}_{112}
+\frac{2}{7}c^{17}_{112}
+\frac{1}{7}c^{19}_{112}
-\frac{1}{7}c^{21}_{112}
-\frac{1}{7}c^{23}_{112}
$,
$-\frac{2}{7}c^{1}_{56}
+\frac{1}{7}c^{3}_{56}
+\frac{1}{7}c^{5}_{56}
-\frac{1}{7}c^{7}_{56}
-\frac{1}{7}c^{9}_{56}
+\frac{1}{7}c^{10}_{56}
+\frac{1}{7}c^{11}_{56}
$;\ \ 
$-\frac{1}{7}c^{2}_{56}
+\frac{1}{7}c^{3}_{56}
-\frac{1}{7}c^{11}_{56}
$,
$-\frac{1}{7}c^{5}_{56}
-\frac{1}{7}c^{6}_{56}
-\frac{1}{7}c^{9}_{56}
$,
$\frac{1}{7}c^{3}_{112}
-\frac{1}{7}c^{9}_{112}
+\frac{1}{7}c^{11}_{112}
+\frac{1}{7}c^{23}_{112}
$;\ \ 
$-\frac{2}{7}c^{1}_{56}
+\frac{1}{7}c^{3}_{56}
+\frac{1}{7}c^{5}_{56}
-\frac{1}{7}c^{7}_{56}
-\frac{1}{7}c^{9}_{56}
+\frac{1}{7}c^{10}_{56}
+\frac{1}{7}c^{11}_{56}
$,
$\frac{1}{7}c^{1}_{112}
+\frac{1}{7}c^{5}_{112}
-\frac{1}{7}c^{15}_{112}
+\frac{1}{7}c^{19}_{112}
$;\ \ 
$\frac{1}{7}c^{5}_{56}
+\frac{1}{7}c^{6}_{56}
+\frac{1}{7}c^{9}_{56}
$)

Fail:
cnd( Tr$_I(\rho(\mathfrak s))$ ) =
56 does not divide
 ord($\rho(\mathfrak t)$) =
14, I = [ 1/14 ]. Prop. B.4 (2)

 \ \color{black}

\noindent 532: (dims,levels) = $(6;14
)$,
irreps = $6_{7,1}^{3}
\hskip -1.5pt \otimes \hskip -1.5pt
1_{2}^{1,0}$,
pord$(\rho_\text{isum}(\mathfrak{t})) = 7$,

\vskip 0.7ex
\hangindent=5.5em \hangafter=1
{\white .}\hskip 1em $\rho_\text{isum}(\mathfrak{t})$ =
 $( \frac{1}{14},
\frac{3}{14},
\frac{5}{14},
\frac{9}{14},
\frac{11}{14},
\frac{13}{14} )
$,

\vskip 0.7ex
\hangindent=5.5em \hangafter=1
{\white .}\hskip 1em $\rho_\text{isum}(\mathfrak{s})$ =
$\mathrm{i}$($\frac{1}{7}c^{5}_{56}
-\frac{1}{7}c^{6}_{56}
+\frac{1}{7}c^{9}_{56}
$,
$\frac{1}{7}c^{3}_{112}
-\frac{1}{7}c^{5}_{112}
-\frac{1}{7}c^{9}_{112}
-\frac{1}{7}c^{11}_{112}
+\frac{2}{7}c^{13}_{112}
+\frac{1}{7}c^{19}_{112}
-\frac{1}{7}c^{21}_{112}
+\frac{1}{7}c^{23}_{112}
$,
$\frac{1}{7}c^{1}_{112}
+\frac{1}{7}c^{5}_{112}
-\frac{1}{7}c^{7}_{112}
-\frac{1}{7}c^{9}_{112}
+\frac{1}{7}c^{15}_{112}
+\frac{2}{7}c^{17}_{112}
+\frac{1}{7}c^{19}_{112}
-\frac{1}{7}c^{23}_{112}
$,
$\frac{2}{7}c^{1}_{56}
-\frac{1}{7}c^{3}_{56}
-\frac{1}{7}c^{5}_{56}
+\frac{1}{7}c^{7}_{56}
+\frac{1}{7}c^{9}_{56}
+\frac{1}{7}c^{10}_{56}
-\frac{1}{7}c^{11}_{56}
$,
$\frac{1}{7}c^{2}_{56}
+\frac{1}{7}c^{3}_{56}
-\frac{1}{7}c^{11}_{56}
$,
$-\frac{1}{7}c^{1}_{112}
+\frac{1}{7}c^{3}_{112}
+\frac{1}{7}c^{11}_{112}
+\frac{1}{7}c^{15}_{112}
$;\ \ 
$-\frac{2}{7}c^{1}_{56}
+\frac{1}{7}c^{3}_{56}
+\frac{1}{7}c^{5}_{56}
-\frac{1}{7}c^{7}_{56}
-\frac{1}{7}c^{9}_{56}
-\frac{1}{7}c^{10}_{56}
+\frac{1}{7}c^{11}_{56}
$,
$\frac{1}{7}c^{5}_{56}
-\frac{1}{7}c^{6}_{56}
+\frac{1}{7}c^{9}_{56}
$,
$\frac{1}{7}c^{1}_{112}
-\frac{1}{7}c^{3}_{112}
-\frac{1}{7}c^{11}_{112}
-\frac{1}{7}c^{15}_{112}
$,
$\frac{1}{7}c^{1}_{112}
+\frac{1}{7}c^{5}_{112}
-\frac{1}{7}c^{7}_{112}
-\frac{1}{7}c^{9}_{112}
+\frac{1}{7}c^{15}_{112}
+\frac{2}{7}c^{17}_{112}
+\frac{1}{7}c^{19}_{112}
-\frac{1}{7}c^{23}_{112}
$,
$-\frac{1}{7}c^{2}_{56}
-\frac{1}{7}c^{3}_{56}
+\frac{1}{7}c^{11}_{56}
$;\ \ 
$\frac{1}{7}c^{2}_{56}
+\frac{1}{7}c^{3}_{56}
-\frac{1}{7}c^{11}_{56}
$,
$-\frac{1}{7}c^{3}_{112}
+\frac{1}{7}c^{5}_{112}
+\frac{1}{7}c^{9}_{112}
+\frac{1}{7}c^{11}_{112}
-\frac{2}{7}c^{13}_{112}
-\frac{1}{7}c^{19}_{112}
+\frac{1}{7}c^{21}_{112}
-\frac{1}{7}c^{23}_{112}
$,
$\frac{1}{7}c^{1}_{112}
-\frac{1}{7}c^{3}_{112}
-\frac{1}{7}c^{11}_{112}
-\frac{1}{7}c^{15}_{112}
$,
$-\frac{2}{7}c^{1}_{56}
+\frac{1}{7}c^{3}_{56}
+\frac{1}{7}c^{5}_{56}
-\frac{1}{7}c^{7}_{56}
-\frac{1}{7}c^{9}_{56}
-\frac{1}{7}c^{10}_{56}
+\frac{1}{7}c^{11}_{56}
$;\ \ 
$-\frac{1}{7}c^{2}_{56}
-\frac{1}{7}c^{3}_{56}
+\frac{1}{7}c^{11}_{56}
$,
$-\frac{1}{7}c^{5}_{56}
+\frac{1}{7}c^{6}_{56}
-\frac{1}{7}c^{9}_{56}
$,
$\frac{1}{7}c^{1}_{112}
+\frac{1}{7}c^{5}_{112}
-\frac{1}{7}c^{7}_{112}
-\frac{1}{7}c^{9}_{112}
+\frac{1}{7}c^{15}_{112}
+\frac{2}{7}c^{17}_{112}
+\frac{1}{7}c^{19}_{112}
-\frac{1}{7}c^{23}_{112}
$;\ \ 
$\frac{2}{7}c^{1}_{56}
-\frac{1}{7}c^{3}_{56}
-\frac{1}{7}c^{5}_{56}
+\frac{1}{7}c^{7}_{56}
+\frac{1}{7}c^{9}_{56}
+\frac{1}{7}c^{10}_{56}
-\frac{1}{7}c^{11}_{56}
$,
$\frac{1}{7}c^{3}_{112}
-\frac{1}{7}c^{5}_{112}
-\frac{1}{7}c^{9}_{112}
-\frac{1}{7}c^{11}_{112}
+\frac{2}{7}c^{13}_{112}
+\frac{1}{7}c^{19}_{112}
-\frac{1}{7}c^{21}_{112}
+\frac{1}{7}c^{23}_{112}
$;\ \ 
$-\frac{1}{7}c^{5}_{56}
+\frac{1}{7}c^{6}_{56}
-\frac{1}{7}c^{9}_{56}
$)

Fail:
cnd( Tr$_I(\rho(\mathfrak s))$ ) =
56 does not divide
 ord($\rho(\mathfrak t)$) =
14, I = [ 1/14 ]. Prop. B.4 (2)

 \ \color{black}

\noindent 533: (dims,levels) = $(6;14
)$,
irreps = $3_{7}^{1}
\hskip -1.5pt \otimes \hskip -1.5pt
2_{2}^{1,0}$,
pord$(\rho_\text{isum}(\mathfrak{t})) = 14$,

\vskip 0.7ex
\hangindent=5.5em \hangafter=1
{\white .}\hskip 1em $\rho_\text{isum}(\mathfrak{t})$ =
 $( \frac{1}{7},
\frac{2}{7},
\frac{4}{7},
\frac{1}{14},
\frac{9}{14},
\frac{11}{14} )
$,

\vskip 0.7ex
\hangindent=5.5em \hangafter=1
{\white .}\hskip 1em $\rho_\text{isum}(\mathfrak{s})$ =
($\frac{1}{2\sqrt{7}}c^{1}_{28}
$,
$-\frac{1}{2\sqrt{7}}c^{3}_{28}
$,
$\frac{1}{2\sqrt{7}}c^{5}_{28}
$,
$-\frac{3}{2\sqrt{21}}c^{5}_{28}
$,
$\frac{3}{2\sqrt{21}}c^{1}_{28}
$,
$\frac{3}{2\sqrt{21}}c^{3}_{28}
$;
$-\frac{1}{2\sqrt{7}}c^{5}_{28}
$,
$\frac{1}{2\sqrt{7}}c^{1}_{28}
$,
$-\frac{3}{2\sqrt{21}}c^{1}_{28}
$,
$-\frac{3}{2\sqrt{21}}c^{3}_{28}
$,
$\frac{3}{2\sqrt{21}}c^{5}_{28}
$;
$\frac{1}{2\sqrt{7}}c^{3}_{28}
$,
$-\frac{3}{2\sqrt{21}}c^{3}_{28}
$,
$\frac{3}{2\sqrt{21}}c^{5}_{28}
$,
$-\frac{3}{2\sqrt{21}}c^{1}_{28}
$;
$-\frac{1}{2\sqrt{7}}c^{3}_{28}
$,
$\frac{1}{2\sqrt{7}}c^{5}_{28}
$,
$-\frac{1}{2\sqrt{7}}c^{1}_{28}
$;
$-\frac{1}{2\sqrt{7}}c^{1}_{28}
$,
$-\frac{1}{2\sqrt{7}}c^{3}_{28}
$;
$\frac{1}{2\sqrt{7}}c^{5}_{28}
$)

Fail:
cnd($\rho(\mathfrak s)_\mathrm{ndeg}$) = 84 does not divide
 ord($\rho(\mathfrak t)$)=14. Prop. B.4 (2)

 \ \color{black}

\noindent 534: (dims,levels) = $(6;14
)$,
irreps = $3_{7}^{3}
\hskip -1.5pt \otimes \hskip -1.5pt
2_{2}^{1,0}$,
pord$(\rho_\text{isum}(\mathfrak{t})) = 14$,

\vskip 0.7ex
\hangindent=5.5em \hangafter=1
{\white .}\hskip 1em $\rho_\text{isum}(\mathfrak{t})$ =
 $( \frac{3}{7},
\frac{5}{7},
\frac{6}{7},
\frac{3}{14},
\frac{5}{14},
\frac{13}{14} )
$,

\vskip 0.7ex
\hangindent=5.5em \hangafter=1
{\white .}\hskip 1em $\rho_\text{isum}(\mathfrak{s})$ =
($\frac{1}{2\sqrt{7}}c^{3}_{28}
$,
$-\frac{1}{2\sqrt{7}}c^{1}_{28}
$,
$\frac{1}{2\sqrt{7}}c^{5}_{28}
$,
$\frac{3}{2\sqrt{21}}c^{1}_{28}
$,
$-\frac{3}{2\sqrt{21}}c^{5}_{28}
$,
$\frac{3}{2\sqrt{21}}c^{3}_{28}
$;
$-\frac{1}{2\sqrt{7}}c^{5}_{28}
$,
$\frac{1}{2\sqrt{7}}c^{3}_{28}
$,
$\frac{3}{2\sqrt{21}}c^{5}_{28}
$,
$-\frac{3}{2\sqrt{21}}c^{3}_{28}
$,
$-\frac{3}{2\sqrt{21}}c^{1}_{28}
$;
$\frac{1}{2\sqrt{7}}c^{1}_{28}
$,
$-\frac{3}{2\sqrt{21}}c^{3}_{28}
$,
$-\frac{3}{2\sqrt{21}}c^{1}_{28}
$,
$\frac{3}{2\sqrt{21}}c^{5}_{28}
$;
$\frac{1}{2\sqrt{7}}c^{5}_{28}
$,
$-\frac{1}{2\sqrt{7}}c^{3}_{28}
$,
$-\frac{1}{2\sqrt{7}}c^{1}_{28}
$;
$-\frac{1}{2\sqrt{7}}c^{1}_{28}
$,
$\frac{1}{2\sqrt{7}}c^{5}_{28}
$;
$-\frac{1}{2\sqrt{7}}c^{3}_{28}
$)

Fail:
cnd($\rho(\mathfrak s)_\mathrm{ndeg}$) = 84 does not divide
 ord($\rho(\mathfrak t)$)=14. Prop. B.4 (2)

 \ \color{black}

\noindent 535: (dims,levels) = $(6;15
)$,
irreps = $3_{5}^{1}
\hskip -1.5pt \otimes \hskip -1.5pt
2_{3}^{1,0}$,
pord$(\rho_\text{isum}(\mathfrak{t})) = 15$,

\vskip 0.7ex
\hangindent=5.5em \hangafter=1
{\white .}\hskip 1em $\rho_\text{isum}(\mathfrak{t})$ =
 $( 0,
\frac{1}{3},
\frac{1}{5},
\frac{4}{5},
\frac{2}{15},
\frac{8}{15} )
$,

\vskip 0.7ex
\hangindent=5.5em \hangafter=1
{\white .}\hskip 1em $\rho_\text{isum}(\mathfrak{s})$ =
$\mathrm{i}$($-\sqrt{\frac{1}{15}}$,
$-\sqrt{\frac{2}{15}}$,
$-\sqrt{\frac{2}{15}}$,
$-\sqrt{\frac{2}{15}}$,
$\sqrt{\frac{4}{15}}$,
$\sqrt{\frac{4}{15}}$;\ \ 
$\sqrt{\frac{1}{15}}$,
$-\sqrt{\frac{4}{15}}$,
$-\sqrt{\frac{4}{15}}$,
$-\sqrt{\frac{2}{15}}$,
$-\sqrt{\frac{2}{15}}$;\ \ 
$\frac{1}{\sqrt{15}\mathrm{i}}s^{3}_{20}
$,
$-\frac{1}{\sqrt{15}}c^{1}_{5}
$,
$\frac{2}{\sqrt{30}}c^{1}_{5}
$,
$-\frac{2}{\sqrt{30}\mathrm{i}}s^{3}_{20}
$;\ \ 
$\frac{1}{\sqrt{15}\mathrm{i}}s^{3}_{20}
$,
$-\frac{2}{\sqrt{30}\mathrm{i}}s^{3}_{20}
$,
$\frac{2}{\sqrt{30}}c^{1}_{5}
$;\ \ 
$-\frac{1}{\sqrt{15}\mathrm{i}}s^{3}_{20}
$,
$\frac{1}{\sqrt{15}}c^{1}_{5}
$;\ \ 
$-\frac{1}{\sqrt{15}\mathrm{i}}s^{3}_{20}
$)

Fail:
cnd($\rho(\mathfrak s)_\mathrm{ndeg}$) = 120 does not divide
 ord($\rho(\mathfrak t)$)=15. Prop. B.4 (2)

 \ \color{black}

\noindent 536: (dims,levels) = $(6;15
)$,
irreps = $3_{5}^{3}
\hskip -1.5pt \otimes \hskip -1.5pt
2_{3}^{1,0}$,
pord$(\rho_\text{isum}(\mathfrak{t})) = 15$,

\vskip 0.7ex
\hangindent=5.5em \hangafter=1
{\white .}\hskip 1em $\rho_\text{isum}(\mathfrak{t})$ =
 $( 0,
\frac{1}{3},
\frac{2}{5},
\frac{3}{5},
\frac{11}{15},
\frac{14}{15} )
$,

\vskip 0.7ex
\hangindent=5.5em \hangafter=1
{\white .}\hskip 1em $\rho_\text{isum}(\mathfrak{s})$ =
$\mathrm{i}$($\sqrt{\frac{1}{15}}$,
$-\sqrt{\frac{2}{15}}$,
$-\sqrt{\frac{2}{15}}$,
$-\sqrt{\frac{2}{15}}$,
$\sqrt{\frac{4}{15}}$,
$\sqrt{\frac{4}{15}}$;\ \ 
$-\sqrt{\frac{1}{15}}$,
$\sqrt{\frac{4}{15}}$,
$\sqrt{\frac{4}{15}}$,
$\sqrt{\frac{2}{15}}$,
$\sqrt{\frac{2}{15}}$;\ \ 
$\frac{1}{\sqrt{15}}c^{1}_{5}
$,
$-\frac{1}{\sqrt{15}\mathrm{i}}s^{3}_{20}
$,
$-\frac{2}{\sqrt{30}}c^{1}_{5}
$,
$\frac{2}{\sqrt{30}\mathrm{i}}s^{3}_{20}
$;\ \ 
$\frac{1}{\sqrt{15}}c^{1}_{5}
$,
$\frac{2}{\sqrt{30}\mathrm{i}}s^{3}_{20}
$,
$-\frac{2}{\sqrt{30}}c^{1}_{5}
$;\ \ 
$-\frac{1}{\sqrt{15}}c^{1}_{5}
$,
$\frac{1}{\sqrt{15}\mathrm{i}}s^{3}_{20}
$;\ \ 
$-\frac{1}{\sqrt{15}}c^{1}_{5}
$)

Fail:
cnd($\rho(\mathfrak s)_\mathrm{ndeg}$) = 120 does not divide
 ord($\rho(\mathfrak t)$)=15. Prop. B.4 (2)

 \ \color{black}

 \color{blue}

\noindent 537: (dims,levels) = $(6;15
)$,
irreps = $3_{3}^{1,0}
\hskip -1.5pt \otimes \hskip -1.5pt
2_{5}^{1}$,
pord$(\rho_\text{isum}(\mathfrak{t})) = 15$,

\vskip 0.7ex
\hangindent=5.5em \hangafter=1
{\white .}\hskip 1em $\rho_\text{isum}(\mathfrak{t})$ =
 $( \frac{1}{5},
\frac{4}{5},
\frac{2}{15},
\frac{7}{15},
\frac{8}{15},
\frac{13}{15} )
$,

\vskip 0.7ex
\hangindent=5.5em \hangafter=1
{\white .}\hskip 1em $\rho_\text{isum}(\mathfrak{s})$ =
$\mathrm{i}$($\frac{1}{3\sqrt{5}}c^{3}_{20}
$,
$\frac{1}{3\sqrt{5}}c^{1}_{20}
$,
$\frac{2}{3\sqrt{5}}c^{1}_{20}
$,
$\frac{2}{3\sqrt{5}}c^{1}_{20}
$,
$\frac{2}{3\sqrt{5}}c^{3}_{20}
$,
$\frac{2}{3\sqrt{5}}c^{3}_{20}
$;\ \ 
$-\frac{1}{3\sqrt{5}}c^{3}_{20}
$,
$-\frac{2}{3\sqrt{5}}c^{3}_{20}
$,
$-\frac{2}{3\sqrt{5}}c^{3}_{20}
$,
$\frac{2}{3\sqrt{5}}c^{1}_{20}
$,
$\frac{2}{3\sqrt{5}}c^{1}_{20}
$;\ \ 
$-\frac{1}{3\sqrt{5}}c^{3}_{20}
$,
$\frac{2}{3\sqrt{5}}c^{3}_{20}
$,
$\frac{1}{3\sqrt{5}}c^{1}_{20}
$,
$-\frac{2}{3\sqrt{5}}c^{1}_{20}
$;\ \ 
$-\frac{1}{3\sqrt{5}}c^{3}_{20}
$,
$-\frac{2}{3\sqrt{5}}c^{1}_{20}
$,
$\frac{1}{3\sqrt{5}}c^{1}_{20}
$;\ \ 
$\frac{1}{3\sqrt{5}}c^{3}_{20}
$,
$-\frac{2}{3\sqrt{5}}c^{3}_{20}
$;\ \ 
$\frac{1}{3\sqrt{5}}c^{3}_{20}
$)

Pass. 

 \ \color{black}

\noindent 538: (dims,levels) = $(6;15
)$,
irreps = $6_{5}^{1}
\hskip -1.5pt \otimes \hskip -1.5pt
1_{3}^{1,0}$,
pord$(\rho_\text{isum}(\mathfrak{t})) = 5$,

\vskip 0.7ex
\hangindent=5.5em \hangafter=1
{\white .}\hskip 1em $\rho_\text{isum}(\mathfrak{t})$ =
 $( \frac{1}{3},
\frac{1}{3},
\frac{2}{15},
\frac{8}{15},
\frac{11}{15},
\frac{14}{15} )
$,

\vskip 0.7ex
\hangindent=5.5em \hangafter=1
{\white .}\hskip 1em $\rho_\text{isum}(\mathfrak{s})$ =
$\mathrm{i}$($-\frac{1}{5}c^{1}_{20}
$,
$\frac{1}{5}c^{3}_{20}
$,
$\frac{4}{5\sqrt{10}}c^{1}_{20}
-\frac{2}{5\sqrt{10}}c^{3}_{20}
$,
$\frac{1}{5}c^{3}_{40}
+\frac{1}{5}c^{7}_{40}
$,
$\frac{1}{5}c^{1}_{20}
+\frac{1}{5}c^{3}_{20}
$,
$\frac{1}{5}c^{1}_{20}
-\frac{1}{5}c^{3}_{20}
$;\ \ 
$\frac{1}{5}c^{1}_{20}
$,
$\frac{1}{5}c^{3}_{40}
+\frac{1}{5}c^{7}_{40}
$,
$-\frac{4}{5\sqrt{10}}c^{1}_{20}
+\frac{2}{5\sqrt{10}}c^{3}_{20}
$,
$\frac{1}{5}c^{1}_{20}
-\frac{1}{5}c^{3}_{20}
$,
$-\frac{1}{5}c^{1}_{20}
-\frac{1}{5}c^{3}_{20}
$;\ \ 
$\frac{1}{5}c^{3}_{20}
$,
$\frac{1}{5}c^{1}_{20}
$,
$-\frac{1}{5}c^{3}_{40}
-\frac{1}{5}c^{7}_{40}
$,
$\frac{4}{5\sqrt{10}}c^{1}_{20}
-\frac{2}{5\sqrt{10}}c^{3}_{20}
$;\ \ 
$-\frac{1}{5}c^{3}_{20}
$,
$\frac{4}{5\sqrt{10}}c^{1}_{20}
-\frac{2}{5\sqrt{10}}c^{3}_{20}
$,
$\frac{1}{5}c^{3}_{40}
+\frac{1}{5}c^{7}_{40}
$;\ \ 
$\frac{1}{5}c^{1}_{20}
$,
$-\frac{1}{5}c^{3}_{20}
$;\ \ 
$-\frac{1}{5}c^{1}_{20}
$)

Fail:
cnd($\rho(\mathfrak s)_\mathrm{ndeg}$) = 40 does not divide
 ord($\rho(\mathfrak t)$)=15. Prop. B.4 (2)

 \ \color{black}

 \color{blue}

\noindent 539: (dims,levels) = $(6;15
)$,
irreps = $3_{3}^{1,0}
\hskip -1.5pt \otimes \hskip -1.5pt
2_{5}^{2}$,
pord$(\rho_\text{isum}(\mathfrak{t})) = 15$,

\vskip 0.7ex
\hangindent=5.5em \hangafter=1
{\white .}\hskip 1em $\rho_\text{isum}(\mathfrak{t})$ =
 $( \frac{2}{5},
\frac{3}{5},
\frac{1}{15},
\frac{4}{15},
\frac{11}{15},
\frac{14}{15} )
$,

\vskip 0.7ex
\hangindent=5.5em \hangafter=1
{\white .}\hskip 1em $\rho_\text{isum}(\mathfrak{s})$ =
$\mathrm{i}$($\frac{1}{3\sqrt{5}}c^{1}_{20}
$,
$\frac{1}{3\sqrt{5}}c^{3}_{20}
$,
$\frac{2}{3\sqrt{5}}c^{1}_{20}
$,
$\frac{2}{3\sqrt{5}}c^{3}_{20}
$,
$\frac{2}{3\sqrt{5}}c^{1}_{20}
$,
$\frac{2}{3\sqrt{5}}c^{3}_{20}
$;\ \ 
$-\frac{1}{3\sqrt{5}}c^{1}_{20}
$,
$\frac{2}{3\sqrt{5}}c^{3}_{20}
$,
$-\frac{2}{3\sqrt{5}}c^{1}_{20}
$,
$\frac{2}{3\sqrt{5}}c^{3}_{20}
$,
$-\frac{2}{3\sqrt{5}}c^{1}_{20}
$;\ \ 
$\frac{1}{3\sqrt{5}}c^{1}_{20}
$,
$\frac{1}{3\sqrt{5}}c^{3}_{20}
$,
$-\frac{2}{3\sqrt{5}}c^{1}_{20}
$,
$-\frac{2}{3\sqrt{5}}c^{3}_{20}
$;\ \ 
$-\frac{1}{3\sqrt{5}}c^{1}_{20}
$,
$-\frac{2}{3\sqrt{5}}c^{3}_{20}
$,
$\frac{2}{3\sqrt{5}}c^{1}_{20}
$;\ \ 
$\frac{1}{3\sqrt{5}}c^{1}_{20}
$,
$\frac{1}{3\sqrt{5}}c^{3}_{20}
$;\ \ 
$-\frac{1}{3\sqrt{5}}c^{1}_{20}
$)

Pass. 

 \ \color{black}

 \color{blue}

\noindent 540: (dims,levels) = $(6;16
)$,
irreps = $6_{16,1}^{1,0}$,
pord$(\rho_\text{isum}(\mathfrak{t})) = 16$,

\vskip 0.7ex
\hangindent=5.5em \hangafter=1
{\white .}\hskip 1em $\rho_\text{isum}(\mathfrak{t})$ =
 $( 0,
\frac{1}{4},
\frac{1}{16},
\frac{5}{16},
\frac{9}{16},
\frac{13}{16} )
$,

\vskip 0.7ex
\hangindent=5.5em \hangafter=1
{\white .}\hskip 1em $\rho_\text{isum}(\mathfrak{s})$ =
$\mathrm{i}$($0$,
$0$,
$\frac{1}{2}$,
$\frac{1}{2}$,
$\frac{1}{2}$,
$\frac{1}{2}$;\ \ 
$0$,
$\frac{1}{2}$,
$-\frac{1}{2}$,
$\frac{1}{2}$,
$-\frac{1}{2}$;\ \ 
$-\sqrt{\frac{1}{8}}$,
$-\sqrt{\frac{1}{8}}$,
$\sqrt{\frac{1}{8}}$,
$\sqrt{\frac{1}{8}}$;\ \ 
$\sqrt{\frac{1}{8}}$,
$\sqrt{\frac{1}{8}}$,
$-\sqrt{\frac{1}{8}}$;\ \ 
$-\sqrt{\frac{1}{8}}$,
$-\sqrt{\frac{1}{8}}$;\ \ 
$\sqrt{\frac{1}{8}}$)

Pass. 

 \ \color{black}

\noindent 541: (dims,levels) = $(6;16
)$,
irreps = $6_{16,2}^{1,0}$,
pord$(\rho_\text{isum}(\mathfrak{t})) = 16$,

\vskip 0.7ex
\hangindent=5.5em \hangafter=1
{\white .}\hskip 1em $\rho_\text{isum}(\mathfrak{t})$ =
 $( 0,
\frac{1}{2},
\frac{1}{16},
\frac{3}{16},
\frac{9}{16},
\frac{11}{16} )
$,

\vskip 0.7ex
\hangindent=5.5em \hangafter=1
{\white .}\hskip 1em $\rho_\text{isum}(\mathfrak{s})$ =
$\mathrm{i}$($0$,
$0$,
$\frac{1}{2}$,
$\frac{1}{2}$,
$\frac{1}{2}$,
$\frac{1}{2}$;\ \ 
$0$,
$\frac{1}{2}$,
$-\frac{1}{2}$,
$\frac{1}{2}$,
$-\frac{1}{2}$;\ \ 
$0$,
$-\frac{1}{2}$,
$0$,
$\frac{1}{2}$;\ \ 
$0$,
$\frac{1}{2}$,
$0$;\ \ 
$0$,
$-\frac{1}{2}$;\ \ 
$0$)

Fail:
all rows of $U \rho(\mathfrak s) U^\dagger$
 contain zero for any block-diagonal $U$. Prop. B.5 (4) eqn. (B.27)

 \ \color{black}

\noindent 542: (dims,levels) = $(6;16
)$,
irreps = $6_{16,3}^{1,0}$,
pord$(\rho_\text{isum}(\mathfrak{t})) = 16$,

\vskip 0.7ex
\hangindent=5.5em \hangafter=1
{\white .}\hskip 1em $\rho_\text{isum}(\mathfrak{t})$ =
 $( 0,
\frac{1}{2},
\frac{1}{16},
\frac{7}{16},
\frac{9}{16},
\frac{15}{16} )
$,

\vskip 0.7ex
\hangindent=5.5em \hangafter=1
{\white .}\hskip 1em $\rho_\text{isum}(\mathfrak{s})$ =
($0$,
$0$,
$\frac{1}{2}$,
$\frac{1}{2}$,
$\frac{1}{2}$,
$\frac{1}{2}$;
$0$,
$\frac{1}{2}$,
$-\frac{1}{2}$,
$\frac{1}{2}$,
$-\frac{1}{2}$;
$0$,
$-\frac{1}{2}$,
$0$,
$\frac{1}{2}$;
$0$,
$\frac{1}{2}$,
$0$;
$0$,
$-\frac{1}{2}$;
$0$)

Fail:
all rows of $U \rho(\mathfrak s) U^\dagger$
 contain zero for any block-diagonal $U$. Prop. B.5 (4) eqn. (B.27)

 \ \color{black}

\noindent 543: (dims,levels) = $(6;16
)$,
irreps = $6_{16,3}^{3,0}$,
pord$(\rho_\text{isum}(\mathfrak{t})) = 16$,

\vskip 0.7ex
\hangindent=5.5em \hangafter=1
{\white .}\hskip 1em $\rho_\text{isum}(\mathfrak{t})$ =
 $( 0,
\frac{1}{2},
\frac{3}{16},
\frac{5}{16},
\frac{11}{16},
\frac{13}{16} )
$,

\vskip 0.7ex
\hangindent=5.5em \hangafter=1
{\white .}\hskip 1em $\rho_\text{isum}(\mathfrak{s})$ =
($0$,
$0$,
$\frac{1}{2}$,
$\frac{1}{2}$,
$\frac{1}{2}$,
$\frac{1}{2}$;
$0$,
$\frac{1}{2}$,
$-\frac{1}{2}$,
$\frac{1}{2}$,
$-\frac{1}{2}$;
$0$,
$-\frac{1}{2}$,
$0$,
$\frac{1}{2}$;
$0$,
$\frac{1}{2}$,
$0$;
$0$,
$-\frac{1}{2}$;
$0$)

Fail:
all rows of $U \rho(\mathfrak s) U^\dagger$
 contain zero for any block-diagonal $U$. Prop. B.5 (4) eqn. (B.27)

 \ \color{black}

\noindent 544: (dims,levels) = $(6;16
)$,
irreps = $6_{16,2}^{5,0}$,
pord$(\rho_\text{isum}(\mathfrak{t})) = 16$,

\vskip 0.7ex
\hangindent=5.5em \hangafter=1
{\white .}\hskip 1em $\rho_\text{isum}(\mathfrak{t})$ =
 $( 0,
\frac{1}{2},
\frac{5}{16},
\frac{7}{16},
\frac{13}{16},
\frac{15}{16} )
$,

\vskip 0.7ex
\hangindent=5.5em \hangafter=1
{\white .}\hskip 1em $\rho_\text{isum}(\mathfrak{s})$ =
$\mathrm{i}$($0$,
$0$,
$\frac{1}{2}$,
$\frac{1}{2}$,
$\frac{1}{2}$,
$\frac{1}{2}$;\ \ 
$0$,
$\frac{1}{2}$,
$-\frac{1}{2}$,
$\frac{1}{2}$,
$-\frac{1}{2}$;\ \ 
$0$,
$\frac{1}{2}$,
$0$,
$-\frac{1}{2}$;\ \ 
$0$,
$-\frac{1}{2}$,
$0$;\ \ 
$0$,
$\frac{1}{2}$;\ \ 
$0$)

Fail:
all rows of $U \rho(\mathfrak s) U^\dagger$
 contain zero for any block-diagonal $U$. Prop. B.5 (4) eqn. (B.27)

 \ \color{black}

 \color{blue}

\noindent 545: (dims,levels) = $(6;16
)$,
irreps = $6_{16,1}^{3,0}$,
pord$(\rho_\text{isum}(\mathfrak{t})) = 16$,

\vskip 0.7ex
\hangindent=5.5em \hangafter=1
{\white .}\hskip 1em $\rho_\text{isum}(\mathfrak{t})$ =
 $( 0,
\frac{3}{4},
\frac{3}{16},
\frac{7}{16},
\frac{11}{16},
\frac{15}{16} )
$,

\vskip 0.7ex
\hangindent=5.5em \hangafter=1
{\white .}\hskip 1em $\rho_\text{isum}(\mathfrak{s})$ =
$\mathrm{i}$($0$,
$0$,
$\frac{1}{2}$,
$\frac{1}{2}$,
$\frac{1}{2}$,
$\frac{1}{2}$;\ \ 
$0$,
$\frac{1}{2}$,
$-\frac{1}{2}$,
$\frac{1}{2}$,
$-\frac{1}{2}$;\ \ 
$-\sqrt{\frac{1}{8}}$,
$\sqrt{\frac{1}{8}}$,
$\sqrt{\frac{1}{8}}$,
$-\sqrt{\frac{1}{8}}$;\ \ 
$\sqrt{\frac{1}{8}}$,
$-\sqrt{\frac{1}{8}}$,
$-\sqrt{\frac{1}{8}}$;\ \ 
$-\sqrt{\frac{1}{8}}$,
$\sqrt{\frac{1}{8}}$;\ \ 
$\sqrt{\frac{1}{8}}$)

Pass. 

 \ \color{black}

\noindent 546: (dims,levels) = $(6;16
)$,
irreps = $6_{16,4}^{1,0}$,
pord$(\rho_\text{isum}(\mathfrak{t})) = 16$,

\vskip 0.7ex
\hangindent=5.5em \hangafter=1
{\white .}\hskip 1em $\rho_\text{isum}(\mathfrak{t})$ =
 $( \frac{1}{8},
\frac{5}{8},
\frac{1}{16},
\frac{5}{16},
\frac{9}{16},
\frac{13}{16} )
$,

\vskip 0.7ex
\hangindent=5.5em \hangafter=1
{\white .}\hskip 1em $\rho_\text{isum}(\mathfrak{s})$ =
($0$,
$0$,
$\frac{1}{2}$,
$\frac{1}{2}$,
$\frac{1}{2}$,
$\frac{1}{2}$;
$0$,
$\frac{1}{2}$,
$-\frac{1}{2}$,
$\frac{1}{2}$,
$-\frac{1}{2}$;
$0$,
$-\frac{1}{2}$,
$0$,
$\frac{1}{2}$;
$0$,
$\frac{1}{2}$,
$0$;
$0$,
$-\frac{1}{2}$;
$0$)

Fail:
all rows of $U \rho(\mathfrak s) U^\dagger$
 contain zero for any block-diagonal $U$. Prop. B.5 (4) eqn. (B.27)

 \ \color{black}

\noindent 547: (dims,levels) = $(6;16
)$,
irreps = $6_{16,4}^{3,0}$,
pord$(\rho_\text{isum}(\mathfrak{t})) = 16$,

\vskip 0.7ex
\hangindent=5.5em \hangafter=1
{\white .}\hskip 1em $\rho_\text{isum}(\mathfrak{t})$ =
 $( \frac{3}{8},
\frac{7}{8},
\frac{3}{16},
\frac{7}{16},
\frac{11}{16},
\frac{15}{16} )
$,

\vskip 0.7ex
\hangindent=5.5em \hangafter=1
{\white .}\hskip 1em $\rho_\text{isum}(\mathfrak{s})$ =
($0$,
$0$,
$\frac{1}{2}$,
$\frac{1}{2}$,
$\frac{1}{2}$,
$\frac{1}{2}$;
$0$,
$\frac{1}{2}$,
$-\frac{1}{2}$,
$\frac{1}{2}$,
$-\frac{1}{2}$;
$0$,
$\frac{1}{2}$,
$0$,
$-\frac{1}{2}$;
$0$,
$-\frac{1}{2}$,
$0$;
$0$,
$\frac{1}{2}$;
$0$)

Fail:
all rows of $U \rho(\mathfrak s) U^\dagger$
 contain zero for any block-diagonal $U$. Prop. B.5 (4) eqn. (B.27)

 \ \color{black}

 \color{blue}

\noindent 548: (dims,levels) = $(6;18
)$,
irreps = $6_{9,3}^{1,0}
\hskip -1.5pt \otimes \hskip -1.5pt
1_{2}^{1,0}$,
pord$(\rho_\text{isum}(\mathfrak{t})) = 9$,

\vskip 0.7ex
\hangindent=5.5em \hangafter=1
{\white .}\hskip 1em $\rho_\text{isum}(\mathfrak{t})$ =
 $( \frac{1}{18},
\frac{5}{18},
\frac{7}{18},
\frac{11}{18},
\frac{13}{18},
\frac{17}{18} )
$,

\vskip 0.7ex
\hangindent=5.5em \hangafter=1
{\white .}\hskip 1em $\rho_\text{isum}(\mathfrak{s})$ =
($-\frac{1}{3}$,
$\frac{1}{3} c_9^4 $,
$\frac{1}{3}$,
$-\frac{1}{3}c^{1}_{9}
$,
$\frac{1}{3}$,
$\frac{1}{3}c^{2}_{9}
$;
$-\frac{1}{3}$,
$-\frac{1}{3}c^{2}_{9}
$,
$\frac{1}{3}$,
$-\frac{1}{3}c^{1}_{9}
$,
$-\frac{1}{3}$;
$-\frac{1}{3}$,
$\frac{1}{3} c_9^4 $,
$-\frac{1}{3}$,
$-\frac{1}{3}c^{1}_{9}
$;
$-\frac{1}{3}$,
$\frac{1}{3}c^{2}_{9}
$,
$\frac{1}{3}$;
$-\frac{1}{3}$,
$-\frac{1}{3} c_9^4 $;
$-\frac{1}{3}$)

Pass. 

 \ \color{black}

\noindent 549: (dims,levels) = $(6;18
)$,
irreps = $6_{9,2}^{1,0}
\hskip -1.5pt \otimes \hskip -1.5pt
1_{2}^{1,0}$,
pord$(\rho_\text{isum}(\mathfrak{t})) = 9$,

\vskip 0.7ex
\hangindent=5.5em \hangafter=1
{\white .}\hskip 1em $\rho_\text{isum}(\mathfrak{t})$ =
 $( \frac{1}{18},
\frac{5}{18},
\frac{7}{18},
\frac{11}{18},
\frac{13}{18},
\frac{17}{18} )
$,

\vskip 0.7ex
\hangindent=5.5em \hangafter=1
{\white .}\hskip 1em $\rho_\text{isum}(\mathfrak{s})$ =
($\frac{1}{3}$,
$-\frac{1}{3}c^{1}_{36}
+\frac{1}{3}c^{5}_{36}
$,
$\frac{1}{3}$,
$\frac{1}{3}c^{5}_{36}
$,
$\frac{1}{3}$,
$\frac{1}{3}c^{1}_{36}
$;
$\frac{1}{3}$,
$\frac{1}{3}c^{1}_{36}
$,
$-\frac{1}{3}$,
$-\frac{1}{3}c^{5}_{36}
$,
$\frac{1}{3}$;
$\frac{1}{3}$,
$\frac{1}{3}c^{1}_{36}
-\frac{1}{3}c^{5}_{36}
$,
$\frac{1}{3}$,
$-\frac{1}{3}c^{5}_{36}
$;
$\frac{1}{3}$,
$-\frac{1}{3}c^{1}_{36}
$,
$-\frac{1}{3}$;
$\frac{1}{3}$,
$-\frac{1}{3}c^{1}_{36}
+\frac{1}{3}c^{5}_{36}
$;
$\frac{1}{3}$)

Fail:
cnd($\rho(\mathfrak s)_\mathrm{ndeg}$) = 36 does not divide
 ord($\rho(\mathfrak t)$)=18. Prop. B.4 (2)

 \ \color{black}

\noindent 550: (dims,levels) = $(6;18
)$,
irreps = $6_{9,1}^{1,0}
\hskip -1.5pt \otimes \hskip -1.5pt
1_{2}^{1,0}$,
pord$(\rho_\text{isum}(\mathfrak{t})) = 9$,

\vskip 0.7ex
\hangindent=5.5em \hangafter=1
{\white .}\hskip 1em $\rho_\text{isum}(\mathfrak{t})$ =
 $( \frac{1}{18},
\frac{5}{18},
\frac{7}{18},
\frac{11}{18},
\frac{13}{18},
\frac{17}{18} )
$,

\vskip 0.7ex
\hangindent=5.5em \hangafter=1
{\white .}\hskip 1em $\rho_\text{isum}(\mathfrak{s})$ =
$\mathrm{i}$($-\frac{1}{3}$,
$-\frac{1}{3}c^{5}_{72}
$,
$\frac{1}{3}$,
$-\frac{1}{3}c^{5}_{72}
+\frac{1}{3}c^{7}_{72}
$,
$\frac{1}{3}$,
$\frac{1}{3}c^{7}_{72}
$;\ \ 
$\frac{1}{3}$,
$-\frac{1}{3}c^{7}_{72}
$,
$-\frac{1}{3}$,
$-\frac{1}{3}c^{5}_{72}
+\frac{1}{3}c^{7}_{72}
$,
$\frac{1}{3}$;\ \ 
$-\frac{1}{3}$,
$-\frac{1}{3}c^{5}_{72}
$,
$-\frac{1}{3}$,
$-\frac{1}{3}c^{5}_{72}
+\frac{1}{3}c^{7}_{72}
$;\ \ 
$\frac{1}{3}$,
$\frac{1}{3}c^{7}_{72}
$,
$-\frac{1}{3}$;\ \ 
$-\frac{1}{3}$,
$\frac{1}{3}c^{5}_{72}
$;\ \ 
$\frac{1}{3}$)

Fail:
cnd( Tr$_I(\rho(\mathfrak s))$ ) =
4 does not divide
 ord($\rho(\mathfrak t)$) =
18, I = [ 1/18 ]. Prop. B.4 (2)

 \ \color{black}

\noindent 551: (dims,levels) = $(6;18
)$,
irreps = $6_{9,1}^{5,0}
\hskip -1.5pt \otimes \hskip -1.5pt
1_{2}^{1,0}$,
pord$(\rho_\text{isum}(\mathfrak{t})) = 9$,

\vskip 0.7ex
\hangindent=5.5em \hangafter=1
{\white .}\hskip 1em $\rho_\text{isum}(\mathfrak{t})$ =
 $( \frac{1}{18},
\frac{5}{18},
\frac{7}{18},
\frac{11}{18},
\frac{13}{18},
\frac{17}{18} )
$,

\vskip 0.7ex
\hangindent=5.5em \hangafter=1
{\white .}\hskip 1em $\rho_\text{isum}(\mathfrak{s})$ =
$\mathrm{i}$($\frac{1}{3}$,
$-\frac{1}{3}c^{1}_{72}
+\frac{1}{3}c^{11}_{72}
$,
$\frac{1}{3}$,
$-\frac{1}{3}c^{1}_{72}
$,
$\frac{1}{3}$,
$-\frac{1}{3}c^{11}_{72}
$;\ \ 
$-\frac{1}{3}$,
$-\frac{1}{3}c^{11}_{72}
$,
$\frac{1}{3}$,
$\frac{1}{3}c^{1}_{72}
$,
$-\frac{1}{3}$;\ \ 
$\frac{1}{3}$,
$\frac{1}{3}c^{1}_{72}
-\frac{1}{3}c^{11}_{72}
$,
$\frac{1}{3}$,
$\frac{1}{3}c^{1}_{72}
$;\ \ 
$-\frac{1}{3}$,
$\frac{1}{3}c^{11}_{72}
$,
$\frac{1}{3}$;\ \ 
$\frac{1}{3}$,
$-\frac{1}{3}c^{1}_{72}
+\frac{1}{3}c^{11}_{72}
$;\ \ 
$-\frac{1}{3}$)

Fail:
cnd( Tr$_I(\rho(\mathfrak s))$ ) =
4 does not divide
 ord($\rho(\mathfrak t)$) =
18, I = [ 1/18 ]. Prop. B.4 (2)

 \ \color{black}

\noindent 552: (dims,levels) = $(6;20
)$,
irreps = $3_{4}^{1,0}
\hskip -1.5pt \otimes \hskip -1.5pt
2_{5}^{1}$,
pord$(\rho_\text{isum}(\mathfrak{t})) = 20$,

\vskip 0.7ex
\hangindent=5.5em \hangafter=1
{\white .}\hskip 1em $\rho_\text{isum}(\mathfrak{t})$ =
 $( \frac{1}{5},
\frac{4}{5},
\frac{1}{20},
\frac{9}{20},
\frac{11}{20},
\frac{19}{20} )
$,

\vskip 0.7ex
\hangindent=5.5em \hangafter=1
{\white .}\hskip 1em $\rho_\text{isum}(\mathfrak{s})$ =
$\mathrm{i}$($0$,
$0$,
$\frac{1}{\sqrt{10}}c^{1}_{20}
$,
$\frac{1}{\sqrt{10}}c^{3}_{20}
$,
$\frac{1}{\sqrt{10}}c^{1}_{20}
$,
$\frac{1}{\sqrt{10}}c^{3}_{20}
$;\ \ 
$0$,
$\frac{1}{\sqrt{10}}c^{3}_{20}
$,
$-\frac{1}{\sqrt{10}}c^{1}_{20}
$,
$\frac{1}{\sqrt{10}}c^{3}_{20}
$,
$-\frac{1}{\sqrt{10}}c^{1}_{20}
$;\ \ 
$-\frac{1}{2\sqrt{5}}c^{3}_{20}
$,
$\frac{1}{2\sqrt{5}}c^{1}_{20}
$,
$\frac{1}{2\sqrt{5}}c^{3}_{20}
$,
$-\frac{1}{2\sqrt{5}}c^{1}_{20}
$;\ \ 
$\frac{1}{2\sqrt{5}}c^{3}_{20}
$,
$-\frac{1}{2\sqrt{5}}c^{1}_{20}
$,
$-\frac{1}{2\sqrt{5}}c^{3}_{20}
$;\ \ 
$-\frac{1}{2\sqrt{5}}c^{3}_{20}
$,
$\frac{1}{2\sqrt{5}}c^{1}_{20}
$;\ \ 
$\frac{1}{2\sqrt{5}}c^{3}_{20}
$)

Fail:
cnd($\rho(\mathfrak s)_\mathrm{ndeg}$) = 40 does not divide
 ord($\rho(\mathfrak t)$)=20. Prop. B.4 (2)

 \ \color{black}

\noindent 553: (dims,levels) = $(6;20
)$,
irreps = $6_{5}^{1}
\hskip -1.5pt \otimes \hskip -1.5pt
1_{4}^{1,0}$,
pord$(\rho_\text{isum}(\mathfrak{t})) = 5$,

\vskip 0.7ex
\hangindent=5.5em \hangafter=1
{\white .}\hskip 1em $\rho_\text{isum}(\mathfrak{t})$ =
 $( \frac{1}{4},
\frac{1}{4},
\frac{1}{20},
\frac{9}{20},
\frac{13}{20},
\frac{17}{20} )
$,

\vskip 0.7ex
\hangindent=5.5em \hangafter=1
{\white .}\hskip 1em $\rho_\text{isum}(\mathfrak{s})$ =
($-\frac{1}{5}c^{1}_{20}
$,
$-\frac{1}{5}c^{3}_{20}
$,
$\frac{1}{5}c^{3}_{40}
+\frac{1}{5}c^{7}_{40}
$,
$-\frac{4}{5\sqrt{10}}c^{1}_{20}
+\frac{2}{5\sqrt{10}}c^{3}_{20}
$,
$\frac{1}{5}c^{1}_{20}
-\frac{1}{5}c^{3}_{20}
$,
$\frac{1}{5}c^{1}_{20}
+\frac{1}{5}c^{3}_{20}
$;
$\frac{1}{5}c^{1}_{20}
$,
$\frac{4}{5\sqrt{10}}c^{1}_{20}
-\frac{2}{5\sqrt{10}}c^{3}_{20}
$,
$\frac{1}{5}c^{3}_{40}
+\frac{1}{5}c^{7}_{40}
$,
$\frac{1}{5}c^{1}_{20}
+\frac{1}{5}c^{3}_{20}
$,
$-\frac{1}{5}c^{1}_{20}
+\frac{1}{5}c^{3}_{20}
$;
$-\frac{1}{5}c^{3}_{20}
$,
$-\frac{1}{5}c^{1}_{20}
$,
$\frac{1}{5}c^{3}_{40}
+\frac{1}{5}c^{7}_{40}
$,
$\frac{4}{5\sqrt{10}}c^{1}_{20}
-\frac{2}{5\sqrt{10}}c^{3}_{20}
$;
$\frac{1}{5}c^{3}_{20}
$,
$-\frac{4}{5\sqrt{10}}c^{1}_{20}
+\frac{2}{5\sqrt{10}}c^{3}_{20}
$,
$\frac{1}{5}c^{3}_{40}
+\frac{1}{5}c^{7}_{40}
$;
$-\frac{1}{5}c^{1}_{20}
$,
$-\frac{1}{5}c^{3}_{20}
$;
$\frac{1}{5}c^{1}_{20}
$)

Fail:
cnd($\rho(\mathfrak s)_\mathrm{ndeg}$) = 40 does not divide
 ord($\rho(\mathfrak t)$)=20. Prop. B.4 (2)

 \ \color{black}

\noindent 554: (dims,levels) = $(6;20
)$,
irreps = $3_{5}^{1}
\hskip -1.5pt \otimes \hskip -1.5pt
2_{4}^{1,0}$,
pord$(\rho_\text{isum}(\mathfrak{t})) = 10$,

\vskip 0.7ex
\hangindent=5.5em \hangafter=1
{\white .}\hskip 1em $\rho_\text{isum}(\mathfrak{t})$ =
 $( \frac{1}{4},
\frac{3}{4},
\frac{1}{20},
\frac{9}{20},
\frac{11}{20},
\frac{19}{20} )
$,

\vskip 0.7ex
\hangindent=5.5em \hangafter=1
{\white .}\hskip 1em $\rho_\text{isum}(\mathfrak{s})$ =
$\mathrm{i}$($-\sqrt{\frac{1}{20}}$,
$\sqrt{\frac{3}{20}}$,
$\sqrt{\frac{1}{10}}$,
$\sqrt{\frac{1}{10}}$,
$-\sqrt{\frac{3}{10}}$,
$-\sqrt{\frac{3}{10}}$;\ \ 
$\sqrt{\frac{1}{20}}$,
$-\sqrt{\frac{3}{10}}$,
$-\sqrt{\frac{3}{10}}$,
$-\sqrt{\frac{1}{10}}$,
$-\sqrt{\frac{1}{10}}$;\ \ 
$\frac{5+\sqrt{5}}{20}$,
$\frac{-5+\sqrt{5}}{20}$,
$-\frac{3}{2\sqrt{15}\mathrm{i}}s^{3}_{20}
$,
$\frac{3}{2\sqrt{15}}c^{1}_{5}
$;\ \ 
$\frac{5+\sqrt{5}}{20}$,
$\frac{3}{2\sqrt{15}}c^{1}_{5}
$,
$-\frac{3}{2\sqrt{15}\mathrm{i}}s^{3}_{20}
$;\ \ 
$-\frac{5+\sqrt{5}}{20}$,
$\frac{5-\sqrt{5}}{20}$;\ \ 
$-\frac{5+\sqrt{5}}{20}$)

Fail:
cnd($\rho(\mathfrak s)_\mathrm{ndeg}$) = 120 does not divide
 ord($\rho(\mathfrak t)$)=20. Prop. B.4 (2)

 \ \color{black}

\noindent 555: (dims,levels) = $(6;20
)$,
irreps = $3_{5}^{3}
\hskip -1.5pt \otimes \hskip -1.5pt
2_{4}^{1,0}$,
pord$(\rho_\text{isum}(\mathfrak{t})) = 10$,

\vskip 0.7ex
\hangindent=5.5em \hangafter=1
{\white .}\hskip 1em $\rho_\text{isum}(\mathfrak{t})$ =
 $( \frac{1}{4},
\frac{3}{4},
\frac{3}{20},
\frac{7}{20},
\frac{13}{20},
\frac{17}{20} )
$,

\vskip 0.7ex
\hangindent=5.5em \hangafter=1
{\white .}\hskip 1em $\rho_\text{isum}(\mathfrak{s})$ =
$\mathrm{i}$($\sqrt{\frac{1}{20}}$,
$\sqrt{\frac{3}{20}}$,
$-\sqrt{\frac{3}{10}}$,
$-\sqrt{\frac{3}{10}}$,
$\sqrt{\frac{1}{10}}$,
$\sqrt{\frac{1}{10}}$;\ \ 
$-\sqrt{\frac{1}{20}}$,
$\sqrt{\frac{1}{10}}$,
$\sqrt{\frac{1}{10}}$,
$\sqrt{\frac{3}{10}}$,
$\sqrt{\frac{3}{10}}$;\ \ 
$\frac{-5+\sqrt{5}}{20}$,
$\frac{5+\sqrt{5}}{20}$,
$-\frac{3}{2\sqrt{15}}c^{1}_{5}
$,
$\frac{3}{2\sqrt{15}\mathrm{i}}s^{3}_{20}
$;\ \ 
$\frac{-5+\sqrt{5}}{20}$,
$\frac{3}{2\sqrt{15}\mathrm{i}}s^{3}_{20}
$,
$-\frac{3}{2\sqrt{15}}c^{1}_{5}
$;\ \ 
$\frac{5-\sqrt{5}}{20}$,
$-\frac{5+\sqrt{5}}{20}$;\ \ 
$\frac{5-\sqrt{5}}{20}$)

Fail:
cnd($\rho(\mathfrak s)_\mathrm{ndeg}$) = 120 does not divide
 ord($\rho(\mathfrak t)$)=20. Prop. B.4 (2)

 \ \color{black}

\noindent 556: (dims,levels) = $(6;20
)$,
irreps = $3_{4}^{1,0}
\hskip -1.5pt \otimes \hskip -1.5pt
2_{5}^{2}$,
pord$(\rho_\text{isum}(\mathfrak{t})) = 20$,

\vskip 0.7ex
\hangindent=5.5em \hangafter=1
{\white .}\hskip 1em $\rho_\text{isum}(\mathfrak{t})$ =
 $( \frac{2}{5},
\frac{3}{5},
\frac{3}{20},
\frac{7}{20},
\frac{13}{20},
\frac{17}{20} )
$,

\vskip 0.7ex
\hangindent=5.5em \hangafter=1
{\white .}\hskip 1em $\rho_\text{isum}(\mathfrak{s})$ =
$\mathrm{i}$($0$,
$0$,
$\frac{1}{\sqrt{10}}c^{1}_{20}
$,
$\frac{1}{\sqrt{10}}c^{3}_{20}
$,
$\frac{1}{\sqrt{10}}c^{1}_{20}
$,
$\frac{1}{\sqrt{10}}c^{3}_{20}
$;\ \ 
$0$,
$\frac{1}{\sqrt{10}}c^{3}_{20}
$,
$-\frac{1}{\sqrt{10}}c^{1}_{20}
$,
$\frac{1}{\sqrt{10}}c^{3}_{20}
$,
$-\frac{1}{\sqrt{10}}c^{1}_{20}
$;\ \ 
$\frac{1}{2\sqrt{5}}c^{1}_{20}
$,
$\frac{1}{2\sqrt{5}}c^{3}_{20}
$,
$-\frac{1}{2\sqrt{5}}c^{1}_{20}
$,
$-\frac{1}{2\sqrt{5}}c^{3}_{20}
$;\ \ 
$-\frac{1}{2\sqrt{5}}c^{1}_{20}
$,
$-\frac{1}{2\sqrt{5}}c^{3}_{20}
$,
$\frac{1}{2\sqrt{5}}c^{1}_{20}
$;\ \ 
$\frac{1}{2\sqrt{5}}c^{1}_{20}
$,
$\frac{1}{2\sqrt{5}}c^{3}_{20}
$;\ \ 
$-\frac{1}{2\sqrt{5}}c^{1}_{20}
$)

Fail:
cnd($\rho(\mathfrak s)_\mathrm{ndeg}$) = 40 does not divide
 ord($\rho(\mathfrak t)$)=20. Prop. B.4 (2)

 \ \color{black}

\noindent 557: (dims,levels) = $(6;21
)$,
irreps = $6_{7,2}^{1}
\hskip -1.5pt \otimes \hskip -1.5pt
1_{3}^{1,0}$,
pord$(\rho_\text{isum}(\mathfrak{t})) = 7$,

\vskip 0.7ex
\hangindent=5.5em \hangafter=1
{\white .}\hskip 1em $\rho_\text{isum}(\mathfrak{t})$ =
 $( \frac{1}{21},
\frac{4}{21},
\frac{10}{21},
\frac{13}{21},
\frac{16}{21},
\frac{19}{21} )
$,

\vskip 0.7ex
\hangindent=5.5em \hangafter=1
{\white .}\hskip 1em $\rho_\text{isum}(\mathfrak{s})$ =
($\frac{3}{7}+\frac{1}{7}c^{1}_{7}
+\frac{1}{7}c^{2}_{7}
$,
$-\frac{2}{7}+\frac{1}{7}c^{1}_{7}
$,
$-\frac{1}{7}c^{3}_{56}
+\frac{1}{7}c^{5}_{56}
-\frac{1}{7}c^{9}_{56}
-\frac{1}{7}c^{11}_{56}
$,
$\frac{1}{7}c^{3}_{56}
+\frac{2}{7}c^{5}_{56}
-\frac{1}{7}c^{7}_{56}
-\frac{2}{7}c^{9}_{56}
+\frac{1}{7}c^{11}_{56}
$,
$-\frac{2}{7}+\frac{1}{7}c^{2}_{7}
$,
$-\frac{2}{7}c^{3}_{56}
-\frac{1}{7}c^{5}_{56}
+\frac{1}{7}c^{7}_{56}
+\frac{1}{7}c^{9}_{56}
-\frac{2}{7}c^{11}_{56}
$;
$\frac{2}{7}-\frac{1}{7}c^{2}_{7}
$,
$-\frac{2}{7}c^{3}_{56}
-\frac{1}{7}c^{5}_{56}
+\frac{1}{7}c^{7}_{56}
+\frac{1}{7}c^{9}_{56}
-\frac{2}{7}c^{11}_{56}
$,
$-\frac{1}{7}c^{3}_{56}
+\frac{1}{7}c^{5}_{56}
-\frac{1}{7}c^{9}_{56}
-\frac{1}{7}c^{11}_{56}
$,
$\frac{3}{7}+\frac{1}{7}c^{1}_{7}
+\frac{1}{7}c^{2}_{7}
$,
$-\frac{1}{7}c^{3}_{56}
-\frac{2}{7}c^{5}_{56}
+\frac{1}{7}c^{7}_{56}
+\frac{2}{7}c^{9}_{56}
-\frac{1}{7}c^{11}_{56}
$;
$\frac{2}{7}-\frac{1}{7}c^{2}_{7}
$,
$-\frac{2}{7}+\frac{1}{7}c^{1}_{7}
$,
$\frac{1}{7}c^{3}_{56}
+\frac{2}{7}c^{5}_{56}
-\frac{1}{7}c^{7}_{56}
-\frac{2}{7}c^{9}_{56}
+\frac{1}{7}c^{11}_{56}
$,
$-\frac{3}{7}-\frac{1}{7}c^{1}_{7}
-\frac{1}{7}c^{2}_{7}
$;
$\frac{3}{7}+\frac{1}{7}c^{1}_{7}
+\frac{1}{7}c^{2}_{7}
$,
$\frac{2}{7}c^{3}_{56}
+\frac{1}{7}c^{5}_{56}
-\frac{1}{7}c^{7}_{56}
-\frac{1}{7}c^{9}_{56}
+\frac{2}{7}c^{11}_{56}
$,
$\frac{2}{7}-\frac{1}{7}c^{2}_{7}
$;
$\frac{2}{7}-\frac{1}{7}c^{1}_{7}
$,
$-\frac{1}{7}c^{3}_{56}
+\frac{1}{7}c^{5}_{56}
-\frac{1}{7}c^{9}_{56}
-\frac{1}{7}c^{11}_{56}
$;
$\frac{2}{7}-\frac{1}{7}c^{1}_{7}
$)

Fail:
cnd($\rho(\mathfrak s)_\mathrm{ndeg}$) = 56 does not divide
 ord($\rho(\mathfrak t)$)=21. Prop. B.4 (2)

 \ \color{black}

\noindent 558: (dims,levels) = $(6;21
)$,
irreps = $6_{7,1}^{1}
\hskip -1.5pt \otimes \hskip -1.5pt
1_{3}^{1,0}$,
pord$(\rho_\text{isum}(\mathfrak{t})) = 7$,

\vskip 0.7ex
\hangindent=5.5em \hangafter=1
{\white .}\hskip 1em $\rho_\text{isum}(\mathfrak{t})$ =
 $( \frac{1}{21},
\frac{4}{21},
\frac{10}{21},
\frac{13}{21},
\frac{16}{21},
\frac{19}{21} )
$,

\vskip 0.7ex
\hangindent=5.5em \hangafter=1
{\white .}\hskip 1em $\rho_\text{isum}(\mathfrak{s})$ =
$\mathrm{i}$($-\frac{2}{7}c^{1}_{56}
+\frac{1}{7}c^{3}_{56}
+\frac{1}{7}c^{5}_{56}
-\frac{1}{7}c^{7}_{56}
-\frac{1}{7}c^{9}_{56}
+\frac{1}{7}c^{10}_{56}
+\frac{1}{7}c^{11}_{56}
$,
$\frac{1}{7}c^{5}_{56}
+\frac{1}{7}c^{6}_{56}
+\frac{1}{7}c^{9}_{56}
$,
$\frac{1}{7}c^{1}_{112}
+\frac{1}{7}c^{3}_{112}
-\frac{1}{7}c^{5}_{112}
-\frac{1}{7}c^{7}_{112}
-\frac{1}{7}c^{9}_{112}
-\frac{1}{7}c^{11}_{112}
+\frac{2}{7}c^{13}_{112}
+\frac{1}{7}c^{15}_{112}
+\frac{2}{7}c^{17}_{112}
+\frac{1}{7}c^{19}_{112}
-\frac{1}{7}c^{21}_{112}
-\frac{1}{7}c^{23}_{112}
$,
$\frac{1}{7}c^{3}_{112}
-\frac{1}{7}c^{9}_{112}
+\frac{1}{7}c^{11}_{112}
+\frac{1}{7}c^{23}_{112}
$,
$\frac{1}{7}c^{2}_{56}
-\frac{1}{7}c^{3}_{56}
+\frac{1}{7}c^{11}_{56}
$,
$\frac{1}{7}c^{1}_{112}
+\frac{1}{7}c^{5}_{112}
-\frac{1}{7}c^{15}_{112}
+\frac{1}{7}c^{19}_{112}
$;\ \ 
$-\frac{1}{7}c^{2}_{56}
+\frac{1}{7}c^{3}_{56}
-\frac{1}{7}c^{11}_{56}
$,
$\frac{1}{7}c^{1}_{112}
+\frac{1}{7}c^{5}_{112}
-\frac{1}{7}c^{15}_{112}
+\frac{1}{7}c^{19}_{112}
$,
$\frac{1}{7}c^{1}_{112}
+\frac{1}{7}c^{3}_{112}
-\frac{1}{7}c^{5}_{112}
-\frac{1}{7}c^{7}_{112}
-\frac{1}{7}c^{9}_{112}
-\frac{1}{7}c^{11}_{112}
+\frac{2}{7}c^{13}_{112}
+\frac{1}{7}c^{15}_{112}
+\frac{2}{7}c^{17}_{112}
+\frac{1}{7}c^{19}_{112}
-\frac{1}{7}c^{21}_{112}
-\frac{1}{7}c^{23}_{112}
$,
$-\frac{2}{7}c^{1}_{56}
+\frac{1}{7}c^{3}_{56}
+\frac{1}{7}c^{5}_{56}
-\frac{1}{7}c^{7}_{56}
-\frac{1}{7}c^{9}_{56}
+\frac{1}{7}c^{10}_{56}
+\frac{1}{7}c^{11}_{56}
$,
$-\frac{1}{7}c^{3}_{112}
+\frac{1}{7}c^{9}_{112}
-\frac{1}{7}c^{11}_{112}
-\frac{1}{7}c^{23}_{112}
$;\ \ 
$\frac{1}{7}c^{2}_{56}
-\frac{1}{7}c^{3}_{56}
+\frac{1}{7}c^{11}_{56}
$,
$-\frac{1}{7}c^{5}_{56}
-\frac{1}{7}c^{6}_{56}
-\frac{1}{7}c^{9}_{56}
$,
$\frac{1}{7}c^{3}_{112}
-\frac{1}{7}c^{9}_{112}
+\frac{1}{7}c^{11}_{112}
+\frac{1}{7}c^{23}_{112}
$,
$-\frac{2}{7}c^{1}_{56}
+\frac{1}{7}c^{3}_{56}
+\frac{1}{7}c^{5}_{56}
-\frac{1}{7}c^{7}_{56}
-\frac{1}{7}c^{9}_{56}
+\frac{1}{7}c^{10}_{56}
+\frac{1}{7}c^{11}_{56}
$;\ \ 
$\frac{2}{7}c^{1}_{56}
-\frac{1}{7}c^{3}_{56}
-\frac{1}{7}c^{5}_{56}
+\frac{1}{7}c^{7}_{56}
+\frac{1}{7}c^{9}_{56}
-\frac{1}{7}c^{10}_{56}
-\frac{1}{7}c^{11}_{56}
$,
$-\frac{1}{7}c^{1}_{112}
-\frac{1}{7}c^{5}_{112}
+\frac{1}{7}c^{15}_{112}
-\frac{1}{7}c^{19}_{112}
$,
$\frac{1}{7}c^{2}_{56}
-\frac{1}{7}c^{3}_{56}
+\frac{1}{7}c^{11}_{56}
$;\ \ 
$-\frac{1}{7}c^{5}_{56}
-\frac{1}{7}c^{6}_{56}
-\frac{1}{7}c^{9}_{56}
$,
$\frac{1}{7}c^{1}_{112}
+\frac{1}{7}c^{3}_{112}
-\frac{1}{7}c^{5}_{112}
-\frac{1}{7}c^{7}_{112}
-\frac{1}{7}c^{9}_{112}
-\frac{1}{7}c^{11}_{112}
+\frac{2}{7}c^{13}_{112}
+\frac{1}{7}c^{15}_{112}
+\frac{2}{7}c^{17}_{112}
+\frac{1}{7}c^{19}_{112}
-\frac{1}{7}c^{21}_{112}
-\frac{1}{7}c^{23}_{112}
$;\ \ 
$\frac{1}{7}c^{5}_{56}
+\frac{1}{7}c^{6}_{56}
+\frac{1}{7}c^{9}_{56}
$)

Fail:
cnd( Tr$_I(\rho(\mathfrak s))$ ) =
56 does not divide
 ord($\rho(\mathfrak t)$) =
21, I = [ 1/21 ]. Prop. B.4 (2)

 \ \color{black}

\noindent 559: (dims,levels) = $(6;21
)$,
irreps = $6_{7,1}^{3}
\hskip -1.5pt \otimes \hskip -1.5pt
1_{3}^{1,0}$,
pord$(\rho_\text{isum}(\mathfrak{t})) = 7$,

\vskip 0.7ex
\hangindent=5.5em \hangafter=1
{\white .}\hskip 1em $\rho_\text{isum}(\mathfrak{t})$ =
 $( \frac{1}{21},
\frac{4}{21},
\frac{10}{21},
\frac{13}{21},
\frac{16}{21},
\frac{19}{21} )
$,

\vskip 0.7ex
\hangindent=5.5em \hangafter=1
{\white .}\hskip 1em $\rho_\text{isum}(\mathfrak{s})$ =
$\mathrm{i}$($\frac{2}{7}c^{1}_{56}
-\frac{1}{7}c^{3}_{56}
-\frac{1}{7}c^{5}_{56}
+\frac{1}{7}c^{7}_{56}
+\frac{1}{7}c^{9}_{56}
+\frac{1}{7}c^{10}_{56}
-\frac{1}{7}c^{11}_{56}
$,
$\frac{1}{7}c^{5}_{56}
-\frac{1}{7}c^{6}_{56}
+\frac{1}{7}c^{9}_{56}
$,
$-\frac{1}{7}c^{1}_{112}
+\frac{1}{7}c^{3}_{112}
+\frac{1}{7}c^{11}_{112}
+\frac{1}{7}c^{15}_{112}
$,
$\frac{1}{7}c^{1}_{112}
+\frac{1}{7}c^{5}_{112}
-\frac{1}{7}c^{7}_{112}
-\frac{1}{7}c^{9}_{112}
+\frac{1}{7}c^{15}_{112}
+\frac{2}{7}c^{17}_{112}
+\frac{1}{7}c^{19}_{112}
-\frac{1}{7}c^{23}_{112}
$,
$\frac{1}{7}c^{2}_{56}
+\frac{1}{7}c^{3}_{56}
-\frac{1}{7}c^{11}_{56}
$,
$\frac{1}{7}c^{3}_{112}
-\frac{1}{7}c^{5}_{112}
-\frac{1}{7}c^{9}_{112}
-\frac{1}{7}c^{11}_{112}
+\frac{2}{7}c^{13}_{112}
+\frac{1}{7}c^{19}_{112}
-\frac{1}{7}c^{21}_{112}
+\frac{1}{7}c^{23}_{112}
$;\ \ 
$-\frac{1}{7}c^{2}_{56}
-\frac{1}{7}c^{3}_{56}
+\frac{1}{7}c^{11}_{56}
$,
$-\frac{1}{7}c^{3}_{112}
+\frac{1}{7}c^{5}_{112}
+\frac{1}{7}c^{9}_{112}
+\frac{1}{7}c^{11}_{112}
-\frac{2}{7}c^{13}_{112}
-\frac{1}{7}c^{19}_{112}
+\frac{1}{7}c^{21}_{112}
-\frac{1}{7}c^{23}_{112}
$,
$-\frac{1}{7}c^{1}_{112}
+\frac{1}{7}c^{3}_{112}
+\frac{1}{7}c^{11}_{112}
+\frac{1}{7}c^{15}_{112}
$,
$-\frac{2}{7}c^{1}_{56}
+\frac{1}{7}c^{3}_{56}
+\frac{1}{7}c^{5}_{56}
-\frac{1}{7}c^{7}_{56}
-\frac{1}{7}c^{9}_{56}
-\frac{1}{7}c^{10}_{56}
+\frac{1}{7}c^{11}_{56}
$,
$-\frac{1}{7}c^{1}_{112}
-\frac{1}{7}c^{5}_{112}
+\frac{1}{7}c^{7}_{112}
+\frac{1}{7}c^{9}_{112}
-\frac{1}{7}c^{15}_{112}
-\frac{2}{7}c^{17}_{112}
-\frac{1}{7}c^{19}_{112}
+\frac{1}{7}c^{23}_{112}
$;\ \ 
$\frac{1}{7}c^{2}_{56}
+\frac{1}{7}c^{3}_{56}
-\frac{1}{7}c^{11}_{56}
$,
$-\frac{1}{7}c^{5}_{56}
+\frac{1}{7}c^{6}_{56}
-\frac{1}{7}c^{9}_{56}
$,
$-\frac{1}{7}c^{1}_{112}
-\frac{1}{7}c^{5}_{112}
+\frac{1}{7}c^{7}_{112}
+\frac{1}{7}c^{9}_{112}
-\frac{1}{7}c^{15}_{112}
-\frac{2}{7}c^{17}_{112}
-\frac{1}{7}c^{19}_{112}
+\frac{1}{7}c^{23}_{112}
$,
$\frac{2}{7}c^{1}_{56}
-\frac{1}{7}c^{3}_{56}
-\frac{1}{7}c^{5}_{56}
+\frac{1}{7}c^{7}_{56}
+\frac{1}{7}c^{9}_{56}
+\frac{1}{7}c^{10}_{56}
-\frac{1}{7}c^{11}_{56}
$;\ \ 
$-\frac{2}{7}c^{1}_{56}
+\frac{1}{7}c^{3}_{56}
+\frac{1}{7}c^{5}_{56}
-\frac{1}{7}c^{7}_{56}
-\frac{1}{7}c^{9}_{56}
-\frac{1}{7}c^{10}_{56}
+\frac{1}{7}c^{11}_{56}
$,
$\frac{1}{7}c^{3}_{112}
-\frac{1}{7}c^{5}_{112}
-\frac{1}{7}c^{9}_{112}
-\frac{1}{7}c^{11}_{112}
+\frac{2}{7}c^{13}_{112}
+\frac{1}{7}c^{19}_{112}
-\frac{1}{7}c^{21}_{112}
+\frac{1}{7}c^{23}_{112}
$,
$-\frac{1}{7}c^{2}_{56}
-\frac{1}{7}c^{3}_{56}
+\frac{1}{7}c^{11}_{56}
$;\ \ 
$\frac{1}{7}c^{5}_{56}
-\frac{1}{7}c^{6}_{56}
+\frac{1}{7}c^{9}_{56}
$,
$-\frac{1}{7}c^{1}_{112}
+\frac{1}{7}c^{3}_{112}
+\frac{1}{7}c^{11}_{112}
+\frac{1}{7}c^{15}_{112}
$;\ \ 
$-\frac{1}{7}c^{5}_{56}
+\frac{1}{7}c^{6}_{56}
-\frac{1}{7}c^{9}_{56}
$)

Fail:
cnd( Tr$_I(\rho(\mathfrak s))$ ) =
56 does not divide
 ord($\rho(\mathfrak t)$) =
21, I = [ 1/21 ]. Prop. B.4 (2)

 \ \color{black}

\noindent 560: (dims,levels) = $(6;21
)$,
irreps = $3_{7}^{1}
\hskip -1.5pt \otimes \hskip -1.5pt
2_{3}^{1,0}$,
pord$(\rho_\text{isum}(\mathfrak{t})) = 21$,

\vskip 0.7ex
\hangindent=5.5em \hangafter=1
{\white .}\hskip 1em $\rho_\text{isum}(\mathfrak{t})$ =
 $( \frac{1}{7},
\frac{2}{7},
\frac{4}{7},
\frac{10}{21},
\frac{13}{21},
\frac{19}{21} )
$,

\vskip 0.7ex
\hangindent=5.5em \hangafter=1
{\white .}\hskip 1em $\rho_\text{isum}(\mathfrak{s})$ =
$\mathrm{i}$($\frac{1}{\sqrt{21}}c^{1}_{28}
$,
$-\frac{1}{\sqrt{21}}c^{3}_{28}
$,
$\frac{1}{\sqrt{21}}c^{5}_{28}
$,
$-\frac{2}{\sqrt{42}}c^{1}_{28}
$,
$-\frac{2}{\sqrt{42}}c^{3}_{28}
$,
$\frac{2}{\sqrt{42}}c^{5}_{28}
$;\ \ 
$-\frac{1}{\sqrt{21}}c^{5}_{28}
$,
$\frac{1}{\sqrt{21}}c^{1}_{28}
$,
$\frac{2}{\sqrt{42}}c^{3}_{28}
$,
$-\frac{2}{\sqrt{42}}c^{5}_{28}
$,
$\frac{2}{\sqrt{42}}c^{1}_{28}
$;\ \ 
$\frac{1}{\sqrt{21}}c^{3}_{28}
$,
$-\frac{2}{\sqrt{42}}c^{5}_{28}
$,
$\frac{2}{\sqrt{42}}c^{1}_{28}
$,
$\frac{2}{\sqrt{42}}c^{3}_{28}
$;\ \ 
$-\frac{1}{\sqrt{21}}c^{1}_{28}
$,
$-\frac{1}{\sqrt{21}}c^{3}_{28}
$,
$\frac{1}{\sqrt{21}}c^{5}_{28}
$;\ \ 
$\frac{1}{\sqrt{21}}c^{5}_{28}
$,
$-\frac{1}{\sqrt{21}}c^{1}_{28}
$;\ \ 
$-\frac{1}{\sqrt{21}}c^{3}_{28}
$)

Fail:
cnd($\rho(\mathfrak s)_\mathrm{ndeg}$) = 168 does not divide
 ord($\rho(\mathfrak t)$)=21. Prop. B.4 (2)

 \ \color{black}

\noindent 561: (dims,levels) = $(6;21
)$,
irreps = $3_{7}^{3}
\hskip -1.5pt \otimes \hskip -1.5pt
2_{3}^{1,0}$,
pord$(\rho_\text{isum}(\mathfrak{t})) = 21$,

\vskip 0.7ex
\hangindent=5.5em \hangafter=1
{\white .}\hskip 1em $\rho_\text{isum}(\mathfrak{t})$ =
 $( \frac{3}{7},
\frac{5}{7},
\frac{6}{7},
\frac{1}{21},
\frac{4}{21},
\frac{16}{21} )
$,

\vskip 0.7ex
\hangindent=5.5em \hangafter=1
{\white .}\hskip 1em $\rho_\text{isum}(\mathfrak{s})$ =
$\mathrm{i}$($\frac{1}{\sqrt{21}}c^{3}_{28}
$,
$-\frac{1}{\sqrt{21}}c^{1}_{28}
$,
$\frac{1}{\sqrt{21}}c^{5}_{28}
$,
$-\frac{2}{\sqrt{42}}c^{1}_{28}
$,
$\frac{2}{\sqrt{42}}c^{5}_{28}
$,
$-\frac{2}{\sqrt{42}}c^{3}_{28}
$;\ \ 
$-\frac{1}{\sqrt{21}}c^{5}_{28}
$,
$\frac{1}{\sqrt{21}}c^{3}_{28}
$,
$-\frac{2}{\sqrt{42}}c^{5}_{28}
$,
$\frac{2}{\sqrt{42}}c^{3}_{28}
$,
$\frac{2}{\sqrt{42}}c^{1}_{28}
$;\ \ 
$\frac{1}{\sqrt{21}}c^{1}_{28}
$,
$\frac{2}{\sqrt{42}}c^{3}_{28}
$,
$\frac{2}{\sqrt{42}}c^{1}_{28}
$,
$-\frac{2}{\sqrt{42}}c^{5}_{28}
$;\ \ 
$\frac{1}{\sqrt{21}}c^{5}_{28}
$,
$-\frac{1}{\sqrt{21}}c^{3}_{28}
$,
$-\frac{1}{\sqrt{21}}c^{1}_{28}
$;\ \ 
$-\frac{1}{\sqrt{21}}c^{1}_{28}
$,
$\frac{1}{\sqrt{21}}c^{5}_{28}
$;\ \ 
$-\frac{1}{\sqrt{21}}c^{3}_{28}
$)

Fail:
cnd($\rho(\mathfrak s)_\mathrm{ndeg}$) = 168 does not divide
 ord($\rho(\mathfrak t)$)=21. Prop. B.4 (2)

 \ \color{black}

\noindent 562: (dims,levels) = $(6;22
)$,
irreps = $6_{11}^{7}
\hskip -1.5pt \otimes \hskip -1.5pt
1_{2}^{1,0}$,
pord$(\rho_\text{isum}(\mathfrak{t})) = 11$,

\vskip 0.7ex
\hangindent=5.5em \hangafter=1
{\white .}\hskip 1em $\rho_\text{isum}(\mathfrak{t})$ =
 $( \frac{1}{2},
\frac{1}{22},
\frac{3}{22},
\frac{5}{22},
\frac{9}{22},
\frac{15}{22} )
$,

\vskip 0.7ex
\hangindent=5.5em \hangafter=1
{\white .}\hskip 1em $\rho_\text{isum}(\mathfrak{s})$ =
$\mathrm{i}$($-\sqrt{\frac{1}{11}}$,
$\sqrt{\frac{2}{11}}$,
$\sqrt{\frac{2}{11}}$,
$\sqrt{\frac{2}{11}}$,
$\sqrt{\frac{2}{11}}$,
$\sqrt{\frac{2}{11}}$;\ \ 
$-\frac{1}{\sqrt{11}}c^{1}_{11}
$,
$\frac{1}{\sqrt{11}\mathrm{i}}s^{9}_{44}
$,
$-\frac{1}{\sqrt{11}}c^{4}_{11}
$,
$-\frac{1}{\sqrt{11}}c^{3}_{11}
$,
$-\frac{1}{\sqrt{11}}c^{2}_{11}
$;\ \ 
$-\frac{1}{\sqrt{11}}c^{3}_{11}
$,
$-\frac{1}{\sqrt{11}}c^{2}_{11}
$,
$-\frac{1}{\sqrt{11}}c^{4}_{11}
$,
$-\frac{1}{\sqrt{11}}c^{1}_{11}
$;\ \ 
$\frac{1}{\sqrt{11}\mathrm{i}}s^{9}_{44}
$,
$-\frac{1}{\sqrt{11}}c^{1}_{11}
$,
$-\frac{1}{\sqrt{11}}c^{3}_{11}
$;\ \ 
$-\frac{1}{\sqrt{11}}c^{2}_{11}
$,
$\frac{1}{\sqrt{11}\mathrm{i}}s^{9}_{44}
$;\ \ 
$-\frac{1}{\sqrt{11}}c^{4}_{11}
$)

Fail:
cnd($\rho(\mathfrak s)_\mathrm{ndeg}$) = 88 does not divide
 ord($\rho(\mathfrak t)$)=22. Prop. B.4 (2)

 \ \color{black}

\noindent 563: (dims,levels) = $(6;22
)$,
irreps = $6_{11}^{1}
\hskip -1.5pt \otimes \hskip -1.5pt
1_{2}^{1,0}$,
pord$(\rho_\text{isum}(\mathfrak{t})) = 11$,

\vskip 0.7ex
\hangindent=5.5em \hangafter=1
{\white .}\hskip 1em $\rho_\text{isum}(\mathfrak{t})$ =
 $( \frac{1}{2},
\frac{7}{22},
\frac{13}{22},
\frac{17}{22},
\frac{19}{22},
\frac{21}{22} )
$,

\vskip 0.7ex
\hangindent=5.5em \hangafter=1
{\white .}\hskip 1em $\rho_\text{isum}(\mathfrak{s})$ =
$\mathrm{i}$($\sqrt{\frac{1}{11}}$,
$\sqrt{\frac{2}{11}}$,
$\sqrt{\frac{2}{11}}$,
$\sqrt{\frac{2}{11}}$,
$\sqrt{\frac{2}{11}}$,
$\sqrt{\frac{2}{11}}$;\ \ 
$\frac{1}{\sqrt{11}}c^{4}_{11}
$,
$-\frac{1}{\sqrt{11}\mathrm{i}}s^{9}_{44}
$,
$\frac{1}{\sqrt{11}}c^{3}_{11}
$,
$\frac{1}{\sqrt{11}}c^{1}_{11}
$,
$\frac{1}{\sqrt{11}}c^{2}_{11}
$;\ \ 
$\frac{1}{\sqrt{11}}c^{2}_{11}
$,
$\frac{1}{\sqrt{11}}c^{1}_{11}
$,
$\frac{1}{\sqrt{11}}c^{4}_{11}
$,
$\frac{1}{\sqrt{11}}c^{3}_{11}
$;\ \ 
$-\frac{1}{\sqrt{11}\mathrm{i}}s^{9}_{44}
$,
$\frac{1}{\sqrt{11}}c^{2}_{11}
$,
$\frac{1}{\sqrt{11}}c^{4}_{11}
$;\ \ 
$\frac{1}{\sqrt{11}}c^{3}_{11}
$,
$-\frac{1}{\sqrt{11}\mathrm{i}}s^{9}_{44}
$;\ \ 
$\frac{1}{\sqrt{11}}c^{1}_{11}
$)

Fail:
cnd($\rho(\mathfrak s)_\mathrm{ndeg}$) = 88 does not divide
 ord($\rho(\mathfrak t)$)=22. Prop. B.4 (2)

 \ \color{black}

\noindent 564: (dims,levels) = $(6;24
)$,
irreps = $3_{8}^{1,0}
\hskip -1.5pt \otimes \hskip -1.5pt
2_{3}^{1,0}$,
pord$(\rho_\text{isum}(\mathfrak{t})) = 24$,

\vskip 0.7ex
\hangindent=5.5em \hangafter=1
{\white .}\hskip 1em $\rho_\text{isum}(\mathfrak{t})$ =
 $( 0,
\frac{1}{3},
\frac{1}{8},
\frac{5}{8},
\frac{11}{24},
\frac{23}{24} )
$,

\vskip 0.7ex
\hangindent=5.5em \hangafter=1
{\white .}\hskip 1em $\rho_\text{isum}(\mathfrak{s})$ =
($0$,
$0$,
$\sqrt{\frac{1}{6}}$,
$\sqrt{\frac{1}{6}}$,
$-\sqrt{\frac{1}{3}}$,
$-\sqrt{\frac{1}{3}}$;
$0$,
$-\sqrt{\frac{1}{3}}$,
$-\sqrt{\frac{1}{3}}$,
$-\sqrt{\frac{1}{6}}$,
$-\sqrt{\frac{1}{6}}$;
$-\sqrt{\frac{1}{12}}$,
$\sqrt{\frac{1}{12}}$,
$\sqrt{\frac{1}{6}}$,
$-\sqrt{\frac{1}{6}}$;
$-\sqrt{\frac{1}{12}}$,
$-\sqrt{\frac{1}{6}}$,
$\sqrt{\frac{1}{6}}$;
$\sqrt{\frac{1}{12}}$,
$-\sqrt{\frac{1}{12}}$;
$\sqrt{\frac{1}{12}}$)

Fail:
$\sigma(\rho(\mathfrak s)_\mathrm{ndeg}) \neq
 (\rho(\mathfrak t)^a \rho(\mathfrak s) \rho(\mathfrak t)^b
 \rho(\mathfrak s) \rho(\mathfrak t)^a)_\mathrm{ndeg}$,
 $\sigma = a$ = 5. Prop. B.5 (3) eqn. (B.25)

 \ \color{black}

\noindent 565: (dims,levels) = $(6;24
)$,
irreps = $3_{8}^{3,0}
\hskip -1.5pt \otimes \hskip -1.5pt
2_{3}^{1,0}$,
pord$(\rho_\text{isum}(\mathfrak{t})) = 24$,

\vskip 0.7ex
\hangindent=5.5em \hangafter=1
{\white .}\hskip 1em $\rho_\text{isum}(\mathfrak{t})$ =
 $( 0,
\frac{1}{3},
\frac{3}{8},
\frac{7}{8},
\frac{5}{24},
\frac{17}{24} )
$,

\vskip 0.7ex
\hangindent=5.5em \hangafter=1
{\white .}\hskip 1em $\rho_\text{isum}(\mathfrak{s})$ =
($0$,
$0$,
$\sqrt{\frac{1}{6}}$,
$\sqrt{\frac{1}{6}}$,
$-\sqrt{\frac{1}{3}}$,
$-\sqrt{\frac{1}{3}}$;
$0$,
$-\sqrt{\frac{1}{3}}$,
$-\sqrt{\frac{1}{3}}$,
$-\sqrt{\frac{1}{6}}$,
$-\sqrt{\frac{1}{6}}$;
$\sqrt{\frac{1}{12}}$,
$-\sqrt{\frac{1}{12}}$,
$\sqrt{\frac{1}{6}}$,
$-\sqrt{\frac{1}{6}}$;
$\sqrt{\frac{1}{12}}$,
$-\sqrt{\frac{1}{6}}$,
$\sqrt{\frac{1}{6}}$;
$-\sqrt{\frac{1}{12}}$,
$\sqrt{\frac{1}{12}}$;
$-\sqrt{\frac{1}{12}}$)

Fail:
$\sigma(\rho(\mathfrak s)_\mathrm{ndeg}) \neq
 (\rho(\mathfrak t)^a \rho(\mathfrak s) \rho(\mathfrak t)^b
 \rho(\mathfrak s) \rho(\mathfrak t)^a)_\mathrm{ndeg}$,
 $\sigma = a$ = 5. Prop. B.5 (3) eqn. (B.25)

 \ \color{black}

\noindent 566: (dims,levels) = $(6;24
)$,
irreps = $3_{3}^{1,0}
\hskip -1.5pt \otimes \hskip -1.5pt
2_{8}^{1,0}$,
pord$(\rho_\text{isum}(\mathfrak{t})) = 12$,

\vskip 0.7ex
\hangindent=5.5em \hangafter=1
{\white .}\hskip 1em $\rho_\text{isum}(\mathfrak{t})$ =
 $( \frac{1}{8},
\frac{3}{8},
\frac{1}{24},
\frac{11}{24},
\frac{17}{24},
\frac{19}{24} )
$,

\vskip 0.7ex
\hangindent=5.5em \hangafter=1
{\white .}\hskip 1em $\rho_\text{isum}(\mathfrak{s})$ =
($\sqrt{\frac{1}{18}}$,
$\sqrt{\frac{1}{18}}$,
$\sqrt{\frac{2}{9}}$,
$\sqrt{\frac{2}{9}}$,
$\sqrt{\frac{2}{9}}$,
$\sqrt{\frac{2}{9}}$;
$-\sqrt{\frac{1}{18}}$,
$-\sqrt{\frac{2}{9}}$,
$\sqrt{\frac{2}{9}}$,
$-\sqrt{\frac{2}{9}}$,
$\sqrt{\frac{2}{9}}$;
$-\sqrt{\frac{1}{18}}$,
$-\sqrt{\frac{2}{9}}$,
$\sqrt{\frac{2}{9}}$,
$\sqrt{\frac{1}{18}}$;
$\sqrt{\frac{1}{18}}$,
$\sqrt{\frac{1}{18}}$,
$-\sqrt{\frac{2}{9}}$;
$-\sqrt{\frac{1}{18}}$,
$-\sqrt{\frac{2}{9}}$;
$\sqrt{\frac{1}{18}}$)

Fail:
Integral: $D_{\rho}(\sigma)_{\theta} \propto $ id,
 for all $\sigma$ and all $\theta$-eigenspaces that can contain unit. Prop. B.5 (6)

 \ \color{black}

\noindent 567: (dims,levels) = $(6;24
)$,
irreps = $6_{8,2}^{1,0}
\hskip -1.5pt \otimes \hskip -1.5pt
1_{3}^{1,0}$,
pord$(\rho_\text{isum}(\mathfrak{t})) = 8$,

\vskip 0.7ex
\hangindent=5.5em \hangafter=1
{\white .}\hskip 1em $\rho_\text{isum}(\mathfrak{t})$ =
 $( \frac{1}{3},
\frac{5}{6},
\frac{1}{12},
\frac{7}{12},
\frac{11}{24},
\frac{17}{24} )
$,

\vskip 0.7ex
\hangindent=5.5em \hangafter=1
{\white .}\hskip 1em $\rho_\text{isum}(\mathfrak{s})$ =
$\mathrm{i}$($-\sqrt{\frac{1}{8}}$,
$\sqrt{\frac{1}{8}}$,
$\sqrt{\frac{1}{8}}$,
$\sqrt{\frac{1}{8}}$,
$\frac{1}{2}$,
$\frac{1}{2}$;\ \ 
$-\sqrt{\frac{1}{8}}$,
$-\sqrt{\frac{1}{8}}$,
$-\sqrt{\frac{1}{8}}$,
$\frac{1}{2}$,
$\frac{1}{2}$;\ \ 
$\sqrt{\frac{1}{8}}$,
$\sqrt{\frac{1}{8}}$,
$\frac{1}{2}$,
$-\frac{1}{2}$;\ \ 
$\sqrt{\frac{1}{8}}$,
$-\frac{1}{2}$,
$\frac{1}{2}$;\ \ 
$0$,
$0$;\ \ 
$0$)

Fail:
$\sigma(\rho(\mathfrak s)_\mathrm{ndeg}) \neq
 (\rho(\mathfrak t)^a \rho(\mathfrak s) \rho(\mathfrak t)^b
 \rho(\mathfrak s) \rho(\mathfrak t)^a)_\mathrm{ndeg}$,
 $\sigma = a$ = 5. Prop. B.5 (3) eqn. (B.25)

 \ \color{black}

\noindent 568: (dims,levels) = $(6;24
)$,
irreps = $6_{8,1}^{1,0}
\hskip -1.5pt \otimes \hskip -1.5pt
1_{3}^{1,0}$,
pord$(\rho_\text{isum}(\mathfrak{t})) = 8$,

\vskip 0.7ex
\hangindent=5.5em \hangafter=1
{\white .}\hskip 1em $\rho_\text{isum}(\mathfrak{t})$ =
 $( \frac{1}{3},
\frac{5}{6},
\frac{5}{24},
\frac{11}{24},
\frac{17}{24},
\frac{23}{24} )
$,

\vskip 0.7ex
\hangindent=5.5em \hangafter=1
{\white .}\hskip 1em $\rho_\text{isum}(\mathfrak{s})$ =
($0$,
$0$,
$\frac{1}{2}$,
$\frac{1}{2}$,
$\frac{1}{2}$,
$\frac{1}{2}$;
$0$,
$\frac{1}{2}$,
$-\frac{1}{2}$,
$\frac{1}{2}$,
$-\frac{1}{2}$;
$0$,
$\frac{1}{2}$,
$0$,
$-\frac{1}{2}$;
$0$,
$-\frac{1}{2}$,
$0$;
$0$,
$\frac{1}{2}$;
$0$)

Fail:
all rows of $U \rho(\mathfrak s) U^\dagger$
 contain zero for any block-diagonal $U$. Prop. B.5 (4) eqn. (B.27)

 \ \color{black}

 \color{blue}

\noindent 569: (dims,levels) = $(6;26
)$,
irreps = $6_{13}^{2}
\hskip -1.5pt \otimes \hskip -1.5pt
1_{2}^{1,0}$,
pord$(\rho_\text{isum}(\mathfrak{t})) = 13$,

\vskip 0.7ex
\hangindent=5.5em \hangafter=1
{\white .}\hskip 1em $\rho_\text{isum}(\mathfrak{t})$ =
 $( \frac{1}{26},
\frac{3}{26},
\frac{9}{26},
\frac{17}{26},
\frac{23}{26},
\frac{25}{26} )
$,

\vskip 0.7ex
\hangindent=5.5em \hangafter=1
{\white .}\hskip 1em $\rho_\text{isum}(\mathfrak{s})$ =
$\mathrm{i}$($-\frac{1}{\sqrt{13}}c^{9}_{52}
$,
$\frac{1}{\sqrt{13}}c^{3}_{52}
$,
$-\frac{1}{\sqrt{13}}c^{1}_{52}
$,
$-\frac{1}{\sqrt{13}}c^{5}_{52}
$,
$\frac{1}{\sqrt{13}}c^{11}_{52}
$,
$\frac{1}{\sqrt{13}}c^{7}_{52}
$;\ \ 
$-\frac{1}{\sqrt{13}}c^{1}_{52}
$,
$-\frac{1}{\sqrt{13}}c^{9}_{52}
$,
$-\frac{1}{\sqrt{13}}c^{7}_{52}
$,
$\frac{1}{\sqrt{13}}c^{5}_{52}
$,
$\frac{1}{\sqrt{13}}c^{11}_{52}
$;\ \ 
$\frac{1}{\sqrt{13}}c^{3}_{52}
$,
$-\frac{1}{\sqrt{13}}c^{11}_{52}
$,
$\frac{1}{\sqrt{13}}c^{7}_{52}
$,
$\frac{1}{\sqrt{13}}c^{5}_{52}
$;\ \ 
$-\frac{1}{\sqrt{13}}c^{3}_{52}
$,
$-\frac{1}{\sqrt{13}}c^{9}_{52}
$,
$-\frac{1}{\sqrt{13}}c^{1}_{52}
$;\ \ 
$\frac{1}{\sqrt{13}}c^{1}_{52}
$,
$-\frac{1}{\sqrt{13}}c^{3}_{52}
$;\ \ 
$\frac{1}{\sqrt{13}}c^{9}_{52}
$)

Pass. 

 \ \color{black}

 \color{blue}

\noindent 570: (dims,levels) = $(6;26
)$,
irreps = $6_{13}^{1}
\hskip -1.5pt \otimes \hskip -1.5pt
1_{2}^{1,0}$,
pord$(\rho_\text{isum}(\mathfrak{t})) = 13$,

\vskip 0.7ex
\hangindent=5.5em \hangafter=1
{\white .}\hskip 1em $\rho_\text{isum}(\mathfrak{t})$ =
 $( \frac{5}{26},
\frac{7}{26},
\frac{11}{26},
\frac{15}{26},
\frac{19}{26},
\frac{21}{26} )
$,

\vskip 0.7ex
\hangindent=5.5em \hangafter=1
{\white .}\hskip 1em $\rho_\text{isum}(\mathfrak{s})$ =
$\mathrm{i}$($\frac{1}{\sqrt{13}}c^{7}_{52}
$,
$-\frac{1}{\sqrt{13}}c^{1}_{52}
$,
$\frac{1}{\sqrt{13}}c^{3}_{52}
$,
$\frac{1}{\sqrt{13}}c^{11}_{52}
$,
$-\frac{1}{\sqrt{13}}c^{5}_{52}
$,
$\frac{1}{\sqrt{13}}c^{9}_{52}
$;\ \ 
$-\frac{1}{\sqrt{13}}c^{11}_{52}
$,
$-\frac{1}{\sqrt{13}}c^{7}_{52}
$,
$-\frac{1}{\sqrt{13}}c^{9}_{52}
$,
$-\frac{1}{\sqrt{13}}c^{3}_{52}
$,
$\frac{1}{\sqrt{13}}c^{5}_{52}
$;\ \ 
$-\frac{1}{\sqrt{13}}c^{5}_{52}
$,
$-\frac{1}{\sqrt{13}}c^{1}_{52}
$,
$\frac{1}{\sqrt{13}}c^{9}_{52}
$,
$\frac{1}{\sqrt{13}}c^{11}_{52}
$;\ \ 
$\frac{1}{\sqrt{13}}c^{5}_{52}
$,
$-\frac{1}{\sqrt{13}}c^{7}_{52}
$,
$-\frac{1}{\sqrt{13}}c^{3}_{52}
$;\ \ 
$\frac{1}{\sqrt{13}}c^{11}_{52}
$,
$-\frac{1}{\sqrt{13}}c^{1}_{52}
$;\ \ 
$-\frac{1}{\sqrt{13}}c^{7}_{52}
$)

Pass. 

 \ \color{black}

\noindent 571: (dims,levels) = $(6;28
)$,
irreps = $3_{7}^{1}
\hskip -1.5pt \otimes \hskip -1.5pt
2_{4}^{1,0}$,
pord$(\rho_\text{isum}(\mathfrak{t})) = 14$,

\vskip 0.7ex
\hangindent=5.5em \hangafter=1
{\white .}\hskip 1em $\rho_\text{isum}(\mathfrak{t})$ =
 $( \frac{1}{28},
\frac{9}{28},
\frac{11}{28},
\frac{15}{28},
\frac{23}{28},
\frac{25}{28} )
$,

\vskip 0.7ex
\hangindent=5.5em \hangafter=1
{\white .}\hskip 1em $\rho_\text{isum}(\mathfrak{s})$ =
$\mathrm{i}$($\frac{1}{2\sqrt{7}}c^{5}_{28}
$,
$-\frac{1}{2\sqrt{7}}c^{1}_{28}
$,
$-\frac{3}{2\sqrt{21}}c^{3}_{28}
$,
$\frac{3}{2\sqrt{21}}c^{5}_{28}
$,
$-\frac{3}{2\sqrt{21}}c^{1}_{28}
$,
$-\frac{1}{2\sqrt{7}}c^{3}_{28}
$;\ \ 
$-\frac{1}{2\sqrt{7}}c^{3}_{28}
$,
$\frac{3}{2\sqrt{21}}c^{5}_{28}
$,
$-\frac{3}{2\sqrt{21}}c^{1}_{28}
$,
$-\frac{3}{2\sqrt{21}}c^{3}_{28}
$,
$\frac{1}{2\sqrt{7}}c^{5}_{28}
$;\ \ 
$\frac{1}{2\sqrt{7}}c^{1}_{28}
$,
$\frac{1}{2\sqrt{7}}c^{3}_{28}
$,
$-\frac{1}{2\sqrt{7}}c^{5}_{28}
$,
$-\frac{3}{2\sqrt{21}}c^{1}_{28}
$;\ \ 
$-\frac{1}{2\sqrt{7}}c^{5}_{28}
$,
$\frac{1}{2\sqrt{7}}c^{1}_{28}
$,
$-\frac{3}{2\sqrt{21}}c^{3}_{28}
$;\ \ 
$\frac{1}{2\sqrt{7}}c^{3}_{28}
$,
$\frac{3}{2\sqrt{21}}c^{5}_{28}
$;\ \ 
$-\frac{1}{2\sqrt{7}}c^{1}_{28}
$)

Fail:
cnd($\rho(\mathfrak s)_\mathrm{ndeg}$) = 84 does not divide
 ord($\rho(\mathfrak t)$)=28. Prop. B.4 (2)

 \ \color{black}

\noindent 572: (dims,levels) = $(6;28
)$,
irreps = $3_{7}^{3}
\hskip -1.5pt \otimes \hskip -1.5pt
2_{4}^{1,0}$,
pord$(\rho_\text{isum}(\mathfrak{t})) = 14$,

\vskip 0.7ex
\hangindent=5.5em \hangafter=1
{\white .}\hskip 1em $\rho_\text{isum}(\mathfrak{t})$ =
 $( \frac{3}{28},
\frac{5}{28},
\frac{13}{28},
\frac{17}{28},
\frac{19}{28},
\frac{27}{28} )
$,

\vskip 0.7ex
\hangindent=5.5em \hangafter=1
{\white .}\hskip 1em $\rho_\text{isum}(\mathfrak{s})$ =
$\mathrm{i}$($\frac{1}{2\sqrt{7}}c^{1}_{28}
$,
$\frac{3}{2\sqrt{21}}c^{5}_{28}
$,
$-\frac{3}{2\sqrt{21}}c^{3}_{28}
$,
$-\frac{3}{2\sqrt{21}}c^{1}_{28}
$,
$\frac{1}{2\sqrt{7}}c^{5}_{28}
$,
$-\frac{1}{2\sqrt{7}}c^{3}_{28}
$;\ \ 
$-\frac{1}{2\sqrt{7}}c^{3}_{28}
$,
$-\frac{1}{2\sqrt{7}}c^{1}_{28}
$,
$\frac{1}{2\sqrt{7}}c^{5}_{28}
$,
$\frac{3}{2\sqrt{21}}c^{3}_{28}
$,
$\frac{3}{2\sqrt{21}}c^{1}_{28}
$;\ \ 
$\frac{1}{2\sqrt{7}}c^{5}_{28}
$,
$-\frac{1}{2\sqrt{7}}c^{3}_{28}
$,
$\frac{3}{2\sqrt{21}}c^{1}_{28}
$,
$-\frac{3}{2\sqrt{21}}c^{5}_{28}
$;\ \ 
$-\frac{1}{2\sqrt{7}}c^{1}_{28}
$,
$-\frac{3}{2\sqrt{21}}c^{5}_{28}
$,
$\frac{3}{2\sqrt{21}}c^{3}_{28}
$;\ \ 
$\frac{1}{2\sqrt{7}}c^{3}_{28}
$,
$\frac{1}{2\sqrt{7}}c^{1}_{28}
$;\ \ 
$-\frac{1}{2\sqrt{7}}c^{5}_{28}
$)

Fail:
cnd($\rho(\mathfrak s)_\mathrm{ndeg}$) = 84 does not divide
 ord($\rho(\mathfrak t)$)=28. Prop. B.4 (2)

 \ \color{black}

\noindent 573: (dims,levels) = $(6;28
)$,
irreps = $6_{7,2}^{1}
\hskip -1.5pt \otimes \hskip -1.5pt
1_{4}^{1,0}$,
pord$(\rho_\text{isum}(\mathfrak{t})) = 7$,

\vskip 0.7ex
\hangindent=5.5em \hangafter=1
{\white .}\hskip 1em $\rho_\text{isum}(\mathfrak{t})$ =
 $( \frac{3}{28},
\frac{11}{28},
\frac{15}{28},
\frac{19}{28},
\frac{23}{28},
\frac{27}{28} )
$,

\vskip 0.7ex
\hangindent=5.5em \hangafter=1
{\white .}\hskip 1em $\rho_\text{isum}(\mathfrak{s})$ =
$\mathrm{i}$($\frac{2}{7}-\frac{1}{7}c^{2}_{7}
$,
$\frac{2}{7}c^{3}_{56}
+\frac{1}{7}c^{5}_{56}
-\frac{1}{7}c^{7}_{56}
-\frac{1}{7}c^{9}_{56}
+\frac{2}{7}c^{11}_{56}
$,
$\frac{1}{7}c^{3}_{56}
-\frac{1}{7}c^{5}_{56}
+\frac{1}{7}c^{9}_{56}
+\frac{1}{7}c^{11}_{56}
$,
$-\frac{3}{7}-\frac{1}{7}c^{1}_{7}
-\frac{1}{7}c^{2}_{7}
$,
$-\frac{1}{7}c^{3}_{56}
-\frac{2}{7}c^{5}_{56}
+\frac{1}{7}c^{7}_{56}
+\frac{2}{7}c^{9}_{56}
-\frac{1}{7}c^{11}_{56}
$,
$-\frac{2}{7}+\frac{1}{7}c^{1}_{7}
$;\ \ 
$\frac{2}{7}-\frac{1}{7}c^{2}_{7}
$,
$-\frac{2}{7}+\frac{1}{7}c^{1}_{7}
$,
$\frac{1}{7}c^{3}_{56}
+\frac{2}{7}c^{5}_{56}
-\frac{1}{7}c^{7}_{56}
-\frac{2}{7}c^{9}_{56}
+\frac{1}{7}c^{11}_{56}
$,
$\frac{3}{7}+\frac{1}{7}c^{1}_{7}
+\frac{1}{7}c^{2}_{7}
$,
$\frac{1}{7}c^{3}_{56}
-\frac{1}{7}c^{5}_{56}
+\frac{1}{7}c^{9}_{56}
+\frac{1}{7}c^{11}_{56}
$;\ \ 
$\frac{3}{7}+\frac{1}{7}c^{1}_{7}
+\frac{1}{7}c^{2}_{7}
$,
$\frac{2}{7}c^{3}_{56}
+\frac{1}{7}c^{5}_{56}
-\frac{1}{7}c^{7}_{56}
-\frac{1}{7}c^{9}_{56}
+\frac{2}{7}c^{11}_{56}
$,
$-\frac{2}{7}+\frac{1}{7}c^{2}_{7}
$,
$-\frac{1}{7}c^{3}_{56}
-\frac{2}{7}c^{5}_{56}
+\frac{1}{7}c^{7}_{56}
+\frac{2}{7}c^{9}_{56}
-\frac{1}{7}c^{11}_{56}
$;\ \ 
$\frac{2}{7}-\frac{1}{7}c^{1}_{7}
$,
$\frac{1}{7}c^{3}_{56}
-\frac{1}{7}c^{5}_{56}
+\frac{1}{7}c^{9}_{56}
+\frac{1}{7}c^{11}_{56}
$,
$\frac{2}{7}-\frac{1}{7}c^{2}_{7}
$;\ \ 
$\frac{2}{7}-\frac{1}{7}c^{1}_{7}
$,
$-\frac{2}{7}c^{3}_{56}
-\frac{1}{7}c^{5}_{56}
+\frac{1}{7}c^{7}_{56}
+\frac{1}{7}c^{9}_{56}
-\frac{2}{7}c^{11}_{56}
$;\ \ 
$\frac{3}{7}+\frac{1}{7}c^{1}_{7}
+\frac{1}{7}c^{2}_{7}
$)

Fail:
cnd($\rho(\mathfrak s)_\mathrm{ndeg}$) = 56 does not divide
 ord($\rho(\mathfrak t)$)=28. Prop. B.4 (2)

 \ \color{black}

\noindent 574: (dims,levels) = $(6;28
)$,
irreps = $6_{7,1}^{3}
\hskip -1.5pt \otimes \hskip -1.5pt
1_{4}^{1,0}$,
pord$(\rho_\text{isum}(\mathfrak{t})) = 7$,

\vskip 0.7ex
\hangindent=5.5em \hangafter=1
{\white .}\hskip 1em $\rho_\text{isum}(\mathfrak{t})$ =
 $( \frac{3}{28},
\frac{11}{28},
\frac{15}{28},
\frac{19}{28},
\frac{23}{28},
\frac{27}{28} )
$,

\vskip 0.7ex
\hangindent=5.5em \hangafter=1
{\white .}\hskip 1em $\rho_\text{isum}(\mathfrak{s})$ =
($\frac{1}{7}c^{2}_{56}
+\frac{1}{7}c^{3}_{56}
-\frac{1}{7}c^{11}_{56}
$,
$\frac{1}{7}c^{3}_{112}
-\frac{1}{7}c^{5}_{112}
-\frac{1}{7}c^{9}_{112}
-\frac{1}{7}c^{11}_{112}
+\frac{2}{7}c^{13}_{112}
+\frac{1}{7}c^{19}_{112}
-\frac{1}{7}c^{21}_{112}
+\frac{1}{7}c^{23}_{112}
$,
$\frac{1}{7}c^{1}_{112}
-\frac{1}{7}c^{3}_{112}
-\frac{1}{7}c^{11}_{112}
-\frac{1}{7}c^{15}_{112}
$,
$\frac{2}{7}c^{1}_{56}
-\frac{1}{7}c^{3}_{56}
-\frac{1}{7}c^{5}_{56}
+\frac{1}{7}c^{7}_{56}
+\frac{1}{7}c^{9}_{56}
+\frac{1}{7}c^{10}_{56}
-\frac{1}{7}c^{11}_{56}
$,
$\frac{1}{7}c^{1}_{112}
+\frac{1}{7}c^{5}_{112}
-\frac{1}{7}c^{7}_{112}
-\frac{1}{7}c^{9}_{112}
+\frac{1}{7}c^{15}_{112}
+\frac{2}{7}c^{17}_{112}
+\frac{1}{7}c^{19}_{112}
-\frac{1}{7}c^{23}_{112}
$,
$\frac{1}{7}c^{5}_{56}
-\frac{1}{7}c^{6}_{56}
+\frac{1}{7}c^{9}_{56}
$;
$-\frac{1}{7}c^{2}_{56}
-\frac{1}{7}c^{3}_{56}
+\frac{1}{7}c^{11}_{56}
$,
$\frac{1}{7}c^{5}_{56}
-\frac{1}{7}c^{6}_{56}
+\frac{1}{7}c^{9}_{56}
$,
$\frac{1}{7}c^{1}_{112}
+\frac{1}{7}c^{5}_{112}
-\frac{1}{7}c^{7}_{112}
-\frac{1}{7}c^{9}_{112}
+\frac{1}{7}c^{15}_{112}
+\frac{2}{7}c^{17}_{112}
+\frac{1}{7}c^{19}_{112}
-\frac{1}{7}c^{23}_{112}
$,
$-\frac{2}{7}c^{1}_{56}
+\frac{1}{7}c^{3}_{56}
+\frac{1}{7}c^{5}_{56}
-\frac{1}{7}c^{7}_{56}
-\frac{1}{7}c^{9}_{56}
-\frac{1}{7}c^{10}_{56}
+\frac{1}{7}c^{11}_{56}
$,
$-\frac{1}{7}c^{1}_{112}
+\frac{1}{7}c^{3}_{112}
+\frac{1}{7}c^{11}_{112}
+\frac{1}{7}c^{15}_{112}
$;
$\frac{2}{7}c^{1}_{56}
-\frac{1}{7}c^{3}_{56}
-\frac{1}{7}c^{5}_{56}
+\frac{1}{7}c^{7}_{56}
+\frac{1}{7}c^{9}_{56}
+\frac{1}{7}c^{10}_{56}
-\frac{1}{7}c^{11}_{56}
$,
$-\frac{1}{7}c^{3}_{112}
+\frac{1}{7}c^{5}_{112}
+\frac{1}{7}c^{9}_{112}
+\frac{1}{7}c^{11}_{112}
-\frac{2}{7}c^{13}_{112}
-\frac{1}{7}c^{19}_{112}
+\frac{1}{7}c^{21}_{112}
-\frac{1}{7}c^{23}_{112}
$,
$\frac{1}{7}c^{2}_{56}
+\frac{1}{7}c^{3}_{56}
-\frac{1}{7}c^{11}_{56}
$,
$\frac{1}{7}c^{1}_{112}
+\frac{1}{7}c^{5}_{112}
-\frac{1}{7}c^{7}_{112}
-\frac{1}{7}c^{9}_{112}
+\frac{1}{7}c^{15}_{112}
+\frac{2}{7}c^{17}_{112}
+\frac{1}{7}c^{19}_{112}
-\frac{1}{7}c^{23}_{112}
$;
$-\frac{1}{7}c^{5}_{56}
+\frac{1}{7}c^{6}_{56}
-\frac{1}{7}c^{9}_{56}
$,
$\frac{1}{7}c^{1}_{112}
-\frac{1}{7}c^{3}_{112}
-\frac{1}{7}c^{11}_{112}
-\frac{1}{7}c^{15}_{112}
$,
$\frac{1}{7}c^{2}_{56}
+\frac{1}{7}c^{3}_{56}
-\frac{1}{7}c^{11}_{56}
$;
$\frac{1}{7}c^{5}_{56}
-\frac{1}{7}c^{6}_{56}
+\frac{1}{7}c^{9}_{56}
$,
$\frac{1}{7}c^{3}_{112}
-\frac{1}{7}c^{5}_{112}
-\frac{1}{7}c^{9}_{112}
-\frac{1}{7}c^{11}_{112}
+\frac{2}{7}c^{13}_{112}
+\frac{1}{7}c^{19}_{112}
-\frac{1}{7}c^{21}_{112}
+\frac{1}{7}c^{23}_{112}
$;
$-\frac{2}{7}c^{1}_{56}
+\frac{1}{7}c^{3}_{56}
+\frac{1}{7}c^{5}_{56}
-\frac{1}{7}c^{7}_{56}
-\frac{1}{7}c^{9}_{56}
-\frac{1}{7}c^{10}_{56}
+\frac{1}{7}c^{11}_{56}
$)

Fail:
cnd( Tr$_I(\rho(\mathfrak s))$ ) =
56 does not divide
 ord($\rho(\mathfrak t)$) =
28, I = [ 3/28 ]. Prop. B.4 (2)

 \ \color{black}

\noindent 575: (dims,levels) = $(6;28
)$,
irreps = $6_{7,1}^{1}
\hskip -1.5pt \otimes \hskip -1.5pt
1_{4}^{1,0}$,
pord$(\rho_\text{isum}(\mathfrak{t})) = 7$,

\vskip 0.7ex
\hangindent=5.5em \hangafter=1
{\white .}\hskip 1em $\rho_\text{isum}(\mathfrak{t})$ =
 $( \frac{3}{28},
\frac{11}{28},
\frac{15}{28},
\frac{19}{28},
\frac{23}{28},
\frac{27}{28} )
$,

\vskip 0.7ex
\hangindent=5.5em \hangafter=1
{\white .}\hskip 1em $\rho_\text{isum}(\mathfrak{s})$ =
($\frac{1}{7}c^{2}_{56}
-\frac{1}{7}c^{3}_{56}
+\frac{1}{7}c^{11}_{56}
$,
$\frac{1}{7}c^{1}_{112}
+\frac{1}{7}c^{5}_{112}
-\frac{1}{7}c^{15}_{112}
+\frac{1}{7}c^{19}_{112}
$,
$\frac{1}{7}c^{1}_{112}
+\frac{1}{7}c^{3}_{112}
-\frac{1}{7}c^{5}_{112}
-\frac{1}{7}c^{7}_{112}
-\frac{1}{7}c^{9}_{112}
-\frac{1}{7}c^{11}_{112}
+\frac{2}{7}c^{13}_{112}
+\frac{1}{7}c^{15}_{112}
+\frac{2}{7}c^{17}_{112}
+\frac{1}{7}c^{19}_{112}
-\frac{1}{7}c^{21}_{112}
-\frac{1}{7}c^{23}_{112}
$,
$\frac{2}{7}c^{1}_{56}
-\frac{1}{7}c^{3}_{56}
-\frac{1}{7}c^{5}_{56}
+\frac{1}{7}c^{7}_{56}
+\frac{1}{7}c^{9}_{56}
-\frac{1}{7}c^{10}_{56}
-\frac{1}{7}c^{11}_{56}
$,
$\frac{1}{7}c^{3}_{112}
-\frac{1}{7}c^{9}_{112}
+\frac{1}{7}c^{11}_{112}
+\frac{1}{7}c^{23}_{112}
$,
$\frac{1}{7}c^{5}_{56}
+\frac{1}{7}c^{6}_{56}
+\frac{1}{7}c^{9}_{56}
$;
$-\frac{1}{7}c^{2}_{56}
+\frac{1}{7}c^{3}_{56}
-\frac{1}{7}c^{11}_{56}
$,
$\frac{1}{7}c^{5}_{56}
+\frac{1}{7}c^{6}_{56}
+\frac{1}{7}c^{9}_{56}
$,
$\frac{1}{7}c^{3}_{112}
-\frac{1}{7}c^{9}_{112}
+\frac{1}{7}c^{11}_{112}
+\frac{1}{7}c^{23}_{112}
$,
$-\frac{2}{7}c^{1}_{56}
+\frac{1}{7}c^{3}_{56}
+\frac{1}{7}c^{5}_{56}
-\frac{1}{7}c^{7}_{56}
-\frac{1}{7}c^{9}_{56}
+\frac{1}{7}c^{10}_{56}
+\frac{1}{7}c^{11}_{56}
$,
$-\frac{1}{7}c^{1}_{112}
-\frac{1}{7}c^{3}_{112}
+\frac{1}{7}c^{5}_{112}
+\frac{1}{7}c^{7}_{112}
+\frac{1}{7}c^{9}_{112}
+\frac{1}{7}c^{11}_{112}
-\frac{2}{7}c^{13}_{112}
-\frac{1}{7}c^{15}_{112}
-\frac{2}{7}c^{17}_{112}
-\frac{1}{7}c^{19}_{112}
+\frac{1}{7}c^{21}_{112}
+\frac{1}{7}c^{23}_{112}
$;
$-\frac{2}{7}c^{1}_{56}
+\frac{1}{7}c^{3}_{56}
+\frac{1}{7}c^{5}_{56}
-\frac{1}{7}c^{7}_{56}
-\frac{1}{7}c^{9}_{56}
+\frac{1}{7}c^{10}_{56}
+\frac{1}{7}c^{11}_{56}
$,
$-\frac{1}{7}c^{1}_{112}
-\frac{1}{7}c^{5}_{112}
+\frac{1}{7}c^{15}_{112}
-\frac{1}{7}c^{19}_{112}
$,
$\frac{1}{7}c^{2}_{56}
-\frac{1}{7}c^{3}_{56}
+\frac{1}{7}c^{11}_{56}
$,
$-\frac{1}{7}c^{3}_{112}
+\frac{1}{7}c^{9}_{112}
-\frac{1}{7}c^{11}_{112}
-\frac{1}{7}c^{23}_{112}
$;
$\frac{1}{7}c^{5}_{56}
+\frac{1}{7}c^{6}_{56}
+\frac{1}{7}c^{9}_{56}
$,
$-\frac{1}{7}c^{1}_{112}
-\frac{1}{7}c^{3}_{112}
+\frac{1}{7}c^{5}_{112}
+\frac{1}{7}c^{7}_{112}
+\frac{1}{7}c^{9}_{112}
+\frac{1}{7}c^{11}_{112}
-\frac{2}{7}c^{13}_{112}
-\frac{1}{7}c^{15}_{112}
-\frac{2}{7}c^{17}_{112}
-\frac{1}{7}c^{19}_{112}
+\frac{1}{7}c^{21}_{112}
+\frac{1}{7}c^{23}_{112}
$,
$\frac{1}{7}c^{2}_{56}
-\frac{1}{7}c^{3}_{56}
+\frac{1}{7}c^{11}_{56}
$;
$-\frac{1}{7}c^{5}_{56}
-\frac{1}{7}c^{6}_{56}
-\frac{1}{7}c^{9}_{56}
$,
$\frac{1}{7}c^{1}_{112}
+\frac{1}{7}c^{5}_{112}
-\frac{1}{7}c^{15}_{112}
+\frac{1}{7}c^{19}_{112}
$;
$\frac{2}{7}c^{1}_{56}
-\frac{1}{7}c^{3}_{56}
-\frac{1}{7}c^{5}_{56}
+\frac{1}{7}c^{7}_{56}
+\frac{1}{7}c^{9}_{56}
-\frac{1}{7}c^{10}_{56}
-\frac{1}{7}c^{11}_{56}
$)

Fail:
cnd( Tr$_I(\rho(\mathfrak s))$ ) =
56 does not divide
 ord($\rho(\mathfrak t)$) =
28, I = [ 3/28 ]. Prop. B.4 (2)

 \ \color{black}

 \color{blue}

\noindent 576: (dims,levels) = $(6;30
)$,
irreps = $3_{3}^{1,0}
\hskip -1.5pt \otimes \hskip -1.5pt
2_{5}^{2}
\hskip -1.5pt \otimes \hskip -1.5pt
1_{2}^{1,0}$,
pord$(\rho_\text{isum}(\mathfrak{t})) = 15$,

\vskip 0.7ex
\hangindent=5.5em \hangafter=1
{\white .}\hskip 1em $\rho_\text{isum}(\mathfrak{t})$ =
 $( \frac{1}{10},
\frac{9}{10},
\frac{7}{30},
\frac{13}{30},
\frac{17}{30},
\frac{23}{30} )
$,

\vskip 0.7ex
\hangindent=5.5em \hangafter=1
{\white .}\hskip 1em $\rho_\text{isum}(\mathfrak{s})$ =
$\mathrm{i}$($\frac{1}{3\sqrt{5}}c^{1}_{20}
$,
$\frac{1}{3\sqrt{5}}c^{3}_{20}
$,
$\frac{2}{3\sqrt{5}}c^{3}_{20}
$,
$\frac{2}{3\sqrt{5}}c^{1}_{20}
$,
$\frac{2}{3\sqrt{5}}c^{3}_{20}
$,
$\frac{2}{3\sqrt{5}}c^{1}_{20}
$;\ \ 
$-\frac{1}{3\sqrt{5}}c^{1}_{20}
$,
$-\frac{2}{3\sqrt{5}}c^{1}_{20}
$,
$\frac{2}{3\sqrt{5}}c^{3}_{20}
$,
$-\frac{2}{3\sqrt{5}}c^{1}_{20}
$,
$\frac{2}{3\sqrt{5}}c^{3}_{20}
$;\ \ 
$-\frac{1}{3\sqrt{5}}c^{1}_{20}
$,
$\frac{1}{3\sqrt{5}}c^{3}_{20}
$,
$\frac{2}{3\sqrt{5}}c^{1}_{20}
$,
$-\frac{2}{3\sqrt{5}}c^{3}_{20}
$;\ \ 
$\frac{1}{3\sqrt{5}}c^{1}_{20}
$,
$-\frac{2}{3\sqrt{5}}c^{3}_{20}
$,
$-\frac{2}{3\sqrt{5}}c^{1}_{20}
$;\ \ 
$-\frac{1}{3\sqrt{5}}c^{1}_{20}
$,
$\frac{1}{3\sqrt{5}}c^{3}_{20}
$;\ \ 
$\frac{1}{3\sqrt{5}}c^{1}_{20}
$)

Pass. 

 \ \color{black}

 \color{blue}

\noindent 577: (dims,levels) = $(6;30
)$,
irreps = $3_{3}^{1,0}
\hskip -1.5pt \otimes \hskip -1.5pt
2_{5}^{1}
\hskip -1.5pt \otimes \hskip -1.5pt
1_{2}^{1,0}$,
pord$(\rho_\text{isum}(\mathfrak{t})) = 15$,

\vskip 0.7ex
\hangindent=5.5em \hangafter=1
{\white .}\hskip 1em $\rho_\text{isum}(\mathfrak{t})$ =
 $( \frac{3}{10},
\frac{7}{10},
\frac{1}{30},
\frac{11}{30},
\frac{19}{30},
\frac{29}{30} )
$,

\vskip 0.7ex
\hangindent=5.5em \hangafter=1
{\white .}\hskip 1em $\rho_\text{isum}(\mathfrak{s})$ =
$\mathrm{i}$($\frac{1}{3\sqrt{5}}c^{3}_{20}
$,
$\frac{1}{3\sqrt{5}}c^{1}_{20}
$,
$\frac{2}{3\sqrt{5}}c^{1}_{20}
$,
$\frac{2}{3\sqrt{5}}c^{1}_{20}
$,
$\frac{2}{3\sqrt{5}}c^{3}_{20}
$,
$\frac{2}{3\sqrt{5}}c^{3}_{20}
$;\ \ 
$-\frac{1}{3\sqrt{5}}c^{3}_{20}
$,
$-\frac{2}{3\sqrt{5}}c^{3}_{20}
$,
$-\frac{2}{3\sqrt{5}}c^{3}_{20}
$,
$\frac{2}{3\sqrt{5}}c^{1}_{20}
$,
$\frac{2}{3\sqrt{5}}c^{1}_{20}
$;\ \ 
$-\frac{1}{3\sqrt{5}}c^{3}_{20}
$,
$\frac{2}{3\sqrt{5}}c^{3}_{20}
$,
$\frac{1}{3\sqrt{5}}c^{1}_{20}
$,
$-\frac{2}{3\sqrt{5}}c^{1}_{20}
$;\ \ 
$-\frac{1}{3\sqrt{5}}c^{3}_{20}
$,
$-\frac{2}{3\sqrt{5}}c^{1}_{20}
$,
$\frac{1}{3\sqrt{5}}c^{1}_{20}
$;\ \ 
$\frac{1}{3\sqrt{5}}c^{3}_{20}
$,
$-\frac{2}{3\sqrt{5}}c^{3}_{20}
$;\ \ 
$\frac{1}{3\sqrt{5}}c^{3}_{20}
$)

Pass. 

 \ \color{black}

\noindent 578: (dims,levels) = $(6;30
)$,
irreps = $3_{5}^{1}
\hskip -1.5pt \otimes \hskip -1.5pt
2_{2}^{1,0}
\hskip -1.5pt \otimes \hskip -1.5pt
1_{3}^{1,0}$,
pord$(\rho_\text{isum}(\mathfrak{t})) = 10$,

\vskip 0.7ex
\hangindent=5.5em \hangafter=1
{\white .}\hskip 1em $\rho_\text{isum}(\mathfrak{t})$ =
 $( \frac{1}{3},
\frac{5}{6},
\frac{2}{15},
\frac{8}{15},
\frac{1}{30},
\frac{19}{30} )
$,

\vskip 0.7ex
\hangindent=5.5em \hangafter=1
{\white .}\hskip 1em $\rho_\text{isum}(\mathfrak{s})$ =
($-\sqrt{\frac{1}{20}}$,
$\sqrt{\frac{3}{20}}$,
$-\sqrt{\frac{1}{10}}$,
$-\sqrt{\frac{1}{10}}$,
$\sqrt{\frac{3}{10}}$,
$\sqrt{\frac{3}{10}}$;
$\sqrt{\frac{1}{20}}$,
$\sqrt{\frac{3}{10}}$,
$\sqrt{\frac{3}{10}}$,
$\sqrt{\frac{1}{10}}$,
$\sqrt{\frac{1}{10}}$;
$\frac{5+\sqrt{5}}{20}$,
$\frac{-5+\sqrt{5}}{20}$,
$\frac{3}{2\sqrt{15}}c^{1}_{5}
$,
$-\frac{3}{2\sqrt{15}\mathrm{i}}s^{3}_{20}
$;
$\frac{5+\sqrt{5}}{20}$,
$-\frac{3}{2\sqrt{15}\mathrm{i}}s^{3}_{20}
$,
$\frac{3}{2\sqrt{15}}c^{1}_{5}
$;
$-\frac{5+\sqrt{5}}{20}$,
$\frac{5-\sqrt{5}}{20}$;
$-\frac{5+\sqrt{5}}{20}$)

Fail:
cnd($\rho(\mathfrak s)_\mathrm{ndeg}$) = 120 does not divide
 ord($\rho(\mathfrak t)$)=30. Prop. B.4 (2)

 \ \color{black}

\noindent 579: (dims,levels) = $(6;30
)$,
irreps = $3_{5}^{3}
\hskip -1.5pt \otimes \hskip -1.5pt
2_{2}^{1,0}
\hskip -1.5pt \otimes \hskip -1.5pt
1_{3}^{1,0}$,
pord$(\rho_\text{isum}(\mathfrak{t})) = 10$,

\vskip 0.7ex
\hangindent=5.5em \hangafter=1
{\white .}\hskip 1em $\rho_\text{isum}(\mathfrak{t})$ =
 $( \frac{1}{3},
\frac{5}{6},
\frac{11}{15},
\frac{14}{15},
\frac{7}{30},
\frac{13}{30} )
$,

\vskip 0.7ex
\hangindent=5.5em \hangafter=1
{\white .}\hskip 1em $\rho_\text{isum}(\mathfrak{s})$ =
($\sqrt{\frac{1}{20}}$,
$\sqrt{\frac{3}{20}}$,
$-\sqrt{\frac{1}{10}}$,
$-\sqrt{\frac{1}{10}}$,
$\sqrt{\frac{3}{10}}$,
$\sqrt{\frac{3}{10}}$;
$-\sqrt{\frac{1}{20}}$,
$-\sqrt{\frac{3}{10}}$,
$-\sqrt{\frac{3}{10}}$,
$-\sqrt{\frac{1}{10}}$,
$-\sqrt{\frac{1}{10}}$;
$\frac{5-\sqrt{5}}{20}$,
$-\frac{5+\sqrt{5}}{20}$,
$-\frac{3}{2\sqrt{15}}c^{1}_{5}
$,
$\frac{3}{2\sqrt{15}\mathrm{i}}s^{3}_{20}
$;
$\frac{5-\sqrt{5}}{20}$,
$\frac{3}{2\sqrt{15}\mathrm{i}}s^{3}_{20}
$,
$-\frac{3}{2\sqrt{15}}c^{1}_{5}
$;
$\frac{-5+\sqrt{5}}{20}$,
$\frac{5+\sqrt{5}}{20}$;
$\frac{-5+\sqrt{5}}{20}$)

Fail:
cnd($\rho(\mathfrak s)_\mathrm{ndeg}$) = 120 does not divide
 ord($\rho(\mathfrak t)$)=30. Prop. B.4 (2)

 \ \color{black}

\noindent 580: (dims,levels) = $(6;30
)$,
irreps = $3_{5}^{3}
\hskip -1.5pt \otimes \hskip -1.5pt
2_{3}^{1,0}
\hskip -1.5pt \otimes \hskip -1.5pt
1_{2}^{1,0}$,
pord$(\rho_\text{isum}(\mathfrak{t})) = 15$,

\vskip 0.7ex
\hangindent=5.5em \hangafter=1
{\white .}\hskip 1em $\rho_\text{isum}(\mathfrak{t})$ =
 $( \frac{1}{2},
\frac{5}{6},
\frac{1}{10},
\frac{9}{10},
\frac{7}{30},
\frac{13}{30} )
$,

\vskip 0.7ex
\hangindent=5.5em \hangafter=1
{\white .}\hskip 1em $\rho_\text{isum}(\mathfrak{s})$ =
$\mathrm{i}$($-\sqrt{\frac{1}{15}}$,
$-\sqrt{\frac{2}{15}}$,
$-\sqrt{\frac{2}{15}}$,
$-\sqrt{\frac{2}{15}}$,
$\sqrt{\frac{4}{15}}$,
$\sqrt{\frac{4}{15}}$;\ \ 
$\sqrt{\frac{1}{15}}$,
$-\sqrt{\frac{4}{15}}$,
$-\sqrt{\frac{4}{15}}$,
$-\sqrt{\frac{2}{15}}$,
$-\sqrt{\frac{2}{15}}$;\ \ 
$-\frac{1}{\sqrt{15}}c^{1}_{5}
$,
$\frac{1}{\sqrt{15}\mathrm{i}}s^{3}_{20}
$,
$-\frac{2}{\sqrt{30}\mathrm{i}}s^{3}_{20}
$,
$\frac{2}{\sqrt{30}}c^{1}_{5}
$;\ \ 
$-\frac{1}{\sqrt{15}}c^{1}_{5}
$,
$\frac{2}{\sqrt{30}}c^{1}_{5}
$,
$-\frac{2}{\sqrt{30}\mathrm{i}}s^{3}_{20}
$;\ \ 
$\frac{1}{\sqrt{15}}c^{1}_{5}
$,
$-\frac{1}{\sqrt{15}\mathrm{i}}s^{3}_{20}
$;\ \ 
$\frac{1}{\sqrt{15}}c^{1}_{5}
$)

Fail:
cnd($\rho(\mathfrak s)_\mathrm{ndeg}$) = 120 does not divide
 ord($\rho(\mathfrak t)$)=30. Prop. B.4 (2)

 \ \color{black}

\noindent 581: (dims,levels) = $(6;30
)$,
irreps = $3_{5}^{1}
\hskip -1.5pt \otimes \hskip -1.5pt
2_{3}^{1,0}
\hskip -1.5pt \otimes \hskip -1.5pt
1_{2}^{1,0}$,
pord$(\rho_\text{isum}(\mathfrak{t})) = 15$,

\vskip 0.7ex
\hangindent=5.5em \hangafter=1
{\white .}\hskip 1em $\rho_\text{isum}(\mathfrak{t})$ =
 $( \frac{1}{2},
\frac{5}{6},
\frac{3}{10},
\frac{7}{10},
\frac{1}{30},
\frac{19}{30} )
$,

\vskip 0.7ex
\hangindent=5.5em \hangafter=1
{\white .}\hskip 1em $\rho_\text{isum}(\mathfrak{s})$ =
$\mathrm{i}$($\sqrt{\frac{1}{15}}$,
$-\sqrt{\frac{2}{15}}$,
$-\sqrt{\frac{2}{15}}$,
$-\sqrt{\frac{2}{15}}$,
$\sqrt{\frac{4}{15}}$,
$\sqrt{\frac{4}{15}}$;\ \ 
$-\sqrt{\frac{1}{15}}$,
$\sqrt{\frac{4}{15}}$,
$\sqrt{\frac{4}{15}}$,
$\sqrt{\frac{2}{15}}$,
$\sqrt{\frac{2}{15}}$;\ \ 
$-\frac{1}{\sqrt{15}\mathrm{i}}s^{3}_{20}
$,
$\frac{1}{\sqrt{15}}c^{1}_{5}
$,
$-\frac{2}{\sqrt{30}}c^{1}_{5}
$,
$\frac{2}{\sqrt{30}\mathrm{i}}s^{3}_{20}
$;\ \ 
$-\frac{1}{\sqrt{15}\mathrm{i}}s^{3}_{20}
$,
$\frac{2}{\sqrt{30}\mathrm{i}}s^{3}_{20}
$,
$-\frac{2}{\sqrt{30}}c^{1}_{5}
$;\ \ 
$\frac{1}{\sqrt{15}\mathrm{i}}s^{3}_{20}
$,
$-\frac{1}{\sqrt{15}}c^{1}_{5}
$;\ \ 
$\frac{1}{\sqrt{15}\mathrm{i}}s^{3}_{20}
$)

Fail:
cnd($\rho(\mathfrak s)_\mathrm{ndeg}$) = 120 does not divide
 ord($\rho(\mathfrak t)$)=30. Prop. B.4 (2)

 \ \color{black}

\noindent 582: (dims,levels) = $(6;30
)$,
irreps = $6_{5}^{1}
\hskip -1.5pt \otimes \hskip -1.5pt
1_{3}^{1,0}
\hskip -1.5pt \otimes \hskip -1.5pt
1_{2}^{1,0}$,
pord$(\rho_\text{isum}(\mathfrak{t})) = 5$,

\vskip 0.7ex
\hangindent=5.5em \hangafter=1
{\white .}\hskip 1em $\rho_\text{isum}(\mathfrak{t})$ =
 $( \frac{5}{6},
\frac{5}{6},
\frac{1}{30},
\frac{7}{30},
\frac{13}{30},
\frac{19}{30} )
$,

\vskip 0.7ex
\hangindent=5.5em \hangafter=1
{\white .}\hskip 1em $\rho_\text{isum}(\mathfrak{s})$ =
$\mathrm{i}$($-\frac{1}{5}c^{1}_{20}
$,
$\frac{1}{5}c^{3}_{20}
$,
$\frac{4}{5\sqrt{10}}c^{1}_{20}
-\frac{2}{5\sqrt{10}}c^{3}_{20}
$,
$\frac{1}{5}c^{1}_{20}
-\frac{1}{5}c^{3}_{20}
$,
$\frac{1}{5}c^{1}_{20}
+\frac{1}{5}c^{3}_{20}
$,
$\frac{1}{5}c^{3}_{40}
+\frac{1}{5}c^{7}_{40}
$;\ \ 
$\frac{1}{5}c^{1}_{20}
$,
$\frac{1}{5}c^{3}_{40}
+\frac{1}{5}c^{7}_{40}
$,
$-\frac{1}{5}c^{1}_{20}
-\frac{1}{5}c^{3}_{20}
$,
$\frac{1}{5}c^{1}_{20}
-\frac{1}{5}c^{3}_{20}
$,
$-\frac{4}{5\sqrt{10}}c^{1}_{20}
+\frac{2}{5\sqrt{10}}c^{3}_{20}
$;\ \ 
$\frac{1}{5}c^{3}_{20}
$,
$\frac{4}{5\sqrt{10}}c^{1}_{20}
-\frac{2}{5\sqrt{10}}c^{3}_{20}
$,
$-\frac{1}{5}c^{3}_{40}
-\frac{1}{5}c^{7}_{40}
$,
$\frac{1}{5}c^{1}_{20}
$;\ \ 
$-\frac{1}{5}c^{1}_{20}
$,
$-\frac{1}{5}c^{3}_{20}
$,
$\frac{1}{5}c^{3}_{40}
+\frac{1}{5}c^{7}_{40}
$;\ \ 
$\frac{1}{5}c^{1}_{20}
$,
$\frac{4}{5\sqrt{10}}c^{1}_{20}
-\frac{2}{5\sqrt{10}}c^{3}_{20}
$;\ \ 
$-\frac{1}{5}c^{3}_{20}
$)

Fail:
cnd($\rho(\mathfrak s)_\mathrm{ndeg}$) = 40 does not divide
 ord($\rho(\mathfrak t)$)=30. Prop. B.4 (2)

 \ \color{black}

\noindent 583: (dims,levels) = $(6;32
)$,
irreps = $6_{32,1}^{1,0}$,
pord$(\rho_\text{isum}(\mathfrak{t})) = 32$,

\vskip 0.7ex
\hangindent=5.5em \hangafter=1
{\white .}\hskip 1em $\rho_\text{isum}(\mathfrak{t})$ =
 $( 0,
\frac{1}{8},
\frac{3}{32},
\frac{11}{32},
\frac{19}{32},
\frac{27}{32} )
$,

\vskip 0.7ex
\hangindent=5.5em \hangafter=1
{\white .}\hskip 1em $\rho_\text{isum}(\mathfrak{s})$ =
$\mathrm{i}$($0$,
$0$,
$\frac{1}{2}$,
$\frac{1}{2}$,
$\frac{1}{2}$,
$\frac{1}{2}$;\ \ 
$0$,
$\frac{1}{2}$,
$-\frac{1}{2}$,
$\frac{1}{2}$,
$-\frac{1}{2}$;\ \ 
$-\frac{1}{4}c^{1}_{16}
$,
$-\frac{1}{4}c^{3}_{16}
$,
$\frac{1}{4}c^{1}_{16}
$,
$\frac{1}{4}c^{3}_{16}
$;\ \ 
$\frac{1}{4}c^{1}_{16}
$,
$\frac{1}{4}c^{3}_{16}
$,
$-\frac{1}{4}c^{1}_{16}
$;\ \ 
$-\frac{1}{4}c^{1}_{16}
$,
$-\frac{1}{4}c^{3}_{16}
$;\ \ 
$\frac{1}{4}c^{1}_{16}
$)

Fail:
$\sigma(\rho(\mathfrak s)_\mathrm{ndeg}) \neq
 (\rho(\mathfrak t)^a \rho(\mathfrak s) \rho(\mathfrak t)^b
 \rho(\mathfrak s) \rho(\mathfrak t)^a)_\mathrm{ndeg}$,
 $\sigma = a$ = 3. Prop. B.5 (3) eqn. (B.25)

 \ \color{black}

\noindent 584: (dims,levels) = $(6;32
)$,
irreps = $6_{32,2}^{1,0}$,
pord$(\rho_\text{isum}(\mathfrak{t})) = 32$,

\vskip 0.7ex
\hangindent=5.5em \hangafter=1
{\white .}\hskip 1em $\rho_\text{isum}(\mathfrak{t})$ =
 $( 0,
\frac{1}{8},
\frac{7}{32},
\frac{15}{32},
\frac{23}{32},
\frac{31}{32} )
$,

\vskip 0.7ex
\hangindent=5.5em \hangafter=1
{\white .}\hskip 1em $\rho_\text{isum}(\mathfrak{s})$ =
($0$,
$0$,
$\frac{1}{2}$,
$\frac{1}{2}$,
$\frac{1}{2}$,
$\frac{1}{2}$;
$0$,
$\frac{1}{2}$,
$-\frac{1}{2}$,
$\frac{1}{2}$,
$-\frac{1}{2}$;
$-\frac{1}{4}c^{1}_{16}
$,
$-\frac{1}{4}c^{3}_{16}
$,
$\frac{1}{4}c^{1}_{16}
$,
$\frac{1}{4}c^{3}_{16}
$;
$\frac{1}{4}c^{1}_{16}
$,
$\frac{1}{4}c^{3}_{16}
$,
$-\frac{1}{4}c^{1}_{16}
$;
$-\frac{1}{4}c^{1}_{16}
$,
$-\frac{1}{4}c^{3}_{16}
$;
$\frac{1}{4}c^{1}_{16}
$)

Fail:
$\sigma(\rho(\mathfrak s)_\mathrm{ndeg}) \neq
 (\rho(\mathfrak t)^a \rho(\mathfrak s) \rho(\mathfrak t)^b
 \rho(\mathfrak s) \rho(\mathfrak t)^a)_\mathrm{ndeg}$,
 $\sigma = a$ = 3. Prop. B.5 (3) eqn. (B.25)

 \ \color{black}

\noindent 585: (dims,levels) = $(6;32
)$,
irreps = $6_{32,1}^{3,0}$,
pord$(\rho_\text{isum}(\mathfrak{t})) = 32$,

\vskip 0.7ex
\hangindent=5.5em \hangafter=1
{\white .}\hskip 1em $\rho_\text{isum}(\mathfrak{t})$ =
 $( 0,
\frac{3}{8},
\frac{1}{32},
\frac{9}{32},
\frac{17}{32},
\frac{25}{32} )
$,

\vskip 0.7ex
\hangindent=5.5em \hangafter=1
{\white .}\hskip 1em $\rho_\text{isum}(\mathfrak{s})$ =
$\mathrm{i}$($0$,
$0$,
$\frac{1}{2}$,
$\frac{1}{2}$,
$\frac{1}{2}$,
$\frac{1}{2}$;\ \ 
$0$,
$\frac{1}{2}$,
$-\frac{1}{2}$,
$\frac{1}{2}$,
$-\frac{1}{2}$;\ \ 
$-\frac{1}{4}c^{3}_{16}
$,
$-\frac{1}{4}c^{1}_{16}
$,
$\frac{1}{4}c^{3}_{16}
$,
$\frac{1}{4}c^{1}_{16}
$;\ \ 
$\frac{1}{4}c^{3}_{16}
$,
$\frac{1}{4}c^{1}_{16}
$,
$-\frac{1}{4}c^{3}_{16}
$;\ \ 
$-\frac{1}{4}c^{3}_{16}
$,
$-\frac{1}{4}c^{1}_{16}
$;\ \ 
$\frac{1}{4}c^{3}_{16}
$)

Fail:
$\sigma(\rho(\mathfrak s)_\mathrm{ndeg}) \neq
 (\rho(\mathfrak t)^a \rho(\mathfrak s) \rho(\mathfrak t)^b
 \rho(\mathfrak s) \rho(\mathfrak t)^a)_\mathrm{ndeg}$,
 $\sigma = a$ = 3. Prop. B.5 (3) eqn. (B.25)

 \ \color{black}

\noindent 586: (dims,levels) = $(6;32
)$,
irreps = $6_{32,2}^{3,0}$,
pord$(\rho_\text{isum}(\mathfrak{t})) = 32$,

\vskip 0.7ex
\hangindent=5.5em \hangafter=1
{\white .}\hskip 1em $\rho_\text{isum}(\mathfrak{t})$ =
 $( 0,
\frac{3}{8},
\frac{5}{32},
\frac{13}{32},
\frac{21}{32},
\frac{29}{32} )
$,

\vskip 0.7ex
\hangindent=5.5em \hangafter=1
{\white .}\hskip 1em $\rho_\text{isum}(\mathfrak{s})$ =
($0$,
$0$,
$\frac{1}{2}$,
$\frac{1}{2}$,
$\frac{1}{2}$,
$\frac{1}{2}$;
$0$,
$\frac{1}{2}$,
$-\frac{1}{2}$,
$\frac{1}{2}$,
$-\frac{1}{2}$;
$-\frac{1}{4}c^{3}_{16}
$,
$-\frac{1}{4}c^{1}_{16}
$,
$\frac{1}{4}c^{3}_{16}
$,
$\frac{1}{4}c^{1}_{16}
$;
$\frac{1}{4}c^{3}_{16}
$,
$\frac{1}{4}c^{1}_{16}
$,
$-\frac{1}{4}c^{3}_{16}
$;
$-\frac{1}{4}c^{3}_{16}
$,
$-\frac{1}{4}c^{1}_{16}
$;
$\frac{1}{4}c^{3}_{16}
$)

Fail:
$\sigma(\rho(\mathfrak s)_\mathrm{ndeg}) \neq
 (\rho(\mathfrak t)^a \rho(\mathfrak s) \rho(\mathfrak t)^b
 \rho(\mathfrak s) \rho(\mathfrak t)^a)_\mathrm{ndeg}$,
 $\sigma = a$ = 3. Prop. B.5 (3) eqn. (B.25)

 \ \color{black}

\noindent 587: (dims,levels) = $(6;32
)$,
irreps = $6_{32,2}^{5,0}$,
pord$(\rho_\text{isum}(\mathfrak{t})) = 32$,

\vskip 0.7ex
\hangindent=5.5em \hangafter=1
{\white .}\hskip 1em $\rho_\text{isum}(\mathfrak{t})$ =
 $( 0,
\frac{5}{8},
\frac{3}{32},
\frac{11}{32},
\frac{19}{32},
\frac{27}{32} )
$,

\vskip 0.7ex
\hangindent=5.5em \hangafter=1
{\white .}\hskip 1em $\rho_\text{isum}(\mathfrak{s})$ =
($0$,
$0$,
$\frac{1}{2}$,
$\frac{1}{2}$,
$\frac{1}{2}$,
$\frac{1}{2}$;
$0$,
$\frac{1}{2}$,
$-\frac{1}{2}$,
$\frac{1}{2}$,
$-\frac{1}{2}$;
$\frac{1}{4}c^{3}_{16}
$,
$-\frac{1}{4}c^{1}_{16}
$,
$-\frac{1}{4}c^{3}_{16}
$,
$\frac{1}{4}c^{1}_{16}
$;
$-\frac{1}{4}c^{3}_{16}
$,
$\frac{1}{4}c^{1}_{16}
$,
$\frac{1}{4}c^{3}_{16}
$;
$\frac{1}{4}c^{3}_{16}
$,
$-\frac{1}{4}c^{1}_{16}
$;
$-\frac{1}{4}c^{3}_{16}
$)

Fail:
$\sigma(\rho(\mathfrak s)_\mathrm{ndeg}) \neq
 (\rho(\mathfrak t)^a \rho(\mathfrak s) \rho(\mathfrak t)^b
 \rho(\mathfrak s) \rho(\mathfrak t)^a)_\mathrm{ndeg}$,
 $\sigma = a$ = 3. Prop. B.5 (3) eqn. (B.25)

 \ \color{black}

\noindent 588: (dims,levels) = $(6;32
)$,
irreps = $6_{32,1}^{5,0}$,
pord$(\rho_\text{isum}(\mathfrak{t})) = 32$,

\vskip 0.7ex
\hangindent=5.5em \hangafter=1
{\white .}\hskip 1em $\rho_\text{isum}(\mathfrak{t})$ =
 $( 0,
\frac{5}{8},
\frac{7}{32},
\frac{15}{32},
\frac{23}{32},
\frac{31}{32} )
$,

\vskip 0.7ex
\hangindent=5.5em \hangafter=1
{\white .}\hskip 1em $\rho_\text{isum}(\mathfrak{s})$ =
$\mathrm{i}$($0$,
$0$,
$\frac{1}{2}$,
$\frac{1}{2}$,
$\frac{1}{2}$,
$\frac{1}{2}$;\ \ 
$0$,
$\frac{1}{2}$,
$-\frac{1}{2}$,
$\frac{1}{2}$,
$-\frac{1}{2}$;\ \ 
$-\frac{1}{4}c^{3}_{16}
$,
$\frac{1}{4}c^{1}_{16}
$,
$\frac{1}{4}c^{3}_{16}
$,
$-\frac{1}{4}c^{1}_{16}
$;\ \ 
$\frac{1}{4}c^{3}_{16}
$,
$-\frac{1}{4}c^{1}_{16}
$,
$-\frac{1}{4}c^{3}_{16}
$;\ \ 
$-\frac{1}{4}c^{3}_{16}
$,
$\frac{1}{4}c^{1}_{16}
$;\ \ 
$\frac{1}{4}c^{3}_{16}
$)

Fail:
$\sigma(\rho(\mathfrak s)_\mathrm{ndeg}) \neq
 (\rho(\mathfrak t)^a \rho(\mathfrak s) \rho(\mathfrak t)^b
 \rho(\mathfrak s) \rho(\mathfrak t)^a)_\mathrm{ndeg}$,
 $\sigma = a$ = 3. Prop. B.5 (3) eqn. (B.25)

 \ \color{black}

\noindent 589: (dims,levels) = $(6;32
)$,
irreps = $6_{32,2}^{7,0}$,
pord$(\rho_\text{isum}(\mathfrak{t})) = 32$,

\vskip 0.7ex
\hangindent=5.5em \hangafter=1
{\white .}\hskip 1em $\rho_\text{isum}(\mathfrak{t})$ =
 $( 0,
\frac{7}{8},
\frac{1}{32},
\frac{9}{32},
\frac{17}{32},
\frac{25}{32} )
$,

\vskip 0.7ex
\hangindent=5.5em \hangafter=1
{\white .}\hskip 1em $\rho_\text{isum}(\mathfrak{s})$ =
($0$,
$0$,
$\frac{1}{2}$,
$\frac{1}{2}$,
$\frac{1}{2}$,
$\frac{1}{2}$;
$0$,
$\frac{1}{2}$,
$-\frac{1}{2}$,
$\frac{1}{2}$,
$-\frac{1}{2}$;
$\frac{1}{4}c^{1}_{16}
$,
$-\frac{1}{4}c^{3}_{16}
$,
$-\frac{1}{4}c^{1}_{16}
$,
$\frac{1}{4}c^{3}_{16}
$;
$-\frac{1}{4}c^{1}_{16}
$,
$\frac{1}{4}c^{3}_{16}
$,
$\frac{1}{4}c^{1}_{16}
$;
$\frac{1}{4}c^{1}_{16}
$,
$-\frac{1}{4}c^{3}_{16}
$;
$-\frac{1}{4}c^{1}_{16}
$)

Fail:
$\sigma(\rho(\mathfrak s)_\mathrm{ndeg}) \neq
 (\rho(\mathfrak t)^a \rho(\mathfrak s) \rho(\mathfrak t)^b
 \rho(\mathfrak s) \rho(\mathfrak t)^a)_\mathrm{ndeg}$,
 $\sigma = a$ = 3. Prop. B.5 (3) eqn. (B.25)

 \ \color{black}

\noindent 590: (dims,levels) = $(6;32
)$,
irreps = $6_{32,1}^{7,0}$,
pord$(\rho_\text{isum}(\mathfrak{t})) = 32$,

\vskip 0.7ex
\hangindent=5.5em \hangafter=1
{\white .}\hskip 1em $\rho_\text{isum}(\mathfrak{t})$ =
 $( 0,
\frac{7}{8},
\frac{5}{32},
\frac{13}{32},
\frac{21}{32},
\frac{29}{32} )
$,

\vskip 0.7ex
\hangindent=5.5em \hangafter=1
{\white .}\hskip 1em $\rho_\text{isum}(\mathfrak{s})$ =
$\mathrm{i}$($0$,
$0$,
$\frac{1}{2}$,
$\frac{1}{2}$,
$\frac{1}{2}$,
$\frac{1}{2}$;\ \ 
$0$,
$\frac{1}{2}$,
$-\frac{1}{2}$,
$\frac{1}{2}$,
$-\frac{1}{2}$;\ \ 
$-\frac{1}{4}c^{1}_{16}
$,
$\frac{1}{4}c^{3}_{16}
$,
$\frac{1}{4}c^{1}_{16}
$,
$-\frac{1}{4}c^{3}_{16}
$;\ \ 
$\frac{1}{4}c^{1}_{16}
$,
$-\frac{1}{4}c^{3}_{16}
$,
$-\frac{1}{4}c^{1}_{16}
$;\ \ 
$-\frac{1}{4}c^{1}_{16}
$,
$\frac{1}{4}c^{3}_{16}
$;\ \ 
$\frac{1}{4}c^{1}_{16}
$)

Fail:
$\sigma(\rho(\mathfrak s)_\mathrm{ndeg}) \neq
 (\rho(\mathfrak t)^a \rho(\mathfrak s) \rho(\mathfrak t)^b
 \rho(\mathfrak s) \rho(\mathfrak t)^a)_\mathrm{ndeg}$,
 $\sigma = a$ = 3. Prop. B.5 (3) eqn. (B.25)

 \ \color{black}

\noindent 591: (dims,levels) = $(6;33
)$,
irreps = $6_{11}^{7}
\hskip -1.5pt \otimes \hskip -1.5pt
1_{3}^{1,0}$,
pord$(\rho_\text{isum}(\mathfrak{t})) = 11$,

\vskip 0.7ex
\hangindent=5.5em \hangafter=1
{\white .}\hskip 1em $\rho_\text{isum}(\mathfrak{t})$ =
 $( \frac{1}{3},
\frac{2}{33},
\frac{8}{33},
\frac{17}{33},
\frac{29}{33},
\frac{32}{33} )
$,

\vskip 0.7ex
\hangindent=5.5em \hangafter=1
{\white .}\hskip 1em $\rho_\text{isum}(\mathfrak{s})$ =
$\mathrm{i}$($\sqrt{\frac{1}{11}}$,
$\sqrt{\frac{2}{11}}$,
$\sqrt{\frac{2}{11}}$,
$\sqrt{\frac{2}{11}}$,
$\sqrt{\frac{2}{11}}$,
$\sqrt{\frac{2}{11}}$;\ \ 
$-\frac{1}{\sqrt{11}\mathrm{i}}s^{9}_{44}
$,
$\frac{1}{\sqrt{11}}c^{1}_{11}
$,
$\frac{1}{\sqrt{11}}c^{3}_{11}
$,
$\frac{1}{\sqrt{11}}c^{4}_{11}
$,
$\frac{1}{\sqrt{11}}c^{2}_{11}
$;\ \ 
$\frac{1}{\sqrt{11}}c^{2}_{11}
$,
$-\frac{1}{\sqrt{11}\mathrm{i}}s^{9}_{44}
$,
$\frac{1}{\sqrt{11}}c^{3}_{11}
$,
$\frac{1}{\sqrt{11}}c^{4}_{11}
$;\ \ 
$\frac{1}{\sqrt{11}}c^{4}_{11}
$,
$\frac{1}{\sqrt{11}}c^{2}_{11}
$,
$\frac{1}{\sqrt{11}}c^{1}_{11}
$;\ \ 
$\frac{1}{\sqrt{11}}c^{1}_{11}
$,
$-\frac{1}{\sqrt{11}\mathrm{i}}s^{9}_{44}
$;\ \ 
$\frac{1}{\sqrt{11}}c^{3}_{11}
$)

Fail:
cnd($\rho(\mathfrak s)_\mathrm{ndeg}$) = 88 does not divide
 ord($\rho(\mathfrak t)$)=33. Prop. B.4 (2)

 \ \color{black}

\noindent 592: (dims,levels) = $(6;33
)$,
irreps = $6_{11}^{1}
\hskip -1.5pt \otimes \hskip -1.5pt
1_{3}^{1,0}$,
pord$(\rho_\text{isum}(\mathfrak{t})) = 11$,

\vskip 0.7ex
\hangindent=5.5em \hangafter=1
{\white .}\hskip 1em $\rho_\text{isum}(\mathfrak{t})$ =
 $( \frac{1}{3},
\frac{5}{33},
\frac{14}{33},
\frac{20}{33},
\frac{23}{33},
\frac{26}{33} )
$,

\vskip 0.7ex
\hangindent=5.5em \hangafter=1
{\white .}\hskip 1em $\rho_\text{isum}(\mathfrak{s})$ =
$\mathrm{i}$($-\sqrt{\frac{1}{11}}$,
$\sqrt{\frac{2}{11}}$,
$\sqrt{\frac{2}{11}}$,
$\sqrt{\frac{2}{11}}$,
$\sqrt{\frac{2}{11}}$,
$\sqrt{\frac{2}{11}}$;\ \ 
$-\frac{1}{\sqrt{11}}c^{4}_{11}
$,
$\frac{1}{\sqrt{11}\mathrm{i}}s^{9}_{44}
$,
$-\frac{1}{\sqrt{11}}c^{3}_{11}
$,
$-\frac{1}{\sqrt{11}}c^{1}_{11}
$,
$-\frac{1}{\sqrt{11}}c^{2}_{11}
$;\ \ 
$-\frac{1}{\sqrt{11}}c^{2}_{11}
$,
$-\frac{1}{\sqrt{11}}c^{1}_{11}
$,
$-\frac{1}{\sqrt{11}}c^{4}_{11}
$,
$-\frac{1}{\sqrt{11}}c^{3}_{11}
$;\ \ 
$\frac{1}{\sqrt{11}\mathrm{i}}s^{9}_{44}
$,
$-\frac{1}{\sqrt{11}}c^{2}_{11}
$,
$-\frac{1}{\sqrt{11}}c^{4}_{11}
$;\ \ 
$-\frac{1}{\sqrt{11}}c^{3}_{11}
$,
$\frac{1}{\sqrt{11}\mathrm{i}}s^{9}_{44}
$;\ \ 
$-\frac{1}{\sqrt{11}}c^{1}_{11}
$)

Fail:
cnd($\rho(\mathfrak s)_\mathrm{ndeg}$) = 88 does not divide
 ord($\rho(\mathfrak t)$)=33. Prop. B.4 (2)

 \ \color{black}

 \color{blue}

\noindent 593: (dims,levels) = $(6;35
)$,
irreps = $3_{7}^{3}
\hskip -1.5pt \otimes \hskip -1.5pt
2_{5}^{2}$,
pord$(\rho_\text{isum}(\mathfrak{t})) = 35$,

\vskip 0.7ex
\hangindent=5.5em \hangafter=1
{\white .}\hskip 1em $\rho_\text{isum}(\mathfrak{t})$ =
 $( \frac{1}{35},
\frac{4}{35},
\frac{9}{35},
\frac{11}{35},
\frac{16}{35},
\frac{29}{35} )
$,

\vskip 0.7ex
\hangindent=5.5em \hangafter=1
{\white .}\hskip 1em $\rho_\text{isum}(\mathfrak{s})$ =
$\mathrm{i}$($-\frac{4}{35}c^{1}_{140}
-\frac{3}{35}c^{3}_{140}
-\frac{1}{7}c^{5}_{140}
+\frac{1}{35}c^{7}_{140}
+\frac{1}{35}c^{9}_{140}
+\frac{4}{35}c^{13}_{140}
+\frac{2}{35}c^{15}_{140}
-\frac{3}{35}c^{17}_{140}
+\frac{9}{35}c^{19}_{140}
-\frac{4}{35}c^{21}_{140}
-\frac{2}{7}c^{23}_{140}
$,
$-\frac{1}{\sqrt{35}}c^{4}_{35}
+\frac{1}{\sqrt{35}}c^{11}_{35}
$,
$\frac{1}{\sqrt{35}}c^{1}_{35}
-\frac{1}{\sqrt{35}}c^{6}_{35}
$,
$\frac{2}{\sqrt{35}}c^{3}_{35}
+\frac{1}{\sqrt{35}}c^{4}_{35}
+\frac{1}{\sqrt{35}}c^{10}_{35}
+\frac{1}{\sqrt{35}}c^{11}_{35}
$,
$-\frac{1}{\sqrt{35}\mathrm{i}}s^{3}_{140}
-\frac{1}{\sqrt{35}\mathrm{i}}s^{17}_{140}
$,
$\frac{2}{35}c^{1}_{140}
-\frac{1}{35}c^{3}_{140}
-\frac{1}{7}c^{5}_{140}
-\frac{3}{35}c^{7}_{140}
+\frac{1}{5}c^{9}_{140}
-\frac{2}{35}c^{13}_{140}
-\frac{1}{35}c^{15}_{140}
-\frac{1}{35}c^{17}_{140}
+\frac{3}{35}c^{19}_{140}
+\frac{2}{35}c^{21}_{140}
-\frac{2}{7}c^{23}_{140}
$;\ \ 
$-\frac{1}{\sqrt{35}\mathrm{i}}s^{3}_{140}
-\frac{1}{\sqrt{35}\mathrm{i}}s^{17}_{140}
$,
$\frac{4}{35}c^{1}_{140}
+\frac{3}{35}c^{3}_{140}
+\frac{1}{7}c^{5}_{140}
-\frac{1}{35}c^{7}_{140}
-\frac{1}{35}c^{9}_{140}
-\frac{4}{35}c^{13}_{140}
-\frac{2}{35}c^{15}_{140}
+\frac{3}{35}c^{17}_{140}
-\frac{9}{35}c^{19}_{140}
+\frac{4}{35}c^{21}_{140}
+\frac{2}{7}c^{23}_{140}
$,
$-\frac{1}{\sqrt{35}}c^{1}_{35}
+\frac{1}{\sqrt{35}}c^{6}_{35}
$,
$-\frac{2}{35}c^{1}_{140}
+\frac{1}{35}c^{3}_{140}
+\frac{1}{7}c^{5}_{140}
+\frac{3}{35}c^{7}_{140}
-\frac{1}{5}c^{9}_{140}
+\frac{2}{35}c^{13}_{140}
+\frac{1}{35}c^{15}_{140}
+\frac{1}{35}c^{17}_{140}
-\frac{3}{35}c^{19}_{140}
-\frac{2}{35}c^{21}_{140}
+\frac{2}{7}c^{23}_{140}
$,
$\frac{2}{\sqrt{35}}c^{3}_{35}
+\frac{1}{\sqrt{35}}c^{4}_{35}
+\frac{1}{\sqrt{35}}c^{10}_{35}
+\frac{1}{\sqrt{35}}c^{11}_{35}
$;\ \ 
$\frac{2}{\sqrt{35}}c^{3}_{35}
+\frac{1}{\sqrt{35}}c^{4}_{35}
+\frac{1}{\sqrt{35}}c^{10}_{35}
+\frac{1}{\sqrt{35}}c^{11}_{35}
$,
$-\frac{2}{35}c^{1}_{140}
+\frac{1}{35}c^{3}_{140}
+\frac{1}{7}c^{5}_{140}
+\frac{3}{35}c^{7}_{140}
-\frac{1}{5}c^{9}_{140}
+\frac{2}{35}c^{13}_{140}
+\frac{1}{35}c^{15}_{140}
+\frac{1}{35}c^{17}_{140}
-\frac{3}{35}c^{19}_{140}
-\frac{2}{35}c^{21}_{140}
+\frac{2}{7}c^{23}_{140}
$,
$\frac{1}{\sqrt{35}}c^{4}_{35}
-\frac{1}{\sqrt{35}}c^{11}_{35}
$,
$-\frac{1}{\sqrt{35}\mathrm{i}}s^{3}_{140}
-\frac{1}{\sqrt{35}\mathrm{i}}s^{17}_{140}
$;\ \ 
$\frac{1}{\sqrt{35}\mathrm{i}}s^{3}_{140}
+\frac{1}{\sqrt{35}\mathrm{i}}s^{17}_{140}
$,
$-\frac{4}{35}c^{1}_{140}
-\frac{3}{35}c^{3}_{140}
-\frac{1}{7}c^{5}_{140}
+\frac{1}{35}c^{7}_{140}
+\frac{1}{35}c^{9}_{140}
+\frac{4}{35}c^{13}_{140}
+\frac{2}{35}c^{15}_{140}
-\frac{3}{35}c^{17}_{140}
+\frac{9}{35}c^{19}_{140}
-\frac{4}{35}c^{21}_{140}
-\frac{2}{7}c^{23}_{140}
$,
$\frac{1}{\sqrt{35}}c^{4}_{35}
-\frac{1}{\sqrt{35}}c^{11}_{35}
$;\ \ 
$-\frac{2}{\sqrt{35}}c^{3}_{35}
-\frac{1}{\sqrt{35}}c^{4}_{35}
-\frac{1}{\sqrt{35}}c^{10}_{35}
-\frac{1}{\sqrt{35}}c^{11}_{35}
$,
$-\frac{1}{\sqrt{35}}c^{1}_{35}
+\frac{1}{\sqrt{35}}c^{6}_{35}
$;\ \ 
$\frac{4}{35}c^{1}_{140}
+\frac{3}{35}c^{3}_{140}
+\frac{1}{7}c^{5}_{140}
-\frac{1}{35}c^{7}_{140}
-\frac{1}{35}c^{9}_{140}
-\frac{4}{35}c^{13}_{140}
-\frac{2}{35}c^{15}_{140}
+\frac{3}{35}c^{17}_{140}
-\frac{9}{35}c^{19}_{140}
+\frac{4}{35}c^{21}_{140}
+\frac{2}{7}c^{23}_{140}
$)

Pass. 

 \ \color{black}

 \color{blue}

\noindent 594: (dims,levels) = $(6;35
)$,
irreps = $3_{7}^{3}
\hskip -1.5pt \otimes \hskip -1.5pt
2_{5}^{1}$,
pord$(\rho_\text{isum}(\mathfrak{t})) = 35$,

\vskip 0.7ex
\hangindent=5.5em \hangafter=1
{\white .}\hskip 1em $\rho_\text{isum}(\mathfrak{t})$ =
 $( \frac{2}{35},
\frac{8}{35},
\frac{18}{35},
\frac{22}{35},
\frac{23}{35},
\frac{32}{35} )
$,

\vskip 0.7ex
\hangindent=5.5em \hangafter=1
{\white .}\hskip 1em $\rho_\text{isum}(\mathfrak{s})$ =
$\mathrm{i}$($\frac{1}{\sqrt{35}}c^{4}_{35}
-\frac{1}{\sqrt{35}}c^{11}_{35}
$,
$-\frac{1}{\sqrt{35}\mathrm{i}}s^{3}_{140}
-\frac{1}{\sqrt{35}\mathrm{i}}s^{17}_{140}
$,
$\frac{4}{35}c^{1}_{140}
+\frac{3}{35}c^{3}_{140}
+\frac{1}{7}c^{5}_{140}
-\frac{1}{35}c^{7}_{140}
-\frac{1}{35}c^{9}_{140}
-\frac{4}{35}c^{13}_{140}
-\frac{2}{35}c^{15}_{140}
+\frac{3}{35}c^{17}_{140}
-\frac{9}{35}c^{19}_{140}
+\frac{4}{35}c^{21}_{140}
+\frac{2}{7}c^{23}_{140}
$,
$\frac{1}{\sqrt{35}}c^{1}_{35}
-\frac{1}{\sqrt{35}}c^{6}_{35}
$,
$\frac{2}{\sqrt{35}}c^{3}_{35}
+\frac{1}{\sqrt{35}}c^{4}_{35}
+\frac{1}{\sqrt{35}}c^{10}_{35}
+\frac{1}{\sqrt{35}}c^{11}_{35}
$,
$\frac{2}{35}c^{1}_{140}
-\frac{1}{35}c^{3}_{140}
-\frac{1}{7}c^{5}_{140}
-\frac{3}{35}c^{7}_{140}
+\frac{1}{5}c^{9}_{140}
-\frac{2}{35}c^{13}_{140}
-\frac{1}{35}c^{15}_{140}
-\frac{1}{35}c^{17}_{140}
+\frac{3}{35}c^{19}_{140}
+\frac{2}{35}c^{21}_{140}
-\frac{2}{7}c^{23}_{140}
$;\ \ 
$\frac{2}{35}c^{1}_{140}
-\frac{1}{35}c^{3}_{140}
-\frac{1}{7}c^{5}_{140}
-\frac{3}{35}c^{7}_{140}
+\frac{1}{5}c^{9}_{140}
-\frac{2}{35}c^{13}_{140}
-\frac{1}{35}c^{15}_{140}
-\frac{1}{35}c^{17}_{140}
+\frac{3}{35}c^{19}_{140}
+\frac{2}{35}c^{21}_{140}
-\frac{2}{7}c^{23}_{140}
$,
$-\frac{1}{\sqrt{35}}c^{4}_{35}
+\frac{1}{\sqrt{35}}c^{11}_{35}
$,
$-\frac{4}{35}c^{1}_{140}
-\frac{3}{35}c^{3}_{140}
-\frac{1}{7}c^{5}_{140}
+\frac{1}{35}c^{7}_{140}
+\frac{1}{35}c^{9}_{140}
+\frac{4}{35}c^{13}_{140}
+\frac{2}{35}c^{15}_{140}
-\frac{3}{35}c^{17}_{140}
+\frac{9}{35}c^{19}_{140}
-\frac{4}{35}c^{21}_{140}
-\frac{2}{7}c^{23}_{140}
$,
$\frac{1}{\sqrt{35}}c^{1}_{35}
-\frac{1}{\sqrt{35}}c^{6}_{35}
$,
$-\frac{2}{\sqrt{35}}c^{3}_{35}
-\frac{1}{\sqrt{35}}c^{4}_{35}
-\frac{1}{\sqrt{35}}c^{10}_{35}
-\frac{1}{\sqrt{35}}c^{11}_{35}
$;\ \ 
$\frac{1}{\sqrt{35}}c^{1}_{35}
-\frac{1}{\sqrt{35}}c^{6}_{35}
$,
$-\frac{2}{\sqrt{35}}c^{3}_{35}
-\frac{1}{\sqrt{35}}c^{4}_{35}
-\frac{1}{\sqrt{35}}c^{10}_{35}
-\frac{1}{\sqrt{35}}c^{11}_{35}
$,
$\frac{2}{35}c^{1}_{140}
-\frac{1}{35}c^{3}_{140}
-\frac{1}{7}c^{5}_{140}
-\frac{3}{35}c^{7}_{140}
+\frac{1}{5}c^{9}_{140}
-\frac{2}{35}c^{13}_{140}
-\frac{1}{35}c^{15}_{140}
-\frac{1}{35}c^{17}_{140}
+\frac{3}{35}c^{19}_{140}
+\frac{2}{35}c^{21}_{140}
-\frac{2}{7}c^{23}_{140}
$,
$\frac{1}{\sqrt{35}\mathrm{i}}s^{3}_{140}
+\frac{1}{\sqrt{35}\mathrm{i}}s^{17}_{140}
$;\ \ 
$-\frac{2}{35}c^{1}_{140}
+\frac{1}{35}c^{3}_{140}
+\frac{1}{7}c^{5}_{140}
+\frac{3}{35}c^{7}_{140}
-\frac{1}{5}c^{9}_{140}
+\frac{2}{35}c^{13}_{140}
+\frac{1}{35}c^{15}_{140}
+\frac{1}{35}c^{17}_{140}
-\frac{3}{35}c^{19}_{140}
-\frac{2}{35}c^{21}_{140}
+\frac{2}{7}c^{23}_{140}
$,
$\frac{1}{\sqrt{35}\mathrm{i}}s^{3}_{140}
+\frac{1}{\sqrt{35}\mathrm{i}}s^{17}_{140}
$,
$\frac{1}{\sqrt{35}}c^{4}_{35}
-\frac{1}{\sqrt{35}}c^{11}_{35}
$;\ \ 
$-\frac{1}{\sqrt{35}}c^{4}_{35}
+\frac{1}{\sqrt{35}}c^{11}_{35}
$,
$-\frac{4}{35}c^{1}_{140}
-\frac{3}{35}c^{3}_{140}
-\frac{1}{7}c^{5}_{140}
+\frac{1}{35}c^{7}_{140}
+\frac{1}{35}c^{9}_{140}
+\frac{4}{35}c^{13}_{140}
+\frac{2}{35}c^{15}_{140}
-\frac{3}{35}c^{17}_{140}
+\frac{9}{35}c^{19}_{140}
-\frac{4}{35}c^{21}_{140}
-\frac{2}{7}c^{23}_{140}
$;\ \ 
$-\frac{1}{\sqrt{35}}c^{1}_{35}
+\frac{1}{\sqrt{35}}c^{6}_{35}
$)

Pass. 

 \ \color{black}

 \color{blue}

\noindent 595: (dims,levels) = $(6;35
)$,
irreps = $3_{7}^{1}
\hskip -1.5pt \otimes \hskip -1.5pt
2_{5}^{1}$,
pord$(\rho_\text{isum}(\mathfrak{t})) = 35$,

\vskip 0.7ex
\hangindent=5.5em \hangafter=1
{\white .}\hskip 1em $\rho_\text{isum}(\mathfrak{t})$ =
 $( \frac{3}{35},
\frac{12}{35},
\frac{13}{35},
\frac{17}{35},
\frac{27}{35},
\frac{33}{35} )
$,

\vskip 0.7ex
\hangindent=5.5em \hangafter=1
{\white .}\hskip 1em $\rho_\text{isum}(\mathfrak{s})$ =
$\mathrm{i}$($\frac{1}{\sqrt{35}}c^{1}_{35}
-\frac{1}{\sqrt{35}}c^{6}_{35}
$,
$\frac{4}{35}c^{1}_{140}
+\frac{3}{35}c^{3}_{140}
+\frac{1}{7}c^{5}_{140}
-\frac{1}{35}c^{7}_{140}
-\frac{1}{35}c^{9}_{140}
-\frac{4}{35}c^{13}_{140}
-\frac{2}{35}c^{15}_{140}
+\frac{3}{35}c^{17}_{140}
-\frac{9}{35}c^{19}_{140}
+\frac{4}{35}c^{21}_{140}
+\frac{2}{7}c^{23}_{140}
$,
$-\frac{1}{\sqrt{35}}c^{4}_{35}
+\frac{1}{\sqrt{35}}c^{11}_{35}
$,
$-\frac{1}{\sqrt{35}\mathrm{i}}s^{3}_{140}
-\frac{1}{\sqrt{35}\mathrm{i}}s^{17}_{140}
$,
$\frac{2}{\sqrt{35}}c^{3}_{35}
+\frac{1}{\sqrt{35}}c^{4}_{35}
+\frac{1}{\sqrt{35}}c^{10}_{35}
+\frac{1}{\sqrt{35}}c^{11}_{35}
$,
$\frac{2}{35}c^{1}_{140}
-\frac{1}{35}c^{3}_{140}
-\frac{1}{7}c^{5}_{140}
-\frac{3}{35}c^{7}_{140}
+\frac{1}{5}c^{9}_{140}
-\frac{2}{35}c^{13}_{140}
-\frac{1}{35}c^{15}_{140}
-\frac{1}{35}c^{17}_{140}
+\frac{3}{35}c^{19}_{140}
+\frac{2}{35}c^{21}_{140}
-\frac{2}{7}c^{23}_{140}
$;\ \ 
$\frac{1}{\sqrt{35}}c^{4}_{35}
-\frac{1}{\sqrt{35}}c^{11}_{35}
$,
$-\frac{1}{\sqrt{35}\mathrm{i}}s^{3}_{140}
-\frac{1}{\sqrt{35}\mathrm{i}}s^{17}_{140}
$,
$-\frac{2}{35}c^{1}_{140}
+\frac{1}{35}c^{3}_{140}
+\frac{1}{7}c^{5}_{140}
+\frac{3}{35}c^{7}_{140}
-\frac{1}{5}c^{9}_{140}
+\frac{2}{35}c^{13}_{140}
+\frac{1}{35}c^{15}_{140}
+\frac{1}{35}c^{17}_{140}
-\frac{3}{35}c^{19}_{140}
-\frac{2}{35}c^{21}_{140}
+\frac{2}{7}c^{23}_{140}
$,
$-\frac{1}{\sqrt{35}}c^{1}_{35}
+\frac{1}{\sqrt{35}}c^{6}_{35}
$,
$\frac{2}{\sqrt{35}}c^{3}_{35}
+\frac{1}{\sqrt{35}}c^{4}_{35}
+\frac{1}{\sqrt{35}}c^{10}_{35}
+\frac{1}{\sqrt{35}}c^{11}_{35}
$;\ \ 
$\frac{2}{35}c^{1}_{140}
-\frac{1}{35}c^{3}_{140}
-\frac{1}{7}c^{5}_{140}
-\frac{3}{35}c^{7}_{140}
+\frac{1}{5}c^{9}_{140}
-\frac{2}{35}c^{13}_{140}
-\frac{1}{35}c^{15}_{140}
-\frac{1}{35}c^{17}_{140}
+\frac{3}{35}c^{19}_{140}
+\frac{2}{35}c^{21}_{140}
-\frac{2}{7}c^{23}_{140}
$,
$\frac{2}{\sqrt{35}}c^{3}_{35}
+\frac{1}{\sqrt{35}}c^{4}_{35}
+\frac{1}{\sqrt{35}}c^{10}_{35}
+\frac{1}{\sqrt{35}}c^{11}_{35}
$,
$\frac{4}{35}c^{1}_{140}
+\frac{3}{35}c^{3}_{140}
+\frac{1}{7}c^{5}_{140}
-\frac{1}{35}c^{7}_{140}
-\frac{1}{35}c^{9}_{140}
-\frac{4}{35}c^{13}_{140}
-\frac{2}{35}c^{15}_{140}
+\frac{3}{35}c^{17}_{140}
-\frac{9}{35}c^{19}_{140}
+\frac{4}{35}c^{21}_{140}
+\frac{2}{7}c^{23}_{140}
$,
$\frac{1}{\sqrt{35}}c^{1}_{35}
-\frac{1}{\sqrt{35}}c^{6}_{35}
$;\ \ 
$-\frac{1}{\sqrt{35}}c^{1}_{35}
+\frac{1}{\sqrt{35}}c^{6}_{35}
$,
$\frac{1}{\sqrt{35}}c^{4}_{35}
-\frac{1}{\sqrt{35}}c^{11}_{35}
$,
$\frac{4}{35}c^{1}_{140}
+\frac{3}{35}c^{3}_{140}
+\frac{1}{7}c^{5}_{140}
-\frac{1}{35}c^{7}_{140}
-\frac{1}{35}c^{9}_{140}
-\frac{4}{35}c^{13}_{140}
-\frac{2}{35}c^{15}_{140}
+\frac{3}{35}c^{17}_{140}
-\frac{9}{35}c^{19}_{140}
+\frac{4}{35}c^{21}_{140}
+\frac{2}{7}c^{23}_{140}
$;\ \ 
$-\frac{2}{35}c^{1}_{140}
+\frac{1}{35}c^{3}_{140}
+\frac{1}{7}c^{5}_{140}
+\frac{3}{35}c^{7}_{140}
-\frac{1}{5}c^{9}_{140}
+\frac{2}{35}c^{13}_{140}
+\frac{1}{35}c^{15}_{140}
+\frac{1}{35}c^{17}_{140}
-\frac{3}{35}c^{19}_{140}
-\frac{2}{35}c^{21}_{140}
+\frac{2}{7}c^{23}_{140}
$,
$-\frac{1}{\sqrt{35}\mathrm{i}}s^{3}_{140}
-\frac{1}{\sqrt{35}\mathrm{i}}s^{17}_{140}
$;\ \ 
$-\frac{1}{\sqrt{35}}c^{4}_{35}
+\frac{1}{\sqrt{35}}c^{11}_{35}
$)

Pass. 

 \ \color{black}

 \color{blue}

\noindent 596: (dims,levels) = $(6;35
)$,
irreps = $3_{7}^{1}
\hskip -1.5pt \otimes \hskip -1.5pt
2_{5}^{2}$,
pord$(\rho_\text{isum}(\mathfrak{t})) = 35$,

\vskip 0.7ex
\hangindent=5.5em \hangafter=1
{\white .}\hskip 1em $\rho_\text{isum}(\mathfrak{t})$ =
 $( \frac{6}{35},
\frac{19}{35},
\frac{24}{35},
\frac{26}{35},
\frac{31}{35},
\frac{34}{35} )
$,

\vskip 0.7ex
\hangindent=5.5em \hangafter=1
{\white .}\hskip 1em $\rho_\text{isum}(\mathfrak{s})$ =
$\mathrm{i}$($-\frac{4}{35}c^{1}_{140}
-\frac{3}{35}c^{3}_{140}
-\frac{1}{7}c^{5}_{140}
+\frac{1}{35}c^{7}_{140}
+\frac{1}{35}c^{9}_{140}
+\frac{4}{35}c^{13}_{140}
+\frac{2}{35}c^{15}_{140}
-\frac{3}{35}c^{17}_{140}
+\frac{9}{35}c^{19}_{140}
-\frac{4}{35}c^{21}_{140}
-\frac{2}{7}c^{23}_{140}
$,
$\frac{1}{\sqrt{35}}c^{1}_{35}
-\frac{1}{\sqrt{35}}c^{6}_{35}
$,
$-\frac{1}{\sqrt{35}}c^{4}_{35}
+\frac{1}{\sqrt{35}}c^{11}_{35}
$,
$-\frac{1}{\sqrt{35}\mathrm{i}}s^{3}_{140}
-\frac{1}{\sqrt{35}\mathrm{i}}s^{17}_{140}
$,
$\frac{2}{\sqrt{35}}c^{3}_{35}
+\frac{1}{\sqrt{35}}c^{4}_{35}
+\frac{1}{\sqrt{35}}c^{10}_{35}
+\frac{1}{\sqrt{35}}c^{11}_{35}
$,
$\frac{2}{35}c^{1}_{140}
-\frac{1}{35}c^{3}_{140}
-\frac{1}{7}c^{5}_{140}
-\frac{3}{35}c^{7}_{140}
+\frac{1}{5}c^{9}_{140}
-\frac{2}{35}c^{13}_{140}
-\frac{1}{35}c^{15}_{140}
-\frac{1}{35}c^{17}_{140}
+\frac{3}{35}c^{19}_{140}
+\frac{2}{35}c^{21}_{140}
-\frac{2}{7}c^{23}_{140}
$;\ \ 
$\frac{2}{\sqrt{35}}c^{3}_{35}
+\frac{1}{\sqrt{35}}c^{4}_{35}
+\frac{1}{\sqrt{35}}c^{10}_{35}
+\frac{1}{\sqrt{35}}c^{11}_{35}
$,
$\frac{4}{35}c^{1}_{140}
+\frac{3}{35}c^{3}_{140}
+\frac{1}{7}c^{5}_{140}
-\frac{1}{35}c^{7}_{140}
-\frac{1}{35}c^{9}_{140}
-\frac{4}{35}c^{13}_{140}
-\frac{2}{35}c^{15}_{140}
+\frac{3}{35}c^{17}_{140}
-\frac{9}{35}c^{19}_{140}
+\frac{4}{35}c^{21}_{140}
+\frac{2}{7}c^{23}_{140}
$,
$\frac{1}{\sqrt{35}}c^{4}_{35}
-\frac{1}{\sqrt{35}}c^{11}_{35}
$,
$-\frac{2}{35}c^{1}_{140}
+\frac{1}{35}c^{3}_{140}
+\frac{1}{7}c^{5}_{140}
+\frac{3}{35}c^{7}_{140}
-\frac{1}{5}c^{9}_{140}
+\frac{2}{35}c^{13}_{140}
+\frac{1}{35}c^{15}_{140}
+\frac{1}{35}c^{17}_{140}
-\frac{3}{35}c^{19}_{140}
-\frac{2}{35}c^{21}_{140}
+\frac{2}{7}c^{23}_{140}
$,
$-\frac{1}{\sqrt{35}\mathrm{i}}s^{3}_{140}
-\frac{1}{\sqrt{35}\mathrm{i}}s^{17}_{140}
$;\ \ 
$-\frac{1}{\sqrt{35}\mathrm{i}}s^{3}_{140}
-\frac{1}{\sqrt{35}\mathrm{i}}s^{17}_{140}
$,
$-\frac{2}{35}c^{1}_{140}
+\frac{1}{35}c^{3}_{140}
+\frac{1}{7}c^{5}_{140}
+\frac{3}{35}c^{7}_{140}
-\frac{1}{5}c^{9}_{140}
+\frac{2}{35}c^{13}_{140}
+\frac{1}{35}c^{15}_{140}
+\frac{1}{35}c^{17}_{140}
-\frac{3}{35}c^{19}_{140}
-\frac{2}{35}c^{21}_{140}
+\frac{2}{7}c^{23}_{140}
$,
$-\frac{1}{\sqrt{35}}c^{1}_{35}
+\frac{1}{\sqrt{35}}c^{6}_{35}
$,
$\frac{2}{\sqrt{35}}c^{3}_{35}
+\frac{1}{\sqrt{35}}c^{4}_{35}
+\frac{1}{\sqrt{35}}c^{10}_{35}
+\frac{1}{\sqrt{35}}c^{11}_{35}
$;\ \ 
$-\frac{2}{\sqrt{35}}c^{3}_{35}
-\frac{1}{\sqrt{35}}c^{4}_{35}
-\frac{1}{\sqrt{35}}c^{10}_{35}
-\frac{1}{\sqrt{35}}c^{11}_{35}
$,
$-\frac{4}{35}c^{1}_{140}
-\frac{3}{35}c^{3}_{140}
-\frac{1}{7}c^{5}_{140}
+\frac{1}{35}c^{7}_{140}
+\frac{1}{35}c^{9}_{140}
+\frac{4}{35}c^{13}_{140}
+\frac{2}{35}c^{15}_{140}
-\frac{3}{35}c^{17}_{140}
+\frac{9}{35}c^{19}_{140}
-\frac{4}{35}c^{21}_{140}
-\frac{2}{7}c^{23}_{140}
$,
$-\frac{1}{\sqrt{35}}c^{1}_{35}
+\frac{1}{\sqrt{35}}c^{6}_{35}
$;\ \ 
$\frac{1}{\sqrt{35}\mathrm{i}}s^{3}_{140}
+\frac{1}{\sqrt{35}\mathrm{i}}s^{17}_{140}
$,
$\frac{1}{\sqrt{35}}c^{4}_{35}
-\frac{1}{\sqrt{35}}c^{11}_{35}
$;\ \ 
$\frac{4}{35}c^{1}_{140}
+\frac{3}{35}c^{3}_{140}
+\frac{1}{7}c^{5}_{140}
-\frac{1}{35}c^{7}_{140}
-\frac{1}{35}c^{9}_{140}
-\frac{4}{35}c^{13}_{140}
-\frac{2}{35}c^{15}_{140}
+\frac{3}{35}c^{17}_{140}
-\frac{9}{35}c^{19}_{140}
+\frac{4}{35}c^{21}_{140}
+\frac{2}{7}c^{23}_{140}
$)

Pass. 

 \ \color{black}

\noindent 597: (dims,levels) = $(6;36
)$,
irreps = $6_{9,1}^{5,0}
\hskip -1.5pt \otimes \hskip -1.5pt
1_{4}^{1,0}$,
pord$(\rho_\text{isum}(\mathfrak{t})) = 9$,

\vskip 0.7ex
\hangindent=5.5em \hangafter=1
{\white .}\hskip 1em $\rho_\text{isum}(\mathfrak{t})$ =
 $( \frac{1}{36},
\frac{5}{36},
\frac{13}{36},
\frac{17}{36},
\frac{25}{36},
\frac{29}{36} )
$,

\vskip 0.7ex
\hangindent=5.5em \hangafter=1
{\white .}\hskip 1em $\rho_\text{isum}(\mathfrak{s})$ =
($-\frac{1}{3}$,
$\frac{1}{3}c^{11}_{72}
$,
$\frac{1}{3}$,
$\frac{1}{3}c^{1}_{72}
$,
$\frac{1}{3}$,
$-\frac{1}{3}c^{1}_{72}
+\frac{1}{3}c^{11}_{72}
$;
$\frac{1}{3}$,
$-\frac{1}{3}c^{1}_{72}
+\frac{1}{3}c^{11}_{72}
$,
$-\frac{1}{3}$,
$\frac{1}{3}c^{1}_{72}
$,
$-\frac{1}{3}$;
$-\frac{1}{3}$,
$\frac{1}{3}c^{11}_{72}
$,
$-\frac{1}{3}$,
$-\frac{1}{3}c^{1}_{72}
$;
$\frac{1}{3}$,
$\frac{1}{3}c^{1}_{72}
-\frac{1}{3}c^{11}_{72}
$,
$\frac{1}{3}$;
$-\frac{1}{3}$,
$\frac{1}{3}c^{11}_{72}
$;
$\frac{1}{3}$)

Fail:
cnd($\rho(\mathfrak s)_\mathrm{ndeg}$) = 72 does not divide
 ord($\rho(\mathfrak t)$)=36. Prop. B.4 (2)

 \ \color{black}

\noindent 598: (dims,levels) = $(6;36
)$,
irreps = $6_{9,1}^{1,0}
\hskip -1.5pt \otimes \hskip -1.5pt
1_{4}^{1,0}$,
pord$(\rho_\text{isum}(\mathfrak{t})) = 9$,

\vskip 0.7ex
\hangindent=5.5em \hangafter=1
{\white .}\hskip 1em $\rho_\text{isum}(\mathfrak{t})$ =
 $( \frac{1}{36},
\frac{5}{36},
\frac{13}{36},
\frac{17}{36},
\frac{25}{36},
\frac{29}{36} )
$,

\vskip 0.7ex
\hangindent=5.5em \hangafter=1
{\white .}\hskip 1em $\rho_\text{isum}(\mathfrak{s})$ =
($\frac{1}{3}$,
$-\frac{1}{3}c^{7}_{72}
$,
$\frac{1}{3}$,
$-\frac{1}{3}c^{5}_{72}
+\frac{1}{3}c^{7}_{72}
$,
$\frac{1}{3}$,
$-\frac{1}{3}c^{5}_{72}
$;
$-\frac{1}{3}$,
$\frac{1}{3}c^{5}_{72}
$,
$-\frac{1}{3}$,
$-\frac{1}{3}c^{5}_{72}
+\frac{1}{3}c^{7}_{72}
$,
$\frac{1}{3}$;
$\frac{1}{3}$,
$-\frac{1}{3}c^{7}_{72}
$,
$\frac{1}{3}$,
$\frac{1}{3}c^{5}_{72}
-\frac{1}{3}c^{7}_{72}
$;
$-\frac{1}{3}$,
$\frac{1}{3}c^{5}_{72}
$,
$\frac{1}{3}$;
$\frac{1}{3}$,
$\frac{1}{3}c^{7}_{72}
$;
$-\frac{1}{3}$)

Fail:
cnd($\rho(\mathfrak s)_\mathrm{ndeg}$) = 72 does not divide
 ord($\rho(\mathfrak t)$)=36. Prop. B.4 (2)

 \ \color{black}

\noindent 599: (dims,levels) = $(6;36
)$,
irreps = $6_{9,2}^{1,0}
\hskip -1.5pt \otimes \hskip -1.5pt
1_{4}^{1,0}$,
pord$(\rho_\text{isum}(\mathfrak{t})) = 9$,

\vskip 0.7ex
\hangindent=5.5em \hangafter=1
{\white .}\hskip 1em $\rho_\text{isum}(\mathfrak{t})$ =
 $( \frac{1}{36},
\frac{5}{36},
\frac{13}{36},
\frac{17}{36},
\frac{25}{36},
\frac{29}{36} )
$,

\vskip 0.7ex
\hangindent=5.5em \hangafter=1
{\white .}\hskip 1em $\rho_\text{isum}(\mathfrak{s})$ =
$\mathrm{i}$($-\frac{1}{3}$,
$\frac{1}{3}c^{1}_{36}
$,
$\frac{1}{3}$,
$\frac{1}{3}c^{5}_{36}
$,
$\frac{1}{3}$,
$\frac{1}{3}c^{1}_{36}
-\frac{1}{3}c^{5}_{36}
$;\ \ 
$-\frac{1}{3}$,
$\frac{1}{3}c^{1}_{36}
-\frac{1}{3}c^{5}_{36}
$,
$\frac{1}{3}$,
$\frac{1}{3}c^{5}_{36}
$,
$\frac{1}{3}$;\ \ 
$-\frac{1}{3}$,
$\frac{1}{3}c^{1}_{36}
$,
$-\frac{1}{3}$,
$-\frac{1}{3}c^{5}_{36}
$;\ \ 
$-\frac{1}{3}$,
$-\frac{1}{3}c^{1}_{36}
+\frac{1}{3}c^{5}_{36}
$,
$-\frac{1}{3}$;\ \ 
$-\frac{1}{3}$,
$\frac{1}{3}c^{1}_{36}
$;\ \ 
$-\frac{1}{3}$)

Fail:
$\sigma(\rho(\mathfrak s)_\mathrm{ndeg}) \neq
 (\rho(\mathfrak t)^a \rho(\mathfrak s) \rho(\mathfrak t)^b
 \rho(\mathfrak s) \rho(\mathfrak t)^a)_\mathrm{ndeg}$,
 $\sigma = a$ = 5. Prop. B.5 (3) eqn. (B.25)

 \ \color{black}

 \color{blue}

\noindent 600: (dims,levels) = $(6;36
)$,
irreps = $6_{9,3}^{1,0}
\hskip -1.5pt \otimes \hskip -1.5pt
1_{4}^{1,0}$,
pord$(\rho_\text{isum}(\mathfrak{t})) = 9$,

\vskip 0.7ex
\hangindent=5.5em \hangafter=1
{\white .}\hskip 1em $\rho_\text{isum}(\mathfrak{t})$ =
 $( \frac{1}{36},
\frac{5}{36},
\frac{13}{36},
\frac{17}{36},
\frac{25}{36},
\frac{29}{36} )
$,

\vskip 0.7ex
\hangindent=5.5em \hangafter=1
{\white .}\hskip 1em $\rho_\text{isum}(\mathfrak{s})$ =
$\mathrm{i}$($\frac{1}{3}$,
$\frac{1}{3}c^{2}_{9}
$,
$\frac{1}{3}$,
$-\frac{1}{3}c^{1}_{9}
$,
$\frac{1}{3}$,
$\frac{1}{3} c_9^4 $;\ \ 
$\frac{1}{3}$,
$\frac{1}{3} c_9^4 $,
$-\frac{1}{3}$,
$\frac{1}{3}c^{1}_{9}
$,
$\frac{1}{3}$;\ \ 
$\frac{1}{3}$,
$-\frac{1}{3}c^{2}_{9}
$,
$\frac{1}{3}$,
$\frac{1}{3}c^{1}_{9}
$;\ \ 
$\frac{1}{3}$,
$-\frac{1}{3} c_9^4 $,
$-\frac{1}{3}$;\ \ 
$\frac{1}{3}$,
$\frac{1}{3}c^{2}_{9}
$;\ \ 
$\frac{1}{3}$)

Pass. 

 \ \color{black}

 \color{blue}

\noindent 601: (dims,levels) = $(6;39
)$,
irreps = $6_{13}^{1}
\hskip -1.5pt \otimes \hskip -1.5pt
1_{3}^{1,0}$,
pord$(\rho_\text{isum}(\mathfrak{t})) = 13$,

\vskip 0.7ex
\hangindent=5.5em \hangafter=1
{\white .}\hskip 1em $\rho_\text{isum}(\mathfrak{t})$ =
 $( \frac{1}{39},
\frac{4}{39},
\frac{10}{39},
\frac{16}{39},
\frac{22}{39},
\frac{25}{39} )
$,

\vskip 0.7ex
\hangindent=5.5em \hangafter=1
{\white .}\hskip 1em $\rho_\text{isum}(\mathfrak{s})$ =
$\mathrm{i}$($-\frac{1}{\sqrt{13}}c^{7}_{52}
$,
$-\frac{1}{\sqrt{13}}c^{1}_{52}
$,
$\frac{1}{\sqrt{13}}c^{3}_{52}
$,
$\frac{1}{\sqrt{13}}c^{11}_{52}
$,
$-\frac{1}{\sqrt{13}}c^{5}_{52}
$,
$\frac{1}{\sqrt{13}}c^{9}_{52}
$;\ \ 
$\frac{1}{\sqrt{13}}c^{11}_{52}
$,
$\frac{1}{\sqrt{13}}c^{7}_{52}
$,
$\frac{1}{\sqrt{13}}c^{9}_{52}
$,
$\frac{1}{\sqrt{13}}c^{3}_{52}
$,
$-\frac{1}{\sqrt{13}}c^{5}_{52}
$;\ \ 
$\frac{1}{\sqrt{13}}c^{5}_{52}
$,
$\frac{1}{\sqrt{13}}c^{1}_{52}
$,
$-\frac{1}{\sqrt{13}}c^{9}_{52}
$,
$-\frac{1}{\sqrt{13}}c^{11}_{52}
$;\ \ 
$-\frac{1}{\sqrt{13}}c^{5}_{52}
$,
$\frac{1}{\sqrt{13}}c^{7}_{52}
$,
$\frac{1}{\sqrt{13}}c^{3}_{52}
$;\ \ 
$-\frac{1}{\sqrt{13}}c^{11}_{52}
$,
$\frac{1}{\sqrt{13}}c^{1}_{52}
$;\ \ 
$\frac{1}{\sqrt{13}}c^{7}_{52}
$)

Pass. 

 \ \color{black}

 \color{blue}

\noindent 602: (dims,levels) = $(6;39
)$,
irreps = $6_{13}^{2}
\hskip -1.5pt \otimes \hskip -1.5pt
1_{3}^{1,0}$,
pord$(\rho_\text{isum}(\mathfrak{t})) = 13$,

\vskip 0.7ex
\hangindent=5.5em \hangafter=1
{\white .}\hskip 1em $\rho_\text{isum}(\mathfrak{t})$ =
 $( \frac{7}{39},
\frac{19}{39},
\frac{28}{39},
\frac{31}{39},
\frac{34}{39},
\frac{37}{39} )
$,

\vskip 0.7ex
\hangindent=5.5em \hangafter=1
{\white .}\hskip 1em $\rho_\text{isum}(\mathfrak{s})$ =
$\mathrm{i}$($-\frac{1}{\sqrt{13}}c^{3}_{52}
$,
$\frac{1}{\sqrt{13}}c^{11}_{52}
$,
$\frac{1}{\sqrt{13}}c^{7}_{52}
$,
$-\frac{1}{\sqrt{13}}c^{5}_{52}
$,
$-\frac{1}{\sqrt{13}}c^{1}_{52}
$,
$\frac{1}{\sqrt{13}}c^{9}_{52}
$;\ \ 
$\frac{1}{\sqrt{13}}c^{3}_{52}
$,
$-\frac{1}{\sqrt{13}}c^{9}_{52}
$,
$\frac{1}{\sqrt{13}}c^{1}_{52}
$,
$-\frac{1}{\sqrt{13}}c^{5}_{52}
$,
$\frac{1}{\sqrt{13}}c^{7}_{52}
$;\ \ 
$-\frac{1}{\sqrt{13}}c^{1}_{52}
$,
$-\frac{1}{\sqrt{13}}c^{3}_{52}
$,
$-\frac{1}{\sqrt{13}}c^{11}_{52}
$,
$\frac{1}{\sqrt{13}}c^{5}_{52}
$;\ \ 
$-\frac{1}{\sqrt{13}}c^{9}_{52}
$,
$\frac{1}{\sqrt{13}}c^{7}_{52}
$,
$-\frac{1}{\sqrt{13}}c^{11}_{52}
$;\ \ 
$\frac{1}{\sqrt{13}}c^{9}_{52}
$,
$\frac{1}{\sqrt{13}}c^{3}_{52}
$;\ \ 
$\frac{1}{\sqrt{13}}c^{1}_{52}
$)

Pass. 

 \ \color{black}

\noindent 603: (dims,levels) = $(6;40
)$,
irreps = $3_{5}^{1}
\hskip -1.5pt \otimes \hskip -1.5pt
2_{8}^{1,0}$,
pord$(\rho_\text{isum}(\mathfrak{t})) = 20$,

\vskip 0.7ex
\hangindent=5.5em \hangafter=1
{\white .}\hskip 1em $\rho_\text{isum}(\mathfrak{t})$ =
 $( \frac{1}{8},
\frac{3}{8},
\frac{7}{40},
\frac{13}{40},
\frac{23}{40},
\frac{37}{40} )
$,

\vskip 0.7ex
\hangindent=5.5em \hangafter=1
{\white .}\hskip 1em $\rho_\text{isum}(\mathfrak{s})$ =
($-\sqrt{\frac{1}{10}}$,
$-\sqrt{\frac{1}{10}}$,
$\sqrt{\frac{1}{5}}$,
$\sqrt{\frac{1}{5}}$,
$\sqrt{\frac{1}{5}}$,
$\sqrt{\frac{1}{5}}$;
$\sqrt{\frac{1}{10}}$,
$-\sqrt{\frac{1}{5}}$,
$\sqrt{\frac{1}{5}}$,
$-\sqrt{\frac{1}{5}}$,
$\sqrt{\frac{1}{5}}$;
$-\frac{1}{\sqrt{10}\mathrm{i}}s^{3}_{20}
$,
$-\frac{1}{\sqrt{10}}c^{1}_{5}
$,
$\frac{1}{\sqrt{10}}c^{1}_{5}
$,
$\frac{1}{\sqrt{10}\mathrm{i}}s^{3}_{20}
$;
$\frac{1}{\sqrt{10}\mathrm{i}}s^{3}_{20}
$,
$\frac{1}{\sqrt{10}\mathrm{i}}s^{3}_{20}
$,
$-\frac{1}{\sqrt{10}}c^{1}_{5}
$;
$-\frac{1}{\sqrt{10}\mathrm{i}}s^{3}_{20}
$,
$-\frac{1}{\sqrt{10}}c^{1}_{5}
$;
$\frac{1}{\sqrt{10}\mathrm{i}}s^{3}_{20}
$)

Fail:
$\sigma(\rho(\mathfrak s)_\mathrm{ndeg}) \neq
 (\rho(\mathfrak t)^a \rho(\mathfrak s) \rho(\mathfrak t)^b
 \rho(\mathfrak s) \rho(\mathfrak t)^a)_\mathrm{ndeg}$,
 $\sigma = a$ = 3. Prop. B.5 (3) eqn. (B.25)

 \ \color{black}

\noindent 604: (dims,levels) = $(6;40
)$,
irreps = $3_{5}^{3}
\hskip -1.5pt \otimes \hskip -1.5pt
2_{8}^{1,0}$,
pord$(\rho_\text{isum}(\mathfrak{t})) = 20$,

\vskip 0.7ex
\hangindent=5.5em \hangafter=1
{\white .}\hskip 1em $\rho_\text{isum}(\mathfrak{t})$ =
 $( \frac{1}{8},
\frac{3}{8},
\frac{21}{40},
\frac{29}{40},
\frac{31}{40},
\frac{39}{40} )
$,

\vskip 0.7ex
\hangindent=5.5em \hangafter=1
{\white .}\hskip 1em $\rho_\text{isum}(\mathfrak{s})$ =
($\sqrt{\frac{1}{10}}$,
$-\sqrt{\frac{1}{10}}$,
$\sqrt{\frac{1}{5}}$,
$\sqrt{\frac{1}{5}}$,
$\sqrt{\frac{1}{5}}$,
$\sqrt{\frac{1}{5}}$;
$-\sqrt{\frac{1}{10}}$,
$-\sqrt{\frac{1}{5}}$,
$-\sqrt{\frac{1}{5}}$,
$\sqrt{\frac{1}{5}}$,
$\sqrt{\frac{1}{5}}$;
$\frac{1}{\sqrt{10}}c^{1}_{5}
$,
$-\frac{1}{\sqrt{10}\mathrm{i}}s^{3}_{20}
$,
$\frac{1}{\sqrt{10}}c^{1}_{5}
$,
$-\frac{1}{\sqrt{10}\mathrm{i}}s^{3}_{20}
$;
$\frac{1}{\sqrt{10}}c^{1}_{5}
$,
$-\frac{1}{\sqrt{10}\mathrm{i}}s^{3}_{20}
$,
$\frac{1}{\sqrt{10}}c^{1}_{5}
$;
$-\frac{1}{\sqrt{10}}c^{1}_{5}
$,
$\frac{1}{\sqrt{10}\mathrm{i}}s^{3}_{20}
$;
$-\frac{1}{\sqrt{10}}c^{1}_{5}
$)

Fail:
$\sigma(\rho(\mathfrak s)_\mathrm{ndeg}) \neq
 (\rho(\mathfrak t)^a \rho(\mathfrak s) \rho(\mathfrak t)^b
 \rho(\mathfrak s) \rho(\mathfrak t)^a)_\mathrm{ndeg}$,
 $\sigma = a$ = 3. Prop. B.5 (3) eqn. (B.25)

 \ \color{black}

\noindent 605: (dims,levels) = $(6;40
)$,
irreps = $3_{8}^{3,0}
\hskip -1.5pt \otimes \hskip -1.5pt
2_{5}^{1}$,
pord$(\rho_\text{isum}(\mathfrak{t})) = 40$,

\vskip 0.7ex
\hangindent=5.5em \hangafter=1
{\white .}\hskip 1em $\rho_\text{isum}(\mathfrak{t})$ =
 $( \frac{1}{5},
\frac{4}{5},
\frac{3}{40},
\frac{7}{40},
\frac{23}{40},
\frac{27}{40} )
$,

\vskip 0.7ex
\hangindent=5.5em \hangafter=1
{\white .}\hskip 1em $\rho_\text{isum}(\mathfrak{s})$ =
($0$,
$0$,
$\frac{1}{\sqrt{10}}c^{3}_{20}
$,
$\frac{1}{\sqrt{10}}c^{1}_{20}
$,
$\frac{1}{\sqrt{10}}c^{3}_{20}
$,
$\frac{1}{\sqrt{10}}c^{1}_{20}
$;
$0$,
$\frac{1}{\sqrt{10}}c^{1}_{20}
$,
$-\frac{1}{\sqrt{10}}c^{3}_{20}
$,
$\frac{1}{\sqrt{10}}c^{1}_{20}
$,
$-\frac{1}{\sqrt{10}}c^{3}_{20}
$;
$\frac{1}{2\sqrt{5}}c^{3}_{20}
$,
$-\frac{1}{2\sqrt{5}}c^{1}_{20}
$,
$-\frac{1}{2\sqrt{5}}c^{3}_{20}
$,
$\frac{1}{2\sqrt{5}}c^{1}_{20}
$;
$-\frac{1}{2\sqrt{5}}c^{3}_{20}
$,
$\frac{1}{2\sqrt{5}}c^{1}_{20}
$,
$\frac{1}{2\sqrt{5}}c^{3}_{20}
$;
$\frac{1}{2\sqrt{5}}c^{3}_{20}
$,
$-\frac{1}{2\sqrt{5}}c^{1}_{20}
$;
$-\frac{1}{2\sqrt{5}}c^{3}_{20}
$)

Fail:
$\sigma(\rho(\mathfrak s)_\mathrm{ndeg}) \neq
 (\rho(\mathfrak t)^a \rho(\mathfrak s) \rho(\mathfrak t)^b
 \rho(\mathfrak s) \rho(\mathfrak t)^a)_\mathrm{ndeg}$,
 $\sigma = a$ = 3. Prop. B.5 (3) eqn. (B.25)

 \ \color{black}

\noindent 606: (dims,levels) = $(6;40
)$,
irreps = $3_{8}^{1,0}
\hskip -1.5pt \otimes \hskip -1.5pt
2_{5}^{1}$,
pord$(\rho_\text{isum}(\mathfrak{t})) = 40$,

\vskip 0.7ex
\hangindent=5.5em \hangafter=1
{\white .}\hskip 1em $\rho_\text{isum}(\mathfrak{t})$ =
 $( \frac{1}{5},
\frac{4}{5},
\frac{13}{40},
\frac{17}{40},
\frac{33}{40},
\frac{37}{40} )
$,

\vskip 0.7ex
\hangindent=5.5em \hangafter=1
{\white .}\hskip 1em $\rho_\text{isum}(\mathfrak{s})$ =
($0$,
$0$,
$\frac{1}{\sqrt{10}}c^{3}_{20}
$,
$\frac{1}{\sqrt{10}}c^{1}_{20}
$,
$\frac{1}{\sqrt{10}}c^{3}_{20}
$,
$\frac{1}{\sqrt{10}}c^{1}_{20}
$;
$0$,
$\frac{1}{\sqrt{10}}c^{1}_{20}
$,
$-\frac{1}{\sqrt{10}}c^{3}_{20}
$,
$\frac{1}{\sqrt{10}}c^{1}_{20}
$,
$-\frac{1}{\sqrt{10}}c^{3}_{20}
$;
$-\frac{1}{2\sqrt{5}}c^{3}_{20}
$,
$\frac{1}{2\sqrt{5}}c^{1}_{20}
$,
$\frac{1}{2\sqrt{5}}c^{3}_{20}
$,
$-\frac{1}{2\sqrt{5}}c^{1}_{20}
$;
$\frac{1}{2\sqrt{5}}c^{3}_{20}
$,
$-\frac{1}{2\sqrt{5}}c^{1}_{20}
$,
$-\frac{1}{2\sqrt{5}}c^{3}_{20}
$;
$-\frac{1}{2\sqrt{5}}c^{3}_{20}
$,
$\frac{1}{2\sqrt{5}}c^{1}_{20}
$;
$\frac{1}{2\sqrt{5}}c^{3}_{20}
$)

Fail:
$\sigma(\rho(\mathfrak s)_\mathrm{ndeg}) \neq
 (\rho(\mathfrak t)^a \rho(\mathfrak s) \rho(\mathfrak t)^b
 \rho(\mathfrak s) \rho(\mathfrak t)^a)_\mathrm{ndeg}$,
 $\sigma = a$ = 3. Prop. B.5 (3) eqn. (B.25)

 \ \color{black}

\noindent 607: (dims,levels) = $(6;40
)$,
irreps = $3_{8}^{1,0}
\hskip -1.5pt \otimes \hskip -1.5pt
2_{5}^{2}$,
pord$(\rho_\text{isum}(\mathfrak{t})) = 40$,

\vskip 0.7ex
\hangindent=5.5em \hangafter=1
{\white .}\hskip 1em $\rho_\text{isum}(\mathfrak{t})$ =
 $( \frac{2}{5},
\frac{3}{5},
\frac{1}{40},
\frac{9}{40},
\frac{21}{40},
\frac{29}{40} )
$,

\vskip 0.7ex
\hangindent=5.5em \hangafter=1
{\white .}\hskip 1em $\rho_\text{isum}(\mathfrak{s})$ =
($0$,
$0$,
$\frac{1}{\sqrt{10}}c^{1}_{20}
$,
$\frac{1}{\sqrt{10}}c^{3}_{20}
$,
$\frac{1}{\sqrt{10}}c^{1}_{20}
$,
$\frac{1}{\sqrt{10}}c^{3}_{20}
$;
$0$,
$\frac{1}{\sqrt{10}}c^{3}_{20}
$,
$-\frac{1}{\sqrt{10}}c^{1}_{20}
$,
$\frac{1}{\sqrt{10}}c^{3}_{20}
$,
$-\frac{1}{\sqrt{10}}c^{1}_{20}
$;
$-\frac{1}{2\sqrt{5}}c^{1}_{20}
$,
$-\frac{1}{2\sqrt{5}}c^{3}_{20}
$,
$\frac{1}{2\sqrt{5}}c^{1}_{20}
$,
$\frac{1}{2\sqrt{5}}c^{3}_{20}
$;
$\frac{1}{2\sqrt{5}}c^{1}_{20}
$,
$\frac{1}{2\sqrt{5}}c^{3}_{20}
$,
$-\frac{1}{2\sqrt{5}}c^{1}_{20}
$;
$-\frac{1}{2\sqrt{5}}c^{1}_{20}
$,
$-\frac{1}{2\sqrt{5}}c^{3}_{20}
$;
$\frac{1}{2\sqrt{5}}c^{1}_{20}
$)

Fail:
$\sigma(\rho(\mathfrak s)_\mathrm{ndeg}) \neq
 (\rho(\mathfrak t)^a \rho(\mathfrak s) \rho(\mathfrak t)^b
 \rho(\mathfrak s) \rho(\mathfrak t)^a)_\mathrm{ndeg}$,
 $\sigma = a$ = 3. Prop. B.5 (3) eqn. (B.25)

 \ \color{black}

\noindent 608: (dims,levels) = $(6;40
)$,
irreps = $3_{8}^{3,0}
\hskip -1.5pt \otimes \hskip -1.5pt
2_{5}^{2}$,
pord$(\rho_\text{isum}(\mathfrak{t})) = 40$,

\vskip 0.7ex
\hangindent=5.5em \hangafter=1
{\white .}\hskip 1em $\rho_\text{isum}(\mathfrak{t})$ =
 $( \frac{2}{5},
\frac{3}{5},
\frac{11}{40},
\frac{19}{40},
\frac{31}{40},
\frac{39}{40} )
$,

\vskip 0.7ex
\hangindent=5.5em \hangafter=1
{\white .}\hskip 1em $\rho_\text{isum}(\mathfrak{s})$ =
($0$,
$0$,
$\frac{1}{\sqrt{10}}c^{1}_{20}
$,
$\frac{1}{\sqrt{10}}c^{3}_{20}
$,
$\frac{1}{\sqrt{10}}c^{1}_{20}
$,
$\frac{1}{\sqrt{10}}c^{3}_{20}
$;
$0$,
$\frac{1}{\sqrt{10}}c^{3}_{20}
$,
$-\frac{1}{\sqrt{10}}c^{1}_{20}
$,
$\frac{1}{\sqrt{10}}c^{3}_{20}
$,
$-\frac{1}{\sqrt{10}}c^{1}_{20}
$;
$\frac{1}{2\sqrt{5}}c^{1}_{20}
$,
$\frac{1}{2\sqrt{5}}c^{3}_{20}
$,
$-\frac{1}{2\sqrt{5}}c^{1}_{20}
$,
$-\frac{1}{2\sqrt{5}}c^{3}_{20}
$;
$-\frac{1}{2\sqrt{5}}c^{1}_{20}
$,
$-\frac{1}{2\sqrt{5}}c^{3}_{20}
$,
$\frac{1}{2\sqrt{5}}c^{1}_{20}
$;
$\frac{1}{2\sqrt{5}}c^{1}_{20}
$,
$\frac{1}{2\sqrt{5}}c^{3}_{20}
$;
$-\frac{1}{2\sqrt{5}}c^{1}_{20}
$)

Fail:
$\sigma(\rho(\mathfrak s)_\mathrm{ndeg}) \neq
 (\rho(\mathfrak t)^a \rho(\mathfrak s) \rho(\mathfrak t)^b
 \rho(\mathfrak s) \rho(\mathfrak t)^a)_\mathrm{ndeg}$,
 $\sigma = a$ = 3. Prop. B.5 (3) eqn. (B.25)

 \ \color{black}

\noindent 609: (dims,levels) = $(6;42
)$,
irreps = $3_{7}^{3}
\hskip -1.5pt \otimes \hskip -1.5pt
2_{2}^{1,0}
\hskip -1.5pt \otimes \hskip -1.5pt
1_{3}^{1,0}$,
pord$(\rho_\text{isum}(\mathfrak{t})) = 14$,

\vskip 0.7ex
\hangindent=5.5em \hangafter=1
{\white .}\hskip 1em $\rho_\text{isum}(\mathfrak{t})$ =
 $( \frac{1}{21},
\frac{4}{21},
\frac{16}{21},
\frac{11}{42},
\frac{23}{42},
\frac{29}{42} )
$,

\vskip 0.7ex
\hangindent=5.5em \hangafter=1
{\white .}\hskip 1em $\rho_\text{isum}(\mathfrak{s})$ =
($-\frac{1}{2\sqrt{7}}c^{5}_{28}
$,
$-\frac{1}{2\sqrt{7}}c^{3}_{28}
$,
$-\frac{1}{2\sqrt{7}}c^{1}_{28}
$,
$\frac{3}{2\sqrt{21}}c^{1}_{28}
$,
$-\frac{3}{2\sqrt{21}}c^{5}_{28}
$,
$\frac{3}{2\sqrt{21}}c^{3}_{28}
$;
$\frac{1}{2\sqrt{7}}c^{1}_{28}
$,
$-\frac{1}{2\sqrt{7}}c^{5}_{28}
$,
$\frac{3}{2\sqrt{21}}c^{5}_{28}
$,
$-\frac{3}{2\sqrt{21}}c^{3}_{28}
$,
$-\frac{3}{2\sqrt{21}}c^{1}_{28}
$;
$\frac{1}{2\sqrt{7}}c^{3}_{28}
$,
$-\frac{3}{2\sqrt{21}}c^{3}_{28}
$,
$-\frac{3}{2\sqrt{21}}c^{1}_{28}
$,
$\frac{3}{2\sqrt{21}}c^{5}_{28}
$;
$-\frac{1}{2\sqrt{7}}c^{3}_{28}
$,
$-\frac{1}{2\sqrt{7}}c^{1}_{28}
$,
$\frac{1}{2\sqrt{7}}c^{5}_{28}
$;
$\frac{1}{2\sqrt{7}}c^{5}_{28}
$,
$-\frac{1}{2\sqrt{7}}c^{3}_{28}
$;
$-\frac{1}{2\sqrt{7}}c^{1}_{28}
$)

Fail:
cnd($\rho(\mathfrak s)_\mathrm{ndeg}$) = 84 does not divide
 ord($\rho(\mathfrak t)$)=42. Prop. B.4 (2)

 \ \color{black}

\noindent 610: (dims,levels) = $(6;42
)$,
irreps = $3_{7}^{1}
\hskip -1.5pt \otimes \hskip -1.5pt
2_{3}^{1,0}
\hskip -1.5pt \otimes \hskip -1.5pt
1_{2}^{1,0}$,
pord$(\rho_\text{isum}(\mathfrak{t})) = 21$,

\vskip 0.7ex
\hangindent=5.5em \hangafter=1
{\white .}\hskip 1em $\rho_\text{isum}(\mathfrak{t})$ =
 $( \frac{1}{14},
\frac{9}{14},
\frac{11}{14},
\frac{5}{42},
\frac{17}{42},
\frac{41}{42} )
$,

\vskip 0.7ex
\hangindent=5.5em \hangafter=1
{\white .}\hskip 1em $\rho_\text{isum}(\mathfrak{s})$ =
$\mathrm{i}$($-\frac{1}{\sqrt{21}}c^{3}_{28}
$,
$\frac{1}{\sqrt{21}}c^{5}_{28}
$,
$-\frac{1}{\sqrt{21}}c^{1}_{28}
$,
$-\frac{2}{\sqrt{42}}c^{1}_{28}
$,
$-\frac{2}{\sqrt{42}}c^{3}_{28}
$,
$\frac{2}{\sqrt{42}}c^{5}_{28}
$;\ \ 
$-\frac{1}{\sqrt{21}}c^{1}_{28}
$,
$-\frac{1}{\sqrt{21}}c^{3}_{28}
$,
$-\frac{2}{\sqrt{42}}c^{3}_{28}
$,
$\frac{2}{\sqrt{42}}c^{5}_{28}
$,
$-\frac{2}{\sqrt{42}}c^{1}_{28}
$;\ \ 
$\frac{1}{\sqrt{21}}c^{5}_{28}
$,
$\frac{2}{\sqrt{42}}c^{5}_{28}
$,
$-\frac{2}{\sqrt{42}}c^{1}_{28}
$,
$-\frac{2}{\sqrt{42}}c^{3}_{28}
$;\ \ 
$-\frac{1}{\sqrt{21}}c^{5}_{28}
$,
$\frac{1}{\sqrt{21}}c^{1}_{28}
$,
$\frac{1}{\sqrt{21}}c^{3}_{28}
$;\ \ 
$\frac{1}{\sqrt{21}}c^{3}_{28}
$,
$-\frac{1}{\sqrt{21}}c^{5}_{28}
$;\ \ 
$\frac{1}{\sqrt{21}}c^{1}_{28}
$)

Fail:
cnd($\rho(\mathfrak s)_\mathrm{ndeg}$) = 168 does not divide
 ord($\rho(\mathfrak t)$)=42. Prop. B.4 (2)

 \ \color{black}

\noindent 611: (dims,levels) = $(6;42
)$,
irreps = $6_{7,2}^{1}
\hskip -1.5pt \otimes \hskip -1.5pt
1_{3}^{1,0}
\hskip -1.5pt \otimes \hskip -1.5pt
1_{2}^{1,0}$,
pord$(\rho_\text{isum}(\mathfrak{t})) = 7$,

\vskip 0.7ex
\hangindent=5.5em \hangafter=1
{\white .}\hskip 1em $\rho_\text{isum}(\mathfrak{t})$ =
 $( \frac{5}{42},
\frac{11}{42},
\frac{17}{42},
\frac{23}{42},
\frac{29}{42},
\frac{41}{42} )
$,

\vskip 0.7ex
\hangindent=5.5em \hangafter=1
{\white .}\hskip 1em $\rho_\text{isum}(\mathfrak{s})$ =
($-\frac{3}{7}-\frac{1}{7}c^{1}_{7}
-\frac{1}{7}c^{2}_{7}
$,
$-\frac{2}{7}c^{3}_{56}
-\frac{1}{7}c^{5}_{56}
+\frac{1}{7}c^{7}_{56}
+\frac{1}{7}c^{9}_{56}
-\frac{2}{7}c^{11}_{56}
$,
$-\frac{2}{7}+\frac{1}{7}c^{2}_{7}
$,
$\frac{1}{7}c^{3}_{56}
+\frac{2}{7}c^{5}_{56}
-\frac{1}{7}c^{7}_{56}
-\frac{2}{7}c^{9}_{56}
+\frac{1}{7}c^{11}_{56}
$,
$-\frac{1}{7}c^{3}_{56}
+\frac{1}{7}c^{5}_{56}
-\frac{1}{7}c^{9}_{56}
-\frac{1}{7}c^{11}_{56}
$,
$-\frac{2}{7}+\frac{1}{7}c^{1}_{7}
$;
$-\frac{2}{7}+\frac{1}{7}c^{1}_{7}
$,
$\frac{1}{7}c^{3}_{56}
-\frac{1}{7}c^{5}_{56}
+\frac{1}{7}c^{9}_{56}
+\frac{1}{7}c^{11}_{56}
$,
$-\frac{2}{7}+\frac{1}{7}c^{2}_{7}
$,
$\frac{3}{7}+\frac{1}{7}c^{1}_{7}
+\frac{1}{7}c^{2}_{7}
$,
$\frac{1}{7}c^{3}_{56}
+\frac{2}{7}c^{5}_{56}
-\frac{1}{7}c^{7}_{56}
-\frac{2}{7}c^{9}_{56}
+\frac{1}{7}c^{11}_{56}
$;
$-\frac{2}{7}+\frac{1}{7}c^{1}_{7}
$,
$-\frac{2}{7}c^{3}_{56}
-\frac{1}{7}c^{5}_{56}
+\frac{1}{7}c^{7}_{56}
+\frac{1}{7}c^{9}_{56}
-\frac{2}{7}c^{11}_{56}
$,
$-\frac{1}{7}c^{3}_{56}
-\frac{2}{7}c^{5}_{56}
+\frac{1}{7}c^{7}_{56}
+\frac{2}{7}c^{9}_{56}
-\frac{1}{7}c^{11}_{56}
$,
$-\frac{3}{7}-\frac{1}{7}c^{1}_{7}
-\frac{1}{7}c^{2}_{7}
$;
$-\frac{3}{7}-\frac{1}{7}c^{1}_{7}
-\frac{1}{7}c^{2}_{7}
$,
$\frac{2}{7}-\frac{1}{7}c^{1}_{7}
$,
$\frac{1}{7}c^{3}_{56}
-\frac{1}{7}c^{5}_{56}
+\frac{1}{7}c^{9}_{56}
+\frac{1}{7}c^{11}_{56}
$;
$-\frac{2}{7}+\frac{1}{7}c^{2}_{7}
$,
$\frac{2}{7}c^{3}_{56}
+\frac{1}{7}c^{5}_{56}
-\frac{1}{7}c^{7}_{56}
-\frac{1}{7}c^{9}_{56}
+\frac{2}{7}c^{11}_{56}
$;
$-\frac{2}{7}+\frac{1}{7}c^{2}_{7}
$)

Fail:
cnd($\rho(\mathfrak s)_\mathrm{ndeg}$) = 56 does not divide
 ord($\rho(\mathfrak t)$)=42. Prop. B.4 (2)

 \ \color{black}

\noindent 612: (dims,levels) = $(6;42
)$,
irreps = $6_{7,1}^{1}
\hskip -1.5pt \otimes \hskip -1.5pt
1_{3}^{1,0}
\hskip -1.5pt \otimes \hskip -1.5pt
1_{2}^{1,0}$,
pord$(\rho_\text{isum}(\mathfrak{t})) = 7$,

\vskip 0.7ex
\hangindent=5.5em \hangafter=1
{\white .}\hskip 1em $\rho_\text{isum}(\mathfrak{t})$ =
 $( \frac{5}{42},
\frac{11}{42},
\frac{17}{42},
\frac{23}{42},
\frac{29}{42},
\frac{41}{42} )
$,

\vskip 0.7ex
\hangindent=5.5em \hangafter=1
{\white .}\hskip 1em $\rho_\text{isum}(\mathfrak{s})$ =
$\mathrm{i}$($-\frac{2}{7}c^{1}_{56}
+\frac{1}{7}c^{3}_{56}
+\frac{1}{7}c^{5}_{56}
-\frac{1}{7}c^{7}_{56}
-\frac{1}{7}c^{9}_{56}
+\frac{1}{7}c^{10}_{56}
+\frac{1}{7}c^{11}_{56}
$,
$\frac{1}{7}c^{1}_{112}
+\frac{1}{7}c^{5}_{112}
-\frac{1}{7}c^{15}_{112}
+\frac{1}{7}c^{19}_{112}
$,
$\frac{1}{7}c^{2}_{56}
-\frac{1}{7}c^{3}_{56}
+\frac{1}{7}c^{11}_{56}
$,
$\frac{1}{7}c^{3}_{112}
-\frac{1}{7}c^{9}_{112}
+\frac{1}{7}c^{11}_{112}
+\frac{1}{7}c^{23}_{112}
$,
$\frac{1}{7}c^{1}_{112}
+\frac{1}{7}c^{3}_{112}
-\frac{1}{7}c^{5}_{112}
-\frac{1}{7}c^{7}_{112}
-\frac{1}{7}c^{9}_{112}
-\frac{1}{7}c^{11}_{112}
+\frac{2}{7}c^{13}_{112}
+\frac{1}{7}c^{15}_{112}
+\frac{2}{7}c^{17}_{112}
+\frac{1}{7}c^{19}_{112}
-\frac{1}{7}c^{21}_{112}
-\frac{1}{7}c^{23}_{112}
$,
$\frac{1}{7}c^{5}_{56}
+\frac{1}{7}c^{6}_{56}
+\frac{1}{7}c^{9}_{56}
$;\ \ 
$\frac{1}{7}c^{5}_{56}
+\frac{1}{7}c^{6}_{56}
+\frac{1}{7}c^{9}_{56}
$,
$\frac{1}{7}c^{1}_{112}
+\frac{1}{7}c^{3}_{112}
-\frac{1}{7}c^{5}_{112}
-\frac{1}{7}c^{7}_{112}
-\frac{1}{7}c^{9}_{112}
-\frac{1}{7}c^{11}_{112}
+\frac{2}{7}c^{13}_{112}
+\frac{1}{7}c^{15}_{112}
+\frac{2}{7}c^{17}_{112}
+\frac{1}{7}c^{19}_{112}
-\frac{1}{7}c^{21}_{112}
-\frac{1}{7}c^{23}_{112}
$,
$\frac{1}{7}c^{2}_{56}
-\frac{1}{7}c^{3}_{56}
+\frac{1}{7}c^{11}_{56}
$,
$-\frac{2}{7}c^{1}_{56}
+\frac{1}{7}c^{3}_{56}
+\frac{1}{7}c^{5}_{56}
-\frac{1}{7}c^{7}_{56}
-\frac{1}{7}c^{9}_{56}
+\frac{1}{7}c^{10}_{56}
+\frac{1}{7}c^{11}_{56}
$,
$-\frac{1}{7}c^{3}_{112}
+\frac{1}{7}c^{9}_{112}
-\frac{1}{7}c^{11}_{112}
-\frac{1}{7}c^{23}_{112}
$;\ \ 
$-\frac{1}{7}c^{5}_{56}
-\frac{1}{7}c^{6}_{56}
-\frac{1}{7}c^{9}_{56}
$,
$-\frac{1}{7}c^{1}_{112}
-\frac{1}{7}c^{5}_{112}
+\frac{1}{7}c^{15}_{112}
-\frac{1}{7}c^{19}_{112}
$,
$\frac{1}{7}c^{3}_{112}
-\frac{1}{7}c^{9}_{112}
+\frac{1}{7}c^{11}_{112}
+\frac{1}{7}c^{23}_{112}
$,
$-\frac{2}{7}c^{1}_{56}
+\frac{1}{7}c^{3}_{56}
+\frac{1}{7}c^{5}_{56}
-\frac{1}{7}c^{7}_{56}
-\frac{1}{7}c^{9}_{56}
+\frac{1}{7}c^{10}_{56}
+\frac{1}{7}c^{11}_{56}
$;\ \ 
$\frac{2}{7}c^{1}_{56}
-\frac{1}{7}c^{3}_{56}
-\frac{1}{7}c^{5}_{56}
+\frac{1}{7}c^{7}_{56}
+\frac{1}{7}c^{9}_{56}
-\frac{1}{7}c^{10}_{56}
-\frac{1}{7}c^{11}_{56}
$,
$-\frac{1}{7}c^{5}_{56}
-\frac{1}{7}c^{6}_{56}
-\frac{1}{7}c^{9}_{56}
$,
$\frac{1}{7}c^{1}_{112}
+\frac{1}{7}c^{3}_{112}
-\frac{1}{7}c^{5}_{112}
-\frac{1}{7}c^{7}_{112}
-\frac{1}{7}c^{9}_{112}
-\frac{1}{7}c^{11}_{112}
+\frac{2}{7}c^{13}_{112}
+\frac{1}{7}c^{15}_{112}
+\frac{2}{7}c^{17}_{112}
+\frac{1}{7}c^{19}_{112}
-\frac{1}{7}c^{21}_{112}
-\frac{1}{7}c^{23}_{112}
$;\ \ 
$\frac{1}{7}c^{2}_{56}
-\frac{1}{7}c^{3}_{56}
+\frac{1}{7}c^{11}_{56}
$,
$\frac{1}{7}c^{1}_{112}
+\frac{1}{7}c^{5}_{112}
-\frac{1}{7}c^{15}_{112}
+\frac{1}{7}c^{19}_{112}
$;\ \ 
$-\frac{1}{7}c^{2}_{56}
+\frac{1}{7}c^{3}_{56}
-\frac{1}{7}c^{11}_{56}
$)

Fail:
cnd( Tr$_I(\rho(\mathfrak s))$ ) =
56 does not divide
 ord($\rho(\mathfrak t)$) =
42, I = [ 5/42 ]. Prop. B.4 (2)

 \ \color{black}

\noindent 613: (dims,levels) = $(6;42
)$,
irreps = $6_{7,1}^{3}
\hskip -1.5pt \otimes \hskip -1.5pt
1_{3}^{1,0}
\hskip -1.5pt \otimes \hskip -1.5pt
1_{2}^{1,0}$,
pord$(\rho_\text{isum}(\mathfrak{t})) = 7$,

\vskip 0.7ex
\hangindent=5.5em \hangafter=1
{\white .}\hskip 1em $\rho_\text{isum}(\mathfrak{t})$ =
 $( \frac{5}{42},
\frac{11}{42},
\frac{17}{42},
\frac{23}{42},
\frac{29}{42},
\frac{41}{42} )
$,

\vskip 0.7ex
\hangindent=5.5em \hangafter=1
{\white .}\hskip 1em $\rho_\text{isum}(\mathfrak{s})$ =
$\mathrm{i}$($\frac{2}{7}c^{1}_{56}
-\frac{1}{7}c^{3}_{56}
-\frac{1}{7}c^{5}_{56}
+\frac{1}{7}c^{7}_{56}
+\frac{1}{7}c^{9}_{56}
+\frac{1}{7}c^{10}_{56}
-\frac{1}{7}c^{11}_{56}
$,
$\frac{1}{7}c^{3}_{112}
-\frac{1}{7}c^{5}_{112}
-\frac{1}{7}c^{9}_{112}
-\frac{1}{7}c^{11}_{112}
+\frac{2}{7}c^{13}_{112}
+\frac{1}{7}c^{19}_{112}
-\frac{1}{7}c^{21}_{112}
+\frac{1}{7}c^{23}_{112}
$,
$\frac{1}{7}c^{2}_{56}
+\frac{1}{7}c^{3}_{56}
-\frac{1}{7}c^{11}_{56}
$,
$\frac{1}{7}c^{1}_{112}
+\frac{1}{7}c^{5}_{112}
-\frac{1}{7}c^{7}_{112}
-\frac{1}{7}c^{9}_{112}
+\frac{1}{7}c^{15}_{112}
+\frac{2}{7}c^{17}_{112}
+\frac{1}{7}c^{19}_{112}
-\frac{1}{7}c^{23}_{112}
$,
$-\frac{1}{7}c^{1}_{112}
+\frac{1}{7}c^{3}_{112}
+\frac{1}{7}c^{11}_{112}
+\frac{1}{7}c^{15}_{112}
$,
$\frac{1}{7}c^{5}_{56}
-\frac{1}{7}c^{6}_{56}
+\frac{1}{7}c^{9}_{56}
$;\ \ 
$-\frac{1}{7}c^{5}_{56}
+\frac{1}{7}c^{6}_{56}
-\frac{1}{7}c^{9}_{56}
$,
$-\frac{1}{7}c^{1}_{112}
+\frac{1}{7}c^{3}_{112}
+\frac{1}{7}c^{11}_{112}
+\frac{1}{7}c^{15}_{112}
$,
$-\frac{1}{7}c^{2}_{56}
-\frac{1}{7}c^{3}_{56}
+\frac{1}{7}c^{11}_{56}
$,
$\frac{2}{7}c^{1}_{56}
-\frac{1}{7}c^{3}_{56}
-\frac{1}{7}c^{5}_{56}
+\frac{1}{7}c^{7}_{56}
+\frac{1}{7}c^{9}_{56}
+\frac{1}{7}c^{10}_{56}
-\frac{1}{7}c^{11}_{56}
$,
$-\frac{1}{7}c^{1}_{112}
-\frac{1}{7}c^{5}_{112}
+\frac{1}{7}c^{7}_{112}
+\frac{1}{7}c^{9}_{112}
-\frac{1}{7}c^{15}_{112}
-\frac{2}{7}c^{17}_{112}
-\frac{1}{7}c^{19}_{112}
+\frac{1}{7}c^{23}_{112}
$;\ \ 
$\frac{1}{7}c^{5}_{56}
-\frac{1}{7}c^{6}_{56}
+\frac{1}{7}c^{9}_{56}
$,
$\frac{1}{7}c^{3}_{112}
-\frac{1}{7}c^{5}_{112}
-\frac{1}{7}c^{9}_{112}
-\frac{1}{7}c^{11}_{112}
+\frac{2}{7}c^{13}_{112}
+\frac{1}{7}c^{19}_{112}
-\frac{1}{7}c^{21}_{112}
+\frac{1}{7}c^{23}_{112}
$,
$-\frac{1}{7}c^{1}_{112}
-\frac{1}{7}c^{5}_{112}
+\frac{1}{7}c^{7}_{112}
+\frac{1}{7}c^{9}_{112}
-\frac{1}{7}c^{15}_{112}
-\frac{2}{7}c^{17}_{112}
-\frac{1}{7}c^{19}_{112}
+\frac{1}{7}c^{23}_{112}
$,
$-\frac{2}{7}c^{1}_{56}
+\frac{1}{7}c^{3}_{56}
+\frac{1}{7}c^{5}_{56}
-\frac{1}{7}c^{7}_{56}
-\frac{1}{7}c^{9}_{56}
-\frac{1}{7}c^{10}_{56}
+\frac{1}{7}c^{11}_{56}
$;\ \ 
$-\frac{2}{7}c^{1}_{56}
+\frac{1}{7}c^{3}_{56}
+\frac{1}{7}c^{5}_{56}
-\frac{1}{7}c^{7}_{56}
-\frac{1}{7}c^{9}_{56}
-\frac{1}{7}c^{10}_{56}
+\frac{1}{7}c^{11}_{56}
$,
$-\frac{1}{7}c^{5}_{56}
+\frac{1}{7}c^{6}_{56}
-\frac{1}{7}c^{9}_{56}
$,
$-\frac{1}{7}c^{1}_{112}
+\frac{1}{7}c^{3}_{112}
+\frac{1}{7}c^{11}_{112}
+\frac{1}{7}c^{15}_{112}
$;\ \ 
$\frac{1}{7}c^{2}_{56}
+\frac{1}{7}c^{3}_{56}
-\frac{1}{7}c^{11}_{56}
$,
$-\frac{1}{7}c^{3}_{112}
+\frac{1}{7}c^{5}_{112}
+\frac{1}{7}c^{9}_{112}
+\frac{1}{7}c^{11}_{112}
-\frac{2}{7}c^{13}_{112}
-\frac{1}{7}c^{19}_{112}
+\frac{1}{7}c^{21}_{112}
-\frac{1}{7}c^{23}_{112}
$;\ \ 
$-\frac{1}{7}c^{2}_{56}
-\frac{1}{7}c^{3}_{56}
+\frac{1}{7}c^{11}_{56}
$)

Fail:
cnd( Tr$_I(\rho(\mathfrak s))$ ) =
56 does not divide
 ord($\rho(\mathfrak t)$) =
42, I = [ 5/42 ]. Prop. B.4 (2)

 \ \color{black}

\noindent 614: (dims,levels) = $(6;42
)$,
irreps = $3_{7}^{3}
\hskip -1.5pt \otimes \hskip -1.5pt
2_{3}^{1,0}
\hskip -1.5pt \otimes \hskip -1.5pt
1_{2}^{1,0}$,
pord$(\rho_\text{isum}(\mathfrak{t})) = 21$,

\vskip 0.7ex
\hangindent=5.5em \hangafter=1
{\white .}\hskip 1em $\rho_\text{isum}(\mathfrak{t})$ =
 $( \frac{3}{14},
\frac{5}{14},
\frac{13}{14},
\frac{11}{42},
\frac{23}{42},
\frac{29}{42} )
$,

\vskip 0.7ex
\hangindent=5.5em \hangafter=1
{\white .}\hskip 1em $\rho_\text{isum}(\mathfrak{s})$ =
$\mathrm{i}$($\frac{1}{\sqrt{21}}c^{5}_{28}
$,
$-\frac{1}{\sqrt{21}}c^{3}_{28}
$,
$-\frac{1}{\sqrt{21}}c^{1}_{28}
$,
$-\frac{2}{\sqrt{42}}c^{1}_{28}
$,
$\frac{2}{\sqrt{42}}c^{5}_{28}
$,
$-\frac{2}{\sqrt{42}}c^{3}_{28}
$;\ \ 
$-\frac{1}{\sqrt{21}}c^{1}_{28}
$,
$\frac{1}{\sqrt{21}}c^{5}_{28}
$,
$\frac{2}{\sqrt{42}}c^{5}_{28}
$,
$-\frac{2}{\sqrt{42}}c^{3}_{28}
$,
$-\frac{2}{\sqrt{42}}c^{1}_{28}
$;\ \ 
$-\frac{1}{\sqrt{21}}c^{3}_{28}
$,
$-\frac{2}{\sqrt{42}}c^{3}_{28}
$,
$-\frac{2}{\sqrt{42}}c^{1}_{28}
$,
$\frac{2}{\sqrt{42}}c^{5}_{28}
$;\ \ 
$\frac{1}{\sqrt{21}}c^{3}_{28}
$,
$\frac{1}{\sqrt{21}}c^{1}_{28}
$,
$-\frac{1}{\sqrt{21}}c^{5}_{28}
$;\ \ 
$-\frac{1}{\sqrt{21}}c^{5}_{28}
$,
$\frac{1}{\sqrt{21}}c^{3}_{28}
$;\ \ 
$\frac{1}{\sqrt{21}}c^{1}_{28}
$)

Fail:
cnd($\rho(\mathfrak s)_\mathrm{ndeg}$) = 168 does not divide
 ord($\rho(\mathfrak t)$)=42. Prop. B.4 (2)

 \ \color{black}

\noindent 615: (dims,levels) = $(6;42
)$,
irreps = $3_{7}^{1}
\hskip -1.5pt \otimes \hskip -1.5pt
2_{2}^{1,0}
\hskip -1.5pt \otimes \hskip -1.5pt
1_{3}^{1,0}$,
pord$(\rho_\text{isum}(\mathfrak{t})) = 14$,

\vskip 0.7ex
\hangindent=5.5em \hangafter=1
{\white .}\hskip 1em $\rho_\text{isum}(\mathfrak{t})$ =
 $( \frac{10}{21},
\frac{13}{21},
\frac{19}{21},
\frac{5}{42},
\frac{17}{42},
\frac{41}{42} )
$,

\vskip 0.7ex
\hangindent=5.5em \hangafter=1
{\white .}\hskip 1em $\rho_\text{isum}(\mathfrak{s})$ =
($\frac{1}{2\sqrt{7}}c^{1}_{28}
$,
$-\frac{1}{2\sqrt{7}}c^{3}_{28}
$,
$\frac{1}{2\sqrt{7}}c^{5}_{28}
$,
$\frac{3}{2\sqrt{21}}c^{3}_{28}
$,
$-\frac{3}{2\sqrt{21}}c^{5}_{28}
$,
$\frac{3}{2\sqrt{21}}c^{1}_{28}
$;
$-\frac{1}{2\sqrt{7}}c^{5}_{28}
$,
$\frac{1}{2\sqrt{7}}c^{1}_{28}
$,
$\frac{3}{2\sqrt{21}}c^{5}_{28}
$,
$-\frac{3}{2\sqrt{21}}c^{1}_{28}
$,
$-\frac{3}{2\sqrt{21}}c^{3}_{28}
$;
$\frac{1}{2\sqrt{7}}c^{3}_{28}
$,
$-\frac{3}{2\sqrt{21}}c^{1}_{28}
$,
$-\frac{3}{2\sqrt{21}}c^{3}_{28}
$,
$\frac{3}{2\sqrt{21}}c^{5}_{28}
$;
$\frac{1}{2\sqrt{7}}c^{5}_{28}
$,
$-\frac{1}{2\sqrt{7}}c^{1}_{28}
$,
$-\frac{1}{2\sqrt{7}}c^{3}_{28}
$;
$-\frac{1}{2\sqrt{7}}c^{3}_{28}
$,
$\frac{1}{2\sqrt{7}}c^{5}_{28}
$;
$-\frac{1}{2\sqrt{7}}c^{1}_{28}
$)

Fail:
cnd($\rho(\mathfrak s)_\mathrm{ndeg}$) = 84 does not divide
 ord($\rho(\mathfrak t)$)=42. Prop. B.4 (2)

 \ \color{black}

\noindent 616: (dims,levels) = $(6;44
)$,
irreps = $6_{11}^{1}
\hskip -1.5pt \otimes \hskip -1.5pt
1_{4}^{1,0}$,
pord$(\rho_\text{isum}(\mathfrak{t})) = 11$,

\vskip 0.7ex
\hangindent=5.5em \hangafter=1
{\white .}\hskip 1em $\rho_\text{isum}(\mathfrak{t})$ =
 $( \frac{1}{4},
\frac{3}{44},
\frac{15}{44},
\frac{23}{44},
\frac{27}{44},
\frac{31}{44} )
$,

\vskip 0.7ex
\hangindent=5.5em \hangafter=1
{\white .}\hskip 1em $\rho_\text{isum}(\mathfrak{s})$ =
($\sqrt{\frac{1}{11}}$,
$\sqrt{\frac{2}{11}}$,
$\sqrt{\frac{2}{11}}$,
$\sqrt{\frac{2}{11}}$,
$\sqrt{\frac{2}{11}}$,
$\sqrt{\frac{2}{11}}$;
$\frac{1}{\sqrt{11}}c^{4}_{11}
$,
$-\frac{1}{\sqrt{11}\mathrm{i}}s^{9}_{44}
$,
$\frac{1}{\sqrt{11}}c^{3}_{11}
$,
$\frac{1}{\sqrt{11}}c^{1}_{11}
$,
$\frac{1}{\sqrt{11}}c^{2}_{11}
$;
$\frac{1}{\sqrt{11}}c^{2}_{11}
$,
$\frac{1}{\sqrt{11}}c^{1}_{11}
$,
$\frac{1}{\sqrt{11}}c^{4}_{11}
$,
$\frac{1}{\sqrt{11}}c^{3}_{11}
$;
$-\frac{1}{\sqrt{11}\mathrm{i}}s^{9}_{44}
$,
$\frac{1}{\sqrt{11}}c^{2}_{11}
$,
$\frac{1}{\sqrt{11}}c^{4}_{11}
$;
$\frac{1}{\sqrt{11}}c^{3}_{11}
$,
$-\frac{1}{\sqrt{11}\mathrm{i}}s^{9}_{44}
$;
$\frac{1}{\sqrt{11}}c^{1}_{11}
$)

Fail:
cnd($\rho(\mathfrak s)_\mathrm{ndeg}$) = 88 does not divide
 ord($\rho(\mathfrak t)$)=44. Prop. B.4 (2)

 \ \color{black}

\noindent 617: (dims,levels) = $(6;44
)$,
irreps = $6_{11}^{7}
\hskip -1.5pt \otimes \hskip -1.5pt
1_{4}^{1,0}$,
pord$(\rho_\text{isum}(\mathfrak{t})) = 11$,

\vskip 0.7ex
\hangindent=5.5em \hangafter=1
{\white .}\hskip 1em $\rho_\text{isum}(\mathfrak{t})$ =
 $( \frac{1}{4},
\frac{7}{44},
\frac{19}{44},
\frac{35}{44},
\frac{39}{44},
\frac{43}{44} )
$,

\vskip 0.7ex
\hangindent=5.5em \hangafter=1
{\white .}\hskip 1em $\rho_\text{isum}(\mathfrak{s})$ =
($-\sqrt{\frac{1}{11}}$,
$\sqrt{\frac{2}{11}}$,
$\sqrt{\frac{2}{11}}$,
$\sqrt{\frac{2}{11}}$,
$\sqrt{\frac{2}{11}}$,
$\sqrt{\frac{2}{11}}$;
$-\frac{1}{\sqrt{11}}c^{2}_{11}
$,
$\frac{1}{\sqrt{11}\mathrm{i}}s^{9}_{44}
$,
$-\frac{1}{\sqrt{11}}c^{3}_{11}
$,
$-\frac{1}{\sqrt{11}}c^{4}_{11}
$,
$-\frac{1}{\sqrt{11}}c^{1}_{11}
$;
$-\frac{1}{\sqrt{11}}c^{4}_{11}
$,
$-\frac{1}{\sqrt{11}}c^{2}_{11}
$,
$-\frac{1}{\sqrt{11}}c^{1}_{11}
$,
$-\frac{1}{\sqrt{11}}c^{3}_{11}
$;
$-\frac{1}{\sqrt{11}}c^{1}_{11}
$,
$\frac{1}{\sqrt{11}\mathrm{i}}s^{9}_{44}
$,
$-\frac{1}{\sqrt{11}}c^{4}_{11}
$;
$-\frac{1}{\sqrt{11}}c^{3}_{11}
$,
$-\frac{1}{\sqrt{11}}c^{2}_{11}
$;
$\frac{1}{\sqrt{11}\mathrm{i}}s^{9}_{44}
$)

Fail:
cnd($\rho(\mathfrak s)_\mathrm{ndeg}$) = 88 does not divide
 ord($\rho(\mathfrak t)$)=44. Prop. B.4 (2)

 \ \color{black}

\noindent 618: (dims,levels) = $(6;48
)$,
irreps = $3_{16}^{1,0}
\hskip -1.5pt \otimes \hskip -1.5pt
2_{3}^{1,0}$,
pord$(\rho_\text{isum}(\mathfrak{t})) = 48$,

\vskip 0.7ex
\hangindent=5.5em \hangafter=1
{\white .}\hskip 1em $\rho_\text{isum}(\mathfrak{t})$ =
 $( \frac{1}{8},
\frac{1}{16},
\frac{9}{16},
\frac{11}{24},
\frac{19}{48},
\frac{43}{48} )
$,

\vskip 0.7ex
\hangindent=5.5em \hangafter=1
{\white .}\hskip 1em $\rho_\text{isum}(\mathfrak{s})$ =
($0$,
$\sqrt{\frac{1}{6}}$,
$\sqrt{\frac{1}{6}}$,
$0$,
$-\sqrt{\frac{1}{3}}$,
$-\sqrt{\frac{1}{3}}$;
$-\sqrt{\frac{1}{12}}$,
$\sqrt{\frac{1}{12}}$,
$-\sqrt{\frac{1}{3}}$,
$\sqrt{\frac{1}{6}}$,
$-\sqrt{\frac{1}{6}}$;
$-\sqrt{\frac{1}{12}}$,
$-\sqrt{\frac{1}{3}}$,
$-\sqrt{\frac{1}{6}}$,
$\sqrt{\frac{1}{6}}$;
$0$,
$-\sqrt{\frac{1}{6}}$,
$-\sqrt{\frac{1}{6}}$;
$\sqrt{\frac{1}{12}}$,
$-\sqrt{\frac{1}{12}}$;
$\sqrt{\frac{1}{12}}$)

Fail:
$\sigma(\rho(\mathfrak s)_\mathrm{ndeg}) \neq
 (\rho(\mathfrak t)^a \rho(\mathfrak s) \rho(\mathfrak t)^b
 \rho(\mathfrak s) \rho(\mathfrak t)^a)_\mathrm{ndeg}$,
 $\sigma = a$ = 5. Prop. B.5 (3) eqn. (B.25)

 \ \color{black}

\noindent 619: (dims,levels) = $(6;48
)$,
irreps = $6_{16,4}^{3,0}
\hskip -1.5pt \otimes \hskip -1.5pt
1_{3}^{1,0}$,
pord$(\rho_\text{isum}(\mathfrak{t})) = 16$,

\vskip 0.7ex
\hangindent=5.5em \hangafter=1
{\white .}\hskip 1em $\rho_\text{isum}(\mathfrak{t})$ =
 $( \frac{5}{24},
\frac{17}{24},
\frac{1}{48},
\frac{13}{48},
\frac{25}{48},
\frac{37}{48} )
$,

\vskip 0.7ex
\hangindent=5.5em \hangafter=1
{\white .}\hskip 1em $\rho_\text{isum}(\mathfrak{s})$ =
($0$,
$0$,
$\frac{1}{2}$,
$\frac{1}{2}$,
$\frac{1}{2}$,
$\frac{1}{2}$;
$0$,
$\frac{1}{2}$,
$-\frac{1}{2}$,
$\frac{1}{2}$,
$-\frac{1}{2}$;
$0$,
$-\frac{1}{2}$,
$0$,
$\frac{1}{2}$;
$0$,
$\frac{1}{2}$,
$0$;
$0$,
$-\frac{1}{2}$;
$0$)

Fail:
all rows of $U \rho(\mathfrak s) U^\dagger$
 contain zero for any block-diagonal $U$. Prop. B.5 (4) eqn. (B.27)

 \ \color{black}

 \color{blue}

\noindent 620: (dims,levels) = $(6;48
)$,
irreps = $6_{16,1}^{3,0}
\hskip -1.5pt \otimes \hskip -1.5pt
1_{3}^{1,0}$,
pord$(\rho_\text{isum}(\mathfrak{t})) = 16$,

\vskip 0.7ex
\hangindent=5.5em \hangafter=1
{\white .}\hskip 1em $\rho_\text{isum}(\mathfrak{t})$ =
 $( \frac{1}{3},
\frac{1}{12},
\frac{1}{48},
\frac{13}{48},
\frac{25}{48},
\frac{37}{48} )
$,

\vskip 0.7ex
\hangindent=5.5em \hangafter=1
{\white .}\hskip 1em $\rho_\text{isum}(\mathfrak{s})$ =
$\mathrm{i}$($0$,
$0$,
$\frac{1}{2}$,
$\frac{1}{2}$,
$\frac{1}{2}$,
$\frac{1}{2}$;\ \ 
$0$,
$\frac{1}{2}$,
$-\frac{1}{2}$,
$\frac{1}{2}$,
$-\frac{1}{2}$;\ \ 
$-\sqrt{\frac{1}{8}}$,
$\sqrt{\frac{1}{8}}$,
$\sqrt{\frac{1}{8}}$,
$-\sqrt{\frac{1}{8}}$;\ \ 
$\sqrt{\frac{1}{8}}$,
$-\sqrt{\frac{1}{8}}$,
$-\sqrt{\frac{1}{8}}$;\ \ 
$-\sqrt{\frac{1}{8}}$,
$\sqrt{\frac{1}{8}}$;\ \ 
$\sqrt{\frac{1}{8}}$)

Pass. 

 \ \color{black}

 \color{blue}

\noindent 621: (dims,levels) = $(6;48
)$,
irreps = $6_{16,1}^{1,0}
\hskip -1.5pt \otimes \hskip -1.5pt
1_{3}^{1,0}$,
pord$(\rho_\text{isum}(\mathfrak{t})) = 16$,

\vskip 0.7ex
\hangindent=5.5em \hangafter=1
{\white .}\hskip 1em $\rho_\text{isum}(\mathfrak{t})$ =
 $( \frac{1}{3},
\frac{7}{12},
\frac{7}{48},
\frac{19}{48},
\frac{31}{48},
\frac{43}{48} )
$,

\vskip 0.7ex
\hangindent=5.5em \hangafter=1
{\white .}\hskip 1em $\rho_\text{isum}(\mathfrak{s})$ =
$\mathrm{i}$($0$,
$0$,
$\frac{1}{2}$,
$\frac{1}{2}$,
$\frac{1}{2}$,
$\frac{1}{2}$;\ \ 
$0$,
$\frac{1}{2}$,
$-\frac{1}{2}$,
$\frac{1}{2}$,
$-\frac{1}{2}$;\ \ 
$\sqrt{\frac{1}{8}}$,
$\sqrt{\frac{1}{8}}$,
$-\sqrt{\frac{1}{8}}$,
$-\sqrt{\frac{1}{8}}$;\ \ 
$-\sqrt{\frac{1}{8}}$,
$-\sqrt{\frac{1}{8}}$,
$\sqrt{\frac{1}{8}}$;\ \ 
$\sqrt{\frac{1}{8}}$,
$\sqrt{\frac{1}{8}}$;\ \ 
$-\sqrt{\frac{1}{8}}$)

Pass. 

 \ \color{black}

\noindent 622: (dims,levels) = $(6;48
)$,
irreps = $6_{16,3}^{3,0}
\hskip -1.5pt \otimes \hskip -1.5pt
1_{3}^{1,0}$,
pord$(\rho_\text{isum}(\mathfrak{t})) = 16$,

\vskip 0.7ex
\hangindent=5.5em \hangafter=1
{\white .}\hskip 1em $\rho_\text{isum}(\mathfrak{t})$ =
 $( \frac{1}{3},
\frac{5}{6},
\frac{1}{48},
\frac{7}{48},
\frac{25}{48},
\frac{31}{48} )
$,

\vskip 0.7ex
\hangindent=5.5em \hangafter=1
{\white .}\hskip 1em $\rho_\text{isum}(\mathfrak{s})$ =
($0$,
$0$,
$\frac{1}{2}$,
$\frac{1}{2}$,
$\frac{1}{2}$,
$\frac{1}{2}$;
$0$,
$\frac{1}{2}$,
$-\frac{1}{2}$,
$\frac{1}{2}$,
$-\frac{1}{2}$;
$0$,
$-\frac{1}{2}$,
$0$,
$\frac{1}{2}$;
$0$,
$\frac{1}{2}$,
$0$;
$0$,
$-\frac{1}{2}$;
$0$)

Fail:
all rows of $U \rho(\mathfrak s) U^\dagger$
 contain zero for any block-diagonal $U$. Prop. B.5 (4) eqn. (B.27)

 \ \color{black}

\noindent 623: (dims,levels) = $(6;48
)$,
irreps = $6_{16,2}^{1,0}
\hskip -1.5pt \otimes \hskip -1.5pt
1_{3}^{1,0}$,
pord$(\rho_\text{isum}(\mathfrak{t})) = 16$,

\vskip 0.7ex
\hangindent=5.5em \hangafter=1
{\white .}\hskip 1em $\rho_\text{isum}(\mathfrak{t})$ =
 $( \frac{1}{3},
\frac{5}{6},
\frac{1}{48},
\frac{19}{48},
\frac{25}{48},
\frac{43}{48} )
$,

\vskip 0.7ex
\hangindent=5.5em \hangafter=1
{\white .}\hskip 1em $\rho_\text{isum}(\mathfrak{s})$ =
$\mathrm{i}$($0$,
$0$,
$\frac{1}{2}$,
$\frac{1}{2}$,
$\frac{1}{2}$,
$\frac{1}{2}$;\ \ 
$0$,
$\frac{1}{2}$,
$-\frac{1}{2}$,
$\frac{1}{2}$,
$-\frac{1}{2}$;\ \ 
$0$,
$\frac{1}{2}$,
$0$,
$-\frac{1}{2}$;\ \ 
$0$,
$-\frac{1}{2}$,
$0$;\ \ 
$0$,
$\frac{1}{2}$;\ \ 
$0$)

Fail:
all rows of $U \rho(\mathfrak s) U^\dagger$
 contain zero for any block-diagonal $U$. Prop. B.5 (4) eqn. (B.27)

 \ \color{black}

\noindent 624: (dims,levels) = $(6;48
)$,
irreps = $6_{16,2}^{5,0}
\hskip -1.5pt \otimes \hskip -1.5pt
1_{3}^{1,0}$,
pord$(\rho_\text{isum}(\mathfrak{t})) = 16$,

\vskip 0.7ex
\hangindent=5.5em \hangafter=1
{\white .}\hskip 1em $\rho_\text{isum}(\mathfrak{t})$ =
 $( \frac{1}{3},
\frac{5}{6},
\frac{7}{48},
\frac{13}{48},
\frac{31}{48},
\frac{37}{48} )
$,

\vskip 0.7ex
\hangindent=5.5em \hangafter=1
{\white .}\hskip 1em $\rho_\text{isum}(\mathfrak{s})$ =
$\mathrm{i}$($0$,
$0$,
$\frac{1}{2}$,
$\frac{1}{2}$,
$\frac{1}{2}$,
$\frac{1}{2}$;\ \ 
$0$,
$\frac{1}{2}$,
$-\frac{1}{2}$,
$\frac{1}{2}$,
$-\frac{1}{2}$;\ \ 
$0$,
$\frac{1}{2}$,
$0$,
$-\frac{1}{2}$;\ \ 
$0$,
$-\frac{1}{2}$,
$0$;\ \ 
$0$,
$\frac{1}{2}$;\ \ 
$0$)

Fail:
all rows of $U \rho(\mathfrak s) U^\dagger$
 contain zero for any block-diagonal $U$. Prop. B.5 (4) eqn. (B.27)

 \ \color{black}

\noindent 625: (dims,levels) = $(6;48
)$,
irreps = $6_{16,3}^{1,0}
\hskip -1.5pt \otimes \hskip -1.5pt
1_{3}^{1,0}$,
pord$(\rho_\text{isum}(\mathfrak{t})) = 16$,

\vskip 0.7ex
\hangindent=5.5em \hangafter=1
{\white .}\hskip 1em $\rho_\text{isum}(\mathfrak{t})$ =
 $( \frac{1}{3},
\frac{5}{6},
\frac{13}{48},
\frac{19}{48},
\frac{37}{48},
\frac{43}{48} )
$,

\vskip 0.7ex
\hangindent=5.5em \hangafter=1
{\white .}\hskip 1em $\rho_\text{isum}(\mathfrak{s})$ =
($0$,
$0$,
$\frac{1}{2}$,
$\frac{1}{2}$,
$\frac{1}{2}$,
$\frac{1}{2}$;
$0$,
$\frac{1}{2}$,
$-\frac{1}{2}$,
$\frac{1}{2}$,
$-\frac{1}{2}$;
$0$,
$\frac{1}{2}$,
$0$,
$-\frac{1}{2}$;
$0$,
$-\frac{1}{2}$,
$0$;
$0$,
$\frac{1}{2}$;
$0$)

Fail:
all rows of $U \rho(\mathfrak s) U^\dagger$
 contain zero for any block-diagonal $U$. Prop. B.5 (4) eqn. (B.27)

 \ \color{black}

\noindent 626: (dims,levels) = $(6;48
)$,
irreps = $3_{16}^{3,0}
\hskip -1.5pt \otimes \hskip -1.5pt
2_{3}^{1,0}$,
pord$(\rho_\text{isum}(\mathfrak{t})) = 48$,

\vskip 0.7ex
\hangindent=5.5em \hangafter=1
{\white .}\hskip 1em $\rho_\text{isum}(\mathfrak{t})$ =
 $( \frac{3}{8},
\frac{3}{16},
\frac{11}{16},
\frac{17}{24},
\frac{1}{48},
\frac{25}{48} )
$,

\vskip 0.7ex
\hangindent=5.5em \hangafter=1
{\white .}\hskip 1em $\rho_\text{isum}(\mathfrak{s})$ =
($0$,
$\sqrt{\frac{1}{6}}$,
$\sqrt{\frac{1}{6}}$,
$0$,
$-\sqrt{\frac{1}{3}}$,
$-\sqrt{\frac{1}{3}}$;
$\sqrt{\frac{1}{12}}$,
$-\sqrt{\frac{1}{12}}$,
$-\sqrt{\frac{1}{3}}$,
$\sqrt{\frac{1}{6}}$,
$-\sqrt{\frac{1}{6}}$;
$\sqrt{\frac{1}{12}}$,
$-\sqrt{\frac{1}{3}}$,
$-\sqrt{\frac{1}{6}}$,
$\sqrt{\frac{1}{6}}$;
$0$,
$-\sqrt{\frac{1}{6}}$,
$-\sqrt{\frac{1}{6}}$;
$-\sqrt{\frac{1}{12}}$,
$\sqrt{\frac{1}{12}}$;
$-\sqrt{\frac{1}{12}}$)

Fail:
$\sigma(\rho(\mathfrak s)_\mathrm{ndeg}) \neq
 (\rho(\mathfrak t)^a \rho(\mathfrak s) \rho(\mathfrak t)^b
 \rho(\mathfrak s) \rho(\mathfrak t)^a)_\mathrm{ndeg}$,
 $\sigma = a$ = 5. Prop. B.5 (3) eqn. (B.25)

 \ \color{black}

\noindent 627: (dims,levels) = $(6;48
)$,
irreps = $6_{16,4}^{1,0}
\hskip -1.5pt \otimes \hskip -1.5pt
1_{3}^{1,0}$,
pord$(\rho_\text{isum}(\mathfrak{t})) = 16$,

\vskip 0.7ex
\hangindent=5.5em \hangafter=1
{\white .}\hskip 1em $\rho_\text{isum}(\mathfrak{t})$ =
 $( \frac{11}{24},
\frac{23}{24},
\frac{7}{48},
\frac{19}{48},
\frac{31}{48},
\frac{43}{48} )
$,

\vskip 0.7ex
\hangindent=5.5em \hangafter=1
{\white .}\hskip 1em $\rho_\text{isum}(\mathfrak{s})$ =
($0$,
$0$,
$\frac{1}{2}$,
$\frac{1}{2}$,
$\frac{1}{2}$,
$\frac{1}{2}$;
$0$,
$\frac{1}{2}$,
$-\frac{1}{2}$,
$\frac{1}{2}$,
$-\frac{1}{2}$;
$0$,
$\frac{1}{2}$,
$0$,
$-\frac{1}{2}$;
$0$,
$-\frac{1}{2}$,
$0$;
$0$,
$\frac{1}{2}$;
$0$)

Fail:
all rows of $U \rho(\mathfrak s) U^\dagger$
 contain zero for any block-diagonal $U$. Prop. B.5 (4) eqn. (B.27)

 \ \color{black}

\noindent 628: (dims,levels) = $(6;48
)$,
irreps = $3_{16}^{5,0}
\hskip -1.5pt \otimes \hskip -1.5pt
2_{3}^{1,0}$,
pord$(\rho_\text{isum}(\mathfrak{t})) = 48$,

\vskip 0.7ex
\hangindent=5.5em \hangafter=1
{\white .}\hskip 1em $\rho_\text{isum}(\mathfrak{t})$ =
 $( \frac{5}{8},
\frac{5}{16},
\frac{13}{16},
\frac{23}{24},
\frac{7}{48},
\frac{31}{48} )
$,

\vskip 0.7ex
\hangindent=5.5em \hangafter=1
{\white .}\hskip 1em $\rho_\text{isum}(\mathfrak{s})$ =
($0$,
$\sqrt{\frac{1}{6}}$,
$\sqrt{\frac{1}{6}}$,
$0$,
$-\sqrt{\frac{1}{3}}$,
$-\sqrt{\frac{1}{3}}$;
$-\sqrt{\frac{1}{12}}$,
$\sqrt{\frac{1}{12}}$,
$-\sqrt{\frac{1}{3}}$,
$-\sqrt{\frac{1}{6}}$,
$\sqrt{\frac{1}{6}}$;
$-\sqrt{\frac{1}{12}}$,
$-\sqrt{\frac{1}{3}}$,
$\sqrt{\frac{1}{6}}$,
$-\sqrt{\frac{1}{6}}$;
$0$,
$-\sqrt{\frac{1}{6}}$,
$-\sqrt{\frac{1}{6}}$;
$\sqrt{\frac{1}{12}}$,
$-\sqrt{\frac{1}{12}}$;
$\sqrt{\frac{1}{12}}$)

Fail:
$\sigma(\rho(\mathfrak s)_\mathrm{ndeg}) \neq
 (\rho(\mathfrak t)^a \rho(\mathfrak s) \rho(\mathfrak t)^b
 \rho(\mathfrak s) \rho(\mathfrak t)^a)_\mathrm{ndeg}$,
 $\sigma = a$ = 5. Prop. B.5 (3) eqn. (B.25)

 \ \color{black}

\noindent 629: (dims,levels) = $(6;48
)$,
irreps = $3_{16}^{7,0}
\hskip -1.5pt \otimes \hskip -1.5pt
2_{3}^{1,0}$,
pord$(\rho_\text{isum}(\mathfrak{t})) = 48$,

\vskip 0.7ex
\hangindent=5.5em \hangafter=1
{\white .}\hskip 1em $\rho_\text{isum}(\mathfrak{t})$ =
 $( \frac{7}{8},
\frac{7}{16},
\frac{15}{16},
\frac{5}{24},
\frac{13}{48},
\frac{37}{48} )
$,

\vskip 0.7ex
\hangindent=5.5em \hangafter=1
{\white .}\hskip 1em $\rho_\text{isum}(\mathfrak{s})$ =
($0$,
$\sqrt{\frac{1}{6}}$,
$\sqrt{\frac{1}{6}}$,
$0$,
$-\sqrt{\frac{1}{3}}$,
$-\sqrt{\frac{1}{3}}$;
$\sqrt{\frac{1}{12}}$,
$-\sqrt{\frac{1}{12}}$,
$-\sqrt{\frac{1}{3}}$,
$\sqrt{\frac{1}{6}}$,
$-\sqrt{\frac{1}{6}}$;
$\sqrt{\frac{1}{12}}$,
$-\sqrt{\frac{1}{3}}$,
$-\sqrt{\frac{1}{6}}$,
$\sqrt{\frac{1}{6}}$;
$0$,
$-\sqrt{\frac{1}{6}}$,
$-\sqrt{\frac{1}{6}}$;
$-\sqrt{\frac{1}{12}}$,
$\sqrt{\frac{1}{12}}$;
$-\sqrt{\frac{1}{12}}$)

Fail:
$\sigma(\rho(\mathfrak s)_\mathrm{ndeg}) \neq
 (\rho(\mathfrak t)^a \rho(\mathfrak s) \rho(\mathfrak t)^b
 \rho(\mathfrak s) \rho(\mathfrak t)^a)_\mathrm{ndeg}$,
 $\sigma = a$ = 5. Prop. B.5 (3) eqn. (B.25)

 \ \color{black}

 \color{blue}

\noindent 630: (dims,levels) = $(6;52
)$,
irreps = $6_{13}^{1}
\hskip -1.5pt \otimes \hskip -1.5pt
1_{4}^{1,0}$,
pord$(\rho_\text{isum}(\mathfrak{t})) = 13$,

\vskip 0.7ex
\hangindent=5.5em \hangafter=1
{\white .}\hskip 1em $\rho_\text{isum}(\mathfrak{t})$ =
 $( \frac{1}{52},
\frac{9}{52},
\frac{17}{52},
\frac{25}{52},
\frac{29}{52},
\frac{49}{52} )
$,

\vskip 0.7ex
\hangindent=5.5em \hangafter=1
{\white .}\hskip 1em $\rho_\text{isum}(\mathfrak{s})$ =
($-\frac{1}{\sqrt{13}}c^{11}_{52}
$,
$\frac{1}{\sqrt{13}}c^{7}_{52}
$,
$-\frac{1}{\sqrt{13}}c^{9}_{52}
$,
$\frac{1}{\sqrt{13}}c^{3}_{52}
$,
$\frac{1}{\sqrt{13}}c^{5}_{52}
$,
$\frac{1}{\sqrt{13}}c^{1}_{52}
$;
$-\frac{1}{\sqrt{13}}c^{5}_{52}
$,
$\frac{1}{\sqrt{13}}c^{1}_{52}
$,
$\frac{1}{\sqrt{13}}c^{9}_{52}
$,
$-\frac{1}{\sqrt{13}}c^{11}_{52}
$,
$\frac{1}{\sqrt{13}}c^{3}_{52}
$;
$\frac{1}{\sqrt{13}}c^{5}_{52}
$,
$\frac{1}{\sqrt{13}}c^{7}_{52}
$,
$-\frac{1}{\sqrt{13}}c^{3}_{52}
$,
$-\frac{1}{\sqrt{13}}c^{11}_{52}
$;
$\frac{1}{\sqrt{13}}c^{11}_{52}
$,
$\frac{1}{\sqrt{13}}c^{1}_{52}
$,
$-\frac{1}{\sqrt{13}}c^{5}_{52}
$;
$-\frac{1}{\sqrt{13}}c^{7}_{52}
$,
$-\frac{1}{\sqrt{13}}c^{9}_{52}
$;
$\frac{1}{\sqrt{13}}c^{7}_{52}
$)

Pass. 

 \ \color{black}

 \color{blue}

\noindent 631: (dims,levels) = $(6;52
)$,
irreps = $6_{13}^{2}
\hskip -1.5pt \otimes \hskip -1.5pt
1_{4}^{1,0}$,
pord$(\rho_\text{isum}(\mathfrak{t})) = 13$,

\vskip 0.7ex
\hangindent=5.5em \hangafter=1
{\white .}\hskip 1em $\rho_\text{isum}(\mathfrak{t})$ =
 $( \frac{5}{52},
\frac{21}{52},
\frac{33}{52},
\frac{37}{52},
\frac{41}{52},
\frac{45}{52} )
$,

\vskip 0.7ex
\hangindent=5.5em \hangafter=1
{\white .}\hskip 1em $\rho_\text{isum}(\mathfrak{s})$ =
($\frac{1}{\sqrt{13}}c^{3}_{52}
$,
$\frac{1}{\sqrt{13}}c^{11}_{52}
$,
$\frac{1}{\sqrt{13}}c^{7}_{52}
$,
$\frac{1}{\sqrt{13}}c^{5}_{52}
$,
$\frac{1}{\sqrt{13}}c^{1}_{52}
$,
$-\frac{1}{\sqrt{13}}c^{9}_{52}
$;
$-\frac{1}{\sqrt{13}}c^{3}_{52}
$,
$\frac{1}{\sqrt{13}}c^{9}_{52}
$,
$\frac{1}{\sqrt{13}}c^{1}_{52}
$,
$-\frac{1}{\sqrt{13}}c^{5}_{52}
$,
$\frac{1}{\sqrt{13}}c^{7}_{52}
$;
$\frac{1}{\sqrt{13}}c^{1}_{52}
$,
$-\frac{1}{\sqrt{13}}c^{3}_{52}
$,
$-\frac{1}{\sqrt{13}}c^{11}_{52}
$,
$\frac{1}{\sqrt{13}}c^{5}_{52}
$;
$\frac{1}{\sqrt{13}}c^{9}_{52}
$,
$-\frac{1}{\sqrt{13}}c^{7}_{52}
$,
$\frac{1}{\sqrt{13}}c^{11}_{52}
$;
$-\frac{1}{\sqrt{13}}c^{9}_{52}
$,
$-\frac{1}{\sqrt{13}}c^{3}_{52}
$;
$-\frac{1}{\sqrt{13}}c^{1}_{52}
$)

Pass. 

 \ \color{black}

 \color{blue}

\noindent 632: (dims,levels) = $(6;56
)$,
irreps = $3_{7}^{3}
\hskip -1.5pt \otimes \hskip -1.5pt
2_{8}^{1,0}$,
pord$(\rho_\text{isum}(\mathfrak{t})) = 28$,

\vskip 0.7ex
\hangindent=5.5em \hangafter=1
{\white .}\hskip 1em $\rho_\text{isum}(\mathfrak{t})$ =
 $( \frac{5}{56},
\frac{13}{56},
\frac{31}{56},
\frac{45}{56},
\frac{47}{56},
\frac{55}{56} )
$,

\vskip 0.7ex
\hangindent=5.5em \hangafter=1
{\white .}\hskip 1em $\rho_\text{isum}(\mathfrak{s})$ =
($\frac{1}{\sqrt{14}}c^{5}_{28}
$,
$\frac{1}{\sqrt{14}}c^{3}_{28}
$,
$\frac{1}{\sqrt{14}}c^{1}_{28}
$,
$\frac{1}{\sqrt{14}}c^{1}_{28}
$,
$-\frac{1}{\sqrt{14}}c^{5}_{28}
$,
$\frac{1}{\sqrt{14}}c^{3}_{28}
$;
$-\frac{1}{\sqrt{14}}c^{1}_{28}
$,
$\frac{1}{\sqrt{14}}c^{5}_{28}
$,
$\frac{1}{\sqrt{14}}c^{5}_{28}
$,
$-\frac{1}{\sqrt{14}}c^{3}_{28}
$,
$-\frac{1}{\sqrt{14}}c^{1}_{28}
$;
$\frac{1}{\sqrt{14}}c^{3}_{28}
$,
$-\frac{1}{\sqrt{14}}c^{3}_{28}
$,
$\frac{1}{\sqrt{14}}c^{1}_{28}
$,
$-\frac{1}{\sqrt{14}}c^{5}_{28}
$;
$-\frac{1}{\sqrt{14}}c^{3}_{28}
$,
$-\frac{1}{\sqrt{14}}c^{1}_{28}
$,
$\frac{1}{\sqrt{14}}c^{5}_{28}
$;
$-\frac{1}{\sqrt{14}}c^{5}_{28}
$,
$\frac{1}{\sqrt{14}}c^{3}_{28}
$;
$\frac{1}{\sqrt{14}}c^{1}_{28}
$)

Pass. 

 \ \color{black}

 \color{blue}

\noindent 633: (dims,levels) = $(6;56
)$,
irreps = $3_{7}^{1}
\hskip -1.5pt \otimes \hskip -1.5pt
2_{8}^{1,0}$,
pord$(\rho_\text{isum}(\mathfrak{t})) = 28$,

\vskip 0.7ex
\hangindent=5.5em \hangafter=1
{\white .}\hskip 1em $\rho_\text{isum}(\mathfrak{t})$ =
 $( \frac{15}{56},
\frac{23}{56},
\frac{29}{56},
\frac{37}{56},
\frac{39}{56},
\frac{53}{56} )
$,

\vskip 0.7ex
\hangindent=5.5em \hangafter=1
{\white .}\hskip 1em $\rho_\text{isum}(\mathfrak{s})$ =
($\frac{1}{\sqrt{14}}c^{1}_{28}
$,
$\frac{1}{\sqrt{14}}c^{3}_{28}
$,
$\frac{1}{\sqrt{14}}c^{1}_{28}
$,
$\frac{1}{\sqrt{14}}c^{3}_{28}
$,
$-\frac{1}{\sqrt{14}}c^{5}_{28}
$,
$-\frac{1}{\sqrt{14}}c^{5}_{28}
$;
$-\frac{1}{\sqrt{14}}c^{5}_{28}
$,
$\frac{1}{\sqrt{14}}c^{3}_{28}
$,
$-\frac{1}{\sqrt{14}}c^{5}_{28}
$,
$\frac{1}{\sqrt{14}}c^{1}_{28}
$,
$\frac{1}{\sqrt{14}}c^{1}_{28}
$;
$-\frac{1}{\sqrt{14}}c^{1}_{28}
$,
$-\frac{1}{\sqrt{14}}c^{3}_{28}
$,
$-\frac{1}{\sqrt{14}}c^{5}_{28}
$,
$\frac{1}{\sqrt{14}}c^{5}_{28}
$;
$\frac{1}{\sqrt{14}}c^{5}_{28}
$,
$\frac{1}{\sqrt{14}}c^{1}_{28}
$,
$-\frac{1}{\sqrt{14}}c^{1}_{28}
$;
$\frac{1}{\sqrt{14}}c^{3}_{28}
$,
$\frac{1}{\sqrt{14}}c^{3}_{28}
$;
$-\frac{1}{\sqrt{14}}c^{3}_{28}
$)

Pass. 

 \ \color{black}

 \color{blue}

\noindent 634: (dims,levels) = $(6;60
)$,
irreps = $3_{3}^{1,0}
\hskip -1.5pt \otimes \hskip -1.5pt
2_{5}^{1}
\hskip -1.5pt \otimes \hskip -1.5pt
1_{4}^{1,0}$,
pord$(\rho_\text{isum}(\mathfrak{t})) = 15$,

\vskip 0.7ex
\hangindent=5.5em \hangafter=1
{\white .}\hskip 1em $\rho_\text{isum}(\mathfrak{t})$ =
 $( \frac{1}{20},
\frac{9}{20},
\frac{7}{60},
\frac{23}{60},
\frac{43}{60},
\frac{47}{60} )
$,

\vskip 0.7ex
\hangindent=5.5em \hangafter=1
{\white .}\hskip 1em $\rho_\text{isum}(\mathfrak{s})$ =
($\frac{1}{3\sqrt{5}}c^{3}_{20}
$,
$\frac{1}{3\sqrt{5}}c^{1}_{20}
$,
$\frac{2}{3\sqrt{5}}c^{1}_{20}
$,
$\frac{2}{3\sqrt{5}}c^{3}_{20}
$,
$\frac{2}{3\sqrt{5}}c^{3}_{20}
$,
$\frac{2}{3\sqrt{5}}c^{1}_{20}
$;
$-\frac{1}{3\sqrt{5}}c^{3}_{20}
$,
$-\frac{2}{3\sqrt{5}}c^{3}_{20}
$,
$\frac{2}{3\sqrt{5}}c^{1}_{20}
$,
$\frac{2}{3\sqrt{5}}c^{1}_{20}
$,
$-\frac{2}{3\sqrt{5}}c^{3}_{20}
$;
$-\frac{1}{3\sqrt{5}}c^{3}_{20}
$,
$-\frac{2}{3\sqrt{5}}c^{1}_{20}
$,
$\frac{1}{3\sqrt{5}}c^{1}_{20}
$,
$\frac{2}{3\sqrt{5}}c^{3}_{20}
$;
$\frac{1}{3\sqrt{5}}c^{3}_{20}
$,
$-\frac{2}{3\sqrt{5}}c^{3}_{20}
$,
$\frac{1}{3\sqrt{5}}c^{1}_{20}
$;
$\frac{1}{3\sqrt{5}}c^{3}_{20}
$,
$-\frac{2}{3\sqrt{5}}c^{1}_{20}
$;
$-\frac{1}{3\sqrt{5}}c^{3}_{20}
$)

Pass. 

 \ \color{black}

\noindent 635: (dims,levels) = $(6;60
)$,
irreps = $3_{5}^{3}
\hskip -1.5pt \otimes \hskip -1.5pt
2_{4}^{1,0}
\hskip -1.5pt \otimes \hskip -1.5pt
1_{3}^{1,0}$,
pord$(\rho_\text{isum}(\mathfrak{t})) = 10$,

\vskip 0.7ex
\hangindent=5.5em \hangafter=1
{\white .}\hskip 1em $\rho_\text{isum}(\mathfrak{t})$ =
 $( \frac{1}{12},
\frac{7}{12},
\frac{11}{60},
\frac{29}{60},
\frac{41}{60},
\frac{59}{60} )
$,

\vskip 0.7ex
\hangindent=5.5em \hangafter=1
{\white .}\hskip 1em $\rho_\text{isum}(\mathfrak{s})$ =
$\mathrm{i}$($-\sqrt{\frac{1}{20}}$,
$\sqrt{\frac{3}{20}}$,
$-\sqrt{\frac{3}{10}}$,
$\sqrt{\frac{1}{10}}$,
$\sqrt{\frac{1}{10}}$,
$-\sqrt{\frac{3}{10}}$;\ \ 
$\sqrt{\frac{1}{20}}$,
$-\sqrt{\frac{1}{10}}$,
$-\sqrt{\frac{3}{10}}$,
$-\sqrt{\frac{3}{10}}$,
$-\sqrt{\frac{1}{10}}$;\ \ 
$\frac{5-\sqrt{5}}{20}$,
$-\frac{3}{2\sqrt{15}\mathrm{i}}s^{3}_{20}
$,
$\frac{3}{2\sqrt{15}}c^{1}_{5}
$,
$-\frac{5+\sqrt{5}}{20}$;\ \ 
$\frac{-5+\sqrt{5}}{20}$,
$\frac{5+\sqrt{5}}{20}$,
$\frac{3}{2\sqrt{15}}c^{1}_{5}
$;\ \ 
$\frac{-5+\sqrt{5}}{20}$,
$-\frac{3}{2\sqrt{15}\mathrm{i}}s^{3}_{20}
$;\ \ 
$\frac{5-\sqrt{5}}{20}$)

Fail:
cnd($\rho(\mathfrak s)_\mathrm{ndeg}$) = 120 does not divide
 ord($\rho(\mathfrak t)$)=60. Prop. B.4 (2)

 \ \color{black}

\noindent 636: (dims,levels) = $(6;60
)$,
irreps = $3_{5}^{1}
\hskip -1.5pt \otimes \hskip -1.5pt
2_{4}^{1,0}
\hskip -1.5pt \otimes \hskip -1.5pt
1_{3}^{1,0}$,
pord$(\rho_\text{isum}(\mathfrak{t})) = 10$,

\vskip 0.7ex
\hangindent=5.5em \hangafter=1
{\white .}\hskip 1em $\rho_\text{isum}(\mathfrak{t})$ =
 $( \frac{1}{12},
\frac{7}{12},
\frac{17}{60},
\frac{23}{60},
\frac{47}{60},
\frac{53}{60} )
$,

\vskip 0.7ex
\hangindent=5.5em \hangafter=1
{\white .}\hskip 1em $\rho_\text{isum}(\mathfrak{s})$ =
$\mathrm{i}$($\sqrt{\frac{1}{20}}$,
$\sqrt{\frac{3}{20}}$,
$\sqrt{\frac{1}{10}}$,
$-\sqrt{\frac{3}{10}}$,
$-\sqrt{\frac{3}{10}}$,
$\sqrt{\frac{1}{10}}$;\ \ 
$-\sqrt{\frac{1}{20}}$,
$\sqrt{\frac{3}{10}}$,
$\sqrt{\frac{1}{10}}$,
$\sqrt{\frac{1}{10}}$,
$\sqrt{\frac{3}{10}}$;\ \ 
$-\frac{5+\sqrt{5}}{20}$,
$-\frac{3}{2\sqrt{15}}c^{1}_{5}
$,
$\frac{3}{2\sqrt{15}\mathrm{i}}s^{3}_{20}
$,
$\frac{5-\sqrt{5}}{20}$;\ \ 
$\frac{5+\sqrt{5}}{20}$,
$\frac{-5+\sqrt{5}}{20}$,
$\frac{3}{2\sqrt{15}\mathrm{i}}s^{3}_{20}
$;\ \ 
$\frac{5+\sqrt{5}}{20}$,
$-\frac{3}{2\sqrt{15}}c^{1}_{5}
$;\ \ 
$-\frac{5+\sqrt{5}}{20}$)

Fail:
cnd($\rho(\mathfrak s)_\mathrm{ndeg}$) = 120 does not divide
 ord($\rho(\mathfrak t)$)=60. Prop. B.4 (2)

 \ \color{black}

\noindent 637: (dims,levels) = $(6;60
)$,
irreps = $3_{4}^{1,0}
\hskip -1.5pt \otimes \hskip -1.5pt
2_{5}^{1}
\hskip -1.5pt \otimes \hskip -1.5pt
1_{3}^{1,0}$,
pord$(\rho_\text{isum}(\mathfrak{t})) = 20$,

\vskip 0.7ex
\hangindent=5.5em \hangafter=1
{\white .}\hskip 1em $\rho_\text{isum}(\mathfrak{t})$ =
 $( \frac{2}{15},
\frac{8}{15},
\frac{17}{60},
\frac{23}{60},
\frac{47}{60},
\frac{53}{60} )
$,

\vskip 0.7ex
\hangindent=5.5em \hangafter=1
{\white .}\hskip 1em $\rho_\text{isum}(\mathfrak{s})$ =
$\mathrm{i}$($0$,
$0$,
$\frac{1}{\sqrt{10}}c^{1}_{20}
$,
$\frac{1}{\sqrt{10}}c^{3}_{20}
$,
$\frac{1}{\sqrt{10}}c^{1}_{20}
$,
$\frac{1}{\sqrt{10}}c^{3}_{20}
$;\ \ 
$0$,
$\frac{1}{\sqrt{10}}c^{3}_{20}
$,
$-\frac{1}{\sqrt{10}}c^{1}_{20}
$,
$\frac{1}{\sqrt{10}}c^{3}_{20}
$,
$-\frac{1}{\sqrt{10}}c^{1}_{20}
$;\ \ 
$\frac{1}{2\sqrt{5}}c^{3}_{20}
$,
$\frac{1}{2\sqrt{5}}c^{1}_{20}
$,
$-\frac{1}{2\sqrt{5}}c^{3}_{20}
$,
$-\frac{1}{2\sqrt{5}}c^{1}_{20}
$;\ \ 
$-\frac{1}{2\sqrt{5}}c^{3}_{20}
$,
$-\frac{1}{2\sqrt{5}}c^{1}_{20}
$,
$\frac{1}{2\sqrt{5}}c^{3}_{20}
$;\ \ 
$\frac{1}{2\sqrt{5}}c^{3}_{20}
$,
$\frac{1}{2\sqrt{5}}c^{1}_{20}
$;\ \ 
$-\frac{1}{2\sqrt{5}}c^{3}_{20}
$)

Fail:
cnd($\rho(\mathfrak s)_\mathrm{ndeg}$) = 40 does not divide
 ord($\rho(\mathfrak t)$)=60. Prop. B.4 (2)

 \ \color{black}

\noindent 638: (dims,levels) = $(6;60
)$,
irreps = $3_{5}^{1}
\hskip -1.5pt \otimes \hskip -1.5pt
2_{3}^{1,0}
\hskip -1.5pt \otimes \hskip -1.5pt
1_{4}^{1,0}$,
pord$(\rho_\text{isum}(\mathfrak{t})) = 15$,

\vskip 0.7ex
\hangindent=5.5em \hangafter=1
{\white .}\hskip 1em $\rho_\text{isum}(\mathfrak{t})$ =
 $( \frac{1}{4},
\frac{7}{12},
\frac{1}{20},
\frac{9}{20},
\frac{23}{60},
\frac{47}{60} )
$,

\vskip 0.7ex
\hangindent=5.5em \hangafter=1
{\white .}\hskip 1em $\rho_\text{isum}(\mathfrak{s})$ =
($\sqrt{\frac{1}{15}}$,
$\sqrt{\frac{2}{15}}$,
$\sqrt{\frac{2}{15}}$,
$\sqrt{\frac{2}{15}}$,
$\sqrt{\frac{4}{15}}$,
$\sqrt{\frac{4}{15}}$;
$-\sqrt{\frac{1}{15}}$,
$\sqrt{\frac{4}{15}}$,
$\sqrt{\frac{4}{15}}$,
$-\sqrt{\frac{2}{15}}$,
$-\sqrt{\frac{2}{15}}$;
$-\frac{1}{\sqrt{15}\mathrm{i}}s^{3}_{20}
$,
$\frac{1}{\sqrt{15}}c^{1}_{5}
$,
$-\frac{2}{\sqrt{30}\mathrm{i}}s^{3}_{20}
$,
$\frac{2}{\sqrt{30}}c^{1}_{5}
$;
$-\frac{1}{\sqrt{15}\mathrm{i}}s^{3}_{20}
$,
$\frac{2}{\sqrt{30}}c^{1}_{5}
$,
$-\frac{2}{\sqrt{30}\mathrm{i}}s^{3}_{20}
$;
$\frac{1}{\sqrt{15}\mathrm{i}}s^{3}_{20}
$,
$-\frac{1}{\sqrt{15}}c^{1}_{5}
$;
$\frac{1}{\sqrt{15}\mathrm{i}}s^{3}_{20}
$)

Fail:
cnd($\rho(\mathfrak s)_\mathrm{ndeg}$) = 120 does not divide
 ord($\rho(\mathfrak t)$)=60. Prop. B.4 (2)

 \ \color{black}

\noindent 639: (dims,levels) = $(6;60
)$,
irreps = $3_{5}^{3}
\hskip -1.5pt \otimes \hskip -1.5pt
2_{3}^{1,0}
\hskip -1.5pt \otimes \hskip -1.5pt
1_{4}^{1,0}$,
pord$(\rho_\text{isum}(\mathfrak{t})) = 15$,

\vskip 0.7ex
\hangindent=5.5em \hangafter=1
{\white .}\hskip 1em $\rho_\text{isum}(\mathfrak{t})$ =
 $( \frac{1}{4},
\frac{7}{12},
\frac{13}{20},
\frac{17}{20},
\frac{11}{60},
\frac{59}{60} )
$,

\vskip 0.7ex
\hangindent=5.5em \hangafter=1
{\white .}\hskip 1em $\rho_\text{isum}(\mathfrak{s})$ =
($-\sqrt{\frac{1}{15}}$,
$\sqrt{\frac{2}{15}}$,
$\sqrt{\frac{2}{15}}$,
$\sqrt{\frac{2}{15}}$,
$\sqrt{\frac{4}{15}}$,
$\sqrt{\frac{4}{15}}$;
$\sqrt{\frac{1}{15}}$,
$-\sqrt{\frac{4}{15}}$,
$-\sqrt{\frac{4}{15}}$,
$\sqrt{\frac{2}{15}}$,
$\sqrt{\frac{2}{15}}$;
$-\frac{1}{\sqrt{15}}c^{1}_{5}
$,
$\frac{1}{\sqrt{15}\mathrm{i}}s^{3}_{20}
$,
$\frac{2}{\sqrt{30}\mathrm{i}}s^{3}_{20}
$,
$-\frac{2}{\sqrt{30}}c^{1}_{5}
$;
$-\frac{1}{\sqrt{15}}c^{1}_{5}
$,
$-\frac{2}{\sqrt{30}}c^{1}_{5}
$,
$\frac{2}{\sqrt{30}\mathrm{i}}s^{3}_{20}
$;
$\frac{1}{\sqrt{15}}c^{1}_{5}
$,
$-\frac{1}{\sqrt{15}\mathrm{i}}s^{3}_{20}
$;
$\frac{1}{\sqrt{15}}c^{1}_{5}
$)

Fail:
cnd($\rho(\mathfrak s)_\mathrm{ndeg}$) = 120 does not divide
 ord($\rho(\mathfrak t)$)=60. Prop. B.4 (2)

 \ \color{black}

\noindent 640: (dims,levels) = $(6;60
)$,
irreps = $6_{5}^{1}
\hskip -1.5pt \otimes \hskip -1.5pt
1_{4}^{1,0}
\hskip -1.5pt \otimes \hskip -1.5pt
1_{3}^{1,0}$,
pord$(\rho_\text{isum}(\mathfrak{t})) = 5$,

\vskip 0.7ex
\hangindent=5.5em \hangafter=1
{\white .}\hskip 1em $\rho_\text{isum}(\mathfrak{t})$ =
 $( \frac{7}{12},
\frac{7}{12},
\frac{11}{60},
\frac{23}{60},
\frac{47}{60},
\frac{59}{60} )
$,

\vskip 0.7ex
\hangindent=5.5em \hangafter=1
{\white .}\hskip 1em $\rho_\text{isum}(\mathfrak{s})$ =
($-\frac{1}{5}c^{1}_{20}
$,
$-\frac{1}{5}c^{3}_{20}
$,
$\frac{1}{5}c^{1}_{20}
+\frac{1}{5}c^{3}_{20}
$,
$\frac{1}{5}c^{3}_{40}
+\frac{1}{5}c^{7}_{40}
$,
$-\frac{4}{5\sqrt{10}}c^{1}_{20}
+\frac{2}{5\sqrt{10}}c^{3}_{20}
$,
$\frac{1}{5}c^{1}_{20}
-\frac{1}{5}c^{3}_{20}
$;
$\frac{1}{5}c^{1}_{20}
$,
$-\frac{1}{5}c^{1}_{20}
+\frac{1}{5}c^{3}_{20}
$,
$\frac{4}{5\sqrt{10}}c^{1}_{20}
-\frac{2}{5\sqrt{10}}c^{3}_{20}
$,
$\frac{1}{5}c^{3}_{40}
+\frac{1}{5}c^{7}_{40}
$,
$\frac{1}{5}c^{1}_{20}
+\frac{1}{5}c^{3}_{20}
$;
$\frac{1}{5}c^{1}_{20}
$,
$\frac{4}{5\sqrt{10}}c^{1}_{20}
-\frac{2}{5\sqrt{10}}c^{3}_{20}
$,
$\frac{1}{5}c^{3}_{40}
+\frac{1}{5}c^{7}_{40}
$,
$-\frac{1}{5}c^{3}_{20}
$;
$-\frac{1}{5}c^{3}_{20}
$,
$-\frac{1}{5}c^{1}_{20}
$,
$\frac{1}{5}c^{3}_{40}
+\frac{1}{5}c^{7}_{40}
$;
$\frac{1}{5}c^{3}_{20}
$,
$-\frac{4}{5\sqrt{10}}c^{1}_{20}
+\frac{2}{5\sqrt{10}}c^{3}_{20}
$;
$-\frac{1}{5}c^{1}_{20}
$)

Fail:
cnd($\rho(\mathfrak s)_\mathrm{ndeg}$) = 40 does not divide
 ord($\rho(\mathfrak t)$)=60. Prop. B.4 (2)

 \ \color{black}

 \color{blue}

\noindent 641: (dims,levels) = $(6;60
)$,
irreps = $3_{3}^{1,0}
\hskip -1.5pt \otimes \hskip -1.5pt
2_{5}^{2}
\hskip -1.5pt \otimes \hskip -1.5pt
1_{4}^{1,0}$,
pord$(\rho_\text{isum}(\mathfrak{t})) = 15$,

\vskip 0.7ex
\hangindent=5.5em \hangafter=1
{\white .}\hskip 1em $\rho_\text{isum}(\mathfrak{t})$ =
 $( \frac{13}{20},
\frac{17}{20},
\frac{11}{60},
\frac{19}{60},
\frac{31}{60},
\frac{59}{60} )
$,

\vskip 0.7ex
\hangindent=5.5em \hangafter=1
{\white .}\hskip 1em $\rho_\text{isum}(\mathfrak{s})$ =
($-\frac{1}{3\sqrt{5}}c^{1}_{20}
$,
$\frac{1}{3\sqrt{5}}c^{3}_{20}
$,
$\frac{2}{3\sqrt{5}}c^{3}_{20}
$,
$\frac{2}{3\sqrt{5}}c^{1}_{20}
$,
$\frac{2}{3\sqrt{5}}c^{3}_{20}
$,
$\frac{2}{3\sqrt{5}}c^{1}_{20}
$;
$\frac{1}{3\sqrt{5}}c^{1}_{20}
$,
$\frac{2}{3\sqrt{5}}c^{1}_{20}
$,
$-\frac{2}{3\sqrt{5}}c^{3}_{20}
$,
$\frac{2}{3\sqrt{5}}c^{1}_{20}
$,
$-\frac{2}{3\sqrt{5}}c^{3}_{20}
$;
$\frac{1}{3\sqrt{5}}c^{1}_{20}
$,
$\frac{2}{3\sqrt{5}}c^{3}_{20}
$,
$-\frac{2}{3\sqrt{5}}c^{1}_{20}
$,
$-\frac{1}{3\sqrt{5}}c^{3}_{20}
$;
$-\frac{1}{3\sqrt{5}}c^{1}_{20}
$,
$-\frac{1}{3\sqrt{5}}c^{3}_{20}
$,
$\frac{2}{3\sqrt{5}}c^{1}_{20}
$;
$\frac{1}{3\sqrt{5}}c^{1}_{20}
$,
$\frac{2}{3\sqrt{5}}c^{3}_{20}
$;
$-\frac{1}{3\sqrt{5}}c^{1}_{20}
$)

Pass. 

 \ \color{black}

\noindent 642: (dims,levels) = $(6;60
)$,
irreps = $3_{4}^{1,0}
\hskip -1.5pt \otimes \hskip -1.5pt
2_{5}^{2}
\hskip -1.5pt \otimes \hskip -1.5pt
1_{3}^{1,0}$,
pord$(\rho_\text{isum}(\mathfrak{t})) = 20$,

\vskip 0.7ex
\hangindent=5.5em \hangafter=1
{\white .}\hskip 1em $\rho_\text{isum}(\mathfrak{t})$ =
 $( \frac{11}{15},
\frac{14}{15},
\frac{11}{60},
\frac{29}{60},
\frac{41}{60},
\frac{59}{60} )
$,

\vskip 0.7ex
\hangindent=5.5em \hangafter=1
{\white .}\hskip 1em $\rho_\text{isum}(\mathfrak{s})$ =
$\mathrm{i}$($0$,
$0$,
$\frac{1}{\sqrt{10}}c^{3}_{20}
$,
$\frac{1}{\sqrt{10}}c^{1}_{20}
$,
$\frac{1}{\sqrt{10}}c^{3}_{20}
$,
$\frac{1}{\sqrt{10}}c^{1}_{20}
$;\ \ 
$0$,
$\frac{1}{\sqrt{10}}c^{1}_{20}
$,
$-\frac{1}{\sqrt{10}}c^{3}_{20}
$,
$\frac{1}{\sqrt{10}}c^{1}_{20}
$,
$-\frac{1}{\sqrt{10}}c^{3}_{20}
$;\ \ 
$-\frac{1}{2\sqrt{5}}c^{1}_{20}
$,
$-\frac{1}{2\sqrt{5}}c^{3}_{20}
$,
$\frac{1}{2\sqrt{5}}c^{1}_{20}
$,
$\frac{1}{2\sqrt{5}}c^{3}_{20}
$;\ \ 
$\frac{1}{2\sqrt{5}}c^{1}_{20}
$,
$\frac{1}{2\sqrt{5}}c^{3}_{20}
$,
$-\frac{1}{2\sqrt{5}}c^{1}_{20}
$;\ \ 
$-\frac{1}{2\sqrt{5}}c^{1}_{20}
$,
$-\frac{1}{2\sqrt{5}}c^{3}_{20}
$;\ \ 
$\frac{1}{2\sqrt{5}}c^{1}_{20}
$)

Fail:
cnd($\rho(\mathfrak s)_\mathrm{ndeg}$) = 40 does not divide
 ord($\rho(\mathfrak t)$)=60. Prop. B.4 (2)

 \ \color{black}

\noindent 643: (dims,levels) = $(6;66
)$,
irreps = $6_{11}^{7}
\hskip -1.5pt \otimes \hskip -1.5pt
1_{3}^{1,0}
\hskip -1.5pt \otimes \hskip -1.5pt
1_{2}^{1,0}$,
pord$(\rho_\text{isum}(\mathfrak{t})) = 11$,

\vskip 0.7ex
\hangindent=5.5em \hangafter=1
{\white .}\hskip 1em $\rho_\text{isum}(\mathfrak{t})$ =
 $( \frac{5}{6},
\frac{1}{66},
\frac{25}{66},
\frac{31}{66},
\frac{37}{66},
\frac{49}{66} )
$,

\vskip 0.7ex
\hangindent=5.5em \hangafter=1
{\white .}\hskip 1em $\rho_\text{isum}(\mathfrak{s})$ =
$\mathrm{i}$($-\sqrt{\frac{1}{11}}$,
$\sqrt{\frac{2}{11}}$,
$\sqrt{\frac{2}{11}}$,
$\sqrt{\frac{2}{11}}$,
$\sqrt{\frac{2}{11}}$,
$\sqrt{\frac{2}{11}}$;\ \ 
$-\frac{1}{\sqrt{11}}c^{4}_{11}
$,
$-\frac{1}{\sqrt{11}}c^{2}_{11}
$,
$-\frac{1}{\sqrt{11}}c^{1}_{11}
$,
$-\frac{1}{\sqrt{11}}c^{3}_{11}
$,
$\frac{1}{\sqrt{11}\mathrm{i}}s^{9}_{44}
$;\ \ 
$-\frac{1}{\sqrt{11}}c^{1}_{11}
$,
$\frac{1}{\sqrt{11}\mathrm{i}}s^{9}_{44}
$,
$-\frac{1}{\sqrt{11}}c^{4}_{11}
$,
$-\frac{1}{\sqrt{11}}c^{3}_{11}
$;\ \ 
$-\frac{1}{\sqrt{11}}c^{3}_{11}
$,
$-\frac{1}{\sqrt{11}}c^{2}_{11}
$,
$-\frac{1}{\sqrt{11}}c^{4}_{11}
$;\ \ 
$\frac{1}{\sqrt{11}\mathrm{i}}s^{9}_{44}
$,
$-\frac{1}{\sqrt{11}}c^{1}_{11}
$;\ \ 
$-\frac{1}{\sqrt{11}}c^{2}_{11}
$)

Fail:
cnd($\rho(\mathfrak s)_\mathrm{ndeg}$) = 88 does not divide
 ord($\rho(\mathfrak t)$)=66. Prop. B.4 (2)

 \ \color{black}

\noindent 644: (dims,levels) = $(6;66
)$,
irreps = $6_{11}^{1}
\hskip -1.5pt \otimes \hskip -1.5pt
1_{3}^{1,0}
\hskip -1.5pt \otimes \hskip -1.5pt
1_{2}^{1,0}$,
pord$(\rho_\text{isum}(\mathfrak{t})) = 11$,

\vskip 0.7ex
\hangindent=5.5em \hangafter=1
{\white .}\hskip 1em $\rho_\text{isum}(\mathfrak{t})$ =
 $( \frac{5}{6},
\frac{7}{66},
\frac{13}{66},
\frac{19}{66},
\frac{43}{66},
\frac{61}{66} )
$,

\vskip 0.7ex
\hangindent=5.5em \hangafter=1
{\white .}\hskip 1em $\rho_\text{isum}(\mathfrak{s})$ =
$\mathrm{i}$($\sqrt{\frac{1}{11}}$,
$\sqrt{\frac{2}{11}}$,
$\sqrt{\frac{2}{11}}$,
$\sqrt{\frac{2}{11}}$,
$\sqrt{\frac{2}{11}}$,
$\sqrt{\frac{2}{11}}$;\ \ 
$-\frac{1}{\sqrt{11}\mathrm{i}}s^{9}_{44}
$,
$\frac{1}{\sqrt{11}}c^{2}_{11}
$,
$\frac{1}{\sqrt{11}}c^{4}_{11}
$,
$\frac{1}{\sqrt{11}}c^{3}_{11}
$,
$\frac{1}{\sqrt{11}}c^{1}_{11}
$;\ \ 
$\frac{1}{\sqrt{11}}c^{3}_{11}
$,
$-\frac{1}{\sqrt{11}\mathrm{i}}s^{9}_{44}
$,
$\frac{1}{\sqrt{11}}c^{1}_{11}
$,
$\frac{1}{\sqrt{11}}c^{4}_{11}
$;\ \ 
$\frac{1}{\sqrt{11}}c^{1}_{11}
$,
$\frac{1}{\sqrt{11}}c^{2}_{11}
$,
$\frac{1}{\sqrt{11}}c^{3}_{11}
$;\ \ 
$\frac{1}{\sqrt{11}}c^{4}_{11}
$,
$-\frac{1}{\sqrt{11}\mathrm{i}}s^{9}_{44}
$;\ \ 
$\frac{1}{\sqrt{11}}c^{2}_{11}
$)

Fail:
cnd($\rho(\mathfrak s)_\mathrm{ndeg}$) = 88 does not divide
 ord($\rho(\mathfrak t)$)=66. Prop. B.4 (2)

 \ \color{black}

 \color{blue}

\noindent 645: (dims,levels) = $(6;70
)$,
irreps = $3_{7}^{3}
\hskip -1.5pt \otimes \hskip -1.5pt
2_{5}^{1}
\hskip -1.5pt \otimes \hskip -1.5pt
1_{2}^{1,0}$,
pord$(\rho_\text{isum}(\mathfrak{t})) = 35$,

\vskip 0.7ex
\hangindent=5.5em \hangafter=1
{\white .}\hskip 1em $\rho_\text{isum}(\mathfrak{t})$ =
 $( \frac{1}{70},
\frac{9}{70},
\frac{11}{70},
\frac{29}{70},
\frac{39}{70},
\frac{51}{70} )
$,

\vskip 0.7ex
\hangindent=5.5em \hangafter=1
{\white .}\hskip 1em $\rho_\text{isum}(\mathfrak{s})$ =
$\mathrm{i}$($-\frac{1}{\sqrt{35}}c^{1}_{35}
+\frac{1}{\sqrt{35}}c^{6}_{35}
$,
$\frac{2}{\sqrt{35}}c^{3}_{35}
+\frac{1}{\sqrt{35}}c^{4}_{35}
+\frac{1}{\sqrt{35}}c^{10}_{35}
+\frac{1}{\sqrt{35}}c^{11}_{35}
$,
$\frac{2}{35}c^{1}_{140}
-\frac{1}{35}c^{3}_{140}
-\frac{1}{7}c^{5}_{140}
-\frac{3}{35}c^{7}_{140}
+\frac{1}{5}c^{9}_{140}
-\frac{2}{35}c^{13}_{140}
-\frac{1}{35}c^{15}_{140}
-\frac{1}{35}c^{17}_{140}
+\frac{3}{35}c^{19}_{140}
+\frac{2}{35}c^{21}_{140}
-\frac{2}{7}c^{23}_{140}
$,
$-\frac{1}{\sqrt{35}\mathrm{i}}s^{3}_{140}
-\frac{1}{\sqrt{35}\mathrm{i}}s^{17}_{140}
$,
$\frac{4}{35}c^{1}_{140}
+\frac{3}{35}c^{3}_{140}
+\frac{1}{7}c^{5}_{140}
-\frac{1}{35}c^{7}_{140}
-\frac{1}{35}c^{9}_{140}
-\frac{4}{35}c^{13}_{140}
-\frac{2}{35}c^{15}_{140}
+\frac{3}{35}c^{17}_{140}
-\frac{9}{35}c^{19}_{140}
+\frac{4}{35}c^{21}_{140}
+\frac{2}{7}c^{23}_{140}
$,
$-\frac{1}{\sqrt{35}}c^{4}_{35}
+\frac{1}{\sqrt{35}}c^{11}_{35}
$;\ \ 
$\frac{2}{35}c^{1}_{140}
-\frac{1}{35}c^{3}_{140}
-\frac{1}{7}c^{5}_{140}
-\frac{3}{35}c^{7}_{140}
+\frac{1}{5}c^{9}_{140}
-\frac{2}{35}c^{13}_{140}
-\frac{1}{35}c^{15}_{140}
-\frac{1}{35}c^{17}_{140}
+\frac{3}{35}c^{19}_{140}
+\frac{2}{35}c^{21}_{140}
-\frac{2}{7}c^{23}_{140}
$,
$\frac{1}{\sqrt{35}\mathrm{i}}s^{3}_{140}
+\frac{1}{\sqrt{35}\mathrm{i}}s^{17}_{140}
$,
$-\frac{1}{\sqrt{35}}c^{4}_{35}
+\frac{1}{\sqrt{35}}c^{11}_{35}
$,
$\frac{1}{\sqrt{35}}c^{1}_{35}
-\frac{1}{\sqrt{35}}c^{6}_{35}
$,
$-\frac{4}{35}c^{1}_{140}
-\frac{3}{35}c^{3}_{140}
-\frac{1}{7}c^{5}_{140}
+\frac{1}{35}c^{7}_{140}
+\frac{1}{35}c^{9}_{140}
+\frac{4}{35}c^{13}_{140}
+\frac{2}{35}c^{15}_{140}
-\frac{3}{35}c^{17}_{140}
+\frac{9}{35}c^{19}_{140}
-\frac{4}{35}c^{21}_{140}
-\frac{2}{7}c^{23}_{140}
$;\ \ 
$\frac{1}{\sqrt{35}}c^{4}_{35}
-\frac{1}{\sqrt{35}}c^{11}_{35}
$,
$-\frac{4}{35}c^{1}_{140}
-\frac{3}{35}c^{3}_{140}
-\frac{1}{7}c^{5}_{140}
+\frac{1}{35}c^{7}_{140}
+\frac{1}{35}c^{9}_{140}
+\frac{4}{35}c^{13}_{140}
+\frac{2}{35}c^{15}_{140}
-\frac{3}{35}c^{17}_{140}
+\frac{9}{35}c^{19}_{140}
-\frac{4}{35}c^{21}_{140}
-\frac{2}{7}c^{23}_{140}
$,
$-\frac{2}{\sqrt{35}}c^{3}_{35}
-\frac{1}{\sqrt{35}}c^{4}_{35}
-\frac{1}{\sqrt{35}}c^{10}_{35}
-\frac{1}{\sqrt{35}}c^{11}_{35}
$,
$-\frac{1}{\sqrt{35}}c^{1}_{35}
+\frac{1}{\sqrt{35}}c^{6}_{35}
$;\ \ 
$\frac{1}{\sqrt{35}}c^{1}_{35}
-\frac{1}{\sqrt{35}}c^{6}_{35}
$,
$\frac{2}{35}c^{1}_{140}
-\frac{1}{35}c^{3}_{140}
-\frac{1}{7}c^{5}_{140}
-\frac{3}{35}c^{7}_{140}
+\frac{1}{5}c^{9}_{140}
-\frac{2}{35}c^{13}_{140}
-\frac{1}{35}c^{15}_{140}
-\frac{1}{35}c^{17}_{140}
+\frac{3}{35}c^{19}_{140}
+\frac{2}{35}c^{21}_{140}
-\frac{2}{7}c^{23}_{140}
$,
$-\frac{2}{\sqrt{35}}c^{3}_{35}
-\frac{1}{\sqrt{35}}c^{4}_{35}
-\frac{1}{\sqrt{35}}c^{10}_{35}
-\frac{1}{\sqrt{35}}c^{11}_{35}
$;\ \ 
$-\frac{1}{\sqrt{35}}c^{4}_{35}
+\frac{1}{\sqrt{35}}c^{11}_{35}
$,
$\frac{1}{\sqrt{35}\mathrm{i}}s^{3}_{140}
+\frac{1}{\sqrt{35}\mathrm{i}}s^{17}_{140}
$;\ \ 
$-\frac{2}{35}c^{1}_{140}
+\frac{1}{35}c^{3}_{140}
+\frac{1}{7}c^{5}_{140}
+\frac{3}{35}c^{7}_{140}
-\frac{1}{5}c^{9}_{140}
+\frac{2}{35}c^{13}_{140}
+\frac{1}{35}c^{15}_{140}
+\frac{1}{35}c^{17}_{140}
-\frac{3}{35}c^{19}_{140}
-\frac{2}{35}c^{21}_{140}
+\frac{2}{7}c^{23}_{140}
$)

Pass. 

 \ \color{black}

 \color{blue}

\noindent 646: (dims,levels) = $(6;70
)$,
irreps = $3_{7}^{1}
\hskip -1.5pt \otimes \hskip -1.5pt
2_{5}^{2}
\hskip -1.5pt \otimes \hskip -1.5pt
1_{2}^{1,0}$,
pord$(\rho_\text{isum}(\mathfrak{t})) = 35$,

\vskip 0.7ex
\hangindent=5.5em \hangafter=1
{\white .}\hskip 1em $\rho_\text{isum}(\mathfrak{t})$ =
 $( \frac{3}{70},
\frac{13}{70},
\frac{17}{70},
\frac{27}{70},
\frac{33}{70},
\frac{47}{70} )
$,

\vskip 0.7ex
\hangindent=5.5em \hangafter=1
{\white .}\hskip 1em $\rho_\text{isum}(\mathfrak{s})$ =
$\mathrm{i}$($-\frac{2}{\sqrt{35}}c^{3}_{35}
-\frac{1}{\sqrt{35}}c^{4}_{35}
-\frac{1}{\sqrt{35}}c^{10}_{35}
-\frac{1}{\sqrt{35}}c^{11}_{35}
$,
$\frac{4}{35}c^{1}_{140}
+\frac{3}{35}c^{3}_{140}
+\frac{1}{7}c^{5}_{140}
-\frac{1}{35}c^{7}_{140}
-\frac{1}{35}c^{9}_{140}
-\frac{4}{35}c^{13}_{140}
-\frac{2}{35}c^{15}_{140}
+\frac{3}{35}c^{17}_{140}
-\frac{9}{35}c^{19}_{140}
+\frac{4}{35}c^{21}_{140}
+\frac{2}{7}c^{23}_{140}
$,
$-\frac{1}{\sqrt{35}}c^{4}_{35}
+\frac{1}{\sqrt{35}}c^{11}_{35}
$,
$\frac{2}{35}c^{1}_{140}
-\frac{1}{35}c^{3}_{140}
-\frac{1}{7}c^{5}_{140}
-\frac{3}{35}c^{7}_{140}
+\frac{1}{5}c^{9}_{140}
-\frac{2}{35}c^{13}_{140}
-\frac{1}{35}c^{15}_{140}
-\frac{1}{35}c^{17}_{140}
+\frac{3}{35}c^{19}_{140}
+\frac{2}{35}c^{21}_{140}
-\frac{2}{7}c^{23}_{140}
$,
$-\frac{1}{\sqrt{35}\mathrm{i}}s^{3}_{140}
-\frac{1}{\sqrt{35}\mathrm{i}}s^{17}_{140}
$,
$\frac{1}{\sqrt{35}}c^{1}_{35}
-\frac{1}{\sqrt{35}}c^{6}_{35}
$;\ \ 
$\frac{1}{\sqrt{35}\mathrm{i}}s^{3}_{140}
+\frac{1}{\sqrt{35}\mathrm{i}}s^{17}_{140}
$,
$-\frac{2}{35}c^{1}_{140}
+\frac{1}{35}c^{3}_{140}
+\frac{1}{7}c^{5}_{140}
+\frac{3}{35}c^{7}_{140}
-\frac{1}{5}c^{9}_{140}
+\frac{2}{35}c^{13}_{140}
+\frac{1}{35}c^{15}_{140}
+\frac{1}{35}c^{17}_{140}
-\frac{3}{35}c^{19}_{140}
-\frac{2}{35}c^{21}_{140}
+\frac{2}{7}c^{23}_{140}
$,
$-\frac{1}{\sqrt{35}}c^{1}_{35}
+\frac{1}{\sqrt{35}}c^{6}_{35}
$,
$-\frac{2}{\sqrt{35}}c^{3}_{35}
-\frac{1}{\sqrt{35}}c^{4}_{35}
-\frac{1}{\sqrt{35}}c^{10}_{35}
-\frac{1}{\sqrt{35}}c^{11}_{35}
$,
$\frac{1}{\sqrt{35}}c^{4}_{35}
-\frac{1}{\sqrt{35}}c^{11}_{35}
$;\ \ 
$\frac{2}{\sqrt{35}}c^{3}_{35}
+\frac{1}{\sqrt{35}}c^{4}_{35}
+\frac{1}{\sqrt{35}}c^{10}_{35}
+\frac{1}{\sqrt{35}}c^{11}_{35}
$,
$\frac{4}{35}c^{1}_{140}
+\frac{3}{35}c^{3}_{140}
+\frac{1}{7}c^{5}_{140}
-\frac{1}{35}c^{7}_{140}
-\frac{1}{35}c^{9}_{140}
-\frac{4}{35}c^{13}_{140}
-\frac{2}{35}c^{15}_{140}
+\frac{3}{35}c^{17}_{140}
-\frac{9}{35}c^{19}_{140}
+\frac{4}{35}c^{21}_{140}
+\frac{2}{7}c^{23}_{140}
$,
$-\frac{1}{\sqrt{35}}c^{1}_{35}
+\frac{1}{\sqrt{35}}c^{6}_{35}
$,
$-\frac{1}{\sqrt{35}\mathrm{i}}s^{3}_{140}
-\frac{1}{\sqrt{35}\mathrm{i}}s^{17}_{140}
$;\ \ 
$-\frac{1}{\sqrt{35}\mathrm{i}}s^{3}_{140}
-\frac{1}{\sqrt{35}\mathrm{i}}s^{17}_{140}
$,
$\frac{1}{\sqrt{35}}c^{4}_{35}
-\frac{1}{\sqrt{35}}c^{11}_{35}
$,
$\frac{2}{\sqrt{35}}c^{3}_{35}
+\frac{1}{\sqrt{35}}c^{4}_{35}
+\frac{1}{\sqrt{35}}c^{10}_{35}
+\frac{1}{\sqrt{35}}c^{11}_{35}
$;\ \ 
$-\frac{4}{35}c^{1}_{140}
-\frac{3}{35}c^{3}_{140}
-\frac{1}{7}c^{5}_{140}
+\frac{1}{35}c^{7}_{140}
+\frac{1}{35}c^{9}_{140}
+\frac{4}{35}c^{13}_{140}
+\frac{2}{35}c^{15}_{140}
-\frac{3}{35}c^{17}_{140}
+\frac{9}{35}c^{19}_{140}
-\frac{4}{35}c^{21}_{140}
-\frac{2}{7}c^{23}_{140}
$,
$-\frac{2}{35}c^{1}_{140}
+\frac{1}{35}c^{3}_{140}
+\frac{1}{7}c^{5}_{140}
+\frac{3}{35}c^{7}_{140}
-\frac{1}{5}c^{9}_{140}
+\frac{2}{35}c^{13}_{140}
+\frac{1}{35}c^{15}_{140}
+\frac{1}{35}c^{17}_{140}
-\frac{3}{35}c^{19}_{140}
-\frac{2}{35}c^{21}_{140}
+\frac{2}{7}c^{23}_{140}
$;\ \ 
$\frac{4}{35}c^{1}_{140}
+\frac{3}{35}c^{3}_{140}
+\frac{1}{7}c^{5}_{140}
-\frac{1}{35}c^{7}_{140}
-\frac{1}{35}c^{9}_{140}
-\frac{4}{35}c^{13}_{140}
-\frac{2}{35}c^{15}_{140}
+\frac{3}{35}c^{17}_{140}
-\frac{9}{35}c^{19}_{140}
+\frac{4}{35}c^{21}_{140}
+\frac{2}{7}c^{23}_{140}
$)

Pass. 

 \ \color{black}

 \color{blue}

\noindent 647: (dims,levels) = $(6;70
)$,
irreps = $3_{7}^{1}
\hskip -1.5pt \otimes \hskip -1.5pt
2_{5}^{1}
\hskip -1.5pt \otimes \hskip -1.5pt
1_{2}^{1,0}$,
pord$(\rho_\text{isum}(\mathfrak{t})) = 35$,

\vskip 0.7ex
\hangindent=5.5em \hangafter=1
{\white .}\hskip 1em $\rho_\text{isum}(\mathfrak{t})$ =
 $( \frac{19}{70},
\frac{31}{70},
\frac{41}{70},
\frac{59}{70},
\frac{61}{70},
\frac{69}{70} )
$,

\vskip 0.7ex
\hangindent=5.5em \hangafter=1
{\white .}\hskip 1em $\rho_\text{isum}(\mathfrak{s})$ =
$\mathrm{i}$($\frac{2}{35}c^{1}_{140}
-\frac{1}{35}c^{3}_{140}
-\frac{1}{7}c^{5}_{140}
-\frac{3}{35}c^{7}_{140}
+\frac{1}{5}c^{9}_{140}
-\frac{2}{35}c^{13}_{140}
-\frac{1}{35}c^{15}_{140}
-\frac{1}{35}c^{17}_{140}
+\frac{3}{35}c^{19}_{140}
+\frac{2}{35}c^{21}_{140}
-\frac{2}{7}c^{23}_{140}
$,
$-\frac{1}{\sqrt{35}\mathrm{i}}s^{3}_{140}
-\frac{1}{\sqrt{35}\mathrm{i}}s^{17}_{140}
$,
$\frac{2}{\sqrt{35}}c^{3}_{35}
+\frac{1}{\sqrt{35}}c^{4}_{35}
+\frac{1}{\sqrt{35}}c^{10}_{35}
+\frac{1}{\sqrt{35}}c^{11}_{35}
$,
$\frac{1}{\sqrt{35}}c^{1}_{35}
-\frac{1}{\sqrt{35}}c^{6}_{35}
$,
$\frac{4}{35}c^{1}_{140}
+\frac{3}{35}c^{3}_{140}
+\frac{1}{7}c^{5}_{140}
-\frac{1}{35}c^{7}_{140}
-\frac{1}{35}c^{9}_{140}
-\frac{4}{35}c^{13}_{140}
-\frac{2}{35}c^{15}_{140}
+\frac{3}{35}c^{17}_{140}
-\frac{9}{35}c^{19}_{140}
+\frac{4}{35}c^{21}_{140}
+\frac{2}{7}c^{23}_{140}
$,
$-\frac{1}{\sqrt{35}}c^{4}_{35}
+\frac{1}{\sqrt{35}}c^{11}_{35}
$;\ \ 
$\frac{1}{\sqrt{35}}c^{4}_{35}
-\frac{1}{\sqrt{35}}c^{11}_{35}
$,
$-\frac{2}{35}c^{1}_{140}
+\frac{1}{35}c^{3}_{140}
+\frac{1}{7}c^{5}_{140}
+\frac{3}{35}c^{7}_{140}
-\frac{1}{5}c^{9}_{140}
+\frac{2}{35}c^{13}_{140}
+\frac{1}{35}c^{15}_{140}
+\frac{1}{35}c^{17}_{140}
-\frac{3}{35}c^{19}_{140}
-\frac{2}{35}c^{21}_{140}
+\frac{2}{7}c^{23}_{140}
$,
$\frac{2}{\sqrt{35}}c^{3}_{35}
+\frac{1}{\sqrt{35}}c^{4}_{35}
+\frac{1}{\sqrt{35}}c^{10}_{35}
+\frac{1}{\sqrt{35}}c^{11}_{35}
$,
$-\frac{1}{\sqrt{35}}c^{1}_{35}
+\frac{1}{\sqrt{35}}c^{6}_{35}
$,
$\frac{4}{35}c^{1}_{140}
+\frac{3}{35}c^{3}_{140}
+\frac{1}{7}c^{5}_{140}
-\frac{1}{35}c^{7}_{140}
-\frac{1}{35}c^{9}_{140}
-\frac{4}{35}c^{13}_{140}
-\frac{2}{35}c^{15}_{140}
+\frac{3}{35}c^{17}_{140}
-\frac{9}{35}c^{19}_{140}
+\frac{4}{35}c^{21}_{140}
+\frac{2}{7}c^{23}_{140}
$;\ \ 
$-\frac{1}{\sqrt{35}}c^{1}_{35}
+\frac{1}{\sqrt{35}}c^{6}_{35}
$,
$\frac{4}{35}c^{1}_{140}
+\frac{3}{35}c^{3}_{140}
+\frac{1}{7}c^{5}_{140}
-\frac{1}{35}c^{7}_{140}
-\frac{1}{35}c^{9}_{140}
-\frac{4}{35}c^{13}_{140}
-\frac{2}{35}c^{15}_{140}
+\frac{3}{35}c^{17}_{140}
-\frac{9}{35}c^{19}_{140}
+\frac{4}{35}c^{21}_{140}
+\frac{2}{7}c^{23}_{140}
$,
$\frac{1}{\sqrt{35}}c^{4}_{35}
-\frac{1}{\sqrt{35}}c^{11}_{35}
$,
$-\frac{1}{\sqrt{35}\mathrm{i}}s^{3}_{140}
-\frac{1}{\sqrt{35}\mathrm{i}}s^{17}_{140}
$;\ \ 
$-\frac{1}{\sqrt{35}}c^{4}_{35}
+\frac{1}{\sqrt{35}}c^{11}_{35}
$,
$-\frac{1}{\sqrt{35}\mathrm{i}}s^{3}_{140}
-\frac{1}{\sqrt{35}\mathrm{i}}s^{17}_{140}
$,
$\frac{2}{35}c^{1}_{140}
-\frac{1}{35}c^{3}_{140}
-\frac{1}{7}c^{5}_{140}
-\frac{3}{35}c^{7}_{140}
+\frac{1}{5}c^{9}_{140}
-\frac{2}{35}c^{13}_{140}
-\frac{1}{35}c^{15}_{140}
-\frac{1}{35}c^{17}_{140}
+\frac{3}{35}c^{19}_{140}
+\frac{2}{35}c^{21}_{140}
-\frac{2}{7}c^{23}_{140}
$;\ \ 
$-\frac{2}{35}c^{1}_{140}
+\frac{1}{35}c^{3}_{140}
+\frac{1}{7}c^{5}_{140}
+\frac{3}{35}c^{7}_{140}
-\frac{1}{5}c^{9}_{140}
+\frac{2}{35}c^{13}_{140}
+\frac{1}{35}c^{15}_{140}
+\frac{1}{35}c^{17}_{140}
-\frac{3}{35}c^{19}_{140}
-\frac{2}{35}c^{21}_{140}
+\frac{2}{7}c^{23}_{140}
$,
$\frac{2}{\sqrt{35}}c^{3}_{35}
+\frac{1}{\sqrt{35}}c^{4}_{35}
+\frac{1}{\sqrt{35}}c^{10}_{35}
+\frac{1}{\sqrt{35}}c^{11}_{35}
$;\ \ 
$\frac{1}{\sqrt{35}}c^{1}_{35}
-\frac{1}{\sqrt{35}}c^{6}_{35}
$)

Pass. 

 \ \color{black}

 \color{blue}

\noindent 648: (dims,levels) = $(6;70
)$,
irreps = $3_{7}^{3}
\hskip -1.5pt \otimes \hskip -1.5pt
2_{5}^{2}
\hskip -1.5pt \otimes \hskip -1.5pt
1_{2}^{1,0}$,
pord$(\rho_\text{isum}(\mathfrak{t})) = 35$,

\vskip 0.7ex
\hangindent=5.5em \hangafter=1
{\white .}\hskip 1em $\rho_\text{isum}(\mathfrak{t})$ =
 $( \frac{23}{70},
\frac{37}{70},
\frac{43}{70},
\frac{53}{70},
\frac{57}{70},
\frac{67}{70} )
$,

\vskip 0.7ex
\hangindent=5.5em \hangafter=1
{\white .}\hskip 1em $\rho_\text{isum}(\mathfrak{s})$ =
$\mathrm{i}$($-\frac{4}{35}c^{1}_{140}
-\frac{3}{35}c^{3}_{140}
-\frac{1}{7}c^{5}_{140}
+\frac{1}{35}c^{7}_{140}
+\frac{1}{35}c^{9}_{140}
+\frac{4}{35}c^{13}_{140}
+\frac{2}{35}c^{15}_{140}
-\frac{3}{35}c^{17}_{140}
+\frac{9}{35}c^{19}_{140}
-\frac{4}{35}c^{21}_{140}
-\frac{2}{7}c^{23}_{140}
$,
$\frac{2}{35}c^{1}_{140}
-\frac{1}{35}c^{3}_{140}
-\frac{1}{7}c^{5}_{140}
-\frac{3}{35}c^{7}_{140}
+\frac{1}{5}c^{9}_{140}
-\frac{2}{35}c^{13}_{140}
-\frac{1}{35}c^{15}_{140}
-\frac{1}{35}c^{17}_{140}
+\frac{3}{35}c^{19}_{140}
+\frac{2}{35}c^{21}_{140}
-\frac{2}{7}c^{23}_{140}
$,
$\frac{2}{\sqrt{35}}c^{3}_{35}
+\frac{1}{\sqrt{35}}c^{4}_{35}
+\frac{1}{\sqrt{35}}c^{10}_{35}
+\frac{1}{\sqrt{35}}c^{11}_{35}
$,
$-\frac{1}{\sqrt{35}\mathrm{i}}s^{3}_{140}
-\frac{1}{\sqrt{35}\mathrm{i}}s^{17}_{140}
$,
$-\frac{1}{\sqrt{35}}c^{4}_{35}
+\frac{1}{\sqrt{35}}c^{11}_{35}
$,
$\frac{1}{\sqrt{35}}c^{1}_{35}
-\frac{1}{\sqrt{35}}c^{6}_{35}
$;\ \ 
$\frac{4}{35}c^{1}_{140}
+\frac{3}{35}c^{3}_{140}
+\frac{1}{7}c^{5}_{140}
-\frac{1}{35}c^{7}_{140}
-\frac{1}{35}c^{9}_{140}
-\frac{4}{35}c^{13}_{140}
-\frac{2}{35}c^{15}_{140}
+\frac{3}{35}c^{17}_{140}
-\frac{9}{35}c^{19}_{140}
+\frac{4}{35}c^{21}_{140}
+\frac{2}{7}c^{23}_{140}
$,
$\frac{1}{\sqrt{35}}c^{4}_{35}
-\frac{1}{\sqrt{35}}c^{11}_{35}
$,
$-\frac{1}{\sqrt{35}}c^{1}_{35}
+\frac{1}{\sqrt{35}}c^{6}_{35}
$,
$\frac{2}{\sqrt{35}}c^{3}_{35}
+\frac{1}{\sqrt{35}}c^{4}_{35}
+\frac{1}{\sqrt{35}}c^{10}_{35}
+\frac{1}{\sqrt{35}}c^{11}_{35}
$,
$-\frac{1}{\sqrt{35}\mathrm{i}}s^{3}_{140}
-\frac{1}{\sqrt{35}\mathrm{i}}s^{17}_{140}
$;\ \ 
$\frac{1}{\sqrt{35}\mathrm{i}}s^{3}_{140}
+\frac{1}{\sqrt{35}\mathrm{i}}s^{17}_{140}
$,
$-\frac{4}{35}c^{1}_{140}
-\frac{3}{35}c^{3}_{140}
-\frac{1}{7}c^{5}_{140}
+\frac{1}{35}c^{7}_{140}
+\frac{1}{35}c^{9}_{140}
+\frac{4}{35}c^{13}_{140}
+\frac{2}{35}c^{15}_{140}
-\frac{3}{35}c^{17}_{140}
+\frac{9}{35}c^{19}_{140}
-\frac{4}{35}c^{21}_{140}
-\frac{2}{7}c^{23}_{140}
$,
$-\frac{1}{\sqrt{35}}c^{1}_{35}
+\frac{1}{\sqrt{35}}c^{6}_{35}
$,
$-\frac{2}{35}c^{1}_{140}
+\frac{1}{35}c^{3}_{140}
+\frac{1}{7}c^{5}_{140}
+\frac{3}{35}c^{7}_{140}
-\frac{1}{5}c^{9}_{140}
+\frac{2}{35}c^{13}_{140}
+\frac{1}{35}c^{15}_{140}
+\frac{1}{35}c^{17}_{140}
-\frac{3}{35}c^{19}_{140}
-\frac{2}{35}c^{21}_{140}
+\frac{2}{7}c^{23}_{140}
$;\ \ 
$-\frac{2}{\sqrt{35}}c^{3}_{35}
-\frac{1}{\sqrt{35}}c^{4}_{35}
-\frac{1}{\sqrt{35}}c^{10}_{35}
-\frac{1}{\sqrt{35}}c^{11}_{35}
$,
$-\frac{2}{35}c^{1}_{140}
+\frac{1}{35}c^{3}_{140}
+\frac{1}{7}c^{5}_{140}
+\frac{3}{35}c^{7}_{140}
-\frac{1}{5}c^{9}_{140}
+\frac{2}{35}c^{13}_{140}
+\frac{1}{35}c^{15}_{140}
+\frac{1}{35}c^{17}_{140}
-\frac{3}{35}c^{19}_{140}
-\frac{2}{35}c^{21}_{140}
+\frac{2}{7}c^{23}_{140}
$,
$\frac{1}{\sqrt{35}}c^{4}_{35}
-\frac{1}{\sqrt{35}}c^{11}_{35}
$;\ \ 
$-\frac{1}{\sqrt{35}\mathrm{i}}s^{3}_{140}
-\frac{1}{\sqrt{35}\mathrm{i}}s^{17}_{140}
$,
$\frac{4}{35}c^{1}_{140}
+\frac{3}{35}c^{3}_{140}
+\frac{1}{7}c^{5}_{140}
-\frac{1}{35}c^{7}_{140}
-\frac{1}{35}c^{9}_{140}
-\frac{4}{35}c^{13}_{140}
-\frac{2}{35}c^{15}_{140}
+\frac{3}{35}c^{17}_{140}
-\frac{9}{35}c^{19}_{140}
+\frac{4}{35}c^{21}_{140}
+\frac{2}{7}c^{23}_{140}
$;\ \ 
$\frac{2}{\sqrt{35}}c^{3}_{35}
+\frac{1}{\sqrt{35}}c^{4}_{35}
+\frac{1}{\sqrt{35}}c^{10}_{35}
+\frac{1}{\sqrt{35}}c^{11}_{35}
$)

Pass. 

 \ \color{black}

 \color{blue}

\noindent 649: (dims,levels) = $(6;78
)$,
irreps = $6_{13}^{1}
\hskip -1.5pt \otimes \hskip -1.5pt
1_{3}^{1,0}
\hskip -1.5pt \otimes \hskip -1.5pt
1_{2}^{1,0}$,
pord$(\rho_\text{isum}(\mathfrak{t})) = 13$,

\vskip 0.7ex
\hangindent=5.5em \hangafter=1
{\white .}\hskip 1em $\rho_\text{isum}(\mathfrak{t})$ =
 $( \frac{5}{78},
\frac{11}{78},
\frac{41}{78},
\frac{47}{78},
\frac{59}{78},
\frac{71}{78} )
$,

\vskip 0.7ex
\hangindent=5.5em \hangafter=1
{\white .}\hskip 1em $\rho_\text{isum}(\mathfrak{s})$ =
$\mathrm{i}$($\frac{1}{\sqrt{13}}c^{11}_{52}
$,
$-\frac{1}{\sqrt{13}}c^{1}_{52}
$,
$-\frac{1}{\sqrt{13}}c^{5}_{52}
$,
$\frac{1}{\sqrt{13}}c^{3}_{52}
$,
$\frac{1}{\sqrt{13}}c^{9}_{52}
$,
$\frac{1}{\sqrt{13}}c^{7}_{52}
$;\ \ 
$-\frac{1}{\sqrt{13}}c^{7}_{52}
$,
$\frac{1}{\sqrt{13}}c^{9}_{52}
$,
$-\frac{1}{\sqrt{13}}c^{5}_{52}
$,
$\frac{1}{\sqrt{13}}c^{11}_{52}
$,
$\frac{1}{\sqrt{13}}c^{3}_{52}
$;\ \ 
$\frac{1}{\sqrt{13}}c^{7}_{52}
$,
$\frac{1}{\sqrt{13}}c^{1}_{52}
$,
$\frac{1}{\sqrt{13}}c^{3}_{52}
$,
$-\frac{1}{\sqrt{13}}c^{11}_{52}
$;\ \ 
$-\frac{1}{\sqrt{13}}c^{11}_{52}
$,
$\frac{1}{\sqrt{13}}c^{7}_{52}
$,
$-\frac{1}{\sqrt{13}}c^{9}_{52}
$;\ \ 
$-\frac{1}{\sqrt{13}}c^{5}_{52}
$,
$\frac{1}{\sqrt{13}}c^{1}_{52}
$;\ \ 
$\frac{1}{\sqrt{13}}c^{5}_{52}
$)

Pass. 

 \ \color{black}

 \color{blue}

\noindent 650: (dims,levels) = $(6;78
)$,
irreps = $6_{13}^{2}
\hskip -1.5pt \otimes \hskip -1.5pt
1_{3}^{1,0}
\hskip -1.5pt \otimes \hskip -1.5pt
1_{2}^{1,0}$,
pord$(\rho_\text{isum}(\mathfrak{t})) = 13$,

\vskip 0.7ex
\hangindent=5.5em \hangafter=1
{\white .}\hskip 1em $\rho_\text{isum}(\mathfrak{t})$ =
 $( \frac{17}{78},
\frac{23}{78},
\frac{29}{78},
\frac{35}{78},
\frac{53}{78},
\frac{77}{78} )
$,

\vskip 0.7ex
\hangindent=5.5em \hangafter=1
{\white .}\hskip 1em $\rho_\text{isum}(\mathfrak{s})$ =
$\mathrm{i}$($\frac{1}{\sqrt{13}}c^{1}_{52}
$,
$\frac{1}{\sqrt{13}}c^{3}_{52}
$,
$\frac{1}{\sqrt{13}}c^{11}_{52}
$,
$-\frac{1}{\sqrt{13}}c^{5}_{52}
$,
$\frac{1}{\sqrt{13}}c^{7}_{52}
$,
$\frac{1}{\sqrt{13}}c^{9}_{52}
$;\ \ 
$\frac{1}{\sqrt{13}}c^{9}_{52}
$,
$-\frac{1}{\sqrt{13}}c^{7}_{52}
$,
$\frac{1}{\sqrt{13}}c^{11}_{52}
$,
$-\frac{1}{\sqrt{13}}c^{5}_{52}
$,
$-\frac{1}{\sqrt{13}}c^{1}_{52}
$;\ \ 
$-\frac{1}{\sqrt{13}}c^{9}_{52}
$,
$-\frac{1}{\sqrt{13}}c^{3}_{52}
$,
$-\frac{1}{\sqrt{13}}c^{1}_{52}
$,
$\frac{1}{\sqrt{13}}c^{5}_{52}
$;\ \ 
$-\frac{1}{\sqrt{13}}c^{1}_{52}
$,
$\frac{1}{\sqrt{13}}c^{9}_{52}
$,
$-\frac{1}{\sqrt{13}}c^{7}_{52}
$;\ \ 
$\frac{1}{\sqrt{13}}c^{3}_{52}
$,
$\frac{1}{\sqrt{13}}c^{11}_{52}
$;\ \ 
$-\frac{1}{\sqrt{13}}c^{3}_{52}
$)

Pass. 

 \ \color{black}

 \color{blue}

\noindent 651: (dims,levels) = $(6;80
)$,
irreps = $3_{16}^{5,0}
\hskip -1.5pt \otimes \hskip -1.5pt
2_{5}^{2}$,
pord$(\rho_\text{isum}(\mathfrak{t})) = 80$,

\vskip 0.7ex
\hangindent=5.5em \hangafter=1
{\white .}\hskip 1em $\rho_\text{isum}(\mathfrak{t})$ =
 $( \frac{1}{40},
\frac{9}{40},
\frac{17}{80},
\frac{33}{80},
\frac{57}{80},
\frac{73}{80} )
$,

\vskip 0.7ex
\hangindent=5.5em \hangafter=1
{\white .}\hskip 1em $\rho_\text{isum}(\mathfrak{s})$ =
($0$,
$0$,
$\frac{1}{\sqrt{10}}c^{1}_{20}
$,
$\frac{1}{\sqrt{10}}c^{3}_{20}
$,
$\frac{1}{\sqrt{10}}c^{1}_{20}
$,
$\frac{1}{\sqrt{10}}c^{3}_{20}
$;
$0$,
$\frac{1}{\sqrt{10}}c^{3}_{20}
$,
$-\frac{1}{\sqrt{10}}c^{1}_{20}
$,
$\frac{1}{\sqrt{10}}c^{3}_{20}
$,
$-\frac{1}{\sqrt{10}}c^{1}_{20}
$;
$-\frac{1}{2\sqrt{5}}c^{1}_{20}
$,
$-\frac{1}{2\sqrt{5}}c^{3}_{20}
$,
$\frac{1}{2\sqrt{5}}c^{1}_{20}
$,
$\frac{1}{2\sqrt{5}}c^{3}_{20}
$;
$\frac{1}{2\sqrt{5}}c^{1}_{20}
$,
$\frac{1}{2\sqrt{5}}c^{3}_{20}
$,
$-\frac{1}{2\sqrt{5}}c^{1}_{20}
$;
$-\frac{1}{2\sqrt{5}}c^{1}_{20}
$,
$-\frac{1}{2\sqrt{5}}c^{3}_{20}
$;
$\frac{1}{2\sqrt{5}}c^{1}_{20}
$)

Pass. 

 \ \color{black}

 \color{blue}

\noindent 652: (dims,levels) = $(6;80
)$,
irreps = $3_{16}^{7,0}
\hskip -1.5pt \otimes \hskip -1.5pt
2_{5}^{1}$,
pord$(\rho_\text{isum}(\mathfrak{t})) = 80$,

\vskip 0.7ex
\hangindent=5.5em \hangafter=1
{\white .}\hskip 1em $\rho_\text{isum}(\mathfrak{t})$ =
 $( \frac{3}{40},
\frac{27}{40},
\frac{11}{80},
\frac{19}{80},
\frac{51}{80},
\frac{59}{80} )
$,

\vskip 0.7ex
\hangindent=5.5em \hangafter=1
{\white .}\hskip 1em $\rho_\text{isum}(\mathfrak{s})$ =
($0$,
$0$,
$\frac{1}{\sqrt{10}}c^{3}_{20}
$,
$\frac{1}{\sqrt{10}}c^{1}_{20}
$,
$\frac{1}{\sqrt{10}}c^{3}_{20}
$,
$\frac{1}{\sqrt{10}}c^{1}_{20}
$;
$0$,
$\frac{1}{\sqrt{10}}c^{1}_{20}
$,
$-\frac{1}{\sqrt{10}}c^{3}_{20}
$,
$\frac{1}{\sqrt{10}}c^{1}_{20}
$,
$-\frac{1}{\sqrt{10}}c^{3}_{20}
$;
$\frac{1}{2\sqrt{5}}c^{3}_{20}
$,
$-\frac{1}{2\sqrt{5}}c^{1}_{20}
$,
$-\frac{1}{2\sqrt{5}}c^{3}_{20}
$,
$\frac{1}{2\sqrt{5}}c^{1}_{20}
$;
$-\frac{1}{2\sqrt{5}}c^{3}_{20}
$,
$\frac{1}{2\sqrt{5}}c^{1}_{20}
$,
$\frac{1}{2\sqrt{5}}c^{3}_{20}
$;
$\frac{1}{2\sqrt{5}}c^{3}_{20}
$,
$-\frac{1}{2\sqrt{5}}c^{1}_{20}
$;
$-\frac{1}{2\sqrt{5}}c^{3}_{20}
$)

Pass. 

 \ \color{black}

 \color{blue}

\noindent 653: (dims,levels) = $(6;80
)$,
irreps = $3_{16}^{3,0}
\hskip -1.5pt \otimes \hskip -1.5pt
2_{5}^{1}$,
pord$(\rho_\text{isum}(\mathfrak{t})) = 80$,

\vskip 0.7ex
\hangindent=5.5em \hangafter=1
{\white .}\hskip 1em $\rho_\text{isum}(\mathfrak{t})$ =
 $( \frac{7}{40},
\frac{23}{40},
\frac{31}{80},
\frac{39}{80},
\frac{71}{80},
\frac{79}{80} )
$,

\vskip 0.7ex
\hangindent=5.5em \hangafter=1
{\white .}\hskip 1em $\rho_\text{isum}(\mathfrak{s})$ =
($0$,
$0$,
$\frac{1}{\sqrt{10}}c^{1}_{20}
$,
$\frac{1}{\sqrt{10}}c^{3}_{20}
$,
$\frac{1}{\sqrt{10}}c^{1}_{20}
$,
$\frac{1}{\sqrt{10}}c^{3}_{20}
$;
$0$,
$\frac{1}{\sqrt{10}}c^{3}_{20}
$,
$-\frac{1}{\sqrt{10}}c^{1}_{20}
$,
$\frac{1}{\sqrt{10}}c^{3}_{20}
$,
$-\frac{1}{\sqrt{10}}c^{1}_{20}
$;
$\frac{1}{2\sqrt{5}}c^{3}_{20}
$,
$\frac{1}{2\sqrt{5}}c^{1}_{20}
$,
$-\frac{1}{2\sqrt{5}}c^{3}_{20}
$,
$-\frac{1}{2\sqrt{5}}c^{1}_{20}
$;
$-\frac{1}{2\sqrt{5}}c^{3}_{20}
$,
$-\frac{1}{2\sqrt{5}}c^{1}_{20}
$,
$\frac{1}{2\sqrt{5}}c^{3}_{20}
$;
$\frac{1}{2\sqrt{5}}c^{3}_{20}
$,
$\frac{1}{2\sqrt{5}}c^{1}_{20}
$;
$-\frac{1}{2\sqrt{5}}c^{3}_{20}
$)

Pass. 

 \ \color{black}

 \color{blue}

\noindent 654: (dims,levels) = $(6;80
)$,
irreps = $3_{16}^{7,0}
\hskip -1.5pt \otimes \hskip -1.5pt
2_{5}^{2}$,
pord$(\rho_\text{isum}(\mathfrak{t})) = 80$,

\vskip 0.7ex
\hangindent=5.5em \hangafter=1
{\white .}\hskip 1em $\rho_\text{isum}(\mathfrak{t})$ =
 $( \frac{11}{40},
\frac{19}{40},
\frac{3}{80},
\frac{27}{80},
\frac{43}{80},
\frac{67}{80} )
$,

\vskip 0.7ex
\hangindent=5.5em \hangafter=1
{\white .}\hskip 1em $\rho_\text{isum}(\mathfrak{s})$ =
($0$,
$0$,
$\frac{1}{\sqrt{10}}c^{3}_{20}
$,
$\frac{1}{\sqrt{10}}c^{1}_{20}
$,
$\frac{1}{\sqrt{10}}c^{3}_{20}
$,
$\frac{1}{\sqrt{10}}c^{1}_{20}
$;
$0$,
$\frac{1}{\sqrt{10}}c^{1}_{20}
$,
$-\frac{1}{\sqrt{10}}c^{3}_{20}
$,
$\frac{1}{\sqrt{10}}c^{1}_{20}
$,
$-\frac{1}{\sqrt{10}}c^{3}_{20}
$;
$-\frac{1}{2\sqrt{5}}c^{1}_{20}
$,
$-\frac{1}{2\sqrt{5}}c^{3}_{20}
$,
$\frac{1}{2\sqrt{5}}c^{1}_{20}
$,
$\frac{1}{2\sqrt{5}}c^{3}_{20}
$;
$\frac{1}{2\sqrt{5}}c^{1}_{20}
$,
$\frac{1}{2\sqrt{5}}c^{3}_{20}
$,
$-\frac{1}{2\sqrt{5}}c^{1}_{20}
$;
$-\frac{1}{2\sqrt{5}}c^{1}_{20}
$,
$-\frac{1}{2\sqrt{5}}c^{3}_{20}
$;
$\frac{1}{2\sqrt{5}}c^{1}_{20}
$)

Pass. 

 \ \color{black}

 \color{blue}

\noindent 655: (dims,levels) = $(6;80
)$,
irreps = $3_{16}^{1,0}
\hskip -1.5pt \otimes \hskip -1.5pt
2_{5}^{1}$,
pord$(\rho_\text{isum}(\mathfrak{t})) = 80$,

\vskip 0.7ex
\hangindent=5.5em \hangafter=1
{\white .}\hskip 1em $\rho_\text{isum}(\mathfrak{t})$ =
 $( \frac{13}{40},
\frac{37}{40},
\frac{21}{80},
\frac{29}{80},
\frac{61}{80},
\frac{69}{80} )
$,

\vskip 0.7ex
\hangindent=5.5em \hangafter=1
{\white .}\hskip 1em $\rho_\text{isum}(\mathfrak{s})$ =
($0$,
$0$,
$\frac{1}{\sqrt{10}}c^{3}_{20}
$,
$\frac{1}{\sqrt{10}}c^{1}_{20}
$,
$\frac{1}{\sqrt{10}}c^{3}_{20}
$,
$\frac{1}{\sqrt{10}}c^{1}_{20}
$;
$0$,
$\frac{1}{\sqrt{10}}c^{1}_{20}
$,
$-\frac{1}{\sqrt{10}}c^{3}_{20}
$,
$\frac{1}{\sqrt{10}}c^{1}_{20}
$,
$-\frac{1}{\sqrt{10}}c^{3}_{20}
$;
$-\frac{1}{2\sqrt{5}}c^{3}_{20}
$,
$\frac{1}{2\sqrt{5}}c^{1}_{20}
$,
$\frac{1}{2\sqrt{5}}c^{3}_{20}
$,
$-\frac{1}{2\sqrt{5}}c^{1}_{20}
$;
$\frac{1}{2\sqrt{5}}c^{3}_{20}
$,
$-\frac{1}{2\sqrt{5}}c^{1}_{20}
$,
$-\frac{1}{2\sqrt{5}}c^{3}_{20}
$;
$-\frac{1}{2\sqrt{5}}c^{3}_{20}
$,
$\frac{1}{2\sqrt{5}}c^{1}_{20}
$;
$\frac{1}{2\sqrt{5}}c^{3}_{20}
$)

Pass. 

 \ \color{black}

 \color{blue}

\noindent 656: (dims,levels) = $(6;80
)$,
irreps = $3_{16}^{5,0}
\hskip -1.5pt \otimes \hskip -1.5pt
2_{5}^{1}$,
pord$(\rho_\text{isum}(\mathfrak{t})) = 80$,

\vskip 0.7ex
\hangindent=5.5em \hangafter=1
{\white .}\hskip 1em $\rho_\text{isum}(\mathfrak{t})$ =
 $( \frac{17}{40},
\frac{33}{40},
\frac{1}{80},
\frac{9}{80},
\frac{41}{80},
\frac{49}{80} )
$,

\vskip 0.7ex
\hangindent=5.5em \hangafter=1
{\white .}\hskip 1em $\rho_\text{isum}(\mathfrak{s})$ =
($0$,
$0$,
$\frac{1}{\sqrt{10}}c^{1}_{20}
$,
$\frac{1}{\sqrt{10}}c^{3}_{20}
$,
$\frac{1}{\sqrt{10}}c^{1}_{20}
$,
$\frac{1}{\sqrt{10}}c^{3}_{20}
$;
$0$,
$\frac{1}{\sqrt{10}}c^{3}_{20}
$,
$-\frac{1}{\sqrt{10}}c^{1}_{20}
$,
$\frac{1}{\sqrt{10}}c^{3}_{20}
$,
$-\frac{1}{\sqrt{10}}c^{1}_{20}
$;
$-\frac{1}{2\sqrt{5}}c^{3}_{20}
$,
$-\frac{1}{2\sqrt{5}}c^{1}_{20}
$,
$\frac{1}{2\sqrt{5}}c^{3}_{20}
$,
$\frac{1}{2\sqrt{5}}c^{1}_{20}
$;
$\frac{1}{2\sqrt{5}}c^{3}_{20}
$,
$\frac{1}{2\sqrt{5}}c^{1}_{20}
$,
$-\frac{1}{2\sqrt{5}}c^{3}_{20}
$;
$-\frac{1}{2\sqrt{5}}c^{3}_{20}
$,
$-\frac{1}{2\sqrt{5}}c^{1}_{20}
$;
$\frac{1}{2\sqrt{5}}c^{3}_{20}
$)

Pass. 

 \ \color{black}

 \color{blue}

\noindent 657: (dims,levels) = $(6;80
)$,
irreps = $3_{16}^{1,0}
\hskip -1.5pt \otimes \hskip -1.5pt
2_{5}^{2}$,
pord$(\rho_\text{isum}(\mathfrak{t})) = 80$,

\vskip 0.7ex
\hangindent=5.5em \hangafter=1
{\white .}\hskip 1em $\rho_\text{isum}(\mathfrak{t})$ =
 $( \frac{21}{40},
\frac{29}{40},
\frac{13}{80},
\frac{37}{80},
\frac{53}{80},
\frac{77}{80} )
$,

\vskip 0.7ex
\hangindent=5.5em \hangafter=1
{\white .}\hskip 1em $\rho_\text{isum}(\mathfrak{s})$ =
($0$,
$0$,
$\frac{1}{\sqrt{10}}c^{3}_{20}
$,
$\frac{1}{\sqrt{10}}c^{1}_{20}
$,
$\frac{1}{\sqrt{10}}c^{3}_{20}
$,
$\frac{1}{\sqrt{10}}c^{1}_{20}
$;
$0$,
$\frac{1}{\sqrt{10}}c^{1}_{20}
$,
$-\frac{1}{\sqrt{10}}c^{3}_{20}
$,
$\frac{1}{\sqrt{10}}c^{1}_{20}
$,
$-\frac{1}{\sqrt{10}}c^{3}_{20}
$;
$\frac{1}{2\sqrt{5}}c^{1}_{20}
$,
$\frac{1}{2\sqrt{5}}c^{3}_{20}
$,
$-\frac{1}{2\sqrt{5}}c^{1}_{20}
$,
$-\frac{1}{2\sqrt{5}}c^{3}_{20}
$;
$-\frac{1}{2\sqrt{5}}c^{1}_{20}
$,
$-\frac{1}{2\sqrt{5}}c^{3}_{20}
$,
$\frac{1}{2\sqrt{5}}c^{1}_{20}
$;
$\frac{1}{2\sqrt{5}}c^{1}_{20}
$,
$\frac{1}{2\sqrt{5}}c^{3}_{20}
$;
$-\frac{1}{2\sqrt{5}}c^{1}_{20}
$)

Pass. 

 \ \color{black}

 \color{blue}

\noindent 658: (dims,levels) = $(6;80
)$,
irreps = $3_{16}^{3,0}
\hskip -1.5pt \otimes \hskip -1.5pt
2_{5}^{2}$,
pord$(\rho_\text{isum}(\mathfrak{t})) = 80$,

\vskip 0.7ex
\hangindent=5.5em \hangafter=1
{\white .}\hskip 1em $\rho_\text{isum}(\mathfrak{t})$ =
 $( \frac{31}{40},
\frac{39}{40},
\frac{7}{80},
\frac{23}{80},
\frac{47}{80},
\frac{63}{80} )
$,

\vskip 0.7ex
\hangindent=5.5em \hangafter=1
{\white .}\hskip 1em $\rho_\text{isum}(\mathfrak{s})$ =
($0$,
$0$,
$\frac{1}{\sqrt{10}}c^{1}_{20}
$,
$\frac{1}{\sqrt{10}}c^{3}_{20}
$,
$\frac{1}{\sqrt{10}}c^{1}_{20}
$,
$\frac{1}{\sqrt{10}}c^{3}_{20}
$;
$0$,
$\frac{1}{\sqrt{10}}c^{3}_{20}
$,
$-\frac{1}{\sqrt{10}}c^{1}_{20}
$,
$\frac{1}{\sqrt{10}}c^{3}_{20}
$,
$-\frac{1}{\sqrt{10}}c^{1}_{20}
$;
$\frac{1}{2\sqrt{5}}c^{1}_{20}
$,
$\frac{1}{2\sqrt{5}}c^{3}_{20}
$,
$-\frac{1}{2\sqrt{5}}c^{1}_{20}
$,
$-\frac{1}{2\sqrt{5}}c^{3}_{20}
$;
$-\frac{1}{2\sqrt{5}}c^{1}_{20}
$,
$-\frac{1}{2\sqrt{5}}c^{3}_{20}
$,
$\frac{1}{2\sqrt{5}}c^{1}_{20}
$;
$\frac{1}{2\sqrt{5}}c^{1}_{20}
$,
$\frac{1}{2\sqrt{5}}c^{3}_{20}
$;
$-\frac{1}{2\sqrt{5}}c^{1}_{20}
$)

Pass. 

 \ \color{black}

\noindent 659: (dims,levels) = $(6;84
)$,
irreps = $6_{7,2}^{1}
\hskip -1.5pt \otimes \hskip -1.5pt
1_{4}^{1,0}
\hskip -1.5pt \otimes \hskip -1.5pt
1_{3}^{1,0}$,
pord$(\rho_\text{isum}(\mathfrak{t})) = 7$,

\vskip 0.7ex
\hangindent=5.5em \hangafter=1
{\white .}\hskip 1em $\rho_\text{isum}(\mathfrak{t})$ =
 $( \frac{1}{84},
\frac{13}{84},
\frac{25}{84},
\frac{37}{84},
\frac{61}{84},
\frac{73}{84} )
$,

\vskip 0.7ex
\hangindent=5.5em \hangafter=1
{\white .}\hskip 1em $\rho_\text{isum}(\mathfrak{s})$ =
$\mathrm{i}$($\frac{2}{7}-\frac{1}{7}c^{1}_{7}
$,
$\frac{1}{7}c^{3}_{56}
-\frac{1}{7}c^{5}_{56}
+\frac{1}{7}c^{9}_{56}
+\frac{1}{7}c^{11}_{56}
$,
$-\frac{2}{7}+\frac{1}{7}c^{2}_{7}
$,
$-\frac{3}{7}-\frac{1}{7}c^{1}_{7}
-\frac{1}{7}c^{2}_{7}
$,
$-\frac{1}{7}c^{3}_{56}
-\frac{2}{7}c^{5}_{56}
+\frac{1}{7}c^{7}_{56}
+\frac{2}{7}c^{9}_{56}
-\frac{1}{7}c^{11}_{56}
$,
$\frac{2}{7}c^{3}_{56}
+\frac{1}{7}c^{5}_{56}
-\frac{1}{7}c^{7}_{56}
-\frac{1}{7}c^{9}_{56}
+\frac{2}{7}c^{11}_{56}
$;\ \ 
$\frac{2}{7}-\frac{1}{7}c^{1}_{7}
$,
$\frac{2}{7}c^{3}_{56}
+\frac{1}{7}c^{5}_{56}
-\frac{1}{7}c^{7}_{56}
-\frac{1}{7}c^{9}_{56}
+\frac{2}{7}c^{11}_{56}
$,
$-\frac{1}{7}c^{3}_{56}
-\frac{2}{7}c^{5}_{56}
+\frac{1}{7}c^{7}_{56}
+\frac{2}{7}c^{9}_{56}
-\frac{1}{7}c^{11}_{56}
$,
$-\frac{3}{7}-\frac{1}{7}c^{1}_{7}
-\frac{1}{7}c^{2}_{7}
$,
$-\frac{2}{7}+\frac{1}{7}c^{2}_{7}
$;\ \ 
$\frac{3}{7}+\frac{1}{7}c^{1}_{7}
+\frac{1}{7}c^{2}_{7}
$,
$\frac{2}{7}-\frac{1}{7}c^{1}_{7}
$,
$\frac{1}{7}c^{3}_{56}
-\frac{1}{7}c^{5}_{56}
+\frac{1}{7}c^{9}_{56}
+\frac{1}{7}c^{11}_{56}
$,
$\frac{1}{7}c^{3}_{56}
+\frac{2}{7}c^{5}_{56}
-\frac{1}{7}c^{7}_{56}
-\frac{2}{7}c^{9}_{56}
+\frac{1}{7}c^{11}_{56}
$;\ \ 
$\frac{2}{7}-\frac{1}{7}c^{2}_{7}
$,
$-\frac{2}{7}c^{3}_{56}
-\frac{1}{7}c^{5}_{56}
+\frac{1}{7}c^{7}_{56}
+\frac{1}{7}c^{9}_{56}
-\frac{2}{7}c^{11}_{56}
$,
$\frac{1}{7}c^{3}_{56}
-\frac{1}{7}c^{5}_{56}
+\frac{1}{7}c^{9}_{56}
+\frac{1}{7}c^{11}_{56}
$;\ \ 
$\frac{2}{7}-\frac{1}{7}c^{2}_{7}
$,
$\frac{2}{7}-\frac{1}{7}c^{1}_{7}
$;\ \ 
$\frac{3}{7}+\frac{1}{7}c^{1}_{7}
+\frac{1}{7}c^{2}_{7}
$)

Fail:
cnd($\rho(\mathfrak s)_\mathrm{ndeg}$) = 56 does not divide
 ord($\rho(\mathfrak t)$)=84. Prop. B.4 (2)

 \ \color{black}

\noindent 660: (dims,levels) = $(6;84
)$,
irreps = $6_{7,1}^{3}
\hskip -1.5pt \otimes \hskip -1.5pt
1_{4}^{1,0}
\hskip -1.5pt \otimes \hskip -1.5pt
1_{3}^{1,0}$,
pord$(\rho_\text{isum}(\mathfrak{t})) = 7$,

\vskip 0.7ex
\hangindent=5.5em \hangafter=1
{\white .}\hskip 1em $\rho_\text{isum}(\mathfrak{t})$ =
 $( \frac{1}{84},
\frac{13}{84},
\frac{25}{84},
\frac{37}{84},
\frac{61}{84},
\frac{73}{84} )
$,

\vskip 0.7ex
\hangindent=5.5em \hangafter=1
{\white .}\hskip 1em $\rho_\text{isum}(\mathfrak{s})$ =
($-\frac{1}{7}c^{5}_{56}
+\frac{1}{7}c^{6}_{56}
-\frac{1}{7}c^{9}_{56}
$,
$\frac{1}{7}c^{1}_{112}
-\frac{1}{7}c^{3}_{112}
-\frac{1}{7}c^{11}_{112}
-\frac{1}{7}c^{15}_{112}
$,
$-\frac{1}{7}c^{2}_{56}
-\frac{1}{7}c^{3}_{56}
+\frac{1}{7}c^{11}_{56}
$,
$\frac{2}{7}c^{1}_{56}
-\frac{1}{7}c^{3}_{56}
-\frac{1}{7}c^{5}_{56}
+\frac{1}{7}c^{7}_{56}
+\frac{1}{7}c^{9}_{56}
+\frac{1}{7}c^{10}_{56}
-\frac{1}{7}c^{11}_{56}
$,
$\frac{1}{7}c^{1}_{112}
+\frac{1}{7}c^{5}_{112}
-\frac{1}{7}c^{7}_{112}
-\frac{1}{7}c^{9}_{112}
+\frac{1}{7}c^{15}_{112}
+\frac{2}{7}c^{17}_{112}
+\frac{1}{7}c^{19}_{112}
-\frac{1}{7}c^{23}_{112}
$,
$\frac{1}{7}c^{3}_{112}
-\frac{1}{7}c^{5}_{112}
-\frac{1}{7}c^{9}_{112}
-\frac{1}{7}c^{11}_{112}
+\frac{2}{7}c^{13}_{112}
+\frac{1}{7}c^{19}_{112}
-\frac{1}{7}c^{21}_{112}
+\frac{1}{7}c^{23}_{112}
$;
$\frac{1}{7}c^{5}_{56}
-\frac{1}{7}c^{6}_{56}
+\frac{1}{7}c^{9}_{56}
$,
$-\frac{1}{7}c^{3}_{112}
+\frac{1}{7}c^{5}_{112}
+\frac{1}{7}c^{9}_{112}
+\frac{1}{7}c^{11}_{112}
-\frac{2}{7}c^{13}_{112}
-\frac{1}{7}c^{19}_{112}
+\frac{1}{7}c^{21}_{112}
-\frac{1}{7}c^{23}_{112}
$,
$\frac{1}{7}c^{1}_{112}
+\frac{1}{7}c^{5}_{112}
-\frac{1}{7}c^{7}_{112}
-\frac{1}{7}c^{9}_{112}
+\frac{1}{7}c^{15}_{112}
+\frac{2}{7}c^{17}_{112}
+\frac{1}{7}c^{19}_{112}
-\frac{1}{7}c^{23}_{112}
$,
$-\frac{2}{7}c^{1}_{56}
+\frac{1}{7}c^{3}_{56}
+\frac{1}{7}c^{5}_{56}
-\frac{1}{7}c^{7}_{56}
-\frac{1}{7}c^{9}_{56}
-\frac{1}{7}c^{10}_{56}
+\frac{1}{7}c^{11}_{56}
$,
$-\frac{1}{7}c^{2}_{56}
-\frac{1}{7}c^{3}_{56}
+\frac{1}{7}c^{11}_{56}
$;
$-\frac{2}{7}c^{1}_{56}
+\frac{1}{7}c^{3}_{56}
+\frac{1}{7}c^{5}_{56}
-\frac{1}{7}c^{7}_{56}
-\frac{1}{7}c^{9}_{56}
-\frac{1}{7}c^{10}_{56}
+\frac{1}{7}c^{11}_{56}
$,
$-\frac{1}{7}c^{5}_{56}
+\frac{1}{7}c^{6}_{56}
-\frac{1}{7}c^{9}_{56}
$,
$\frac{1}{7}c^{1}_{112}
-\frac{1}{7}c^{3}_{112}
-\frac{1}{7}c^{11}_{112}
-\frac{1}{7}c^{15}_{112}
$,
$\frac{1}{7}c^{1}_{112}
+\frac{1}{7}c^{5}_{112}
-\frac{1}{7}c^{7}_{112}
-\frac{1}{7}c^{9}_{112}
+\frac{1}{7}c^{15}_{112}
+\frac{2}{7}c^{17}_{112}
+\frac{1}{7}c^{19}_{112}
-\frac{1}{7}c^{23}_{112}
$;
$\frac{1}{7}c^{2}_{56}
+\frac{1}{7}c^{3}_{56}
-\frac{1}{7}c^{11}_{56}
$,
$\frac{1}{7}c^{3}_{112}
-\frac{1}{7}c^{5}_{112}
-\frac{1}{7}c^{9}_{112}
-\frac{1}{7}c^{11}_{112}
+\frac{2}{7}c^{13}_{112}
+\frac{1}{7}c^{19}_{112}
-\frac{1}{7}c^{21}_{112}
+\frac{1}{7}c^{23}_{112}
$,
$-\frac{1}{7}c^{1}_{112}
+\frac{1}{7}c^{3}_{112}
+\frac{1}{7}c^{11}_{112}
+\frac{1}{7}c^{15}_{112}
$;
$-\frac{1}{7}c^{2}_{56}
-\frac{1}{7}c^{3}_{56}
+\frac{1}{7}c^{11}_{56}
$,
$-\frac{1}{7}c^{5}_{56}
+\frac{1}{7}c^{6}_{56}
-\frac{1}{7}c^{9}_{56}
$;
$\frac{2}{7}c^{1}_{56}
-\frac{1}{7}c^{3}_{56}
-\frac{1}{7}c^{5}_{56}
+\frac{1}{7}c^{7}_{56}
+\frac{1}{7}c^{9}_{56}
+\frac{1}{7}c^{10}_{56}
-\frac{1}{7}c^{11}_{56}
$)

Fail:
cnd( Tr$_I(\rho(\mathfrak s))$ ) =
56 does not divide
 ord($\rho(\mathfrak t)$) =
84, I = [ 1/84 ]. Prop. B.4 (2)

 \ \color{black}

\noindent 661: (dims,levels) = $(6;84
)$,
irreps = $6_{7,1}^{1}
\hskip -1.5pt \otimes \hskip -1.5pt
1_{4}^{1,0}
\hskip -1.5pt \otimes \hskip -1.5pt
1_{3}^{1,0}$,
pord$(\rho_\text{isum}(\mathfrak{t})) = 7$,

\vskip 0.7ex
\hangindent=5.5em \hangafter=1
{\white .}\hskip 1em $\rho_\text{isum}(\mathfrak{t})$ =
 $( \frac{1}{84},
\frac{13}{84},
\frac{25}{84},
\frac{37}{84},
\frac{61}{84},
\frac{73}{84} )
$,

\vskip 0.7ex
\hangindent=5.5em \hangafter=1
{\white .}\hskip 1em $\rho_\text{isum}(\mathfrak{s})$ =
($\frac{1}{7}c^{5}_{56}
+\frac{1}{7}c^{6}_{56}
+\frac{1}{7}c^{9}_{56}
$,
$\frac{1}{7}c^{1}_{112}
+\frac{1}{7}c^{3}_{112}
-\frac{1}{7}c^{5}_{112}
-\frac{1}{7}c^{7}_{112}
-\frac{1}{7}c^{9}_{112}
-\frac{1}{7}c^{11}_{112}
+\frac{2}{7}c^{13}_{112}
+\frac{1}{7}c^{15}_{112}
+\frac{2}{7}c^{17}_{112}
+\frac{1}{7}c^{19}_{112}
-\frac{1}{7}c^{21}_{112}
-\frac{1}{7}c^{23}_{112}
$,
$-\frac{1}{7}c^{2}_{56}
+\frac{1}{7}c^{3}_{56}
-\frac{1}{7}c^{11}_{56}
$,
$\frac{2}{7}c^{1}_{56}
-\frac{1}{7}c^{3}_{56}
-\frac{1}{7}c^{5}_{56}
+\frac{1}{7}c^{7}_{56}
+\frac{1}{7}c^{9}_{56}
-\frac{1}{7}c^{10}_{56}
-\frac{1}{7}c^{11}_{56}
$,
$\frac{1}{7}c^{3}_{112}
-\frac{1}{7}c^{9}_{112}
+\frac{1}{7}c^{11}_{112}
+\frac{1}{7}c^{23}_{112}
$,
$\frac{1}{7}c^{1}_{112}
+\frac{1}{7}c^{5}_{112}
-\frac{1}{7}c^{15}_{112}
+\frac{1}{7}c^{19}_{112}
$;
$-\frac{1}{7}c^{5}_{56}
-\frac{1}{7}c^{6}_{56}
-\frac{1}{7}c^{9}_{56}
$,
$\frac{1}{7}c^{1}_{112}
+\frac{1}{7}c^{5}_{112}
-\frac{1}{7}c^{15}_{112}
+\frac{1}{7}c^{19}_{112}
$,
$-\frac{1}{7}c^{3}_{112}
+\frac{1}{7}c^{9}_{112}
-\frac{1}{7}c^{11}_{112}
-\frac{1}{7}c^{23}_{112}
$,
$\frac{2}{7}c^{1}_{56}
-\frac{1}{7}c^{3}_{56}
-\frac{1}{7}c^{5}_{56}
+\frac{1}{7}c^{7}_{56}
+\frac{1}{7}c^{9}_{56}
-\frac{1}{7}c^{10}_{56}
-\frac{1}{7}c^{11}_{56}
$,
$\frac{1}{7}c^{2}_{56}
-\frac{1}{7}c^{3}_{56}
+\frac{1}{7}c^{11}_{56}
$;
$\frac{2}{7}c^{1}_{56}
-\frac{1}{7}c^{3}_{56}
-\frac{1}{7}c^{5}_{56}
+\frac{1}{7}c^{7}_{56}
+\frac{1}{7}c^{9}_{56}
-\frac{1}{7}c^{10}_{56}
-\frac{1}{7}c^{11}_{56}
$,
$-\frac{1}{7}c^{5}_{56}
-\frac{1}{7}c^{6}_{56}
-\frac{1}{7}c^{9}_{56}
$,
$\frac{1}{7}c^{1}_{112}
+\frac{1}{7}c^{3}_{112}
-\frac{1}{7}c^{5}_{112}
-\frac{1}{7}c^{7}_{112}
-\frac{1}{7}c^{9}_{112}
-\frac{1}{7}c^{11}_{112}
+\frac{2}{7}c^{13}_{112}
+\frac{1}{7}c^{15}_{112}
+\frac{2}{7}c^{17}_{112}
+\frac{1}{7}c^{19}_{112}
-\frac{1}{7}c^{21}_{112}
-\frac{1}{7}c^{23}_{112}
$,
$-\frac{1}{7}c^{3}_{112}
+\frac{1}{7}c^{9}_{112}
-\frac{1}{7}c^{11}_{112}
-\frac{1}{7}c^{23}_{112}
$;
$\frac{1}{7}c^{2}_{56}
-\frac{1}{7}c^{3}_{56}
+\frac{1}{7}c^{11}_{56}
$,
$\frac{1}{7}c^{1}_{112}
+\frac{1}{7}c^{5}_{112}
-\frac{1}{7}c^{15}_{112}
+\frac{1}{7}c^{19}_{112}
$,
$-\frac{1}{7}c^{1}_{112}
-\frac{1}{7}c^{3}_{112}
+\frac{1}{7}c^{5}_{112}
+\frac{1}{7}c^{7}_{112}
+\frac{1}{7}c^{9}_{112}
+\frac{1}{7}c^{11}_{112}
-\frac{2}{7}c^{13}_{112}
-\frac{1}{7}c^{15}_{112}
-\frac{2}{7}c^{17}_{112}
-\frac{1}{7}c^{19}_{112}
+\frac{1}{7}c^{21}_{112}
+\frac{1}{7}c^{23}_{112}
$;
$-\frac{1}{7}c^{2}_{56}
+\frac{1}{7}c^{3}_{56}
-\frac{1}{7}c^{11}_{56}
$,
$-\frac{1}{7}c^{5}_{56}
-\frac{1}{7}c^{6}_{56}
-\frac{1}{7}c^{9}_{56}
$;
$-\frac{2}{7}c^{1}_{56}
+\frac{1}{7}c^{3}_{56}
+\frac{1}{7}c^{5}_{56}
-\frac{1}{7}c^{7}_{56}
-\frac{1}{7}c^{9}_{56}
+\frac{1}{7}c^{10}_{56}
+\frac{1}{7}c^{11}_{56}
$)

Fail:
cnd( Tr$_I(\rho(\mathfrak s))$ ) =
56 does not divide
 ord($\rho(\mathfrak t)$) =
84, I = [ 1/84 ]. Prop. B.4 (2)

 \ \color{black}

\noindent 662: (dims,levels) = $(6;84
)$,
irreps = $3_{7}^{3}
\hskip -1.5pt \otimes \hskip -1.5pt
2_{4}^{1,0}
\hskip -1.5pt \otimes \hskip -1.5pt
1_{3}^{1,0}$,
pord$(\rho_\text{isum}(\mathfrak{t})) = 14$,

\vskip 0.7ex
\hangindent=5.5em \hangafter=1
{\white .}\hskip 1em $\rho_\text{isum}(\mathfrak{t})$ =
 $( \frac{1}{84},
\frac{25}{84},
\frac{37}{84},
\frac{43}{84},
\frac{67}{84},
\frac{79}{84} )
$,

\vskip 0.7ex
\hangindent=5.5em \hangafter=1
{\white .}\hskip 1em $\rho_\text{isum}(\mathfrak{s})$ =
$\mathrm{i}$($\frac{1}{2\sqrt{7}}c^{3}_{28}
$,
$-\frac{1}{2\sqrt{7}}c^{1}_{28}
$,
$\frac{1}{2\sqrt{7}}c^{5}_{28}
$,
$-\frac{3}{2\sqrt{21}}c^{3}_{28}
$,
$-\frac{3}{2\sqrt{21}}c^{1}_{28}
$,
$\frac{3}{2\sqrt{21}}c^{5}_{28}
$;\ \ 
$-\frac{1}{2\sqrt{7}}c^{5}_{28}
$,
$\frac{1}{2\sqrt{7}}c^{3}_{28}
$,
$\frac{3}{2\sqrt{21}}c^{1}_{28}
$,
$-\frac{3}{2\sqrt{21}}c^{5}_{28}
$,
$\frac{3}{2\sqrt{21}}c^{3}_{28}
$;\ \ 
$\frac{1}{2\sqrt{7}}c^{1}_{28}
$,
$-\frac{3}{2\sqrt{21}}c^{5}_{28}
$,
$\frac{3}{2\sqrt{21}}c^{3}_{28}
$,
$\frac{3}{2\sqrt{21}}c^{1}_{28}
$;\ \ 
$-\frac{1}{2\sqrt{7}}c^{3}_{28}
$,
$-\frac{1}{2\sqrt{7}}c^{1}_{28}
$,
$\frac{1}{2\sqrt{7}}c^{5}_{28}
$;\ \ 
$\frac{1}{2\sqrt{7}}c^{5}_{28}
$,
$-\frac{1}{2\sqrt{7}}c^{3}_{28}
$;\ \ 
$-\frac{1}{2\sqrt{7}}c^{1}_{28}
$)

Fail:
$\sigma(\rho(\mathfrak s)_\mathrm{ndeg}) \neq
 (\rho(\mathfrak t)^a \rho(\mathfrak s) \rho(\mathfrak t)^b
 \rho(\mathfrak s) \rho(\mathfrak t)^a)_\mathrm{ndeg}$,
 $\sigma = a$ = 5. Prop. B.5 (3) eqn. (B.25)

 \ \color{black}

\noindent 663: (dims,levels) = $(6;84
)$,
irreps = $3_{7}^{3}
\hskip -1.5pt \otimes \hskip -1.5pt
2_{3}^{1,0}
\hskip -1.5pt \otimes \hskip -1.5pt
1_{4}^{1,0}$,
pord$(\rho_\text{isum}(\mathfrak{t})) = 21$,

\vskip 0.7ex
\hangindent=5.5em \hangafter=1
{\white .}\hskip 1em $\rho_\text{isum}(\mathfrak{t})$ =
 $( \frac{3}{28},
\frac{19}{28},
\frac{27}{28},
\frac{1}{84},
\frac{25}{84},
\frac{37}{84} )
$,

\vskip 0.7ex
\hangindent=5.5em \hangafter=1
{\white .}\hskip 1em $\rho_\text{isum}(\mathfrak{s})$ =
($-\frac{1}{\sqrt{21}}c^{1}_{28}
$,
$-\frac{1}{\sqrt{21}}c^{5}_{28}
$,
$\frac{1}{\sqrt{21}}c^{3}_{28}
$,
$-\frac{2}{\sqrt{42}}c^{5}_{28}
$,
$\frac{2}{\sqrt{42}}c^{3}_{28}
$,
$\frac{2}{\sqrt{42}}c^{1}_{28}
$;
$-\frac{1}{\sqrt{21}}c^{3}_{28}
$,
$-\frac{1}{\sqrt{21}}c^{1}_{28}
$,
$-\frac{2}{\sqrt{42}}c^{3}_{28}
$,
$-\frac{2}{\sqrt{42}}c^{1}_{28}
$,
$\frac{2}{\sqrt{42}}c^{5}_{28}
$;
$\frac{1}{\sqrt{21}}c^{5}_{28}
$,
$-\frac{2}{\sqrt{42}}c^{1}_{28}
$,
$\frac{2}{\sqrt{42}}c^{5}_{28}
$,
$-\frac{2}{\sqrt{42}}c^{3}_{28}
$;
$\frac{1}{\sqrt{21}}c^{3}_{28}
$,
$\frac{1}{\sqrt{21}}c^{1}_{28}
$,
$-\frac{1}{\sqrt{21}}c^{5}_{28}
$;
$-\frac{1}{\sqrt{21}}c^{5}_{28}
$,
$\frac{1}{\sqrt{21}}c^{3}_{28}
$;
$\frac{1}{\sqrt{21}}c^{1}_{28}
$)

Fail:
cnd($\rho(\mathfrak s)_\mathrm{ndeg}$) = 168 does not divide
 ord($\rho(\mathfrak t)$)=84. Prop. B.4 (2)

 \ \color{black}

\noindent 664: (dims,levels) = $(6;84
)$,
irreps = $3_{7}^{1}
\hskip -1.5pt \otimes \hskip -1.5pt
2_{4}^{1,0}
\hskip -1.5pt \otimes \hskip -1.5pt
1_{3}^{1,0}$,
pord$(\rho_\text{isum}(\mathfrak{t})) = 14$,

\vskip 0.7ex
\hangindent=5.5em \hangafter=1
{\white .}\hskip 1em $\rho_\text{isum}(\mathfrak{t})$ =
 $( \frac{13}{84},
\frac{19}{84},
\frac{31}{84},
\frac{55}{84},
\frac{61}{84},
\frac{73}{84} )
$,

\vskip 0.7ex
\hangindent=5.5em \hangafter=1
{\white .}\hskip 1em $\rho_\text{isum}(\mathfrak{s})$ =
$\mathrm{i}$($\frac{1}{2\sqrt{7}}c^{3}_{28}
$,
$\frac{3}{2\sqrt{21}}c^{5}_{28}
$,
$-\frac{3}{2\sqrt{21}}c^{1}_{28}
$,
$-\frac{3}{2\sqrt{21}}c^{3}_{28}
$,
$\frac{1}{2\sqrt{7}}c^{5}_{28}
$,
$-\frac{1}{2\sqrt{7}}c^{1}_{28}
$;\ \ 
$-\frac{1}{2\sqrt{7}}c^{1}_{28}
$,
$-\frac{1}{2\sqrt{7}}c^{3}_{28}
$,
$\frac{1}{2\sqrt{7}}c^{5}_{28}
$,
$\frac{3}{2\sqrt{21}}c^{1}_{28}
$,
$\frac{3}{2\sqrt{21}}c^{3}_{28}
$;\ \ 
$\frac{1}{2\sqrt{7}}c^{5}_{28}
$,
$-\frac{1}{2\sqrt{7}}c^{1}_{28}
$,
$\frac{3}{2\sqrt{21}}c^{3}_{28}
$,
$-\frac{3}{2\sqrt{21}}c^{5}_{28}
$;\ \ 
$-\frac{1}{2\sqrt{7}}c^{3}_{28}
$,
$-\frac{3}{2\sqrt{21}}c^{5}_{28}
$,
$\frac{3}{2\sqrt{21}}c^{1}_{28}
$;\ \ 
$\frac{1}{2\sqrt{7}}c^{1}_{28}
$,
$\frac{1}{2\sqrt{7}}c^{3}_{28}
$;\ \ 
$-\frac{1}{2\sqrt{7}}c^{5}_{28}
$)

Fail:
$\sigma(\rho(\mathfrak s)_\mathrm{ndeg}) \neq
 (\rho(\mathfrak t)^a \rho(\mathfrak s) \rho(\mathfrak t)^b
 \rho(\mathfrak s) \rho(\mathfrak t)^a)_\mathrm{ndeg}$,
 $\sigma = a$ = 5. Prop. B.5 (3) eqn. (B.25)

 \ \color{black}

\noindent 665: (dims,levels) = $(6;84
)$,
irreps = $3_{7}^{1}
\hskip -1.5pt \otimes \hskip -1.5pt
2_{3}^{1,0}
\hskip -1.5pt \otimes \hskip -1.5pt
1_{4}^{1,0}$,
pord$(\rho_\text{isum}(\mathfrak{t})) = 21$,

\vskip 0.7ex
\hangindent=5.5em \hangafter=1
{\white .}\hskip 1em $\rho_\text{isum}(\mathfrak{t})$ =
 $( \frac{11}{28},
\frac{15}{28},
\frac{23}{28},
\frac{13}{84},
\frac{61}{84},
\frac{73}{84} )
$,

\vskip 0.7ex
\hangindent=5.5em \hangafter=1
{\white .}\hskip 1em $\rho_\text{isum}(\mathfrak{s})$ =
($-\frac{1}{\sqrt{21}}c^{1}_{28}
$,
$\frac{1}{\sqrt{21}}c^{3}_{28}
$,
$-\frac{1}{\sqrt{21}}c^{5}_{28}
$,
$-\frac{2}{\sqrt{42}}c^{5}_{28}
$,
$\frac{2}{\sqrt{42}}c^{1}_{28}
$,
$\frac{2}{\sqrt{42}}c^{3}_{28}
$;
$\frac{1}{\sqrt{21}}c^{5}_{28}
$,
$-\frac{1}{\sqrt{21}}c^{1}_{28}
$,
$-\frac{2}{\sqrt{42}}c^{1}_{28}
$,
$-\frac{2}{\sqrt{42}}c^{3}_{28}
$,
$\frac{2}{\sqrt{42}}c^{5}_{28}
$;
$-\frac{1}{\sqrt{21}}c^{3}_{28}
$,
$-\frac{2}{\sqrt{42}}c^{3}_{28}
$,
$\frac{2}{\sqrt{42}}c^{5}_{28}
$,
$-\frac{2}{\sqrt{42}}c^{1}_{28}
$;
$\frac{1}{\sqrt{21}}c^{3}_{28}
$,
$-\frac{1}{\sqrt{21}}c^{5}_{28}
$,
$\frac{1}{\sqrt{21}}c^{1}_{28}
$;
$\frac{1}{\sqrt{21}}c^{1}_{28}
$,
$\frac{1}{\sqrt{21}}c^{3}_{28}
$;
$-\frac{1}{\sqrt{21}}c^{5}_{28}
$)

Fail:
cnd($\rho(\mathfrak s)_\mathrm{ndeg}$) = 168 does not divide
 ord($\rho(\mathfrak t)$)=84. Prop. B.4 (2)

 \ \color{black}

\noindent 666: (dims,levels) = $(6;96
)$,
irreps = $6_{32,2}^{7,0}
\hskip -1.5pt \otimes \hskip -1.5pt
1_{3}^{1,0}$,
pord$(\rho_\text{isum}(\mathfrak{t})) = 32$,

\vskip 0.7ex
\hangindent=5.5em \hangafter=1
{\white .}\hskip 1em $\rho_\text{isum}(\mathfrak{t})$ =
 $( \frac{1}{3},
\frac{5}{24},
\frac{11}{96},
\frac{35}{96},
\frac{59}{96},
\frac{83}{96} )
$,

\vskip 0.7ex
\hangindent=5.5em \hangafter=1
{\white .}\hskip 1em $\rho_\text{isum}(\mathfrak{s})$ =
($0$,
$0$,
$\frac{1}{2}$,
$\frac{1}{2}$,
$\frac{1}{2}$,
$\frac{1}{2}$;
$0$,
$\frac{1}{2}$,
$-\frac{1}{2}$,
$\frac{1}{2}$,
$-\frac{1}{2}$;
$-\frac{1}{4}c^{1}_{16}
$,
$\frac{1}{4}c^{3}_{16}
$,
$\frac{1}{4}c^{1}_{16}
$,
$-\frac{1}{4}c^{3}_{16}
$;
$\frac{1}{4}c^{1}_{16}
$,
$-\frac{1}{4}c^{3}_{16}
$,
$-\frac{1}{4}c^{1}_{16}
$;
$-\frac{1}{4}c^{1}_{16}
$,
$\frac{1}{4}c^{3}_{16}
$;
$\frac{1}{4}c^{1}_{16}
$)

Fail:
$\sigma(\rho(\mathfrak s)_\mathrm{ndeg}) \neq
 (\rho(\mathfrak t)^a \rho(\mathfrak s) \rho(\mathfrak t)^b
 \rho(\mathfrak s) \rho(\mathfrak t)^a)_\mathrm{ndeg}$,
 $\sigma = a$ = 5. Prop. B.5 (3) eqn. (B.25)

 \ \color{black}

\noindent 667: (dims,levels) = $(6;96
)$,
irreps = $6_{32,1}^{7,0}
\hskip -1.5pt \otimes \hskip -1.5pt
1_{3}^{1,0}$,
pord$(\rho_\text{isum}(\mathfrak{t})) = 32$,

\vskip 0.7ex
\hangindent=5.5em \hangafter=1
{\white .}\hskip 1em $\rho_\text{isum}(\mathfrak{t})$ =
 $( \frac{1}{3},
\frac{5}{24},
\frac{23}{96},
\frac{47}{96},
\frac{71}{96},
\frac{95}{96} )
$,

\vskip 0.7ex
\hangindent=5.5em \hangafter=1
{\white .}\hskip 1em $\rho_\text{isum}(\mathfrak{s})$ =
$\mathrm{i}$($0$,
$0$,
$\frac{1}{2}$,
$\frac{1}{2}$,
$\frac{1}{2}$,
$\frac{1}{2}$;\ \ 
$0$,
$\frac{1}{2}$,
$-\frac{1}{2}$,
$\frac{1}{2}$,
$-\frac{1}{2}$;\ \ 
$\frac{1}{4}c^{1}_{16}
$,
$-\frac{1}{4}c^{3}_{16}
$,
$-\frac{1}{4}c^{1}_{16}
$,
$\frac{1}{4}c^{3}_{16}
$;\ \ 
$-\frac{1}{4}c^{1}_{16}
$,
$\frac{1}{4}c^{3}_{16}
$,
$\frac{1}{4}c^{1}_{16}
$;\ \ 
$\frac{1}{4}c^{1}_{16}
$,
$-\frac{1}{4}c^{3}_{16}
$;\ \ 
$-\frac{1}{4}c^{1}_{16}
$)

Fail:
$\sigma(\rho(\mathfrak s)_\mathrm{ndeg}) \neq
 (\rho(\mathfrak t)^a \rho(\mathfrak s) \rho(\mathfrak t)^b
 \rho(\mathfrak s) \rho(\mathfrak t)^a)_\mathrm{ndeg}$,
 $\sigma = a$ = 5. Prop. B.5 (3) eqn. (B.25)

 \ \color{black}

\noindent 668: (dims,levels) = $(6;96
)$,
irreps = $6_{32,2}^{1,0}
\hskip -1.5pt \otimes \hskip -1.5pt
1_{3}^{1,0}$,
pord$(\rho_\text{isum}(\mathfrak{t})) = 32$,

\vskip 0.7ex
\hangindent=5.5em \hangafter=1
{\white .}\hskip 1em $\rho_\text{isum}(\mathfrak{t})$ =
 $( \frac{1}{3},
\frac{11}{24},
\frac{5}{96},
\frac{29}{96},
\frac{53}{96},
\frac{77}{96} )
$,

\vskip 0.7ex
\hangindent=5.5em \hangafter=1
{\white .}\hskip 1em $\rho_\text{isum}(\mathfrak{s})$ =
($0$,
$0$,
$\frac{1}{2}$,
$\frac{1}{2}$,
$\frac{1}{2}$,
$\frac{1}{2}$;
$0$,
$\frac{1}{2}$,
$-\frac{1}{2}$,
$\frac{1}{2}$,
$-\frac{1}{2}$;
$-\frac{1}{4}c^{1}_{16}
$,
$-\frac{1}{4}c^{3}_{16}
$,
$\frac{1}{4}c^{1}_{16}
$,
$\frac{1}{4}c^{3}_{16}
$;
$\frac{1}{4}c^{1}_{16}
$,
$\frac{1}{4}c^{3}_{16}
$,
$-\frac{1}{4}c^{1}_{16}
$;
$-\frac{1}{4}c^{1}_{16}
$,
$-\frac{1}{4}c^{3}_{16}
$;
$\frac{1}{4}c^{1}_{16}
$)

Fail:
$\sigma(\rho(\mathfrak s)_\mathrm{ndeg}) \neq
 (\rho(\mathfrak t)^a \rho(\mathfrak s) \rho(\mathfrak t)^b
 \rho(\mathfrak s) \rho(\mathfrak t)^a)_\mathrm{ndeg}$,
 $\sigma = a$ = 5. Prop. B.5 (3) eqn. (B.25)

 \ \color{black}

\noindent 669: (dims,levels) = $(6;96
)$,
irreps = $6_{32,1}^{1,0}
\hskip -1.5pt \otimes \hskip -1.5pt
1_{3}^{1,0}$,
pord$(\rho_\text{isum}(\mathfrak{t})) = 32$,

\vskip 0.7ex
\hangindent=5.5em \hangafter=1
{\white .}\hskip 1em $\rho_\text{isum}(\mathfrak{t})$ =
 $( \frac{1}{3},
\frac{11}{24},
\frac{17}{96},
\frac{41}{96},
\frac{65}{96},
\frac{89}{96} )
$,

\vskip 0.7ex
\hangindent=5.5em \hangafter=1
{\white .}\hskip 1em $\rho_\text{isum}(\mathfrak{s})$ =
$\mathrm{i}$($0$,
$0$,
$\frac{1}{2}$,
$\frac{1}{2}$,
$\frac{1}{2}$,
$\frac{1}{2}$;\ \ 
$0$,
$\frac{1}{2}$,
$-\frac{1}{2}$,
$\frac{1}{2}$,
$-\frac{1}{2}$;\ \ 
$\frac{1}{4}c^{1}_{16}
$,
$\frac{1}{4}c^{3}_{16}
$,
$-\frac{1}{4}c^{1}_{16}
$,
$-\frac{1}{4}c^{3}_{16}
$;\ \ 
$-\frac{1}{4}c^{1}_{16}
$,
$-\frac{1}{4}c^{3}_{16}
$,
$\frac{1}{4}c^{1}_{16}
$;\ \ 
$\frac{1}{4}c^{1}_{16}
$,
$\frac{1}{4}c^{3}_{16}
$;\ \ 
$-\frac{1}{4}c^{1}_{16}
$)

Fail:
$\sigma(\rho(\mathfrak s)_\mathrm{ndeg}) \neq
 (\rho(\mathfrak t)^a \rho(\mathfrak s) \rho(\mathfrak t)^b
 \rho(\mathfrak s) \rho(\mathfrak t)^a)_\mathrm{ndeg}$,
 $\sigma = a$ = 5. Prop. B.5 (3) eqn. (B.25)

 \ \color{black}

\noindent 670: (dims,levels) = $(6;96
)$,
irreps = $6_{32,1}^{3,0}
\hskip -1.5pt \otimes \hskip -1.5pt
1_{3}^{1,0}$,
pord$(\rho_\text{isum}(\mathfrak{t})) = 32$,

\vskip 0.7ex
\hangindent=5.5em \hangafter=1
{\white .}\hskip 1em $\rho_\text{isum}(\mathfrak{t})$ =
 $( \frac{1}{3},
\frac{17}{24},
\frac{11}{96},
\frac{35}{96},
\frac{59}{96},
\frac{83}{96} )
$,

\vskip 0.7ex
\hangindent=5.5em \hangafter=1
{\white .}\hskip 1em $\rho_\text{isum}(\mathfrak{s})$ =
$\mathrm{i}$($0$,
$0$,
$\frac{1}{2}$,
$\frac{1}{2}$,
$\frac{1}{2}$,
$\frac{1}{2}$;\ \ 
$0$,
$\frac{1}{2}$,
$-\frac{1}{2}$,
$\frac{1}{2}$,
$-\frac{1}{2}$;\ \ 
$\frac{1}{4}c^{3}_{16}
$,
$\frac{1}{4}c^{1}_{16}
$,
$-\frac{1}{4}c^{3}_{16}
$,
$-\frac{1}{4}c^{1}_{16}
$;\ \ 
$-\frac{1}{4}c^{3}_{16}
$,
$-\frac{1}{4}c^{1}_{16}
$,
$\frac{1}{4}c^{3}_{16}
$;\ \ 
$\frac{1}{4}c^{3}_{16}
$,
$\frac{1}{4}c^{1}_{16}
$;\ \ 
$-\frac{1}{4}c^{3}_{16}
$)

Fail:
$\sigma(\rho(\mathfrak s)_\mathrm{ndeg}) \neq
 (\rho(\mathfrak t)^a \rho(\mathfrak s) \rho(\mathfrak t)^b
 \rho(\mathfrak s) \rho(\mathfrak t)^a)_\mathrm{ndeg}$,
 $\sigma = a$ = 5. Prop. B.5 (3) eqn. (B.25)

 \ \color{black}

\noindent 671: (dims,levels) = $(6;96
)$,
irreps = $6_{32,2}^{3,0}
\hskip -1.5pt \otimes \hskip -1.5pt
1_{3}^{1,0}$,
pord$(\rho_\text{isum}(\mathfrak{t})) = 32$,

\vskip 0.7ex
\hangindent=5.5em \hangafter=1
{\white .}\hskip 1em $\rho_\text{isum}(\mathfrak{t})$ =
 $( \frac{1}{3},
\frac{17}{24},
\frac{23}{96},
\frac{47}{96},
\frac{71}{96},
\frac{95}{96} )
$,

\vskip 0.7ex
\hangindent=5.5em \hangafter=1
{\white .}\hskip 1em $\rho_\text{isum}(\mathfrak{s})$ =
($0$,
$0$,
$\frac{1}{2}$,
$\frac{1}{2}$,
$\frac{1}{2}$,
$\frac{1}{2}$;
$0$,
$\frac{1}{2}$,
$-\frac{1}{2}$,
$\frac{1}{2}$,
$-\frac{1}{2}$;
$\frac{1}{4}c^{3}_{16}
$,
$\frac{1}{4}c^{1}_{16}
$,
$-\frac{1}{4}c^{3}_{16}
$,
$-\frac{1}{4}c^{1}_{16}
$;
$-\frac{1}{4}c^{3}_{16}
$,
$-\frac{1}{4}c^{1}_{16}
$,
$\frac{1}{4}c^{3}_{16}
$;
$\frac{1}{4}c^{3}_{16}
$,
$\frac{1}{4}c^{1}_{16}
$;
$-\frac{1}{4}c^{3}_{16}
$)

Fail:
$\sigma(\rho(\mathfrak s)_\mathrm{ndeg}) \neq
 (\rho(\mathfrak t)^a \rho(\mathfrak s) \rho(\mathfrak t)^b
 \rho(\mathfrak s) \rho(\mathfrak t)^a)_\mathrm{ndeg}$,
 $\sigma = a$ = 5. Prop. B.5 (3) eqn. (B.25)

 \ \color{black}

\noindent 672: (dims,levels) = $(6;96
)$,
irreps = $6_{32,1}^{5,0}
\hskip -1.5pt \otimes \hskip -1.5pt
1_{3}^{1,0}$,
pord$(\rho_\text{isum}(\mathfrak{t})) = 32$,

\vskip 0.7ex
\hangindent=5.5em \hangafter=1
{\white .}\hskip 1em $\rho_\text{isum}(\mathfrak{t})$ =
 $( \frac{1}{3},
\frac{23}{24},
\frac{5}{96},
\frac{29}{96},
\frac{53}{96},
\frac{77}{96} )
$,

\vskip 0.7ex
\hangindent=5.5em \hangafter=1
{\white .}\hskip 1em $\rho_\text{isum}(\mathfrak{s})$ =
$\mathrm{i}$($0$,
$0$,
$\frac{1}{2}$,
$\frac{1}{2}$,
$\frac{1}{2}$,
$\frac{1}{2}$;\ \ 
$0$,
$\frac{1}{2}$,
$-\frac{1}{2}$,
$\frac{1}{2}$,
$-\frac{1}{2}$;\ \ 
$-\frac{1}{4}c^{3}_{16}
$,
$\frac{1}{4}c^{1}_{16}
$,
$\frac{1}{4}c^{3}_{16}
$,
$-\frac{1}{4}c^{1}_{16}
$;\ \ 
$\frac{1}{4}c^{3}_{16}
$,
$-\frac{1}{4}c^{1}_{16}
$,
$-\frac{1}{4}c^{3}_{16}
$;\ \ 
$-\frac{1}{4}c^{3}_{16}
$,
$\frac{1}{4}c^{1}_{16}
$;\ \ 
$\frac{1}{4}c^{3}_{16}
$)

Fail:
$\sigma(\rho(\mathfrak s)_\mathrm{ndeg}) \neq
 (\rho(\mathfrak t)^a \rho(\mathfrak s) \rho(\mathfrak t)^b
 \rho(\mathfrak s) \rho(\mathfrak t)^a)_\mathrm{ndeg}$,
 $\sigma = a$ = 5. Prop. B.5 (3) eqn. (B.25)

 \ \color{black}

\noindent 673: (dims,levels) = $(6;96
)$,
irreps = $6_{32,2}^{5,0}
\hskip -1.5pt \otimes \hskip -1.5pt
1_{3}^{1,0}$,
pord$(\rho_\text{isum}(\mathfrak{t})) = 32$,

\vskip 0.7ex
\hangindent=5.5em \hangafter=1
{\white .}\hskip 1em $\rho_\text{isum}(\mathfrak{t})$ =
 $( \frac{1}{3},
\frac{23}{24},
\frac{17}{96},
\frac{41}{96},
\frac{65}{96},
\frac{89}{96} )
$,

\vskip 0.7ex
\hangindent=5.5em \hangafter=1
{\white .}\hskip 1em $\rho_\text{isum}(\mathfrak{s})$ =
($0$,
$0$,
$\frac{1}{2}$,
$\frac{1}{2}$,
$\frac{1}{2}$,
$\frac{1}{2}$;
$0$,
$\frac{1}{2}$,
$-\frac{1}{2}$,
$\frac{1}{2}$,
$-\frac{1}{2}$;
$-\frac{1}{4}c^{3}_{16}
$,
$\frac{1}{4}c^{1}_{16}
$,
$\frac{1}{4}c^{3}_{16}
$,
$-\frac{1}{4}c^{1}_{16}
$;
$\frac{1}{4}c^{3}_{16}
$,
$-\frac{1}{4}c^{1}_{16}
$,
$-\frac{1}{4}c^{3}_{16}
$;
$-\frac{1}{4}c^{3}_{16}
$,
$\frac{1}{4}c^{1}_{16}
$;
$\frac{1}{4}c^{3}_{16}
$)

Fail:
$\sigma(\rho(\mathfrak s)_\mathrm{ndeg}) \neq
 (\rho(\mathfrak t)^a \rho(\mathfrak s) \rho(\mathfrak t)^b
 \rho(\mathfrak s) \rho(\mathfrak t)^a)_\mathrm{ndeg}$,
 $\sigma = a$ = 5. Prop. B.5 (3) eqn. (B.25)

 \ \color{black}

 \color{blue}

\noindent 674: (dims,levels) = $(6;105
)$,
irreps = $3_{7}^{1}
\hskip -1.5pt \otimes \hskip -1.5pt
2_{5}^{2}
\hskip -1.5pt \otimes \hskip -1.5pt
1_{3}^{1,0}$,
pord$(\rho_\text{isum}(\mathfrak{t})) = 35$,

\vskip 0.7ex
\hangindent=5.5em \hangafter=1
{\white .}\hskip 1em $\rho_\text{isum}(\mathfrak{t})$ =
 $( \frac{2}{105},
\frac{8}{105},
\frac{23}{105},
\frac{32}{105},
\frac{53}{105},
\frac{92}{105} )
$,

\vskip 0.7ex
\hangindent=5.5em \hangafter=1
{\white .}\hskip 1em $\rho_\text{isum}(\mathfrak{s})$ =
$\mathrm{i}$($-\frac{1}{\sqrt{35}\mathrm{i}}s^{3}_{140}
-\frac{1}{\sqrt{35}\mathrm{i}}s^{17}_{140}
$,
$\frac{2}{35}c^{1}_{140}
-\frac{1}{35}c^{3}_{140}
-\frac{1}{7}c^{5}_{140}
-\frac{3}{35}c^{7}_{140}
+\frac{1}{5}c^{9}_{140}
-\frac{2}{35}c^{13}_{140}
-\frac{1}{35}c^{15}_{140}
-\frac{1}{35}c^{17}_{140}
+\frac{3}{35}c^{19}_{140}
+\frac{2}{35}c^{21}_{140}
-\frac{2}{7}c^{23}_{140}
$,
$\frac{1}{\sqrt{35}}c^{1}_{35}
-\frac{1}{\sqrt{35}}c^{6}_{35}
$,
$\frac{2}{\sqrt{35}}c^{3}_{35}
+\frac{1}{\sqrt{35}}c^{4}_{35}
+\frac{1}{\sqrt{35}}c^{10}_{35}
+\frac{1}{\sqrt{35}}c^{11}_{35}
$,
$-\frac{1}{\sqrt{35}}c^{4}_{35}
+\frac{1}{\sqrt{35}}c^{11}_{35}
$,
$\frac{4}{35}c^{1}_{140}
+\frac{3}{35}c^{3}_{140}
+\frac{1}{7}c^{5}_{140}
-\frac{1}{35}c^{7}_{140}
-\frac{1}{35}c^{9}_{140}
-\frac{4}{35}c^{13}_{140}
-\frac{2}{35}c^{15}_{140}
+\frac{3}{35}c^{17}_{140}
-\frac{9}{35}c^{19}_{140}
+\frac{4}{35}c^{21}_{140}
+\frac{2}{7}c^{23}_{140}
$;\ \ 
$-\frac{2}{\sqrt{35}}c^{3}_{35}
-\frac{1}{\sqrt{35}}c^{4}_{35}
-\frac{1}{\sqrt{35}}c^{10}_{35}
-\frac{1}{\sqrt{35}}c^{11}_{35}
$,
$-\frac{4}{35}c^{1}_{140}
-\frac{3}{35}c^{3}_{140}
-\frac{1}{7}c^{5}_{140}
+\frac{1}{35}c^{7}_{140}
+\frac{1}{35}c^{9}_{140}
+\frac{4}{35}c^{13}_{140}
+\frac{2}{35}c^{15}_{140}
-\frac{3}{35}c^{17}_{140}
+\frac{9}{35}c^{19}_{140}
-\frac{4}{35}c^{21}_{140}
-\frac{2}{7}c^{23}_{140}
$,
$\frac{1}{\sqrt{35}}c^{1}_{35}
-\frac{1}{\sqrt{35}}c^{6}_{35}
$,
$\frac{1}{\sqrt{35}\mathrm{i}}s^{3}_{140}
+\frac{1}{\sqrt{35}\mathrm{i}}s^{17}_{140}
$,
$-\frac{1}{\sqrt{35}}c^{4}_{35}
+\frac{1}{\sqrt{35}}c^{11}_{35}
$;\ \ 
$\frac{1}{\sqrt{35}\mathrm{i}}s^{3}_{140}
+\frac{1}{\sqrt{35}\mathrm{i}}s^{17}_{140}
$,
$-\frac{1}{\sqrt{35}}c^{4}_{35}
+\frac{1}{\sqrt{35}}c^{11}_{35}
$,
$-\frac{2}{\sqrt{35}}c^{3}_{35}
-\frac{1}{\sqrt{35}}c^{4}_{35}
-\frac{1}{\sqrt{35}}c^{10}_{35}
-\frac{1}{\sqrt{35}}c^{11}_{35}
$,
$\frac{2}{35}c^{1}_{140}
-\frac{1}{35}c^{3}_{140}
-\frac{1}{7}c^{5}_{140}
-\frac{3}{35}c^{7}_{140}
+\frac{1}{5}c^{9}_{140}
-\frac{2}{35}c^{13}_{140}
-\frac{1}{35}c^{15}_{140}
-\frac{1}{35}c^{17}_{140}
+\frac{3}{35}c^{19}_{140}
+\frac{2}{35}c^{21}_{140}
-\frac{2}{7}c^{23}_{140}
$;\ \ 
$\frac{4}{35}c^{1}_{140}
+\frac{3}{35}c^{3}_{140}
+\frac{1}{7}c^{5}_{140}
-\frac{1}{35}c^{7}_{140}
-\frac{1}{35}c^{9}_{140}
-\frac{4}{35}c^{13}_{140}
-\frac{2}{35}c^{15}_{140}
+\frac{3}{35}c^{17}_{140}
-\frac{9}{35}c^{19}_{140}
+\frac{4}{35}c^{21}_{140}
+\frac{2}{7}c^{23}_{140}
$,
$\frac{2}{35}c^{1}_{140}
-\frac{1}{35}c^{3}_{140}
-\frac{1}{7}c^{5}_{140}
-\frac{3}{35}c^{7}_{140}
+\frac{1}{5}c^{9}_{140}
-\frac{2}{35}c^{13}_{140}
-\frac{1}{35}c^{15}_{140}
-\frac{1}{35}c^{17}_{140}
+\frac{3}{35}c^{19}_{140}
+\frac{2}{35}c^{21}_{140}
-\frac{2}{7}c^{23}_{140}
$,
$-\frac{1}{\sqrt{35}\mathrm{i}}s^{3}_{140}
-\frac{1}{\sqrt{35}\mathrm{i}}s^{17}_{140}
$;\ \ 
$-\frac{4}{35}c^{1}_{140}
-\frac{3}{35}c^{3}_{140}
-\frac{1}{7}c^{5}_{140}
+\frac{1}{35}c^{7}_{140}
+\frac{1}{35}c^{9}_{140}
+\frac{4}{35}c^{13}_{140}
+\frac{2}{35}c^{15}_{140}
-\frac{3}{35}c^{17}_{140}
+\frac{9}{35}c^{19}_{140}
-\frac{4}{35}c^{21}_{140}
-\frac{2}{7}c^{23}_{140}
$,
$\frac{1}{\sqrt{35}}c^{1}_{35}
-\frac{1}{\sqrt{35}}c^{6}_{35}
$;\ \ 
$\frac{2}{\sqrt{35}}c^{3}_{35}
+\frac{1}{\sqrt{35}}c^{4}_{35}
+\frac{1}{\sqrt{35}}c^{10}_{35}
+\frac{1}{\sqrt{35}}c^{11}_{35}
$)

Pass. 

 \ \color{black}

 \color{blue}

\noindent 675: (dims,levels) = $(6;105
)$,
irreps = $3_{7}^{1}
\hskip -1.5pt \otimes \hskip -1.5pt
2_{5}^{1}
\hskip -1.5pt \otimes \hskip -1.5pt
1_{3}^{1,0}$,
pord$(\rho_\text{isum}(\mathfrak{t})) = 35$,

\vskip 0.7ex
\hangindent=5.5em \hangafter=1
{\white .}\hskip 1em $\rho_\text{isum}(\mathfrak{t})$ =
 $( \frac{11}{105},
\frac{29}{105},
\frac{44}{105},
\frac{71}{105},
\frac{74}{105},
\frac{86}{105} )
$,

\vskip 0.7ex
\hangindent=5.5em \hangafter=1
{\white .}\hskip 1em $\rho_\text{isum}(\mathfrak{s})$ =
$\mathrm{i}$($-\frac{2}{35}c^{1}_{140}
+\frac{1}{35}c^{3}_{140}
+\frac{1}{7}c^{5}_{140}
+\frac{3}{35}c^{7}_{140}
-\frac{1}{5}c^{9}_{140}
+\frac{2}{35}c^{13}_{140}
+\frac{1}{35}c^{15}_{140}
+\frac{1}{35}c^{17}_{140}
-\frac{3}{35}c^{19}_{140}
-\frac{2}{35}c^{21}_{140}
+\frac{2}{7}c^{23}_{140}
$,
$-\frac{1}{\sqrt{35}\mathrm{i}}s^{3}_{140}
-\frac{1}{\sqrt{35}\mathrm{i}}s^{17}_{140}
$,
$\frac{2}{\sqrt{35}}c^{3}_{35}
+\frac{1}{\sqrt{35}}c^{4}_{35}
+\frac{1}{\sqrt{35}}c^{10}_{35}
+\frac{1}{\sqrt{35}}c^{11}_{35}
$,
$\frac{1}{\sqrt{35}}c^{1}_{35}
-\frac{1}{\sqrt{35}}c^{6}_{35}
$,
$\frac{4}{35}c^{1}_{140}
+\frac{3}{35}c^{3}_{140}
+\frac{1}{7}c^{5}_{140}
-\frac{1}{35}c^{7}_{140}
-\frac{1}{35}c^{9}_{140}
-\frac{4}{35}c^{13}_{140}
-\frac{2}{35}c^{15}_{140}
+\frac{3}{35}c^{17}_{140}
-\frac{9}{35}c^{19}_{140}
+\frac{4}{35}c^{21}_{140}
+\frac{2}{7}c^{23}_{140}
$,
$-\frac{1}{\sqrt{35}}c^{4}_{35}
+\frac{1}{\sqrt{35}}c^{11}_{35}
$;\ \ 
$-\frac{1}{\sqrt{35}}c^{4}_{35}
+\frac{1}{\sqrt{35}}c^{11}_{35}
$,
$\frac{2}{35}c^{1}_{140}
-\frac{1}{35}c^{3}_{140}
-\frac{1}{7}c^{5}_{140}
-\frac{3}{35}c^{7}_{140}
+\frac{1}{5}c^{9}_{140}
-\frac{2}{35}c^{13}_{140}
-\frac{1}{35}c^{15}_{140}
-\frac{1}{35}c^{17}_{140}
+\frac{3}{35}c^{19}_{140}
+\frac{2}{35}c^{21}_{140}
-\frac{2}{7}c^{23}_{140}
$,
$-\frac{2}{\sqrt{35}}c^{3}_{35}
-\frac{1}{\sqrt{35}}c^{4}_{35}
-\frac{1}{\sqrt{35}}c^{10}_{35}
-\frac{1}{\sqrt{35}}c^{11}_{35}
$,
$\frac{1}{\sqrt{35}}c^{1}_{35}
-\frac{1}{\sqrt{35}}c^{6}_{35}
$,
$-\frac{4}{35}c^{1}_{140}
-\frac{3}{35}c^{3}_{140}
-\frac{1}{7}c^{5}_{140}
+\frac{1}{35}c^{7}_{140}
+\frac{1}{35}c^{9}_{140}
+\frac{4}{35}c^{13}_{140}
+\frac{2}{35}c^{15}_{140}
-\frac{3}{35}c^{17}_{140}
+\frac{9}{35}c^{19}_{140}
-\frac{4}{35}c^{21}_{140}
-\frac{2}{7}c^{23}_{140}
$;\ \ 
$\frac{1}{\sqrt{35}}c^{1}_{35}
-\frac{1}{\sqrt{35}}c^{6}_{35}
$,
$-\frac{4}{35}c^{1}_{140}
-\frac{3}{35}c^{3}_{140}
-\frac{1}{7}c^{5}_{140}
+\frac{1}{35}c^{7}_{140}
+\frac{1}{35}c^{9}_{140}
+\frac{4}{35}c^{13}_{140}
+\frac{2}{35}c^{15}_{140}
-\frac{3}{35}c^{17}_{140}
+\frac{9}{35}c^{19}_{140}
-\frac{4}{35}c^{21}_{140}
-\frac{2}{7}c^{23}_{140}
$,
$-\frac{1}{\sqrt{35}}c^{4}_{35}
+\frac{1}{\sqrt{35}}c^{11}_{35}
$,
$\frac{1}{\sqrt{35}\mathrm{i}}s^{3}_{140}
+\frac{1}{\sqrt{35}\mathrm{i}}s^{17}_{140}
$;\ \ 
$\frac{1}{\sqrt{35}}c^{4}_{35}
-\frac{1}{\sqrt{35}}c^{11}_{35}
$,
$\frac{1}{\sqrt{35}\mathrm{i}}s^{3}_{140}
+\frac{1}{\sqrt{35}\mathrm{i}}s^{17}_{140}
$,
$-\frac{2}{35}c^{1}_{140}
+\frac{1}{35}c^{3}_{140}
+\frac{1}{7}c^{5}_{140}
+\frac{3}{35}c^{7}_{140}
-\frac{1}{5}c^{9}_{140}
+\frac{2}{35}c^{13}_{140}
+\frac{1}{35}c^{15}_{140}
+\frac{1}{35}c^{17}_{140}
-\frac{3}{35}c^{19}_{140}
-\frac{2}{35}c^{21}_{140}
+\frac{2}{7}c^{23}_{140}
$;\ \ 
$\frac{2}{35}c^{1}_{140}
-\frac{1}{35}c^{3}_{140}
-\frac{1}{7}c^{5}_{140}
-\frac{3}{35}c^{7}_{140}
+\frac{1}{5}c^{9}_{140}
-\frac{2}{35}c^{13}_{140}
-\frac{1}{35}c^{15}_{140}
-\frac{1}{35}c^{17}_{140}
+\frac{3}{35}c^{19}_{140}
+\frac{2}{35}c^{21}_{140}
-\frac{2}{7}c^{23}_{140}
$,
$-\frac{2}{\sqrt{35}}c^{3}_{35}
-\frac{1}{\sqrt{35}}c^{4}_{35}
-\frac{1}{\sqrt{35}}c^{10}_{35}
-\frac{1}{\sqrt{35}}c^{11}_{35}
$;\ \ 
$-\frac{1}{\sqrt{35}}c^{1}_{35}
+\frac{1}{\sqrt{35}}c^{6}_{35}
$)

Pass. 

 \ \color{black}

 \color{blue}

\noindent 676: (dims,levels) = $(6;105
)$,
irreps = $3_{7}^{3}
\hskip -1.5pt \otimes \hskip -1.5pt
2_{5}^{2}
\hskip -1.5pt \otimes \hskip -1.5pt
1_{3}^{1,0}$,
pord$(\rho_\text{isum}(\mathfrak{t})) = 35$,

\vskip 0.7ex
\hangindent=5.5em \hangafter=1
{\white .}\hskip 1em $\rho_\text{isum}(\mathfrak{t})$ =
 $( \frac{17}{105},
\frac{38}{105},
\frac{47}{105},
\frac{62}{105},
\frac{68}{105},
\frac{83}{105} )
$,

\vskip 0.7ex
\hangindent=5.5em \hangafter=1
{\white .}\hskip 1em $\rho_\text{isum}(\mathfrak{s})$ =
$\mathrm{i}$($\frac{4}{35}c^{1}_{140}
+\frac{3}{35}c^{3}_{140}
+\frac{1}{7}c^{5}_{140}
-\frac{1}{35}c^{7}_{140}
-\frac{1}{35}c^{9}_{140}
-\frac{4}{35}c^{13}_{140}
-\frac{2}{35}c^{15}_{140}
+\frac{3}{35}c^{17}_{140}
-\frac{9}{35}c^{19}_{140}
+\frac{4}{35}c^{21}_{140}
+\frac{2}{7}c^{23}_{140}
$,
$\frac{2}{35}c^{1}_{140}
-\frac{1}{35}c^{3}_{140}
-\frac{1}{7}c^{5}_{140}
-\frac{3}{35}c^{7}_{140}
+\frac{1}{5}c^{9}_{140}
-\frac{2}{35}c^{13}_{140}
-\frac{1}{35}c^{15}_{140}
-\frac{1}{35}c^{17}_{140}
+\frac{3}{35}c^{19}_{140}
+\frac{2}{35}c^{21}_{140}
-\frac{2}{7}c^{23}_{140}
$,
$\frac{2}{\sqrt{35}}c^{3}_{35}
+\frac{1}{\sqrt{35}}c^{4}_{35}
+\frac{1}{\sqrt{35}}c^{10}_{35}
+\frac{1}{\sqrt{35}}c^{11}_{35}
$,
$-\frac{1}{\sqrt{35}\mathrm{i}}s^{3}_{140}
-\frac{1}{\sqrt{35}\mathrm{i}}s^{17}_{140}
$,
$-\frac{1}{\sqrt{35}}c^{4}_{35}
+\frac{1}{\sqrt{35}}c^{11}_{35}
$,
$\frac{1}{\sqrt{35}}c^{1}_{35}
-\frac{1}{\sqrt{35}}c^{6}_{35}
$;\ \ 
$-\frac{4}{35}c^{1}_{140}
-\frac{3}{35}c^{3}_{140}
-\frac{1}{7}c^{5}_{140}
+\frac{1}{35}c^{7}_{140}
+\frac{1}{35}c^{9}_{140}
+\frac{4}{35}c^{13}_{140}
+\frac{2}{35}c^{15}_{140}
-\frac{3}{35}c^{17}_{140}
+\frac{9}{35}c^{19}_{140}
-\frac{4}{35}c^{21}_{140}
-\frac{2}{7}c^{23}_{140}
$,
$-\frac{1}{\sqrt{35}}c^{4}_{35}
+\frac{1}{\sqrt{35}}c^{11}_{35}
$,
$\frac{1}{\sqrt{35}}c^{1}_{35}
-\frac{1}{\sqrt{35}}c^{6}_{35}
$,
$-\frac{2}{\sqrt{35}}c^{3}_{35}
-\frac{1}{\sqrt{35}}c^{4}_{35}
-\frac{1}{\sqrt{35}}c^{10}_{35}
-\frac{1}{\sqrt{35}}c^{11}_{35}
$,
$\frac{1}{\sqrt{35}\mathrm{i}}s^{3}_{140}
+\frac{1}{\sqrt{35}\mathrm{i}}s^{17}_{140}
$;\ \ 
$-\frac{1}{\sqrt{35}\mathrm{i}}s^{3}_{140}
-\frac{1}{\sqrt{35}\mathrm{i}}s^{17}_{140}
$,
$\frac{4}{35}c^{1}_{140}
+\frac{3}{35}c^{3}_{140}
+\frac{1}{7}c^{5}_{140}
-\frac{1}{35}c^{7}_{140}
-\frac{1}{35}c^{9}_{140}
-\frac{4}{35}c^{13}_{140}
-\frac{2}{35}c^{15}_{140}
+\frac{3}{35}c^{17}_{140}
-\frac{9}{35}c^{19}_{140}
+\frac{4}{35}c^{21}_{140}
+\frac{2}{7}c^{23}_{140}
$,
$\frac{1}{\sqrt{35}}c^{1}_{35}
-\frac{1}{\sqrt{35}}c^{6}_{35}
$,
$\frac{2}{35}c^{1}_{140}
-\frac{1}{35}c^{3}_{140}
-\frac{1}{7}c^{5}_{140}
-\frac{3}{35}c^{7}_{140}
+\frac{1}{5}c^{9}_{140}
-\frac{2}{35}c^{13}_{140}
-\frac{1}{35}c^{15}_{140}
-\frac{1}{35}c^{17}_{140}
+\frac{3}{35}c^{19}_{140}
+\frac{2}{35}c^{21}_{140}
-\frac{2}{7}c^{23}_{140}
$;\ \ 
$\frac{2}{\sqrt{35}}c^{3}_{35}
+\frac{1}{\sqrt{35}}c^{4}_{35}
+\frac{1}{\sqrt{35}}c^{10}_{35}
+\frac{1}{\sqrt{35}}c^{11}_{35}
$,
$\frac{2}{35}c^{1}_{140}
-\frac{1}{35}c^{3}_{140}
-\frac{1}{7}c^{5}_{140}
-\frac{3}{35}c^{7}_{140}
+\frac{1}{5}c^{9}_{140}
-\frac{2}{35}c^{13}_{140}
-\frac{1}{35}c^{15}_{140}
-\frac{1}{35}c^{17}_{140}
+\frac{3}{35}c^{19}_{140}
+\frac{2}{35}c^{21}_{140}
-\frac{2}{7}c^{23}_{140}
$,
$-\frac{1}{\sqrt{35}}c^{4}_{35}
+\frac{1}{\sqrt{35}}c^{11}_{35}
$;\ \ 
$\frac{1}{\sqrt{35}\mathrm{i}}s^{3}_{140}
+\frac{1}{\sqrt{35}\mathrm{i}}s^{17}_{140}
$,
$-\frac{4}{35}c^{1}_{140}
-\frac{3}{35}c^{3}_{140}
-\frac{1}{7}c^{5}_{140}
+\frac{1}{35}c^{7}_{140}
+\frac{1}{35}c^{9}_{140}
+\frac{4}{35}c^{13}_{140}
+\frac{2}{35}c^{15}_{140}
-\frac{3}{35}c^{17}_{140}
+\frac{9}{35}c^{19}_{140}
-\frac{4}{35}c^{21}_{140}
-\frac{2}{7}c^{23}_{140}
$;\ \ 
$-\frac{2}{\sqrt{35}}c^{3}_{35}
-\frac{1}{\sqrt{35}}c^{4}_{35}
-\frac{1}{\sqrt{35}}c^{10}_{35}
-\frac{1}{\sqrt{35}}c^{11}_{35}
$)

Pass. 

 \ \color{black}

 \color{blue}

\noindent 677: (dims,levels) = $(6;105
)$,
irreps = $3_{7}^{3}
\hskip -1.5pt \otimes \hskip -1.5pt
2_{5}^{1}
\hskip -1.5pt \otimes \hskip -1.5pt
1_{3}^{1,0}$,
pord$(\rho_\text{isum}(\mathfrak{t})) = 35$,

\vskip 0.7ex
\hangindent=5.5em \hangafter=1
{\white .}\hskip 1em $\rho_\text{isum}(\mathfrak{t})$ =
 $( \frac{26}{105},
\frac{41}{105},
\frac{59}{105},
\frac{89}{105},
\frac{101}{105},
\frac{104}{105} )
$,

\vskip 0.7ex
\hangindent=5.5em \hangafter=1
{\white .}\hskip 1em $\rho_\text{isum}(\mathfrak{s})$ =
$\mathrm{i}$($-\frac{1}{\sqrt{35}}c^{1}_{35}
+\frac{1}{\sqrt{35}}c^{6}_{35}
$,
$\frac{2}{35}c^{1}_{140}
-\frac{1}{35}c^{3}_{140}
-\frac{1}{7}c^{5}_{140}
-\frac{3}{35}c^{7}_{140}
+\frac{1}{5}c^{9}_{140}
-\frac{2}{35}c^{13}_{140}
-\frac{1}{35}c^{15}_{140}
-\frac{1}{35}c^{17}_{140}
+\frac{3}{35}c^{19}_{140}
+\frac{2}{35}c^{21}_{140}
-\frac{2}{7}c^{23}_{140}
$,
$\frac{2}{\sqrt{35}}c^{3}_{35}
+\frac{1}{\sqrt{35}}c^{4}_{35}
+\frac{1}{\sqrt{35}}c^{10}_{35}
+\frac{1}{\sqrt{35}}c^{11}_{35}
$,
$-\frac{1}{\sqrt{35}\mathrm{i}}s^{3}_{140}
-\frac{1}{\sqrt{35}\mathrm{i}}s^{17}_{140}
$,
$-\frac{1}{\sqrt{35}}c^{4}_{35}
+\frac{1}{\sqrt{35}}c^{11}_{35}
$,
$\frac{4}{35}c^{1}_{140}
+\frac{3}{35}c^{3}_{140}
+\frac{1}{7}c^{5}_{140}
-\frac{1}{35}c^{7}_{140}
-\frac{1}{35}c^{9}_{140}
-\frac{4}{35}c^{13}_{140}
-\frac{2}{35}c^{15}_{140}
+\frac{3}{35}c^{17}_{140}
-\frac{9}{35}c^{19}_{140}
+\frac{4}{35}c^{21}_{140}
+\frac{2}{7}c^{23}_{140}
$;\ \ 
$\frac{1}{\sqrt{35}}c^{4}_{35}
-\frac{1}{\sqrt{35}}c^{11}_{35}
$,
$\frac{1}{\sqrt{35}\mathrm{i}}s^{3}_{140}
+\frac{1}{\sqrt{35}\mathrm{i}}s^{17}_{140}
$,
$-\frac{4}{35}c^{1}_{140}
-\frac{3}{35}c^{3}_{140}
-\frac{1}{7}c^{5}_{140}
+\frac{1}{35}c^{7}_{140}
+\frac{1}{35}c^{9}_{140}
+\frac{4}{35}c^{13}_{140}
+\frac{2}{35}c^{15}_{140}
-\frac{3}{35}c^{17}_{140}
+\frac{9}{35}c^{19}_{140}
-\frac{4}{35}c^{21}_{140}
-\frac{2}{7}c^{23}_{140}
$,
$-\frac{1}{\sqrt{35}}c^{1}_{35}
+\frac{1}{\sqrt{35}}c^{6}_{35}
$,
$-\frac{2}{\sqrt{35}}c^{3}_{35}
-\frac{1}{\sqrt{35}}c^{4}_{35}
-\frac{1}{\sqrt{35}}c^{10}_{35}
-\frac{1}{\sqrt{35}}c^{11}_{35}
$;\ \ 
$\frac{2}{35}c^{1}_{140}
-\frac{1}{35}c^{3}_{140}
-\frac{1}{7}c^{5}_{140}
-\frac{3}{35}c^{7}_{140}
+\frac{1}{5}c^{9}_{140}
-\frac{2}{35}c^{13}_{140}
-\frac{1}{35}c^{15}_{140}
-\frac{1}{35}c^{17}_{140}
+\frac{3}{35}c^{19}_{140}
+\frac{2}{35}c^{21}_{140}
-\frac{2}{7}c^{23}_{140}
$,
$-\frac{1}{\sqrt{35}}c^{4}_{35}
+\frac{1}{\sqrt{35}}c^{11}_{35}
$,
$-\frac{4}{35}c^{1}_{140}
-\frac{3}{35}c^{3}_{140}
-\frac{1}{7}c^{5}_{140}
+\frac{1}{35}c^{7}_{140}
+\frac{1}{35}c^{9}_{140}
+\frac{4}{35}c^{13}_{140}
+\frac{2}{35}c^{15}_{140}
-\frac{3}{35}c^{17}_{140}
+\frac{9}{35}c^{19}_{140}
-\frac{4}{35}c^{21}_{140}
-\frac{2}{7}c^{23}_{140}
$,
$\frac{1}{\sqrt{35}}c^{1}_{35}
-\frac{1}{\sqrt{35}}c^{6}_{35}
$;\ \ 
$\frac{1}{\sqrt{35}}c^{1}_{35}
-\frac{1}{\sqrt{35}}c^{6}_{35}
$,
$-\frac{2}{\sqrt{35}}c^{3}_{35}
-\frac{1}{\sqrt{35}}c^{4}_{35}
-\frac{1}{\sqrt{35}}c^{10}_{35}
-\frac{1}{\sqrt{35}}c^{11}_{35}
$,
$\frac{2}{35}c^{1}_{140}
-\frac{1}{35}c^{3}_{140}
-\frac{1}{7}c^{5}_{140}
-\frac{3}{35}c^{7}_{140}
+\frac{1}{5}c^{9}_{140}
-\frac{2}{35}c^{13}_{140}
-\frac{1}{35}c^{15}_{140}
-\frac{1}{35}c^{17}_{140}
+\frac{3}{35}c^{19}_{140}
+\frac{2}{35}c^{21}_{140}
-\frac{2}{7}c^{23}_{140}
$;\ \ 
$-\frac{2}{35}c^{1}_{140}
+\frac{1}{35}c^{3}_{140}
+\frac{1}{7}c^{5}_{140}
+\frac{3}{35}c^{7}_{140}
-\frac{1}{5}c^{9}_{140}
+\frac{2}{35}c^{13}_{140}
+\frac{1}{35}c^{15}_{140}
+\frac{1}{35}c^{17}_{140}
-\frac{3}{35}c^{19}_{140}
-\frac{2}{35}c^{21}_{140}
+\frac{2}{7}c^{23}_{140}
$,
$\frac{1}{\sqrt{35}\mathrm{i}}s^{3}_{140}
+\frac{1}{\sqrt{35}\mathrm{i}}s^{17}_{140}
$;\ \ 
$-\frac{1}{\sqrt{35}}c^{4}_{35}
+\frac{1}{\sqrt{35}}c^{11}_{35}
$)

Pass. 

 \ \color{black}

\noindent 678: (dims,levels) = $(6;120
)$,
irreps = $3_{8}^{3,0}
\hskip -1.5pt \otimes \hskip -1.5pt
2_{5}^{1}
\hskip -1.5pt \otimes \hskip -1.5pt
1_{3}^{1,0}$,
pord$(\rho_\text{isum}(\mathfrak{t})) = 40$,

\vskip 0.7ex
\hangindent=5.5em \hangafter=1
{\white .}\hskip 1em $\rho_\text{isum}(\mathfrak{t})$ =
 $( \frac{2}{15},
\frac{8}{15},
\frac{1}{120},
\frac{49}{120},
\frac{61}{120},
\frac{109}{120} )
$,

\vskip 0.7ex
\hangindent=5.5em \hangafter=1
{\white .}\hskip 1em $\rho_\text{isum}(\mathfrak{s})$ =
($0$,
$0$,
$\frac{1}{\sqrt{10}}c^{3}_{20}
$,
$\frac{1}{\sqrt{10}}c^{1}_{20}
$,
$\frac{1}{\sqrt{10}}c^{3}_{20}
$,
$\frac{1}{\sqrt{10}}c^{1}_{20}
$;
$0$,
$\frac{1}{\sqrt{10}}c^{1}_{20}
$,
$-\frac{1}{\sqrt{10}}c^{3}_{20}
$,
$\frac{1}{\sqrt{10}}c^{1}_{20}
$,
$-\frac{1}{\sqrt{10}}c^{3}_{20}
$;
$-\frac{1}{2\sqrt{5}}c^{3}_{20}
$,
$-\frac{1}{2\sqrt{5}}c^{1}_{20}
$,
$\frac{1}{2\sqrt{5}}c^{3}_{20}
$,
$\frac{1}{2\sqrt{5}}c^{1}_{20}
$;
$\frac{1}{2\sqrt{5}}c^{3}_{20}
$,
$\frac{1}{2\sqrt{5}}c^{1}_{20}
$,
$-\frac{1}{2\sqrt{5}}c^{3}_{20}
$;
$-\frac{1}{2\sqrt{5}}c^{3}_{20}
$,
$-\frac{1}{2\sqrt{5}}c^{1}_{20}
$;
$\frac{1}{2\sqrt{5}}c^{3}_{20}
$)

Fail:
$\sigma(\rho(\mathfrak s)_\mathrm{ndeg}) \neq
 (\rho(\mathfrak t)^a \rho(\mathfrak s) \rho(\mathfrak t)^b
 \rho(\mathfrak s) \rho(\mathfrak t)^a)_\mathrm{ndeg}$,
 $\sigma = a$ = 11. Prop. B.5 (3) eqn. (B.25)

 \ \color{black}

\noindent 679: (dims,levels) = $(6;120
)$,
irreps = $3_{8}^{1,0}
\hskip -1.5pt \otimes \hskip -1.5pt
2_{5}^{1}
\hskip -1.5pt \otimes \hskip -1.5pt
1_{3}^{1,0}$,
pord$(\rho_\text{isum}(\mathfrak{t})) = 40$,

\vskip 0.7ex
\hangindent=5.5em \hangafter=1
{\white .}\hskip 1em $\rho_\text{isum}(\mathfrak{t})$ =
 $( \frac{2}{15},
\frac{8}{15},
\frac{19}{120},
\frac{31}{120},
\frac{79}{120},
\frac{91}{120} )
$,

\vskip 0.7ex
\hangindent=5.5em \hangafter=1
{\white .}\hskip 1em $\rho_\text{isum}(\mathfrak{s})$ =
($0$,
$0$,
$\frac{1}{\sqrt{10}}c^{1}_{20}
$,
$\frac{1}{\sqrt{10}}c^{3}_{20}
$,
$\frac{1}{\sqrt{10}}c^{1}_{20}
$,
$\frac{1}{\sqrt{10}}c^{3}_{20}
$;
$0$,
$\frac{1}{\sqrt{10}}c^{3}_{20}
$,
$-\frac{1}{\sqrt{10}}c^{1}_{20}
$,
$\frac{1}{\sqrt{10}}c^{3}_{20}
$,
$-\frac{1}{\sqrt{10}}c^{1}_{20}
$;
$-\frac{1}{2\sqrt{5}}c^{3}_{20}
$,
$-\frac{1}{2\sqrt{5}}c^{1}_{20}
$,
$\frac{1}{2\sqrt{5}}c^{3}_{20}
$,
$\frac{1}{2\sqrt{5}}c^{1}_{20}
$;
$\frac{1}{2\sqrt{5}}c^{3}_{20}
$,
$\frac{1}{2\sqrt{5}}c^{1}_{20}
$,
$-\frac{1}{2\sqrt{5}}c^{3}_{20}
$;
$-\frac{1}{2\sqrt{5}}c^{3}_{20}
$,
$-\frac{1}{2\sqrt{5}}c^{1}_{20}
$;
$\frac{1}{2\sqrt{5}}c^{3}_{20}
$)

Fail:
$\sigma(\rho(\mathfrak s)_\mathrm{ndeg}) \neq
 (\rho(\mathfrak t)^a \rho(\mathfrak s) \rho(\mathfrak t)^b
 \rho(\mathfrak s) \rho(\mathfrak t)^a)_\mathrm{ndeg}$,
 $\sigma = a$ = 11. Prop. B.5 (3) eqn. (B.25)

 \ \color{black}

\noindent 680: (dims,levels) = $(6;120
)$,
irreps = $3_{5}^{3}
\hskip -1.5pt \otimes \hskip -1.5pt
2_{8}^{1,0}
\hskip -1.5pt \otimes \hskip -1.5pt
1_{3}^{1,0}$,
pord$(\rho_\text{isum}(\mathfrak{t})) = 20$,

\vskip 0.7ex
\hangindent=5.5em \hangafter=1
{\white .}\hskip 1em $\rho_\text{isum}(\mathfrak{t})$ =
 $( \frac{11}{24},
\frac{17}{24},
\frac{7}{120},
\frac{13}{120},
\frac{37}{120},
\frac{103}{120} )
$,

\vskip 0.7ex
\hangindent=5.5em \hangafter=1
{\white .}\hskip 1em $\rho_\text{isum}(\mathfrak{s})$ =
($\sqrt{\frac{1}{10}}$,
$-\sqrt{\frac{1}{10}}$,
$\sqrt{\frac{1}{5}}$,
$\sqrt{\frac{1}{5}}$,
$\sqrt{\frac{1}{5}}$,
$\sqrt{\frac{1}{5}}$;
$-\sqrt{\frac{1}{10}}$,
$-\sqrt{\frac{1}{5}}$,
$\sqrt{\frac{1}{5}}$,
$\sqrt{\frac{1}{5}}$,
$-\sqrt{\frac{1}{5}}$;
$\frac{1}{\sqrt{10}}c^{1}_{5}
$,
$-\frac{1}{\sqrt{10}\mathrm{i}}s^{3}_{20}
$,
$\frac{1}{\sqrt{10}}c^{1}_{5}
$,
$-\frac{1}{\sqrt{10}\mathrm{i}}s^{3}_{20}
$;
$-\frac{1}{\sqrt{10}}c^{1}_{5}
$,
$\frac{1}{\sqrt{10}\mathrm{i}}s^{3}_{20}
$,
$\frac{1}{\sqrt{10}}c^{1}_{5}
$;
$-\frac{1}{\sqrt{10}}c^{1}_{5}
$,
$-\frac{1}{\sqrt{10}\mathrm{i}}s^{3}_{20}
$;
$\frac{1}{\sqrt{10}}c^{1}_{5}
$)

Fail:
$\sigma(\rho(\mathfrak s)_\mathrm{ndeg}) \neq
 (\rho(\mathfrak t)^a \rho(\mathfrak s) \rho(\mathfrak t)^b
 \rho(\mathfrak s) \rho(\mathfrak t)^a)_\mathrm{ndeg}$,
 $\sigma = a$ = 11. Prop. B.5 (3) eqn. (B.25)

 \ \color{black}

\noindent 681: (dims,levels) = $(6;120
)$,
irreps = $3_{5}^{1}
\hskip -1.5pt \otimes \hskip -1.5pt
2_{8}^{1,0}
\hskip -1.5pt \otimes \hskip -1.5pt
1_{3}^{1,0}$,
pord$(\rho_\text{isum}(\mathfrak{t})) = 20$,

\vskip 0.7ex
\hangindent=5.5em \hangafter=1
{\white .}\hskip 1em $\rho_\text{isum}(\mathfrak{t})$ =
 $( \frac{11}{24},
\frac{17}{24},
\frac{31}{120},
\frac{61}{120},
\frac{79}{120},
\frac{109}{120} )
$,

\vskip 0.7ex
\hangindent=5.5em \hangafter=1
{\white .}\hskip 1em $\rho_\text{isum}(\mathfrak{s})$ =
($-\sqrt{\frac{1}{10}}$,
$-\sqrt{\frac{1}{10}}$,
$\sqrt{\frac{1}{5}}$,
$\sqrt{\frac{1}{5}}$,
$\sqrt{\frac{1}{5}}$,
$\sqrt{\frac{1}{5}}$;
$\sqrt{\frac{1}{10}}$,
$\sqrt{\frac{1}{5}}$,
$-\sqrt{\frac{1}{5}}$,
$\sqrt{\frac{1}{5}}$,
$-\sqrt{\frac{1}{5}}$;
$\frac{1}{\sqrt{10}\mathrm{i}}s^{3}_{20}
$,
$\frac{1}{\sqrt{10}\mathrm{i}}s^{3}_{20}
$,
$-\frac{1}{\sqrt{10}}c^{1}_{5}
$,
$-\frac{1}{\sqrt{10}}c^{1}_{5}
$;
$-\frac{1}{\sqrt{10}\mathrm{i}}s^{3}_{20}
$,
$-\frac{1}{\sqrt{10}}c^{1}_{5}
$,
$\frac{1}{\sqrt{10}}c^{1}_{5}
$;
$\frac{1}{\sqrt{10}\mathrm{i}}s^{3}_{20}
$,
$\frac{1}{\sqrt{10}\mathrm{i}}s^{3}_{20}
$;
$-\frac{1}{\sqrt{10}\mathrm{i}}s^{3}_{20}
$)

Fail:
$\sigma(\rho(\mathfrak s)_\mathrm{ndeg}) \neq
 (\rho(\mathfrak t)^a \rho(\mathfrak s) \rho(\mathfrak t)^b
 \rho(\mathfrak s) \rho(\mathfrak t)^a)_\mathrm{ndeg}$,
 $\sigma = a$ = 11. Prop. B.5 (3) eqn. (B.25)

 \ \color{black}

\noindent 682: (dims,levels) = $(6;120
)$,
irreps = $3_{8}^{1,0}
\hskip -1.5pt \otimes \hskip -1.5pt
2_{5}^{2}
\hskip -1.5pt \otimes \hskip -1.5pt
1_{3}^{1,0}$,
pord$(\rho_\text{isum}(\mathfrak{t})) = 40$,

\vskip 0.7ex
\hangindent=5.5em \hangafter=1
{\white .}\hskip 1em $\rho_\text{isum}(\mathfrak{t})$ =
 $( \frac{11}{15},
\frac{14}{15},
\frac{7}{120},
\frac{43}{120},
\frac{67}{120},
\frac{103}{120} )
$,

\vskip 0.7ex
\hangindent=5.5em \hangafter=1
{\white .}\hskip 1em $\rho_\text{isum}(\mathfrak{s})$ =
($0$,
$0$,
$\frac{1}{\sqrt{10}}c^{3}_{20}
$,
$\frac{1}{\sqrt{10}}c^{1}_{20}
$,
$\frac{1}{\sqrt{10}}c^{3}_{20}
$,
$\frac{1}{\sqrt{10}}c^{1}_{20}
$;
$0$,
$\frac{1}{\sqrt{10}}c^{1}_{20}
$,
$-\frac{1}{\sqrt{10}}c^{3}_{20}
$,
$\frac{1}{\sqrt{10}}c^{1}_{20}
$,
$-\frac{1}{\sqrt{10}}c^{3}_{20}
$;
$\frac{1}{2\sqrt{5}}c^{1}_{20}
$,
$\frac{1}{2\sqrt{5}}c^{3}_{20}
$,
$-\frac{1}{2\sqrt{5}}c^{1}_{20}
$,
$-\frac{1}{2\sqrt{5}}c^{3}_{20}
$;
$-\frac{1}{2\sqrt{5}}c^{1}_{20}
$,
$-\frac{1}{2\sqrt{5}}c^{3}_{20}
$,
$\frac{1}{2\sqrt{5}}c^{1}_{20}
$;
$\frac{1}{2\sqrt{5}}c^{1}_{20}
$,
$\frac{1}{2\sqrt{5}}c^{3}_{20}
$;
$-\frac{1}{2\sqrt{5}}c^{1}_{20}
$)

Fail:
$\sigma(\rho(\mathfrak s)_\mathrm{ndeg}) \neq
 (\rho(\mathfrak t)^a \rho(\mathfrak s) \rho(\mathfrak t)^b
 \rho(\mathfrak s) \rho(\mathfrak t)^a)_\mathrm{ndeg}$,
 $\sigma = a$ = 11. Prop. B.5 (3) eqn. (B.25)

 \ \color{black}

\noindent 683: (dims,levels) = $(6;120
)$,
irreps = $3_{8}^{3,0}
\hskip -1.5pt \otimes \hskip -1.5pt
2_{5}^{2}
\hskip -1.5pt \otimes \hskip -1.5pt
1_{3}^{1,0}$,
pord$(\rho_\text{isum}(\mathfrak{t})) = 40$,

\vskip 0.7ex
\hangindent=5.5em \hangafter=1
{\white .}\hskip 1em $\rho_\text{isum}(\mathfrak{t})$ =
 $( \frac{11}{15},
\frac{14}{15},
\frac{13}{120},
\frac{37}{120},
\frac{73}{120},
\frac{97}{120} )
$,

\vskip 0.7ex
\hangindent=5.5em \hangafter=1
{\white .}\hskip 1em $\rho_\text{isum}(\mathfrak{s})$ =
($0$,
$0$,
$\frac{1}{\sqrt{10}}c^{1}_{20}
$,
$\frac{1}{\sqrt{10}}c^{3}_{20}
$,
$\frac{1}{\sqrt{10}}c^{1}_{20}
$,
$\frac{1}{\sqrt{10}}c^{3}_{20}
$;
$0$,
$\frac{1}{\sqrt{10}}c^{3}_{20}
$,
$-\frac{1}{\sqrt{10}}c^{1}_{20}
$,
$\frac{1}{\sqrt{10}}c^{3}_{20}
$,
$-\frac{1}{\sqrt{10}}c^{1}_{20}
$;
$\frac{1}{2\sqrt{5}}c^{1}_{20}
$,
$\frac{1}{2\sqrt{5}}c^{3}_{20}
$,
$-\frac{1}{2\sqrt{5}}c^{1}_{20}
$,
$-\frac{1}{2\sqrt{5}}c^{3}_{20}
$;
$-\frac{1}{2\sqrt{5}}c^{1}_{20}
$,
$-\frac{1}{2\sqrt{5}}c^{3}_{20}
$,
$\frac{1}{2\sqrt{5}}c^{1}_{20}
$;
$\frac{1}{2\sqrt{5}}c^{1}_{20}
$,
$\frac{1}{2\sqrt{5}}c^{3}_{20}
$;
$-\frac{1}{2\sqrt{5}}c^{1}_{20}
$)

Fail:
$\sigma(\rho(\mathfrak s)_\mathrm{ndeg}) \neq
 (\rho(\mathfrak t)^a \rho(\mathfrak s) \rho(\mathfrak t)^b
 \rho(\mathfrak s) \rho(\mathfrak t)^a)_\mathrm{ndeg}$,
 $\sigma = a$ = 11. Prop. B.5 (3) eqn. (B.25)

 \ \color{black}

\noindent 684: (dims,levels) = $(6;132
)$,
irreps = $6_{11}^{1}
\hskip -1.5pt \otimes \hskip -1.5pt
1_{4}^{1,0}
\hskip -1.5pt \otimes \hskip -1.5pt
1_{3}^{1,0}$,
pord$(\rho_\text{isum}(\mathfrak{t})) = 11$,

\vskip 0.7ex
\hangindent=5.5em \hangafter=1
{\white .}\hskip 1em $\rho_\text{isum}(\mathfrak{t})$ =
 $( \frac{7}{12},
\frac{5}{132},
\frac{53}{132},
\frac{89}{132},
\frac{113}{132},
\frac{125}{132} )
$,

\vskip 0.7ex
\hangindent=5.5em \hangafter=1
{\white .}\hskip 1em $\rho_\text{isum}(\mathfrak{s})$ =
($\sqrt{\frac{1}{11}}$,
$\sqrt{\frac{2}{11}}$,
$\sqrt{\frac{2}{11}}$,
$\sqrt{\frac{2}{11}}$,
$\sqrt{\frac{2}{11}}$,
$\sqrt{\frac{2}{11}}$;
$\frac{1}{\sqrt{11}}c^{1}_{11}
$,
$\frac{1}{\sqrt{11}}c^{2}_{11}
$,
$\frac{1}{\sqrt{11}}c^{3}_{11}
$,
$\frac{1}{\sqrt{11}}c^{4}_{11}
$,
$-\frac{1}{\sqrt{11}\mathrm{i}}s^{9}_{44}
$;
$\frac{1}{\sqrt{11}}c^{4}_{11}
$,
$-\frac{1}{\sqrt{11}\mathrm{i}}s^{9}_{44}
$,
$\frac{1}{\sqrt{11}}c^{3}_{11}
$,
$\frac{1}{\sqrt{11}}c^{1}_{11}
$;
$\frac{1}{\sqrt{11}}c^{2}_{11}
$,
$\frac{1}{\sqrt{11}}c^{1}_{11}
$,
$\frac{1}{\sqrt{11}}c^{4}_{11}
$;
$-\frac{1}{\sqrt{11}\mathrm{i}}s^{9}_{44}
$,
$\frac{1}{\sqrt{11}}c^{2}_{11}
$;
$\frac{1}{\sqrt{11}}c^{3}_{11}
$)

Fail:
cnd($\rho(\mathfrak s)_\mathrm{ndeg}$) = 88 does not divide
 ord($\rho(\mathfrak t)$)=132. Prop. B.4 (2)

 \ \color{black}

\noindent 685: (dims,levels) = $(6;132
)$,
irreps = $6_{11}^{7}
\hskip -1.5pt \otimes \hskip -1.5pt
1_{4}^{1,0}
\hskip -1.5pt \otimes \hskip -1.5pt
1_{3}^{1,0}$,
pord$(\rho_\text{isum}(\mathfrak{t})) = 11$,

\vskip 0.7ex
\hangindent=5.5em \hangafter=1
{\white .}\hskip 1em $\rho_\text{isum}(\mathfrak{t})$ =
 $( \frac{7}{12},
\frac{17}{132},
\frac{29}{132},
\frac{41}{132},
\frac{65}{132},
\frac{101}{132} )
$,

\vskip 0.7ex
\hangindent=5.5em \hangafter=1
{\white .}\hskip 1em $\rho_\text{isum}(\mathfrak{s})$ =
($-\sqrt{\frac{1}{11}}$,
$\sqrt{\frac{2}{11}}$,
$\sqrt{\frac{2}{11}}$,
$\sqrt{\frac{2}{11}}$,
$\sqrt{\frac{2}{11}}$,
$\sqrt{\frac{2}{11}}$;
$-\frac{1}{\sqrt{11}}c^{1}_{11}
$,
$\frac{1}{\sqrt{11}\mathrm{i}}s^{9}_{44}
$,
$-\frac{1}{\sqrt{11}}c^{4}_{11}
$,
$-\frac{1}{\sqrt{11}}c^{3}_{11}
$,
$-\frac{1}{\sqrt{11}}c^{2}_{11}
$;
$-\frac{1}{\sqrt{11}}c^{3}_{11}
$,
$-\frac{1}{\sqrt{11}}c^{2}_{11}
$,
$-\frac{1}{\sqrt{11}}c^{4}_{11}
$,
$-\frac{1}{\sqrt{11}}c^{1}_{11}
$;
$\frac{1}{\sqrt{11}\mathrm{i}}s^{9}_{44}
$,
$-\frac{1}{\sqrt{11}}c^{1}_{11}
$,
$-\frac{1}{\sqrt{11}}c^{3}_{11}
$;
$-\frac{1}{\sqrt{11}}c^{2}_{11}
$,
$\frac{1}{\sqrt{11}\mathrm{i}}s^{9}_{44}
$;
$-\frac{1}{\sqrt{11}}c^{4}_{11}
$)

Fail:
cnd($\rho(\mathfrak s)_\mathrm{ndeg}$) = 88 does not divide
 ord($\rho(\mathfrak t)$)=132. Prop. B.4 (2)

 \ \color{black}

 \color{blue}

\noindent 686: (dims,levels) = $(6;140
)$,
irreps = $3_{7}^{1}
\hskip -1.5pt \otimes \hskip -1.5pt
2_{5}^{1}
\hskip -1.5pt \otimes \hskip -1.5pt
1_{4}^{1,0}$,
pord$(\rho_\text{isum}(\mathfrak{t})) = 35$,

\vskip 0.7ex
\hangindent=5.5em \hangafter=1
{\white .}\hskip 1em $\rho_\text{isum}(\mathfrak{t})$ =
 $( \frac{3}{140},
\frac{27}{140},
\frac{47}{140},
\frac{83}{140},
\frac{87}{140},
\frac{103}{140} )
$,

\vskip 0.7ex
\hangindent=5.5em \hangafter=1
{\white .}\hskip 1em $\rho_\text{isum}(\mathfrak{s})$ =
($\frac{2}{35}c^{1}_{140}
-\frac{1}{35}c^{3}_{140}
-\frac{1}{7}c^{5}_{140}
-\frac{3}{35}c^{7}_{140}
+\frac{1}{5}c^{9}_{140}
-\frac{2}{35}c^{13}_{140}
-\frac{1}{35}c^{15}_{140}
-\frac{1}{35}c^{17}_{140}
+\frac{3}{35}c^{19}_{140}
+\frac{2}{35}c^{21}_{140}
-\frac{2}{7}c^{23}_{140}
$,
$-\frac{1}{\sqrt{35}\mathrm{i}}s^{3}_{140}
-\frac{1}{\sqrt{35}\mathrm{i}}s^{17}_{140}
$,
$\frac{2}{\sqrt{35}}c^{3}_{35}
+\frac{1}{\sqrt{35}}c^{4}_{35}
+\frac{1}{\sqrt{35}}c^{10}_{35}
+\frac{1}{\sqrt{35}}c^{11}_{35}
$,
$-\frac{1}{\sqrt{35}}c^{1}_{35}
+\frac{1}{\sqrt{35}}c^{6}_{35}
$,
$\frac{4}{35}c^{1}_{140}
+\frac{3}{35}c^{3}_{140}
+\frac{1}{7}c^{5}_{140}
-\frac{1}{35}c^{7}_{140}
-\frac{1}{35}c^{9}_{140}
-\frac{4}{35}c^{13}_{140}
-\frac{2}{35}c^{15}_{140}
+\frac{3}{35}c^{17}_{140}
-\frac{9}{35}c^{19}_{140}
+\frac{4}{35}c^{21}_{140}
+\frac{2}{7}c^{23}_{140}
$,
$\frac{1}{\sqrt{35}}c^{4}_{35}
-\frac{1}{\sqrt{35}}c^{11}_{35}
$;
$\frac{1}{\sqrt{35}}c^{4}_{35}
-\frac{1}{\sqrt{35}}c^{11}_{35}
$,
$-\frac{2}{35}c^{1}_{140}
+\frac{1}{35}c^{3}_{140}
+\frac{1}{7}c^{5}_{140}
+\frac{3}{35}c^{7}_{140}
-\frac{1}{5}c^{9}_{140}
+\frac{2}{35}c^{13}_{140}
+\frac{1}{35}c^{15}_{140}
+\frac{1}{35}c^{17}_{140}
-\frac{3}{35}c^{19}_{140}
-\frac{2}{35}c^{21}_{140}
+\frac{2}{7}c^{23}_{140}
$,
$-\frac{2}{\sqrt{35}}c^{3}_{35}
-\frac{1}{\sqrt{35}}c^{4}_{35}
-\frac{1}{\sqrt{35}}c^{10}_{35}
-\frac{1}{\sqrt{35}}c^{11}_{35}
$,
$-\frac{1}{\sqrt{35}}c^{1}_{35}
+\frac{1}{\sqrt{35}}c^{6}_{35}
$,
$-\frac{4}{35}c^{1}_{140}
-\frac{3}{35}c^{3}_{140}
-\frac{1}{7}c^{5}_{140}
+\frac{1}{35}c^{7}_{140}
+\frac{1}{35}c^{9}_{140}
+\frac{4}{35}c^{13}_{140}
+\frac{2}{35}c^{15}_{140}
-\frac{3}{35}c^{17}_{140}
+\frac{9}{35}c^{19}_{140}
-\frac{4}{35}c^{21}_{140}
-\frac{2}{7}c^{23}_{140}
$;
$-\frac{1}{\sqrt{35}}c^{1}_{35}
+\frac{1}{\sqrt{35}}c^{6}_{35}
$,
$-\frac{4}{35}c^{1}_{140}
-\frac{3}{35}c^{3}_{140}
-\frac{1}{7}c^{5}_{140}
+\frac{1}{35}c^{7}_{140}
+\frac{1}{35}c^{9}_{140}
+\frac{4}{35}c^{13}_{140}
+\frac{2}{35}c^{15}_{140}
-\frac{3}{35}c^{17}_{140}
+\frac{9}{35}c^{19}_{140}
-\frac{4}{35}c^{21}_{140}
-\frac{2}{7}c^{23}_{140}
$,
$\frac{1}{\sqrt{35}}c^{4}_{35}
-\frac{1}{\sqrt{35}}c^{11}_{35}
$,
$\frac{1}{\sqrt{35}\mathrm{i}}s^{3}_{140}
+\frac{1}{\sqrt{35}\mathrm{i}}s^{17}_{140}
$;
$-\frac{1}{\sqrt{35}}c^{4}_{35}
+\frac{1}{\sqrt{35}}c^{11}_{35}
$,
$\frac{1}{\sqrt{35}\mathrm{i}}s^{3}_{140}
+\frac{1}{\sqrt{35}\mathrm{i}}s^{17}_{140}
$,
$\frac{2}{35}c^{1}_{140}
-\frac{1}{35}c^{3}_{140}
-\frac{1}{7}c^{5}_{140}
-\frac{3}{35}c^{7}_{140}
+\frac{1}{5}c^{9}_{140}
-\frac{2}{35}c^{13}_{140}
-\frac{1}{35}c^{15}_{140}
-\frac{1}{35}c^{17}_{140}
+\frac{3}{35}c^{19}_{140}
+\frac{2}{35}c^{21}_{140}
-\frac{2}{7}c^{23}_{140}
$;
$-\frac{2}{35}c^{1}_{140}
+\frac{1}{35}c^{3}_{140}
+\frac{1}{7}c^{5}_{140}
+\frac{3}{35}c^{7}_{140}
-\frac{1}{5}c^{9}_{140}
+\frac{2}{35}c^{13}_{140}
+\frac{1}{35}c^{15}_{140}
+\frac{1}{35}c^{17}_{140}
-\frac{3}{35}c^{19}_{140}
-\frac{2}{35}c^{21}_{140}
+\frac{2}{7}c^{23}_{140}
$,
$-\frac{2}{\sqrt{35}}c^{3}_{35}
-\frac{1}{\sqrt{35}}c^{4}_{35}
-\frac{1}{\sqrt{35}}c^{10}_{35}
-\frac{1}{\sqrt{35}}c^{11}_{35}
$;
$\frac{1}{\sqrt{35}}c^{1}_{35}
-\frac{1}{\sqrt{35}}c^{6}_{35}
$)

Pass. 

 \ \color{black}

 \color{blue}

\noindent 687: (dims,levels) = $(6;140
)$,
irreps = $3_{7}^{3}
\hskip -1.5pt \otimes \hskip -1.5pt
2_{5}^{2}
\hskip -1.5pt \otimes \hskip -1.5pt
1_{4}^{1,0}$,
pord$(\rho_\text{isum}(\mathfrak{t})) = 35$,

\vskip 0.7ex
\hangindent=5.5em \hangafter=1
{\white .}\hskip 1em $\rho_\text{isum}(\mathfrak{t})$ =
 $( \frac{11}{140},
\frac{39}{140},
\frac{51}{140},
\frac{71}{140},
\frac{79}{140},
\frac{99}{140} )
$,

\vskip 0.7ex
\hangindent=5.5em \hangafter=1
{\white .}\hskip 1em $\rho_\text{isum}(\mathfrak{s})$ =
($-\frac{4}{35}c^{1}_{140}
-\frac{3}{35}c^{3}_{140}
-\frac{1}{7}c^{5}_{140}
+\frac{1}{35}c^{7}_{140}
+\frac{1}{35}c^{9}_{140}
+\frac{4}{35}c^{13}_{140}
+\frac{2}{35}c^{15}_{140}
-\frac{3}{35}c^{17}_{140}
+\frac{9}{35}c^{19}_{140}
-\frac{4}{35}c^{21}_{140}
-\frac{2}{7}c^{23}_{140}
$,
$-\frac{2}{35}c^{1}_{140}
+\frac{1}{35}c^{3}_{140}
+\frac{1}{7}c^{5}_{140}
+\frac{3}{35}c^{7}_{140}
-\frac{1}{5}c^{9}_{140}
+\frac{2}{35}c^{13}_{140}
+\frac{1}{35}c^{15}_{140}
+\frac{1}{35}c^{17}_{140}
-\frac{3}{35}c^{19}_{140}
-\frac{2}{35}c^{21}_{140}
+\frac{2}{7}c^{23}_{140}
$,
$\frac{2}{\sqrt{35}}c^{3}_{35}
+\frac{1}{\sqrt{35}}c^{4}_{35}
+\frac{1}{\sqrt{35}}c^{10}_{35}
+\frac{1}{\sqrt{35}}c^{11}_{35}
$,
$-\frac{1}{\sqrt{35}\mathrm{i}}s^{3}_{140}
-\frac{1}{\sqrt{35}\mathrm{i}}s^{17}_{140}
$,
$\frac{1}{\sqrt{35}}c^{4}_{35}
-\frac{1}{\sqrt{35}}c^{11}_{35}
$,
$-\frac{1}{\sqrt{35}}c^{1}_{35}
+\frac{1}{\sqrt{35}}c^{6}_{35}
$;
$\frac{4}{35}c^{1}_{140}
+\frac{3}{35}c^{3}_{140}
+\frac{1}{7}c^{5}_{140}
-\frac{1}{35}c^{7}_{140}
-\frac{1}{35}c^{9}_{140}
-\frac{4}{35}c^{13}_{140}
-\frac{2}{35}c^{15}_{140}
+\frac{3}{35}c^{17}_{140}
-\frac{9}{35}c^{19}_{140}
+\frac{4}{35}c^{21}_{140}
+\frac{2}{7}c^{23}_{140}
$,
$-\frac{1}{\sqrt{35}}c^{4}_{35}
+\frac{1}{\sqrt{35}}c^{11}_{35}
$,
$\frac{1}{\sqrt{35}}c^{1}_{35}
-\frac{1}{\sqrt{35}}c^{6}_{35}
$,
$\frac{2}{\sqrt{35}}c^{3}_{35}
+\frac{1}{\sqrt{35}}c^{4}_{35}
+\frac{1}{\sqrt{35}}c^{10}_{35}
+\frac{1}{\sqrt{35}}c^{11}_{35}
$,
$-\frac{1}{\sqrt{35}\mathrm{i}}s^{3}_{140}
-\frac{1}{\sqrt{35}\mathrm{i}}s^{17}_{140}
$;
$\frac{1}{\sqrt{35}\mathrm{i}}s^{3}_{140}
+\frac{1}{\sqrt{35}\mathrm{i}}s^{17}_{140}
$,
$-\frac{4}{35}c^{1}_{140}
-\frac{3}{35}c^{3}_{140}
-\frac{1}{7}c^{5}_{140}
+\frac{1}{35}c^{7}_{140}
+\frac{1}{35}c^{9}_{140}
+\frac{4}{35}c^{13}_{140}
+\frac{2}{35}c^{15}_{140}
-\frac{3}{35}c^{17}_{140}
+\frac{9}{35}c^{19}_{140}
-\frac{4}{35}c^{21}_{140}
-\frac{2}{7}c^{23}_{140}
$,
$\frac{1}{\sqrt{35}}c^{1}_{35}
-\frac{1}{\sqrt{35}}c^{6}_{35}
$,
$\frac{2}{35}c^{1}_{140}
-\frac{1}{35}c^{3}_{140}
-\frac{1}{7}c^{5}_{140}
-\frac{3}{35}c^{7}_{140}
+\frac{1}{5}c^{9}_{140}
-\frac{2}{35}c^{13}_{140}
-\frac{1}{35}c^{15}_{140}
-\frac{1}{35}c^{17}_{140}
+\frac{3}{35}c^{19}_{140}
+\frac{2}{35}c^{21}_{140}
-\frac{2}{7}c^{23}_{140}
$;
$-\frac{2}{\sqrt{35}}c^{3}_{35}
-\frac{1}{\sqrt{35}}c^{4}_{35}
-\frac{1}{\sqrt{35}}c^{10}_{35}
-\frac{1}{\sqrt{35}}c^{11}_{35}
$,
$\frac{2}{35}c^{1}_{140}
-\frac{1}{35}c^{3}_{140}
-\frac{1}{7}c^{5}_{140}
-\frac{3}{35}c^{7}_{140}
+\frac{1}{5}c^{9}_{140}
-\frac{2}{35}c^{13}_{140}
-\frac{1}{35}c^{15}_{140}
-\frac{1}{35}c^{17}_{140}
+\frac{3}{35}c^{19}_{140}
+\frac{2}{35}c^{21}_{140}
-\frac{2}{7}c^{23}_{140}
$,
$-\frac{1}{\sqrt{35}}c^{4}_{35}
+\frac{1}{\sqrt{35}}c^{11}_{35}
$;
$-\frac{1}{\sqrt{35}\mathrm{i}}s^{3}_{140}
-\frac{1}{\sqrt{35}\mathrm{i}}s^{17}_{140}
$,
$\frac{4}{35}c^{1}_{140}
+\frac{3}{35}c^{3}_{140}
+\frac{1}{7}c^{5}_{140}
-\frac{1}{35}c^{7}_{140}
-\frac{1}{35}c^{9}_{140}
-\frac{4}{35}c^{13}_{140}
-\frac{2}{35}c^{15}_{140}
+\frac{3}{35}c^{17}_{140}
-\frac{9}{35}c^{19}_{140}
+\frac{4}{35}c^{21}_{140}
+\frac{2}{7}c^{23}_{140}
$;
$\frac{2}{\sqrt{35}}c^{3}_{35}
+\frac{1}{\sqrt{35}}c^{4}_{35}
+\frac{1}{\sqrt{35}}c^{10}_{35}
+\frac{1}{\sqrt{35}}c^{11}_{35}
$)

Pass. 

 \ \color{black}

 \color{blue}

\noindent 688: (dims,levels) = $(6;140
)$,
irreps = $3_{7}^{1}
\hskip -1.5pt \otimes \hskip -1.5pt
2_{5}^{2}
\hskip -1.5pt \otimes \hskip -1.5pt
1_{4}^{1,0}$,
pord$(\rho_\text{isum}(\mathfrak{t})) = 35$,

\vskip 0.7ex
\hangindent=5.5em \hangafter=1
{\white .}\hskip 1em $\rho_\text{isum}(\mathfrak{t})$ =
 $( \frac{19}{140},
\frac{31}{140},
\frac{59}{140},
\frac{111}{140},
\frac{131}{140},
\frac{139}{140} )
$,

\vskip 0.7ex
\hangindent=5.5em \hangafter=1
{\white .}\hskip 1em $\rho_\text{isum}(\mathfrak{s})$ =
($-\frac{1}{\sqrt{35}\mathrm{i}}s^{3}_{140}
-\frac{1}{\sqrt{35}\mathrm{i}}s^{17}_{140}
$,
$\frac{1}{\sqrt{35}}c^{4}_{35}
-\frac{1}{\sqrt{35}}c^{11}_{35}
$,
$\frac{2}{\sqrt{35}}c^{3}_{35}
+\frac{1}{\sqrt{35}}c^{4}_{35}
+\frac{1}{\sqrt{35}}c^{10}_{35}
+\frac{1}{\sqrt{35}}c^{11}_{35}
$,
$-\frac{2}{35}c^{1}_{140}
+\frac{1}{35}c^{3}_{140}
+\frac{1}{7}c^{5}_{140}
+\frac{3}{35}c^{7}_{140}
-\frac{1}{5}c^{9}_{140}
+\frac{2}{35}c^{13}_{140}
+\frac{1}{35}c^{15}_{140}
+\frac{1}{35}c^{17}_{140}
-\frac{3}{35}c^{19}_{140}
-\frac{2}{35}c^{21}_{140}
+\frac{2}{7}c^{23}_{140}
$,
$-\frac{1}{\sqrt{35}}c^{1}_{35}
+\frac{1}{\sqrt{35}}c^{6}_{35}
$,
$\frac{4}{35}c^{1}_{140}
+\frac{3}{35}c^{3}_{140}
+\frac{1}{7}c^{5}_{140}
-\frac{1}{35}c^{7}_{140}
-\frac{1}{35}c^{9}_{140}
-\frac{4}{35}c^{13}_{140}
-\frac{2}{35}c^{15}_{140}
+\frac{3}{35}c^{17}_{140}
-\frac{9}{35}c^{19}_{140}
+\frac{4}{35}c^{21}_{140}
+\frac{2}{7}c^{23}_{140}
$;
$-\frac{4}{35}c^{1}_{140}
-\frac{3}{35}c^{3}_{140}
-\frac{1}{7}c^{5}_{140}
+\frac{1}{35}c^{7}_{140}
+\frac{1}{35}c^{9}_{140}
+\frac{4}{35}c^{13}_{140}
+\frac{2}{35}c^{15}_{140}
-\frac{3}{35}c^{17}_{140}
+\frac{9}{35}c^{19}_{140}
-\frac{4}{35}c^{21}_{140}
-\frac{2}{7}c^{23}_{140}
$,
$-\frac{2}{35}c^{1}_{140}
+\frac{1}{35}c^{3}_{140}
+\frac{1}{7}c^{5}_{140}
+\frac{3}{35}c^{7}_{140}
-\frac{1}{5}c^{9}_{140}
+\frac{2}{35}c^{13}_{140}
+\frac{1}{35}c^{15}_{140}
+\frac{1}{35}c^{17}_{140}
-\frac{3}{35}c^{19}_{140}
-\frac{2}{35}c^{21}_{140}
+\frac{2}{7}c^{23}_{140}
$,
$\frac{1}{\sqrt{35}\mathrm{i}}s^{3}_{140}
+\frac{1}{\sqrt{35}\mathrm{i}}s^{17}_{140}
$,
$-\frac{2}{\sqrt{35}}c^{3}_{35}
-\frac{1}{\sqrt{35}}c^{4}_{35}
-\frac{1}{\sqrt{35}}c^{10}_{35}
-\frac{1}{\sqrt{35}}c^{11}_{35}
$,
$-\frac{1}{\sqrt{35}}c^{1}_{35}
+\frac{1}{\sqrt{35}}c^{6}_{35}
$;
$\frac{4}{35}c^{1}_{140}
+\frac{3}{35}c^{3}_{140}
+\frac{1}{7}c^{5}_{140}
-\frac{1}{35}c^{7}_{140}
-\frac{1}{35}c^{9}_{140}
-\frac{4}{35}c^{13}_{140}
-\frac{2}{35}c^{15}_{140}
+\frac{3}{35}c^{17}_{140}
-\frac{9}{35}c^{19}_{140}
+\frac{4}{35}c^{21}_{140}
+\frac{2}{7}c^{23}_{140}
$,
$-\frac{1}{\sqrt{35}}c^{1}_{35}
+\frac{1}{\sqrt{35}}c^{6}_{35}
$,
$\frac{1}{\sqrt{35}}c^{4}_{35}
-\frac{1}{\sqrt{35}}c^{11}_{35}
$,
$-\frac{1}{\sqrt{35}\mathrm{i}}s^{3}_{140}
-\frac{1}{\sqrt{35}\mathrm{i}}s^{17}_{140}
$;
$-\frac{2}{\sqrt{35}}c^{3}_{35}
-\frac{1}{\sqrt{35}}c^{4}_{35}
-\frac{1}{\sqrt{35}}c^{10}_{35}
-\frac{1}{\sqrt{35}}c^{11}_{35}
$,
$-\frac{4}{35}c^{1}_{140}
-\frac{3}{35}c^{3}_{140}
-\frac{1}{7}c^{5}_{140}
+\frac{1}{35}c^{7}_{140}
+\frac{1}{35}c^{9}_{140}
+\frac{4}{35}c^{13}_{140}
+\frac{2}{35}c^{15}_{140}
-\frac{3}{35}c^{17}_{140}
+\frac{9}{35}c^{19}_{140}
-\frac{4}{35}c^{21}_{140}
-\frac{2}{7}c^{23}_{140}
$,
$\frac{1}{\sqrt{35}}c^{4}_{35}
-\frac{1}{\sqrt{35}}c^{11}_{35}
$;
$\frac{1}{\sqrt{35}\mathrm{i}}s^{3}_{140}
+\frac{1}{\sqrt{35}\mathrm{i}}s^{17}_{140}
$,
$-\frac{2}{35}c^{1}_{140}
+\frac{1}{35}c^{3}_{140}
+\frac{1}{7}c^{5}_{140}
+\frac{3}{35}c^{7}_{140}
-\frac{1}{5}c^{9}_{140}
+\frac{2}{35}c^{13}_{140}
+\frac{1}{35}c^{15}_{140}
+\frac{1}{35}c^{17}_{140}
-\frac{3}{35}c^{19}_{140}
-\frac{2}{35}c^{21}_{140}
+\frac{2}{7}c^{23}_{140}
$;
$\frac{2}{\sqrt{35}}c^{3}_{35}
+\frac{1}{\sqrt{35}}c^{4}_{35}
+\frac{1}{\sqrt{35}}c^{10}_{35}
+\frac{1}{\sqrt{35}}c^{11}_{35}
$)

Pass. 

 \ \color{black}

 \color{blue}

\noindent 689: (dims,levels) = $(6;140
)$,
irreps = $3_{7}^{3}
\hskip -1.5pt \otimes \hskip -1.5pt
2_{5}^{1}
\hskip -1.5pt \otimes \hskip -1.5pt
1_{4}^{1,0}$,
pord$(\rho_\text{isum}(\mathfrak{t})) = 35$,

\vskip 0.7ex
\hangindent=5.5em \hangafter=1
{\white .}\hskip 1em $\rho_\text{isum}(\mathfrak{t})$ =
 $( \frac{23}{140},
\frac{43}{140},
\frac{67}{140},
\frac{107}{140},
\frac{123}{140},
\frac{127}{140} )
$,

\vskip 0.7ex
\hangindent=5.5em \hangafter=1
{\white .}\hskip 1em $\rho_\text{isum}(\mathfrak{s})$ =
($\frac{1}{\sqrt{35}}c^{1}_{35}
-\frac{1}{\sqrt{35}}c^{6}_{35}
$,
$-\frac{2}{35}c^{1}_{140}
+\frac{1}{35}c^{3}_{140}
+\frac{1}{7}c^{5}_{140}
+\frac{3}{35}c^{7}_{140}
-\frac{1}{5}c^{9}_{140}
+\frac{2}{35}c^{13}_{140}
+\frac{1}{35}c^{15}_{140}
+\frac{1}{35}c^{17}_{140}
-\frac{3}{35}c^{19}_{140}
-\frac{2}{35}c^{21}_{140}
+\frac{2}{7}c^{23}_{140}
$,
$\frac{2}{\sqrt{35}}c^{3}_{35}
+\frac{1}{\sqrt{35}}c^{4}_{35}
+\frac{1}{\sqrt{35}}c^{10}_{35}
+\frac{1}{\sqrt{35}}c^{11}_{35}
$,
$-\frac{1}{\sqrt{35}\mathrm{i}}s^{3}_{140}
-\frac{1}{\sqrt{35}\mathrm{i}}s^{17}_{140}
$,
$\frac{1}{\sqrt{35}}c^{4}_{35}
-\frac{1}{\sqrt{35}}c^{11}_{35}
$,
$\frac{4}{35}c^{1}_{140}
+\frac{3}{35}c^{3}_{140}
+\frac{1}{7}c^{5}_{140}
-\frac{1}{35}c^{7}_{140}
-\frac{1}{35}c^{9}_{140}
-\frac{4}{35}c^{13}_{140}
-\frac{2}{35}c^{15}_{140}
+\frac{3}{35}c^{17}_{140}
-\frac{9}{35}c^{19}_{140}
+\frac{4}{35}c^{21}_{140}
+\frac{2}{7}c^{23}_{140}
$;
$-\frac{1}{\sqrt{35}}c^{4}_{35}
+\frac{1}{\sqrt{35}}c^{11}_{35}
$,
$\frac{1}{\sqrt{35}\mathrm{i}}s^{3}_{140}
+\frac{1}{\sqrt{35}\mathrm{i}}s^{17}_{140}
$,
$-\frac{4}{35}c^{1}_{140}
-\frac{3}{35}c^{3}_{140}
-\frac{1}{7}c^{5}_{140}
+\frac{1}{35}c^{7}_{140}
+\frac{1}{35}c^{9}_{140}
+\frac{4}{35}c^{13}_{140}
+\frac{2}{35}c^{15}_{140}
-\frac{3}{35}c^{17}_{140}
+\frac{9}{35}c^{19}_{140}
-\frac{4}{35}c^{21}_{140}
-\frac{2}{7}c^{23}_{140}
$,
$\frac{1}{\sqrt{35}}c^{1}_{35}
-\frac{1}{\sqrt{35}}c^{6}_{35}
$,
$-\frac{2}{\sqrt{35}}c^{3}_{35}
-\frac{1}{\sqrt{35}}c^{4}_{35}
-\frac{1}{\sqrt{35}}c^{10}_{35}
-\frac{1}{\sqrt{35}}c^{11}_{35}
$;
$-\frac{2}{35}c^{1}_{140}
+\frac{1}{35}c^{3}_{140}
+\frac{1}{7}c^{5}_{140}
+\frac{3}{35}c^{7}_{140}
-\frac{1}{5}c^{9}_{140}
+\frac{2}{35}c^{13}_{140}
+\frac{1}{35}c^{15}_{140}
+\frac{1}{35}c^{17}_{140}
-\frac{3}{35}c^{19}_{140}
-\frac{2}{35}c^{21}_{140}
+\frac{2}{7}c^{23}_{140}
$,
$\frac{1}{\sqrt{35}}c^{4}_{35}
-\frac{1}{\sqrt{35}}c^{11}_{35}
$,
$-\frac{4}{35}c^{1}_{140}
-\frac{3}{35}c^{3}_{140}
-\frac{1}{7}c^{5}_{140}
+\frac{1}{35}c^{7}_{140}
+\frac{1}{35}c^{9}_{140}
+\frac{4}{35}c^{13}_{140}
+\frac{2}{35}c^{15}_{140}
-\frac{3}{35}c^{17}_{140}
+\frac{9}{35}c^{19}_{140}
-\frac{4}{35}c^{21}_{140}
-\frac{2}{7}c^{23}_{140}
$,
$-\frac{1}{\sqrt{35}}c^{1}_{35}
+\frac{1}{\sqrt{35}}c^{6}_{35}
$;
$-\frac{1}{\sqrt{35}}c^{1}_{35}
+\frac{1}{\sqrt{35}}c^{6}_{35}
$,
$-\frac{2}{\sqrt{35}}c^{3}_{35}
-\frac{1}{\sqrt{35}}c^{4}_{35}
-\frac{1}{\sqrt{35}}c^{10}_{35}
-\frac{1}{\sqrt{35}}c^{11}_{35}
$,
$-\frac{2}{35}c^{1}_{140}
+\frac{1}{35}c^{3}_{140}
+\frac{1}{7}c^{5}_{140}
+\frac{3}{35}c^{7}_{140}
-\frac{1}{5}c^{9}_{140}
+\frac{2}{35}c^{13}_{140}
+\frac{1}{35}c^{15}_{140}
+\frac{1}{35}c^{17}_{140}
-\frac{3}{35}c^{19}_{140}
-\frac{2}{35}c^{21}_{140}
+\frac{2}{7}c^{23}_{140}
$;
$\frac{2}{35}c^{1}_{140}
-\frac{1}{35}c^{3}_{140}
-\frac{1}{7}c^{5}_{140}
-\frac{3}{35}c^{7}_{140}
+\frac{1}{5}c^{9}_{140}
-\frac{2}{35}c^{13}_{140}
-\frac{1}{35}c^{15}_{140}
-\frac{1}{35}c^{17}_{140}
+\frac{3}{35}c^{19}_{140}
+\frac{2}{35}c^{21}_{140}
-\frac{2}{7}c^{23}_{140}
$,
$\frac{1}{\sqrt{35}\mathrm{i}}s^{3}_{140}
+\frac{1}{\sqrt{35}\mathrm{i}}s^{17}_{140}
$;
$\frac{1}{\sqrt{35}}c^{4}_{35}
-\frac{1}{\sqrt{35}}c^{11}_{35}
$)

Pass. 

 \ \color{black}

 \color{blue}

\noindent 690: (dims,levels) = $(6;156
)$,
irreps = $6_{13}^{2}
\hskip -1.5pt \otimes \hskip -1.5pt
1_{4}^{1,0}
\hskip -1.5pt \otimes \hskip -1.5pt
1_{3}^{1,0}$,
pord$(\rho_\text{isum}(\mathfrak{t})) = 13$,

\vskip 0.7ex
\hangindent=5.5em \hangafter=1
{\white .}\hskip 1em $\rho_\text{isum}(\mathfrak{t})$ =
 $( \frac{7}{156},
\frac{19}{156},
\frac{31}{156},
\frac{67}{156},
\frac{115}{156},
\frac{151}{156} )
$,

\vskip 0.7ex
\hangindent=5.5em \hangafter=1
{\white .}\hskip 1em $\rho_\text{isum}(\mathfrak{s})$ =
($\frac{1}{\sqrt{13}}c^{9}_{52}
$,
$\frac{1}{\sqrt{13}}c^{7}_{52}
$,
$\frac{1}{\sqrt{13}}c^{11}_{52}
$,
$\frac{1}{\sqrt{13}}c^{5}_{52}
$,
$\frac{1}{\sqrt{13}}c^{1}_{52}
$,
$\frac{1}{\sqrt{13}}c^{3}_{52}
$;
$-\frac{1}{\sqrt{13}}c^{9}_{52}
$,
$\frac{1}{\sqrt{13}}c^{3}_{52}
$,
$-\frac{1}{\sqrt{13}}c^{1}_{52}
$,
$\frac{1}{\sqrt{13}}c^{5}_{52}
$,
$-\frac{1}{\sqrt{13}}c^{11}_{52}
$;
$-\frac{1}{\sqrt{13}}c^{1}_{52}
$,
$-\frac{1}{\sqrt{13}}c^{9}_{52}
$,
$\frac{1}{\sqrt{13}}c^{7}_{52}
$,
$-\frac{1}{\sqrt{13}}c^{5}_{52}
$;
$\frac{1}{\sqrt{13}}c^{3}_{52}
$,
$\frac{1}{\sqrt{13}}c^{11}_{52}
$,
$-\frac{1}{\sqrt{13}}c^{7}_{52}
$;
$-\frac{1}{\sqrt{13}}c^{3}_{52}
$,
$-\frac{1}{\sqrt{13}}c^{9}_{52}
$;
$\frac{1}{\sqrt{13}}c^{1}_{52}
$)

Pass. 

 \ \color{black}

 \color{blue}

\noindent 691: (dims,levels) = $(6;156
)$,
irreps = $6_{13}^{1}
\hskip -1.5pt \otimes \hskip -1.5pt
1_{4}^{1,0}
\hskip -1.5pt \otimes \hskip -1.5pt
1_{3}^{1,0}$,
pord$(\rho_\text{isum}(\mathfrak{t})) = 13$,

\vskip 0.7ex
\hangindent=5.5em \hangafter=1
{\white .}\hskip 1em $\rho_\text{isum}(\mathfrak{t})$ =
 $( \frac{43}{156},
\frac{55}{156},
\frac{79}{156},
\frac{103}{156},
\frac{127}{156},
\frac{139}{156} )
$,

\vskip 0.7ex
\hangindent=5.5em \hangafter=1
{\white .}\hskip 1em $\rho_\text{isum}(\mathfrak{s})$ =
($\frac{1}{\sqrt{13}}c^{7}_{52}
$,
$\frac{1}{\sqrt{13}}c^{1}_{52}
$,
$\frac{1}{\sqrt{13}}c^{3}_{52}
$,
$\frac{1}{\sqrt{13}}c^{11}_{52}
$,
$\frac{1}{\sqrt{13}}c^{5}_{52}
$,
$-\frac{1}{\sqrt{13}}c^{9}_{52}
$;
$-\frac{1}{\sqrt{13}}c^{11}_{52}
$,
$\frac{1}{\sqrt{13}}c^{7}_{52}
$,
$\frac{1}{\sqrt{13}}c^{9}_{52}
$,
$-\frac{1}{\sqrt{13}}c^{3}_{52}
$,
$\frac{1}{\sqrt{13}}c^{5}_{52}
$;
$-\frac{1}{\sqrt{13}}c^{5}_{52}
$,
$-\frac{1}{\sqrt{13}}c^{1}_{52}
$,
$-\frac{1}{\sqrt{13}}c^{9}_{52}
$,
$-\frac{1}{\sqrt{13}}c^{11}_{52}
$;
$\frac{1}{\sqrt{13}}c^{5}_{52}
$,
$\frac{1}{\sqrt{13}}c^{7}_{52}
$,
$\frac{1}{\sqrt{13}}c^{3}_{52}
$;
$\frac{1}{\sqrt{13}}c^{11}_{52}
$,
$-\frac{1}{\sqrt{13}}c^{1}_{52}
$;
$-\frac{1}{\sqrt{13}}c^{7}_{52}
$)

Pass. 

 \ \color{black}

 \color{blue}

\noindent 692: (dims,levels) = $(6;168
)$,
irreps = $3_{7}^{1}
\hskip -1.5pt \otimes \hskip -1.5pt
2_{8}^{1,0}
\hskip -1.5pt \otimes \hskip -1.5pt
1_{3}^{1,0}$,
pord$(\rho_\text{isum}(\mathfrak{t})) = 28$,

\vskip 0.7ex
\hangindent=5.5em \hangafter=1
{\white .}\hskip 1em $\rho_\text{isum}(\mathfrak{t})$ =
 $( \frac{5}{168},
\frac{47}{168},
\frac{101}{168},
\frac{125}{168},
\frac{143}{168},
\frac{167}{168} )
$,

\vskip 0.7ex
\hangindent=5.5em \hangafter=1
{\white .}\hskip 1em $\rho_\text{isum}(\mathfrak{s})$ =
($\frac{1}{\sqrt{14}}c^{3}_{28}
$,
$\frac{1}{\sqrt{14}}c^{3}_{28}
$,
$-\frac{1}{\sqrt{14}}c^{5}_{28}
$,
$\frac{1}{\sqrt{14}}c^{1}_{28}
$,
$-\frac{1}{\sqrt{14}}c^{5}_{28}
$,
$\frac{1}{\sqrt{14}}c^{1}_{28}
$;
$-\frac{1}{\sqrt{14}}c^{3}_{28}
$,
$-\frac{1}{\sqrt{14}}c^{5}_{28}
$,
$\frac{1}{\sqrt{14}}c^{1}_{28}
$,
$\frac{1}{\sqrt{14}}c^{5}_{28}
$,
$-\frac{1}{\sqrt{14}}c^{1}_{28}
$;
$\frac{1}{\sqrt{14}}c^{1}_{28}
$,
$\frac{1}{\sqrt{14}}c^{3}_{28}
$,
$\frac{1}{\sqrt{14}}c^{1}_{28}
$,
$\frac{1}{\sqrt{14}}c^{3}_{28}
$;
$-\frac{1}{\sqrt{14}}c^{5}_{28}
$,
$\frac{1}{\sqrt{14}}c^{3}_{28}
$,
$-\frac{1}{\sqrt{14}}c^{5}_{28}
$;
$-\frac{1}{\sqrt{14}}c^{1}_{28}
$,
$-\frac{1}{\sqrt{14}}c^{3}_{28}
$;
$\frac{1}{\sqrt{14}}c^{5}_{28}
$)

Pass. 

 \ \color{black}

 \color{blue}

\noindent 693: (dims,levels) = $(6;168
)$,
irreps = $3_{7}^{3}
\hskip -1.5pt \otimes \hskip -1.5pt
2_{8}^{1,0}
\hskip -1.5pt \otimes \hskip -1.5pt
1_{3}^{1,0}$,
pord$(\rho_\text{isum}(\mathfrak{t})) = 28$,

\vskip 0.7ex
\hangindent=5.5em \hangafter=1
{\white .}\hskip 1em $\rho_\text{isum}(\mathfrak{t})$ =
 $( \frac{23}{168},
\frac{29}{168},
\frac{53}{168},
\frac{71}{168},
\frac{95}{168},
\frac{149}{168} )
$,

\vskip 0.7ex
\hangindent=5.5em \hangafter=1
{\white .}\hskip 1em $\rho_\text{isum}(\mathfrak{s})$ =
($-\frac{1}{\sqrt{14}}c^{3}_{28}
$,
$\frac{1}{\sqrt{14}}c^{1}_{28}
$,
$-\frac{1}{\sqrt{14}}c^{5}_{28}
$,
$\frac{1}{\sqrt{14}}c^{1}_{28}
$,
$-\frac{1}{\sqrt{14}}c^{5}_{28}
$,
$\frac{1}{\sqrt{14}}c^{3}_{28}
$;
$-\frac{1}{\sqrt{14}}c^{5}_{28}
$,
$\frac{1}{\sqrt{14}}c^{3}_{28}
$,
$\frac{1}{\sqrt{14}}c^{5}_{28}
$,
$-\frac{1}{\sqrt{14}}c^{3}_{28}
$,
$\frac{1}{\sqrt{14}}c^{1}_{28}
$;
$\frac{1}{\sqrt{14}}c^{1}_{28}
$,
$-\frac{1}{\sqrt{14}}c^{3}_{28}
$,
$-\frac{1}{\sqrt{14}}c^{1}_{28}
$,
$-\frac{1}{\sqrt{14}}c^{5}_{28}
$;
$\frac{1}{\sqrt{14}}c^{5}_{28}
$,
$-\frac{1}{\sqrt{14}}c^{3}_{28}
$,
$-\frac{1}{\sqrt{14}}c^{1}_{28}
$;
$-\frac{1}{\sqrt{14}}c^{1}_{28}
$,
$\frac{1}{\sqrt{14}}c^{5}_{28}
$;
$\frac{1}{\sqrt{14}}c^{3}_{28}
$)

Pass. 

 \ \color{black}

 \color{blue}

\noindent 694: (dims,levels) = $(6;210
)$,
irreps = $3_{7}^{1}
\hskip -1.5pt \otimes \hskip -1.5pt
2_{5}^{2}
\hskip -1.5pt \otimes \hskip -1.5pt
1_{3}^{1,0}
\hskip -1.5pt \otimes \hskip -1.5pt
1_{2}^{1,0}$,
pord$(\rho_\text{isum}(\mathfrak{t})) = 35$,

\vskip 0.7ex
\hangindent=5.5em \hangafter=1
{\white .}\hskip 1em $\rho_\text{isum}(\mathfrak{t})$ =
 $( \frac{1}{210},
\frac{79}{210},
\frac{109}{210},
\frac{121}{210},
\frac{151}{210},
\frac{169}{210} )
$,

\vskip 0.7ex
\hangindent=5.5em \hangafter=1
{\white .}\hskip 1em $\rho_\text{isum}(\mathfrak{s})$ =
$\mathrm{i}$($\frac{4}{35}c^{1}_{140}
+\frac{3}{35}c^{3}_{140}
+\frac{1}{7}c^{5}_{140}
-\frac{1}{35}c^{7}_{140}
-\frac{1}{35}c^{9}_{140}
-\frac{4}{35}c^{13}_{140}
-\frac{2}{35}c^{15}_{140}
+\frac{3}{35}c^{17}_{140}
-\frac{9}{35}c^{19}_{140}
+\frac{4}{35}c^{21}_{140}
+\frac{2}{7}c^{23}_{140}
$,
$\frac{1}{\sqrt{35}}c^{1}_{35}
-\frac{1}{\sqrt{35}}c^{6}_{35}
$,
$-\frac{1}{\sqrt{35}}c^{4}_{35}
+\frac{1}{\sqrt{35}}c^{11}_{35}
$,
$-\frac{1}{\sqrt{35}\mathrm{i}}s^{3}_{140}
-\frac{1}{\sqrt{35}\mathrm{i}}s^{17}_{140}
$,
$\frac{2}{\sqrt{35}}c^{3}_{35}
+\frac{1}{\sqrt{35}}c^{4}_{35}
+\frac{1}{\sqrt{35}}c^{10}_{35}
+\frac{1}{\sqrt{35}}c^{11}_{35}
$,
$\frac{2}{35}c^{1}_{140}
-\frac{1}{35}c^{3}_{140}
-\frac{1}{7}c^{5}_{140}
-\frac{3}{35}c^{7}_{140}
+\frac{1}{5}c^{9}_{140}
-\frac{2}{35}c^{13}_{140}
-\frac{1}{35}c^{15}_{140}
-\frac{1}{35}c^{17}_{140}
+\frac{3}{35}c^{19}_{140}
+\frac{2}{35}c^{21}_{140}
-\frac{2}{7}c^{23}_{140}
$;\ \ 
$-\frac{2}{\sqrt{35}}c^{3}_{35}
-\frac{1}{\sqrt{35}}c^{4}_{35}
-\frac{1}{\sqrt{35}}c^{10}_{35}
-\frac{1}{\sqrt{35}}c^{11}_{35}
$,
$-\frac{4}{35}c^{1}_{140}
-\frac{3}{35}c^{3}_{140}
-\frac{1}{7}c^{5}_{140}
+\frac{1}{35}c^{7}_{140}
+\frac{1}{35}c^{9}_{140}
+\frac{4}{35}c^{13}_{140}
+\frac{2}{35}c^{15}_{140}
-\frac{3}{35}c^{17}_{140}
+\frac{9}{35}c^{19}_{140}
-\frac{4}{35}c^{21}_{140}
-\frac{2}{7}c^{23}_{140}
$,
$-\frac{1}{\sqrt{35}}c^{4}_{35}
+\frac{1}{\sqrt{35}}c^{11}_{35}
$,
$\frac{2}{35}c^{1}_{140}
-\frac{1}{35}c^{3}_{140}
-\frac{1}{7}c^{5}_{140}
-\frac{3}{35}c^{7}_{140}
+\frac{1}{5}c^{9}_{140}
-\frac{2}{35}c^{13}_{140}
-\frac{1}{35}c^{15}_{140}
-\frac{1}{35}c^{17}_{140}
+\frac{3}{35}c^{19}_{140}
+\frac{2}{35}c^{21}_{140}
-\frac{2}{7}c^{23}_{140}
$,
$\frac{1}{\sqrt{35}\mathrm{i}}s^{3}_{140}
+\frac{1}{\sqrt{35}\mathrm{i}}s^{17}_{140}
$;\ \ 
$\frac{1}{\sqrt{35}\mathrm{i}}s^{3}_{140}
+\frac{1}{\sqrt{35}\mathrm{i}}s^{17}_{140}
$,
$\frac{2}{35}c^{1}_{140}
-\frac{1}{35}c^{3}_{140}
-\frac{1}{7}c^{5}_{140}
-\frac{3}{35}c^{7}_{140}
+\frac{1}{5}c^{9}_{140}
-\frac{2}{35}c^{13}_{140}
-\frac{1}{35}c^{15}_{140}
-\frac{1}{35}c^{17}_{140}
+\frac{3}{35}c^{19}_{140}
+\frac{2}{35}c^{21}_{140}
-\frac{2}{7}c^{23}_{140}
$,
$\frac{1}{\sqrt{35}}c^{1}_{35}
-\frac{1}{\sqrt{35}}c^{6}_{35}
$,
$-\frac{2}{\sqrt{35}}c^{3}_{35}
-\frac{1}{\sqrt{35}}c^{4}_{35}
-\frac{1}{\sqrt{35}}c^{10}_{35}
-\frac{1}{\sqrt{35}}c^{11}_{35}
$;\ \ 
$\frac{2}{\sqrt{35}}c^{3}_{35}
+\frac{1}{\sqrt{35}}c^{4}_{35}
+\frac{1}{\sqrt{35}}c^{10}_{35}
+\frac{1}{\sqrt{35}}c^{11}_{35}
$,
$\frac{4}{35}c^{1}_{140}
+\frac{3}{35}c^{3}_{140}
+\frac{1}{7}c^{5}_{140}
-\frac{1}{35}c^{7}_{140}
-\frac{1}{35}c^{9}_{140}
-\frac{4}{35}c^{13}_{140}
-\frac{2}{35}c^{15}_{140}
+\frac{3}{35}c^{17}_{140}
-\frac{9}{35}c^{19}_{140}
+\frac{4}{35}c^{21}_{140}
+\frac{2}{7}c^{23}_{140}
$,
$\frac{1}{\sqrt{35}}c^{1}_{35}
-\frac{1}{\sqrt{35}}c^{6}_{35}
$;\ \ 
$-\frac{1}{\sqrt{35}\mathrm{i}}s^{3}_{140}
-\frac{1}{\sqrt{35}\mathrm{i}}s^{17}_{140}
$,
$-\frac{1}{\sqrt{35}}c^{4}_{35}
+\frac{1}{\sqrt{35}}c^{11}_{35}
$;\ \ 
$-\frac{4}{35}c^{1}_{140}
-\frac{3}{35}c^{3}_{140}
-\frac{1}{7}c^{5}_{140}
+\frac{1}{35}c^{7}_{140}
+\frac{1}{35}c^{9}_{140}
+\frac{4}{35}c^{13}_{140}
+\frac{2}{35}c^{15}_{140}
-\frac{3}{35}c^{17}_{140}
+\frac{9}{35}c^{19}_{140}
-\frac{4}{35}c^{21}_{140}
-\frac{2}{7}c^{23}_{140}
$)

Pass. 

 \ \color{black}

 \color{blue}

\noindent 695: (dims,levels) = $(6;210
)$,
irreps = $3_{7}^{3}
\hskip -1.5pt \otimes \hskip -1.5pt
2_{5}^{1}
\hskip -1.5pt \otimes \hskip -1.5pt
1_{3}^{1,0}
\hskip -1.5pt \otimes \hskip -1.5pt
1_{2}^{1,0}$,
pord$(\rho_\text{isum}(\mathfrak{t})) = 35$,

\vskip 0.7ex
\hangindent=5.5em \hangafter=1
{\white .}\hskip 1em $\rho_\text{isum}(\mathfrak{t})$ =
 $( \frac{13}{210},
\frac{73}{210},
\frac{97}{210},
\frac{103}{210},
\frac{157}{210},
\frac{187}{210} )
$,

\vskip 0.7ex
\hangindent=5.5em \hangafter=1
{\white .}\hskip 1em $\rho_\text{isum}(\mathfrak{s})$ =
$\mathrm{i}$($-\frac{2}{35}c^{1}_{140}
+\frac{1}{35}c^{3}_{140}
+\frac{1}{7}c^{5}_{140}
+\frac{3}{35}c^{7}_{140}
-\frac{1}{5}c^{9}_{140}
+\frac{2}{35}c^{13}_{140}
+\frac{1}{35}c^{15}_{140}
+\frac{1}{35}c^{17}_{140}
-\frac{3}{35}c^{19}_{140}
-\frac{2}{35}c^{21}_{140}
+\frac{2}{7}c^{23}_{140}
$,
$-\frac{1}{\sqrt{35}}c^{4}_{35}
+\frac{1}{\sqrt{35}}c^{11}_{35}
$,
$\frac{4}{35}c^{1}_{140}
+\frac{3}{35}c^{3}_{140}
+\frac{1}{7}c^{5}_{140}
-\frac{1}{35}c^{7}_{140}
-\frac{1}{35}c^{9}_{140}
-\frac{4}{35}c^{13}_{140}
-\frac{2}{35}c^{15}_{140}
+\frac{3}{35}c^{17}_{140}
-\frac{9}{35}c^{19}_{140}
+\frac{4}{35}c^{21}_{140}
+\frac{2}{7}c^{23}_{140}
$,
$\frac{1}{\sqrt{35}}c^{1}_{35}
-\frac{1}{\sqrt{35}}c^{6}_{35}
$,
$\frac{2}{\sqrt{35}}c^{3}_{35}
+\frac{1}{\sqrt{35}}c^{4}_{35}
+\frac{1}{\sqrt{35}}c^{10}_{35}
+\frac{1}{\sqrt{35}}c^{11}_{35}
$,
$-\frac{1}{\sqrt{35}\mathrm{i}}s^{3}_{140}
-\frac{1}{\sqrt{35}\mathrm{i}}s^{17}_{140}
$;\ \ 
$-\frac{1}{\sqrt{35}}c^{1}_{35}
+\frac{1}{\sqrt{35}}c^{6}_{35}
$,
$-\frac{2}{\sqrt{35}}c^{3}_{35}
-\frac{1}{\sqrt{35}}c^{4}_{35}
-\frac{1}{\sqrt{35}}c^{10}_{35}
-\frac{1}{\sqrt{35}}c^{11}_{35}
$,
$-\frac{2}{35}c^{1}_{140}
+\frac{1}{35}c^{3}_{140}
+\frac{1}{7}c^{5}_{140}
+\frac{3}{35}c^{7}_{140}
-\frac{1}{5}c^{9}_{140}
+\frac{2}{35}c^{13}_{140}
+\frac{1}{35}c^{15}_{140}
+\frac{1}{35}c^{17}_{140}
-\frac{3}{35}c^{19}_{140}
-\frac{2}{35}c^{21}_{140}
+\frac{2}{7}c^{23}_{140}
$,
$\frac{1}{\sqrt{35}\mathrm{i}}s^{3}_{140}
+\frac{1}{\sqrt{35}\mathrm{i}}s^{17}_{140}
$,
$-\frac{4}{35}c^{1}_{140}
-\frac{3}{35}c^{3}_{140}
-\frac{1}{7}c^{5}_{140}
+\frac{1}{35}c^{7}_{140}
+\frac{1}{35}c^{9}_{140}
+\frac{4}{35}c^{13}_{140}
+\frac{2}{35}c^{15}_{140}
-\frac{3}{35}c^{17}_{140}
+\frac{9}{35}c^{19}_{140}
-\frac{4}{35}c^{21}_{140}
-\frac{2}{7}c^{23}_{140}
$;\ \ 
$\frac{2}{35}c^{1}_{140}
-\frac{1}{35}c^{3}_{140}
-\frac{1}{7}c^{5}_{140}
-\frac{3}{35}c^{7}_{140}
+\frac{1}{5}c^{9}_{140}
-\frac{2}{35}c^{13}_{140}
-\frac{1}{35}c^{15}_{140}
-\frac{1}{35}c^{17}_{140}
+\frac{3}{35}c^{19}_{140}
+\frac{2}{35}c^{21}_{140}
-\frac{2}{7}c^{23}_{140}
$,
$\frac{1}{\sqrt{35}\mathrm{i}}s^{3}_{140}
+\frac{1}{\sqrt{35}\mathrm{i}}s^{17}_{140}
$,
$-\frac{1}{\sqrt{35}}c^{4}_{35}
+\frac{1}{\sqrt{35}}c^{11}_{35}
$,
$\frac{1}{\sqrt{35}}c^{1}_{35}
-\frac{1}{\sqrt{35}}c^{6}_{35}
$;\ \ 
$\frac{1}{\sqrt{35}}c^{4}_{35}
-\frac{1}{\sqrt{35}}c^{11}_{35}
$,
$-\frac{4}{35}c^{1}_{140}
-\frac{3}{35}c^{3}_{140}
-\frac{1}{7}c^{5}_{140}
+\frac{1}{35}c^{7}_{140}
+\frac{1}{35}c^{9}_{140}
+\frac{4}{35}c^{13}_{140}
+\frac{2}{35}c^{15}_{140}
-\frac{3}{35}c^{17}_{140}
+\frac{9}{35}c^{19}_{140}
-\frac{4}{35}c^{21}_{140}
-\frac{2}{7}c^{23}_{140}
$,
$-\frac{2}{\sqrt{35}}c^{3}_{35}
-\frac{1}{\sqrt{35}}c^{4}_{35}
-\frac{1}{\sqrt{35}}c^{10}_{35}
-\frac{1}{\sqrt{35}}c^{11}_{35}
$;\ \ 
$\frac{1}{\sqrt{35}}c^{1}_{35}
-\frac{1}{\sqrt{35}}c^{6}_{35}
$,
$\frac{2}{35}c^{1}_{140}
-\frac{1}{35}c^{3}_{140}
-\frac{1}{7}c^{5}_{140}
-\frac{3}{35}c^{7}_{140}
+\frac{1}{5}c^{9}_{140}
-\frac{2}{35}c^{13}_{140}
-\frac{1}{35}c^{15}_{140}
-\frac{1}{35}c^{17}_{140}
+\frac{3}{35}c^{19}_{140}
+\frac{2}{35}c^{21}_{140}
-\frac{2}{7}c^{23}_{140}
$;\ \ 
$-\frac{1}{\sqrt{35}}c^{4}_{35}
+\frac{1}{\sqrt{35}}c^{11}_{35}
$)

Pass. 

 \ \color{black}

 \color{blue}

\noindent 696: (dims,levels) = $(6;210
)$,
irreps = $3_{7}^{3}
\hskip -1.5pt \otimes \hskip -1.5pt
2_{5}^{2}
\hskip -1.5pt \otimes \hskip -1.5pt
1_{3}^{1,0}
\hskip -1.5pt \otimes \hskip -1.5pt
1_{2}^{1,0}$,
pord$(\rho_\text{isum}(\mathfrak{t})) = 35$,

\vskip 0.7ex
\hangindent=5.5em \hangafter=1
{\white .}\hskip 1em $\rho_\text{isum}(\mathfrak{t})$ =
 $( \frac{19}{210},
\frac{31}{210},
\frac{61}{210},
\frac{139}{210},
\frac{181}{210},
\frac{199}{210} )
$,

\vskip 0.7ex
\hangindent=5.5em \hangafter=1
{\white .}\hskip 1em $\rho_\text{isum}(\mathfrak{s})$ =
$\mathrm{i}$($-\frac{2}{\sqrt{35}}c^{3}_{35}
-\frac{1}{\sqrt{35}}c^{4}_{35}
-\frac{1}{\sqrt{35}}c^{10}_{35}
-\frac{1}{\sqrt{35}}c^{11}_{35}
$,
$\frac{2}{35}c^{1}_{140}
-\frac{1}{35}c^{3}_{140}
-\frac{1}{7}c^{5}_{140}
-\frac{3}{35}c^{7}_{140}
+\frac{1}{5}c^{9}_{140}
-\frac{2}{35}c^{13}_{140}
-\frac{1}{35}c^{15}_{140}
-\frac{1}{35}c^{17}_{140}
+\frac{3}{35}c^{19}_{140}
+\frac{2}{35}c^{21}_{140}
-\frac{2}{7}c^{23}_{140}
$,
$-\frac{1}{\sqrt{35}}c^{4}_{35}
+\frac{1}{\sqrt{35}}c^{11}_{35}
$,
$-\frac{1}{\sqrt{35}\mathrm{i}}s^{3}_{140}
-\frac{1}{\sqrt{35}\mathrm{i}}s^{17}_{140}
$,
$\frac{1}{\sqrt{35}}c^{1}_{35}
-\frac{1}{\sqrt{35}}c^{6}_{35}
$,
$\frac{4}{35}c^{1}_{140}
+\frac{3}{35}c^{3}_{140}
+\frac{1}{7}c^{5}_{140}
-\frac{1}{35}c^{7}_{140}
-\frac{1}{35}c^{9}_{140}
-\frac{4}{35}c^{13}_{140}
-\frac{2}{35}c^{15}_{140}
+\frac{3}{35}c^{17}_{140}
-\frac{9}{35}c^{19}_{140}
+\frac{4}{35}c^{21}_{140}
+\frac{2}{7}c^{23}_{140}
$;\ \ 
$-\frac{1}{\sqrt{35}\mathrm{i}}s^{3}_{140}
-\frac{1}{\sqrt{35}\mathrm{i}}s^{17}_{140}
$,
$\frac{4}{35}c^{1}_{140}
+\frac{3}{35}c^{3}_{140}
+\frac{1}{7}c^{5}_{140}
-\frac{1}{35}c^{7}_{140}
-\frac{1}{35}c^{9}_{140}
-\frac{4}{35}c^{13}_{140}
-\frac{2}{35}c^{15}_{140}
+\frac{3}{35}c^{17}_{140}
-\frac{9}{35}c^{19}_{140}
+\frac{4}{35}c^{21}_{140}
+\frac{2}{7}c^{23}_{140}
$,
$\frac{1}{\sqrt{35}}c^{4}_{35}
-\frac{1}{\sqrt{35}}c^{11}_{35}
$,
$\frac{2}{\sqrt{35}}c^{3}_{35}
+\frac{1}{\sqrt{35}}c^{4}_{35}
+\frac{1}{\sqrt{35}}c^{10}_{35}
+\frac{1}{\sqrt{35}}c^{11}_{35}
$,
$-\frac{1}{\sqrt{35}}c^{1}_{35}
+\frac{1}{\sqrt{35}}c^{6}_{35}
$;\ \ 
$\frac{2}{\sqrt{35}}c^{3}_{35}
+\frac{1}{\sqrt{35}}c^{4}_{35}
+\frac{1}{\sqrt{35}}c^{10}_{35}
+\frac{1}{\sqrt{35}}c^{11}_{35}
$,
$-\frac{1}{\sqrt{35}}c^{1}_{35}
+\frac{1}{\sqrt{35}}c^{6}_{35}
$,
$-\frac{1}{\sqrt{35}\mathrm{i}}s^{3}_{140}
-\frac{1}{\sqrt{35}\mathrm{i}}s^{17}_{140}
$,
$-\frac{2}{35}c^{1}_{140}
+\frac{1}{35}c^{3}_{140}
+\frac{1}{7}c^{5}_{140}
+\frac{3}{35}c^{7}_{140}
-\frac{1}{5}c^{9}_{140}
+\frac{2}{35}c^{13}_{140}
+\frac{1}{35}c^{15}_{140}
+\frac{1}{35}c^{17}_{140}
-\frac{3}{35}c^{19}_{140}
-\frac{2}{35}c^{21}_{140}
+\frac{2}{7}c^{23}_{140}
$;\ \ 
$-\frac{4}{35}c^{1}_{140}
-\frac{3}{35}c^{3}_{140}
-\frac{1}{7}c^{5}_{140}
+\frac{1}{35}c^{7}_{140}
+\frac{1}{35}c^{9}_{140}
+\frac{4}{35}c^{13}_{140}
+\frac{2}{35}c^{15}_{140}
-\frac{3}{35}c^{17}_{140}
+\frac{9}{35}c^{19}_{140}
-\frac{4}{35}c^{21}_{140}
-\frac{2}{7}c^{23}_{140}
$,
$-\frac{2}{35}c^{1}_{140}
+\frac{1}{35}c^{3}_{140}
+\frac{1}{7}c^{5}_{140}
+\frac{3}{35}c^{7}_{140}
-\frac{1}{5}c^{9}_{140}
+\frac{2}{35}c^{13}_{140}
+\frac{1}{35}c^{15}_{140}
+\frac{1}{35}c^{17}_{140}
-\frac{3}{35}c^{19}_{140}
-\frac{2}{35}c^{21}_{140}
+\frac{2}{7}c^{23}_{140}
$,
$-\frac{2}{\sqrt{35}}c^{3}_{35}
-\frac{1}{\sqrt{35}}c^{4}_{35}
-\frac{1}{\sqrt{35}}c^{10}_{35}
-\frac{1}{\sqrt{35}}c^{11}_{35}
$;\ \ 
$\frac{4}{35}c^{1}_{140}
+\frac{3}{35}c^{3}_{140}
+\frac{1}{7}c^{5}_{140}
-\frac{1}{35}c^{7}_{140}
-\frac{1}{35}c^{9}_{140}
-\frac{4}{35}c^{13}_{140}
-\frac{2}{35}c^{15}_{140}
+\frac{3}{35}c^{17}_{140}
-\frac{9}{35}c^{19}_{140}
+\frac{4}{35}c^{21}_{140}
+\frac{2}{7}c^{23}_{140}
$,
$\frac{1}{\sqrt{35}}c^{4}_{35}
-\frac{1}{\sqrt{35}}c^{11}_{35}
$;\ \ 
$\frac{1}{\sqrt{35}\mathrm{i}}s^{3}_{140}
+\frac{1}{\sqrt{35}\mathrm{i}}s^{17}_{140}
$)

Pass. 

 \ \color{black}

 \color{blue}

\noindent 697: (dims,levels) = $(6;210
)$,
irreps = $3_{7}^{1}
\hskip -1.5pt \otimes \hskip -1.5pt
2_{5}^{1}
\hskip -1.5pt \otimes \hskip -1.5pt
1_{3}^{1,0}
\hskip -1.5pt \otimes \hskip -1.5pt
1_{2}^{1,0}$,
pord$(\rho_\text{isum}(\mathfrak{t})) = 35$,

\vskip 0.7ex
\hangindent=5.5em \hangafter=1
{\white .}\hskip 1em $\rho_\text{isum}(\mathfrak{t})$ =
 $( \frac{37}{210},
\frac{43}{210},
\frac{67}{210},
\frac{127}{210},
\frac{163}{210},
\frac{193}{210} )
$,

\vskip 0.7ex
\hangindent=5.5em \hangafter=1
{\white .}\hskip 1em $\rho_\text{isum}(\mathfrak{s})$ =
$\mathrm{i}$($-\frac{1}{\sqrt{35}}c^{4}_{35}
+\frac{1}{\sqrt{35}}c^{11}_{35}
$,
$-\frac{1}{\sqrt{35}\mathrm{i}}s^{3}_{140}
-\frac{1}{\sqrt{35}\mathrm{i}}s^{17}_{140}
$,
$\frac{2}{35}c^{1}_{140}
-\frac{1}{35}c^{3}_{140}
-\frac{1}{7}c^{5}_{140}
-\frac{3}{35}c^{7}_{140}
+\frac{1}{5}c^{9}_{140}
-\frac{2}{35}c^{13}_{140}
-\frac{1}{35}c^{15}_{140}
-\frac{1}{35}c^{17}_{140}
+\frac{3}{35}c^{19}_{140}
+\frac{2}{35}c^{21}_{140}
-\frac{2}{7}c^{23}_{140}
$,
$\frac{1}{\sqrt{35}}c^{1}_{35}
-\frac{1}{\sqrt{35}}c^{6}_{35}
$,
$\frac{2}{\sqrt{35}}c^{3}_{35}
+\frac{1}{\sqrt{35}}c^{4}_{35}
+\frac{1}{\sqrt{35}}c^{10}_{35}
+\frac{1}{\sqrt{35}}c^{11}_{35}
$,
$\frac{4}{35}c^{1}_{140}
+\frac{3}{35}c^{3}_{140}
+\frac{1}{7}c^{5}_{140}
-\frac{1}{35}c^{7}_{140}
-\frac{1}{35}c^{9}_{140}
-\frac{4}{35}c^{13}_{140}
-\frac{2}{35}c^{15}_{140}
+\frac{3}{35}c^{17}_{140}
-\frac{9}{35}c^{19}_{140}
+\frac{4}{35}c^{21}_{140}
+\frac{2}{7}c^{23}_{140}
$;\ \ 
$-\frac{2}{35}c^{1}_{140}
+\frac{1}{35}c^{3}_{140}
+\frac{1}{7}c^{5}_{140}
+\frac{3}{35}c^{7}_{140}
-\frac{1}{5}c^{9}_{140}
+\frac{2}{35}c^{13}_{140}
+\frac{1}{35}c^{15}_{140}
+\frac{1}{35}c^{17}_{140}
-\frac{3}{35}c^{19}_{140}
-\frac{2}{35}c^{21}_{140}
+\frac{2}{7}c^{23}_{140}
$,
$\frac{2}{\sqrt{35}}c^{3}_{35}
+\frac{1}{\sqrt{35}}c^{4}_{35}
+\frac{1}{\sqrt{35}}c^{10}_{35}
+\frac{1}{\sqrt{35}}c^{11}_{35}
$,
$\frac{4}{35}c^{1}_{140}
+\frac{3}{35}c^{3}_{140}
+\frac{1}{7}c^{5}_{140}
-\frac{1}{35}c^{7}_{140}
-\frac{1}{35}c^{9}_{140}
-\frac{4}{35}c^{13}_{140}
-\frac{2}{35}c^{15}_{140}
+\frac{3}{35}c^{17}_{140}
-\frac{9}{35}c^{19}_{140}
+\frac{4}{35}c^{21}_{140}
+\frac{2}{7}c^{23}_{140}
$,
$-\frac{1}{\sqrt{35}}c^{1}_{35}
+\frac{1}{\sqrt{35}}c^{6}_{35}
$,
$\frac{1}{\sqrt{35}}c^{4}_{35}
-\frac{1}{\sqrt{35}}c^{11}_{35}
$;\ \ 
$\frac{1}{\sqrt{35}}c^{1}_{35}
-\frac{1}{\sqrt{35}}c^{6}_{35}
$,
$-\frac{1}{\sqrt{35}}c^{4}_{35}
+\frac{1}{\sqrt{35}}c^{11}_{35}
$,
$\frac{4}{35}c^{1}_{140}
+\frac{3}{35}c^{3}_{140}
+\frac{1}{7}c^{5}_{140}
-\frac{1}{35}c^{7}_{140}
-\frac{1}{35}c^{9}_{140}
-\frac{4}{35}c^{13}_{140}
-\frac{2}{35}c^{15}_{140}
+\frac{3}{35}c^{17}_{140}
-\frac{9}{35}c^{19}_{140}
+\frac{4}{35}c^{21}_{140}
+\frac{2}{7}c^{23}_{140}
$,
$-\frac{1}{\sqrt{35}\mathrm{i}}s^{3}_{140}
-\frac{1}{\sqrt{35}\mathrm{i}}s^{17}_{140}
$;\ \ 
$\frac{2}{35}c^{1}_{140}
-\frac{1}{35}c^{3}_{140}
-\frac{1}{7}c^{5}_{140}
-\frac{3}{35}c^{7}_{140}
+\frac{1}{5}c^{9}_{140}
-\frac{2}{35}c^{13}_{140}
-\frac{1}{35}c^{15}_{140}
-\frac{1}{35}c^{17}_{140}
+\frac{3}{35}c^{19}_{140}
+\frac{2}{35}c^{21}_{140}
-\frac{2}{7}c^{23}_{140}
$,
$-\frac{1}{\sqrt{35}\mathrm{i}}s^{3}_{140}
-\frac{1}{\sqrt{35}\mathrm{i}}s^{17}_{140}
$,
$\frac{2}{\sqrt{35}}c^{3}_{35}
+\frac{1}{\sqrt{35}}c^{4}_{35}
+\frac{1}{\sqrt{35}}c^{10}_{35}
+\frac{1}{\sqrt{35}}c^{11}_{35}
$;\ \ 
$\frac{1}{\sqrt{35}}c^{4}_{35}
-\frac{1}{\sqrt{35}}c^{11}_{35}
$,
$-\frac{2}{35}c^{1}_{140}
+\frac{1}{35}c^{3}_{140}
+\frac{1}{7}c^{5}_{140}
+\frac{3}{35}c^{7}_{140}
-\frac{1}{5}c^{9}_{140}
+\frac{2}{35}c^{13}_{140}
+\frac{1}{35}c^{15}_{140}
+\frac{1}{35}c^{17}_{140}
-\frac{3}{35}c^{19}_{140}
-\frac{2}{35}c^{21}_{140}
+\frac{2}{7}c^{23}_{140}
$;\ \ 
$-\frac{1}{\sqrt{35}}c^{1}_{35}
+\frac{1}{\sqrt{35}}c^{6}_{35}
$)

Pass. 

 \ \color{black}

 \color{blue}

\noindent 698: (dims,levels) = $(6;240
)$,
irreps = $3_{16}^{7,0}
\hskip -1.5pt \otimes \hskip -1.5pt
2_{5}^{1}
\hskip -1.5pt \otimes \hskip -1.5pt
1_{3}^{1,0}$,
pord$(\rho_\text{isum}(\mathfrak{t})) = 80$,

\vskip 0.7ex
\hangindent=5.5em \hangafter=1
{\white .}\hskip 1em $\rho_\text{isum}(\mathfrak{t})$ =
 $( \frac{1}{120},
\frac{49}{120},
\frac{17}{240},
\frac{113}{240},
\frac{137}{240},
\frac{233}{240} )
$,

\vskip 0.7ex
\hangindent=5.5em \hangafter=1
{\white .}\hskip 1em $\rho_\text{isum}(\mathfrak{s})$ =
($0$,
$0$,
$\frac{1}{\sqrt{10}}c^{3}_{20}
$,
$\frac{1}{\sqrt{10}}c^{1}_{20}
$,
$\frac{1}{\sqrt{10}}c^{3}_{20}
$,
$\frac{1}{\sqrt{10}}c^{1}_{20}
$;
$0$,
$\frac{1}{\sqrt{10}}c^{1}_{20}
$,
$-\frac{1}{\sqrt{10}}c^{3}_{20}
$,
$\frac{1}{\sqrt{10}}c^{1}_{20}
$,
$-\frac{1}{\sqrt{10}}c^{3}_{20}
$;
$-\frac{1}{2\sqrt{5}}c^{3}_{20}
$,
$-\frac{1}{2\sqrt{5}}c^{1}_{20}
$,
$\frac{1}{2\sqrt{5}}c^{3}_{20}
$,
$\frac{1}{2\sqrt{5}}c^{1}_{20}
$;
$\frac{1}{2\sqrt{5}}c^{3}_{20}
$,
$\frac{1}{2\sqrt{5}}c^{1}_{20}
$,
$-\frac{1}{2\sqrt{5}}c^{3}_{20}
$;
$-\frac{1}{2\sqrt{5}}c^{3}_{20}
$,
$-\frac{1}{2\sqrt{5}}c^{1}_{20}
$;
$\frac{1}{2\sqrt{5}}c^{3}_{20}
$)

Pass. 

 \ \color{black}

 \color{blue}

\noindent 699: (dims,levels) = $(6;240
)$,
irreps = $3_{16}^{1,0}
\hskip -1.5pt \otimes \hskip -1.5pt
2_{5}^{2}
\hskip -1.5pt \otimes \hskip -1.5pt
1_{3}^{1,0}$,
pord$(\rho_\text{isum}(\mathfrak{t})) = 80$,

\vskip 0.7ex
\hangindent=5.5em \hangafter=1
{\white .}\hskip 1em $\rho_\text{isum}(\mathfrak{t})$ =
 $( \frac{7}{120},
\frac{103}{120},
\frac{71}{240},
\frac{119}{240},
\frac{191}{240},
\frac{239}{240} )
$,

\vskip 0.7ex
\hangindent=5.5em \hangafter=1
{\white .}\hskip 1em $\rho_\text{isum}(\mathfrak{s})$ =
($0$,
$0$,
$\frac{1}{\sqrt{10}}c^{3}_{20}
$,
$\frac{1}{\sqrt{10}}c^{1}_{20}
$,
$\frac{1}{\sqrt{10}}c^{3}_{20}
$,
$\frac{1}{\sqrt{10}}c^{1}_{20}
$;
$0$,
$\frac{1}{\sqrt{10}}c^{1}_{20}
$,
$-\frac{1}{\sqrt{10}}c^{3}_{20}
$,
$\frac{1}{\sqrt{10}}c^{1}_{20}
$,
$-\frac{1}{\sqrt{10}}c^{3}_{20}
$;
$-\frac{1}{2\sqrt{5}}c^{1}_{20}
$,
$\frac{1}{2\sqrt{5}}c^{3}_{20}
$,
$\frac{1}{2\sqrt{5}}c^{1}_{20}
$,
$-\frac{1}{2\sqrt{5}}c^{3}_{20}
$;
$\frac{1}{2\sqrt{5}}c^{1}_{20}
$,
$-\frac{1}{2\sqrt{5}}c^{3}_{20}
$,
$-\frac{1}{2\sqrt{5}}c^{1}_{20}
$;
$-\frac{1}{2\sqrt{5}}c^{1}_{20}
$,
$\frac{1}{2\sqrt{5}}c^{3}_{20}
$;
$\frac{1}{2\sqrt{5}}c^{1}_{20}
$)

Pass. 

 \ \color{black}

 \color{blue}

\noindent 700: (dims,levels) = $(6;240
)$,
irreps = $3_{16}^{3,0}
\hskip -1.5pt \otimes \hskip -1.5pt
2_{5}^{2}
\hskip -1.5pt \otimes \hskip -1.5pt
1_{3}^{1,0}$,
pord$(\rho_\text{isum}(\mathfrak{t})) = 80$,

\vskip 0.7ex
\hangindent=5.5em \hangafter=1
{\white .}\hskip 1em $\rho_\text{isum}(\mathfrak{t})$ =
 $( \frac{13}{120},
\frac{37}{120},
\frac{29}{240},
\frac{101}{240},
\frac{149}{240},
\frac{221}{240} )
$,

\vskip 0.7ex
\hangindent=5.5em \hangafter=1
{\white .}\hskip 1em $\rho_\text{isum}(\mathfrak{s})$ =
($0$,
$0$,
$\frac{1}{\sqrt{10}}c^{3}_{20}
$,
$\frac{1}{\sqrt{10}}c^{1}_{20}
$,
$\frac{1}{\sqrt{10}}c^{3}_{20}
$,
$\frac{1}{\sqrt{10}}c^{1}_{20}
$;
$0$,
$\frac{1}{\sqrt{10}}c^{1}_{20}
$,
$-\frac{1}{\sqrt{10}}c^{3}_{20}
$,
$\frac{1}{\sqrt{10}}c^{1}_{20}
$,
$-\frac{1}{\sqrt{10}}c^{3}_{20}
$;
$-\frac{1}{2\sqrt{5}}c^{1}_{20}
$,
$-\frac{1}{2\sqrt{5}}c^{3}_{20}
$,
$\frac{1}{2\sqrt{5}}c^{1}_{20}
$,
$\frac{1}{2\sqrt{5}}c^{3}_{20}
$;
$\frac{1}{2\sqrt{5}}c^{1}_{20}
$,
$\frac{1}{2\sqrt{5}}c^{3}_{20}
$,
$-\frac{1}{2\sqrt{5}}c^{1}_{20}
$;
$-\frac{1}{2\sqrt{5}}c^{1}_{20}
$,
$-\frac{1}{2\sqrt{5}}c^{3}_{20}
$;
$\frac{1}{2\sqrt{5}}c^{1}_{20}
$)

Pass. 

 \ \color{black}

 \color{blue}

\noindent 701: (dims,levels) = $(6;240
)$,
irreps = $3_{16}^{5,0}
\hskip -1.5pt \otimes \hskip -1.5pt
2_{5}^{1}
\hskip -1.5pt \otimes \hskip -1.5pt
1_{3}^{1,0}$,
pord$(\rho_\text{isum}(\mathfrak{t})) = 80$,

\vskip 0.7ex
\hangindent=5.5em \hangafter=1
{\white .}\hskip 1em $\rho_\text{isum}(\mathfrak{t})$ =
 $( \frac{19}{120},
\frac{91}{120},
\frac{83}{240},
\frac{107}{240},
\frac{203}{240},
\frac{227}{240} )
$,

\vskip 0.7ex
\hangindent=5.5em \hangafter=1
{\white .}\hskip 1em $\rho_\text{isum}(\mathfrak{s})$ =
($0$,
$0$,
$\frac{1}{\sqrt{10}}c^{3}_{20}
$,
$\frac{1}{\sqrt{10}}c^{1}_{20}
$,
$\frac{1}{\sqrt{10}}c^{3}_{20}
$,
$\frac{1}{\sqrt{10}}c^{1}_{20}
$;
$0$,
$\frac{1}{\sqrt{10}}c^{1}_{20}
$,
$-\frac{1}{\sqrt{10}}c^{3}_{20}
$,
$\frac{1}{\sqrt{10}}c^{1}_{20}
$,
$-\frac{1}{\sqrt{10}}c^{3}_{20}
$;
$-\frac{1}{2\sqrt{5}}c^{3}_{20}
$,
$\frac{1}{2\sqrt{5}}c^{1}_{20}
$,
$\frac{1}{2\sqrt{5}}c^{3}_{20}
$,
$-\frac{1}{2\sqrt{5}}c^{1}_{20}
$;
$\frac{1}{2\sqrt{5}}c^{3}_{20}
$,
$-\frac{1}{2\sqrt{5}}c^{1}_{20}
$,
$-\frac{1}{2\sqrt{5}}c^{3}_{20}
$;
$-\frac{1}{2\sqrt{5}}c^{3}_{20}
$,
$\frac{1}{2\sqrt{5}}c^{1}_{20}
$;
$\frac{1}{2\sqrt{5}}c^{3}_{20}
$)

Pass. 

 \ \color{black}

 \color{blue}

\noindent 702: (dims,levels) = $(6;240
)$,
irreps = $3_{16}^{1,0}
\hskip -1.5pt \otimes \hskip -1.5pt
2_{5}^{1}
\hskip -1.5pt \otimes \hskip -1.5pt
1_{3}^{1,0}$,
pord$(\rho_\text{isum}(\mathfrak{t})) = 80$,

\vskip 0.7ex
\hangindent=5.5em \hangafter=1
{\white .}\hskip 1em $\rho_\text{isum}(\mathfrak{t})$ =
 $( \frac{31}{120},
\frac{79}{120},
\frac{23}{240},
\frac{47}{240},
\frac{143}{240},
\frac{167}{240} )
$,

\vskip 0.7ex
\hangindent=5.5em \hangafter=1
{\white .}\hskip 1em $\rho_\text{isum}(\mathfrak{s})$ =
($0$,
$0$,
$\frac{1}{\sqrt{10}}c^{1}_{20}
$,
$\frac{1}{\sqrt{10}}c^{3}_{20}
$,
$\frac{1}{\sqrt{10}}c^{1}_{20}
$,
$\frac{1}{\sqrt{10}}c^{3}_{20}
$;
$0$,
$\frac{1}{\sqrt{10}}c^{3}_{20}
$,
$-\frac{1}{\sqrt{10}}c^{1}_{20}
$,
$\frac{1}{\sqrt{10}}c^{3}_{20}
$,
$-\frac{1}{\sqrt{10}}c^{1}_{20}
$;
$-\frac{1}{2\sqrt{5}}c^{3}_{20}
$,
$-\frac{1}{2\sqrt{5}}c^{1}_{20}
$,
$\frac{1}{2\sqrt{5}}c^{3}_{20}
$,
$\frac{1}{2\sqrt{5}}c^{1}_{20}
$;
$\frac{1}{2\sqrt{5}}c^{3}_{20}
$,
$\frac{1}{2\sqrt{5}}c^{1}_{20}
$,
$-\frac{1}{2\sqrt{5}}c^{3}_{20}
$;
$-\frac{1}{2\sqrt{5}}c^{3}_{20}
$,
$-\frac{1}{2\sqrt{5}}c^{1}_{20}
$;
$\frac{1}{2\sqrt{5}}c^{3}_{20}
$)

Pass. 

 \ \color{black}

 \color{blue}

\noindent 703: (dims,levels) = $(6;240
)$,
irreps = $3_{16}^{5,0}
\hskip -1.5pt \otimes \hskip -1.5pt
2_{5}^{2}
\hskip -1.5pt \otimes \hskip -1.5pt
1_{3}^{1,0}$,
pord$(\rho_\text{isum}(\mathfrak{t})) = 80$,

\vskip 0.7ex
\hangindent=5.5em \hangafter=1
{\white .}\hskip 1em $\rho_\text{isum}(\mathfrak{t})$ =
 $( \frac{43}{120},
\frac{67}{120},
\frac{11}{240},
\frac{59}{240},
\frac{131}{240},
\frac{179}{240} )
$,

\vskip 0.7ex
\hangindent=5.5em \hangafter=1
{\white .}\hskip 1em $\rho_\text{isum}(\mathfrak{s})$ =
($0$,
$0$,
$\frac{1}{\sqrt{10}}c^{1}_{20}
$,
$\frac{1}{\sqrt{10}}c^{3}_{20}
$,
$\frac{1}{\sqrt{10}}c^{1}_{20}
$,
$\frac{1}{\sqrt{10}}c^{3}_{20}
$;
$0$,
$\frac{1}{\sqrt{10}}c^{3}_{20}
$,
$-\frac{1}{\sqrt{10}}c^{1}_{20}
$,
$\frac{1}{\sqrt{10}}c^{3}_{20}
$,
$-\frac{1}{\sqrt{10}}c^{1}_{20}
$;
$-\frac{1}{2\sqrt{5}}c^{1}_{20}
$,
$-\frac{1}{2\sqrt{5}}c^{3}_{20}
$,
$\frac{1}{2\sqrt{5}}c^{1}_{20}
$,
$\frac{1}{2\sqrt{5}}c^{3}_{20}
$;
$\frac{1}{2\sqrt{5}}c^{1}_{20}
$,
$\frac{1}{2\sqrt{5}}c^{3}_{20}
$,
$-\frac{1}{2\sqrt{5}}c^{1}_{20}
$;
$-\frac{1}{2\sqrt{5}}c^{1}_{20}
$,
$-\frac{1}{2\sqrt{5}}c^{3}_{20}
$;
$\frac{1}{2\sqrt{5}}c^{1}_{20}
$)

Pass. 

 \ \color{black}

 \color{blue}

\noindent 704: (dims,levels) = $(6;240
)$,
irreps = $3_{16}^{3,0}
\hskip -1.5pt \otimes \hskip -1.5pt
2_{5}^{1}
\hskip -1.5pt \otimes \hskip -1.5pt
1_{3}^{1,0}$,
pord$(\rho_\text{isum}(\mathfrak{t})) = 80$,

\vskip 0.7ex
\hangindent=5.5em \hangafter=1
{\white .}\hskip 1em $\rho_\text{isum}(\mathfrak{t})$ =
 $( \frac{61}{120},
\frac{109}{120},
\frac{53}{240},
\frac{77}{240},
\frac{173}{240},
\frac{197}{240} )
$,

\vskip 0.7ex
\hangindent=5.5em \hangafter=1
{\white .}\hskip 1em $\rho_\text{isum}(\mathfrak{s})$ =
($0$,
$0$,
$\frac{1}{\sqrt{10}}c^{1}_{20}
$,
$\frac{1}{\sqrt{10}}c^{3}_{20}
$,
$\frac{1}{\sqrt{10}}c^{1}_{20}
$,
$\frac{1}{\sqrt{10}}c^{3}_{20}
$;
$0$,
$\frac{1}{\sqrt{10}}c^{3}_{20}
$,
$-\frac{1}{\sqrt{10}}c^{1}_{20}
$,
$\frac{1}{\sqrt{10}}c^{3}_{20}
$,
$-\frac{1}{\sqrt{10}}c^{1}_{20}
$;
$\frac{1}{2\sqrt{5}}c^{3}_{20}
$,
$\frac{1}{2\sqrt{5}}c^{1}_{20}
$,
$-\frac{1}{2\sqrt{5}}c^{3}_{20}
$,
$-\frac{1}{2\sqrt{5}}c^{1}_{20}
$;
$-\frac{1}{2\sqrt{5}}c^{3}_{20}
$,
$-\frac{1}{2\sqrt{5}}c^{1}_{20}
$,
$\frac{1}{2\sqrt{5}}c^{3}_{20}
$;
$\frac{1}{2\sqrt{5}}c^{3}_{20}
$,
$\frac{1}{2\sqrt{5}}c^{1}_{20}
$;
$-\frac{1}{2\sqrt{5}}c^{3}_{20}
$)

Pass. 

 \ \color{black}

 \color{blue}

\noindent 705: (dims,levels) = $(6;240
)$,
irreps = $3_{16}^{7,0}
\hskip -1.5pt \otimes \hskip -1.5pt
2_{5}^{2}
\hskip -1.5pt \otimes \hskip -1.5pt
1_{3}^{1,0}$,
pord$(\rho_\text{isum}(\mathfrak{t})) = 80$,

\vskip 0.7ex
\hangindent=5.5em \hangafter=1
{\white .}\hskip 1em $\rho_\text{isum}(\mathfrak{t})$ =
 $( \frac{73}{120},
\frac{97}{120},
\frac{41}{240},
\frac{89}{240},
\frac{161}{240},
\frac{209}{240} )
$,

\vskip 0.7ex
\hangindent=5.5em \hangafter=1
{\white .}\hskip 1em $\rho_\text{isum}(\mathfrak{s})$ =
($0$,
$0$,
$\frac{1}{\sqrt{10}}c^{1}_{20}
$,
$\frac{1}{\sqrt{10}}c^{3}_{20}
$,
$\frac{1}{\sqrt{10}}c^{1}_{20}
$,
$\frac{1}{\sqrt{10}}c^{3}_{20}
$;
$0$,
$\frac{1}{\sqrt{10}}c^{3}_{20}
$,
$-\frac{1}{\sqrt{10}}c^{1}_{20}
$,
$\frac{1}{\sqrt{10}}c^{3}_{20}
$,
$-\frac{1}{\sqrt{10}}c^{1}_{20}
$;
$\frac{1}{2\sqrt{5}}c^{1}_{20}
$,
$\frac{1}{2\sqrt{5}}c^{3}_{20}
$,
$-\frac{1}{2\sqrt{5}}c^{1}_{20}
$,
$-\frac{1}{2\sqrt{5}}c^{3}_{20}
$;
$-\frac{1}{2\sqrt{5}}c^{1}_{20}
$,
$-\frac{1}{2\sqrt{5}}c^{3}_{20}
$,
$\frac{1}{2\sqrt{5}}c^{1}_{20}
$;
$\frac{1}{2\sqrt{5}}c^{1}_{20}
$,
$\frac{1}{2\sqrt{5}}c^{3}_{20}
$;
$-\frac{1}{2\sqrt{5}}c^{1}_{20}
$)

Pass. 

 \ \color{black}

 \color{blue}

\noindent 706: (dims,levels) = $(6;420
)$,
irreps = $3_{7}^{3}
\hskip -1.5pt \otimes \hskip -1.5pt
2_{5}^{2}
\hskip -1.5pt \otimes \hskip -1.5pt
1_{4}^{1,0}
\hskip -1.5pt \otimes \hskip -1.5pt
1_{3}^{1,0}$,
pord$(\rho_\text{isum}(\mathfrak{t})) = 35$,

\vskip 0.7ex
\hangindent=5.5em \hangafter=1
{\white .}\hskip 1em $\rho_\text{isum}(\mathfrak{t})$ =
 $( \frac{17}{420},
\frac{173}{420},
\frac{257}{420},
\frac{293}{420},
\frac{353}{420},
\frac{377}{420} )
$,

\vskip 0.7ex
\hangindent=5.5em \hangafter=1
{\white .}\hskip 1em $\rho_\text{isum}(\mathfrak{s})$ =
($\frac{2}{\sqrt{35}}c^{3}_{35}
+\frac{1}{\sqrt{35}}c^{4}_{35}
+\frac{1}{\sqrt{35}}c^{10}_{35}
+\frac{1}{\sqrt{35}}c^{11}_{35}
$,
$-\frac{1}{\sqrt{35}}c^{1}_{35}
+\frac{1}{\sqrt{35}}c^{6}_{35}
$,
$-\frac{1}{\sqrt{35}\mathrm{i}}s^{3}_{140}
-\frac{1}{\sqrt{35}\mathrm{i}}s^{17}_{140}
$,
$-\frac{2}{35}c^{1}_{140}
+\frac{1}{35}c^{3}_{140}
+\frac{1}{7}c^{5}_{140}
+\frac{3}{35}c^{7}_{140}
-\frac{1}{5}c^{9}_{140}
+\frac{2}{35}c^{13}_{140}
+\frac{1}{35}c^{15}_{140}
+\frac{1}{35}c^{17}_{140}
-\frac{3}{35}c^{19}_{140}
-\frac{2}{35}c^{21}_{140}
+\frac{2}{7}c^{23}_{140}
$,
$\frac{1}{\sqrt{35}}c^{4}_{35}
-\frac{1}{\sqrt{35}}c^{11}_{35}
$,
$\frac{4}{35}c^{1}_{140}
+\frac{3}{35}c^{3}_{140}
+\frac{1}{7}c^{5}_{140}
-\frac{1}{35}c^{7}_{140}
-\frac{1}{35}c^{9}_{140}
-\frac{4}{35}c^{13}_{140}
-\frac{2}{35}c^{15}_{140}
+\frac{3}{35}c^{17}_{140}
-\frac{9}{35}c^{19}_{140}
+\frac{4}{35}c^{21}_{140}
+\frac{2}{7}c^{23}_{140}
$;
$-\frac{4}{35}c^{1}_{140}
-\frac{3}{35}c^{3}_{140}
-\frac{1}{7}c^{5}_{140}
+\frac{1}{35}c^{7}_{140}
+\frac{1}{35}c^{9}_{140}
+\frac{4}{35}c^{13}_{140}
+\frac{2}{35}c^{15}_{140}
-\frac{3}{35}c^{17}_{140}
+\frac{9}{35}c^{19}_{140}
-\frac{4}{35}c^{21}_{140}
-\frac{2}{7}c^{23}_{140}
$,
$-\frac{2}{35}c^{1}_{140}
+\frac{1}{35}c^{3}_{140}
+\frac{1}{7}c^{5}_{140}
+\frac{3}{35}c^{7}_{140}
-\frac{1}{5}c^{9}_{140}
+\frac{2}{35}c^{13}_{140}
+\frac{1}{35}c^{15}_{140}
+\frac{1}{35}c^{17}_{140}
-\frac{3}{35}c^{19}_{140}
-\frac{2}{35}c^{21}_{140}
+\frac{2}{7}c^{23}_{140}
$,
$-\frac{2}{\sqrt{35}}c^{3}_{35}
-\frac{1}{\sqrt{35}}c^{4}_{35}
-\frac{1}{\sqrt{35}}c^{10}_{35}
-\frac{1}{\sqrt{35}}c^{11}_{35}
$,
$\frac{1}{\sqrt{35}\mathrm{i}}s^{3}_{140}
+\frac{1}{\sqrt{35}\mathrm{i}}s^{17}_{140}
$,
$\frac{1}{\sqrt{35}}c^{4}_{35}
-\frac{1}{\sqrt{35}}c^{11}_{35}
$;
$\frac{4}{35}c^{1}_{140}
+\frac{3}{35}c^{3}_{140}
+\frac{1}{7}c^{5}_{140}
-\frac{1}{35}c^{7}_{140}
-\frac{1}{35}c^{9}_{140}
-\frac{4}{35}c^{13}_{140}
-\frac{2}{35}c^{15}_{140}
+\frac{3}{35}c^{17}_{140}
-\frac{9}{35}c^{19}_{140}
+\frac{4}{35}c^{21}_{140}
+\frac{2}{7}c^{23}_{140}
$,
$\frac{1}{\sqrt{35}}c^{4}_{35}
-\frac{1}{\sqrt{35}}c^{11}_{35}
$,
$-\frac{1}{\sqrt{35}}c^{1}_{35}
+\frac{1}{\sqrt{35}}c^{6}_{35}
$,
$\frac{2}{\sqrt{35}}c^{3}_{35}
+\frac{1}{\sqrt{35}}c^{4}_{35}
+\frac{1}{\sqrt{35}}c^{10}_{35}
+\frac{1}{\sqrt{35}}c^{11}_{35}
$;
$\frac{1}{\sqrt{35}\mathrm{i}}s^{3}_{140}
+\frac{1}{\sqrt{35}\mathrm{i}}s^{17}_{140}
$,
$-\frac{4}{35}c^{1}_{140}
-\frac{3}{35}c^{3}_{140}
-\frac{1}{7}c^{5}_{140}
+\frac{1}{35}c^{7}_{140}
+\frac{1}{35}c^{9}_{140}
+\frac{4}{35}c^{13}_{140}
+\frac{2}{35}c^{15}_{140}
-\frac{3}{35}c^{17}_{140}
+\frac{9}{35}c^{19}_{140}
-\frac{4}{35}c^{21}_{140}
-\frac{2}{7}c^{23}_{140}
$,
$-\frac{1}{\sqrt{35}}c^{1}_{35}
+\frac{1}{\sqrt{35}}c^{6}_{35}
$;
$-\frac{2}{\sqrt{35}}c^{3}_{35}
-\frac{1}{\sqrt{35}}c^{4}_{35}
-\frac{1}{\sqrt{35}}c^{10}_{35}
-\frac{1}{\sqrt{35}}c^{11}_{35}
$,
$-\frac{2}{35}c^{1}_{140}
+\frac{1}{35}c^{3}_{140}
+\frac{1}{7}c^{5}_{140}
+\frac{3}{35}c^{7}_{140}
-\frac{1}{5}c^{9}_{140}
+\frac{2}{35}c^{13}_{140}
+\frac{1}{35}c^{15}_{140}
+\frac{1}{35}c^{17}_{140}
-\frac{3}{35}c^{19}_{140}
-\frac{2}{35}c^{21}_{140}
+\frac{2}{7}c^{23}_{140}
$;
$-\frac{1}{\sqrt{35}\mathrm{i}}s^{3}_{140}
-\frac{1}{\sqrt{35}\mathrm{i}}s^{17}_{140}
$)

Pass. 

 \ \color{black}

 \color{blue}

\noindent 707: (dims,levels) = $(6;420
)$,
irreps = $3_{7}^{1}
\hskip -1.5pt \otimes \hskip -1.5pt
2_{5}^{1}
\hskip -1.5pt \otimes \hskip -1.5pt
1_{4}^{1,0}
\hskip -1.5pt \otimes \hskip -1.5pt
1_{3}^{1,0}$,
pord$(\rho_\text{isum}(\mathfrak{t})) = 35$,

\vskip 0.7ex
\hangindent=5.5em \hangafter=1
{\white .}\hskip 1em $\rho_\text{isum}(\mathfrak{t})$ =
 $( \frac{29}{420},
\frac{149}{420},
\frac{221}{420},
\frac{281}{420},
\frac{389}{420},
\frac{401}{420} )
$,

\vskip 0.7ex
\hangindent=5.5em \hangafter=1
{\white .}\hskip 1em $\rho_\text{isum}(\mathfrak{s})$ =
($\frac{1}{\sqrt{35}}c^{1}_{35}
-\frac{1}{\sqrt{35}}c^{6}_{35}
$,
$\frac{1}{\sqrt{35}}c^{4}_{35}
-\frac{1}{\sqrt{35}}c^{11}_{35}
$,
$\frac{4}{35}c^{1}_{140}
+\frac{3}{35}c^{3}_{140}
+\frac{1}{7}c^{5}_{140}
-\frac{1}{35}c^{7}_{140}
-\frac{1}{35}c^{9}_{140}
-\frac{4}{35}c^{13}_{140}
-\frac{2}{35}c^{15}_{140}
+\frac{3}{35}c^{17}_{140}
-\frac{9}{35}c^{19}_{140}
+\frac{4}{35}c^{21}_{140}
+\frac{2}{7}c^{23}_{140}
$,
$-\frac{1}{\sqrt{35}\mathrm{i}}s^{3}_{140}
-\frac{1}{\sqrt{35}\mathrm{i}}s^{17}_{140}
$,
$-\frac{2}{35}c^{1}_{140}
+\frac{1}{35}c^{3}_{140}
+\frac{1}{7}c^{5}_{140}
+\frac{3}{35}c^{7}_{140}
-\frac{1}{5}c^{9}_{140}
+\frac{2}{35}c^{13}_{140}
+\frac{1}{35}c^{15}_{140}
+\frac{1}{35}c^{17}_{140}
-\frac{3}{35}c^{19}_{140}
-\frac{2}{35}c^{21}_{140}
+\frac{2}{7}c^{23}_{140}
$,
$\frac{2}{\sqrt{35}}c^{3}_{35}
+\frac{1}{\sqrt{35}}c^{4}_{35}
+\frac{1}{\sqrt{35}}c^{10}_{35}
+\frac{1}{\sqrt{35}}c^{11}_{35}
$;
$\frac{2}{35}c^{1}_{140}
-\frac{1}{35}c^{3}_{140}
-\frac{1}{7}c^{5}_{140}
-\frac{3}{35}c^{7}_{140}
+\frac{1}{5}c^{9}_{140}
-\frac{2}{35}c^{13}_{140}
-\frac{1}{35}c^{15}_{140}
-\frac{1}{35}c^{17}_{140}
+\frac{3}{35}c^{19}_{140}
+\frac{2}{35}c^{21}_{140}
-\frac{2}{7}c^{23}_{140}
$,
$\frac{1}{\sqrt{35}\mathrm{i}}s^{3}_{140}
+\frac{1}{\sqrt{35}\mathrm{i}}s^{17}_{140}
$,
$-\frac{2}{\sqrt{35}}c^{3}_{35}
-\frac{1}{\sqrt{35}}c^{4}_{35}
-\frac{1}{\sqrt{35}}c^{10}_{35}
-\frac{1}{\sqrt{35}}c^{11}_{35}
$,
$\frac{1}{\sqrt{35}}c^{1}_{35}
-\frac{1}{\sqrt{35}}c^{6}_{35}
$,
$-\frac{4}{35}c^{1}_{140}
-\frac{3}{35}c^{3}_{140}
-\frac{1}{7}c^{5}_{140}
+\frac{1}{35}c^{7}_{140}
+\frac{1}{35}c^{9}_{140}
+\frac{4}{35}c^{13}_{140}
+\frac{2}{35}c^{15}_{140}
-\frac{3}{35}c^{17}_{140}
+\frac{9}{35}c^{19}_{140}
-\frac{4}{35}c^{21}_{140}
-\frac{2}{7}c^{23}_{140}
$;
$\frac{1}{\sqrt{35}}c^{4}_{35}
-\frac{1}{\sqrt{35}}c^{11}_{35}
$,
$-\frac{2}{35}c^{1}_{140}
+\frac{1}{35}c^{3}_{140}
+\frac{1}{7}c^{5}_{140}
+\frac{3}{35}c^{7}_{140}
-\frac{1}{5}c^{9}_{140}
+\frac{2}{35}c^{13}_{140}
+\frac{1}{35}c^{15}_{140}
+\frac{1}{35}c^{17}_{140}
-\frac{3}{35}c^{19}_{140}
-\frac{2}{35}c^{21}_{140}
+\frac{2}{7}c^{23}_{140}
$,
$-\frac{2}{\sqrt{35}}c^{3}_{35}
-\frac{1}{\sqrt{35}}c^{4}_{35}
-\frac{1}{\sqrt{35}}c^{10}_{35}
-\frac{1}{\sqrt{35}}c^{11}_{35}
$,
$-\frac{1}{\sqrt{35}}c^{1}_{35}
+\frac{1}{\sqrt{35}}c^{6}_{35}
$;
$-\frac{1}{\sqrt{35}}c^{1}_{35}
+\frac{1}{\sqrt{35}}c^{6}_{35}
$,
$-\frac{4}{35}c^{1}_{140}
-\frac{3}{35}c^{3}_{140}
-\frac{1}{7}c^{5}_{140}
+\frac{1}{35}c^{7}_{140}
+\frac{1}{35}c^{9}_{140}
+\frac{4}{35}c^{13}_{140}
+\frac{2}{35}c^{15}_{140}
-\frac{3}{35}c^{17}_{140}
+\frac{9}{35}c^{19}_{140}
-\frac{4}{35}c^{21}_{140}
-\frac{2}{7}c^{23}_{140}
$,
$\frac{1}{\sqrt{35}}c^{4}_{35}
-\frac{1}{\sqrt{35}}c^{11}_{35}
$;
$-\frac{1}{\sqrt{35}}c^{4}_{35}
+\frac{1}{\sqrt{35}}c^{11}_{35}
$,
$\frac{1}{\sqrt{35}\mathrm{i}}s^{3}_{140}
+\frac{1}{\sqrt{35}\mathrm{i}}s^{17}_{140}
$;
$-\frac{2}{35}c^{1}_{140}
+\frac{1}{35}c^{3}_{140}
+\frac{1}{7}c^{5}_{140}
+\frac{3}{35}c^{7}_{140}
-\frac{1}{5}c^{9}_{140}
+\frac{2}{35}c^{13}_{140}
+\frac{1}{35}c^{15}_{140}
+\frac{1}{35}c^{17}_{140}
-\frac{3}{35}c^{19}_{140}
-\frac{2}{35}c^{21}_{140}
+\frac{2}{7}c^{23}_{140}
$)

Pass. 

 \ \color{black}

 \color{blue}

\noindent 708: (dims,levels) = $(6;420
)$,
irreps = $3_{7}^{3}
\hskip -1.5pt \otimes \hskip -1.5pt
2_{5}^{1}
\hskip -1.5pt \otimes \hskip -1.5pt
1_{4}^{1,0}
\hskip -1.5pt \otimes \hskip -1.5pt
1_{3}^{1,0}$,
pord$(\rho_\text{isum}(\mathfrak{t})) = 35$,

\vskip 0.7ex
\hangindent=5.5em \hangafter=1
{\white .}\hskip 1em $\rho_\text{isum}(\mathfrak{t})$ =
 $( \frac{41}{420},
\frac{89}{420},
\frac{101}{420},
\frac{209}{420},
\frac{269}{420},
\frac{341}{420} )
$,

\vskip 0.7ex
\hangindent=5.5em \hangafter=1
{\white .}\hskip 1em $\rho_\text{isum}(\mathfrak{s})$ =
($-\frac{1}{\sqrt{35}}c^{1}_{35}
+\frac{1}{\sqrt{35}}c^{6}_{35}
$,
$\frac{2}{\sqrt{35}}c^{3}_{35}
+\frac{1}{\sqrt{35}}c^{4}_{35}
+\frac{1}{\sqrt{35}}c^{10}_{35}
+\frac{1}{\sqrt{35}}c^{11}_{35}
$,
$-\frac{2}{35}c^{1}_{140}
+\frac{1}{35}c^{3}_{140}
+\frac{1}{7}c^{5}_{140}
+\frac{3}{35}c^{7}_{140}
-\frac{1}{5}c^{9}_{140}
+\frac{2}{35}c^{13}_{140}
+\frac{1}{35}c^{15}_{140}
+\frac{1}{35}c^{17}_{140}
-\frac{3}{35}c^{19}_{140}
-\frac{2}{35}c^{21}_{140}
+\frac{2}{7}c^{23}_{140}
$,
$-\frac{1}{\sqrt{35}\mathrm{i}}s^{3}_{140}
-\frac{1}{\sqrt{35}\mathrm{i}}s^{17}_{140}
$,
$\frac{4}{35}c^{1}_{140}
+\frac{3}{35}c^{3}_{140}
+\frac{1}{7}c^{5}_{140}
-\frac{1}{35}c^{7}_{140}
-\frac{1}{35}c^{9}_{140}
-\frac{4}{35}c^{13}_{140}
-\frac{2}{35}c^{15}_{140}
+\frac{3}{35}c^{17}_{140}
-\frac{9}{35}c^{19}_{140}
+\frac{4}{35}c^{21}_{140}
+\frac{2}{7}c^{23}_{140}
$,
$\frac{1}{\sqrt{35}}c^{4}_{35}
-\frac{1}{\sqrt{35}}c^{11}_{35}
$;
$\frac{2}{35}c^{1}_{140}
-\frac{1}{35}c^{3}_{140}
-\frac{1}{7}c^{5}_{140}
-\frac{3}{35}c^{7}_{140}
+\frac{1}{5}c^{9}_{140}
-\frac{2}{35}c^{13}_{140}
-\frac{1}{35}c^{15}_{140}
-\frac{1}{35}c^{17}_{140}
+\frac{3}{35}c^{19}_{140}
+\frac{2}{35}c^{21}_{140}
-\frac{2}{7}c^{23}_{140}
$,
$-\frac{1}{\sqrt{35}\mathrm{i}}s^{3}_{140}
-\frac{1}{\sqrt{35}\mathrm{i}}s^{17}_{140}
$,
$-\frac{1}{\sqrt{35}}c^{4}_{35}
+\frac{1}{\sqrt{35}}c^{11}_{35}
$,
$\frac{1}{\sqrt{35}}c^{1}_{35}
-\frac{1}{\sqrt{35}}c^{6}_{35}
$,
$\frac{4}{35}c^{1}_{140}
+\frac{3}{35}c^{3}_{140}
+\frac{1}{7}c^{5}_{140}
-\frac{1}{35}c^{7}_{140}
-\frac{1}{35}c^{9}_{140}
-\frac{4}{35}c^{13}_{140}
-\frac{2}{35}c^{15}_{140}
+\frac{3}{35}c^{17}_{140}
-\frac{9}{35}c^{19}_{140}
+\frac{4}{35}c^{21}_{140}
+\frac{2}{7}c^{23}_{140}
$;
$\frac{1}{\sqrt{35}}c^{4}_{35}
-\frac{1}{\sqrt{35}}c^{11}_{35}
$,
$\frac{4}{35}c^{1}_{140}
+\frac{3}{35}c^{3}_{140}
+\frac{1}{7}c^{5}_{140}
-\frac{1}{35}c^{7}_{140}
-\frac{1}{35}c^{9}_{140}
-\frac{4}{35}c^{13}_{140}
-\frac{2}{35}c^{15}_{140}
+\frac{3}{35}c^{17}_{140}
-\frac{9}{35}c^{19}_{140}
+\frac{4}{35}c^{21}_{140}
+\frac{2}{7}c^{23}_{140}
$,
$\frac{2}{\sqrt{35}}c^{3}_{35}
+\frac{1}{\sqrt{35}}c^{4}_{35}
+\frac{1}{\sqrt{35}}c^{10}_{35}
+\frac{1}{\sqrt{35}}c^{11}_{35}
$,
$-\frac{1}{\sqrt{35}}c^{1}_{35}
+\frac{1}{\sqrt{35}}c^{6}_{35}
$;
$\frac{1}{\sqrt{35}}c^{1}_{35}
-\frac{1}{\sqrt{35}}c^{6}_{35}
$,
$\frac{2}{35}c^{1}_{140}
-\frac{1}{35}c^{3}_{140}
-\frac{1}{7}c^{5}_{140}
-\frac{3}{35}c^{7}_{140}
+\frac{1}{5}c^{9}_{140}
-\frac{2}{35}c^{13}_{140}
-\frac{1}{35}c^{15}_{140}
-\frac{1}{35}c^{17}_{140}
+\frac{3}{35}c^{19}_{140}
+\frac{2}{35}c^{21}_{140}
-\frac{2}{7}c^{23}_{140}
$,
$\frac{2}{\sqrt{35}}c^{3}_{35}
+\frac{1}{\sqrt{35}}c^{4}_{35}
+\frac{1}{\sqrt{35}}c^{10}_{35}
+\frac{1}{\sqrt{35}}c^{11}_{35}
$;
$-\frac{1}{\sqrt{35}}c^{4}_{35}
+\frac{1}{\sqrt{35}}c^{11}_{35}
$,
$-\frac{1}{\sqrt{35}\mathrm{i}}s^{3}_{140}
-\frac{1}{\sqrt{35}\mathrm{i}}s^{17}_{140}
$;
$-\frac{2}{35}c^{1}_{140}
+\frac{1}{35}c^{3}_{140}
+\frac{1}{7}c^{5}_{140}
+\frac{3}{35}c^{7}_{140}
-\frac{1}{5}c^{9}_{140}
+\frac{2}{35}c^{13}_{140}
+\frac{1}{35}c^{15}_{140}
+\frac{1}{35}c^{17}_{140}
-\frac{3}{35}c^{19}_{140}
-\frac{2}{35}c^{21}_{140}
+\frac{2}{7}c^{23}_{140}
$)

Pass. 

 \ \color{black}

 \color{blue}

\noindent 709: (dims,levels) = $(6;420
)$,
irreps = $3_{7}^{1}
\hskip -1.5pt \otimes \hskip -1.5pt
2_{5}^{2}
\hskip -1.5pt \otimes \hskip -1.5pt
1_{4}^{1,0}
\hskip -1.5pt \otimes \hskip -1.5pt
1_{3}^{1,0}$,
pord$(\rho_\text{isum}(\mathfrak{t})) = 35$,

\vskip 0.7ex
\hangindent=5.5em \hangafter=1
{\white .}\hskip 1em $\rho_\text{isum}(\mathfrak{t})$ =
 $( \frac{53}{420},
\frac{113}{420},
\frac{137}{420},
\frac{197}{420},
\frac{233}{420},
\frac{317}{420} )
$,

\vskip 0.7ex
\hangindent=5.5em \hangafter=1
{\white .}\hskip 1em $\rho_\text{isum}(\mathfrak{s})$ =
($-\frac{2}{\sqrt{35}}c^{3}_{35}
-\frac{1}{\sqrt{35}}c^{4}_{35}
-\frac{1}{\sqrt{35}}c^{10}_{35}
-\frac{1}{\sqrt{35}}c^{11}_{35}
$,
$\frac{4}{35}c^{1}_{140}
+\frac{3}{35}c^{3}_{140}
+\frac{1}{7}c^{5}_{140}
-\frac{1}{35}c^{7}_{140}
-\frac{1}{35}c^{9}_{140}
-\frac{4}{35}c^{13}_{140}
-\frac{2}{35}c^{15}_{140}
+\frac{3}{35}c^{17}_{140}
-\frac{9}{35}c^{19}_{140}
+\frac{4}{35}c^{21}_{140}
+\frac{2}{7}c^{23}_{140}
$,
$\frac{1}{\sqrt{35}}c^{4}_{35}
-\frac{1}{\sqrt{35}}c^{11}_{35}
$,
$-\frac{2}{35}c^{1}_{140}
+\frac{1}{35}c^{3}_{140}
+\frac{1}{7}c^{5}_{140}
+\frac{3}{35}c^{7}_{140}
-\frac{1}{5}c^{9}_{140}
+\frac{2}{35}c^{13}_{140}
+\frac{1}{35}c^{15}_{140}
+\frac{1}{35}c^{17}_{140}
-\frac{3}{35}c^{19}_{140}
-\frac{2}{35}c^{21}_{140}
+\frac{2}{7}c^{23}_{140}
$,
$-\frac{1}{\sqrt{35}\mathrm{i}}s^{3}_{140}
-\frac{1}{\sqrt{35}\mathrm{i}}s^{17}_{140}
$,
$-\frac{1}{\sqrt{35}}c^{1}_{35}
+\frac{1}{\sqrt{35}}c^{6}_{35}
$;
$\frac{1}{\sqrt{35}\mathrm{i}}s^{3}_{140}
+\frac{1}{\sqrt{35}\mathrm{i}}s^{17}_{140}
$,
$\frac{2}{35}c^{1}_{140}
-\frac{1}{35}c^{3}_{140}
-\frac{1}{7}c^{5}_{140}
-\frac{3}{35}c^{7}_{140}
+\frac{1}{5}c^{9}_{140}
-\frac{2}{35}c^{13}_{140}
-\frac{1}{35}c^{15}_{140}
-\frac{1}{35}c^{17}_{140}
+\frac{3}{35}c^{19}_{140}
+\frac{2}{35}c^{21}_{140}
-\frac{2}{7}c^{23}_{140}
$,
$\frac{1}{\sqrt{35}}c^{1}_{35}
-\frac{1}{\sqrt{35}}c^{6}_{35}
$,
$-\frac{2}{\sqrt{35}}c^{3}_{35}
-\frac{1}{\sqrt{35}}c^{4}_{35}
-\frac{1}{\sqrt{35}}c^{10}_{35}
-\frac{1}{\sqrt{35}}c^{11}_{35}
$,
$-\frac{1}{\sqrt{35}}c^{4}_{35}
+\frac{1}{\sqrt{35}}c^{11}_{35}
$;
$\frac{2}{\sqrt{35}}c^{3}_{35}
+\frac{1}{\sqrt{35}}c^{4}_{35}
+\frac{1}{\sqrt{35}}c^{10}_{35}
+\frac{1}{\sqrt{35}}c^{11}_{35}
$,
$\frac{4}{35}c^{1}_{140}
+\frac{3}{35}c^{3}_{140}
+\frac{1}{7}c^{5}_{140}
-\frac{1}{35}c^{7}_{140}
-\frac{1}{35}c^{9}_{140}
-\frac{4}{35}c^{13}_{140}
-\frac{2}{35}c^{15}_{140}
+\frac{3}{35}c^{17}_{140}
-\frac{9}{35}c^{19}_{140}
+\frac{4}{35}c^{21}_{140}
+\frac{2}{7}c^{23}_{140}
$,
$\frac{1}{\sqrt{35}}c^{1}_{35}
-\frac{1}{\sqrt{35}}c^{6}_{35}
$,
$-\frac{1}{\sqrt{35}\mathrm{i}}s^{3}_{140}
-\frac{1}{\sqrt{35}\mathrm{i}}s^{17}_{140}
$;
$-\frac{1}{\sqrt{35}\mathrm{i}}s^{3}_{140}
-\frac{1}{\sqrt{35}\mathrm{i}}s^{17}_{140}
$,
$-\frac{1}{\sqrt{35}}c^{4}_{35}
+\frac{1}{\sqrt{35}}c^{11}_{35}
$,
$\frac{2}{\sqrt{35}}c^{3}_{35}
+\frac{1}{\sqrt{35}}c^{4}_{35}
+\frac{1}{\sqrt{35}}c^{10}_{35}
+\frac{1}{\sqrt{35}}c^{11}_{35}
$;
$-\frac{4}{35}c^{1}_{140}
-\frac{3}{35}c^{3}_{140}
-\frac{1}{7}c^{5}_{140}
+\frac{1}{35}c^{7}_{140}
+\frac{1}{35}c^{9}_{140}
+\frac{4}{35}c^{13}_{140}
+\frac{2}{35}c^{15}_{140}
-\frac{3}{35}c^{17}_{140}
+\frac{9}{35}c^{19}_{140}
-\frac{4}{35}c^{21}_{140}
-\frac{2}{7}c^{23}_{140}
$,
$\frac{2}{35}c^{1}_{140}
-\frac{1}{35}c^{3}_{140}
-\frac{1}{7}c^{5}_{140}
-\frac{3}{35}c^{7}_{140}
+\frac{1}{5}c^{9}_{140}
-\frac{2}{35}c^{13}_{140}
-\frac{1}{35}c^{15}_{140}
-\frac{1}{35}c^{17}_{140}
+\frac{3}{35}c^{19}_{140}
+\frac{2}{35}c^{21}_{140}
-\frac{2}{7}c^{23}_{140}
$;
$\frac{4}{35}c^{1}_{140}
+\frac{3}{35}c^{3}_{140}
+\frac{1}{7}c^{5}_{140}
-\frac{1}{35}c^{7}_{140}
-\frac{1}{35}c^{9}_{140}
-\frac{4}{35}c^{13}_{140}
-\frac{2}{35}c^{15}_{140}
+\frac{3}{35}c^{17}_{140}
-\frac{9}{35}c^{19}_{140}
+\frac{4}{35}c^{21}_{140}
+\frac{2}{7}c^{23}_{140}
$)

Pass. 

 \ \color{black}

\

\section{Rank-2 modular data}
\label{Section5}

\subsection{A list of 2-dimensional irrep-sum $\SL$ representations}

The following is a list of 2-dimensional irrep-sum $\SL$ representations, with
3 types of representations omitted and one type skipped.  For details and notations, see Section
3 of this file.

\

 \color{blue}

\noindent 1: (dims,levels) = $(2;5
)$,
irreps = $2_{5}^{1}$,
pord$(\rho_\text{isum}(\mathfrak{t})) = 5$,

\vskip 0.7ex
\hangindent=5.5em \hangafter=1
{\white .}\hskip 1em $\rho_\text{isum}(\mathfrak{t})$ =
 $( \frac{1}{5},
\frac{4}{5} )
$,

\vskip 0.7ex
\hangindent=5.5em \hangafter=1
{\white .}\hskip 1em $\rho_\text{isum}(\mathfrak{s})$ =
$\mathrm{i}$($-\frac{1}{\sqrt{5}}c^{3}_{20}
$,
$\frac{1}{\sqrt{5}}c^{1}_{20}
$;\ \ 
$\frac{1}{\sqrt{5}}c^{3}_{20}
$)

Pass. 

 \ \color{black}

 \color{blue}

\noindent 2: (dims,levels) = $(2;5
)$,
irreps = $2_{5}^{2}$,
pord$(\rho_\text{isum}(\mathfrak{t})) = 5$,

\vskip 0.7ex
\hangindent=5.5em \hangafter=1
{\white .}\hskip 1em $\rho_\text{isum}(\mathfrak{t})$ =
 $( \frac{2}{5},
\frac{3}{5} )
$,

\vskip 0.7ex
\hangindent=5.5em \hangafter=1
{\white .}\hskip 1em $\rho_\text{isum}(\mathfrak{s})$ =
$\mathrm{i}$($-\frac{1}{\sqrt{5}}c^{1}_{20}
$,
$\frac{1}{\sqrt{5}}c^{3}_{20}
$;\ \ 
$\frac{1}{\sqrt{5}}c^{1}_{20}
$)

Pass. 

 \ \color{black}

 \color{blue}

\noindent 3: (dims,levels) = $(2;10
)$,
irreps = $2_{5}^{2}
\hskip -1.5pt \otimes \hskip -1.5pt
1_{2}^{1,0}$,
pord$(\rho_\text{isum}(\mathfrak{t})) = 5$,

\vskip 0.7ex
\hangindent=5.5em \hangafter=1
{\white .}\hskip 1em $\rho_\text{isum}(\mathfrak{t})$ =
 $( \frac{1}{10},
\frac{9}{10} )
$,

\vskip 0.7ex
\hangindent=5.5em \hangafter=1
{\white .}\hskip 1em $\rho_\text{isum}(\mathfrak{s})$ =
$\mathrm{i}$($-\frac{1}{\sqrt{5}}c^{1}_{20}
$,
$\frac{1}{\sqrt{5}}c^{3}_{20}
$;\ \ 
$\frac{1}{\sqrt{5}}c^{1}_{20}
$)

Pass. 

 \ \color{black}

 \color{blue}

\noindent 4: (dims,levels) = $(2;10
)$,
irreps = $2_{5}^{1}
\hskip -1.5pt \otimes \hskip -1.5pt
1_{2}^{1,0}$,
pord$(\rho_\text{isum}(\mathfrak{t})) = 5$,

\vskip 0.7ex
\hangindent=5.5em \hangafter=1
{\white .}\hskip 1em $\rho_\text{isum}(\mathfrak{t})$ =
 $( \frac{3}{10},
\frac{7}{10} )
$,

\vskip 0.7ex
\hangindent=5.5em \hangafter=1
{\white .}\hskip 1em $\rho_\text{isum}(\mathfrak{s})$ =
$\mathrm{i}$($-\frac{1}{\sqrt{5}}c^{3}_{20}
$,
$\frac{1}{\sqrt{5}}c^{1}_{20}
$;\ \ 
$\frac{1}{\sqrt{5}}c^{3}_{20}
$)

Pass. 

 \ \color{black}

 \color{blue}

\noindent 5: (dims,levels) = $(2;15
)$,
irreps = $2_{5}^{1}
\hskip -1.5pt \otimes \hskip -1.5pt
1_{3}^{1,0}$,
pord$(\rho_\text{isum}(\mathfrak{t})) = 5$,

\vskip 0.7ex
\hangindent=5.5em \hangafter=1
{\white .}\hskip 1em $\rho_\text{isum}(\mathfrak{t})$ =
 $( \frac{2}{15},
\frac{8}{15} )
$,

\vskip 0.7ex
\hangindent=5.5em \hangafter=1
{\white .}\hskip 1em $\rho_\text{isum}(\mathfrak{s})$ =
$\mathrm{i}$($\frac{1}{\sqrt{5}}c^{3}_{20}
$,
$\frac{1}{\sqrt{5}}c^{1}_{20}
$;\ \ 
$-\frac{1}{\sqrt{5}}c^{3}_{20}
$)

Pass. 

 \ \color{black}

 \color{blue}

\noindent 6: (dims,levels) = $(2;15
)$,
irreps = $2_{5}^{2}
\hskip -1.5pt \otimes \hskip -1.5pt
1_{3}^{1,0}$,
pord$(\rho_\text{isum}(\mathfrak{t})) = 5$,

\vskip 0.7ex
\hangindent=5.5em \hangafter=1
{\white .}\hskip 1em $\rho_\text{isum}(\mathfrak{t})$ =
 $( \frac{11}{15},
\frac{14}{15} )
$,

\vskip 0.7ex
\hangindent=5.5em \hangafter=1
{\white .}\hskip 1em $\rho_\text{isum}(\mathfrak{s})$ =
$\mathrm{i}$($-\frac{1}{\sqrt{5}}c^{1}_{20}
$,
$\frac{1}{\sqrt{5}}c^{3}_{20}
$;\ \ 
$\frac{1}{\sqrt{5}}c^{1}_{20}
$)

Pass. 

 \ \color{black}

 \color{blue}

\noindent 7: (dims,levels) = $(2;20
)$,
irreps = $2_{5}^{1}
\hskip -1.5pt \otimes \hskip -1.5pt
1_{4}^{1,0}$,
pord$(\rho_\text{isum}(\mathfrak{t})) = 5$,

\vskip 0.7ex
\hangindent=5.5em \hangafter=1
{\white .}\hskip 1em $\rho_\text{isum}(\mathfrak{t})$ =
 $( \frac{1}{20},
\frac{9}{20} )
$,

\vskip 0.7ex
\hangindent=5.5em \hangafter=1
{\white .}\hskip 1em $\rho_\text{isum}(\mathfrak{s})$ =
($-\frac{1}{\sqrt{5}}c^{3}_{20}
$,
$\frac{1}{\sqrt{5}}c^{1}_{20}
$;
$\frac{1}{\sqrt{5}}c^{3}_{20}
$)

Pass. 

 \ \color{black}

 \color{blue}

\noindent 8: (dims,levels) = $(2;20
)$,
irreps = $2_{5}^{2}
\hskip -1.5pt \otimes \hskip -1.5pt
1_{4}^{1,0}$,
pord$(\rho_\text{isum}(\mathfrak{t})) = 5$,

\vskip 0.7ex
\hangindent=5.5em \hangafter=1
{\white .}\hskip 1em $\rho_\text{isum}(\mathfrak{t})$ =
 $( \frac{13}{20},
\frac{17}{20} )
$,

\vskip 0.7ex
\hangindent=5.5em \hangafter=1
{\white .}\hskip 1em $\rho_\text{isum}(\mathfrak{s})$ =
($\frac{1}{\sqrt{5}}c^{1}_{20}
$,
$\frac{1}{\sqrt{5}}c^{3}_{20}
$;
$-\frac{1}{\sqrt{5}}c^{1}_{20}
$)

Pass. 

 \ \color{black}

 \color{blue}

\noindent 9: (dims,levels) = $(2;30
)$,
irreps = $2_{5}^{1}
\hskip -1.5pt \otimes \hskip -1.5pt
1_{3}^{1,0}
\hskip -1.5pt \otimes \hskip -1.5pt
1_{2}^{1,0}$,
pord$(\rho_\text{isum}(\mathfrak{t})) = 5$,

\vskip 0.7ex
\hangindent=5.5em \hangafter=1
{\white .}\hskip 1em $\rho_\text{isum}(\mathfrak{t})$ =
 $( \frac{1}{30},
\frac{19}{30} )
$,

\vskip 0.7ex
\hangindent=5.5em \hangafter=1
{\white .}\hskip 1em $\rho_\text{isum}(\mathfrak{s})$ =
$\mathrm{i}$($\frac{1}{\sqrt{5}}c^{3}_{20}
$,
$\frac{1}{\sqrt{5}}c^{1}_{20}
$;\ \ 
$-\frac{1}{\sqrt{5}}c^{3}_{20}
$)

Pass. 

 \ \color{black}

 \color{blue}

\noindent 10: (dims,levels) = $(2;30
)$,
irreps = $2_{5}^{2}
\hskip -1.5pt \otimes \hskip -1.5pt
1_{3}^{1,0}
\hskip -1.5pt \otimes \hskip -1.5pt
1_{2}^{1,0}$,
pord$(\rho_\text{isum}(\mathfrak{t})) = 5$,

\vskip 0.7ex
\hangindent=5.5em \hangafter=1
{\white .}\hskip 1em $\rho_\text{isum}(\mathfrak{t})$ =
 $( \frac{7}{30},
\frac{13}{30} )
$,

\vskip 0.7ex
\hangindent=5.5em \hangafter=1
{\white .}\hskip 1em $\rho_\text{isum}(\mathfrak{s})$ =
$\mathrm{i}$($\frac{1}{\sqrt{5}}c^{1}_{20}
$,
$\frac{1}{\sqrt{5}}c^{3}_{20}
$;\ \ 
$-\frac{1}{\sqrt{5}}c^{1}_{20}
$)

Pass. 

 \ \color{black}

 \color{blue}

\noindent 11: (dims,levels) = $(2;60
)$,
irreps = $2_{5}^{2}
\hskip -1.5pt \otimes \hskip -1.5pt
1_{4}^{1,0}
\hskip -1.5pt \otimes \hskip -1.5pt
1_{3}^{1,0}$,
pord$(\rho_\text{isum}(\mathfrak{t})) = 5$,

\vskip 0.7ex
\hangindent=5.5em \hangafter=1
{\white .}\hskip 1em $\rho_\text{isum}(\mathfrak{t})$ =
 $( \frac{11}{60},
\frac{59}{60} )
$,

\vskip 0.7ex
\hangindent=5.5em \hangafter=1
{\white .}\hskip 1em $\rho_\text{isum}(\mathfrak{s})$ =
($-\frac{1}{\sqrt{5}}c^{1}_{20}
$,
$\frac{1}{\sqrt{5}}c^{3}_{20}
$;
$\frac{1}{\sqrt{5}}c^{1}_{20}
$)

Pass. 

 \ \color{black}

 \color{blue}

\noindent 12: (dims,levels) = $(2;60
)$,
irreps = $2_{5}^{1}
\hskip -1.5pt \otimes \hskip -1.5pt
1_{4}^{1,0}
\hskip -1.5pt \otimes \hskip -1.5pt
1_{3}^{1,0}$,
pord$(\rho_\text{isum}(\mathfrak{t})) = 5$,

\vskip 0.7ex
\hangindent=5.5em \hangafter=1
{\white .}\hskip 1em $\rho_\text{isum}(\mathfrak{t})$ =
 $( \frac{23}{60},
\frac{47}{60} )
$,

\vskip 0.7ex
\hangindent=5.5em \hangafter=1
{\white .}\hskip 1em $\rho_\text{isum}(\mathfrak{s})$ =
($-\frac{1}{\sqrt{5}}c^{3}_{20}
$,
$\frac{1}{\sqrt{5}}c^{1}_{20}
$;
$\frac{1}{\sqrt{5}}c^{3}_{20}
$)

Pass. 

 \ \color{black}

\

\subsection{A list of passing GT orbits}

The above passing representations can be grouped into GT orbits.  The following
list displays one representative representation for each GT orbit.  For details
and notations, see Appendix B.2.

\

\noindent1. (dims;levels) =$(2;5
)$,
irreps = $2_{5}^{1}$,
pord$(\rho_\text{isum}(\mathfrak{t})) = 5$,

\vskip 0.7ex
\hangindent=4em \hangafter=1
 $\rho_\text{isum}(\mathfrak{t})$ =
 $( \frac{1}{5},
\frac{4}{5} )
$,

\vskip 0.7ex
\hangindent=4em \hangafter=1
 $\rho_\text{isum}(\mathfrak{s})$ =
$\mathrm{i}$($-\frac{1}{\sqrt{5}}c^{3}_{20}
$,
$\frac{1}{\sqrt{5}}c^{1}_{20}
$;\ \ 
$\frac{1}{\sqrt{5}}c^{3}_{20}
$)

Resolved. Number of valid $(S,T)$ pairs = 1.

\vskip 2ex

\

\subsection{A list of rank-2 $S,T$ matrices from resolved representations}

From the representative representation in each GT orbit, if it is revolved, we
can compute all the $S,T$ matrices coming from such a representation.  The
computation steps are displayed below.  For details and notations, see Section
1 of this file.

\

\noindent1. (dims;levels) =$(2;5
)$,
irreps = $2_{5}^{1}$,
pord$(\tilde\rho(\mathfrak{t})) = 5$,

\vskip 0.7ex
\hangindent=4em \hangafter=1
 $\tilde\rho(\mathfrak{t})$ =
 $( \frac{1}{5},
\frac{4}{5} )
$,

\vskip 0.7ex
\hangindent=4em \hangafter=1
 $\tilde\rho(\mathfrak{s})$ =
$\mathrm{i}$($-\frac{1}{\sqrt{5}}c^{3}_{20}
$,
$\frac{1}{\sqrt{5}}c^{1}_{20}
$;\ \ 
$\frac{1}{\sqrt{5}}c^{3}_{20}
$)

 \vskip 1ex \setlength{\leftskip}{2em}

\grey{Try $U_0$ =
$\begin{pmatrix}
1 \\ 
\end{pmatrix}
$ $\oplus
\begin{pmatrix}
1 \\ 
\end{pmatrix}
$:}\ \ \ \ \ 
\grey{$U_0\tilde\rho(\mathfrak{s})U_0^\dagger$ =}

\grey{$\begin{pmatrix}
-\frac{1}{\sqrt{5}}c^{3}_{20}
\mathrm{i},
& \frac{1}{\sqrt{5}}c^{1}_{20}
\mathrm{i} \\ 
\frac{1}{\sqrt{5}}c^{1}_{20}
\mathrm{i},
& \frac{1}{\sqrt{5}}c^{3}_{20}
\mathrm{i} \\ 
\end{pmatrix}
$}

\grey{Try different $u$'s and signed diagonal matrix $V_\mathrm{sd}$'s:}

 \grey{
\begin{tabular}{|r|l|l|}
\hline
$2_{5}^{1}:\ u$ 
 & 0 & 1\\ 
 \hline
$D_\rho$ conditions 
 & 0 & 0\\ 
 \hline
$[\rho(\mathfrak{s})\rho(\mathfrak{t})]^3
 = \rho^2(\mathfrak{s}) = \tilde C$ 
 & 0 & 0\\ 
 \hline
$\rho(\mathfrak{s})_{iu}\rho^*(\mathfrak{s})_{ju} \in \mathbb{R}$ 
 & 0 & 0\\ 
 \hline
$\rho(\mathfrak{s})_{i u} \neq 0$  
 & 0 & 0\\ 
 \hline
$\mathrm{cnd}(S)$, $\mathrm{cnd}(\rho(\mathfrak{s}))$ 
 & 0 & 0\\ 
 \hline
$\mathrm{norm}(D^2)$ factors
 & 0 & 0\\ 
 \hline
$1/\rho(\mathfrak{s})_{iu} = $ cyc-int 
 & 0 & 0\\ 
 \hline
norm$(1/\rho(\mathfrak{s})_{iu})$ factors
 & 0 & 0\\ 
 \hline
$\frac{S_{ij}}{S_{uj}} = $ cyc-int
 & 0 & 0\\ 
 \hline
$N^{ij}_k \in \mathbb{N}$
 & 0 & 0\\ 
 \hline
$\exists\ j \text{ that } \frac{S_{ij}}{S_{uj}} \geq 1 $
 & 0 & 0\\ 
 \hline
FS indicator
 & 0 & 0\\ 
 \hline
$C = $ perm-mat
 & 0 & 0\\ 
 \hline
\end{tabular}

Number of valid $(S,T)$ pairs: 1 \vskip 2ex }%grey

Total number of valid $(S,T)$ pairs: 1

 \vskip 4ex

\ \setlength{\leftskip}{0em}

\

The $S,T$ matrices obtained above (the black or the blue entries below), plus
their Galois conjugations (the grey entries below), form the following list of
rank-2 $S,T$ matrices.  For details and notations, see Section 2 of this file.

\

\noindent1. ind = $(2;5
)_{1}^{1}$:\ \ 
$d_i$ = ($1.0$,
$1.618$) 

\vskip 0.7ex
\hangindent=3em \hangafter=1
$D^2=$ 3.618 = 
 $\frac{5+\sqrt{5}}{2}$

\vskip 0.7ex
\hangindent=3em \hangafter=1
$T = ( 0,
\frac{2}{5} )
$,

\vskip 0.7ex
\hangindent=3em \hangafter=1
$S$ = ($ 1$,
$ \frac{1+\sqrt{5}}{2}$;\ \ 
$ -1$)

\vskip 1ex 
\color{grey}

\noindent2. ind = $(2;5
)_{1}^{4}$:\ \ 
$d_i$ = ($1.0$,
$1.618$) 

\vskip 0.7ex
\hangindent=3em \hangafter=1
$D^2=$ 3.618 = 
 $\frac{5+\sqrt{5}}{2}$

\vskip 0.7ex
\hangindent=3em \hangafter=1
$T = ( 0,
\frac{3}{5} )
$,

\vskip 0.7ex
\hangindent=3em \hangafter=1
$S$ = ($ 1$,
$ \frac{1+\sqrt{5}}{2}$;\ \ 
$ -1$)

\vskip 1ex 
\color{grey}

\noindent3. ind = $(2;5
)_{1}^{3}$:\ \ 
$d_i$ = ($1.0$,
$-0.618$) 

\vskip 0.7ex
\hangindent=3em \hangafter=1
$D^2=$ 1.381 = 
 $\frac{5-\sqrt{5}}{2}$

\vskip 0.7ex
\hangindent=3em \hangafter=1
$T = ( 0,
\frac{1}{5} )
$,

\vskip 0.7ex
\hangindent=3em \hangafter=1
$S$ = ($ 1$,
$ \frac{1-\sqrt{5}}{2}$;\ \ 
$ -1$)

Not pseudo-unitary. 

\vskip 1ex 
\color{grey}

\noindent4. ind = $(2;5
)_{1}^{2}$:\ \ 
$d_i$ = ($1.0$,
$-0.618$) 

\vskip 0.7ex
\hangindent=3em \hangafter=1
$D^2=$ 1.381 = 
 $\frac{5-\sqrt{5}}{2}$

\vskip 0.7ex
\hangindent=3em \hangafter=1
$T = ( 0,
\frac{4}{5} )
$,

\vskip 0.7ex
\hangindent=3em \hangafter=1
$S$ = ($ 1$,
$ \frac{1-\sqrt{5}}{2}$;\ \ 
$ -1$)

Not pseudo-unitary. 

\vskip 1ex 

 \color{black} \vskip 2ex

\

The above list includes all modular data (unitary or non-unitary) from resolved
$\SL$ representations and non-integral MTCs.
Since for rank 2, all the passing $\SL$ representations are resolved, the list
actually includes all modular data from non-integral MTCs.

\section{Rank-3 modular data}
\label{Section6}

\subsection{A list of 3-dimensional irrep-sum $\SL$ representations}

The following is a list of 3-dimensional irrep-sum $\SL$ representations, with
3 types of representations omitted and one type skipped.  For details and notations, see Section
3 of this file.

\

\noindent 1: (dims,levels) = $(3;5
)$,
irreps = $3_{5}^{1}$,
pord$(\rho_\text{isum}(\mathfrak{t})) = 5$,

\vskip 0.7ex
\hangindent=5.5em \hangafter=1
{\white .}\hskip 1em $\rho_\text{isum}(\mathfrak{t})$ =
 $( 0,
\frac{1}{5},
\frac{4}{5} )
$,

\vskip 0.7ex
\hangindent=5.5em \hangafter=1
{\white .}\hskip 1em $\rho_\text{isum}(\mathfrak{s})$ =
($\sqrt{\frac{1}{5}}$,
$-\sqrt{\frac{2}{5}}$,
$-\sqrt{\frac{2}{5}}$;
$-\frac{5+\sqrt{5}}{10}$,
$\frac{5-\sqrt{5}}{10}$;
$-\frac{5+\sqrt{5}}{10}$)

Fail:
cnd($\rho(\mathfrak s)_\mathrm{ndeg}$) = 40 does not divide
 ord($\rho(\mathfrak t)$)=5. Prop. B.4 (2)

 \ \color{black}

\noindent 2: (dims,levels) = $(3;5
)$,
irreps = $3_{5}^{3}$,
pord$(\rho_\text{isum}(\mathfrak{t})) = 5$,

\vskip 0.7ex
\hangindent=5.5em \hangafter=1
{\white .}\hskip 1em $\rho_\text{isum}(\mathfrak{t})$ =
 $( 0,
\frac{2}{5},
\frac{3}{5} )
$,

\vskip 0.7ex
\hangindent=5.5em \hangafter=1
{\white .}\hskip 1em $\rho_\text{isum}(\mathfrak{s})$ =
($-\sqrt{\frac{1}{5}}$,
$-\sqrt{\frac{2}{5}}$,
$-\sqrt{\frac{2}{5}}$;
$\frac{-5+\sqrt{5}}{10}$,
$\frac{5+\sqrt{5}}{10}$;
$\frac{-5+\sqrt{5}}{10}$)

Fail:
cnd($\rho(\mathfrak s)_\mathrm{ndeg}$) = 40 does not divide
 ord($\rho(\mathfrak t)$)=5. Prop. B.4 (2)

 \ \color{black}

 \color{blue}

\noindent 3: (dims,levels) = $(3;7
)$,
irreps = $3_{7}^{1}$,
pord$(\rho_\text{isum}(\mathfrak{t})) = 7$,

\vskip 0.7ex
\hangindent=5.5em \hangafter=1
{\white .}\hskip 1em $\rho_\text{isum}(\mathfrak{t})$ =
 $( \frac{1}{7},
\frac{2}{7},
\frac{4}{7} )
$,

\vskip 0.7ex
\hangindent=5.5em \hangafter=1
{\white .}\hskip 1em $\rho_\text{isum}(\mathfrak{s})$ =
($-\frac{1}{\sqrt{7}}c^{1}_{28}
$,
$-\frac{1}{\sqrt{7}}c^{3}_{28}
$,
$\frac{1}{\sqrt{7}}c^{5}_{28}
$;
$\frac{1}{\sqrt{7}}c^{5}_{28}
$,
$-\frac{1}{\sqrt{7}}c^{1}_{28}
$;
$-\frac{1}{\sqrt{7}}c^{3}_{28}
$)

Pass. 

 \ \color{black}

 \color{blue}

\noindent 4: (dims,levels) = $(3;7
)$,
irreps = $3_{7}^{3}$,
pord$(\rho_\text{isum}(\mathfrak{t})) = 7$,

\vskip 0.7ex
\hangindent=5.5em \hangafter=1
{\white .}\hskip 1em $\rho_\text{isum}(\mathfrak{t})$ =
 $( \frac{3}{7},
\frac{5}{7},
\frac{6}{7} )
$,

\vskip 0.7ex
\hangindent=5.5em \hangafter=1
{\white .}\hskip 1em $\rho_\text{isum}(\mathfrak{s})$ =
($-\frac{1}{\sqrt{7}}c^{3}_{28}
$,
$-\frac{1}{\sqrt{7}}c^{1}_{28}
$,
$\frac{1}{\sqrt{7}}c^{5}_{28}
$;
$\frac{1}{\sqrt{7}}c^{5}_{28}
$,
$-\frac{1}{\sqrt{7}}c^{3}_{28}
$;
$-\frac{1}{\sqrt{7}}c^{1}_{28}
$)

Pass. 

 \ \color{black}

\noindent 5: (dims,levels) = $(3;8
)$,
irreps = $3_{8}^{1,0}$,
pord$(\rho_\text{isum}(\mathfrak{t})) = 8$,

\vskip 0.7ex
\hangindent=5.5em \hangafter=1
{\white .}\hskip 1em $\rho_\text{isum}(\mathfrak{t})$ =
 $( 0,
\frac{1}{8},
\frac{5}{8} )
$,

\vskip 0.7ex
\hangindent=5.5em \hangafter=1
{\white .}\hskip 1em $\rho_\text{isum}(\mathfrak{s})$ =
$\mathrm{i}$($0$,
$\sqrt{\frac{1}{2}}$,
$\sqrt{\frac{1}{2}}$;\ \ 
$-\frac{1}{2}$,
$\frac{1}{2}$;\ \ 
$-\frac{1}{2}$)

Fail:
$\sigma(\rho(\mathfrak s)_\mathrm{ndeg}) \neq
 (\rho(\mathfrak t)^a \rho(\mathfrak s) \rho(\mathfrak t)^b
 \rho(\mathfrak s) \rho(\mathfrak t)^a)_\mathrm{ndeg}$,
 $\sigma = a$ = 3. Prop. B.5 (3) eqn. (B.25)

 \ \color{black}

\noindent 6: (dims,levels) = $(3;8
)$,
irreps = $3_{8}^{3,0}$,
pord$(\rho_\text{isum}(\mathfrak{t})) = 8$,

\vskip 0.7ex
\hangindent=5.5em \hangafter=1
{\white .}\hskip 1em $\rho_\text{isum}(\mathfrak{t})$ =
 $( 0,
\frac{3}{8},
\frac{7}{8} )
$,

\vskip 0.7ex
\hangindent=5.5em \hangafter=1
{\white .}\hskip 1em $\rho_\text{isum}(\mathfrak{s})$ =
$\mathrm{i}$($0$,
$\sqrt{\frac{1}{2}}$,
$\sqrt{\frac{1}{2}}$;\ \ 
$\frac{1}{2}$,
$-\frac{1}{2}$;\ \ 
$\frac{1}{2}$)

Fail:
$\sigma(\rho(\mathfrak s)_\mathrm{ndeg}) \neq
 (\rho(\mathfrak t)^a \rho(\mathfrak s) \rho(\mathfrak t)^b
 \rho(\mathfrak s) \rho(\mathfrak t)^a)_\mathrm{ndeg}$,
 $\sigma = a$ = 3. Prop. B.5 (3) eqn. (B.25)

 \ \color{black}

\noindent 7: (dims,levels) = $(3;10
)$,
irreps = $3_{5}^{3}
\hskip -1.5pt \otimes \hskip -1.5pt
1_{2}^{1,0}$,
pord$(\rho_\text{isum}(\mathfrak{t})) = 5$,

\vskip 0.7ex
\hangindent=5.5em \hangafter=1
{\white .}\hskip 1em $\rho_\text{isum}(\mathfrak{t})$ =
 $( \frac{1}{2},
\frac{1}{10},
\frac{9}{10} )
$,

\vskip 0.7ex
\hangindent=5.5em \hangafter=1
{\white .}\hskip 1em $\rho_\text{isum}(\mathfrak{s})$ =
($\sqrt{\frac{1}{5}}$,
$-\sqrt{\frac{2}{5}}$,
$-\sqrt{\frac{2}{5}}$;
$\frac{5-\sqrt{5}}{10}$,
$-\frac{5+\sqrt{5}}{10}$;
$\frac{5-\sqrt{5}}{10}$)

Fail:
cnd($\rho(\mathfrak s)_\mathrm{ndeg}$) = 40 does not divide
 ord($\rho(\mathfrak t)$)=10. Prop. B.4 (2)

 \ \color{black}

\noindent 8: (dims,levels) = $(3;10
)$,
irreps = $3_{5}^{1}
\hskip -1.5pt \otimes \hskip -1.5pt
1_{2}^{1,0}$,
pord$(\rho_\text{isum}(\mathfrak{t})) = 5$,

\vskip 0.7ex
\hangindent=5.5em \hangafter=1
{\white .}\hskip 1em $\rho_\text{isum}(\mathfrak{t})$ =
 $( \frac{1}{2},
\frac{3}{10},
\frac{7}{10} )
$,

\vskip 0.7ex
\hangindent=5.5em \hangafter=1
{\white .}\hskip 1em $\rho_\text{isum}(\mathfrak{s})$ =
($-\sqrt{\frac{1}{5}}$,
$-\sqrt{\frac{2}{5}}$,
$-\sqrt{\frac{2}{5}}$;
$\frac{5+\sqrt{5}}{10}$,
$\frac{-5+\sqrt{5}}{10}$;
$\frac{5+\sqrt{5}}{10}$)

Fail:
cnd($\rho(\mathfrak s)_\mathrm{ndeg}$) = 40 does not divide
 ord($\rho(\mathfrak t)$)=10. Prop. B.4 (2)

 \ \color{black}

 \color{blue}

\noindent 9: (dims,levels) = $(3;14
)$,
irreps = $3_{7}^{1}
\hskip -1.5pt \otimes \hskip -1.5pt
1_{2}^{1,0}$,
pord$(\rho_\text{isum}(\mathfrak{t})) = 7$,

\vskip 0.7ex
\hangindent=5.5em \hangafter=1
{\white .}\hskip 1em $\rho_\text{isum}(\mathfrak{t})$ =
 $( \frac{1}{14},
\frac{9}{14},
\frac{11}{14} )
$,

\vskip 0.7ex
\hangindent=5.5em \hangafter=1
{\white .}\hskip 1em $\rho_\text{isum}(\mathfrak{s})$ =
($\frac{1}{\sqrt{7}}c^{3}_{28}
$,
$\frac{1}{\sqrt{7}}c^{5}_{28}
$,
$-\frac{1}{\sqrt{7}}c^{1}_{28}
$;
$\frac{1}{\sqrt{7}}c^{1}_{28}
$,
$\frac{1}{\sqrt{7}}c^{3}_{28}
$;
$-\frac{1}{\sqrt{7}}c^{5}_{28}
$)

Pass. 

 \ \color{black}

 \color{blue}

\noindent 10: (dims,levels) = $(3;14
)$,
irreps = $3_{7}^{3}
\hskip -1.5pt \otimes \hskip -1.5pt
1_{2}^{1,0}$,
pord$(\rho_\text{isum}(\mathfrak{t})) = 7$,

\vskip 0.7ex
\hangindent=5.5em \hangafter=1
{\white .}\hskip 1em $\rho_\text{isum}(\mathfrak{t})$ =
 $( \frac{3}{14},
\frac{5}{14},
\frac{13}{14} )
$,

\vskip 0.7ex
\hangindent=5.5em \hangafter=1
{\white .}\hskip 1em $\rho_\text{isum}(\mathfrak{s})$ =
($-\frac{1}{\sqrt{7}}c^{5}_{28}
$,
$-\frac{1}{\sqrt{7}}c^{3}_{28}
$,
$-\frac{1}{\sqrt{7}}c^{1}_{28}
$;
$\frac{1}{\sqrt{7}}c^{1}_{28}
$,
$-\frac{1}{\sqrt{7}}c^{5}_{28}
$;
$\frac{1}{\sqrt{7}}c^{3}_{28}
$)

Pass. 

 \ \color{black}

\noindent 11: (dims,levels) = $(3;15
)$,
irreps = $3_{5}^{1}
\hskip -1.5pt \otimes \hskip -1.5pt
1_{3}^{1,0}$,
pord$(\rho_\text{isum}(\mathfrak{t})) = 5$,

\vskip 0.7ex
\hangindent=5.5em \hangafter=1
{\white .}\hskip 1em $\rho_\text{isum}(\mathfrak{t})$ =
 $( \frac{1}{3},
\frac{2}{15},
\frac{8}{15} )
$,

\vskip 0.7ex
\hangindent=5.5em \hangafter=1
{\white .}\hskip 1em $\rho_\text{isum}(\mathfrak{s})$ =
($\sqrt{\frac{1}{5}}$,
$-\sqrt{\frac{2}{5}}$,
$-\sqrt{\frac{2}{5}}$;
$-\frac{5+\sqrt{5}}{10}$,
$\frac{5-\sqrt{5}}{10}$;
$-\frac{5+\sqrt{5}}{10}$)

Fail:
cnd($\rho(\mathfrak s)_\mathrm{ndeg}$) = 40 does not divide
 ord($\rho(\mathfrak t)$)=15. Prop. B.4 (2)

 \ \color{black}

\noindent 12: (dims,levels) = $(3;15
)$,
irreps = $3_{5}^{3}
\hskip -1.5pt \otimes \hskip -1.5pt
1_{3}^{1,0}$,
pord$(\rho_\text{isum}(\mathfrak{t})) = 5$,

\vskip 0.7ex
\hangindent=5.5em \hangafter=1
{\white .}\hskip 1em $\rho_\text{isum}(\mathfrak{t})$ =
 $( \frac{1}{3},
\frac{11}{15},
\frac{14}{15} )
$,

\vskip 0.7ex
\hangindent=5.5em \hangafter=1
{\white .}\hskip 1em $\rho_\text{isum}(\mathfrak{s})$ =
($-\sqrt{\frac{1}{5}}$,
$-\sqrt{\frac{2}{5}}$,
$-\sqrt{\frac{2}{5}}$;
$\frac{-5+\sqrt{5}}{10}$,
$\frac{5+\sqrt{5}}{10}$;
$\frac{-5+\sqrt{5}}{10}$)

Fail:
cnd($\rho(\mathfrak s)_\mathrm{ndeg}$) = 40 does not divide
 ord($\rho(\mathfrak t)$)=15. Prop. B.4 (2)

 \ \color{black}

 \color{blue}

\noindent 13: (dims,levels) = $(3;16
)$,
irreps = $3_{16}^{1,0}$,
pord$(\rho_\text{isum}(\mathfrak{t})) = 16$,

\vskip 0.7ex
\hangindent=5.5em \hangafter=1
{\white .}\hskip 1em $\rho_\text{isum}(\mathfrak{t})$ =
 $( \frac{1}{8},
\frac{1}{16},
\frac{9}{16} )
$,

\vskip 0.7ex
\hangindent=5.5em \hangafter=1
{\white .}\hskip 1em $\rho_\text{isum}(\mathfrak{s})$ =
$\mathrm{i}$($0$,
$\sqrt{\frac{1}{2}}$,
$\sqrt{\frac{1}{2}}$;\ \ 
$-\frac{1}{2}$,
$\frac{1}{2}$;\ \ 
$-\frac{1}{2}$)

Pass. 

 \ \color{black}

 \color{blue}

\noindent 14: (dims,levels) = $(3;16
)$,
irreps = $3_{16}^{3,0}$,
pord$(\rho_\text{isum}(\mathfrak{t})) = 16$,

\vskip 0.7ex
\hangindent=5.5em \hangafter=1
{\white .}\hskip 1em $\rho_\text{isum}(\mathfrak{t})$ =
 $( \frac{3}{8},
\frac{3}{16},
\frac{11}{16} )
$,

\vskip 0.7ex
\hangindent=5.5em \hangafter=1
{\white .}\hskip 1em $\rho_\text{isum}(\mathfrak{s})$ =
$\mathrm{i}$($0$,
$\sqrt{\frac{1}{2}}$,
$\sqrt{\frac{1}{2}}$;\ \ 
$\frac{1}{2}$,
$-\frac{1}{2}$;\ \ 
$\frac{1}{2}$)

Pass. 

 \ \color{black}

 \color{blue}

\noindent 15: (dims,levels) = $(3;16
)$,
irreps = $3_{16}^{5,0}$,
pord$(\rho_\text{isum}(\mathfrak{t})) = 16$,

\vskip 0.7ex
\hangindent=5.5em \hangafter=1
{\white .}\hskip 1em $\rho_\text{isum}(\mathfrak{t})$ =
 $( \frac{5}{8},
\frac{5}{16},
\frac{13}{16} )
$,

\vskip 0.7ex
\hangindent=5.5em \hangafter=1
{\white .}\hskip 1em $\rho_\text{isum}(\mathfrak{s})$ =
$\mathrm{i}$($0$,
$\sqrt{\frac{1}{2}}$,
$\sqrt{\frac{1}{2}}$;\ \ 
$-\frac{1}{2}$,
$\frac{1}{2}$;\ \ 
$-\frac{1}{2}$)

Pass. 

 \ \color{black}

 \color{blue}

\noindent 16: (dims,levels) = $(3;16
)$,
irreps = $3_{16}^{7,0}$,
pord$(\rho_\text{isum}(\mathfrak{t})) = 16$,

\vskip 0.7ex
\hangindent=5.5em \hangafter=1
{\white .}\hskip 1em $\rho_\text{isum}(\mathfrak{t})$ =
 $( \frac{7}{8},
\frac{7}{16},
\frac{15}{16} )
$,

\vskip 0.7ex
\hangindent=5.5em \hangafter=1
{\white .}\hskip 1em $\rho_\text{isum}(\mathfrak{s})$ =
$\mathrm{i}$($0$,
$\sqrt{\frac{1}{2}}$,
$\sqrt{\frac{1}{2}}$;\ \ 
$\frac{1}{2}$,
$-\frac{1}{2}$;\ \ 
$\frac{1}{2}$)

Pass. 

 \ \color{black}

\noindent 17: (dims,levels) = $(3;20
)$,
irreps = $3_{5}^{1}
\hskip -1.5pt \otimes \hskip -1.5pt
1_{4}^{1,0}$,
pord$(\rho_\text{isum}(\mathfrak{t})) = 5$,

\vskip 0.7ex
\hangindent=5.5em \hangafter=1
{\white .}\hskip 1em $\rho_\text{isum}(\mathfrak{t})$ =
 $( \frac{1}{4},
\frac{1}{20},
\frac{9}{20} )
$,

\vskip 0.7ex
\hangindent=5.5em \hangafter=1
{\white .}\hskip 1em $\rho_\text{isum}(\mathfrak{s})$ =
$\mathrm{i}$($\sqrt{\frac{1}{5}}$,
$\sqrt{\frac{2}{5}}$,
$\sqrt{\frac{2}{5}}$;\ \ 
$-\frac{5+\sqrt{5}}{10}$,
$\frac{5-\sqrt{5}}{10}$;\ \ 
$-\frac{5+\sqrt{5}}{10}$)

Fail:
cnd($\rho(\mathfrak s)_\mathrm{ndeg}$) = 40 does not divide
 ord($\rho(\mathfrak t)$)=20. Prop. B.4 (2)

 \ \color{black}

\noindent 18: (dims,levels) = $(3;20
)$,
irreps = $3_{5}^{3}
\hskip -1.5pt \otimes \hskip -1.5pt
1_{4}^{1,0}$,
pord$(\rho_\text{isum}(\mathfrak{t})) = 5$,

\vskip 0.7ex
\hangindent=5.5em \hangafter=1
{\white .}\hskip 1em $\rho_\text{isum}(\mathfrak{t})$ =
 $( \frac{1}{4},
\frac{13}{20},
\frac{17}{20} )
$,

\vskip 0.7ex
\hangindent=5.5em \hangafter=1
{\white .}\hskip 1em $\rho_\text{isum}(\mathfrak{s})$ =
$\mathrm{i}$($-\sqrt{\frac{1}{5}}$,
$\sqrt{\frac{2}{5}}$,
$\sqrt{\frac{2}{5}}$;\ \ 
$\frac{-5+\sqrt{5}}{10}$,
$\frac{5+\sqrt{5}}{10}$;\ \ 
$\frac{-5+\sqrt{5}}{10}$)

Fail:
cnd($\rho(\mathfrak s)_\mathrm{ndeg}$) = 40 does not divide
 ord($\rho(\mathfrak t)$)=20. Prop. B.4 (2)

 \ \color{black}

 \color{blue}

\noindent 19: (dims,levels) = $(3;21
)$,
irreps = $3_{7}^{3}
\hskip -1.5pt \otimes \hskip -1.5pt
1_{3}^{1,0}$,
pord$(\rho_\text{isum}(\mathfrak{t})) = 7$,

\vskip 0.7ex
\hangindent=5.5em \hangafter=1
{\white .}\hskip 1em $\rho_\text{isum}(\mathfrak{t})$ =
 $( \frac{1}{21},
\frac{4}{21},
\frac{16}{21} )
$,

\vskip 0.7ex
\hangindent=5.5em \hangafter=1
{\white .}\hskip 1em $\rho_\text{isum}(\mathfrak{s})$ =
($\frac{1}{\sqrt{7}}c^{5}_{28}
$,
$-\frac{1}{\sqrt{7}}c^{3}_{28}
$,
$-\frac{1}{\sqrt{7}}c^{1}_{28}
$;
$-\frac{1}{\sqrt{7}}c^{1}_{28}
$,
$\frac{1}{\sqrt{7}}c^{5}_{28}
$;
$-\frac{1}{\sqrt{7}}c^{3}_{28}
$)

Pass. 

 \ \color{black}

 \color{blue}

\noindent 20: (dims,levels) = $(3;21
)$,
irreps = $3_{7}^{1}
\hskip -1.5pt \otimes \hskip -1.5pt
1_{3}^{1,0}$,
pord$(\rho_\text{isum}(\mathfrak{t})) = 7$,

\vskip 0.7ex
\hangindent=5.5em \hangafter=1
{\white .}\hskip 1em $\rho_\text{isum}(\mathfrak{t})$ =
 $( \frac{10}{21},
\frac{13}{21},
\frac{19}{21} )
$,

\vskip 0.7ex
\hangindent=5.5em \hangafter=1
{\white .}\hskip 1em $\rho_\text{isum}(\mathfrak{s})$ =
($-\frac{1}{\sqrt{7}}c^{1}_{28}
$,
$-\frac{1}{\sqrt{7}}c^{3}_{28}
$,
$\frac{1}{\sqrt{7}}c^{5}_{28}
$;
$\frac{1}{\sqrt{7}}c^{5}_{28}
$,
$-\frac{1}{\sqrt{7}}c^{1}_{28}
$;
$-\frac{1}{\sqrt{7}}c^{3}_{28}
$)

Pass. 

 \ \color{black}

\noindent 21: (dims,levels) = $(3;24
)$,
irreps = $3_{8}^{3,0}
\hskip -1.5pt \otimes \hskip -1.5pt
1_{3}^{1,0}$,
pord$(\rho_\text{isum}(\mathfrak{t})) = 8$,

\vskip 0.7ex
\hangindent=5.5em \hangafter=1
{\white .}\hskip 1em $\rho_\text{isum}(\mathfrak{t})$ =
 $( \frac{1}{3},
\frac{5}{24},
\frac{17}{24} )
$,

\vskip 0.7ex
\hangindent=5.5em \hangafter=1
{\white .}\hskip 1em $\rho_\text{isum}(\mathfrak{s})$ =
$\mathrm{i}$($0$,
$\sqrt{\frac{1}{2}}$,
$\sqrt{\frac{1}{2}}$;\ \ 
$\frac{1}{2}$,
$-\frac{1}{2}$;\ \ 
$\frac{1}{2}$)

Fail:
$\sigma(\rho(\mathfrak s)_\mathrm{ndeg}) \neq
 (\rho(\mathfrak t)^a \rho(\mathfrak s) \rho(\mathfrak t)^b
 \rho(\mathfrak s) \rho(\mathfrak t)^a)_\mathrm{ndeg}$,
 $\sigma = a$ = 5. Prop. B.5 (3) eqn. (B.25)

 \ \color{black}

\noindent 22: (dims,levels) = $(3;24
)$,
irreps = $3_{8}^{1,0}
\hskip -1.5pt \otimes \hskip -1.5pt
1_{3}^{1,0}$,
pord$(\rho_\text{isum}(\mathfrak{t})) = 8$,

\vskip 0.7ex
\hangindent=5.5em \hangafter=1
{\white .}\hskip 1em $\rho_\text{isum}(\mathfrak{t})$ =
 $( \frac{1}{3},
\frac{11}{24},
\frac{23}{24} )
$,

\vskip 0.7ex
\hangindent=5.5em \hangafter=1
{\white .}\hskip 1em $\rho_\text{isum}(\mathfrak{s})$ =
$\mathrm{i}$($0$,
$\sqrt{\frac{1}{2}}$,
$\sqrt{\frac{1}{2}}$;\ \ 
$-\frac{1}{2}$,
$\frac{1}{2}$;\ \ 
$-\frac{1}{2}$)

Fail:
$\sigma(\rho(\mathfrak s)_\mathrm{ndeg}) \neq
 (\rho(\mathfrak t)^a \rho(\mathfrak s) \rho(\mathfrak t)^b
 \rho(\mathfrak s) \rho(\mathfrak t)^a)_\mathrm{ndeg}$,
 $\sigma = a$ = 5. Prop. B.5 (3) eqn. (B.25)

 \ \color{black}

 \color{blue}

\noindent 23: (dims,levels) = $(3;28
)$,
irreps = $3_{7}^{3}
\hskip -1.5pt \otimes \hskip -1.5pt
1_{4}^{1,0}$,
pord$(\rho_\text{isum}(\mathfrak{t})) = 7$,

\vskip 0.7ex
\hangindent=5.5em \hangafter=1
{\white .}\hskip 1em $\rho_\text{isum}(\mathfrak{t})$ =
 $( \frac{3}{28},
\frac{19}{28},
\frac{27}{28} )
$,

\vskip 0.7ex
\hangindent=5.5em \hangafter=1
{\white .}\hskip 1em $\rho_\text{isum}(\mathfrak{s})$ =
$\mathrm{i}$($-\frac{1}{\sqrt{7}}c^{1}_{28}
$,
$\frac{1}{\sqrt{7}}c^{5}_{28}
$,
$-\frac{1}{\sqrt{7}}c^{3}_{28}
$;\ \ 
$-\frac{1}{\sqrt{7}}c^{3}_{28}
$,
$-\frac{1}{\sqrt{7}}c^{1}_{28}
$;\ \ 
$\frac{1}{\sqrt{7}}c^{5}_{28}
$)

Pass. 

 \ \color{black}

 \color{blue}

\noindent 24: (dims,levels) = $(3;28
)$,
irreps = $3_{7}^{1}
\hskip -1.5pt \otimes \hskip -1.5pt
1_{4}^{1,0}$,
pord$(\rho_\text{isum}(\mathfrak{t})) = 7$,

\vskip 0.7ex
\hangindent=5.5em \hangafter=1
{\white .}\hskip 1em $\rho_\text{isum}(\mathfrak{t})$ =
 $( \frac{11}{28},
\frac{15}{28},
\frac{23}{28} )
$,

\vskip 0.7ex
\hangindent=5.5em \hangafter=1
{\white .}\hskip 1em $\rho_\text{isum}(\mathfrak{s})$ =
$\mathrm{i}$($-\frac{1}{\sqrt{7}}c^{1}_{28}
$,
$-\frac{1}{\sqrt{7}}c^{3}_{28}
$,
$\frac{1}{\sqrt{7}}c^{5}_{28}
$;\ \ 
$\frac{1}{\sqrt{7}}c^{5}_{28}
$,
$-\frac{1}{\sqrt{7}}c^{1}_{28}
$;\ \ 
$-\frac{1}{\sqrt{7}}c^{3}_{28}
$)

Pass. 

 \ \color{black}

\noindent 25: (dims,levels) = $(3;30
)$,
irreps = $3_{5}^{1}
\hskip -1.5pt \otimes \hskip -1.5pt
1_{3}^{1,0}
\hskip -1.5pt \otimes \hskip -1.5pt
1_{2}^{1,0}$,
pord$(\rho_\text{isum}(\mathfrak{t})) = 5$,

\vskip 0.7ex
\hangindent=5.5em \hangafter=1
{\white .}\hskip 1em $\rho_\text{isum}(\mathfrak{t})$ =
 $( \frac{5}{6},
\frac{1}{30},
\frac{19}{30} )
$,

\vskip 0.7ex
\hangindent=5.5em \hangafter=1
{\white .}\hskip 1em $\rho_\text{isum}(\mathfrak{s})$ =
($-\sqrt{\frac{1}{5}}$,
$-\sqrt{\frac{2}{5}}$,
$-\sqrt{\frac{2}{5}}$;
$\frac{5+\sqrt{5}}{10}$,
$\frac{-5+\sqrt{5}}{10}$;
$\frac{5+\sqrt{5}}{10}$)

Fail:
cnd($\rho(\mathfrak s)_\mathrm{ndeg}$) = 40 does not divide
 ord($\rho(\mathfrak t)$)=30. Prop. B.4 (2)

 \ \color{black}

\noindent 26: (dims,levels) = $(3;30
)$,
irreps = $3_{5}^{3}
\hskip -1.5pt \otimes \hskip -1.5pt
1_{3}^{1,0}
\hskip -1.5pt \otimes \hskip -1.5pt
1_{2}^{1,0}$,
pord$(\rho_\text{isum}(\mathfrak{t})) = 5$,

\vskip 0.7ex
\hangindent=5.5em \hangafter=1
{\white .}\hskip 1em $\rho_\text{isum}(\mathfrak{t})$ =
 $( \frac{5}{6},
\frac{7}{30},
\frac{13}{30} )
$,

\vskip 0.7ex
\hangindent=5.5em \hangafter=1
{\white .}\hskip 1em $\rho_\text{isum}(\mathfrak{s})$ =
($\sqrt{\frac{1}{5}}$,
$-\sqrt{\frac{2}{5}}$,
$-\sqrt{\frac{2}{5}}$;
$\frac{5-\sqrt{5}}{10}$,
$-\frac{5+\sqrt{5}}{10}$;
$\frac{5-\sqrt{5}}{10}$)

Fail:
cnd($\rho(\mathfrak s)_\mathrm{ndeg}$) = 40 does not divide
 ord($\rho(\mathfrak t)$)=30. Prop. B.4 (2)

 \ \color{black}

 \color{blue}

\noindent 27: (dims,levels) = $(3;42
)$,
irreps = $3_{7}^{1}
\hskip -1.5pt \otimes \hskip -1.5pt
1_{3}^{1,0}
\hskip -1.5pt \otimes \hskip -1.5pt
1_{2}^{1,0}$,
pord$(\rho_\text{isum}(\mathfrak{t})) = 7$,

\vskip 0.7ex
\hangindent=5.5em \hangafter=1
{\white .}\hskip 1em $\rho_\text{isum}(\mathfrak{t})$ =
 $( \frac{5}{42},
\frac{17}{42},
\frac{41}{42} )
$,

\vskip 0.7ex
\hangindent=5.5em \hangafter=1
{\white .}\hskip 1em $\rho_\text{isum}(\mathfrak{s})$ =
($-\frac{1}{\sqrt{7}}c^{5}_{28}
$,
$-\frac{1}{\sqrt{7}}c^{1}_{28}
$,
$-\frac{1}{\sqrt{7}}c^{3}_{28}
$;
$\frac{1}{\sqrt{7}}c^{3}_{28}
$,
$-\frac{1}{\sqrt{7}}c^{5}_{28}
$;
$\frac{1}{\sqrt{7}}c^{1}_{28}
$)

Pass. 

 \ \color{black}

 \color{blue}

\noindent 28: (dims,levels) = $(3;42
)$,
irreps = $3_{7}^{3}
\hskip -1.5pt \otimes \hskip -1.5pt
1_{3}^{1,0}
\hskip -1.5pt \otimes \hskip -1.5pt
1_{2}^{1,0}$,
pord$(\rho_\text{isum}(\mathfrak{t})) = 7$,

\vskip 0.7ex
\hangindent=5.5em \hangafter=1
{\white .}\hskip 1em $\rho_\text{isum}(\mathfrak{t})$ =
 $( \frac{11}{42},
\frac{23}{42},
\frac{29}{42} )
$,

\vskip 0.7ex
\hangindent=5.5em \hangafter=1
{\white .}\hskip 1em $\rho_\text{isum}(\mathfrak{s})$ =
($\frac{1}{\sqrt{7}}c^{3}_{28}
$,
$-\frac{1}{\sqrt{7}}c^{1}_{28}
$,
$\frac{1}{\sqrt{7}}c^{5}_{28}
$;
$-\frac{1}{\sqrt{7}}c^{5}_{28}
$,
$\frac{1}{\sqrt{7}}c^{3}_{28}
$;
$\frac{1}{\sqrt{7}}c^{1}_{28}
$)

Pass. 

 \ \color{black}

 \color{blue}

\noindent 29: (dims,levels) = $(3;48
)$,
irreps = $3_{16}^{7,0}
\hskip -1.5pt \otimes \hskip -1.5pt
1_{3}^{1,0}$,
pord$(\rho_\text{isum}(\mathfrak{t})) = 16$,

\vskip 0.7ex
\hangindent=5.5em \hangafter=1
{\white .}\hskip 1em $\rho_\text{isum}(\mathfrak{t})$ =
 $( \frac{5}{24},
\frac{13}{48},
\frac{37}{48} )
$,

\vskip 0.7ex
\hangindent=5.5em \hangafter=1
{\white .}\hskip 1em $\rho_\text{isum}(\mathfrak{s})$ =
$\mathrm{i}$($0$,
$\sqrt{\frac{1}{2}}$,
$\sqrt{\frac{1}{2}}$;\ \ 
$\frac{1}{2}$,
$-\frac{1}{2}$;\ \ 
$\frac{1}{2}$)

Pass. 

 \ \color{black}

 \color{blue}

\noindent 30: (dims,levels) = $(3;48
)$,
irreps = $3_{16}^{1,0}
\hskip -1.5pt \otimes \hskip -1.5pt
1_{3}^{1,0}$,
pord$(\rho_\text{isum}(\mathfrak{t})) = 16$,

\vskip 0.7ex
\hangindent=5.5em \hangafter=1
{\white .}\hskip 1em $\rho_\text{isum}(\mathfrak{t})$ =
 $( \frac{11}{24},
\frac{19}{48},
\frac{43}{48} )
$,

\vskip 0.7ex
\hangindent=5.5em \hangafter=1
{\white .}\hskip 1em $\rho_\text{isum}(\mathfrak{s})$ =
$\mathrm{i}$($0$,
$\sqrt{\frac{1}{2}}$,
$\sqrt{\frac{1}{2}}$;\ \ 
$-\frac{1}{2}$,
$\frac{1}{2}$;\ \ 
$-\frac{1}{2}$)

Pass. 

 \ \color{black}

 \color{blue}

\noindent 31: (dims,levels) = $(3;48
)$,
irreps = $3_{16}^{3,0}
\hskip -1.5pt \otimes \hskip -1.5pt
1_{3}^{1,0}$,
pord$(\rho_\text{isum}(\mathfrak{t})) = 16$,

\vskip 0.7ex
\hangindent=5.5em \hangafter=1
{\white .}\hskip 1em $\rho_\text{isum}(\mathfrak{t})$ =
 $( \frac{17}{24},
\frac{1}{48},
\frac{25}{48} )
$,

\vskip 0.7ex
\hangindent=5.5em \hangafter=1
{\white .}\hskip 1em $\rho_\text{isum}(\mathfrak{s})$ =
$\mathrm{i}$($0$,
$\sqrt{\frac{1}{2}}$,
$\sqrt{\frac{1}{2}}$;\ \ 
$\frac{1}{2}$,
$-\frac{1}{2}$;\ \ 
$\frac{1}{2}$)

Pass. 

 \ \color{black}

 \color{blue}

\noindent 32: (dims,levels) = $(3;48
)$,
irreps = $3_{16}^{5,0}
\hskip -1.5pt \otimes \hskip -1.5pt
1_{3}^{1,0}$,
pord$(\rho_\text{isum}(\mathfrak{t})) = 16$,

\vskip 0.7ex
\hangindent=5.5em \hangafter=1
{\white .}\hskip 1em $\rho_\text{isum}(\mathfrak{t})$ =
 $( \frac{23}{24},
\frac{7}{48},
\frac{31}{48} )
$,

\vskip 0.7ex
\hangindent=5.5em \hangafter=1
{\white .}\hskip 1em $\rho_\text{isum}(\mathfrak{s})$ =
$\mathrm{i}$($0$,
$\sqrt{\frac{1}{2}}$,
$\sqrt{\frac{1}{2}}$;\ \ 
$-\frac{1}{2}$,
$\frac{1}{2}$;\ \ 
$-\frac{1}{2}$)

Pass. 

 \ \color{black}

\noindent 33: (dims,levels) = $(3;60
)$,
irreps = $3_{5}^{3}
\hskip -1.5pt \otimes \hskip -1.5pt
1_{4}^{1,0}
\hskip -1.5pt \otimes \hskip -1.5pt
1_{3}^{1,0}$,
pord$(\rho_\text{isum}(\mathfrak{t})) = 5$,

\vskip 0.7ex
\hangindent=5.5em \hangafter=1
{\white .}\hskip 1em $\rho_\text{isum}(\mathfrak{t})$ =
 $( \frac{7}{12},
\frac{11}{60},
\frac{59}{60} )
$,

\vskip 0.7ex
\hangindent=5.5em \hangafter=1
{\white .}\hskip 1em $\rho_\text{isum}(\mathfrak{s})$ =
$\mathrm{i}$($-\sqrt{\frac{1}{5}}$,
$\sqrt{\frac{2}{5}}$,
$\sqrt{\frac{2}{5}}$;\ \ 
$\frac{-5+\sqrt{5}}{10}$,
$\frac{5+\sqrt{5}}{10}$;\ \ 
$\frac{-5+\sqrt{5}}{10}$)

Fail:
cnd($\rho(\mathfrak s)_\mathrm{ndeg}$) = 40 does not divide
 ord($\rho(\mathfrak t)$)=60. Prop. B.4 (2)

 \ \color{black}

\noindent 34: (dims,levels) = $(3;60
)$,
irreps = $3_{5}^{1}
\hskip -1.5pt \otimes \hskip -1.5pt
1_{4}^{1,0}
\hskip -1.5pt \otimes \hskip -1.5pt
1_{3}^{1,0}$,
pord$(\rho_\text{isum}(\mathfrak{t})) = 5$,

\vskip 0.7ex
\hangindent=5.5em \hangafter=1
{\white .}\hskip 1em $\rho_\text{isum}(\mathfrak{t})$ =
 $( \frac{7}{12},
\frac{23}{60},
\frac{47}{60} )
$,

\vskip 0.7ex
\hangindent=5.5em \hangafter=1
{\white .}\hskip 1em $\rho_\text{isum}(\mathfrak{s})$ =
$\mathrm{i}$($\sqrt{\frac{1}{5}}$,
$\sqrt{\frac{2}{5}}$,
$\sqrt{\frac{2}{5}}$;\ \ 
$-\frac{5+\sqrt{5}}{10}$,
$\frac{5-\sqrt{5}}{10}$;\ \ 
$-\frac{5+\sqrt{5}}{10}$)

Fail:
cnd($\rho(\mathfrak s)_\mathrm{ndeg}$) = 40 does not divide
 ord($\rho(\mathfrak t)$)=60. Prop. B.4 (2)

 \ \color{black}

 \color{blue}

\noindent 35: (dims,levels) = $(3;84
)$,
irreps = $3_{7}^{3}
\hskip -1.5pt \otimes \hskip -1.5pt
1_{4}^{1,0}
\hskip -1.5pt \otimes \hskip -1.5pt
1_{3}^{1,0}$,
pord$(\rho_\text{isum}(\mathfrak{t})) = 7$,

\vskip 0.7ex
\hangindent=5.5em \hangafter=1
{\white .}\hskip 1em $\rho_\text{isum}(\mathfrak{t})$ =
 $( \frac{1}{84},
\frac{25}{84},
\frac{37}{84} )
$,

\vskip 0.7ex
\hangindent=5.5em \hangafter=1
{\white .}\hskip 1em $\rho_\text{isum}(\mathfrak{s})$ =
$\mathrm{i}$($-\frac{1}{\sqrt{7}}c^{3}_{28}
$,
$-\frac{1}{\sqrt{7}}c^{1}_{28}
$,
$\frac{1}{\sqrt{7}}c^{5}_{28}
$;\ \ 
$\frac{1}{\sqrt{7}}c^{5}_{28}
$,
$-\frac{1}{\sqrt{7}}c^{3}_{28}
$;\ \ 
$-\frac{1}{\sqrt{7}}c^{1}_{28}
$)

Pass. 

 \ \color{black}

 \color{blue}

\noindent 36: (dims,levels) = $(3;84
)$,
irreps = $3_{7}^{1}
\hskip -1.5pt \otimes \hskip -1.5pt
1_{4}^{1,0}
\hskip -1.5pt \otimes \hskip -1.5pt
1_{3}^{1,0}$,
pord$(\rho_\text{isum}(\mathfrak{t})) = 7$,

\vskip 0.7ex
\hangindent=5.5em \hangafter=1
{\white .}\hskip 1em $\rho_\text{isum}(\mathfrak{t})$ =
 $( \frac{13}{84},
\frac{61}{84},
\frac{73}{84} )
$,

\vskip 0.7ex
\hangindent=5.5em \hangafter=1
{\white .}\hskip 1em $\rho_\text{isum}(\mathfrak{s})$ =
$\mathrm{i}$($-\frac{1}{\sqrt{7}}c^{3}_{28}
$,
$\frac{1}{\sqrt{7}}c^{5}_{28}
$,
$-\frac{1}{\sqrt{7}}c^{1}_{28}
$;\ \ 
$-\frac{1}{\sqrt{7}}c^{1}_{28}
$,
$-\frac{1}{\sqrt{7}}c^{3}_{28}
$;\ \ 
$\frac{1}{\sqrt{7}}c^{5}_{28}
$)

Pass. 

 \ \color{black}

\

\subsection{A list of passing GT orbits}

The above passing representations can be grouped into GT orbits.  The following
list displays one representative representation for each GT orbit.  For details
and notations, see Appendix B.2.

\

\noindent1. (dims;levels) =$(3;7
)$,
irreps = $3_{7}^{1}$,
pord$(\rho_\text{isum}(\mathfrak{t})) = 7$,

\vskip 0.7ex
\hangindent=4em \hangafter=1
 $\rho_\text{isum}(\mathfrak{t})$ =
 $( \frac{1}{7},
\frac{2}{7},
\frac{4}{7} )
$,

\vskip 0.7ex
\hangindent=4em \hangafter=1
 $\rho_\text{isum}(\mathfrak{s})$ =
($-\frac{1}{\sqrt{7}}c^{1}_{28}
$,
$-\frac{1}{\sqrt{7}}c^{3}_{28}
$,
$\frac{1}{\sqrt{7}}c^{5}_{28}
$;
$\frac{1}{\sqrt{7}}c^{5}_{28}
$,
$-\frac{1}{\sqrt{7}}c^{1}_{28}
$;
$-\frac{1}{\sqrt{7}}c^{3}_{28}
$)

Resolved. Number of valid $(S,T)$ pairs = 1.

\vskip 2ex

 \noindent2. (dims;levels) =$(3;16
)$,
irreps = $3_{16}^{1,0}$,
pord$(\rho_\text{isum}(\mathfrak{t})) = 16$,

\vskip 0.7ex
\hangindent=4em \hangafter=1
 $\rho_\text{isum}(\mathfrak{t})$ =
 $( \frac{1}{8},
\frac{1}{16},
\frac{9}{16} )
$,

\vskip 0.7ex
\hangindent=4em \hangafter=1
 $\rho_\text{isum}(\mathfrak{s})$ =
$\mathrm{i}$($0$,
$\sqrt{\frac{1}{2}}$,
$\sqrt{\frac{1}{2}}$;\ \ 
$-\frac{1}{2}$,
$\frac{1}{2}$;\ \ 
$-\frac{1}{2}$)

Resolved. Number of valid $(S,T)$ pairs = 2.

\vskip 2ex

\

\subsection{A list of rank-3 $S,T$ matrices from resolved representations}

From the representative representation in each GT orbit, if it is revolved, we
can compute all the $S,T$ matrices coming from such a representation.  The
computation steps are displayed below.  For details and notations, see Section
1 of this file.

\

\noindent1. (dims;levels) =$(3;7
)$,
irreps = $3_{7}^{1}$,
pord$(\tilde\rho(\mathfrak{t})) = 7$,

\vskip 0.7ex
\hangindent=4em \hangafter=1
 $\tilde\rho(\mathfrak{t})$ =
 $( \frac{1}{7},
\frac{2}{7},
\frac{4}{7} )
$,

\vskip 0.7ex
\hangindent=4em \hangafter=1
 $\tilde\rho(\mathfrak{s})$ =
($-\frac{1}{\sqrt{7}}c^{1}_{28}
$,
$-\frac{1}{\sqrt{7}}c^{3}_{28}
$,
$\frac{1}{\sqrt{7}}c^{5}_{28}
$;
$\frac{1}{\sqrt{7}}c^{5}_{28}
$,
$-\frac{1}{\sqrt{7}}c^{1}_{28}
$;
$-\frac{1}{\sqrt{7}}c^{3}_{28}
$)

 \vskip 1ex \setlength{\leftskip}{2em}

\grey{Try $U_0$ =
$\begin{pmatrix}
1 \\ 
\end{pmatrix}
$ $\oplus
\begin{pmatrix}
1 \\ 
\end{pmatrix}
$ $\oplus
\begin{pmatrix}
1 \\ 
\end{pmatrix}
$:}\ \ \ \ \ 
\grey{$U_0\tilde\rho(\mathfrak{s})U_0^\dagger$ =}

\grey{$\begin{pmatrix}
-\frac{1}{\sqrt{7}}c^{1}_{28}
,
& -\frac{1}{\sqrt{7}}c^{3}_{28}
,
& \frac{1}{\sqrt{7}}c^{5}_{28}
 \\ 
-\frac{1}{\sqrt{7}}c^{3}_{28}
,
& \frac{1}{\sqrt{7}}c^{5}_{28}
,
& -\frac{1}{\sqrt{7}}c^{1}_{28}
 \\ 
\frac{1}{\sqrt{7}}c^{5}_{28}
,
& -\frac{1}{\sqrt{7}}c^{1}_{28}
,
& -\frac{1}{\sqrt{7}}c^{3}_{28}
 \\ 
\end{pmatrix}
$}

\grey{Try different $u$'s and signed diagonal matrix $V_\mathrm{sd}$'s:}

 \grey{
\begin{tabular}{|r|l|l|l|}
\hline
$3_{7}^{1}:\ u$ 
 & 0 & 1 & 2\\ 
 \hline
$D_\rho$ conditions 
 & 0 & 0 & 0\\ 
 \hline
$[\rho(\mathfrak{s})\rho(\mathfrak{t})]^3
 = \rho^2(\mathfrak{s}) = \tilde C$ 
 & 0 & 0 & 0\\ 
 \hline
$\rho(\mathfrak{s})_{iu}\rho^*(\mathfrak{s})_{ju} \in \mathbb{R}$ 
 & 0 & 0 & 0\\ 
 \hline
$\rho(\mathfrak{s})_{i u} \neq 0$  
 & 0 & 0 & 0\\ 
 \hline
$\mathrm{cnd}(S)$, $\mathrm{cnd}(\rho(\mathfrak{s}))$ 
 & 0 & 0 & 0\\ 
 \hline
$\mathrm{norm}(D^2)$ factors
 & 0 & 0 & 0\\ 
 \hline
$1/\rho(\mathfrak{s})_{iu} = $ cyc-int 
 & 0 & 0 & 0\\ 
 \hline
norm$(1/\rho(\mathfrak{s})_{iu})$ factors
 & 0 & 0 & 0\\ 
 \hline
$\frac{S_{ij}}{S_{uj}} = $ cyc-int
 & 0 & 0 & 0\\ 
 \hline
$N^{ij}_k \in \mathbb{N}$
 & 0 & 0 & 0\\ 
 \hline
$\exists\ j \text{ that } \frac{S_{ij}}{S_{uj}} \geq 1 $
 & 0 & 0 & 0\\ 
 \hline
FS indicator
 & 0 & 0 & 0\\ 
 \hline
$C = $ perm-mat
 & 0 & 0 & 0\\ 
 \hline
\end{tabular}

Number of valid $(S,T)$ pairs: 1 \vskip 2ex }%grey

Total number of valid $(S,T)$ pairs: 1

 \vskip 4ex

\ \setlength{\leftskip}{0em} 

\noindent2. (dims;levels) =$(3;16
)$,
irreps = $3_{16}^{1,0}$,
pord$(\tilde\rho(\mathfrak{t})) = 16$,

\vskip 0.7ex
\hangindent=4em \hangafter=1
 $\tilde\rho(\mathfrak{t})$ =
 $( \frac{1}{16},
\frac{1}{8},
\frac{9}{16} )
$,

\vskip 0.7ex
\hangindent=4em \hangafter=1
 $\tilde\rho(\mathfrak{s})$ =
$\mathrm{i}$($-\frac{1}{2}$,
$\sqrt{\frac{1}{2}}$,
$\frac{1}{2}$;\ \ 
$0$,
$\sqrt{\frac{1}{2}}$;\ \ 
$-\frac{1}{2}$)

 \vskip 1ex \setlength{\leftskip}{2em}

\grey{Try $U_0$ =
$\begin{pmatrix}
1 \\ 
\end{pmatrix}
$ $\oplus
\begin{pmatrix}
1 \\ 
\end{pmatrix}
$ $\oplus
\begin{pmatrix}
1 \\ 
\end{pmatrix}
$:}\ \ \ \ \ 
\grey{$U_0\tilde\rho(\mathfrak{s})U_0^\dagger$ =}

\grey{$\begin{pmatrix}
-\frac{1}{2}\mathrm{i},
& \sqrt{\frac{1}{2}}\mathrm{i},
& \frac{1}{2}\mathrm{i} \\ 
\sqrt{\frac{1}{2}}\mathrm{i},
& 0,
& \sqrt{\frac{1}{2}}\mathrm{i} \\ 
\frac{1}{2}\mathrm{i},
& \sqrt{\frac{1}{2}}\mathrm{i},
& -\frac{1}{2}\mathrm{i} \\ 
\end{pmatrix}
$}

\grey{Try different $u$'s and signed diagonal matrix $V_\mathrm{sd}$'s:}

 \grey{
\begin{tabular}{|r|l|l|l|}
\hline
$3_{16}^{1,0}:\ u$ 
 & 0 & 1 & 2\\ 
 \hline
$D_\rho$ conditions 
 & 0 & 0 & 0\\ 
 \hline
$[\rho(\mathfrak{s})\rho(\mathfrak{t})]^3
 = \rho^2(\mathfrak{s}) = \tilde C$ 
 & 0 & 0 & 0\\ 
 \hline
$\rho(\mathfrak{s})_{iu}\rho^*(\mathfrak{s})_{ju} \in \mathbb{R}$ 
 & 0 & 0 & 0\\ 
 \hline
$\rho(\mathfrak{s})_{i u} \neq 0$  
 & 0 & 2 & 0\\ 
 \hline
$\mathrm{cnd}(S)$, $\mathrm{cnd}(\rho(\mathfrak{s}))$ 
 & 0 & - & 0\\ 
 \hline
$\mathrm{norm}(D^2)$ factors
 & 0 & - & 0\\ 
 \hline
$1/\rho(\mathfrak{s})_{iu} = $ cyc-int 
 & 0 & - & 0\\ 
 \hline
norm$(1/\rho(\mathfrak{s})_{iu})$ factors
 & 0 & - & 0\\ 
 \hline
$\frac{S_{ij}}{S_{uj}} = $ cyc-int
 & 0 & - & 0\\ 
 \hline
$N^{ij}_k \in \mathbb{N}$
 & 0 & - & 0\\ 
 \hline
$\exists\ j \text{ that } \frac{S_{ij}}{S_{uj}} \geq 1 $
 & 0 & - & 0\\ 
 \hline
FS indicator
 & 0 & - & 0\\ 
 \hline
$C = $ perm-mat
 & 0 & - & 0\\ 
 \hline
\end{tabular}

Number of valid $(S,T)$ pairs: 2 \vskip 2ex }%grey

Total number of valid $(S,T)$ pairs: 2

 \vskip 4ex

\ \setlength{\leftskip}{0em}

\

The $S,T$ matrices obtained above (the black or the blue entries below), plus
their Galois conjugations (the grey entries below), form the following list of
rank-3 $S,T$ matrices.  For details and notations, see Section 2 of this file.

\

\noindent1. ind = $(3;7
)_{1}^{1}$:\ \ 
$d_i$ = ($1.0$,
$1.801$,
$2.246$) 

\vskip 0.7ex
\hangindent=3em \hangafter=1
$D^2=$ 9.295 = 
 $6+3c^{1}_{7}
+c^{2}_{7}
$

\vskip 0.7ex
\hangindent=3em \hangafter=1
$T = ( 0,
\frac{1}{7},
\frac{5}{7} )
$,

\vskip 0.7ex
\hangindent=3em \hangafter=1
$S$ = ($ 1$,
$ \xi_{7}^{2,1}$,
$ \xi_{7}^{3,1}$;\ \ 
$ -\xi_{7}^{3,1}$,
$ 1$;\ \ 
$ -\xi_{7}^{2,1}$)

\vskip 1ex 
\color{grey}

\noindent2. ind = $(3;7
)_{1}^{6}$:\ \ 
$d_i$ = ($1.0$,
$1.801$,
$2.246$) 

\vskip 0.7ex
\hangindent=3em \hangafter=1
$D^2=$ 9.295 = 
 $6+3c^{1}_{7}
+c^{2}_{7}
$

\vskip 0.7ex
\hangindent=3em \hangafter=1
$T = ( 0,
\frac{6}{7},
\frac{2}{7} )
$,

\vskip 0.7ex
\hangindent=3em \hangafter=1
$S$ = ($ 1$,
$ \xi_{7}^{2,1}$,
$ \xi_{7}^{3,1}$;\ \ 
$ -\xi_{7}^{3,1}$,
$ 1$;\ \ 
$ -\xi_{7}^{2,1}$)

\vskip 1ex 
\color{grey}

\noindent3. ind = $(3;7
)_{1}^{3}$:\ \ 
$d_i$ = ($1.0$,
$0.445$,
$-0.801$) 

\vskip 0.7ex
\hangindent=3em \hangafter=1
$D^2=$ 1.841 = 
 $3-2c^{1}_{7}
-3c^{2}_{7}
$

\vskip 0.7ex
\hangindent=3em \hangafter=1
$T = ( 0,
\frac{3}{7},
\frac{1}{7} )
$,

\vskip 0.7ex
\hangindent=3em \hangafter=1
$S$ = ($ 1$,
$ -c^{2}_{7}
$,
$ -\xi_{7}^{2,3}$;\ \ 
$ \xi_{7}^{2,3}$,
$ 1$;\ \ 
$ c^{2}_{7}
$)

Not pseudo-unitary. 

\vskip 1ex 
\color{grey}

\noindent4. ind = $(3;7
)_{1}^{2}$:\ \ 
$d_i$ = ($1.0$,
$0.554$,
$-1.246$) 

\vskip 0.7ex
\hangindent=3em \hangafter=1
$D^2=$ 2.862 = 
 $5-c^{1}_{7}
+2c^{2}_{7}
$

\vskip 0.7ex
\hangindent=3em \hangafter=1
$T = ( 0,
\frac{3}{7},
\frac{2}{7} )
$,

\vskip 0.7ex
\hangindent=3em \hangafter=1
$S$ = ($ 1$,
$ \xi_{7}^{1,2}$,
$ -c^{1}_{7}
$;\ \ 
$ c^{1}_{7}
$,
$ 1$;\ \ 
$ -\xi_{7}^{1,2}$)

Not pseudo-unitary. 

\vskip 1ex 
\color{grey}

\noindent5. ind = $(3;7
)_{1}^{5}$:\ \ 
$d_i$ = ($1.0$,
$0.554$,
$-1.246$) 

\vskip 0.7ex
\hangindent=3em \hangafter=1
$D^2=$ 2.862 = 
 $5-c^{1}_{7}
+2c^{2}_{7}
$

\vskip 0.7ex
\hangindent=3em \hangafter=1
$T = ( 0,
\frac{4}{7},
\frac{5}{7} )
$,

\vskip 0.7ex
\hangindent=3em \hangafter=1
$S$ = ($ 1$,
$ \xi_{7}^{1,2}$,
$ -c^{1}_{7}
$;\ \ 
$ c^{1}_{7}
$,
$ 1$;\ \ 
$ -\xi_{7}^{1,2}$)

Not pseudo-unitary. 

\vskip 1ex 
\color{grey}

\noindent6. ind = $(3;7
)_{1}^{4}$:\ \ 
$d_i$ = ($1.0$,
$0.445$,
$-0.801$) 

\vskip 0.7ex
\hangindent=3em \hangafter=1
$D^2=$ 1.841 = 
 $3-2c^{1}_{7}
-3c^{2}_{7}
$

\vskip 0.7ex
\hangindent=3em \hangafter=1
$T = ( 0,
\frac{4}{7},
\frac{6}{7} )
$,

\vskip 0.7ex
\hangindent=3em \hangafter=1
$S$ = ($ 1$,
$ -c^{2}_{7}
$,
$ -\xi_{7}^{2,3}$;\ \ 
$ \xi_{7}^{2,3}$,
$ 1$;\ \ 
$ c^{2}_{7}
$)

Not pseudo-unitary. 

\vskip 1ex 

 \color{black} \vskip 2ex

\noindent7. ind = $(3;16
)_{1}^{1}$:\ \ 
$d_i$ = ($1.0$,
$1.0$,
$1.414$) 

\vskip 0.7ex
\hangindent=3em \hangafter=1
$D^2=$ 4.0 = 
 $4$

\vskip 0.7ex
\hangindent=3em \hangafter=1
$T = ( 0,
\frac{1}{2},
\frac{1}{16} )
$,

\vskip 0.7ex
\hangindent=3em \hangafter=1
$S$ = ($ 1$,
$ 1$,
$ \sqrt{2}$;\ \ 
$ 1$,
$ -\sqrt{2}$;\ \ 
$0$)

\vskip 1ex 
\color{grey}

\noindent8. ind = $(3;16
)_{1}^{7}$:\ \ 
$d_i$ = ($1.0$,
$1.0$,
$1.414$) 

\vskip 0.7ex
\hangindent=3em \hangafter=1
$D^2=$ 4.0 = 
 $4$

\vskip 0.7ex
\hangindent=3em \hangafter=1
$T = ( 0,
\frac{1}{2},
\frac{7}{16} )
$,

\vskip 0.7ex
\hangindent=3em \hangafter=1
$S$ = ($ 1$,
$ 1$,
$ \sqrt{2}$;\ \ 
$ 1$,
$ -\sqrt{2}$;\ \ 
$0$)

\vskip 1ex 
\color{grey}

\noindent9. ind = $(3;16
)_{1}^{9}$:\ \ 
$d_i$ = ($1.0$,
$1.0$,
$1.414$) 

\vskip 0.7ex
\hangindent=3em \hangafter=1
$D^2=$ 4.0 = 
 $4$

\vskip 0.7ex
\hangindent=3em \hangafter=1
$T = ( 0,
\frac{1}{2},
\frac{9}{16} )
$,

\vskip 0.7ex
\hangindent=3em \hangafter=1
$S$ = ($ 1$,
$ 1$,
$ \sqrt{2}$;\ \ 
$ 1$,
$ -\sqrt{2}$;\ \ 
$0$)

\vskip 1ex 
\color{grey}

\noindent10. ind = $(3;16
)_{1}^{15}$:\ \ 
$d_i$ = ($1.0$,
$1.0$,
$1.414$) 

\vskip 0.7ex
\hangindent=3em \hangafter=1
$D^2=$ 4.0 = 
 $4$

\vskip 0.7ex
\hangindent=3em \hangafter=1
$T = ( 0,
\frac{1}{2},
\frac{15}{16} )
$,

\vskip 0.7ex
\hangindent=3em \hangafter=1
$S$ = ($ 1$,
$ 1$,
$ \sqrt{2}$;\ \ 
$ 1$,
$ -\sqrt{2}$;\ \ 
$0$)

\vskip 1ex 
\color{grey}

\noindent11. ind = $(3;16
)_{1}^{3}$:\ \ 
$d_i$ = ($1.0$,
$1.0$,
$-1.414$) 

\vskip 0.7ex
\hangindent=3em \hangafter=1
$D^2=$ 4.0 = 
 $4$

\vskip 0.7ex
\hangindent=3em \hangafter=1
$T = ( 0,
\frac{1}{2},
\frac{3}{16} )
$,

\vskip 0.7ex
\hangindent=3em \hangafter=1
$S$ = ($ 1$,
$ 1$,
$ -\sqrt{2}$;\ \ 
$ 1$,
$ \sqrt{2}$;\ \ 
$0$)

Pseudo-unitary $\sim$  
$(3;16
)_{2}^{1}$

\vskip 1ex 
\color{grey}

\noindent12. ind = $(3;16
)_{1}^{5}$:\ \ 
$d_i$ = ($1.0$,
$1.0$,
$-1.414$) 

\vskip 0.7ex
\hangindent=3em \hangafter=1
$D^2=$ 4.0 = 
 $4$

\vskip 0.7ex
\hangindent=3em \hangafter=1
$T = ( 0,
\frac{1}{2},
\frac{5}{16} )
$,

\vskip 0.7ex
\hangindent=3em \hangafter=1
$S$ = ($ 1$,
$ 1$,
$ -\sqrt{2}$;\ \ 
$ 1$,
$ \sqrt{2}$;\ \ 
$0$)

Pseudo-unitary $\sim$  
$(3;16
)_{2}^{7}$

\vskip 1ex 
\color{grey}

\noindent13. ind = $(3;16
)_{1}^{11}$:\ \ 
$d_i$ = ($1.0$,
$1.0$,
$-1.414$) 

\vskip 0.7ex
\hangindent=3em \hangafter=1
$D^2=$ 4.0 = 
 $4$

\vskip 0.7ex
\hangindent=3em \hangafter=1
$T = ( 0,
\frac{1}{2},
\frac{11}{16} )
$,

\vskip 0.7ex
\hangindent=3em \hangafter=1
$S$ = ($ 1$,
$ 1$,
$ -\sqrt{2}$;\ \ 
$ 1$,
$ \sqrt{2}$;\ \ 
$0$)

Pseudo-unitary $\sim$  
$(3;16
)_{2}^{9}$

\vskip 1ex 
\color{grey}

\noindent14. ind = $(3;16
)_{1}^{13}$:\ \ 
$d_i$ = ($1.0$,
$1.0$,
$-1.414$) 

\vskip 0.7ex
\hangindent=3em \hangafter=1
$D^2=$ 4.0 = 
 $4$

\vskip 0.7ex
\hangindent=3em \hangafter=1
$T = ( 0,
\frac{1}{2},
\frac{13}{16} )
$,

\vskip 0.7ex
\hangindent=3em \hangafter=1
$S$ = ($ 1$,
$ 1$,
$ -\sqrt{2}$;\ \ 
$ 1$,
$ \sqrt{2}$;\ \ 
$0$)

Pseudo-unitary $\sim$  
$(3;16
)_{2}^{15}$

\vskip 1ex 

 \color{black} \vskip 2ex

\noindent15. ind = $(3;16
)_{2}^{1}$:\ \ 
$d_i$ = ($1.0$,
$1.0$,
$1.414$) 

\vskip 0.7ex
\hangindent=3em \hangafter=1
$D^2=$ 4.0 = 
 $4$

\vskip 0.7ex
\hangindent=3em \hangafter=1
$T = ( 0,
\frac{1}{2},
\frac{3}{16} )
$,

\vskip 0.7ex
\hangindent=3em \hangafter=1
$S$ = ($ 1$,
$ 1$,
$ \sqrt{2}$;\ \ 
$ 1$,
$ -\sqrt{2}$;\ \ 
$0$)

\vskip 1ex 
\color{grey}

\noindent16. ind = $(3;16
)_{2}^{7}$:\ \ 
$d_i$ = ($1.0$,
$1.0$,
$1.414$) 

\vskip 0.7ex
\hangindent=3em \hangafter=1
$D^2=$ 4.0 = 
 $4$

\vskip 0.7ex
\hangindent=3em \hangafter=1
$T = ( 0,
\frac{1}{2},
\frac{5}{16} )
$,

\vskip 0.7ex
\hangindent=3em \hangafter=1
$S$ = ($ 1$,
$ 1$,
$ \sqrt{2}$;\ \ 
$ 1$,
$ -\sqrt{2}$;\ \ 
$0$)

\vskip 1ex 
\color{grey}

\noindent17. ind = $(3;16
)_{2}^{9}$:\ \ 
$d_i$ = ($1.0$,
$1.0$,
$1.414$) 

\vskip 0.7ex
\hangindent=3em \hangafter=1
$D^2=$ 4.0 = 
 $4$

\vskip 0.7ex
\hangindent=3em \hangafter=1
$T = ( 0,
\frac{1}{2},
\frac{11}{16} )
$,

\vskip 0.7ex
\hangindent=3em \hangafter=1
$S$ = ($ 1$,
$ 1$,
$ \sqrt{2}$;\ \ 
$ 1$,
$ -\sqrt{2}$;\ \ 
$0$)

\vskip 1ex 
\color{grey}

\noindent18. ind = $(3;16
)_{2}^{15}$:\ \ 
$d_i$ = ($1.0$,
$1.0$,
$1.414$) 

\vskip 0.7ex
\hangindent=3em \hangafter=1
$D^2=$ 4.0 = 
 $4$

\vskip 0.7ex
\hangindent=3em \hangafter=1
$T = ( 0,
\frac{1}{2},
\frac{13}{16} )
$,

\vskip 0.7ex
\hangindent=3em \hangafter=1
$S$ = ($ 1$,
$ 1$,
$ \sqrt{2}$;\ \ 
$ 1$,
$ -\sqrt{2}$;\ \ 
$0$)

\vskip 1ex 
\color{grey}

\noindent19. ind = $(3;16
)_{2}^{11}$:\ \ 
$d_i$ = ($1.0$,
$1.0$,
$-1.414$) 

\vskip 0.7ex
\hangindent=3em \hangafter=1
$D^2=$ 4.0 = 
 $4$

\vskip 0.7ex
\hangindent=3em \hangafter=1
$T = ( 0,
\frac{1}{2},
\frac{1}{16} )
$,

\vskip 0.7ex
\hangindent=3em \hangafter=1
$S$ = ($ 1$,
$ 1$,
$ -\sqrt{2}$;\ \ 
$ 1$,
$ \sqrt{2}$;\ \ 
$0$)

Pseudo-unitary $\sim$  
$(3;16
)_{1}^{1}$

\vskip 1ex 
\color{grey}

\noindent20. ind = $(3;16
)_{2}^{13}$:\ \ 
$d_i$ = ($1.0$,
$1.0$,
$-1.414$) 

\vskip 0.7ex
\hangindent=3em \hangafter=1
$D^2=$ 4.0 = 
 $4$

\vskip 0.7ex
\hangindent=3em \hangafter=1
$T = ( 0,
\frac{1}{2},
\frac{7}{16} )
$,

\vskip 0.7ex
\hangindent=3em \hangafter=1
$S$ = ($ 1$,
$ 1$,
$ -\sqrt{2}$;\ \ 
$ 1$,
$ \sqrt{2}$;\ \ 
$0$)

Pseudo-unitary $\sim$  
$(3;16
)_{1}^{7}$

\vskip 1ex 
\color{grey}

\noindent21. ind = $(3;16
)_{2}^{3}$:\ \ 
$d_i$ = ($1.0$,
$1.0$,
$-1.414$) 

\vskip 0.7ex
\hangindent=3em \hangafter=1
$D^2=$ 4.0 = 
 $4$

\vskip 0.7ex
\hangindent=3em \hangafter=1
$T = ( 0,
\frac{1}{2},
\frac{9}{16} )
$,

\vskip 0.7ex
\hangindent=3em \hangafter=1
$S$ = ($ 1$,
$ 1$,
$ -\sqrt{2}$;\ \ 
$ 1$,
$ \sqrt{2}$;\ \ 
$0$)

Pseudo-unitary $\sim$  
$(3;16
)_{1}^{9}$

\vskip 1ex 
\color{grey}

\noindent22. ind = $(3;16
)_{2}^{5}$:\ \ 
$d_i$ = ($1.0$,
$1.0$,
$-1.414$) 

\vskip 0.7ex
\hangindent=3em \hangafter=1
$D^2=$ 4.0 = 
 $4$

\vskip 0.7ex
\hangindent=3em \hangafter=1
$T = ( 0,
\frac{1}{2},
\frac{15}{16} )
$,

\vskip 0.7ex
\hangindent=3em \hangafter=1
$S$ = ($ 1$,
$ 1$,
$ -\sqrt{2}$;\ \ 
$ 1$,
$ \sqrt{2}$;\ \ 
$0$)

Pseudo-unitary $\sim$  
$(3;16
)_{1}^{15}$

\vskip 1ex 

 \color{black} \vskip 2ex

\

The above list includes all modular data (unitary or non-unitary) from resolved
$\SL$ representations and non-integral MTCs.
Since for rank 3, all the passing $\SL$ representations are resolved, the list
actually includes all modular data from non-integral MTCs.

\section{Rank-4 modular data}
\label{Section6}

\subsection{A list of 4-dimensional irrep-sum $\SL$ representations}

The following is a list of 4-dimensional irrep-sum $\SL$ representations, with
3 types of representations omitted and one type skipped.  For details and notations, see Section
3 of this file.

\

 \color{blue}

\noindent 1: (dims,levels) = $(3\oplus
1;5,
1
)$,
irreps = $3_{5}^{1}\oplus
1_{1}^{1}$,
pord$(\rho_\text{isum}(\mathfrak{t})) = 5$,

\vskip 0.7ex
\hangindent=5.5em \hangafter=1
{\white .}\hskip 1em $\rho_\text{isum}(\mathfrak{t})$ =
 $( 0,
\frac{1}{5},
\frac{4}{5} )
\oplus
( 0 )
$,

\vskip 0.7ex
\hangindent=5.5em \hangafter=1
{\white .}\hskip 1em $\rho_\text{isum}(\mathfrak{s})$ =
($\sqrt{\frac{1}{5}}$,
$-\sqrt{\frac{2}{5}}$,
$-\sqrt{\frac{2}{5}}$;
$-\frac{5+\sqrt{5}}{10}$,
$\frac{5-\sqrt{5}}{10}$;
$-\frac{5+\sqrt{5}}{10}$)
 $\oplus$
($1$)

Pass. 

 \ \color{black}

 \color{blue}

\noindent 2: (dims,levels) = $(3\oplus
1;5,
1
)$,
irreps = $3_{5}^{3}\oplus
1_{1}^{1}$,
pord$(\rho_\text{isum}(\mathfrak{t})) = 5$,

\vskip 0.7ex
\hangindent=5.5em \hangafter=1
{\white .}\hskip 1em $\rho_\text{isum}(\mathfrak{t})$ =
 $( 0,
\frac{2}{5},
\frac{3}{5} )
\oplus
( 0 )
$,

\vskip 0.7ex
\hangindent=5.5em \hangafter=1
{\white .}\hskip 1em $\rho_\text{isum}(\mathfrak{s})$ =
($-\sqrt{\frac{1}{5}}$,
$-\sqrt{\frac{2}{5}}$,
$-\sqrt{\frac{2}{5}}$;
$\frac{-5+\sqrt{5}}{10}$,
$\frac{5+\sqrt{5}}{10}$;
$\frac{-5+\sqrt{5}}{10}$)
 $\oplus$
($1$)

Pass. 

 \ \color{black}

\noindent 3: (dims,levels) = $(3\oplus
1;8,
1
)$,
irreps = $3_{8}^{1,0}\oplus
1_{1}^{1}$,
pord$(\rho_\text{isum}(\mathfrak{t})) = 8$,

\vskip 0.7ex
\hangindent=5.5em \hangafter=1
{\white .}\hskip 1em $\rho_\text{isum}(\mathfrak{t})$ =
 $( 0,
\frac{1}{8},
\frac{5}{8} )
\oplus
( 0 )
$,

\vskip 0.7ex
\hangindent=5.5em \hangafter=1
{\white .}\hskip 1em $\rho_\text{isum}(\mathfrak{s})$ =
$\mathrm{i}$($0$,
$\sqrt{\frac{1}{2}}$,
$\sqrt{\frac{1}{2}}$;\ \ 
$-\frac{1}{2}$,
$\frac{1}{2}$;\ \ 
$-\frac{1}{2}$)
 $\oplus$
($1$)

Fail:
Integral: $D_{\rho}(\sigma)_{\theta} \propto $ id,
 for all $\sigma$ and all $\theta$-eigenspaces that can contain unit. Prop. B.5 (6)

 \ \color{black}

\noindent 4: (dims,levels) = $(3\oplus
1;8,
1
)$,
irreps = $3_{8}^{3,0}\oplus
1_{1}^{1}$,
pord$(\rho_\text{isum}(\mathfrak{t})) = 8$,

\vskip 0.7ex
\hangindent=5.5em \hangafter=1
{\white .}\hskip 1em $\rho_\text{isum}(\mathfrak{t})$ =
 $( 0,
\frac{3}{8},
\frac{7}{8} )
\oplus
( 0 )
$,

\vskip 0.7ex
\hangindent=5.5em \hangafter=1
{\white .}\hskip 1em $\rho_\text{isum}(\mathfrak{s})$ =
$\mathrm{i}$($0$,
$\sqrt{\frac{1}{2}}$,
$\sqrt{\frac{1}{2}}$;\ \ 
$\frac{1}{2}$,
$-\frac{1}{2}$;\ \ 
$\frac{1}{2}$)
 $\oplus$
($1$)

Fail:
Integral: $D_{\rho}(\sigma)_{\theta} \propto $ id,
 for all $\sigma$ and all $\theta$-eigenspaces that can contain unit. Prop. B.5 (6)

 \ \color{black}

 \color{blue}

\noindent 5: (dims,levels) = $(3\oplus
1;10,
2
)$,
irreps = $3_{5}^{3}
\hskip -1.5pt \otimes \hskip -1.5pt
1_{2}^{1,0}\oplus
1_{2}^{1,0}$,
pord$(\rho_\text{isum}(\mathfrak{t})) = 5$,

\vskip 0.7ex
\hangindent=5.5em \hangafter=1
{\white .}\hskip 1em $\rho_\text{isum}(\mathfrak{t})$ =
 $( \frac{1}{2},
\frac{1}{10},
\frac{9}{10} )
\oplus
( \frac{1}{2} )
$,

\vskip 0.7ex
\hangindent=5.5em \hangafter=1
{\white .}\hskip 1em $\rho_\text{isum}(\mathfrak{s})$ =
($\sqrt{\frac{1}{5}}$,
$-\sqrt{\frac{2}{5}}$,
$-\sqrt{\frac{2}{5}}$;
$\frac{5-\sqrt{5}}{10}$,
$-\frac{5+\sqrt{5}}{10}$;
$\frac{5-\sqrt{5}}{10}$)
 $\oplus$
($-1$)

Pass. 

 \ \color{black}

 \color{blue}

\noindent 6: (dims,levels) = $(3\oplus
1;10,
2
)$,
irreps = $3_{5}^{1}
\hskip -1.5pt \otimes \hskip -1.5pt
1_{2}^{1,0}\oplus
1_{2}^{1,0}$,
pord$(\rho_\text{isum}(\mathfrak{t})) = 5$,

\vskip 0.7ex
\hangindent=5.5em \hangafter=1
{\white .}\hskip 1em $\rho_\text{isum}(\mathfrak{t})$ =
 $( \frac{1}{2},
\frac{3}{10},
\frac{7}{10} )
\oplus
( \frac{1}{2} )
$,

\vskip 0.7ex
\hangindent=5.5em \hangafter=1
{\white .}\hskip 1em $\rho_\text{isum}(\mathfrak{s})$ =
($-\sqrt{\frac{1}{5}}$,
$-\sqrt{\frac{2}{5}}$,
$-\sqrt{\frac{2}{5}}$;
$\frac{5+\sqrt{5}}{10}$,
$\frac{-5+\sqrt{5}}{10}$;
$\frac{5+\sqrt{5}}{10}$)
 $\oplus$
($-1$)

Pass. 

 \ \color{black}

 \color{blue}

\noindent 7: (dims,levels) = $(3\oplus
1;15,
3
)$,
irreps = $3_{5}^{1}
\hskip -1.5pt \otimes \hskip -1.5pt
1_{3}^{1,0}\oplus
1_{3}^{1,0}$,
pord$(\rho_\text{isum}(\mathfrak{t})) = 5$,

\vskip 0.7ex
\hangindent=5.5em \hangafter=1
{\white .}\hskip 1em $\rho_\text{isum}(\mathfrak{t})$ =
 $( \frac{1}{3},
\frac{2}{15},
\frac{8}{15} )
\oplus
( \frac{1}{3} )
$,

\vskip 0.7ex
\hangindent=5.5em \hangafter=1
{\white .}\hskip 1em $\rho_\text{isum}(\mathfrak{s})$ =
($\sqrt{\frac{1}{5}}$,
$-\sqrt{\frac{2}{5}}$,
$-\sqrt{\frac{2}{5}}$;
$-\frac{5+\sqrt{5}}{10}$,
$\frac{5-\sqrt{5}}{10}$;
$-\frac{5+\sqrt{5}}{10}$)
 $\oplus$
($1$)

Pass. 

 \ \color{black}

 \color{blue}

\noindent 8: (dims,levels) = $(3\oplus
1;15,
3
)$,
irreps = $3_{5}^{3}
\hskip -1.5pt \otimes \hskip -1.5pt
1_{3}^{1,0}\oplus
1_{3}^{1,0}$,
pord$(\rho_\text{isum}(\mathfrak{t})) = 5$,

\vskip 0.7ex
\hangindent=5.5em \hangafter=1
{\white .}\hskip 1em $\rho_\text{isum}(\mathfrak{t})$ =
 $( \frac{1}{3},
\frac{11}{15},
\frac{14}{15} )
\oplus
( \frac{1}{3} )
$,

\vskip 0.7ex
\hangindent=5.5em \hangafter=1
{\white .}\hskip 1em $\rho_\text{isum}(\mathfrak{s})$ =
($-\sqrt{\frac{1}{5}}$,
$-\sqrt{\frac{2}{5}}$,
$-\sqrt{\frac{2}{5}}$;
$\frac{-5+\sqrt{5}}{10}$,
$\frac{5+\sqrt{5}}{10}$;
$\frac{-5+\sqrt{5}}{10}$)
 $\oplus$
($1$)

Pass. 

 \ \color{black}

 \color{blue}

\noindent 9: (dims,levels) = $(3\oplus
1;20,
4
)$,
irreps = $3_{5}^{1}
\hskip -1.5pt \otimes \hskip -1.5pt
1_{4}^{1,0}\oplus
1_{4}^{1,0}$,
pord$(\rho_\text{isum}(\mathfrak{t})) = 5$,

\vskip 0.7ex
\hangindent=5.5em \hangafter=1
{\white .}\hskip 1em $\rho_\text{isum}(\mathfrak{t})$ =
 $( \frac{1}{4},
\frac{1}{20},
\frac{9}{20} )
\oplus
( \frac{1}{4} )
$,

\vskip 0.7ex
\hangindent=5.5em \hangafter=1
{\white .}\hskip 1em $\rho_\text{isum}(\mathfrak{s})$ =
$\mathrm{i}$($\sqrt{\frac{1}{5}}$,
$\sqrt{\frac{2}{5}}$,
$\sqrt{\frac{2}{5}}$;\ \ 
$-\frac{5+\sqrt{5}}{10}$,
$\frac{5-\sqrt{5}}{10}$;\ \ 
$-\frac{5+\sqrt{5}}{10}$)
 $\oplus$
$\mathrm{i}$($1$)

Pass. 

 \ \color{black}

 \color{blue}

\noindent 10: (dims,levels) = $(3\oplus
1;20,
4
)$,
irreps = $3_{5}^{3}
\hskip -1.5pt \otimes \hskip -1.5pt
1_{4}^{1,0}\oplus
1_{4}^{1,0}$,
pord$(\rho_\text{isum}(\mathfrak{t})) = 5$,

\vskip 0.7ex
\hangindent=5.5em \hangafter=1
{\white .}\hskip 1em $\rho_\text{isum}(\mathfrak{t})$ =
 $( \frac{1}{4},
\frac{13}{20},
\frac{17}{20} )
\oplus
( \frac{1}{4} )
$,

\vskip 0.7ex
\hangindent=5.5em \hangafter=1
{\white .}\hskip 1em $\rho_\text{isum}(\mathfrak{s})$ =
$\mathrm{i}$($-\sqrt{\frac{1}{5}}$,
$\sqrt{\frac{2}{5}}$,
$\sqrt{\frac{2}{5}}$;\ \ 
$\frac{-5+\sqrt{5}}{10}$,
$\frac{5+\sqrt{5}}{10}$;\ \ 
$\frac{-5+\sqrt{5}}{10}$)
 $\oplus$
$\mathrm{i}$($1$)

Pass. 

 \ \color{black}

\noindent 11: (dims,levels) = $(3\oplus
1;24,
3
)$,
irreps = $3_{8}^{3,0}
\hskip -1.5pt \otimes \hskip -1.5pt
1_{3}^{1,0}\oplus
1_{3}^{1,0}$,
pord$(\rho_\text{isum}(\mathfrak{t})) = 8$,

\vskip 0.7ex
\hangindent=5.5em \hangafter=1
{\white .}\hskip 1em $\rho_\text{isum}(\mathfrak{t})$ =
 $( \frac{1}{3},
\frac{5}{24},
\frac{17}{24} )
\oplus
( \frac{1}{3} )
$,

\vskip 0.7ex
\hangindent=5.5em \hangafter=1
{\white .}\hskip 1em $\rho_\text{isum}(\mathfrak{s})$ =
$\mathrm{i}$($0$,
$\sqrt{\frac{1}{2}}$,
$\sqrt{\frac{1}{2}}$;\ \ 
$\frac{1}{2}$,
$-\frac{1}{2}$;\ \ 
$\frac{1}{2}$)
 $\oplus$
($1$)

Fail:
Integral: $D_{\rho}(\sigma)_{\theta} \propto $ id,
 for all $\sigma$ and all $\theta$-eigenspaces that can contain unit. Prop. B.5 (6)

 \ \color{black}

\noindent 12: (dims,levels) = $(3\oplus
1;24,
3
)$,
irreps = $3_{8}^{1,0}
\hskip -1.5pt \otimes \hskip -1.5pt
1_{3}^{1,0}\oplus
1_{3}^{1,0}$,
pord$(\rho_\text{isum}(\mathfrak{t})) = 8$,

\vskip 0.7ex
\hangindent=5.5em \hangafter=1
{\white .}\hskip 1em $\rho_\text{isum}(\mathfrak{t})$ =
 $( \frac{1}{3},
\frac{11}{24},
\frac{23}{24} )
\oplus
( \frac{1}{3} )
$,

\vskip 0.7ex
\hangindent=5.5em \hangafter=1
{\white .}\hskip 1em $\rho_\text{isum}(\mathfrak{s})$ =
$\mathrm{i}$($0$,
$\sqrt{\frac{1}{2}}$,
$\sqrt{\frac{1}{2}}$;\ \ 
$-\frac{1}{2}$,
$\frac{1}{2}$;\ \ 
$-\frac{1}{2}$)
 $\oplus$
($1$)

Fail:
Integral: $D_{\rho}(\sigma)_{\theta} \propto $ id,
 for all $\sigma$ and all $\theta$-eigenspaces that can contain unit. Prop. B.5 (6)

 \ \color{black}

 \color{blue}

\noindent 13: (dims,levels) = $(3\oplus
1;30,
6
)$,
irreps = $3_{5}^{1}
\hskip -1.5pt \otimes \hskip -1.5pt
1_{3}^{1,0}
\hskip -1.5pt \otimes \hskip -1.5pt
1_{2}^{1,0}\oplus
1_{3}^{1,0}
\hskip -1.5pt \otimes \hskip -1.5pt
1_{2}^{1,0}$,
pord$(\rho_\text{isum}(\mathfrak{t})) = 5$,

\vskip 0.7ex
\hangindent=5.5em \hangafter=1
{\white .}\hskip 1em $\rho_\text{isum}(\mathfrak{t})$ =
 $( \frac{5}{6},
\frac{1}{30},
\frac{19}{30} )
\oplus
( \frac{5}{6} )
$,

\vskip 0.7ex
\hangindent=5.5em \hangafter=1
{\white .}\hskip 1em $\rho_\text{isum}(\mathfrak{s})$ =
($-\sqrt{\frac{1}{5}}$,
$-\sqrt{\frac{2}{5}}$,
$-\sqrt{\frac{2}{5}}$;
$\frac{5+\sqrt{5}}{10}$,
$\frac{-5+\sqrt{5}}{10}$;
$\frac{5+\sqrt{5}}{10}$)
 $\oplus$
($-1$)

Pass. 

 \ \color{black}

 \color{blue}

\noindent 14: (dims,levels) = $(3\oplus
1;30,
6
)$,
irreps = $3_{5}^{3}
\hskip -1.5pt \otimes \hskip -1.5pt
1_{3}^{1,0}
\hskip -1.5pt \otimes \hskip -1.5pt
1_{2}^{1,0}\oplus
1_{3}^{1,0}
\hskip -1.5pt \otimes \hskip -1.5pt
1_{2}^{1,0}$,
pord$(\rho_\text{isum}(\mathfrak{t})) = 5$,

\vskip 0.7ex
\hangindent=5.5em \hangafter=1
{\white .}\hskip 1em $\rho_\text{isum}(\mathfrak{t})$ =
 $( \frac{5}{6},
\frac{7}{30},
\frac{13}{30} )
\oplus
( \frac{5}{6} )
$,

\vskip 0.7ex
\hangindent=5.5em \hangafter=1
{\white .}\hskip 1em $\rho_\text{isum}(\mathfrak{s})$ =
($\sqrt{\frac{1}{5}}$,
$-\sqrt{\frac{2}{5}}$,
$-\sqrt{\frac{2}{5}}$;
$\frac{5-\sqrt{5}}{10}$,
$-\frac{5+\sqrt{5}}{10}$;
$\frac{5-\sqrt{5}}{10}$)
 $\oplus$
($-1$)

Pass. 

 \ \color{black}

 \color{blue}

\noindent 15: (dims,levels) = $(3\oplus
1;60,
12
)$,
irreps = $3_{5}^{3}
\hskip -1.5pt \otimes \hskip -1.5pt
1_{4}^{1,0}
\hskip -1.5pt \otimes \hskip -1.5pt
1_{3}^{1,0}\oplus
1_{4}^{1,0}
\hskip -1.5pt \otimes \hskip -1.5pt
1_{3}^{1,0}$,
pord$(\rho_\text{isum}(\mathfrak{t})) = 5$,

\vskip 0.7ex
\hangindent=5.5em \hangafter=1
{\white .}\hskip 1em $\rho_\text{isum}(\mathfrak{t})$ =
 $( \frac{7}{12},
\frac{11}{60},
\frac{59}{60} )
\oplus
( \frac{7}{12} )
$,

\vskip 0.7ex
\hangindent=5.5em \hangafter=1
{\white .}\hskip 1em $\rho_\text{isum}(\mathfrak{s})$ =
$\mathrm{i}$($-\sqrt{\frac{1}{5}}$,
$\sqrt{\frac{2}{5}}$,
$\sqrt{\frac{2}{5}}$;\ \ 
$\frac{-5+\sqrt{5}}{10}$,
$\frac{5+\sqrt{5}}{10}$;\ \ 
$\frac{-5+\sqrt{5}}{10}$)
 $\oplus$
$\mathrm{i}$($1$)

Pass. 

 \ \color{black}

 \color{blue}

\noindent 16: (dims,levels) = $(3\oplus
1;60,
12
)$,
irreps = $3_{5}^{1}
\hskip -1.5pt \otimes \hskip -1.5pt
1_{4}^{1,0}
\hskip -1.5pt \otimes \hskip -1.5pt
1_{3}^{1,0}\oplus
1_{4}^{1,0}
\hskip -1.5pt \otimes \hskip -1.5pt
1_{3}^{1,0}$,
pord$(\rho_\text{isum}(\mathfrak{t})) = 5$,

\vskip 0.7ex
\hangindent=5.5em \hangafter=1
{\white .}\hskip 1em $\rho_\text{isum}(\mathfrak{t})$ =
 $( \frac{7}{12},
\frac{23}{60},
\frac{47}{60} )
\oplus
( \frac{7}{12} )
$,

\vskip 0.7ex
\hangindent=5.5em \hangafter=1
{\white .}\hskip 1em $\rho_\text{isum}(\mathfrak{s})$ =
$\mathrm{i}$($\sqrt{\frac{1}{5}}$,
$\sqrt{\frac{2}{5}}$,
$\sqrt{\frac{2}{5}}$;\ \ 
$-\frac{5+\sqrt{5}}{10}$,
$\frac{5-\sqrt{5}}{10}$;\ \ 
$-\frac{5+\sqrt{5}}{10}$)
 $\oplus$
$\mathrm{i}$($1$)

Pass. 

 \ \color{black}

 \color{blue}

\noindent 17: (dims,levels) = $(4;5
)$,
irreps = $4_{5,2}^{1}$,
pord$(\rho_\text{isum}(\mathfrak{t})) = 5$,

\vskip 0.7ex
\hangindent=5.5em \hangafter=1
{\white .}\hskip 1em $\rho_\text{isum}(\mathfrak{t})$ =
 $( \frac{1}{5},
\frac{2}{5},
\frac{3}{5},
\frac{4}{5} )
$,

\vskip 0.7ex
\hangindent=5.5em \hangafter=1
{\white .}\hskip 1em $\rho_\text{isum}(\mathfrak{s})$ =
($\sqrt{\frac{1}{5}}$,
$\frac{-5+\sqrt{5}}{10}$,
$-\frac{5+\sqrt{5}}{10}$,
$\sqrt{\frac{1}{5}}$;
$-\sqrt{\frac{1}{5}}$,
$\sqrt{\frac{1}{5}}$,
$\frac{5+\sqrt{5}}{10}$;
$-\sqrt{\frac{1}{5}}$,
$\frac{5-\sqrt{5}}{10}$;
$\sqrt{\frac{1}{5}}$)

Pass. 

 \ \color{black}

\noindent 18: (dims,levels) = $(4;5
)$,
irreps = $4_{5,1}^{1}$,
pord$(\rho_\text{isum}(\mathfrak{t})) = 5$,

\vskip 0.7ex
\hangindent=5.5em \hangafter=1
{\white .}\hskip 1em $\rho_\text{isum}(\mathfrak{t})$ =
 $( \frac{1}{5},
\frac{2}{5},
\frac{3}{5},
\frac{4}{5} )
$,

\vskip 0.7ex
\hangindent=5.5em \hangafter=1
{\white .}\hskip 1em $\rho_\text{isum}(\mathfrak{s})$ =
$\mathrm{i}$($\frac{1}{5}c^{1}_{20}
+\frac{1}{5}c^{3}_{20}
$,
$\frac{2}{5}c^{2}_{15}
+\frac{1}{5}c^{3}_{15}
$,
$-\frac{1}{5}+\frac{2}{5}c^{1}_{15}
-\frac{1}{5}c^{3}_{15}
$,
$\frac{1}{5}c^{1}_{20}
-\frac{1}{5}c^{3}_{20}
$;\ \ 
$-\frac{1}{5}c^{1}_{20}
+\frac{1}{5}c^{3}_{20}
$,
$-\frac{1}{5}c^{1}_{20}
-\frac{1}{5}c^{3}_{20}
$,
$\frac{1}{5}-\frac{2}{5}c^{1}_{15}
+\frac{1}{5}c^{3}_{15}
$;\ \ 
$\frac{1}{5}c^{1}_{20}
-\frac{1}{5}c^{3}_{20}
$,
$\frac{2}{5}c^{2}_{15}
+\frac{1}{5}c^{3}_{15}
$;\ \ 
$-\frac{1}{5}c^{1}_{20}
-\frac{1}{5}c^{3}_{20}
$)

Fail:
cnd($\rho(\mathfrak s)_\mathrm{ndeg}$) = 60 does not divide
 ord($\rho(\mathfrak t)$)=5. Prop. B.4 (2)

 \ \color{black}

\noindent 19: (dims,levels) = $(4;7
)$,
irreps = $4_{7}^{1}$,
pord$(\rho_\text{isum}(\mathfrak{t})) = 7$,

\vskip 0.7ex
\hangindent=5.5em \hangafter=1
{\white .}\hskip 1em $\rho_\text{isum}(\mathfrak{t})$ =
 $( 0,
\frac{1}{7},
\frac{2}{7},
\frac{4}{7} )
$,

\vskip 0.7ex
\hangindent=5.5em \hangafter=1
{\white .}\hskip 1em $\rho_\text{isum}(\mathfrak{s})$ =
$\mathrm{i}$($-\sqrt{\frac{1}{7}}$,
$\sqrt{\frac{2}{7}}$,
$\sqrt{\frac{2}{7}}$,
$\sqrt{\frac{2}{7}}$;\ \ 
$-\frac{1}{\sqrt{7}}c^{2}_{7}
$,
$-\frac{1}{\sqrt{7}}c^{1}_{7}
$,
$\frac{1}{\sqrt{7}\mathrm{i}}s^{5}_{28}
$;\ \ 
$\frac{1}{\sqrt{7}\mathrm{i}}s^{5}_{28}
$,
$-\frac{1}{\sqrt{7}}c^{2}_{7}
$;\ \ 
$-\frac{1}{\sqrt{7}}c^{1}_{7}
$)

Fail:
cnd($\rho(\mathfrak s)_\mathrm{ndeg}$) = 56 does not divide
 ord($\rho(\mathfrak t)$)=7. Prop. B.4 (2)

 \ \color{black}

\noindent 20: (dims,levels) = $(4;7
)$,
irreps = $4_{7}^{3}$,
pord$(\rho_\text{isum}(\mathfrak{t})) = 7$,

\vskip 0.7ex
\hangindent=5.5em \hangafter=1
{\white .}\hskip 1em $\rho_\text{isum}(\mathfrak{t})$ =
 $( 0,
\frac{3}{7},
\frac{5}{7},
\frac{6}{7} )
$,

\vskip 0.7ex
\hangindent=5.5em \hangafter=1
{\white .}\hskip 1em $\rho_\text{isum}(\mathfrak{s})$ =
$\mathrm{i}$($\sqrt{\frac{1}{7}}$,
$\sqrt{\frac{2}{7}}$,
$\sqrt{\frac{2}{7}}$,
$\sqrt{\frac{2}{7}}$;\ \ 
$\frac{1}{\sqrt{7}}c^{1}_{7}
$,
$\frac{1}{\sqrt{7}}c^{2}_{7}
$,
$-\frac{1}{\sqrt{7}\mathrm{i}}s^{5}_{28}
$;\ \ 
$-\frac{1}{\sqrt{7}\mathrm{i}}s^{5}_{28}
$,
$\frac{1}{\sqrt{7}}c^{1}_{7}
$;\ \ 
$\frac{1}{\sqrt{7}}c^{2}_{7}
$)

Fail:
cnd($\rho(\mathfrak s)_\mathrm{ndeg}$) = 56 does not divide
 ord($\rho(\mathfrak t)$)=7. Prop. B.4 (2)

 \ \color{black}

 \color{blue}

\noindent 21: (dims,levels) = $(4;9
)$,
irreps = $4_{9,1}^{1,0}$,
pord$(\rho_\text{isum}(\mathfrak{t})) = 9$,

\vskip 0.7ex
\hangindent=5.5em \hangafter=1
{\white .}\hskip 1em $\rho_\text{isum}(\mathfrak{t})$ =
 $( 0,
\frac{1}{9},
\frac{4}{9},
\frac{7}{9} )
$,

\vskip 0.7ex
\hangindent=5.5em \hangafter=1
{\white .}\hskip 1em $\rho_\text{isum}(\mathfrak{s})$ =
$\mathrm{i}$($0$,
$\sqrt{\frac{1}{3}}$,
$\sqrt{\frac{1}{3}}$,
$\sqrt{\frac{1}{3}}$;\ \ 
$-\frac{1}{3}c^{1}_{36}
$,
$\frac{1}{3}c^{1}_{36}
-\frac{1}{3}c^{5}_{36}
$,
$\frac{1}{3}c^{5}_{36}
$;\ \ 
$\frac{1}{3}c^{5}_{36}
$,
$-\frac{1}{3}c^{1}_{36}
$;\ \ 
$\frac{1}{3}c^{1}_{36}
-\frac{1}{3}c^{5}_{36}
$)

Pass. 

 \ \color{black}

\noindent 22: (dims,levels) = $(4;9
)$,
irreps = $4_{9,2}^{1,0}$,
pord$(\rho_\text{isum}(\mathfrak{t})) = 9$,

\vskip 0.7ex
\hangindent=5.5em \hangafter=1
{\white .}\hskip 1em $\rho_\text{isum}(\mathfrak{t})$ =
 $( 0,
\frac{1}{9},
\frac{4}{9},
\frac{7}{9} )
$,

\vskip 0.7ex
\hangindent=5.5em \hangafter=1
{\white .}\hskip 1em $\rho_\text{isum}(\mathfrak{s})$ =
($0$,
$-\sqrt{\frac{1}{3}}$,
$-\sqrt{\frac{1}{3}}$,
$-\sqrt{\frac{1}{3}}$;
$\frac{1}{3}c^{2}_{9}
$,
$\frac{1}{3} c_9^4 $,
$\frac{1}{3}c^{1}_{9}
$;
$\frac{1}{3}c^{1}_{9}
$,
$\frac{1}{3}c^{2}_{9}
$;
$\frac{1}{3} c_9^4 $)

Fail:
cnd($\rho(\mathfrak s)_\mathrm{ndeg}$) = 36 does not divide
 ord($\rho(\mathfrak t)$)=9. Prop. B.4 (2)

 \ \color{black}

 \color{blue}

\noindent 23: (dims,levels) = $(4;9
)$,
irreps = $4_{9,1}^{2,0}$,
pord$(\rho_\text{isum}(\mathfrak{t})) = 9$,

\vskip 0.7ex
\hangindent=5.5em \hangafter=1
{\white .}\hskip 1em $\rho_\text{isum}(\mathfrak{t})$ =
 $( 0,
\frac{2}{9},
\frac{5}{9},
\frac{8}{9} )
$,

\vskip 0.7ex
\hangindent=5.5em \hangafter=1
{\white .}\hskip 1em $\rho_\text{isum}(\mathfrak{s})$ =
$\mathrm{i}$($0$,
$\sqrt{\frac{1}{3}}$,
$\sqrt{\frac{1}{3}}$,
$\sqrt{\frac{1}{3}}$;\ \ 
$-\frac{1}{3}c^{1}_{36}
+\frac{1}{3}c^{5}_{36}
$,
$\frac{1}{3}c^{1}_{36}
$,
$-\frac{1}{3}c^{5}_{36}
$;\ \ 
$-\frac{1}{3}c^{5}_{36}
$,
$-\frac{1}{3}c^{1}_{36}
+\frac{1}{3}c^{5}_{36}
$;\ \ 
$\frac{1}{3}c^{1}_{36}
$)

Pass. 

 \ \color{black}

\noindent 24: (dims,levels) = $(4;9
)$,
irreps = $4_{9,2}^{5,0}$,
pord$(\rho_\text{isum}(\mathfrak{t})) = 9$,

\vskip 0.7ex
\hangindent=5.5em \hangafter=1
{\white .}\hskip 1em $\rho_\text{isum}(\mathfrak{t})$ =
 $( 0,
\frac{2}{9},
\frac{5}{9},
\frac{8}{9} )
$,

\vskip 0.7ex
\hangindent=5.5em \hangafter=1
{\white .}\hskip 1em $\rho_\text{isum}(\mathfrak{s})$ =
($0$,
$-\sqrt{\frac{1}{3}}$,
$-\sqrt{\frac{1}{3}}$,
$-\sqrt{\frac{1}{3}}$;
$\frac{1}{3} c_9^4 $,
$\frac{1}{3}c^{2}_{9}
$,
$\frac{1}{3}c^{1}_{9}
$;
$\frac{1}{3}c^{1}_{9}
$,
$\frac{1}{3} c_9^4 $;
$\frac{1}{3}c^{2}_{9}
$)

Fail:
cnd($\rho(\mathfrak s)_\mathrm{ndeg}$) = 36 does not divide
 ord($\rho(\mathfrak t)$)=9. Prop. B.4 (2)

 \ \color{black}

\noindent 25: (dims,levels) = $(4;10
)$,
irreps = $4_{5,1}^{1}
\hskip -1.5pt \otimes \hskip -1.5pt
1_{2}^{1,0}$,
pord$(\rho_\text{isum}(\mathfrak{t})) = 5$,

\vskip 0.7ex
\hangindent=5.5em \hangafter=1
{\white .}\hskip 1em $\rho_\text{isum}(\mathfrak{t})$ =
 $( \frac{1}{10},
\frac{3}{10},
\frac{7}{10},
\frac{9}{10} )
$,

\vskip 0.7ex
\hangindent=5.5em \hangafter=1
{\white .}\hskip 1em $\rho_\text{isum}(\mathfrak{s})$ =
$\mathrm{i}$($-\frac{1}{5}c^{1}_{20}
+\frac{1}{5}c^{3}_{20}
$,
$\frac{2}{5}c^{2}_{15}
+\frac{1}{5}c^{3}_{15}
$,
$-\frac{1}{5}+\frac{2}{5}c^{1}_{15}
-\frac{1}{5}c^{3}_{15}
$,
$\frac{1}{5}c^{1}_{20}
+\frac{1}{5}c^{3}_{20}
$;\ \ 
$\frac{1}{5}c^{1}_{20}
+\frac{1}{5}c^{3}_{20}
$,
$-\frac{1}{5}c^{1}_{20}
+\frac{1}{5}c^{3}_{20}
$,
$\frac{1}{5}-\frac{2}{5}c^{1}_{15}
+\frac{1}{5}c^{3}_{15}
$;\ \ 
$-\frac{1}{5}c^{1}_{20}
-\frac{1}{5}c^{3}_{20}
$,
$\frac{2}{5}c^{2}_{15}
+\frac{1}{5}c^{3}_{15}
$;\ \ 
$\frac{1}{5}c^{1}_{20}
-\frac{1}{5}c^{3}_{20}
$)

Fail:
cnd($\rho(\mathfrak s)_\mathrm{ndeg}$) = 60 does not divide
 ord($\rho(\mathfrak t)$)=10. Prop. B.4 (2)

 \ \color{black}

 \color{blue}

\noindent 26: (dims,levels) = $(4;10
)$,
irreps = $4_{5,2}^{1}
\hskip -1.5pt \otimes \hskip -1.5pt
1_{2}^{1,0}$,
pord$(\rho_\text{isum}(\mathfrak{t})) = 5$,

\vskip 0.7ex
\hangindent=5.5em \hangafter=1
{\white .}\hskip 1em $\rho_\text{isum}(\mathfrak{t})$ =
 $( \frac{1}{10},
\frac{3}{10},
\frac{7}{10},
\frac{9}{10} )
$,

\vskip 0.7ex
\hangindent=5.5em \hangafter=1
{\white .}\hskip 1em $\rho_\text{isum}(\mathfrak{s})$ =
($\sqrt{\frac{1}{5}}$,
$\frac{-5+\sqrt{5}}{10}$,
$-\frac{5+\sqrt{5}}{10}$,
$\sqrt{\frac{1}{5}}$;
$-\sqrt{\frac{1}{5}}$,
$\sqrt{\frac{1}{5}}$,
$\frac{5+\sqrt{5}}{10}$;
$-\sqrt{\frac{1}{5}}$,
$\frac{5-\sqrt{5}}{10}$;
$\sqrt{\frac{1}{5}}$)

Pass. 

 \ \color{black}

\noindent 27: (dims,levels) = $(4;10
)$,
irreps = $2_{5}^{1}
\hskip -1.5pt \otimes \hskip -1.5pt
2_{2}^{1,0}$,
pord$(\rho_\text{isum}(\mathfrak{t})) = 10$,

\vskip 0.7ex
\hangindent=5.5em \hangafter=1
{\white .}\hskip 1em $\rho_\text{isum}(\mathfrak{t})$ =
 $( \frac{1}{5},
\frac{4}{5},
\frac{3}{10},
\frac{7}{10} )
$,

\vskip 0.7ex
\hangindent=5.5em \hangafter=1
{\white .}\hskip 1em $\rho_\text{isum}(\mathfrak{s})$ =
$\mathrm{i}$($\frac{1}{2\sqrt{5}}c^{3}_{20}
$,
$\frac{1}{2\sqrt{5}}c^{1}_{20}
$,
$\frac{3}{2\sqrt{15}}c^{1}_{20}
$,
$\frac{3}{2\sqrt{15}}c^{3}_{20}
$;\ \ 
$-\frac{1}{2\sqrt{5}}c^{3}_{20}
$,
$-\frac{3}{2\sqrt{15}}c^{3}_{20}
$,
$\frac{3}{2\sqrt{15}}c^{1}_{20}
$;\ \ 
$\frac{1}{2\sqrt{5}}c^{3}_{20}
$,
$-\frac{1}{2\sqrt{5}}c^{1}_{20}
$;\ \ 
$-\frac{1}{2\sqrt{5}}c^{3}_{20}
$)

Fail:
cnd($\rho(\mathfrak s)_\mathrm{ndeg}$) = 60 does not divide
 ord($\rho(\mathfrak t)$)=10. Prop. B.4 (2)

 \ \color{black}

\noindent 28: (dims,levels) = $(4;10
)$,
irreps = $2_{5}^{2}
\hskip -1.5pt \otimes \hskip -1.5pt
2_{2}^{1,0}$,
pord$(\rho_\text{isum}(\mathfrak{t})) = 10$,

\vskip 0.7ex
\hangindent=5.5em \hangafter=1
{\white .}\hskip 1em $\rho_\text{isum}(\mathfrak{t})$ =
 $( \frac{2}{5},
\frac{3}{5},
\frac{1}{10},
\frac{9}{10} )
$,

\vskip 0.7ex
\hangindent=5.5em \hangafter=1
{\white .}\hskip 1em $\rho_\text{isum}(\mathfrak{s})$ =
$\mathrm{i}$($\frac{1}{2\sqrt{5}}c^{1}_{20}
$,
$\frac{1}{2\sqrt{5}}c^{3}_{20}
$,
$\frac{3}{2\sqrt{15}}c^{3}_{20}
$,
$\frac{3}{2\sqrt{15}}c^{1}_{20}
$;\ \ 
$-\frac{1}{2\sqrt{5}}c^{1}_{20}
$,
$-\frac{3}{2\sqrt{15}}c^{1}_{20}
$,
$\frac{3}{2\sqrt{15}}c^{3}_{20}
$;\ \ 
$\frac{1}{2\sqrt{5}}c^{1}_{20}
$,
$-\frac{1}{2\sqrt{5}}c^{3}_{20}
$;\ \ 
$-\frac{1}{2\sqrt{5}}c^{1}_{20}
$)

Fail:
cnd($\rho(\mathfrak s)_\mathrm{ndeg}$) = 60 does not divide
 ord($\rho(\mathfrak t)$)=10. Prop. B.4 (2)

 \ \color{black}

\noindent 29: (dims,levels) = $(4;14
)$,
irreps = $4_{7}^{1}
\hskip -1.5pt \otimes \hskip -1.5pt
1_{2}^{1,0}$,
pord$(\rho_\text{isum}(\mathfrak{t})) = 7$,

\vskip 0.7ex
\hangindent=5.5em \hangafter=1
{\white .}\hskip 1em $\rho_\text{isum}(\mathfrak{t})$ =
 $( \frac{1}{2},
\frac{1}{14},
\frac{9}{14},
\frac{11}{14} )
$,

\vskip 0.7ex
\hangindent=5.5em \hangafter=1
{\white .}\hskip 1em $\rho_\text{isum}(\mathfrak{s})$ =
$\mathrm{i}$($\sqrt{\frac{1}{7}}$,
$\sqrt{\frac{2}{7}}$,
$\sqrt{\frac{2}{7}}$,
$\sqrt{\frac{2}{7}}$;\ \ 
$\frac{1}{\sqrt{7}}c^{1}_{7}
$,
$-\frac{1}{\sqrt{7}\mathrm{i}}s^{5}_{28}
$,
$\frac{1}{\sqrt{7}}c^{2}_{7}
$;\ \ 
$\frac{1}{\sqrt{7}}c^{2}_{7}
$,
$\frac{1}{\sqrt{7}}c^{1}_{7}
$;\ \ 
$-\frac{1}{\sqrt{7}\mathrm{i}}s^{5}_{28}
$)

Fail:
cnd($\rho(\mathfrak s)_\mathrm{ndeg}$) = 56 does not divide
 ord($\rho(\mathfrak t)$)=14. Prop. B.4 (2)

 \ \color{black}

\noindent 30: (dims,levels) = $(4;14
)$,
irreps = $4_{7}^{3}
\hskip -1.5pt \otimes \hskip -1.5pt
1_{2}^{1,0}$,
pord$(\rho_\text{isum}(\mathfrak{t})) = 7$,

\vskip 0.7ex
\hangindent=5.5em \hangafter=1
{\white .}\hskip 1em $\rho_\text{isum}(\mathfrak{t})$ =
 $( \frac{1}{2},
\frac{3}{14},
\frac{5}{14},
\frac{13}{14} )
$,

\vskip 0.7ex
\hangindent=5.5em \hangafter=1
{\white .}\hskip 1em $\rho_\text{isum}(\mathfrak{s})$ =
$\mathrm{i}$($-\sqrt{\frac{1}{7}}$,
$\sqrt{\frac{2}{7}}$,
$\sqrt{\frac{2}{7}}$,
$\sqrt{\frac{2}{7}}$;\ \ 
$\frac{1}{\sqrt{7}\mathrm{i}}s^{5}_{28}
$,
$-\frac{1}{\sqrt{7}}c^{1}_{7}
$,
$-\frac{1}{\sqrt{7}}c^{2}_{7}
$;\ \ 
$-\frac{1}{\sqrt{7}}c^{2}_{7}
$,
$\frac{1}{\sqrt{7}\mathrm{i}}s^{5}_{28}
$;\ \ 
$-\frac{1}{\sqrt{7}}c^{1}_{7}
$)

Fail:
cnd($\rho(\mathfrak s)_\mathrm{ndeg}$) = 56 does not divide
 ord($\rho(\mathfrak t)$)=14. Prop. B.4 (2)

 \ \color{black}

\noindent 31: (dims,levels) = $(4;15
)$,
irreps = $4_{5,1}^{1}
\hskip -1.5pt \otimes \hskip -1.5pt
1_{3}^{1,0}$,
pord$(\rho_\text{isum}(\mathfrak{t})) = 5$,

\vskip 0.7ex
\hangindent=5.5em \hangafter=1
{\white .}\hskip 1em $\rho_\text{isum}(\mathfrak{t})$ =
 $( \frac{2}{15},
\frac{8}{15},
\frac{11}{15},
\frac{14}{15} )
$,

\vskip 0.7ex
\hangindent=5.5em \hangafter=1
{\white .}\hskip 1em $\rho_\text{isum}(\mathfrak{s})$ =
$\mathrm{i}$($-\frac{1}{5}c^{1}_{20}
-\frac{1}{5}c^{3}_{20}
$,
$\frac{1}{5}c^{1}_{20}
-\frac{1}{5}c^{3}_{20}
$,
$-\frac{1}{5}+\frac{2}{5}c^{1}_{15}
-\frac{1}{5}c^{3}_{15}
$,
$\frac{2}{5}c^{2}_{15}
+\frac{1}{5}c^{3}_{15}
$;\ \ 
$\frac{1}{5}c^{1}_{20}
+\frac{1}{5}c^{3}_{20}
$,
$-\frac{2}{5}c^{2}_{15}
-\frac{1}{5}c^{3}_{15}
$,
$-\frac{1}{5}+\frac{2}{5}c^{1}_{15}
-\frac{1}{5}c^{3}_{15}
$;\ \ 
$-\frac{1}{5}c^{1}_{20}
+\frac{1}{5}c^{3}_{20}
$,
$\frac{1}{5}c^{1}_{20}
+\frac{1}{5}c^{3}_{20}
$;\ \ 
$\frac{1}{5}c^{1}_{20}
-\frac{1}{5}c^{3}_{20}
$)

Fail:
cnd($\rho(\mathfrak s)_\mathrm{ndeg}$) = 60 does not divide
 ord($\rho(\mathfrak t)$)=15. Prop. B.4 (2)

 \ \color{black}

 \color{blue}

\noindent 32: (dims,levels) = $(4;15
)$,
irreps = $4_{5,2}^{1}
\hskip -1.5pt \otimes \hskip -1.5pt
1_{3}^{1,0}$,
pord$(\rho_\text{isum}(\mathfrak{t})) = 5$,

\vskip 0.7ex
\hangindent=5.5em \hangafter=1
{\white .}\hskip 1em $\rho_\text{isum}(\mathfrak{t})$ =
 $( \frac{2}{15},
\frac{8}{15},
\frac{11}{15},
\frac{14}{15} )
$,

\vskip 0.7ex
\hangindent=5.5em \hangafter=1
{\white .}\hskip 1em $\rho_\text{isum}(\mathfrak{s})$ =
($\sqrt{\frac{1}{5}}$,
$\sqrt{\frac{1}{5}}$,
$-\frac{5+\sqrt{5}}{10}$,
$\frac{-5+\sqrt{5}}{10}$;
$\sqrt{\frac{1}{5}}$,
$\frac{5-\sqrt{5}}{10}$,
$\frac{5+\sqrt{5}}{10}$;
$-\sqrt{\frac{1}{5}}$,
$\sqrt{\frac{1}{5}}$;
$-\sqrt{\frac{1}{5}}$)

Pass. 

 \ \color{black}

\noindent 33: (dims,levels) = $(4;15
)$,
irreps = $2_{5}^{1}
\hskip -1.5pt \otimes \hskip -1.5pt
2_{3}^{1,0}$,
pord$(\rho_\text{isum}(\mathfrak{t})) = 15$,

\vskip 0.7ex
\hangindent=5.5em \hangafter=1
{\white .}\hskip 1em $\rho_\text{isum}(\mathfrak{t})$ =
 $( \frac{1}{5},
\frac{4}{5},
\frac{2}{15},
\frac{8}{15} )
$,

\vskip 0.7ex
\hangindent=5.5em \hangafter=1
{\white .}\hskip 1em $\rho_\text{isum}(\mathfrak{s})$ =
($-\frac{1}{\sqrt{15}}c^{3}_{20}
$,
$\frac{1}{\sqrt{15}}c^{1}_{20}
$,
$\frac{2}{\sqrt{30}}c^{1}_{20}
$,
$-\frac{2}{\sqrt{30}}c^{3}_{20}
$;
$\frac{1}{\sqrt{15}}c^{3}_{20}
$,
$\frac{2}{\sqrt{30}}c^{3}_{20}
$,
$\frac{2}{\sqrt{30}}c^{1}_{20}
$;
$-\frac{1}{\sqrt{15}}c^{3}_{20}
$,
$-\frac{1}{\sqrt{15}}c^{1}_{20}
$;
$\frac{1}{\sqrt{15}}c^{3}_{20}
$)

Fail:
cnd($\rho(\mathfrak s)_\mathrm{ndeg}$) = 120 does not divide
 ord($\rho(\mathfrak t)$)=15. Prop. B.4 (2)

 \ \color{black}

\noindent 34: (dims,levels) = $(4;15
)$,
irreps = $2_{5}^{2}
\hskip -1.5pt \otimes \hskip -1.5pt
2_{3}^{1,0}$,
pord$(\rho_\text{isum}(\mathfrak{t})) = 15$,

\vskip 0.7ex
\hangindent=5.5em \hangafter=1
{\white .}\hskip 1em $\rho_\text{isum}(\mathfrak{t})$ =
 $( \frac{2}{5},
\frac{3}{5},
\frac{11}{15},
\frac{14}{15} )
$,

\vskip 0.7ex
\hangindent=5.5em \hangafter=1
{\white .}\hskip 1em $\rho_\text{isum}(\mathfrak{s})$ =
($-\frac{1}{\sqrt{15}}c^{1}_{20}
$,
$-\frac{1}{\sqrt{15}}c^{3}_{20}
$,
$\frac{2}{\sqrt{30}}c^{1}_{20}
$,
$-\frac{2}{\sqrt{30}}c^{3}_{20}
$;
$\frac{1}{\sqrt{15}}c^{1}_{20}
$,
$\frac{2}{\sqrt{30}}c^{3}_{20}
$,
$\frac{2}{\sqrt{30}}c^{1}_{20}
$;
$\frac{1}{\sqrt{15}}c^{1}_{20}
$,
$-\frac{1}{\sqrt{15}}c^{3}_{20}
$;
$-\frac{1}{\sqrt{15}}c^{1}_{20}
$)

Fail:
cnd($\rho(\mathfrak s)_\mathrm{ndeg}$) = 120 does not divide
 ord($\rho(\mathfrak t)$)=15. Prop. B.4 (2)

 \ \color{black}

 \color{blue}

\noindent 35: (dims,levels) = $(4;18
)$,
irreps = $4_{9,1}^{2,0}
\hskip -1.5pt \otimes \hskip -1.5pt
1_{2}^{1,0}$,
pord$(\rho_\text{isum}(\mathfrak{t})) = 9$,

\vskip 0.7ex
\hangindent=5.5em \hangafter=1
{\white .}\hskip 1em $\rho_\text{isum}(\mathfrak{t})$ =
 $( \frac{1}{2},
\frac{1}{18},
\frac{7}{18},
\frac{13}{18} )
$,

\vskip 0.7ex
\hangindent=5.5em \hangafter=1
{\white .}\hskip 1em $\rho_\text{isum}(\mathfrak{s})$ =
$\mathrm{i}$($0$,
$\sqrt{\frac{1}{3}}$,
$\sqrt{\frac{1}{3}}$,
$\sqrt{\frac{1}{3}}$;\ \ 
$\frac{1}{3}c^{5}_{36}
$,
$\frac{1}{3}c^{1}_{36}
-\frac{1}{3}c^{5}_{36}
$,
$-\frac{1}{3}c^{1}_{36}
$;\ \ 
$-\frac{1}{3}c^{1}_{36}
$,
$\frac{1}{3}c^{5}_{36}
$;\ \ 
$\frac{1}{3}c^{1}_{36}
-\frac{1}{3}c^{5}_{36}
$)

Pass. 

 \ \color{black}

\noindent 36: (dims,levels) = $(4;18
)$,
irreps = $4_{9,2}^{5,0}
\hskip -1.5pt \otimes \hskip -1.5pt
1_{2}^{1,0}$,
pord$(\rho_\text{isum}(\mathfrak{t})) = 9$,

\vskip 0.7ex
\hangindent=5.5em \hangafter=1
{\white .}\hskip 1em $\rho_\text{isum}(\mathfrak{t})$ =
 $( \frac{1}{2},
\frac{1}{18},
\frac{7}{18},
\frac{13}{18} )
$,

\vskip 0.7ex
\hangindent=5.5em \hangafter=1
{\white .}\hskip 1em $\rho_\text{isum}(\mathfrak{s})$ =
($0$,
$-\sqrt{\frac{1}{3}}$,
$-\sqrt{\frac{1}{3}}$,
$-\sqrt{\frac{1}{3}}$;
$-\frac{1}{3}c^{1}_{9}
$,
$-\frac{1}{3} c_9^4 $,
$-\frac{1}{3}c^{2}_{9}
$;
$-\frac{1}{3}c^{2}_{9}
$,
$-\frac{1}{3}c^{1}_{9}
$;
$-\frac{1}{3} c_9^4 $)

Fail:
cnd($\rho(\mathfrak s)_\mathrm{ndeg}$) = 36 does not divide
 ord($\rho(\mathfrak t)$)=18. Prop. B.4 (2)

 \ \color{black}

 \color{blue}

\noindent 37: (dims,levels) = $(4;18
)$,
irreps = $4_{9,1}^{1,0}
\hskip -1.5pt \otimes \hskip -1.5pt
1_{2}^{1,0}$,
pord$(\rho_\text{isum}(\mathfrak{t})) = 9$,

\vskip 0.7ex
\hangindent=5.5em \hangafter=1
{\white .}\hskip 1em $\rho_\text{isum}(\mathfrak{t})$ =
 $( \frac{1}{2},
\frac{5}{18},
\frac{11}{18},
\frac{17}{18} )
$,

\vskip 0.7ex
\hangindent=5.5em \hangafter=1
{\white .}\hskip 1em $\rho_\text{isum}(\mathfrak{s})$ =
$\mathrm{i}$($0$,
$\sqrt{\frac{1}{3}}$,
$\sqrt{\frac{1}{3}}$,
$\sqrt{\frac{1}{3}}$;\ \ 
$-\frac{1}{3}c^{1}_{36}
+\frac{1}{3}c^{5}_{36}
$,
$-\frac{1}{3}c^{5}_{36}
$,
$\frac{1}{3}c^{1}_{36}
$;\ \ 
$\frac{1}{3}c^{1}_{36}
$,
$-\frac{1}{3}c^{1}_{36}
+\frac{1}{3}c^{5}_{36}
$;\ \ 
$-\frac{1}{3}c^{5}_{36}
$)

Pass. 

 \ \color{black}

\noindent 38: (dims,levels) = $(4;18
)$,
irreps = $4_{9,2}^{1,0}
\hskip -1.5pt \otimes \hskip -1.5pt
1_{2}^{1,0}$,
pord$(\rho_\text{isum}(\mathfrak{t})) = 9$,

\vskip 0.7ex
\hangindent=5.5em \hangafter=1
{\white .}\hskip 1em $\rho_\text{isum}(\mathfrak{t})$ =
 $( \frac{1}{2},
\frac{5}{18},
\frac{11}{18},
\frac{17}{18} )
$,

\vskip 0.7ex
\hangindent=5.5em \hangafter=1
{\white .}\hskip 1em $\rho_\text{isum}(\mathfrak{s})$ =
($0$,
$-\sqrt{\frac{1}{3}}$,
$-\sqrt{\frac{1}{3}}$,
$-\sqrt{\frac{1}{3}}$;
$-\frac{1}{3} c_9^4 $,
$-\frac{1}{3}c^{1}_{9}
$,
$-\frac{1}{3}c^{2}_{9}
$;
$-\frac{1}{3}c^{2}_{9}
$,
$-\frac{1}{3} c_9^4 $;
$-\frac{1}{3}c^{1}_{9}
$)

Fail:
cnd($\rho(\mathfrak s)_\mathrm{ndeg}$) = 36 does not divide
 ord($\rho(\mathfrak t)$)=18. Prop. B.4 (2)

 \ \color{black}

\noindent 39: (dims,levels) = $(4;20
)$,
irreps = $2_{5}^{1}
\hskip -1.5pt \otimes \hskip -1.5pt
2_{4}^{1,0}$,
pord$(\rho_\text{isum}(\mathfrak{t})) = 10$,

\vskip 0.7ex
\hangindent=5.5em \hangafter=1
{\white .}\hskip 1em $\rho_\text{isum}(\mathfrak{t})$ =
 $( \frac{1}{20},
\frac{9}{20},
\frac{11}{20},
\frac{19}{20} )
$,

\vskip 0.7ex
\hangindent=5.5em \hangafter=1
{\white .}\hskip 1em $\rho_\text{isum}(\mathfrak{s})$ =
($\frac{1}{2\sqrt{5}}c^{3}_{20}
$,
$\frac{1}{2\sqrt{5}}c^{1}_{20}
$,
$-\frac{3}{2\sqrt{15}}c^{3}_{20}
$,
$\frac{3}{2\sqrt{15}}c^{1}_{20}
$;
$-\frac{1}{2\sqrt{5}}c^{3}_{20}
$,
$-\frac{3}{2\sqrt{15}}c^{1}_{20}
$,
$-\frac{3}{2\sqrt{15}}c^{3}_{20}
$;
$-\frac{1}{2\sqrt{5}}c^{3}_{20}
$,
$\frac{1}{2\sqrt{5}}c^{1}_{20}
$;
$\frac{1}{2\sqrt{5}}c^{3}_{20}
$)

Fail:
cnd($\rho(\mathfrak s)_\mathrm{ndeg}$) = 60 does not divide
 ord($\rho(\mathfrak t)$)=20. Prop. B.4 (2)

 \ \color{black}

\noindent 40: (dims,levels) = $(4;20
)$,
irreps = $4_{5,1}^{1}
\hskip -1.5pt \otimes \hskip -1.5pt
1_{4}^{1,0}$,
pord$(\rho_\text{isum}(\mathfrak{t})) = 5$,

\vskip 0.7ex
\hangindent=5.5em \hangafter=1
{\white .}\hskip 1em $\rho_\text{isum}(\mathfrak{t})$ =
 $( \frac{1}{20},
\frac{9}{20},
\frac{13}{20},
\frac{17}{20} )
$,

\vskip 0.7ex
\hangindent=5.5em \hangafter=1
{\white .}\hskip 1em $\rho_\text{isum}(\mathfrak{s})$ =
($\frac{1}{5}c^{1}_{20}
+\frac{1}{5}c^{3}_{20}
$,
$\frac{1}{5}c^{1}_{20}
-\frac{1}{5}c^{3}_{20}
$,
$-\frac{1}{5}+\frac{2}{5}c^{1}_{15}
-\frac{1}{5}c^{3}_{15}
$,
$\frac{2}{5}c^{2}_{15}
+\frac{1}{5}c^{3}_{15}
$;
$-\frac{1}{5}c^{1}_{20}
-\frac{1}{5}c^{3}_{20}
$,
$\frac{2}{5}c^{2}_{15}
+\frac{1}{5}c^{3}_{15}
$,
$\frac{1}{5}-\frac{2}{5}c^{1}_{15}
+\frac{1}{5}c^{3}_{15}
$;
$\frac{1}{5}c^{1}_{20}
-\frac{1}{5}c^{3}_{20}
$,
$-\frac{1}{5}c^{1}_{20}
-\frac{1}{5}c^{3}_{20}
$;
$-\frac{1}{5}c^{1}_{20}
+\frac{1}{5}c^{3}_{20}
$)

Fail:
cnd($\rho(\mathfrak s)_\mathrm{ndeg}$) = 60 does not divide
 ord($\rho(\mathfrak t)$)=20. Prop. B.4 (2)

 \ \color{black}

 \color{blue}

\noindent 41: (dims,levels) = $(4;20
)$,
irreps = $4_{5,2}^{1}
\hskip -1.5pt \otimes \hskip -1.5pt
1_{4}^{1,0}$,
pord$(\rho_\text{isum}(\mathfrak{t})) = 5$,

\vskip 0.7ex
\hangindent=5.5em \hangafter=1
{\white .}\hskip 1em $\rho_\text{isum}(\mathfrak{t})$ =
 $( \frac{1}{20},
\frac{9}{20},
\frac{13}{20},
\frac{17}{20} )
$,

\vskip 0.7ex
\hangindent=5.5em \hangafter=1
{\white .}\hskip 1em $\rho_\text{isum}(\mathfrak{s})$ =
$\mathrm{i}$($\sqrt{\frac{1}{5}}$,
$\sqrt{\frac{1}{5}}$,
$-\frac{5+\sqrt{5}}{10}$,
$\frac{-5+\sqrt{5}}{10}$;\ \ 
$\sqrt{\frac{1}{5}}$,
$\frac{5-\sqrt{5}}{10}$,
$\frac{5+\sqrt{5}}{10}$;\ \ 
$-\sqrt{\frac{1}{5}}$,
$\sqrt{\frac{1}{5}}$;\ \ 
$-\sqrt{\frac{1}{5}}$)

Pass. 

 \ \color{black}

\noindent 42: (dims,levels) = $(4;20
)$,
irreps = $2_{5}^{2}
\hskip -1.5pt \otimes \hskip -1.5pt
2_{4}^{1,0}$,
pord$(\rho_\text{isum}(\mathfrak{t})) = 10$,

\vskip 0.7ex
\hangindent=5.5em \hangafter=1
{\white .}\hskip 1em $\rho_\text{isum}(\mathfrak{t})$ =
 $( \frac{3}{20},
\frac{7}{20},
\frac{13}{20},
\frac{17}{20} )
$,

\vskip 0.7ex
\hangindent=5.5em \hangafter=1
{\white .}\hskip 1em $\rho_\text{isum}(\mathfrak{s})$ =
($\frac{1}{2\sqrt{5}}c^{1}_{20}
$,
$\frac{1}{2\sqrt{5}}c^{3}_{20}
$,
$\frac{3}{2\sqrt{15}}c^{1}_{20}
$,
$-\frac{3}{2\sqrt{15}}c^{3}_{20}
$;
$-\frac{1}{2\sqrt{5}}c^{1}_{20}
$,
$\frac{3}{2\sqrt{15}}c^{3}_{20}
$,
$\frac{3}{2\sqrt{15}}c^{1}_{20}
$;
$-\frac{1}{2\sqrt{5}}c^{1}_{20}
$,
$\frac{1}{2\sqrt{5}}c^{3}_{20}
$;
$\frac{1}{2\sqrt{5}}c^{1}_{20}
$)

Fail:
cnd($\rho(\mathfrak s)_\mathrm{ndeg}$) = 60 does not divide
 ord($\rho(\mathfrak t)$)=20. Prop. B.4 (2)

 \ \color{black}

\noindent 43: (dims,levels) = $(4;21
)$,
irreps = $4_{7}^{3}
\hskip -1.5pt \otimes \hskip -1.5pt
1_{3}^{1,0}$,
pord$(\rho_\text{isum}(\mathfrak{t})) = 7$,

\vskip 0.7ex
\hangindent=5.5em \hangafter=1
{\white .}\hskip 1em $\rho_\text{isum}(\mathfrak{t})$ =
 $( \frac{1}{3},
\frac{1}{21},
\frac{4}{21},
\frac{16}{21} )
$,

\vskip 0.7ex
\hangindent=5.5em \hangafter=1
{\white .}\hskip 1em $\rho_\text{isum}(\mathfrak{s})$ =
$\mathrm{i}$($\sqrt{\frac{1}{7}}$,
$\sqrt{\frac{2}{7}}$,
$\sqrt{\frac{2}{7}}$,
$\sqrt{\frac{2}{7}}$;\ \ 
$-\frac{1}{\sqrt{7}\mathrm{i}}s^{5}_{28}
$,
$\frac{1}{\sqrt{7}}c^{1}_{7}
$,
$\frac{1}{\sqrt{7}}c^{2}_{7}
$;\ \ 
$\frac{1}{\sqrt{7}}c^{2}_{7}
$,
$-\frac{1}{\sqrt{7}\mathrm{i}}s^{5}_{28}
$;\ \ 
$\frac{1}{\sqrt{7}}c^{1}_{7}
$)

Fail:
cnd($\rho(\mathfrak s)_\mathrm{ndeg}$) = 56 does not divide
 ord($\rho(\mathfrak t)$)=21. Prop. B.4 (2)

 \ \color{black}

\noindent 44: (dims,levels) = $(4;21
)$,
irreps = $4_{7}^{1}
\hskip -1.5pt \otimes \hskip -1.5pt
1_{3}^{1,0}$,
pord$(\rho_\text{isum}(\mathfrak{t})) = 7$,

\vskip 0.7ex
\hangindent=5.5em \hangafter=1
{\white .}\hskip 1em $\rho_\text{isum}(\mathfrak{t})$ =
 $( \frac{1}{3},
\frac{10}{21},
\frac{13}{21},
\frac{19}{21} )
$,

\vskip 0.7ex
\hangindent=5.5em \hangafter=1
{\white .}\hskip 1em $\rho_\text{isum}(\mathfrak{s})$ =
$\mathrm{i}$($-\sqrt{\frac{1}{7}}$,
$\sqrt{\frac{2}{7}}$,
$\sqrt{\frac{2}{7}}$,
$\sqrt{\frac{2}{7}}$;\ \ 
$-\frac{1}{\sqrt{7}}c^{2}_{7}
$,
$-\frac{1}{\sqrt{7}}c^{1}_{7}
$,
$\frac{1}{\sqrt{7}\mathrm{i}}s^{5}_{28}
$;\ \ 
$\frac{1}{\sqrt{7}\mathrm{i}}s^{5}_{28}
$,
$-\frac{1}{\sqrt{7}}c^{2}_{7}
$;\ \ 
$-\frac{1}{\sqrt{7}}c^{1}_{7}
$)

Fail:
cnd($\rho(\mathfrak s)_\mathrm{ndeg}$) = 56 does not divide
 ord($\rho(\mathfrak t)$)=21. Prop. B.4 (2)

 \ \color{black}

\noindent 45: (dims,levels) = $(4;24
)$,
irreps = $2_{8}^{1,0}
\hskip -1.5pt \otimes \hskip -1.5pt
2_{3}^{1,0}$,
pord$(\rho_\text{isum}(\mathfrak{t})) = 12$,

\vskip 0.7ex
\hangindent=5.5em \hangafter=1
{\white .}\hskip 1em $\rho_\text{isum}(\mathfrak{t})$ =
 $( \frac{1}{8},
\frac{3}{8},
\frac{11}{24},
\frac{17}{24} )
$,

\vskip 0.7ex
\hangindent=5.5em \hangafter=1
{\white .}\hskip 1em $\rho_\text{isum}(\mathfrak{s})$ =
$\mathrm{i}$($\sqrt{\frac{1}{6}}$,
$\sqrt{\frac{1}{6}}$,
$\sqrt{\frac{1}{3}}$,
$\sqrt{\frac{1}{3}}$;\ \ 
$-\sqrt{\frac{1}{6}}$,
$\sqrt{\frac{1}{3}}$,
$-\sqrt{\frac{1}{3}}$;\ \ 
$-\sqrt{\frac{1}{6}}$,
$-\sqrt{\frac{1}{6}}$;\ \ 
$\sqrt{\frac{1}{6}}$)

Fail:
$\sigma(\rho(\mathfrak s)_\mathrm{ndeg}) \neq
 (\rho(\mathfrak t)^a \rho(\mathfrak s) \rho(\mathfrak t)^b
 \rho(\mathfrak s) \rho(\mathfrak t)^a)_\mathrm{ndeg}$,
 $\sigma = a$ = 5. Prop. B.5 (3) eqn. (B.25)

 \ \color{black}

\noindent 46: (dims,levels) = $(4;28
)$,
irreps = $4_{7}^{3}
\hskip -1.5pt \otimes \hskip -1.5pt
1_{4}^{1,0}$,
pord$(\rho_\text{isum}(\mathfrak{t})) = 7$,

\vskip 0.7ex
\hangindent=5.5em \hangafter=1
{\white .}\hskip 1em $\rho_\text{isum}(\mathfrak{t})$ =
 $( \frac{1}{4},
\frac{3}{28},
\frac{19}{28},
\frac{27}{28} )
$,

\vskip 0.7ex
\hangindent=5.5em \hangafter=1
{\white .}\hskip 1em $\rho_\text{isum}(\mathfrak{s})$ =
($-\sqrt{\frac{1}{7}}$,
$\sqrt{\frac{2}{7}}$,
$\sqrt{\frac{2}{7}}$,
$\sqrt{\frac{2}{7}}$;
$-\frac{1}{\sqrt{7}}c^{2}_{7}
$,
$\frac{1}{\sqrt{7}\mathrm{i}}s^{5}_{28}
$,
$-\frac{1}{\sqrt{7}}c^{1}_{7}
$;
$-\frac{1}{\sqrt{7}}c^{1}_{7}
$,
$-\frac{1}{\sqrt{7}}c^{2}_{7}
$;
$\frac{1}{\sqrt{7}\mathrm{i}}s^{5}_{28}
$)

Fail:
cnd($\rho(\mathfrak s)_\mathrm{ndeg}$) = 56 does not divide
 ord($\rho(\mathfrak t)$)=28. Prop. B.4 (2)

 \ \color{black}

\noindent 47: (dims,levels) = $(4;28
)$,
irreps = $4_{7}^{1}
\hskip -1.5pt \otimes \hskip -1.5pt
1_{4}^{1,0}$,
pord$(\rho_\text{isum}(\mathfrak{t})) = 7$,

\vskip 0.7ex
\hangindent=5.5em \hangafter=1
{\white .}\hskip 1em $\rho_\text{isum}(\mathfrak{t})$ =
 $( \frac{1}{4},
\frac{11}{28},
\frac{15}{28},
\frac{23}{28} )
$,

\vskip 0.7ex
\hangindent=5.5em \hangafter=1
{\white .}\hskip 1em $\rho_\text{isum}(\mathfrak{s})$ =
($\sqrt{\frac{1}{7}}$,
$\sqrt{\frac{2}{7}}$,
$\sqrt{\frac{2}{7}}$,
$\sqrt{\frac{2}{7}}$;
$\frac{1}{\sqrt{7}}c^{2}_{7}
$,
$\frac{1}{\sqrt{7}}c^{1}_{7}
$,
$-\frac{1}{\sqrt{7}\mathrm{i}}s^{5}_{28}
$;
$-\frac{1}{\sqrt{7}\mathrm{i}}s^{5}_{28}
$,
$\frac{1}{\sqrt{7}}c^{2}_{7}
$;
$\frac{1}{\sqrt{7}}c^{1}_{7}
$)

Fail:
cnd($\rho(\mathfrak s)_\mathrm{ndeg}$) = 56 does not divide
 ord($\rho(\mathfrak t)$)=28. Prop. B.4 (2)

 \ \color{black}

\noindent 48: (dims,levels) = $(4;30
)$,
irreps = $4_{5,1}^{1}
\hskip -1.5pt \otimes \hskip -1.5pt
1_{3}^{1,0}
\hskip -1.5pt \otimes \hskip -1.5pt
1_{2}^{1,0}$,
pord$(\rho_\text{isum}(\mathfrak{t})) = 5$,

\vskip 0.7ex
\hangindent=5.5em \hangafter=1
{\white .}\hskip 1em $\rho_\text{isum}(\mathfrak{t})$ =
 $( \frac{1}{30},
\frac{7}{30},
\frac{13}{30},
\frac{19}{30} )
$,

\vskip 0.7ex
\hangindent=5.5em \hangafter=1
{\white .}\hskip 1em $\rho_\text{isum}(\mathfrak{s})$ =
$\mathrm{i}$($-\frac{1}{5}c^{1}_{20}
-\frac{1}{5}c^{3}_{20}
$,
$\frac{2}{5}c^{2}_{15}
+\frac{1}{5}c^{3}_{15}
$,
$-\frac{1}{5}+\frac{2}{5}c^{1}_{15}
-\frac{1}{5}c^{3}_{15}
$,
$\frac{1}{5}c^{1}_{20}
-\frac{1}{5}c^{3}_{20}
$;\ \ 
$\frac{1}{5}c^{1}_{20}
-\frac{1}{5}c^{3}_{20}
$,
$\frac{1}{5}c^{1}_{20}
+\frac{1}{5}c^{3}_{20}
$,
$-\frac{1}{5}+\frac{2}{5}c^{1}_{15}
-\frac{1}{5}c^{3}_{15}
$;\ \ 
$-\frac{1}{5}c^{1}_{20}
+\frac{1}{5}c^{3}_{20}
$,
$-\frac{2}{5}c^{2}_{15}
-\frac{1}{5}c^{3}_{15}
$;\ \ 
$\frac{1}{5}c^{1}_{20}
+\frac{1}{5}c^{3}_{20}
$)

Fail:
cnd($\rho(\mathfrak s)_\mathrm{ndeg}$) = 60 does not divide
 ord($\rho(\mathfrak t)$)=30. Prop. B.4 (2)

 \ \color{black}

 \color{blue}

\noindent 49: (dims,levels) = $(4;30
)$,
irreps = $4_{5,2}^{1}
\hskip -1.5pt \otimes \hskip -1.5pt
1_{3}^{1,0}
\hskip -1.5pt \otimes \hskip -1.5pt
1_{2}^{1,0}$,
pord$(\rho_\text{isum}(\mathfrak{t})) = 5$,

\vskip 0.7ex
\hangindent=5.5em \hangafter=1
{\white .}\hskip 1em $\rho_\text{isum}(\mathfrak{t})$ =
 $( \frac{1}{30},
\frac{7}{30},
\frac{13}{30},
\frac{19}{30} )
$,

\vskip 0.7ex
\hangindent=5.5em \hangafter=1
{\white .}\hskip 1em $\rho_\text{isum}(\mathfrak{s})$ =
($-\sqrt{\frac{1}{5}}$,
$\frac{-5+\sqrt{5}}{10}$,
$-\frac{5+\sqrt{5}}{10}$,
$\sqrt{\frac{1}{5}}$;
$\sqrt{\frac{1}{5}}$,
$-\sqrt{\frac{1}{5}}$,
$-\frac{5+\sqrt{5}}{10}$;
$\sqrt{\frac{1}{5}}$,
$\frac{-5+\sqrt{5}}{10}$;
$-\sqrt{\frac{1}{5}}$)

Pass. 

 \ \color{black}

\noindent 50: (dims,levels) = $(4;30
)$,
irreps = $2_{5}^{2}
\hskip -1.5pt \otimes \hskip -1.5pt
2_{3}^{1,0}
\hskip -1.5pt \otimes \hskip -1.5pt
1_{2}^{1,0}$,
pord$(\rho_\text{isum}(\mathfrak{t})) = 15$,

\vskip 0.7ex
\hangindent=5.5em \hangafter=1
{\white .}\hskip 1em $\rho_\text{isum}(\mathfrak{t})$ =
 $( \frac{1}{10},
\frac{9}{10},
\frac{7}{30},
\frac{13}{30} )
$,

\vskip 0.7ex
\hangindent=5.5em \hangafter=1
{\white .}\hskip 1em $\rho_\text{isum}(\mathfrak{s})$ =
($-\frac{1}{\sqrt{15}}c^{1}_{20}
$,
$-\frac{1}{\sqrt{15}}c^{3}_{20}
$,
$-\frac{2}{\sqrt{30}}c^{3}_{20}
$,
$\frac{2}{\sqrt{30}}c^{1}_{20}
$;
$\frac{1}{\sqrt{15}}c^{1}_{20}
$,
$\frac{2}{\sqrt{30}}c^{1}_{20}
$,
$\frac{2}{\sqrt{30}}c^{3}_{20}
$;
$-\frac{1}{\sqrt{15}}c^{1}_{20}
$,
$-\frac{1}{\sqrt{15}}c^{3}_{20}
$;
$\frac{1}{\sqrt{15}}c^{1}_{20}
$)

Fail:
cnd($\rho(\mathfrak s)_\mathrm{ndeg}$) = 120 does not divide
 ord($\rho(\mathfrak t)$)=30. Prop. B.4 (2)

 \ \color{black}

\noindent 51: (dims,levels) = $(4;30
)$,
irreps = $2_{5}^{1}
\hskip -1.5pt \otimes \hskip -1.5pt
2_{2}^{1,0}
\hskip -1.5pt \otimes \hskip -1.5pt
1_{3}^{1,0}$,
pord$(\rho_\text{isum}(\mathfrak{t})) = 10$,

\vskip 0.7ex
\hangindent=5.5em \hangafter=1
{\white .}\hskip 1em $\rho_\text{isum}(\mathfrak{t})$ =
 $( \frac{2}{15},
\frac{8}{15},
\frac{1}{30},
\frac{19}{30} )
$,

\vskip 0.7ex
\hangindent=5.5em \hangafter=1
{\white .}\hskip 1em $\rho_\text{isum}(\mathfrak{s})$ =
$\mathrm{i}$($-\frac{1}{2\sqrt{5}}c^{3}_{20}
$,
$\frac{1}{2\sqrt{5}}c^{1}_{20}
$,
$\frac{3}{2\sqrt{15}}c^{1}_{20}
$,
$\frac{3}{2\sqrt{15}}c^{3}_{20}
$;\ \ 
$\frac{1}{2\sqrt{5}}c^{3}_{20}
$,
$\frac{3}{2\sqrt{15}}c^{3}_{20}
$,
$-\frac{3}{2\sqrt{15}}c^{1}_{20}
$;\ \ 
$-\frac{1}{2\sqrt{5}}c^{3}_{20}
$,
$\frac{1}{2\sqrt{5}}c^{1}_{20}
$;\ \ 
$\frac{1}{2\sqrt{5}}c^{3}_{20}
$)

Fail:
cnd($\rho(\mathfrak s)_\mathrm{ndeg}$) = 60 does not divide
 ord($\rho(\mathfrak t)$)=30. Prop. B.4 (2)

 \ \color{black}

\noindent 52: (dims,levels) = $(4;30
)$,
irreps = $2_{5}^{1}
\hskip -1.5pt \otimes \hskip -1.5pt
2_{3}^{1,0}
\hskip -1.5pt \otimes \hskip -1.5pt
1_{2}^{1,0}$,
pord$(\rho_\text{isum}(\mathfrak{t})) = 15$,

\vskip 0.7ex
\hangindent=5.5em \hangafter=1
{\white .}\hskip 1em $\rho_\text{isum}(\mathfrak{t})$ =
 $( \frac{3}{10},
\frac{7}{10},
\frac{1}{30},
\frac{19}{30} )
$,

\vskip 0.7ex
\hangindent=5.5em \hangafter=1
{\white .}\hskip 1em $\rho_\text{isum}(\mathfrak{s})$ =
($-\frac{1}{\sqrt{15}}c^{3}_{20}
$,
$\frac{1}{\sqrt{15}}c^{1}_{20}
$,
$\frac{2}{\sqrt{30}}c^{1}_{20}
$,
$-\frac{2}{\sqrt{30}}c^{3}_{20}
$;
$\frac{1}{\sqrt{15}}c^{3}_{20}
$,
$\frac{2}{\sqrt{30}}c^{3}_{20}
$,
$\frac{2}{\sqrt{30}}c^{1}_{20}
$;
$-\frac{1}{\sqrt{15}}c^{3}_{20}
$,
$-\frac{1}{\sqrt{15}}c^{1}_{20}
$;
$\frac{1}{\sqrt{15}}c^{3}_{20}
$)

Fail:
cnd($\rho(\mathfrak s)_\mathrm{ndeg}$) = 120 does not divide
 ord($\rho(\mathfrak t)$)=30. Prop. B.4 (2)

 \ \color{black}

\noindent 53: (dims,levels) = $(4;30
)$,
irreps = $2_{5}^{2}
\hskip -1.5pt \otimes \hskip -1.5pt
2_{2}^{1,0}
\hskip -1.5pt \otimes \hskip -1.5pt
1_{3}^{1,0}$,
pord$(\rho_\text{isum}(\mathfrak{t})) = 10$,

\vskip 0.7ex
\hangindent=5.5em \hangafter=1
{\white .}\hskip 1em $\rho_\text{isum}(\mathfrak{t})$ =
 $( \frac{11}{15},
\frac{14}{15},
\frac{7}{30},
\frac{13}{30} )
$,

\vskip 0.7ex
\hangindent=5.5em \hangafter=1
{\white .}\hskip 1em $\rho_\text{isum}(\mathfrak{s})$ =
$\mathrm{i}$($\frac{1}{2\sqrt{5}}c^{1}_{20}
$,
$\frac{1}{2\sqrt{5}}c^{3}_{20}
$,
$\frac{3}{2\sqrt{15}}c^{1}_{20}
$,
$\frac{3}{2\sqrt{15}}c^{3}_{20}
$;\ \ 
$-\frac{1}{2\sqrt{5}}c^{1}_{20}
$,
$\frac{3}{2\sqrt{15}}c^{3}_{20}
$,
$-\frac{3}{2\sqrt{15}}c^{1}_{20}
$;\ \ 
$-\frac{1}{2\sqrt{5}}c^{1}_{20}
$,
$-\frac{1}{2\sqrt{5}}c^{3}_{20}
$;\ \ 
$\frac{1}{2\sqrt{5}}c^{1}_{20}
$)

Fail:
cnd($\rho(\mathfrak s)_\mathrm{ndeg}$) = 60 does not divide
 ord($\rho(\mathfrak t)$)=30. Prop. B.4 (2)

 \ \color{black}

\noindent 54: (dims,levels) = $(4;36
)$,
irreps = $4_{9,2}^{1,0}
\hskip -1.5pt \otimes \hskip -1.5pt
1_{4}^{1,0}$,
pord$(\rho_\text{isum}(\mathfrak{t})) = 9$,

\vskip 0.7ex
\hangindent=5.5em \hangafter=1
{\white .}\hskip 1em $\rho_\text{isum}(\mathfrak{t})$ =
 $( \frac{1}{4},
\frac{1}{36},
\frac{13}{36},
\frac{25}{36} )
$,

\vskip 0.7ex
\hangindent=5.5em \hangafter=1
{\white .}\hskip 1em $\rho_\text{isum}(\mathfrak{s})$ =
$\mathrm{i}$($0$,
$\sqrt{\frac{1}{3}}$,
$\sqrt{\frac{1}{3}}$,
$\sqrt{\frac{1}{3}}$;\ \ 
$\frac{1}{3} c_9^4 $,
$\frac{1}{3}c^{1}_{9}
$,
$\frac{1}{3}c^{2}_{9}
$;\ \ 
$\frac{1}{3}c^{2}_{9}
$,
$\frac{1}{3} c_9^4 $;\ \ 
$\frac{1}{3}c^{1}_{9}
$)

Fail:
$\sigma(\rho(\mathfrak s)_\mathrm{ndeg}) \neq
 (\rho(\mathfrak t)^a \rho(\mathfrak s) \rho(\mathfrak t)^b
 \rho(\mathfrak s) \rho(\mathfrak t)^a)_\mathrm{ndeg}$,
 $\sigma = a$ = 5. Prop. B.5 (3) eqn. (B.25)

 \ \color{black}

 \color{blue}

\noindent 55: (dims,levels) = $(4;36
)$,
irreps = $4_{9,1}^{1,0}
\hskip -1.5pt \otimes \hskip -1.5pt
1_{4}^{1,0}$,
pord$(\rho_\text{isum}(\mathfrak{t})) = 9$,

\vskip 0.7ex
\hangindent=5.5em \hangafter=1
{\white .}\hskip 1em $\rho_\text{isum}(\mathfrak{t})$ =
 $( \frac{1}{4},
\frac{1}{36},
\frac{13}{36},
\frac{25}{36} )
$,

\vskip 0.7ex
\hangindent=5.5em \hangafter=1
{\white .}\hskip 1em $\rho_\text{isum}(\mathfrak{s})$ =
($0$,
$-\sqrt{\frac{1}{3}}$,
$-\sqrt{\frac{1}{3}}$,
$-\sqrt{\frac{1}{3}}$;
$-\frac{1}{3}c^{1}_{36}
+\frac{1}{3}c^{5}_{36}
$,
$-\frac{1}{3}c^{5}_{36}
$,
$\frac{1}{3}c^{1}_{36}
$;
$\frac{1}{3}c^{1}_{36}
$,
$-\frac{1}{3}c^{1}_{36}
+\frac{1}{3}c^{5}_{36}
$;
$-\frac{1}{3}c^{5}_{36}
$)

Pass. 

 \ \color{black}

\noindent 56: (dims,levels) = $(4;36
)$,
irreps = $4_{9,2}^{5,0}
\hskip -1.5pt \otimes \hskip -1.5pt
1_{4}^{1,0}$,
pord$(\rho_\text{isum}(\mathfrak{t})) = 9$,

\vskip 0.7ex
\hangindent=5.5em \hangafter=1
{\white .}\hskip 1em $\rho_\text{isum}(\mathfrak{t})$ =
 $( \frac{1}{4},
\frac{5}{36},
\frac{17}{36},
\frac{29}{36} )
$,

\vskip 0.7ex
\hangindent=5.5em \hangafter=1
{\white .}\hskip 1em $\rho_\text{isum}(\mathfrak{s})$ =
$\mathrm{i}$($0$,
$\sqrt{\frac{1}{3}}$,
$\sqrt{\frac{1}{3}}$,
$\sqrt{\frac{1}{3}}$;\ \ 
$\frac{1}{3}c^{2}_{9}
$,
$\frac{1}{3}c^{1}_{9}
$,
$\frac{1}{3} c_9^4 $;\ \ 
$\frac{1}{3} c_9^4 $,
$\frac{1}{3}c^{2}_{9}
$;\ \ 
$\frac{1}{3}c^{1}_{9}
$)

Fail:
$\sigma(\rho(\mathfrak s)_\mathrm{ndeg}) \neq
 (\rho(\mathfrak t)^a \rho(\mathfrak s) \rho(\mathfrak t)^b
 \rho(\mathfrak s) \rho(\mathfrak t)^a)_\mathrm{ndeg}$,
 $\sigma = a$ = 5. Prop. B.5 (3) eqn. (B.25)

 \ \color{black}

 \color{blue}

\noindent 57: (dims,levels) = $(4;36
)$,
irreps = $4_{9,1}^{2,0}
\hskip -1.5pt \otimes \hskip -1.5pt
1_{4}^{1,0}$,
pord$(\rho_\text{isum}(\mathfrak{t})) = 9$,

\vskip 0.7ex
\hangindent=5.5em \hangafter=1
{\white .}\hskip 1em $\rho_\text{isum}(\mathfrak{t})$ =
 $( \frac{1}{4},
\frac{5}{36},
\frac{17}{36},
\frac{29}{36} )
$,

\vskip 0.7ex
\hangindent=5.5em \hangafter=1
{\white .}\hskip 1em $\rho_\text{isum}(\mathfrak{s})$ =
($0$,
$-\sqrt{\frac{1}{3}}$,
$-\sqrt{\frac{1}{3}}$,
$-\sqrt{\frac{1}{3}}$;
$-\frac{1}{3}c^{1}_{36}
$,
$\frac{1}{3}c^{5}_{36}
$,
$\frac{1}{3}c^{1}_{36}
-\frac{1}{3}c^{5}_{36}
$;
$\frac{1}{3}c^{1}_{36}
-\frac{1}{3}c^{5}_{36}
$,
$-\frac{1}{3}c^{1}_{36}
$;
$\frac{1}{3}c^{5}_{36}
$)

Pass. 

 \ \color{black}

 \color{blue}

\noindent 58: (dims,levels) = $(4;40
)$,
irreps = $2_{8}^{1,0}
\hskip -1.5pt \otimes \hskip -1.5pt
2_{5}^{1}$,
pord$(\rho_\text{isum}(\mathfrak{t})) = 20$,

\vskip 0.7ex
\hangindent=5.5em \hangafter=1
{\white .}\hskip 1em $\rho_\text{isum}(\mathfrak{t})$ =
 $( \frac{7}{40},
\frac{13}{40},
\frac{23}{40},
\frac{37}{40} )
$,

\vskip 0.7ex
\hangindent=5.5em \hangafter=1
{\white .}\hskip 1em $\rho_\text{isum}(\mathfrak{s})$ =
$\mathrm{i}$($\frac{1}{\sqrt{10}}c^{3}_{20}
$,
$\frac{1}{\sqrt{10}}c^{1}_{20}
$,
$\frac{1}{\sqrt{10}}c^{1}_{20}
$,
$\frac{1}{\sqrt{10}}c^{3}_{20}
$;\ \ 
$\frac{1}{\sqrt{10}}c^{3}_{20}
$,
$-\frac{1}{\sqrt{10}}c^{3}_{20}
$,
$-\frac{1}{\sqrt{10}}c^{1}_{20}
$;\ \ 
$-\frac{1}{\sqrt{10}}c^{3}_{20}
$,
$\frac{1}{\sqrt{10}}c^{1}_{20}
$;\ \ 
$-\frac{1}{\sqrt{10}}c^{3}_{20}
$)

Pass. 

 \ \color{black}

 \color{blue}

\noindent 59: (dims,levels) = $(4;40
)$,
irreps = $2_{8}^{1,0}
\hskip -1.5pt \otimes \hskip -1.5pt
2_{5}^{2}$,
pord$(\rho_\text{isum}(\mathfrak{t})) = 20$,

\vskip 0.7ex
\hangindent=5.5em \hangafter=1
{\white .}\hskip 1em $\rho_\text{isum}(\mathfrak{t})$ =
 $( \frac{21}{40},
\frac{29}{40},
\frac{31}{40},
\frac{39}{40} )
$,

\vskip 0.7ex
\hangindent=5.5em \hangafter=1
{\white .}\hskip 1em $\rho_\text{isum}(\mathfrak{s})$ =
$\mathrm{i}$($\frac{1}{\sqrt{10}}c^{1}_{20}
$,
$\frac{1}{\sqrt{10}}c^{3}_{20}
$,
$\frac{1}{\sqrt{10}}c^{1}_{20}
$,
$\frac{1}{\sqrt{10}}c^{3}_{20}
$;\ \ 
$-\frac{1}{\sqrt{10}}c^{1}_{20}
$,
$\frac{1}{\sqrt{10}}c^{3}_{20}
$,
$-\frac{1}{\sqrt{10}}c^{1}_{20}
$;\ \ 
$-\frac{1}{\sqrt{10}}c^{1}_{20}
$,
$-\frac{1}{\sqrt{10}}c^{3}_{20}
$;\ \ 
$\frac{1}{\sqrt{10}}c^{1}_{20}
$)

Pass. 

 \ \color{black}

\noindent 60: (dims,levels) = $(4;42
)$,
irreps = $4_{7}^{1}
\hskip -1.5pt \otimes \hskip -1.5pt
1_{3}^{1,0}
\hskip -1.5pt \otimes \hskip -1.5pt
1_{2}^{1,0}$,
pord$(\rho_\text{isum}(\mathfrak{t})) = 7$,

\vskip 0.7ex
\hangindent=5.5em \hangafter=1
{\white .}\hskip 1em $\rho_\text{isum}(\mathfrak{t})$ =
 $( \frac{5}{6},
\frac{5}{42},
\frac{17}{42},
\frac{41}{42} )
$,

\vskip 0.7ex
\hangindent=5.5em \hangafter=1
{\white .}\hskip 1em $\rho_\text{isum}(\mathfrak{s})$ =
$\mathrm{i}$($\sqrt{\frac{1}{7}}$,
$\sqrt{\frac{2}{7}}$,
$\sqrt{\frac{2}{7}}$,
$\sqrt{\frac{2}{7}}$;\ \ 
$-\frac{1}{\sqrt{7}\mathrm{i}}s^{5}_{28}
$,
$\frac{1}{\sqrt{7}}c^{2}_{7}
$,
$\frac{1}{\sqrt{7}}c^{1}_{7}
$;\ \ 
$\frac{1}{\sqrt{7}}c^{1}_{7}
$,
$-\frac{1}{\sqrt{7}\mathrm{i}}s^{5}_{28}
$;\ \ 
$\frac{1}{\sqrt{7}}c^{2}_{7}
$)

Fail:
cnd($\rho(\mathfrak s)_\mathrm{ndeg}$) = 56 does not divide
 ord($\rho(\mathfrak t)$)=42. Prop. B.4 (2)

 \ \color{black}

\noindent 61: (dims,levels) = $(4;42
)$,
irreps = $4_{7}^{3}
\hskip -1.5pt \otimes \hskip -1.5pt
1_{3}^{1,0}
\hskip -1.5pt \otimes \hskip -1.5pt
1_{2}^{1,0}$,
pord$(\rho_\text{isum}(\mathfrak{t})) = 7$,

\vskip 0.7ex
\hangindent=5.5em \hangafter=1
{\white .}\hskip 1em $\rho_\text{isum}(\mathfrak{t})$ =
 $( \frac{5}{6},
\frac{11}{42},
\frac{23}{42},
\frac{29}{42} )
$,

\vskip 0.7ex
\hangindent=5.5em \hangafter=1
{\white .}\hskip 1em $\rho_\text{isum}(\mathfrak{s})$ =
$\mathrm{i}$($-\sqrt{\frac{1}{7}}$,
$\sqrt{\frac{2}{7}}$,
$\sqrt{\frac{2}{7}}$,
$\sqrt{\frac{2}{7}}$;\ \ 
$-\frac{1}{\sqrt{7}}c^{1}_{7}
$,
$-\frac{1}{\sqrt{7}}c^{2}_{7}
$,
$\frac{1}{\sqrt{7}\mathrm{i}}s^{5}_{28}
$;\ \ 
$\frac{1}{\sqrt{7}\mathrm{i}}s^{5}_{28}
$,
$-\frac{1}{\sqrt{7}}c^{1}_{7}
$;\ \ 
$-\frac{1}{\sqrt{7}}c^{2}_{7}
$)

Fail:
cnd($\rho(\mathfrak s)_\mathrm{ndeg}$) = 56 does not divide
 ord($\rho(\mathfrak t)$)=42. Prop. B.4 (2)

 \ \color{black}

\noindent 62: (dims,levels) = $(4;60
)$,
irreps = $2_{5}^{1}
\hskip -1.5pt \otimes \hskip -1.5pt
2_{3}^{1,0}
\hskip -1.5pt \otimes \hskip -1.5pt
1_{4}^{1,0}$,
pord$(\rho_\text{isum}(\mathfrak{t})) = 15$,

\vskip 0.7ex
\hangindent=5.5em \hangafter=1
{\white .}\hskip 1em $\rho_\text{isum}(\mathfrak{t})$ =
 $( \frac{1}{20},
\frac{9}{20},
\frac{23}{60},
\frac{47}{60} )
$,

\vskip 0.7ex
\hangindent=5.5em \hangafter=1
{\white .}\hskip 1em $\rho_\text{isum}(\mathfrak{s})$ =
$\mathrm{i}$($\frac{1}{\sqrt{15}}c^{3}_{20}
$,
$\frac{1}{\sqrt{15}}c^{1}_{20}
$,
$-\frac{2}{\sqrt{30}}c^{3}_{20}
$,
$\frac{2}{\sqrt{30}}c^{1}_{20}
$;\ \ 
$-\frac{1}{\sqrt{15}}c^{3}_{20}
$,
$-\frac{2}{\sqrt{30}}c^{1}_{20}
$,
$-\frac{2}{\sqrt{30}}c^{3}_{20}
$;\ \ 
$-\frac{1}{\sqrt{15}}c^{3}_{20}
$,
$\frac{1}{\sqrt{15}}c^{1}_{20}
$;\ \ 
$\frac{1}{\sqrt{15}}c^{3}_{20}
$)

Fail:
cnd($\rho(\mathfrak s)_\mathrm{ndeg}$) = 120 does not divide
 ord($\rho(\mathfrak t)$)=60. Prop. B.4 (2)

 \ \color{black}

 \color{blue}

\noindent 63: (dims,levels) = $(4;60
)$,
irreps = $4_{5,2}^{1}
\hskip -1.5pt \otimes \hskip -1.5pt
1_{4}^{1,0}
\hskip -1.5pt \otimes \hskip -1.5pt
1_{3}^{1,0}$,
pord$(\rho_\text{isum}(\mathfrak{t})) = 5$,

\vskip 0.7ex
\hangindent=5.5em \hangafter=1
{\white .}\hskip 1em $\rho_\text{isum}(\mathfrak{t})$ =
 $( \frac{11}{60},
\frac{23}{60},
\frac{47}{60},
\frac{59}{60} )
$,

\vskip 0.7ex
\hangindent=5.5em \hangafter=1
{\white .}\hskip 1em $\rho_\text{isum}(\mathfrak{s})$ =
$\mathrm{i}$($-\sqrt{\frac{1}{5}}$,
$\frac{-5+\sqrt{5}}{10}$,
$-\frac{5+\sqrt{5}}{10}$,
$\sqrt{\frac{1}{5}}$;\ \ 
$\sqrt{\frac{1}{5}}$,
$-\sqrt{\frac{1}{5}}$,
$-\frac{5+\sqrt{5}}{10}$;\ \ 
$\sqrt{\frac{1}{5}}$,
$\frac{-5+\sqrt{5}}{10}$;\ \ 
$-\sqrt{\frac{1}{5}}$)

Pass. 

 \ \color{black}

\noindent 64: (dims,levels) = $(4;60
)$,
irreps = $4_{5,1}^{1}
\hskip -1.5pt \otimes \hskip -1.5pt
1_{4}^{1,0}
\hskip -1.5pt \otimes \hskip -1.5pt
1_{3}^{1,0}$,
pord$(\rho_\text{isum}(\mathfrak{t})) = 5$,

\vskip 0.7ex
\hangindent=5.5em \hangafter=1
{\white .}\hskip 1em $\rho_\text{isum}(\mathfrak{t})$ =
 $( \frac{11}{60},
\frac{23}{60},
\frac{47}{60},
\frac{59}{60} )
$,

\vskip 0.7ex
\hangindent=5.5em \hangafter=1
{\white .}\hskip 1em $\rho_\text{isum}(\mathfrak{s})$ =
($-\frac{1}{5}c^{1}_{20}
+\frac{1}{5}c^{3}_{20}
$,
$\frac{2}{5}c^{2}_{15}
+\frac{1}{5}c^{3}_{15}
$,
$-\frac{1}{5}+\frac{2}{5}c^{1}_{15}
-\frac{1}{5}c^{3}_{15}
$,
$\frac{1}{5}c^{1}_{20}
+\frac{1}{5}c^{3}_{20}
$;
$\frac{1}{5}c^{1}_{20}
+\frac{1}{5}c^{3}_{20}
$,
$-\frac{1}{5}c^{1}_{20}
+\frac{1}{5}c^{3}_{20}
$,
$\frac{1}{5}-\frac{2}{5}c^{1}_{15}
+\frac{1}{5}c^{3}_{15}
$;
$-\frac{1}{5}c^{1}_{20}
-\frac{1}{5}c^{3}_{20}
$,
$\frac{2}{5}c^{2}_{15}
+\frac{1}{5}c^{3}_{15}
$;
$\frac{1}{5}c^{1}_{20}
-\frac{1}{5}c^{3}_{20}
$)

Fail:
$\sigma(\rho(\mathfrak s)_\mathrm{ndeg}) \neq
 (\rho(\mathfrak t)^a \rho(\mathfrak s) \rho(\mathfrak t)^b
 \rho(\mathfrak s) \rho(\mathfrak t)^a)_\mathrm{ndeg}$,
 $\sigma = a$ = 7. Prop. B.5 (3) eqn. (B.25)

 \ \color{black}

\noindent 65: (dims,levels) = $(4;60
)$,
irreps = $2_{5}^{2}
\hskip -1.5pt \otimes \hskip -1.5pt
2_{4}^{1,0}
\hskip -1.5pt \otimes \hskip -1.5pt
1_{3}^{1,0}$,
pord$(\rho_\text{isum}(\mathfrak{t})) = 10$,

\vskip 0.7ex
\hangindent=5.5em \hangafter=1
{\white .}\hskip 1em $\rho_\text{isum}(\mathfrak{t})$ =
 $( \frac{11}{60},
\frac{29}{60},
\frac{41}{60},
\frac{59}{60} )
$,

\vskip 0.7ex
\hangindent=5.5em \hangafter=1
{\white .}\hskip 1em $\rho_\text{isum}(\mathfrak{s})$ =
($\frac{1}{2\sqrt{5}}c^{1}_{20}
$,
$-\frac{3}{2\sqrt{15}}c^{3}_{20}
$,
$\frac{3}{2\sqrt{15}}c^{1}_{20}
$,
$\frac{1}{2\sqrt{5}}c^{3}_{20}
$;
$\frac{1}{2\sqrt{5}}c^{1}_{20}
$,
$\frac{1}{2\sqrt{5}}c^{3}_{20}
$,
$\frac{3}{2\sqrt{15}}c^{1}_{20}
$;
$-\frac{1}{2\sqrt{5}}c^{1}_{20}
$,
$\frac{3}{2\sqrt{15}}c^{3}_{20}
$;
$-\frac{1}{2\sqrt{5}}c^{1}_{20}
$)

Fail:
$\sigma(\rho(\mathfrak s)_\mathrm{ndeg}) \neq
 (\rho(\mathfrak t)^a \rho(\mathfrak s) \rho(\mathfrak t)^b
 \rho(\mathfrak s) \rho(\mathfrak t)^a)_\mathrm{ndeg}$,
 $\sigma = a$ = 7. Prop. B.5 (3) eqn. (B.25)

 \ \color{black}

\noindent 66: (dims,levels) = $(4;60
)$,
irreps = $2_{5}^{1}
\hskip -1.5pt \otimes \hskip -1.5pt
2_{4}^{1,0}
\hskip -1.5pt \otimes \hskip -1.5pt
1_{3}^{1,0}$,
pord$(\rho_\text{isum}(\mathfrak{t})) = 10$,

\vskip 0.7ex
\hangindent=5.5em \hangafter=1
{\white .}\hskip 1em $\rho_\text{isum}(\mathfrak{t})$ =
 $( \frac{17}{60},
\frac{23}{60},
\frac{47}{60},
\frac{53}{60} )
$,

\vskip 0.7ex
\hangindent=5.5em \hangafter=1
{\white .}\hskip 1em $\rho_\text{isum}(\mathfrak{s})$ =
($\frac{1}{2\sqrt{5}}c^{3}_{20}
$,
$\frac{3}{2\sqrt{15}}c^{1}_{20}
$,
$-\frac{3}{2\sqrt{15}}c^{3}_{20}
$,
$\frac{1}{2\sqrt{5}}c^{1}_{20}
$;
$\frac{1}{2\sqrt{5}}c^{3}_{20}
$,
$\frac{1}{2\sqrt{5}}c^{1}_{20}
$,
$-\frac{3}{2\sqrt{15}}c^{3}_{20}
$;
$-\frac{1}{2\sqrt{5}}c^{3}_{20}
$,
$-\frac{3}{2\sqrt{15}}c^{1}_{20}
$;
$-\frac{1}{2\sqrt{5}}c^{3}_{20}
$)

Fail:
$\sigma(\rho(\mathfrak s)_\mathrm{ndeg}) \neq
 (\rho(\mathfrak t)^a \rho(\mathfrak s) \rho(\mathfrak t)^b
 \rho(\mathfrak s) \rho(\mathfrak t)^a)_\mathrm{ndeg}$,
 $\sigma = a$ = 7. Prop. B.5 (3) eqn. (B.25)

 \ \color{black}

\noindent 67: (dims,levels) = $(4;60
)$,
irreps = $2_{5}^{2}
\hskip -1.5pt \otimes \hskip -1.5pt
2_{3}^{1,0}
\hskip -1.5pt \otimes \hskip -1.5pt
1_{4}^{1,0}$,
pord$(\rho_\text{isum}(\mathfrak{t})) = 15$,

\vskip 0.7ex
\hangindent=5.5em \hangafter=1
{\white .}\hskip 1em $\rho_\text{isum}(\mathfrak{t})$ =
 $( \frac{13}{20},
\frac{17}{20},
\frac{11}{60},
\frac{59}{60} )
$,

\vskip 0.7ex
\hangindent=5.5em \hangafter=1
{\white .}\hskip 1em $\rho_\text{isum}(\mathfrak{s})$ =
$\mathrm{i}$($-\frac{1}{\sqrt{15}}c^{1}_{20}
$,
$\frac{1}{\sqrt{15}}c^{3}_{20}
$,
$-\frac{2}{\sqrt{30}}c^{3}_{20}
$,
$\frac{2}{\sqrt{30}}c^{1}_{20}
$;\ \ 
$\frac{1}{\sqrt{15}}c^{1}_{20}
$,
$-\frac{2}{\sqrt{30}}c^{1}_{20}
$,
$-\frac{2}{\sqrt{30}}c^{3}_{20}
$;\ \ 
$-\frac{1}{\sqrt{15}}c^{1}_{20}
$,
$-\frac{1}{\sqrt{15}}c^{3}_{20}
$;\ \ 
$\frac{1}{\sqrt{15}}c^{1}_{20}
$)

Fail:
cnd($\rho(\mathfrak s)_\mathrm{ndeg}$) = 120 does not divide
 ord($\rho(\mathfrak t)$)=60. Prop. B.4 (2)

 \ \color{black}

\noindent 68: (dims,levels) = $(4;84
)$,
irreps = $4_{7}^{3}
\hskip -1.5pt \otimes \hskip -1.5pt
1_{4}^{1,0}
\hskip -1.5pt \otimes \hskip -1.5pt
1_{3}^{1,0}$,
pord$(\rho_\text{isum}(\mathfrak{t})) = 7$,

\vskip 0.7ex
\hangindent=5.5em \hangafter=1
{\white .}\hskip 1em $\rho_\text{isum}(\mathfrak{t})$ =
 $( \frac{7}{12},
\frac{1}{84},
\frac{25}{84},
\frac{37}{84} )
$,

\vskip 0.7ex
\hangindent=5.5em \hangafter=1
{\white .}\hskip 1em $\rho_\text{isum}(\mathfrak{s})$ =
($-\sqrt{\frac{1}{7}}$,
$\sqrt{\frac{2}{7}}$,
$\sqrt{\frac{2}{7}}$,
$\sqrt{\frac{2}{7}}$;
$-\frac{1}{\sqrt{7}}c^{1}_{7}
$,
$-\frac{1}{\sqrt{7}}c^{2}_{7}
$,
$\frac{1}{\sqrt{7}\mathrm{i}}s^{5}_{28}
$;
$\frac{1}{\sqrt{7}\mathrm{i}}s^{5}_{28}
$,
$-\frac{1}{\sqrt{7}}c^{1}_{7}
$;
$-\frac{1}{\sqrt{7}}c^{2}_{7}
$)

Fail:
cnd($\rho(\mathfrak s)_\mathrm{ndeg}$) = 56 does not divide
 ord($\rho(\mathfrak t)$)=84. Prop. B.4 (2)

 \ \color{black}

\noindent 69: (dims,levels) = $(4;84
)$,
irreps = $4_{7}^{1}
\hskip -1.5pt \otimes \hskip -1.5pt
1_{4}^{1,0}
\hskip -1.5pt \otimes \hskip -1.5pt
1_{3}^{1,0}$,
pord$(\rho_\text{isum}(\mathfrak{t})) = 7$,

\vskip 0.7ex
\hangindent=5.5em \hangafter=1
{\white .}\hskip 1em $\rho_\text{isum}(\mathfrak{t})$ =
 $( \frac{7}{12},
\frac{13}{84},
\frac{61}{84},
\frac{73}{84} )
$,

\vskip 0.7ex
\hangindent=5.5em \hangafter=1
{\white .}\hskip 1em $\rho_\text{isum}(\mathfrak{s})$ =
($\sqrt{\frac{1}{7}}$,
$\sqrt{\frac{2}{7}}$,
$\sqrt{\frac{2}{7}}$,
$\sqrt{\frac{2}{7}}$;
$\frac{1}{\sqrt{7}}c^{1}_{7}
$,
$-\frac{1}{\sqrt{7}\mathrm{i}}s^{5}_{28}
$,
$\frac{1}{\sqrt{7}}c^{2}_{7}
$;
$\frac{1}{\sqrt{7}}c^{2}_{7}
$,
$\frac{1}{\sqrt{7}}c^{1}_{7}
$;
$-\frac{1}{\sqrt{7}\mathrm{i}}s^{5}_{28}
$)

Fail:
cnd($\rho(\mathfrak s)_\mathrm{ndeg}$) = 56 does not divide
 ord($\rho(\mathfrak t)$)=84. Prop. B.4 (2)

 \ \color{black}

 \color{blue}

\noindent 70: (dims,levels) = $(4;120
)$,
irreps = $2_{8}^{1,0}
\hskip -1.5pt \otimes \hskip -1.5pt
2_{5}^{2}
\hskip -1.5pt \otimes \hskip -1.5pt
1_{3}^{1,0}$,
pord$(\rho_\text{isum}(\mathfrak{t})) = 20$,

\vskip 0.7ex
\hangindent=5.5em \hangafter=1
{\white .}\hskip 1em $\rho_\text{isum}(\mathfrak{t})$ =
 $( \frac{7}{120},
\frac{13}{120},
\frac{37}{120},
\frac{103}{120} )
$,

\vskip 0.7ex
\hangindent=5.5em \hangafter=1
{\white .}\hskip 1em $\rho_\text{isum}(\mathfrak{s})$ =
$\mathrm{i}$($-\frac{1}{\sqrt{10}}c^{1}_{20}
$,
$\frac{1}{\sqrt{10}}c^{3}_{20}
$,
$\frac{1}{\sqrt{10}}c^{1}_{20}
$,
$\frac{1}{\sqrt{10}}c^{3}_{20}
$;\ \ 
$-\frac{1}{\sqrt{10}}c^{1}_{20}
$,
$\frac{1}{\sqrt{10}}c^{3}_{20}
$,
$\frac{1}{\sqrt{10}}c^{1}_{20}
$;\ \ 
$\frac{1}{\sqrt{10}}c^{1}_{20}
$,
$-\frac{1}{\sqrt{10}}c^{3}_{20}
$;\ \ 
$\frac{1}{\sqrt{10}}c^{1}_{20}
$)

Pass. 

 \ \color{black}

 \color{blue}

\noindent 71: (dims,levels) = $(4;120
)$,
irreps = $2_{8}^{1,0}
\hskip -1.5pt \otimes \hskip -1.5pt
2_{5}^{1}
\hskip -1.5pt \otimes \hskip -1.5pt
1_{3}^{1,0}$,
pord$(\rho_\text{isum}(\mathfrak{t})) = 20$,

\vskip 0.7ex
\hangindent=5.5em \hangafter=1
{\white .}\hskip 1em $\rho_\text{isum}(\mathfrak{t})$ =
 $( \frac{31}{120},
\frac{61}{120},
\frac{79}{120},
\frac{109}{120} )
$,

\vskip 0.7ex
\hangindent=5.5em \hangafter=1
{\white .}\hskip 1em $\rho_\text{isum}(\mathfrak{s})$ =
$\mathrm{i}$($-\frac{1}{\sqrt{10}}c^{3}_{20}
$,
$\frac{1}{\sqrt{10}}c^{3}_{20}
$,
$\frac{1}{\sqrt{10}}c^{1}_{20}
$,
$\frac{1}{\sqrt{10}}c^{1}_{20}
$;\ \ 
$\frac{1}{\sqrt{10}}c^{3}_{20}
$,
$-\frac{1}{\sqrt{10}}c^{1}_{20}
$,
$\frac{1}{\sqrt{10}}c^{1}_{20}
$;\ \ 
$\frac{1}{\sqrt{10}}c^{3}_{20}
$,
$\frac{1}{\sqrt{10}}c^{3}_{20}
$;\ \ 
$-\frac{1}{\sqrt{10}}c^{3}_{20}
$)

Pass. 

 \ \color{black}

\

\subsection{A list of passing GT orbits}

The above passing representations can be grouped into GT orbits.  The following
list displays one representative representation for each GT orbit.  For details
and notations, see Appendix B.2.

\

\noindent1. (dims;levels) =$(3\oplus
1;5,
1
)$,
irreps = $3_{5}^{1}\oplus
1_{1}^{1}$,
pord$(\rho_\text{isum}(\mathfrak{t})) = 5$,

\vskip 0.7ex
\hangindent=4em \hangafter=1
 $\rho_\text{isum}(\mathfrak{t})$ =
 $( 0,
\frac{1}{5},
\frac{4}{5} )
\oplus
( 0 )
$,

\vskip 0.7ex
\hangindent=4em \hangafter=1
 $\rho_\text{isum}(\mathfrak{s})$ =
($\sqrt{\frac{1}{5}}$,
$-\sqrt{\frac{2}{5}}$,
$-\sqrt{\frac{2}{5}}$;
$-\frac{5+\sqrt{5}}{10}$,
$\frac{5-\sqrt{5}}{10}$;
$-\frac{5+\sqrt{5}}{10}$)
 $\oplus$
($1$)

Resolved. Number of valid $(S,T)$ pairs = 2.

\vskip 2ex

 \noindent2. (dims;levels) =$(4;5
)$,
irreps = $4_{5,2}^{1}$,
pord$(\rho_\text{isum}(\mathfrak{t})) = 5$,

\vskip 0.7ex
\hangindent=4em \hangafter=1
 $\rho_\text{isum}(\mathfrak{t})$ =
 $( \frac{1}{5},
\frac{2}{5},
\frac{3}{5},
\frac{4}{5} )
$,

\vskip 0.7ex
\hangindent=4em \hangafter=1
 $\rho_\text{isum}(\mathfrak{s})$ =
($\sqrt{\frac{1}{5}}$,
$\frac{-5+\sqrt{5}}{10}$,
$-\frac{5+\sqrt{5}}{10}$,
$\sqrt{\frac{1}{5}}$;
$-\sqrt{\frac{1}{5}}$,
$\sqrt{\frac{1}{5}}$,
$\frac{5+\sqrt{5}}{10}$;
$-\sqrt{\frac{1}{5}}$,
$\frac{5-\sqrt{5}}{10}$;
$\sqrt{\frac{1}{5}}$)

Resolved. Number of valid $(S,T)$ pairs = 1.

\vskip 2ex

 \noindent3. (dims;levels) =$(4;9
)$,
irreps = $4_{9,1}^{1,0}$,
pord$(\rho_\text{isum}(\mathfrak{t})) = 9$,

\vskip 0.7ex
\hangindent=4em \hangafter=1
 $\rho_\text{isum}(\mathfrak{t})$ =
 $( 0,
\frac{1}{9},
\frac{4}{9},
\frac{7}{9} )
$,

\vskip 0.7ex
\hangindent=4em \hangafter=1
 $\rho_\text{isum}(\mathfrak{s})$ =
$\mathrm{i}$($0$,
$\sqrt{\frac{1}{3}}$,
$\sqrt{\frac{1}{3}}$,
$\sqrt{\frac{1}{3}}$;\ \ 
$-\frac{1}{3}c^{1}_{36}
$,
$\frac{1}{3}c^{1}_{36}
-\frac{1}{3}c^{5}_{36}
$,
$\frac{1}{3}c^{5}_{36}
$;\ \ 
$\frac{1}{3}c^{5}_{36}
$,
$-\frac{1}{3}c^{1}_{36}
$;\ \ 
$\frac{1}{3}c^{1}_{36}
-\frac{1}{3}c^{5}_{36}
$)

Resolved. Number of valid $(S,T)$ pairs = 1.

\vskip 2ex

 \noindent4. (dims;levels) =$(4;40
)$,
irreps = $2_{8}^{1,6}
\hskip -1.5pt \otimes \hskip -1.5pt
2_{5}^{2}$,
pord$(\rho_\text{isum}(\mathfrak{t})) = 20$,

\vskip 0.7ex
\hangindent=4em \hangafter=1
 $\rho_\text{isum}(\mathfrak{t})$ =
 $( \frac{1}{40},
\frac{9}{40},
\frac{11}{40},
\frac{19}{40} )
$,

\vskip 0.7ex
\hangindent=4em \hangafter=1
 $\rho_\text{isum}(\mathfrak{s})$ =
$\mathrm{i}$($-\frac{1}{\sqrt{10}}c^{1}_{20}
$,
$\frac{1}{\sqrt{10}}c^{3}_{20}
$,
$\frac{1}{\sqrt{10}}c^{1}_{20}
$,
$\frac{1}{\sqrt{10}}c^{3}_{20}
$;\ \ 
$\frac{1}{\sqrt{10}}c^{1}_{20}
$,
$-\frac{1}{\sqrt{10}}c^{3}_{20}
$,
$\frac{1}{\sqrt{10}}c^{1}_{20}
$;\ \ 
$\frac{1}{\sqrt{10}}c^{1}_{20}
$,
$\frac{1}{\sqrt{10}}c^{3}_{20}
$;\ \ 
$-\frac{1}{\sqrt{10}}c^{1}_{20}
$)

Resolved. Number of valid $(S,T)$ pairs = 2.

\vskip 2ex

\

\subsection{A list of rank-4 $S,T$ matrices from resolved representations}

From the representative representation in each GT orbit, if it is revolved, we
can compute all the $S,T$ matrices coming from such a representation.  The
computation steps are displayed below.  For details and notations, see Section
1 of this file.

\

\noindent1. (dims;levels) =$(3 , 
1;5,
1
)$,
irreps = $3_{5}^{1}\oplus
1_{1}^{1}$,
pord$(\tilde\rho(\mathfrak{t})) = 5$,

\vskip 0.7ex
\hangindent=4em \hangafter=1
 $\tilde\rho(\mathfrak{t})$ =
 $( 0,
0,
\frac{1}{5},
\frac{4}{5} )
$,

\vskip 0.7ex
\hangindent=4em \hangafter=1
 $\tilde\rho(\mathfrak{s})$ =
($\sqrt{\frac{1}{5}}$,
$0$,
$-\sqrt{\frac{2}{5}}$,
$-\sqrt{\frac{2}{5}}$;
$1$,
$0$,
$0$;
$-\frac{5+\sqrt{5}}{10}$,
$\frac{5-\sqrt{5}}{10}$;
$-\frac{5+\sqrt{5}}{10}$)

 \vskip 1ex \setlength{\leftskip}{2em}

\grey{Try $U_0$ =
$\begin{pmatrix}
1,
& 0 \\ 
0,
& 1 \\ 
\end{pmatrix}
$ $\oplus
\begin{pmatrix}
1 \\ 
\end{pmatrix}
$ $\oplus
\begin{pmatrix}
1 \\ 
\end{pmatrix}
$:}\ \ \ \ \ 
\grey{$U_0\tilde\rho(\mathfrak{s})U_0^\dagger$ =}

\grey{$\begin{pmatrix}
\sqrt{\frac{1}{5}},
& 0,
& -\sqrt{\frac{2}{5}},
& -\sqrt{\frac{2}{5}} \\ 
0,
& 1,
& 0,
& 0 \\ 
-\sqrt{\frac{2}{5}},
& 0,
& -\frac{5+\sqrt{5}}{10},
& \frac{5-\sqrt{5}}{10} \\ 
-\sqrt{\frac{2}{5}},
& 0,
& \frac{5-\sqrt{5}}{10},
& -\frac{5+\sqrt{5}}{10} \\ 
\end{pmatrix}
$}

\grey{Try different $u$'s and signed diagonal matrix $V_\mathrm{sd}$'s:}

 \grey{
\begin{tabular}{|r|l|l|l|l|}
\hline
$3_{5}^{1}\oplus
1_{1}^{1}:\ u$ 
 & 0 & 1 & 2 & 3\\ 
 \hline
$D_\rho$ conditions 
 & 1 & 1 & 1 & 1\\ 
 \hline
$[\rho(\mathfrak{s})\rho(\mathfrak{t})]^3
 = \rho^2(\mathfrak{s}) = \tilde C$ 
 & 0 & 0 & 0 & 0\\ 
 \hline
$\rho(\mathfrak{s})_{iu}\rho^*(\mathfrak{s})_{ju} \in \mathbb{R}$ 
 & 0 & 0 & 0 & 0\\ 
 \hline
$\rho(\mathfrak{s})_{i u} \neq 0$  
 & 1 & 1 & 1 & 1\\ 
 \hline
$\mathrm{cnd}(S)$, $\mathrm{cnd}(\rho(\mathfrak{s}))$ 
 & - & - & - & -\\ 
 \hline
$\mathrm{norm}(D^2)$ factors
 & - & - & - & -\\ 
 \hline
$1/\rho(\mathfrak{s})_{iu} = $ cyc-int 
 & - & - & - & -\\ 
 \hline
norm$(1/\rho(\mathfrak{s})_{iu})$ factors
 & - & - & - & -\\ 
 \hline
$\frac{S_{ij}}{S_{uj}} = $ cyc-int
 & - & - & - & -\\ 
 \hline
$N^{ij}_k \in \mathbb{N}$
 & - & - & - & -\\ 
 \hline
$\exists\ j \text{ that } \frac{S_{ij}}{S_{uj}} \geq 1 $
 & - & - & - & -\\ 
 \hline
FS indicator
 & - & - & - & -\\ 
 \hline
$C = $ perm-mat
 & - & - & - & -\\ 
 \hline
\end{tabular}

Number of valid $(S,T)$ pairs: 0 \vskip 2ex }%grey

\grey{Try $U_0$ =
$\begin{pmatrix}
\sqrt{\frac{1}{2}},
& \sqrt{\frac{1}{2}} \\ 
\sqrt{\frac{1}{2}},
& -\sqrt{\frac{1}{2}} \\ 
\end{pmatrix}
$ $\oplus
\begin{pmatrix}
1 \\ 
\end{pmatrix}
$ $\oplus
\begin{pmatrix}
1 \\ 
\end{pmatrix}
$:}\ \ \ \ \ 
\grey{$U_0\tilde\rho(\mathfrak{s})U_0^\dagger$ =}

\grey{$\begin{pmatrix}
\frac{5+\sqrt{5}}{10},
& \frac{-5+\sqrt{5}}{10},
& -\sqrt{\frac{1}{5}},
& -\sqrt{\frac{1}{5}} \\ 
\frac{-5+\sqrt{5}}{10},
& \frac{5+\sqrt{5}}{10},
& -\sqrt{\frac{1}{5}},
& -\sqrt{\frac{1}{5}} \\ 
-\sqrt{\frac{1}{5}},
& -\sqrt{\frac{1}{5}},
& -\frac{5+\sqrt{5}}{10},
& \frac{5-\sqrt{5}}{10} \\ 
-\sqrt{\frac{1}{5}},
& -\sqrt{\frac{1}{5}},
& \frac{5-\sqrt{5}}{10},
& -\frac{5+\sqrt{5}}{10} \\ 
\end{pmatrix}
$}

\grey{Try different $u$'s and signed diagonal matrix $V_\mathrm{sd}$'s:}

 \grey{
\begin{tabular}{|r|l|l|l|l|}
\hline
$3_{5}^{1}\oplus
1_{1}^{1}:\ u$ 
 & 0 & 1 & 2 & 3\\ 
 \hline
$D_\rho$ conditions 
 & 0 & 0 & 0 & 0\\ 
 \hline
$[\rho(\mathfrak{s})\rho(\mathfrak{t})]^3
 = \rho^2(\mathfrak{s}) = \tilde C$ 
 & 0 & 0 & 0 & 0\\ 
 \hline
$\rho(\mathfrak{s})_{iu}\rho^*(\mathfrak{s})_{ju} \in \mathbb{R}$ 
 & 0 & 0 & 0 & 0\\ 
 \hline
$\rho(\mathfrak{s})_{i u} \neq 0$  
 & 0 & 0 & 0 & 0\\ 
 \hline
$\mathrm{cnd}(S)$, $\mathrm{cnd}(\rho(\mathfrak{s}))$ 
 & 0 & 0 & 0 & 0\\ 
 \hline
$\mathrm{norm}(D^2)$ factors
 & 0 & 0 & 0 & 0\\ 
 \hline
$1/\rho(\mathfrak{s})_{iu} = $ cyc-int 
 & 0 & 0 & 0 & 0\\ 
 \hline
norm$(1/\rho(\mathfrak{s})_{iu})$ factors
 & 0 & 0 & 0 & 0\\ 
 \hline
$\frac{S_{ij}}{S_{uj}} = $ cyc-int
 & 0 & 0 & 0 & 0\\ 
 \hline
$N^{ij}_k \in \mathbb{N}$
 & 0 & 0 & 0 & 0\\ 
 \hline
$\exists\ j \text{ that } \frac{S_{ij}}{S_{uj}} \geq 1 $
 & 0 & 0 & 0 & 0\\ 
 \hline
FS indicator
 & 0 & 0 & 0 & 0\\ 
 \hline
$C = $ perm-mat
 & 0 & 0 & 0 & 0\\ 
 \hline
\end{tabular}

Number of valid $(S,T)$ pairs: 2 \vskip 2ex }%grey

\grey{Try $U_0$ =
$\begin{pmatrix}
\sqrt{\frac{1}{2}},
& -\sqrt{\frac{1}{2}} \\ 
-\sqrt{\frac{1}{2}},
& -\sqrt{\frac{1}{2}} \\ 
\end{pmatrix}
$ $\oplus
\begin{pmatrix}
1 \\ 
\end{pmatrix}
$ $\oplus
\begin{pmatrix}
1 \\ 
\end{pmatrix}
$:}\ \ \ \ \ 
\grey{$U_0\tilde\rho(\mathfrak{s})U_0^\dagger$ =}

\grey{$\begin{pmatrix}
\frac{5+\sqrt{5}}{10},
& \frac{5-\sqrt{5}}{10},
& -\sqrt{\frac{1}{5}},
& -\sqrt{\frac{1}{5}} \\ 
\frac{5-\sqrt{5}}{10},
& \frac{5+\sqrt{5}}{10},
& \sqrt{\frac{1}{5}},
& \sqrt{\frac{1}{5}} \\ 
-\sqrt{\frac{1}{5}},
& \sqrt{\frac{1}{5}},
& -\frac{5+\sqrt{5}}{10},
& \frac{5-\sqrt{5}}{10} \\ 
-\sqrt{\frac{1}{5}},
& \sqrt{\frac{1}{5}},
& \frac{5-\sqrt{5}}{10},
& -\frac{5+\sqrt{5}}{10} \\ 
\end{pmatrix}
$}

\grey{Try different $u$'s and signed diagonal matrix $V_\mathrm{sd}$'s:}

 \grey{
\begin{tabular}{|r|l|l|l|l|}
\hline
$3_{5}^{1}\oplus
1_{1}^{1}:\ u$ 
 & 0 & 1 & 2 & 3\\ 
 \hline
$D_\rho$ conditions 
 & 0 & 0 & 0 & 0\\ 
 \hline
$[\rho(\mathfrak{s})\rho(\mathfrak{t})]^3
 = \rho^2(\mathfrak{s}) = \tilde C$ 
 & 0 & 0 & 0 & 0\\ 
 \hline
$\rho(\mathfrak{s})_{iu}\rho^*(\mathfrak{s})_{ju} \in \mathbb{R}$ 
 & 0 & 0 & 0 & 0\\ 
 \hline
$\rho(\mathfrak{s})_{i u} \neq 0$  
 & 0 & 0 & 0 & 0\\ 
 \hline
$\mathrm{cnd}(S)$, $\mathrm{cnd}(\rho(\mathfrak{s}))$ 
 & 0 & 0 & 0 & 0\\ 
 \hline
$\mathrm{norm}(D^2)$ factors
 & 0 & 0 & 0 & 0\\ 
 \hline
$1/\rho(\mathfrak{s})_{iu} = $ cyc-int 
 & 0 & 0 & 0 & 0\\ 
 \hline
norm$(1/\rho(\mathfrak{s})_{iu})$ factors
 & 0 & 0 & 0 & 0\\ 
 \hline
$\frac{S_{ij}}{S_{uj}} = $ cyc-int
 & 0 & 0 & 0 & 0\\ 
 \hline
$N^{ij}_k \in \mathbb{N}$
 & 0 & 0 & 0 & 0\\ 
 \hline
$\exists\ j \text{ that } \frac{S_{ij}}{S_{uj}} \geq 1 $
 & 0 & 0 & 0 & 0\\ 
 \hline
FS indicator
 & 0 & 0 & 0 & 0\\ 
 \hline
$C = $ perm-mat
 & 0 & 0 & 0 & 0\\ 
 \hline
\end{tabular}

Number of valid $(S,T)$ pairs: 2 \vskip 2ex }%grey

Total number of valid $(S,T)$ pairs: 2

 \vskip 4ex

\ \setlength{\leftskip}{0em} 

\noindent2. (dims;levels) =$(4;5
)$,
irreps = $4_{5,2}^{1}$,
pord$(\tilde\rho(\mathfrak{t})) = 5$,

\vskip 0.7ex
\hangindent=4em \hangafter=1
 $\tilde\rho(\mathfrak{t})$ =
 $( \frac{1}{5},
\frac{2}{5},
\frac{3}{5},
\frac{4}{5} )
$,

\vskip 0.7ex
\hangindent=4em \hangafter=1
 $\tilde\rho(\mathfrak{s})$ =
($\sqrt{\frac{1}{5}}$,
$\frac{-5+\sqrt{5}}{10}$,
$-\frac{5+\sqrt{5}}{10}$,
$\sqrt{\frac{1}{5}}$;
$-\sqrt{\frac{1}{5}}$,
$\sqrt{\frac{1}{5}}$,
$\frac{5+\sqrt{5}}{10}$;
$-\sqrt{\frac{1}{5}}$,
$\frac{5-\sqrt{5}}{10}$;
$\sqrt{\frac{1}{5}}$)

 \vskip 1ex \setlength{\leftskip}{2em}

\grey{Try $U_0$ =
$\begin{pmatrix}
1 \\ 
\end{pmatrix}
$ $\oplus
\begin{pmatrix}
1 \\ 
\end{pmatrix}
$ $\oplus
\begin{pmatrix}
1 \\ 
\end{pmatrix}
$ $\oplus
\begin{pmatrix}
1 \\ 
\end{pmatrix}
$:}\ \ \ \ \ 
\grey{$U_0\tilde\rho(\mathfrak{s})U_0^\dagger$ =}

\grey{$\begin{pmatrix}
\sqrt{\frac{1}{5}},
& \frac{-5+\sqrt{5}}{10},
& -\frac{5+\sqrt{5}}{10},
& \sqrt{\frac{1}{5}} \\ 
\frac{-5+\sqrt{5}}{10},
& -\sqrt{\frac{1}{5}},
& \sqrt{\frac{1}{5}},
& \frac{5+\sqrt{5}}{10} \\ 
-\frac{5+\sqrt{5}}{10},
& \sqrt{\frac{1}{5}},
& -\sqrt{\frac{1}{5}},
& \frac{5-\sqrt{5}}{10} \\ 
\sqrt{\frac{1}{5}},
& \frac{5+\sqrt{5}}{10},
& \frac{5-\sqrt{5}}{10},
& \sqrt{\frac{1}{5}} \\ 
\end{pmatrix}
$}

\grey{Try different $u$'s and signed diagonal matrix $V_\mathrm{sd}$'s:}

 \grey{
\begin{tabular}{|r|l|l|l|l|}
\hline
$4_{5,2}^{1}:\ u$ 
 & 0 & 1 & 2 & 3\\ 
 \hline
$D_\rho$ conditions 
 & 0 & 0 & 0 & 0\\ 
 \hline
$[\rho(\mathfrak{s})\rho(\mathfrak{t})]^3
 = \rho^2(\mathfrak{s}) = \tilde C$ 
 & 0 & 0 & 0 & 0\\ 
 \hline
$\rho(\mathfrak{s})_{iu}\rho^*(\mathfrak{s})_{ju} \in \mathbb{R}$ 
 & 0 & 0 & 0 & 0\\ 
 \hline
$\rho(\mathfrak{s})_{i u} \neq 0$  
 & 0 & 0 & 0 & 0\\ 
 \hline
$\mathrm{cnd}(S)$, $\mathrm{cnd}(\rho(\mathfrak{s}))$ 
 & 0 & 0 & 0 & 0\\ 
 \hline
$\mathrm{norm}(D^2)$ factors
 & 0 & 0 & 0 & 0\\ 
 \hline
$1/\rho(\mathfrak{s})_{iu} = $ cyc-int 
 & 0 & 0 & 0 & 0\\ 
 \hline
norm$(1/\rho(\mathfrak{s})_{iu})$ factors
 & 0 & 0 & 0 & 0\\ 
 \hline
$\frac{S_{ij}}{S_{uj}} = $ cyc-int
 & 0 & 0 & 0 & 0\\ 
 \hline
$N^{ij}_k \in \mathbb{N}$
 & 0 & 0 & 0 & 0\\ 
 \hline
$\exists\ j \text{ that } \frac{S_{ij}}{S_{uj}} \geq 1 $
 & 0 & 0 & 0 & 0\\ 
 \hline
FS indicator
 & 0 & 0 & 0 & 0\\ 
 \hline
$C = $ perm-mat
 & 0 & 0 & 0 & 0\\ 
 \hline
\end{tabular}

Number of valid $(S,T)$ pairs: 1 \vskip 2ex }%grey

Total number of valid $(S,T)$ pairs: 1

 \vskip 4ex

\ \setlength{\leftskip}{0em} 

\noindent3. (dims;levels) =$(4;9
)$,
irreps = $4_{9,1}^{1,0}$,
pord$(\tilde\rho(\mathfrak{t})) = 9$,

\vskip 0.7ex
\hangindent=4em \hangafter=1
 $\tilde\rho(\mathfrak{t})$ =
 $( 0,
\frac{1}{9},
\frac{4}{9},
\frac{7}{9} )
$,

\vskip 0.7ex
\hangindent=4em \hangafter=1
 $\tilde\rho(\mathfrak{s})$ =
$\mathrm{i}$($0$,
$\sqrt{\frac{1}{3}}$,
$\sqrt{\frac{1}{3}}$,
$\sqrt{\frac{1}{3}}$;\ \ 
$-\frac{1}{3}c^{1}_{36}
$,
$\frac{1}{3}c^{1}_{36}
-\frac{1}{3}c^{5}_{36}
$,
$\frac{1}{3}c^{5}_{36}
$;\ \ 
$\frac{1}{3}c^{5}_{36}
$,
$-\frac{1}{3}c^{1}_{36}
$;\ \ 
$\frac{1}{3}c^{1}_{36}
-\frac{1}{3}c^{5}_{36}
$)

 \vskip 1ex \setlength{\leftskip}{2em}

\grey{Try $U_0$ =
$\begin{pmatrix}
1 \\ 
\end{pmatrix}
$ $\oplus
\begin{pmatrix}
1 \\ 
\end{pmatrix}
$ $\oplus
\begin{pmatrix}
1 \\ 
\end{pmatrix}
$ $\oplus
\begin{pmatrix}
1 \\ 
\end{pmatrix}
$:}\ \ \ \ \ 
\grey{$U_0\tilde\rho(\mathfrak{s})U_0^\dagger$ =}

\grey{$\begin{pmatrix}
0,
& \sqrt{\frac{1}{3}}\mathrm{i},
& \sqrt{\frac{1}{3}}\mathrm{i},
& \sqrt{\frac{1}{3}}\mathrm{i} \\ 
\sqrt{\frac{1}{3}}\mathrm{i},
& -\frac{1}{3}c^{1}_{36}
\mathrm{i},
& (\frac{1}{3}c^{1}_{36}
-\frac{1}{3}c^{5}_{36}
)\mathrm{i},
& \frac{1}{3}c^{5}_{36}
\mathrm{i} \\ 
\sqrt{\frac{1}{3}}\mathrm{i},
& (\frac{1}{3}c^{1}_{36}
-\frac{1}{3}c^{5}_{36}
)\mathrm{i},
& \frac{1}{3}c^{5}_{36}
\mathrm{i},
& -\frac{1}{3}c^{1}_{36}
\mathrm{i} \\ 
\sqrt{\frac{1}{3}}\mathrm{i},
& \frac{1}{3}c^{5}_{36}
\mathrm{i},
& -\frac{1}{3}c^{1}_{36}
\mathrm{i},
& (\frac{1}{3}c^{1}_{36}
-\frac{1}{3}c^{5}_{36}
)\mathrm{i} \\ 
\end{pmatrix}
$}

\grey{Try different $u$'s and signed diagonal matrix $V_\mathrm{sd}$'s:}

 \grey{
\begin{tabular}{|r|l|l|l|l|}
\hline
$4_{9,1}^{1,0}:\ u$ 
 & 0 & 1 & 2 & 3\\ 
 \hline
$D_\rho$ conditions 
 & 0 & 0 & 0 & 0\\ 
 \hline
$[\rho(\mathfrak{s})\rho(\mathfrak{t})]^3
 = \rho^2(\mathfrak{s}) = \tilde C$ 
 & 0 & 0 & 0 & 0\\ 
 \hline
$\rho(\mathfrak{s})_{iu}\rho^*(\mathfrak{s})_{ju} \in \mathbb{R}$ 
 & 0 & 0 & 0 & 0\\ 
 \hline
$\rho(\mathfrak{s})_{i u} \neq 0$  
 & 2 & 0 & 0 & 0\\ 
 \hline
$\mathrm{cnd}(S)$, $\mathrm{cnd}(\rho(\mathfrak{s}))$ 
 & - & 0 & 0 & 0\\ 
 \hline
$\mathrm{norm}(D^2)$ factors
 & - & 0 & 0 & 0\\ 
 \hline
$1/\rho(\mathfrak{s})_{iu} = $ cyc-int 
 & - & 0 & 0 & 0\\ 
 \hline
norm$(1/\rho(\mathfrak{s})_{iu})$ factors
 & - & 0 & 0 & 0\\ 
 \hline
$\frac{S_{ij}}{S_{uj}} = $ cyc-int
 & - & 0 & 0 & 0\\ 
 \hline
$N^{ij}_k \in \mathbb{N}$
 & - & 0 & 0 & 0\\ 
 \hline
$\exists\ j \text{ that } \frac{S_{ij}}{S_{uj}} \geq 1 $
 & - & 0 & 0 & 0\\ 
 \hline
FS indicator
 & - & 0 & 0 & 0\\ 
 \hline
$C = $ perm-mat
 & - & 0 & 0 & 0\\ 
 \hline
\end{tabular}

Number of valid $(S,T)$ pairs: 1 \vskip 2ex }%grey

Total number of valid $(S,T)$ pairs: 1

 \vskip 4ex

\ \setlength{\leftskip}{0em} 

\noindent4. (dims;levels) =$(4;40
)$,
irreps = $2_{8}^{1,6}
\hskip -1.5pt \otimes \hskip -1.5pt
2_{5}^{2}$,
pord$(\tilde\rho(\mathfrak{t})) = 20$,

\vskip 0.7ex
\hangindent=4em \hangafter=1
 $\tilde\rho(\mathfrak{t})$ =
 $( \frac{1}{40},
\frac{9}{40},
\frac{11}{40},
\frac{19}{40} )
$,

\vskip 0.7ex
\hangindent=4em \hangafter=1
 $\tilde\rho(\mathfrak{s})$ =
$\mathrm{i}$($-\frac{1}{\sqrt{10}}c^{1}_{20}
$,
$\frac{1}{\sqrt{10}}c^{3}_{20}
$,
$\frac{1}{\sqrt{10}}c^{1}_{20}
$,
$\frac{1}{\sqrt{10}}c^{3}_{20}
$;\ \ 
$\frac{1}{\sqrt{10}}c^{1}_{20}
$,
$-\frac{1}{\sqrt{10}}c^{3}_{20}
$,
$\frac{1}{\sqrt{10}}c^{1}_{20}
$;\ \ 
$\frac{1}{\sqrt{10}}c^{1}_{20}
$,
$\frac{1}{\sqrt{10}}c^{3}_{20}
$;\ \ 
$-\frac{1}{\sqrt{10}}c^{1}_{20}
$)

 \vskip 1ex \setlength{\leftskip}{2em}

\grey{Try $U_0$ =
$\begin{pmatrix}
1 \\ 
\end{pmatrix}
$ $\oplus
\begin{pmatrix}
1 \\ 
\end{pmatrix}
$ $\oplus
\begin{pmatrix}
1 \\ 
\end{pmatrix}
$ $\oplus
\begin{pmatrix}
1 \\ 
\end{pmatrix}
$:}\ \ \ \ \ 
\grey{$U_0\tilde\rho(\mathfrak{s})U_0^\dagger$ =}

\grey{$\begin{pmatrix}
-\frac{1}{\sqrt{10}}c^{1}_{20}
\mathrm{i},
& \frac{1}{\sqrt{10}}c^{3}_{20}
\mathrm{i},
& \frac{1}{\sqrt{10}}c^{1}_{20}
\mathrm{i},
& \frac{1}{\sqrt{10}}c^{3}_{20}
\mathrm{i} \\ 
\frac{1}{\sqrt{10}}c^{3}_{20}
\mathrm{i},
& \frac{1}{\sqrt{10}}c^{1}_{20}
\mathrm{i},
& -\frac{1}{\sqrt{10}}c^{3}_{20}
\mathrm{i},
& \frac{1}{\sqrt{10}}c^{1}_{20}
\mathrm{i} \\ 
\frac{1}{\sqrt{10}}c^{1}_{20}
\mathrm{i},
& -\frac{1}{\sqrt{10}}c^{3}_{20}
\mathrm{i},
& \frac{1}{\sqrt{10}}c^{1}_{20}
\mathrm{i},
& \frac{1}{\sqrt{10}}c^{3}_{20}
\mathrm{i} \\ 
\frac{1}{\sqrt{10}}c^{3}_{20}
\mathrm{i},
& \frac{1}{\sqrt{10}}c^{1}_{20}
\mathrm{i},
& \frac{1}{\sqrt{10}}c^{3}_{20}
\mathrm{i},
& -\frac{1}{\sqrt{10}}c^{1}_{20}
\mathrm{i} \\ 
\end{pmatrix}
$}

\grey{Try different $u$'s and signed diagonal matrix $V_\mathrm{sd}$'s:}

 \grey{
\begin{tabular}{|r|l|l|l|l|}
\hline
$2_{8}^{1,6}
\hskip -1.5pt \otimes \hskip -1.5pt
2_{5}^{2}:\ u$ 
 & 0 & 1 & 2 & 3\\ 
 \hline
$D_\rho$ conditions 
 & 0 & 0 & 0 & 0\\ 
 \hline
$[\rho(\mathfrak{s})\rho(\mathfrak{t})]^3
 = \rho^2(\mathfrak{s}) = \tilde C$ 
 & 0 & 0 & 0 & 0\\ 
 \hline
$\rho(\mathfrak{s})_{iu}\rho^*(\mathfrak{s})_{ju} \in \mathbb{R}$ 
 & 0 & 0 & 0 & 0\\ 
 \hline
$\rho(\mathfrak{s})_{i u} \neq 0$  
 & 0 & 0 & 0 & 0\\ 
 \hline
$\mathrm{cnd}(S)$, $\mathrm{cnd}(\rho(\mathfrak{s}))$ 
 & 0 & 0 & 0 & 0\\ 
 \hline
$\mathrm{norm}(D^2)$ factors
 & 0 & 0 & 0 & 0\\ 
 \hline
$1/\rho(\mathfrak{s})_{iu} = $ cyc-int 
 & 0 & 0 & 0 & 0\\ 
 \hline
norm$(1/\rho(\mathfrak{s})_{iu})$ factors
 & 0 & 0 & 0 & 0\\ 
 \hline
$\frac{S_{ij}}{S_{uj}} = $ cyc-int
 & 0 & 0 & 0 & 0\\ 
 \hline
$N^{ij}_k \in \mathbb{N}$
 & 0 & 0 & 0 & 0\\ 
 \hline
$\exists\ j \text{ that } \frac{S_{ij}}{S_{uj}} \geq 1 $
 & 0 & 0 & 0 & 0\\ 
 \hline
FS indicator
 & 0 & 0 & 0 & 0\\ 
 \hline
$C = $ perm-mat
 & 0 & 0 & 0 & 0\\ 
 \hline
\end{tabular}

Number of valid $(S,T)$ pairs: 2 \vskip 2ex }%grey

Total number of valid $(S,T)$ pairs: 2

 \vskip 4ex

\ \setlength{\leftskip}{0em}

\

The $S,T$ matrices obtained above (the black or the blue entries below), plus
their Galois conjugations (the grey entries below), form the following list of
rank-4 $S,T$ matrices.  For details and notations, see Section 2 of this file.

\

\noindent1. ind = $(3 , 
1;5,
1
)_{1}^{1}$:\ \ 
$d_i$ = ($1.0$,
$1.618$,
$1.618$,
$2.618$) 

\vskip 0.7ex
\hangindent=3em \hangafter=1
$D^2=$ 13.90 = 
 $\frac{15+5\sqrt{5}}{2}$

\vskip 0.7ex
\hangindent=3em \hangafter=1
$T = ( 0,
\frac{2}{5},
\frac{2}{5},
\frac{4}{5} )
$,

\vskip 0.7ex
\hangindent=3em \hangafter=1
$S$ = ($ 1$,
$ \frac{1+\sqrt{5}}{2}$,
$ \frac{1+\sqrt{5}}{2}$,
$ \frac{3+\sqrt{5}}{2}$;\ \ 
$ -1$,
$ \frac{3+\sqrt{5}}{2}$,
$ -\frac{1+\sqrt{5}}{2}$;\ \ 
$ -1$,
$ -\frac{1+\sqrt{5}}{2}$;\ \ 
$ 1$)

\vskip 1ex 
\color{grey}

\noindent2. ind = $(3 , 
1;5,
1
)_{1}^{4}$:\ \ 
$d_i$ = ($1.0$,
$1.618$,
$1.618$,
$2.618$) 

\vskip 0.7ex
\hangindent=3em \hangafter=1
$D^2=$ 13.90 = 
 $\frac{15+5\sqrt{5}}{2}$

\vskip 0.7ex
\hangindent=3em \hangafter=1
$T = ( 0,
\frac{3}{5},
\frac{3}{5},
\frac{1}{5} )
$,

\vskip 0.7ex
\hangindent=3em \hangafter=1
$S$ = ($ 1$,
$ \frac{1+\sqrt{5}}{2}$,
$ \frac{1+\sqrt{5}}{2}$,
$ \frac{3+\sqrt{5}}{2}$;\ \ 
$ -1$,
$ \frac{3+\sqrt{5}}{2}$,
$ -\frac{1+\sqrt{5}}{2}$;\ \ 
$ -1$,
$ -\frac{1+\sqrt{5}}{2}$;\ \ 
$ 1$)

\vskip 1ex 
\color{grey}

\noindent3. ind = $(3 , 
1;5,
1
)_{1}^{3}$:\ \ 
$d_i$ = ($1.0$,
$0.381$,
$-0.618$,
$-0.618$) 

\vskip 0.7ex
\hangindent=3em \hangafter=1
$D^2=$ 1.909 = 
 $\frac{15-5\sqrt{5}}{2}$

\vskip 0.7ex
\hangindent=3em \hangafter=1
$T = ( 0,
\frac{2}{5},
\frac{1}{5},
\frac{1}{5} )
$,

\vskip 0.7ex
\hangindent=3em \hangafter=1
$S$ = ($ 1$,
$ \frac{3-\sqrt{5}}{2}$,
$ \frac{1-\sqrt{5}}{2}$,
$ \frac{1-\sqrt{5}}{2}$;\ \ 
$ 1$,
$ \frac{-1+\sqrt{5}}{2}$,
$ \frac{-1+\sqrt{5}}{2}$;\ \ 
$ -1$,
$ \frac{3-\sqrt{5}}{2}$;\ \ 
$ -1$)

Not pseudo-unitary. 

\vskip 1ex 
\color{grey}

\noindent4. ind = $(3 , 
1;5,
1
)_{1}^{2}$:\ \ 
$d_i$ = ($1.0$,
$0.381$,
$-0.618$,
$-0.618$) 

\vskip 0.7ex
\hangindent=3em \hangafter=1
$D^2=$ 1.909 = 
 $\frac{15-5\sqrt{5}}{2}$

\vskip 0.7ex
\hangindent=3em \hangafter=1
$T = ( 0,
\frac{3}{5},
\frac{4}{5},
\frac{4}{5} )
$,

\vskip 0.7ex
\hangindent=3em \hangafter=1
$S$ = ($ 1$,
$ \frac{3-\sqrt{5}}{2}$,
$ \frac{1-\sqrt{5}}{2}$,
$ \frac{1-\sqrt{5}}{2}$;\ \ 
$ 1$,
$ \frac{-1+\sqrt{5}}{2}$,
$ \frac{-1+\sqrt{5}}{2}$;\ \ 
$ -1$,
$ \frac{3-\sqrt{5}}{2}$;\ \ 
$ -1$)

Not pseudo-unitary. 

\vskip 1ex 

 \color{black} \vskip 2ex

\noindent5. ind = $(3 , 
1;5,
1
)_{2}^{1}$:\ \ 
$d_i$ = ($1.0$,
$1.618$,
$1.618$,
$2.618$) 

\vskip 0.7ex
\hangindent=3em \hangafter=1
$D^2=$ 13.90 = 
 $\frac{15+5\sqrt{5}}{2}$

\vskip 0.7ex
\hangindent=3em \hangafter=1
$T = ( 0,
\frac{2}{5},
\frac{3}{5},
0 )
$,

\vskip 0.7ex
\hangindent=3em \hangafter=1
$S$ = ($ 1$,
$ \frac{1+\sqrt{5}}{2}$,
$ \frac{1+\sqrt{5}}{2}$,
$ \frac{3+\sqrt{5}}{2}$;\ \ 
$ -1$,
$ \frac{3+\sqrt{5}}{2}$,
$ -\frac{1+\sqrt{5}}{2}$;\ \ 
$ -1$,
$ -\frac{1+\sqrt{5}}{2}$;\ \ 
$ 1$)

\vskip 1ex 
\color{grey}

\noindent6. ind = $(3 , 
1;5,
1
)_{2}^{2}$:\ \ 
$d_i$ = ($1.0$,
$0.381$,
$-0.618$,
$-0.618$) 

\vskip 0.7ex
\hangindent=3em \hangafter=1
$D^2=$ 1.909 = 
 $\frac{15-5\sqrt{5}}{2}$

\vskip 0.7ex
\hangindent=3em \hangafter=1
$T = ( 0,
0,
\frac{1}{5},
\frac{4}{5} )
$,

\vskip 0.7ex
\hangindent=3em \hangafter=1
$S$ = ($ 1$,
$ \frac{3-\sqrt{5}}{2}$,
$ \frac{1-\sqrt{5}}{2}$,
$ \frac{1-\sqrt{5}}{2}$;\ \ 
$ 1$,
$ \frac{-1+\sqrt{5}}{2}$,
$ \frac{-1+\sqrt{5}}{2}$;\ \ 
$ -1$,
$ \frac{3-\sqrt{5}}{2}$;\ \ 
$ -1$)

Not pseudo-unitary. 

\vskip 1ex 

 \color{black} \vskip 2ex 
\color{blue}

\noindent7. ind = $(4;5
)_{1}^{1}$:\ \ 
$d_i$ = ($1.0$,
$1.618$,
$-0.618$,
$-1.0$) 

\vskip 0.7ex
\hangindent=3em \hangafter=1
$D^2=$ 5.0 = 
 $5$

\vskip 0.7ex
\hangindent=3em \hangafter=1
$T = ( 0,
\frac{2}{5},
\frac{1}{5},
\frac{3}{5} )
$,

\vskip 0.7ex
\hangindent=3em \hangafter=1
$S$ = ($ 1$,
$ \frac{1+\sqrt{5}}{2}$,
$ \frac{1-\sqrt{5}}{2}$,
$ -1$;\ \ 
$ -1$,
$ -1$,
$ \frac{-1+\sqrt{5}}{2}$;\ \ 
$ -1$,
$ -\frac{1+\sqrt{5}}{2}$;\ \ 
$ 1$)

Not pseudo-unitary. 

\vskip 1ex 
\color{grey}

\noindent8. ind = $(4;5
)_{1}^{2}$:\ \ 
$d_i$ = ($1.0$,
$1.618$,
$-0.618$,
$-1.0$) 

\vskip 0.7ex
\hangindent=3em \hangafter=1
$D^2=$ 5.0 = 
 $5$

\vskip 0.7ex
\hangindent=3em \hangafter=1
$T = ( 0,
\frac{2}{5},
\frac{4}{5},
\frac{1}{5} )
$,

\vskip 0.7ex
\hangindent=3em \hangafter=1
$S$ = ($ 1$,
$ \frac{1+\sqrt{5}}{2}$,
$ \frac{1-\sqrt{5}}{2}$,
$ -1$;\ \ 
$ -1$,
$ -1$,
$ \frac{-1+\sqrt{5}}{2}$;\ \ 
$ -1$,
$ -\frac{1+\sqrt{5}}{2}$;\ \ 
$ 1$)

Not pseudo-unitary. 

\vskip 1ex 
\color{grey}

\noindent9. ind = $(4;5
)_{1}^{3}$:\ \ 
$d_i$ = ($1.0$,
$1.618$,
$-0.618$,
$-1.0$) 

\vskip 0.7ex
\hangindent=3em \hangafter=1
$D^2=$ 5.0 = 
 $5$

\vskip 0.7ex
\hangindent=3em \hangafter=1
$T = ( 0,
\frac{3}{5},
\frac{1}{5},
\frac{4}{5} )
$,

\vskip 0.7ex
\hangindent=3em \hangafter=1
$S$ = ($ 1$,
$ \frac{1+\sqrt{5}}{2}$,
$ \frac{1-\sqrt{5}}{2}$,
$ -1$;\ \ 
$ -1$,
$ -1$,
$ \frac{-1+\sqrt{5}}{2}$;\ \ 
$ -1$,
$ -\frac{1+\sqrt{5}}{2}$;\ \ 
$ 1$)

Not pseudo-unitary. 

\vskip 1ex 
\color{grey}

\noindent10. ind = $(4;5
)_{1}^{4}$:\ \ 
$d_i$ = ($1.0$,
$1.618$,
$-0.618$,
$-1.0$) 

\vskip 0.7ex
\hangindent=3em \hangafter=1
$D^2=$ 5.0 = 
 $5$

\vskip 0.7ex
\hangindent=3em \hangafter=1
$T = ( 0,
\frac{3}{5},
\frac{4}{5},
\frac{2}{5} )
$,

\vskip 0.7ex
\hangindent=3em \hangafter=1
$S$ = ($ 1$,
$ \frac{1+\sqrt{5}}{2}$,
$ \frac{1-\sqrt{5}}{2}$,
$ -1$;\ \ 
$ -1$,
$ -1$,
$ \frac{-1+\sqrt{5}}{2}$;\ \ 
$ -1$,
$ -\frac{1+\sqrt{5}}{2}$;\ \ 
$ 1$)

Not pseudo-unitary. 

\vskip 1ex 

 \color{black} \vskip 2ex

\noindent11. ind = $(4;9
)_{1}^{1}$:\ \ 
$d_i$ = ($1.0$,
$1.879$,
$2.532$,
$2.879$) 

\vskip 0.7ex
\hangindent=3em \hangafter=1
$D^2=$ 19.234 = 
 $9+6c^{1}_{9}
+3c^{2}_{9}
$

\vskip 0.7ex
\hangindent=3em \hangafter=1
$T = ( 0,
\frac{1}{3},
\frac{2}{9},
\frac{2}{3} )
$,

\vskip 0.7ex
\hangindent=3em \hangafter=1
$S$ = ($ 1$,
$  -c_9^4 $,
$ 1+c^{1}_{9}
$,
$ \xi_{9}^{4}$;\ \ 
$ -\xi_{9}^{4}$,
$ 1+c^{1}_{9}
$,
$ -1$;\ \ 
$0$,
$ -1-c^{1}_{9}
$;\ \ 
$  -c_9^4 $)

\vskip 1ex 
\color{grey}

\noindent12. ind = $(4;9
)_{1}^{8}$:\ \ 
$d_i$ = ($1.0$,
$1.879$,
$2.532$,
$2.879$) 

\vskip 0.7ex
\hangindent=3em \hangafter=1
$D^2=$ 19.234 = 
 $9+6c^{1}_{9}
+3c^{2}_{9}
$

\vskip 0.7ex
\hangindent=3em \hangafter=1
$T = ( 0,
\frac{2}{3},
\frac{7}{9},
\frac{1}{3} )
$,

\vskip 0.7ex
\hangindent=3em \hangafter=1
$S$ = ($ 1$,
$  -c_9^4 $,
$ 1+c^{1}_{9}
$,
$ \xi_{9}^{4}$;\ \ 
$ -\xi_{9}^{4}$,
$ 1+c^{1}_{9}
$,
$ -1$;\ \ 
$0$,
$ -1-c^{1}_{9}
$;\ \ 
$  -c_9^4 $)

\vskip 1ex 
\color{grey}

\noindent13. ind = $(4;9
)_{1}^{5}$:\ \ 
$d_i$ = ($1.0$,
$0.652$,
$-0.347$,
$-0.879$) 

\vskip 0.7ex
\hangindent=3em \hangafter=1
$D^2=$ 2.319 = 
 $9-3c^{1}_{9}
-6c^{2}_{9}
$

\vskip 0.7ex
\hangindent=3em \hangafter=1
$T = ( 0,
\frac{1}{3},
\frac{2}{3},
\frac{1}{9} )
$,

\vskip 0.7ex
\hangindent=3em \hangafter=1
$S$ = ($ 1$,
$ \xi_{9}^{2,4}$,
$ -c^{2}_{9}
$,
$ 1-c^{1}_{9}
-c^{2}_{9}
$;\ \ 
$ -c^{2}_{9}
$,
$ -1$,
$ -1+c^{1}_{9}
+c^{2}_{9}
$;\ \ 
$ -\xi_{9}^{2,4}$,
$ 1-c^{1}_{9}
-c^{2}_{9}
$;\ \ 
$0$)

Not pseudo-unitary. 

\vskip 1ex 
\color{grey}

\noindent14. ind = $(4;9
)_{1}^{2}$:\ \ 
$d_i$ = ($1.0$,
$1.347$,
$-0.532$,
$-1.532$) 

\vskip 0.7ex
\hangindent=3em \hangafter=1
$D^2=$ 5.445 = 
 $9-3c^{1}_{9}
+3c^{2}_{9}
$

\vskip 0.7ex
\hangindent=3em \hangafter=1
$T = ( 0,
\frac{4}{9},
\frac{1}{3},
\frac{2}{3} )
$,

\vskip 0.7ex
\hangindent=3em \hangafter=1
$S$ = ($ 1$,
$ 1+c^{2}_{9}
$,
$ -\xi_{9}^{1,2}$,
$ -c^{1}_{9}
$;\ \ 
$0$,
$ -1-c^{2}_{9}
$,
$ 1+c^{2}_{9}
$;\ \ 
$ -c^{1}_{9}
$,
$ -1$;\ \ 
$ \xi_{9}^{1,2}$)

Not pseudo-unitary. 

\vskip 1ex 
\color{grey}

\noindent15. ind = $(4;9
)_{1}^{7}$:\ \ 
$d_i$ = ($1.0$,
$1.347$,
$-0.532$,
$-1.532$) 

\vskip 0.7ex
\hangindent=3em \hangafter=1
$D^2=$ 5.445 = 
 $9-3c^{1}_{9}
+3c^{2}_{9}
$

\vskip 0.7ex
\hangindent=3em \hangafter=1
$T = ( 0,
\frac{5}{9},
\frac{2}{3},
\frac{1}{3} )
$,

\vskip 0.7ex
\hangindent=3em \hangafter=1
$S$ = ($ 1$,
$ 1+c^{2}_{9}
$,
$ -\xi_{9}^{1,2}$,
$ -c^{1}_{9}
$;\ \ 
$0$,
$ -1-c^{2}_{9}
$,
$ 1+c^{2}_{9}
$;\ \ 
$ -c^{1}_{9}
$,
$ -1$;\ \ 
$ \xi_{9}^{1,2}$)

Not pseudo-unitary. 

\vskip 1ex 
\color{grey}

\noindent16. ind = $(4;9
)_{1}^{4}$:\ \ 
$d_i$ = ($1.0$,
$0.652$,
$-0.347$,
$-0.879$) 

\vskip 0.7ex
\hangindent=3em \hangafter=1
$D^2=$ 2.319 = 
 $9-3c^{1}_{9}
-6c^{2}_{9}
$

\vskip 0.7ex
\hangindent=3em \hangafter=1
$T = ( 0,
\frac{2}{3},
\frac{1}{3},
\frac{8}{9} )
$,

\vskip 0.7ex
\hangindent=3em \hangafter=1
$S$ = ($ 1$,
$ \xi_{9}^{2,4}$,
$ -c^{2}_{9}
$,
$ 1-c^{1}_{9}
-c^{2}_{9}
$;\ \ 
$ -c^{2}_{9}
$,
$ -1$,
$ -1+c^{1}_{9}
+c^{2}_{9}
$;\ \ 
$ -\xi_{9}^{2,4}$,
$ 1-c^{1}_{9}
-c^{2}_{9}
$;\ \ 
$0$)

Not pseudo-unitary. 

\vskip 1ex 

 \color{black} \vskip 2ex

\noindent17. ind = $(4;40
)_{1}^{1}$:\ \ 
$d_i$ = ($1.0$,
$1.0$,
$1.618$,
$1.618$) 

\vskip 0.7ex
\hangindent=3em \hangafter=1
$D^2=$ 7.236 = 
 $5+\sqrt{5}$

\vskip 0.7ex
\hangindent=3em \hangafter=1
$T = ( 0,
\frac{1}{4},
\frac{2}{5},
\frac{13}{20} )
$,

\vskip 0.7ex
\hangindent=3em \hangafter=1
$S$ = ($ 1$,
$ 1$,
$ \frac{1+\sqrt{5}}{2}$,
$ \frac{1+\sqrt{5}}{2}$;\ \ 
$ -1$,
$ \frac{1+\sqrt{5}}{2}$,
$ -\frac{1+\sqrt{5}}{2}$;\ \ 
$ -1$,
$ -1$;\ \ 
$ 1$)

\vskip 1ex 
\color{grey}

\noindent18. ind = $(4;40
)_{1}^{9}$:\ \ 
$d_i$ = ($1.0$,
$1.0$,
$1.618$,
$1.618$) 

\vskip 0.7ex
\hangindent=3em \hangafter=1
$D^2=$ 7.236 = 
 $5+\sqrt{5}$

\vskip 0.7ex
\hangindent=3em \hangafter=1
$T = ( 0,
\frac{1}{4},
\frac{3}{5},
\frac{17}{20} )
$,

\vskip 0.7ex
\hangindent=3em \hangafter=1
$S$ = ($ 1$,
$ 1$,
$ \frac{1+\sqrt{5}}{2}$,
$ \frac{1+\sqrt{5}}{2}$;\ \ 
$ -1$,
$ \frac{1+\sqrt{5}}{2}$,
$ -\frac{1+\sqrt{5}}{2}$;\ \ 
$ -1$,
$ -1$;\ \ 
$ 1$)

\vskip 1ex 
\color{grey}

\noindent19. ind = $(4;40
)_{1}^{11}$:\ \ 
$d_i$ = ($1.0$,
$1.0$,
$1.618$,
$1.618$) 

\vskip 0.7ex
\hangindent=3em \hangafter=1
$D^2=$ 7.236 = 
 $5+\sqrt{5}$

\vskip 0.7ex
\hangindent=3em \hangafter=1
$T = ( 0,
\frac{3}{4},
\frac{2}{5},
\frac{3}{20} )
$,

\vskip 0.7ex
\hangindent=3em \hangafter=1
$S$ = ($ 1$,
$ 1$,
$ \frac{1+\sqrt{5}}{2}$,
$ \frac{1+\sqrt{5}}{2}$;\ \ 
$ -1$,
$ \frac{1+\sqrt{5}}{2}$,
$ -\frac{1+\sqrt{5}}{2}$;\ \ 
$ -1$,
$ -1$;\ \ 
$ 1$)

\vskip 1ex 
\color{grey}

\noindent20. ind = $(4;40
)_{1}^{19}$:\ \ 
$d_i$ = ($1.0$,
$1.0$,
$1.618$,
$1.618$) 

\vskip 0.7ex
\hangindent=3em \hangafter=1
$D^2=$ 7.236 = 
 $5+\sqrt{5}$

\vskip 0.7ex
\hangindent=3em \hangafter=1
$T = ( 0,
\frac{3}{4},
\frac{3}{5},
\frac{7}{20} )
$,

\vskip 0.7ex
\hangindent=3em \hangafter=1
$S$ = ($ 1$,
$ 1$,
$ \frac{1+\sqrt{5}}{2}$,
$ \frac{1+\sqrt{5}}{2}$;\ \ 
$ -1$,
$ \frac{1+\sqrt{5}}{2}$,
$ -\frac{1+\sqrt{5}}{2}$;\ \ 
$ -1$,
$ -1$;\ \ 
$ 1$)

\vskip 1ex 
\color{grey}

\noindent21. ind = $(4;40
)_{1}^{13}$:\ \ 
$d_i$ = ($1.0$,
$1.0$,
$-0.618$,
$-0.618$) 

\vskip 0.7ex
\hangindent=3em \hangafter=1
$D^2=$ 2.763 = 
 $5-\sqrt{5}$

\vskip 0.7ex
\hangindent=3em \hangafter=1
$T = ( 0,
\frac{1}{4},
\frac{1}{5},
\frac{9}{20} )
$,

\vskip 0.7ex
\hangindent=3em \hangafter=1
$S$ = ($ 1$,
$ 1$,
$ \frac{1-\sqrt{5}}{2}$,
$ \frac{1-\sqrt{5}}{2}$;\ \ 
$ -1$,
$ \frac{1-\sqrt{5}}{2}$,
$ \frac{-1+\sqrt{5}}{2}$;\ \ 
$ -1$,
$ -1$;\ \ 
$ 1$)

Not pseudo-unitary. 

\vskip 1ex 
\color{grey}

\noindent22. ind = $(4;40
)_{1}^{17}$:\ \ 
$d_i$ = ($1.0$,
$1.0$,
$-0.618$,
$-0.618$) 

\vskip 0.7ex
\hangindent=3em \hangafter=1
$D^2=$ 2.763 = 
 $5-\sqrt{5}$

\vskip 0.7ex
\hangindent=3em \hangafter=1
$T = ( 0,
\frac{1}{4},
\frac{4}{5},
\frac{1}{20} )
$,

\vskip 0.7ex
\hangindent=3em \hangafter=1
$S$ = ($ 1$,
$ 1$,
$ \frac{1-\sqrt{5}}{2}$,
$ \frac{1-\sqrt{5}}{2}$;\ \ 
$ -1$,
$ \frac{1-\sqrt{5}}{2}$,
$ \frac{-1+\sqrt{5}}{2}$;\ \ 
$ -1$,
$ -1$;\ \ 
$ 1$)

Not pseudo-unitary. 

\vskip 1ex 
\color{grey}

\noindent23. ind = $(4;40
)_{1}^{3}$:\ \ 
$d_i$ = ($1.0$,
$1.0$,
$-0.618$,
$-0.618$) 

\vskip 0.7ex
\hangindent=3em \hangafter=1
$D^2=$ 2.763 = 
 $5-\sqrt{5}$

\vskip 0.7ex
\hangindent=3em \hangafter=1
$T = ( 0,
\frac{3}{4},
\frac{1}{5},
\frac{19}{20} )
$,

\vskip 0.7ex
\hangindent=3em \hangafter=1
$S$ = ($ 1$,
$ 1$,
$ \frac{1-\sqrt{5}}{2}$,
$ \frac{1-\sqrt{5}}{2}$;\ \ 
$ -1$,
$ \frac{1-\sqrt{5}}{2}$,
$ \frac{-1+\sqrt{5}}{2}$;\ \ 
$ -1$,
$ -1$;\ \ 
$ 1$)

Not pseudo-unitary. 

\vskip 1ex 
\color{grey}

\noindent24. ind = $(4;40
)_{1}^{7}$:\ \ 
$d_i$ = ($1.0$,
$1.0$,
$-0.618$,
$-0.618$) 

\vskip 0.7ex
\hangindent=3em \hangafter=1
$D^2=$ 2.763 = 
 $5-\sqrt{5}$

\vskip 0.7ex
\hangindent=3em \hangafter=1
$T = ( 0,
\frac{3}{4},
\frac{4}{5},
\frac{11}{20} )
$,

\vskip 0.7ex
\hangindent=3em \hangafter=1
$S$ = ($ 1$,
$ 1$,
$ \frac{1-\sqrt{5}}{2}$,
$ \frac{1-\sqrt{5}}{2}$;\ \ 
$ -1$,
$ \frac{1-\sqrt{5}}{2}$,
$ \frac{-1+\sqrt{5}}{2}$;\ \ 
$ -1$,
$ -1$;\ \ 
$ 1$)

Not pseudo-unitary. 

\vskip 1ex 

 \color{black} \vskip 2ex 
\color{blue}

\noindent25. ind = $(4;40
)_{2}^{1}$:\ \ 
$d_i$ = ($1.0$,
$0.618$,
$-0.618$,
$-1.0$) 

\vskip 0.7ex
\hangindent=3em \hangafter=1
$D^2=$ 2.763 = 
 $5-\sqrt{5}$

\vskip 0.7ex
\hangindent=3em \hangafter=1
$T = ( 0,
\frac{1}{20},
\frac{4}{5},
\frac{1}{4} )
$,

\vskip 0.7ex
\hangindent=3em \hangafter=1
$S$ = ($ 1$,
$ \frac{-1+\sqrt{5}}{2}$,
$ \frac{1-\sqrt{5}}{2}$,
$ -1$;\ \ 
$ 1$,
$ 1$,
$ \frac{-1+\sqrt{5}}{2}$;\ \ 
$ -1$,
$ \frac{-1+\sqrt{5}}{2}$;\ \ 
$ -1$)

Not pseudo-unitary. 

\vskip 1ex 
\color{grey}

\noindent26. ind = $(4;40
)_{2}^{13}$:\ \ 
$d_i$ = ($1.0$,
$1.618$,
$-1.0$,
$-1.618$) 

\vskip 0.7ex
\hangindent=3em \hangafter=1
$D^2=$ 7.236 = 
 $5+\sqrt{5}$

\vskip 0.7ex
\hangindent=3em \hangafter=1
$T = ( 0,
\frac{2}{5},
\frac{1}{4},
\frac{13}{20} )
$,

\vskip 0.7ex
\hangindent=3em \hangafter=1
$S$ = ($ 1$,
$ \frac{1+\sqrt{5}}{2}$,
$ -1$,
$ -\frac{1+\sqrt{5}}{2}$;\ \ 
$ -1$,
$ -\frac{1+\sqrt{5}}{2}$,
$ 1$;\ \ 
$ -1$,
$ -\frac{1+\sqrt{5}}{2}$;\ \ 
$ 1$)

Pseudo-unitary $\sim$  
$(4;40
)_{1}^{1}$

\vskip 1ex 
\color{grey}

\noindent27. ind = $(4;40
)_{2}^{3}$:\ \ 
$d_i$ = ($1.0$,
$1.618$,
$-1.0$,
$-1.618$) 

\vskip 0.7ex
\hangindent=3em \hangafter=1
$D^2=$ 7.236 = 
 $5+\sqrt{5}$

\vskip 0.7ex
\hangindent=3em \hangafter=1
$T = ( 0,
\frac{2}{5},
\frac{3}{4},
\frac{3}{20} )
$,

\vskip 0.7ex
\hangindent=3em \hangafter=1
$S$ = ($ 1$,
$ \frac{1+\sqrt{5}}{2}$,
$ -1$,
$ -\frac{1+\sqrt{5}}{2}$;\ \ 
$ -1$,
$ -\frac{1+\sqrt{5}}{2}$,
$ 1$;\ \ 
$ -1$,
$ -\frac{1+\sqrt{5}}{2}$;\ \ 
$ 1$)

Pseudo-unitary $\sim$  
$(4;40
)_{1}^{11}$

\vskip 1ex 
\color{grey}

\noindent28. ind = $(4;40
)_{2}^{9}$:\ \ 
$d_i$ = ($1.0$,
$0.618$,
$-0.618$,
$-1.0$) 

\vskip 0.7ex
\hangindent=3em \hangafter=1
$D^2=$ 2.763 = 
 $5-\sqrt{5}$

\vskip 0.7ex
\hangindent=3em \hangafter=1
$T = ( 0,
\frac{9}{20},
\frac{1}{5},
\frac{1}{4} )
$,

\vskip 0.7ex
\hangindent=3em \hangafter=1
$S$ = ($ 1$,
$ \frac{-1+\sqrt{5}}{2}$,
$ \frac{1-\sqrt{5}}{2}$,
$ -1$;\ \ 
$ 1$,
$ 1$,
$ \frac{-1+\sqrt{5}}{2}$;\ \ 
$ -1$,
$ \frac{-1+\sqrt{5}}{2}$;\ \ 
$ -1$)

Not pseudo-unitary. 

\vskip 1ex 
\color{grey}

\noindent29. ind = $(4;40
)_{2}^{11}$:\ \ 
$d_i$ = ($1.0$,
$0.618$,
$-0.618$,
$-1.0$) 

\vskip 0.7ex
\hangindent=3em \hangafter=1
$D^2=$ 2.763 = 
 $5-\sqrt{5}$

\vskip 0.7ex
\hangindent=3em \hangafter=1
$T = ( 0,
\frac{11}{20},
\frac{4}{5},
\frac{3}{4} )
$,

\vskip 0.7ex
\hangindent=3em \hangafter=1
$S$ = ($ 1$,
$ \frac{-1+\sqrt{5}}{2}$,
$ \frac{1-\sqrt{5}}{2}$,
$ -1$;\ \ 
$ 1$,
$ 1$,
$ \frac{-1+\sqrt{5}}{2}$;\ \ 
$ -1$,
$ \frac{-1+\sqrt{5}}{2}$;\ \ 
$ -1$)

Not pseudo-unitary. 

\vskip 1ex 
\color{grey}

\noindent30. ind = $(4;40
)_{2}^{17}$:\ \ 
$d_i$ = ($1.0$,
$1.618$,
$-1.0$,
$-1.618$) 

\vskip 0.7ex
\hangindent=3em \hangafter=1
$D^2=$ 7.236 = 
 $5+\sqrt{5}$

\vskip 0.7ex
\hangindent=3em \hangafter=1
$T = ( 0,
\frac{3}{5},
\frac{1}{4},
\frac{17}{20} )
$,

\vskip 0.7ex
\hangindent=3em \hangafter=1
$S$ = ($ 1$,
$ \frac{1+\sqrt{5}}{2}$,
$ -1$,
$ -\frac{1+\sqrt{5}}{2}$;\ \ 
$ -1$,
$ -\frac{1+\sqrt{5}}{2}$,
$ 1$;\ \ 
$ -1$,
$ -\frac{1+\sqrt{5}}{2}$;\ \ 
$ 1$)

Pseudo-unitary $\sim$  
$(4;40
)_{1}^{9}$

\vskip 1ex 
\color{grey}

\noindent31. ind = $(4;40
)_{2}^{7}$:\ \ 
$d_i$ = ($1.0$,
$1.618$,
$-1.0$,
$-1.618$) 

\vskip 0.7ex
\hangindent=3em \hangafter=1
$D^2=$ 7.236 = 
 $5+\sqrt{5}$

\vskip 0.7ex
\hangindent=3em \hangafter=1
$T = ( 0,
\frac{3}{5},
\frac{3}{4},
\frac{7}{20} )
$,

\vskip 0.7ex
\hangindent=3em \hangafter=1
$S$ = ($ 1$,
$ \frac{1+\sqrt{5}}{2}$,
$ -1$,
$ -\frac{1+\sqrt{5}}{2}$;\ \ 
$ -1$,
$ -\frac{1+\sqrt{5}}{2}$,
$ 1$;\ \ 
$ -1$,
$ -\frac{1+\sqrt{5}}{2}$;\ \ 
$ 1$)

Pseudo-unitary $\sim$  
$(4;40
)_{1}^{19}$

\vskip 1ex 
\color{grey}

\noindent32. ind = $(4;40
)_{2}^{19}$:\ \ 
$d_i$ = ($1.0$,
$0.618$,
$-0.618$,
$-1.0$) 

\vskip 0.7ex
\hangindent=3em \hangafter=1
$D^2=$ 2.763 = 
 $5-\sqrt{5}$

\vskip 0.7ex
\hangindent=3em \hangafter=1
$T = ( 0,
\frac{19}{20},
\frac{1}{5},
\frac{3}{4} )
$,

\vskip 0.7ex
\hangindent=3em \hangafter=1
$S$ = ($ 1$,
$ \frac{-1+\sqrt{5}}{2}$,
$ \frac{1-\sqrt{5}}{2}$,
$ -1$;\ \ 
$ 1$,
$ 1$,
$ \frac{-1+\sqrt{5}}{2}$;\ \ 
$ -1$,
$ \frac{-1+\sqrt{5}}{2}$;\ \ 
$ -1$)

Not pseudo-unitary. 

\vskip 1ex 

 \color{black} \vskip 2ex

\

The above list includes all modular data (unitary or non-unitary) from resolved
$\SL$ representations and non-integral MTCs.
Since for rank 4, all the passing $\SL$ representations are resolved, the list
actually includes all modular data from non-integral MTCs.

\section{Rank-5 modular data}
\label{Section8}

\subsection{A list of 5-dimensional irrep-sum $\SL$ representations}

The following is a list of 5-dimensional irrep-sum $\SL$ representations, with
3 types of representations omitted and one type skipped.  For details and
notations, see Section 3 of this file.

\

\noindent 1: (dims,levels) = $(3\oplus
1\oplus
1;5,
1,
1
)$,
irreps = $3_{5}^{1}\oplus
1_{1}^{1}\oplus
1_{1}^{1}$,
pord$(\rho_\text{isum}(\mathfrak{t})) = 5$,

\vskip 0.7ex
\hangindent=5.5em \hangafter=1
{\white .}\hskip 1em $\rho_\text{isum}(\mathfrak{t})$ =
 $( 0,
\frac{1}{5},
\frac{4}{5} )
\oplus
( 0 )
\oplus
( 0 )
$,

\vskip 0.7ex
\hangindent=5.5em \hangafter=1
{\white .}\hskip 1em $\rho_\text{isum}(\mathfrak{s})$ =
($\sqrt{\frac{1}{5}}$,
$-\sqrt{\frac{2}{5}}$,
$-\sqrt{\frac{2}{5}}$;
$-\frac{5+\sqrt{5}}{10}$,
$\frac{5-\sqrt{5}}{10}$;
$-\frac{5+\sqrt{5}}{10}$)
 $\oplus$
($1$)
 $\oplus$
($1$)

Fail:
dims is $(d,1,1,\cdots)$, $d=(\mathrm{ord}(T)+1)/2$, self-dual, ...
 Prop. B.3 (3)

 \ \color{black}

\noindent 2: (dims,levels) = $(3\oplus
1\oplus
1;5,
1,
1
)$,
irreps = $3_{5}^{3}\oplus
1_{1}^{1}\oplus
1_{1}^{1}$,
pord$(\rho_\text{isum}(\mathfrak{t})) = 5$,

\vskip 0.7ex
\hangindent=5.5em \hangafter=1
{\white .}\hskip 1em $\rho_\text{isum}(\mathfrak{t})$ =
 $( 0,
\frac{2}{5},
\frac{3}{5} )
\oplus
( 0 )
\oplus
( 0 )
$,

\vskip 0.7ex
\hangindent=5.5em \hangafter=1
{\white .}\hskip 1em $\rho_\text{isum}(\mathfrak{s})$ =
($-\sqrt{\frac{1}{5}}$,
$-\sqrt{\frac{2}{5}}$,
$-\sqrt{\frac{2}{5}}$;
$\frac{-5+\sqrt{5}}{10}$,
$\frac{5+\sqrt{5}}{10}$;
$\frac{-5+\sqrt{5}}{10}$)
 $\oplus$
($1$)
 $\oplus$
($1$)

Fail:
dims is $(d,1,1,\cdots)$, $d=(\mathrm{ord}(T)+1)/2$, self-dual, ...
 Prop. B.3 (3)

 \ \color{black}

\noindent 3: (dims,levels) = $(3\oplus
1\oplus
1;8,
1,
1
)$,
irreps = $3_{8}^{1,0}\oplus
1_{1}^{1}\oplus
1_{1}^{1}$,
pord$(\rho_\text{isum}(\mathfrak{t})) = 8$,

\vskip 0.7ex
\hangindent=5.5em \hangafter=1
{\white .}\hskip 1em $\rho_\text{isum}(\mathfrak{t})$ =
 $( 0,
\frac{1}{8},
\frac{5}{8} )
\oplus
( 0 )
\oplus
( 0 )
$,

\vskip 0.7ex
\hangindent=5.5em \hangafter=1
{\white .}\hskip 1em $\rho_\text{isum}(\mathfrak{s})$ =
$\mathrm{i}$($0$,
$\sqrt{\frac{1}{2}}$,
$\sqrt{\frac{1}{2}}$;\ \ 
$-\frac{1}{2}$,
$\frac{1}{2}$;\ \ 
$-\frac{1}{2}$)
 $\oplus$
($1$)
 $\oplus$
($1$)

Fail:
for $\rho = \rho_1+l\chi, ...,
 (\rho_1(\mathfrak s)/\chi(\mathfrak s))^2\neq$ id. Prop. B.3 (2)

 \ \color{black}

\noindent 4: (dims,levels) = $(3\oplus
1\oplus
1;8,
1,
1
)$,
irreps = $3_{8}^{3,0}\oplus
1_{1}^{1}\oplus
1_{1}^{1}$,
pord$(\rho_\text{isum}(\mathfrak{t})) = 8$,

\vskip 0.7ex
\hangindent=5.5em \hangafter=1
{\white .}\hskip 1em $\rho_\text{isum}(\mathfrak{t})$ =
 $( 0,
\frac{3}{8},
\frac{7}{8} )
\oplus
( 0 )
\oplus
( 0 )
$,

\vskip 0.7ex
\hangindent=5.5em \hangafter=1
{\white .}\hskip 1em $\rho_\text{isum}(\mathfrak{s})$ =
$\mathrm{i}$($0$,
$\sqrt{\frac{1}{2}}$,
$\sqrt{\frac{1}{2}}$;\ \ 
$\frac{1}{2}$,
$-\frac{1}{2}$;\ \ 
$\frac{1}{2}$)
 $\oplus$
($1$)
 $\oplus$
($1$)

Fail:
for $\rho = \rho_1+l\chi, ...,
 (\rho_1(\mathfrak s)/\chi(\mathfrak s))^2\neq$ id. Prop. B.3 (2)

 \ \color{black}

\noindent 5: (dims,levels) = $(3\oplus
1\oplus
1;10,
2,
2
)$,
irreps = $3_{5}^{3}
\hskip -1.5pt \otimes \hskip -1.5pt
1_{2}^{1,0}\oplus
1_{2}^{1,0}\oplus
1_{2}^{1,0}$,
pord$(\rho_\text{isum}(\mathfrak{t})) = 5$,

\vskip 0.7ex
\hangindent=5.5em \hangafter=1
{\white .}\hskip 1em $\rho_\text{isum}(\mathfrak{t})$ =
 $( \frac{1}{2},
\frac{1}{10},
\frac{9}{10} )
\oplus
( \frac{1}{2} )
\oplus
( \frac{1}{2} )
$,

\vskip 0.7ex
\hangindent=5.5em \hangafter=1
{\white .}\hskip 1em $\rho_\text{isum}(\mathfrak{s})$ =
($\sqrt{\frac{1}{5}}$,
$-\sqrt{\frac{2}{5}}$,
$-\sqrt{\frac{2}{5}}$;
$\frac{5-\sqrt{5}}{10}$,
$-\frac{5+\sqrt{5}}{10}$;
$\frac{5-\sqrt{5}}{10}$)
 $\oplus$
($-1$)
 $\oplus$
($-1$)

Fail:
dims is $(d,1,1,\cdots)$, $d=(\mathrm{ord}(T)+1)/2$, self-dual, ...
 Prop. B.3 (3)

 \ \color{black}

\noindent 6: (dims,levels) = $(3\oplus
1\oplus
1;10,
2,
2
)$,
irreps = $3_{5}^{1}
\hskip -1.5pt \otimes \hskip -1.5pt
1_{2}^{1,0}\oplus
1_{2}^{1,0}\oplus
1_{2}^{1,0}$,
pord$(\rho_\text{isum}(\mathfrak{t})) = 5$,

\vskip 0.7ex
\hangindent=5.5em \hangafter=1
{\white .}\hskip 1em $\rho_\text{isum}(\mathfrak{t})$ =
 $( \frac{1}{2},
\frac{3}{10},
\frac{7}{10} )
\oplus
( \frac{1}{2} )
\oplus
( \frac{1}{2} )
$,

\vskip 0.7ex
\hangindent=5.5em \hangafter=1
{\white .}\hskip 1em $\rho_\text{isum}(\mathfrak{s})$ =
($-\sqrt{\frac{1}{5}}$,
$-\sqrt{\frac{2}{5}}$,
$-\sqrt{\frac{2}{5}}$;
$\frac{5+\sqrt{5}}{10}$,
$\frac{-5+\sqrt{5}}{10}$;
$\frac{5+\sqrt{5}}{10}$)
 $\oplus$
($-1$)
 $\oplus$
($-1$)

Fail:
dims is $(d,1,1,\cdots)$, $d=(\mathrm{ord}(T)+1)/2$, self-dual, ...
 Prop. B.3 (3)

 \ \color{black}

\noindent 7: (dims,levels) = $(3\oplus
1\oplus
1;15,
3,
3
)$,
irreps = $3_{5}^{1}
\hskip -1.5pt \otimes \hskip -1.5pt
1_{3}^{1,0}\oplus
1_{3}^{1,0}\oplus
1_{3}^{1,0}$,
pord$(\rho_\text{isum}(\mathfrak{t})) = 5$,

\vskip 0.7ex
\hangindent=5.5em \hangafter=1
{\white .}\hskip 1em $\rho_\text{isum}(\mathfrak{t})$ =
 $( \frac{1}{3},
\frac{2}{15},
\frac{8}{15} )
\oplus
( \frac{1}{3} )
\oplus
( \frac{1}{3} )
$,

\vskip 0.7ex
\hangindent=5.5em \hangafter=1
{\white .}\hskip 1em $\rho_\text{isum}(\mathfrak{s})$ =
($\sqrt{\frac{1}{5}}$,
$-\sqrt{\frac{2}{5}}$,
$-\sqrt{\frac{2}{5}}$;
$-\frac{5+\sqrt{5}}{10}$,
$\frac{5-\sqrt{5}}{10}$;
$-\frac{5+\sqrt{5}}{10}$)
 $\oplus$
($1$)
 $\oplus$
($1$)

Fail:
dims is $(d,1,1,\cdots)$, $d=(\mathrm{ord}(T)+1)/2$, self-dual, ...
 Prop. B.3 (3)

 \ \color{black}

\noindent 8: (dims,levels) = $(3\oplus
1\oplus
1;15,
3,
3
)$,
irreps = $3_{5}^{3}
\hskip -1.5pt \otimes \hskip -1.5pt
1_{3}^{1,0}\oplus
1_{3}^{1,0}\oplus
1_{3}^{1,0}$,
pord$(\rho_\text{isum}(\mathfrak{t})) = 5$,

\vskip 0.7ex
\hangindent=5.5em \hangafter=1
{\white .}\hskip 1em $\rho_\text{isum}(\mathfrak{t})$ =
 $( \frac{1}{3},
\frac{11}{15},
\frac{14}{15} )
\oplus
( \frac{1}{3} )
\oplus
( \frac{1}{3} )
$,

\vskip 0.7ex
\hangindent=5.5em \hangafter=1
{\white .}\hskip 1em $\rho_\text{isum}(\mathfrak{s})$ =
($-\sqrt{\frac{1}{5}}$,
$-\sqrt{\frac{2}{5}}$,
$-\sqrt{\frac{2}{5}}$;
$\frac{-5+\sqrt{5}}{10}$,
$\frac{5+\sqrt{5}}{10}$;
$\frac{-5+\sqrt{5}}{10}$)
 $\oplus$
($1$)
 $\oplus$
($1$)

Fail:
dims is $(d,1,1,\cdots)$, $d=(\mathrm{ord}(T)+1)/2$, self-dual, ...
 Prop. B.3 (3)

 \ \color{black}

\noindent 9: (dims,levels) = $(3\oplus
1\oplus
1;20,
4,
4
)$,
irreps = $3_{5}^{1}
\hskip -1.5pt \otimes \hskip -1.5pt
1_{4}^{1,0}\oplus
1_{4}^{1,0}\oplus
1_{4}^{1,0}$,
pord$(\rho_\text{isum}(\mathfrak{t})) = 5$,

\vskip 0.7ex
\hangindent=5.5em \hangafter=1
{\white .}\hskip 1em $\rho_\text{isum}(\mathfrak{t})$ =
 $( \frac{1}{4},
\frac{1}{20},
\frac{9}{20} )
\oplus
( \frac{1}{4} )
\oplus
( \frac{1}{4} )
$,

\vskip 0.7ex
\hangindent=5.5em \hangafter=1
{\white .}\hskip 1em $\rho_\text{isum}(\mathfrak{s})$ =
$\mathrm{i}$($\sqrt{\frac{1}{5}}$,
$\sqrt{\frac{2}{5}}$,
$\sqrt{\frac{2}{5}}$;\ \ 
$-\frac{5+\sqrt{5}}{10}$,
$\frac{5-\sqrt{5}}{10}$;\ \ 
$-\frac{5+\sqrt{5}}{10}$)
 $\oplus$
$\mathrm{i}$($1$)
 $\oplus$
$\mathrm{i}$($1$)

Fail:
dims is $(d,1,1,\cdots)$, $d=(\mathrm{ord}(T)+1)/2$, self-dual, ...
 Prop. B.3 (3)

 \ \color{black}

\noindent 10: (dims,levels) = $(3\oplus
1\oplus
1;20,
4,
4
)$,
irreps = $3_{5}^{3}
\hskip -1.5pt \otimes \hskip -1.5pt
1_{4}^{1,0}\oplus
1_{4}^{1,0}\oplus
1_{4}^{1,0}$,
pord$(\rho_\text{isum}(\mathfrak{t})) = 5$,

\vskip 0.7ex
\hangindent=5.5em \hangafter=1
{\white .}\hskip 1em $\rho_\text{isum}(\mathfrak{t})$ =
 $( \frac{1}{4},
\frac{13}{20},
\frac{17}{20} )
\oplus
( \frac{1}{4} )
\oplus
( \frac{1}{4} )
$,

\vskip 0.7ex
\hangindent=5.5em \hangafter=1
{\white .}\hskip 1em $\rho_\text{isum}(\mathfrak{s})$ =
$\mathrm{i}$($-\sqrt{\frac{1}{5}}$,
$\sqrt{\frac{2}{5}}$,
$\sqrt{\frac{2}{5}}$;\ \ 
$\frac{-5+\sqrt{5}}{10}$,
$\frac{5+\sqrt{5}}{10}$;\ \ 
$\frac{-5+\sqrt{5}}{10}$)
 $\oplus$
$\mathrm{i}$($1$)
 $\oplus$
$\mathrm{i}$($1$)

Fail:
dims is $(d,1,1,\cdots)$, $d=(\mathrm{ord}(T)+1)/2$, self-dual, ...
 Prop. B.3 (3)

 \ \color{black}

\noindent 11: (dims,levels) = $(3\oplus
1\oplus
1;24,
3,
3
)$,
irreps = $3_{8}^{3,0}
\hskip -1.5pt \otimes \hskip -1.5pt
1_{3}^{1,0}\oplus
1_{3}^{1,0}\oplus
1_{3}^{1,0}$,
pord$(\rho_\text{isum}(\mathfrak{t})) = 8$,

\vskip 0.7ex
\hangindent=5.5em \hangafter=1
{\white .}\hskip 1em $\rho_\text{isum}(\mathfrak{t})$ =
 $( \frac{1}{3},
\frac{5}{24},
\frac{17}{24} )
\oplus
( \frac{1}{3} )
\oplus
( \frac{1}{3} )
$,

\vskip 0.7ex
\hangindent=5.5em \hangafter=1
{\white .}\hskip 1em $\rho_\text{isum}(\mathfrak{s})$ =
$\mathrm{i}$($0$,
$\sqrt{\frac{1}{2}}$,
$\sqrt{\frac{1}{2}}$;\ \ 
$\frac{1}{2}$,
$-\frac{1}{2}$;\ \ 
$\frac{1}{2}$)
 $\oplus$
($1$)
 $\oplus$
($1$)

Fail:
for $\rho = \rho_1+l\chi, ...,
 (\rho_1(\mathfrak s)/\chi(\mathfrak s))^2\neq$ id. Prop. B.3 (2)

 \ \color{black}

\noindent 12: (dims,levels) = $(3\oplus
1\oplus
1;24,
3,
3
)$,
irreps = $3_{8}^{1,0}
\hskip -1.5pt \otimes \hskip -1.5pt
1_{3}^{1,0}\oplus
1_{3}^{1,0}\oplus
1_{3}^{1,0}$,
pord$(\rho_\text{isum}(\mathfrak{t})) = 8$,

\vskip 0.7ex
\hangindent=5.5em \hangafter=1
{\white .}\hskip 1em $\rho_\text{isum}(\mathfrak{t})$ =
 $( \frac{1}{3},
\frac{11}{24},
\frac{23}{24} )
\oplus
( \frac{1}{3} )
\oplus
( \frac{1}{3} )
$,

\vskip 0.7ex
\hangindent=5.5em \hangafter=1
{\white .}\hskip 1em $\rho_\text{isum}(\mathfrak{s})$ =
$\mathrm{i}$($0$,
$\sqrt{\frac{1}{2}}$,
$\sqrt{\frac{1}{2}}$;\ \ 
$-\frac{1}{2}$,
$\frac{1}{2}$;\ \ 
$-\frac{1}{2}$)
 $\oplus$
($1$)
 $\oplus$
($1$)

Fail:
for $\rho = \rho_1+l\chi, ...,
 (\rho_1(\mathfrak s)/\chi(\mathfrak s))^2\neq$ id. Prop. B.3 (2)

 \ \color{black}

\noindent 13: (dims,levels) = $(3\oplus
1\oplus
1;30,
6,
6
)$,
irreps = $3_{5}^{1}
\hskip -1.5pt \otimes \hskip -1.5pt
1_{3}^{1,0}
\hskip -1.5pt \otimes \hskip -1.5pt
1_{2}^{1,0}\oplus
1_{3}^{1,0}
\hskip -1.5pt \otimes \hskip -1.5pt
1_{2}^{1,0}\oplus
1_{3}^{1,0}
\hskip -1.5pt \otimes \hskip -1.5pt
1_{2}^{1,0}$,
pord$(\rho_\text{isum}(\mathfrak{t})) = 5$,

\vskip 0.7ex
\hangindent=5.5em \hangafter=1
{\white .}\hskip 1em $\rho_\text{isum}(\mathfrak{t})$ =
 $( \frac{5}{6},
\frac{1}{30},
\frac{19}{30} )
\oplus
( \frac{5}{6} )
\oplus
( \frac{5}{6} )
$,

\vskip 0.7ex
\hangindent=5.5em \hangafter=1
{\white .}\hskip 1em $\rho_\text{isum}(\mathfrak{s})$ =
($-\sqrt{\frac{1}{5}}$,
$-\sqrt{\frac{2}{5}}$,
$-\sqrt{\frac{2}{5}}$;
$\frac{5+\sqrt{5}}{10}$,
$\frac{-5+\sqrt{5}}{10}$;
$\frac{5+\sqrt{5}}{10}$)
 $\oplus$
($-1$)
 $\oplus$
($-1$)

Fail:
dims is $(d,1,1,\cdots)$, $d=(\mathrm{ord}(T)+1)/2$, self-dual, ...
 Prop. B.3 (3)

 \ \color{black}

\noindent 14: (dims,levels) = $(3\oplus
1\oplus
1;30,
6,
6
)$,
irreps = $3_{5}^{3}
\hskip -1.5pt \otimes \hskip -1.5pt
1_{3}^{1,0}
\hskip -1.5pt \otimes \hskip -1.5pt
1_{2}^{1,0}\oplus
1_{3}^{1,0}
\hskip -1.5pt \otimes \hskip -1.5pt
1_{2}^{1,0}\oplus
1_{3}^{1,0}
\hskip -1.5pt \otimes \hskip -1.5pt
1_{2}^{1,0}$,
pord$(\rho_\text{isum}(\mathfrak{t})) = 5$,

\vskip 0.7ex
\hangindent=5.5em \hangafter=1
{\white .}\hskip 1em $\rho_\text{isum}(\mathfrak{t})$ =
 $( \frac{5}{6},
\frac{7}{30},
\frac{13}{30} )
\oplus
( \frac{5}{6} )
\oplus
( \frac{5}{6} )
$,

\vskip 0.7ex
\hangindent=5.5em \hangafter=1
{\white .}\hskip 1em $\rho_\text{isum}(\mathfrak{s})$ =
($\sqrt{\frac{1}{5}}$,
$-\sqrt{\frac{2}{5}}$,
$-\sqrt{\frac{2}{5}}$;
$\frac{5-\sqrt{5}}{10}$,
$-\frac{5+\sqrt{5}}{10}$;
$\frac{5-\sqrt{5}}{10}$)
 $\oplus$
($-1$)
 $\oplus$
($-1$)

Fail:
dims is $(d,1,1,\cdots)$, $d=(\mathrm{ord}(T)+1)/2$, self-dual, ...
 Prop. B.3 (3)

 \ \color{black}

\noindent 15: (dims,levels) = $(3\oplus
1\oplus
1;60,
12,
12
)$,
irreps = $3_{5}^{3}
\hskip -1.5pt \otimes \hskip -1.5pt
1_{4}^{1,0}
\hskip -1.5pt \otimes \hskip -1.5pt
1_{3}^{1,0}\oplus
1_{4}^{1,0}
\hskip -1.5pt \otimes \hskip -1.5pt
1_{3}^{1,0}\oplus
1_{4}^{1,0}
\hskip -1.5pt \otimes \hskip -1.5pt
1_{3}^{1,0}$,
pord$(\rho_\text{isum}(\mathfrak{t})) = 5$,

\vskip 0.7ex
\hangindent=5.5em \hangafter=1
{\white .}\hskip 1em $\rho_\text{isum}(\mathfrak{t})$ =
 $( \frac{7}{12},
\frac{11}{60},
\frac{59}{60} )
\oplus
( \frac{7}{12} )
\oplus
( \frac{7}{12} )
$,

\vskip 0.7ex
\hangindent=5.5em \hangafter=1
{\white .}\hskip 1em $\rho_\text{isum}(\mathfrak{s})$ =
$\mathrm{i}$($-\sqrt{\frac{1}{5}}$,
$\sqrt{\frac{2}{5}}$,
$\sqrt{\frac{2}{5}}$;\ \ 
$\frac{-5+\sqrt{5}}{10}$,
$\frac{5+\sqrt{5}}{10}$;\ \ 
$\frac{-5+\sqrt{5}}{10}$)
 $\oplus$
$\mathrm{i}$($1$)
 $\oplus$
$\mathrm{i}$($1$)

Fail:
dims is $(d,1,1,\cdots)$, $d=(\mathrm{ord}(T)+1)/2$, self-dual, ...
 Prop. B.3 (3)

 \ \color{black}

\noindent 16: (dims,levels) = $(3\oplus
1\oplus
1;60,
12,
12
)$,
irreps = $3_{5}^{1}
\hskip -1.5pt \otimes \hskip -1.5pt
1_{4}^{1,0}
\hskip -1.5pt \otimes \hskip -1.5pt
1_{3}^{1,0}\oplus
1_{4}^{1,0}
\hskip -1.5pt \otimes \hskip -1.5pt
1_{3}^{1,0}\oplus
1_{4}^{1,0}
\hskip -1.5pt \otimes \hskip -1.5pt
1_{3}^{1,0}$,
pord$(\rho_\text{isum}(\mathfrak{t})) = 5$,

\vskip 0.7ex
\hangindent=5.5em \hangafter=1
{\white .}\hskip 1em $\rho_\text{isum}(\mathfrak{t})$ =
 $( \frac{7}{12},
\frac{23}{60},
\frac{47}{60} )
\oplus
( \frac{7}{12} )
\oplus
( \frac{7}{12} )
$,

\vskip 0.7ex
\hangindent=5.5em \hangafter=1
{\white .}\hskip 1em $\rho_\text{isum}(\mathfrak{s})$ =
$\mathrm{i}$($\sqrt{\frac{1}{5}}$,
$\sqrt{\frac{2}{5}}$,
$\sqrt{\frac{2}{5}}$;\ \ 
$-\frac{5+\sqrt{5}}{10}$,
$\frac{5-\sqrt{5}}{10}$;\ \ 
$-\frac{5+\sqrt{5}}{10}$)
 $\oplus$
$\mathrm{i}$($1$)
 $\oplus$
$\mathrm{i}$($1$)

Fail:
dims is $(d,1,1,\cdots)$, $d=(\mathrm{ord}(T)+1)/2$, self-dual, ...
 Prop. B.3 (3)

 \ \color{black}

\noindent 17: (dims,levels) = $(3\oplus
2;4,
3
)$,
irreps = $3_{4}^{1,0}\oplus
2_{3}^{1,0}$,
pord$(\rho_\text{isum}(\mathfrak{t})) = 12$,

\vskip 0.7ex
\hangindent=5.5em \hangafter=1
{\white .}\hskip 1em $\rho_\text{isum}(\mathfrak{t})$ =
 $( 0,
\frac{1}{4},
\frac{3}{4} )
\oplus
( 0,
\frac{1}{3} )
$,

\vskip 0.7ex
\hangindent=5.5em \hangafter=1
{\white .}\hskip 1em $\rho_\text{isum}(\mathfrak{s})$ =
($0$,
$\sqrt{\frac{1}{2}}$,
$\sqrt{\frac{1}{2}}$;
$-\frac{1}{2}$,
$\frac{1}{2}$;
$-\frac{1}{2}$)
 $\oplus$
$\mathrm{i}$($-\sqrt{\frac{1}{3}}$,
$\sqrt{\frac{2}{3}}$;\ \ 
$\sqrt{\frac{1}{3}}$)

Fail:
Tr$_I(C) = -1 <$  0 for I = [ 1/3 ]. Prop. B.4 (1) eqn. (B.18)

 \ \color{black}

\noindent 18: (dims,levels) = $(3\oplus
2;4,
3
)$,
irreps = $3_{4}^{1,0}\oplus
2_{3}^{1,8}$,
pord$(\rho_\text{isum}(\mathfrak{t})) = 12$,

\vskip 0.7ex
\hangindent=5.5em \hangafter=1
{\white .}\hskip 1em $\rho_\text{isum}(\mathfrak{t})$ =
 $( 0,
\frac{1}{4},
\frac{3}{4} )
\oplus
( 0,
\frac{2}{3} )
$,

\vskip 0.7ex
\hangindent=5.5em \hangafter=1
{\white .}\hskip 1em $\rho_\text{isum}(\mathfrak{s})$ =
($0$,
$\sqrt{\frac{1}{2}}$,
$\sqrt{\frac{1}{2}}$;
$-\frac{1}{2}$,
$\frac{1}{2}$;
$-\frac{1}{2}$)
 $\oplus$
$\mathrm{i}$($\sqrt{\frac{1}{3}}$,
$\sqrt{\frac{2}{3}}$;\ \ 
$-\sqrt{\frac{1}{3}}$)

Fail:
Tr$_I(C) = -1 <$  0 for I = [ 2/3 ]. Prop. B.4 (1) eqn. (B.18)

 \ \color{black}

\noindent 19: (dims,levels) = $(3\oplus
2;4,
12
)$,
irreps = $3_{4}^{1,0}\oplus
2_{3}^{1,0}
\hskip -1.5pt \otimes \hskip -1.5pt
1_{4}^{1,0}$,
pord$(\rho_\text{isum}(\mathfrak{t})) = 12$,

\vskip 0.7ex
\hangindent=5.5em \hangafter=1
{\white .}\hskip 1em $\rho_\text{isum}(\mathfrak{t})$ =
 $( 0,
\frac{1}{4},
\frac{3}{4} )
\oplus
( \frac{1}{4},
\frac{7}{12} )
$,

\vskip 0.7ex
\hangindent=5.5em \hangafter=1
{\white .}\hskip 1em $\rho_\text{isum}(\mathfrak{s})$ =
($0$,
$\sqrt{\frac{1}{2}}$,
$\sqrt{\frac{1}{2}}$;
$-\frac{1}{2}$,
$\frac{1}{2}$;
$-\frac{1}{2}$)
 $\oplus$
($\sqrt{\frac{1}{3}}$,
$\sqrt{\frac{2}{3}}$;
$-\sqrt{\frac{1}{3}}$)

Fail:
cnd($\rho(\mathfrak s)_\mathrm{ndeg}$) = 24 does not divide
 ord($\rho(\mathfrak t)$)=12. Prop. B.4 (2)

 \ \color{black}

\noindent 20: (dims,levels) = $(3\oplus
2;4,
12
)$,
irreps = $3_{4}^{1,0}\oplus
2_{3}^{1,8}
\hskip -1.5pt \otimes \hskip -1.5pt
1_{4}^{1,0}$,
pord$(\rho_\text{isum}(\mathfrak{t})) = 12$,

\vskip 0.7ex
\hangindent=5.5em \hangafter=1
{\white .}\hskip 1em $\rho_\text{isum}(\mathfrak{t})$ =
 $( 0,
\frac{1}{4},
\frac{3}{4} )
\oplus
( \frac{1}{4},
\frac{11}{12} )
$,

\vskip 0.7ex
\hangindent=5.5em \hangafter=1
{\white .}\hskip 1em $\rho_\text{isum}(\mathfrak{s})$ =
($0$,
$\sqrt{\frac{1}{2}}$,
$\sqrt{\frac{1}{2}}$;
$-\frac{1}{2}$,
$\frac{1}{2}$;
$-\frac{1}{2}$)
 $\oplus$
($-\sqrt{\frac{1}{3}}$,
$\sqrt{\frac{2}{3}}$;
$\sqrt{\frac{1}{3}}$)

Fail:
cnd($\rho(\mathfrak s)_\mathrm{ndeg}$) = 24 does not divide
 ord($\rho(\mathfrak t)$)=12. Prop. B.4 (2)

 \ \color{black}

\noindent 21: (dims,levels) = $(3\oplus
2;4,
12
)$,
irreps = $3_{4}^{1,0}\oplus
2_{3}^{1,0}
\hskip -1.5pt \otimes \hskip -1.5pt
1_{4}^{1,6}$,
pord$(\rho_\text{isum}(\mathfrak{t})) = 12$,

\vskip 0.7ex
\hangindent=5.5em \hangafter=1
{\white .}\hskip 1em $\rho_\text{isum}(\mathfrak{t})$ =
 $( 0,
\frac{1}{4},
\frac{3}{4} )
\oplus
( \frac{3}{4},
\frac{1}{12} )
$,

\vskip 0.7ex
\hangindent=5.5em \hangafter=1
{\white .}\hskip 1em $\rho_\text{isum}(\mathfrak{s})$ =
($0$,
$\sqrt{\frac{1}{2}}$,
$\sqrt{\frac{1}{2}}$;
$-\frac{1}{2}$,
$\frac{1}{2}$;
$-\frac{1}{2}$)
 $\oplus$
($-\sqrt{\frac{1}{3}}$,
$\sqrt{\frac{2}{3}}$;
$\sqrt{\frac{1}{3}}$)

Fail:
cnd($\rho(\mathfrak s)_\mathrm{ndeg}$) = 24 does not divide
 ord($\rho(\mathfrak t)$)=12. Prop. B.4 (2)

 \ \color{black}

\noindent 22: (dims,levels) = $(3\oplus
2;4,
12
)$,
irreps = $3_{4}^{1,0}\oplus
2_{3}^{1,8}
\hskip -1.5pt \otimes \hskip -1.5pt
1_{4}^{1,6}$,
pord$(\rho_\text{isum}(\mathfrak{t})) = 12$,

\vskip 0.7ex
\hangindent=5.5em \hangafter=1
{\white .}\hskip 1em $\rho_\text{isum}(\mathfrak{t})$ =
 $( 0,
\frac{1}{4},
\frac{3}{4} )
\oplus
( \frac{3}{4},
\frac{5}{12} )
$,

\vskip 0.7ex
\hangindent=5.5em \hangafter=1
{\white .}\hskip 1em $\rho_\text{isum}(\mathfrak{s})$ =
($0$,
$\sqrt{\frac{1}{2}}$,
$\sqrt{\frac{1}{2}}$;
$-\frac{1}{2}$,
$\frac{1}{2}$;
$-\frac{1}{2}$)
 $\oplus$
($\sqrt{\frac{1}{3}}$,
$\sqrt{\frac{2}{3}}$;
$-\sqrt{\frac{1}{3}}$)

Fail:
cnd($\rho(\mathfrak s)_\mathrm{ndeg}$) = 24 does not divide
 ord($\rho(\mathfrak t)$)=12. Prop. B.4 (2)

 \ \color{black}

 \color{blue}

\noindent 23: (dims,levels) = $(3\oplus
2;5,
2
)$,
irreps = $3_{5}^{1}\oplus
2_{2}^{1,0}$,
pord$(\rho_\text{isum}(\mathfrak{t})) = 10$,

\vskip 0.7ex
\hangindent=5.5em \hangafter=1
{\white .}\hskip 1em $\rho_\text{isum}(\mathfrak{t})$ =
 $( 0,
\frac{1}{5},
\frac{4}{5} )
\oplus
( 0,
\frac{1}{2} )
$,

\vskip 0.7ex
\hangindent=5.5em \hangafter=1
{\white .}\hskip 1em $\rho_\text{isum}(\mathfrak{s})$ =
($\sqrt{\frac{1}{5}}$,
$-\sqrt{\frac{2}{5}}$,
$-\sqrt{\frac{2}{5}}$;
$-\frac{5+\sqrt{5}}{10}$,
$\frac{5-\sqrt{5}}{10}$;
$-\frac{5+\sqrt{5}}{10}$)
 $\oplus$
($-\frac{1}{2}$,
$-\sqrt{\frac{3}{4}}$;
$\frac{1}{2}$)

Pass. 

 \ \color{black}

 \color{blue}

\noindent 24: (dims,levels) = $(3\oplus
2;5,
2
)$,
irreps = $3_{5}^{3}\oplus
2_{2}^{1,0}$,
pord$(\rho_\text{isum}(\mathfrak{t})) = 10$,

\vskip 0.7ex
\hangindent=5.5em \hangafter=1
{\white .}\hskip 1em $\rho_\text{isum}(\mathfrak{t})$ =
 $( 0,
\frac{2}{5},
\frac{3}{5} )
\oplus
( 0,
\frac{1}{2} )
$,

\vskip 0.7ex
\hangindent=5.5em \hangafter=1
{\white .}\hskip 1em $\rho_\text{isum}(\mathfrak{s})$ =
($-\sqrt{\frac{1}{5}}$,
$-\sqrt{\frac{2}{5}}$,
$-\sqrt{\frac{2}{5}}$;
$\frac{-5+\sqrt{5}}{10}$,
$\frac{5+\sqrt{5}}{10}$;
$\frac{-5+\sqrt{5}}{10}$)
 $\oplus$
($-\frac{1}{2}$,
$-\sqrt{\frac{3}{4}}$;
$\frac{1}{2}$)

Pass. 

 \ \color{black}

\noindent 25: (dims,levels) = $(3\oplus
2;5,
3
)$,
irreps = $3_{5}^{1}\oplus
2_{3}^{1,0}$,
pord$(\rho_\text{isum}(\mathfrak{t})) = 15$,

\vskip 0.7ex
\hangindent=5.5em \hangafter=1
{\white .}\hskip 1em $\rho_\text{isum}(\mathfrak{t})$ =
 $( 0,
\frac{1}{5},
\frac{4}{5} )
\oplus
( 0,
\frac{1}{3} )
$,

\vskip 0.7ex
\hangindent=5.5em \hangafter=1
{\white .}\hskip 1em $\rho_\text{isum}(\mathfrak{s})$ =
($\sqrt{\frac{1}{5}}$,
$-\sqrt{\frac{2}{5}}$,
$-\sqrt{\frac{2}{5}}$;
$-\frac{5+\sqrt{5}}{10}$,
$\frac{5-\sqrt{5}}{10}$;
$-\frac{5+\sqrt{5}}{10}$)
 $\oplus$
$\mathrm{i}$($-\sqrt{\frac{1}{3}}$,
$\sqrt{\frac{2}{3}}$;\ \ 
$\sqrt{\frac{1}{3}}$)

Fail:
Tr$_I(C) = -1 <$  0 for I = [ 1/3 ]. Prop. B.4 (1) eqn. (B.18)

 \ \color{black}

\noindent 26: (dims,levels) = $(3\oplus
2;5,
3
)$,
irreps = $3_{5}^{1}\oplus
2_{3}^{1,8}$,
pord$(\rho_\text{isum}(\mathfrak{t})) = 15$,

\vskip 0.7ex
\hangindent=5.5em \hangafter=1
{\white .}\hskip 1em $\rho_\text{isum}(\mathfrak{t})$ =
 $( 0,
\frac{1}{5},
\frac{4}{5} )
\oplus
( 0,
\frac{2}{3} )
$,

\vskip 0.7ex
\hangindent=5.5em \hangafter=1
{\white .}\hskip 1em $\rho_\text{isum}(\mathfrak{s})$ =
($\sqrt{\frac{1}{5}}$,
$-\sqrt{\frac{2}{5}}$,
$-\sqrt{\frac{2}{5}}$;
$-\frac{5+\sqrt{5}}{10}$,
$\frac{5-\sqrt{5}}{10}$;
$-\frac{5+\sqrt{5}}{10}$)
 $\oplus$
$\mathrm{i}$($\sqrt{\frac{1}{3}}$,
$\sqrt{\frac{2}{3}}$;\ \ 
$-\sqrt{\frac{1}{3}}$)

Fail:
Tr$_I(C) = -1 <$  0 for I = [ 2/3 ]. Prop. B.4 (1) eqn. (B.18)

 \ \color{black}

\noindent 27: (dims,levels) = $(3\oplus
2;5,
3
)$,
irreps = $3_{5}^{3}\oplus
2_{3}^{1,0}$,
pord$(\rho_\text{isum}(\mathfrak{t})) = 15$,

\vskip 0.7ex
\hangindent=5.5em \hangafter=1
{\white .}\hskip 1em $\rho_\text{isum}(\mathfrak{t})$ =
 $( 0,
\frac{2}{5},
\frac{3}{5} )
\oplus
( 0,
\frac{1}{3} )
$,

\vskip 0.7ex
\hangindent=5.5em \hangafter=1
{\white .}\hskip 1em $\rho_\text{isum}(\mathfrak{s})$ =
($-\sqrt{\frac{1}{5}}$,
$-\sqrt{\frac{2}{5}}$,
$-\sqrt{\frac{2}{5}}$;
$\frac{-5+\sqrt{5}}{10}$,
$\frac{5+\sqrt{5}}{10}$;
$\frac{-5+\sqrt{5}}{10}$)
 $\oplus$
$\mathrm{i}$($-\sqrt{\frac{1}{3}}$,
$\sqrt{\frac{2}{3}}$;\ \ 
$\sqrt{\frac{1}{3}}$)

Fail:
Tr$_I(C) = -1 <$  0 for I = [ 1/3 ]. Prop. B.4 (1) eqn. (B.18)

 \ \color{black}

\noindent 28: (dims,levels) = $(3\oplus
2;5,
3
)$,
irreps = $3_{5}^{3}\oplus
2_{3}^{1,8}$,
pord$(\rho_\text{isum}(\mathfrak{t})) = 15$,

\vskip 0.7ex
\hangindent=5.5em \hangafter=1
{\white .}\hskip 1em $\rho_\text{isum}(\mathfrak{t})$ =
 $( 0,
\frac{2}{5},
\frac{3}{5} )
\oplus
( 0,
\frac{2}{3} )
$,

\vskip 0.7ex
\hangindent=5.5em \hangafter=1
{\white .}\hskip 1em $\rho_\text{isum}(\mathfrak{s})$ =
($-\sqrt{\frac{1}{5}}$,
$-\sqrt{\frac{2}{5}}$,
$-\sqrt{\frac{2}{5}}$;
$\frac{-5+\sqrt{5}}{10}$,
$\frac{5+\sqrt{5}}{10}$;
$\frac{-5+\sqrt{5}}{10}$)
 $\oplus$
$\mathrm{i}$($\sqrt{\frac{1}{3}}$,
$\sqrt{\frac{2}{3}}$;\ \ 
$-\sqrt{\frac{1}{3}}$)

Fail:
Tr$_I(C) = -1 <$  0 for I = [ 2/3 ]. Prop. B.4 (1) eqn. (B.18)

 \ \color{black}

\noindent 29: (dims,levels) = $(3\oplus
2;5,
5
)$,
irreps = $3_{5}^{1}\oplus
2_{5}^{1}$,
pord$(\rho_\text{isum}(\mathfrak{t})) = 5$,

\vskip 0.7ex
\hangindent=5.5em \hangafter=1
{\white .}\hskip 1em $\rho_\text{isum}(\mathfrak{t})$ =
 $( 0,
\frac{1}{5},
\frac{4}{5} )
\oplus
( \frac{1}{5},
\frac{4}{5} )
$,

\vskip 0.7ex
\hangindent=5.5em \hangafter=1
{\white .}\hskip 1em $\rho_\text{isum}(\mathfrak{s})$ =
($\sqrt{\frac{1}{5}}$,
$-\sqrt{\frac{2}{5}}$,
$-\sqrt{\frac{2}{5}}$;
$-\frac{5+\sqrt{5}}{10}$,
$\frac{5-\sqrt{5}}{10}$;
$-\frac{5+\sqrt{5}}{10}$)
 $\oplus$
$\mathrm{i}$($-\frac{1}{\sqrt{5}}c^{3}_{20}
$,
$\frac{1}{\sqrt{5}}c^{1}_{20}
$;\ \ 
$\frac{1}{\sqrt{5}}c^{3}_{20}
$)

Fail:
Integral: $D_{\rho}(\sigma)_{\theta} \propto $ id,
 for all $\sigma$ and all $\theta$-eigenspaces that can contain unit. Prop. B.5 (6)

 \ \color{black}

\noindent 30: (dims,levels) = $(3\oplus
2;5,
5
)$,
irreps = $3_{5}^{3}\oplus
2_{5}^{2}$,
pord$(\rho_\text{isum}(\mathfrak{t})) = 5$,

\vskip 0.7ex
\hangindent=5.5em \hangafter=1
{\white .}\hskip 1em $\rho_\text{isum}(\mathfrak{t})$ =
 $( 0,
\frac{2}{5},
\frac{3}{5} )
\oplus
( \frac{2}{5},
\frac{3}{5} )
$,

\vskip 0.7ex
\hangindent=5.5em \hangafter=1
{\white .}\hskip 1em $\rho_\text{isum}(\mathfrak{s})$ =
($-\sqrt{\frac{1}{5}}$,
$-\sqrt{\frac{2}{5}}$,
$-\sqrt{\frac{2}{5}}$;
$\frac{-5+\sqrt{5}}{10}$,
$\frac{5+\sqrt{5}}{10}$;
$\frac{-5+\sqrt{5}}{10}$)
 $\oplus$
$\mathrm{i}$($-\frac{1}{\sqrt{5}}c^{1}_{20}
$,
$\frac{1}{\sqrt{5}}c^{3}_{20}
$;\ \ 
$\frac{1}{\sqrt{5}}c^{1}_{20}
$)

Fail:
Integral: $D_{\rho}(\sigma)_{\theta} \propto $ id,
 for all $\sigma$ and all $\theta$-eigenspaces that can contain unit. Prop. B.5 (6)

 \ \color{black}

\noindent 31: (dims,levels) = $(3\oplus
2;8,
2
)$,
irreps = $3_{8}^{1,0}\oplus
2_{2}^{1,0}$,
pord$(\rho_\text{isum}(\mathfrak{t})) = 8$,

\vskip 0.7ex
\hangindent=5.5em \hangafter=1
{\white .}\hskip 1em $\rho_\text{isum}(\mathfrak{t})$ =
 $( 0,
\frac{1}{8},
\frac{5}{8} )
\oplus
( 0,
\frac{1}{2} )
$,

\vskip 0.7ex
\hangindent=5.5em \hangafter=1
{\white .}\hskip 1em $\rho_\text{isum}(\mathfrak{s})$ =
$\mathrm{i}$($0$,
$\sqrt{\frac{1}{2}}$,
$\sqrt{\frac{1}{2}}$;\ \ 
$-\frac{1}{2}$,
$\frac{1}{2}$;\ \ 
$-\frac{1}{2}$)
 $\oplus$
($-\frac{1}{2}$,
$-\sqrt{\frac{3}{4}}$;
$\frac{1}{2}$)

Fail:
Tr$_I(C) = -1 <$  0 for I = [ 1/2 ]. Prop. B.4 (1) eqn. (B.18)

 \ \color{black}

\noindent 32: (dims,levels) = $(3\oplus
2;8,
2
)$,
irreps = $3_{8}^{3,0}\oplus
2_{2}^{1,0}$,
pord$(\rho_\text{isum}(\mathfrak{t})) = 8$,

\vskip 0.7ex
\hangindent=5.5em \hangafter=1
{\white .}\hskip 1em $\rho_\text{isum}(\mathfrak{t})$ =
 $( 0,
\frac{3}{8},
\frac{7}{8} )
\oplus
( 0,
\frac{1}{2} )
$,

\vskip 0.7ex
\hangindent=5.5em \hangafter=1
{\white .}\hskip 1em $\rho_\text{isum}(\mathfrak{s})$ =
$\mathrm{i}$($0$,
$\sqrt{\frac{1}{2}}$,
$\sqrt{\frac{1}{2}}$;\ \ 
$\frac{1}{2}$,
$-\frac{1}{2}$;\ \ 
$\frac{1}{2}$)
 $\oplus$
($-\frac{1}{2}$,
$-\sqrt{\frac{3}{4}}$;
$\frac{1}{2}$)

Fail:
Tr$_I(C) = -1 <$  0 for I = [ 1/2 ]. Prop. B.4 (1) eqn. (B.18)

 \ \color{black}

 \color{blue}

\noindent 33: (dims,levels) = $(3\oplus
2;8,
3
)$,
irreps = $3_{8}^{1,0}\oplus
2_{3}^{1,0}$,
pord$(\rho_\text{isum}(\mathfrak{t})) = 24$,

\vskip 0.7ex
\hangindent=5.5em \hangafter=1
{\white .}\hskip 1em $\rho_\text{isum}(\mathfrak{t})$ =
 $( 0,
\frac{1}{8},
\frac{5}{8} )
\oplus
( 0,
\frac{1}{3} )
$,

\vskip 0.7ex
\hangindent=5.5em \hangafter=1
{\white .}\hskip 1em $\rho_\text{isum}(\mathfrak{s})$ =
$\mathrm{i}$($0$,
$\sqrt{\frac{1}{2}}$,
$\sqrt{\frac{1}{2}}$;\ \ 
$-\frac{1}{2}$,
$\frac{1}{2}$;\ \ 
$-\frac{1}{2}$)
 $\oplus$
$\mathrm{i}$($-\sqrt{\frac{1}{3}}$,
$\sqrt{\frac{2}{3}}$;\ \ 
$\sqrt{\frac{1}{3}}$)

Pass. 

 \ \color{black}

 \color{blue}

\noindent 34: (dims,levels) = $(3\oplus
2;8,
3
)$,
irreps = $3_{8}^{1,0}\oplus
2_{3}^{1,8}$,
pord$(\rho_\text{isum}(\mathfrak{t})) = 24$,

\vskip 0.7ex
\hangindent=5.5em \hangafter=1
{\white .}\hskip 1em $\rho_\text{isum}(\mathfrak{t})$ =
 $( 0,
\frac{1}{8},
\frac{5}{8} )
\oplus
( 0,
\frac{2}{3} )
$,

\vskip 0.7ex
\hangindent=5.5em \hangafter=1
{\white .}\hskip 1em $\rho_\text{isum}(\mathfrak{s})$ =
$\mathrm{i}$($0$,
$\sqrt{\frac{1}{2}}$,
$\sqrt{\frac{1}{2}}$;\ \ 
$-\frac{1}{2}$,
$\frac{1}{2}$;\ \ 
$-\frac{1}{2}$)
 $\oplus$
$\mathrm{i}$($\sqrt{\frac{1}{3}}$,
$\sqrt{\frac{2}{3}}$;\ \ 
$-\sqrt{\frac{1}{3}}$)

Pass. 

 \ \color{black}

 \color{blue}

\noindent 35: (dims,levels) = $(3\oplus
2;8,
3
)$,
irreps = $3_{8}^{3,0}\oplus
2_{3}^{1,0}$,
pord$(\rho_\text{isum}(\mathfrak{t})) = 24$,

\vskip 0.7ex
\hangindent=5.5em \hangafter=1
{\white .}\hskip 1em $\rho_\text{isum}(\mathfrak{t})$ =
 $( 0,
\frac{3}{8},
\frac{7}{8} )
\oplus
( 0,
\frac{1}{3} )
$,

\vskip 0.7ex
\hangindent=5.5em \hangafter=1
{\white .}\hskip 1em $\rho_\text{isum}(\mathfrak{s})$ =
$\mathrm{i}$($0$,
$\sqrt{\frac{1}{2}}$,
$\sqrt{\frac{1}{2}}$;\ \ 
$\frac{1}{2}$,
$-\frac{1}{2}$;\ \ 
$\frac{1}{2}$)
 $\oplus$
$\mathrm{i}$($-\sqrt{\frac{1}{3}}$,
$\sqrt{\frac{2}{3}}$;\ \ 
$\sqrt{\frac{1}{3}}$)

Pass. 

 \ \color{black}

 \color{blue}

\noindent 36: (dims,levels) = $(3\oplus
2;8,
3
)$,
irreps = $3_{8}^{3,0}\oplus
2_{3}^{1,8}$,
pord$(\rho_\text{isum}(\mathfrak{t})) = 24$,

\vskip 0.7ex
\hangindent=5.5em \hangafter=1
{\white .}\hskip 1em $\rho_\text{isum}(\mathfrak{t})$ =
 $( 0,
\frac{3}{8},
\frac{7}{8} )
\oplus
( 0,
\frac{2}{3} )
$,

\vskip 0.7ex
\hangindent=5.5em \hangafter=1
{\white .}\hskip 1em $\rho_\text{isum}(\mathfrak{s})$ =
$\mathrm{i}$($0$,
$\sqrt{\frac{1}{2}}$,
$\sqrt{\frac{1}{2}}$;\ \ 
$\frac{1}{2}$,
$-\frac{1}{2}$;\ \ 
$\frac{1}{2}$)
 $\oplus$
$\mathrm{i}$($\sqrt{\frac{1}{3}}$,
$\sqrt{\frac{2}{3}}$;\ \ 
$-\sqrt{\frac{1}{3}}$)

Pass. 

 \ \color{black}

\noindent 37: (dims,levels) = $(3\oplus
2;8,
8
)$,
irreps = $3_{8}^{1,0}\oplus
2_{8}^{1,0}$,
pord$(\rho_\text{isum}(\mathfrak{t})) = 8$,

\vskip 0.7ex
\hangindent=5.5em \hangafter=1
{\white .}\hskip 1em $\rho_\text{isum}(\mathfrak{t})$ =
 $( 0,
\frac{1}{8},
\frac{5}{8} )
\oplus
( \frac{1}{8},
\frac{3}{8} )
$,

\vskip 0.7ex
\hangindent=5.5em \hangafter=1
{\white .}\hskip 1em $\rho_\text{isum}(\mathfrak{s})$ =
$\mathrm{i}$($0$,
$\sqrt{\frac{1}{2}}$,
$\sqrt{\frac{1}{2}}$;\ \ 
$-\frac{1}{2}$,
$\frac{1}{2}$;\ \ 
$-\frac{1}{2}$)
 $\oplus$
($-\sqrt{\frac{1}{2}}$,
$\sqrt{\frac{1}{2}}$;
$\sqrt{\frac{1}{2}}$)

Fail:
Tr$_I(C) = -1 <$  0 for I = [ 3/8 ]. Prop. B.4 (1) eqn. (B.18)

 \ \color{black}

\noindent 38: (dims,levels) = $(3\oplus
2;8,
8
)$,
irreps = $3_{8}^{1,0}\oplus
2_{8}^{1,9}$,
pord$(\rho_\text{isum}(\mathfrak{t})) = 8$,

\vskip 0.7ex
\hangindent=5.5em \hangafter=1
{\white .}\hskip 1em $\rho_\text{isum}(\mathfrak{t})$ =
 $( 0,
\frac{1}{8},
\frac{5}{8} )
\oplus
( \frac{1}{8},
\frac{7}{8} )
$,

\vskip 0.7ex
\hangindent=5.5em \hangafter=1
{\white .}\hskip 1em $\rho_\text{isum}(\mathfrak{s})$ =
$\mathrm{i}$($0$,
$\sqrt{\frac{1}{2}}$,
$\sqrt{\frac{1}{2}}$;\ \ 
$-\frac{1}{2}$,
$\frac{1}{2}$;\ \ 
$-\frac{1}{2}$)
 $\oplus$
$\mathrm{i}$($-\sqrt{\frac{1}{2}}$,
$\sqrt{\frac{1}{2}}$;\ \ 
$\sqrt{\frac{1}{2}}$)

Fail:
$\sigma(\rho(\mathfrak s)_\mathrm{ndeg}) \neq
 (\rho(\mathfrak t)^a \rho(\mathfrak s) \rho(\mathfrak t)^b
 \rho(\mathfrak s) \rho(\mathfrak t)^a)_\mathrm{ndeg}$,
 $\sigma = a$ = 3. Prop. B.5 (3) eqn. (B.25)

 \ \color{black}

\noindent 39: (dims,levels) = $(3\oplus
2;8,
8
)$,
irreps = $3_{8}^{1,0}\oplus
2_{8}^{1,3}$,
pord$(\rho_\text{isum}(\mathfrak{t})) = 8$,

\vskip 0.7ex
\hangindent=5.5em \hangafter=1
{\white .}\hskip 1em $\rho_\text{isum}(\mathfrak{t})$ =
 $( 0,
\frac{1}{8},
\frac{5}{8} )
\oplus
( \frac{3}{8},
\frac{5}{8} )
$,

\vskip 0.7ex
\hangindent=5.5em \hangafter=1
{\white .}\hskip 1em $\rho_\text{isum}(\mathfrak{s})$ =
$\mathrm{i}$($0$,
$\sqrt{\frac{1}{2}}$,
$\sqrt{\frac{1}{2}}$;\ \ 
$-\frac{1}{2}$,
$\frac{1}{2}$;\ \ 
$-\frac{1}{2}$)
 $\oplus$
$\mathrm{i}$($-\sqrt{\frac{1}{2}}$,
$\sqrt{\frac{1}{2}}$;\ \ 
$\sqrt{\frac{1}{2}}$)

Fail:
$\sigma(\rho(\mathfrak s)_\mathrm{ndeg}) \neq
 (\rho(\mathfrak t)^a \rho(\mathfrak s) \rho(\mathfrak t)^b
 \rho(\mathfrak s) \rho(\mathfrak t)^a)_\mathrm{ndeg}$,
 $\sigma = a$ = 3. Prop. B.5 (3) eqn. (B.25)

 \ \color{black}

\noindent 40: (dims,levels) = $(3\oplus
2;8,
8
)$,
irreps = $3_{8}^{1,0}\oplus
2_{8}^{1,6}$,
pord$(\rho_\text{isum}(\mathfrak{t})) = 8$,

\vskip 0.7ex
\hangindent=5.5em \hangafter=1
{\white .}\hskip 1em $\rho_\text{isum}(\mathfrak{t})$ =
 $( 0,
\frac{1}{8},
\frac{5}{8} )
\oplus
( \frac{5}{8},
\frac{7}{8} )
$,

\vskip 0.7ex
\hangindent=5.5em \hangafter=1
{\white .}\hskip 1em $\rho_\text{isum}(\mathfrak{s})$ =
$\mathrm{i}$($0$,
$\sqrt{\frac{1}{2}}$,
$\sqrt{\frac{1}{2}}$;\ \ 
$-\frac{1}{2}$,
$\frac{1}{2}$;\ \ 
$-\frac{1}{2}$)
 $\oplus$
($\sqrt{\frac{1}{2}}$,
$\sqrt{\frac{1}{2}}$;
$-\sqrt{\frac{1}{2}}$)

Fail:
Tr$_I(C) = -1 <$  0 for I = [ 7/8 ]. Prop. B.4 (1) eqn. (B.18)

 \ \color{black}

\noindent 41: (dims,levels) = $(3\oplus
2;8,
8
)$,
irreps = $3_{8}^{3,0}\oplus
2_{8}^{1,0}$,
pord$(\rho_\text{isum}(\mathfrak{t})) = 8$,

\vskip 0.7ex
\hangindent=5.5em \hangafter=1
{\white .}\hskip 1em $\rho_\text{isum}(\mathfrak{t})$ =
 $( 0,
\frac{3}{8},
\frac{7}{8} )
\oplus
( \frac{1}{8},
\frac{3}{8} )
$,

\vskip 0.7ex
\hangindent=5.5em \hangafter=1
{\white .}\hskip 1em $\rho_\text{isum}(\mathfrak{s})$ =
$\mathrm{i}$($0$,
$\sqrt{\frac{1}{2}}$,
$\sqrt{\frac{1}{2}}$;\ \ 
$\frac{1}{2}$,
$-\frac{1}{2}$;\ \ 
$\frac{1}{2}$)
 $\oplus$
($-\sqrt{\frac{1}{2}}$,
$\sqrt{\frac{1}{2}}$;
$\sqrt{\frac{1}{2}}$)

Fail:
Tr$_I(C) = -1 <$  0 for I = [ 1/8 ]. Prop. B.4 (1) eqn. (B.18)

 \ \color{black}

\noindent 42: (dims,levels) = $(3\oplus
2;8,
8
)$,
irreps = $3_{8}^{3,0}\oplus
2_{8}^{1,9}$,
pord$(\rho_\text{isum}(\mathfrak{t})) = 8$,

\vskip 0.7ex
\hangindent=5.5em \hangafter=1
{\white .}\hskip 1em $\rho_\text{isum}(\mathfrak{t})$ =
 $( 0,
\frac{3}{8},
\frac{7}{8} )
\oplus
( \frac{1}{8},
\frac{7}{8} )
$,

\vskip 0.7ex
\hangindent=5.5em \hangafter=1
{\white .}\hskip 1em $\rho_\text{isum}(\mathfrak{s})$ =
$\mathrm{i}$($0$,
$\sqrt{\frac{1}{2}}$,
$\sqrt{\frac{1}{2}}$;\ \ 
$\frac{1}{2}$,
$-\frac{1}{2}$;\ \ 
$\frac{1}{2}$)
 $\oplus$
$\mathrm{i}$($-\sqrt{\frac{1}{2}}$,
$\sqrt{\frac{1}{2}}$;\ \ 
$\sqrt{\frac{1}{2}}$)

Fail:
$\sigma(\rho(\mathfrak s)_\mathrm{ndeg}) \neq
 (\rho(\mathfrak t)^a \rho(\mathfrak s) \rho(\mathfrak t)^b
 \rho(\mathfrak s) \rho(\mathfrak t)^a)_\mathrm{ndeg}$,
 $\sigma = a$ = 3. Prop. B.5 (3) eqn. (B.25)

 \ \color{black}

\noindent 43: (dims,levels) = $(3\oplus
2;8,
8
)$,
irreps = $3_{8}^{3,0}\oplus
2_{8}^{1,3}$,
pord$(\rho_\text{isum}(\mathfrak{t})) = 8$,

\vskip 0.7ex
\hangindent=5.5em \hangafter=1
{\white .}\hskip 1em $\rho_\text{isum}(\mathfrak{t})$ =
 $( 0,
\frac{3}{8},
\frac{7}{8} )
\oplus
( \frac{3}{8},
\frac{5}{8} )
$,

\vskip 0.7ex
\hangindent=5.5em \hangafter=1
{\white .}\hskip 1em $\rho_\text{isum}(\mathfrak{s})$ =
$\mathrm{i}$($0$,
$\sqrt{\frac{1}{2}}$,
$\sqrt{\frac{1}{2}}$;\ \ 
$\frac{1}{2}$,
$-\frac{1}{2}$;\ \ 
$\frac{1}{2}$)
 $\oplus$
$\mathrm{i}$($-\sqrt{\frac{1}{2}}$,
$\sqrt{\frac{1}{2}}$;\ \ 
$\sqrt{\frac{1}{2}}$)

Fail:
$\sigma(\rho(\mathfrak s)_\mathrm{ndeg}) \neq
 (\rho(\mathfrak t)^a \rho(\mathfrak s) \rho(\mathfrak t)^b
 \rho(\mathfrak s) \rho(\mathfrak t)^a)_\mathrm{ndeg}$,
 $\sigma = a$ = 3. Prop. B.5 (3) eqn. (B.25)

 \ \color{black}

\noindent 44: (dims,levels) = $(3\oplus
2;8,
8
)$,
irreps = $3_{8}^{3,0}\oplus
2_{8}^{1,6}$,
pord$(\rho_\text{isum}(\mathfrak{t})) = 8$,

\vskip 0.7ex
\hangindent=5.5em \hangafter=1
{\white .}\hskip 1em $\rho_\text{isum}(\mathfrak{t})$ =
 $( 0,
\frac{3}{8},
\frac{7}{8} )
\oplus
( \frac{5}{8},
\frac{7}{8} )
$,

\vskip 0.7ex
\hangindent=5.5em \hangafter=1
{\white .}\hskip 1em $\rho_\text{isum}(\mathfrak{s})$ =
$\mathrm{i}$($0$,
$\sqrt{\frac{1}{2}}$,
$\sqrt{\frac{1}{2}}$;\ \ 
$\frac{1}{2}$,
$-\frac{1}{2}$;\ \ 
$\frac{1}{2}$)
 $\oplus$
($\sqrt{\frac{1}{2}}$,
$\sqrt{\frac{1}{2}}$;
$-\sqrt{\frac{1}{2}}$)

Fail:
Tr$_I(C) = -1 <$  0 for I = [ 5/8 ]. Prop. B.4 (1) eqn. (B.18)

 \ \color{black}

 \color{blue}

\noindent 45: (dims,levels) = $(3\oplus
2;10,
2
)$,
irreps = $3_{5}^{3}
\hskip -1.5pt \otimes \hskip -1.5pt
1_{2}^{1,0}\oplus
2_{2}^{1,0}$,
pord$(\rho_\text{isum}(\mathfrak{t})) = 10$,

\vskip 0.7ex
\hangindent=5.5em \hangafter=1
{\white .}\hskip 1em $\rho_\text{isum}(\mathfrak{t})$ =
 $( \frac{1}{2},
\frac{1}{10},
\frac{9}{10} )
\oplus
( 0,
\frac{1}{2} )
$,

\vskip 0.7ex
\hangindent=5.5em \hangafter=1
{\white .}\hskip 1em $\rho_\text{isum}(\mathfrak{s})$ =
($\sqrt{\frac{1}{5}}$,
$-\sqrt{\frac{2}{5}}$,
$-\sqrt{\frac{2}{5}}$;
$\frac{5-\sqrt{5}}{10}$,
$-\frac{5+\sqrt{5}}{10}$;
$\frac{5-\sqrt{5}}{10}$)
 $\oplus$
($-\frac{1}{2}$,
$-\sqrt{\frac{3}{4}}$;
$\frac{1}{2}$)

Pass. 

 \ \color{black}

 \color{blue}

\noindent 46: (dims,levels) = $(3\oplus
2;10,
2
)$,
irreps = $3_{5}^{1}
\hskip -1.5pt \otimes \hskip -1.5pt
1_{2}^{1,0}\oplus
2_{2}^{1,0}$,
pord$(\rho_\text{isum}(\mathfrak{t})) = 10$,

\vskip 0.7ex
\hangindent=5.5em \hangafter=1
{\white .}\hskip 1em $\rho_\text{isum}(\mathfrak{t})$ =
 $( \frac{1}{2},
\frac{3}{10},
\frac{7}{10} )
\oplus
( 0,
\frac{1}{2} )
$,

\vskip 0.7ex
\hangindent=5.5em \hangafter=1
{\white .}\hskip 1em $\rho_\text{isum}(\mathfrak{s})$ =
($-\sqrt{\frac{1}{5}}$,
$-\sqrt{\frac{2}{5}}$,
$-\sqrt{\frac{2}{5}}$;
$\frac{5+\sqrt{5}}{10}$,
$\frac{-5+\sqrt{5}}{10}$;
$\frac{5+\sqrt{5}}{10}$)
 $\oplus$
($-\frac{1}{2}$,
$-\sqrt{\frac{3}{4}}$;
$\frac{1}{2}$)

Pass. 

 \ \color{black}

\noindent 47: (dims,levels) = $(3\oplus
2;10,
6
)$,
irreps = $3_{5}^{3}
\hskip -1.5pt \otimes \hskip -1.5pt
1_{2}^{1,0}\oplus
2_{3}^{1,8}
\hskip -1.5pt \otimes \hskip -1.5pt
1_{2}^{1,0}$,
pord$(\rho_\text{isum}(\mathfrak{t})) = 15$,

\vskip 0.7ex
\hangindent=5.5em \hangafter=1
{\white .}\hskip 1em $\rho_\text{isum}(\mathfrak{t})$ =
 $( \frac{1}{2},
\frac{1}{10},
\frac{9}{10} )
\oplus
( \frac{1}{2},
\frac{1}{6} )
$,

\vskip 0.7ex
\hangindent=5.5em \hangafter=1
{\white .}\hskip 1em $\rho_\text{isum}(\mathfrak{s})$ =
($\sqrt{\frac{1}{5}}$,
$-\sqrt{\frac{2}{5}}$,
$-\sqrt{\frac{2}{5}}$;
$\frac{5-\sqrt{5}}{10}$,
$-\frac{5+\sqrt{5}}{10}$;
$\frac{5-\sqrt{5}}{10}$)
 $\oplus$
$\mathrm{i}$($-\sqrt{\frac{1}{3}}$,
$\sqrt{\frac{2}{3}}$;\ \ 
$\sqrt{\frac{1}{3}}$)

Fail:
Tr$_I(C) = -1 <$  0 for I = [ 1/6 ]. Prop. B.4 (1) eqn. (B.18)

 \ \color{black}

\noindent 48: (dims,levels) = $(3\oplus
2;10,
6
)$,
irreps = $3_{5}^{3}
\hskip -1.5pt \otimes \hskip -1.5pt
1_{2}^{1,0}\oplus
2_{3}^{1,0}
\hskip -1.5pt \otimes \hskip -1.5pt
1_{2}^{1,0}$,
pord$(\rho_\text{isum}(\mathfrak{t})) = 15$,

\vskip 0.7ex
\hangindent=5.5em \hangafter=1
{\white .}\hskip 1em $\rho_\text{isum}(\mathfrak{t})$ =
 $( \frac{1}{2},
\frac{1}{10},
\frac{9}{10} )
\oplus
( \frac{1}{2},
\frac{5}{6} )
$,

\vskip 0.7ex
\hangindent=5.5em \hangafter=1
{\white .}\hskip 1em $\rho_\text{isum}(\mathfrak{s})$ =
($\sqrt{\frac{1}{5}}$,
$-\sqrt{\frac{2}{5}}$,
$-\sqrt{\frac{2}{5}}$;
$\frac{5-\sqrt{5}}{10}$,
$-\frac{5+\sqrt{5}}{10}$;
$\frac{5-\sqrt{5}}{10}$)
 $\oplus$
$\mathrm{i}$($\sqrt{\frac{1}{3}}$,
$\sqrt{\frac{2}{3}}$;\ \ 
$-\sqrt{\frac{1}{3}}$)

Fail:
Tr$_I(C) = -1 <$  0 for I = [ 5/6 ]. Prop. B.4 (1) eqn. (B.18)

 \ \color{black}

\noindent 49: (dims,levels) = $(3\oplus
2;10,
6
)$,
irreps = $3_{5}^{1}
\hskip -1.5pt \otimes \hskip -1.5pt
1_{2}^{1,0}\oplus
2_{3}^{1,8}
\hskip -1.5pt \otimes \hskip -1.5pt
1_{2}^{1,0}$,
pord$(\rho_\text{isum}(\mathfrak{t})) = 15$,

\vskip 0.7ex
\hangindent=5.5em \hangafter=1
{\white .}\hskip 1em $\rho_\text{isum}(\mathfrak{t})$ =
 $( \frac{1}{2},
\frac{3}{10},
\frac{7}{10} )
\oplus
( \frac{1}{2},
\frac{1}{6} )
$,

\vskip 0.7ex
\hangindent=5.5em \hangafter=1
{\white .}\hskip 1em $\rho_\text{isum}(\mathfrak{s})$ =
($-\sqrt{\frac{1}{5}}$,
$-\sqrt{\frac{2}{5}}$,
$-\sqrt{\frac{2}{5}}$;
$\frac{5+\sqrt{5}}{10}$,
$\frac{-5+\sqrt{5}}{10}$;
$\frac{5+\sqrt{5}}{10}$)
 $\oplus$
$\mathrm{i}$($-\sqrt{\frac{1}{3}}$,
$\sqrt{\frac{2}{3}}$;\ \ 
$\sqrt{\frac{1}{3}}$)

Fail:
Tr$_I(C) = -1 <$  0 for I = [ 1/6 ]. Prop. B.4 (1) eqn. (B.18)

 \ \color{black}

\noindent 50: (dims,levels) = $(3\oplus
2;10,
6
)$,
irreps = $3_{5}^{1}
\hskip -1.5pt \otimes \hskip -1.5pt
1_{2}^{1,0}\oplus
2_{3}^{1,0}
\hskip -1.5pt \otimes \hskip -1.5pt
1_{2}^{1,0}$,
pord$(\rho_\text{isum}(\mathfrak{t})) = 15$,

\vskip 0.7ex
\hangindent=5.5em \hangafter=1
{\white .}\hskip 1em $\rho_\text{isum}(\mathfrak{t})$ =
 $( \frac{1}{2},
\frac{3}{10},
\frac{7}{10} )
\oplus
( \frac{1}{2},
\frac{5}{6} )
$,

\vskip 0.7ex
\hangindent=5.5em \hangafter=1
{\white .}\hskip 1em $\rho_\text{isum}(\mathfrak{s})$ =
($-\sqrt{\frac{1}{5}}$,
$-\sqrt{\frac{2}{5}}$,
$-\sqrt{\frac{2}{5}}$;
$\frac{5+\sqrt{5}}{10}$,
$\frac{-5+\sqrt{5}}{10}$;
$\frac{5+\sqrt{5}}{10}$)
 $\oplus$
$\mathrm{i}$($\sqrt{\frac{1}{3}}$,
$\sqrt{\frac{2}{3}}$;\ \ 
$-\sqrt{\frac{1}{3}}$)

Fail:
Tr$_I(C) = -1 <$  0 for I = [ 5/6 ]. Prop. B.4 (1) eqn. (B.18)

 \ \color{black}

\noindent 51: (dims,levels) = $(3\oplus
2;10,
10
)$,
irreps = $3_{5}^{3}
\hskip -1.5pt \otimes \hskip -1.5pt
1_{2}^{1,0}\oplus
2_{5}^{2}
\hskip -1.5pt \otimes \hskip -1.5pt
1_{2}^{1,0}$,
pord$(\rho_\text{isum}(\mathfrak{t})) = 5$,

\vskip 0.7ex
\hangindent=5.5em \hangafter=1
{\white .}\hskip 1em $\rho_\text{isum}(\mathfrak{t})$ =
 $( \frac{1}{2},
\frac{1}{10},
\frac{9}{10} )
\oplus
( \frac{1}{10},
\frac{9}{10} )
$,

\vskip 0.7ex
\hangindent=5.5em \hangafter=1
{\white .}\hskip 1em $\rho_\text{isum}(\mathfrak{s})$ =
($\sqrt{\frac{1}{5}}$,
$-\sqrt{\frac{2}{5}}$,
$-\sqrt{\frac{2}{5}}$;
$\frac{5-\sqrt{5}}{10}$,
$-\frac{5+\sqrt{5}}{10}$;
$\frac{5-\sqrt{5}}{10}$)
 $\oplus$
$\mathrm{i}$($-\frac{1}{\sqrt{5}}c^{1}_{20}
$,
$\frac{1}{\sqrt{5}}c^{3}_{20}
$;\ \ 
$\frac{1}{\sqrt{5}}c^{1}_{20}
$)

Fail:
Integral: $D_{\rho}(\sigma)_{\theta} \propto $ id,
 for all $\sigma$ and all $\theta$-eigenspaces that can contain unit. Prop. B.5 (6)

 \ \color{black}

\noindent 52: (dims,levels) = $(3\oplus
2;10,
10
)$,
irreps = $3_{5}^{1}
\hskip -1.5pt \otimes \hskip -1.5pt
1_{2}^{1,0}\oplus
2_{5}^{1}
\hskip -1.5pt \otimes \hskip -1.5pt
1_{2}^{1,0}$,
pord$(\rho_\text{isum}(\mathfrak{t})) = 5$,

\vskip 0.7ex
\hangindent=5.5em \hangafter=1
{\white .}\hskip 1em $\rho_\text{isum}(\mathfrak{t})$ =
 $( \frac{1}{2},
\frac{3}{10},
\frac{7}{10} )
\oplus
( \frac{3}{10},
\frac{7}{10} )
$,

\vskip 0.7ex
\hangindent=5.5em \hangafter=1
{\white .}\hskip 1em $\rho_\text{isum}(\mathfrak{s})$ =
($-\sqrt{\frac{1}{5}}$,
$-\sqrt{\frac{2}{5}}$,
$-\sqrt{\frac{2}{5}}$;
$\frac{5+\sqrt{5}}{10}$,
$\frac{-5+\sqrt{5}}{10}$;
$\frac{5+\sqrt{5}}{10}$)
 $\oplus$
$\mathrm{i}$($-\frac{1}{\sqrt{5}}c^{3}_{20}
$,
$\frac{1}{\sqrt{5}}c^{1}_{20}
$;\ \ 
$\frac{1}{\sqrt{5}}c^{3}_{20}
$)

Fail:
Integral: $D_{\rho}(\sigma)_{\theta} \propto $ id,
 for all $\sigma$ and all $\theta$-eigenspaces that can contain unit. Prop. B.5 (6)

 \ \color{black}

\noindent 53: (dims,levels) = $(3\oplus
2;12,
3
)$,
irreps = $3_{4}^{1,0}
\hskip -1.5pt \otimes \hskip -1.5pt
1_{3}^{1,0}\oplus
2_{3}^{1,0}$,
pord$(\rho_\text{isum}(\mathfrak{t})) = 12$,

\vskip 0.7ex
\hangindent=5.5em \hangafter=1
{\white .}\hskip 1em $\rho_\text{isum}(\mathfrak{t})$ =
 $( \frac{1}{3},
\frac{1}{12},
\frac{7}{12} )
\oplus
( 0,
\frac{1}{3} )
$,

\vskip 0.7ex
\hangindent=5.5em \hangafter=1
{\white .}\hskip 1em $\rho_\text{isum}(\mathfrak{s})$ =
($0$,
$\sqrt{\frac{1}{2}}$,
$\sqrt{\frac{1}{2}}$;
$-\frac{1}{2}$,
$\frac{1}{2}$;
$-\frac{1}{2}$)
 $\oplus$
$\mathrm{i}$($-\sqrt{\frac{1}{3}}$,
$\sqrt{\frac{2}{3}}$;\ \ 
$\sqrt{\frac{1}{3}}$)

Fail:
Tr$_I(C) = -1 <$  0 for I = [ 0 ]. Prop. B.4 (1) eqn. (B.18)

 \ \color{black}

\noindent 54: (dims,levels) = $(3\oplus
2;12,
3
)$,
irreps = $3_{4}^{1,0}
\hskip -1.5pt \otimes \hskip -1.5pt
1_{3}^{1,0}\oplus
2_{3}^{1,4}$,
pord$(\rho_\text{isum}(\mathfrak{t})) = 12$,

\vskip 0.7ex
\hangindent=5.5em \hangafter=1
{\white .}\hskip 1em $\rho_\text{isum}(\mathfrak{t})$ =
 $( \frac{1}{3},
\frac{1}{12},
\frac{7}{12} )
\oplus
( \frac{1}{3},
\frac{2}{3} )
$,

\vskip 0.7ex
\hangindent=5.5em \hangafter=1
{\white .}\hskip 1em $\rho_\text{isum}(\mathfrak{s})$ =
($0$,
$\sqrt{\frac{1}{2}}$,
$\sqrt{\frac{1}{2}}$;
$-\frac{1}{2}$,
$\frac{1}{2}$;
$-\frac{1}{2}$)
 $\oplus$
$\mathrm{i}$($-\sqrt{\frac{1}{3}}$,
$\sqrt{\frac{2}{3}}$;\ \ 
$\sqrt{\frac{1}{3}}$)

Fail:
Tr$_I(C) = -1 <$  0 for I = [ 2/3 ]. Prop. B.4 (1) eqn. (B.18)

 \ \color{black}

\noindent 55: (dims,levels) = $(3\oplus
2;12,
12
)$,
irreps = $3_{4}^{1,0}
\hskip -1.5pt \otimes \hskip -1.5pt
1_{3}^{1,0}\oplus
2_{3}^{1,4}
\hskip -1.5pt \otimes \hskip -1.5pt
1_{4}^{1,6}$,
pord$(\rho_\text{isum}(\mathfrak{t})) = 12$,

\vskip 0.7ex
\hangindent=5.5em \hangafter=1
{\white .}\hskip 1em $\rho_\text{isum}(\mathfrak{t})$ =
 $( \frac{1}{3},
\frac{1}{12},
\frac{7}{12} )
\oplus
( \frac{1}{12},
\frac{5}{12} )
$,

\vskip 0.7ex
\hangindent=5.5em \hangafter=1
{\white .}\hskip 1em $\rho_\text{isum}(\mathfrak{s})$ =
($0$,
$\sqrt{\frac{1}{2}}$,
$\sqrt{\frac{1}{2}}$;
$-\frac{1}{2}$,
$\frac{1}{2}$;
$-\frac{1}{2}$)
 $\oplus$
($-\sqrt{\frac{1}{3}}$,
$\sqrt{\frac{2}{3}}$;
$\sqrt{\frac{1}{3}}$)

Fail:
cnd($\rho(\mathfrak s)_\mathrm{ndeg}$) = 24 does not divide
 ord($\rho(\mathfrak t)$)=12. Prop. B.4 (2)

 \ \color{black}

\noindent 56: (dims,levels) = $(3\oplus
2;12,
12
)$,
irreps = $3_{4}^{1,0}
\hskip -1.5pt \otimes \hskip -1.5pt
1_{3}^{1,0}\oplus
2_{3}^{1,0}
\hskip -1.5pt \otimes \hskip -1.5pt
1_{4}^{1,0}$,
pord$(\rho_\text{isum}(\mathfrak{t})) = 12$,

\vskip 0.7ex
\hangindent=5.5em \hangafter=1
{\white .}\hskip 1em $\rho_\text{isum}(\mathfrak{t})$ =
 $( \frac{1}{3},
\frac{1}{12},
\frac{7}{12} )
\oplus
( \frac{1}{4},
\frac{7}{12} )
$,

\vskip 0.7ex
\hangindent=5.5em \hangafter=1
{\white .}\hskip 1em $\rho_\text{isum}(\mathfrak{s})$ =
($0$,
$\sqrt{\frac{1}{2}}$,
$\sqrt{\frac{1}{2}}$;
$-\frac{1}{2}$,
$\frac{1}{2}$;
$-\frac{1}{2}$)
 $\oplus$
($\sqrt{\frac{1}{3}}$,
$\sqrt{\frac{2}{3}}$;
$-\sqrt{\frac{1}{3}}$)

Fail:
cnd($\rho(\mathfrak s)_\mathrm{ndeg}$) = 24 does not divide
 ord($\rho(\mathfrak t)$)=12. Prop. B.4 (2)

 \ \color{black}

\noindent 57: (dims,levels) = $(3\oplus
2;12,
12
)$,
irreps = $3_{4}^{1,0}
\hskip -1.5pt \otimes \hskip -1.5pt
1_{3}^{1,0}\oplus
2_{3}^{1,4}
\hskip -1.5pt \otimes \hskip -1.5pt
1_{4}^{1,0}$,
pord$(\rho_\text{isum}(\mathfrak{t})) = 12$,

\vskip 0.7ex
\hangindent=5.5em \hangafter=1
{\white .}\hskip 1em $\rho_\text{isum}(\mathfrak{t})$ =
 $( \frac{1}{3},
\frac{1}{12},
\frac{7}{12} )
\oplus
( \frac{7}{12},
\frac{11}{12} )
$,

\vskip 0.7ex
\hangindent=5.5em \hangafter=1
{\white .}\hskip 1em $\rho_\text{isum}(\mathfrak{s})$ =
($0$,
$\sqrt{\frac{1}{2}}$,
$\sqrt{\frac{1}{2}}$;
$-\frac{1}{2}$,
$\frac{1}{2}$;
$-\frac{1}{2}$)
 $\oplus$
($\sqrt{\frac{1}{3}}$,
$\sqrt{\frac{2}{3}}$;
$-\sqrt{\frac{1}{3}}$)

Fail:
cnd($\rho(\mathfrak s)_\mathrm{ndeg}$) = 24 does not divide
 ord($\rho(\mathfrak t)$)=12. Prop. B.4 (2)

 \ \color{black}

\noindent 58: (dims,levels) = $(3\oplus
2;12,
12
)$,
irreps = $3_{4}^{1,0}
\hskip -1.5pt \otimes \hskip -1.5pt
1_{3}^{1,0}\oplus
2_{3}^{1,0}
\hskip -1.5pt \otimes \hskip -1.5pt
1_{4}^{1,6}$,
pord$(\rho_\text{isum}(\mathfrak{t})) = 12$,

\vskip 0.7ex
\hangindent=5.5em \hangafter=1
{\white .}\hskip 1em $\rho_\text{isum}(\mathfrak{t})$ =
 $( \frac{1}{3},
\frac{1}{12},
\frac{7}{12} )
\oplus
( \frac{3}{4},
\frac{1}{12} )
$,

\vskip 0.7ex
\hangindent=5.5em \hangafter=1
{\white .}\hskip 1em $\rho_\text{isum}(\mathfrak{s})$ =
($0$,
$\sqrt{\frac{1}{2}}$,
$\sqrt{\frac{1}{2}}$;
$-\frac{1}{2}$,
$\frac{1}{2}$;
$-\frac{1}{2}$)
 $\oplus$
($-\sqrt{\frac{1}{3}}$,
$\sqrt{\frac{2}{3}}$;
$\sqrt{\frac{1}{3}}$)

Fail:
cnd($\rho(\mathfrak s)_\mathrm{ndeg}$) = 24 does not divide
 ord($\rho(\mathfrak t)$)=12. Prop. B.4 (2)

 \ \color{black}

\noindent 59: (dims,levels) = $(3\oplus
2;15,
3
)$,
irreps = $3_{5}^{1}
\hskip -1.5pt \otimes \hskip -1.5pt
1_{3}^{1,0}\oplus
2_{3}^{1,0}$,
pord$(\rho_\text{isum}(\mathfrak{t})) = 15$,

\vskip 0.7ex
\hangindent=5.5em \hangafter=1
{\white .}\hskip 1em $\rho_\text{isum}(\mathfrak{t})$ =
 $( \frac{1}{3},
\frac{2}{15},
\frac{8}{15} )
\oplus
( 0,
\frac{1}{3} )
$,

\vskip 0.7ex
\hangindent=5.5em \hangafter=1
{\white .}\hskip 1em $\rho_\text{isum}(\mathfrak{s})$ =
($\sqrt{\frac{1}{5}}$,
$-\sqrt{\frac{2}{5}}$,
$-\sqrt{\frac{2}{5}}$;
$-\frac{5+\sqrt{5}}{10}$,
$\frac{5-\sqrt{5}}{10}$;
$-\frac{5+\sqrt{5}}{10}$)
 $\oplus$
$\mathrm{i}$($-\sqrt{\frac{1}{3}}$,
$\sqrt{\frac{2}{3}}$;\ \ 
$\sqrt{\frac{1}{3}}$)

Fail:
Tr$_I(C) = -1 <$  0 for I = [ 0 ]. Prop. B.4 (1) eqn. (B.18)

 \ \color{black}

\noindent 60: (dims,levels) = $(3\oplus
2;15,
3
)$,
irreps = $3_{5}^{1}
\hskip -1.5pt \otimes \hskip -1.5pt
1_{3}^{1,0}\oplus
2_{3}^{1,4}$,
pord$(\rho_\text{isum}(\mathfrak{t})) = 15$,

\vskip 0.7ex
\hangindent=5.5em \hangafter=1
{\white .}\hskip 1em $\rho_\text{isum}(\mathfrak{t})$ =
 $( \frac{1}{3},
\frac{2}{15},
\frac{8}{15} )
\oplus
( \frac{1}{3},
\frac{2}{3} )
$,

\vskip 0.7ex
\hangindent=5.5em \hangafter=1
{\white .}\hskip 1em $\rho_\text{isum}(\mathfrak{s})$ =
($\sqrt{\frac{1}{5}}$,
$-\sqrt{\frac{2}{5}}$,
$-\sqrt{\frac{2}{5}}$;
$-\frac{5+\sqrt{5}}{10}$,
$\frac{5-\sqrt{5}}{10}$;
$-\frac{5+\sqrt{5}}{10}$)
 $\oplus$
$\mathrm{i}$($-\sqrt{\frac{1}{3}}$,
$\sqrt{\frac{2}{3}}$;\ \ 
$\sqrt{\frac{1}{3}}$)

Fail:
Tr$_I(C) = -1 <$  0 for I = [ 2/3 ]. Prop. B.4 (1) eqn. (B.18)

 \ \color{black}

\noindent 61: (dims,levels) = $(3\oplus
2;15,
3
)$,
irreps = $3_{5}^{3}
\hskip -1.5pt \otimes \hskip -1.5pt
1_{3}^{1,0}\oplus
2_{3}^{1,0}$,
pord$(\rho_\text{isum}(\mathfrak{t})) = 15$,

\vskip 0.7ex
\hangindent=5.5em \hangafter=1
{\white .}\hskip 1em $\rho_\text{isum}(\mathfrak{t})$ =
 $( \frac{1}{3},
\frac{11}{15},
\frac{14}{15} )
\oplus
( 0,
\frac{1}{3} )
$,

\vskip 0.7ex
\hangindent=5.5em \hangafter=1
{\white .}\hskip 1em $\rho_\text{isum}(\mathfrak{s})$ =
($-\sqrt{\frac{1}{5}}$,
$-\sqrt{\frac{2}{5}}$,
$-\sqrt{\frac{2}{5}}$;
$\frac{-5+\sqrt{5}}{10}$,
$\frac{5+\sqrt{5}}{10}$;
$\frac{-5+\sqrt{5}}{10}$)
 $\oplus$
$\mathrm{i}$($-\sqrt{\frac{1}{3}}$,
$\sqrt{\frac{2}{3}}$;\ \ 
$\sqrt{\frac{1}{3}}$)

Fail:
Tr$_I(C) = -1 <$  0 for I = [ 0 ]. Prop. B.4 (1) eqn. (B.18)

 \ \color{black}

\noindent 62: (dims,levels) = $(3\oplus
2;15,
3
)$,
irreps = $3_{5}^{3}
\hskip -1.5pt \otimes \hskip -1.5pt
1_{3}^{1,0}\oplus
2_{3}^{1,4}$,
pord$(\rho_\text{isum}(\mathfrak{t})) = 15$,

\vskip 0.7ex
\hangindent=5.5em \hangafter=1
{\white .}\hskip 1em $\rho_\text{isum}(\mathfrak{t})$ =
 $( \frac{1}{3},
\frac{11}{15},
\frac{14}{15} )
\oplus
( \frac{1}{3},
\frac{2}{3} )
$,

\vskip 0.7ex
\hangindent=5.5em \hangafter=1
{\white .}\hskip 1em $\rho_\text{isum}(\mathfrak{s})$ =
($-\sqrt{\frac{1}{5}}$,
$-\sqrt{\frac{2}{5}}$,
$-\sqrt{\frac{2}{5}}$;
$\frac{-5+\sqrt{5}}{10}$,
$\frac{5+\sqrt{5}}{10}$;
$\frac{-5+\sqrt{5}}{10}$)
 $\oplus$
$\mathrm{i}$($-\sqrt{\frac{1}{3}}$,
$\sqrt{\frac{2}{3}}$;\ \ 
$\sqrt{\frac{1}{3}}$)

Fail:
Tr$_I(C) = -1 <$  0 for I = [ 2/3 ]. Prop. B.4 (1) eqn. (B.18)

 \ \color{black}

 \color{blue}

\noindent 63: (dims,levels) = $(3\oplus
2;15,
6
)$,
irreps = $3_{5}^{1}
\hskip -1.5pt \otimes \hskip -1.5pt
1_{3}^{1,0}\oplus
2_{2}^{1,0}
\hskip -1.5pt \otimes \hskip -1.5pt
1_{3}^{1,0}$,
pord$(\rho_\text{isum}(\mathfrak{t})) = 10$,

\vskip 0.7ex
\hangindent=5.5em \hangafter=1
{\white .}\hskip 1em $\rho_\text{isum}(\mathfrak{t})$ =
 $( \frac{1}{3},
\frac{2}{15},
\frac{8}{15} )
\oplus
( \frac{1}{3},
\frac{5}{6} )
$,

\vskip 0.7ex
\hangindent=5.5em \hangafter=1
{\white .}\hskip 1em $\rho_\text{isum}(\mathfrak{s})$ =
($\sqrt{\frac{1}{5}}$,
$-\sqrt{\frac{2}{5}}$,
$-\sqrt{\frac{2}{5}}$;
$-\frac{5+\sqrt{5}}{10}$,
$\frac{5-\sqrt{5}}{10}$;
$-\frac{5+\sqrt{5}}{10}$)
 $\oplus$
($-\frac{1}{2}$,
$-\sqrt{\frac{3}{4}}$;
$\frac{1}{2}$)

Pass. 

 \ \color{black}

 \color{blue}

\noindent 64: (dims,levels) = $(3\oplus
2;15,
6
)$,
irreps = $3_{5}^{3}
\hskip -1.5pt \otimes \hskip -1.5pt
1_{3}^{1,0}\oplus
2_{2}^{1,0}
\hskip -1.5pt \otimes \hskip -1.5pt
1_{3}^{1,0}$,
pord$(\rho_\text{isum}(\mathfrak{t})) = 10$,

\vskip 0.7ex
\hangindent=5.5em \hangafter=1
{\white .}\hskip 1em $\rho_\text{isum}(\mathfrak{t})$ =
 $( \frac{1}{3},
\frac{11}{15},
\frac{14}{15} )
\oplus
( \frac{1}{3},
\frac{5}{6} )
$,

\vskip 0.7ex
\hangindent=5.5em \hangafter=1
{\white .}\hskip 1em $\rho_\text{isum}(\mathfrak{s})$ =
($-\sqrt{\frac{1}{5}}$,
$-\sqrt{\frac{2}{5}}$,
$-\sqrt{\frac{2}{5}}$;
$\frac{-5+\sqrt{5}}{10}$,
$\frac{5+\sqrt{5}}{10}$;
$\frac{-5+\sqrt{5}}{10}$)
 $\oplus$
($-\frac{1}{2}$,
$-\sqrt{\frac{3}{4}}$;
$\frac{1}{2}$)

Pass. 

 \ \color{black}

\noindent 65: (dims,levels) = $(3\oplus
2;15,
15
)$,
irreps = $3_{5}^{1}
\hskip -1.5pt \otimes \hskip -1.5pt
1_{3}^{1,0}\oplus
2_{5}^{1}
\hskip -1.5pt \otimes \hskip -1.5pt
1_{3}^{1,0}$,
pord$(\rho_\text{isum}(\mathfrak{t})) = 5$,

\vskip 0.7ex
\hangindent=5.5em \hangafter=1
{\white .}\hskip 1em $\rho_\text{isum}(\mathfrak{t})$ =
 $( \frac{1}{3},
\frac{2}{15},
\frac{8}{15} )
\oplus
( \frac{2}{15},
\frac{8}{15} )
$,

\vskip 0.7ex
\hangindent=5.5em \hangafter=1
{\white .}\hskip 1em $\rho_\text{isum}(\mathfrak{s})$ =
($\sqrt{\frac{1}{5}}$,
$-\sqrt{\frac{2}{5}}$,
$-\sqrt{\frac{2}{5}}$;
$-\frac{5+\sqrt{5}}{10}$,
$\frac{5-\sqrt{5}}{10}$;
$-\frac{5+\sqrt{5}}{10}$)
 $\oplus$
$\mathrm{i}$($\frac{1}{\sqrt{5}}c^{3}_{20}
$,
$\frac{1}{\sqrt{5}}c^{1}_{20}
$;\ \ 
$-\frac{1}{\sqrt{5}}c^{3}_{20}
$)

Fail:
Integral: $D_{\rho}(\sigma)_{\theta} \propto $ id,
 for all $\sigma$ and all $\theta$-eigenspaces that can contain unit. Prop. B.5 (6)

 \ \color{black}

\noindent 66: (dims,levels) = $(3\oplus
2;15,
15
)$,
irreps = $3_{5}^{3}
\hskip -1.5pt \otimes \hskip -1.5pt
1_{3}^{1,0}\oplus
2_{5}^{2}
\hskip -1.5pt \otimes \hskip -1.5pt
1_{3}^{1,0}$,
pord$(\rho_\text{isum}(\mathfrak{t})) = 5$,

\vskip 0.7ex
\hangindent=5.5em \hangafter=1
{\white .}\hskip 1em $\rho_\text{isum}(\mathfrak{t})$ =
 $( \frac{1}{3},
\frac{11}{15},
\frac{14}{15} )
\oplus
( \frac{11}{15},
\frac{14}{15} )
$,

\vskip 0.7ex
\hangindent=5.5em \hangafter=1
{\white .}\hskip 1em $\rho_\text{isum}(\mathfrak{s})$ =
($-\sqrt{\frac{1}{5}}$,
$-\sqrt{\frac{2}{5}}$,
$-\sqrt{\frac{2}{5}}$;
$\frac{-5+\sqrt{5}}{10}$,
$\frac{5+\sqrt{5}}{10}$;
$\frac{-5+\sqrt{5}}{10}$)
 $\oplus$
$\mathrm{i}$($-\frac{1}{\sqrt{5}}c^{1}_{20}
$,
$\frac{1}{\sqrt{5}}c^{3}_{20}
$;\ \ 
$\frac{1}{\sqrt{5}}c^{1}_{20}
$)

Fail:
Integral: $D_{\rho}(\sigma)_{\theta} \propto $ id,
 for all $\sigma$ and all $\theta$-eigenspaces that can contain unit. Prop. B.5 (6)

 \ \color{black}

\noindent 67: (dims,levels) = $(3\oplus
2;16,
8
)$,
irreps = $3_{16}^{1,0}\oplus
2_{8}^{1,0}$,
pord$(\rho_\text{isum}(\mathfrak{t})) = 16$,

\vskip 0.7ex
\hangindent=5.5em \hangafter=1
{\white .}\hskip 1em $\rho_\text{isum}(\mathfrak{t})$ =
 $( \frac{1}{8},
\frac{1}{16},
\frac{9}{16} )
\oplus
( \frac{1}{8},
\frac{3}{8} )
$,

\vskip 0.7ex
\hangindent=5.5em \hangafter=1
{\white .}\hskip 1em $\rho_\text{isum}(\mathfrak{s})$ =
$\mathrm{i}$($0$,
$\sqrt{\frac{1}{2}}$,
$\sqrt{\frac{1}{2}}$;\ \ 
$-\frac{1}{2}$,
$\frac{1}{2}$;\ \ 
$-\frac{1}{2}$)
 $\oplus$
($-\sqrt{\frac{1}{2}}$,
$\sqrt{\frac{1}{2}}$;
$\sqrt{\frac{1}{2}}$)

Fail:
Tr$_I(C) = -1 <$  0 for I = [ 3/8 ]. Prop. B.4 (1) eqn. (B.18)

 \ \color{black}

\noindent 68: (dims,levels) = $(3\oplus
2;16,
8
)$,
irreps = $3_{16}^{1,0}\oplus
2_{8}^{1,9}$,
pord$(\rho_\text{isum}(\mathfrak{t})) = 16$,

\vskip 0.7ex
\hangindent=5.5em \hangafter=1
{\white .}\hskip 1em $\rho_\text{isum}(\mathfrak{t})$ =
 $( \frac{1}{8},
\frac{1}{16},
\frac{9}{16} )
\oplus
( \frac{1}{8},
\frac{7}{8} )
$,

\vskip 0.7ex
\hangindent=5.5em \hangafter=1
{\white .}\hskip 1em $\rho_\text{isum}(\mathfrak{s})$ =
$\mathrm{i}$($0$,
$\sqrt{\frac{1}{2}}$,
$\sqrt{\frac{1}{2}}$;\ \ 
$-\frac{1}{2}$,
$\frac{1}{2}$;\ \ 
$-\frac{1}{2}$)
 $\oplus$
$\mathrm{i}$($-\sqrt{\frac{1}{2}}$,
$\sqrt{\frac{1}{2}}$;\ \ 
$\sqrt{\frac{1}{2}}$)

Fail:
Integral: $D_{\rho}(\sigma)_{\theta} \propto $ id,
 for all $\sigma$ and all $\theta$-eigenspaces that can contain unit. Prop. B.5 (6)

 \ \color{black}

\noindent 69: (dims,levels) = $(3\oplus
2;16,
8
)$,
irreps = $3_{16}^{3,0}\oplus
2_{8}^{1,0}$,
pord$(\rho_\text{isum}(\mathfrak{t})) = 16$,

\vskip 0.7ex
\hangindent=5.5em \hangafter=1
{\white .}\hskip 1em $\rho_\text{isum}(\mathfrak{t})$ =
 $( \frac{3}{8},
\frac{3}{16},
\frac{11}{16} )
\oplus
( \frac{1}{8},
\frac{3}{8} )
$,

\vskip 0.7ex
\hangindent=5.5em \hangafter=1
{\white .}\hskip 1em $\rho_\text{isum}(\mathfrak{s})$ =
$\mathrm{i}$($0$,
$\sqrt{\frac{1}{2}}$,
$\sqrt{\frac{1}{2}}$;\ \ 
$\frac{1}{2}$,
$-\frac{1}{2}$;\ \ 
$\frac{1}{2}$)
 $\oplus$
($-\sqrt{\frac{1}{2}}$,
$\sqrt{\frac{1}{2}}$;
$\sqrt{\frac{1}{2}}$)

Fail:
Tr$_I(C) = -1 <$  0 for I = [ 1/8 ]. Prop. B.4 (1) eqn. (B.18)

 \ \color{black}

\noindent 70: (dims,levels) = $(3\oplus
2;16,
8
)$,
irreps = $3_{16}^{3,0}\oplus
2_{8}^{1,3}$,
pord$(\rho_\text{isum}(\mathfrak{t})) = 16$,

\vskip 0.7ex
\hangindent=5.5em \hangafter=1
{\white .}\hskip 1em $\rho_\text{isum}(\mathfrak{t})$ =
 $( \frac{3}{8},
\frac{3}{16},
\frac{11}{16} )
\oplus
( \frac{3}{8},
\frac{5}{8} )
$,

\vskip 0.7ex
\hangindent=5.5em \hangafter=1
{\white .}\hskip 1em $\rho_\text{isum}(\mathfrak{s})$ =
$\mathrm{i}$($0$,
$\sqrt{\frac{1}{2}}$,
$\sqrt{\frac{1}{2}}$;\ \ 
$\frac{1}{2}$,
$-\frac{1}{2}$;\ \ 
$\frac{1}{2}$)
 $\oplus$
$\mathrm{i}$($-\sqrt{\frac{1}{2}}$,
$\sqrt{\frac{1}{2}}$;\ \ 
$\sqrt{\frac{1}{2}}$)

Fail:
Integral: $D_{\rho}(\sigma)_{\theta} \propto $ id,
 for all $\sigma$ and all $\theta$-eigenspaces that can contain unit. Prop. B.5 (6)

 \ \color{black}

\noindent 71: (dims,levels) = $(3\oplus
2;16,
8
)$,
irreps = $3_{16}^{5,0}\oplus
2_{8}^{1,3}$,
pord$(\rho_\text{isum}(\mathfrak{t})) = 16$,

\vskip 0.7ex
\hangindent=5.5em \hangafter=1
{\white .}\hskip 1em $\rho_\text{isum}(\mathfrak{t})$ =
 $( \frac{5}{8},
\frac{5}{16},
\frac{13}{16} )
\oplus
( \frac{3}{8},
\frac{5}{8} )
$,

\vskip 0.7ex
\hangindent=5.5em \hangafter=1
{\white .}\hskip 1em $\rho_\text{isum}(\mathfrak{s})$ =
$\mathrm{i}$($0$,
$\sqrt{\frac{1}{2}}$,
$\sqrt{\frac{1}{2}}$;\ \ 
$-\frac{1}{2}$,
$\frac{1}{2}$;\ \ 
$-\frac{1}{2}$)
 $\oplus$
$\mathrm{i}$($-\sqrt{\frac{1}{2}}$,
$\sqrt{\frac{1}{2}}$;\ \ 
$\sqrt{\frac{1}{2}}$)

Fail:
Integral: $D_{\rho}(\sigma)_{\theta} \propto $ id,
 for all $\sigma$ and all $\theta$-eigenspaces that can contain unit. Prop. B.5 (6)

 \ \color{black}

\noindent 72: (dims,levels) = $(3\oplus
2;16,
8
)$,
irreps = $3_{16}^{5,0}\oplus
2_{8}^{1,6}$,
pord$(\rho_\text{isum}(\mathfrak{t})) = 16$,

\vskip 0.7ex
\hangindent=5.5em \hangafter=1
{\white .}\hskip 1em $\rho_\text{isum}(\mathfrak{t})$ =
 $( \frac{5}{8},
\frac{5}{16},
\frac{13}{16} )
\oplus
( \frac{5}{8},
\frac{7}{8} )
$,

\vskip 0.7ex
\hangindent=5.5em \hangafter=1
{\white .}\hskip 1em $\rho_\text{isum}(\mathfrak{s})$ =
$\mathrm{i}$($0$,
$\sqrt{\frac{1}{2}}$,
$\sqrt{\frac{1}{2}}$;\ \ 
$-\frac{1}{2}$,
$\frac{1}{2}$;\ \ 
$-\frac{1}{2}$)
 $\oplus$
($\sqrt{\frac{1}{2}}$,
$\sqrt{\frac{1}{2}}$;
$-\sqrt{\frac{1}{2}}$)

Fail:
Tr$_I(C) = -1 <$  0 for I = [ 7/8 ]. Prop. B.4 (1) eqn. (B.18)

 \ \color{black}

\noindent 73: (dims,levels) = $(3\oplus
2;16,
8
)$,
irreps = $3_{16}^{7,0}\oplus
2_{8}^{1,9}$,
pord$(\rho_\text{isum}(\mathfrak{t})) = 16$,

\vskip 0.7ex
\hangindent=5.5em \hangafter=1
{\white .}\hskip 1em $\rho_\text{isum}(\mathfrak{t})$ =
 $( \frac{7}{8},
\frac{7}{16},
\frac{15}{16} )
\oplus
( \frac{1}{8},
\frac{7}{8} )
$,

\vskip 0.7ex
\hangindent=5.5em \hangafter=1
{\white .}\hskip 1em $\rho_\text{isum}(\mathfrak{s})$ =
$\mathrm{i}$($0$,
$\sqrt{\frac{1}{2}}$,
$\sqrt{\frac{1}{2}}$;\ \ 
$\frac{1}{2}$,
$-\frac{1}{2}$;\ \ 
$\frac{1}{2}$)
 $\oplus$
$\mathrm{i}$($-\sqrt{\frac{1}{2}}$,
$\sqrt{\frac{1}{2}}$;\ \ 
$\sqrt{\frac{1}{2}}$)

Fail:
Integral: $D_{\rho}(\sigma)_{\theta} \propto $ id,
 for all $\sigma$ and all $\theta$-eigenspaces that can contain unit. Prop. B.5 (6)

 \ \color{black}

\noindent 74: (dims,levels) = $(3\oplus
2;16,
8
)$,
irreps = $3_{16}^{7,0}\oplus
2_{8}^{1,6}$,
pord$(\rho_\text{isum}(\mathfrak{t})) = 16$,

\vskip 0.7ex
\hangindent=5.5em \hangafter=1
{\white .}\hskip 1em $\rho_\text{isum}(\mathfrak{t})$ =
 $( \frac{7}{8},
\frac{7}{16},
\frac{15}{16} )
\oplus
( \frac{5}{8},
\frac{7}{8} )
$,

\vskip 0.7ex
\hangindent=5.5em \hangafter=1
{\white .}\hskip 1em $\rho_\text{isum}(\mathfrak{s})$ =
$\mathrm{i}$($0$,
$\sqrt{\frac{1}{2}}$,
$\sqrt{\frac{1}{2}}$;\ \ 
$\frac{1}{2}$,
$-\frac{1}{2}$;\ \ 
$\frac{1}{2}$)
 $\oplus$
($\sqrt{\frac{1}{2}}$,
$\sqrt{\frac{1}{2}}$;
$-\sqrt{\frac{1}{2}}$)

Fail:
Tr$_I(C) = -1 <$  0 for I = [ 5/8 ]. Prop. B.4 (1) eqn. (B.18)

 \ \color{black}

 \color{blue}

\noindent 75: (dims,levels) = $(3\oplus
2;20,
4
)$,
irreps = $3_{5}^{1}
\hskip -1.5pt \otimes \hskip -1.5pt
1_{4}^{1,0}\oplus
2_{4}^{1,0}$,
pord$(\rho_\text{isum}(\mathfrak{t})) = 10$,

\vskip 0.7ex
\hangindent=5.5em \hangafter=1
{\white .}\hskip 1em $\rho_\text{isum}(\mathfrak{t})$ =
 $( \frac{1}{4},
\frac{1}{20},
\frac{9}{20} )
\oplus
( \frac{1}{4},
\frac{3}{4} )
$,

\vskip 0.7ex
\hangindent=5.5em \hangafter=1
{\white .}\hskip 1em $\rho_\text{isum}(\mathfrak{s})$ =
$\mathrm{i}$($\sqrt{\frac{1}{5}}$,
$\sqrt{\frac{2}{5}}$,
$\sqrt{\frac{2}{5}}$;\ \ 
$-\frac{5+\sqrt{5}}{10}$,
$\frac{5-\sqrt{5}}{10}$;\ \ 
$-\frac{5+\sqrt{5}}{10}$)
 $\oplus$
$\mathrm{i}$($-\frac{1}{2}$,
$\sqrt{\frac{3}{4}}$;\ \ 
$\frac{1}{2}$)

Pass. 

 \ \color{black}

 \color{blue}

\noindent 76: (dims,levels) = $(3\oplus
2;20,
4
)$,
irreps = $3_{5}^{3}
\hskip -1.5pt \otimes \hskip -1.5pt
1_{4}^{1,0}\oplus
2_{4}^{1,0}$,
pord$(\rho_\text{isum}(\mathfrak{t})) = 10$,

\vskip 0.7ex
\hangindent=5.5em \hangafter=1
{\white .}\hskip 1em $\rho_\text{isum}(\mathfrak{t})$ =
 $( \frac{1}{4},
\frac{13}{20},
\frac{17}{20} )
\oplus
( \frac{1}{4},
\frac{3}{4} )
$,

\vskip 0.7ex
\hangindent=5.5em \hangafter=1
{\white .}\hskip 1em $\rho_\text{isum}(\mathfrak{s})$ =
$\mathrm{i}$($-\sqrt{\frac{1}{5}}$,
$\sqrt{\frac{2}{5}}$,
$\sqrt{\frac{2}{5}}$;\ \ 
$\frac{-5+\sqrt{5}}{10}$,
$\frac{5+\sqrt{5}}{10}$;\ \ 
$\frac{-5+\sqrt{5}}{10}$)
 $\oplus$
$\mathrm{i}$($-\frac{1}{2}$,
$\sqrt{\frac{3}{4}}$;\ \ 
$\frac{1}{2}$)

Pass. 

 \ \color{black}

\noindent 77: (dims,levels) = $(3\oplus
2;20,
12
)$,
irreps = $3_{5}^{1}
\hskip -1.5pt \otimes \hskip -1.5pt
1_{4}^{1,0}\oplus
2_{3}^{1,0}
\hskip -1.5pt \otimes \hskip -1.5pt
1_{4}^{1,0}$,
pord$(\rho_\text{isum}(\mathfrak{t})) = 15$,

\vskip 0.7ex
\hangindent=5.5em \hangafter=1
{\white .}\hskip 1em $\rho_\text{isum}(\mathfrak{t})$ =
 $( \frac{1}{4},
\frac{1}{20},
\frac{9}{20} )
\oplus
( \frac{1}{4},
\frac{7}{12} )
$,

\vskip 0.7ex
\hangindent=5.5em \hangafter=1
{\white .}\hskip 1em $\rho_\text{isum}(\mathfrak{s})$ =
$\mathrm{i}$($\sqrt{\frac{1}{5}}$,
$\sqrt{\frac{2}{5}}$,
$\sqrt{\frac{2}{5}}$;\ \ 
$-\frac{5+\sqrt{5}}{10}$,
$\frac{5-\sqrt{5}}{10}$;\ \ 
$-\frac{5+\sqrt{5}}{10}$)
 $\oplus$
($\sqrt{\frac{1}{3}}$,
$\sqrt{\frac{2}{3}}$;
$-\sqrt{\frac{1}{3}}$)

Fail:
Tr$_I(C) = -1 <$  0 for I = [ 7/12 ]. Prop. B.4 (1) eqn. (B.18)

 \ \color{black}

\noindent 78: (dims,levels) = $(3\oplus
2;20,
12
)$,
irreps = $3_{5}^{1}
\hskip -1.5pt \otimes \hskip -1.5pt
1_{4}^{1,0}\oplus
2_{3}^{1,8}
\hskip -1.5pt \otimes \hskip -1.5pt
1_{4}^{1,0}$,
pord$(\rho_\text{isum}(\mathfrak{t})) = 15$,

\vskip 0.7ex
\hangindent=5.5em \hangafter=1
{\white .}\hskip 1em $\rho_\text{isum}(\mathfrak{t})$ =
 $( \frac{1}{4},
\frac{1}{20},
\frac{9}{20} )
\oplus
( \frac{1}{4},
\frac{11}{12} )
$,

\vskip 0.7ex
\hangindent=5.5em \hangafter=1
{\white .}\hskip 1em $\rho_\text{isum}(\mathfrak{s})$ =
$\mathrm{i}$($\sqrt{\frac{1}{5}}$,
$\sqrt{\frac{2}{5}}$,
$\sqrt{\frac{2}{5}}$;\ \ 
$-\frac{5+\sqrt{5}}{10}$,
$\frac{5-\sqrt{5}}{10}$;\ \ 
$-\frac{5+\sqrt{5}}{10}$)
 $\oplus$
($-\sqrt{\frac{1}{3}}$,
$\sqrt{\frac{2}{3}}$;
$\sqrt{\frac{1}{3}}$)

Fail:
Tr$_I(C) = -1 <$  0 for I = [ 11/12 ]. Prop. B.4 (1) eqn. (B.18)

 \ \color{black}

\noindent 79: (dims,levels) = $(3\oplus
2;20,
12
)$,
irreps = $3_{5}^{3}
\hskip -1.5pt \otimes \hskip -1.5pt
1_{4}^{1,0}\oplus
2_{3}^{1,0}
\hskip -1.5pt \otimes \hskip -1.5pt
1_{4}^{1,0}$,
pord$(\rho_\text{isum}(\mathfrak{t})) = 15$,

\vskip 0.7ex
\hangindent=5.5em \hangafter=1
{\white .}\hskip 1em $\rho_\text{isum}(\mathfrak{t})$ =
 $( \frac{1}{4},
\frac{13}{20},
\frac{17}{20} )
\oplus
( \frac{1}{4},
\frac{7}{12} )
$,

\vskip 0.7ex
\hangindent=5.5em \hangafter=1
{\white .}\hskip 1em $\rho_\text{isum}(\mathfrak{s})$ =
$\mathrm{i}$($-\sqrt{\frac{1}{5}}$,
$\sqrt{\frac{2}{5}}$,
$\sqrt{\frac{2}{5}}$;\ \ 
$\frac{-5+\sqrt{5}}{10}$,
$\frac{5+\sqrt{5}}{10}$;\ \ 
$\frac{-5+\sqrt{5}}{10}$)
 $\oplus$
($\sqrt{\frac{1}{3}}$,
$\sqrt{\frac{2}{3}}$;
$-\sqrt{\frac{1}{3}}$)

Fail:
Tr$_I(C) = -1 <$  0 for I = [ 7/12 ]. Prop. B.4 (1) eqn. (B.18)

 \ \color{black}

\noindent 80: (dims,levels) = $(3\oplus
2;20,
12
)$,
irreps = $3_{5}^{3}
\hskip -1.5pt \otimes \hskip -1.5pt
1_{4}^{1,0}\oplus
2_{3}^{1,8}
\hskip -1.5pt \otimes \hskip -1.5pt
1_{4}^{1,0}$,
pord$(\rho_\text{isum}(\mathfrak{t})) = 15$,

\vskip 0.7ex
\hangindent=5.5em \hangafter=1
{\white .}\hskip 1em $\rho_\text{isum}(\mathfrak{t})$ =
 $( \frac{1}{4},
\frac{13}{20},
\frac{17}{20} )
\oplus
( \frac{1}{4},
\frac{11}{12} )
$,

\vskip 0.7ex
\hangindent=5.5em \hangafter=1
{\white .}\hskip 1em $\rho_\text{isum}(\mathfrak{s})$ =
$\mathrm{i}$($-\sqrt{\frac{1}{5}}$,
$\sqrt{\frac{2}{5}}$,
$\sqrt{\frac{2}{5}}$;\ \ 
$\frac{-5+\sqrt{5}}{10}$,
$\frac{5+\sqrt{5}}{10}$;\ \ 
$\frac{-5+\sqrt{5}}{10}$)
 $\oplus$
($-\sqrt{\frac{1}{3}}$,
$\sqrt{\frac{2}{3}}$;
$\sqrt{\frac{1}{3}}$)

Fail:
Tr$_I(C) = -1 <$  0 for I = [ 11/12 ]. Prop. B.4 (1) eqn. (B.18)

 \ \color{black}

\noindent 81: (dims,levels) = $(3\oplus
2;20,
20
)$,
irreps = $3_{5}^{1}
\hskip -1.5pt \otimes \hskip -1.5pt
1_{4}^{1,0}\oplus
2_{5}^{1}
\hskip -1.5pt \otimes \hskip -1.5pt
1_{4}^{1,0}$,
pord$(\rho_\text{isum}(\mathfrak{t})) = 5$,

\vskip 0.7ex
\hangindent=5.5em \hangafter=1
{\white .}\hskip 1em $\rho_\text{isum}(\mathfrak{t})$ =
 $( \frac{1}{4},
\frac{1}{20},
\frac{9}{20} )
\oplus
( \frac{1}{20},
\frac{9}{20} )
$,

\vskip 0.7ex
\hangindent=5.5em \hangafter=1
{\white .}\hskip 1em $\rho_\text{isum}(\mathfrak{s})$ =
$\mathrm{i}$($\sqrt{\frac{1}{5}}$,
$\sqrt{\frac{2}{5}}$,
$\sqrt{\frac{2}{5}}$;\ \ 
$-\frac{5+\sqrt{5}}{10}$,
$\frac{5-\sqrt{5}}{10}$;\ \ 
$-\frac{5+\sqrt{5}}{10}$)
 $\oplus$
($-\frac{1}{\sqrt{5}}c^{3}_{20}
$,
$\frac{1}{\sqrt{5}}c^{1}_{20}
$;
$\frac{1}{\sqrt{5}}c^{3}_{20}
$)

Fail:
Integral: $D_{\rho}(\sigma)_{\theta} \propto $ id,
 for all $\sigma$ and all $\theta$-eigenspaces that can contain unit. Prop. B.5 (6)

 \ \color{black}

\noindent 82: (dims,levels) = $(3\oplus
2;20,
20
)$,
irreps = $3_{5}^{3}
\hskip -1.5pt \otimes \hskip -1.5pt
1_{4}^{1,0}\oplus
2_{5}^{2}
\hskip -1.5pt \otimes \hskip -1.5pt
1_{4}^{1,0}$,
pord$(\rho_\text{isum}(\mathfrak{t})) = 5$,

\vskip 0.7ex
\hangindent=5.5em \hangafter=1
{\white .}\hskip 1em $\rho_\text{isum}(\mathfrak{t})$ =
 $( \frac{1}{4},
\frac{13}{20},
\frac{17}{20} )
\oplus
( \frac{13}{20},
\frac{17}{20} )
$,

\vskip 0.7ex
\hangindent=5.5em \hangafter=1
{\white .}\hskip 1em $\rho_\text{isum}(\mathfrak{s})$ =
$\mathrm{i}$($-\sqrt{\frac{1}{5}}$,
$\sqrt{\frac{2}{5}}$,
$\sqrt{\frac{2}{5}}$;\ \ 
$\frac{-5+\sqrt{5}}{10}$,
$\frac{5+\sqrt{5}}{10}$;\ \ 
$\frac{-5+\sqrt{5}}{10}$)
 $\oplus$
($\frac{1}{\sqrt{5}}c^{1}_{20}
$,
$\frac{1}{\sqrt{5}}c^{3}_{20}
$;
$-\frac{1}{\sqrt{5}}c^{1}_{20}
$)

Fail:
Integral: $D_{\rho}(\sigma)_{\theta} \propto $ id,
 for all $\sigma$ and all $\theta$-eigenspaces that can contain unit. Prop. B.5 (6)

 \ \color{black}

 \color{blue}

\noindent 83: (dims,levels) = $(3\oplus
2;24,
3
)$,
irreps = $3_{8}^{3,0}
\hskip -1.5pt \otimes \hskip -1.5pt
1_{3}^{1,0}\oplus
2_{3}^{1,0}$,
pord$(\rho_\text{isum}(\mathfrak{t})) = 24$,

\vskip 0.7ex
\hangindent=5.5em \hangafter=1
{\white .}\hskip 1em $\rho_\text{isum}(\mathfrak{t})$ =
 $( \frac{1}{3},
\frac{5}{24},
\frac{17}{24} )
\oplus
( 0,
\frac{1}{3} )
$,

\vskip 0.7ex
\hangindent=5.5em \hangafter=1
{\white .}\hskip 1em $\rho_\text{isum}(\mathfrak{s})$ =
$\mathrm{i}$($0$,
$\sqrt{\frac{1}{2}}$,
$\sqrt{\frac{1}{2}}$;\ \ 
$\frac{1}{2}$,
$-\frac{1}{2}$;\ \ 
$\frac{1}{2}$)
 $\oplus$
$\mathrm{i}$($-\sqrt{\frac{1}{3}}$,
$\sqrt{\frac{2}{3}}$;\ \ 
$\sqrt{\frac{1}{3}}$)

Pass. 

 \ \color{black}

 \color{blue}

\noindent 84: (dims,levels) = $(3\oplus
2;24,
3
)$,
irreps = $3_{8}^{3,0}
\hskip -1.5pt \otimes \hskip -1.5pt
1_{3}^{1,0}\oplus
2_{3}^{1,4}$,
pord$(\rho_\text{isum}(\mathfrak{t})) = 24$,

\vskip 0.7ex
\hangindent=5.5em \hangafter=1
{\white .}\hskip 1em $\rho_\text{isum}(\mathfrak{t})$ =
 $( \frac{1}{3},
\frac{5}{24},
\frac{17}{24} )
\oplus
( \frac{1}{3},
\frac{2}{3} )
$,

\vskip 0.7ex
\hangindent=5.5em \hangafter=1
{\white .}\hskip 1em $\rho_\text{isum}(\mathfrak{s})$ =
$\mathrm{i}$($0$,
$\sqrt{\frac{1}{2}}$,
$\sqrt{\frac{1}{2}}$;\ \ 
$\frac{1}{2}$,
$-\frac{1}{2}$;\ \ 
$\frac{1}{2}$)
 $\oplus$
$\mathrm{i}$($-\sqrt{\frac{1}{3}}$,
$\sqrt{\frac{2}{3}}$;\ \ 
$\sqrt{\frac{1}{3}}$)

Pass. 

 \ \color{black}

 \color{blue}

\noindent 85: (dims,levels) = $(3\oplus
2;24,
3
)$,
irreps = $3_{8}^{1,0}
\hskip -1.5pt \otimes \hskip -1.5pt
1_{3}^{1,0}\oplus
2_{3}^{1,0}$,
pord$(\rho_\text{isum}(\mathfrak{t})) = 24$,

\vskip 0.7ex
\hangindent=5.5em \hangafter=1
{\white .}\hskip 1em $\rho_\text{isum}(\mathfrak{t})$ =
 $( \frac{1}{3},
\frac{11}{24},
\frac{23}{24} )
\oplus
( 0,
\frac{1}{3} )
$,

\vskip 0.7ex
\hangindent=5.5em \hangafter=1
{\white .}\hskip 1em $\rho_\text{isum}(\mathfrak{s})$ =
$\mathrm{i}$($0$,
$\sqrt{\frac{1}{2}}$,
$\sqrt{\frac{1}{2}}$;\ \ 
$-\frac{1}{2}$,
$\frac{1}{2}$;\ \ 
$-\frac{1}{2}$)
 $\oplus$
$\mathrm{i}$($-\sqrt{\frac{1}{3}}$,
$\sqrt{\frac{2}{3}}$;\ \ 
$\sqrt{\frac{1}{3}}$)

Pass. 

 \ \color{black}

 \color{blue}

\noindent 86: (dims,levels) = $(3\oplus
2;24,
3
)$,
irreps = $3_{8}^{1,0}
\hskip -1.5pt \otimes \hskip -1.5pt
1_{3}^{1,0}\oplus
2_{3}^{1,4}$,
pord$(\rho_\text{isum}(\mathfrak{t})) = 24$,

\vskip 0.7ex
\hangindent=5.5em \hangafter=1
{\white .}\hskip 1em $\rho_\text{isum}(\mathfrak{t})$ =
 $( \frac{1}{3},
\frac{11}{24},
\frac{23}{24} )
\oplus
( \frac{1}{3},
\frac{2}{3} )
$,

\vskip 0.7ex
\hangindent=5.5em \hangafter=1
{\white .}\hskip 1em $\rho_\text{isum}(\mathfrak{s})$ =
$\mathrm{i}$($0$,
$\sqrt{\frac{1}{2}}$,
$\sqrt{\frac{1}{2}}$;\ \ 
$-\frac{1}{2}$,
$\frac{1}{2}$;\ \ 
$-\frac{1}{2}$)
 $\oplus$
$\mathrm{i}$($-\sqrt{\frac{1}{3}}$,
$\sqrt{\frac{2}{3}}$;\ \ 
$\sqrt{\frac{1}{3}}$)

Pass. 

 \ \color{black}

\noindent 87: (dims,levels) = $(3\oplus
2;24,
6
)$,
irreps = $3_{8}^{3,0}
\hskip -1.5pt \otimes \hskip -1.5pt
1_{3}^{1,0}\oplus
2_{2}^{1,0}
\hskip -1.5pt \otimes \hskip -1.5pt
1_{3}^{1,0}$,
pord$(\rho_\text{isum}(\mathfrak{t})) = 8$,

\vskip 0.7ex
\hangindent=5.5em \hangafter=1
{\white .}\hskip 1em $\rho_\text{isum}(\mathfrak{t})$ =
 $( \frac{1}{3},
\frac{5}{24},
\frac{17}{24} )
\oplus
( \frac{1}{3},
\frac{5}{6} )
$,

\vskip 0.7ex
\hangindent=5.5em \hangafter=1
{\white .}\hskip 1em $\rho_\text{isum}(\mathfrak{s})$ =
$\mathrm{i}$($0$,
$\sqrt{\frac{1}{2}}$,
$\sqrt{\frac{1}{2}}$;\ \ 
$\frac{1}{2}$,
$-\frac{1}{2}$;\ \ 
$\frac{1}{2}$)
 $\oplus$
($-\frac{1}{2}$,
$-\sqrt{\frac{3}{4}}$;
$\frac{1}{2}$)

Fail:
Tr$_I(C) = -1 <$  0 for I = [ 5/6 ]. Prop. B.4 (1) eqn. (B.18)

 \ \color{black}

\noindent 88: (dims,levels) = $(3\oplus
2;24,
6
)$,
irreps = $3_{8}^{1,0}
\hskip -1.5pt \otimes \hskip -1.5pt
1_{3}^{1,0}\oplus
2_{2}^{1,0}
\hskip -1.5pt \otimes \hskip -1.5pt
1_{3}^{1,0}$,
pord$(\rho_\text{isum}(\mathfrak{t})) = 8$,

\vskip 0.7ex
\hangindent=5.5em \hangafter=1
{\white .}\hskip 1em $\rho_\text{isum}(\mathfrak{t})$ =
 $( \frac{1}{3},
\frac{11}{24},
\frac{23}{24} )
\oplus
( \frac{1}{3},
\frac{5}{6} )
$,

\vskip 0.7ex
\hangindent=5.5em \hangafter=1
{\white .}\hskip 1em $\rho_\text{isum}(\mathfrak{s})$ =
$\mathrm{i}$($0$,
$\sqrt{\frac{1}{2}}$,
$\sqrt{\frac{1}{2}}$;\ \ 
$-\frac{1}{2}$,
$\frac{1}{2}$;\ \ 
$-\frac{1}{2}$)
 $\oplus$
($-\frac{1}{2}$,
$-\sqrt{\frac{3}{4}}$;
$\frac{1}{2}$)

Fail:
Tr$_I(C) = -1 <$  0 for I = [ 5/6 ]. Prop. B.4 (1) eqn. (B.18)

 \ \color{black}

\noindent 89: (dims,levels) = $(3\oplus
2;24,
24
)$,
irreps = $3_{8}^{3,0}
\hskip -1.5pt \otimes \hskip -1.5pt
1_{3}^{1,0}\oplus
2_{8}^{1,9}
\hskip -1.5pt \otimes \hskip -1.5pt
1_{3}^{1,0}$,
pord$(\rho_\text{isum}(\mathfrak{t})) = 8$,

\vskip 0.7ex
\hangindent=5.5em \hangafter=1
{\white .}\hskip 1em $\rho_\text{isum}(\mathfrak{t})$ =
 $( \frac{1}{3},
\frac{5}{24},
\frac{17}{24} )
\oplus
( \frac{5}{24},
\frac{11}{24} )
$,

\vskip 0.7ex
\hangindent=5.5em \hangafter=1
{\white .}\hskip 1em $\rho_\text{isum}(\mathfrak{s})$ =
$\mathrm{i}$($0$,
$\sqrt{\frac{1}{2}}$,
$\sqrt{\frac{1}{2}}$;\ \ 
$\frac{1}{2}$,
$-\frac{1}{2}$;\ \ 
$\frac{1}{2}$)
 $\oplus$
$\mathrm{i}$($\sqrt{\frac{1}{2}}$,
$\sqrt{\frac{1}{2}}$;\ \ 
$-\sqrt{\frac{1}{2}}$)

Fail:
$\sigma(\rho(\mathfrak s)_\mathrm{ndeg}) \neq
 (\rho(\mathfrak t)^a \rho(\mathfrak s) \rho(\mathfrak t)^b
 \rho(\mathfrak s) \rho(\mathfrak t)^a)_\mathrm{ndeg}$,
 $\sigma = a$ = 5. Prop. B.5 (3) eqn. (B.25)

 \ \color{black}

\noindent 90: (dims,levels) = $(3\oplus
2;24,
24
)$,
irreps = $3_{8}^{3,0}
\hskip -1.5pt \otimes \hskip -1.5pt
1_{3}^{1,0}\oplus
2_{8}^{1,6}
\hskip -1.5pt \otimes \hskip -1.5pt
1_{3}^{1,0}$,
pord$(\rho_\text{isum}(\mathfrak{t})) = 8$,

\vskip 0.7ex
\hangindent=5.5em \hangafter=1
{\white .}\hskip 1em $\rho_\text{isum}(\mathfrak{t})$ =
 $( \frac{1}{3},
\frac{5}{24},
\frac{17}{24} )
\oplus
( \frac{5}{24},
\frac{23}{24} )
$,

\vskip 0.7ex
\hangindent=5.5em \hangafter=1
{\white .}\hskip 1em $\rho_\text{isum}(\mathfrak{s})$ =
$\mathrm{i}$($0$,
$\sqrt{\frac{1}{2}}$,
$\sqrt{\frac{1}{2}}$;\ \ 
$\frac{1}{2}$,
$-\frac{1}{2}$;\ \ 
$\frac{1}{2}$)
 $\oplus$
($-\sqrt{\frac{1}{2}}$,
$\sqrt{\frac{1}{2}}$;
$\sqrt{\frac{1}{2}}$)

Fail:
Tr$_I(C) = -1 <$  0 for I = [ 23/24 ]. Prop. B.4 (1) eqn. (B.18)

 \ \color{black}

\noindent 91: (dims,levels) = $(3\oplus
2;24,
24
)$,
irreps = $3_{8}^{3,0}
\hskip -1.5pt \otimes \hskip -1.5pt
1_{3}^{1,0}\oplus
2_{8}^{1,0}
\hskip -1.5pt \otimes \hskip -1.5pt
1_{3}^{1,0}$,
pord$(\rho_\text{isum}(\mathfrak{t})) = 8$,

\vskip 0.7ex
\hangindent=5.5em \hangafter=1
{\white .}\hskip 1em $\rho_\text{isum}(\mathfrak{t})$ =
 $( \frac{1}{3},
\frac{5}{24},
\frac{17}{24} )
\oplus
( \frac{11}{24},
\frac{17}{24} )
$,

\vskip 0.7ex
\hangindent=5.5em \hangafter=1
{\white .}\hskip 1em $\rho_\text{isum}(\mathfrak{s})$ =
$\mathrm{i}$($0$,
$\sqrt{\frac{1}{2}}$,
$\sqrt{\frac{1}{2}}$;\ \ 
$\frac{1}{2}$,
$-\frac{1}{2}$;\ \ 
$\frac{1}{2}$)
 $\oplus$
($-\sqrt{\frac{1}{2}}$,
$\sqrt{\frac{1}{2}}$;
$\sqrt{\frac{1}{2}}$)

Fail:
Tr$_I(C) = -1 <$  0 for I = [ 11/24 ]. Prop. B.4 (1) eqn. (B.18)

 \ \color{black}

\noindent 92: (dims,levels) = $(3\oplus
2;24,
24
)$,
irreps = $3_{8}^{3,0}
\hskip -1.5pt \otimes \hskip -1.5pt
1_{3}^{1,0}\oplus
2_{8}^{1,3}
\hskip -1.5pt \otimes \hskip -1.5pt
1_{3}^{1,0}$,
pord$(\rho_\text{isum}(\mathfrak{t})) = 8$,

\vskip 0.7ex
\hangindent=5.5em \hangafter=1
{\white .}\hskip 1em $\rho_\text{isum}(\mathfrak{t})$ =
 $( \frac{1}{3},
\frac{5}{24},
\frac{17}{24} )
\oplus
( \frac{17}{24},
\frac{23}{24} )
$,

\vskip 0.7ex
\hangindent=5.5em \hangafter=1
{\white .}\hskip 1em $\rho_\text{isum}(\mathfrak{s})$ =
$\mathrm{i}$($0$,
$\sqrt{\frac{1}{2}}$,
$\sqrt{\frac{1}{2}}$;\ \ 
$\frac{1}{2}$,
$-\frac{1}{2}$;\ \ 
$\frac{1}{2}$)
 $\oplus$
$\mathrm{i}$($-\sqrt{\frac{1}{2}}$,
$\sqrt{\frac{1}{2}}$;\ \ 
$\sqrt{\frac{1}{2}}$)

Fail:
$\sigma(\rho(\mathfrak s)_\mathrm{ndeg}) \neq
 (\rho(\mathfrak t)^a \rho(\mathfrak s) \rho(\mathfrak t)^b
 \rho(\mathfrak s) \rho(\mathfrak t)^a)_\mathrm{ndeg}$,
 $\sigma = a$ = 5. Prop. B.5 (3) eqn. (B.25)

 \ \color{black}

\noindent 93: (dims,levels) = $(3\oplus
2;24,
24
)$,
irreps = $3_{8}^{1,0}
\hskip -1.5pt \otimes \hskip -1.5pt
1_{3}^{1,0}\oplus
2_{8}^{1,9}
\hskip -1.5pt \otimes \hskip -1.5pt
1_{3}^{1,0}$,
pord$(\rho_\text{isum}(\mathfrak{t})) = 8$,

\vskip 0.7ex
\hangindent=5.5em \hangafter=1
{\white .}\hskip 1em $\rho_\text{isum}(\mathfrak{t})$ =
 $( \frac{1}{3},
\frac{11}{24},
\frac{23}{24} )
\oplus
( \frac{5}{24},
\frac{11}{24} )
$,

\vskip 0.7ex
\hangindent=5.5em \hangafter=1
{\white .}\hskip 1em $\rho_\text{isum}(\mathfrak{s})$ =
$\mathrm{i}$($0$,
$\sqrt{\frac{1}{2}}$,
$\sqrt{\frac{1}{2}}$;\ \ 
$-\frac{1}{2}$,
$\frac{1}{2}$;\ \ 
$-\frac{1}{2}$)
 $\oplus$
$\mathrm{i}$($\sqrt{\frac{1}{2}}$,
$\sqrt{\frac{1}{2}}$;\ \ 
$-\sqrt{\frac{1}{2}}$)

Fail:
$\sigma(\rho(\mathfrak s)_\mathrm{ndeg}) \neq
 (\rho(\mathfrak t)^a \rho(\mathfrak s) \rho(\mathfrak t)^b
 \rho(\mathfrak s) \rho(\mathfrak t)^a)_\mathrm{ndeg}$,
 $\sigma = a$ = 5. Prop. B.5 (3) eqn. (B.25)

 \ \color{black}

\noindent 94: (dims,levels) = $(3\oplus
2;24,
24
)$,
irreps = $3_{8}^{1,0}
\hskip -1.5pt \otimes \hskip -1.5pt
1_{3}^{1,0}\oplus
2_{8}^{1,6}
\hskip -1.5pt \otimes \hskip -1.5pt
1_{3}^{1,0}$,
pord$(\rho_\text{isum}(\mathfrak{t})) = 8$,

\vskip 0.7ex
\hangindent=5.5em \hangafter=1
{\white .}\hskip 1em $\rho_\text{isum}(\mathfrak{t})$ =
 $( \frac{1}{3},
\frac{11}{24},
\frac{23}{24} )
\oplus
( \frac{5}{24},
\frac{23}{24} )
$,

\vskip 0.7ex
\hangindent=5.5em \hangafter=1
{\white .}\hskip 1em $\rho_\text{isum}(\mathfrak{s})$ =
$\mathrm{i}$($0$,
$\sqrt{\frac{1}{2}}$,
$\sqrt{\frac{1}{2}}$;\ \ 
$-\frac{1}{2}$,
$\frac{1}{2}$;\ \ 
$-\frac{1}{2}$)
 $\oplus$
($-\sqrt{\frac{1}{2}}$,
$\sqrt{\frac{1}{2}}$;
$\sqrt{\frac{1}{2}}$)

Fail:
Tr$_I(C) = -1 <$  0 for I = [ 5/24 ]. Prop. B.4 (1) eqn. (B.18)

 \ \color{black}

\noindent 95: (dims,levels) = $(3\oplus
2;24,
24
)$,
irreps = $3_{8}^{1,0}
\hskip -1.5pt \otimes \hskip -1.5pt
1_{3}^{1,0}\oplus
2_{8}^{1,0}
\hskip -1.5pt \otimes \hskip -1.5pt
1_{3}^{1,0}$,
pord$(\rho_\text{isum}(\mathfrak{t})) = 8$,

\vskip 0.7ex
\hangindent=5.5em \hangafter=1
{\white .}\hskip 1em $\rho_\text{isum}(\mathfrak{t})$ =
 $( \frac{1}{3},
\frac{11}{24},
\frac{23}{24} )
\oplus
( \frac{11}{24},
\frac{17}{24} )
$,

\vskip 0.7ex
\hangindent=5.5em \hangafter=1
{\white .}\hskip 1em $\rho_\text{isum}(\mathfrak{s})$ =
$\mathrm{i}$($0$,
$\sqrt{\frac{1}{2}}$,
$\sqrt{\frac{1}{2}}$;\ \ 
$-\frac{1}{2}$,
$\frac{1}{2}$;\ \ 
$-\frac{1}{2}$)
 $\oplus$
($-\sqrt{\frac{1}{2}}$,
$\sqrt{\frac{1}{2}}$;
$\sqrt{\frac{1}{2}}$)

Fail:
Tr$_I(C) = -1 <$  0 for I = [ 17/24 ]. Prop. B.4 (1) eqn. (B.18)

 \ \color{black}

\noindent 96: (dims,levels) = $(3\oplus
2;24,
24
)$,
irreps = $3_{8}^{1,0}
\hskip -1.5pt \otimes \hskip -1.5pt
1_{3}^{1,0}\oplus
2_{8}^{1,3}
\hskip -1.5pt \otimes \hskip -1.5pt
1_{3}^{1,0}$,
pord$(\rho_\text{isum}(\mathfrak{t})) = 8$,

\vskip 0.7ex
\hangindent=5.5em \hangafter=1
{\white .}\hskip 1em $\rho_\text{isum}(\mathfrak{t})$ =
 $( \frac{1}{3},
\frac{11}{24},
\frac{23}{24} )
\oplus
( \frac{17}{24},
\frac{23}{24} )
$,

\vskip 0.7ex
\hangindent=5.5em \hangafter=1
{\white .}\hskip 1em $\rho_\text{isum}(\mathfrak{s})$ =
$\mathrm{i}$($0$,
$\sqrt{\frac{1}{2}}$,
$\sqrt{\frac{1}{2}}$;\ \ 
$-\frac{1}{2}$,
$\frac{1}{2}$;\ \ 
$-\frac{1}{2}$)
 $\oplus$
$\mathrm{i}$($-\sqrt{\frac{1}{2}}$,
$\sqrt{\frac{1}{2}}$;\ \ 
$\sqrt{\frac{1}{2}}$)

Fail:
$\sigma(\rho(\mathfrak s)_\mathrm{ndeg}) \neq
 (\rho(\mathfrak t)^a \rho(\mathfrak s) \rho(\mathfrak t)^b
 \rho(\mathfrak s) \rho(\mathfrak t)^a)_\mathrm{ndeg}$,
 $\sigma = a$ = 5. Prop. B.5 (3) eqn. (B.25)

 \ \color{black}

\noindent 97: (dims,levels) = $(3\oplus
2;30,
6
)$,
irreps = $3_{5}^{1}
\hskip -1.5pt \otimes \hskip -1.5pt
1_{3}^{1,0}
\hskip -1.5pt \otimes \hskip -1.5pt
1_{2}^{1,0}\oplus
2_{3}^{1,4}
\hskip -1.5pt \otimes \hskip -1.5pt
1_{2}^{1,0}$,
pord$(\rho_\text{isum}(\mathfrak{t})) = 15$,

\vskip 0.7ex
\hangindent=5.5em \hangafter=1
{\white .}\hskip 1em $\rho_\text{isum}(\mathfrak{t})$ =
 $( \frac{5}{6},
\frac{1}{30},
\frac{19}{30} )
\oplus
( \frac{1}{6},
\frac{5}{6} )
$,

\vskip 0.7ex
\hangindent=5.5em \hangafter=1
{\white .}\hskip 1em $\rho_\text{isum}(\mathfrak{s})$ =
($-\sqrt{\frac{1}{5}}$,
$-\sqrt{\frac{2}{5}}$,
$-\sqrt{\frac{2}{5}}$;
$\frac{5+\sqrt{5}}{10}$,
$\frac{-5+\sqrt{5}}{10}$;
$\frac{5+\sqrt{5}}{10}$)
 $\oplus$
$\mathrm{i}$($-\sqrt{\frac{1}{3}}$,
$\sqrt{\frac{2}{3}}$;\ \ 
$\sqrt{\frac{1}{3}}$)

Fail:
Tr$_I(C) = -1 <$  0 for I = [ 1/6 ]. Prop. B.4 (1) eqn. (B.18)

 \ \color{black}

 \color{blue}

\noindent 98: (dims,levels) = $(3\oplus
2;30,
6
)$,
irreps = $3_{5}^{1}
\hskip -1.5pt \otimes \hskip -1.5pt
1_{3}^{1,0}
\hskip -1.5pt \otimes \hskip -1.5pt
1_{2}^{1,0}\oplus
2_{2}^{1,0}
\hskip -1.5pt \otimes \hskip -1.5pt
1_{3}^{1,0}$,
pord$(\rho_\text{isum}(\mathfrak{t})) = 10$,

\vskip 0.7ex
\hangindent=5.5em \hangafter=1
{\white .}\hskip 1em $\rho_\text{isum}(\mathfrak{t})$ =
 $( \frac{5}{6},
\frac{1}{30},
\frac{19}{30} )
\oplus
( \frac{1}{3},
\frac{5}{6} )
$,

\vskip 0.7ex
\hangindent=5.5em \hangafter=1
{\white .}\hskip 1em $\rho_\text{isum}(\mathfrak{s})$ =
($-\sqrt{\frac{1}{5}}$,
$-\sqrt{\frac{2}{5}}$,
$-\sqrt{\frac{2}{5}}$;
$\frac{5+\sqrt{5}}{10}$,
$\frac{-5+\sqrt{5}}{10}$;
$\frac{5+\sqrt{5}}{10}$)
 $\oplus$
($-\frac{1}{2}$,
$-\sqrt{\frac{3}{4}}$;
$\frac{1}{2}$)

Pass. 

 \ \color{black}

\noindent 99: (dims,levels) = $(3\oplus
2;30,
6
)$,
irreps = $3_{5}^{1}
\hskip -1.5pt \otimes \hskip -1.5pt
1_{3}^{1,0}
\hskip -1.5pt \otimes \hskip -1.5pt
1_{2}^{1,0}\oplus
2_{3}^{1,0}
\hskip -1.5pt \otimes \hskip -1.5pt
1_{2}^{1,0}$,
pord$(\rho_\text{isum}(\mathfrak{t})) = 15$,

\vskip 0.7ex
\hangindent=5.5em \hangafter=1
{\white .}\hskip 1em $\rho_\text{isum}(\mathfrak{t})$ =
 $( \frac{5}{6},
\frac{1}{30},
\frac{19}{30} )
\oplus
( \frac{1}{2},
\frac{5}{6} )
$,

\vskip 0.7ex
\hangindent=5.5em \hangafter=1
{\white .}\hskip 1em $\rho_\text{isum}(\mathfrak{s})$ =
($-\sqrt{\frac{1}{5}}$,
$-\sqrt{\frac{2}{5}}$,
$-\sqrt{\frac{2}{5}}$;
$\frac{5+\sqrt{5}}{10}$,
$\frac{-5+\sqrt{5}}{10}$;
$\frac{5+\sqrt{5}}{10}$)
 $\oplus$
$\mathrm{i}$($\sqrt{\frac{1}{3}}$,
$\sqrt{\frac{2}{3}}$;\ \ 
$-\sqrt{\frac{1}{3}}$)

Fail:
Tr$_I(C) = -1 <$  0 for I = [ 1/2 ]. Prop. B.4 (1) eqn. (B.18)

 \ \color{black}

\noindent 100: (dims,levels) = $(3\oplus
2;30,
6
)$,
irreps = $3_{5}^{3}
\hskip -1.5pt \otimes \hskip -1.5pt
1_{3}^{1,0}
\hskip -1.5pt \otimes \hskip -1.5pt
1_{2}^{1,0}\oplus
2_{3}^{1,4}
\hskip -1.5pt \otimes \hskip -1.5pt
1_{2}^{1,0}$,
pord$(\rho_\text{isum}(\mathfrak{t})) = 15$,

\vskip 0.7ex
\hangindent=5.5em \hangafter=1
{\white .}\hskip 1em $\rho_\text{isum}(\mathfrak{t})$ =
 $( \frac{5}{6},
\frac{7}{30},
\frac{13}{30} )
\oplus
( \frac{1}{6},
\frac{5}{6} )
$,

\vskip 0.7ex
\hangindent=5.5em \hangafter=1
{\white .}\hskip 1em $\rho_\text{isum}(\mathfrak{s})$ =
($\sqrt{\frac{1}{5}}$,
$-\sqrt{\frac{2}{5}}$,
$-\sqrt{\frac{2}{5}}$;
$\frac{5-\sqrt{5}}{10}$,
$-\frac{5+\sqrt{5}}{10}$;
$\frac{5-\sqrt{5}}{10}$)
 $\oplus$
$\mathrm{i}$($-\sqrt{\frac{1}{3}}$,
$\sqrt{\frac{2}{3}}$;\ \ 
$\sqrt{\frac{1}{3}}$)

Fail:
Tr$_I(C) = -1 <$  0 for I = [ 1/6 ]. Prop. B.4 (1) eqn. (B.18)

 \ \color{black}

 \color{blue}

\noindent 101: (dims,levels) = $(3\oplus
2;30,
6
)$,
irreps = $3_{5}^{3}
\hskip -1.5pt \otimes \hskip -1.5pt
1_{3}^{1,0}
\hskip -1.5pt \otimes \hskip -1.5pt
1_{2}^{1,0}\oplus
2_{2}^{1,0}
\hskip -1.5pt \otimes \hskip -1.5pt
1_{3}^{1,0}$,
pord$(\rho_\text{isum}(\mathfrak{t})) = 10$,

\vskip 0.7ex
\hangindent=5.5em \hangafter=1
{\white .}\hskip 1em $\rho_\text{isum}(\mathfrak{t})$ =
 $( \frac{5}{6},
\frac{7}{30},
\frac{13}{30} )
\oplus
( \frac{1}{3},
\frac{5}{6} )
$,

\vskip 0.7ex
\hangindent=5.5em \hangafter=1
{\white .}\hskip 1em $\rho_\text{isum}(\mathfrak{s})$ =
($\sqrt{\frac{1}{5}}$,
$-\sqrt{\frac{2}{5}}$,
$-\sqrt{\frac{2}{5}}$;
$\frac{5-\sqrt{5}}{10}$,
$-\frac{5+\sqrt{5}}{10}$;
$\frac{5-\sqrt{5}}{10}$)
 $\oplus$
($-\frac{1}{2}$,
$-\sqrt{\frac{3}{4}}$;
$\frac{1}{2}$)

Pass. 

 \ \color{black}

\noindent 102: (dims,levels) = $(3\oplus
2;30,
6
)$,
irreps = $3_{5}^{3}
\hskip -1.5pt \otimes \hskip -1.5pt
1_{3}^{1,0}
\hskip -1.5pt \otimes \hskip -1.5pt
1_{2}^{1,0}\oplus
2_{3}^{1,0}
\hskip -1.5pt \otimes \hskip -1.5pt
1_{2}^{1,0}$,
pord$(\rho_\text{isum}(\mathfrak{t})) = 15$,

\vskip 0.7ex
\hangindent=5.5em \hangafter=1
{\white .}\hskip 1em $\rho_\text{isum}(\mathfrak{t})$ =
 $( \frac{5}{6},
\frac{7}{30},
\frac{13}{30} )
\oplus
( \frac{1}{2},
\frac{5}{6} )
$,

\vskip 0.7ex
\hangindent=5.5em \hangafter=1
{\white .}\hskip 1em $\rho_\text{isum}(\mathfrak{s})$ =
($\sqrt{\frac{1}{5}}$,
$-\sqrt{\frac{2}{5}}$,
$-\sqrt{\frac{2}{5}}$;
$\frac{5-\sqrt{5}}{10}$,
$-\frac{5+\sqrt{5}}{10}$;
$\frac{5-\sqrt{5}}{10}$)
 $\oplus$
$\mathrm{i}$($\sqrt{\frac{1}{3}}$,
$\sqrt{\frac{2}{3}}$;\ \ 
$-\sqrt{\frac{1}{3}}$)

Fail:
Tr$_I(C) = -1 <$  0 for I = [ 1/2 ]. Prop. B.4 (1) eqn. (B.18)

 \ \color{black}

\noindent 103: (dims,levels) = $(3\oplus
2;30,
30
)$,
irreps = $3_{5}^{1}
\hskip -1.5pt \otimes \hskip -1.5pt
1_{3}^{1,0}
\hskip -1.5pt \otimes \hskip -1.5pt
1_{2}^{1,0}\oplus
2_{5}^{1}
\hskip -1.5pt \otimes \hskip -1.5pt
1_{3}^{1,0}
\hskip -1.5pt \otimes \hskip -1.5pt
1_{2}^{1,0}$,
pord$(\rho_\text{isum}(\mathfrak{t})) = 5$,

\vskip 0.7ex
\hangindent=5.5em \hangafter=1
{\white .}\hskip 1em $\rho_\text{isum}(\mathfrak{t})$ =
 $( \frac{5}{6},
\frac{1}{30},
\frac{19}{30} )
\oplus
( \frac{1}{30},
\frac{19}{30} )
$,

\vskip 0.7ex
\hangindent=5.5em \hangafter=1
{\white .}\hskip 1em $\rho_\text{isum}(\mathfrak{s})$ =
($-\sqrt{\frac{1}{5}}$,
$-\sqrt{\frac{2}{5}}$,
$-\sqrt{\frac{2}{5}}$;
$\frac{5+\sqrt{5}}{10}$,
$\frac{-5+\sqrt{5}}{10}$;
$\frac{5+\sqrt{5}}{10}$)
 $\oplus$
$\mathrm{i}$($\frac{1}{\sqrt{5}}c^{3}_{20}
$,
$\frac{1}{\sqrt{5}}c^{1}_{20}
$;\ \ 
$-\frac{1}{\sqrt{5}}c^{3}_{20}
$)

Fail:
Integral: $D_{\rho}(\sigma)_{\theta} \propto $ id,
 for all $\sigma$ and all $\theta$-eigenspaces that can contain unit. Prop. B.5 (6)

 \ \color{black}

\noindent 104: (dims,levels) = $(3\oplus
2;30,
30
)$,
irreps = $3_{5}^{3}
\hskip -1.5pt \otimes \hskip -1.5pt
1_{3}^{1,0}
\hskip -1.5pt \otimes \hskip -1.5pt
1_{2}^{1,0}\oplus
2_{5}^{2}
\hskip -1.5pt \otimes \hskip -1.5pt
1_{3}^{1,0}
\hskip -1.5pt \otimes \hskip -1.5pt
1_{2}^{1,0}$,
pord$(\rho_\text{isum}(\mathfrak{t})) = 5$,

\vskip 0.7ex
\hangindent=5.5em \hangafter=1
{\white .}\hskip 1em $\rho_\text{isum}(\mathfrak{t})$ =
 $( \frac{5}{6},
\frac{7}{30},
\frac{13}{30} )
\oplus
( \frac{7}{30},
\frac{13}{30} )
$,

\vskip 0.7ex
\hangindent=5.5em \hangafter=1
{\white .}\hskip 1em $\rho_\text{isum}(\mathfrak{s})$ =
($\sqrt{\frac{1}{5}}$,
$-\sqrt{\frac{2}{5}}$,
$-\sqrt{\frac{2}{5}}$;
$\frac{5-\sqrt{5}}{10}$,
$-\frac{5+\sqrt{5}}{10}$;
$\frac{5-\sqrt{5}}{10}$)
 $\oplus$
$\mathrm{i}$($\frac{1}{\sqrt{5}}c^{1}_{20}
$,
$\frac{1}{\sqrt{5}}c^{3}_{20}
$;\ \ 
$-\frac{1}{\sqrt{5}}c^{1}_{20}
$)

Fail:
Integral: $D_{\rho}(\sigma)_{\theta} \propto $ id,
 for all $\sigma$ and all $\theta$-eigenspaces that can contain unit. Prop. B.5 (6)

 \ \color{black}

\noindent 105: (dims,levels) = $(3\oplus
2;48,
24
)$,
irreps = $3_{16}^{7,0}
\hskip -1.5pt \otimes \hskip -1.5pt
1_{3}^{1,0}\oplus
2_{8}^{1,9}
\hskip -1.5pt \otimes \hskip -1.5pt
1_{3}^{1,0}$,
pord$(\rho_\text{isum}(\mathfrak{t})) = 16$,

\vskip 0.7ex
\hangindent=5.5em \hangafter=1
{\white .}\hskip 1em $\rho_\text{isum}(\mathfrak{t})$ =
 $( \frac{5}{24},
\frac{13}{48},
\frac{37}{48} )
\oplus
( \frac{5}{24},
\frac{11}{24} )
$,

\vskip 0.7ex
\hangindent=5.5em \hangafter=1
{\white .}\hskip 1em $\rho_\text{isum}(\mathfrak{s})$ =
$\mathrm{i}$($0$,
$\sqrt{\frac{1}{2}}$,
$\sqrt{\frac{1}{2}}$;\ \ 
$\frac{1}{2}$,
$-\frac{1}{2}$;\ \ 
$\frac{1}{2}$)
 $\oplus$
$\mathrm{i}$($\sqrt{\frac{1}{2}}$,
$\sqrt{\frac{1}{2}}$;\ \ 
$-\sqrt{\frac{1}{2}}$)

Fail:
Integral: $D_{\rho}(\sigma)_{\theta} \propto $ id,
 for all $\sigma$ and all $\theta$-eigenspaces that can contain unit. Prop. B.5 (6)

 \ \color{black}

\noindent 106: (dims,levels) = $(3\oplus
2;48,
24
)$,
irreps = $3_{16}^{7,0}
\hskip -1.5pt \otimes \hskip -1.5pt
1_{3}^{1,0}\oplus
2_{8}^{1,6}
\hskip -1.5pt \otimes \hskip -1.5pt
1_{3}^{1,0}$,
pord$(\rho_\text{isum}(\mathfrak{t})) = 16$,

\vskip 0.7ex
\hangindent=5.5em \hangafter=1
{\white .}\hskip 1em $\rho_\text{isum}(\mathfrak{t})$ =
 $( \frac{5}{24},
\frac{13}{48},
\frac{37}{48} )
\oplus
( \frac{5}{24},
\frac{23}{24} )
$,

\vskip 0.7ex
\hangindent=5.5em \hangafter=1
{\white .}\hskip 1em $\rho_\text{isum}(\mathfrak{s})$ =
$\mathrm{i}$($0$,
$\sqrt{\frac{1}{2}}$,
$\sqrt{\frac{1}{2}}$;\ \ 
$\frac{1}{2}$,
$-\frac{1}{2}$;\ \ 
$\frac{1}{2}$)
 $\oplus$
($-\sqrt{\frac{1}{2}}$,
$\sqrt{\frac{1}{2}}$;
$\sqrt{\frac{1}{2}}$)

Fail:
Tr$_I(C) = -1 <$  0 for I = [ 23/24 ]. Prop. B.4 (1) eqn. (B.18)

 \ \color{black}

\noindent 107: (dims,levels) = $(3\oplus
2;48,
24
)$,
irreps = $3_{16}^{1,0}
\hskip -1.5pt \otimes \hskip -1.5pt
1_{3}^{1,0}\oplus
2_{8}^{1,9}
\hskip -1.5pt \otimes \hskip -1.5pt
1_{3}^{1,0}$,
pord$(\rho_\text{isum}(\mathfrak{t})) = 16$,

\vskip 0.7ex
\hangindent=5.5em \hangafter=1
{\white .}\hskip 1em $\rho_\text{isum}(\mathfrak{t})$ =
 $( \frac{11}{24},
\frac{19}{48},
\frac{43}{48} )
\oplus
( \frac{5}{24},
\frac{11}{24} )
$,

\vskip 0.7ex
\hangindent=5.5em \hangafter=1
{\white .}\hskip 1em $\rho_\text{isum}(\mathfrak{s})$ =
$\mathrm{i}$($0$,
$\sqrt{\frac{1}{2}}$,
$\sqrt{\frac{1}{2}}$;\ \ 
$-\frac{1}{2}$,
$\frac{1}{2}$;\ \ 
$-\frac{1}{2}$)
 $\oplus$
$\mathrm{i}$($\sqrt{\frac{1}{2}}$,
$\sqrt{\frac{1}{2}}$;\ \ 
$-\sqrt{\frac{1}{2}}$)

Fail:
Integral: $D_{\rho}(\sigma)_{\theta} \propto $ id,
 for all $\sigma$ and all $\theta$-eigenspaces that can contain unit. Prop. B.5 (6)

 \ \color{black}

\noindent 108: (dims,levels) = $(3\oplus
2;48,
24
)$,
irreps = $3_{16}^{1,0}
\hskip -1.5pt \otimes \hskip -1.5pt
1_{3}^{1,0}\oplus
2_{8}^{1,0}
\hskip -1.5pt \otimes \hskip -1.5pt
1_{3}^{1,0}$,
pord$(\rho_\text{isum}(\mathfrak{t})) = 16$,

\vskip 0.7ex
\hangindent=5.5em \hangafter=1
{\white .}\hskip 1em $\rho_\text{isum}(\mathfrak{t})$ =
 $( \frac{11}{24},
\frac{19}{48},
\frac{43}{48} )
\oplus
( \frac{11}{24},
\frac{17}{24} )
$,

\vskip 0.7ex
\hangindent=5.5em \hangafter=1
{\white .}\hskip 1em $\rho_\text{isum}(\mathfrak{s})$ =
$\mathrm{i}$($0$,
$\sqrt{\frac{1}{2}}$,
$\sqrt{\frac{1}{2}}$;\ \ 
$-\frac{1}{2}$,
$\frac{1}{2}$;\ \ 
$-\frac{1}{2}$)
 $\oplus$
($-\sqrt{\frac{1}{2}}$,
$\sqrt{\frac{1}{2}}$;
$\sqrt{\frac{1}{2}}$)

Fail:
Tr$_I(C) = -1 <$  0 for I = [ 17/24 ]. Prop. B.4 (1) eqn. (B.18)

 \ \color{black}

\noindent 109: (dims,levels) = $(3\oplus
2;48,
24
)$,
irreps = $3_{16}^{3,0}
\hskip -1.5pt \otimes \hskip -1.5pt
1_{3}^{1,0}\oplus
2_{8}^{1,0}
\hskip -1.5pt \otimes \hskip -1.5pt
1_{3}^{1,0}$,
pord$(\rho_\text{isum}(\mathfrak{t})) = 16$,

\vskip 0.7ex
\hangindent=5.5em \hangafter=1
{\white .}\hskip 1em $\rho_\text{isum}(\mathfrak{t})$ =
 $( \frac{17}{24},
\frac{1}{48},
\frac{25}{48} )
\oplus
( \frac{11}{24},
\frac{17}{24} )
$,

\vskip 0.7ex
\hangindent=5.5em \hangafter=1
{\white .}\hskip 1em $\rho_\text{isum}(\mathfrak{s})$ =
$\mathrm{i}$($0$,
$\sqrt{\frac{1}{2}}$,
$\sqrt{\frac{1}{2}}$;\ \ 
$\frac{1}{2}$,
$-\frac{1}{2}$;\ \ 
$\frac{1}{2}$)
 $\oplus$
($-\sqrt{\frac{1}{2}}$,
$\sqrt{\frac{1}{2}}$;
$\sqrt{\frac{1}{2}}$)

Fail:
Tr$_I(C) = -1 <$  0 for I = [ 11/24 ]. Prop. B.4 (1) eqn. (B.18)

 \ \color{black}

\noindent 110: (dims,levels) = $(3\oplus
2;48,
24
)$,
irreps = $3_{16}^{3,0}
\hskip -1.5pt \otimes \hskip -1.5pt
1_{3}^{1,0}\oplus
2_{8}^{1,3}
\hskip -1.5pt \otimes \hskip -1.5pt
1_{3}^{1,0}$,
pord$(\rho_\text{isum}(\mathfrak{t})) = 16$,

\vskip 0.7ex
\hangindent=5.5em \hangafter=1
{\white .}\hskip 1em $\rho_\text{isum}(\mathfrak{t})$ =
 $( \frac{17}{24},
\frac{1}{48},
\frac{25}{48} )
\oplus
( \frac{17}{24},
\frac{23}{24} )
$,

\vskip 0.7ex
\hangindent=5.5em \hangafter=1
{\white .}\hskip 1em $\rho_\text{isum}(\mathfrak{s})$ =
$\mathrm{i}$($0$,
$\sqrt{\frac{1}{2}}$,
$\sqrt{\frac{1}{2}}$;\ \ 
$\frac{1}{2}$,
$-\frac{1}{2}$;\ \ 
$\frac{1}{2}$)
 $\oplus$
$\mathrm{i}$($-\sqrt{\frac{1}{2}}$,
$\sqrt{\frac{1}{2}}$;\ \ 
$\sqrt{\frac{1}{2}}$)

Fail:
Integral: $D_{\rho}(\sigma)_{\theta} \propto $ id,
 for all $\sigma$ and all $\theta$-eigenspaces that can contain unit. Prop. B.5 (6)

 \ \color{black}

\noindent 111: (dims,levels) = $(3\oplus
2;48,
24
)$,
irreps = $3_{16}^{5,0}
\hskip -1.5pt \otimes \hskip -1.5pt
1_{3}^{1,0}\oplus
2_{8}^{1,6}
\hskip -1.5pt \otimes \hskip -1.5pt
1_{3}^{1,0}$,
pord$(\rho_\text{isum}(\mathfrak{t})) = 16$,

\vskip 0.7ex
\hangindent=5.5em \hangafter=1
{\white .}\hskip 1em $\rho_\text{isum}(\mathfrak{t})$ =
 $( \frac{23}{24},
\frac{7}{48},
\frac{31}{48} )
\oplus
( \frac{5}{24},
\frac{23}{24} )
$,

\vskip 0.7ex
\hangindent=5.5em \hangafter=1
{\white .}\hskip 1em $\rho_\text{isum}(\mathfrak{s})$ =
$\mathrm{i}$($0$,
$\sqrt{\frac{1}{2}}$,
$\sqrt{\frac{1}{2}}$;\ \ 
$-\frac{1}{2}$,
$\frac{1}{2}$;\ \ 
$-\frac{1}{2}$)
 $\oplus$
($-\sqrt{\frac{1}{2}}$,
$\sqrt{\frac{1}{2}}$;
$\sqrt{\frac{1}{2}}$)

Fail:
Tr$_I(C) = -1 <$  0 for I = [ 5/24 ]. Prop. B.4 (1) eqn. (B.18)

 \ \color{black}

\noindent 112: (dims,levels) = $(3\oplus
2;48,
24
)$,
irreps = $3_{16}^{5,0}
\hskip -1.5pt \otimes \hskip -1.5pt
1_{3}^{1,0}\oplus
2_{8}^{1,3}
\hskip -1.5pt \otimes \hskip -1.5pt
1_{3}^{1,0}$,
pord$(\rho_\text{isum}(\mathfrak{t})) = 16$,

\vskip 0.7ex
\hangindent=5.5em \hangafter=1
{\white .}\hskip 1em $\rho_\text{isum}(\mathfrak{t})$ =
 $( \frac{23}{24},
\frac{7}{48},
\frac{31}{48} )
\oplus
( \frac{17}{24},
\frac{23}{24} )
$,

\vskip 0.7ex
\hangindent=5.5em \hangafter=1
{\white .}\hskip 1em $\rho_\text{isum}(\mathfrak{s})$ =
$\mathrm{i}$($0$,
$\sqrt{\frac{1}{2}}$,
$\sqrt{\frac{1}{2}}$;\ \ 
$-\frac{1}{2}$,
$\frac{1}{2}$;\ \ 
$-\frac{1}{2}$)
 $\oplus$
$\mathrm{i}$($-\sqrt{\frac{1}{2}}$,
$\sqrt{\frac{1}{2}}$;\ \ 
$\sqrt{\frac{1}{2}}$)

Fail:
Integral: $D_{\rho}(\sigma)_{\theta} \propto $ id,
 for all $\sigma$ and all $\theta$-eigenspaces that can contain unit. Prop. B.5 (6)

 \ \color{black}

 \color{blue}

\noindent 113: (dims,levels) = $(3\oplus
2;60,
12
)$,
irreps = $3_{5}^{3}
\hskip -1.5pt \otimes \hskip -1.5pt
1_{4}^{1,0}
\hskip -1.5pt \otimes \hskip -1.5pt
1_{3}^{1,0}\oplus
2_{4}^{1,0}
\hskip -1.5pt \otimes \hskip -1.5pt
1_{3}^{1,0}$,
pord$(\rho_\text{isum}(\mathfrak{t})) = 10$,

\vskip 0.7ex
\hangindent=5.5em \hangafter=1
{\white .}\hskip 1em $\rho_\text{isum}(\mathfrak{t})$ =
 $( \frac{7}{12},
\frac{11}{60},
\frac{59}{60} )
\oplus
( \frac{1}{12},
\frac{7}{12} )
$,

\vskip 0.7ex
\hangindent=5.5em \hangafter=1
{\white .}\hskip 1em $\rho_\text{isum}(\mathfrak{s})$ =
$\mathrm{i}$($-\sqrt{\frac{1}{5}}$,
$\sqrt{\frac{2}{5}}$,
$\sqrt{\frac{2}{5}}$;\ \ 
$\frac{-5+\sqrt{5}}{10}$,
$\frac{5+\sqrt{5}}{10}$;\ \ 
$\frac{-5+\sqrt{5}}{10}$)
 $\oplus$
$\mathrm{i}$($\frac{1}{2}$,
$\sqrt{\frac{3}{4}}$;\ \ 
$-\frac{1}{2}$)

Pass. 

 \ \color{black}

\noindent 114: (dims,levels) = $(3\oplus
2;60,
12
)$,
irreps = $3_{5}^{3}
\hskip -1.5pt \otimes \hskip -1.5pt
1_{4}^{1,0}
\hskip -1.5pt \otimes \hskip -1.5pt
1_{3}^{1,0}\oplus
2_{3}^{1,0}
\hskip -1.5pt \otimes \hskip -1.5pt
1_{4}^{1,0}$,
pord$(\rho_\text{isum}(\mathfrak{t})) = 15$,

\vskip 0.7ex
\hangindent=5.5em \hangafter=1
{\white .}\hskip 1em $\rho_\text{isum}(\mathfrak{t})$ =
 $( \frac{7}{12},
\frac{11}{60},
\frac{59}{60} )
\oplus
( \frac{1}{4},
\frac{7}{12} )
$,

\vskip 0.7ex
\hangindent=5.5em \hangafter=1
{\white .}\hskip 1em $\rho_\text{isum}(\mathfrak{s})$ =
$\mathrm{i}$($-\sqrt{\frac{1}{5}}$,
$\sqrt{\frac{2}{5}}$,
$\sqrt{\frac{2}{5}}$;\ \ 
$\frac{-5+\sqrt{5}}{10}$,
$\frac{5+\sqrt{5}}{10}$;\ \ 
$\frac{-5+\sqrt{5}}{10}$)
 $\oplus$
($\sqrt{\frac{1}{3}}$,
$\sqrt{\frac{2}{3}}$;
$-\sqrt{\frac{1}{3}}$)

Fail:
Tr$_I(C) = -1 <$  0 for I = [ 1/4 ]. Prop. B.4 (1) eqn. (B.18)

 \ \color{black}

\noindent 115: (dims,levels) = $(3\oplus
2;60,
12
)$,
irreps = $3_{5}^{3}
\hskip -1.5pt \otimes \hskip -1.5pt
1_{4}^{1,0}
\hskip -1.5pt \otimes \hskip -1.5pt
1_{3}^{1,0}\oplus
2_{3}^{1,4}
\hskip -1.5pt \otimes \hskip -1.5pt
1_{4}^{1,0}$,
pord$(\rho_\text{isum}(\mathfrak{t})) = 15$,

\vskip 0.7ex
\hangindent=5.5em \hangafter=1
{\white .}\hskip 1em $\rho_\text{isum}(\mathfrak{t})$ =
 $( \frac{7}{12},
\frac{11}{60},
\frac{59}{60} )
\oplus
( \frac{7}{12},
\frac{11}{12} )
$,

\vskip 0.7ex
\hangindent=5.5em \hangafter=1
{\white .}\hskip 1em $\rho_\text{isum}(\mathfrak{s})$ =
$\mathrm{i}$($-\sqrt{\frac{1}{5}}$,
$\sqrt{\frac{2}{5}}$,
$\sqrt{\frac{2}{5}}$;\ \ 
$\frac{-5+\sqrt{5}}{10}$,
$\frac{5+\sqrt{5}}{10}$;\ \ 
$\frac{-5+\sqrt{5}}{10}$)
 $\oplus$
($\sqrt{\frac{1}{3}}$,
$\sqrt{\frac{2}{3}}$;
$-\sqrt{\frac{1}{3}}$)

Fail:
Tr$_I(C) = -1 <$  0 for I = [ 11/12 ]. Prop. B.4 (1) eqn. (B.18)

 \ \color{black}

 \color{blue}

\noindent 116: (dims,levels) = $(3\oplus
2;60,
12
)$,
irreps = $3_{5}^{1}
\hskip -1.5pt \otimes \hskip -1.5pt
1_{4}^{1,0}
\hskip -1.5pt \otimes \hskip -1.5pt
1_{3}^{1,0}\oplus
2_{4}^{1,0}
\hskip -1.5pt \otimes \hskip -1.5pt
1_{3}^{1,0}$,
pord$(\rho_\text{isum}(\mathfrak{t})) = 10$,

\vskip 0.7ex
\hangindent=5.5em \hangafter=1
{\white .}\hskip 1em $\rho_\text{isum}(\mathfrak{t})$ =
 $( \frac{7}{12},
\frac{23}{60},
\frac{47}{60} )
\oplus
( \frac{1}{12},
\frac{7}{12} )
$,

\vskip 0.7ex
\hangindent=5.5em \hangafter=1
{\white .}\hskip 1em $\rho_\text{isum}(\mathfrak{s})$ =
$\mathrm{i}$($\sqrt{\frac{1}{5}}$,
$\sqrt{\frac{2}{5}}$,
$\sqrt{\frac{2}{5}}$;\ \ 
$-\frac{5+\sqrt{5}}{10}$,
$\frac{5-\sqrt{5}}{10}$;\ \ 
$-\frac{5+\sqrt{5}}{10}$)
 $\oplus$
$\mathrm{i}$($\frac{1}{2}$,
$\sqrt{\frac{3}{4}}$;\ \ 
$-\frac{1}{2}$)

Pass. 

 \ \color{black}

\noindent 117: (dims,levels) = $(3\oplus
2;60,
12
)$,
irreps = $3_{5}^{1}
\hskip -1.5pt \otimes \hskip -1.5pt
1_{4}^{1,0}
\hskip -1.5pt \otimes \hskip -1.5pt
1_{3}^{1,0}\oplus
2_{3}^{1,0}
\hskip -1.5pt \otimes \hskip -1.5pt
1_{4}^{1,0}$,
pord$(\rho_\text{isum}(\mathfrak{t})) = 15$,

\vskip 0.7ex
\hangindent=5.5em \hangafter=1
{\white .}\hskip 1em $\rho_\text{isum}(\mathfrak{t})$ =
 $( \frac{7}{12},
\frac{23}{60},
\frac{47}{60} )
\oplus
( \frac{1}{4},
\frac{7}{12} )
$,

\vskip 0.7ex
\hangindent=5.5em \hangafter=1
{\white .}\hskip 1em $\rho_\text{isum}(\mathfrak{s})$ =
$\mathrm{i}$($\sqrt{\frac{1}{5}}$,
$\sqrt{\frac{2}{5}}$,
$\sqrt{\frac{2}{5}}$;\ \ 
$-\frac{5+\sqrt{5}}{10}$,
$\frac{5-\sqrt{5}}{10}$;\ \ 
$-\frac{5+\sqrt{5}}{10}$)
 $\oplus$
($\sqrt{\frac{1}{3}}$,
$\sqrt{\frac{2}{3}}$;
$-\sqrt{\frac{1}{3}}$)

Fail:
Tr$_I(C) = -1 <$  0 for I = [ 1/4 ]. Prop. B.4 (1) eqn. (B.18)

 \ \color{black}

\noindent 118: (dims,levels) = $(3\oplus
2;60,
12
)$,
irreps = $3_{5}^{1}
\hskip -1.5pt \otimes \hskip -1.5pt
1_{4}^{1,0}
\hskip -1.5pt \otimes \hskip -1.5pt
1_{3}^{1,0}\oplus
2_{3}^{1,4}
\hskip -1.5pt \otimes \hskip -1.5pt
1_{4}^{1,0}$,
pord$(\rho_\text{isum}(\mathfrak{t})) = 15$,

\vskip 0.7ex
\hangindent=5.5em \hangafter=1
{\white .}\hskip 1em $\rho_\text{isum}(\mathfrak{t})$ =
 $( \frac{7}{12},
\frac{23}{60},
\frac{47}{60} )
\oplus
( \frac{7}{12},
\frac{11}{12} )
$,

\vskip 0.7ex
\hangindent=5.5em \hangafter=1
{\white .}\hskip 1em $\rho_\text{isum}(\mathfrak{s})$ =
$\mathrm{i}$($\sqrt{\frac{1}{5}}$,
$\sqrt{\frac{2}{5}}$,
$\sqrt{\frac{2}{5}}$;\ \ 
$-\frac{5+\sqrt{5}}{10}$,
$\frac{5-\sqrt{5}}{10}$;\ \ 
$-\frac{5+\sqrt{5}}{10}$)
 $\oplus$
($\sqrt{\frac{1}{3}}$,
$\sqrt{\frac{2}{3}}$;
$-\sqrt{\frac{1}{3}}$)

Fail:
Tr$_I(C) = -1 <$  0 for I = [ 11/12 ]. Prop. B.4 (1) eqn. (B.18)

 \ \color{black}

\noindent 119: (dims,levels) = $(3\oplus
2;60,
60
)$,
irreps = $3_{5}^{3}
\hskip -1.5pt \otimes \hskip -1.5pt
1_{4}^{1,0}
\hskip -1.5pt \otimes \hskip -1.5pt
1_{3}^{1,0}\oplus
2_{5}^{2}
\hskip -1.5pt \otimes \hskip -1.5pt
1_{4}^{1,0}
\hskip -1.5pt \otimes \hskip -1.5pt
1_{3}^{1,0}$,
pord$(\rho_\text{isum}(\mathfrak{t})) = 5$,

\vskip 0.7ex
\hangindent=5.5em \hangafter=1
{\white .}\hskip 1em $\rho_\text{isum}(\mathfrak{t})$ =
 $( \frac{7}{12},
\frac{11}{60},
\frac{59}{60} )
\oplus
( \frac{11}{60},
\frac{59}{60} )
$,

\vskip 0.7ex
\hangindent=5.5em \hangafter=1
{\white .}\hskip 1em $\rho_\text{isum}(\mathfrak{s})$ =
$\mathrm{i}$($-\sqrt{\frac{1}{5}}$,
$\sqrt{\frac{2}{5}}$,
$\sqrt{\frac{2}{5}}$;\ \ 
$\frac{-5+\sqrt{5}}{10}$,
$\frac{5+\sqrt{5}}{10}$;\ \ 
$\frac{-5+\sqrt{5}}{10}$)
 $\oplus$
($-\frac{1}{\sqrt{5}}c^{1}_{20}
$,
$\frac{1}{\sqrt{5}}c^{3}_{20}
$;
$\frac{1}{\sqrt{5}}c^{1}_{20}
$)

Fail:
Integral: $D_{\rho}(\sigma)_{\theta} \propto $ id,
 for all $\sigma$ and all $\theta$-eigenspaces that can contain unit. Prop. B.5 (6)

 \ \color{black}

\noindent 120: (dims,levels) = $(3\oplus
2;60,
60
)$,
irreps = $3_{5}^{1}
\hskip -1.5pt \otimes \hskip -1.5pt
1_{4}^{1,0}
\hskip -1.5pt \otimes \hskip -1.5pt
1_{3}^{1,0}\oplus
2_{5}^{1}
\hskip -1.5pt \otimes \hskip -1.5pt
1_{4}^{1,0}
\hskip -1.5pt \otimes \hskip -1.5pt
1_{3}^{1,0}$,
pord$(\rho_\text{isum}(\mathfrak{t})) = 5$,

\vskip 0.7ex
\hangindent=5.5em \hangafter=1
{\white .}\hskip 1em $\rho_\text{isum}(\mathfrak{t})$ =
 $( \frac{7}{12},
\frac{23}{60},
\frac{47}{60} )
\oplus
( \frac{23}{60},
\frac{47}{60} )
$,

\vskip 0.7ex
\hangindent=5.5em \hangafter=1
{\white .}\hskip 1em $\rho_\text{isum}(\mathfrak{s})$ =
$\mathrm{i}$($\sqrt{\frac{1}{5}}$,
$\sqrt{\frac{2}{5}}$,
$\sqrt{\frac{2}{5}}$;\ \ 
$-\frac{5+\sqrt{5}}{10}$,
$\frac{5-\sqrt{5}}{10}$;\ \ 
$-\frac{5+\sqrt{5}}{10}$)
 $\oplus$
($-\frac{1}{\sqrt{5}}c^{3}_{20}
$,
$\frac{1}{\sqrt{5}}c^{1}_{20}
$;
$\frac{1}{\sqrt{5}}c^{3}_{20}
$)

Fail:
Integral: $D_{\rho}(\sigma)_{\theta} \propto $ id,
 for all $\sigma$ and all $\theta$-eigenspaces that can contain unit. Prop. B.5 (6)

 \ \color{black}

 \color{blue}

\noindent 121: (dims,levels) = $(4\oplus
1;7,
1
)$,
irreps = $4_{7}^{1}\oplus
1_{1}^{1}$,
pord$(\rho_\text{isum}(\mathfrak{t})) = 7$,

\vskip 0.7ex
\hangindent=5.5em \hangafter=1
{\white .}\hskip 1em $\rho_\text{isum}(\mathfrak{t})$ =
 $( 0,
\frac{1}{7},
\frac{2}{7},
\frac{4}{7} )
\oplus
( 0 )
$,

\vskip 0.7ex
\hangindent=5.5em \hangafter=1
{\white .}\hskip 1em $\rho_\text{isum}(\mathfrak{s})$ =
$\mathrm{i}$($-\sqrt{\frac{1}{7}}$,
$\sqrt{\frac{2}{7}}$,
$\sqrt{\frac{2}{7}}$,
$\sqrt{\frac{2}{7}}$;\ \ 
$-\frac{1}{\sqrt{7}}c^{2}_{7}
$,
$-\frac{1}{\sqrt{7}}c^{1}_{7}
$,
$\frac{1}{\sqrt{7}\mathrm{i}}s^{5}_{28}
$;\ \ 
$\frac{1}{\sqrt{7}\mathrm{i}}s^{5}_{28}
$,
$-\frac{1}{\sqrt{7}}c^{2}_{7}
$;\ \ 
$-\frac{1}{\sqrt{7}}c^{1}_{7}
$)
 $\oplus$
($1$)

Pass. 

 \ \color{black}

 \color{blue}

\noindent 122: (dims,levels) = $(4\oplus
1;7,
1
)$,
irreps = $4_{7}^{3}\oplus
1_{1}^{1}$,
pord$(\rho_\text{isum}(\mathfrak{t})) = 7$,

\vskip 0.7ex
\hangindent=5.5em \hangafter=1
{\white .}\hskip 1em $\rho_\text{isum}(\mathfrak{t})$ =
 $( 0,
\frac{3}{7},
\frac{5}{7},
\frac{6}{7} )
\oplus
( 0 )
$,

\vskip 0.7ex
\hangindent=5.5em \hangafter=1
{\white .}\hskip 1em $\rho_\text{isum}(\mathfrak{s})$ =
$\mathrm{i}$($\sqrt{\frac{1}{7}}$,
$\sqrt{\frac{2}{7}}$,
$\sqrt{\frac{2}{7}}$,
$\sqrt{\frac{2}{7}}$;\ \ 
$\frac{1}{\sqrt{7}}c^{1}_{7}
$,
$\frac{1}{\sqrt{7}}c^{2}_{7}
$,
$-\frac{1}{\sqrt{7}\mathrm{i}}s^{5}_{28}
$;\ \ 
$-\frac{1}{\sqrt{7}\mathrm{i}}s^{5}_{28}
$,
$\frac{1}{\sqrt{7}}c^{1}_{7}
$;\ \ 
$\frac{1}{\sqrt{7}}c^{2}_{7}
$)
 $\oplus$
($1$)

Pass. 

 \ \color{black}

 \color{blue}

\noindent 123: (dims,levels) = $(4\oplus
1;9,
1
)$,
irreps = $4_{9,1}^{1,0}\oplus
1_{1}^{1}$,
pord$(\rho_\text{isum}(\mathfrak{t})) = 9$,

\vskip 0.7ex
\hangindent=5.5em \hangafter=1
{\white .}\hskip 1em $\rho_\text{isum}(\mathfrak{t})$ =
 $( 0,
\frac{1}{9},
\frac{4}{9},
\frac{7}{9} )
\oplus
( 0 )
$,

\vskip 0.7ex
\hangindent=5.5em \hangafter=1
{\white .}\hskip 1em $\rho_\text{isum}(\mathfrak{s})$ =
$\mathrm{i}$($0$,
$\sqrt{\frac{1}{3}}$,
$\sqrt{\frac{1}{3}}$,
$\sqrt{\frac{1}{3}}$;\ \ 
$-\frac{1}{3}c^{1}_{36}
$,
$\frac{1}{3}c^{1}_{36}
-\frac{1}{3}c^{5}_{36}
$,
$\frac{1}{3}c^{5}_{36}
$;\ \ 
$\frac{1}{3}c^{5}_{36}
$,
$-\frac{1}{3}c^{1}_{36}
$;\ \ 
$\frac{1}{3}c^{1}_{36}
-\frac{1}{3}c^{5}_{36}
$)
 $\oplus$
($1$)

Pass. 

 \ \color{black}

 \color{blue}

\noindent 124: (dims,levels) = $(4\oplus
1;9,
1
)$,
irreps = $4_{9,2}^{1,0}\oplus
1_{1}^{1}$,
pord$(\rho_\text{isum}(\mathfrak{t})) = 9$,

\vskip 0.7ex
\hangindent=5.5em \hangafter=1
{\white .}\hskip 1em $\rho_\text{isum}(\mathfrak{t})$ =
 $( 0,
\frac{1}{9},
\frac{4}{9},
\frac{7}{9} )
\oplus
( 0 )
$,

\vskip 0.7ex
\hangindent=5.5em \hangafter=1
{\white .}\hskip 1em $\rho_\text{isum}(\mathfrak{s})$ =
($0$,
$-\sqrt{\frac{1}{3}}$,
$-\sqrt{\frac{1}{3}}$,
$-\sqrt{\frac{1}{3}}$;
$\frac{1}{3}c^{2}_{9}
$,
$\frac{1}{3} c_9^4 $,
$\frac{1}{3}c^{1}_{9}
$;
$\frac{1}{3}c^{1}_{9}
$,
$\frac{1}{3}c^{2}_{9}
$;
$\frac{1}{3} c_9^4 $)
 $\oplus$
($1$)

Pass. 

 \ \color{black}

 \color{blue}

\noindent 125: (dims,levels) = $(4\oplus
1;9,
1
)$,
irreps = $4_{9,1}^{2,0}\oplus
1_{1}^{1}$,
pord$(\rho_\text{isum}(\mathfrak{t})) = 9$,

\vskip 0.7ex
\hangindent=5.5em \hangafter=1
{\white .}\hskip 1em $\rho_\text{isum}(\mathfrak{t})$ =
 $( 0,
\frac{2}{9},
\frac{5}{9},
\frac{8}{9} )
\oplus
( 0 )
$,

\vskip 0.7ex
\hangindent=5.5em \hangafter=1
{\white .}\hskip 1em $\rho_\text{isum}(\mathfrak{s})$ =
$\mathrm{i}$($0$,
$\sqrt{\frac{1}{3}}$,
$\sqrt{\frac{1}{3}}$,
$\sqrt{\frac{1}{3}}$;\ \ 
$-\frac{1}{3}c^{1}_{36}
+\frac{1}{3}c^{5}_{36}
$,
$\frac{1}{3}c^{1}_{36}
$,
$-\frac{1}{3}c^{5}_{36}
$;\ \ 
$-\frac{1}{3}c^{5}_{36}
$,
$-\frac{1}{3}c^{1}_{36}
+\frac{1}{3}c^{5}_{36}
$;\ \ 
$\frac{1}{3}c^{1}_{36}
$)
 $\oplus$
($1$)

Pass. 

 \ \color{black}

 \color{blue}

\noindent 126: (dims,levels) = $(4\oplus
1;9,
1
)$,
irreps = $4_{9,2}^{5,0}\oplus
1_{1}^{1}$,
pord$(\rho_\text{isum}(\mathfrak{t})) = 9$,

\vskip 0.7ex
\hangindent=5.5em \hangafter=1
{\white .}\hskip 1em $\rho_\text{isum}(\mathfrak{t})$ =
 $( 0,
\frac{2}{9},
\frac{5}{9},
\frac{8}{9} )
\oplus
( 0 )
$,

\vskip 0.7ex
\hangindent=5.5em \hangafter=1
{\white .}\hskip 1em $\rho_\text{isum}(\mathfrak{s})$ =
($0$,
$-\sqrt{\frac{1}{3}}$,
$-\sqrt{\frac{1}{3}}$,
$-\sqrt{\frac{1}{3}}$;
$\frac{1}{3} c_9^4 $,
$\frac{1}{3}c^{2}_{9}
$,
$\frac{1}{3}c^{1}_{9}
$;
$\frac{1}{3}c^{1}_{9}
$,
$\frac{1}{3} c_9^4 $;
$\frac{1}{3}c^{2}_{9}
$)
 $\oplus$
($1$)

Pass. 

 \ \color{black}

 \color{blue}

\noindent 127: (dims,levels) = $(4\oplus
1;14,
2
)$,
irreps = $4_{7}^{1}
\hskip -1.5pt \otimes \hskip -1.5pt
1_{2}^{1,0}\oplus
1_{2}^{1,0}$,
pord$(\rho_\text{isum}(\mathfrak{t})) = 7$,

\vskip 0.7ex
\hangindent=5.5em \hangafter=1
{\white .}\hskip 1em $\rho_\text{isum}(\mathfrak{t})$ =
 $( \frac{1}{2},
\frac{1}{14},
\frac{9}{14},
\frac{11}{14} )
\oplus
( \frac{1}{2} )
$,

\vskip 0.7ex
\hangindent=5.5em \hangafter=1
{\white .}\hskip 1em $\rho_\text{isum}(\mathfrak{s})$ =
$\mathrm{i}$($\sqrt{\frac{1}{7}}$,
$\sqrt{\frac{2}{7}}$,
$\sqrt{\frac{2}{7}}$,
$\sqrt{\frac{2}{7}}$;\ \ 
$\frac{1}{\sqrt{7}}c^{1}_{7}
$,
$-\frac{1}{\sqrt{7}\mathrm{i}}s^{5}_{28}
$,
$\frac{1}{\sqrt{7}}c^{2}_{7}
$;\ \ 
$\frac{1}{\sqrt{7}}c^{2}_{7}
$,
$\frac{1}{\sqrt{7}}c^{1}_{7}
$;\ \ 
$-\frac{1}{\sqrt{7}\mathrm{i}}s^{5}_{28}
$)
 $\oplus$
($-1$)

Pass. 

 \ \color{black}

 \color{blue}

\noindent 128: (dims,levels) = $(4\oplus
1;14,
2
)$,
irreps = $4_{7}^{3}
\hskip -1.5pt \otimes \hskip -1.5pt
1_{2}^{1,0}\oplus
1_{2}^{1,0}$,
pord$(\rho_\text{isum}(\mathfrak{t})) = 7$,

\vskip 0.7ex
\hangindent=5.5em \hangafter=1
{\white .}\hskip 1em $\rho_\text{isum}(\mathfrak{t})$ =
 $( \frac{1}{2},
\frac{3}{14},
\frac{5}{14},
\frac{13}{14} )
\oplus
( \frac{1}{2} )
$,

\vskip 0.7ex
\hangindent=5.5em \hangafter=1
{\white .}\hskip 1em $\rho_\text{isum}(\mathfrak{s})$ =
$\mathrm{i}$($-\sqrt{\frac{1}{7}}$,
$\sqrt{\frac{2}{7}}$,
$\sqrt{\frac{2}{7}}$,
$\sqrt{\frac{2}{7}}$;\ \ 
$\frac{1}{\sqrt{7}\mathrm{i}}s^{5}_{28}
$,
$-\frac{1}{\sqrt{7}}c^{1}_{7}
$,
$-\frac{1}{\sqrt{7}}c^{2}_{7}
$;\ \ 
$-\frac{1}{\sqrt{7}}c^{2}_{7}
$,
$\frac{1}{\sqrt{7}\mathrm{i}}s^{5}_{28}
$;\ \ 
$-\frac{1}{\sqrt{7}}c^{1}_{7}
$)
 $\oplus$
($-1$)

Pass. 

 \ \color{black}

 \color{blue}

\noindent 129: (dims,levels) = $(4\oplus
1;18,
2
)$,
irreps = $4_{9,1}^{2,0}
\hskip -1.5pt \otimes \hskip -1.5pt
1_{2}^{1,0}\oplus
1_{2}^{1,0}$,
pord$(\rho_\text{isum}(\mathfrak{t})) = 9$,

\vskip 0.7ex
\hangindent=5.5em \hangafter=1
{\white .}\hskip 1em $\rho_\text{isum}(\mathfrak{t})$ =
 $( \frac{1}{2},
\frac{1}{18},
\frac{7}{18},
\frac{13}{18} )
\oplus
( \frac{1}{2} )
$,

\vskip 0.7ex
\hangindent=5.5em \hangafter=1
{\white .}\hskip 1em $\rho_\text{isum}(\mathfrak{s})$ =
$\mathrm{i}$($0$,
$\sqrt{\frac{1}{3}}$,
$\sqrt{\frac{1}{3}}$,
$\sqrt{\frac{1}{3}}$;\ \ 
$\frac{1}{3}c^{5}_{36}
$,
$\frac{1}{3}c^{1}_{36}
-\frac{1}{3}c^{5}_{36}
$,
$-\frac{1}{3}c^{1}_{36}
$;\ \ 
$-\frac{1}{3}c^{1}_{36}
$,
$\frac{1}{3}c^{5}_{36}
$;\ \ 
$\frac{1}{3}c^{1}_{36}
-\frac{1}{3}c^{5}_{36}
$)
 $\oplus$
($-1$)

Pass. 

 \ \color{black}

 \color{blue}

\noindent 130: (dims,levels) = $(4\oplus
1;18,
2
)$,
irreps = $4_{9,2}^{5,0}
\hskip -1.5pt \otimes \hskip -1.5pt
1_{2}^{1,0}\oplus
1_{2}^{1,0}$,
pord$(\rho_\text{isum}(\mathfrak{t})) = 9$,

\vskip 0.7ex
\hangindent=5.5em \hangafter=1
{\white .}\hskip 1em $\rho_\text{isum}(\mathfrak{t})$ =
 $( \frac{1}{2},
\frac{1}{18},
\frac{7}{18},
\frac{13}{18} )
\oplus
( \frac{1}{2} )
$,

\vskip 0.7ex
\hangindent=5.5em \hangafter=1
{\white .}\hskip 1em $\rho_\text{isum}(\mathfrak{s})$ =
($0$,
$-\sqrt{\frac{1}{3}}$,
$-\sqrt{\frac{1}{3}}$,
$-\sqrt{\frac{1}{3}}$;
$-\frac{1}{3}c^{1}_{9}
$,
$-\frac{1}{3} c_9^4 $,
$-\frac{1}{3}c^{2}_{9}
$;
$-\frac{1}{3}c^{2}_{9}
$,
$-\frac{1}{3}c^{1}_{9}
$;
$-\frac{1}{3} c_9^4 $)
 $\oplus$
($-1$)

Pass. 

 \ \color{black}

 \color{blue}

\noindent 131: (dims,levels) = $(4\oplus
1;18,
2
)$,
irreps = $4_{9,1}^{1,0}
\hskip -1.5pt \otimes \hskip -1.5pt
1_{2}^{1,0}\oplus
1_{2}^{1,0}$,
pord$(\rho_\text{isum}(\mathfrak{t})) = 9$,

\vskip 0.7ex
\hangindent=5.5em \hangafter=1
{\white .}\hskip 1em $\rho_\text{isum}(\mathfrak{t})$ =
 $( \frac{1}{2},
\frac{5}{18},
\frac{11}{18},
\frac{17}{18} )
\oplus
( \frac{1}{2} )
$,

\vskip 0.7ex
\hangindent=5.5em \hangafter=1
{\white .}\hskip 1em $\rho_\text{isum}(\mathfrak{s})$ =
$\mathrm{i}$($0$,
$\sqrt{\frac{1}{3}}$,
$\sqrt{\frac{1}{3}}$,
$\sqrt{\frac{1}{3}}$;\ \ 
$-\frac{1}{3}c^{1}_{36}
+\frac{1}{3}c^{5}_{36}
$,
$-\frac{1}{3}c^{5}_{36}
$,
$\frac{1}{3}c^{1}_{36}
$;\ \ 
$\frac{1}{3}c^{1}_{36}
$,
$-\frac{1}{3}c^{1}_{36}
+\frac{1}{3}c^{5}_{36}
$;\ \ 
$-\frac{1}{3}c^{5}_{36}
$)
 $\oplus$
($-1$)

Pass. 

 \ \color{black}

 \color{blue}

\noindent 132: (dims,levels) = $(4\oplus
1;18,
2
)$,
irreps = $4_{9,2}^{1,0}
\hskip -1.5pt \otimes \hskip -1.5pt
1_{2}^{1,0}\oplus
1_{2}^{1,0}$,
pord$(\rho_\text{isum}(\mathfrak{t})) = 9$,

\vskip 0.7ex
\hangindent=5.5em \hangafter=1
{\white .}\hskip 1em $\rho_\text{isum}(\mathfrak{t})$ =
 $( \frac{1}{2},
\frac{5}{18},
\frac{11}{18},
\frac{17}{18} )
\oplus
( \frac{1}{2} )
$,

\vskip 0.7ex
\hangindent=5.5em \hangafter=1
{\white .}\hskip 1em $\rho_\text{isum}(\mathfrak{s})$ =
($0$,
$-\sqrt{\frac{1}{3}}$,
$-\sqrt{\frac{1}{3}}$,
$-\sqrt{\frac{1}{3}}$;
$-\frac{1}{3} c_9^4 $,
$-\frac{1}{3}c^{1}_{9}
$,
$-\frac{1}{3}c^{2}_{9}
$;
$-\frac{1}{3}c^{2}_{9}
$,
$-\frac{1}{3} c_9^4 $;
$-\frac{1}{3}c^{1}_{9}
$)
 $\oplus$
($-1$)

Pass. 

 \ \color{black}

 \color{blue}

\noindent 133: (dims,levels) = $(4\oplus
1;21,
3
)$,
irreps = $4_{7}^{3}
\hskip -1.5pt \otimes \hskip -1.5pt
1_{3}^{1,0}\oplus
1_{3}^{1,0}$,
pord$(\rho_\text{isum}(\mathfrak{t})) = 7$,

\vskip 0.7ex
\hangindent=5.5em \hangafter=1
{\white .}\hskip 1em $\rho_\text{isum}(\mathfrak{t})$ =
 $( \frac{1}{3},
\frac{1}{21},
\frac{4}{21},
\frac{16}{21} )
\oplus
( \frac{1}{3} )
$,

\vskip 0.7ex
\hangindent=5.5em \hangafter=1
{\white .}\hskip 1em $\rho_\text{isum}(\mathfrak{s})$ =
$\mathrm{i}$($\sqrt{\frac{1}{7}}$,
$\sqrt{\frac{2}{7}}$,
$\sqrt{\frac{2}{7}}$,
$\sqrt{\frac{2}{7}}$;\ \ 
$-\frac{1}{\sqrt{7}\mathrm{i}}s^{5}_{28}
$,
$\frac{1}{\sqrt{7}}c^{1}_{7}
$,
$\frac{1}{\sqrt{7}}c^{2}_{7}
$;\ \ 
$\frac{1}{\sqrt{7}}c^{2}_{7}
$,
$-\frac{1}{\sqrt{7}\mathrm{i}}s^{5}_{28}
$;\ \ 
$\frac{1}{\sqrt{7}}c^{1}_{7}
$)
 $\oplus$
($1$)

Pass. 

 \ \color{black}

 \color{blue}

\noindent 134: (dims,levels) = $(4\oplus
1;21,
3
)$,
irreps = $4_{7}^{1}
\hskip -1.5pt \otimes \hskip -1.5pt
1_{3}^{1,0}\oplus
1_{3}^{1,0}$,
pord$(\rho_\text{isum}(\mathfrak{t})) = 7$,

\vskip 0.7ex
\hangindent=5.5em \hangafter=1
{\white .}\hskip 1em $\rho_\text{isum}(\mathfrak{t})$ =
 $( \frac{1}{3},
\frac{10}{21},
\frac{13}{21},
\frac{19}{21} )
\oplus
( \frac{1}{3} )
$,

\vskip 0.7ex
\hangindent=5.5em \hangafter=1
{\white .}\hskip 1em $\rho_\text{isum}(\mathfrak{s})$ =
$\mathrm{i}$($-\sqrt{\frac{1}{7}}$,
$\sqrt{\frac{2}{7}}$,
$\sqrt{\frac{2}{7}}$,
$\sqrt{\frac{2}{7}}$;\ \ 
$-\frac{1}{\sqrt{7}}c^{2}_{7}
$,
$-\frac{1}{\sqrt{7}}c^{1}_{7}
$,
$\frac{1}{\sqrt{7}\mathrm{i}}s^{5}_{28}
$;\ \ 
$\frac{1}{\sqrt{7}\mathrm{i}}s^{5}_{28}
$,
$-\frac{1}{\sqrt{7}}c^{2}_{7}
$;\ \ 
$-\frac{1}{\sqrt{7}}c^{1}_{7}
$)
 $\oplus$
($1$)

Pass. 

 \ \color{black}

 \color{blue}

\noindent 135: (dims,levels) = $(4\oplus
1;28,
4
)$,
irreps = $4_{7}^{3}
\hskip -1.5pt \otimes \hskip -1.5pt
1_{4}^{1,0}\oplus
1_{4}^{1,0}$,
pord$(\rho_\text{isum}(\mathfrak{t})) = 7$,

\vskip 0.7ex
\hangindent=5.5em \hangafter=1
{\white .}\hskip 1em $\rho_\text{isum}(\mathfrak{t})$ =
 $( \frac{1}{4},
\frac{3}{28},
\frac{19}{28},
\frac{27}{28} )
\oplus
( \frac{1}{4} )
$,

\vskip 0.7ex
\hangindent=5.5em \hangafter=1
{\white .}\hskip 1em $\rho_\text{isum}(\mathfrak{s})$ =
($-\sqrt{\frac{1}{7}}$,
$\sqrt{\frac{2}{7}}$,
$\sqrt{\frac{2}{7}}$,
$\sqrt{\frac{2}{7}}$;
$-\frac{1}{\sqrt{7}}c^{2}_{7}
$,
$\frac{1}{\sqrt{7}\mathrm{i}}s^{5}_{28}
$,
$-\frac{1}{\sqrt{7}}c^{1}_{7}
$;
$-\frac{1}{\sqrt{7}}c^{1}_{7}
$,
$-\frac{1}{\sqrt{7}}c^{2}_{7}
$;
$\frac{1}{\sqrt{7}\mathrm{i}}s^{5}_{28}
$)
 $\oplus$
$\mathrm{i}$($1$)

Pass. 

 \ \color{black}

 \color{blue}

\noindent 136: (dims,levels) = $(4\oplus
1;28,
4
)$,
irreps = $4_{7}^{1}
\hskip -1.5pt \otimes \hskip -1.5pt
1_{4}^{1,0}\oplus
1_{4}^{1,0}$,
pord$(\rho_\text{isum}(\mathfrak{t})) = 7$,

\vskip 0.7ex
\hangindent=5.5em \hangafter=1
{\white .}\hskip 1em $\rho_\text{isum}(\mathfrak{t})$ =
 $( \frac{1}{4},
\frac{11}{28},
\frac{15}{28},
\frac{23}{28} )
\oplus
( \frac{1}{4} )
$,

\vskip 0.7ex
\hangindent=5.5em \hangafter=1
{\white .}\hskip 1em $\rho_\text{isum}(\mathfrak{s})$ =
($\sqrt{\frac{1}{7}}$,
$\sqrt{\frac{2}{7}}$,
$\sqrt{\frac{2}{7}}$,
$\sqrt{\frac{2}{7}}$;
$\frac{1}{\sqrt{7}}c^{2}_{7}
$,
$\frac{1}{\sqrt{7}}c^{1}_{7}
$,
$-\frac{1}{\sqrt{7}\mathrm{i}}s^{5}_{28}
$;
$-\frac{1}{\sqrt{7}\mathrm{i}}s^{5}_{28}
$,
$\frac{1}{\sqrt{7}}c^{2}_{7}
$;
$\frac{1}{\sqrt{7}}c^{1}_{7}
$)
 $\oplus$
$\mathrm{i}$($1$)

Pass. 

 \ \color{black}

 \color{blue}

\noindent 137: (dims,levels) = $(4\oplus
1;36,
4
)$,
irreps = $4_{9,2}^{1,0}
\hskip -1.5pt \otimes \hskip -1.5pt
1_{4}^{1,0}\oplus
1_{4}^{1,0}$,
pord$(\rho_\text{isum}(\mathfrak{t})) = 9$,

\vskip 0.7ex
\hangindent=5.5em \hangafter=1
{\white .}\hskip 1em $\rho_\text{isum}(\mathfrak{t})$ =
 $( \frac{1}{4},
\frac{1}{36},
\frac{13}{36},
\frac{25}{36} )
\oplus
( \frac{1}{4} )
$,

\vskip 0.7ex
\hangindent=5.5em \hangafter=1
{\white .}\hskip 1em $\rho_\text{isum}(\mathfrak{s})$ =
$\mathrm{i}$($0$,
$\sqrt{\frac{1}{3}}$,
$\sqrt{\frac{1}{3}}$,
$\sqrt{\frac{1}{3}}$;\ \ 
$\frac{1}{3} c_9^4 $,
$\frac{1}{3}c^{1}_{9}
$,
$\frac{1}{3}c^{2}_{9}
$;\ \ 
$\frac{1}{3}c^{2}_{9}
$,
$\frac{1}{3} c_9^4 $;\ \ 
$\frac{1}{3}c^{1}_{9}
$)
 $\oplus$
$\mathrm{i}$($1$)

Pass. 

 \ \color{black}

 \color{blue}

\noindent 138: (dims,levels) = $(4\oplus
1;36,
4
)$,
irreps = $4_{9,1}^{1,0}
\hskip -1.5pt \otimes \hskip -1.5pt
1_{4}^{1,0}\oplus
1_{4}^{1,0}$,
pord$(\rho_\text{isum}(\mathfrak{t})) = 9$,

\vskip 0.7ex
\hangindent=5.5em \hangafter=1
{\white .}\hskip 1em $\rho_\text{isum}(\mathfrak{t})$ =
 $( \frac{1}{4},
\frac{1}{36},
\frac{13}{36},
\frac{25}{36} )
\oplus
( \frac{1}{4} )
$,

\vskip 0.7ex
\hangindent=5.5em \hangafter=1
{\white .}\hskip 1em $\rho_\text{isum}(\mathfrak{s})$ =
($0$,
$-\sqrt{\frac{1}{3}}$,
$-\sqrt{\frac{1}{3}}$,
$-\sqrt{\frac{1}{3}}$;
$-\frac{1}{3}c^{1}_{36}
+\frac{1}{3}c^{5}_{36}
$,
$-\frac{1}{3}c^{5}_{36}
$,
$\frac{1}{3}c^{1}_{36}
$;
$\frac{1}{3}c^{1}_{36}
$,
$-\frac{1}{3}c^{1}_{36}
+\frac{1}{3}c^{5}_{36}
$;
$-\frac{1}{3}c^{5}_{36}
$)
 $\oplus$
$\mathrm{i}$($1$)

Pass. 

 \ \color{black}

 \color{blue}

\noindent 139: (dims,levels) = $(4\oplus
1;36,
4
)$,
irreps = $4_{9,2}^{5,0}
\hskip -1.5pt \otimes \hskip -1.5pt
1_{4}^{1,0}\oplus
1_{4}^{1,0}$,
pord$(\rho_\text{isum}(\mathfrak{t})) = 9$,

\vskip 0.7ex
\hangindent=5.5em \hangafter=1
{\white .}\hskip 1em $\rho_\text{isum}(\mathfrak{t})$ =
 $( \frac{1}{4},
\frac{5}{36},
\frac{17}{36},
\frac{29}{36} )
\oplus
( \frac{1}{4} )
$,

\vskip 0.7ex
\hangindent=5.5em \hangafter=1
{\white .}\hskip 1em $\rho_\text{isum}(\mathfrak{s})$ =
$\mathrm{i}$($0$,
$\sqrt{\frac{1}{3}}$,
$\sqrt{\frac{1}{3}}$,
$\sqrt{\frac{1}{3}}$;\ \ 
$\frac{1}{3}c^{2}_{9}
$,
$\frac{1}{3}c^{1}_{9}
$,
$\frac{1}{3} c_9^4 $;\ \ 
$\frac{1}{3} c_9^4 $,
$\frac{1}{3}c^{2}_{9}
$;\ \ 
$\frac{1}{3}c^{1}_{9}
$)
 $\oplus$
$\mathrm{i}$($1$)

Pass. 

 \ \color{black}

 \color{blue}

\noindent 140: (dims,levels) = $(4\oplus
1;36,
4
)$,
irreps = $4_{9,1}^{2,0}
\hskip -1.5pt \otimes \hskip -1.5pt
1_{4}^{1,0}\oplus
1_{4}^{1,0}$,
pord$(\rho_\text{isum}(\mathfrak{t})) = 9$,

\vskip 0.7ex
\hangindent=5.5em \hangafter=1
{\white .}\hskip 1em $\rho_\text{isum}(\mathfrak{t})$ =
 $( \frac{1}{4},
\frac{5}{36},
\frac{17}{36},
\frac{29}{36} )
\oplus
( \frac{1}{4} )
$,

\vskip 0.7ex
\hangindent=5.5em \hangafter=1
{\white .}\hskip 1em $\rho_\text{isum}(\mathfrak{s})$ =
($0$,
$-\sqrt{\frac{1}{3}}$,
$-\sqrt{\frac{1}{3}}$,
$-\sqrt{\frac{1}{3}}$;
$-\frac{1}{3}c^{1}_{36}
$,
$\frac{1}{3}c^{5}_{36}
$,
$\frac{1}{3}c^{1}_{36}
-\frac{1}{3}c^{5}_{36}
$;
$\frac{1}{3}c^{1}_{36}
-\frac{1}{3}c^{5}_{36}
$,
$-\frac{1}{3}c^{1}_{36}
$;
$\frac{1}{3}c^{5}_{36}
$)
 $\oplus$
$\mathrm{i}$($1$)

Pass. 

 \ \color{black}

 \color{blue}

\noindent 141: (dims,levels) = $(4\oplus
1;42,
6
)$,
irreps = $4_{7}^{1}
\hskip -1.5pt \otimes \hskip -1.5pt
1_{3}^{1,0}
\hskip -1.5pt \otimes \hskip -1.5pt
1_{2}^{1,0}\oplus
1_{3}^{1,0}
\hskip -1.5pt \otimes \hskip -1.5pt
1_{2}^{1,0}$,
pord$(\rho_\text{isum}(\mathfrak{t})) = 7$,

\vskip 0.7ex
\hangindent=5.5em \hangafter=1
{\white .}\hskip 1em $\rho_\text{isum}(\mathfrak{t})$ =
 $( \frac{5}{6},
\frac{5}{42},
\frac{17}{42},
\frac{41}{42} )
\oplus
( \frac{5}{6} )
$,

\vskip 0.7ex
\hangindent=5.5em \hangafter=1
{\white .}\hskip 1em $\rho_\text{isum}(\mathfrak{s})$ =
$\mathrm{i}$($\sqrt{\frac{1}{7}}$,
$\sqrt{\frac{2}{7}}$,
$\sqrt{\frac{2}{7}}$,
$\sqrt{\frac{2}{7}}$;\ \ 
$-\frac{1}{\sqrt{7}\mathrm{i}}s^{5}_{28}
$,
$\frac{1}{\sqrt{7}}c^{2}_{7}
$,
$\frac{1}{\sqrt{7}}c^{1}_{7}
$;\ \ 
$\frac{1}{\sqrt{7}}c^{1}_{7}
$,
$-\frac{1}{\sqrt{7}\mathrm{i}}s^{5}_{28}
$;\ \ 
$\frac{1}{\sqrt{7}}c^{2}_{7}
$)
 $\oplus$
($-1$)

Pass. 

 \ \color{black}

 \color{blue}

\noindent 142: (dims,levels) = $(4\oplus
1;42,
6
)$,
irreps = $4_{7}^{3}
\hskip -1.5pt \otimes \hskip -1.5pt
1_{3}^{1,0}
\hskip -1.5pt \otimes \hskip -1.5pt
1_{2}^{1,0}\oplus
1_{3}^{1,0}
\hskip -1.5pt \otimes \hskip -1.5pt
1_{2}^{1,0}$,
pord$(\rho_\text{isum}(\mathfrak{t})) = 7$,

\vskip 0.7ex
\hangindent=5.5em \hangafter=1
{\white .}\hskip 1em $\rho_\text{isum}(\mathfrak{t})$ =
 $( \frac{5}{6},
\frac{11}{42},
\frac{23}{42},
\frac{29}{42} )
\oplus
( \frac{5}{6} )
$,

\vskip 0.7ex
\hangindent=5.5em \hangafter=1
{\white .}\hskip 1em $\rho_\text{isum}(\mathfrak{s})$ =
$\mathrm{i}$($-\sqrt{\frac{1}{7}}$,
$\sqrt{\frac{2}{7}}$,
$\sqrt{\frac{2}{7}}$,
$\sqrt{\frac{2}{7}}$;\ \ 
$-\frac{1}{\sqrt{7}}c^{1}_{7}
$,
$-\frac{1}{\sqrt{7}}c^{2}_{7}
$,
$\frac{1}{\sqrt{7}\mathrm{i}}s^{5}_{28}
$;\ \ 
$\frac{1}{\sqrt{7}\mathrm{i}}s^{5}_{28}
$,
$-\frac{1}{\sqrt{7}}c^{1}_{7}
$;\ \ 
$-\frac{1}{\sqrt{7}}c^{2}_{7}
$)
 $\oplus$
($-1$)

Pass. 

 \ \color{black}

 \color{blue}

\noindent 143: (dims,levels) = $(4\oplus
1;84,
12
)$,
irreps = $4_{7}^{3}
\hskip -1.5pt \otimes \hskip -1.5pt
1_{4}^{1,0}
\hskip -1.5pt \otimes \hskip -1.5pt
1_{3}^{1,0}\oplus
1_{4}^{1,0}
\hskip -1.5pt \otimes \hskip -1.5pt
1_{3}^{1,0}$,
pord$(\rho_\text{isum}(\mathfrak{t})) = 7$,

\vskip 0.7ex
\hangindent=5.5em \hangafter=1
{\white .}\hskip 1em $\rho_\text{isum}(\mathfrak{t})$ =
 $( \frac{7}{12},
\frac{1}{84},
\frac{25}{84},
\frac{37}{84} )
\oplus
( \frac{7}{12} )
$,

\vskip 0.7ex
\hangindent=5.5em \hangafter=1
{\white .}\hskip 1em $\rho_\text{isum}(\mathfrak{s})$ =
($-\sqrt{\frac{1}{7}}$,
$\sqrt{\frac{2}{7}}$,
$\sqrt{\frac{2}{7}}$,
$\sqrt{\frac{2}{7}}$;
$-\frac{1}{\sqrt{7}}c^{1}_{7}
$,
$-\frac{1}{\sqrt{7}}c^{2}_{7}
$,
$\frac{1}{\sqrt{7}\mathrm{i}}s^{5}_{28}
$;
$\frac{1}{\sqrt{7}\mathrm{i}}s^{5}_{28}
$,
$-\frac{1}{\sqrt{7}}c^{1}_{7}
$;
$-\frac{1}{\sqrt{7}}c^{2}_{7}
$)
 $\oplus$
$\mathrm{i}$($1$)

Pass. 

 \ \color{black}

 \color{blue}

\noindent 144: (dims,levels) = $(4\oplus
1;84,
12
)$,
irreps = $4_{7}^{1}
\hskip -1.5pt \otimes \hskip -1.5pt
1_{4}^{1,0}
\hskip -1.5pt \otimes \hskip -1.5pt
1_{3}^{1,0}\oplus
1_{4}^{1,0}
\hskip -1.5pt \otimes \hskip -1.5pt
1_{3}^{1,0}$,
pord$(\rho_\text{isum}(\mathfrak{t})) = 7$,

\vskip 0.7ex
\hangindent=5.5em \hangafter=1
{\white .}\hskip 1em $\rho_\text{isum}(\mathfrak{t})$ =
 $( \frac{7}{12},
\frac{13}{84},
\frac{61}{84},
\frac{73}{84} )
\oplus
( \frac{7}{12} )
$,

\vskip 0.7ex
\hangindent=5.5em \hangafter=1
{\white .}\hskip 1em $\rho_\text{isum}(\mathfrak{s})$ =
($\sqrt{\frac{1}{7}}$,
$\sqrt{\frac{2}{7}}$,
$\sqrt{\frac{2}{7}}$,
$\sqrt{\frac{2}{7}}$;
$\frac{1}{\sqrt{7}}c^{1}_{7}
$,
$-\frac{1}{\sqrt{7}\mathrm{i}}s^{5}_{28}
$,
$\frac{1}{\sqrt{7}}c^{2}_{7}
$;
$\frac{1}{\sqrt{7}}c^{2}_{7}
$,
$\frac{1}{\sqrt{7}}c^{1}_{7}
$;
$-\frac{1}{\sqrt{7}\mathrm{i}}s^{5}_{28}
$)
 $\oplus$
$\mathrm{i}$($1$)

Pass. 

 \ \color{black}

\noindent 145: (dims,levels) = $(5;5
)$,
irreps = $5_{5}^{1}$,
pord$(\rho_\text{isum}(\mathfrak{t})) = 5$,

\vskip 0.7ex
\hangindent=5.5em \hangafter=1
{\white .}\hskip 1em $\rho_\text{isum}(\mathfrak{t})$ =
 $( 0,
\frac{1}{5},
\frac{2}{5},
\frac{3}{5},
\frac{4}{5} )
$,

\vskip 0.7ex
\hangindent=5.5em \hangafter=1
{\white .}\hskip 1em $\rho_\text{isum}(\mathfrak{s})$ =
($-\frac{1}{5}$,
$\sqrt{\frac{6}{25}}$,
$\sqrt{\frac{6}{25}}$,
$\sqrt{\frac{6}{25}}$,
$\sqrt{\frac{6}{25}}$;
$\frac{3-\sqrt{5}}{10}$,
$-\frac{1+\sqrt{5}}{5}$,
$\frac{-1+\sqrt{5}}{5}$,
$\frac{3+\sqrt{5}}{10}$;
$\frac{3+\sqrt{5}}{10}$,
$\frac{3-\sqrt{5}}{10}$,
$\frac{-1+\sqrt{5}}{5}$;
$\frac{3+\sqrt{5}}{10}$,
$-\frac{1+\sqrt{5}}{5}$;
$\frac{3-\sqrt{5}}{10}$)

Fail:
cnd($\rho(\mathfrak s)_\mathrm{ndeg}$) = 120 does not divide
 ord($\rho(\mathfrak t)$)=5. Prop. B.4 (2)

 \ \color{black}

\noindent 146: (dims,levels) = $(5;10
)$,
irreps = $5_{5}^{1}
\hskip -1.5pt \otimes \hskip -1.5pt
1_{2}^{1,0}$,
pord$(\rho_\text{isum}(\mathfrak{t})) = 5$,

\vskip 0.7ex
\hangindent=5.5em \hangafter=1
{\white .}\hskip 1em $\rho_\text{isum}(\mathfrak{t})$ =
 $( \frac{1}{2},
\frac{1}{10},
\frac{3}{10},
\frac{7}{10},
\frac{9}{10} )
$,

\vskip 0.7ex
\hangindent=5.5em \hangafter=1
{\white .}\hskip 1em $\rho_\text{isum}(\mathfrak{s})$ =
($\frac{1}{5}$,
$\sqrt{\frac{6}{25}}$,
$\sqrt{\frac{6}{25}}$,
$\sqrt{\frac{6}{25}}$,
$\sqrt{\frac{6}{25}}$;
$-\frac{3+\sqrt{5}}{10}$,
$\frac{1+\sqrt{5}}{5}$,
$\frac{1-\sqrt{5}}{5}$,
$\frac{-3+\sqrt{5}}{10}$;
$\frac{-3+\sqrt{5}}{10}$,
$-\frac{3+\sqrt{5}}{10}$,
$\frac{1-\sqrt{5}}{5}$;
$\frac{-3+\sqrt{5}}{10}$,
$\frac{1+\sqrt{5}}{5}$;
$-\frac{3+\sqrt{5}}{10}$)

Fail:
cnd($\rho(\mathfrak s)_\mathrm{ndeg}$) = 120 does not divide
 ord($\rho(\mathfrak t)$)=10. Prop. B.4 (2)

 \ \color{black}

 \color{blue}

\noindent 147: (dims,levels) = $(5;11
)$,
irreps = $5_{11}^{1}$,
pord$(\rho_\text{isum}(\mathfrak{t})) = 11$,

\vskip 0.7ex
\hangindent=5.5em \hangafter=1
{\white .}\hskip 1em $\rho_\text{isum}(\mathfrak{t})$ =
 $( \frac{1}{11},
\frac{3}{11},
\frac{4}{11},
\frac{5}{11},
\frac{9}{11} )
$,

\vskip 0.7ex
\hangindent=5.5em \hangafter=1
{\white .}\hskip 1em $\rho_\text{isum}(\mathfrak{s})$ =
($-\frac{1}{\sqrt{11}}c^{3}_{44}
$,
$-\frac{1}{\sqrt{11}}c^{7}_{44}
$,
$-\frac{1}{\sqrt{11}}c^{5}_{44}
$,
$-\frac{1}{\sqrt{11}}c^{1}_{44}
$,
$-\frac{1}{\sqrt{11}}c^{9}_{44}
$;
$\frac{1}{\sqrt{11}}c^{9}_{44}
$,
$-\frac{1}{\sqrt{11}}c^{3}_{44}
$,
$\frac{1}{\sqrt{11}}c^{5}_{44}
$,
$\frac{1}{\sqrt{11}}c^{1}_{44}
$;
$\frac{1}{\sqrt{11}}c^{1}_{44}
$,
$\frac{1}{\sqrt{11}}c^{9}_{44}
$,
$\frac{1}{\sqrt{11}}c^{7}_{44}
$;
$\frac{1}{\sqrt{11}}c^{7}_{44}
$,
$-\frac{1}{\sqrt{11}}c^{3}_{44}
$;
$\frac{1}{\sqrt{11}}c^{5}_{44}
$)

Pass. 

 \ \color{black}

 \color{blue}

\noindent 148: (dims,levels) = $(5;11
)$,
irreps = $5_{11}^{2}$,
pord$(\rho_\text{isum}(\mathfrak{t})) = 11$,

\vskip 0.7ex
\hangindent=5.5em \hangafter=1
{\white .}\hskip 1em $\rho_\text{isum}(\mathfrak{t})$ =
 $( \frac{2}{11},
\frac{6}{11},
\frac{7}{11},
\frac{8}{11},
\frac{10}{11} )
$,

\vskip 0.7ex
\hangindent=5.5em \hangafter=1
{\white .}\hskip 1em $\rho_\text{isum}(\mathfrak{s})$ =
($\frac{1}{\sqrt{11}}c^{5}_{44}
$,
$\frac{1}{\sqrt{11}}c^{3}_{44}
$,
$-\frac{1}{\sqrt{11}}c^{7}_{44}
$,
$-\frac{1}{\sqrt{11}}c^{1}_{44}
$,
$-\frac{1}{\sqrt{11}}c^{9}_{44}
$;
$\frac{1}{\sqrt{11}}c^{7}_{44}
$,
$\frac{1}{\sqrt{11}}c^{9}_{44}
$,
$\frac{1}{\sqrt{11}}c^{5}_{44}
$,
$\frac{1}{\sqrt{11}}c^{1}_{44}
$;
$\frac{1}{\sqrt{11}}c^{1}_{44}
$,
$-\frac{1}{\sqrt{11}}c^{3}_{44}
$,
$\frac{1}{\sqrt{11}}c^{5}_{44}
$;
$\frac{1}{\sqrt{11}}c^{9}_{44}
$,
$\frac{1}{\sqrt{11}}c^{7}_{44}
$;
$-\frac{1}{\sqrt{11}}c^{3}_{44}
$)

Pass. 

 \ \color{black}

\noindent 149: (dims,levels) = $(5;15
)$,
irreps = $5_{5}^{1}
\hskip -1.5pt \otimes \hskip -1.5pt
1_{3}^{1,0}$,
pord$(\rho_\text{isum}(\mathfrak{t})) = 5$,

\vskip 0.7ex
\hangindent=5.5em \hangafter=1
{\white .}\hskip 1em $\rho_\text{isum}(\mathfrak{t})$ =
 $( \frac{1}{3},
\frac{2}{15},
\frac{8}{15},
\frac{11}{15},
\frac{14}{15} )
$,

\vskip 0.7ex
\hangindent=5.5em \hangafter=1
{\white .}\hskip 1em $\rho_\text{isum}(\mathfrak{s})$ =
($-\frac{1}{5}$,
$\sqrt{\frac{6}{25}}$,
$\sqrt{\frac{6}{25}}$,
$\sqrt{\frac{6}{25}}$,
$\sqrt{\frac{6}{25}}$;
$\frac{3-\sqrt{5}}{10}$,
$\frac{3+\sqrt{5}}{10}$,
$\frac{-1+\sqrt{5}}{5}$,
$-\frac{1+\sqrt{5}}{5}$;
$\frac{3-\sqrt{5}}{10}$,
$-\frac{1+\sqrt{5}}{5}$,
$\frac{-1+\sqrt{5}}{5}$;
$\frac{3+\sqrt{5}}{10}$,
$\frac{3-\sqrt{5}}{10}$;
$\frac{3+\sqrt{5}}{10}$)

Fail:
cnd($\rho(\mathfrak s)_\mathrm{ndeg}$) = 120 does not divide
 ord($\rho(\mathfrak t)$)=15. Prop. B.4 (2)

 \ \color{black}

\noindent 150: (dims,levels) = $(5;20
)$,
irreps = $5_{5}^{1}
\hskip -1.5pt \otimes \hskip -1.5pt
1_{4}^{1,0}$,
pord$(\rho_\text{isum}(\mathfrak{t})) = 5$,

\vskip 0.7ex
\hangindent=5.5em \hangafter=1
{\white .}\hskip 1em $\rho_\text{isum}(\mathfrak{t})$ =
 $( \frac{1}{4},
\frac{1}{20},
\frac{9}{20},
\frac{13}{20},
\frac{17}{20} )
$,

\vskip 0.7ex
\hangindent=5.5em \hangafter=1
{\white .}\hskip 1em $\rho_\text{isum}(\mathfrak{s})$ =
$\mathrm{i}$($-\frac{1}{5}$,
$\sqrt{\frac{6}{25}}$,
$\sqrt{\frac{6}{25}}$,
$\sqrt{\frac{6}{25}}$,
$\sqrt{\frac{6}{25}}$;\ \ 
$\frac{3-\sqrt{5}}{10}$,
$\frac{3+\sqrt{5}}{10}$,
$\frac{-1+\sqrt{5}}{5}$,
$-\frac{1+\sqrt{5}}{5}$;\ \ 
$\frac{3-\sqrt{5}}{10}$,
$-\frac{1+\sqrt{5}}{5}$,
$\frac{-1+\sqrt{5}}{5}$;\ \ 
$\frac{3+\sqrt{5}}{10}$,
$\frac{3-\sqrt{5}}{10}$;\ \ 
$\frac{3+\sqrt{5}}{10}$)

Fail:
cnd($\rho(\mathfrak s)_\mathrm{ndeg}$) = 120 does not divide
 ord($\rho(\mathfrak t)$)=20. Prop. B.4 (2)

 \ \color{black}

 \color{blue}

\noindent 151: (dims,levels) = $(5;22
)$,
irreps = $5_{11}^{2}
\hskip -1.5pt \otimes \hskip -1.5pt
1_{2}^{1,0}$,
pord$(\rho_\text{isum}(\mathfrak{t})) = 11$,

\vskip 0.7ex
\hangindent=5.5em \hangafter=1
{\white .}\hskip 1em $\rho_\text{isum}(\mathfrak{t})$ =
 $( \frac{1}{22},
\frac{3}{22},
\frac{5}{22},
\frac{9}{22},
\frac{15}{22} )
$,

\vskip 0.7ex
\hangindent=5.5em \hangafter=1
{\white .}\hskip 1em $\rho_\text{isum}(\mathfrak{s})$ =
($-\frac{1}{\sqrt{11}}c^{7}_{44}
$,
$-\frac{1}{\sqrt{11}}c^{9}_{44}
$,
$-\frac{1}{\sqrt{11}}c^{5}_{44}
$,
$-\frac{1}{\sqrt{11}}c^{1}_{44}
$,
$\frac{1}{\sqrt{11}}c^{3}_{44}
$;
$-\frac{1}{\sqrt{11}}c^{1}_{44}
$,
$\frac{1}{\sqrt{11}}c^{3}_{44}
$,
$-\frac{1}{\sqrt{11}}c^{5}_{44}
$,
$-\frac{1}{\sqrt{11}}c^{7}_{44}
$;
$-\frac{1}{\sqrt{11}}c^{9}_{44}
$,
$-\frac{1}{\sqrt{11}}c^{7}_{44}
$,
$-\frac{1}{\sqrt{11}}c^{1}_{44}
$;
$\frac{1}{\sqrt{11}}c^{3}_{44}
$,
$-\frac{1}{\sqrt{11}}c^{9}_{44}
$;
$-\frac{1}{\sqrt{11}}c^{5}_{44}
$)

Pass. 

 \ \color{black}

 \color{blue}

\noindent 152: (dims,levels) = $(5;22
)$,
irreps = $5_{11}^{1}
\hskip -1.5pt \otimes \hskip -1.5pt
1_{2}^{1,0}$,
pord$(\rho_\text{isum}(\mathfrak{t})) = 11$,

\vskip 0.7ex
\hangindent=5.5em \hangafter=1
{\white .}\hskip 1em $\rho_\text{isum}(\mathfrak{t})$ =
 $( \frac{7}{22},
\frac{13}{22},
\frac{17}{22},
\frac{19}{22},
\frac{21}{22} )
$,

\vskip 0.7ex
\hangindent=5.5em \hangafter=1
{\white .}\hskip 1em $\rho_\text{isum}(\mathfrak{s})$ =
($-\frac{1}{\sqrt{11}}c^{5}_{44}
$,
$-\frac{1}{\sqrt{11}}c^{9}_{44}
$,
$-\frac{1}{\sqrt{11}}c^{1}_{44}
$,
$-\frac{1}{\sqrt{11}}c^{7}_{44}
$,
$\frac{1}{\sqrt{11}}c^{3}_{44}
$;
$\frac{1}{\sqrt{11}}c^{3}_{44}
$,
$-\frac{1}{\sqrt{11}}c^{7}_{44}
$,
$-\frac{1}{\sqrt{11}}c^{5}_{44}
$,
$-\frac{1}{\sqrt{11}}c^{1}_{44}
$;
$-\frac{1}{\sqrt{11}}c^{9}_{44}
$,
$\frac{1}{\sqrt{11}}c^{3}_{44}
$,
$-\frac{1}{\sqrt{11}}c^{5}_{44}
$;
$-\frac{1}{\sqrt{11}}c^{1}_{44}
$,
$-\frac{1}{\sqrt{11}}c^{9}_{44}
$;
$-\frac{1}{\sqrt{11}}c^{7}_{44}
$)

Pass. 

 \ \color{black}

\noindent 153: (dims,levels) = $(5;30
)$,
irreps = $5_{5}^{1}
\hskip -1.5pt \otimes \hskip -1.5pt
1_{3}^{1,0}
\hskip -1.5pt \otimes \hskip -1.5pt
1_{2}^{1,0}$,
pord$(\rho_\text{isum}(\mathfrak{t})) = 5$,

\vskip 0.7ex
\hangindent=5.5em \hangafter=1
{\white .}\hskip 1em $\rho_\text{isum}(\mathfrak{t})$ =
 $( \frac{5}{6},
\frac{1}{30},
\frac{7}{30},
\frac{13}{30},
\frac{19}{30} )
$,

\vskip 0.7ex
\hangindent=5.5em \hangafter=1
{\white .}\hskip 1em $\rho_\text{isum}(\mathfrak{s})$ =
($\frac{1}{5}$,
$\sqrt{\frac{6}{25}}$,
$\sqrt{\frac{6}{25}}$,
$\sqrt{\frac{6}{25}}$,
$\sqrt{\frac{6}{25}}$;
$\frac{-3+\sqrt{5}}{10}$,
$\frac{1+\sqrt{5}}{5}$,
$\frac{1-\sqrt{5}}{5}$,
$-\frac{3+\sqrt{5}}{10}$;
$-\frac{3+\sqrt{5}}{10}$,
$\frac{-3+\sqrt{5}}{10}$,
$\frac{1-\sqrt{5}}{5}$;
$-\frac{3+\sqrt{5}}{10}$,
$\frac{1+\sqrt{5}}{5}$;
$\frac{-3+\sqrt{5}}{10}$)

Fail:
cnd($\rho(\mathfrak s)_\mathrm{ndeg}$) = 120 does not divide
 ord($\rho(\mathfrak t)$)=30. Prop. B.4 (2)

 \ \color{black}

 \color{blue}

\noindent 154: (dims,levels) = $(5;33
)$,
irreps = $5_{11}^{2}
\hskip -1.5pt \otimes \hskip -1.5pt
1_{3}^{1,0}$,
pord$(\rho_\text{isum}(\mathfrak{t})) = 11$,

\vskip 0.7ex
\hangindent=5.5em \hangafter=1
{\white .}\hskip 1em $\rho_\text{isum}(\mathfrak{t})$ =
 $( \frac{2}{33},
\frac{8}{33},
\frac{17}{33},
\frac{29}{33},
\frac{32}{33} )
$,

\vskip 0.7ex
\hangindent=5.5em \hangafter=1
{\white .}\hskip 1em $\rho_\text{isum}(\mathfrak{s})$ =
($\frac{1}{\sqrt{11}}c^{9}_{44}
$,
$-\frac{1}{\sqrt{11}}c^{7}_{44}
$,
$-\frac{1}{\sqrt{11}}c^{1}_{44}
$,
$-\frac{1}{\sqrt{11}}c^{5}_{44}
$,
$\frac{1}{\sqrt{11}}c^{3}_{44}
$;
$-\frac{1}{\sqrt{11}}c^{3}_{44}
$,
$\frac{1}{\sqrt{11}}c^{9}_{44}
$,
$\frac{1}{\sqrt{11}}c^{1}_{44}
$,
$\frac{1}{\sqrt{11}}c^{5}_{44}
$;
$\frac{1}{\sqrt{11}}c^{5}_{44}
$,
$-\frac{1}{\sqrt{11}}c^{3}_{44}
$,
$\frac{1}{\sqrt{11}}c^{7}_{44}
$;
$\frac{1}{\sqrt{11}}c^{7}_{44}
$,
$\frac{1}{\sqrt{11}}c^{9}_{44}
$;
$\frac{1}{\sqrt{11}}c^{1}_{44}
$)

Pass. 

 \ \color{black}

 \color{blue}

\noindent 155: (dims,levels) = $(5;33
)$,
irreps = $5_{11}^{1}
\hskip -1.5pt \otimes \hskip -1.5pt
1_{3}^{1,0}$,
pord$(\rho_\text{isum}(\mathfrak{t})) = 11$,

\vskip 0.7ex
\hangindent=5.5em \hangafter=1
{\white .}\hskip 1em $\rho_\text{isum}(\mathfrak{t})$ =
 $( \frac{5}{33},
\frac{14}{33},
\frac{20}{33},
\frac{23}{33},
\frac{26}{33} )
$,

\vskip 0.7ex
\hangindent=5.5em \hangafter=1
{\white .}\hskip 1em $\rho_\text{isum}(\mathfrak{s})$ =
($\frac{1}{\sqrt{11}}c^{5}_{44}
$,
$-\frac{1}{\sqrt{11}}c^{9}_{44}
$,
$-\frac{1}{\sqrt{11}}c^{1}_{44}
$,
$-\frac{1}{\sqrt{11}}c^{7}_{44}
$,
$\frac{1}{\sqrt{11}}c^{3}_{44}
$;
$-\frac{1}{\sqrt{11}}c^{3}_{44}
$,
$\frac{1}{\sqrt{11}}c^{7}_{44}
$,
$\frac{1}{\sqrt{11}}c^{5}_{44}
$,
$\frac{1}{\sqrt{11}}c^{1}_{44}
$;
$\frac{1}{\sqrt{11}}c^{9}_{44}
$,
$-\frac{1}{\sqrt{11}}c^{3}_{44}
$,
$\frac{1}{\sqrt{11}}c^{5}_{44}
$;
$\frac{1}{\sqrt{11}}c^{1}_{44}
$,
$\frac{1}{\sqrt{11}}c^{9}_{44}
$;
$\frac{1}{\sqrt{11}}c^{7}_{44}
$)

Pass. 

 \ \color{black}

 \color{blue}

\noindent 156: (dims,levels) = $(5;44
)$,
irreps = $5_{11}^{1}
\hskip -1.5pt \otimes \hskip -1.5pt
1_{4}^{1,0}$,
pord$(\rho_\text{isum}(\mathfrak{t})) = 11$,

\vskip 0.7ex
\hangindent=5.5em \hangafter=1
{\white .}\hskip 1em $\rho_\text{isum}(\mathfrak{t})$ =
 $( \frac{3}{44},
\frac{15}{44},
\frac{23}{44},
\frac{27}{44},
\frac{31}{44} )
$,

\vskip 0.7ex
\hangindent=5.5em \hangafter=1
{\white .}\hskip 1em $\rho_\text{isum}(\mathfrak{s})$ =
$\mathrm{i}$($\frac{1}{\sqrt{11}}c^{5}_{44}
$,
$-\frac{1}{\sqrt{11}}c^{9}_{44}
$,
$-\frac{1}{\sqrt{11}}c^{1}_{44}
$,
$-\frac{1}{\sqrt{11}}c^{7}_{44}
$,
$\frac{1}{\sqrt{11}}c^{3}_{44}
$;\ \ 
$-\frac{1}{\sqrt{11}}c^{3}_{44}
$,
$\frac{1}{\sqrt{11}}c^{7}_{44}
$,
$\frac{1}{\sqrt{11}}c^{5}_{44}
$,
$\frac{1}{\sqrt{11}}c^{1}_{44}
$;\ \ 
$\frac{1}{\sqrt{11}}c^{9}_{44}
$,
$-\frac{1}{\sqrt{11}}c^{3}_{44}
$,
$\frac{1}{\sqrt{11}}c^{5}_{44}
$;\ \ 
$\frac{1}{\sqrt{11}}c^{1}_{44}
$,
$\frac{1}{\sqrt{11}}c^{9}_{44}
$;\ \ 
$\frac{1}{\sqrt{11}}c^{7}_{44}
$)

Pass. 

 \ \color{black}

 \color{blue}

\noindent 157: (dims,levels) = $(5;44
)$,
irreps = $5_{11}^{2}
\hskip -1.5pt \otimes \hskip -1.5pt
1_{4}^{1,0}$,
pord$(\rho_\text{isum}(\mathfrak{t})) = 11$,

\vskip 0.7ex
\hangindent=5.5em \hangafter=1
{\white .}\hskip 1em $\rho_\text{isum}(\mathfrak{t})$ =
 $( \frac{7}{44},
\frac{19}{44},
\frac{35}{44},
\frac{39}{44},
\frac{43}{44} )
$,

\vskip 0.7ex
\hangindent=5.5em \hangafter=1
{\white .}\hskip 1em $\rho_\text{isum}(\mathfrak{s})$ =
$\mathrm{i}$($-\frac{1}{\sqrt{11}}c^{3}_{44}
$,
$-\frac{1}{\sqrt{11}}c^{9}_{44}
$,
$-\frac{1}{\sqrt{11}}c^{1}_{44}
$,
$-\frac{1}{\sqrt{11}}c^{5}_{44}
$,
$-\frac{1}{\sqrt{11}}c^{7}_{44}
$;\ \ 
$\frac{1}{\sqrt{11}}c^{5}_{44}
$,
$-\frac{1}{\sqrt{11}}c^{3}_{44}
$,
$\frac{1}{\sqrt{11}}c^{7}_{44}
$,
$\frac{1}{\sqrt{11}}c^{1}_{44}
$;\ \ 
$\frac{1}{\sqrt{11}}c^{7}_{44}
$,
$\frac{1}{\sqrt{11}}c^{9}_{44}
$,
$\frac{1}{\sqrt{11}}c^{5}_{44}
$;\ \ 
$\frac{1}{\sqrt{11}}c^{1}_{44}
$,
$-\frac{1}{\sqrt{11}}c^{3}_{44}
$;\ \ 
$\frac{1}{\sqrt{11}}c^{9}_{44}
$)

Pass. 

 \ \color{black}

\noindent 158: (dims,levels) = $(5;60
)$,
irreps = $5_{5}^{1}
\hskip -1.5pt \otimes \hskip -1.5pt
1_{4}^{1,0}
\hskip -1.5pt \otimes \hskip -1.5pt
1_{3}^{1,0}$,
pord$(\rho_\text{isum}(\mathfrak{t})) = 5$,

\vskip 0.7ex
\hangindent=5.5em \hangafter=1
{\white .}\hskip 1em $\rho_\text{isum}(\mathfrak{t})$ =
 $( \frac{7}{12},
\frac{11}{60},
\frac{23}{60},
\frac{47}{60},
\frac{59}{60} )
$,

\vskip 0.7ex
\hangindent=5.5em \hangafter=1
{\white .}\hskip 1em $\rho_\text{isum}(\mathfrak{s})$ =
$\mathrm{i}$($-\frac{1}{5}$,
$\sqrt{\frac{6}{25}}$,
$\sqrt{\frac{6}{25}}$,
$\sqrt{\frac{6}{25}}$,
$\sqrt{\frac{6}{25}}$;\ \ 
$\frac{3+\sqrt{5}}{10}$,
$-\frac{1+\sqrt{5}}{5}$,
$\frac{-1+\sqrt{5}}{5}$,
$\frac{3-\sqrt{5}}{10}$;\ \ 
$\frac{3-\sqrt{5}}{10}$,
$\frac{3+\sqrt{5}}{10}$,
$\frac{-1+\sqrt{5}}{5}$;\ \ 
$\frac{3-\sqrt{5}}{10}$,
$-\frac{1+\sqrt{5}}{5}$;\ \ 
$\frac{3+\sqrt{5}}{10}$)

Fail:
cnd($\rho(\mathfrak s)_\mathrm{ndeg}$) = 120 does not divide
 ord($\rho(\mathfrak t)$)=60. Prop. B.4 (2)

 \ \color{black}

 \color{blue}

\noindent 159: (dims,levels) = $(5;66
)$,
irreps = $5_{11}^{2}
\hskip -1.5pt \otimes \hskip -1.5pt
1_{3}^{1,0}
\hskip -1.5pt \otimes \hskip -1.5pt
1_{2}^{1,0}$,
pord$(\rho_\text{isum}(\mathfrak{t})) = 11$,

\vskip 0.7ex
\hangindent=5.5em \hangafter=1
{\white .}\hskip 1em $\rho_\text{isum}(\mathfrak{t})$ =
 $( \frac{1}{66},
\frac{25}{66},
\frac{31}{66},
\frac{37}{66},
\frac{49}{66} )
$,

\vskip 0.7ex
\hangindent=5.5em \hangafter=1
{\white .}\hskip 1em $\rho_\text{isum}(\mathfrak{s})$ =
($-\frac{1}{\sqrt{11}}c^{5}_{44}
$,
$\frac{1}{\sqrt{11}}c^{3}_{44}
$,
$-\frac{1}{\sqrt{11}}c^{7}_{44}
$,
$-\frac{1}{\sqrt{11}}c^{1}_{44}
$,
$-\frac{1}{\sqrt{11}}c^{9}_{44}
$;
$-\frac{1}{\sqrt{11}}c^{7}_{44}
$,
$-\frac{1}{\sqrt{11}}c^{9}_{44}
$,
$-\frac{1}{\sqrt{11}}c^{5}_{44}
$,
$-\frac{1}{\sqrt{11}}c^{1}_{44}
$;
$-\frac{1}{\sqrt{11}}c^{1}_{44}
$,
$\frac{1}{\sqrt{11}}c^{3}_{44}
$,
$-\frac{1}{\sqrt{11}}c^{5}_{44}
$;
$-\frac{1}{\sqrt{11}}c^{9}_{44}
$,
$-\frac{1}{\sqrt{11}}c^{7}_{44}
$;
$\frac{1}{\sqrt{11}}c^{3}_{44}
$)

Pass. 

 \ \color{black}

 \color{blue}

\noindent 160: (dims,levels) = $(5;66
)$,
irreps = $5_{11}^{1}
\hskip -1.5pt \otimes \hskip -1.5pt
1_{3}^{1,0}
\hskip -1.5pt \otimes \hskip -1.5pt
1_{2}^{1,0}$,
pord$(\rho_\text{isum}(\mathfrak{t})) = 11$,

\vskip 0.7ex
\hangindent=5.5em \hangafter=1
{\white .}\hskip 1em $\rho_\text{isum}(\mathfrak{t})$ =
 $( \frac{7}{66},
\frac{13}{66},
\frac{19}{66},
\frac{43}{66},
\frac{61}{66} )
$,

\vskip 0.7ex
\hangindent=5.5em \hangafter=1
{\white .}\hskip 1em $\rho_\text{isum}(\mathfrak{s})$ =
($-\frac{1}{\sqrt{11}}c^{9}_{44}
$,
$\frac{1}{\sqrt{11}}c^{3}_{44}
$,
$-\frac{1}{\sqrt{11}}c^{5}_{44}
$,
$-\frac{1}{\sqrt{11}}c^{1}_{44}
$,
$-\frac{1}{\sqrt{11}}c^{7}_{44}
$;
$-\frac{1}{\sqrt{11}}c^{1}_{44}
$,
$-\frac{1}{\sqrt{11}}c^{9}_{44}
$,
$-\frac{1}{\sqrt{11}}c^{7}_{44}
$,
$-\frac{1}{\sqrt{11}}c^{5}_{44}
$;
$-\frac{1}{\sqrt{11}}c^{7}_{44}
$,
$\frac{1}{\sqrt{11}}c^{3}_{44}
$,
$-\frac{1}{\sqrt{11}}c^{1}_{44}
$;
$-\frac{1}{\sqrt{11}}c^{5}_{44}
$,
$-\frac{1}{\sqrt{11}}c^{9}_{44}
$;
$\frac{1}{\sqrt{11}}c^{3}_{44}
$)

Pass. 

 \ \color{black}

 \color{blue}

\noindent 161: (dims,levels) = $(5;132
)$,
irreps = $5_{11}^{1}
\hskip -1.5pt \otimes \hskip -1.5pt
1_{4}^{1,0}
\hskip -1.5pt \otimes \hskip -1.5pt
1_{3}^{1,0}$,
pord$(\rho_\text{isum}(\mathfrak{t})) = 11$,

\vskip 0.7ex
\hangindent=5.5em \hangafter=1
{\white .}\hskip 1em $\rho_\text{isum}(\mathfrak{t})$ =
 $( \frac{5}{132},
\frac{53}{132},
\frac{89}{132},
\frac{113}{132},
\frac{125}{132} )
$,

\vskip 0.7ex
\hangindent=5.5em \hangafter=1
{\white .}\hskip 1em $\rho_\text{isum}(\mathfrak{s})$ =
$\mathrm{i}$($\frac{1}{\sqrt{11}}c^{7}_{44}
$,
$\frac{1}{\sqrt{11}}c^{3}_{44}
$,
$-\frac{1}{\sqrt{11}}c^{1}_{44}
$,
$-\frac{1}{\sqrt{11}}c^{5}_{44}
$,
$-\frac{1}{\sqrt{11}}c^{9}_{44}
$;\ \ 
$\frac{1}{\sqrt{11}}c^{5}_{44}
$,
$\frac{1}{\sqrt{11}}c^{9}_{44}
$,
$\frac{1}{\sqrt{11}}c^{1}_{44}
$,
$\frac{1}{\sqrt{11}}c^{7}_{44}
$;\ \ 
$-\frac{1}{\sqrt{11}}c^{3}_{44}
$,
$\frac{1}{\sqrt{11}}c^{7}_{44}
$,
$\frac{1}{\sqrt{11}}c^{5}_{44}
$;\ \ 
$\frac{1}{\sqrt{11}}c^{9}_{44}
$,
$-\frac{1}{\sqrt{11}}c^{3}_{44}
$;\ \ 
$\frac{1}{\sqrt{11}}c^{1}_{44}
$)

Pass. 

 \ \color{black}

 \color{blue}

\noindent 162: (dims,levels) = $(5;132
)$,
irreps = $5_{11}^{2}
\hskip -1.5pt \otimes \hskip -1.5pt
1_{4}^{1,0}
\hskip -1.5pt \otimes \hskip -1.5pt
1_{3}^{1,0}$,
pord$(\rho_\text{isum}(\mathfrak{t})) = 11$,

\vskip 0.7ex
\hangindent=5.5em \hangafter=1
{\white .}\hskip 1em $\rho_\text{isum}(\mathfrak{t})$ =
 $( \frac{17}{132},
\frac{29}{132},
\frac{41}{132},
\frac{65}{132},
\frac{101}{132} )
$,

\vskip 0.7ex
\hangindent=5.5em \hangafter=1
{\white .}\hskip 1em $\rho_\text{isum}(\mathfrak{s})$ =
$\mathrm{i}$($\frac{1}{\sqrt{11}}c^{7}_{44}
$,
$-\frac{1}{\sqrt{11}}c^{9}_{44}
$,
$-\frac{1}{\sqrt{11}}c^{5}_{44}
$,
$-\frac{1}{\sqrt{11}}c^{1}_{44}
$,
$\frac{1}{\sqrt{11}}c^{3}_{44}
$;\ \ 
$\frac{1}{\sqrt{11}}c^{1}_{44}
$,
$-\frac{1}{\sqrt{11}}c^{3}_{44}
$,
$\frac{1}{\sqrt{11}}c^{5}_{44}
$,
$\frac{1}{\sqrt{11}}c^{7}_{44}
$;\ \ 
$\frac{1}{\sqrt{11}}c^{9}_{44}
$,
$\frac{1}{\sqrt{11}}c^{7}_{44}
$,
$\frac{1}{\sqrt{11}}c^{1}_{44}
$;\ \ 
$-\frac{1}{\sqrt{11}}c^{3}_{44}
$,
$\frac{1}{\sqrt{11}}c^{9}_{44}
$;\ \ 
$\frac{1}{\sqrt{11}}c^{5}_{44}
$)

Pass. 

 \ \color{black}

\

\subsection{A list of passing GT orbits}

The above passing representations can be grouped into GT orbits.  The following
list displays one representative representation for each GT orbit.  For details
and notations, see Appendix B.2.

\

\noindent1. (dims;levels) =$(3\oplus
2;5,
2
)$,
irreps = $3_{5}^{1}\oplus
2_{2}^{1,0}$,
pord$(\rho_\text{isum}(\mathfrak{t})) = 10$,

\vskip 0.7ex
\hangindent=4em \hangafter=1
 $\rho_\text{isum}(\mathfrak{t})$ =
 $( 0,
\frac{1}{5},
\frac{4}{5} )
\oplus
( 0,
\frac{1}{2} )
$,

\vskip 0.7ex
\hangindent=4em \hangafter=1
 $\rho_\text{isum}(\mathfrak{s})$ =
($\sqrt{\frac{1}{5}}$,
$-\sqrt{\frac{2}{5}}$,
$-\sqrt{\frac{2}{5}}$;
$-\frac{5+\sqrt{5}}{10}$,
$\frac{5-\sqrt{5}}{10}$;
$-\frac{5+\sqrt{5}}{10}$)
 $\oplus$
($-\frac{1}{2}$,
$-\sqrt{\frac{3}{4}}$;
$\frac{1}{2}$)

Resolved. Number of valid $(S,T)$ pairs = 0.

\vskip 2ex

 \noindent2. (dims;levels) =$(3\oplus
2;8,
3
)$,
irreps = $3_{8}^{1,0}\oplus
2_{3}^{1,0}$,
pord$(\rho_\text{isum}(\mathfrak{t})) = 24$,

\vskip 0.7ex
\hangindent=4em \hangafter=1
 $\rho_\text{isum}(\mathfrak{t})$ =
 $( 0,
\frac{1}{8},
\frac{5}{8} )
\oplus
( 0,
\frac{1}{3} )
$,

\vskip 0.7ex
\hangindent=4em \hangafter=1
 $\rho_\text{isum}(\mathfrak{s})$ =
$\mathrm{i}$($0$,
$\sqrt{\frac{1}{2}}$,
$\sqrt{\frac{1}{2}}$;\ \ 
$-\frac{1}{2}$,
$\frac{1}{2}$;\ \ 
$-\frac{1}{2}$)
 $\oplus$
$\mathrm{i}$($-\sqrt{\frac{1}{3}}$,
$\sqrt{\frac{2}{3}}$;\ \ 
$\sqrt{\frac{1}{3}}$)

Resolved. Number of valid $(S,T)$ pairs = 2.

\vskip 2ex

 \noindent3. (dims;levels) =$(4\oplus
1;7,
1
)$,
irreps = $4_{7}^{1}\oplus
1_{1}^{1}$,
pord$(\rho_\text{isum}(\mathfrak{t})) = 7$,

\vskip 0.7ex
\hangindent=4em \hangafter=1
 $\rho_\text{isum}(\mathfrak{t})$ =
 $( 0,
\frac{1}{7},
\frac{2}{7},
\frac{4}{7} )
\oplus
( 0 )
$,

\vskip 0.7ex
\hangindent=4em \hangafter=1
 $\rho_\text{isum}(\mathfrak{s})$ =
$\mathrm{i}$($-\sqrt{\frac{1}{7}}$,
$\sqrt{\frac{2}{7}}$,
$\sqrt{\frac{2}{7}}$,
$\sqrt{\frac{2}{7}}$;\ \ 
$-\frac{1}{\sqrt{7}}c^{2}_{7}
$,
$-\frac{1}{\sqrt{7}}c^{1}_{7}
$,
$\frac{1}{\sqrt{7}\mathrm{i}}s^{5}_{28}
$;\ \ 
$\frac{1}{\sqrt{7}\mathrm{i}}s^{5}_{28}
$,
$-\frac{1}{\sqrt{7}}c^{2}_{7}
$;\ \ 
$-\frac{1}{\sqrt{7}}c^{1}_{7}
$)
 $\oplus$
($1$)

Resolved. Number of valid $(S,T)$ pairs = 1.

\vskip 2ex

 \noindent4. (dims;levels) =$(4\oplus
1;9,
1
)$,
irreps = $4_{9,1}^{1,0}\oplus
1_{1}^{1}$,
pord$(\rho_\text{isum}(\mathfrak{t})) = 9$,

\vskip 0.7ex
\hangindent=4em \hangafter=1
 $\rho_\text{isum}(\mathfrak{t})$ =
 $( 0,
\frac{1}{9},
\frac{4}{9},
\frac{7}{9} )
\oplus
( 0 )
$,

\vskip 0.7ex
\hangindent=4em \hangafter=1
 $\rho_\text{isum}(\mathfrak{s})$ =
$\mathrm{i}$($0$,
$\sqrt{\frac{1}{3}}$,
$\sqrt{\frac{1}{3}}$,
$\sqrt{\frac{1}{3}}$;\ \ 
$-\frac{1}{3}c^{1}_{36}
$,
$\frac{1}{3}c^{1}_{36}
-\frac{1}{3}c^{5}_{36}
$,
$\frac{1}{3}c^{5}_{36}
$;\ \ 
$\frac{1}{3}c^{5}_{36}
$,
$-\frac{1}{3}c^{1}_{36}
$;\ \ 
$\frac{1}{3}c^{1}_{36}
-\frac{1}{3}c^{5}_{36}
$)
 $\oplus$
($1$)

Resolved. Number of valid $(S,T)$ pairs = 0.

\vskip 2ex

 \noindent5. (dims;levels) =$(4\oplus
1;9,
1
;a)$,
irreps = $4_{9,2}^{1,0}\oplus
1_{1}^{1}$,
pord$(\rho_\text{isum}(\mathfrak{t})) = 9$,

\vskip 0.7ex
\hangindent=4em \hangafter=1
 $\rho_\text{isum}(\mathfrak{t})$ =
 $( 0,
\frac{1}{9},
\frac{4}{9},
\frac{7}{9} )
\oplus
( 0 )
$,

\vskip 0.7ex
\hangindent=4em \hangafter=1
 $\rho_\text{isum}(\mathfrak{s})$ =
($0$,
$-\sqrt{\frac{1}{3}}$,
$-\sqrt{\frac{1}{3}}$,
$-\sqrt{\frac{1}{3}}$;
$\frac{1}{3}c^{2}_{9}
$,
$\frac{1}{3} c_9^4 $,
$\frac{1}{3}c^{1}_{9}
$;
$\frac{1}{3}c^{1}_{9}
$,
$\frac{1}{3}c^{2}_{9}
$;
$\frac{1}{3} c_9^4 $)
 $\oplus$
($1$)

Unresolved. 

\vskip 2ex

\noindent6. (dims;levels) =$(5;11
)$,
irreps = $5_{11}^{1}$,
pord$(\rho_\text{isum}(\mathfrak{t})) = 11$,

\vskip 0.7ex
\hangindent=4em \hangafter=1
 $\rho_\text{isum}(\mathfrak{t})$ =
 $( \frac{1}{11},
\frac{3}{11},
\frac{4}{11},
\frac{5}{11},
\frac{9}{11} )
$,

\vskip 0.7ex
\hangindent=4em \hangafter=1
 $\rho_\text{isum}(\mathfrak{s})$ =
($-\frac{1}{\sqrt{11}}c^{3}_{44}
$,
$-\frac{1}{\sqrt{11}}c^{7}_{44}
$,
$-\frac{1}{\sqrt{11}}c^{5}_{44}
$,
$-\frac{1}{\sqrt{11}}c^{1}_{44}
$,
$-\frac{1}{\sqrt{11}}c^{9}_{44}
$;
$\frac{1}{\sqrt{11}}c^{9}_{44}
$,
$-\frac{1}{\sqrt{11}}c^{3}_{44}
$,
$\frac{1}{\sqrt{11}}c^{5}_{44}
$,
$\frac{1}{\sqrt{11}}c^{1}_{44}
$;
$\frac{1}{\sqrt{11}}c^{1}_{44}
$,
$\frac{1}{\sqrt{11}}c^{9}_{44}
$,
$\frac{1}{\sqrt{11}}c^{7}_{44}
$;
$\frac{1}{\sqrt{11}}c^{7}_{44}
$,
$-\frac{1}{\sqrt{11}}c^{3}_{44}
$;
$\frac{1}{\sqrt{11}}c^{5}_{44}
$)

Resolved. Number of valid $(S,T)$ pairs = 1.

\vskip 2ex

\

\subsection{A list of rank-5 $S,T$ matrices from resolved representations}

From the representative representation in each GT orbit, if it is revolved, we
can compute all the $S,T$ matrices coming from such a representation.  The
computation steps are displayed below.  For details and notations, see Section
1 of this file.

\

\noindent1. (dims;levels) =$(3 , 
2;5,
2
)$,
irreps = $3_{5}^{1}\oplus
2_{2}^{1,0}$,
pord$(\tilde\rho(\mathfrak{t})) = 10$,

\vskip 0.7ex
\hangindent=4em \hangafter=1
 $\tilde\rho(\mathfrak{t})$ =
 $( 0,
0,
\frac{1}{5},
\frac{1}{2},
\frac{4}{5} )
$,

\vskip 0.7ex
\hangindent=4em \hangafter=1
 $\tilde\rho(\mathfrak{s})$ =
($\sqrt{\frac{1}{5}}$,
$0$,
$-\sqrt{\frac{2}{5}}$,
$0$,
$-\sqrt{\frac{2}{5}}$;
$-\frac{1}{2}$,
$0$,
$-\sqrt{\frac{3}{4}}$,
$0$;
$-\frac{5+\sqrt{5}}{10}$,
$0$,
$\frac{5-\sqrt{5}}{10}$;
$\frac{1}{2}$,
$0$;
$-\frac{5+\sqrt{5}}{10}$)

 \vskip 1ex \setlength{\leftskip}{2em}

\grey{Try $U_0$ =
$\begin{pmatrix}
1,
& 0 \\ 
0,
& 1 \\ 
\end{pmatrix}
$ $\oplus
\begin{pmatrix}
1 \\ 
\end{pmatrix}
$ $\oplus
\begin{pmatrix}
1 \\ 
\end{pmatrix}
$ $\oplus
\begin{pmatrix}
1 \\ 
\end{pmatrix}
$:}\ \ \ \ \ 
\grey{$U_0\tilde\rho(\mathfrak{s})U_0^\dagger$ =}

\grey{$\begin{pmatrix}
\sqrt{\frac{1}{5}},
& 0,
& -\sqrt{\frac{2}{5}},
& 0,
& -\sqrt{\frac{2}{5}} \\ 
0,
& -\frac{1}{2},
& 0,
& -\sqrt{\frac{3}{4}},
& 0 \\ 
-\sqrt{\frac{2}{5}},
& 0,
& -\frac{5+\sqrt{5}}{10},
& 0,
& \frac{5-\sqrt{5}}{10} \\ 
0,
& -\sqrt{\frac{3}{4}},
& 0,
& \frac{1}{2},
& 0 \\ 
-\sqrt{\frac{2}{5}},
& 0,
& \frac{5-\sqrt{5}}{10},
& 0,
& -\frac{5+\sqrt{5}}{10} \\ 
\end{pmatrix}
$}

\grey{Try different $u$'s and signed diagonal matrix $V_\mathrm{sd}$'s:}

 \grey{
\begin{tabular}{|r|l|l|l|l|l|}
\hline
$3_{5}^{1}\oplus
2_{2}^{1,0}:\ u$ 
 & 0 & 1 & 2 & 3 & 4\\ 
 \hline
$D_\rho$ conditions 
 & 1 & 1 & 1 & 1 & 1\\ 
 \hline
$[\rho(\mathfrak{s})\rho(\mathfrak{t})]^3
 = \rho^2(\mathfrak{s}) = \tilde C$ 
 & 0 & 0 & 0 & 0 & 0\\ 
 \hline
$\rho(\mathfrak{s})_{iu}\rho^*(\mathfrak{s})_{ju} \in \mathbb{R}$ 
 & 0 & 0 & 0 & 0 & 0\\ 
 \hline
$\rho(\mathfrak{s})_{i u} \neq 0$  
 & 1 & 1 & 1 & 1 & 1\\ 
 \hline
$\mathrm{cnd}(S)$, $\mathrm{cnd}(\rho(\mathfrak{s}))$ 
 & - & - & - & - & -\\ 
 \hline
$\mathrm{norm}(D^2)$ factors
 & - & - & - & - & -\\ 
 \hline
$1/\rho(\mathfrak{s})_{iu} = $ cyc-int 
 & - & - & - & - & -\\ 
 \hline
norm$(1/\rho(\mathfrak{s})_{iu})$ factors
 & - & - & - & - & -\\ 
 \hline
$\frac{S_{ij}}{S_{uj}} = $ cyc-int
 & - & - & - & - & -\\ 
 \hline
$N^{ij}_k \in \mathbb{N}$
 & - & - & - & - & -\\ 
 \hline
$\exists\ j \text{ that } \frac{S_{ij}}{S_{uj}} \geq 1 $
 & - & - & - & - & -\\ 
 \hline
FS indicator
 & - & - & - & - & -\\ 
 \hline
$C = $ perm-mat
 & - & - & - & - & -\\ 
 \hline
\end{tabular}

Number of valid $(S,T)$ pairs: 0 \vskip 2ex }%grey

\grey{Try $U_0$ =
$\begin{pmatrix}
\sqrt{\frac{1}{2}},
& \sqrt{\frac{1}{2}} \\ 
\sqrt{\frac{1}{2}},
& -\sqrt{\frac{1}{2}} \\ 
\end{pmatrix}
$ $\oplus
\begin{pmatrix}
1 \\ 
\end{pmatrix}
$ $\oplus
\begin{pmatrix}
1 \\ 
\end{pmatrix}
$ $\oplus
\begin{pmatrix}
1 \\ 
\end{pmatrix}
$:}\ \ \ \ \ 
\grey{$U_0\tilde\rho(\mathfrak{s})U_0^\dagger$ =}

\grey{$\begin{pmatrix}
\frac{-5+2\sqrt{5}}{20},
& \frac{5+2\sqrt{5}}{20},
& -\sqrt{\frac{1}{5}},
& -\sqrt{\frac{3}{8}},
& -\sqrt{\frac{1}{5}} \\ 
\frac{5+2\sqrt{5}}{20},
& \frac{-5+2\sqrt{5}}{20},
& -\sqrt{\frac{1}{5}},
& \sqrt{\frac{3}{8}},
& -\sqrt{\frac{1}{5}} \\ 
-\sqrt{\frac{1}{5}},
& -\sqrt{\frac{1}{5}},
& -\frac{5+\sqrt{5}}{10},
& 0,
& \frac{5-\sqrt{5}}{10} \\ 
-\sqrt{\frac{3}{8}},
& \sqrt{\frac{3}{8}},
& 0,
& \frac{1}{2},
& 0 \\ 
-\sqrt{\frac{1}{5}},
& -\sqrt{\frac{1}{5}},
& \frac{5-\sqrt{5}}{10},
& 0,
& -\frac{5+\sqrt{5}}{10} \\ 
\end{pmatrix}
$}

\grey{Try different $u$'s and signed diagonal matrix $V_\mathrm{sd}$'s:}

 \grey{
\begin{tabular}{|r|l|l|l|l|l|}
\hline
$3_{5}^{1}\oplus
2_{2}^{1,0}:\ u$ 
 & 0 & 1 & 2 & 3 & 4\\ 
 \hline
$D_\rho$ conditions 
 & 1 & 1 & 1 & 1 & 1\\ 
 \hline
$[\rho(\mathfrak{s})\rho(\mathfrak{t})]^3
 = \rho^2(\mathfrak{s}) = \tilde C$ 
 & 0 & 0 & 0 & 0 & 0\\ 
 \hline
$\rho(\mathfrak{s})_{iu}\rho^*(\mathfrak{s})_{ju} \in \mathbb{R}$ 
 & 0 & 0 & 0 & 0 & 0\\ 
 \hline
$\rho(\mathfrak{s})_{i u} \neq 0$  
 & 0 & 0 & 1 & 1 & 1\\ 
 \hline
$\mathrm{cnd}(S)$, $\mathrm{cnd}(\rho(\mathfrak{s}))$ 
 & 3 & 3 & - & - & -\\ 
 \hline
$\mathrm{norm}(D^2)$ factors
 & 0 & 0 & - & - & -\\ 
 \hline
$1/\rho(\mathfrak{s})_{iu} = $ cyc-int 
 & 1 & 1 & - & - & -\\ 
 \hline
norm$(1/\rho(\mathfrak{s})_{iu})$ factors
 & 2 & 2 & - & - & -\\ 
 \hline
$\frac{S_{ij}}{S_{uj}} = $ cyc-int
 & 1 & 1 & - & - & -\\ 
 \hline
$N^{ij}_k \in \mathbb{N}$
 & 4 & 4 & - & - & -\\ 
 \hline
$\exists\ j \text{ that } \frac{S_{ij}}{S_{uj}} \geq 1 $
 & 1 & 1 & - & - & -\\ 
 \hline
FS indicator
 & 2 & 2 & - & - & -\\ 
 \hline
$C = $ perm-mat
 & 0 & 0 & - & - & -\\ 
 \hline
\end{tabular}

Number of valid $(S,T)$ pairs: 0 \vskip 2ex }%grey

\grey{Try $U_0$ =
$\begin{pmatrix}
\sqrt{\frac{1}{2}},
& -\sqrt{\frac{1}{2}} \\ 
-\sqrt{\frac{1}{2}},
& -\sqrt{\frac{1}{2}} \\ 
\end{pmatrix}
$ $\oplus
\begin{pmatrix}
1 \\ 
\end{pmatrix}
$ $\oplus
\begin{pmatrix}
1 \\ 
\end{pmatrix}
$ $\oplus
\begin{pmatrix}
1 \\ 
\end{pmatrix}
$:}\ \ \ \ \ 
\grey{$U_0\tilde\rho(\mathfrak{s})U_0^\dagger$ =}

\grey{$\begin{pmatrix}
\frac{-5+2\sqrt{5}}{20},
& -\frac{5+2\sqrt{5}}{20},
& -\sqrt{\frac{1}{5}},
& \sqrt{\frac{3}{8}},
& -\sqrt{\frac{1}{5}} \\ 
-\frac{5+2\sqrt{5}}{20},
& \frac{-5+2\sqrt{5}}{20},
& \sqrt{\frac{1}{5}},
& \sqrt{\frac{3}{8}},
& \sqrt{\frac{1}{5}} \\ 
-\sqrt{\frac{1}{5}},
& \sqrt{\frac{1}{5}},
& -\frac{5+\sqrt{5}}{10},
& 0,
& \frac{5-\sqrt{5}}{10} \\ 
\sqrt{\frac{3}{8}},
& \sqrt{\frac{3}{8}},
& 0,
& \frac{1}{2},
& 0 \\ 
-\sqrt{\frac{1}{5}},
& \sqrt{\frac{1}{5}},
& \frac{5-\sqrt{5}}{10},
& 0,
& -\frac{5+\sqrt{5}}{10} \\ 
\end{pmatrix}
$}

\grey{Try different $u$'s and signed diagonal matrix $V_\mathrm{sd}$'s:}

 \grey{
\begin{tabular}{|r|l|l|l|l|l|}
\hline
$3_{5}^{1}\oplus
2_{2}^{1,0}:\ u$ 
 & 0 & 1 & 2 & 3 & 4\\ 
 \hline
$D_\rho$ conditions 
 & 1 & 1 & 1 & 1 & 1\\ 
 \hline
$[\rho(\mathfrak{s})\rho(\mathfrak{t})]^3
 = \rho^2(\mathfrak{s}) = \tilde C$ 
 & 0 & 0 & 0 & 0 & 0\\ 
 \hline
$\rho(\mathfrak{s})_{iu}\rho^*(\mathfrak{s})_{ju} \in \mathbb{R}$ 
 & 0 & 0 & 0 & 0 & 0\\ 
 \hline
$\rho(\mathfrak{s})_{i u} \neq 0$  
 & 0 & 0 & 1 & 1 & 1\\ 
 \hline
$\mathrm{cnd}(S)$, $\mathrm{cnd}(\rho(\mathfrak{s}))$ 
 & 3 & 3 & - & - & -\\ 
 \hline
$\mathrm{norm}(D^2)$ factors
 & 0 & 0 & - & - & -\\ 
 \hline
$1/\rho(\mathfrak{s})_{iu} = $ cyc-int 
 & 1 & 1 & - & - & -\\ 
 \hline
norm$(1/\rho(\mathfrak{s})_{iu})$ factors
 & 2 & 2 & - & - & -\\ 
 \hline
$\frac{S_{ij}}{S_{uj}} = $ cyc-int
 & 1 & 1 & - & - & -\\ 
 \hline
$N^{ij}_k \in \mathbb{N}$
 & 4 & 4 & - & - & -\\ 
 \hline
$\exists\ j \text{ that } \frac{S_{ij}}{S_{uj}} \geq 1 $
 & 1 & 1 & - & - & -\\ 
 \hline
FS indicator
 & 2 & 2 & - & - & -\\ 
 \hline
$C = $ perm-mat
 & 0 & 0 & - & - & -\\ 
 \hline
\end{tabular}

Number of valid $(S,T)$ pairs: 0 \vskip 2ex }%grey

Total number of valid $(S,T)$ pairs: 0

 \vskip 4ex

\ \setlength{\leftskip}{0em}

\ \setlength{\leftskip}{0em} 

\noindent2. (dims;levels) =$(3 , 
2;8,
3
)$,
irreps = $3_{8}^{1,0}\oplus
2_{3}^{1,0}$,
pord$(\tilde\rho(\mathfrak{t})) = 24$,

\vskip 0.7ex
\hangindent=4em \hangafter=1
 $\tilde\rho(\mathfrak{t})$ =
 $( 0,
0,
\frac{1}{8},
\frac{1}{3},
\frac{5}{8} )
$,

\vskip 0.7ex
\hangindent=4em \hangafter=1
 $\tilde\rho(\mathfrak{s})$ =
$\mathrm{i}$($0$,
$0$,
$\sqrt{\frac{1}{2}}$,
$0$,
$\sqrt{\frac{1}{2}}$;\ \ 
$-\sqrt{\frac{1}{3}}$,
$0$,
$\sqrt{\frac{2}{3}}$,
$0$;\ \ 
$-\frac{1}{2}$,
$0$,
$\frac{1}{2}$;\ \ 
$\sqrt{\frac{1}{3}}$,
$0$;\ \ 
$-\frac{1}{2}$)

 \vskip 1ex \setlength{\leftskip}{2em}

\grey{Try $U_0$ =
$\begin{pmatrix}
1,
& 0 \\ 
0,
& 1 \\ 
\end{pmatrix}
$ $\oplus
\begin{pmatrix}
1 \\ 
\end{pmatrix}
$ $\oplus
\begin{pmatrix}
1 \\ 
\end{pmatrix}
$ $\oplus
\begin{pmatrix}
1 \\ 
\end{pmatrix}
$:}\ \ \ \ \ 
\grey{$U_0\tilde\rho(\mathfrak{s})U_0^\dagger$ =}

\grey{$\begin{pmatrix}
0,
& 0,
& \sqrt{\frac{1}{2}}\mathrm{i},
& 0,
& \sqrt{\frac{1}{2}}\mathrm{i} \\ 
0,
& -\sqrt{\frac{1}{3}}\mathrm{i},
& 0,
& \sqrt{\frac{2}{3}}\mathrm{i},
& 0 \\ 
\sqrt{\frac{1}{2}}\mathrm{i},
& 0,
& -\frac{1}{2}\mathrm{i},
& 0,
& \frac{1}{2}\mathrm{i} \\ 
0,
& \sqrt{\frac{2}{3}}\mathrm{i},
& 0,
& \sqrt{\frac{1}{3}}\mathrm{i},
& 0 \\ 
\sqrt{\frac{1}{2}}\mathrm{i},
& 0,
& \frac{1}{2}\mathrm{i},
& 0,
& -\frac{1}{2}\mathrm{i} \\ 
\end{pmatrix}
$}

\grey{Try different $u$'s and signed diagonal matrix $V_\mathrm{sd}$'s:}

 \grey{
\begin{tabular}{|r|l|l|l|l|l|}
\hline
$3_{8}^{1,0}\oplus
2_{3}^{1,0}:\ u$ 
 & 0 & 1 & 2 & 3 & 4\\ 
 \hline
$D_\rho$ conditions 
 & 1 & 1 & 1 & 1 & 1\\ 
 \hline
$[\rho(\mathfrak{s})\rho(\mathfrak{t})]^3
 = \rho^2(\mathfrak{s}) = \tilde C$ 
 & 0 & 0 & 0 & 0 & 0\\ 
 \hline
$\rho(\mathfrak{s})_{iu}\rho^*(\mathfrak{s})_{ju} \in \mathbb{R}$ 
 & 0 & 0 & 0 & 0 & 0\\ 
 \hline
$\rho(\mathfrak{s})_{i u} \neq 0$  
 & 2 & 1 & 1 & 1 & 1\\ 
 \hline
$\mathrm{cnd}(S)$, $\mathrm{cnd}(\rho(\mathfrak{s}))$ 
 & - & - & - & - & -\\ 
 \hline
$\mathrm{norm}(D^2)$ factors
 & - & - & - & - & -\\ 
 \hline
$1/\rho(\mathfrak{s})_{iu} = $ cyc-int 
 & - & - & - & - & -\\ 
 \hline
norm$(1/\rho(\mathfrak{s})_{iu})$ factors
 & - & - & - & - & -\\ 
 \hline
$\frac{S_{ij}}{S_{uj}} = $ cyc-int
 & - & - & - & - & -\\ 
 \hline
$N^{ij}_k \in \mathbb{N}$
 & - & - & - & - & -\\ 
 \hline
$\exists\ j \text{ that } \frac{S_{ij}}{S_{uj}} \geq 1 $
 & - & - & - & - & -\\ 
 \hline
FS indicator
 & - & - & - & - & -\\ 
 \hline
$C = $ perm-mat
 & - & - & - & - & -\\ 
 \hline
\end{tabular}

Number of valid $(S,T)$ pairs: 0 \vskip 2ex }%grey

\grey{Try $U_0$ =
$\begin{pmatrix}
\sqrt{\frac{1}{2}},
& \sqrt{\frac{1}{2}} \\ 
\sqrt{\frac{1}{2}},
& -\sqrt{\frac{1}{2}} \\ 
\end{pmatrix}
$ $\oplus
\begin{pmatrix}
1 \\ 
\end{pmatrix}
$ $\oplus
\begin{pmatrix}
1 \\ 
\end{pmatrix}
$ $\oplus
\begin{pmatrix}
1 \\ 
\end{pmatrix}
$:}\ \ \ \ \ 
\grey{$U_0\tilde\rho(\mathfrak{s})U_0^\dagger$ =}

\grey{$\begin{pmatrix}
-\sqrt{\frac{1}{12}}\mathrm{i},
& \sqrt{\frac{1}{12}}\mathrm{i},
& \frac{1}{2}\mathrm{i},
& \sqrt{\frac{1}{3}}\mathrm{i},
& \frac{1}{2}\mathrm{i} \\ 
\sqrt{\frac{1}{12}}\mathrm{i},
& -\sqrt{\frac{1}{12}}\mathrm{i},
& \frac{1}{2}\mathrm{i},
& -\sqrt{\frac{1}{3}}\mathrm{i},
& \frac{1}{2}\mathrm{i} \\ 
\frac{1}{2}\mathrm{i},
& \frac{1}{2}\mathrm{i},
& -\frac{1}{2}\mathrm{i},
& 0,
& \frac{1}{2}\mathrm{i} \\ 
\sqrt{\frac{1}{3}}\mathrm{i},
& -\sqrt{\frac{1}{3}}\mathrm{i},
& 0,
& \sqrt{\frac{1}{3}}\mathrm{i},
& 0 \\ 
\frac{1}{2}\mathrm{i},
& \frac{1}{2}\mathrm{i},
& \frac{1}{2}\mathrm{i},
& 0,
& -\frac{1}{2}\mathrm{i} \\ 
\end{pmatrix}
$}

\grey{Try different $u$'s and signed diagonal matrix $V_\mathrm{sd}$'s:}

 \grey{
\begin{tabular}{|r|l|l|l|l|l|}
\hline
$3_{8}^{1,0}\oplus
2_{3}^{1,0}:\ u$ 
 & 0 & 1 & 2 & 3 & 4\\ 
 \hline
$D_\rho$ conditions 
 & 0 & 0 & 0 & 0 & 0\\ 
 \hline
$[\rho(\mathfrak{s})\rho(\mathfrak{t})]^3
 = \rho^2(\mathfrak{s}) = \tilde C$ 
 & 0 & 0 & 0 & 0 & 0\\ 
 \hline
$\rho(\mathfrak{s})_{iu}\rho^*(\mathfrak{s})_{ju} \in \mathbb{R}$ 
 & 0 & 0 & 0 & 0 & 0\\ 
 \hline
$\rho(\mathfrak{s})_{i u} \neq 0$  
 & 0 & 0 & 1 & 1 & 1\\ 
 \hline
$\mathrm{cnd}(S)$, $\mathrm{cnd}(\rho(\mathfrak{s}))$ 
 & 0 & 0 & - & - & -\\ 
 \hline
$\mathrm{norm}(D^2)$ factors
 & 0 & 0 & - & - & -\\ 
 \hline
$1/\rho(\mathfrak{s})_{iu} = $ cyc-int 
 & 0 & 0 & - & - & -\\ 
 \hline
norm$(1/\rho(\mathfrak{s})_{iu})$ factors
 & 0 & 0 & - & - & -\\ 
 \hline
$\frac{S_{ij}}{S_{uj}} = $ cyc-int
 & 0 & 0 & - & - & -\\ 
 \hline
$N^{ij}_k \in \mathbb{N}$
 & 0 & 0 & - & - & -\\ 
 \hline
$\exists\ j \text{ that } \frac{S_{ij}}{S_{uj}} \geq 1 $
 & 0 & 0 & - & - & -\\ 
 \hline
FS indicator
 & 0 & 0 & - & - & -\\ 
 \hline
$C = $ perm-mat
 & 0 & 0 & - & - & -\\ 
 \hline
\end{tabular}

Number of valid $(S,T)$ pairs: 2 \vskip 2ex }%grey

\grey{Try $U_0$ =
$\begin{pmatrix}
\sqrt{\frac{1}{2}},
& -\sqrt{\frac{1}{2}} \\ 
-\sqrt{\frac{1}{2}},
& -\sqrt{\frac{1}{2}} \\ 
\end{pmatrix}
$ $\oplus
\begin{pmatrix}
1 \\ 
\end{pmatrix}
$ $\oplus
\begin{pmatrix}
1 \\ 
\end{pmatrix}
$ $\oplus
\begin{pmatrix}
1 \\ 
\end{pmatrix}
$:}\ \ \ \ \ 
\grey{$U_0\tilde\rho(\mathfrak{s})U_0^\dagger$ =}

\grey{$\begin{pmatrix}
-\sqrt{\frac{1}{12}}\mathrm{i},
& -\sqrt{\frac{1}{12}}\mathrm{i},
& \frac{1}{2}\mathrm{i},
& -\sqrt{\frac{1}{3}}\mathrm{i},
& \frac{1}{2}\mathrm{i} \\ 
-\sqrt{\frac{1}{12}}\mathrm{i},
& -\sqrt{\frac{1}{12}}\mathrm{i},
& -\frac{1}{2}\mathrm{i},
& -\sqrt{\frac{1}{3}}\mathrm{i},
& -\frac{1}{2}\mathrm{i} \\ 
\frac{1}{2}\mathrm{i},
& -\frac{1}{2}\mathrm{i},
& -\frac{1}{2}\mathrm{i},
& 0,
& \frac{1}{2}\mathrm{i} \\ 
-\sqrt{\frac{1}{3}}\mathrm{i},
& -\sqrt{\frac{1}{3}}\mathrm{i},
& 0,
& \sqrt{\frac{1}{3}}\mathrm{i},
& 0 \\ 
\frac{1}{2}\mathrm{i},
& -\frac{1}{2}\mathrm{i},
& \frac{1}{2}\mathrm{i},
& 0,
& -\frac{1}{2}\mathrm{i} \\ 
\end{pmatrix}
$}

\grey{Try different $u$'s and signed diagonal matrix $V_\mathrm{sd}$'s:}

 \grey{
\begin{tabular}{|r|l|l|l|l|l|}
\hline
$3_{8}^{1,0}\oplus
2_{3}^{1,0}:\ u$ 
 & 0 & 1 & 2 & 3 & 4\\ 
 \hline
$D_\rho$ conditions 
 & 0 & 0 & 0 & 0 & 0\\ 
 \hline
$[\rho(\mathfrak{s})\rho(\mathfrak{t})]^3
 = \rho^2(\mathfrak{s}) = \tilde C$ 
 & 0 & 0 & 0 & 0 & 0\\ 
 \hline
$\rho(\mathfrak{s})_{iu}\rho^*(\mathfrak{s})_{ju} \in \mathbb{R}$ 
 & 0 & 0 & 0 & 0 & 0\\ 
 \hline
$\rho(\mathfrak{s})_{i u} \neq 0$  
 & 0 & 0 & 1 & 1 & 1\\ 
 \hline
$\mathrm{cnd}(S)$, $\mathrm{cnd}(\rho(\mathfrak{s}))$ 
 & 0 & 0 & - & - & -\\ 
 \hline
$\mathrm{norm}(D^2)$ factors
 & 0 & 0 & - & - & -\\ 
 \hline
$1/\rho(\mathfrak{s})_{iu} = $ cyc-int 
 & 0 & 0 & - & - & -\\ 
 \hline
norm$(1/\rho(\mathfrak{s})_{iu})$ factors
 & 0 & 0 & - & - & -\\ 
 \hline
$\frac{S_{ij}}{S_{uj}} = $ cyc-int
 & 0 & 0 & - & - & -\\ 
 \hline
$N^{ij}_k \in \mathbb{N}$
 & 0 & 0 & - & - & -\\ 
 \hline
$\exists\ j \text{ that } \frac{S_{ij}}{S_{uj}} \geq 1 $
 & 0 & 0 & - & - & -\\ 
 \hline
FS indicator
 & 0 & 0 & - & - & -\\ 
 \hline
$C = $ perm-mat
 & 0 & 0 & - & - & -\\ 
 \hline
\end{tabular}

Number of valid $(S,T)$ pairs: 2 \vskip 2ex }%grey

Total number of valid $(S,T)$ pairs: 2

 \vskip 4ex

\ \setlength{\leftskip}{0em} 

\noindent3. (dims;levels) =$(4 , 
1;7,
1
)$,
irreps = $4_{7}^{1}\oplus
1_{1}^{1}$,
pord$(\tilde\rho(\mathfrak{t})) = 7$,

\vskip 0.7ex
\hangindent=4em \hangafter=1
 $\tilde\rho(\mathfrak{t})$ =
 $( 0,
0,
\frac{1}{7},
\frac{2}{7},
\frac{4}{7} )
$,

\vskip 0.7ex
\hangindent=4em \hangafter=1
 $\tilde\rho(\mathfrak{s})$ =
($-\sqrt{\frac{1}{7}}\mathrm{i}$,
$0$,
$\sqrt{\frac{2}{7}}\mathrm{i}$,
$\sqrt{\frac{2}{7}}\mathrm{i}$,
$\sqrt{\frac{2}{7}}\mathrm{i}$;
$1$,
$0$,
$0$,
$0$;
$-\frac{1}{\sqrt{7}}c^{2}_{7}
\mathrm{i}$,
$-\frac{1}{\sqrt{7}}c^{1}_{7}
\mathrm{i}$,
$\frac{1}{\sqrt{7}\mathrm{i}}s^{5}_{28}
\mathrm{i}$;
$\frac{1}{\sqrt{7}\mathrm{i}}s^{5}_{28}
\mathrm{i}$,
$-\frac{1}{\sqrt{7}}c^{2}_{7}
\mathrm{i}$;
$-\frac{1}{\sqrt{7}}c^{1}_{7}
\mathrm{i}$)

 \vskip 1ex \setlength{\leftskip}{2em}

\grey{Try $U_0$ =
$\begin{pmatrix}
1,
& 0 \\ 
0,
& 1 \\ 
\end{pmatrix}
$ $\oplus
\begin{pmatrix}
1 \\ 
\end{pmatrix}
$ $\oplus
\begin{pmatrix}
1 \\ 
\end{pmatrix}
$ $\oplus
\begin{pmatrix}
1 \\ 
\end{pmatrix}
$:}\ \ \ \ \ 
\grey{$U_0\tilde\rho(\mathfrak{s})U_0^\dagger$ =}

\grey{$\begin{pmatrix}
-\sqrt{\frac{1}{7}}\mathrm{i},
& 0,
& \sqrt{\frac{2}{7}}\mathrm{i},
& \sqrt{\frac{2}{7}}\mathrm{i},
& \sqrt{\frac{2}{7}}\mathrm{i} \\ 
0,
& 1,
& 0,
& 0,
& 0 \\ 
\sqrt{\frac{2}{7}}\mathrm{i},
& 0,
& -\frac{1}{\sqrt{7}}c^{2}_{7}
\mathrm{i},
& -\frac{1}{\sqrt{7}}c^{1}_{7}
\mathrm{i},
& \frac{1}{\sqrt{7}\mathrm{i}}s^{5}_{28}
\mathrm{i} \\ 
\sqrt{\frac{2}{7}}\mathrm{i},
& 0,
& -\frac{1}{\sqrt{7}}c^{1}_{7}
\mathrm{i},
& \frac{1}{\sqrt{7}\mathrm{i}}s^{5}_{28}
\mathrm{i},
& -\frac{1}{\sqrt{7}}c^{2}_{7}
\mathrm{i} \\ 
\sqrt{\frac{2}{7}}\mathrm{i},
& 0,
& \frac{1}{\sqrt{7}\mathrm{i}}s^{5}_{28}
\mathrm{i},
& -\frac{1}{\sqrt{7}}c^{2}_{7}
\mathrm{i},
& -\frac{1}{\sqrt{7}}c^{1}_{7}
\mathrm{i} \\ 
\end{pmatrix}
$}

\grey{Try different $u$'s and signed diagonal matrix $V_\mathrm{sd}$'s:}

 \grey{
\begin{tabular}{|r|l|l|l|l|l|}
\hline
$4_{7}^{1}\oplus
1_{1}^{1}:\ u$ 
 & 0 & 1 & 2 & 3 & 4\\ 
 \hline
$D_\rho$ conditions 
 & 1 & 1 & 1 & 1 & 1\\ 
 \hline
$[\rho(\mathfrak{s})\rho(\mathfrak{t})]^3
 = \rho^2(\mathfrak{s}) = \tilde C$ 
 & 0 & 0 & 0 & 0 & 0\\ 
 \hline
$\rho(\mathfrak{s})_{iu}\rho^*(\mathfrak{s})_{ju} \in \mathbb{R}$ 
 & 0 & 0 & 0 & 0 & 0\\ 
 \hline
$\rho(\mathfrak{s})_{i u} \neq 0$  
 & 1 & 1 & 1 & 1 & 1\\ 
 \hline
$\mathrm{cnd}(S)$, $\mathrm{cnd}(\rho(\mathfrak{s}))$ 
 & - & - & - & - & -\\ 
 \hline
$\mathrm{norm}(D^2)$ factors
 & - & - & - & - & -\\ 
 \hline
$1/\rho(\mathfrak{s})_{iu} = $ cyc-int 
 & - & - & - & - & -\\ 
 \hline
norm$(1/\rho(\mathfrak{s})_{iu})$ factors
 & - & - & - & - & -\\ 
 \hline
$\frac{S_{ij}}{S_{uj}} = $ cyc-int
 & - & - & - & - & -\\ 
 \hline
$N^{ij}_k \in \mathbb{N}$
 & - & - & - & - & -\\ 
 \hline
$\exists\ j \text{ that } \frac{S_{ij}}{S_{uj}} \geq 1 $
 & - & - & - & - & -\\ 
 \hline
FS indicator
 & - & - & - & - & -\\ 
 \hline
$C = $ perm-mat
 & - & - & - & - & -\\ 
 \hline
\end{tabular}

Number of valid $(S,T)$ pairs: 0 \vskip 2ex }%grey

\grey{Try $U_0$ =
$\begin{pmatrix}
\sqrt{\frac{1}{2}},
& \sqrt{\frac{1}{2}} \\ 
\sqrt{\frac{1}{2}},
& -\sqrt{\frac{1}{2}} \\ 
\end{pmatrix}
$ $\oplus
\begin{pmatrix}
1 \\ 
\end{pmatrix}
$ $\oplus
\begin{pmatrix}
1 \\ 
\end{pmatrix}
$ $\oplus
\begin{pmatrix}
1 \\ 
\end{pmatrix}
$:}\ \ \ \ \ 
\grey{$U_0\tilde\rho(\mathfrak{s})U_0^\dagger$ =}

\grey{$\begin{pmatrix}
\frac{1}{\sqrt{7}}\zeta^{-1}_{28}
+\frac{1}{\sqrt{7}}\zeta^{3}_{28}
-\frac{1}{\sqrt{7}}\zeta^{5}_{28}
,
& -\frac{3}{7}+\frac{1}{7}\zeta^{-1}_{7}
+\frac{1}{7}\zeta^{-2}_{7}
+\frac{1}{7}\zeta^{3}_{7}
,
& \sqrt{\frac{1}{7}}\mathrm{i},
& \sqrt{\frac{1}{7}}\mathrm{i},
& \sqrt{\frac{1}{7}}\mathrm{i} \\ 
-\frac{3}{7}+\frac{1}{7}\zeta^{-1}_{7}
+\frac{1}{7}\zeta^{-2}_{7}
+\frac{1}{7}\zeta^{3}_{7}
,
& \frac{1}{\sqrt{7}}\zeta^{-1}_{28}
+\frac{1}{\sqrt{7}}\zeta^{3}_{28}
-\frac{1}{\sqrt{7}}\zeta^{5}_{28}
,
& \sqrt{\frac{1}{7}}\mathrm{i},
& \sqrt{\frac{1}{7}}\mathrm{i},
& \sqrt{\frac{1}{7}}\mathrm{i} \\ 
\sqrt{\frac{1}{7}}\mathrm{i},
& \sqrt{\frac{1}{7}}\mathrm{i},
& -\frac{1}{\sqrt{7}}c^{2}_{7}
\mathrm{i},
& -\frac{1}{\sqrt{7}}c^{1}_{7}
\mathrm{i},
& \frac{1}{\sqrt{7}\mathrm{i}}s^{5}_{28}
\mathrm{i} \\ 
\sqrt{\frac{1}{7}}\mathrm{i},
& \sqrt{\frac{1}{7}}\mathrm{i},
& -\frac{1}{\sqrt{7}}c^{1}_{7}
\mathrm{i},
& \frac{1}{\sqrt{7}\mathrm{i}}s^{5}_{28}
\mathrm{i},
& -\frac{1}{\sqrt{7}}c^{2}_{7}
\mathrm{i} \\ 
\sqrt{\frac{1}{7}}\mathrm{i},
& \sqrt{\frac{1}{7}}\mathrm{i},
& \frac{1}{\sqrt{7}\mathrm{i}}s^{5}_{28}
\mathrm{i},
& -\frac{1}{\sqrt{7}}c^{2}_{7}
\mathrm{i},
& -\frac{1}{\sqrt{7}}c^{1}_{7}
\mathrm{i} \\ 
\end{pmatrix}
$}

\grey{Try different $u$'s and signed diagonal matrix $V_\mathrm{sd}$'s:}

 \grey{
\begin{tabular}{|r|l|l|l|l|l|}
\hline
$4_{7}^{1}\oplus
1_{1}^{1}:\ u$ 
 & 0 & 1 & 2 & 3 & 4\\ 
 \hline
$D_\rho$ conditions 
 & 0 & 0 & 0 & 0 & 0\\ 
 \hline
$[\rho(\mathfrak{s})\rho(\mathfrak{t})]^3
 = \rho^2(\mathfrak{s}) = \tilde C$ 
 & 0 & 0 & 0 & 0 & 0\\ 
 \hline
$\rho(\mathfrak{s})_{iu}\rho^*(\mathfrak{s})_{ju} \in \mathbb{R}$ 
 & 1 & 1 & 0 & 0 & 0\\ 
 \hline
$\rho(\mathfrak{s})_{i u} \neq 0$  
 & 0 & 0 & 0 & 0 & 0\\ 
 \hline
$\mathrm{cnd}(S)$, $\mathrm{cnd}(\rho(\mathfrak{s}))$ 
 & 0 & 0 & 0 & 0 & 0\\ 
 \hline
$\mathrm{norm}(D^2)$ factors
 & 2 & 2 & 0 & 0 & 0\\ 
 \hline
$1/\rho(\mathfrak{s})_{iu} = $ cyc-int 
 & 1 & 1 & 0 & 0 & 0\\ 
 \hline
norm$(1/\rho(\mathfrak{s})_{iu})$ factors
 & 8 & 8 & 0 & 0 & 0\\ 
 \hline
$\frac{S_{ij}}{S_{uj}} = $ cyc-int
 & 1 & 1 & 0 & 0 & 0\\ 
 \hline
$N^{ij}_k \in \mathbb{N}$
 & 3 & 3 & 0 & 0 & 0\\ 
 \hline
$\exists\ j \text{ that } \frac{S_{ij}}{S_{uj}} \geq 1 $
 & 2 & 2 & 0 & 0 & 0\\ 
 \hline
FS indicator
 & 0 & 0 & 0 & 0 & 0\\ 
 \hline
$C = $ perm-mat
 & 1 & 1 & 0 & 0 & 0\\ 
 \hline
\end{tabular}

Number of valid $(S,T)$ pairs: 1 \vskip 2ex }%grey

\grey{Try $U_0$ =
$\begin{pmatrix}
\sqrt{\frac{1}{2}},
& -\sqrt{\frac{1}{2}} \\ 
-\sqrt{\frac{1}{2}},
& -\sqrt{\frac{1}{2}} \\ 
\end{pmatrix}
$ $\oplus
\begin{pmatrix}
1 \\ 
\end{pmatrix}
$ $\oplus
\begin{pmatrix}
1 \\ 
\end{pmatrix}
$ $\oplus
\begin{pmatrix}
1 \\ 
\end{pmatrix}
$:}\ \ \ \ \ 
\grey{$U_0\tilde\rho(\mathfrak{s})U_0^\dagger$ =}

\grey{$\begin{pmatrix}
\frac{1}{\sqrt{7}}\zeta^{-1}_{28}
+\frac{1}{\sqrt{7}}\zeta^{3}_{28}
-\frac{1}{\sqrt{7}}\zeta^{5}_{28}
,
& \frac{1}{\sqrt{7}}\zeta^{1}_{28}
+\frac{1}{\sqrt{7}}\zeta^{-3}_{28}
-\frac{1}{\sqrt{7}}\zeta^{-5}_{28}
,
& \sqrt{\frac{1}{7}}\mathrm{i},
& \sqrt{\frac{1}{7}}\mathrm{i},
& \sqrt{\frac{1}{7}}\mathrm{i} \\ 
\frac{1}{\sqrt{7}}\zeta^{1}_{28}
+\frac{1}{\sqrt{7}}\zeta^{-3}_{28}
-\frac{1}{\sqrt{7}}\zeta^{-5}_{28}
,
& \frac{1}{\sqrt{7}}\zeta^{-1}_{28}
+\frac{1}{\sqrt{7}}\zeta^{3}_{28}
-\frac{1}{\sqrt{7}}\zeta^{5}_{28}
,
& -\sqrt{\frac{1}{7}}\mathrm{i},
& -\sqrt{\frac{1}{7}}\mathrm{i},
& -\sqrt{\frac{1}{7}}\mathrm{i} \\ 
\sqrt{\frac{1}{7}}\mathrm{i},
& -\sqrt{\frac{1}{7}}\mathrm{i},
& -\frac{1}{\sqrt{7}}c^{2}_{7}
\mathrm{i},
& -\frac{1}{\sqrt{7}}c^{1}_{7}
\mathrm{i},
& \frac{1}{\sqrt{7}\mathrm{i}}s^{5}_{28}
\mathrm{i} \\ 
\sqrt{\frac{1}{7}}\mathrm{i},
& -\sqrt{\frac{1}{7}}\mathrm{i},
& -\frac{1}{\sqrt{7}}c^{1}_{7}
\mathrm{i},
& \frac{1}{\sqrt{7}\mathrm{i}}s^{5}_{28}
\mathrm{i},
& -\frac{1}{\sqrt{7}}c^{2}_{7}
\mathrm{i} \\ 
\sqrt{\frac{1}{7}}\mathrm{i},
& -\sqrt{\frac{1}{7}}\mathrm{i},
& \frac{1}{\sqrt{7}\mathrm{i}}s^{5}_{28}
\mathrm{i},
& -\frac{1}{\sqrt{7}}c^{2}_{7}
\mathrm{i},
& -\frac{1}{\sqrt{7}}c^{1}_{7}
\mathrm{i} \\ 
\end{pmatrix}
$}

\grey{Try different $u$'s and signed diagonal matrix $V_\mathrm{sd}$'s:}

 \grey{
\begin{tabular}{|r|l|l|l|l|l|}
\hline
$4_{7}^{1}\oplus
1_{1}^{1}:\ u$ 
 & 0 & 1 & 2 & 3 & 4\\ 
 \hline
$D_\rho$ conditions 
 & 0 & 0 & 0 & 0 & 0\\ 
 \hline
$[\rho(\mathfrak{s})\rho(\mathfrak{t})]^3
 = \rho^2(\mathfrak{s}) = \tilde C$ 
 & 0 & 0 & 0 & 0 & 0\\ 
 \hline
$\rho(\mathfrak{s})_{iu}\rho^*(\mathfrak{s})_{ju} \in \mathbb{R}$ 
 & 1 & 1 & 0 & 0 & 0\\ 
 \hline
$\rho(\mathfrak{s})_{i u} \neq 0$  
 & 0 & 0 & 0 & 0 & 0\\ 
 \hline
$\mathrm{cnd}(S)$, $\mathrm{cnd}(\rho(\mathfrak{s}))$ 
 & 0 & 0 & 0 & 0 & 0\\ 
 \hline
$\mathrm{norm}(D^2)$ factors
 & 2 & 2 & 0 & 0 & 0\\ 
 \hline
$1/\rho(\mathfrak{s})_{iu} = $ cyc-int 
 & 1 & 1 & 0 & 0 & 0\\ 
 \hline
norm$(1/\rho(\mathfrak{s})_{iu})$ factors
 & 8 & 8 & 0 & 0 & 0\\ 
 \hline
$\frac{S_{ij}}{S_{uj}} = $ cyc-int
 & 1 & 1 & 0 & 0 & 0\\ 
 \hline
$N^{ij}_k \in \mathbb{N}$
 & 3 & 3 & 0 & 0 & 0\\ 
 \hline
$\exists\ j \text{ that } \frac{S_{ij}}{S_{uj}} \geq 1 $
 & 2 & 2 & 0 & 0 & 0\\ 
 \hline
FS indicator
 & 0 & 0 & 0 & 0 & 0\\ 
 \hline
$C = $ perm-mat
 & 1 & 1 & 0 & 0 & 0\\ 
 \hline
\end{tabular}

Number of valid $(S,T)$ pairs: 1 \vskip 2ex }%grey

Total number of valid $(S,T)$ pairs: 1

 \vskip 4ex

\ \setlength{\leftskip}{0em} 

\noindent4. (dims;levels) =$(4 , 
1;9,
1
)$,
irreps = $4_{9,1}^{1,0}\oplus
1_{1}^{1}$,
pord$(\tilde\rho(\mathfrak{t})) = 9$,

\vskip 0.7ex
\hangindent=4em \hangafter=1
 $\tilde\rho(\mathfrak{t})$ =
 $( 0,
0,
\frac{1}{9},
\frac{4}{9},
\frac{7}{9} )
$,

\vskip 0.7ex
\hangindent=4em \hangafter=1
 $\tilde\rho(\mathfrak{s})$ =
($0$,
$0$,
$\sqrt{\frac{1}{3}}\mathrm{i}$,
$\sqrt{\frac{1}{3}}\mathrm{i}$,
$\sqrt{\frac{1}{3}}\mathrm{i}$;
$1$,
$0$,
$0$,
$0$;
$-\frac{1}{3}c^{1}_{36}
\mathrm{i}$,
$(\frac{1}{3}c^{1}_{36}
-\frac{1}{3}c^{5}_{36}
)\mathrm{i}$,
$\frac{1}{3}c^{5}_{36}
\mathrm{i}$;
$\frac{1}{3}c^{5}_{36}
\mathrm{i}$,
$-\frac{1}{3}c^{1}_{36}
\mathrm{i}$;
$(\frac{1}{3}c^{1}_{36}
-\frac{1}{3}c^{5}_{36}
)\mathrm{i}$)

 \vskip 1ex \setlength{\leftskip}{2em}

\grey{Try $U_0$ =
$\begin{pmatrix}
1,
& 0 \\ 
0,
& 1 \\ 
\end{pmatrix}
$ $\oplus
\begin{pmatrix}
1 \\ 
\end{pmatrix}
$ $\oplus
\begin{pmatrix}
1 \\ 
\end{pmatrix}
$ $\oplus
\begin{pmatrix}
1 \\ 
\end{pmatrix}
$:}\ \ \ \ \ 
\grey{$U_0\tilde\rho(\mathfrak{s})U_0^\dagger$ =}

\grey{$\begin{pmatrix}
0,
& 0,
& \sqrt{\frac{1}{3}}\mathrm{i},
& \sqrt{\frac{1}{3}}\mathrm{i},
& \sqrt{\frac{1}{3}}\mathrm{i} \\ 
0,
& 1,
& 0,
& 0,
& 0 \\ 
\sqrt{\frac{1}{3}}\mathrm{i},
& 0,
& -\frac{1}{3}c^{1}_{36}
\mathrm{i},
& (\frac{1}{3}c^{1}_{36}
-\frac{1}{3}c^{5}_{36}
)\mathrm{i},
& \frac{1}{3}c^{5}_{36}
\mathrm{i} \\ 
\sqrt{\frac{1}{3}}\mathrm{i},
& 0,
& (\frac{1}{3}c^{1}_{36}
-\frac{1}{3}c^{5}_{36}
)\mathrm{i},
& \frac{1}{3}c^{5}_{36}
\mathrm{i},
& -\frac{1}{3}c^{1}_{36}
\mathrm{i} \\ 
\sqrt{\frac{1}{3}}\mathrm{i},
& 0,
& \frac{1}{3}c^{5}_{36}
\mathrm{i},
& -\frac{1}{3}c^{1}_{36}
\mathrm{i},
& (\frac{1}{3}c^{1}_{36}
-\frac{1}{3}c^{5}_{36}
)\mathrm{i} \\ 
\end{pmatrix}
$}

\grey{Try different $u$'s and signed diagonal matrix $V_\mathrm{sd}$'s:}

 \grey{
\begin{tabular}{|r|l|l|l|l|l|}
\hline
$4_{9,1}^{1,0}\oplus
1_{1}^{1}:\ u$ 
 & 0 & 1 & 2 & 3 & 4\\ 
 \hline
$D_\rho$ conditions 
 & 0 & 0 & 0 & 0 & 0\\ 
 \hline
$[\rho(\mathfrak{s})\rho(\mathfrak{t})]^3
 = \rho^2(\mathfrak{s}) = \tilde C$ 
 & 0 & 0 & 0 & 0 & 0\\ 
 \hline
$\rho(\mathfrak{s})_{iu}\rho^*(\mathfrak{s})_{ju} \in \mathbb{R}$ 
 & 0 & 0 & 0 & 0 & 0\\ 
 \hline
$\rho(\mathfrak{s})_{i u} \neq 0$  
 & 2 & 1 & 1 & 1 & 1\\ 
 \hline
$\mathrm{cnd}(S)$, $\mathrm{cnd}(\rho(\mathfrak{s}))$ 
 & - & - & - & - & -\\ 
 \hline
$\mathrm{norm}(D^2)$ factors
 & - & - & - & - & -\\ 
 \hline
$1/\rho(\mathfrak{s})_{iu} = $ cyc-int 
 & - & - & - & - & -\\ 
 \hline
norm$(1/\rho(\mathfrak{s})_{iu})$ factors
 & - & - & - & - & -\\ 
 \hline
$\frac{S_{ij}}{S_{uj}} = $ cyc-int
 & - & - & - & - & -\\ 
 \hline
$N^{ij}_k \in \mathbb{N}$
 & - & - & - & - & -\\ 
 \hline
$\exists\ j \text{ that } \frac{S_{ij}}{S_{uj}} \geq 1 $
 & - & - & - & - & -\\ 
 \hline
FS indicator
 & - & - & - & - & -\\ 
 \hline
$C = $ perm-mat
 & - & - & - & - & -\\ 
 \hline
\end{tabular}

Number of valid $(S,T)$ pairs: 0 \vskip 2ex }%grey

\grey{Try $U_0$ =
$\begin{pmatrix}
\sqrt{\frac{1}{2}},
& \sqrt{\frac{1}{2}} \\ 
\sqrt{\frac{1}{2}},
& -\sqrt{\frac{1}{2}} \\ 
\end{pmatrix}
$ $\oplus
\begin{pmatrix}
1 \\ 
\end{pmatrix}
$ $\oplus
\begin{pmatrix}
1 \\ 
\end{pmatrix}
$ $\oplus
\begin{pmatrix}
1 \\ 
\end{pmatrix}
$:}\ \ \ \ \ 
\grey{$U_0\tilde\rho(\mathfrak{s})U_0^\dagger$ =}

\grey{$\begin{pmatrix}
\frac{1}{2},
& -\frac{1}{2},
& \sqrt{\frac{1}{6}}\mathrm{i},
& \sqrt{\frac{1}{6}}\mathrm{i},
& \sqrt{\frac{1}{6}}\mathrm{i} \\ 
-\frac{1}{2},
& \frac{1}{2},
& \sqrt{\frac{1}{6}}\mathrm{i},
& \sqrt{\frac{1}{6}}\mathrm{i},
& \sqrt{\frac{1}{6}}\mathrm{i} \\ 
\sqrt{\frac{1}{6}}\mathrm{i},
& \sqrt{\frac{1}{6}}\mathrm{i},
& -\frac{1}{3}c^{1}_{36}
\mathrm{i},
& (\frac{1}{3}c^{1}_{36}
-\frac{1}{3}c^{5}_{36}
)\mathrm{i},
& \frac{1}{3}c^{5}_{36}
\mathrm{i} \\ 
\sqrt{\frac{1}{6}}\mathrm{i},
& \sqrt{\frac{1}{6}}\mathrm{i},
& (\frac{1}{3}c^{1}_{36}
-\frac{1}{3}c^{5}_{36}
)\mathrm{i},
& \frac{1}{3}c^{5}_{36}
\mathrm{i},
& -\frac{1}{3}c^{1}_{36}
\mathrm{i} \\ 
\sqrt{\frac{1}{6}}\mathrm{i},
& \sqrt{\frac{1}{6}}\mathrm{i},
& \frac{1}{3}c^{5}_{36}
\mathrm{i},
& -\frac{1}{3}c^{1}_{36}
\mathrm{i},
& (\frac{1}{3}c^{1}_{36}
-\frac{1}{3}c^{5}_{36}
)\mathrm{i} \\ 
\end{pmatrix}
$}

\grey{Try different $u$'s and signed diagonal matrix $V_\mathrm{sd}$'s:}

 \grey{
\begin{tabular}{|r|l|l|l|l|l|}
\hline
$4_{9,1}^{1,0}\oplus
1_{1}^{1}:\ u$ 
 & 0 & 1 & 2 & 3 & 4\\ 
 \hline
$D_\rho$ conditions 
 & 1 & 1 & 1 & 1 & 1\\ 
 \hline
$[\rho(\mathfrak{s})\rho(\mathfrak{t})]^3
 = \rho^2(\mathfrak{s}) = \tilde C$ 
 & 0 & 0 & 0 & 0 & 0\\ 
 \hline
$\rho(\mathfrak{s})_{iu}\rho^*(\mathfrak{s})_{ju} \in \mathbb{R}$ 
 & 1 & 1 & 0 & 0 & 0\\ 
 \hline
$\rho(\mathfrak{s})_{i u} \neq 0$  
 & 0 & 0 & 0 & 0 & 0\\ 
 \hline
$\mathrm{cnd}(S)$, $\mathrm{cnd}(\rho(\mathfrak{s}))$ 
 & 3 & 3 & 3 & 3 & 3\\ 
 \hline
$\mathrm{norm}(D^2)$ factors
 & 1 & 1 & 0 & 0 & 0\\ 
 \hline
$1/\rho(\mathfrak{s})_{iu} = $ cyc-int 
 & 0 & 0 & 0 & 0 & 0\\ 
 \hline
norm$(1/\rho(\mathfrak{s})_{iu})$ factors
 & 4 & 4 & 1 & 1 & 1\\ 
 \hline
$\frac{S_{ij}}{S_{uj}} = $ cyc-int
 & 1 & 1 & 1 & 1 & 1\\ 
 \hline
$N^{ij}_k \in \mathbb{N}$
 & 4 & 4 & 4 & 4 & 4\\ 
 \hline
$\exists\ j \text{ that } \frac{S_{ij}}{S_{uj}} \geq 1 $
 & 2 & 2 & 1 & 1 & 1\\ 
 \hline
FS indicator
 & 2 & 2 & 2 & 2 & 2\\ 
 \hline
$C = $ perm-mat
 & 1 & 1 & 0 & 0 & 0\\ 
 \hline
\end{tabular}

Number of valid $(S,T)$ pairs: 0 \vskip 2ex }%grey

\grey{Try $U_0$ =
$\begin{pmatrix}
\sqrt{\frac{1}{2}},
& -\sqrt{\frac{1}{2}} \\ 
-\sqrt{\frac{1}{2}},
& -\sqrt{\frac{1}{2}} \\ 
\end{pmatrix}
$ $\oplus
\begin{pmatrix}
1 \\ 
\end{pmatrix}
$ $\oplus
\begin{pmatrix}
1 \\ 
\end{pmatrix}
$ $\oplus
\begin{pmatrix}
1 \\ 
\end{pmatrix}
$:}\ \ \ \ \ 
\grey{$U_0\tilde\rho(\mathfrak{s})U_0^\dagger$ =}

\grey{$\begin{pmatrix}
\frac{1}{2},
& \frac{1}{2},
& \sqrt{\frac{1}{6}}\mathrm{i},
& \sqrt{\frac{1}{6}}\mathrm{i},
& \sqrt{\frac{1}{6}}\mathrm{i} \\ 
\frac{1}{2},
& \frac{1}{2},
& -\sqrt{\frac{1}{6}}\mathrm{i},
& -\sqrt{\frac{1}{6}}\mathrm{i},
& -\sqrt{\frac{1}{6}}\mathrm{i} \\ 
\sqrt{\frac{1}{6}}\mathrm{i},
& -\sqrt{\frac{1}{6}}\mathrm{i},
& -\frac{1}{3}c^{1}_{36}
\mathrm{i},
& (\frac{1}{3}c^{1}_{36}
-\frac{1}{3}c^{5}_{36}
)\mathrm{i},
& \frac{1}{3}c^{5}_{36}
\mathrm{i} \\ 
\sqrt{\frac{1}{6}}\mathrm{i},
& -\sqrt{\frac{1}{6}}\mathrm{i},
& (\frac{1}{3}c^{1}_{36}
-\frac{1}{3}c^{5}_{36}
)\mathrm{i},
& \frac{1}{3}c^{5}_{36}
\mathrm{i},
& -\frac{1}{3}c^{1}_{36}
\mathrm{i} \\ 
\sqrt{\frac{1}{6}}\mathrm{i},
& -\sqrt{\frac{1}{6}}\mathrm{i},
& \frac{1}{3}c^{5}_{36}
\mathrm{i},
& -\frac{1}{3}c^{1}_{36}
\mathrm{i},
& (\frac{1}{3}c^{1}_{36}
-\frac{1}{3}c^{5}_{36}
)\mathrm{i} \\ 
\end{pmatrix}
$}

\grey{Try different $u$'s and signed diagonal matrix $V_\mathrm{sd}$'s:}

 \grey{
\begin{tabular}{|r|l|l|l|l|l|}
\hline
$4_{9,1}^{1,0}\oplus
1_{1}^{1}:\ u$ 
 & 0 & 1 & 2 & 3 & 4\\ 
 \hline
$D_\rho$ conditions 
 & 1 & 1 & 1 & 1 & 1\\ 
 \hline
$[\rho(\mathfrak{s})\rho(\mathfrak{t})]^3
 = \rho^2(\mathfrak{s}) = \tilde C$ 
 & 0 & 0 & 0 & 0 & 0\\ 
 \hline
$\rho(\mathfrak{s})_{iu}\rho^*(\mathfrak{s})_{ju} \in \mathbb{R}$ 
 & 1 & 1 & 0 & 0 & 0\\ 
 \hline
$\rho(\mathfrak{s})_{i u} \neq 0$  
 & 0 & 0 & 0 & 0 & 0\\ 
 \hline
$\mathrm{cnd}(S)$, $\mathrm{cnd}(\rho(\mathfrak{s}))$ 
 & 3 & 3 & 3 & 3 & 3\\ 
 \hline
$\mathrm{norm}(D^2)$ factors
 & 1 & 1 & 0 & 0 & 0\\ 
 \hline
$1/\rho(\mathfrak{s})_{iu} = $ cyc-int 
 & 0 & 0 & 0 & 0 & 0\\ 
 \hline
norm$(1/\rho(\mathfrak{s})_{iu})$ factors
 & 4 & 4 & 1 & 1 & 1\\ 
 \hline
$\frac{S_{ij}}{S_{uj}} = $ cyc-int
 & 1 & 1 & 1 & 1 & 1\\ 
 \hline
$N^{ij}_k \in \mathbb{N}$
 & 4 & 4 & 4 & 4 & 4\\ 
 \hline
$\exists\ j \text{ that } \frac{S_{ij}}{S_{uj}} \geq 1 $
 & 2 & 2 & 1 & 1 & 1\\ 
 \hline
FS indicator
 & 2 & 2 & 2 & 2 & 2\\ 
 \hline
$C = $ perm-mat
 & 1 & 1 & 1 & 1 & 1\\ 
 \hline
\end{tabular}

Number of valid $(S,T)$ pairs: 0 \vskip 2ex }%grey

Total number of valid $(S,T)$ pairs: 0

 \vskip 4ex

\ \setlength{\leftskip}{0em}

\ \setlength{\leftskip}{0em} 

\noindent5. (dims;levels) =$(4 , 
1;9,
1
;a)$,
irreps = $4_{9,2}^{1,0}\oplus
1_{1}^{1}$,
pord$(\tilde\rho(\mathfrak{t})) = 9$,

\vskip 0.7ex
\hangindent=4em \hangafter=1
 $\tilde\rho(\mathfrak{t})$ =
 $( 0,
0,
\frac{1}{9},
\frac{4}{9},
\frac{7}{9} )
$,

\vskip 0.7ex
\hangindent=4em \hangafter=1
 $\tilde\rho(\mathfrak{s})$ =
($0$,
$0$,
$-\sqrt{\frac{1}{3}}$,
$-\sqrt{\frac{1}{3}}$,
$-\sqrt{\frac{1}{3}}$;
$1$,
$0$,
$0$,
$0$;
$\frac{1}{3}c^{2}_{9}
$,
$\frac{1}{3} c_9^4 $,
$\frac{1}{3}c^{1}_{9}
$;
$\frac{1}{3}c^{1}_{9}
$,
$\frac{1}{3}c^{2}_{9}
$;
$\frac{1}{3} c_9^4 $)

 \vskip 1ex \setlength{\leftskip}{2em}

Unresolved representation

\ \setlength{\leftskip}{0em} 

\noindent6. (dims;levels) =$(5;11
)$,
irreps = $5_{11}^{1}$,
pord$(\tilde\rho(\mathfrak{t})) = 11$,

\vskip 0.7ex
\hangindent=4em \hangafter=1
 $\tilde\rho(\mathfrak{t})$ =
 $( \frac{1}{11},
\frac{3}{11},
\frac{4}{11},
\frac{5}{11},
\frac{9}{11} )
$,

\vskip 0.7ex
\hangindent=4em \hangafter=1
 $\tilde\rho(\mathfrak{s})$ =
($-\frac{1}{\sqrt{11}}c^{3}_{44}
$,
$-\frac{1}{\sqrt{11}}c^{7}_{44}
$,
$-\frac{1}{\sqrt{11}}c^{5}_{44}
$,
$-\frac{1}{\sqrt{11}}c^{1}_{44}
$,
$-\frac{1}{\sqrt{11}}c^{9}_{44}
$;
$\frac{1}{\sqrt{11}}c^{9}_{44}
$,
$-\frac{1}{\sqrt{11}}c^{3}_{44}
$,
$\frac{1}{\sqrt{11}}c^{5}_{44}
$,
$\frac{1}{\sqrt{11}}c^{1}_{44}
$;
$\frac{1}{\sqrt{11}}c^{1}_{44}
$,
$\frac{1}{\sqrt{11}}c^{9}_{44}
$,
$\frac{1}{\sqrt{11}}c^{7}_{44}
$;
$\frac{1}{\sqrt{11}}c^{7}_{44}
$,
$-\frac{1}{\sqrt{11}}c^{3}_{44}
$;
$\frac{1}{\sqrt{11}}c^{5}_{44}
$)

 \vskip 1ex \setlength{\leftskip}{2em}

\grey{Try $U_0$ =
$\begin{pmatrix}
1 \\ 
\end{pmatrix}
$ $\oplus
\begin{pmatrix}
1 \\ 
\end{pmatrix}
$ $\oplus
\begin{pmatrix}
1 \\ 
\end{pmatrix}
$ $\oplus
\begin{pmatrix}
1 \\ 
\end{pmatrix}
$ $\oplus
\begin{pmatrix}
1 \\ 
\end{pmatrix}
$:}\ \ \ \ \ 
\grey{$U_0\tilde\rho(\mathfrak{s})U_0^\dagger$ =}

\grey{$\begin{pmatrix}
-\frac{1}{\sqrt{11}}c^{3}_{44}
,
& -\frac{1}{\sqrt{11}}c^{7}_{44}
,
& -\frac{1}{\sqrt{11}}c^{5}_{44}
,
& -\frac{1}{\sqrt{11}}c^{1}_{44}
,
& -\frac{1}{\sqrt{11}}c^{9}_{44}
 \\ 
-\frac{1}{\sqrt{11}}c^{7}_{44}
,
& \frac{1}{\sqrt{11}}c^{9}_{44}
,
& -\frac{1}{\sqrt{11}}c^{3}_{44}
,
& \frac{1}{\sqrt{11}}c^{5}_{44}
,
& \frac{1}{\sqrt{11}}c^{1}_{44}
 \\ 
-\frac{1}{\sqrt{11}}c^{5}_{44}
,
& -\frac{1}{\sqrt{11}}c^{3}_{44}
,
& \frac{1}{\sqrt{11}}c^{1}_{44}
,
& \frac{1}{\sqrt{11}}c^{9}_{44}
,
& \frac{1}{\sqrt{11}}c^{7}_{44}
 \\ 
-\frac{1}{\sqrt{11}}c^{1}_{44}
,
& \frac{1}{\sqrt{11}}c^{5}_{44}
,
& \frac{1}{\sqrt{11}}c^{9}_{44}
,
& \frac{1}{\sqrt{11}}c^{7}_{44}
,
& -\frac{1}{\sqrt{11}}c^{3}_{44}
 \\ 
-\frac{1}{\sqrt{11}}c^{9}_{44}
,
& \frac{1}{\sqrt{11}}c^{1}_{44}
,
& \frac{1}{\sqrt{11}}c^{7}_{44}
,
& -\frac{1}{\sqrt{11}}c^{3}_{44}
,
& \frac{1}{\sqrt{11}}c^{5}_{44}
 \\ 
\end{pmatrix}
$}

\grey{Try different $u$'s and signed diagonal matrix $V_\mathrm{sd}$'s:}

 \grey{
\begin{tabular}{|r|l|l|l|l|l|}
\hline
$5_{11}^{1}:\ u$ 
 & 0 & 1 & 2 & 3 & 4\\ 
 \hline
$D_\rho$ conditions 
 & 0 & 0 & 0 & 0 & 0\\ 
 \hline
$[\rho(\mathfrak{s})\rho(\mathfrak{t})]^3
 = \rho^2(\mathfrak{s}) = \tilde C$ 
 & 0 & 0 & 0 & 0 & 0\\ 
 \hline
$\rho(\mathfrak{s})_{iu}\rho^*(\mathfrak{s})_{ju} \in \mathbb{R}$ 
 & 0 & 0 & 0 & 0 & 0\\ 
 \hline
$\rho(\mathfrak{s})_{i u} \neq 0$  
 & 0 & 0 & 0 & 0 & 0\\ 
 \hline
$\mathrm{cnd}(S)$, $\mathrm{cnd}(\rho(\mathfrak{s}))$ 
 & 0 & 0 & 0 & 0 & 0\\ 
 \hline
$\mathrm{norm}(D^2)$ factors
 & 0 & 0 & 0 & 0 & 0\\ 
 \hline
$1/\rho(\mathfrak{s})_{iu} = $ cyc-int 
 & 0 & 0 & 0 & 0 & 0\\ 
 \hline
norm$(1/\rho(\mathfrak{s})_{iu})$ factors
 & 0 & 0 & 0 & 0 & 0\\ 
 \hline
$\frac{S_{ij}}{S_{uj}} = $ cyc-int
 & 0 & 0 & 0 & 0 & 0\\ 
 \hline
$N^{ij}_k \in \mathbb{N}$
 & 0 & 0 & 0 & 0 & 0\\ 
 \hline
$\exists\ j \text{ that } \frac{S_{ij}}{S_{uj}} \geq 1 $
 & 0 & 0 & 0 & 0 & 0\\ 
 \hline
FS indicator
 & 0 & 0 & 0 & 0 & 0\\ 
 \hline
$C = $ perm-mat
 & 0 & 0 & 0 & 0 & 0\\ 
 \hline
\end{tabular}

Number of valid $(S,T)$ pairs: 1 \vskip 2ex }%grey

Total number of valid $(S,T)$ pairs: 1

 \vskip 4ex

\ \setlength{\leftskip}{0em}

\

The $S,T$ matrices obtained above (the black or the blue entries below), plus
their Galois conjugations (the grey entries below), form the following list of
rank-5 $S,T$ matrices.  For details and notations, see Section 2 of this file.

\

\noindent1. ind = $(3 , 
2;8,
3
)_{1}^{1}$:\ \ 
$d_i$ = ($1.0$,
$1.0$,
$1.732$,
$1.732$,
$2.0$) 

\vskip 0.7ex
\hangindent=3em \hangafter=1
$D^2=$ 12.0 = 
 $12$

\vskip 0.7ex
\hangindent=3em \hangafter=1
$T = ( 0,
0,
\frac{1}{8},
\frac{5}{8},
\frac{1}{3} )
$,

\vskip 0.7ex
\hangindent=3em \hangafter=1
$S$ = ($ 1$,
$ 1$,
$ \sqrt{3}$,
$ \sqrt{3}$,
$ 2$;\ \ 
$ 1$,
$ -\sqrt{3}$,
$ -\sqrt{3}$,
$ 2$;\ \ 
$ \sqrt{3}$,
$ -\sqrt{3}$,
$0$;\ \ 
$ \sqrt{3}$,
$0$;\ \ 
$ -2$)

\vskip 1ex 
\color{grey}

\noindent2. ind = $(3 , 
2;8,
3
)_{1}^{11}$:\ \ 
$d_i$ = ($1.0$,
$1.0$,
$1.732$,
$1.732$,
$2.0$) 

\vskip 0.7ex
\hangindent=3em \hangafter=1
$D^2=$ 12.0 = 
 $12$

\vskip 0.7ex
\hangindent=3em \hangafter=1
$T = ( 0,
0,
\frac{3}{8},
\frac{7}{8},
\frac{2}{3} )
$,

\vskip 0.7ex
\hangindent=3em \hangafter=1
$S$ = ($ 1$,
$ 1$,
$ \sqrt{3}$,
$ \sqrt{3}$,
$ 2$;\ \ 
$ 1$,
$ -\sqrt{3}$,
$ -\sqrt{3}$,
$ 2$;\ \ 
$ \sqrt{3}$,
$ -\sqrt{3}$,
$0$;\ \ 
$ \sqrt{3}$,
$0$;\ \ 
$ -2$)

\vskip 1ex 
\color{grey}

\noindent3. ind = $(3 , 
2;8,
3
)_{1}^{7}$:\ \ 
$d_i$ = ($1.0$,
$1.0$,
$2.0$,
$-1.732$,
$-1.732$) 

\vskip 0.7ex
\hangindent=3em \hangafter=1
$D^2=$ 12.0 = 
 $12$

\vskip 0.7ex
\hangindent=3em \hangafter=1
$T = ( 0,
0,
\frac{1}{3},
\frac{3}{8},
\frac{7}{8} )
$,

\vskip 0.7ex
\hangindent=3em \hangafter=1
$S$ = ($ 1$,
$ 1$,
$ 2$,
$ -\sqrt{3}$,
$ -\sqrt{3}$;\ \ 
$ 1$,
$ 2$,
$ \sqrt{3}$,
$ \sqrt{3}$;\ \ 
$ -2$,
$0$,
$0$;\ \ 
$ -\sqrt{3}$,
$ \sqrt{3}$;\ \ 
$ -\sqrt{3}$)

Pseudo-unitary $\sim$  
$(3 , 
2;8,
3
)_{2}^{11}$

\vskip 1ex 
\color{grey}

\noindent4. ind = $(3 , 
2;8,
3
)_{1}^{5}$:\ \ 
$d_i$ = ($1.0$,
$1.0$,
$2.0$,
$-1.732$,
$-1.732$) 

\vskip 0.7ex
\hangindent=3em \hangafter=1
$D^2=$ 12.0 = 
 $12$

\vskip 0.7ex
\hangindent=3em \hangafter=1
$T = ( 0,
0,
\frac{2}{3},
\frac{1}{8},
\frac{5}{8} )
$,

\vskip 0.7ex
\hangindent=3em \hangafter=1
$S$ = ($ 1$,
$ 1$,
$ 2$,
$ -\sqrt{3}$,
$ -\sqrt{3}$;\ \ 
$ 1$,
$ 2$,
$ \sqrt{3}$,
$ \sqrt{3}$;\ \ 
$ -2$,
$0$,
$0$;\ \ 
$ -\sqrt{3}$,
$ \sqrt{3}$;\ \ 
$ -\sqrt{3}$)

Pseudo-unitary $\sim$  
$(3 , 
2;8,
3
)_{2}^{1}$

\vskip 1ex 

 \color{black} \vskip 2ex

\noindent5. ind = $(3 , 
2;8,
3
)_{2}^{1}$:\ \ 
$d_i$ = ($1.0$,
$1.0$,
$1.732$,
$1.732$,
$2.0$) 

\vskip 0.7ex
\hangindent=3em \hangafter=1
$D^2=$ 12.0 = 
 $12$

\vskip 0.7ex
\hangindent=3em \hangafter=1
$T = ( 0,
0,
\frac{1}{8},
\frac{5}{8},
\frac{2}{3} )
$,

\vskip 0.7ex
\hangindent=3em \hangafter=1
$S$ = ($ 1$,
$ 1$,
$ \sqrt{3}$,
$ \sqrt{3}$,
$ 2$;\ \ 
$ 1$,
$ -\sqrt{3}$,
$ -\sqrt{3}$,
$ 2$;\ \ 
$ -\sqrt{3}$,
$ \sqrt{3}$,
$0$;\ \ 
$ -\sqrt{3}$,
$0$;\ \ 
$ -2$)

\vskip 1ex 
\color{grey}

\noindent6. ind = $(3 , 
2;8,
3
)_{2}^{11}$:\ \ 
$d_i$ = ($1.0$,
$1.0$,
$1.732$,
$1.732$,
$2.0$) 

\vskip 0.7ex
\hangindent=3em \hangafter=1
$D^2=$ 12.0 = 
 $12$

\vskip 0.7ex
\hangindent=3em \hangafter=1
$T = ( 0,
0,
\frac{3}{8},
\frac{7}{8},
\frac{1}{3} )
$,

\vskip 0.7ex
\hangindent=3em \hangafter=1
$S$ = ($ 1$,
$ 1$,
$ \sqrt{3}$,
$ \sqrt{3}$,
$ 2$;\ \ 
$ 1$,
$ -\sqrt{3}$,
$ -\sqrt{3}$,
$ 2$;\ \ 
$ -\sqrt{3}$,
$ \sqrt{3}$,
$0$;\ \ 
$ -\sqrt{3}$,
$0$;\ \ 
$ -2$)

\vskip 1ex 
\color{grey}

\noindent7. ind = $(3 , 
2;8,
3
)_{2}^{5}$:\ \ 
$d_i$ = ($1.0$,
$1.0$,
$2.0$,
$-1.732$,
$-1.732$) 

\vskip 0.7ex
\hangindent=3em \hangafter=1
$D^2=$ 12.0 = 
 $12$

\vskip 0.7ex
\hangindent=3em \hangafter=1
$T = ( 0,
0,
\frac{1}{3},
\frac{1}{8},
\frac{5}{8} )
$,

\vskip 0.7ex
\hangindent=3em \hangafter=1
$S$ = ($ 1$,
$ 1$,
$ 2$,
$ -\sqrt{3}$,
$ -\sqrt{3}$;\ \ 
$ 1$,
$ 2$,
$ \sqrt{3}$,
$ \sqrt{3}$;\ \ 
$ -2$,
$0$,
$0$;\ \ 
$ \sqrt{3}$,
$ -\sqrt{3}$;\ \ 
$ \sqrt{3}$)

Pseudo-unitary $\sim$  
$(3 , 
2;8,
3
)_{1}^{1}$

\vskip 1ex 
\color{grey}

\noindent8. ind = $(3 , 
2;8,
3
)_{2}^{7}$:\ \ 
$d_i$ = ($1.0$,
$1.0$,
$2.0$,
$-1.732$,
$-1.732$) 

\vskip 0.7ex
\hangindent=3em \hangafter=1
$D^2=$ 12.0 = 
 $12$

\vskip 0.7ex
\hangindent=3em \hangafter=1
$T = ( 0,
0,
\frac{2}{3},
\frac{3}{8},
\frac{7}{8} )
$,

\vskip 0.7ex
\hangindent=3em \hangafter=1
$S$ = ($ 1$,
$ 1$,
$ 2$,
$ -\sqrt{3}$,
$ -\sqrt{3}$;\ \ 
$ 1$,
$ 2$,
$ \sqrt{3}$,
$ \sqrt{3}$;\ \ 
$ -2$,
$0$,
$0$;\ \ 
$ \sqrt{3}$,
$ -\sqrt{3}$;\ \ 
$ \sqrt{3}$)

Pseudo-unitary $\sim$  
$(3 , 
2;8,
3
)_{1}^{11}$

\vskip 1ex 

 \color{black} \vskip 2ex

\noindent9. ind = $(4 , 
1;7,
1
)_{1}^{1}$:\ \ 
$d_i$ = ($1.0$,
$2.246$,
$2.246$,
$2.801$,
$4.48$) 

\vskip 0.7ex
\hangindent=3em \hangafter=1
$D^2=$ 35.342 = 
 $21+14c^{1}_{7}
+7c^{2}_{7}
$

\vskip 0.7ex
\hangindent=3em \hangafter=1
$T = ( 0,
\frac{1}{7},
\frac{1}{7},
\frac{6}{7},
\frac{4}{7} )
$,

\vskip 0.7ex
\hangindent=3em \hangafter=1
$S$ = ($ 1$,
$ \xi_{7}^{3}$,
$ \xi_{7}^{3}$,
$ \xi_{14}^{3}$,
$ \xi_{14}^{5}$;\ \ 
$ s^{1}_{7}
+\zeta^{2}_{7}
+\zeta^{3}_{7}
$,
$ -1-2\zeta^{1}_{7}
-\zeta^{2}_{7}
-\zeta^{3}_{7}
$,
$ -\xi_{7}^{3}$,
$ \xi_{7}^{3}$;\ \ 
$ s^{1}_{7}
+\zeta^{2}_{7}
+\zeta^{3}_{7}
$,
$ -\xi_{7}^{3}$,
$ \xi_{7}^{3}$;\ \ 
$ \xi_{14}^{5}$,
$ -1$;\ \ 
$ -\xi_{14}^{3}$)

\vskip 1ex 
\color{grey}

\noindent10. ind = $(4 , 
1;7,
1
)_{1}^{6}$:\ \ 
$d_i$ = ($1.0$,
$2.246$,
$2.246$,
$2.801$,
$4.48$) 

\vskip 0.7ex
\hangindent=3em \hangafter=1
$D^2=$ 35.342 = 
 $21+14c^{1}_{7}
+7c^{2}_{7}
$

\vskip 0.7ex
\hangindent=3em \hangafter=1
$T = ( 0,
\frac{6}{7},
\frac{6}{7},
\frac{1}{7},
\frac{3}{7} )
$,

\vskip 0.7ex
\hangindent=3em \hangafter=1
$S$ = ($ 1$,
$ \xi_{7}^{3}$,
$ \xi_{7}^{3}$,
$ \xi_{14}^{3}$,
$ \xi_{14}^{5}$;\ \ 
$ -1-2\zeta^{1}_{7}
-\zeta^{2}_{7}
-\zeta^{3}_{7}
$,
$ s^{1}_{7}
+\zeta^{2}_{7}
+\zeta^{3}_{7}
$,
$ -\xi_{7}^{3}$,
$ \xi_{7}^{3}$;\ \ 
$ -1-2\zeta^{1}_{7}
-\zeta^{2}_{7}
-\zeta^{3}_{7}
$,
$ -\xi_{7}^{3}$,
$ \xi_{7}^{3}$;\ \ 
$ \xi_{14}^{5}$,
$ -1$;\ \ 
$ -\xi_{14}^{3}$)

\vskip 1ex 
\color{grey}

\noindent11. ind = $(4 , 
1;7,
1
)_{1}^{2}$:\ \ 
$d_i$ = ($1.0$,
$0.554$,
$0.554$,
$-0.246$,
$-0.692$) 

\vskip 0.7ex
\hangindent=3em \hangafter=1
$D^2=$ 2.155 = 
 $14-7c^{1}_{7}
+7c^{2}_{7}
$

\vskip 0.7ex
\hangindent=3em \hangafter=1
$T = ( 0,
\frac{2}{7},
\frac{2}{7},
\frac{5}{7},
\frac{1}{7} )
$,

\vskip 0.7ex
\hangindent=3em \hangafter=1
$S$ = ($ 1$,
$ \xi_{7}^{1,2}$,
$ \xi_{7}^{1,2}$,
$ -\xi_{14}^{1,5}$,
$ -\xi_{14}^{3,5}$;\ \ 
$ -1-\zeta^{1}_{7}
-2\zeta^{-2}_{7}
-\zeta^{3}_{7}
$,
$ \zeta^{1}_{7}
-s^{2}_{7}
+\zeta^{3}_{7}
$,
$ -\xi_{7}^{1,2}$,
$ \xi_{7}^{1,2}$;\ \ 
$ -1-\zeta^{1}_{7}
-2\zeta^{-2}_{7}
-\zeta^{3}_{7}
$,
$ -\xi_{7}^{1,2}$,
$ \xi_{7}^{1,2}$;\ \ 
$ -\xi_{14}^{3,5}$,
$ -1$;\ \ 
$ \xi_{14}^{1,5}$)

Not pseudo-unitary. 

\vskip 1ex 
\color{grey}

\noindent12. ind = $(4 , 
1;7,
1
)_{1}^{4}$:\ \ 
$d_i$ = ($1.0$,
$1.445$,
$-0.356$,
$-0.801$,
$-0.801$) 

\vskip 0.7ex
\hangindent=3em \hangafter=1
$D^2=$ 4.501 = 
 $7-7c^{1}_{7}
-14c^{2}_{7}
$

\vskip 0.7ex
\hangindent=3em \hangafter=1
$T = ( 0,
\frac{3}{7},
\frac{2}{7},
\frac{4}{7},
\frac{4}{7} )
$,

\vskip 0.7ex
\hangindent=3em \hangafter=1
$S$ = ($ 1$,
$ \xi_{14}^{5,3}$,
$ -\xi_{14}^{1,3}$,
$ -\xi_{7}^{2,3}$,
$ -\xi_{7}^{2,3}$;\ \ 
$ -\xi_{14}^{1,3}$,
$ -1$,
$ \xi_{7}^{2,3}$,
$ \xi_{7}^{2,3}$;\ \ 
$ -\xi_{14}^{5,3}$,
$ -\xi_{7}^{2,3}$,
$ -\xi_{7}^{2,3}$;\ \ 
$ -1-\zeta^{-1}_{7}
-\zeta^{2}_{7}
-2\zeta^{3}_{7}
$,
$ 1+\zeta^{1}_{7}
+2\zeta^{-1}_{7}
+2\zeta^{2}_{7}
+\zeta^{-2}_{7}
+2\zeta^{3}_{7}
$;\ \ 
$ -1-\zeta^{-1}_{7}
-\zeta^{2}_{7}
-2\zeta^{3}_{7}
$)

Not pseudo-unitary. 

\vskip 1ex 
\color{grey}

\noindent13. ind = $(4 , 
1;7,
1
)_{1}^{3}$:\ \ 
$d_i$ = ($1.0$,
$1.445$,
$-0.356$,
$-0.801$,
$-0.801$) 

\vskip 0.7ex
\hangindent=3em \hangafter=1
$D^2=$ 4.501 = 
 $7-7c^{1}_{7}
-14c^{2}_{7}
$

\vskip 0.7ex
\hangindent=3em \hangafter=1
$T = ( 0,
\frac{4}{7},
\frac{5}{7},
\frac{3}{7},
\frac{3}{7} )
$,

\vskip 0.7ex
\hangindent=3em \hangafter=1
$S$ = ($ 1$,
$ \xi_{14}^{5,3}$,
$ -\xi_{14}^{1,3}$,
$ -\xi_{7}^{2,3}$,
$ -\xi_{7}^{2,3}$;\ \ 
$ -\xi_{14}^{1,3}$,
$ -1$,
$ \xi_{7}^{2,3}$,
$ \xi_{7}^{2,3}$;\ \ 
$ -\xi_{14}^{5,3}$,
$ -\xi_{7}^{2,3}$,
$ -\xi_{7}^{2,3}$;\ \ 
$ 1+\zeta^{1}_{7}
+2\zeta^{-1}_{7}
+2\zeta^{2}_{7}
+\zeta^{-2}_{7}
+2\zeta^{3}_{7}
$,
$ -1-\zeta^{-1}_{7}
-\zeta^{2}_{7}
-2\zeta^{3}_{7}
$;\ \ 
$ 1+\zeta^{1}_{7}
+2\zeta^{-1}_{7}
+2\zeta^{2}_{7}
+\zeta^{-2}_{7}
+2\zeta^{3}_{7}
$)

Not pseudo-unitary. 

\vskip 1ex 
\color{grey}

\noindent14. ind = $(4 , 
1;7,
1
)_{1}^{5}$:\ \ 
$d_i$ = ($1.0$,
$0.554$,
$0.554$,
$-0.246$,
$-0.692$) 

\vskip 0.7ex
\hangindent=3em \hangafter=1
$D^2=$ 2.155 = 
 $14-7c^{1}_{7}
+7c^{2}_{7}
$

\vskip 0.7ex
\hangindent=3em \hangafter=1
$T = ( 0,
\frac{5}{7},
\frac{5}{7},
\frac{2}{7},
\frac{6}{7} )
$,

\vskip 0.7ex
\hangindent=3em \hangafter=1
$S$ = ($ 1$,
$ \xi_{7}^{1,2}$,
$ \xi_{7}^{1,2}$,
$ -\xi_{14}^{1,5}$,
$ -\xi_{14}^{3,5}$;\ \ 
$ \zeta^{1}_{7}
-s^{2}_{7}
+\zeta^{3}_{7}
$,
$ -1-\zeta^{1}_{7}
-2\zeta^{-2}_{7}
-\zeta^{3}_{7}
$,
$ -\xi_{7}^{1,2}$,
$ \xi_{7}^{1,2}$;\ \ 
$ \zeta^{1}_{7}
-s^{2}_{7}
+\zeta^{3}_{7}
$,
$ -\xi_{7}^{1,2}$,
$ \xi_{7}^{1,2}$;\ \ 
$ -\xi_{14}^{3,5}$,
$ -1$;\ \ 
$ \xi_{14}^{1,5}$)

Not pseudo-unitary. 

\vskip 1ex 

 \color{black} \vskip 2ex

\noindent15. ind = $(5;11
)_{1}^{1}$:\ \ 
$d_i$ = ($1.0$,
$1.918$,
$2.682$,
$3.228$,
$3.513$) 

\vskip 0.7ex
\hangindent=3em \hangafter=1
$D^2=$ 34.646 = 
 $15+10c^{1}_{11}
+6c^{2}_{11}
+3c^{3}_{11}
+c^{4}_{11}
$

\vskip 0.7ex
\hangindent=3em \hangafter=1
$T = ( 0,
\frac{2}{11},
\frac{9}{11},
\frac{10}{11},
\frac{5}{11} )
$,

\vskip 0.7ex
\hangindent=3em \hangafter=1
$S$ = ($ 1$,
$ \xi_{11}^{2}$,
$ \xi_{11}^{3}$,
$ \xi_{11}^{4}$,
$ \xi_{11}^{5}$;\ \ 
$ -\xi_{11}^{4}$,
$ \xi_{11}^{5}$,
$ -\xi_{11}^{3}$,
$ 1$;\ \ 
$ \xi_{11}^{2}$,
$ -1$,
$ -\xi_{11}^{4}$;\ \ 
$ \xi_{11}^{5}$,
$ -\xi_{11}^{2}$;\ \ 
$ \xi_{11}^{3}$)

\vskip 1ex 
\color{grey}

\noindent16. ind = $(5;11
)_{1}^{10}$:\ \ 
$d_i$ = ($1.0$,
$1.918$,
$2.682$,
$3.228$,
$3.513$) 

\vskip 0.7ex
\hangindent=3em \hangafter=1
$D^2=$ 34.646 = 
 $15+10c^{1}_{11}
+6c^{2}_{11}
+3c^{3}_{11}
+c^{4}_{11}
$

\vskip 0.7ex
\hangindent=3em \hangafter=1
$T = ( 0,
\frac{9}{11},
\frac{2}{11},
\frac{1}{11},
\frac{6}{11} )
$,

\vskip 0.7ex
\hangindent=3em \hangafter=1
$S$ = ($ 1$,
$ \xi_{11}^{2}$,
$ \xi_{11}^{3}$,
$ \xi_{11}^{4}$,
$ \xi_{11}^{5}$;\ \ 
$ -\xi_{11}^{4}$,
$ \xi_{11}^{5}$,
$ -\xi_{11}^{3}$,
$ 1$;\ \ 
$ \xi_{11}^{2}$,
$ -1$,
$ -\xi_{11}^{4}$;\ \ 
$ \xi_{11}^{5}$,
$ -\xi_{11}^{2}$;\ \ 
$ \xi_{11}^{3}$)

\vskip 1ex 
\color{grey}

\noindent17. ind = $(5;11
)_{1}^{9}$:\ \ 
$d_i$ = ($1.0$,
$0.521$,
$1.830$,
$-1.397$,
$-1.682$) 

\vskip 0.7ex
\hangindent=3em \hangafter=1
$D^2=$ 9.408 = 
 $12-3c^{1}_{11}
+7c^{2}_{11}
-2c^{3}_{11}
+3c^{4}_{11}
$

\vskip 0.7ex
\hangindent=3em \hangafter=1
$T = ( 0,
\frac{1}{11},
\frac{4}{11},
\frac{2}{11},
\frac{7}{11} )
$,

\vskip 0.7ex
\hangindent=3em \hangafter=1
$S$ = ($ 1$,
$ \xi_{11}^{1,2}$,
$ \xi_{11}^{5,2}$,
$ -\xi_{11}^{3,2}$,
$ -c^{1}_{11}
$;\ \ 
$ \xi_{11}^{5,2}$,
$ \xi_{11}^{3,2}$,
$ c^{1}_{11}
$,
$ 1$;\ \ 
$ -c^{1}_{11}
$,
$ -1$,
$ \xi_{11}^{1,2}$;\ \ 
$ \xi_{11}^{1,2}$,
$ -\xi_{11}^{5,2}$;\ \ 
$ \xi_{11}^{3,2}$)

Not pseudo-unitary. 

\vskip 1ex 
\color{grey}

\noindent18. ind = $(5;11
)_{1}^{6}$:\ \ 
$d_i$ = ($1.0$,
$0.284$,
$0.763$,
$-0.546$,
$-0.918$) 

\vskip 0.7ex
\hangindent=3em \hangafter=1
$D^2=$ 2.806 = 
 $5-4c^{1}_{11}
-9c^{2}_{11}
-10c^{3}_{11}
-7c^{4}_{11}
$

\vskip 0.7ex
\hangindent=3em \hangafter=1
$T = ( 0,
\frac{1}{11},
\frac{8}{11},
\frac{5}{11},
\frac{10}{11} )
$,

\vskip 0.7ex
\hangindent=3em \hangafter=1
$S$ = ($ 1$,
$ -c^{3}_{11}
$,
$ \xi_{11}^{3,5}$,
$ -\xi_{11}^{2,5}$,
$ -\xi_{11}^{4,5}$;\ \ 
$ \xi_{11}^{2,5}$,
$ 1$,
$ \xi_{11}^{4,5}$,
$ \xi_{11}^{3,5}$;\ \ 
$ -\xi_{11}^{4,5}$,
$ c^{3}_{11}
$,
$ \xi_{11}^{2,5}$;\ \ 
$ \xi_{11}^{3,5}$,
$ -1$;\ \ 
$ -c^{3}_{11}
$)

Not pseudo-unitary. 

\vskip 1ex 
\color{grey}

\noindent19. ind = $(5;11
)_{1}^{7}$:\ \ 
$d_i$ = ($1.0$,
$1.88$,
$-0.309$,
$-0.594$,
$-0.830$) 

\vskip 0.7ex
\hangindent=3em \hangafter=1
$D^2=$ 3.323 = 
 $14+2c^{1}_{11}
-c^{2}_{11}
+5c^{3}_{11}
+9c^{4}_{11}
$

\vskip 0.7ex
\hangindent=3em \hangafter=1
$T = ( 0,
\frac{4}{11},
\frac{8}{11},
\frac{2}{11},
\frac{3}{11} )
$,

\vskip 0.7ex
\hangindent=3em \hangafter=1
$S$ = ($ 1$,
$ \xi_{11}^{5,4}$,
$ -\xi_{11}^{1,4}$,
$ -\xi_{11}^{2,4}$,
$ -c^{2}_{11}
$;\ \ 
$ -\xi_{11}^{2,4}$,
$ -1$,
$ c^{2}_{11}
$,
$ \xi_{11}^{1,4}$;\ \ 
$ -c^{2}_{11}
$,
$ -\xi_{11}^{5,4}$,
$ -\xi_{11}^{2,4}$;\ \ 
$ -\xi_{11}^{1,4}$,
$ 1$;\ \ 
$ -\xi_{11}^{5,4}$)

Not pseudo-unitary. 

\vskip 1ex 
\color{grey}

\noindent20. ind = $(5;11
)_{1}^{3}$:\ \ 
$d_i$ = ($1.0$,
$0.715$,
$1.309$,
$-0.372$,
$-1.203$) 

\vskip 0.7ex
\hangindent=3em \hangafter=1
$D^2=$ 4.814 = 
 $9-5c^{1}_{11}
-3c^{2}_{11}
+4c^{3}_{11}
-6c^{4}_{11}
$

\vskip 0.7ex
\hangindent=3em \hangafter=1
$T = ( 0,
\frac{5}{11},
\frac{6}{11},
\frac{8}{11},
\frac{4}{11} )
$,

\vskip 0.7ex
\hangindent=3em \hangafter=1
$S$ = ($ 1$,
$ \xi_{11}^{2,3}$,
$ -c^{4}_{11}
$,
$ -\xi_{11}^{1,3}$,
$ -\xi_{11}^{4,3}$;\ \ 
$ -c^{4}_{11}
$,
$ -\xi_{11}^{4,3}$,
$ -1$,
$ \xi_{11}^{1,3}$;\ \ 
$ \xi_{11}^{1,3}$,
$ -\xi_{11}^{2,3}$,
$ 1$;\ \ 
$ -\xi_{11}^{4,3}$,
$ c^{4}_{11}
$;\ \ 
$ \xi_{11}^{2,3}$)

Not pseudo-unitary. 

\vskip 1ex 
\color{grey}

\noindent21. ind = $(5;11
)_{1}^{8}$:\ \ 
$d_i$ = ($1.0$,
$0.715$,
$1.309$,
$-0.372$,
$-1.203$) 

\vskip 0.7ex
\hangindent=3em \hangafter=1
$D^2=$ 4.814 = 
 $9-5c^{1}_{11}
-3c^{2}_{11}
+4c^{3}_{11}
-6c^{4}_{11}
$

\vskip 0.7ex
\hangindent=3em \hangafter=1
$T = ( 0,
\frac{6}{11},
\frac{5}{11},
\frac{3}{11},
\frac{7}{11} )
$,

\vskip 0.7ex
\hangindent=3em \hangafter=1
$S$ = ($ 1$,
$ \xi_{11}^{2,3}$,
$ -c^{4}_{11}
$,
$ -\xi_{11}^{1,3}$,
$ -\xi_{11}^{4,3}$;\ \ 
$ -c^{4}_{11}
$,
$ -\xi_{11}^{4,3}$,
$ -1$,
$ \xi_{11}^{1,3}$;\ \ 
$ \xi_{11}^{1,3}$,
$ -\xi_{11}^{2,3}$,
$ 1$;\ \ 
$ -\xi_{11}^{4,3}$,
$ c^{4}_{11}
$;\ \ 
$ \xi_{11}^{2,3}$)

Not pseudo-unitary. 

\vskip 1ex 
\color{grey}

\noindent22. ind = $(5;11
)_{1}^{4}$:\ \ 
$d_i$ = ($1.0$,
$1.88$,
$-0.309$,
$-0.594$,
$-0.830$) 

\vskip 0.7ex
\hangindent=3em \hangafter=1
$D^2=$ 3.323 = 
 $14+2c^{1}_{11}
-c^{2}_{11}
+5c^{3}_{11}
+9c^{4}_{11}
$

\vskip 0.7ex
\hangindent=3em \hangafter=1
$T = ( 0,
\frac{7}{11},
\frac{3}{11},
\frac{9}{11},
\frac{8}{11} )
$,

\vskip 0.7ex
\hangindent=3em \hangafter=1
$S$ = ($ 1$,
$ \xi_{11}^{5,4}$,
$ -\xi_{11}^{1,4}$,
$ -\xi_{11}^{2,4}$,
$ -c^{2}_{11}
$;\ \ 
$ -\xi_{11}^{2,4}$,
$ -1$,
$ c^{2}_{11}
$,
$ \xi_{11}^{1,4}$;\ \ 
$ -c^{2}_{11}
$,
$ -\xi_{11}^{5,4}$,
$ -\xi_{11}^{2,4}$;\ \ 
$ -\xi_{11}^{1,4}$,
$ 1$;\ \ 
$ -\xi_{11}^{5,4}$)

Not pseudo-unitary. 

\vskip 1ex 
\color{grey}

\noindent23. ind = $(5;11
)_{1}^{5}$:\ \ 
$d_i$ = ($1.0$,
$0.284$,
$0.763$,
$-0.546$,
$-0.918$) 

\vskip 0.7ex
\hangindent=3em \hangafter=1
$D^2=$ 2.806 = 
 $5-4c^{1}_{11}
-9c^{2}_{11}
-10c^{3}_{11}
-7c^{4}_{11}
$

\vskip 0.7ex
\hangindent=3em \hangafter=1
$T = ( 0,
\frac{10}{11},
\frac{3}{11},
\frac{6}{11},
\frac{1}{11} )
$,

\vskip 0.7ex
\hangindent=3em \hangafter=1
$S$ = ($ 1$,
$ -c^{3}_{11}
$,
$ \xi_{11}^{3,5}$,
$ -\xi_{11}^{2,5}$,
$ -\xi_{11}^{4,5}$;\ \ 
$ \xi_{11}^{2,5}$,
$ 1$,
$ \xi_{11}^{4,5}$,
$ \xi_{11}^{3,5}$;\ \ 
$ -\xi_{11}^{4,5}$,
$ c^{3}_{11}
$,
$ \xi_{11}^{2,5}$;\ \ 
$ \xi_{11}^{3,5}$,
$ -1$;\ \ 
$ -c^{3}_{11}
$)

Not pseudo-unitary. 

\vskip 1ex 
\color{grey}

\noindent24. ind = $(5;11
)_{1}^{2}$:\ \ 
$d_i$ = ($1.0$,
$0.521$,
$1.830$,
$-1.397$,
$-1.682$) 

\vskip 0.7ex
\hangindent=3em \hangafter=1
$D^2=$ 9.408 = 
 $12-3c^{1}_{11}
+7c^{2}_{11}
-2c^{3}_{11}
+3c^{4}_{11}
$

\vskip 0.7ex
\hangindent=3em \hangafter=1
$T = ( 0,
\frac{10}{11},
\frac{7}{11},
\frac{9}{11},
\frac{4}{11} )
$,

\vskip 0.7ex
\hangindent=3em \hangafter=1
$S$ = ($ 1$,
$ \xi_{11}^{1,2}$,
$ \xi_{11}^{5,2}$,
$ -\xi_{11}^{3,2}$,
$ -c^{1}_{11}
$;\ \ 
$ \xi_{11}^{5,2}$,
$ \xi_{11}^{3,2}$,
$ c^{1}_{11}
$,
$ 1$;\ \ 
$ -c^{1}_{11}
$,
$ -1$,
$ \xi_{11}^{1,2}$;\ \ 
$ \xi_{11}^{1,2}$,
$ -\xi_{11}^{5,2}$;\ \ 
$ \xi_{11}^{3,2}$)

Not pseudo-unitary. 

\vskip 1ex 

 \color{black} \vskip 2ex

\

The above list includes all modular data (unitary or non-unitary) from resolved
$\SL$ representations and non-integral MTCs.
For rank 5, there is only one GT orbit, $(4\oplus 1;9,1;a)$, that is unresolved.
Through a further calculation, we find that such a  GT orbit does not produce
any valid $S,T$ matrices.  Thus the list actually includes all modular data
from non-integral MTCs.

\end{document}